\def\ifundefined#1#2{\expandafter\ifx\csname#1\endcsname\relax\input #2\fi}
\input amssym.def

\def\RR{{\Bbb R}}
\def\CC{{\Bbb C}}
\def\NN{{\Bbb N}}
\def\ZZ{{\Bbb Z}}
\def\FF{{\Bbb F}}
\def\QQ{{\Bbb Q}}
\def\AA{{\Bbb A}}

\def\DD{{\Bbb D}}
\def\EE{{\Bbb E}}
\def\GG{{\Bbb G}}
\def\HH{{\Bbb H}}

\def\MM{{\Bbb M}}

\def\PP{{\Bbb P}}
\def\SS{{\Bbb S}}
\def\TT{{\Bbb T}}
\def\UU{{\Bbb U}}

\def\Ma{{\cal\char'101}}
\def\Mb{{\cal\char'102}}
\def\Mc{{\cal\char'103}}
\def\Md{{\cal\char'104}}
\def\Me{{\cal\char'105}}
\def\Mf{{\cal\char'106}}
\def\Mg{{\cal\char'107}}
\def\Mh{{\cal\char'110}}
\def\Mi{{\cal\char'111}}
\def\Mj{{\cal\char'112}}
\def\Mk{{\cal\char'113}}
\def\Ml{{\cal\char'114}}
\def\Mm{{\cal\char'115}}
\def\Mn{{\cal\char'116}}
\def\Mo{{\cal\char'117}}
\def\Mp{{\cal\char'120}}

\def\Mr{{\cal\char'122}}
\def\Ms{{\cal\char'123}}
\def\Mt{{\cal\char'124}}
\def\Mu{{\cal\char'125}}
\def\Mv{{\cal\char'126}}
\def\Mw{{\cal\char'127}}
\def\Mx{{\cal\char'130}}
\def\My{{\cal\char'131}}
\def\Mz{{\cal\char'132}}

\def\pmb#1{\setbox0=\hbox{$#1$}       
     \kern-.025em\copy0\kern-\wd0
     \kern.05em\copy0\kern-\wd0
     \kern-.025em\box0}



\def\bsl{{\backslash}}

\def\lra{{\longrightarrow}}

\def\abs#1{\left\vert #1\right\vert}    

\def\Cross{\bigm| \kern-5.5pt \not \ \, }
\def\cross{\mid \kern-5.0pt \not \ \, }             
\def\notto{\hbox{$~\rightarrow~\kern-1.5em\hbox{/}\ \ $}}

\def\al{\alpha}
\def\be{\beta}
\def\ga{\gamma}

\def\om{\omega}
\def\Om{\Omega}

\def\dl{\delta}
\def\Dl{\Delta}
\def\vph{\varphi}
\def\vep{\varepsilon}

\def\sg{\sigma}

\hyphenation{math-ema-ticians}
\hyphenation{pa-ra-meters}
\hyphenation{pa-ra-meter}
\hyphenation{lem-ma}
\hyphenation{lem-mas}
\hyphenation{to-po-logy}
\hyphenation{to-po-logies}
\hyphenation{homo-logy}
\hyphenation{homo-mor-phy}

\def\nSigma{\Sigma \kern-8.3pt \bigm|\,}

\def\got#1{\hbox{\teneuler #1}}

\font\teneufm=eufm10
\font\eighteufm=eufm8
\font\fiveeufm=eufm5

\newfam\eufam
\textfont\eufam=\teneufm
\scriptfont\eufam=\eighteufm
\scriptscriptfont\eufam=\fiveeufm

\def\got{\fam=\eufam\teneufm}

\def\boxit#1{\vbox{\hrule\hbox{\vrule\kern2.0pt
       \vbox{\kern2.0pt#1\kern2.0pt}\kern2.0pt\vrule}\hrule}}

\def\vlra#1{\hbox{\kern-1pt
       \hbox{\raise2.38pt\hbox{\vbox{\hrule width#1 height0.26pt}}}
       \kern-4.0pt$\rightarrow$}}

\def\vlla#1{\hbox{$\leftarrow$\kern-1.0pt
       \hbox{\raise2.38pt\hbox{\vbox{\hrule width#1 height0.26pt}}}}}

\def\vlda#1{\hbox{$\leftarrow$\kern-1.0pt
       \hbox{\raise2.38pt\hbox{\vbox{\hrule width#1 height0.26pt}}}
       \kern-4.0pt$\rightarrow$}}

\def\longra#1#2#3{\,\lower3pt\hbox{${\buildrel\mathop{#2}
\over{{\vlra{#1}}\atop{#3}}}$}\,}

\def\longla#1#2#3{\,\lower3pt\hbox{${\buildrel\mathop{#2}
\over{{\vlla{#1}}\atop{#3}}}$}\,}

\def\longda#1#2#3{\,\lower3pt\hbox{${\buildrel\mathop{#2}
\over{{\vlda{#1}}\atop{#3}}}$}\,}

\def\overrightharpoonup#1{\vbox{\ialign{##\crcr
	$\rightharpoonup$\crcr\noalign{\kern-1pt\nointerlineskip}
	$\hfil\displaystyle{#1}\hfil$\crcr}}}

\catcode`@=11
\def\@@dalembert#1#2{\setbox0\hbox{$#1\rm I$}
  \vrule height.90\ht0 depth.1\ht0 width.04\ht0
  \rlap{\vrule height.90\ht0 depth-.86\ht0 width.8\ht0}
  \vrule height0\ht0 depth.1\ht0 width.8\ht0
  \vrule height.9\ht0 depth.1\ht0 width.1\ht0 }
\def\dalembert{\mathord{\mkern2mu\mathpalette\@@dalembert{}\mkern2mu}}

\def\@@varcirc#1#2{\mathord{\lower#1ex\hbox{\m@th${#2\mathchar\hex0017 }$}}}
\def\varcirc{\mathchoice
  {\@@varcirc{.91}\displaystyle}{\@@varcirc{.91}\textstyle}
{\@@varcirc{.45}\scriptscriptstyle}}
\catcode`@=12

\font\tensf=cmss10 \font\sevensf=cmss8 at 7pt
\newfam\sffam\def\sf{\fam\sffam\tensf}
\textfont\sffam=\tensf\scriptfont\sffam=\sevensf
\def\M{{\mathord{\sf M}}}
\input amssym.def
\input amssym
\magnification=1200
\font\bigsl=cmsl10 scaled\magstep1
\font\bigsll=cmsl10 scaled\magstep3
\tolerance=500
\overfullrule=0pt
\hsize=6.40 true in
\hoffset=.10 true in
\voffset=0.1 true in
\vsize=8.70 true in
\null
\vskip 0.2 in
\centerline{\bigsll Points of integral canonical models of Shimura varieties}
\smallskip
\centerline{\bigsll of preabelian type, p-divisible groups, and applications}
\bigskip
\centerline{\sl Adrian Vasiu, U of Arizona}
\medskip
\centerline{submitted for publ. in J. Ast\'erisque on 4/9/01; up dated version 8/16/01}
\smallskip
\centerline{MSC 2000: Primary 11G10, 11G18, 11G25, 11G35, 14F30, 14G35, 14K10}
\medskip
\centerline{\it (The first part out of three)}
\bigskip
{\bigsl Contents of the first part}
\bigskip
\item{\S 1} Introduction \dotfill\ \ \ 5
\medskip
{1.0.} Quick overview \dotfill\ \ \ 5
\smallskip
{1.1.} Introducing Shimura (filtered) $\sg$-crystals \dotfill\ \ \ 6
\smallskip
{1.2.} The core: some basic results \dotfill\ \ \ 7
\smallskip
{1.3.} Standard Hodge situations \dotfill\ \ \ 19
\smallskip
{1.4.} Main geometric concepts \dotfill\ \ \ 21
\smallskip
{1.5.} Shimura (filtered) Lie $\sg$-crystals attached to points \dotfill\ \ \ 23
\smallskip
{1.6.} Stratifications \dotfill\ \ \ 24
\smallskip
{1.7.} On the contents of 4.3 and 4.6 \dotfill\ \ \ 30
\smallskip
{1.8.} Crystalline coordinates \dotfill\ \ \ 31
\smallskip
{1.9.} A Galois property \dotfill\ \ \ 32
\smallskip
{1.10.} On the passage from the Hodge type to the preabelian type \dotfill\ \ \ 33
\smallskip
{1.11.} Dense Hecke orbits and functorial aspects \dotfill\ \ \ 34
\smallskip
{1.12.} Some important applications \dotfill\ \ \ 35
\smallskip
{1.13.} The geometric picture of the case $p=2$ \dotfill\ \ \ 39
\smallskip
{1.14.} Afterthoughts, the Second Main Corollary and advices \dotfill\ \ \ 42
\smallskip
{1.15.} On the contents of \S5-14 \dotfill\ \ \ 45
\smallskip
{1.16.} Some history \dotfill\ \ \ 49
\medskip
\item{\S 2} Preliminaries \dotfill\ \ \ 51
\medskip
{2.1.} Notations, conventions and some definitions \dotfill\ \ \ 51
\smallskip
{2.2.} Fontaine categories and Shimura crystals \dotfill\ \ \ 60
\smallskip
{2.3.} The standard Hodge situation \dotfill\ \ \ 137
\smallskip
{2.4.} Some identifications \dotfill\ \ \ 179
\medskip
\item{\S 3} The basic results \dotfill\ \ 183
\medskip
{3.0.} The starting setting \dotfill\ \ \ 183
\smallskip
{3.1.} The first group of basic results \dotfill\ \ \ 184
\smallskip
{3.2.} The outline of the proof of 3.1.0 \dotfill\ \ \ 187
\smallskip
{3.3.} About the proof of 3.2.5 \dotfill\ \ \ 191
\smallskip
{3.4.} Split Shimura Lie $\sg$-crystals \dotfill\ \ \ 191
\smallskip
{3.5.} The proof of 3.4.6\dotfill\ \ \ 200
\smallskip
{3.6.} The second group of basic results: global deformations
\item{} $\,\,\,\,$ $\,\,$ $\,$ and some principles for different Fontaine categories \dotfill\ \ \ 209
\smallskip
{3.7.} The proofs of 3.2.6-7 \dotfill\ \ \ 314
\smallskip
{3.8.} The non quasi-split context \dotfill\ \ \ 315
\smallskip
{3.9.} Supplements to 3.1 and applications \dotfill\ \ \ 316
\smallskip
{3.10.} Terminology and formulas \dotfill\ \ \ 324
\smallskip
{3.11.} The case $k=\bar k$ and applications \dotfill\ \ \ 332
\smallskip
{3.12.} Shimura $F$-crystals and generic points \dotfill\ \ \ 342
\smallskip
{3.13.} Deviations of Shimura (filtered) (Lie) $\sg$-crystals \dotfill\ \ \ 342
\smallskip
{3.14.} The case $p=2$ \dotfill\ \ \ 373
\smallskip
{3.15.} Some conclusions\dotfill\ \ \ 380
\medskip
\item{\S 4} Applications of the basic results to integral canonical models of
Shimura
\item{}     varieties of preabelian type, to $p$-divisible groups, and to abelian varieties \dotfill\ \ 415
\medskip
{4.1.} The Shimura-ordinary type \dotfill\ \ \ 415
\smallskip
{4.2.} $G$-ordinary points \dotfill\ \ \ 422
\smallskip
{4.3.} Cyclic factors and (refined) Lie stable $p$-ranks \dotfill\ \ \ 426
\smallskip
{4.4.} $G$-canonical lifts \dotfill\ \ \ 433
\smallskip
{4.5.} Different stratifications \dotfill\ \ \ 444
\smallskip
{4.6.} Examples and main properties \dotfill\ \ \ 466
\smallskip
{4.7.} Crystalline coordinates \dotfill\ \ \ 482
\smallskip
{4.8.} The Galois property of $G$-ordinary points \dotfill\ \ \ 502
\smallskip
{4.9.} The passage from Shimura varieties of Hodge type 
\item{} $\,\,\,\,$ $\,\,$ $\,$ to Shimura varieties of preabelian type \dotfill\ \ \ 505
\smallskip
{4.10.} Hecke orbits of $G$-ordinary points \dotfill\ \ \ 520
\smallskip
{4.11.} The functorial behavior of $G$-ordinary points\dotfill\ \ \ 521
\smallskip
{4.12.} Integral and generalized Manin problems and applications \dotfill\ \ \ 525
\smallskip
{4.13.} The local form of the invariance principle \dotfill\ \ \ 551
\smallskip
{4.14.} Final remarks \dotfill\ \ \ 552
\medskip\noindent
$\,$ AE Addenda and errata to [Va2] \dotfill\ \ 569
\medskip\noindent
$\,$ Appendix \dotfill\ \ 583
\medskip\noindent
$\,$ References \dotfill\ \ 589
\medskip\medskip
{\bigsl Contents of the second part}
\bigskip
\item{\S 5} Shimura $p$-divisible groups, the proof of a conjecture 
\item{}     of Milne, and automorphic vector bundles...\dotfill\ \ 00
\medskip
\item{\S 6} Integral canonical models of Shimura varieties of preabelian type
\item{}     (resp. of Hodge type) with respect to primes dividing 3 (resp. with
\item{}     respect to primes dividing a prime $p\ge 3$ and non-hyperspecial subgroups)...\dotfill\ \ 00
\medskip
\item{\S 7} Twisted Shimura varieties \dotfill\ \ 00
\medskip
\item{\S 8} Shimura varieties and the Mumford--Tate conjecture, part III \dotfill\ \ 00
\medskip
\item{\S 9} The null strata and $G$-supersingular Shimura $\sg$-crystals \dotfill\ \ 00
\medskip
\item{\S 10} The number and dimensions of  strata of the standard list 
\item{}      of stratifications and Grothendieck's specialization categories \dotfill\ \ 00
\medskip\medskip
{\bigsl Contents of the third part}
\bigskip
\item{\S 11} An isogeny property \dotfill\ \ 00
\medskip
\item{\S 12} Points with values in finite fields \dotfill\ \ 00
\medskip
\item{\S 13} The Langlands--Rapoport conjecture \dotfill\ \ 00
\medskip
\item{\S 14} Crystalline applications to the proof of the conjecture \dotfill\ \ 00
\vfill\eject
\medskip
\centerline{}
\vfill\eject
\centerline{}
\bigskip
\bigskip
\centerline{\bigsll {\bf \S 1 Introduction}}
\bigskip\bigskip
\medskip
{\bf 1.0. Quick overview.} Recently we managed to prove the existence of integral canonical models of Shimura varieties of preabelian type with respect to primes $p\ge 3$ (cf. [Va2, 6.4.1 and defs. of 3.2.6]; see [Va2, 6.8] and either 4.6.4 or [Va3] for the case of Shimura varieties of preabelian type which are neither of compact type nor of abelian type, and see \S 6 and [Va5] for the case $p=3$). For the general theory of such integral models we refer to [Va2, \S 3]; for different types of Shimura varieties we refer to [Va2, 2.5]. With this paper, which is a natural continuation of [Va2], we start an extensive program of proving that many well known results for integral canonical models of Siegel modular varieties remain true (under proper formulation) for all integral canonical models of Shimura varieties of preabelian type. Here we are concerned with points of these models with values in regular, formally smooth schemes over rings of Witt vectors of a given length (finite or infinite) of a perfect field. Special attention is paid to points with values in (the Witt ring of) a perfect field. 
\smallskip
For the sake of generality and motivated by applications to many other geometric contexts, most of the things are stated in the context and language of Fontaine categories (see 1.2.4 and 1.2.9). So this paper (together with the ones continuing it) represents:
\medskip
{\bf F1.} {\it A foundation for integral aspects of Shimura varieties, and}
\smallskip
{\bf F2.} {\it A foundation for Fontaine categories.} 
\medskip
For a faster publication of our main results pertaining to the program, we split them in three parts, to be published separately. The first part is the present paper and so contains \S 1-4, the second part will contain \S 5-10 and the third part will contain \S 11-14. The results mentioned in [Va2, 1.6-7] are proved in \S2-14.
In \S1-4 we concentrate (this was the initial goal) in reobtaining in a new manner: 
\medskip
\item{$\bullet$}
the theory (due to Serre--Tate; see [Me, Appendix] for its first features) of ordinary abelian varieties and $p$-divisible groups (over perfect fields of characteristic $p$) and
of their canonical lifts, 
\smallskip
\item{$\bullet$}
the theory (due to Serre--Tate and Dwork; see [Ka1] and [Ka3-4]) of canonical coordinates of the formal moduli scheme of deformations of an ordinary abelian variety over an algebraically closed field, 
\smallskip
\item{$\bullet$}
a proof of Manin's problem (first solved in [Ta, p. 98]),
\medskip\noindent
and in extending these theories to the geometric context of integral canonical models of Shimura varieties of preabelian type and even more generally, to the abstract context (see 1.1) of Shimura $\sg$-crystals. The case $p=2$ is treated in maximal generality (like no restriction in the abstract context and, in the geometric context, we logically restrict to integral canonical models of Shimura varieties of abelian type in mixed characteristic $(0,2)$ whose existence is proved here).  
\smallskip
On the way we obtain many other results pertaining to $p$-divisible groups, to abelian varieties, to the  mentioned integral canonical models and to Fontaine categories. 
\smallskip
With $p\ge 2$, two things are worth being mentioned from the very beginning. First, we show that the number of isomorphism classes of $p$-divisible groups of a given rank over an algebraically closed field $\bar k$ of characteristic $p$ which are definable over a fixed finite field $\FF_{p^q}$ (with $p^q$ elements; here $q\in\NN$), is finite and does not depend on $\bar k$ (see 1.6.4). Second, we construct versal global deformations, i.e. versal deformations over the $p$-adic completion of (a particular type of) pro-\'etale schemes over smooth, affine schemes over Witt rings of perfect fields (of characteristic $p$), of $p$-divisible groups over such fields, cf. Theorem 1 of 3.15.1. Theorem 1 of 3.15.1 complements:
\medskip
-- Grothendieck's theorem of [Il, 4.4] (we get stronger forms than what the combination of loc. cit. and Artin's approximation theorem can produce);
\smallskip
-- Raynaud's theorem of [BBM, 3.1.1] (it is effective, i.e. it has --cf. its proof-- an explicit, concrete aspect and it offers in many regular contexts --see the application of its proof to the Theorem of 3.15.2-- more information on embeddings of finite, flat, commutative group schemes of $p$-power order into $p$-divisible groups);
\smallskip
-- Faltings's theorem of [Fa1, 7.1] (it is effective and works for $p=2$ as well);
\smallskip
-- de Jong's theorem of [dJ2, first main result] (it is effective; moreover it handles the case of truncations as well: for instance, see its application in the proof of Theorem of 3.15.2).
\medskip
Over algebraically closed fields of characteristic $p$, the words ``$p$-adic completion" and ``pro-\'etale" can be dropped (cf. Theorem 13 of 1.12). 
\medskip\smallskip
{\bf 1.1. Introducing Shimura (filtered) $\sg$-crystals.} In what follows $k$ is a perfect field of positive characteristic $p$, $W(k)$ is the ring of (infinite) Witt vectors of $k$ and $\sg$ is the Frobenius automorphism of $W(k)$ or of $B(k):=W(k)[{1\over p}]$. We denote by $\bar E$ the algebraic closure of a field $E$. The proof of the basic result 5.1 of [Va2] leads to the study of $\sg$-crystals over $k$ endowed with an extra structure (cf. [Va2, 5.6.5]).
Moving from the geometric context of loc. cit. to a more general and abstract context, we start considering a $\sg$-crystal $(M,\vph)$ over $k$ (so $M$ is a free $W(k)$-module of finite rank and $\vph$ is a $\sg$-linear endomorphism of $M$ producing a $\sg$-linear automorphism of $M[{1\over p}]$) and a quasi-split reductive subgroup $G$ of $GL(M)$ satisfying the following axioms:
\medskip
\item{\bf i)} There is a family of tensors
$(t_\al)_{\al\in\Mj}$ in spaces of the form 
$M^{\otimes n}\otimes_{W(k)} M^{\ast\otimes m}[{1\over p}]$ (with $M^\ast$ the dual of $M$ and with $n,m\in\NN\cup \{0\}$), such that $\vph(t_{\al})=t_{\al}$, $\forall\al\in\Mj$, and $G$ is the Zariski closure of the subgroup of $GL(M[{1\over p}])$ fixing $t_\al$, $\forall\al\in\Mj$, in $GL(M)$;
\smallskip
\item{\bf ii)} There is a cocharacter $\mu:\GG _m\to G$ producing a direct sum decomposition
$M=F^1\oplus F^0$ such that $\be\in\GG_m(W(k))$ acts through $\mu$ on $F^i$ as the multiplication with $\be^{-i}$,
$i=\overline{0,1}$. Moreover $F^1$ is a proper direct summand of $M$ (this last condition is inserted just to avoid trivial cases, cf. 2.2.9 11)); 
\smallskip
\item{\bf iii)} The triple $(M,F^1,\vph)$ is a filtered $\sg$-crystal (i.e. we have $\vph({1\over p}F^1+M)=M$).
\medskip
The triple $(M,\vph,G)$ (resp. quadruple $(M,F^1,\vph,G)$) is called a Shimura (resp. a Shimura filtered) $\sg$-crystal. The pair $({\rm Lie}(G),\vph)$ is called the Shimura Lie $\sg$-crystal attached to $(M,\vph,G)$. Here ${\rm Lie}(G)$ is the Lie algebra of $G$ and we still denote by $\vph$ the resulting $\sg$-linear Lie automorphism of ${\rm Lie}(G)[{1\over p}]$ obtained via the logical inclusion ${\rm Lie}(G)\subset {\rm End}(M)$ and the natural action of $\vph$ on ${\rm End}(M)[{1\over p}]$ that takes $m\in {\rm End}(M)[{1\over p}]$ into $\vph\circ m\circ\vph^{-1}\in {\rm End}(M)[{1\over p}]$. Axiom i) guarantees that $\vph$ normalizes ${\rm Lie}(G)[{1\over p}]$; we have $\vph(p{\rm Lie}(G))\subset {\rm Lie}(G)$, cf. axiom iii). See 2.1 for the (standard) way $\vph$ acts on the spaces of i).
\smallskip
Fixing the group $G$, the family of tensors $(t_\al)_{\al\in\Mj}$ and the $G(W(k))$-conjugacy class $[\mu]$ of the cocharacter $\mu$, every $\sg$-linear endomorphism $\vph_1$ of $M$ corresponding to any other $\sg$-crystal $(M,\vph_1)$ satisfying the logical analogues of axioms i) to iii), up to $G(W(k))$-conjugation, can be put in the form $\vph_1=g\vph$, with $g\in G(W(k))$ (cf. the second paragraph of \S 3). 
\smallskip
If $G=GL(M)$, then we can take the set $\Mj$ of indices to be the empty set and so we regain the classical context of Serre--Tate's (ordinary) theory.
\medskip\smallskip
{\bf 1.2. The core: some basic results.} The case $p=2$ is treated as asides. So, in all that follows, without a special reference we assume $p>2$. Our basic results are grouped in 3.1 and 3.6. 
\medskip
{\bf 1.2.0. On the first group of basic results.} In 3.1 we present the first group of basic results. For our (standard) conventions on Newton polygons see 2.1. We have: 
\medskip
{\bf Theorem 1.} {\it {\bf a)} Among all  Shimura $\sg$-crystals $(M,g\vph,G)$ (with
$g\in G(W(k))$), the ones which have the smallest Newton polygon (in the sense that all others have a Newton polygon strictly above it) are precisely the ones $(M,g_1\vph,G)$ (called $G$-ordinary $\sg$-crystals) for which there is a cocharacter $\mu:\GG_m\to G$ (as in the above axioms), such that 
the parabolic Lie subalgebra ${\got p}_1$ of ${\rm Lie}(G)$ corresponding to non-negative slopes of the Shimura Lie $\sg$-crystal $\bigl({\rm Lie}(G),g_1\vph)$ (see 2.2.3 3)) is contained in the Lie algebra $F^0({\rm Lie}(G))$ of the parabolic subgroup of $G$ normalizing the direct summand $F^1$ of $M$ defined (as in 1.1 ii)) by $\mu$. 
\smallskip
{\bf b)} If $(M,g_1\vph,G)$ is a $G$-ordinary $\sg$-crystal, then there is a unique Shimura filtered $\sg$-crystal $(M,F^1,g_1\vph,G)$ (called the $G$-canonical lift of $(M,g_1\vph,G)$) such that for each $x\in {\got p}_1$, we have $x(F^1)\subset F^1$.
\smallskip
{\bf c)} Among all Shimura $\sg$-crystals $(M,g\vph,G)$ (with $g\in G(W(k))$), the ones which have the smallest Newton polygon of their attached Shimura Lie $\sg$-crystals $({\rm Lie}(G),g\vph)$, are precisely the $G$-ordinary $\sg$-crystals.
\smallskip
{\bf d)} Any Shimura $\sg$-crystal $(M,g\vph,G)$ is the specialization of a $G_{W(k_1)}$-ordinary $\sg_1$-crystal $(M\otimes_{W(k)} W(k_1),g_1(\vph\otimes\sg_1),G_{W(k_1)})$, with $g_1\in G(W(k_1))$ (here $\sg_1$ is the Frobenius automorphism of the Witt ring of a perfect field $k_1$ containing $k$; we can take $k_1$ to be the field of fractions of the perfection of $k[[X]]$, where $X$ is an independent variable).}
\medskip
Referring to a), the quadruple 
$$({\rm Lie}(G),g_1\vph,F^0({\rm Lie}(G)),F^1({\rm Lie}(G))),$$ 
with $F^1({\rm Lie}(G))$ as the unipotent radical of $F^0({\rm Lie}(G))$, is called the Shimura filtered Lie $\sg$-crystal attached to $(M,F^1,\vph,G)$ and the inclusion ${\got p}_1\subset F^0({\rm Lie}(G))$ is interpreted as (cf. def. 2.2.12 a)): it is of parabolic type. See 3.1.1.2 for the interpretation of the mentioned inclusion (resp. of the uniqueness part of b)) in terms of (see def. of 2.2.8 7) and of 2.2.9 6)) $G$-endomorphisms (resp. of automorphisms); so 3.1.1.2 represents the generalization of the classical result that any endomorphism of an ordinary $p$-divisible group over $k$ lifts to its canonical lift.
\medskip
{\bf 1.2.1. On the proof of Theorem 1.} The proof of Theorem 1 is achieved by first (an idea suggested by G. Faltings) proving it for 
Shimura Lie $\sg$-crystals (${\rm Lie}(G),g\vph$). To achieve this and to complete the proof of Theorem 1 we construct (see 3.6.7 for the exact form needed) global deformations of such (as in 1.1) Shimura filtered $\sg$-crystals; their goal is: to connect any
Shimura $\sg$-crystal $(M,g\vph,G)$ (with $g\in G(W(k))$) with (the extension to $\bar k$ of) a $G$-ordinary $\sg$-crystal $(M,g_1\vph,G)$. Here by global we mean: deformations over the $p$-adic completion of pro-\'etale schemes over smooth, affine $W(k)$-schemes, whose special fibres are connected and have a dense set of $\bar k$-valued points. This takes care of (a stronger form of) d) of Theorem 1, as it can be checked easily starting from Grothendieck--Katz specialization theorem (see [Ka2, 2.3.1-2]) and from 1.2.1.0 2) and 3) below.
\medskip
{\bf 1.2.1.0. Tools from the Lie context.} The proof of Theorem 1 for Shimura Lie $\sg$-crystals $({\rm Lie}(G),g\vph)$ (with $g\in G(W(k))$), means three things:
\medskip\noindent
\item{\bf 1)} there is $g_1\in G(W(k))$ such that for any $x\in {\got p}_1$ (defined in a) of Theorem 1) we have $x(F^1)\subset F^1$ (cf. 3.2.3);
\smallskip\noindent
\item{\bf 2)} for any $g\in G(W(k))$, the Newton polygon of $({\rm Lie}(G),g\vph)$ is strictly above or equal to the Newton polygon of $({\rm Lie}(G),g_1\vph)$; 
\smallskip\noindent
\item{\bf 3)}  if in 2)  we have equality of Newton polygons, then an explicit expression of $g$ can be determined and the  $\sg$-crystal $(M,g\vph)$ has the same Newton polygon as the $\sg$-crystal $(M,g_1\vph)$ (cf. 3.4.11 and 3.4.13).
\medskip 
Part 1) is easy (as $G$ is a quasi-split reductive group), while 3) is implicitly obtained on the way --see below-- of proving 2) (cf. 3.4.8-11 and 3.7.4). From 1) to 3) and 1.2.1 we get  a) and c) of Theorem 1 (see 3.2.4). Based on the fact that $[\mu]$ is uniquely determined by $(M,\vph,G)$ (cf. Fact 2 of 2.2.9 3)), the uniqueness part of b) of Theorem 1 is just a simple property of parabolic subgroups of $G$ (cf. 3.2.8).
\smallskip
For the proof of 2) we first show that for any $g\in G(W(k))$, the Newton polygon of $({\rm Lie}(G),g\vph)$ can not be strictly below the Newton polygon of $({\rm Lie}(G),g_1\vph)$ (cf. 3.4.12). This is achieved (cf. the proof of 3.4.6 needed to get 3.4.8, and so needed to get 3.4.12) by using Lie stable $p$-ranks of Shimura $\sg$-crystals, via their interpretation (see 3.4.5.1 B; see also the proof of 3.9.2) in terms of multiplicities of (suitable --i.e. minimal in some relative sense-- non-positive) slopes of Newton polygons of attached Shimura Lie $\sg$-crystals. 
\smallskip
The Lie stable $p$-rank of a Shimura $\sg$-crystal $(M,\vph,G)$ is a non-negative integer associated to $({\rm Lie}(G),\vph)$: though it is defined just in 3.9.1 and 3.9.1.1, the very essence of its definition is already captured in 3.4.5 and 3.4.5.1. Here is a sample of how it can be computed. If the adjoint group $G^{\rm ad}$ of $G$ is an absolutely simple group and $\mu$ of 1.1 ii) does not factor through the center of $G$, then it is the dimension over $k$ of the stable image (i.e. of the image of the ``sufficiently high" iterate) of the $\sg$-linear endomorphism $FSHW$ of ${\rm Lie}(G_k)={\rm Lie}(G)/p{\rm Lie}(G)$ defined by the reduction mod $p$ of the $\sg$-linear endomorphism $p\vph$ of ${\rm Lie}(G)$ (warning: here $\vph$ is viewed as a $\sg$-linear automorphism of ${\rm End}(M)[{1\over p}]$ and so $p\vph$ denotes its multiplication by $p$). As $p\vph$ normalizes the Lie algebras of the derived group $G^{\rm der}$ of $G$ as well as of the adjoint group $G^{\rm ad}$ of $G$ (see 2.2.13 for the adjoint context), this Lie stable $p$-rank can be computed as well (see 3.9.6) using the resulting $\sg$-linear endomorphism $FSHW^{\rm der}$ (resp. $FSHW^{\rm ad}$) of ${\rm Lie}(G^{\rm der})/p{\rm Lie}(G^{\rm der})$ (resp. of ${\rm Lie}(G^{\rm ad})/p{\rm Lie}(G^{\rm ad})$). 
\smallskip
3.4.6 is a theorem on such Lie stable $p$-ranks capturing the very essence behind 2), 3) and Theorem 2 below. Its proof is carried out in 3.5. It involves standard techniques pertaining to (Lie algebras of) parabolic subgroups of $G_{\bar k}$: the well known classification of Shimura varieties of Hodge type (see [De2]; see also [Sa] or the abstract context of [Se2, \S3]), allows us to restate 3.4.6 in terms of
\medskip
-- possible dimensions of suitable intersections of Lie algebras of unipotent radicals of such parabolic subgroups, and of 
\smallskip
-- suitable inclusions between such intersections. 
\medskip\noindent
Accordingly, the language of root systems occupies a central stage in 3.4-5.
\smallskip
To conclude the proof of 2) we use standard specialization arguments in the context of the mentioned global deformations (cf. 1.2.1; see 3.7.1 and 3.7.6).
\medskip
{\bf 1.2.1.1. Variants and complements.} For a simpler way of proving a) and c) of Theorem 1 see 3.4.14. See also the last sentence of 1.4.1 for a third approach. For the sake of completeness and of future references, 3.4-5 handle as well with full details the cases of Shimura Lie $\sg$-crystals (see def. 2.2.11) which are not attachable to some Shimura $\sg$-crystals (like the ones involving the $E_6$ and $E_7$ Lie types, etc.).
\smallskip
The first group of basic results and their proofs are presented in 3.1-5 and 3.7, while 3.8 shows that Theorem 1 remains true even if $G$ is just a reductive group (i.e. it is not necessarily quasi-split). The construction of the global deformations mentioned in 1.2.1 is carried out in 3.6, being in fact part of the second group of basic results.
\medskip
{\bf 1.2.1.2. On the organization of 2.2.} In 2.2 we present basic facts pertaining to Shimura (adjoint) (filtered) (Lie) ($\sg$- or $F$-) crystals; most of them are indispensable for the reading of \S3-4 (and so of 3.1-5). Special attention is paid to the Lie context (see the first paragraph of 2.2.3 4) for a justification): see 2.2.2-3, 2.2.4 B, 2.2.8 b), 2.2.11-13, 2.2.16.5, 2.2.23 B, etc. The notion of Shimura $p$-divisible groups (these are $p$-divisible groups endowed with a family of crystalline --de Rham-- tensors and satisfying an axiom inspired from the axioms of 1.1) is introduced in 2.2.20. See 2.2.22 (resp. 2.2.23) for some complements on Shimura (Lie) $\sg$-crystals (resp. on quasi-polarizations); in particular, 2.2.22 3) deals with different classification problems. 
\medskip
{\bf 1.2.2. A supplement to 1.2.0.} We have (cf. 3.9.2) the following important supplement to 1.2.0:
\medskip
{\bf Theorem 2.} {\it A Shimura $\sg$-crystal $(M,\vph,G)$ is a $G$-ordinary $\sg$-crystal iff for any $g\in G(W(k))$, the Lie stable $p$-rank of $(M,g\vph,G)$ is smaller than or equal to the Lie stable $p$-rank of $(M,\vph,G)$.} 
\medskip
In particular, any Shimura $\sg$-crystal $(M,g\vph, G)$ having the same Lie stable $p$-rank as a $G$-ordinary $\sg$-crystal $(M,g_1\vph,G)$, is a $G$-ordinary $\sg$-crystal itself. Theorem 2 represents the generalization of the characterization of ordinary $p$-divisible groups over $k$ in terms of $p$-ranks of Hasse--Witt maps or matrices: we get Faltings--Shimura--Hasse--Witt (reduced or adjoint) maps and matrices ``encoding" the values of these Lie stable $p$-ranks (cf. 3.4.5, 3.9.1, 3.9.6, 4.3.7 and 4.3.9; see also 3.13.7.6 for their non-compact variant). The essence of its proof is 3.5. The $\sg$-linear endomorphism $FSHW^{\rm der}$ (resp. $FSHW^{\rm ad}$) of 1.2.1.0, is the most elementary sample of a Faltings--Shimura--Hasse--Witt map (resp. adjoint map); the dimension over $k$ of its stable image is nothing else but the multiplicity of the slope $-1$ of the isocrystal $({\rm Lie}(G)[{1\over p}],\vph)$. Such maps can be defined as well in global contexts (like when $k$ gets replaced by some regular, formally smooth $k$-scheme $S_k$ and we deal with a suitable $F$-crystal with tensors --see 2.2.10 for samples-- over $S_k$). 
\smallskip
We can refer to Lie stable $p$-ranks as well as Faltings--Shimura--Hasse--Witt invariants. As we speak about Faltings--Shimura--Hasse--Witt maps as well as (see 1.6.1) stratifications, and as in more general contexts (like of the class $\Mc\Ml_3$ of 1.2.9), there is a demarcation between the two (see 3.13.7.9 for some explanations), to be short and to avoid confusion, in \S1-14, in connection to Shimura Lie $\sg$-crystals we refer just to Lie stable $p$-ranks of them. 
\medskip
{\bf 1.2.2.1. A mod $p$ interpretation.} Theorem 2 implies: the fact that $(M,g\vph,G)$ is or is not $G$-ordinary depends only on the expression of $g$ mod $p$ (cf. 3.9.3). See 3.13.7.1.2 for a more of a principle formulation of this, closer in spirit to the classical context of Serre--Tate's (ordinary) theory; so 3.13.7.1.2 represents the generalization of the following well known result: the fact that a $p$-divisible group $D$ over $k$ is ordinary or not can be read out from its maximal subgroup $D[p]$ annihilated by $p$.
\medskip
{\bf 1.2.3. The list of principles.} 
In section 3.6 we start the presentation of twelve fundamental principles of the crystalline theory. These principles are: 
\medskip
{\bf PR1} the $\nabla$ principle (cf. 3.6.1.3, 3.6.18.4 and 3.6.18.6 a); see also 3.6.18.4.4 1) and 3.6.18.9); 
\smallskip
{\bf PR2} the moduli principle (cf. 3.6.1.3 1) and 3.6.18.4.2); 
\smallskip
{\bf PR3} the surjectivity principle (cf. 3.6.1.3 5), 3.6.8.1.2 a), 3.6.8.2, 3.6.18.4.5 and 3.6.18.6 b));
\smallskip
{\bf PR4} the uniqueness principle (cf. 3.6.1.2, 3.6.18.8 and 3.6.18.8.1); 
\smallskip
{\bf PR5} the (very) (weak) gluing principle (cf. 3.6.2.2 and 3.6.19); 
\smallskip
{\bf PR6} the integrability principle (cf. 3.6.18.4.1 and 3.6.18.4.1.1); 
\smallskip
{\bf PR7} the deformation principle (cf. 3.6.14 and 3.6.14.1-4);
\smallskip
{\bf PR8} the boundedness principle (cf. 3.15.7; see also 1.6.3);
\smallskip
{\bf PR9} the purity principle (cf. 3.6.8.1.4, 3.15.10 and 4.5.16.1);
\smallskip
{\bf PR10} the (local) integral Manin problems (cf. a great part of 4.12);
\smallskip  
{\bf PR11} the invariance principle (cf. 4.13 and \S 9-10); 
\smallskip
{\bf PR12} the slice principle (cf. \S 7).
\medskip
These principles are accompanied by eight other useful properties:
\medskip
$\bullet$ the constructibility property (cf. 3.6.8.1.2 b));
\smallskip
$\bullet$ the constancy property (cf. 3.6.8.9); 
\smallskip
$\bullet$ the touching  property (cf. 3.6.18.4.3); 
\smallskip
$\bullet$ the inducing property (cf. 3.6.18.5); 
\smallskip
$\bullet$ the liftability property (cf. 3.6.18.5.2, 3.6.18.5.4 1) and 3.6.18.5.9); 
\smallskip
$\bullet$ the inductive property (cf. 3.9.7.2);
\smallskip
$\bullet$ the homomorphism property (cf. Corollary of 3.15.8);
\smallskip
$\bullet$ the rigidity property (cf. \S 7; see also 3.6.1.4 5) and 3.13.5.4). 
\medskip
3.6.18.4-6 can be extended to a relative situation (i.e. to a situation where a smooth group --not necessarily reductive-- is involved, as in the context of 1.1): see the whole of 3.6.18.7. It would be too long to fully describe these principles (or properties) here. So we concentrate in outlying some of the main features.
\medskip
{\bf 1.2.4. The first three principles in a convenient Shimura $\sg$-crystal context.} For the definition of Fontaine categories $\Mm\Mf_{[0,1]}(X)$ and $\Mm\Mf_{[0,1]}^\nabla(X)$, of their $p$-divisible objects and for how ($p$-divisible) objects are pulled back through pro-\'etale, affine morphisms see 2.1, 2.2.1 c) and respectively 3.6.1.1.1 1), 2) and 5). Here $X$ stands for a regular, formally smooth, affine $W(k)$-scheme such that either itself or its $p$-adic completion is endowed with an arbitrary (but fixed) Frobenius lift compatible (see 2.1) with $\sg$; we always assume that the sheaf $\Om_{X_k/k}$ on $X_k$ of relative differentials is locally free of finite rank. Though in this introduction we deal just with such an $X$, it is worth pointing out that in many situations $X$ can be non-affine. For instance, when we are dealing with $\Mm\Mf_{[0,1]}^\nabla(X)$, $X$ can be non-affine and moreover we do not need to assume that a Frobenius lift is specified (cf. 2.2.1 c)). 
\smallskip
We come back to 1.1. Let $U={\rm Spec}(R)$ be an open, affine subscheme of $G$ such that the identity element $a_0$ of $G(W(k))$ belongs to $U(W(k)$). We assume the existence of an \'etale morphism $U\to {\rm Spec}(W(k)[x_1,...,x_{d(G)}])$, where $d(G)$ is the relative dimension of $G$, such that $a_0$ factors through the closed subscheme of $U$ defined by $x_i=0$, $i=\overline{1,d(G)}$; here we identify naturally each $x_i$ with an element of $R$. Let $R^\wedge$ be the $p$-adic completion of $R$ and let $\Phi_R$ be the Frobenius lift of $U^\wedge:={\rm Spec}(R^\wedge)$ which at the level of rings takes $x_i$ (identified naturally with an element of $R$) into $x_i^p$, $\forall i\in\{1,...,d(G)\}$, and is compatible with $\sg$ (for the case of a general type Frobenius lift $\Phi_R$ compatible with $\sg$ see 3.6.18.8.3). Let $M_U$ be the $p$-divisible object of $\Mm\Mf_{[0,1]}(U)$ defined by the triple
$$
(M\otimes_{W(k)} R^\wedge,F^1\otimes_{W(k)} R^\wedge,\Phi_{M_U}).
$$ 
Here $\Phi_{M_U}$ is the $\Phi_R$-linear endomorphism of $M\otimes_{W(k)} R^\wedge$ which takes an element $m\in M$ into $u(\vph(m))$, with $u\in G(R^\wedge)$ as the universal element naturally defined by the inclusion $U\hookrightarrow G$. So $\Phi_{M_U}(t_{\al})=t_{\al}$, $\forall\al\in\Mj$. $m$ and $t_{\alpha}$'s are viewed naturally as tensors of the tensor algebra of $(M\oplus M^*)\otimes_{W(k)} R^\wedge[{1\over p}]$ and $\Phi_{M_U}$ acts in the standard way (see 2.1) on this tensor algebra.  
\smallskip
The $\nabla$, the moduli and the surjectivity  principles in their simplified and convenient form involving Shimura $\sg$-crystals can be formulated as follows (cf. 3.6.1.3):
\medskip
{\bf Theorem 3.} {\it There is a pro-\'etale, affine morphism $\ell:U_1={\rm Spec}(R_1)\to U$, with the special fibre of $U_1$ a geometrically connected $k$-scheme, such that:}
\medskip\noindent
\item{a)} {\it The pull back $M_{U_1}$ of $M_U$ to (the $p$-adic completion of) $U_1$, is a $p$-divisible object of $\Mm\Mf_{[0,1]}^\nabla(U_1)$. In other words (cf. 3.6.1.1.1 2)), there is a unique connection 
$$
\nabla:M\otimes_{W(k)} R_1^\wedge\to M\otimes_{W(k)}\Om_{R/W(k)}\otimes_R R_1^\wedge
$$ 
for which the Frobenius endomorphism $\Phi_{M_{U_1}}$ of $M\otimes_{W(k)} R_1^\wedge$ (defined by $\Phi_{M_U}$; the Frobenius lift of $R_1^\wedge$ being the one naturally induced from $\Phi_R$ via $\ell$) is $\nabla$-parallel (in the usual sense); $\nabla$ is integrable and nilpotent mod $p$. Moreover, $\nabla$ respects the $G$-action (i.e. we have $\nabla(t_{\al})=0$, $\forall\al\in\Mj$);}
\smallskip\noindent
\item{b)} {\it There is a closed subscheme $\Mb$ of the special fibre $U_k$ of $U$ not containing the origin of $G_k$, with the property that the fibres of $\ell$ above points of $U_k\setminus\Mb$ are non-empty.}
\medskip
From Theorem 3 and [Fa1, 7.1] we deduce the existence of a $p$-divisible group $D_1$ over ${\rm Spec}(R_1^\wedge)$ corresponding (via the antiequivalence of loc. cit. applied in the context of a suitable projective system) to this $p$-divisible object of $\Mm\Mf_{[0,1]}^\nabla(U_1)$ (cf. 3.6.2.0; see also b) of the Fact of 3.6.19). As $\nabla$ annihilates $t_{\al}$, $\forall\alpha\in\Mj$, the pair $\Md_1:=(D_1,(t_{\al})_{\al\in\Mj})$ is a Shimura $p$-divisible group.
\smallskip 
In their full form, the above three principles in the context of $M_U$ and $\Phi_R$ describe as well how we can choose $\ell$ to get an universal property and what extra properties it enjoys (subject to this special choice). 
\medskip
{\bf 1.2.5. The essence of the second group of basic results: the first six principles.} For the sake of generality, we now describe PR1-6 in wider contexts (then the one of 1.2.4). Let $X$ be as in 1.2.4.
The $\nabla$ principle for an arbitrary object ${\got C}$ of $\Mm\Mf_{[0,1]}(X)$ shows the existence, in the \'etale topology of $X$ (i.e. when we pass to an \'etale, affine morphism $X^\prime\to X$, surjective mod $p$), of connections on the pull back ${\got C}_{X^\prime}$ of ${\got C}$ to $X^\prime$; such connections satisfy all equations except the ones referring to the integrability part (see 3.6.1.1.1 1), 2) and 4)) needed to be satisfied in order that ${\got C}_{X^\prime}$ is obtainable from an object of $\Mm\Mf_{[0,1]}^{\nabla}(X^\prime)$ by forgetting its connection. It also computes the number $\nabla({\got C})$ of connections on ${\got C}$ (in the case when the special fibre of $X$ is a local strictly henselian scheme), cf. 3.6.18.4.4 1), and expresses upper bounds for this number (in the general case), cf. 3.6.18.9. 
\smallskip
The $\nabla$ principle for a $p$-divisible object ${\got C}$ of $\Mm\Mf_{[0,1]}(X)$ shows the existence of such connections in the pro-\'etale topology of $X$ (i.e. when we pass to a pro-\'etale, affine morphism $X^{\prime}\to X$, surjective mod $p$) and (under some conditions on $X$; see b) of Fact of 3.6.19 and see 3.15.4 for some extensions of it) asserts (via the integrability principle to be mentioned below) the existence of a $p$-divisible group $D_{X^\prime}$ over the $p$-adic completion $X^{\prime\wedge}$ of $X^\prime$ whose associated filtered $F$-crystal, when viewed just as a $p$-divisible object of $\Mm\Mf_{[0,1]}(X^\prime)$ (i.e. when viewed without connection), is the pull back of ${\got C}$ to $X^\prime$; to be consistent with 1.2.4, $X^{\prime\wedge}$ is viewed here as an affine scheme. In many cases the uniqueness of $D_{X^\prime}$ is implied by [Fa1, 7.1] (see b) of Fact of 3.6.19); when loc. cit. does not apply, we rely on the fully faithfulness results of [BM, \S 4] and the deformation theory of [Me, ch. 4-5] in order to get more general uniqueness statements (see 3.15.2-4).
\smallskip
Regardless of ${\got C}$ being an object or a $p$-divisible object, the $\nabla$ principle also describes how such an $X$-scheme $X^\prime$ can be obtained: always $X^\prime_k$ is an affine $X_k$-scheme of a very particular type. To detail this we work locally in the Zariski topology of $X_k$. So we consider the context in which ${\got C}$ is a $p$-divisible object, $\Om_{X_k/k}$ is a free sheaf on $X_k$ and the underlying sheaf of ${\got C}$ (see defs. of 2.2.1 c)) is as well free. We have (cf. 3.6.18.4.2 b) for the mod $p$ part):
\medskip
{\bf Key Property.} {\it $X^\prime_k$ is the projective limit of schemes $X_m$, $m\in\NN$, where $X_m$ is obtained from $X_{m-1}$ (with $X_0:=X_k$) by using a quasi Artin--Schreier system of equations (in a number of variables which does not depend on $m$), while $X^\prime$ itself is obtained as a projective limit of schemes by just lifting what we get mod $p$.} 
\medskip
These systems of equation (cf. the terminology of 3.6.18.4.6 A; they are also called of first type, cf. 3.6.8.9) are systems of $n$ equations in $n$ variables $x_1$,..., $x_n$ with coefficients in an $\FF_p$-algebra $R(p)$, of the form
$$x_i=L_i(x_1^p,...,x_n^p),\leqno (1)$$
$i\in\{1,...,n\}$, where $L_i$'s are linear forms (not necessarily homogeneous). Here $n\in\NN$. Any such system of equations defines naturally an \'etale, affine ${\rm Spec}(R(p))$-scheme (see 3.6.8.1.2 a)).
\smallskip
The uniqueness principle tells us when the  number $\nabla({\got C})$ is 1. There are two forms of it: one in terms only of the Frobenius lift of $X$ (or of its $p$-adic completion) (see 3.6.18.8) and another one in terms only of the Frobenius endomorphism of (the underlying sheaf of) ${\got C}$ (see 3.6.18.8.1 a)). The integrability principle asserts that all connections, whose existence is asserted by the $\nabla$ principle, are in fact integrable. The ideas involved in obtaining it rely on a natural algebraization process (which allows us --via Chinese Remainder Theorem-- to perform suitable global deformations to a generic --i.e. ordinary-- context) and on the touching property; they can be easily grasped from 3.6.18.4.1. These last two principles together with the inducing property of 3.6.18.5 form a natural extension of [Fa2, th. 10].  
\smallskip
The moduli principle shows the existence of a moduli formal scheme of connections (satisfying the equations mentioned in the first paragraph of 1.2.5) on (pull backs of) an arbitrary object or $p$-divisible object ${\got C}$ of $\Mm\Mf_{[0,1]}(X)$. One form of the surjectivity principle asserts that, in the case of a $p$-divisible object ${\got C}$, any connected component of the special fibre of a such resulting moduli formal scheme, maps onto an open, dense subscheme of a connected component of the special fibre of $X$ (cf. 3.6.18.4.5). 
\smallskip
The gluing principle has many variants (cf. 3.6.19). One such variant shows that in many cases two filtered $F$-crystals (which are ``modeled" on Shimura $p$-divisible groups) over the (assumed to be connected) special fibre $Y_k$ of a separated, smooth $W(k)$-scheme $Y$, under some conditions which make them to look similar ``around" a $k$-valued point $y$ of $Y$ and modulo logical changes of filtrations, are isomorphic over a connected $Y_k$-scheme which is the special fibre of a pro-\'etale, affine morphism $Y_1\to Y$ whose image contains an open, dense subscheme of $Y_k$ and to which $y$ lifts. 
\medskip
{\bf 1.2.5.1. Some ideas.} The guiding main new idea for obtaining these principles: we work inductively mod $p^n$, $n\in\NN$ (even in the case of a regular ring of formal power series over $W(k)$). What we gain by doing this: the systems of equations we get (as in 1.2.5, by working locally in the Zariski topology of $X_k$) are of a very simple and practical form (for instance, see 3.6.8; see also the simple, independent Fundamental Lemma 3.6.8.1, whose proof is based on intersection theory). 
\smallskip
Using a subclass of quasi Artin--Schreier systems of equations (i.e. the class of the so called Artin--Schreier systems of equations) with coefficients in $\FF_p$-algebras, we define Artin--Schreier fundamental groups (of different levels) of $\FF_p$-algebras whose spectra are connected (see 3.6.18.4.6 A and B). See 3.6.18.4.6 C and D for the Fundamental Lemma of the $\pi_1$-theory in positive characteristic. 
\smallskip
The idea behind the surjectivity principle is the constancy property: the quasi Artin--Schreier systems of equations we get working with different powers of $p$ are very close (in shape) one to another, i.e. they differ, cf. 3.6.8.9, only by adding some constants to the linear forms involved (so the abstract surjectivity principle of 3.6.1.8.2 a) applies). For an application of 3.6.8.1 to difference systems of equations see 3.6.8.1.1. We view 3.6.8.1.2 c) and part of 3.6.8.1 as the purely algebraic analogue of the well known theorem of S. Lang on connected algebraic groups over finite fields.
\medskip
{\bf 1.2.6. Some details on the contents of 3.9-10 and 3.12-14.} In 3.9.4-5 we present the notion of refined Lie stable $p$-rank of a Shimura $\sg$-crystal $(M,\vph,G)$ and a theorem for it similar to Theorem 2. The idea behind this is: for each factor of the adjoint group $G^{\rm ad}$ of $G$ whose generic fibre has a Lie algebra normalized by $\vph$, we can define a Lie stable $p$-rank; so the ``refined" part refers to the fact that we keep track of all these Lie stable $p$-ranks and not only of a particular sum of them. Strictly speaking, in 3.4-5 we just deal with Lie stable $p$-ranks associated to a fixed such factor of $G^{\rm ad}$ (see their parts referring to $I_0$ of 3.4.0).
\smallskip
In 3.10 we introduce the terminology (and notion)
of cyclic Shimura (adjoint) Lie $\sg$-crystals to be constantly used, as well as logical formulas (see 3.10.6-7) for them. 
\smallskip
In 3.12 we show that for many Shimura filtered $F$-crystals $(M,F^1,\vph,G,H,\tilde f)$ (defined in 2.2.10; so $H$ is a smooth subscheme of $G$ such that $H(W(k))$ contains $a_0$ of 1.2.4), the Shimura $\sg_{\nu}$-crystal over the algebraic closure $\nu$ of the field of fractions of the completion of the special fibre of $H$ in the special fibre of $a_0$, is a $G_{W(\nu)}$-ordinary $\sg_{\nu}$-crystal (here $\sg_{\nu}$ is the Frobenius automorphism of $W(\nu)$). 3.12.1 is an immediate consequence of d) of Theorem 1 and of (restatement 2.2.21 UP of) [Fa2, rm. iii) after th. 10].
\smallskip
By a lift of a Shimura $\sg$-crystal $(M,\vph,G)$ we mean a direct summand $F^1$ of $M$ such that the quadruple $(M,F^1,\vph,G)$ is a Shimura filtered $\sg$-crystal. Our hopes for a good (understandable and useful) theory of lifts of Shimura $\sg$-crystals $(M,\vph,G)$ which are not $G$-ordinary (or are not potentially cyclic diagonalizable; see defs. 2.2.1 d) and 2.2.22 1)) are gathered in 3.13: we define different (too many to be itemized) deviations of Shimura (filtered) (Lie) $\sg$-crystals. Here deviations are thought as measures of how far (or close) a Shimura $\sg$-crystal is from being $G$-ordinary or potentially cyclic diagonalizable. We draw attention to Conjecture 3.13.3.1 3) and to Problem 3.13.4 7).
\smallskip
The case $p=2$ is treated in 3.14: with very few exceptions (all of them are explicitly itemized), everything we get for $p\ge 3$ remains true for $p=2$ without any significant modification. In particular, Theorems 1, 2 and 3 remain true for $p=2$. With $p\ge 3$, 3.5.4 a) is a well known statement on standard ${\got s}{\got l}_2$-subalgebras of the Lie algebra of a semisimple group $SG$ over $\bar k$ which has a simple adjoint; it plays an important role in the proofs of 3.5. Unfortunately it is not true in general for $p=2$. However, for $p=2$ one gets around this problem by using any one (the choice depends on the context and so on $SG$) of the following three things (see 3.14 C; see also Step 2 of 3.13.7.3 and 3.14 J):
\medskip
{\bf a)} a slightly more involved argument pertaining to root systems;
\smallskip
{\bf b)} reduced structures of some normalizing group schemes over $\bar k$ in order to make them smooth;
\smallskip
{\bf c)} the classification of $SG$'s for which 3.5.4 a) does not hold.
\medskip
{\bf 1.2.7. Some conclusions and Dieudonn\'e's theories.} Some conclusions (with $p\ge 2$) obtained by combining the basic results of 3.1 and 3.6, are included in 3.15 (see also 3.1.8); they are detailed below and in 1.6.3-4. Some of them are refined in \S4 (like Theorem 1 of 3.15.1 is refined --see Theorem 13 of 1.12-- in 4.12.12 and 4.12.12.2).
\smallskip
In different parts of \S 3, we reobtain many well known results of crystalline Dieudonn\'e theories over (special fibres of) regular, formally smooth $W(k)$-schemes; moreover, in many situations we weaken --see below-- the often used restriction on such special fibres to be locally of finite type, as well as in many situations we work in a context involving finite, flat, commutative group schemes of $p$-power rank instead of one of (truncations of) $p$-divisible groups. For samples, see 3.6.18.5.1, 3.6.18.5.3, 3.6.18.5.7-8, Fact of 3.6.19, 3.6.20 3), 3.14 B and 3.15.2-3. To our knowledge, some results are new; like 3.15.2 (and its variants referred to in 3.15.3 3) and 6) and in 3.15.4). The idea behind 3.15.2 is the moduli principle; so, though a major result, it is not stated explicitly in this introduction. Here, in connection to these theories we just mention:
\medskip
$\bullet$ In 3.6.20 3) and 3.15.3 4) we show that (for $p\ge 3$) the results [Fa1, 7.1] and [Fa2, th. 10] can be deduced one from another (to our knowledge this is the first time when global Dieudonn\'e theories --like [Fa1, 7.1]-- are used to get local ones --like [Fa2, th. 10]--; usually the implications go the other way round, see [dJ2] and the proof of Theorem 2 of 3.15.1);
\smallskip
$\bullet$ In 3.6.20 3), in 3.14 B6-8, in the proofs of the two Theorems of 3.15.1 and in 3.15.3 3) and 6) we make the connection with [dJ1-2];
\smallskip
$\bullet$ The new ideas we introduce are based on (a more systematic use of) the theory of projective systems of (quasi-coherent sheaves on) schemes (see [EGA IV, \S 8]) and on some of the first six principles of 1.2.3. So some of the results are proved in wider contexts, like the ones involving pro-\'etale morphisms or finite, flat, commutative group schemes of $p$-power rank; moreover, we are stating everything in terms (i.e. in the language) of Fontaine categories (even when we deal with the special fibre of a regular, formally smooth $W(k)$-scheme; see 3.14 B6 for a sample): see 2.2.1.0 for Dieudonn\'e's crystalline, contravariant functor in the language of such categories. 
\medskip
Many of the ideas presented in \S1-4, can be easily adapted to wider contexts, which are either not necessarily regular, formally smooth (cf. 3.6.20 4)) or are involving a (non-perfect) field having a finite $p$-basis (cf. 3.15.3 6) and 3.15.4.1); as our work on $p$-divisible groups emerged from the smooth context of a standard Hodge situation (see 1.3), these extensions are just very partially exploited here.
\medskip
{\bf 1.2.8. On the context of 3.11.} The special features of the case $k=\bar k$ and some of their applications are presented in 3.11. Using defs. 2.2.1 d), 2.2.22 1) and 2.2.12, we have (cf. also 2.2.9 6)):
\medskip
{\bf Theorem 4.} {\it If $k=\bar k$, then all $G$-canonical lifts $(M,F^1,g\vph,G)$ ($g\in G(W(k))$) are isomorphic to each other (under isomorphisms defined by elements of $G(W(k))$), are cyclic diagonalizable and their Shimura filtered Lie $\sg$-crystals are cyclic diagonalizable and of Borel and parabolic type.}
\medskip
In particular, if $k=\bar k$ there is a $W(k)$-basis $\{e_i|i\in B(M)\}$ of $M$ and a permutation $\pi$ of the set $B(M)$, such that 
$$g\vph(e_i)=p^{\vep_i}e_{\pi(i)},$$ 
where $\vep_i$ is $1$ or $0$ depending on the fact that $e_i$ is or is not in $F^1$. All proofs of Theorem 4 presented (or mentioned) in 3.11.1 are in essence based on 1.2.1.0 3) and on properties of parabolic subgroups of $G$ which allow us to use the fact (see 3.6.8.1) that quasi Artin--Schreier systems of equations with coefficients in $k$ have always solutions. If $G=GL(M)$ then we can assume $\pi$ is the trivial permutation (cf. 3.1.1.1). 
\smallskip
Using Theorem 4, we get (see 3.1.4 and 3.11.2 A) that for any $G$-ordinary $\sg$-crystal $(M,\vph,G)$, there is a unique cocharacter $\mu:\GG_m\to G$ as in 1.1 ii) and such that the elements of the parabolic Lie subalgebra of ${\rm Lie}(G)$ corresponding to non-negative (resp. to non-positive) slopes of $({\rm Lie}(G),\vph)$ take $F^1$ (resp.  $F^0$) into itself; it is referred as the canonical split of $(M,\vph,G)$. It is the inverse of the canonical split cocharacter of $(M,F^1,\vph,G)$ as defined in [Wi], cf. 3.1.5. 
Useful direct sum decompositions (in the \'etale context with $\ZZ_p$-coefficients or in the integral crystalline context) can be found in 3.11.8.1 (see also 3.11.9).
\medskip
{\bf 1.2.9. General Fontaine categories.} References to Fontaine categories $\Mm\Mf_{[a,b]}(*)$ and $\Mm\Mf_{[a,b]}^\nabla(*)$, with $a,b\in\ZZ$, $b\ge a+2$, are made (besides some foundational material of parts of \S 2 and of 3.6.1 and some isolated remarks) in 3.6.8.9, 3.6.18.5.5, 3.6.18.7.1 b) and c), 3.6.18.7.1.1 b), 3.6.18.7.3, 3.6.18.7.4 4), 3.6.18.8, 3.6.18.8.1 b), end of 3.6.18.9, 3.13.7.8-9, 3.15.4-10, 4.5.4, 4.5.15.2.4 5), 4.5.18.1, 4.14.5 and in Appendix. Such categories were first used in [FL] (resp. in [Fa1]) for the case when $X={\rm Spec}(W(k))$ (resp. when $X$ is smooth over $W(k)$ and either it is affine or we are in situations in which, based on gluings arguments as in [Fa1, proof of 2.3], we do not need to specify a Frobenius lift). Here we just point out four things.
\smallskip
First, we define a wide (and wild) variety of (subcategories of) Fontaine categories (of objects of or $p$-divisible objects): see 2.2.1.7 for a glimpse. Different tannakian considerations are gathered in 2.2.4. In particular, we draw attention to the list of $\pi_1$-groupoids of 2.2.4 J. They are of interest among themselves. Moreover, they also form the first attempt to:
\medskip
-- define (quotients of) crystalline fundamental groupoids (and so indirectly groups) with integral coefficients, and to
\smallskip
-- (simultaneously) work as well in other contexts. 
\medskip
We recall that previous works of Deligne (see [De6]) and others (more recently see [Shi1-2] and [HK]) were centered on a rational, unipotent context. Here, by ``other contexts" we mean contexts which do pay (at least partially) attention to reductive (group) aspects and by integral coefficients we mean $\ZZ_p$ or $W(k)$ coefficients (versus $\QQ_p$ or $B(k)$). We go as far as: 
\medskip
{\bf i)} to present some tools (which hopefully will be used later on and) which pertain to Fontaine categories (see all of 2.2.4);
\smallskip
{\bf ii)} to define, without assuming that $X$ or its $p$-adic completion has a Frobenius lift, a category $\Mm\Mf^{\nabla(p+tens)}(X)$ which has most of the features one could think of (the still to be defined) category ``$\Mm\Mf^\nabla(X)$" (see 2.2.4 C and D);
\smallskip
{\bf iii)} to speak about tannakian categories and almost fibre functors over $W_n(k)$, $n\in\NN$ (see ii) and 2.2.4 F);
\smallskip
{\bf iv)} to define Fontaine's crystalline fundamental $W(k)$-groupoid of a proper, smooth, geometrically connected $W(k)$-scheme (see 2.2.4 K); 
\smallskip
{\bf v)} to define the crystalline fundamental $W(k)$-groupoid in crystals of the special fibre $X_k$ of $X$ (see 2.2.4.1 1));
\smallskip
{\bf vi)} and to include general affine versions of iv).  
\medskip
Second, most of the properties and principles we mentioned for $\Mm\Mf_{[0,1]}(*)$ in 1.2.4-5, are either not true (in general) or are true in a much weaker form for $\Mm\Mf_{[a,b]}(*)$; this is so due to the fact that working as usual modulo powers of $p$ and following the approach of 1.2.5, we get systems of equations in $n$ variables which besides the equations of (1) of 1.2.5 do involve in general some extra equations (cf. 3.6.8.9; we refer to such systems of equations of being of third type). 
\smallskip
Third, some of the results have to be stated not in terms of multiplicities of the slopes $0$ and $1$ or of the slope $-1$ (as we usually do in the context of $\Mm\Mf_{[0,1]}(*)$: see 3.6.18.4, 3.6.18.7.0, etc.) but in terms of suitable pseudo-multiplicities of the slope $-1$, as defined in 3.6.18.5.5 B. The reason is: the liftability property referred to in 1.2.3 does not hold (in general) in the context of $\Mm\Mf_{[a,b]}(*)$.
\smallskip
Fourth, in the study of (truncations of) $p$-divisible objects of $\Mm\Mf(W(k))$ (see 2.2.1 c) for defs.) endowed with some extra structures, we distinguish five major classes. Working (for simplicity of language) in a non-filtered context, they are listed below in such a way that the next ones contain the previous ones:
\medskip
{\bf 1)} The class $\Mc\Ml_1$ of Shimura $\sg$-crystals;
\smallskip
{\bf 2)} The class $\Mc\Ml_2$ of generalized Shimura $p$-divisible objects over $k$ (see 2.2.8 3) and 4));
\smallskip
{\bf 3)} The class $\Mc\Ml_3$ of $p$-divisible objects with a reductive structure over $k$ (see 2.2.8 3a) and 4a));
\smallskip
{\bf 4)} The class $\Mc\Ml_4$ of $\sg$-$\Ms$-crystals (see 3.6.1.5);
\smallskip
{\bf 5)} The most general context (involving a non-necessarily smooth group scheme).
\medskip
This division is based on degrees of complexities. A great part of the results of 1.2.1-8 (pertaining to $\Mc\Ml_1$) extend automatically to $\Mc\Ml_2$ and even to $\Mc\Ml_3$ or to $\Mc\Ml_4$. For instance, this paper ``handles" entirely $\Mc\Ml_2$ (i.e. everything we get here for $\Mc\Ml_1$ is obtained as well for $\Mc\Ml_2$). For the (local and global) deformation theory in the context of $\Mc\Ml_2$ we refer to 3.6.18.7.1 c), 3.6.18.7.3 C and 3.15.6. One might wonder (based on the above use of the word entirely), why we do not treat $\Mc\Ml_2$ and $\Mc\Ml_1$ as just one class. We have two reasons for this:
\medskip
-- presently, in general, we do not have geometric interpretations for representatives of $\Mc\Ml_2$ (like $p$-divisible groups, Verschiebung maps of them or the existence of integral canonical models of Shimura varieties of preabelian type);
\smallskip
-- occasionally we prove results pertaining to $\Mc\Ml_2$ by reduction to the context of $\Mc\Ml_1$ (like the integrability principle; see 3.15.6 D).
\medskip
Most of the ideas apply as well directly to $\Mc\Ml_3$; however, in Step 1 of Appendix we state the things only as far as all details are obviously the same as for $\Mc\Ml_2$. One think is worth mentioning: we do not present a deformation theory for $\Mc\Ml_3$; however, such a theory can be substituted from many points of view (like ``smallest Newton polygons") either by 3.6.6, 3.6.6.0 and Grothendieck--Katz' specialization theorem or by [RR] (see Step 1 of Appendix). For avoiding repetitions, some results pertaining to $\Mc\Ml_2$ are stated just in Step 1 of Appendix, as refinements of results pertaining to $\Mc\Ml_3$. In connection to $\Mc\Ml_4$ (resp. to $\Mc\Ml_5\setminus\Mc\Ml_4$) we refer to 3.6.1.6 and Appendix (resp. to 3.9.9). In 3.9.9 we construct quasi-versal (in the sense of 3.6.19 B) global deformations (over $W(k)$-schemes as in 1.2.1) in contexts provided by suitable quintuples $(M,F^1,\vph,H,G)$, with $(M,F^1,\vph,G)$ as in 1.1 and with $H$ a closed (but not necessarily smooth), integral subgroup of $G$; in particular, we get such global deformations for quasi-polarized $p$-divisible groups over $W(k)$.
\medskip\smallskip
{\bf 1.3. Standard Hodge situations.} For a better presentation of geometric applications (see \S 4) of the basic results of \S 3, we introduce in 2.3 the standard Hodge situation; it generalizes [Ko2, ch. 5].
\smallskip
We start with an injective map $f:{\rm Sh}(G,X)\hookrightarrow {\rm Sh}\bigl({\rm GSp}(W,\psi),S\bigr)$ of Shimura pairs. Here the pair $({\rm GSp}(W,\psi),S)$ defines a Siegel modular variety, cf. [Va2, 2.5]. Let $\ZZ_{(p)}$ be the localization of $\ZZ$ with respect to the rational prime $p\ge 3$. 
\smallskip
Let $v$ be a prime of the reflex field $E(G,X)$ dividing $p$ and let $O_{(v)}$ be the localization of the ring of integers of $E(G,X)$ with respect to $v$.
We assume the existence of a $\ZZ_{(p)}$-lattice $L_{(p)}$ of $W$ which is
crystalline well positioned for the map $f$ with respect to $v$ (cf. def. 2.3.4; see also below). The triple $\bigl(f,L_{(p)},v\bigr)$  is called a standard Hodge situation.
\smallskip
Let $K_p:=\{g\in {\rm GSp}(W,\psi)(\QQ_p)\bigm| g(L_{(p)}\otimes_{\ZZ_{(p)}} \ZZ_p)=L_{(p)}\otimes_{\ZZ_{(p)}} \ZZ_p\}$.
The hypothesis on $L_{(p)}$ first implies: $L_{(p)}$ is a good $\ZZ_{(p)}$-lattice with respect to $f$. This means (see def. [Va1, 5.8.3]) that
$\psi$ induces a perfect form $L_{(p)}\otimes_{\ZZ_{(p)}} L_{(p)}\to\ZZ_{(p)}$ (i.e. the induced map from $L_{(p)}$ into its dual $L^\ast_{(p)}$ is an isomorphism) and the Zariski closure
$G_{\ZZ_{(p)}}$ of $G$ in ${\rm GSp}(L_{(p)},\psi)$ is a reductive group over $\ZZ_{(p)}$. So $K_p$ is a hyperspecial subgroup of ${\rm GSp}(W,\psi)(\QQ_p)$ and the intersection 
$$H:={\rm GSp}(W,\psi)(\QQ_p)\cap K_p$$ 
is a hyperspecial subgroup of $G(\QQ_p)$.
\smallskip
Let $e\in\NN$ be such that $2e=\dim_\QQ(W)$. We choose a
$\ZZ$-lattice $L_\ZZ$ of $W$ such that $\psi$ induces a perfect form 
$\psi:L_\ZZ\otimes_{\ZZ} L_\ZZ\to\ZZ$ and $L_{(p)}=L_\ZZ\otimes_{\ZZ} \ZZ_{(p)}$. Let
$(v_\al)_{\al\in\Mj}$ be a family of tensors in spaces of the form $W^{*\otimes n}\otimes_{\QQ} W^{\otimes n}$, with $n\in\NN$ and with $W^*$ as the dual of $W$, such that $G$ is the subgroup of ${\rm GSp}(W,\psi)$ fixing $v_{\al}$, $\forall\al\in\Mj$. The choice of the
lattice $L_\ZZ$ and of the family $(v_\al)_{\al\in\Mj}$ allows the interpretation of
${\rm Sh}(G,X)(\CC)$ as the set of isomorphism classes of principally polarized abelian
varieties over $\CC$ of dimension $e$, carrying a family of Hodge
cycles indexed by the set $\Mj$ and satisfying some additional conditions and level structures (cf. [Va2, 4.1]).
\smallskip
It is known (for instance, see [Va2, 3.2.9 and 4.1]) that the $\ZZ_{(p)}$-scheme $\Mm$ parameterizing isomorphism classes
of principally polarized abelian schemes of dimension $e$ over $\ZZ_{(p)}$-schemes, having level-$N$ symplectic similitude structure for any $N\in\NN$ relatively prime to $p$, together
with the canonical action of ${\rm GSp}(W,\psi)(\AA^p_f)$ on it, is an integral canonical model of
${\rm Sh}\bigl({\rm GSp}(W,\psi),S\bigr)/K_p$. These structures and this action are defined naturally --see [Va2, 4.1]-- via $L$. See 2.1 for the $\QQ$--algebras $\AA_f^p$ and $\AA_f$. 
\smallskip
The normalization $\Mn$ of the Zariski closure of ${\rm Sh}(G,X)/H$ in $\Mm_{O_{(v)}}$ is (cf. 2.3.3) the integral canonical model of ${\rm Sh}(G,X)/H$ (or of Shimura quadruple $(G,X,H,v)$; see [Va2, 3.2.6] for the definition of these quadruples). Let $k(v)$ be the residue field of $v$. The universal principally polarized abelian scheme over $\Mm$ gives birth,
by pull back, to a principally polarized abelian scheme $(\Ma,\Mp_{\Ma})$ over $\Mn$, which we call special. 
\medskip
{\bf 1.3.1. A basic assumption.} The
crystalline well positioned property of $L_{(p)}$ for the map $f$ with respect to $v$ also says: any
point $z\in\Mn(W(k))$ gives birth to a principally polarized abelian scheme $(A,p_A)$ over $W(k)$ (obtained from $(\Ma,\Mp_{\Ma})$ by pull back) endowed naturally with a family $(w_\al)_{\al\in\Mj}$ of Hodge cycles such that the filtered $\sg$-crystal $(M,F^1,\vph)$ (with $M:=H^1_{\rm crys}(A/W(k))$, with $\vph$ as the  $\sg$-linear endomorphism of $M$ and with $F^1$ as the Hodge filtration of $M$ defined by $A$) gets a natural structure of a ``not necessarily quasi-split" Shimura filtered $\sg$-crystal $(M,F^1,\vph,\tilde G_{W(k)})$. In other words:
\medskip
\item{(*)} {\it The subgroup of $GL(M)$ obtained by taking the Zariski closure of the subgroup of $GL\bigl(M[{1\over p}]\bigr)$ fixing the perfect alternating form $p_A:M\otimes_{W(k)} M\to W(k)(1)$ (induced by the polarization $p_A$ of the abelian variety $A$ and still denoted by $p_A$) and the de Rham component $t_\al$ of $w_\al$ (this component is a 
tensor in the tensor algebra of $(M\oplus M^*)[{1\over p}]$ fixed by $\vph$, cf. the Corollary of 2.3.10), $\forall\al\in\Mj$, is a reductive group (not a priori quasi-split) $\tilde G_{W(k)}$ over $W(k)$.}
\medskip
There is a finite, nilpotent Galois cover ${\rm Spec}(W(k_1))$ of ${\rm Spec}(W(k))$ of index of nilpotency at most two and of a very particular type (see 2.3.9 b) and d)) such that $\tilde G_{W(k_1)}$ is isomorphic to $G_{W(k_1)}$. Moreover, for many standard Hodge situations we can show (see 2.3.9 c)) that $\tilde G_{W(k)}$ is in fact automatically isomorphic to $G_{W(k)}$. So, just for the sake of convenience, we make the convention (see 2.3.9.2) that, without any extra reference, we consider from now on only points ${\rm Spec}(k)\to\Mn$ with the property that for a (any) lift of it ${\rm Spec}(W(k))\to\Mn$, the reductive group over $W(k)$ we get as in (*), is (isomorphic to) $G_{W(k)}$ and so is quasi-split. The proof of 2.3.9 uses Fontaine's comparison theory, properties of reductive groups as well as a great part of [Va2, 6.4-6].
Based on this convention, the triple $(M,\vph,G_{W(k)})$ (resp. the quadruple $(M,F^1,\vph,G_{W(k)})$) is called the Shimura  $\sg$-crystal (resp. the Shimura filtered $\sg$-crystal) attached to the point $y\in\Mn(k)$ defined by $z$ (resp. attached to $z$). When we want to emphasize $p_A$ we also speak about the principally quasi-polarized Shimura (resp. Shimura filtered) $\sg$-crystal attached to $y$ (resp. to $z$). Notation: $(M,\vph,G_{W(k)},p_A)$ (resp. $(M,F^1,\vph,G_{W(k)},p_A)$). 
\medskip
{\bf 1.3.1.1. Remark.} (*) and the local deformation theory of [Va2, 5.4] imply (see 2.3.11): for any regular, formally smooth $W(k)$-scheme $\tilde U$ and for every morphism $\tilde U\to \Mn$, we get naturally (via pull back) a principally quasi-polarized Shimura $p$-divisible group over $\tilde U$.
\medskip
{\bf 1.3.2. The deformation principle.} This principle in its essence says (cf. 3.6.14):
\medskip
{\bf Theorem 5.} {\it We consider a point $y:{\rm Spec}(k)\to\Mn_{W(k)}$. We assume that a starting assumption (see SA1 of 3.6.14) is satisfied for it. Then there is a smooth morphism $\Mn(y)\to\Mn_{W(k)}$ and a pro-\'etale morphism $m^1:\Mn(y)^1\to\Mn(y)$, with $\Mn(y)^1$ and affine scheme having a geometrically connected special fibre $\Mn(y)^1_k$, such that:
\medskip
{\bf a)} ${\rm Im}(m^1)$ contains an open, 
dense subscheme of $\Mn(y)_k$ and $y\in m^1(\Mn(y)^1(k))$, and
\smallskip
{\bf b)} the principally quasi-polarized Shimura $p$-divisible group $\Md^1$ over the $p$-adic completion of $\Mn(y)^1$, obtained by natural pull back (see 1.3.1.1), has a very nice (down to earth) expression. 
\medskip
What we mean by b): $\Md^1$ is obtained by pulling back (through the $p$-adic completion of a pro-\'etale morphism) the principally quasi-polarized Shimura $p$-divisible group $\Md_1$ over ${\rm Spec}(R_1^\wedge)$ produced entirely as in Theorem 3, but starting from the principally quasi-polarized Shimura $\sg$-crystal attached to $y$ (so here $U$ is a suitable open subscheme of the subgroup of $G_{W(k)}$ fixing the perfect alternating form $p_A$ on $M$).} 
\medskip
We get this principle by combining standard arguments with our approach of working inductively modulo powers of $p$; though lengthy, the proof of Theorem 5 is very instructive. The starting assumption refers to the fact that, working mod $p$, the variant of Theorem 5, obtained by working with henselizations of localizations of smooth $\Mn_k$-schemes in $k$-valued points lifting $y$, is a priori satisfied. In practice it is very hard to check if a starting condition is satisfied or not. However the proof of Theorem 5 involves some very important techniques and its variants modulo powers of $p$ (see 3.6.14.4) involving no assumption, are also very useful for different local computations. Moreover, the fact that a starting condition is satisfied or not, is at the very root of many extremely important problems (see \S 5-8; see also 3.6.19 and the idea of 4.10.5).
\smallskip
There are variants of this principle: we work with a point $y_{H_0}:{\rm Spec}(k)\to\Mn_{W(k)}/H_0$ (here $H_0$ is a compact, open subgroup of $G(\AA_f^p)$ such that the quotient morphism $\Mn\to\Mn/H_0$ is a pro-\'etale cover and ($\Ma,\Mp_{\Ma})$ descends to $\Mn/H_0$) or with an \'etale morphism $\Mn(y)\to\Mn_{W(k)}$, cf. 3.6.14.1-3.
\medskip\smallskip
{\bf 1.4. Main geometric concepts.} Starting from Theorem 1, we are able to
\medskip
-- reobtain (see 4.6.1 3)) the well
known results (pertaining to integral canonical models of Siegel modular varieties and to universal principally polarized abelian schemes over them) concerning the existence of an ordinary type and the existence of the canonical lift of a (principally polarized) ordinary abelian variety, 
\smallskip
-- as well as to extend them to any integral canonical model $\Mn$ of a Shimura variety ${\rm Sh}(G,X)$ of Hodge type and to every special principally polarized abelian scheme ($\Ma,\Mp_{\Ma})$ over $\Mn$, arising (as in the end of 1.3) from a standard Hodge situation $(f,L_{(p)},v)$. 
\medskip
For instance, with the notations of 1.3 we get the following three things.
\medskip
{\bf a)} First we get a Shimura-ordinary type $\tau$ (cf. 4.1). It depends only on $f$ and $v$ and not on the choice of $L_{(p)}$. It can be described as follows. 
\smallskip
The isomorphism $L_{(p)}\tilde\to L^\ast_{(p)}$ induced by $\psi$ allows us to identify $G_{\ZZ_{(p)}}$
with a subgroup of $GL(L^\ast_{(p)})$. Let $L_p^\ast:=L_{(p)}^\ast\otimes_{\ZZ_{(p)}} \ZZ_p$; so $G_{\ZZ_p}$ is naturally a subgroup of $GL(L_p^\ast)$. Let $T$ be a torus of $G_{\ZZ_p}$ such that:
\medskip
\item{1)} there is a cocharacter $\mu:\GG_m\to T_{W(k(v))}$ which over an embedding of $W(k(v))$ in $\CC$ (extending the composite inclusion $O_{(v)}\subset E(G,X)\subset\CC$) is
$G(\CC)$-conjugate to the cocharacters $\mu_x^\ast:\GG_m\to G_\CC$, $x\in X$, which are the cocharacters $\mu_x$ mentioned in [Va2, 2.2], but viewed in the dual context (so they are the cocharacters of $G_{\CC}$ defining the Hodge $\QQ$--structures on $W^\ast$ associated to points $x\in X$);
\smallskip
\item{2)} there is a Borel subgroup $B$ of $G_{\ZZ_p}$ containing $T$ and whose Lie algebra is such that its elements take the $F^1$-filtration $F^1$ of $L^\ast_p\otimes_{\ZZ_p} W(k(v))$ defined by $\mu$ into itself. 
\medskip
We recall that $\mu$ gives birth to a direct sum decomposition $L^\ast_p\otimes_{\ZZ_p} W(k(v))=F^1\oplus F^0$, with
$\be\in \GG_m(W(k(v)))$ acting through $\mu$ on $F^i$ as the multiplication with $\be^{-i}$, $i=\overline{0,1}$; so 2) can be restated as: $b(F^1)\subset F^1$, $\forall b\in {\rm Lie}(B)$.
\smallskip
$\tau$ is the formal isogeny type associated to the $\sg$-crystal over $k(v)$ 
$$
(L^\ast_p\otimes_{\ZZ_p} W(k(v)),\vph)
$$
defined by $\vph:=\sg\circ\mu({1\over p})$. Here $\sg$ acts as identity on $L^\ast_p$ and as
Frobenius automorphism on $W((k(v))$, while 
${1\over p}\in\GG_m\bigl(W\bigl(k(v)\bigl[{1\over p}\bigr]\bigr)\bigr)$. Due to 1) and the expression of $\vph$, the slopes of $\tau$ have as denominators divisors of the degree $[k(v):\FF_p]$ (see 4.6 P4).
\medskip
{\bf b)} Second we get $G$-ordinary (or Shimura-ordinary) points of $\Mn_{k(v)}$ (cf. 4.2). One way to define them: they are those points, with values in fields, which have the property that the abelian varieties, obtained from $\Ma$ by pull backs via them, have $\tau$ as their formal isogeny type. They are dense in $\Mn_{k(v)}$.
\medskip
{\bf c)} Third we get $G$-canonical (or Shimura-canonical) lifts of $G$-ordinary points of $\Mn_{k(v)}$ with values in perfect fields (cf. 4.4); they are uniquely determined (see 1.5). These $G$-canonical lifts are points of $\Mn$ with values in Witt rings over perfect fields.
\medskip
{\bf 1.4.1. On the passage from abstract to geometric contexts.} 1.4 b) and c) are a consequence of Theorems 1 and 2, via a natural algebraization process (see 2.3.15-16) of the local deformation theory of [Va2, 5.4]; this process is supported by:
\medskip
-- the natural interpretation of loc. cit. in terms of filtrations (see 2.3.17; see also 2.4 for an abstract presentation in the right context of PD-hulls);
\smallskip
-- Fontaine's comparison theory (it is the essence behind 4.2.3, for instance, cf. [Va2, 5.2.17.2]);
\smallskip
-- and by a general Lemma (see 3.6.6) on Shimura $\sg$-crystals. 
\medskip
$\tau$ is a usual ordinary type (i.e. it is the formal isogeny type $e(1,0)+e(0,1)$; see 2.1) if and only if $k(v)$ is $\FF_p$ (cf. 4.6 P1: it is an immediate consequence of how $\tau$ is defined). If $k(v)=\FF_p$, then the abelian variety over $W(k)$ obtained from $\Ma$ by pull back through any $W(k)$-valued $G$-canonical lift of $\Mn$, is the canonical lift of an ordinary abelian
variety (cf. 4.6 P2: it is a consequence of 4.6 P1 and of b) of Theorem 1; cf. also 2.3.17).
\smallskip
The proof of 3.6.6 relies heavily on the language of root systems and of the classification mentioned in 1.2.1.0. In fact 3.6.6 (in the form of the Exercise 3.6.6.0) is from some points of view equivalent to a) and c) of Theorem 1: see Step 1 of Appendix.
\medskip\smallskip
{\bf 1.5. Shimura (filtered) Lie $\sg$-crystals attached to points.} Let $y:{\rm Spec}(k)\to\Mn_{k(v)}$. The Shimura Lie $\sg$-crystal $({\got g},\vph)$ attached to the Shimura $\sg$-crystal attached to $y$, is also referred as the Shimura Lie $\sg$-crystal attached to $y$ (cf. 2.3.10; see also 1.3.1); here ${\got g}:={\rm Lie}(G_{W(k)})$. Similarly (see 4.3.1), we speak about the Faltings--Shimura--Hasse--Witt (adjoint) map attached to $y$. 
\smallskip
Any lift $z:{\rm Spec}(W(k))\to\Mn$ of $y$, makes $({\got g},\vph)$ to be a filtered Lie $\sg$-crystal, i.e. it makes it to be a $p$-divisible object of
$\Mm\Mf_{[-1,1]}(W(k))$ endowed naturally with a Lie structure. So we get a filtration
$$0=F^2({\got g})\subset F^1({\got g})\subset F^0({\got g})\subset F^{-1}({\got g})={\got g}$$ 
such that
$$\vph\bigl({1\over p}F^1({\got g})+F^0({\got g})+p{\got g}\bigr)={\got g};$$ 
$F^0({\got g})$ is a parabolic Lie subalgebra of {\got g} having $F^1({\got g})$ as its nilpotent radical, $F^1({\got g})$ is abelian and $[{\got g},F^1({\got g})]\subset F^0({\got g})$. We refer to the quadruple $({\got g},\vph,F^0({\got g}),F^1({\got g}))$ as the Shimura filtered Lie $\sg$-crystal attached to $z$. a) of Theorem 1 implies (via 2.3.17): $y$ is a $G$-ordinary point iff there is a $W(k)$-valued point $z$ of $\Mn$ lifting it and whose attached Shimura filtered Lie $\sg$-crystal is of parabolic type. Such a lift $z$, when exists, is unique (cf. b) of Theorem 1 and 2.3.17) and defines the $G$-canonical lift of $y$.
\smallskip
The $G$-ordinary points of $\Mn_{k(v)}$ (resp. the $G$-canonical lifts of $G$-ordinary points) with values in perfect fields have all usual properties of ordinary points (resp. of canonical lifts of ordinary points) of $\Mm$ with values in perfect fields, which can be read out from Shimura (resp. Shimura filtered) Lie $\sg$-crystals attached to points of $\Mm_{\FF_p}$ with values in perfect fields (resp. in Witt rings of perfect fields). As samples of this philosophy see 1.7-9 and 4.4.4. Of course, these properties have to be reformulated accordingly (if $k(v)\neq\FF_p$).
\medskip
{\bf 1.5.1. A new, more general concept: $U$-ordinariness.} A natural question arises: for which other $k$-valued point of $\Mn_{k(v)}$ we can (identify and) define naturally a lift of it to a $W(k)$-valued point of $\Mn$, which is uniquely determined by some conditions? In 4.4.13, as a digression, we answer partially this question: we define $U$-ordinary points and $U$-canonical lifts. The $U$-ordinary points of $\Mn_{k(v)}$ with values in $k$ are those points whose attached Shimura $\sg$-crystals have a unique lift of quasi CM type in the sense of 2.2.17 (and so in particular are potentially cyclic diagonalizable): these unique lifts define the $U$-canonical lifts. Any $G$-ordinary point is a $U$-ordinary point and any $G$-canonical lift is a $U$-canonical lift (see 4.4.13.3); but there are tremendously many examples of standard Hodge situations for which we have (plenty of) $U$-ordinary points which are not $G$-ordinary. These examples do not show up in the classical setting of $\Mm$: they typically show up when $G^{\rm ad}_{\QQ_p}$ has certain simple factors which are not absolutely simple. See 4.4.13.3.1 for first examples, in the abstract context of Shimura (Lie) $\sg$-crystals. The simplest concrete geometric examples can be obtained starting from 4.12.12.6.6 3). 
\smallskip
The passage from $G$-ordinary points to $U$-ordinary points is achieved via (see 4.4.13.2): we move from keeping track of parabolic Lie subalgebras of ${\rm Lie}(G_{W(k)})$ corresponding to non-negative slopes, to keeping track of Lie subalgebras (they are automatically reductive in the potentially cyclic diagonalizable context, cf. 2.2.19.1-2) of ${\rm Lie}(G_{W(k)})$ corresponding only to the slope $0$; here slopes are with respect to some $\vph$ as in 1.3.1. When we are in the context of $U$-ordinary points and these last Lie subalgebras are abelian, we get a subclass of $U$-ordinary points which we call $T$-ordinary points. The $G$-ordinary points are just exceptionally $T$-ordinary points (see 4.6 P7).
\smallskip
Our hopes for good (i.e. unique in some sense) lifts for other $k$-valued points of $\Mn_{k(v)}$ are expressed in question $Q_3$) of 4.5.15.
\medskip\smallskip
{\bf 1.6. Stratifications.} The Newton polygons of (isocrystals defined by) Shimura Lie $\sg$-crystals attached to points of $\Mn_{k(v)}$ with values in perfect fields, achieve (cf. 4.5) a stratification 
(called the canonical Lie stratification) of $\Mn_{k(v)}$ into $G(\AA^p_f)$-invariant, reduced, locally closed subschemes. It is similar (in nature and properties) to the well known Newton polygon stratification of $\Mm_{\FF_p}$. 
\smallskip
More generally, for any representation 
$$\rho: G_{\ZZ_p}\to GL(N),$$ 
with $N$ a free $\ZZ_p$-module of finite rank, we define (via Newton polygons, see 4.5.4) a $\rho$-stratification of $\Mn_{k(v)}$ in $G(\AA_f^p)$-invariant, reduced, locally closed subschemes. The intersection of all these $\rho$-stratifications defines the absolute stratification. This absolute stratification is the most refined stratification of $\Mn_{k(v)}$ which can be obtained naturally (i.e. in an adequate Shimura context, involving representations of $G_{\ZZ_p}$) with the help of Newton polygons. The $\rho$-stratification depends only on the restriction of $\rho$ to $G_{\QQ_p}^{\rm der}$ (cf. 4.5.6 9)); however, there is a lot of data lost if we consider just such restrictions (see end of 4.5.4 for a glimpse and see 1.15.1 for a general principle).
\smallskip
In particular, we also define (see 4.5.2) the refined
canonical Lie stratification of $\Mn_{k(v)}$: it is the intersection of those $\rho$-stratifications which have $\rho$ as a subrepresentation of the adjoint representation of $G_{\ZZ_p}$ on ${\rm Lie}(G^{\rm ad}_{\ZZ_p})$. If $(G^{\rm ad},X^{\rm ad})$ has all its simple factors of $C_n$ type or of some particular $B_n$ type, with $n\in\NN$, then this stratification coincides (see 4.5.6.1) with the absolute stratification; but this is not necessarily true in general (see 4.5.6.2). The $G$-ordinary points of $\Mn_{k(v)}$ are points of the generic (i.e. dense) open stratum (it does not matter which stratification we use, cf. 4.5.6 8)).
\medskip
{\bf 1.6.1. The (quasi-) ultra stratification.} Warning: till 1.9.1 we assume $k=\bar k$. A long list of other types of stratifications of $\Mn_{k(v)}$ is presented in 4.5.9 and 4.5.15. To explain them let us start mentioning some main points.
\medskip
$\bullet$ Many of them are refinements of the absolute stratification of $\Mn_{k(v)}$ (see question $Q_2$) of 4.5.15).
\smallskip
$\bullet$ The most useful ones are the pseudo-ultra, the quasi-ultra (or the Faltings--Shimura--Hasse--Witt) and the ultra stratification (see 4.5.15.1).
\smallskip
$\bullet$ The quasi-ultra stratification is defined purely in terms of inner isomorphism classes (over algebraically closed fields) of Faltings--Shimura--Hasse--Witt adjoint maps attached to points of $\Mn_{k(v)}$ with values in such fields, while the ultra stratification is defined as the intersection of the quasi-ultra and the absolute stratifications (it is 4.2.10 which allows us to speak about inner isomorphism classes; for instance, referring to the map $FSHW^{\rm ad}$ and notations of 1.2.1.0, we allow only those isomorphisms between $FSHW^{\rm ad}$ and $i_g\circ FSHW^{\rm ad}$, with $i_g$ as the inner automorphism of ${\rm Lie}(G^{\rm ad})/p{\rm Lie}(G^{\rm ad})$ defined by $g\in G^{\rm ad}_k(k)$, which are defined as well by such inner automorphisms).
\smallskip
$\bullet$ In order to define these last two stratifications we either assume that the set of such inner isomorphism classes is finite (we expect --see 4.5.15.1; see also 1.6.2-- that this set is always finite, under no restrictions) or allow stratifications (they are formalized in 2.1) with potentially an infinite number of strata.
\smallskip
$\bullet$ The quasi-ultra stratification can be studied without using any Verschiebung map; its main properties can be listed as follows (see 3.13.7.1-7, 4.5.15.2 and 4.5.15.2.1-5 for more details):
\medskip
{\bf a)} its strata are in one-to-one correspondence to a subset of the set of orbits of some algebraic group action $\Mg\Ma$ over the algebraic closure $\FF$ of $\FF_p$: $\Mg\Ma$ is constructed purely in terms of $G_{\FF_p}^{\rm ad}$ and of the special fibre of the cocharacter $\mu$ of 1.4 1) (in the abstract context, see $\TT^0$ of 3.13.7.1); more precisely, the number of its strata is conjectured to be the number of elements of the quotient set of the Weyl group $WG$ of $G^{\rm ad}_{\FF}$ which parameterizes orbits of $\Mg\Ma$ defined naturally by elements of $WG$ (this can be grasped from 1.6.2; see 4.12.12.6 for why all elements of $WG$ do show up);
\smallskip
{\bf b)} which stratum specializes to which stratum can be read out from $\Mg\Ma$ (see end of 4.5.15.2.1);
\smallskip
{\bf c)} all its strata are regular (this is a consequence of a) and of 2.3.15);
\smallskip
{\bf d)} the dimensions of its strata can be read out from $\Mg\Ma$; in particular, their codimensions are given by dimensions of automorphism groups of Faltings--Shimura--Hasse--Witt adjoint maps attached to $\FF$-valued points of $\Mn_{k(v)}$ (again this is a consequence of 2.3.15);
\smallskip
{\bf e)} if $f$ is an isomorphism (i.e. if $\Mn=\Mm$), it coincides with the canonical stratification as defined (for suitable quotients of $\Mm_{\FF_p}$) in [EO] (see also its published first preliminary version [Oo3]); 
\smallskip
{\bf f)} the quasi-affineness property of loc. cit. implies that all its strata are quasi-affine.  
\medskip
So the philosophy of the quasi-ultra stratification is:
\medskip
{\bf Ph.} {\it Faltings--Shimura--Hasse--Witt adjoint maps are the adjoint Lie analogue (in any relative context pertaining to Shimura Lie $\sg$-crystals) of truncations mod $p$ of $p$-divisible groups over $\FF_p$-schemes.}
\medskip
The only thing we use in \S1-14 from [Oo3] (cf. the quasi-affineness part of f)), can be in fact entirely avoided, cf. 4.9.17.0.0. 
\medskip
{\bf 1.6.2. A main expectation.} The main expectation on the precise number of orbits of $\Mg\Ma$ can be formulated (with the def. of 2.2.14) in an abstract context as follows (the notations are as in 1.1-2):
\medskip
{\bf Expectation (the CM level one property).} {\it For any Shimura $\sg$-crystal $(M,\vph,G)$ over $k$, there is a lift $F^1$ of it such that the truncation mod $p$ of $(M,F^1,\vph,G)$ is isomorphic, under an isomorphism defined by an element of $G(k)$, to the truncation mod $p$ of $(M,F^1,g\vph,G)$, where $g\in G(W(k))$ is such that, for a maximal torus $T$ of the parabolic subgroup of $G$ normalizing $F^1$, the quadruple $(M,F^1,g\vph,T)$ is a Shimura filtered $\sg$-crystal.}
\medskip
This should be provable by slightly refining the proofs of 3.4.6 and 3.6.6, cf. 3.5.3-4 and the way the group actions $\TT$ and $\TT^0$ of 3.13.7.1 are constructed. 3.13.7.2 makes the connection between this Expectation and inner isomorphism classes of Faltings--Shimura--Hasse--Witt adjoint maps. In 3.13.7.1 we develop some steps towards the proof of this Expectation: it gets  reduced to checking that a Bruhat decomposition (naturally associated to $(M,\vph,G)$) holds for $G(k)$, see 3.13.7.1.3-4. We do believe that such decompositions can be also checked abstractly, in an adequate $\sg$- or $F$-context pertaining to Tits systems: see 3) of 3.13.7.4 D for a precise statement on such decompositions in the $F$-context for reductive groups. Here ``$F$-" is used as in the Steinberg's reformulation of Lang's theorem (see [Hu2, 8.3]). 
\smallskip
For some extra reduction steps towards the proof of the existence of such decompositions (and so towards the proof of this Expectation), see 3.13.7.3-4, 3.13.7.6 and 4.6.7-8. For the case when we are dealing with a reductive group whose adjoint has all its simple factors of $A_1$ Lie type, such decompositions are proved to exist in 3.13.7.5. See also 3.13.7.6.0 for how they can be used in connection to inner isomorphism classes of Faltings--Shimura--Hasse--Witt adjoint maps which are $\sg$-linear endomorphisms of the Lie algebra of an adjoint group over $k$ of whose simple factors are all of $A_1$ Lie type. In turns out that the mentioned $A_1$ Lie type case is just a particular case of the so called essentially Borel subgroup case: see Corollary 2 of 3.13.7.5 and 3.13.7.8.2 for precise statements and a description of such a case.
\medskip
{\bf 1.6.2.1. Some generalizations.} A great part of the theories of 1.6.1-2 can be redone in terms of (inner isomorphism classes of) truncations mod $p$ of Shimura $\sg$-crystals (see 2.2.14 for defs.): see 4.5.15.2.4 3) for (geometric) justifications of why we still prefer the adjoint context. Despite this and the philosophy of 1.6.1, we view the refined canonical stratifications hinted at in 4.5.9 as the right generalization of the canonical stratification of [EO]. If one wants, the quasi-ultra stratifications are the generalization of the generalized canonical stratifications. Motivation: they provide the right context which can be immediately extended to any geometric context involving $p$-divisible objects with a reductive structure over $k$. 
\smallskip
In other words, the abstract Expectation of 3) of 3.13.7.4 D has often significant and extremely important (adaptable to geometric) interpretations (similar to the ones of 1.6.1 and referring to isomorphism classes) in the much wider context of $\Mc\Ml_3$: see 3.13.7.8 for details; see 3.13.7.8.1 for a very concrete example. See also [Va12] for geometric examples involving other classes of polarized varieties (like polarized K3 surfaces or cubic fourfolds, etc.). Here we will just mention that in 2.2.14.2 we introduce Fontaine truncations: they provide the right context for such extensions and interpretations, leading to the most general possible Faltings--Shimura--Hasse--Witt stratifications. See 3.13.7.9 for the general form of Faltings--Shimura--Hasse--Witt (all types of) maps. See also 3.13.7.8.3-4, 4.5.15.2.4 3), 4.14.5 and Appendix for some extra useful aspects and comments. 
\medskip
{\bf 1.6.3. The boundedness principle and some new stratifications.} In all that follows we speak only about stratifications of $\Mn_{k(v)}$ and we allow them to have an infinite number of strata. Accordingly, till the end of 1.6.3 we take $k=\FF$ and we consider (see the conventions of 2.1) for the simplicity of language, only those strata which are locally closed subschemes of $\Mn_{k(v)}$. In 4.5.16 and 4.9.9.1 we start a deep study of the ultra stratification of $\Mn_{k(v)}$. We introduce two extra types of stratifications: the Faltings--Shimura--Dieudonn\'e and the ultimate stratifications. Each one of these types has three variants: adjoint, principal or standard; warning: most common we drop the word standard. 
\smallskip
The Faltings--Shimura--Dieudonn\'e principal stratification (resp. adjoint stratification) is the stratification of $\Mn_{k(v)}$ defined by (warning!) inner isomorphism classes of principally quasi-polarized Shimura $\sg$-crystals (resp. of Shimura adjoint Lie $\sg$-crystals) attached to $k$-valued points of $\Mn_{k(v)}$. Related to the use of the word inner we mention two things. First, in the adjoint context it is 4.2.10 which allows us to speak about inner isomorphisms between Shimura adjoint Lie $\sg$-crystals attached to $k$-valued points of $\Mn_{k(v)}$. Second, in the principally quasi-polarized Shimura $\sg$-crystals, the fact that the subgroup of $G_{W(k)}$ fixing $p_A$ of 1.3.1 is connected (see the part of 2.3.3 referring to $G^0_{\ZZ_p}$), we can speak about isomorphism classes of quadruples $(M,\vph,(t_{\al})_{\al\in\Mj},p_A)$ (of elements) as in 1.3.1; they are referred as inner isomorphism classes of principally quasi-polarized Shimura $\sg$-crystals attached to $k$-valued points of $\Mn_{k(v)}$. The Faltings--Shimura--Dieudonn\'e (standard) stratification is defined similarly to the Faltings--Shimura--Dieudonn\'e principal stratification but without keeping track of principal quasi-polarizations. The Faltings--Shimura--Dieudonn\'e (standard or adjoint) stratification is $G(\AA_f^p)$-invariant (see Fact 6 of 2.3.11) and (cf. 4.9.9) invariant under the automorphism group ${\rm Aut}(G,X,H)$ (see 2.3.5.4 for a review on it). In many cases (like when the center of $G$ is a torus and the connected components of $\Mn_{k(v)}$ are permuted transitively by $G(\AA_f^p)$), the Faltings--Shimura--Dieudonn\'e stratification and the Faltings--Shimura--Dieudonn\'e adjoint stratification coincide but we are not able to list all cases when this happens. Also, if these last two stratifications coincide and if $({{-1}\over p})=-1$, then they also coincide with the Faltings--Shimura--Dieudonn\'e principal stratification; we do not know what to say if $({{-1}\over p})=1$. See 4.5.16.0 2) for the last two sentences. Warning: the principal quasi-polarizations (as $p_A$ at the end of 1.3.1) are viewed here as fixed cycles and not up to $\GG_m(\ZZ_p)$-multiples (of such fixed cycles).
\smallskip
The ultimate stratification of $\Mn_{k(v)}$ is the $G(\AA_f^p)$-invariant, ${\rm Aut}(G,X,H)$-invariant maximal refinement of the Faltings--Shimura--Dieudonn\'e stratification obtainable by decomposing its strata into open closed subschemes. One can define as well the ultimate principal stratification of $\Mn_{k(v)}$: we just have to replace $G(\AA_f^p)$ and ${\rm Aut}(G,X,H)$ by their maximal subgroups leaving invariant the Faltings--Shimura--Dieudonn\'e principal stratification; as we do not know how to compute these groups, in \S 4 we do not mention it. For getting nice functorial purposes (see end of 1.10.1), the ultimate adjoint stratification of $\Mn_{k(v)}$ is obtained similarly (but slightly different) from the  Faltings--Shimura--Dieudonn\'e adjoint stratification of $\Mn_{k(v)}$; see 4.9.9.1.
\smallskip
So the ultimate adjoint stratification is a refinement of the Faltings--Shimura--Dieudonn\'e adjoint stratification, which is a refinement of the ultra stratification, which at its turn is a refinement of the quasi-ultra or of the absolute stratification.  Also, there is nothing we can refine to the ultimate stratification, from any point of view pertaining to Shimura $p$-divisible groups attached to points of $\Mn_{k(v)}$. These new two types of stratifications have one disadvantage: in general their number of strata is infinite; however, see rm. 4.9.9.3. We have (cf. 4.5.16 and 4.5.16.1):
\medskip
{\bf Theorem 6 (the third form of the purity principle and a geometric form of the boundedness principle).} {\it The Faltings--Shimura--Dieudonn\'e (principal) stratification is well defined, i.e. its strata are indeed reduced, locally closed subschemes of (suitable pull backs of) $\Mn_{k(v)}$. It satisfies the purity property (as defined in 2.1) and each stratum of it is regular, quasi-affine and all its connected components have the same dimension.}
\medskip
We have versions of this (see 4.5.16, 4.5.16.1, 4.9.9 and 4.9.9.1) for the Faltings--Shimura--Dieudonn\'e adjoint stratification or for the ultimate types of stratifications. The idea of the proof of Theorem 6 is very simple: its very essence is the Fundamental Lemma 3.6.15 B. The simplest form of 3.6.15 B asserts that for any Shimura $\sg$-crystal $(M,\vph,G)$ and for any $h\in G(W(k))$, there is $n\in\NN$ such that $(M,h\vph,G)$ and $(M,gh\vph,G)$ are isomorphic under an isomorphism defined by an element $e(h,g)$ of $G(W(k))$, $\forall g\in G(W(k))$ congruent to the identity mod $p^n$. Its proof involves some fundamental techniques: in general we can not choose $e(h,g)$ to be congruent to the identity mod $p^n$; so in order to be able to use the fact that any quasi Artin--Schreier system of equations with coefficients in $k$ has a solution in $k$, we use the stairs method, i.e. we consider suitable $W(k)$-submodules of ${\rm Lie}(G)$ which (warning!) do not pay at all attention to the Lie structure but are obtained using Dieudonn\'e's classification (see [Di] or [Man]) of isocrystals over $k$ (so 3.6.15 B is true even in the context of $\Mc\Ml_5$; see rm. 1) of 3.15.7 B). Simple estimates (see 3.15.7 C) show that $n$ can be taken to be independent on $h$. This is the main ingredient needed to get (via standard techniques of algebraic geometry) that the Faltings--Shimura--Dieudonn\'e (principal) stratification of $\Mn_{k(v)}$ is well defined. The estimates of 3.15.7 C handle the most general context of (see 2.1) latticed isocrystals over $k$. 
\smallskip
The essence of the above paragraph is presented as the boundedness principle in 3.15.7. See 3.15.8 for the homomorphisms form of 3.6.15 B; it is based on the fact that we can always choose $e(h,g)$ to be congruent to the identity mod $p^{n-m}$, where (for $n$ bounded below as in 3.6.15 B) $m\in\NN$ does not depend as well on $h$. Using it, we get (see 3.15.8) a new proof of ``Katz' open part" of the Grothendieck--Katz' specialization theorem. 
\smallskip
The third form of the purity principle can be checked immediately using schemes of isomorphisms between two finite, flat, commutative group schemes over the same base. Even more, in 3.15.10 we show how the homomorphism form of 3.6.15 B can be used to get the second form of the purity principle, i.e. to get that all stratifications by Newton polygons defined by (see 2.1 for conventions on them) $F$-crystals in locally free sheaves of locally finite ranks over suitable $\FF_p$-schemes $S$, satisfy the purity property. This was obtained previously under the hypothesis that $S$ is locally noetherian in [dJO, 4.1]. The philosophy (which among other things allows a weakening of this local noetherian assumption; see 3.15.10 and 3.15.10.1 1) and 2)) is: 
\medskip
{\bf Ph.} {\it Newton polygons are entirely ``encoded" in the existence of suitable morphisms between truncations modulo effectively computable (positive, integral and high) powers of $p$ of $F$-crystals in locally free sheaves of locally finite ranks.} 
\medskip
Our interest in Newton polygons is considerably reduced (cf. the Real Problem of 1.6.5); so this second form of the purity principle is presented at the very end of \S 3. It is worth pointing out that in 3.15.10 we use only tools which were already available in print in 1979 (i.e. we use besides 3.6.15 B --its proof uses the simplest type of quasi Artin--Schreier systems of equations with coefficients in $k$ which obviously have solutions over $k$-- just [Ka2, 2.7.4]). One might wonder why it took 21 years to get [dJ0, 4.1] (and that by using strong tools --like alterations, overconvergent $F$-crystals, etc.). We think the reason is: too little mathematical effort was put in connection to the Real Problem (according to us, most of what has been done previously from the rational point of view --like Newton polygons, etc.-- has to be redone entirely --unfortunately, often starting from the very scratch-- from the perspective of the Real Problem). 
\smallskip
We came across Theorem 6 being motivated by the integral Manin problem (see 1.12), cf. 3.6.15 A. There is extremely few literature in connection to it. The most we can say: [LO] has some material which can be used to have more computations performed (unfortunately, very little of loc. cit. can be adapted to the context of Shimura $p$-divisible groups; moreover, loc. cit. can not be used in connection to $\Mc\Ml_2$, provided we are in a context involving the Shimura Lie $\sg$-crystals mentioned in 1.2.1.1). We refer to \S 9-10 for our computations involving different stratifications we mentioned above. Warning: we presently have no reason to think that in general the ultimate (adjoint) and the Faltings--Shimura--Dieudonn\'e (adjoint) stratifications of $\Mn_{k(v)}$ coincide; [IKO, 2.10] is a nice result pertaining to isomorphism classes of principal polarizations on products of supersingular elliptic curves over $k$, which tells us to be cautious in trying to identify these last two (adjoint) stratifications.
\medskip
{\bf 1.6.4. The First Main Corollary.} Theorem 6 and its abstract analogues represent, in our opinion, a turn in the study of (Shimura) $p$-divisible groups: we need to build up their ``libraries" (i.e to list them). For instance, the following consequence of 3.6.15 B tells us that it is possible to build an ``atlas" for isomorphism classes of Shimura $p$-divisible groups over $\FF$ which are defined over a fixed finite field, similar to the well known ``atlas" for finite groups or to [Cr]. We have:
\medskip
{\bf First Main Corollary.} {\it Let $q\in\NN$. Let $\FF_{p^q}$ be the field with $p^q$ elements. The number $D(r,d,q)$ of isomorphism classes of $p$-divisible groups of fixed rank $r\in\NN$ and dimension $d\in\{0,...,r\}$ over $k$ which are definable over $\FF_{p^q}$ is finite and does not depend on $k$. Similarly, for any Shimura $\sg$-crystal over $k$, the number $D(Cl(M,\vph,G),q)$ of elements of the class $Cl(M,\vph,G)$ (see 2.2.22 2); it is a set formed by isomorphism classes) definable over $\FF_{p^q}$ is finite and does not depend on $k$.}
\medskip
See 3.15.7 G for concrete (though very gross) upper bounds for $D(r,q)$ or for $D(Cl(M,\vph,G),q)$. For instance, we get $D(r,q)\le p^{qr^2h(r)}$, where $h:\NN\to\NN$ is an effectively (and easily) computable function. 
\medskip
{\bf 1.6.5. Grothendieck's specialization categories.} Presently we have full control on the Expectation of 1.6.2 for Shimura varieties of $A_n$, $C_n$ or $D_n^{\HH}$ Lie types (cf. 4.6.7-8). So, based on this and the last sentence of 1.6.4, in 4.9.9.2 we already introduce Grothendieck's (adjoint) specialization category of $\Mn_k$ (or of $\Mn_{k(v)}$) which keeps track of specializations between strata of the ultimate (adjoint) stratification of $\Mn_{k(v)}$; some of the most natural invariants associated to it are mentioned in 4.12.12.8. 
The philosophy behind this category can be roughly summarized as follows.
\medskip
{\bf The Real Problem.} {\it Let $p\ge 2$. The real problem of the theory of Shimura $p$-divisible groups over some convenient perfect field $\tilde k$ (usually finite or algebraically closed) of characteristic $p$, is to list (often using families) their isomorphism classes and decide which such isomorphism classes specialize to which other. The same applies to the theory of $p$-divisible objects with tensors over $\tilde k$ (we mostly have in mind the context of $\Mc\Ml_3$ but not only).} 
\medskip
Not to make this paper too long, we include here extremely few applications of Theorem 6 (see \S 9-11 for many such applications, cf. also 1.15.6). A last think: we have a result similar to Theorem 6 in the abstract context of $\Mc\Ml_2$ (see [Va6]; see also 3.15.6 D and E).
\medskip
{\bf 1.6.6. (Refined) Artin--Schreier stratifications.} Each quasi Artin--Schreier system of equations with coefficients in a reduced $\FF_p$-algebra $R$, defines naturally (based on 3.6.8.1.2 a)) an Artin--Schreier and a refined Artin--Schreier stratification of ${\rm Spec}(R)$ in a finite number of strata, cf. 3.6.8.1.3. In 4.5.18, we extend these notions of types of stratifications to an arbitrary reduced $\FF_p$-scheme: examples (abstract as well as concrete ones pertaining to $\Mn_{k(v)}$) and comments are included. We draw attention to the Conjecture of 4.5.18.1. It is related to the first form (pertaining to quasi Artin--Schreier systems of equations; see 3.6.8.1.4) of the purity principle: it conjectures that the second and third form of this principle are particular cases of its first form. 
\medskip
{\bf 1.6.7. Toric points.} In 4.5.11 we enlarge the class of $U$-ordinary points: we introduce toric points; the $k$-valued toric points of $\Mn_{k(v)}$ are defined by requiring that their attached Shimura $\sg$-crystals have lifts of quasi CM type, without imposing the uniqueness of such lifts. 
\medskip\smallskip
{\bf 1.7. On the contents of 4.3 and 4.6.} In 4.3 we restate 1.2.2 in the context of a standard Hodge situation as in 1.3. We use refined Lie stable $p$-ranks: each simple factor of $G^{\rm ad}_{\ZZ_p}$ can be used (cf. 4.2.10) to define (independent) Lie stable $p$-ranks of its corresponding cyclic adjoint factors (see 3.10.1 for their definition) of Shimura $\sg$-crystals attached to points of $\Mn_{k(v)}$ with values in $k$.
\smallskip
In 4.6 we include some examples to illustrate the general theory presented in 4.1-5 and to point out the complexities which are arising. One sample: in 4.6.1 1) we point out that it can happen that the formal isogeny type $\tau$ introduced in 1.4 does not have integral slopes. We also list the main properties of the general theory; some particular situations are also referred to (for instance, see 4.6 P11 for the case when the centralizer of $G$ in $GSp(W,\psi)$ is a 1 dimensional torus). As a simple and straightforward application of 1.4 b) and of [Va2, 5.8.6] we get a general principle concerning non-ordinary reductions of an abelian variety $A$ over a number field $E$, with respect to primes of $E$:
\medskip
{\bf Theorem 7 (the non-ordinary reduction criterion).} {\it If the reflex field of the Shimura variety defined naturally (see 2.1) by the Mumford--Tate group of $A$, is not $\QQ$, then $A$ has a non-ordinary reduction with respect to an infinite number of primes of $E$. Moreover, there are very precise recipes of how ``to cook" such primes (see 4.6.2.1; see also 4.6.2.2).}
\medskip
4.6.5-8 contain different complements meant to emphasize how complexities or simplifications arise depending on the type of $f$, on $k(v)$ and of the types of the simple factors of $(G^{\rm ad},X^{\rm ad})$. In connection to 4.6.4 see 1.14.3. 
\medskip\smallskip
{\bf 1.8. Crystalline coordinates.} We consider a standard Hodge situation $\bigl(f,L_{(p)},v\bigr)$ as in 1.3. Let $v^{\rm ad}$ be the prime of $E(G^{\rm ad},X^{\rm ad})$ divided by $v$ and let $d$ be the dimension of $X$ as a complex manifold. Let $y:{\rm Spec}(k)\hookrightarrow\Mn_{W(k)}$ be defined by a $G$-ordinary point $y_0\colon {\rm Spec}(k)\to\Mn_{k(v)}$ and let $(A_y,p_{A_y}):=y_0^*(\Ma,\Mp_{\Ma})$. The principally quasi-polarized filtered $F$-crystal (associated to the natural pull back of ($\Ma,\Mp_{\Ma})$) over the completion $k[[x_1,\ldots,x_d]]$ of the local ring of $y$ in $\Mn_{k}$, can be put in a very practical form (cf. 4.7.2). Using this we get:
\medskip
{\bf Theorem 8.} {\it The moduli formal scheme $\Mm_y$ of $G$-deformations (i.e. of deformations factoring through $\Mn$) of $(A_y,p_{A_y})$ over spectra of artinian local $W(k)$-algebras having $k$ as their residue field, is naturally isomorphic to the formal torus of the completion $\Mt_d$ of $\GG_m^d$ (viewed over $W(k)$) in the special fibre of its identity element section. Under such an isomorphism, the principally polarized abelian scheme over ${\rm Spec}(W(k))$ obtained by pull back via the $G$-canonical lift of $y_0$, corresponds to the identity element of $\Mt_d(W(k))$.} 
\medskip
Theorem 8 and its logical (abstract) variants in contexts without (principal) polarizations (see 4.7.11 6)) represents a generalization of the following three things:
\medskip
$\bullet$ the main result of [No1] (in the less general case of a standard Hodge situation);
\smallskip
$\bullet$ [De3, 2.1.3];
\smallskip
$\bullet$ [Ka3, 3.7.1-3 and 4.3.1-2] (see 4.7.11 8); for $p=2$ cf. also with 4.14.3 E). 
\medskip
So we obtain $G$-multiplicative (crystalline) coordinates of $\Mm_y$. When $k(v^{\rm ad})=\FF_p$, the situation is very close in spirit to the one in [De3, 1.4] (cf. 4.7.5-29; see also its abstract generalization of 4.7.11 4) and 8)): we obtain $G$-canonical multiplicative (crystalline) coordinates. In order to make the connection with previous works, the proofs of 4.7.5-22 follow the pattern of [De2, 1.4.2 and 1.4.7] (cf. 4.7.10). 
\smallskip
Theorem 8 is a very particular case of 4.7.11 6) and so is part of a general theory of crystalline coordinates for Shimura $p$-divisible groups over $k$ presented in 4.7.11. It is the inducing property of 3.6.18.5 --and its variants of 3.6.18.7 in the relative context-- which allows significant simplifications and advancements in such a theory: our approach in 4.7.11 is ``cook and pick". By this we mean: we consider a priori a $p$-divisible object with tensors of some (logical) Fontaine category $\Mm\Mf_{[0,1]}(W(k)[[x_1,...,x_m]])$, $m\in\NN$, and then we pick up a connection (see 3.6.18.4 and 3.6.18.7.0) which makes it to be viewed as a $p$-divisible object with tensors of $\Mm\Mf_{[0,1]}^\nabla(W(k)[[x_1,...,x_m]])$ and whose Kodaira--Spencer map is injective (cf. a) and b) of 4.7.11 2)); the ``cooking" is done in such a way that the ``picking" is possible (it is enough to be able to ``perform the picking" mod $p$, cf. the lifting property). This approach, in most cases (like the ones of the classical situation of [De2] and [Ka3-4] or like the ones pertaining to the case $k(v^{\rm ad})=\FF_p$), avoids entirely any computation (even for $p=2$, cf. 4.14.3 E).
\smallskip
As a particular case of this general theory we have (see 4.7.14.1 for arguments):
\medskip
{\bf Corollary.} {\it If $A_k$ is an abelian variety over $k$ of dimension $d_0$ having a lift $A$ (over $W(k)$) of quasi CM type (see def. 2.2.17), then the moduli formal scheme of deformations of $A_k$ is naturally isomorphic to the formal torus $\Mt_{d^2_0}$ obtained by completing $\GG_m^{d^2_0}$ (viewed over $W(k)$) in its identity element section, with the abelian variety $A$ corresponding to the identity element of $\Mt_{d^2_0}(W(k))$.} 
\medskip
This is the most general (abelian varieties) context in which we can (presently) prove the existence of ``good" (i.e. useful and natural) crystalline coordinates. For variants of Theorem 8 and of this Corollary, including the additive case, see 4.7.14.2 and 4.7.11 5) and 6); here, for instance, in the context of Theorem 8, by the additive case we mean that the Frobenius lift of $W(k)[[x_1,...,x_d]]$ is as in 4.7.1 and we try to get even more practical forms of 4.7.2. The main merit of the Theorem of 4.7.11 6) is: its simple form works out in all cyclic diagonalizable cases.
\medskip
{\bf 1.8.1. Complements.} For a list of the new types of equations (differential or not) showing up in the theory of crystalline coordinates when we move from a classical ordinary context to a cyclic diagonalizable one, see end of 4.7.11 1). For a digression on crystalline coordinates in the context of an ordinary K3 surface over $k$ (including the case $p=2$) see 4.7.11 2A). For a sample of the lifting process (to be fully elaborated in \S 7 and [Va12]) which allows the passage (of many properties) from a Shimura (filtered) $\sg$-crystal context to the attached Shimura (filtered) Lie $\sg$-crystal context, see 4.7.11 2B). In 4.14.2 we point out that all of 4.7.11 can be performed in the context of $\Mc\Ml_2$, without any reference (even if possible via a lifting process) to $p$-divisible groups. 
\medskip\smallskip
{\bf 1.9. A Galois property.} In 4.8 we present the following Galois property of $G$-ordinary points of $\Mn_{k(v)}$; it answers a question of
C. -L. Chai. We consider a standard Hodge situation $(f,L_{(p)},v)$. With the notations of 1.3, we consider a $G$-ordinary point ${\rm Spec}(k)\to\Mn_{k(v)}$ and
a lift of it $z:{\rm Spec}(V)\to\Mn$, with $V$ a DVR which is a finite, flat extension of $W(k)$. Let 
$K:=V\bigl[{1\over p}\bigr]$ and let $A_V:=z^*(\Ma)$. We have (see 4.8.2):
\medskip
{\bf Theorem 9.}
{\it There is a finite field extension $K_1$ of $K$ such that the Galois representation 
$$\rho:{\rm Gal}(\bar K/K_1)\to GL(H^1_{\acute et}(A_{\bar K},\ZZ_p))(\ZZ_p)=
GL(L_{(p)}^*\otimes_{\ZZ_{(p)}} \ZZ_p)(\ZZ_p)$$ 
 factors through the group of $\ZZ_p$-valued points of an integral, solvable subgroup of $G_{\ZZ_p}$.}
\medskip
The  identification $H^1_{\acute et}(A_{\bar K},\ZZ_p)=
L_{(p)}^*\otimes_{\ZZ_{(p)}} \ZZ_p$ used here is not canonical, cf. [Va2, top of p. 473]. The proof of Theorem 9, as R. Pink pointed  out to us, is a direct consequence of Theorem 4 and of 4.7.2 (and so of the proof of Theorem 8).
For variants of this Galois property (including the case of a finite field) see 4.8.3 d) and e). 
\medskip
{\bf 1.9.1. The abstract generalized Serre--Tate (ordinary) theory.} a) to c) of Theorem 1, Theorems 2 and 4, the abstract form of Theorem 9 (see 4.8.3 e)) as well as 3.1.1.2, 3.1.4-5, 3.9.3, 3.9.4 and 3.13.7.1.2 form the abstract generalization of the Serre--Tate's (ordinary) theory for Shimura $\sg$-crystals (or Shimura $p$-divisible groups) over $k$; see 1.13 for the case $p=2$. See Step 1 of Appendix for the context of $\Mc\Ml_2$.
\medskip\smallskip
{\bf 1.10. On the passage from the Hodge type to the preabelian type.} The phenomena of 1.4-8 involving $\Mn$ are an intrinsic ``property" of $\Mn$ (and of
$\Mn_{k(v)}$), i.e. do not depend on the special principally polarized abelian scheme $(\Ma,\Mp_{\Ma})$ we got over $\Mn$ (only $\tau$ of 1.4 a) does depend on it); so they ``hold" for (i.e. can be transferred to) special fibres of integral canonical models of Shimura varieties of preabelian type. For instance, we have (cf. 4.9.8 and 4.9.8.2) the following combined version of the density property of Shimura-ordinary strata and of the existence of refined canonical Lie stratifications and of Shimura-canonical lifts.
\medskip
{\bf Theorem 10.}
{\it Let $\Mn^1_{k(v^1)}$ be the special fibre of an integral canonical
model $\Mn^1$ of a Shimura variety ${\rm Sh}(G^1,X^1)$ of preabelian type with respect to a
prime $v^1$ of $E(G^1,X^1)$ dividing a rational prime $p>2$. Then $\Mn^1_{k(v^1)}$ has a refined canonical Lie stratification in $G^1(\AA^p_f)$-invariant, reduced, locally closed subschemes; it has precisely 1 open stratum (if $\dim_{\CC}(X^1)\ge 1$ then its complement is of pure codimension 1). In the case of a standard Hodge situation we regain the refined canonical Lie stratification mentioned in 1.6. The points (with values in fields) of the open stratum are called
$G^1$-ordinary (or Shimura-ordinary) points of $\Mn^1_{k(v^1)}$; any such point with values in 
$k$, has a uniquely determined (in the same --parabolic type-- way as in 1.5 but expressed over $\bar k$) $G^1$-canonical (also called Shimura-canonical) lift ${\rm Spec}(W(k))\to\Mn^1$.
\smallskip
The refined canonical Lie stratifications and the Shimura-canonical lifts are functorial with respect to finite maps (see def. [Va2, 2.4 and 3.2.6]) between Shimura quadruples.}
\medskip
In 4.9 we present a proof of this for primes $p\ge 5$ (for the case $p=3$ see \S 6).
The proof relies on the following result. 
\medskip
{\bf The Existence Property.} {\it We consider a Shimura quadruple $(G_0,X_0,H_0,v_0)$, with $(G_0,X_0)$
defining a Shimura variety of adjoint, abelian type, with $v_0$ a prime of $E(G_0,X_0)$
dividing a rational prime $p\ge 3$ such that $G_0$ is unramified over $\QQ_p$, and with $H_0$ a hyperspecial subgroup of $G_0(\QQ_p)$. Then 
there is a standard Hodge situation $(f^0,L^0_{(p)},v^0)$, with 
$$f^0:{\rm Sh}(G^0,X^0)\hookrightarrow {\rm Sh}({\rm GSp}(W,\psi),S)$$ 
an injective map such that $(G^{0{\rm ad}},X^{0{\rm ad}})=
(G_0,X_0)$ and $v^0$ divides $v_0$, and with $G^{0\rm der}$ as the maximal isogeny cover of $G^{0\rm ad}$ allowed by the abelian type (for $p\ge 5$, cf. [Va2, 6.4.2] and the translation of [Va2, 5.1, 5.6.5 and 5.6.9] --see 2.3.6-- in the language of standard Hodge situations; for $p=3$, cf. \S 6). Moreover, if we define a hyperspecial subgroup $H^0$ of $G^0(\QQ_p)$ as in 1.3, we can assume (cf. [Va2, 3.2.7.1]) that $H^{0\rm ad}=H_0$.}
\medskip
In order to get Theorem 10 from this result, we need to use a diagonal trick (see 4.9.2): its goal is to ``combine" (i.e. put together) two standard Hodge situations involving two Shimura quadruples having the same adjoint, in order to show the independence on which such standard Hodge situation we use to define (or get) some properties. The diagonal trick itself is supported:
\medskip
a) by the density part of 1.4 b);
\smallskip
b) by a variant (see 4.9.2.0 A) of 2.3.13 in a non a priori reductive context;
\smallskip
c) by the direct sum decompositions of 3.11.8.1;
\smallskip
d) and by Fontaine's comparison theory (used --for $p=3$-- in a way similar to the detailed one of the review --with $p=2$-- 2.3.18.1 E but in the adjoint context). 
\medskip
For the $G^1(\AA^p_f)$-invariant part of Theorem 10 we use the strong approximation theorem for adjoint $\QQ$--groups (see 4.9.7.1 for a second possible approach).
\medskip
{\bf 1.10.1. Extension of terminology.} Strictly speaking 1.10 a) to d) are not needed to prove Theorem 10, as in connection to Newton polygons we can work just rationally. However, they are indispensable for integral purposes. In other words, using 1.10 a) to d) and standard arguments on right translations by Hecke operators, we define (see 4.9.17) Shimura adjoint (resp. Shimura adjoint filtered) Lie $\sg$-crystals attached to points of $\Mn^1_{k(v^1)}$ with values in algebraically closed fields (resp. with values in Witt rings of such fields); so we define also toric points (in particular, $U$-ordinary and $T$-ordinary points) of $\Mn^1_{k(v^1)}$ (see 4.9.17.2 and 4.9.17.5), as well as $U$-canonical lifts of $U$-ordinary points with values in perfect fields. Moreover, Theorem 10 has variants (see 4.9.9) in other contexts, like the one of the Faltings--Shimura--Dieudonn\'e adjoint (or the quasi-ultra, or the ultimate adjoint, etc.) stratifications; again, for these variants 1.10 a) to d) are indispensable. 
\medskip\smallskip
{\bf 1.11. Dense Hecke orbits and functorial aspects.} We do expect (cf. 4.10.1; cf. also [Ch2, p. 441]) that the open stratum of $\Mn^1_{k(v^1)}$ is the smallest $G^1(\AA^p_f)$-invariant, open subscheme of $\Mn^1_{k(v^1)}$. This is proved for $\Mm_{\FF_p}$ in [Ch2]. Using Theorem 10, we prove it for Shimura curves  (cf. 4.10.3). For a more substantial approach to the general case see 4.10.5.
\smallskip
In 4.11 we prove the functorial behavior of $G$-ordinary points and of $G$-canonical lifts of them.
Let $f_0:(G,X,H,v)\hookrightarrow (G_1,X_1,H_1,v_1)$ be an injective map between two Shimura quadruples of Hodge type, having integral canonical models $\Mn$ and respectively $\Mn_1$. We assume $v$ is relatively prime to 2. We have a natural morphism (still denoted by) $f_0:\Mn\to\Mn_1$ (cf. [Va2, 3.2.7 4)]). We also assume the existence of an injective map 
$$f_1:(G_1,X_1,H_1,v_1)\hookrightarrow\bigl({\rm GSp}(W,\psi),S,{\rm GSp}(W,\psi)(L_{(p)}\otimes_{\ZZ_{(p)}} \ZZ_p),p\bigr)$$
 (with $L_{(p)}$ a $\ZZ_{(p)}$-lattice of $W$ for which we have a perfect form  $\psi:L_{(p)}\otimes_{\ZZ_{(p)}} L_{(p)}\to\ZZ_{(p)}$)  such that $(f_1,L_{(p)},v_1)$ and $(f_1\circ f_0,L_{(p)},v)$ are standard Hodge situations. So the Zariski closure of $G$ (resp. of $G_1$) in ${\rm GSp}(L_{(p)},\psi)$ is a reductive group $G_{\ZZ_{(p)}}$ (resp. ${G_1}_{\ZZ_{(p)}}$) over $\ZZ_{(p)}$.
Let $k(v)$ and $k(v_1)$ be the residue field of $v$ and respectively of $v_1$. We have (cf. 4.11.1.1 and 4.11.2):
\medskip
{\bf Theorem 11.} {\it A $G$-ordinary
point of $\Mn_{k(v)}$ is mapped by $f_0$ into a $G_1$-ordinary point of $\Mn_{1k(v_1)}$ iff there is a torus $T$ of $G_{\ZZ_p}$ such that:
\medskip
1) there is a cocharacter
$\mu:\GG_m\hookrightarrow T_{W(k(v))}$, which over an embedding of $W(k(v))$ in $\CC$, extending the composite inclusion $O_{(v)}\subset E(G,X)\subset\CC$, is $G(\CC)$-conjugate to the cocharacters
$\mu^\ast_x$, $x\in X$, defined as in 1.4 (see also 2.3.1);
\smallskip
2) there is a parabolic subgroup $P_1$ of ${G_1}_{\ZZ_p}$ containing $T$ and a Borel subgroup of $G_{\ZZ_p}$ and whose Lie algebra has the property that its elements take the $F^1$-filtration of $L^\ast_{(p)}\otimes_{\ZZ_p} W(k(v))$ defined by $\mu$ (as in 1.4 a)) into itself.
\medskip
Moreover, if these conditions hold, then a $G$-canonical lift of $\Mn$ is mapped by $f_0$ into a $G_1$-canonical lift of $\Mn_1$.} 
\medskip
The equivalent conditions of Theorem 11 imply $k(v)=k(v_1)$ but the converse of this does not hold (see 4.11.4-6). The proof of Theorem 11 relies on 1.4 a), 4.4.13.2 and on standard arguments on reductive groups over $\ZZ_p$ or $W(k)$.
\medskip\smallskip
{\bf 1.12. Some important applications.} Till 1.13 we assume $k=\bar k$. In 4.12-13 we present two important applications inspired from the deformation principle of 3.6.14. First, in 4.12 we formulate the integral Manin problem for a standard Hodge situation as well as for integral canonical models of Shimura varieties of preabelian type with respect to primes relatively prime to 2. For a standard Hodge situation $(f,L_{(p)},v)$, with the notations of 1.3, the problem is to determine all principally quasi-polarized Shimura $\sg$-crystals attached to $k$-valued points of (a connected component of) $\Mn_k$. A solution (cf. \S 11; see the proof of 4.12.12 for some particular cases) to this problem is provided by Milne's conjecture (see d) of 4.4.1 3) and 1.15.1) and by the completion property. The completion property is stated in 3.6.15 A as an expectation and will be proved in \S 11, starting from the isogeny property to be referred in 1.15.7 and from the strong form of the generalized Manin problem (see B) of Theorem 14); see the proof of 4.12.12 and 4.12.12.0 for samples (in contexts in which the isogeny property can be checked immediately). For the initial motivations behind the completion property, see 3.6.15 B and 3.6.14.4.1. In the language of 1.6.5, the completion property (over $\FF$) says: the Grothendieck's specialization category of $\Mn_k$ has as many objects as we can think of. 
\smallskip
There is a second approach to the completion property of 1.12: see 4.12.12.7.  
\smallskip
For different local (resp. global) integral Manin problems see 4.12.14-19 (resp. 4.12.21). For a new (though very close in spirit to [Ta]) and simple solution to the original Manin problem we refer to 4.12.12.4.
\smallskip
Second, in 4.13 we present (partially; to be continued in \S 9-10, cf. 1.15.6) the invariance principle. Let $p\ge 3$ be a rational prime. Let $(G_i,X_i,H_i,v_i)$, $i=\overline{1,2}$, be two Shimura quadruples of preabelian type, with $v_i$ dividing $p$. Let $\Mn_i$ be the integral canonical model of $(G_i,X_i,H_i,v_i)$, $i=\overline{1,2}$. Let $H_{0i}\subset G_i(\AA_f^p)$ be a subgroup such that the subgroup $H_{0i}\times H_i$ of $G_i(\AA_f)$ is smooth for the pair $(G_i,X_i)$ (see def. [Va2, 2.11]) and an extra mild condition holds (satisfied for instance, if the image of $H_{0i}$ in $G^{\rm ad}_i(\QQ_l)$, for some prime $l$ different from $p$, has no pro-$p$ subgroup), $i=\overline{1,2}$. Let $y_i:{\rm Spec}(k)\to\Mn_{ik}/H_{0i}$, $i=\overline{1,2}$, be points such that the non-trivial part of the Shimura adjoint Lie $\sg$-crystal attached to $y_1$ is isomorphic to the non-trivial part of the Shimura adjoint Lie $\sg$-crystal attached to $y_2$ (see 2.3.10 and 3.10.1 for definitions). Let $O_{y_i}$ be the local ring of $y_i$ (in $\Mn_{iW(k)}$) and let $O_{y_i}^{h\wedge}$ be the $p$-adic completion of its henselization, $i=\overline{1,2}$. The refined canonical Lie stratification of $\Mn_{ik(v_i)}$ makes ${\rm Spec}(O(p)_{y_i}^h):={\rm Spec}(O_{y_i}^{h\wedge}/pO_{y_i}^{h\wedge})$ to be a stratified scheme. One variant of the local form of the invariance principle is (see 4.13.1; see 4.13.2 and 4.13.3 3) for other variants): 
\medskip
{\bf Theorem 12.} {\it Let $n\in\NN$. There is an isomorphism 
$$\rho:{\rm Spec}(O_{y_1}^{h\wedge})\tilde\to {\rm Spec}(O_{y_2}^{h\wedge})$$
inducing an isomorphism 
$$\rho(p):{\rm Spec}(O(p)_{y_1}^h)\tilde\to {\rm Spec}(O(p)_{y_2}^h)$$
of stratified schemes, such that the pull back through $\rho$ of the non-trivial part of the Shimura adjoint filtered Lie $F$-crystal over ${\rm Spec}(O(p)_{y_2}^h)$ defined by the natural morphism ${\rm Spec}(O_{y_2}^{h\wedge})\to\Mn_2$, when taken mod $p^n$, is isomorphic to the non-trivial part taken mod $p^n$, of the Shimura adjoint filtered Lie $F$-crystal over ${\rm Spec}(O(p)_{y_1}^h)$ defined by the natural morphism ${\rm Spec}(O_{y_1}^{h\wedge})\to\Mn_1$ (see 4.9.20 and 4.13 for the definition of these Lie $F$-crystals).}
\medskip 
Theorem 12 is a consequence of variants modulo powers of $p$ of Theorem 5.
\medskip
As a direct application of what the completion property says in some particular cases in which $G^{\rm ad}$ is of $A_n$ type, for future references, we state here the combined version of 4.12.12, 4.12.12.2 and 4.12.13 2). We have:
\medskip
{\bf Theorem 13.} {\it We assume $p\ge 2$. Let $(d,r-d)$ be a pair of non-negative integers. There is $q\in\NN$, a geometrically connected, projective, smooth scheme $\Ml_{(r,d-r)}$ over $W(\FF_{p^q})$ and a $p$-divisible group $D_{\Ml_{(r,d-r)}}$ over $\Ml_{(r,d-r)}$ of rank $r$ and dimension $d$ such that:
\medskip
{\bf a)} $D_{\Ml_{(r,d-r)}}$ is versal in each $k$-valued point of $\Ml_{(r,d-r)}$;
\smallskip
{\bf b)} the relative dimension of $\Ml_{(r,d-r)}$ over $W(\FF_{p^q})$ is exactly $d(r-d)$;
\smallskip
{\bf c)} any $p$-divisible group over $W(k)$ of rank $r$ and dimension $d$ is obtained by pulling back $D_{\Ml_{(r,d-r)}}$ via a $W(k)$-valued point of $\Ml_{(r,d-r)}$. 
\medskip
Moreover $\Ml_{(r,d-r)}$ is of general type.}
\medskip
Theorem 13 represents as well a second (totally independent) proof of the ordinary specialization result of 3.1.8.1. We refer to $\Ml_{(r,d-r)}$ as a complete Shimura envelope (or library). The word library is a more recent terminology suggested by [LO, \S 4] and [dJO, \S 5] which speak about catalogues of $p$-divisible groups. Warning: the catalogues of loc. cit. are ``complete" in their specialized class of $p$-divisible groups (i.e. from the point of view c)) but not from the point of view of the projectiveness of $\Ml_{(r,d-r)}$. 
\medskip
{\bf 1.12.1. Null strata and the generalized Manin problem.} For some more applications of Theorem 13 and its proof we refer to:
\medskip
-- 4.12.12.5 and 4.12.13 4) (they represent a new --independent-- approach to the deformation theory of $p$-divisible groups, in the complete, local, Noetherian case having a perfect residue field), and to 
\smallskip
-- 4.12.12.6 (it says that each integral canonical model $\Mn^1_{k(v^1)}$ mentioned in Theorem 10 has $k$-valued toric points of all possible, i.e. expected, types; this is the second part of B) below). 
\medskip
4.12.12.6 has itself many applications (see 4.12.12.6.1-3). Here we will just mention three of them. First, when combined with 3.13.7.6.0 it allows us (mainly see 4.12.12.6.2 A) to list all the main properties of the $A_1$ Lie type case of the quasi-ultra stratifications mentioned in 1.10.1; this (together with its $p=2$ version, see 4.14.3 J) represents the generalization of [GO]. Some properties are completely new, like the simple formula (see end of 4.12.12.6.6 3)) computing the number of strata of such a quasi ultra-stratification (in the $A_1$ Lie type case) having $U$-ordinary points. 
\smallskip
Second (see 4.12.12.6.3-4 for proofs), we have: 
\medskip
{\bf Theorem 14.} {\it {\bf A)} The null stratum of the special fibre $\Mn^1_{k(v^1)}$ of an integral canonical model as in Theorem 10 has toric points (and so it is always non-empty).
\smallskip
{\bf B) (the strong --Weyl-- form of the generalized Manin problem)} Any stratum of the absolute stratification of $\Mn^1_{k(v^1)}$ has a toric point and all expected toric points do show up.
\smallskip
{\bf C)} The number of strata of the absolute stratification of $\Mn^1_{k(v^1)}$ is at most equal to the number of elements of the Weyl group of $G^1_{\CC}$.}
\medskip
The null stratum of $\Mn^1_{k(v^1)}$ is the
stratum of the refined canonical Lie (or absolute) stratification, which is ``opposite" to the
open, dense stratum and is the natural generalization of the supersingular stratum of $\Mm_{\FF_p}$. With the notations of 1.3-5, a point $y:{\rm Spec}(k)\to\Mn_{k(v)}$ is a point of the null stratum iff the Shimura Lie $\sg$-crystal attached to it has all slopes $0$. Obviously, the null stratum is closed. In the classical situation of $\Mm_{\FF_p}$ one shows that its supersingular stratum is non-empty by taking products of supersingular elliptic curves. In the general context, for the proof of A) (as well as of 4.12.12.6) we use the following two things:
\medskip
{\bf a)} the proof of Theorem 10, to reduce the situation to the context of a standard Hodge situation (so in b) to e) below we use the notations of 1.3-5);
\smallskip
{\bf b)} [Har, 5.5.3] and standard properties of Shimura subpairs of $(G,X)$ (see [Va2, 3.2.11 and 3.3.3]; see also 4.12.12.6.0 2)); based on them we show the existence of suitable maximal $\ZZ_{(p)}$-tori of $G_{\ZZ_{(p)}}$. 
\medskip
See 4.12.12.6.4.3 for abstract (general) forms of A) (and so implicitly of b)). B) is a consequence of the following three things:
\medskip
{\bf c)} the logical abstract form of A) in the context of Shimura $\sg$-crystals;
\smallskip
{\bf d)} 4.12.12.6;
\smallskip
{\bf e)} the (abstract) inductive property of 3.9.7.2. 
\medskip
The idea of 3.9.7.2 is: if $k=\bar k$ and if $(M,\vph,G)$ is a Shimura $\sg$-crystal with the property that not all slopes of $({\rm Lie}(G)[{1\over p}],\vph)$ are zero, then from many points of view (including Newton polygons), the study of $(M,\vph,G)$ gets reduced to the study of a Shimura $\sg$-crystal $(M,\vph_1,G_1)$ with $G_1$ a reductive subgroup of $G$ different from $G$. This inductive approach (of dropping the relative dimension of $G$) applies to $\Mc\Ml_2$ as well; however, the arguments of 3.9.7.2 do not apply (in general) to $\Mc\Ml_3$, though one could hope to mimic them in many situations. 
\smallskip
To our knowledge, the strong --Weyl-- form of the generalized Manin problem was not known before even for Siegel modular varieties; implicitly, the Fact of 4.12.9 1) was not known previously. A weaker form of B) can be roughly formulated as: all Newton polygons ``showing up" are the ones obtained using specializations of (see [Va2, 2.10]) special points of $\Mn$ with values in fields of fractions of DVR's of mixed characteristic; this was predicted implicitly by the Langlands--Rapoport conjecture (see the definition of admissible homomorphisms in [Mi1, \S 4]; see also [Ch3, rm. after Theorem 1] for a down to earth translation). 
\smallskip
Third, the study of the orbit map of the group action $\Mg\Ma$ of 1.6.1 a) splits up in terms of double cosets of suitable parabolic subgroups of $G^{\rm ad}_{\FF}$ (see Fact of Step 3 of 3.13.7.3 and 3.13.7.6.1 for the abstract context). 4.12.12.6 allows us to list out all cases when there is such a double coset consisting just of 1 orbit and corresponding to a stratum of the quasi-ultra stratification of $\Mn_{k(v)}$ of dimension $0$ (see 3.13.7.6.1-2 for the abstract context and see 4.12.12.6.2 B for the geometric context of special fibres of the integral canonical models of Theorem 10). We have:
\medskip
{\bf Theorem 15.} {\it In the listed cases, there is a unique stratum of the quasi-ultra stratification of $\Mn_{k(v)}$ of dimension $0$. All its points are toric points, called pivotal. Even more, the Shimura adjoint Lie $\sg$-crystals attached to $k$-valued points of it are all inner isomorphic (it is 4.2.10 which allows us to speak about inner; restricting to a connected component of $\Mn_{k(v)}$, we have a form of this in terms of Shimura $\sg$-crystals attached to $k$-valued points of it).}
\medskip
The toric part of Theorem 15 is an application of the reduction steps of 3.13.7.4 B and 3.13.7.6.0 and of a (significantly simplified form of) the inductive proof of 3.6.15 B (see 3.13.7.6.3 for the abstract context). Theorem 15 represents the generalization of the following classical result: if $k=\bar k$, then any $\sg$-crystal $(M,\vph)$ over $k$ whose truncation mod $p$ is isomorphic to the truncation mod $p$ of the $\sg$-crystal $(M,\vph_1)$ of a product of supersingular elliptic curves over $k$, is in fact isomorphic to $(M,\vph_1)$. In the classical situation, one uses the fact that $\vph_1^2$ is $p$ times a $\sg^2$-linear automorphism of $M$. This elementary trick does not work in the general context and this explains why we need all the ``Weyl staff" of 3.13.7.1-6 for proving Theorem 15. 
\smallskip
In the remaining (i.e. not listed) cases (involving most of Shimura varieties of $A_n$ or $D_n^{\HH}$ type) it can be checked using relatively simple computations the existence of strata of the quasi-ultra stratification of $\Mn_{k(v)}$  having dimension 0 and toric points: we refer to \S 9-10 for these computations (see also Corollary 2 of 4.6.7). 
\smallskip
It is worth mentioning that 3.13.7.6.1 and 3.13.6.7.3 apply as well to the class $\Mc\Ml_2$, cf. 3.13.7.7.
\medskip
{\bf 1.12.2. Degrees of definition.} One of the things which comes again and again is: degrees of definition. To detail what we mean by this we refer to 1.6.2. 3.6.15 B implies: if $k=\FF$, then $(M,\vph,G)$ is obtained from a Shimura $\sg$-crystal over a finite field by extension of scalars. So roughly speaking the degree of definition of $(M,\vph,G)$ is the smallest possible degree over $\FF_p$ of such a finite field. As this definition is too technical for practical computations, in most of the cases we allow $k$ to be arbitrary and moreover:
\medskip
-- we either restrict (see 3.11.3) to specific Shimura $\sg$-crystals for which we can define (using some logical lifts) immediately its degree of definition (later on we will check that the two definitions coincide; for instance, this is done implicitly in 3.11.8 for the Shimura-ordinary context), or
\smallskip
-- we work with potentially cyclic diagonalizable Shimura filtered $\sg$-crystals and we define (see 2.2.16.4 and 2.2.22 1)) their degrees of definition using their canonical split cocharacters (of [Wi]). 
\medskip
Similarly, in adjoint Shimura contexts we define $A$-degrees of definition. For samples, see 2.2.16.4-5, 3.11.3, etc.
\medskip\smallskip
{\bf 1.13. The geometric picture of the case $p=2$.} For the case of a $p=2$ standard Hodge situation introduced in 2.3.18, see 4.12.12.1 and 4.14.3. The needed well known facts from crystalline cohomology pertaining to the mixed characteristic $(0,2)$ are reviewed (parts of it with new personal proofs) in 2.3.18.1. For instance, based on them and on 3.14, we get that all of 1.4-6 holds for a $p=2$ standard Hodge situation $(f,L_{(2)},v)$ modulo two modifications (cf. 4.14.3). 
\smallskip
First, we do not prove that the $U$-canonical lifts of 1.5.1 (and so implicitly the $G$-canonical lifts of 1.4 c)) can be always uniquely identified, i.e. determined (we recall that to locally lift a $2$-divisible group over an algebraically closed field of characteristic $2$ it is not enough to lift the filtration in the crystalline cohomology context; see 2.3.18.1). So presently, it is convenient to allow a finite number of $U$-canonical lifts of a $k$-valued $U$-ordinary point $y$ of $\Mn_{k(v)}$: a $W(k)$-valued point of  $\Mn$ lifting $y$ is a $U$-canonical lift iff its attached Shimura filtered $\sg$-crystal is a $U$-canonical lift. Accordingly, 3.1.1.2 in geometric situations has to be slightly modified for $p=2$ (for instance, as in 4.14.3 B4). On the other hand, 4.14.3 B3 shows that we can define uniquely a $U$-canonical lift of $y$, provided the Shimura adjoint Lie $\sg$-crystal attached to $y$ does not have slope $-1$. 
\smallskip
Second, the part of 1.12.1 referring to 4.12.12.6 has to be weaken as suggested by 4.14.3 J (however, 4.14.3.2.5 points out variants where such a weakening can be considerably eliminated). 
\medskip
{\bf 1.13.1. The $D$ case of PEL type embeddings for $p=2$.} Warning: as in [Mi1, 1.1], in the definition of Shimura varieties of PEL type --and so of PEL type embeddings-- we allow $\QQ$--semisimple algebras which are not necessarily simple. For the case $p=2$ of a PEL type embedding involving the $D$ case (see [Ko2, ch. 5] for the meaning of this), left aside in loc. cit. and in [LR], see 2.3.18 B: it represents the first new instances of $p=2$ standard Hodge situations (while reobtaining all previous ones). 2.3.18 B finally completes the proof of the existence and uniqueness of integral canonical models in the integral context of PEL type embeddings. This proof was started in the original edition (dated 1965) of [MFK] and continued --as a direct application of loc. cit. and of Serre--Tate and Grothendieck--Messing deformation theories (see [Me])-- in [LR] and [Ko2]. We formulate 2.3.18 B here, in a slightly different (and weaker) form, more convenient to make the passage from [Ko2] (without using too much of the language of this paper). Let $f$, $L_{(2)}$ and $v$ (resp. $\FF$) have the same significance as in 1.3 (resp. 1.6.1) but for $p=2$. We assume that $\psi$ induces a perfect form $L_{(2)}\otimes_{\ZZ_{(2)}} L_{(2)}\to\ZZ_{(2)}$, that the Zariski closure of $G$ in $GL(L_{(2)})$ is a reductive group $G_{\ZZ_{(2)}}$ over $\ZZ_{(2)}$ and that all simple factors of $G^{\rm ad}$ are of some $D_n$ Lie type. Let $\Mb$ be a $\ZZ_{(2)}$-subalgebra of ${\rm End}(L_{(2)})$ such that:
\medskip
-- it is self dual with respect to $\psi$;
\smallskip
-- over $W(\FF)$ is a product of matrix $W(\FF)$-algebras;
\smallskip
-- $G$ is the connected component of the origin of the subgroup of $GSp(W,\psi)$ fixing its elements.
\medskip
Let $H$, $K_2$, $\Mm$ and $O_{(v)}$ be defined as in 1.3. We have:
\medskip
{\bf Theorem 16.} {\it The normalization $\Mn$ of the Zariski closure of ${\rm Sh}(G,X)/H$ in $\Mm_{O_{(v)}}$ is the integral canonical model of the Shimura quadruple $(G,X,H,v)$.} 
\medskip
The new ideas needed to handle this left aside case (i.e. to prove Theorem 16), besides the ones of [Ko2, ch. 5] and [Va2, \S 5], consist in:
\medskip
{\bf a)} a slight refinement of parts of [Ko2] (for instance, we have in mind 7.1-2 of loc. cit.) in the $p=2$ context of the D case (see 2.3.18 B2-3);
\smallskip
{\bf b)} a slight refinement of [Fa2, th. 10] (see iii) to v) of 2.3.18 B9);
\smallskip
{\bf c)} performing a greater part of [Va2, 5.2-3] in the context of the $W(\FF)$-algebra $Re$ instead of the $W(\FF)$-algebra $\tilde Re$ of loc. cit.
\medskip
See below for why $\Mn$ is uniquely determined as an integral canonical model.
\medskip
{\bf 1.13.2. Uniqueness of integral canonical models in mixed characteristic $(0,2)$.} 2.2.1.4 fulfills [Va2, 3.2.1.4 5)]. The new thing it brings: any regular, formally smooth scheme over a DVR of mixed characteristic $(0,2)$ and index of ramification $1$ is healthy in the sense of [Va2, 3.2.1 2)]. 2.2.1.4 implies (see 2.2.1.5) the uniqueness of (local) integral canonical models of Shimura varieties in mixed characteristic $(0,2)$. This last result was not known previously, even for the Siegel modular variety of dimension $3$ or $6$ (i.e. of genus $2$ or $3$). What 2.2.1.4 brings new to Step B of [Va2, 3.2.17]: we alternate using an ``analysis" mod $p$ (with $p\ge 2$) with the approach of loc. cit. pertaining to mixed characteristic; so we can appeal to the language of 2.2.1.0 and so to the results of [Fa1, p. 31-3]. 
\smallskip
See 2.2.1.5.2 for a second approach to the uniqueness of projective (local) integral canonical models of Shimura varieties of compact type in arbitrary mixed characteristic.
\medskip
{\bf 1.13.3. The passage from the Hodge type to the abelian type in mixed characteristic $(0,2)$.} 4.14.3.2 achieves the passage of the existence of integral canonical models in mixed characteristic $(0,2)$ from the Hodge type to the abelian type. [Va2, 6.2.2 a)] has a very restricted application in mixed characteristic $(0,2)$. So in order to achieve the mentioned passage, i.e. to get a $p=2$ version of [Va2, 6.2.2 b)], as we have 2.2.1.5, we just need to get $p=2$ variants of the corrected versions (see AE.4 and AE.4.1-2) of [Va2, 6.2.2.1]. They are developed along different parts of this paper (for instance, see 4.14.3.2.1-2) and involve mild conditions. For instance, the main assumption we need in order to achieve the mentioned passage refers to the fact that we start from a $p=2$ standard Hodge situation $(f,L_{(2)},v)$ whose associated integral canonical model $\Mn$ is such that each connected component of $\Mn_{k(v)}$ has $k$-valued quasi-pivotal points (i.e. toric points whose attached Shimura adjoint Lie $\sg$-crystals do not have slope $-1$). The $p=2$ version of 4.12.12.6 implies that this holds, provided the connected components of $\Mn_{k(v)}$ are permuted transitively by $G(\AA_f^2)$ (for instance, if the Shimura pair $(G,X)$ involving $f$ is of compact type, cf. 2.3.3.1 and its $p=2$ version). This assumption is related to the fact (it is just a natural extension of the part of 2.3.18.1 C referring to no slope $0$ or no slope $1$) that for $k$-valued quasi-pivotal points, the lifts of it to $W(k)$-valued points of $\Mn$ are in one-to-one correspondence to (logical) lifts of filtrations in the crystalline cohomology context, cf. 4.14.3 B3. 
\smallskip
Warning: we still do not know how to show that [Va2, 6.8.2] applies as in [Va2, 6.8.3] for $p=2$, i.e. we still do not know how to achieve in mixed characteristic $(0,2)$ the passage from the abelian type to the preabelian type.
\smallskip
See 2.4.1-2 for abstract variants of 2.3.17 for $p=2$; the approach to them is based on [Fa2, th. 10] and the remarks after.  
\medskip
{\bf 1.13.3.1. Generalized Serre lemma.} 4.14.3.2.1-2 lead to a form of the generalized Serre lemma of [Va2, 6.4.6 6)] (see also AE.4.2) in mixed characteristic $(0,2)$, see 4.14.4.1. Using this we regain the classical Serre lemma of [Mu, p. 207]: see 4.14.4.2. To our knowledge, this is the first proof of loc. cit. in a $p$-adic (and not $l$-adic) context.
\medskip
{\bf 1.13.4. Some samples of results with $p=2$.} It is worth mentioning, as samples, some of the things we get in connection to Theorem 16. For instance, we have (all these can be read out from different parts of 4.14.3 --see 4.14.3.3--, cf. also the above part of 1.13):
\medskip
{\bf a)} The special fibre $\Mn_{k(v)}$ of $\Mn$ has a Zariski dense set of Shimura-ordinary points. For all such points with values in a perfect field $k$ of characteristic $2$ we identify a finite set of Shimura-canonical lifts via: their attached Shimura adjoint filtered Lie $\sg$-crystals are of parabolic type. The $W(\FF)$-valued Shimura-canonical lifts give (by pulling back $\Ma$) to abelian varieties over $W(\FF)$ having complex multiplication.
\smallskip
{\bf b)} If $(G,X)$ is of compact type, then any Shimura quadruple of abelian type having the same adjoint as $(G,X,H,v)$, has a uniquely determined integral canonical model. 
\smallskip
{\bf c)} The quasi-ultra stratification of $\Mn_{k(v)}$ has regular, quasi-affine strata whose codimensions can be easily computed. Its number of strata is at least as expected (i.e. predicted in terms of a  quotient set of the Weyl group of $G_{\CC}$).
\smallskip
{\bf d)} If $G^{\rm ad}_{\ZZ_2}$ is absolutely simple of $D_n$ Lie type and if either $n$ is even and $G^{\rm ad}_{\ZZ_2}$ is split or $n$ is odd and $G^{\rm ad}_{\ZZ_2}$ is non-split, then $\Mn_{k(v)}$ has pivotal points to which the $p=2$ version of Theorem 15 applies. 
\smallskip
{\bf e)} If $k=\bar k$ and the residue field of the prime of $E(G^{\rm ad},X^{\rm ad})$ divided by $v$ is $\FF_2$, then the completion of $\Mn$ in any $k$-valued Shimura-ordinary point $y$ has a canonical structure of a formal torus, with the Shimura-canonical lifts of $y$ corresponding to its $2$-torsion points.
\medskip\smallskip
{\bf 1.14. Afterthoughts, the Second Main Corollary and advices.} Here we include the Second Main Corollary and different specifications meant to ease the reading.
\medskip
{\bf 1.14.1. List of invariants.} 
4.14.6 contains the complete list of invariants of a given $(p=2)$ standard Hodge situation, introduced along different parts of \S2-4; the main class invariants (in the sense of 2.2.22 2)) are listed in 4.14.6.1, while 4.14.6.2 contains references to the types of points (with values in fields) defined in different parts of \S 4. 
\medskip
{\bf 1.14.2. More on the overall organization.}  Following \S4, there is an addenda and errata to [Va2]. We refer to it as AE. The reading of AE.4 is indispensable to different parts of 2.3.3-11 and 4.4-14. To AE.0-3 and AE.5-6, the reader ought to look only if needed. Some afterthoughts are presented in 4.14 and in Appendix. Appendix makes the connection and comparison with [RR] and so in particular contains some general forms of the results of \S 2-4 (see also the last paragraph of 1.2.9).
\smallskip
Two proofs of Theorem 1 are presented. The first one, less than 24 pages, starts in 3.1 and ends in 3.4.14, via 3.5; it completely avoids the use of 3.6. 
The second one is just a more elaborated form of the first one, in the sense that it leads to a much stronger form of d) of Theorem 1 (cf. Theorem 2 of 3.15.1). It starts in 3.1 and ends in 3.7.4; it uses from 3.6 just the part till 3.6.8 inclusive. Also for 4.1-8 and 4.11 which form a sequence continuing 2.3, we need from 3.6 just 3.6.6: 3.6.6 is independent of the other parts of 3.6 except of the first sentence of 3.6.0. 
\smallskip
The parts 3.6.13-14 are needed just in 3.6.19, in 4.3.6 (as an exemplification) and in 4.12-14 (for $p=2$ we need to add 3.14 E and 4.14.3 K). They are independent of what follows after them in 3.9 and 3.11-12; so we do not hesitate to use for their proofs results stated in 3.9 and 3.11-12 and so, their reading should proceed only after the reading of 3.9 and 3.11-12. We felt that it is appropriate to include them as part of 3.6, for having the deformation theory part (i.e. one of the most beautiful parts) of \S1-4 presented fluently. 3.6.15-16 (resp. 3.6.17) are first used in 3.15.7-10 (resp. in 3.11.1). 
\smallskip
The sequence 3.6.17-20 is very little related to what follows after it; so it can be read at any (needed) time after 3.6.1-14. To make \S 4 as little dependent as possible on \S 3 and to avoid repetitions, in few places of \S 2-3 we refer to \S 4: as a typical situation is the proof of 2.2.16; these references are made either as possible variants or in such a way that the reader will have no difficulty at all in looking in advance at one or two paragraphs of \S 4 for some independent arguments.
\smallskip
For the present work we were motivated by very concrete problems, emerging from the theory of Shimura varieties. Whenever we were aware from the very beginning that things work out in a more general context, we worked in such a general context; when not or when pedagogical reasons required, we stated (or presented) the things in the more general context as remarks following the ``simpler" context. We tried to present the ideas in the most natural order. But we did not hesitate (in order to bring deeper understanding) to point out different connections between them, which somehow were going against this order. 
\medskip
{\bf 1.14.3. Starred results.} There are two facts stated (and referred) here whose proofs are postponed to future papers: 4.2.8.1 and d) of 4.4.1 3); they are labeled with a $*$. Whenever we use them we present a second proof entirely independent of them. The only two exceptions are:
\medskip
-- the starred Claim of the proof of 4.11.2, as its use is entirely avoidable, and
\smallskip
-- the starred answer of 4.12.5, as it is easy to restate it in a way which does not rely on d) of 4.4.1 3) (see the paragraph following it); also we would like to point out that for most purposes (like ones involving Newton polygons), the $\ZZ_p$-structures of 2.2.9 8) are ``equally good" (though not as precise) as d) of 4.4.1 3). 
\medskip
Occasionally we mention other results to be proved in \S5-14: they are never used here and so they are not starred. So \S1-4 are self contained modulo the existence of integral canonical models of Shimura quadruples $(G,X,H,v)$ of preabelian type with $v$ dividing $p=3$; warning: the limitations of [Va2, 6.8.6] are removed in 4.6.4. 
\smallskip
In other words, in 4.6.4 we prove that the criterion [Va2, 6.8.2 b)] holds (with $p\ge 3$); the proof relies on the simple fact that, in the context of a standard PEL situation, the abelian varieties obtained by natural pull backs via $G$-canonical lifts of $G$-ordinary points with values in $\FF$, are having complex multiplication (cf. 4.4.5; see also 1.15.8). As a conclusion we get (cf. 4.6.4 B):
\medskip
{\bf Second Main Corollary.} {\it Everywhere in [Va2] the single stars (i.e. all stars except the two ones of [Va2, 5.6.5 h)]) can be removed.} 
\medskip
 So the reader, if desires, can everywhere we are dealing with:
\medskip
-- a standard Hodge situation $(f,L_{(p)},v)$ such that $v$ divides $3$ and not all simple factors of $G^{\rm ad}$ are of some $A_n$ Lie type (cf. 2.3.5.1), can add the assumption that $\Mn$ a priori exists;
\smallskip
-- a Shimura quadruple $(G,X,H,v)$ such that $v$ divides $3$ and not all simple factors of $G^{\rm ad}$ are of some $A_n$ Lie type, 
can add the assumption that there is a standard Hodge situation $(f^1,L^1_{(p)},v^1)$ such that $(G^{1{\rm ad}},X^{1{\rm ad}},H^{1{\rm ad}},v^{1{\rm ad}})=(G^{\rm ad},X^{\rm ad},H^{\rm ad},v^{\rm ad})$ and $G^{1{\rm der}}$ is maximal among isogeny covers of $G^{\rm ad}$ allowed by the abelian type (i.e. can add that the analogue of [Va2, 6.4.2] holds; cf. also the Existence Property of 1.10).
\medskip
{\bf 1.14.4. Promises.} 2.2.1.4, the Second Main Corollary and 4.9.8.2 are among the first places meant to fulfill (in due time) all promises and expectations of [Va2]. 
\medskip
{\bf 1.14.5. Some extra tips.} Due to the fact that the results of \S1-4 have been largely detailed above, most common we do no stop at the beginning of the other chapters \S 2-4 (or of their sections), to outline their contents. However, wherever the title of a section does not make its main goal entirely obvious, this main goal is stated explicitly in the opening paragraph. 
\smallskip
At the first lecture, the reader can skip the case $p=2$, 3.6.8.3 and Step 2 of the second way (solution) of 4.6.4 A. Moreover, the reader will benefit by having a general look on 3.6.8 from the very beginning; to avoid repetitions, in some more general contexts we still refer to it (even if it uses very concrete situations) (cf. also 3.6.8.4 1) to 3)). \S1-4 is organized as a foundation: inevitably, in some places (very few in number), it contains material which seems isolated from the rest of it but which is used quite a lot in \S5-14. For instance, many of the concepts of 2.2.4, 2.2.22 2), 3.13.1-6 and 3.13.8-9 are just very partially exploited (or referred) here. As a very wide variety of terminology is introduced and used, often, for convenience, we point backwards to (i.e. we remind) some definitions.
\smallskip
For us, techniques used for proving some results are equally important as these results. So thus, many results in this paper are getting two or more proofs. Moreover, we do not hesitate to include results which are particular cases of our future work, as long as their proofs are involving interesting techniques or are much simpler; as a sample: 2.3.18 B is a very particular case of the $p=2$ theory of Shimura varieties of abelian type to be fully developed in \S6, [Va5] and [Va3]. We strongly encourage the reader to peruse all techniques introduced in this paper. Motivation: we consistently strove to capture the very essence of the mathematical problems we dealt with here.
\smallskip
Warning: many (relatively well known) results pertaining to reductive groups are mentioned with proofs or appropriate references at different parts of \S1-4; we felt this offers pedagogical advantages (in comparison with a gathering of them --in tens of pages-- in one section) even if few and small repetitions are made.
\smallskip
Warning: based on the Second Main Corollary, in 4.9-14 we refer to results of [Va2] having a single star, without mentioning each time that all such stars have been removed by 4.6.4 B. 
\smallskip
Besides the results mentioned in 1.2-15, this paper contains many other ideas, exercises and examples. The ideas we did not have time to deal with, are presented in the form of problems and questions or as remarks. 63 (4.5.6 itself counts as 10) instructive and simple exercises, almost all with full hints or solutions, are included in \S 2-4 and AE.4.2. 
\medskip\smallskip
{\bf 1.15. On the contents of \S5-14.} We present now in a very brief manner the contents of the other two parts continuing this paper. We recall that without stating the things otherwise we assume $p\ge 3$.
\medskip
{\bf 1.15.1. On \S5.} In \S 5 we elaborate on the notion of Shimura $p$-divisible group introduced in 2.2.20 and prove the conjecture of Milne stated in [Va2, 5.6.6], in the wider context of Shimura $p$-divisible groups; so we also prove d) of 4.4.1 3) and remove the two stars of [Va1, 5.6.5 h)]. 
\smallskip
As a corollary of Milne's conjecture we get the integrality theorem of Shimura $p$-divisible groups and so also the integrality theorem of Shimura varieties of Hodge type. Here we formulate it only for a $W(k)$-valued point of an integral canonical model of a Shimura variety of Hodge type (using 2.3.13.1 it can be stated for many other valued points; see also 4.4.12). 
\smallskip
If $(f,L_{(p)},v)$ is a standard Hodge situation then, with the notation of 1.3, the $p$-component of the \'etale component of $w_\al$ is integral (i.e. it is an element of the tensor algebra of
$H^1_{\acute et}\bigl(A_{\overline{K_0}},\ZZ_p\bigr)\oplus (H^1_{\acute et}\bigl(A_{\overline{K_0}},\ZZ_p\bigr))^*$, where $\overline{K_0}$ is the algebraic closure of
$B(k)$) iff $t_\al$ is integral (i.e. it is an element of the tensor algebra of $M\oplus M^*$). What is worth remarking:
there is no restriction on the degree of $w_\al$, $\al\in\Mj$, and it also holds  for $p=3$. 
\smallskip
We use Milne's conjecture to show how, to a standard Hodge situation $(f,L_{(p)},v)$ and to a representation $\rho:G_{\ZZ_p}\to GL(N)$ (with $N$ a finitely generated $\ZZ_p$-module, not necessarily free) it is naturally associated and object (if $N$ as a set is finite) or a non-trivial almost $p$-divisible object (otherwise) ${\got C}_\rho$ of $\Mm\Mf^{\nabla(big-p+tens)}(\Mn)$ (with $\Mm\Mf^{\nabla(big-p+tens)}(\Mn)$ as in 2.2.4.1 5); so the necessary precautions are taken in cases we do not get something to which the gluing arguments of [Fa1, proof of 2.3] apply). For almost $p$-divisible objects we refer to 2.2.1.7 1). If $N$ is a free $\ZZ_p$-module of finite rank, then ${\got C}_{\rho}$ is a filtered $F$-crystal in locally free sheaves on $\Mn_{k(v)}$ (cf. the conjecturally expected result mentioned in [Ch2, end of $(Q1)_D$]). The association $\rho\to{\got C}_\rho$ is functorial and respects the usual operations (tensor products, identity elements, duals when defined, etc.). It is uniquely determined by requesting two other natural properties to be satisfied.
\smallskip
In the last part of \S 5 we introduce the notion of automorphic vector bundle on an integral canonical model of a Shimura variety. We use the existence of ${\got C}_\rho$'s to prove the existence and uniqueness of automorphic vector bundles on $\Mn$ (associated to representations of $G_{\ZZ_{(p)}}$ on $\ZZ_{(p)}$-free modules of finite rank or to representations of $G_{O_{(v)}}$ on $O_{(v)}$-free modules of finite rank) or on $\Mn\times_{O_{(v)}} O_{(v_1)}$ (with $O_{(v_1)}$ a local ring which is an \'etale $O_{(v)}$-algebra) (associated to representations of $G_{O_{(v_1)}}$ on $O_{(v_1)}$-free modules of finite rank). We also extend these results to 
Shimura varieties of preabelian type, under slight conditions (as in [Mi4, ch. III]) on the representations involved.
\medskip
{\bf 1.15.2. On  \S6.} In \S 6 we prove that for any injective map $f:(G,X)\hookrightarrow
\bigl({\rm GSp}(W,\psi),S\bigr)$, every $\ZZ_{(p)}$-lattice $L_{(p)}$ (of $W$) good with respect to $f$, is automatically crystalline well positioned for the map $f$ (with respect to any prime of $E(G,X)$ dividing $p$). We also prove the existence of integral canonical models of Shimura varieties of preabelian type with respect to primes dividing 3 and the existence of integral canonical models of Shimura varieties of Hodge type with respect to primes dividing $p\ge 3$ and (plenty of) non-hyperspecial subgroups. Related to proofs we just mention: Milne's conjecture of 1.15.1 is the main new ingredient besides [Va2] needed to overcome the limitations of [Va2, 5.1] for $p=3$. 
\medskip
{\bf 1.15.3. On  \S7.} In \S 7 we introduce a method by which, from particular Shimura $\sg$-crystals, we construct integral, normal, faithfully flat schemes over local rings which are DVR's formally \'etale over $\ZZ_{(p)}$, with $p\ge 2$, which either have fibres having a dense set of formally smooth points or have the extension property (of [Va2, 3.2.3 3)]). By this method we recover all connected components of integral canonical models of Shimura varieties of Hodge type. We (conjecturally) expect to get some other interesting integral models, even smooth ones; this expectation is a heuristic extrapolation of the existence of generalized Shimura $p$-divisible objects $(M,\vph,G)$ over $k$
for which the pair $({\rm Lie}(G^{\rm ad}),\vph)$ has factors of one of the special $E_6$, $E_7$ or $D_n^{\rm mixed}$ types
(see 3.10.5-6 for the definition of these types). We call these normal schemes twisted connected Shimura schemes.
These schemes and their variants (obtained using generalized Shimura $p$-divisible objects) are supposed to lead to the construction of integral canonical models of Shimura varieties of special type.
Also, they either lead to very beautiful mathematics or provide a negative answer (which we do not expect) to Oort's conjecture [Oo1, p. 6].
\smallskip
\S 7 contains as well many other complements on Shimura (Lie) $\sg$-crystals (like the slice principle, the rigidity property, etc.). See also 2.2.6 4).
\medskip\smallskip
{\bf 1.15.4. On  \S8.} In \S 8 (as a continuation of [Va10]), we apply some of the ideas of \S 3, \S 5  and \S 7 to prove the Mumford--Tate conjecture for many abelian varieties over number fields.
We also present a simple strategy meant to prove the following general criterion (it is known that it implies this conjecture):
\medskip
{\bf Expectation (the split criterion).} {\it Let $A$ be an abelian variety over a number field $E$. There is $N(A)\in\NN$ (very easily expressible in terms of $A$) such that,
if $\ell$ is a prime greater or equal to $N(A)$ and the reductive group $G_{\QQ_\ell}$ over $\QQ_\ell$, defined as the connected component of the origin of the algebraic envelope of the Galois representation $\rho:{\rm Gal}(\bar E/E)\to GL(H_{\acute et}^1(A_{\bar E},\QQ_\ell))(\QQ_\ell)$, is a split reductive group, then the Mumford--Tate conjecture is true for $A$ with respect to $\ell$, i.e. we have $G_{\QQ_\ell}=H_{A\QQ_\ell}$. Here $H_A$ is the Mumford--Tate group of $A$.}
\medskip
{\bf 1.15.5. On  \S9.} In \S 9 we deal in great detail with the null stratum of the special fibre of an integral canonical model of a Shimura variety of preabelian type: influenced by properties (see [LO, 0.2]) of supersingular strata of special fibres of integral canonical models of Siegel modular varieties, we formulate our expectations about the properties enjoyed by these null strata and prove some of these expectations. In particular, we get a non-supersingular criterion similar to the non-ordinary one of 1.7. Warning: referring to 1.6, the null stratum of $\Mn_{k(v)}$ is not always the pull back of the supersingular stratum of $\Mm_{k(v)}$ via the morphism $\Mn\to\Mm$ (see 1.3), as the toric part $(G^{\rm ab},X^{\rm ab})$ of $(G,X)$ plays a role in the definition of this morphism.
\medskip
{\bf 1.15.6. On  \S10.} The first goal of \S 10 is to show that the number $N(G,X,v)$ of strata of the refined canonical Lie stratification of the special fibre of an integral canonical model (of a Shimura variety ${\rm Sh}(G,X)$ of preabelian type with respect to a prime $v$ of $E(G,X)$ dividing $p>2$) and the dimensions of the geometrically connected components of these strata are easily and effectively (algorithmically) computable. To achieve this first goal, the first philosophy is: 
\medskip
{\bf Ph.} {\it The number $N(G,X,v)$, the mentioned
dimensions, as well as the local geometry of strata (like the type of singularities in the \'etale topology of a given 
stratum, like which stratum specializes to which stratum, etc.) depend only on the pair
$(G^{\rm adnc}_{\QQ_p},\mu_g^{\rm ad})$, where $G_{\QQ_p}^{\rm adnc}$ is the smallest direct product factor of $G_{\QQ_p}^{\rm ad}$ with the property that the generic fibre $\mu_g^{\rm ad}$ of the cocharacter $\mu^{\rm ad}:\GG_m\to G^{\rm ad}_{W(k(v^{\rm ad}))}$ (with
$k(v^{\rm ad})$ as the residue field of the prime $v^{\rm ad}$ of $E(G^{\rm ad},X^{\rm ad})$ divided by $v$) factors through the extension of $G_{\QQ_p}^{\rm adnc}$ to $B(k(v^{\rm ad}))$. Here $\mu^{\rm ad}$ is a cocharacter such that $\mu_g^{\rm ad}$ becomes (under an embedding of $B(k(v^{\rm ad}))$ in $\CC$ extending the inclusion
$E(G^{\rm ad},X^{\rm ad})\subset\CC$) $G^{\rm ad}(\CC)$-conjugate to the cocharacters $\mu_x$, $x\in X^{\rm ad}$, of [Va2, 2.2].}
\medskip
This forms the global form of the invariance principle.  
The second philosophy, meant to express very concretely what we mean by ``easily and effectively", is ``captured" by 3.6.14.4, by the formulas of 4.5.15.2 and 4.5.15.2.1 computing the dimensions of the strata of the quasi-ultra stratifications, by the CM level one (expected to hold) property (of 1.6.2) and by the completion property (or by its weaker --intermediary-- form 4.12.12.6). 
\smallskip
We have versions of this Ph for the other stratifications of 1.10.1 as well as for the (presently abstract) context of $\Mc\Ml_2$.
\smallskip
The second goal of \S 10 is to continue the study of the ultimate types of stratifications and to get progress with respect to the Real Problem of 1.6.5. Implicitly we include an extensive study of Grothendieck's (adjoint) specialization categories and of the invariants introduced in 4.12.12.8. In particular, starting from the methods of \S1-4, we generalize [dJO, \S 5] from four points of view:
\medskip
{\bf 1)} we include global forms of different isogeny properties (like of [dJO, 2.17]);
\smallskip
{\bf 2)} we work in a reductive context with cycles (warning: [dJO, \S 5] can not be adapted to such a context);
\smallskip
{\bf 3)} we do not restrict to the isocline context (i.e. to the context of Shimura $p$-divisible groups having only one slope);
\smallskip
{\bf 4)} we go beyond [dJO, 5.13] by being, in many situations, very specific on which $p$-divisible groups specializes to which $p$-divisible group (i.e. we introduce and solve in many situations Grothendieck's integral problem).
\medskip
These generalizations are possible due to the fact that 3.6.15 B has logical global variants; for the constant Newton polygon context they are essentially obvious (cf. also 3.15.10.1 3)) and this takes care of 1) and 2) for this context. The main tools for 3) and 4) are: Theorems 3 and 6 and 3.15.7-9 (specifically, see 3.15.7 BP2). See [Va6] for extra generalizations of [dJO, \S 5] in the context of $\Mc\Ml_2$.
\medskip
{\bf 1.15.7. On  \S11.} The main purpose of \S 11 is to prove the isogeny property mentioned in [Va2, 1.7]. The idea is: we follow entirely the pattern of [FC, 4.3, p. 264] but in the crystalline context provided by Fontaine categories; this approach is supported by the fact that the isogeny property is trivial for $G$-ordinary points. The hard part is when we deal with points of the null strata; for such points we need to combine loc. cit. with ideas of [Mi1].
\medskip
{\bf 1.15.8. On  \S12.} In \S 12 (and \S 6) we prove the following result. Let $(G,X,H,v)$ be a Shimura quadruple, with $(G,X)$ a Shimura pair of adjoint, abelian type and with $(v,2)=1$. Then there is an injective map
$f:(G_0,X_0,H_0,v_0)\hookrightarrow  (G_1,X_1,H_1,v_1)$ such that:
\medskip
{\bf 1)} $(G_0^{\rm ad},X_0^{\rm ad})=(G,X)$ and (under the natural inclusion $E(G,X)\subset E(G_0,X_0)$) $v_0$ divides $v$;
\smallskip 
{\bf 2)} $(G_1,X_1,H_1,v_1)$ is the quadruple of a Shimura variety of
PEL type, with all simple factors of $G^{\rm ad}$ of some $A_\ell$ Lie type ($\ell\in\NN$ depends on the factor);
\smallskip
{\bf 3)} provided $(G,X)$ has no simple factor of $D_n^{\HH}$ type ($n\in\NN$, $n>3$), the residue field $k(v_0)$ of $v_0$ is the same as the residue field $k(v_1)$ of $v_1$;
\smallskip
{\bf 4)} there is an injective map $f_1:(G_1,X_1)\hookrightarrow ({\rm GSp}(W,\psi),S)$ and a $\ZZ_{(p)}$-lattice $L_{(p)}$ of $W$ (with $p$ the rational prime divided by $v$), with the property that the triples $(f_1,L_{(p)},v_1)$
and $(f_1\circ f,L_{(p)},v_0)$ are standard Hodge situations;
\smallskip
{\bf 5)} if $(G,X)$ has no simple factor of $D_n^{\HH}$ type ($n\in\NN$, $n>3$) (so we have $k(v_0)=k(v_1)$) or if some other (simple to check) conditions hold, then we are in a context in which Theorem 11 applies.
\medskip
 As a consequence of this and its proof we prove the result concerning special points mentioned in [Va2, 1.6.1] as well as the following result.
\medskip
{\bf 1.15.8.1. A result.} Let $A$ be an abelian scheme over the ring $W(k)$ of Witt vectors of a finite field $k$ of characteristic $p\ge 3$. Let $(s_{\be})_{\be\in I}$ be a family of Hodge cycles of $A$, including a polarization of degree relatively prime to $p$. We assume the Zariski closure of the subgroup $G_p$ of $GL(L_p[{1\over p}])$ fixing the $p$-components of the \'etale components of the Hodge cycles of the family $(s_{\be})_{\be\in I}$ in $GL(L_p)$, with $L_p:=H^1_{\acute et}(A_{\overline{B(k)}},\ZZ_p)$, is a reductive group over $\ZZ_p$. Let $\pi$ be the Frobenius endomorphism of $A_k$. Then a positive, integral power of $\pi$, as an automorphism of $H^1_{\rm crys}(A_k/W(k))[{1\over p}]$, fixes the de Rham component of $s_{\be}$, $\forall\be\in I$.
\medskip
{\bf 1.15.9. On  \S13-14.} In \S 13-14, we prove the results (pertaining to Langlands--Rapoport conjecture) stated in [Va2, 1.7].
\medskip
{\bf 1.16. Some history.} A small part of the results presented in \S1-4 (essentially Theorems 1 and 2 in the context of a standard Hodge situation and their applications, presented in a) and b) of 1.4), is a revised and improved version of the second part of our thesis [Va1] (see also [Va2, 1.6]). Strictly speaking b) of Theorem 1 was not stated in [Va1]; however all ingredients needed for getting it can be found in loc. cit. In particular, the starting idea for getting a) and c) of Theorem 1, of using Faltings--Shimura--Hasse--Witt shifts (see 3.4.5), which is the essence of all types of computations (for instance, like the ones involving some types of Faltings--Shimura--Hasse--Witt maps), was first introduced in [Va1, 5.5.4] based on a suggestion of G. Faltings (and this explains our terminology). To our knowledge, Shimura-ordinary types were first defined in [Va1, 5.1.3] and the use of Weyl elements in connection to the truncation mod $p$ of Shimura (Lie) $\sg$-crystals was initiated in [Va1, 5.6.2-3]. Implicitly, cyclic factors over (truncation mod $p$) of Lie $\sg$-crystals (see 3.10 for elaborated definitions) were first used in [Va1, \S 5]; see also [Va2, p. 495] and the part of AE.4 referring to it. 
\smallskip
[Va13] was the initial draft (version) of this paper. Proposition 2.3.15, which is very useful for different local computations, was already part of [Va13] (being just the logical globalization of part of [Va2, 5.4]).
\smallskip
We came across quasi Artin--Schreier systems of equations (resp. integral Manin problems) on Oct. 1997 (resp. Feb. 1996): they were part of the version of this paper released on Dec. 1997. 
\medskip\smallskip
{\bf 1.16.1. Supplementary references.} Many results of \S1-4 are completely new; so, there is relatively very few literature to be mentioned in connection to most of them. The literature which influenced us and is directly related to the present work, is mentioned (directly or indirectly) at different moments along the manuscript. The most we can add to what was mentioned in 1.0-15 is as follows:
\medskip
{\bf a)} we do not know if weaker forms of Theorem 3 for ordinary elliptic curves are traceable in the literature;
\smallskip
{\bf b)} Theorems 2 and 5 and the First and Second Main Corollaries are untraceable in the literature;
\smallskip
{\bf c)} tiny bits of 1.9.1 and of 1.4 were known before (for instance, see references in 3.1.8.1 and 4.6.1 3));
\smallskip
{\bf d)} in connection to Theorem 7 we refer to 4.6.2.5;
\smallskip
{\bf e)} if $k(v)=\FF_p$, then Theorem 9 is a consequence of Serre--Tate and Dwork theories of coordinates mentioned in 1.0;
\smallskip
{\bf f)} except some very particular cases (basically pertaining to curves), there is no literature in connection to Theorem 10;
\smallskip
{\bf g)} outside of the case $k(v_1)=k(v_2)=\FF_p$ (which is trivial), there is no literature in connection to Theorem 11;
\smallskip
{\bf h)} in connection to Theorem 12 we refer to 4.13.3 1);
\smallskip
{\bf i)} Theorem 14 was not known before outside of the case of Siegel modular varieties and of trivial cases (like of Shimura curves). 
\medskip
Also we do not hesitate to make connections to the literature which did not influence us, even if this literature appeared later than the first versions of this paper.
\medskip
{\bf 1.16.2. Acknowledgements.} We would like to thank Prof. Gerd Faltings for his encouragement to approach the topics of 1.2, 1.4 and 1.7, for his deep insights and for his suggestion mentioned in 1.2 and in 1.16. Here, as well as in [Va2], we are very much obliged to [Fa1, ch. II and VII] and [Fa2, \S 7]. Though most of the results from the mentioned two loc. cit. are restated in different parts of this paper (see 2.2.1.1 2), 2.2.21, etc.), we still often (as a homage) refer directly to them. Directly (for instance, see [Fo], [FL], [FI] and [CF]) or indirectly (via [Fa1-3] and [Wi]) we are very much obliged to the work of Prof. Jean-Marc Fontaine. We would like to thank Prof. Richard Pink for pointing out the thing mentioned in 1.9, and Prof. Arthur Ogus for some comments and questions, for the suggestion of 3.6.1.4 2) and for explaining us things in connection to the case $p=2$ for ordinary elliptic curves (these things motivated us to include 2.3.18.1). We benefitted from two seminar talks of Prof. Minhyong Kim on [HK] and [Shi1-2]: for them and for some answers we thank him (2.2.4 I to K had them as the very starting point). Parts of AE.6 are based on discussions with Prof. James Stuart Milne: we thank him for them. We would like to thank Ms. Rahel Boller for the typing of [Va13]. We would like to thank FIM, ETH Z\"urich, UC at Berkeley, Univ. of Utah and Univ. of Arizona for providing us with very good conditions for the writing of this manuscript. This work was partially supported by the NSF grant DMS 97-05376.
\vfill\eject
\centerline{}
\bigskip
\bigskip
\centerline{\bigsll {\bf \S 2 Preliminaries}}
\bigskip\bigskip
\medskip
{\bf 2.1. Notations, conventions and some definitions.} Reductive groups over fields are always assumed to be connected. So the fibres of a reductive group over a scheme are connected. For a reductive group $G$ over a scheme $X$ we denote by $G^{\rm der}$, $Z(G)$, $G^{\rm ab}$ and $G^{\rm ad}$, respectively, the derived group of $G$, the center of $G$, the maximal abelian quotient of $G$ and the adjoint group of $G$ (see [SGA3, Vol. III, p. 257] for derived groups and abelianizations). 
The group scheme over $X$ parameterizing automorphisms of $G$ is denoted by $Aut(G)$ (see [SGA3, Vol. III, p. 328]). $G$ is called quasi-split, if it has a Borel subgroup.
\smallskip
We consider now an arbitrary affine group scheme $G$ over $X$. The identity element of the group $G(X)$ is often referred as the origin; if $X$ is connected, the connected component of $G$ containing the origin is referred as the connected component of the origin.
We denote by $Lie(G)$ the Lie algebra of $G$. For an approach to it via representable contravariant functors from the category of $X$-schemes into the category of Lie algebras over different rings of global sections, see [SGA3, Vol. I, Exp. II]. If $G$ is smooth or if $X$ is integral and the tangent space of the origin of any geometric fibre of $G$ has the same dimension, then we always identify $Lie(G)$ with a quasi-coherent locally free sheaf ${\rm Lie}(G)$ on $X$ endowed with a Lie structure; if moreover $X={\rm Spec}(R)$ is affine, we identify ${\rm Lie}(G)$ with a Lie algebra over $R$ whose underlying $R$-module is locally free of locally finite rank. Often we denote ${\rm Lie}(G)$ by ${\got g}$. For any cocharacter $\mu:\GG_m\to G$ and for any affine $X$-scheme $X^1$ we denote by $[\mu_{X^1}]$ the $G(X^1)$-conjugacy class of cocharacters defined by the extension $\mu_{X^1}$ of $\mu$ to $X^1$; if $X^1={\rm Spec}(R^1)$ is affine, we denote $[\mu_{X^1}]$ also by $[\mu_{R^1}]$. 
\smallskip
Let $p$ be a rational prime. We say that a reductive group $G$ over $\QQ$ is unramified over $\QQ_p$ if $G_{\QQ_p}$ is unramified over $\QQ_p$, i.e. if $G_{\QQ_p}$ is quasi-split and splits over the maximal unramified algebraic field extension $\QQ_p^{\rm un}$ of $\QQ_p$. The set of rational primes $l$ such that $G$ is unramified over $\QQ_l$ is denoted by $\Mu(G)$ (cf. [Va2, 6.7]). If $G$ is a reductive group over $S$, where $S$ is either a DVR of mixed characteristic $(0,p)$ or the field of fractions of such a DVR, we denote by $t(G)\in\NN$ the biggest non-negative, integral power of $p$ dividing the order of an element of $G(S)$ of finite order. If $G$ is a reductive group over $\QQ$, its small torsion number $\tilde t(G)\in\NN$ is defined by
$$\tilde t(G):=\prod_{p\in\Mu(G)} t(G_{\ZZ_p}) \prod_{p\notin\Mu(G),\, p\, {\rm a}\, {\rm prime}} t(G_{\QQ_p}),$$
where, for $p\in\Mu(G)$, $G_{\ZZ_p}$ is a reductive group ([Ti2] implies that it does not matter which one) over $\ZZ_p$ having $G_{\QQ_p}$ as its generic fibre (the existence of such numbers $t(*)$ can be easily deduced from the part of AE.4 referring to $n\in\NN$ and the Hints of the Exercise of AE.4.2). Here and in what follows we refer to the chapter ``Addenda and Errata to [Va2]" as AE. We use freely the Definitions 1 and 2 of AE.4.
\smallskip
If $X$ is a set endowed with an equivalence relation $R\subset X\times X$, then for any element $x\in X$, we denote by $[x]\in X/R$ the equivalence class of $x$. Let $n$, $m\in\NN$. An $n+m$-tuple whose entries contain (in the proper order) the entries of an $n$-tuple $t_n$, is referred as an extension of $t_n$. We always look at an $m$-tuple $(n_1,...,n_m)$ whose entries are (denoted by a lower right index to the same small letter and are) in $\ZZ$ from the circular point of view; so for $q\in\ZZ$ we refer, without extra comment, to $n_q:=n_{q_1}$, where $q_1\in\{1,...,m\}$ is such that $m|q-q_1$. We denote by $\abs{B}$ the number of elements of a finite set $B$. For a hermitian symmetric domain $D$ which is a disjoint union of the same connected hermitian symmetric domain $D^0$, we denote by ${\rm mid}(D)$ the maximal (irreducible) dimension (as a complex manifold) of an irreducible factor of $D^0$, cf. [He, p. 381]. 
\smallskip
Our main source of references on Shimura varieties is [Va2].
The expression $(G,X)$ always denotes a pair defining a Shimura variety, $E(G,X)$ denotes its attached reflex field and ${\rm Sh}(G,X)$ denotes the Shimura variety defined by it,
identified with the canonical model of the complex variety. For an arbitrary compact subgroup $K$ of $G(\AA_f)$, we denote by ${\rm Sh}_K(G,X)$ or ${\rm Sh}(G,X)/K$, the quotient of ${\rm Sh}(G,X)$ by $K$. Any
$x\in X$ and $\al\in G(\AA_f)$ define a complex point $[x,\al]$ of ${\rm Sh}_K(G,X)$. As in this paper we encounter several quadruples and triples, the quadruples (resp. triples) introduced in
[Va2, 3.2.6] are called Shimura quadruples (resp. triples). Warning: we use freely the notations and terminology of [Va2, 2.4 and 3.2.6]. For generalities on Shimura varieties of PEL type, we refer to [Ko2, ch. 5] and to [Mi1, p. 161]; warning in last loc. cit. one needs to add that the axiom [Mi1, (SV2)] holds. The embeddings in Siegel modular varieties defining Shimura varieties of PEL type used in [Mi1, p. 161] are referred as PEL type embeddings.
\smallskip
If $A$ is a set and $f:A\to A$ is an endomorphism of it, we denote by $f^m$ the $m$-times composite of $f$ with itself; the identity bijection of $A$ is denoted by $1_A$. If $f:A\to B$ and $g:B\to C$ are morphisms in some category we refer to $g\circ f$ as the composite of $f$ with $g$. If $0\to A\to B\to C\to 0$ is any short exact sequence, we refer to $B$ as the extension of $C$ by $A$.   
\smallskip
If $k$ is a field we denote by $\bar k$ its algebraic closure. Let ${\rm Gal}(k)$ the group of automorphisms of $\bar k$ acting identically on (the perfection of) $k$. For a field $k$ of characteristic $p$,
$W(k)$ is the DVR of Witt vectors of $k$, while $B(k):=W(k)[{1\over p}]$ is the field of fractions of $W(k)$. We denote ${\rm Gal}(B(k))$ by $\Gamma_k$. ${\bf \al}_p$ is the finite, flat, commutative group scheme over $k$ of rank $p$, defined as the kernel of the Frobenius endomorphism $F:\GG_a\to\GG_a$. For $D$ a $p$-divisible group or a finite, flat, commutative group scheme over $k$, we denote by $a(D):=\dim_k({\rm Hom}({\bf \al}_p,D))$ its $a$-number. Occasionally, we use the Witt ring $W(R)$, with $R$ an arbitrary $k$-algebra, and $W_n(R):=W(R)/p^nW(R)$. 
\smallskip
In all that follows we assume $k$ is perfect. We denote by $\sg$ (resp. $\overline{\sg}$) the Frobenius automorphism of $W(k)$ and $B(k)$ (resp. of $W(\bar k)$ and $B(\bar k)$). Often, to avoid repetitions and for the sake of preciseness, we also denote $\sg$ by $\sg_k$. 
\smallskip
By a $\sg$-crystal (over $k$) we mean a pair $(M,\vph)$, with $M$ a free $W(k)$-module of finite rank and with $\vph$ a $\sg$-linear endomorphism of $M$ which is a $\sg$-linear automorphism of $M[{1\over p}]$. We refer to $(M,\vph^n)$ as a $\sg^n$-crystal. Warning: by a filtered $\sg$-crystal (over $k$) we mean a triple $(M,F^1,\vph)$, with $(M,\vph)$ a $\sg$-crystal and with $F^1$ a $W(k)$-submodule of $M$ which is a direct summand of $M$ such that $\vph({1\over p}F^1+M)=M$; this convention is not standard but, as we often deal just with contexts involving $p$-divisible groups, it is very practical. We always denote by $F^0$ a direct supplement of $F^1$ in $M$. By a principal quasi-polarization of a $\sg$-crystal $(M,\vph)$ we mean a perfect alternating form $\psi:M\otimes_{W(k)} M\to W(k)$ such that $\psi(\vph(x),\vph(y))=p\sg(\psi(x,y))$, $\forall x,y\in M$. We often write it as $\psi:M\otimes_{W(k)} M\to W(k)(1)$. 
\smallskip
By a $\sg^n$-isocrystal (over $k$) we mean a pair $(M,\vph)$, with $M$ a finite dimensional $B(k)$-vector space and with $\vph$ a $\sg^n$-linear automorphism of $M$. We still denote by $\vph$ the Frobenius endomorphism of any free $W(k)$-submodule or $B(k)$-vector subspace $N$ of different tensor products (taken in any order) of a given number of copies of $M$ with a given number of copies of its dual $M^*$, which is taken by $\vph$ into itself; here $\vph$ takes $f\in M^*$ into $\sg^n\circ f\circ\vph^{-1}\in M^*$ and acts in the usual tensor way on all such tensor products. Warning: if $\vph$ takes $N$ into $N$, then $\forall m\in\NN$ we make the convention that $p^m\vph$ acts on $N$ via the rule: $n\in N$ is mapped into $p^m\vph(n)$; so $p^m\vph$ is an abbreviation for $(p^m1_N)\circ\vph$. In particular, $\vph$ takes $e\in {\rm End}(M)=M\otimes_{B(k)} M^*$ into $\vph\circ e\circ\vph^{-1}$, while $p\vph$ takes $e$ into $p(\vph\circ e\circ\vph^{-1})$. By a latticed $\sg^n$-isocrystal (over $k$) we mean a pair $(M,\vph)$, with $M$ a free $W(k)$-module of finite rank, such that the pair $(M[{1\over p}],\vph)$ is a $\sg^n$-isocrystal. A morphism between two latticed $\sg^n$-isocrystals $(M_1,\vph_1)\to (M_2,\vph_2)$ is a $W(k)$-linear map $f:M_1\to M_2$ such that by inverting $p$ it is a morphism between $\sg^n$-isocrystals. When $n=1$, we speak simply about (latticed) isocrystals (over $k$). 
\smallskip
Always $V_0$ denotes a Witt ring over an algebraically closed
field of positive characteristic and then $K_0$ automatically denotes its field of fractions. If $v$ is a prime of a number field $E$, we denote by $k(v)$ the residue field of it and by $O_{(v)}$ the localization of the ring of integers of $E$ with respect to it. The normalization of $\ZZ_{(p)}$ in $E$ is denoted by $E_{(p)}$. If $O$ is a local ring, we denote by $O^{\rm h}$, $O^{\rm sh}$ and $\widehat O$ respectively its henselization, its strict henselization and its completion.
\smallskip
We usually write $\ZZ_{(p)}$ instead of $O_{(p)}$. The ring of finite ad\`eles
$\hat\ZZ\otimes_\ZZ\QQ$ is denoted by $\AA_f$ and the ring of finite ad\`eles with the $p$-component omitted is denoted by $\AA^p_f$. So we have $\AA_f=\AA_f^p\times\QQ_p$. We use freely different Tate-twists:
$\QQ(1),\QQ_p(1),\ZZ_p(1),\AA_f(1)$ etc. We denote by $\FF_{p^d}$ the field with $p^d$ elements ($d\in \NN$) and by $\FF$ its algebraic closure. If $G$ is a reductive group over $\QQ$, then its torsion number $t(G)\in\NN$ is defined by
$$t(G):=\prod_{p\in\Mu(G)} t(G_{W(\FF)}) \prod_{p\notin\Mu(G),\, p\, {\rm a}\, {\rm prime}} t(G_{B(\FF)}).$$
We have: $\tilde t(G)$ divides $t(G)$. An integral, separated $\FF_p$-scheme is called perfect if it has an open, affine cover by spectra of perfect $\FF_p$-algebras.
\smallskip
If $R$ is a faithfully flat $\ZZ_{(p)}$-algebra, we denote by $R^\wedge$ its $p$-adic completion. Similarly, if $X={\rm Spec}(R)$, with $R$ as before, we denote $X^\wedge:={\rm Spec}(R^\wedge)$. 
\smallskip
By a $p$-adic formal scheme $X$ we mean a formal scheme which locally can be obtained by completing a faithfully flat $\ZZ_{(p)}$-scheme along its special fibre. By a morphism from $X$ into another $p$-adic formal scheme $Y$ or into a scheme $Y$, we mean a morphism $f:X\to Y$ of locally ringed spaces; we refer to $f$ as a $p$-adic formal scheme over $Y$. So we have a natural morphism $X\to {\rm Spec}(\ZZ_p)$. The fibre $X_{\ZZ_p/p^s\ZZ_p}$ is a scheme, $\forall s\in\NN$. If $X$ is a flat $\ZZ_p$-scheme, we denote by $X^\wedge$ the $p$-adic formal scheme it defines naturally. Warning: if $X={\rm Spec}(R)$, then without a special reference $X^\wedge$ denotes ${\rm Spec}(R^\wedge)$; accordingly, without a special reference, the $p$-adic completion of an affine scheme (resp. of a morphism $f$ between affine schemes) is viewed as an affine scheme (resp. as a morphism $f^\wedge$ between affine schemes).
The sheaf of local rings implicit in the definition of a scheme or of a formal scheme $X$, is denoted by $\Mo_X$. The category of locally free $\Mo_X$-sheaves of locally finite ranks (resp. of coherent $\Mo_X$-sheaves) is denoted by $LF(X)$ (resp. by $(COH(X)$). By a geometric point of a scheme we mean a point of it with values in an algebraically closed field. 
\smallskip
The category of $p$-divisible groups (resp. of finite, flat, commutative group schemes of $p$ power order) over a (formal) scheme $S$ is denoted by $p-DG(S)$ (resp. by $p-FF(S)$). The reduced scheme of a scheme $S$ is denoted $S_{\rm red}$. If $S$ is an $\FF_p$-scheme, then without a special reference on their type, all crystals and $F$-crystals on $S$ are in locally free sheaves which locally have finite ranks and, if they are filtered, the filtrations are given by direct summands. All $F$-crystals are such that by pulling them back through geometric points of $S$ we get isocrystals over algebraically closed fields as defined above. If $S$ is a $\ZZ_{(p)}$-scheme (resp. $W(k)$-scheme) on which $p$ is locally nilpotent, the crystalline site we use $CRIS(S/{\rm Spec}(\ZZ_p))$ (resp. $CRIS(S/{\rm Spec}(W(k)))$) is Berthelot's crystalline site as used in [Be, ch. III, \S 4] and [BM]. The PD-structure on the ideal $(p)$ of $\ZZ_p$ (resp. of $W(k)$)) is the usual one from characteristic 0.
\smallskip
A morphism $X_1\to X$ of formal schemes is called \'etale, if locally in their topology it is defined by a morphism ${\rm Spf}(A_1)\to {\rm Spf}(A)$, where $A$ (resp. $A_1$) is complete with respect to the topology defined by an ideal $I$ of itself (resp. with respect to the topology defined by its ideal $IA_1$), such that the morphism ${\rm Spec}(A_1/I^sA_1)\to {\rm Spec}(A/I^s)$ is \'etale, $\forall s\in\NN$. Similarly we define formally \'etale morphisms between formal schemes.
\smallskip
If $G$ is an affine group scheme over $W(k)$, the $W(k)$-valued point of the $p$-adic completion of $G$ or of the completion of $G$ with respect to its origin, defined by the origin of $G$, is still referred as origin.
\smallskip
Let $K$ be an arbitrary field. By a stratification $\Ms$ of a reduced $K$-scheme $S$ in potentially an infinite number of strata, we mean that:
\medskip
a) for any field $L$ which is either $K$ or an algebraically closed field containing $\bar K$ but different from $\bar K$, a set $\Ms_L$ of disjoint reduced, locally closed subschemes of $S_L$ is given such that each $\bar L$-valued point of $S$ factors through an element of it;
\smallskip
b) if $i_{12}:L_1\hookrightarrow L_2$ is an embedding between two fields as in a), then the pull back to $L_2$ of any member of $\Ms_{L_1}$, is an element of $\Ms_{L_2}$ (so we have a natural transition injective map $\Ms(i_{12}):\Ms_{L_1}\hookrightarrow\Ms_{L_2}$). 
\medskip
Most of the time, we just mention how $\Ms_K$ is constructed (as the constructions of $\Ms_L$'s and of $\Ms(i_{12})$'s are automatic). The inductive limit of all maps $\Ms(i_{12})$ is a class which in general is not a set. We refer to it as the class of $\Ms$. Each element of some $\Ms_L$ is referred as a stratum of $\Ms$. $\Ms$ is said to satisfy the purity property if for each $L$ as above, any element of $\Ms_L$ is an affine $S$-scheme; this is a more practical and general definition then any other one relying on codimension 1 statements on complements (of course, if $\Ms$ satisfies the purity property and $S_L$ is locally noetherian, then the Zariski closure $T_1$ of any irreducible component $T_0$ of an element of $\Ms_L$ in $S_L$ is such that $T_1\setminus T_0$ is either empty or of pure codimension 1). Without a special reference, we always deal with stratifications of $S$ such that all transition injective maps $\Ms(i_{12})$ are bijections and $\Ms_K$ is a finite set. 
\smallskip
By a Frobenius lift of a $\ZZ_{(p)}$-scheme $X$ (resp. of a $\ZZ_{(p)}$-algebra $R$) we mean an endomorphism of $X$ (resp. of $R$) which mod $p$ is the usual Frobenius endomorphism (for instance, for $R/pR$, it takes $x\in R/pR$ into $x^p$). If $R$ and $S$ are two arbitrary $\ZZ_{(p)}$-algebras having Frobenius lifts $\Phi_R:R\to R$ and respectively $\Phi_S:S\to S$, if $M$ is a projective $R$-module endowed with a $\Phi_R$-linear map $\vph:M\to M$, and if $b:R\to S$ is a $\ZZ_{(p)}$-homomorphism compatible with the mentioned Frobenius lifts, we denote by $\vph\otimes 1$ the $\Phi_S$-linear map on $M\otimes_R S$ induced from $\vph$ via $b$. As above, if $R$ is a flat $\ZZ_{(p)}$-algebra and if $\vph$ becomes an isomorphism by inverting $p$, we still denote by $\vph$ the $\Phi_R$-linear endomorphism of any $R$-submodule or $R[{1\over p}]$-submodule of a tensor product of a number of copies of $M[{1\over p}]$ with a number of copies of its dual, taken by $\vph$ into itself. Moreover, we use the same convention on $p^m\vph$, $m\in\NN$, as above (in the case $R=W(k)$).
\smallskip
Let $a$, $b\in\ZZ$, $b\ge a$. Let $S(a,b):=\{a,a+1,...,b-1,b\}$ and let $SS(a,b):=\{a-b,a-b+1,...,b-a-1,b-a\}$. We denote by $\Mm\Mf_{[a,b]}(*)$,
$\Mm\Mf(*)$, $\Mm\Mf^\nabla_{[a,b]}(*)$ and $\Mm\Mf^\nabla(*)$ the categories introduced in [Fa1, ch. II]; here $*\in\{R,X\}$. Whenever required (see loc. cit.), the $p$-adic completion of $*$ is equipped with a Frobenius lift. We refer to them as Fontaine categories (of objects); an arbitrary such category is often denoted as $FC$. We allow $R$ (resp. $X$) to be (the $p$-adic completion of) a regular, formally smooth $W(k)$-algebra (resp. $W(k)$-scheme or $p$-adic formal scheme over $W(k)$). The objects of $FC$ are denoted using tuples. For instance, fixing a Frobenius lift $\Phi_R$ of $R^\wedge$, an object of $\Mm\Mf_{[a,b]}^\nabla(R)$ is denoted as a quadruple $(M,(F^i(M))_{i\in S(a,b)},(\vph_i)_{i\in S(a,b)},\nabla)$: $M$ is an $R$-module such that each element of it is annihilated by some power of $p$, $(F^i(M))_{i\in S(a,b)}$ is a decreasing filtration of $M$, $\vph_i:F^i(M)\to M$ is a $\Phi_R$-linear map, and $\nabla$ is a connection on $M$ (the axioms being as in loc. cit.). Warning: in 2.2.4 we adopt the (naive) point of view that all Fontaine categories are small, i.e. its morphisms are forming a set. This can be fully formalized using the notion of universe as in [SGA4]; however, we prefer to keep the requirements to the minimum and so we never mention such a universe (i.e. we think that from the point of view of 2.2.4 J, the naive approach  is ``enough", as we can work as well with skeletons of categories).  
\smallskip
We always assume, for the sake of simplicity of presentation, that the $R/pR$-module $\Om_{R/pR/k}$ (resp. the $\Mo_{X_k}$-sheaf $\Om_{X_k/k}$) of relative differentials is projective (resp. is locally free) of locally finite rank in each point of ${\rm Spec}(R/pR)$ (resp. of $X_k$). In what follows we refer to [Fa1, 2.1], with $R$ and $X$ as in the above paragraph, without any other further comment: the proofs of loc. cit. apply entirely to the context we have allowed (working in the faithfully flat topology we can assume $R$ is a complete, local ring having an algebraically closed residue field; for this context the inductive argument on dimensions of loc. cit. applies entirely). We refer to [Fa1, 2.3] as the standard gluing process (or arguments) in the context of Fontaine categories involving connections (the word ``gluing" is used for the first time in the last sentence of the paragraph before [Fa1, e) of p. 35]). 
\smallskip
We can speak as well about Fontaine categories for regular, formally smooth, affine formal schemes (when allowed just formal schemes) over $W(k)$, as well as we can replace $W(k)$ by $W_n(k)$ (in this last case we deal only with objects annihilated by $p^n$ in the logical sense recalled in 2.2.1 c) below); the context of formal schemes (which are not necessarily $p$-adic formal schemes) shows up just in 3.15.4. 
\smallskip
As a policy, besides some definitions, we always start with some Fontaine category involving a scheme $X$ and, as the constructions force us, we move to some Fontaine categories involving a $p$-adic formal scheme over $X$. Just one exception: 3.6.18.4.2 is stated simultaneously in the context of schemes as well as of $p$-adic formal schemes.  
\smallskip
We often refer to $\Phi_R$ as a Frobenius lift of ${\rm Spec}(R)$ and so we write as well $\Phi_R:{\rm Spec}(R)\to {\rm Spec}(R)$. If $r:W(k)\to R$ is a $\ZZ_{(p)}$-homomorphism and $\Phi_R(r(x))=r(\sg(x))$, $\forall x\in W(k)$, we say $\Phi_R$ is compatible with $\sg$. In what follows, we always assume that the Frobenius lifts of $W(k)$-algebras are compatible with $\sg$. The same applies to the context of regular, formally smooth ($p$-adic formal) schemes over $W(k)$.\smallskip
We use freely [Be, 1.6.5, p. 247] as well as its filtered version; accordingly, we have a natural functor $\Mf_{[0,b]}$ from the category $\Mm\Mf_{[0,b]}^\nabla(X)$ to the category of filtered $F$-crystals on $X_k$ on coherent sheaves; here $b\in\NN\cup\{0\}$ or $b\in S(0,p-1)$ depending on the fact that $X^\wedge$ is or is not equipped with a Frobenius lift. For instance, if $X={\rm Spec}(R^\wedge)$ and a Frobenius lift of $X$ is fixed, then $\Mf_{[0,b]}$ takes an object $(M,(F^i(M))_{i\in S(0,b)},(\vph_i)_{i\in S(0,b)},\nabla)$ of $\Mm\Mf_{[0,b]}^\nabla(X)$ into the filtered $F$-crystal on $X_k$ whose evaluation at $X$ is the quadruple $(M,(F^i(M))_{i\in S(0,b)},\vph_0,\nabla)$ (the action of $\Mf_{[0,b]}$ on morphisms being the logical one).
\smallskip
For fitting the convention of [De2, 1.1.6] on Hodge structures, the canonical split (cocharacters) of objects of $p-\Mm\Mf(W(k))$ (see 2.2.1 c) for the definition of this category) to be used are the inverses of the ones defined in [Wi, p. 512]. So whenever we refer to them we mean the inverses of the cocharacters of loc. cit., though they are still referred as canonical splits. We use freely (the functorial aspect of) loc. cit. in the formed mentioned.
\smallskip
The zero element of any module or ring is denoted by $0$. For any finitely generated projective module $M$ over a commutative ring $R$, we denote by $M^\ast$ its dual. The $R$-span (i.e. the $R$-submodule of $M$) generated by a finite number of elements $e_1$,..., $e_s$ of $M$, is denoted by $<e_1,...,e_s>$. We denote by $M^{\otimes s}\otimes_R M^{*\otimes t}$, with $s,t\in \NN\cup\{0\}$, the tensor product of $s$-copies of $M$ with $t$-copies of $M^*$ taken in this order. If  $v\in M^{\otimes s}\otimes_R M^{*\otimes t}$, we denote by ${\rm deg}(v):=s+t$ its degree, and we call it a homogeneous tensor of degree $s+t$. By the essential tensor algebra of $M\oplus  M^\ast$ we mean 
$$
\Mt(M):=\oplus_{s,t\in\NN\cup\{0\}} M^{\otimes s}\otimes_R M^{*\otimes t}.
$$ 
We do not use the natural structure of $\Mt(M)$ as an $R$-algebra (it is non-commutative if the rank of $M$ is at least 2 in some point of ${\rm Spec}(R)$). We always identify ${\rm End}(M)$ with $M\otimes_R M^*$. So ${\rm End}({\rm End}(M))={\rm End}(M\otimes_R M^*)=M\otimes_R M^*\otimes_R M^*\otimes_R M$ is always identified (by changing the order) with $M^{\otimes 2}\otimes_R M^{*\otimes 2}$. A family of tensors of $\Mt(M)$ is denoted in the form $(v_{\alpha})_{\alpha\in\Mj}$, with $\Mj$ as the set of indices; if all its tensors are homogeneous of degree not greater than $u\in\NN\cup\{0\}$, we say it is of degree at most $u$. Let $M_1$ be another finitely generated, projective $R$-module. Any isomorphism $f:M\tilde\to M_1$ extends naturally to an isomorphism $\Mt(M)\tilde\to\Mt(M)$ and so we often speak about $f$ taking $v_j$ into some specific element of $\Mt(M_1)$.
\smallskip
If $M$ is endowed with a decreasing filtration by direct summands $F^i(M)$, $i\in\ZZ$, then $\Mt(M)$ gets naturally equipped with a filtration (called the tensor product filtration) by direct summands; moreover, whenever we have a direct summand $N$ of $\Mt(M)$ such that its intersection with any direct summand $F^i(\Mt(M))$ of $\Mt(M)$ which is part of its tensor product filtration, is a direct summand of $N$, we endow it with this intersection (also called induced) filtration; so $F^i(N):=N\cap F^i(\Mt(M))$, $\forall i\in\ZZ$. We refer to $F^i(\Mt(M))$ (resp. to $F^i(M)$) as the $F^i$-filtration of $\Mt(M)$ (resp. of $M$). In particular, we often refer to the $F^0$-filtration of $\Mt(M)$ or of ${\rm End}(M)$ or to the $F^1$-filtration of ${\rm End}(M)$. If $F^1$ is a direct summand of $M$, by the $F^0$-filtration of $\Mt(M)$ defined by it, we mean the $F^0$-filtration of $\Mt(M)$  defined by the filtration of $M$ defined by: $F^2(M)=\{0\}$, $F^1(M)=F^1$ and $F^0(M)=M$. 
\smallskip
All filtrations of $M$ to be considered are decreasing and such that $F^i(M)=M$ for $i<<0$ and $F^j(M)=\{0\}$ for $j>>0$; so often we just mention those direct summands $F^i(M)$'s which are different from $\{0\}$ and $M$. Warning: this does not apply to Fontaine's rings of 2.3.18.1 E and of 4.9.2.0. By a range of a filtration of $M$ we mean an interval $[a,b]$, with $a$, $b\in\ZZ$, such that $F^a(M)=M$ and $F^{b+1}(M)=\{0\}$. If $F^2(M)=\{0\}$ and $F^0(M)=M$, then we just mention the direct summand $F^1(M)$, and we usually denote it (in a simpler way) by $F^1$; moreover, a direct supplement of $F^1$ in $M$ is denoted by $F^0$. If we use some lower right index $j$ to denote the filtration of $M$ (i.e. if we work with $F^i_j(M)$, $i\in\ZZ$), then we automatically use the same type of index for all other filtrations introduced in the above paragraph: we write $F^i_j(\Mt(M))$, $F^i_j(N)$, etc. Warning: we still refer to $F^i_j(\Mt(M))$ as the $F^i$-filtration of $\Mt(M)$, etc. Similarly, the right lower indices are transfered in the context of $F^1$ and $F^0$: we often consider $F^1_j$ and a supplement of it $F^0_j$ in $M$. 
\smallskip
Similarly, we define the essential tensor algebra $\Mt(\Mf)$ of $\Mf\oplus\Mf^*$, for a locally free $\Mo_X$-sheaf $\Mf$ of locally finite rank on a ringed space $(X,\Mo_X)$. If this ringed space is a scheme, then $GL(\Mf)$ is the reductive group scheme over it of $\Mo_X$-linear automorphisms of $\Mf$. 
\smallskip
A bilinear form on $M$ is called perfect if it induces an isomorphism
from $M$ into its dual $M^\ast$. A pair $(M,\psi)$ with $M$ as above and $\psi:M\otimes_R M\to R$ a perfect alternating form on it, is called a symplectic space over $R$. $GL(M)$, $SL(M)$, $Sp(M,\psi)$, $GSp(M,\psi)$, etc., are viewed as group schemes over $R$; so, for an $R$-algebra $R_1$, $GL(M)(R_1)$ denotes the group of $R_1$-valued points of $GL(M)$ (i.e. the group of $R_1$-linear automorphisms of $M\otimes_R R_1$), etc. Warning: if $G$ is a closed subgroup of $GL(M)$, we often identify a morphism ${\rm Spec}(R_1)\to G$, with a matrix automorphism of $M\otimes_R R_1$; in particular, we often refer to $1_M$ as the identity element of $G(R)$. If $f_1$ and $f_2$ are two $\ZZ$-endomorphisms or automorphisms of $M$ (or of the localization of $M$ with respect to some element of $R$), most common we denote $f_1\circ f_2$ by $f_1f_2$. We usually use the same notation for two perfect bilinear forms if they are obtained one from another by an extension of scalars; the same applies to tensors of two essential tensor algebras obtained one from another by an extension of scalars. If $(W,\psi)$ is a symplectic space over $\QQ$, we denote by
${\rm Sh}\bigl({\rm GSp}(W,\psi),S\bigr)$ the Siegel modular variety associated to it (for instance, see [Va2, Example 2 of 2.5]). 
\smallskip
For a finite field extension $k_0\hookrightarrow k$ and for a reductive group $G$ over $k$,
${\rm Res}_{k/k_0}G$ denotes the group scheme over $k_0$ obtained from $G$ by Weil restriction of scalars (see [BLR, 7.6]). If $f:G\to G_1$ is a homomorphism of group schemes which is a closed embedding, then we identify $G$ with its image (through $f$) in $G_1$. 
\smallskip
We  use freely the terminology of Hodge cycles of abelian varieties over a field of characteristic $0$ used in [De4] and, by natural extension, this terminology is also used for abelian schemes over a reduced $\QQ$--scheme $S$. Even more: if $S$ is a reduced, flat $\ZZ_{(p)}$-scheme and $A$ is an abelian scheme over $S$, we speak about Hodge cycles of $A$, though strictly speaking, we should speak about Hodge cycles of the generic fibre of $A$. Any polarization of an arbitrary abelian scheme $A$ is denoted by $p_A$ (sometimes $\Mp_A$) and, by abuse of notation, we also denote by $p_A$ the different maps on the cohomologies (homologies) of $A$ induced by it. Pairs of the form  $(A,p_A)$  always denote a polarized abelian scheme. If $\Ms$ is a set of rational primes and if $\ZZ[{1\over {\Ms}}]$ is the subring of $\QQ$ generated by $\ZZ$ and the inverses of all primes of $\Ms$, then by a $\ZZ[{1\over {\Ms}}]$-isogeny between two abelian schemes $A_1$ and $A_2$ defined over a scheme, we mean a $\QQ$--isomorphism which induces an isomorphism $A_1[N]\tilde\to A_2[N]$, $\forall N\in\NN$ which is relatively prime to all primes in $\Ms$. 
\smallskip
We consider an abelian scheme or a $p$-divisible group $A$ over a scheme $S$. We denote by $A^t$ its dual and by $A[p^n]$ the kernel of its multiplication by $p^n$. We now assume $S$ is the spectrum of a perfect field $k_1$ of positive characteristic $p$. By the $\sg_{k_1}$-crystal (resp. by the isocrystal) of (or attached to) $A$ we mean the $\sg_{k_1}$-crystal defined by the first crystalline cohomology group of $A$ with coefficients in $W(k_1)$ (resp. we mean the isocrystal over $k_1$ defined by such a $\sg_{k_1}$-crystal). Similarly, if $S={\rm Spec}(W(k_1))$, we speak about the filtered $\sg_{k_1}$-crystal of (or attached to) $A$. If $S={\rm Spec}(k_1)$, with $k_1$ arbitrary, by the formal isogeny type of $A$, we mean the formal isogeny type (see [Man]) defined by the isocrystal of $A_{\overline{k_1}}$. We use the notation $(n,m)$ instead of the notation $G_{n,m}$ of loc. cit; so $(1,1)$ is the formal isogeny type of (the isocrystal of) a supersingular elliptic curve, etc.  
\smallskip
If $S$ is the spectrum of a field of characteristic zero and if $A$ is an abelian variety, then we speak about the Shimura variety ${\rm Sh}(G,X)$ attached to $A$: choosing an arbitrary polarization $p_A$ of $A$, the Shimura pair $(G,X)$ is defined as in [Va2, 2.12 3)]; it is independent on the choice of $p_A$.
\smallskip
For Artin's approximation theorem we refer to [BLR, p. 91]. We use it freely as well as its intermediary form of [BLR, th. 12 of p. 83]. For classical Dieudonn\'e theories over $k$ or $W(k)$ we refer to [Fo]; mainly we need [Fo, p. 152 and 160 and \S 5 of ch. 4]. For Lang's theorem on (torsors of) connected algebraic groups over finite fields, we refer to [Bo2, 16.3-4 and 16.6]. For standard properties and terminology of involutions on semisimple algebras over fields we refer to [KMRT, ch. 1]. For the specialization theorem (of Grothendieck--Katz) we refer to [Ka2, 2.3.1-2]. For Grothendieck's algebraization theorem we refer to [EGA III, ch. 5]. For Serre--Tate's deformation theory we refer to [Me, ch. 5] and [Ka3, ch. 1]. For Zariski's Main Theorem we refer to [Ra, ch. IV]. For the ordinary conjecture concerning the ordinary reductions of an abelian variety over a number field, we refer to [Oo1, p. 11]. We use freely the fact that any regular, $p$-adically complete, formally smooth $W(k)$-algebra has Frobenius lifts. 
\smallskip
All limits and colimits to be used are assumed to be filtered. Projective systems of (formal) schemes have affine transition morphisms. A projective limit of a projective system indexed in the standard manner by elements of $\NN$, is referred as an $\NN$-projective limit, while the projective system itself is called an $\NN$-projective system; similarly we speak about an $\NN$-inductive limit or system or about an $\NN$-pro-\'etale scheme over another scheme. As the transition morphisms are affine, the projective limit of a projective system of formal scheme with \'etale transition morphisms, does exist as as formal scheme. So similarly, we speak about an ($\NN$-) pro-\'etale formal scheme over another formal scheme.
\smallskip 
Occasionally, we refer to the $\NN$-pro-\'etale topology of a scheme: this can be defined entirely similar to the \'etale topology. What we need: a property $\Mp$ pertaining to an ``object" $O$ over a scheme $S$ is satisfied locally in the $\NN$-pro-\'etale topology of $S$, if for each point $y$ of $S$ with values in a field there is an $\NN$-pro-\'etale morphism $S_1\to S$ whose image contains $y$, such that $\Mp$ holds for the (assumed to make sense) pull back of $O$ to $S_1$. The same applies to the context of formal schemes.
\smallskip
We often refer to a Newton polygon being below, above, strictly below or strictly above another one, in the usual sense (so the Newton polygon of a supersingular elliptic curve over $k$ is strictly above the Newton polygon of an ordinary elliptic curve over $k$). The Newton polygon of a set of Newton polygons which is below all of them, is called the smallest Newton polygon of the set. We use freely the specialization theorem to define stratifications using Newton polygons, even in a context where we are not dealing with an $F$-crystal over an $\FF_p$-scheme $S$ but with ``modifications" of such an $F$-crystal, where the Frobenius endomorphisms of modules over flat $\ZZ_p$-algebras associated naturally to it are divided by the same positive, integral power of $p$. This applies in particular to the context of Lie $F$-crystals as to be defined below (for instance, see 2.2.2). In particular, we allow Newton polygons which have negative slopes and so, we deviate from the initial approach suggested in [Va2, 1.6.2] of using Tate twists in order to get only non-negative slopes. By the Newton polygon of a latticed isocrystal $(M,\vph)$ we mean the Newton polygon of $(M[{1\over p}],\vph)$.
\smallskip
The set (often viewed as a matrix algebra over $R$) of $m\times m$ matrices with coefficients in $R$ is denoted by $M_m(R)$; the group of invertible such matrices is denoted by $GL_m(R)$. If $pR=\{0\}$ and $B\in M_m(R)$, then $B^{[p]}$ denotes the matrix whose entries are $p$-powers of entries of $B$.
\smallskip
Most common we avoid using ``an index to another index"; for instance, if $G_1$ (resp. $G^1$) is a reductive $\QQ$-group, then its extension to $\QQ_p$ (resp. its derived subgroup) is usually denoted by $G_{1\QQ_p}$ (resp. by $G^{1\rm der}$) and not by ${G_1}_{\QQ_p}$ (resp. by ${G^1}^{\rm der}$). 
\smallskip
For the theories of healthy schemes and of different extension properties (like $EP$, $EEP$, $SEP$, etc.) see [Va2, 3.2.1-3]. The expression ``with respect to" is abbreviated as w.r.t. For minimal weights we refer to [Bou3, ch. 8, \S 7.3], to [Sa] and to [De2, p. 261]. For the theory of well positioned families of tensors we refer to [Va2, 4.3.4-17] (see also AE.3 for the correction of [Va2, 4.3.6 3)]); we use freely the terminology of [Va2, 4.3.4]. 
\medskip\smallskip
{\bf 2.2. Fontaine categories and Shimura crystals.} In what follows we introduce a language to be consistently used in \S2-14 as well as in the subsequent papers (like [Va3-12]). It is an ``interplay" between the general, abstract contexts of Fontaine categories and the particular, concrete contexts which have direct applications to the study of Shimura varieties.
Let $k$ be an arbitrary perfect field of characteristic $p>0$. 
\medskip
{\bf 2.2.1. Basic definitions and few comments. a)} Let $r\in\QQ$. An isocrystal $(M,\vph)$ over $k$ is said to be $r$-symmetric, if $\forall\al\in\QQ$ such that $r-\al$ is a slope of it, $r+\al$ is also a slope of it of the same multiplicity. 
\medskip
{\bf b)} We say that a $k$-scheme $X_k$:
\medskip
-- has the $ALP$ property, if the only open subscheme of it containing all its $\bar k$-valued points is $X_k$ itself; 
\smallskip
-- is an almost global (to be abbreviated as $AG$) $k$-scheme, if for any connected component $\Mc$ of it there is a partition of the reduced scheme $\Mc_{\rm red}$ in a finite number of reduced, locally closed subschemes, each stratum being an $\NN$-pro-\'etale cover of a $k$-scheme of finite type.
\medskip
$ALP$ stands for ``algebraic points". The $ALP$ property is stable under pro-\'etale covers; so any $AG$ $k$-scheme has the $ALP$ property.
\medskip
{\bf c)} Let $R$ be a regular, formally smooth $W(k)$-algebra. We have in mind mainly the case when ${\rm Spec}(R/pR)$ has the $ALP$ property (like the case when ${\rm Spec}(R/pR)$ is a local or an $AG$ $k$-scheme); but for the sake of generality we do not assume we are dealing just with such a case. Let $\Phi_R$ be a Frobenius lift of $R$ (resp. of $R^\wedge$).  
\smallskip
A finitely generated projective $R$-module
(resp. $R^\wedge$-module) $M$, together with a decreasing filtration $(F^i(M))_{i\in \ZZ}$ of it by direct summands and a $\Phi_R$-linear endomorphism 
$$\vph:M[{1\over p}]\to M[{1\over p}],$$
is said to be a $p$-divisible object of $\Mm\Mf(R)$ if for any $i\in\ZZ$, $\vph$ restricted to $F^i(M)$ is of the form $p^i\vph_i$, with 
$$\vph_i:F^i(M)\to M,$$ 
and the $R$-linear (resp. $R^\wedge$-linear) map 
$$\vph_{\ZZ}:\oplus_{i\in\ZZ} \tilde F^i(M)\otimes_R\, _{\Phi_R}R\to M$$ 
(resp. $\vph_{\ZZ}:\oplus_{i\in\ZZ} \tilde F^i(M)\otimes_{R^\wedge}\, _{\Phi_R} R^\wedge\to M$), whose restriction to the summand $\tilde F^i(M)\otimes_R {}_{\Phi_R}R$ (resp. to $\tilde F^i(M)\otimes_{R^\wedge} {}_{\Phi_R} R^\wedge$) acts as $\vph_i\otimes 1$, is an isomorphism; here, for $i\in\ZZ$, $\tilde F^i(M)$ is a direct supplement of $F^{i+1}(M)$ in $F^i(M)$. We denote such a $p$-divisible object by a triple 
$${\got C}:=(M,(F^i(M))_{i\in\ZZ},\vph). 
$$
\indent
For any $n\in \NN$, $M/p^nM$ together with the filtration
$(F^i(M)/p^nF^i(M))_{i\in \ZZ}$ of it and the induced family of $\Phi_R$-linear maps $(\vph_i)_{i\in\ZZ}$ is an object of $\Mm\Mf(R)$ (we recall, cf. 2.1, that $\tilde F^i(M)=\{0\}$ for $i$ very big or very small); here we still denote by $\vph_i$ its reduction mod $p^n$. We denote this object by ${\got C}/p^n{\got C}$ and we call it the truncation mod $p^n$ of ${\got C}$; often, for the sake of convenience, as $\vph$ and $F^i(M)$'s determine $\vph_i$'s, we denote it also by the triple 
$$(M/p^nM,(F^i(M)/p^nF^i(M))_{i\in\ZZ},\vph)$$ 
instead of the triple $(M/p^nM,(F^i(M)/p^nF^i(M))_{i\in\ZZ},(\vph_i)_{i\in\ZZ})$; warning: here we do not think of $\vph$ as some reduction modulo some power of $p$. We refer to $M$ as the $R$-module (resp. the $R^\wedge$-module) underlying ${\got C}$ and to $M/p^nM$ as the $R/p^nR$-module underlying ${\got C}/p^n{\got C}$. We refer to $(F^i(M))_{i\in\ZZ}$ as the filtration underlying ${\got C}$, etc. Similarly, we speak about the (filtered) module or sheaf (of modules) underlying an object of any Fontaine category and about the truncation mod $p^n$ of any such object; if an object is isomorphic to its truncation mod $p^n$, we say it is annihilated by $p^n$. If $F^0(M)=M$, then $\vph$ takes $M$ into $M$ and we refer to it as a Frobenius endomorphism of $M$ and we write $\vph:M\to M$. If $F^0(M)\neq M$, then we refer to $\vph$ as the Frobenius endomorphism of $M[{1\over p}]$.
\smallskip
The notations $\Mm\Mf(R)$ and $\Mm\Mf(R^\wedge)$ always denote the same category. Accordingly, if $M$ is a projective $R$-module, we identify this $p$-divisible object ${\got C}$ with the one obtained from it by the extension of scalars from $R$ to $R^\wedge$; in other words, we identify ${\got C}$ with the $p$-divisible object ${\got C}_{R^\wedge}:=(M\otimes_R R^\wedge,(F^i(M)\otimes_R R^\wedge)_{i\in\ZZ},\vph\otimes 1)$, the Frobenius lift of $R^\wedge$ being the $p$-adic completion of $\Phi_R$. However, as we do not want to lose the information on ${\got C}$ of being ``algebraizable", we always prefer to work with ${\got C}$ (resp. with ${\got C}/p^n{\got C}$) instead of ${\got C}_{R^\wedge}$ (resp. of ${\got C}_{R^\wedge}/p^n{\got C}_{R^\wedge}$).
\smallskip
If ${\got C}_1=(M_1,(F^i(M_1))_{i\in\ZZ},\vph^1)$ is another such $p$-divisible object of $\Mm\Mf(R)$ with $M_1$ an $R$-module (resp. an $R^\wedge$-module), then by a morphism from ${\got C}$ to ${\got C}_1$, we mean an $R$-linear map $q:M\to M_1\otimes_R R^\wedge$ (resp., when $M$ and $M_1$ are $R^\wedge$-modules, we mean an $R^\wedge$-linear map $q:M\to M_1$) respecting the filtrations and Frobenius endomorphisms. Notation: 
$${\got C}\to {\got C}_1.$$ 
The $\ZZ_p$-module of morphisms from ${\got C}$ to ${\got C}_1$ is denoted by ${\rm Hom}({\got C},{\got C}_1)$.
We get the category $p-\Mm\Mf(R)$ of $p$-divisible objects of $\Mm\Mf(R)$; it is a $\ZZ_p$-linear category. 
\smallskip
Warning: a $p$-divisible object of $\Mm\Mf(R)$ is not a particular type of object of $\Mm\Mf(R)$ (of course, it is an object of $p-\Mm\Mf(R)$); it is (i.e. it can be identified with) a suitable $\NN$-projective system of objects of $\Mm\Mf(R)$. In other words, we can identify (cf. definitions) ${\got C}$ with the naturally obtained $\NN$-projective system defined by ${\got C}/p^n{\got C}$, $n\in\NN$. So the category $p-\Mm\Mf(R)$ is equivalent to a full subcategory of the category $\NN-pro-\Mm\Mf(R)$ of $\NN$-projective systems of objects of $\Mm\Mf(R)$; we rarely use this equivalence, as we think the notations involving projective $R$-modules (resp. $R^\wedge$-modules) are much simpler than the ones involving $\NN$-projective systems. 
\smallskip
If moreover we fix a family of tensors $(t_{\al})_{\al\in\Mj}$ of the $F^0$-filtration of $\Mt(M[{1\over p}])$ fixed by $\vph$, then we speak about a $p$-divisible object with tensors of $\Mm\Mf(R)$. Notation: a quadruple $(M,(F^i(M))_{i\in\ZZ},\vph,(t_{\al})_{\al\in\Mj})$ or a pair $({\got C},(t_{\al})_{\al\in\Mj})$.
\smallskip
We also speak about $p$-divisible objects (with tensors) of $\Mm\Mf_{[a,b]}(R)$, of $\Mm\Mf_{[a,b]}^\nabla(R)$ or of $\Mm\Mf^\nabla(R)$, or about $p$-divisible objects (with tensors) of $\Mm\Mf_{[a,b]}(X)$, etc., with $a,b\in\ZZ$, $b\ge a$, and with $X$ a regular, formally smooth $W(k)$-scheme (resp. a regular, formally smooth $p$-adic formal scheme over $W(k)$). Warning: in this last case (pertaining to $X$) the tensors are global sections of sheaves of the form $\Mt(\Mf)\otimes_{\ZZ_{(p)}} \QQ$, with $\Mf$ a locally free $\Mo_{X^\wedge}$-sheaf (resp. $\Mo_X$-sheaf) of locally finite rank. Without specific references, we always assume that either $X$ itself or its $p$-adic completion is equipped with a Frobenius lift or we are in a situation where the standard gluing process of [Fa1, 2.3] applies, i.e. we are in the context of $\Mm\Mf_{[a,b]}^\nabla(X)$, with $b-a\le p-1$ (however, see 2.2.4 C and D for natural extensions). Warning: as the connections involving objects of $\Mm\Mf^\nabla(X)$ are $p$-nilpotent (see [Fa1, p. 34]), the restriction $p>2$ in (the proof of) [Fa1, 2.3] is not needed; so in what follows, we refer to loc. cit. and its proof without extra comment even for $p=2$. We also speak about $\ZZ_p$-linear categories defined by such $p$-divisible objects (resp. by $\NN$-projective systems of objects): $p-\Mm\Mf_{[a,b]}(R)$ (resp. $\NN-pro-\Mm\Mf_{[a,b]}(R)$), etc. When we speak about a $p$-divisible object with tensors of $\Mm\Mf^\nabla(R)$, etc., we require the tensors to be (as well) horizontal, i.e. to be annihilated by the (natural extension to $\Mt(M)$ of the) connection involved. We speak as well about truncations of  $p$-divisible objects of $\Mm\Mf_{[a,b]}(R)$, or of $\Mm\Mf^\nabla(R)$, etc.
\smallskip
We recall from 2.1, that whenever we deal with some (object or $p$-divisible object of a) Fontaine category involving a regular, formally smooth ($p$-adic formal) scheme $Y$ over $W(k)$, we always assume $\Om_{Y_k/k}$ is a locally free $\Mo_{Y_k}$-sheaf of locally finite rank. Warning: as we allowed $R$ (or $X$) to be regular, formally smooth (and not only smooth) over $W(k)$, some precautions are in order. In the case of a $p$-divisible object of $\Mm\Mf^\nabla(R)$, we work with the $p$-adic completion $\Om_{R/W(k)}^\wedge$ of the $R$-module of relative differentials $\Om_{R/W(k)}$ of $R$ over $W(k)$; the same applies in the context of $X$. For instance, if $R=W(k)[[z_1,...,z_m]]$, $m\in\NN$, and if $M$ is a free $R$-module of finite rank, then all connections 
$$
\nabla:M\to M\otimes_{W(k)} \Om_{R/W(k)}
$$ 
to be considered, factor through $M\otimes_R \bar\Om_{R/W(k)}$, where $\bar\Om_{R/W(k)}$ is the free $R$-submodule of $\Om_{R/W(k)}$ generated by $dz_1$,..., $dz_m$; we have a natural identification $\bar\Om_{R/W(k)}=\Om_{R/W(k)}^\wedge$. 
\smallskip
A $p$-divisible object of $\Mm\Mf^\nabla(R)$ is denoted as a quadruple $(M,(F^i(M))_{i\in\ZZ},\vph,\nabla)$ or as a pair $({\got C},\nabla)$, where ${\got C}$ is as above; warning: due to the conventions of the previous paragraph, we always consider $M$ to be an $R^\wedge$-module.
\smallskip
We use freely tensor products of ($p$-divisible) objects of some Fontaine category, whenever they are defined (when we are dealing with a Fontaine category of objects involving filtrations in a fixed range, like $\Mm\Mf_{[a,b]}(R)$, then such tensor products might be defined just in a ``larger" Fontaine category, like $\Mm\Mf(R)$). For instance, if above $M$ and $M_1$ are $R^\wedge$-modules, ${\got C}\otimes_{R^\wedge} {\got C}_1$ is the $p$-divisible object of $\Mm\Mf(R)$ whose underlying filtered $R^\wedge$-module is $M\otimes_{R^\wedge} M_1$ endowed with the tensor product filtration, the Frobenius endomorphism of $M[{1\over p}]\otimes_{R^\wedge[{1\over p}]} M_1[{1\over p}]=M\otimes_{R^\wedge} M_1[{1\over p}]$ being $\vph\otimes\vph^1$.
\smallskip
We speak about a principal quasi-polarization of a $p$-divisible object of $\Mm\Mf_{[0,1]}(R)$ or of $\Mm\Mf_{[0,1]}^\nabla(R)$. For instance, if ${\got C}=(N,F^1,\vph)$ is a $p$-divisible object of $\Mm\Mf_{[0,1]}(R)$ with $N$ an $R^\wedge$-module, the dual ${\got C}^*$ of ${\got C}$ is a $p$-divisible object of $\Mm\Mf_{[-1,0]}(R)$ and so ${\got C}^*(1)$ is a $p$-divisible object of $\Mm\Mf_{[0,1]}(R)$; so by a principal quasi-polarization of ${\got C}$, we mean an isomorphism 
$$i_{\got C}:{\got C}\tilde\to{\got C}^*(1)
$$ 
producing a perfect alternating form
$$
p_{\got C}:{\got C}\otimes_{R^\wedge} {\got C}\to R^\wedge(1)
$$
($p_{\got C}$ is a morphism in the category $p-\Mm\Mf_{[0,2]}(R)$). Here, for $m\in\ZZ$, 
$$
R^\wedge(m):=R^\wedge\otimes_{W(k)} W(k)(m)
$$ 
is the $p$-divisible object of $\Mm\Mf_{[m,m]}(R)$ defined by the $R^\wedge$-module $R^\wedge$ and the $\Phi_R$-linear map of $R^\wedge[{1\over p}]$ which takes $1$ into $p^m$. Usually we identify canonically $i_{\got C}$ with $p_{\got C}$; we refer to the pair $({\got C},p_{\got C})$ or to the quadruple $(N,F^1,\vph,p_N)$ (with $p_N:N\otimes_{R^\wedge} N\to R^{\wedge}$ defining $p_{\got C}$) as a principally quasi-polarized $p$-divisible object of $\Mm\Mf_{[0,1]}(R)$. We often write $p_N:N\otimes_{R^\wedge} N\to R^{\wedge}(1)$. The truncation mod $p^n$ of $R^\wedge(m)$ is also denoted as $R/p^nR(m)$. 
\smallskip
Similarly, we speak about principally quasi-polarized objects of $\Mm\Mf_{[0,1]}(R)$ or of $\Mm\Mf_{[0,1]}^\nabla(R)$; most commonly we get such principally quasi-polarized objects by taking truncations (modulo some integral positive power of $p$) of principally quasi-polarized $p$-divisible objects (of $\Mm\Mf_{[0,1]}(R)$ or of $\Mm\Mf_{[0,1]}^\nabla(R)$). The last two paragraphs apply as well in the context of $X$.
\medskip
{\bf d)} A $p$-divisible object ${\got C}=(M,(F^i(M))_{i\in S(a,b)},\vph)$ of $\Mm\Mf_{[a,b]}(W(k))$ is said to be cyclic diagonalizable if $M\neq\{0\}$ and there is a $W(k)$-basis $\Mb=\{e_i|i\in A\}$ of $M$, with $A:=S(1,\dim_{W(k)}(M))$, such that the following two conditions hold:
\medskip
\item{i)} $\forall i\in S(a,b)$, a subset of $\Mb$ is a $W(k)$-basis of $F^i(M)$;
\smallskip
\item{ii)} there is a permutation $\pi$ of $A$ with the property that $\vph(e_i)=p^{n(i)}e_{\pi(i)}$, $\forall i\in A$, where $n(i)\in S(a,b)$ is uniquely determined by requesting $e_i\in F^{n(i)}(M)\setminus F^{n(i)+1}(M)$.
\medskip
When $\pi=1_A$ we drop the word cyclic. If $\pi$ is an $\abs{A}$-cycle, then we often replace the word cyclic by circular. 
\smallskip
Similarly, we speak about the truncation mod $p^n$ of ${\got C}$ of being cyclic (or circular) diagonalizable; keeping i) we just need to replace in ii) the equality $\vph(e_i)=p^{n(i)}e_{\pi(i)}$ by the requirement that $\vph_{n(i)}(e_i)$ mod $p^n$ is $e_{\pi(i)}$ mod $p^n$; here $\vph_i$'s are obtained from $\vph$ as in c).
\medskip
{\bf e)} Let $O$ be a DVR and let $u$ be a uniformizer of it. If $S$ is an $O$-scheme, then an $\Mo_S$-sheaf is said to have the DC property (here DC stands for direct cyclic) if it is a finite direct sum of $\Mo_S$-sheaves of the form $\Mo_S/u^n\Mo_S$, with $n\in\NN\cup\{0\}$. Similarly, if $R$ is an $O$-algebra, then an $R$-module is said to have the DC property if as an $\Mo_{{\rm Spec}(R)}$-sheaf it has the DC property. If we allow $n$ to be $\infty$ as well (with $\Mo_S/u^{\infty}\Mo_S:=\Mo_S$), then we speak about the EDC property, with $E$ standing for extended.
\medskip
{\bf f)} We refer to the morphism ${\got C}\to {\got C}_1$ of c). We denote it by ${\got q}$. We assume $M$ and $M_1$ are $R^\wedge$-modules. ${\got q}$ is called an epimorphism or a monomorphism if $q$ is so. ${\got q}$ is called an isogeny if $q$ becomes an isomorphism by inverting $p$ and moreover there is $m\in\NN$ such that $q(M)$ contains $p^mM_1$. If $q$ after inverting $p$ is an isomorphism, then we define ${\rm Coker}({\got q})$ as follows. We can assume $X_k$ is connected. Let $n\in\NN$. $M_2:=M_1/q(M)$ is $p$-adically complete. Moreover, $M_2/p^nM_2$ locally has the DC property (being the underlying $R^\wedge$-module of ${\rm Coker}({\got C}/p^n{\got C}\to {\got C}_1/p^n{\got C}_1)$, cf. [Fa1, p. 31-33]). So, as $M_2[{1\over p}]=\{0\}$ and as $X_k$ is connected, using a Teichm\"uller lift of ${\rm Spec}(R^\wedge)$ with values in the Witt ring of a perfect field, we get that we can assume $p^nM_2=\{0\}$. So we take:
$${\rm Coker}({\got q}):={\rm Coker}({\got C}/p^n{\got C}\to {\got C}_1/p^n{\got C}_1);$$
it is well defined independently of the choice of $n$ (subject to $p^nM_2=\{0\}$). The same applies to all other Fontaine categories of $p$-divisible objects. 
\smallskip
Similarly, if ${\got C}$ (resp. ${\got C}_1$) is a $p$-divisible object (resp. an object annihilated by $p^n$) of the same Fontaine category $FC$, then the torsion $\ZZ_p$-module 
$${\rm Hom}({\got C},{\got C}_1):={\rm Hom}_{FC}({\got C}/p^n{\got C},{\got C}_1)$$ 
is well defined independently of the choice of $n$ (subject to the mentioned property) and any element ${\got q}$ of it is referred as a morphism ${\got C}\to {\got C}_1$. If at the level of underlying modules ${\got q}$ is an epimorphism, we say ${\got q}$ is an epimorphism and we refer to ${\got C}$ as a lift of ${\got C}_1$; if moreover the resulting epimorphism ${\got C}/p{\got C}\to {\got C}_1$ is an isomorphism, we refer to ${\got C}$ as a strict lift of ${\got C}_1$.
\medskip
{\bf 2.2.1.0. The crystalline, contravariant Dieudonn\'e functor in the language of Fontaine categories.} The initial motivation for using ($p$-divisible objects of) Fontaine categories, stems from the existence of the crystalline, contravariant Dieudonn\'e functor $\DD$. Let $X$, $R$ and $\Phi_R:R^\wedge\to R^\wedge$ be as in 2.2.1 c).
\smallskip
We have a contravariant, $\ZZ_p$-linear functor 
$$p-DG(X)\to p-\Mm\Mf_{[0,1]}^\nabla(X).$$ 
Its existence is a particular case of Grothendieck--Messing theory (see [Me, ch. 4-5]), via a natural limit process; in other words, the functor associates to a $p$-divisible group $D$ over $X$ the $p$-divisible object of $\Mm\Mf_{[0,1]}^\nabla(X)$ defined naturally by the $\NN$-projective system of objects of $\Mm\Mf_{[0,1]}^\nabla(X)$, whose object corresponding to $n\in\NN$ is annihilated by $p^n$ and is naturally defined (via reduction mod $p^n$) by the dual of the Lie algebra of the universal vector extension of the pull back (warning!) of $D[p^{n+1}]$ to $X\times_{W(k)} W_{n+1}(k)$. Such a functor is also constructed in [BBM, ch. 3]. In [BM, ch. 3] it is shown the compatibility of the two constructions (warning: [BM, 3.2.11] applies here even for $p=2$ due to the fact that we are not dealing with crystals on ${\rm CRIS}(X_{W_{n+1}(k)}/{\rm Spec}(W(k)))$ but just with locally free $\Mo_{X_{W_{n+1}(k)}}$-sheaves obtained by evaluated them at $X_{W_{n+1}(k)}$).  Following [BBM, ch. 3], this functor is denoted by $\DD$. 
\smallskip
We recall briefly the construction of $\DD$, in the wider (in fact intermediary) context of $p-FF(X)$, i.e. of finite, flat, commutative group schemes of $p$-power order over $X$.
Let $G$ and $\tilde G$ be such group schemes. We first assume that $X={\rm Spec}(R^\wedge)$, that $G$ (resp. $\tilde G$) is a closed subgroup of an abelian scheme $A$ (resp. $\tilde A$) over $R^\wedge$ and that there is $m\in\NN$ such that both $G$ and $\tilde G$ are annihilated by $p^m$. Let $A_1:=A/G$. Let ${\got C}$ (resp. ${\got C}_1$) be the $p$-divisible object of $\Mm\Mf_{[0,1]}^\nabla(R)$ naturally defined by $A$ (resp. by $A_1$). We recall that the underlying $R^\wedge$-module of ${\got C}$ is $H^1_{\rm crys}(A/R^\wedge)$ and is naturally equipped with a $\Phi_R$-linear endomorphism, that the filtration of $H^1_{\rm crys}(A/R^\wedge)$ is the Hodge filtration, and that the connection on $H^1_{\rm crys}(A/R^\wedge)$ is the $p$-adic completion of the Gauss--Manin connection. Similarly, we get ${\got C}_1$. 
\smallskip
Now (cf. [Fo] for the particular case $R=W(k)$ and cf. [BBM, 3.1.2] for the general case; see also [BM, p. 190-1]) 
$$\DD(G)=(M,F^1,\Phi_0,\Phi_1,\nabla)$$ 
is defined as the cokernel (see 2.2.1 f)) of the isogeny ${\got C}_1\hookrightarrow {\got C}$ associated to the natural isogeny $A\twoheadrightarrow A_1$. $M$ is annihilated by $p^m$. 
\smallskip
Similarly we define $\DD(\tilde G)$. If $f_G:G\to \tilde G$ is a homomorphism, then the homomorphism $(1_G,f_G):G\to G\times\tilde G$ is a closed emebedding. So $\DD(f_G):\DD(\tilde G)\to\DD(G)$ is defined naturally via the epimorphism $$\DD(G)\oplus \DD(\tilde G)=\DD(G\times\tilde G)\twoheadrightarrow\DD(G)$$
 associated naturally to $(1_G,f_G)$; here $\DD(G\times\tilde G)$ is computed via the product embedding of $G\times\tilde G$ in $A\times\tilde A_1$. One checks that $\DD(G)$ and $\DD(f_G)$ are well defined (i.e. do not depend on the choices made) in the standard diagonal way. For instance, if $G$ is a closed subgroup of another abelian scheme $B$ over $R^\wedge$, then by embedding $G$ diagonally in $A\times_{R^\wedge} B$, we get (via standard arguments involving the snake lemma in the context of any one of the two projections of $A\times_{R^\wedge} B$ onto its factors) that $\DD(G)$ defined via $A\times_{R^\wedge} B$ is canonically isomorphic to $\DD(G)$ defined via $A$ or $B$.
\smallskip
To define $\DD(G)$ in general, we use Raynaud's theorem of [BBM, 3.1.1]: locally in the Zariski topology of $X$, $G$ is a subgroup of an abelian scheme and is annihilated by $p^m$ for some $m\in\NN$; so $\DD(G)$ is defined as above, via standard arguments of gluing (based on the previous paragraph). 
\smallskip
Warning: in [BBM], [BM] or [dJ1-2] the language of Fontaine categories is not used. Referring to the case $X={\rm Spec}(R^\wedge)$, the mentioned places (with [BBM, ch. 3] as the main reference) are using instead of the $\Phi_R$-linear map $\Phi_1:F^1\to M$, a Verschiebung $R^\wedge$-linear map $V:M\to M\otimes_{R^\wedge}\, _{\Phi_R}R^\wedge$; identifying $\Phi_1$ with an $R^\wedge$-linear map $F^1\otimes_R\, _{\Phi_R}R^\wedge\to M$, we have 
$$V\circ\Phi_1(x)=x,\leqno (VPHIONE)$$
$\forall x\in F^1\otimes_{R^\wedge}\, _{\Phi_R}R^\wedge$. Similarly, identifying $\Phi_0$ with an $R^\wedge$-linear map $M\otimes_R\, _{\Phi_R}R^\wedge\to M$, we have
$$V\circ\Phi_0(x)=px,\leqno (VPHIZERO)$$
$\forall x\in M\otimes_{R^\wedge}\, _{\Phi_R}R^\wedge$. As $M$ is generated by the images of $\Phi_1$ and $\Phi_0$, we get:
\medskip
{\bf Fact 1.} {\it $V$ is uniquely determined by $\Phi_0$ and $\Phi_1$.}
\medskip
So we feel it is appropriate to denote $(M,F^1,\Phi_0,\Phi_1,\nabla)$ by $\DD(G)$; so in what follows, (the functor) $\DD$ is used for $p-FF(X)$ as well as for $p-DG(X)$.
\smallskip
If now $D$ is a $p$-divisible group over $X$, then $\DD(D)$ is the $\NN$-projective system defined naturally by $\DD(D[p^n])$, $n\in\NN$: one checks that the underlying $\Mo_X/p^n\Mo_X$-sheaf of $\DD(D[p^n])$ is projective of rank equal to the rank of $D$ and $\DD(i_n)$ is an epimorphism, $\forall n\in\NN$; here $i_n:D[p^n]\hookrightarrow D[p^{n+1}]$ is the natural inclusion. Argument: based on [Fa1, 2.1 ii)], using Teichm\"uller lifts (of open, affine subschemes of $X$ endowed with a Frobenius lift) it is enough to consider the case when $X={\rm Spec}(W(k))$. But this case is a consequence of the classical Dieudonn\'e theory. We get the first sentence of:
\medskip
{\bf Fact 2.} {\it $\DD(D)/p^n\DD(D)$ is naturally identifiable with $\DD(D[p^n])$. It depends only on the pull back of $D[p^{n+1}]$ to $X_{W_{n+1}(k)}$.}
\medskip
The second sentence follows from the above part referring to [Me] or to the connection between $V$ and $\Phi_1$. 
\smallskip
The language to be used in what follows: to any $p$-divisible group $D$ over $X$ (resp. over ${\rm Spec}(R)$ or ${\rm Spec}(R^\wedge$)) it corresponds (or it is associated) a $p$-divisible object $\DD(D)$ of $\Mm\Mf_{[0,1]}^\nabla(X)$ (resp. of $\Mm\Mf_{[0,1]}^\nabla(R)$); similarly, we say: to a finite, flat, commutative group scheme $G$ of $p$-power order over $X$ (resp. over ${\rm Spec}(R)$ or ${\rm Spec}(R^\wedge$)) it corresponds (or it is associated) an object $\DD(G)$ of $\Mm\Mf_{[0,1]}^\nabla(X)$ (resp. of $\Mm\Mf_{[0,1]}^\nabla(R)$). 
\smallskip
Whenever $D$ is uniquely determined by $\DD(D)$, we also say $D$ is associated to $\DD(D)$; [BM, \S 4] and [Me, ch. 4-5] can be combined to give plenty of situations, with $p>2$, in which $D$ is uniquely determined by $\DD(D)$.
\medskip
{\bf 2.2.1.1. Remarks.} {\bf 1)} One might wonder if the terminology $p$-divisible objects is the right one. Our justification is: the below philosophies of 2) and of 3.6.18.5.1. We do think that:
\medskip
-- the terminology ``Barsotti--Tate objects" would have complicated very much the terminology (for instance, later on --see 2.2.8 3)-- we introduce Shimura $p$-divisible objects);
\smallskip
-- the terminology ``formal objects" would have forced us to always be in a formal context and moreover it would have not allowed us to speak about $p$-adic formal schemes as a way of emphasizing the type of formal schemes we are mostly dealing with;
\smallskip
-- the terminology ``torsion free objects" offers some advantages (i.e. in some sense it is more accurate) but using it the results below would have been significantly less uniformly stated;
\smallskip
-- the terminology of [FL, \S 1] used for the case of $W(k)$ is neither practical nor accurate for the general cases of $R$'s and $X$'s as in 2.2.1 c) (that is why [Fa1] speaks just about objects).   
\medskip
{\bf 2)} If $p>2$ and $X$ is the $p$-adic completion of a smooth $W(k)$-scheme, then  [Fa1, 7.1] can be restated as follows. The contravariant, $\ZZ_p$-linear functor
$$\DD:p-FF(X)\to\Mm\Mf_{[0,1]}^\nabla(X)$$ 
is an antiequivalence of categories. This implies: $\DD:p-DG(X)\to p-\Mm\Mf_{[0,1]}^\nabla(X)$ is an antiequivalence of categories. So occasionally, we refer to $\DD^{-1}$.
\smallskip
{\bf 3)} In general, the categories of objects or $p$-divisible objects introduced in 2.2.1 c) do depend on the choice of Frobenius lifts. However, the category $\Mm\Mf_{[0,p-1]}^\nabla(X)$ does not depend (up to canonical isomorphisms) on such choices, cf. [Fa1, 2.3]. One of the main advantages of working with $p$-divisible objects instead of objects is:
\medskip
{\bf Fact.} {\it $p-\Mm\Mf_{[a,b]}^\nabla(X)$ (as well as $p-\Mm\Mf^\nabla(X)$) does not depend (up to canonical isomorphism) on the choice of a Frobenius lift of $X$.}
\medskip
This is a consequence of the fact that, when we have no $p$-torsion involved, canonical isomorphisms can be constructed as in loc. cit. Argument: we can assume $a=0$ and so, to check that some gluing maps (constructed as in [De3, 1.1.3.4] using Taylor series) do define isomorphisms, we can work locally in the faithfully flat topology; for the case of $X={\rm Spec}(W(k)[[x_1,...,x_m]])$ we refer to loc. cit.  
\smallskip
We always identify $p-\Mm\Mf_{[0,b]}^\nabla(X)$ with a full subcategory of the category of filtered $F$-crystals on $X_k$ (to be compared with the case of $\Mm\Mf_{[0,b]}^\nabla(X)$ of 2.1); here $b\in\NN\cup\{0\}$.
\smallskip
{\bf 4)} Any cyclic diagonalizable $p$-divisible object is a direct sum of circular ones. Moreover, any circular diagonalizable $p$-divisible object has a unique slope.
\smallskip
{\bf 5)} $\Mm\Mf(X)$ and $\Mm\Mf(X^\wedge)$ denote the same category of objects; the same applies to all other Fontaine categories of objects or of $p$-divisible objects. Moreover, we always identify $\Mm\Mf_{[a,b]}(X)$ (resp. $\Mm\Mf_{[a,b]}^\nabla(X)$) with a full subcategory of $\Mm\Mf_{[a^\prime,b^\prime]}(X)$ (resp. of $\Mm\Mf_{[a^\prime,b^\prime]}^\nabla(X)$) and of $\Mm\Mf(X)$ (resp. of $\Mm\Mf^\nabla(X)$), where $a^\prime\le a\le b\le b^\prime$ are integers. Also, we always identify $\Mm\Mf(R)$ with $\Mm\Mf({\rm Spec}(R))$, etc. Warning: the natural functor
$$FORG:\Mm\Mf^\nabla_{[a,b]}(X)\to\Mm\Mf_{[a,b]}(X)$$ 
which ``forgets connections" is not full in general. 
\smallskip
{\bf 6)} All Fontaine categories of 2.2.1 c) involving objects are abelian, cf. [Fa1, Cor. of p. 33] (it extends automatically to the context involving connections). This is not true for the context of $p$-divisible objects (the problem being with cokernels). However, it is trivial to see that all Fontaine categories of 2.2.1 c) involving $p$-divisible objects have kernels, images and arbitrary finite fibre products. Moreover, if $a$, $b$, $R$ and $\Phi_R:R^\wedge\to R^\wedge$ are as in 2.2.1 c), we have:
\medskip
{\bf Fact.} {\it Locally in the Zariski topology of $R$, any object of $\Mm\Mf_{[a,b]}(R)$ is annihilated by $p^n$ for some $n\in\NN$ and is the cokernel of an isogeny between two $p$-divisible objects of $\Mm\Mf_{[a,b]}(R)$.}  
\medskip
{\bf Proof:} Let ${\got C}_1:=(M,(F^i(M))_{i\in S(a,b)},(\vph_i)_{i\in S(a,b)})$ be an object of $\Mm\Mf_{[a,b]}(R)$. We choose a direct summand $\tilde F^i(M)$ of $F^{i+1}(M)$ in $F^i(M)$, cf. [Fa1, 2.1 i)], $i\in S(a,b)$. Localizing, we can assume $M=\oplus_{j\in S(1,m)} R/p^{n_j}R$, with $m\in\NN\cup\{0\}$ and all $n_j$'s in $\NN$, cf. [Fa1, 2.1 ii)]. So we can take $n:={\rm max}\{n_j|j\in S(1,m)\}$. 
\smallskip
Localizing further on we can assume all $\tilde F^i(M)$'s have the DC property. The argument goes as follows. $\tilde F^i(M)/p\tilde F^i(M)$'s are projective $R/pR$-modules and so localizing we can assume they are free. We work with a fixed $i\in S(a,b)$. By induction on $n$, we can assume $p\tilde F^i(M)$ has the DC property. Let $\{e_i^1,...,e_i^{m_i}\}$, with $m_i\in\NN\cup\{0\}$, be a set of elements of $\tilde F^i(M)$ such that $p\tilde F^i(M)$ is a direct sum of the cyclic $R$-modules generated by $pe_i^j$, $j\in S(1,m_i)$, and the annihilator of each such $pe_i^j$ is an ideal of $R$ generated by a positive, integral power of $p$. This implies that the $R/pR$-submodule of $\tilde F^i(M)/p\tilde F^i(M)$ generated by the images of $e_i^j$'s is free of rank $m_i$; it is easy to see that is as well a direct summand. So let $\{e_i^{m_i+1},...,e_i^{q_i}\}$ be a set of elements of $\tilde F^i(M)$ such that the images of $e_i^j$, $j\in S(1,q_i)$, in $\tilde F^i(M)/p\tilde F^i(M)$ are forming an $R/pR$-basis of this free $R/pR$-module. By replacing each $e_i^j$, with $j\in S(m_i+1,q_i)$, by a linear combination of $e_i^s$'s, with $s\in S(1,m_i)$, we can assume $pe_i^j=0$, $\forall j\in S(m_i+1,q_i)$. Let $s(i,j)\in\NN$ be such that $p^{s(i,j)}R$ is the annihilator of $e_i^j$, $j\in S(1,q_i)$. We get: $\tilde F^i(M)$ is (isomorphic to) the direct sum $\oplus_{j\in S(1,q_i)} (R/p^{s(i,j)}R)e_i^j$. 
\smallskip
We now consider the free $R^\wedge$-module $M_R$ having as an $R^\wedge$-basis a set formed by elements $f_i^j$, $i\in S(a,b)$ and $j\in S(1,q_i)$. We consider the $R^\wedge$-epimorphism
$$q:M_R\twoheadrightarrow M$$
that takes each $f_i^j$ into $e_i^j$.
For $i\in S(a,b)$, let $F^i(M_R)$ be the $R^\wedge$-submodule of $M_R$ generated by those $f_{i^\prime}^j$, with $i^\prime\in S(i,b)$ and $j\in S(1,q_{i^\prime})$. $F^i(M_R)\cap {\rm Ker}(q)$ is a free $R^\wedge$-module having $\{p^{s(i^\prime,j)}f_{i^\prime}^j|i^\prime\in S(i,b)\, {\rm and}\, j\in S(1,q_{i^\prime})\}$ as an $R^\wedge$-basis. We consider a $\Phi_R$-linear endomorphism $\vph$ of $M_R[{1\over p}]$ that takes $f_i^j$ into $p^ig_i^j$, where $g_i^j\in M_R$ is an arbitrary element with the property that $q(g_i^j)=\vph_i(e_i^j)$. The elements $g_i^j$, $i\in S(a,b)$ and $j\in S(1,q_i)$, generate $M_R$, as this is so mod $p$ ($q$ mod $p$ being an isomorphism). So the triple $(M_R,(F^i(M_R))_{i\in S(a,b)},\vph)$ is a $p$-divisible object ${\got C}$ of $\Mm\Mf_{[a,b]}(R)$. Moreover, from the very construction $q$ defines a morphism ${\got C}\to {\got C}_1$ and so (cf. also on the above part on ${\rm Ker}(q)$) ${\rm Ker}(q)$ is the underlying $R^\wedge$-module of a $p$-divisible subobject ${\got C}_2$ of ${\got C}$. Obviously ${\got C}_1$ is the cokernel of the natural inclusion ${\got C}_2\hookrightarrow {\got C}$. This proves the Fact.
\medskip
From the part of the proof referring to $\tilde F^i(M)$ having the DC property, we get:
\medskip
{\bf Corollary.} {\it An $R$-module has the DC property iff $R/pR$ and $pR$ have it.}
\medskip
{\bf Exercise.} Show that the kernel of a morphism from a $p$-divisible object into an object of the same Fontaine category FC is well defined and is a $p$-divisible object of FC. Hint: use the fact that FC is an abelian category and the above proof.   
\medskip
{\bf 7)} If $FC$ is a Fontaine category of objects, we denote by $FC[p^n]$ (resp. by $p-FC$) its full subcategory formed by objects annihilated by $p^n$ (resp. the category of $p$-divisible objects of $FC$). We refer to it as a Fontaine category of objects annihilated by $p^n$ or as the truncation mod $p^n$ of $FC$ (resp. as the Fontaine category of $p$-divisible objects of $FC$). The operation of taking truncations mod $p^n$ can be viewed as a functor 
$$TR_n(FC):p-FC\to FC[p^n].$$ 
\indent
{\bf Corollary.}  {\it If $FC=\Mm\Mf(R)$ with $R$ local, then $TR_1(FC)$ is essentially surjective.} 
\medskip
{\bf Proof:} The proof of 6) implies $TR_1(FC)$ is surjective on objects; but its arguments can be used as well for morphisms. This ends the proof. 
\medskip
We assume $X_k$ has a finite number of connected components. Let $FC$ be any one of the following four categories: $\Mm\Mf(X)$, $\Mm\Mf^\nabla(X)$, $\Mm\Mf_{[a,b]}(X)$ or $\Mm\Mf_{[a,b]}^\nabla(X)$. Let ${\got q}:{\got C}\hookrightarrow {\got C}_1$ be a monomorphism of $p-FC$. We have:
\medskip
{\bf Fact.} {\it There is a canonical epimorphism ${\got e}:{\got C}_1\twoheadrightarrow {\got C}_2$ of $p-FC$ such that ${\got e}\circ {\got q}=0$ and the factorization ${\got C}\hookrightarrow {\rm Ker}({\got e})$ is an isogeny.}
\medskip
{\bf Proof:} Based on the canonical aspect of the Fact, we can work locally and so we can assume $X={\rm Spec}(R)$, with $R$ local and as in 2.2.1 c). The proof in the context involving connections is entirely the same. So not to introduce extra notations we assume we are in the context of ${\got q}$ of 2.2.1 c), with $M$ and $M_1$ as $R^\wedge$-modules  and with $q$ a monomorphism. Let $M_1^\prime:=\{x\in M_1|p^nx\in q(M) \, {\rm for}\, {\rm some}\, n\in\NN\}$. The quotient module 
$$M_2:=M_1/M_1^\prime$$ 
is a locally free $R^\wedge$-module. The argument for this goes as follows. $M_1/q(M)$ is $p$-adically complete and so is the $\NN$-projective limit of $M_1/p^nM_1+q(M)$'s, $n\in\NN$. $M_1/p^nM_1+q(M)$ is the underlying $R$-module of the kernel of the truncation mod $p^n$ of ${\got q}$ and so it has the DC property. We write it as the direct sum of a free $R/p^nR$-module $M_1(n)$ and of an $R/p^{n-1}R$-module $N_1(n)$. There is $n_0\in\NN$ such that $N_1(n)=N_1(n_0)$, $\forall n\ge n_0$: this can be seen using Teichm\"uller lifts (one for each connected component of ${\rm Spec}(R^\wedge)$; the Fact is trivial for Witt rings as we can argue at the level of lattices). So for $n\ge n_0$, we have $M_1(n+1)/p^nM_1(n+1)=M_1(n)$. So $M_2$ is the projective limit of $M_1(n)$'s and so is a free $R^\wedge$-module. We also get that we can identify $M_1/q(M)$ with $N_1(n_0)\oplus M_2$. The same arguments apply for filtrations: the image of $F^i(M_1)$ in $M_1/q(M)=N_1(n_0)\oplus M_2$ is a direct sum of a direct summand of $N_1(n_0)$ and of a direct summand $F^i(M_2)$ of $M_2$, $\forall i\in\ZZ$.  
\smallskip
For $i\in\ZZ$, the $\Phi_R$-linear endomorphism of the image of $F^i(M)$ in $M_1/q(M)$ naturally defined by $p^{-i}\vph^1$ takes $F^i(M_2)$ into $M_2$ (as $\vph^1$ becomes an isomorphism after inverting $p$). It is easy to see that this implies that $M_2$ together with $(F^i(M_2))_{i\in\ZZ}$ and the $\Phi_R$-linear endomorphism of $M_2[{1\over p}]$ defined naturally by $\vph^1$ is an object ${\got C}_2$ of $p-FC$ and moreover the $R^\wedge$-epimorphism $M_1\twoheadrightarrow M_2$ defines an epimorphism ${\got e}:{\got C}_1\twoheadrightarrow {\got C}_2$ of $p-FC$ having the required properties. This proves the Fact.
\medskip
From this Fact and 6) we deduce that the category 
$$isog-p-FC:=p-FC\otimes_{\ZZ_p} \QQ_p$$ 
of isogeny classes of $p$-divisible objects of $FC$ is an abelian, $\QQ_p$-linear category.   
\medskip
{\bf 2.2.1.2. A result of Wintenberger.} We recall a simplified form of the main result of [Wi]. If $(M,(F^i(M))_{i\in\ZZ},\vph,(t_{\alpha})_{\alpha\in\Mj})$ is a $p$-divisible object with tensors of $\Mm\Mf(W(k))$, then (cf. loc. cit. and the convention of 2.1 pertaining to it) there is a cocharacter $\mu:\GG_m\to GL(M)$ such that:
\medskip
-- the image of $\mu$ fixes $t_{\alpha}$, $\forall\alpha\in\Mj$;
\smallskip
-- it produces a direct sum decomposition $M=\oplus_{i\in\ZZ} \tilde F^i$, with $\beta\in\GG_m(W(k))$ acting through it on $\tilde F^i$ as the multiplication with $\beta^{-i}$, and with $\tilde F^i$ as a direct supplement of $F^{i+1}(M)$ in $F^i(M)$.
\medskip
{\bf 2.2.1.3. Pull backs.} We consider a $W(k)$-morphism $m:X_1\to X$, where each $X$ or $X_1$ is either a regular, formally smooth $W(k)$-scheme or is a regular, formally smooth, $p$-adic formal scheme over $W(k)$; warning: if $X$ is a formal scheme, then $X_1$ is also (automatically) a formal scheme. 
\smallskip
Let $a\in\ZZ$. For any object (resp. $p$-divisible object) ${\got C}$ of $\Mm\Mf_{[a,a+p-1]}^\nabla(X)$, we can define its pull back $m^*({\got C})$ (to $X_1$ or through $m$). The construction of the object (resp. $p$-divisible object) $m^*({\got C})$ of $\Mm\Mf^\nabla(X_1)$ can be obtained by entirely following the pattern of [Fa1, p. 34-35] (for the case of $p$-divisible objects it is somehow more convenient to follow the pattern of the top of [Fa2, p. 135] referring to $\tilde M$). 
\smallskip
To be shorter, we briefly recall the construction of $m^*({\got C})$ under the extra assumption (automatically satisfied if $m$ is a formally closed embedding or is formally smooth, or if it factors as a composite of such morphisms), that locally in the Zariski topology of $X_k$ and ${X_1}_k$, we can choose Frobenius lifts of the $p$-adic completions of $X$ and $X_1$ which are compatible with $m$ (these are the situations we need in what follows). So, due to [Fa1, 2.3], we can assume that $X={\rm Spec}(R)$ and $X_1={\rm Spec}(R_1)$ are affine schemes, $p$-adically complete and endowed with Frobenius lifts $\Phi_R$ and respectively $\Phi_{R_1}$ which are compatible with the natural $W(k)$-homomorphism (still denoted by $m$) $m:R\to R_1$. We also assume ${\got C}$ is a $p$-divisible object $(M,(F^i(M))_{i\in S(a,a+p-1)},\vph,\nabla)$: the case of an object is entirely the same. We now take 
$$m^*({\got C}):=(M\otimes_R R_1,(F^i(M)\otimes_R R_1)_{i\in S(a,a+p-1)},\vph\otimes 1,\nabla_1),$$
where $\nabla_1$ is the connection on $M\otimes_R R_1$ naturally induced by $\nabla$.
\smallskip
In case $X$ and $X_1$ are equipped with Frobenius lifts compatible with $m$, then we can define in the same manner the pull back $m^*({\got C})$ (to $X_1$ or through $m$) of any object or $p$-divisible of $\Mm\Mf(X)$ or of $\Mm\Mf^\nabla(X)$.
\smallskip
In case $X={\rm Spec}(R)$ and $X_1={\rm Spec}(R_1)$ are affine $W(k)$-schemes, we also refer to such an object or $p$-divisible object $m^*({\got C})$ as the extension of ${\got C}$ to $R_1$ or to $X_1$ (through $m$) or as the extension of ${\got C}$ through $m$. We also say: $m^*({\got C})$ is induced from ${\got C}$ via the $W(k)$-homomorphism $m:R\to R_1$. 
\medskip
{\bf 2.2.1.4. An application.} Let $O$ be a DVR of mixed characteristic $(0,p)$ having index of ramification $1$. Referring to the end of 2.1, we have:
\medskip
{\bf Proposition.} {\it Any regular, formally smooth $O$-scheme is healthy and $p$-healthy.} 
\medskip
{\bf Proof:} The case $p\ge 3$ is well known, see [Va2, 3.2.2 1) and 3.2.17]: loc. cit. proves even more for such primes; but the new proof below, under above restrictions on $O$, works for all primes. Following the Steps of [Va2, 3.2.17], only Step B of loc. cit. needs to be redone differently. We just need to show that any $p$-divisible group $D$ over $U$ extends to a $p$-divisible group over $Y:={\rm Spec}(W(\bar k)[[T]])$, where $U$ is the only open subscheme of $Y$ different from $Y$ and such that $(Y,U)$ is an extensible pair (see def. [Va2, 3.2.1 1)]). Let $\Mo$ be the local ring of $Y$ which is a DVR faithfully flat over $W(\bar k)$. Let $\Mo_1$ be a complete DVR which has an algebraically closed residue field and index of ramification $1$ and which is a faithfully flat $\Mo$-algebra. From [FC, 6.2, p. 181] we get that for any $n\in\NN$, $D[p^n]$ extends uniquely to a finite, flat group scheme $G_n$ over $Y$. We just need to show that for any $n,m\in\NN$ the sequence
$$0\to G_n\to G_{n+m}\to G_m\to 0\leqno (SEQ)$$
naturally defined by the standard short exact sequence
$$0\to D[p^n]\to D[p^{n+m}]\to D[p^m]\to 0,$$
is in fact a short exact sequence.
Let 
$$0\to\DD(G_m)\to\DD(G_{n+m})\to\DD(G_n)\to 0\leqno (DSEQ)$$
be the complex of $\Mm\Mf_{[0,1]}^\nabla(Y)$ induced by (SEQ) via 2.2.1.0. Let $M_m$, $M_{m+n}$ and $M_n$ be the underlying $W(\bar k)[[T]]$-modules of $\DD(G_m)$, $\DD(G_{n+m})$ and respectively of $\DD(G_n)$. From the proof of 2.2.1.1 6) we get: they have the DC property. We claim that the complex 
$$0\to M_m\to M_{n+m}\to M_n\to 0\leqno (MOD)$$
of $W(\bar k)[[T]]$-modules defining (DSEQ) is a short exact sequence. To prove this, due to [Fa1, 2.1], it is enough to check that the sequence (MODT) obtained from (MOD) by inverting $T$ is a short exact sequence. But (MODT) over $\Mo_1$ is nothing else but the short exact sequence defined by pulling back (SEQ) to $\Mo_1$ and applying $\DD$ (for instance, cf. the classical Dieudonn\'e theory over the residue field of $\Mo_1$). Another way to see this is to use the arguments of 2.2.1.0 on the independence of $\DD(G)$ on the embedding of $G$ in an abelian scheme; over $\Mo_1$ we can use as well embeddings into $p$-divisible groups over $\Mo_1$ (cf. also [BM, p. 190-1]). 
\smallskip
We conclude (DSEQ) is a short exact sequence. Using the DC property and simple arguments involving lengths of artinian $W(\bar k)$-modules, we get that (DSEQ) mod $T$ is as well a short exact sequence. So from the classical Dieudonn\'e theory we get that:
\medskip
-- if $p>2$, (SEQ) mod $T$ is a short exact sequence;
\smallskip
-- if $p=2$, (SEQ) mod $(2,T)$ is a sort exact sequence.
\medskip
Regardless of the parity of $p$, using Nakayama's lemma we get that $G_n$ is a closed subgroup of $G_{n+m}$. This implies: $G_m=G_{n+m}/G_n$. So (SEQ) is a short exact sequence. This ends the proof. 
\medskip
{\bf 2.2.1.4.1. Corollary.} {\it Let $(Y,U)$ be an extensible pair, with $Y$ a regular, formally smooth $O$-scheme of dimension two. Then any short exact sequence of finite, flat, commutative group schemes of $p$-power order over $U$ extends to a short exact sequence of finite, flat group schemes over $Y$.}
\medskip
{\bf Proof:} We can assume $Y$ is as in the proof of 2.2.1.4 and so the Corollary follows from this proof. 
\medskip
{\bf 2.2.1.4.2. Corollary.} {\it Any regular, formally smooth $\ZZ$-scheme is healthy.}
\medskip
The case $p=2$ of 2.1.4.1-2 answers a question of P. Deligne.
\medskip
{\bf 2.2.1.4.3. Remark.} One might wonder if in the proof of 2.2.1.4, the condition $O$ is of index of ramification $1$ is needed or not. A great part of the proof of 2.2.1.4 can be adapted (even for $p=2$) for an arbitrary index of ramification: however, (presently) there is one obstacle. We sketch what this ``obstacle" is (see \S 6 for precise definitions and more details). We assume $k=\bar k$. If $V$ is a finite, flat, totally ramified DVR extension of $V_0:=W(k)$, then following the pattern of [Fa2, \S 3] we define a category $p-\widetilde{\Mm\Mf}_{[0,1]}^\nabla(R_V[[T]])$ of $p$-divisible objects, whose underlying modules are free $R_V[[T]]$-modules of finite ranks; here $R_V$ is a $W(k)$-algebra obtained from $V$ in the same way as $Re$ was constructed in [Va2, 5.2.1], starting from the choice of a uniformizer $\pi_V$ of $V$. Here the tilde on the top of $\Mm\Mf$ refers to the fact that $F^1$-filtrations of underlying $R_V[[T]]$-modules are not defined by direct summands. There is a logical passage from $p$-divisible objects to objects (so we get the category $\widetilde{\Mm\Mf}_{[0,1]}^\nabla(R_V[[T]])$ of objects) and a similar functor $\DD$ as in 2.2.1.0. So the ``obstacle" shows up as a question: 
\medskip
{\bf Q.} {\it For which morphisms of $\widetilde{\Mm\Mf}_{[0,1]}^\nabla(R_V[[T]])$ can [Fa1, 2.1] be adapted?} 
\medskip
One approach towards answering this question is to try to associate to any morphism of $\widetilde{\Mm\Mf}_{[0,1]}^\nabla(R_V[[T]])$ ``a pull back" of it (via the natural $W(k)$-epimorphisms $R_V[[T]]\twoheadrightarrow V[[T]]\twoheadrightarrow k[[T]]$) to a morphism between two objects of $\Mm\Mf_{[0,1]}^\nabla(W(k)[[T]])[p]$. Such ``pull backs" can be constructed in many cases: it is easy to see that any object of $p-\widetilde{\Mm\Mf}_{[0,1]}^\nabla(R_V[[T]])$ is obtained by extension of scalars from an object of $p-\widetilde{\Mm\Mf}_{[0,1]}^\nabla(W(k)[[T_1,T]])$, where $W(k)[[T_1]]$ is the ring of formal power series ``defining" $R_V$ starting (see [Va2, 5.2.1]) from the $W(k)$-epimorphism $W(k)[[T_1]]\twoheadrightarrow V$, which takes $T_1$ into the chosen uniformizer of $V$.
\smallskip
The same method could be be used in connection to other regular $O$-schemes whose special fibres are having some properties (like are irreducible, etc.). Within these lines, we hope to achieve the classification of $S$-healthy schemes ``tentatively defined" in [Va2, top of p. 430].
\medskip
From 2.2.1.4, as in [Va2, 3.2.4], we get:
\medskip
{\bf 2.2.1.5. Corollary.} {\it Any (local) integral canonical model $\Mn$ of a Shimura quadruple $(G,X,H,v)$, with $v$ dividing $2$, has the SEP and so it is uniquely determined.}
\medskip
This result was known previously just for zero dimensional Shimura varieties (see [Va2, 3.2.8]). However, if ${\rm Sh}(G,X)$ is the elliptic modular curve, then it can be easily checked that $(G,X,H,v)$ has a uniquely determined smooth integral model having the SEP. If ${\rm Sh}(G,X)$ is a Siegel modular variety of genus at least $2$, it was not known previously that $(G,X,H,v)$ has a uniquely determined smooth integral model having either the SEP or the EP.   
\medskip
{\bf 2.2.1.5.1. Corollary.} {\it In [Va2, 3.2.7 4)], the conditions $(v,2)=1$ and $(p,2)=1$ can be removed. In particular, referring to 2.2.1.5, we get that the group ${\rm Aut}(G,X,H)$ of automorphisms of $(G,X,H)$ (see [Va2, 3.2.7 9)]) acts naturally on $\Mn$. Moreover in [Va2, 3.2.12] the condition $e<p-1$ can be replaced by $e\le max\{1,p-2\}$.} 
\medskip
{\bf 2.2.1.5.2. Uniqueness of projective integral canonical models: a second approach.} It is worth pointing out that there is a second approach (besides the one based on 2.2.1.4) to the uniqueness of a projective (local) integral canonical model $\Mn$ of a Shimura quadruple $(G,X,H,v)$ of compact type, with $v$ dividing an arbitrary prime $p$. To present it, we work in a slightly more general context and as the arguments for the local context (over $O_v$) are the same, we concentrate just on integral canonical models. Let $\Mn_1$ be another integral canonical model of $(G,X,H,v)$ which is quasi-projective. Let $H_0$ be a compact, open subgroup of $G(\AA_f^p)$ such that $\Mn_1$ (resp. $\Mn$) is a pro-\'etale cover of the $O_{(v)}$-scheme $\Mn_1/H_0$ (resp. $\Mn/H_0$) of finite type. We can assume $H_0$ is small enough so that the connected components of ${\rm Sh}_{H_0\times H}(G,X)$ are of general type. 
\smallskip
To show that $\Mn=\Mn_1$, it is enough to show that $\Mn/H_0=\Mn_1/H_0$. It is enough to show that $\Mn_1/H_0$ is projective over $O_{(v)}$ (cf. [MM, Theorem 2]). Using Nagata's embedding theorem (see [Vo]) we deduce that each connected component $\Mc_0$ of $\Mn_1/H_0$ is embeddable as an open subscheme into an integral, proper $O_{(v)}$-scheme $\Mc$. We can assume $\Mc_{E(G,X)}=\Mc_{0E(G,X)}$. If $\Mn_1/H_0$ is not projective over $O_{(v)}$, then $\Mc_0$ is not $\Mc$. So let $m_{\Mc}:{\rm Spec}(V)\to\Mc$ be a morphism, with $V$ a strictly henselian complete DVR of mixed characteristic $(0,p)$, such that the special point of ${\rm Spec}(V)$ maps into a point of $\Mc\setminus\Mc_0$. Let $K_V:=V[{1\over p}]$. Let $m_V:={\rm Spec}(K_V)\to{\rm Sh}_{H_0\times H}(G,X)$ be the resulting morphism. It extends to a morphism ${\rm Spec}(V)\to\Mn/H_0$ (as $\Mn/H_0$ is proper over $O_{(v)}$). It lifts to a morphism ${\rm Spec}(V)\to\Mn$. As $\Mn_1$ has the EP, the generic fibre of this last morphism extends to a morphism ${\rm Spec}(V)\to\Mn_1$. This contradicts the choice of $m_{\Mc}$. We conclude: $\Mn_1$ is canonically isomorphic to $\Mn$.
\medskip
{\bf 2.2.1.6. The non-filtered context.} Occasionally we speak about non-filtered Fontaine categories of $p$-divisible objects: $p-\Mm(*)$, $p-\Mm_{[a,b]}(*)$, $p-\Mm^\nabla(*)$ and $p-\Mm_{[a,b]}^\nabla(*)$, with $*\in\{R,X\}$, where $R$ and $X$ are as in 2.2.1 c), as well as about objects of them with tensors. Not to be long, as a variation, we detail the definition of an object with tensors of $p-\Mm_{[a,b]}^\nabla(R)$. A triple $(M,\vph,(t_{\al})_{\al\in\Mj})$ is called an object with tensors of $p-\Mm_{[a,b]}^\nabla(R)$, if locally in the Zariski topology of ${\rm Spec}(R)$, can be extended to a quadruple $(M,(F^i(M))_{i\in S(a,b)},\vph,(t_{\al})_{\al\in\Mj})$ defining a $p$-divisible object with tensors of $\Mm\Mf_{[a,b]}^\nabla(R)$. The morphisms in non-filtered Fontaine categories of $p$-divisible objects are the logical ones (for instance, in the context without connections, they are morphisms between underlying sheaves respecting --after inverting p-- Frobenius endomorphisms). Warning: if $a<0$ then $\vph(M)$ is not necessarily included in $M$.
\smallskip
Moreover, we speak about pull backs of objects of $p-\Mm^\nabla(X)$ through points of $X_k$ with values in perfect fields: we can assume $X$ is local and that a Frobenius lift of $X^\wedge$ is fixed; considering Teichm\"uller lifts of $X^\wedge$ with values in the Witt rings of such fields, the pull backs are defined as in 2.2.1.3 (by just ignoring filtrations).
\smallskip
We could define as well non-filtered Fontaine categories of objects: as we do not use them at all in this paper, except if $(a,b)=(0,1)$ or if we are in the context of $k$-schemes $S$ for which we do not know that $\Om_{S/k}$ is free or of finite type, these excepted contexts are detailed just at appropriate places (see b) of 2.2.4 B, 3.15.2 and 3.15.8-10 below).
However, we always identify $p-\Mm_{[0,b]}^\nabla(X)$ with a full subcategory of the category of $F$-crystals on $X_k$; so we speak about morphisms between truncations of objects of $p-\Mm_{[0,b]}^\nabla(X)$: we view them as morphisms between $F$-crystals on $X_k$ in coherent sheaves.
\medskip
{\bf 2.2.1.7. Complements.} {\bf 1)} We refer to the category $\NN-pro-\Mm\Mf(R)$ of 2.2.1 c). Each object ${\got C}$ of $\Mm\Mf(R)$ can be identified with the projective limit of ${\got C}/p^n{\got C}$, $n\in\NN$. So we view $\Mm\Mf(R)$ as a full subcategory of $\NN-pro-\Mm\Mf(R)$. By an almost $p$-divisible object of $\Mm\Mf(R)$ we mean an object ${\got C}$ of $\NN-pro-\Mm\Mf(R)$ which is the extension of a $p$-divisible object ${\got C}_1$ of $\Mm\Mf(R)$ by an object ${\got C}_2$ of $\Mm\Mf(X)$; when both ${\got C}_1$ and ${\got C}_2$ are not the zero object of $\NN-pro-\Mm\Mf(R)$, we speak about a non-trivial almost $p$-divisible object $\Mc$ of $\Mm\Mf(R)$. The usefulness of this notion springs from the fact that often the cokernel of a morphism of $p-\Mm\Mf(R)$ is an almost $p$-divisible object of $\Mm\Mf(R)$. Starting from 2.2.1.1 6) and 7), it is easy to see that the full subcategory $a-p-\Mm\Mf(R)$ of $\NN-pro-\Mm\Mf(R)$ whose objects are almost $p$-divisible objects of $\Mm\Mf(X)$ is abelian and $\ZZ_p$-linear. Its objects can be interpreted in a similar way to 2.2.1 c) (i.e. using finitely generated $R^\wedge$-modules which are direct sums of projective ones and of ones which locally in the Zariski topology have the DC property).
\smallskip
Similarly we define $a-p-FC$ for any Fontaine category $FC$ of objects.
\smallskip
{\bf 2)} We always identify $\Mm\Mf(W(k))$ with $\Mm\Mf^\nabla(W(k))$. An object of $\Mm\Mf(R)$ or of $\Mm\Mf^\nabla(R)$ is called trivial, if it is either the zero object or it is the extension of scalars of some  object of $\Mm\Mf(W(k))$ of the form $W_q(k)(m)$, $q\in\NN$, $m\in\ZZ$. If $m=0$ then we speak about a totally trivial object. A short exact sequence $0\to\Mo_1\hookrightarrow\Mo_2\twoheadrightarrow\Mo_3\to 0$ in $\Mm\Mf(R)$ or in $\Mm\Mf^\nabla(R)$ is said to be of direct summand type, if the underlying $R$-module of $\Mo_1$ is a direct summand of the underlying $R$-module of $\Mo_2$. 
\smallskip
An object of $\Mm\Mf(R)$ or of $\Mm\Mf^\nabla(R)$ is called solvable (resp. topologically solvable) if it is obtained from trivial ones via a finite number of short exact sequences of direct summand type (resp. via a finite number of short exact sequences). An object of $\Mm\Mf(R)$ or of $\Mm\Mf^\nabla(R)$ is called unipotent (resp. topologically unipotent) if it is obtained from totally trivial ones via a finite number of short exact sequences of direct summand type (resp. via a finite number of short exact sequences). An object of $\Mm\Mf(R)$ or of $\Mm\Mf^\nabla(R)$ is called quasi-solvable (resp. quasi-unipotent), if locally in the $\NN$-pro-\'etale topology of ${\rm Spec}(R)$ is solvable (resp. unipotent). Here, if ${\rm Spec}(Q)\to {\rm Spec}(R)$ is an \'etale morphism defined by a homomorphism $m:R\to Q$, then we endow $Q^\wedge$ with the unique Frobenius lift $\Phi_Q$ such that $\Phi_Q\circ m^\wedge=m^\wedge\circ\Phi_R$ and the pull back operations of objects of $\Mm\Mf(R)$ or of $\Mm\Mf^\nabla(R)$ are as in 2.2.1.3. Similarly, an object of $\Mm\Mf(R)$ or of $\Mm\Mf^\nabla(R)$ is called topologically quasi-solvable (resp. topologically quasi-unipotent), if locally in the $\NN$-pro-\'etale topology of ${\rm Spec}(R)$ is topologically solvable (resp. topologically unipotent).
\smallskip
Similarly we define (topologically) (quasi-) solvable and (quasi-) unipotent objects (resp. $p$-divisible objects) for any Fontaine category $FC$ (resp. $p-FC$) of objects (resp. of $p$-divisible objects). We get the full subcategories $top-q-solv-FC$, $top-solv-FC$, $top-q-unip-FC$, $top-unip-FC$, $q-solv-FC$, $solv-FC$, $q-unip-FC$ and respectively $unip-FC$ of $FC$ whose objects are topologically quasi-solvable, topologically solvable, topologically quasi-unipotent, topologically unipotent, quasi-solvable, solvable, quasi-unipotent and respectively unipotent. Similarly, we get $top-q-solv-p-FC$, $top-solv-p-FC$, $top-q-unip-p-FC$, $top-unip-p-FC$, $q-solv-p-FC$, $solv-p-FC$, $q-unip-p-FC$ and $unip-p-FC$. Of course, 1) and 2) can be combined. Moreover, the same applies to $FC[p^n]$ ($n\in\NN$): we speak about $top-q-solv-FC[p^n]$, $top-solv-FC[p^n]$, $q-solv-FC[p^n]$, $solv-FC[p^n]$, etc.; here, for instance, $solv-FC[p^n]$ is the full subcategory of $solv-FC$ formed by objects annihilated by $p^n$. An object of $FC[p]$ is solvable (resp. unipotent) iff it is topologically solvable (resp. topologically unipotent).  
\smallskip
{\bf 3)} We summarize different connections between most of the categories we have defined so far. Let $FC$ be a Fontaine category of objects. Denoting a full subcategory by a hook right arrow and a subcategory by the inclusion sign, we have
$$FC[p]\hookrightarrow FC[p^2]\hookrightarrow...\hookrightarrow FC[p^n]\hookrightarrow...\hookrightarrow p-FC\hookrightarrow a-p-FC\hookrightarrow \NN-pro-FC$$
$$unip-FC\hookrightarrow solv-FC\hookrightarrow q-solv-FC\hookrightarrow top-q-solv-FC$$
$$unip-FC\hookrightarrow q-unip-FC\hookrightarrow top-q-unip-FC\hookrightarrow top-q-solv-FC,$$
and, with $X$ as in 2.2.1 c), we have
$$isog-p-\Mm\Mf(X)\supset p-\Mm\Mf(X)\subset p-\Mm(X)$$
$$isog-p-\Mm\Mf(X)\supset p-\Mm\Mf^\nabla(X)\subset p-\Mm^\nabla(X).$$
\indent
{\bf 4)} Let $X$ (resp. $R$) be a flat scheme or algebra over $W(k)$. We assume a Frobenius lift of $X^\wedge$ (resp. of $R^\wedge$) is chosen. As in 2.2.1 c) (resp. in 2.2.1.6) we define $\Mm\Mf(X)$, $\Mm\Mf^\nabla(X)$, $p-\Mm\Mf(X)$, $p-\Mm\Mf^\nabla(X)$, $p-\Mm(X)$, $p-\Mm^\nabla(X)$, etc. Warning: we do not claim that all properties of the smooth context extend as well. For instance, we do not know when $\Mm\Mf(X)$ is abelian.  
\smallskip
{\bf 5)} Let $X_k$ be a regular, formally smooth $k$-scheme. We consider the full subcategory $UNIP(X_k)$ (resp. $TOP-UNIP(X_k)$) of the category of crystals on $X_k$ on coherent sheaves with the property that each object $Ob$ of it has the following three properties:
\medskip
i) there is $n\in\NN$ (depending on $Ob$) such that all sheaves of modules we get by evaluating $Ob$ at different thickenings of $X_k$-schemes are annihilated by $p^n$;
\smallskip
ii) if $\tilde X$ is an arbitrary regular, formally smooth $W(k)$-scheme lifting an open subscheme $\tilde X_k$ of $X_k$, the $\Mo_{\tilde X_{W_n(k)}}$-sheaf $\Mf$ obtained by evaluating $Ob$ at $\tilde X_{W_n(k)}$ (cf. also i)), locally in the Zariski topology of $\tilde X_{W_n(k)}$ has the DC property;
\smallskip
iii) the pair $(\Mf,\nabla)$, with $\nabla$ the connection on $\Mf$ defined by $Ob$, is unipotent, i.e. is obtained via short exact sequences of direct summand type in the same sense as of 2) (resp. via short exact sequences) from pull backs of $W_n(k)$-modules of finite type (and endowed with trivial connections).
\medskip
We refer to $UNIP(X_k)$ (resp. to $TOP-UNIP(X_k)$) as the category of unipotent (resp. of topologically unipotent) objects on $X_k$. As in 1), 2.2.1 c) and 2.2.1.1 7) we define the following 5 categories $p-UNIP(X_k)$, $isog-p-UNIP(X_k)$, $\NN-pro-UNIP(X_k)$, $UNIP(X_k)[p^n]$, $a-p-UNIP(X_k)$ (resp. $p-TOP-UNIP(X_k)$, $isog-p-TOP-UNIP(X_k)$, $\NN-pro-TOP-UNIP(X_k)$, $TOP-UNIP(X_k)[p^n]$, $a-p-TOP-UNIP(X_k)$) as well as truncations mod $p^m$, $m\in\NN$, of their objects. $UNIP(X_k)$ is a full subcategory of $TOP-UNIP(X_k)$ and the same remains true in the context of the mentioned 5 categories. To check ii) and iii) it is enough to deal with specific lifts of the members of an arbitrary open, affine cover of $X_k$.
\medskip
{\bf 2.2.1.8. Samples for 2.2.1 a) and d).} The isocrystal of an abelian variety over $k$ is 1/2-symmetric. If $k=\bar k$ then the filtered $\sg$-crystal of the canonical lift of an ordinary abelian
variety over $k$ is diagonalizable and a $p$-divisible object of
$\Mm\Mf_{[0,1]}(W(k))$.
\medskip
{\bf 2.2.2. Definitions. 1)} Let $a\in\ZZ$. An $a$-Lie $\sg$-crystal (resp. a Lie isocrystal over $k$) is a pair $({\got g},\vph)$, with {\got g}
a Lie algebra over $W(k)$ (resp. over $B(k)$), which as a $W(k)$-module (resp. as a $B(k)$-vector space) is free of finite rank (resp. of finite dimension), and
with $\vph$ a $\sg$-linear Lie automorphism of ${\got g}\otimes_{W(k)} B(k)$ such that
$\vph(p^a{\got g})\subset{\got g}$ (resp. with $\vph$ a $\sg$-linear Lie automorphism of $\got g$). If $a=1$ we drop it.
\smallskip
{\bf 2)} A Lie $\sg$-crystal is said to be ordinary or of ordinary type, if the
slopes of the isocrystal $({\got g}\otimes_{W(k)} B(k),\vph)$ are (warning!) precisely $-1$, 0 and 1, with the
multiplicity of the slope 1 equal to the multiplicity of the slope $-1$. 
\smallskip
{\bf 3)} A filtered Lie $\sg$-crystal is a
quadruple (${\got g},\vph,F^0({\got g}),F^1({\got g})$) such that:
\medskip
-- the pair $({\got g},\vph)$ is a Lie
$\sg$-crystal; 
\smallskip
-- $F^0({\got g})$ and $F^1({\got g})$ are Lie subalgebras of {\got g}, which as $W(k)$-submodules are direct summands of ${\got g}$;
\smallskip
-- $F^1({\got g})$ is abelian and $[{\got g},F^1({\got g})]\subset F^0({\got g})$; 
\smallskip
-- the quadruple (${\got g},F^0({\got g}),F^1({\got g}),\vph$) is a $p$-divisible object of $\Mm\Mf_{[-1,1]}(W(k))$, i.e. $\vph\bigl({1\over p}F^1({\got g})+
F^0({\got g})+p{\got g}\bigr)={\got g}$ (we have ${\got g}=F^{-1}({\got g})$ and $F^2({\got g})=\{0\}$).
\medskip
{\bf 4)} Let $R$ and $\Phi_R:R^\wedge\to R^\wedge$ be as in 2.2.1 c). By a Lie $p$-divisible object  of $\Mm\Mf_{[-1,1]}(R)$ (resp. of $\Mm\Mf_{[-1,1]}^\nabla(R)$) we mean a $p$-divisible object $\Ml$ of $\Mm\Mf_{[-1,1]}(R)$ (resp. of $\Mm\Mf_{[-1,1]}^\nabla(R)$) together with a morphism 
$$m_{\Ml}:=\Ml\otimes_{R^\wedge} \Ml\to \Ml$$ 
between $p$-divisible objects of $\Mm\Mf_{[-2,2]}(R)$ (resp. of $\Mm\Mf_{[-2,2]}^\nabla(R)$) satisfying the usual axioms of a Lie algebra. In practice the Lie algebra structure of $\Ml$ is obvious (i.e. is in-built in the underlying $R^\wedge$-module ${\got g}$ of $\Ml$) and so we always omit mentioning $m_{\Ml}$; so the notations are as in 3) above. We also refer to a Lie $p$-divisible object  of $\Mm\Mf_{[-1,1]}^\nabla(R)$ as a filtered Lie $F$-crystal over $R/pR$ (this matches to 3) above in the case $R=W(k)$); when we do not want to mention its filtration we refer to it as a Lie $F$-crystal over $R/pR$. 
\smallskip
All these apply to the context when we are dealing with an $X$ as in 2.2.1 c) instead of $R$ or when we are dealing with objects instead of $p$-divisible objects.
\medskip
{\bf 2.2.3. Remarks. 1)}  The slopes of a Lie (resp. filtered Lie) $\sg$-crystal are rational numbers of the interval
$[-1,\infty)$ (resp. of the interval $[-1,1]$). Most common:
\medskip
-- for Lie $\sg$-crystals $({\got g},\vph)$, {\got g} is the Lie algebra of a reductive group and so $\bigl([{\got g},{\got g}]\otimes_{W(k)} B(k),\vph\bigr)$ is a $0$-symmetric isocrystal (i.e. the multiplicity of a slope $\al$ is the same as the multiplicity of the slope $-\al$; this is a consequence of the fact that the Killing form on $[{\got g},{\got g}]\otimes_{W(k)} B(k)$ is non-degenerate and --as $\vph$ preserves the Lie bracket-- fixed by $\vph$; see also 3) below) whose slopes belong to the interval [-1,1];
\smallskip
-- for filtered Lie $\sg$-crystals (${\got g},\vph,F^0({\got g}),F^1({\got g})$), {\got g} is the Lie algebra of a reductive group and $F^0({\got g})$ is a parabolic Lie subalgebra of {\got g} having $F^1({\got g})$ as its nilpotent radical.
\medskip
We sometimes need to use Lie $\sg$-crystals $({\got g},\vph)$, with ${\got g}$ as a Lie subalgebra of the Lie algebra of a reductive group over $W(k)$.
\smallskip
{\bf 2)} We introduce similarly to 2.2.2 3) the notion of an $[a,b]$-filtered Lie $\sg$-crystal $\Ml:=\bigl({\got g},\vph,F^{a+1}({\got g}),\ldots,F^b({\got g})\bigr)$, where $a,b\in\ZZ$, $b\ge a$, by requiring that it is a $p$-divisible object of $\Mm\Mf_{[a,b]}(W(k))$ and ${\got g}$ has a natural Lie algebra structure such that the Lie bracket morphism ${\got g}\otimes_{W(k)} {\got g}\to {\got g}$ defines a morphism $\Ml\otimes_{W(k)}\Ml\to\Ml$ between objects of $p-\Mm\Mf_{[2a,2b]}(W(k))$. 
So an $[-1,1]$-filtered Lie $\sg$-crystal is the same as a filtered Lie $\sg$-crystal. Forgetting the Lie structure, the pair $({\got g},\vph)$ becomes a latticed isocrystal (if $a\notin\NN$ then it is a $\sg$-crystal). As in 2.2.3 4), we define (filtered) $[a,b]$-Lie $F$-crystals. We always drop $[-1,1]$. 
\smallskip
{\bf 3)} We recall the following well known result:
\medskip
{\bf Lemma.} {\it Any isocrystal $(V,\vph_V)$ over $k$ is a direct sum of isocrystals having only one slope.}
\medskip
{\bf Proof (slightly sketched):} Dieudonn\'e's classification (see [Di] or [Man]) of isocrystals over $\bar k$ implies that this holds over $\bar k$ and moreover $V\otimes_{B(k)} B(\bar k)$ has a nice $W(\bar k)$-lattice $L$. The word ``nice" is used here in the following sense: there is $n$ such that $p^n(\vph_V\otimes 1)(L)\subset L$ and $L$ is a direct sum of $W(\bar k)$-lattices taken by $p^{n}(\vph_V\otimes 1)$ into themselves in such a way that the Hodge numbers of the restriction of $p^{n}(\vph_V\otimes 1)$ to any member of this direct sum, are all the same. But any such $W(\bar k)$-lattice $L$ is defined over the strict henselization of $W(k)$ and so over $W(k_1)$, with $k_1$ a finite field extension of $k$. Using (here is the sketched part) standard matrix computations over $W(k_1)$ (similar to the ones defining the $F^1$-filtrations of canonical lifts of ordinary $p$-divisible groups over $k_1$) the statement of the Lemma is true over $k_1$. Using standard Galois descent, the Lemma follows. 
\medskip
So if $\al\in\RR$ and if $(N,\vph_N)$ is an isocrystal (resp. a latticed isocrystal) over $k$, then we can speak about the $B(k)$-vector subspace (resp. $W(k)$-submodule which is a direct summand) $W_{\al}(N,\vph_N)$ of $N$ corresponding to slopes of $(N,\vph_N)$ greater or equal to $\al$: $W_{\al}(N,\vph_N)$ is maximal with the property that it is taken by $\vph_N$ (resp. by $\vph_N$ times some positive, integral power of $p$) into itself and all slopes of $(W_{\al}(N),\vph_N)$ are greater or equal to $\al$. Similarly, we define $W^{\al}(N,\vph_N)$, by replacing ``greater" by ``smaller", and $W(\al)(N,\vph_N)$ by deleting ``greater or". So 
$$W(\al)(N,\vph_N)\supset W_{\al}(N,\vph_N)\cap W^{\al}(N,\vph_N)$$
and by inverting $p$ we have an isomorphism. 
\smallskip
Let $({\got g},\vph)$ be a Lie isocrystal over $k$. Then ${\got p}_{\ge 0}:=W_{0}({\got g},\vph)$ is a Lie subalgebra of ${\got g}$.
We call the (Lie) isocrystal $({\got p}_{\ge 0},\vph)$ (resp. the Lie algebra ${\got p}_{\ge 0}$) the (Lie) subisocrystal of $({\got g},\vph)$ (resp. the Lie subalgebra of $\got g$) corresponding to non-negative slopes of $({\got g},\vph)$. Similarly we define things for non-positive (or positive, or negative) slopes, or for the slope $0$, or for slopes greater (resp. smaller) than some $r\in [0,\infty)$ (resp. than some $r\in (-\infty,0]$). In particular, ${\got p}_{\le 0}:=W^{0}({\got g},\vph)$ and ${\got p}_{=0}:=W(0)({\got g},\vph)$ are often used. Similarly, we denote by ${\got p}_{>0}$ (resp. by ${\got p}_{<0}$) the Lie subalgebra of $\got g$ corresponding to positive (resp. negative) slopes of $({\got g},\vph)$. Warning: similarly we define and denote things for a Lie $\sg$-crystal. 
\medskip
{\bf Claim.} {\it If $\got g$ is the Lie algebra of a semisimple group over $B(k)$, then ${\got p}_{\ge 0}$ is a parabolic Lie subalgebra of $\got g$.}
\medskip
{\bf Proof:} As this well known Claim is used a lot in what follows, we include a proof of it. Let (cf. Lemma)
$${\got g}=\oplus_{\al\in SL} {\got c}_{\al}$$ 
be the slope decomposition of ${\got g}$ achieved by $\vph$. So the isocrystal $({\got c}_{\al},\vph)$ has only one slope $\al$, $\forall\al\in SL$. Let $G_{B(k)}$ be the adjoint group over $B(k)$ having ${\got g}$ as its Lie algebra. ${\got c}_{\al}$ is perpendicular on ${\got c}_{\be}$ w.r.t. the Killing form $KIL$ on ${\got g}$, provided $\be\in SL$ is different from $-\al\in SL$. Denoting ${\got C}_{s}=\{0\}$, if $s\notin SL$, this is a consequence of the following Lie bracket property
$$[{\got c}_{\al},{\got c}_{\be}]\subset {\got c}_{\al+\be},\leqno (LBP)$$
$\forall\al$, $\be\in SL$.
\smallskip
As the center of ${\got g}$ is trivial, ${\got p}_{=0}$ is the Lie subalgebra of ${\got g}$ normalizing ${\got c}_{\al}$, $\forall\al\in SL$. So (see [Bo2, 7.1 and 7.4]) it is the Lie algebra of a connected subgroup $P_{=0}$ of $G_{B(k)}$. We know that the restriction of $KIL$ to ${\got p}_{=0}$ is non-degenerate. From [Bou1, ch. 1, th. 4] we get that $P_{=0}$ is a reductive group.
\smallskip
We can assume $k=\bar k$. So we can assume there is a maximal torus $T$ of $P_{=0}$ whose Lie algebra is stable under $\vph$. The centralizer of $T$ in $G_{B(k)}$ is a reductive group $C(T)$ (see [Bo2, 11.12]) whose Lie algebra is normalized by $\vph$. The restriction of $KIL$ to this Lie algebra is still non-degenerate and $T=C(T)\cap P_{=0}$; so to check that $T$ is a maximal torus of $G_{B(k)}$ we just need to check that $C(T)^{\rm der}$ is trivial. It is enough to show: if $C(T)^{\rm der}$ is non-trivial, then $({\rm Lie}(C(T)^{\rm der}),\vph)$ has slope $0$ with positive multiplicity. Shifting our attention from $C(T)^{\rm der}$ to $G_{B(k)}$ (in order to use previous notations), we need to show $T$ is non-trivial. If $T$ is trivial then $P_{=0}$ is trivial and then, for $\al\in SL\cap (0,\infty]$ and for the greatest element $\be\in SL$ we get $[{\got c}_{\be},{\got c}_{-\be}]=\{0\}$, $[{\got c}_{\be},{\got c}_{\al}]=\{0\}$ and $[{\got c}_{-\be},{\got c}_{-\al}]=\{0\}$; so (cf. Jacobi's identity) ${\got c}_{\be}$ and ${\got c}_{-\be}$ are perpendicular w.r.t. $KIL$ and so are included in the annihilator of $KIL$. Contradiction. So $T$ is a maximal torus of $G_{B(k)}$.
\smallskip
Let $\Phi$ be the set of roots of the action of $T_{\overline{B(k)}}$ via inner conjugation on ${\got g}\otimes_{B(k)} \overline{B(k)}$, and let $\Phi^{\ge 0}$ be its subset corresponding to the action of $T_{\overline{B(k)}}$ on ${\got p}_{\ge 0}$. $\Phi^{\ge 0}$ is stable under addition inside $\Phi$ and moreover $\Phi=-\Phi^{\ge 0}\cup\Phi^{\ge 0}$. This last equality is a consequence of the fact that ${\got c}_{\al}$ is the dual of ${\got c}_{-\al}$ w.r.t. $KIL$, $\forall\al\in SL$. It is easy to see that this implies $\Phi^{\ge 0}$ contains a base of $\Phi$. So the Claim follows from [Bo2, 14.17-18] applied over $\overline{B(k)}$.  
\medskip
{\bf Comments.} As ${\got p}_{>0}$ is perpendicular on ${\got p}_{\ge 0}$ w.r.t. $KIL$ and as ${\got p}_{=0}$ is reductive, from [Bou1, Prop. 6 b) of p. 81] we get that ${\got p}_{>0}$ is the nilpotent radical of ${\got p}_{\ge 0}$ and so it is included in a Borel Lie subalgebra of ${\got g}\otimes_{B(k)} \overline{B(k)}$ containing ${\rm Lie}(T)\otimes_{B(k)} \overline{B(k)}$. So, as $T$ normalizes ${\got c}_{\al}$, from [Bo2, p. 184] applied over $\overline{B(k)}$, we get that ${\got c}_{\al}$ is the Lie algebra of a unipotent subgroup of $G_{B(k)}$, $\forall\al\in SL$. From [Bo2, 7.1] we get that all these hold for any (perfect) $k$.
\smallskip
We have a duality of language: positive versus negative, non-negative versus non-positive, etc. So the Claim implies as well that ${\got p}_{\le 0}$ is the Lie algebra of a parabolic subgroup of $G_{B(k)}$. However, often in problems pertaining to Lie $\sg$-crystals it is more convenient to work with a specific choice (like with ${\got p}_{\ge 0}$ instead of ${\got p}_{\le 0}$, etc.); warning: the choice depends on the problem. 
\medskip
The above Claim still makes sense over $W(k)$ due to the following general simple Fact:
\medskip
{\bf Fact.} {\it Let $D$ be an integral Dedekind scheme such that all its closed points are of positive characteristic. Let $G_D$ be a reductive group scheme over $D$. Then the Zariski closure in $G_D$ of any parabolic (resp. Borel) subgroup of the generic fibre of $G_D$, is a parabolic (resp. Borel) subgroup $P_D$ (resp. $B_D$) of $G_D$. In particular, $G_D$ is a quasi-split reductive group iff its generic fibre is.}
\medskip
{\bf Proof:} As the operation of taking the Zariski closure is compatible with flat morphisms, we can assume $D$ is the spectrum of a complete DVR $V$ having an algebraically closed residue field. Let $K_V$ be the field of fractions of $V$. We can restrict just to the part involving $P_D$. The special shape of $D$ implies $G_D$ is split; we deduce the existence of a parabolic subgroup $P^1_D$ of $G_D$ whose generic fibre is $G_D(K_V)$-conjugate to the generic fibre of $P_D$ (cf. also [Bo2, 21.11] applied over $K_V$). Using Iwasawa's decomposition (see [Ti2, 3.3.2]), we get that $P_D$ and $P_D^1$ are $G_D(D)$-conjugate and so the Fact follows. 
\smallskip
{\bf 4)} One may ask: why Lie $\sg$-crystals? There are three reasons. First, this is the best way to keep track of the relative situation which involves tensors (independently of some fixed ``embedding": see 2.2.8 for a sample; see also 4.9.8 below). Second, we have natural equivalences between different categories $\Mm\Mf_{[a,b]}(*)$ and $\Mm\Mf_{[a+c,b+c]}(*)$ (or $\Mm\Mf_{[a,b]}^\nabla(*)$ and $\Mm\Mf_{[a+c,b+c]}^\nabla(*)$), with $a,b,c\in\ZZ$, $b\ge a$. When we take $End$'s of their objects (see 2.2.4 B below), $c$ goes automatically away; so we have a better way to express our results (as samples see 3.6.18.7.0, 3.6.18.8.1 b), 4.7.11 4), etc.) and moreover, we do not want to lose track of the extra Lie structure of such $End$'s objects. Third, we can read very easily different versal conditions (involving deformations of --to be defined in 2.2.8 4a)-- generalized Shimura $p$-divisible objects over $k$) by looking at attached Lie $F$-crystals.
\smallskip
In some sense, \S1-4 grew out of a systematic exploitation of Lie $\sg$-crystals. Accordingly, to always have a good way of keeping track (without always emphasizing the extra Lie structures) of when we are dealing with filtered Lie $\sg$-crystals instead of just some $p$-divisible objects of $\Mm\Mf_{[-1,1]}(W(k))$, most of the time we use the ordering $({\got g},\vph,F^0({\got g}),F^1({\got g}))$ instead of the standard one $({\got g},F^0({\got g}),F^1({\got g}),\vph)$; the same applies to the situation when $W(k)$ is replaced by some $R$ or $X$ as in 2.2.1 c).
\medskip
{\bf 2.2.4. Operations with ($p$-divisible) objects.} In this section we deal with different ``tannakian" considerations. We mostly concentrate on what we need; however, some material which is not needed in the rest of this paper is as well included. In what follows we use the tannakian language as of [De5, ch. 2] (see also [Mi4, p. 288-9]) except that we allow the base scheme not to be necessarily the spectrum of a field. Till end of 2.2.4.1 we have $n\in\NN$.
\medskip
{\bf A. Duals.} No Fontaine category $FC$ of objects has identity object with respect to tensor products. However, $FC[p^n]$, $p-FC$ and $isog-p-FC$ do have such identity objects. For instance, referring to 2.2.1 c), $R/p^nR(0)$ is the identity element of $\Mm\Mf(R)[p^n]$. So we can define duals of objects of $FC[p^n]$, of $p-FC$ and of $isog-p-FC$. For instance, if ${\got C}$ is an object of $FC[p^n]$, then its dual 
$${\got C}^*:=Hom_{FC}({\got C},R/p^nR(0))$$ 
is computed as follows. Let us assume $FC=\Mm\Mf_{[a,b]}(X)$ (the $\nabla$ context is the same). If ${\got C}=(M,(F^i(M))_{i\in S(a,b)},(\vph_i)_{i\in S(a,b)})$, then 
$$Hom_{FC}({\got C},R/p^nR(0)):=({\rm Hom}_R(M,R/p^nR),F^i({\rm Hom}_R(M,R/p^nR))_{i\in S^-},(\tilde\vph_i)_{i\in S^-}),$$
with $S^-:=S(-b,-a)$, with
$$F^i({\rm Hom}_R(M,R/p^nR)):={\rm Ker}({\rm Hom}_R(M,R/p^nR),{\rm Hom}_R(F^{1-i}(M),R/p^nR))$$ 
and with $\tilde\vph_i$ defined by the formula
$$\tilde\vph_i(f_i)(\vph_j(x_j)):=p^{-i-j}f_i(x_j),$$ 
where $i\in S(-b,-a)$, $j\in S(a,b)$,  $x_j\in F^j(M)$ and $f_i\in F^i({\rm Hom}_R(M,R/p^nR))$. This last equality makes sense as for $i+j>0$ we have $f_i(x_j)=0$.
It does not depend on the choice of $n$ such that ${\got C}$ is an object of $FC[p^n]$: the natural monomorphism $R/p^{n}R(0)\hookrightarrow R/p^{n+1}R(0)$ defined by the multiplication (inside the ring $R/p^{n+1}R$) with $p$, achieves a natural isomorphism 
$$Hom_{FC}({\got C},R/p^nR(0))\tilde\to Hom_{FC}({\got C},R/p^{n+1}R(0)).$$ 
If $X_k$ has a finite number of connected components, then 
$$FC=\cup_{n\in\NN} FC[p^n]$$
(cf. [Fa1, 2.1] ii)). So, as the dual objects can be defined restricting to each connected component of $X^\wedge$, we get that, in general, each object of $FC$ has a uniquely determined (up to isomorphism) dual.
\smallskip
Similarly, referring to 2.2.1.7 5) we have:
\medskip
{\bf Lemma.}  {\it a) $TOP-UNIP(X_k)[p^n]$ and $UNIP(X_k)[p^n]$ are $\ZZ/p^n\ZZ$-linear, abelian categories endowed with tensors products and duals.
\smallskip
b) We assume that any object of $UNIP(X_k)[p]$ (resp. of $TOP-UNIP(X_k)[p]$) is the cokernel of an isogeny between two objects of $p-UNIP(X_k)$ (resp. of $p-TOP-UNIP(X_k)$). Then any object of $UNIP(X_k)$ (resp. of $TOP-UNIP(X_k)$) is the cokernel of an isogeny between two objects of $p-UNIP(X_k)$ (resp. of $p-TOP-UNIP(X_k)$).}
\medskip
{\bf Proof:} We can assume $X_k$ is connected. a) is just a variant of [Fa1, p. 31-3]: only the abelian part is not entirely trivial. As we have duals (this is as above), we just need to check that $UNIP(X_k)[p^n]$ and $TOP-UNIP(X_k)[p^n]$ have kernels and images. We first consider the case of $TOP-UNIP(X_k)[p^n]$. Its object $(\Mo_{X_k},d)$, where $d$ is just the $d$-operator viewed as a connection on $\Mo_{X_k}$, is simple (i.e. has no proper subobjects) and moreover, each object of $TOP-UNIP(X_k)[p^n]$ is obtained from it through a finite sequence of short exact sequences. Let $m:\Mu_1\to\Mu_2$ be a morphism of $TOP-UNIP(X_k)[p^n]$. Using induction on the length of the underlying $W(k_1)$-module of the pull back of $\Mu_1\oplus\Mu_2$ via a $W(k)$-morphism $z_1:{\rm Spec}(W(k_1))\to X$, with $k_1$ a perfect field containing $k$, based on the fact that $\Mo_{X_k}$ is the only simple object of $TOP-UNIP(X_k)[p^n]$, we can assume $\Mu_1=\Mu_2=(\Mo_{X_k},d)$. In this case $m$ is either an isomorphism or it is the zero morphism. The passage from $TOP-UNIP(X_k)[p^n]$ to $UNIP(X_k)[p^n]$ is an immediate consequence of the following Sublemma, whose easy proof is left as an exercise. 
\medskip
{\bf Sublemma.} {\it Let $O$ and $\pi$ be as in 2.2.1 e). Let $R_O$ be a flat $O$-algebra. Let $0\to M_1\hookrightarrow M\twoheadrightarrow M_2\to 0$ be a short exact sequence of direct summand type, with $M$, $M_1$ and $M_2$ as $R_O$-modules having the DC property. Let $N$ be a $R_O$-submodule of $M$ such that $N$, $N\cap M_1$ and $N/N\cap M_1$ have the DC property. Then $N\cap M_1$ is a direct summand of $N$.}
\medskip
For b) we can restrict to the case $UNIP(X_k)$ (as the case of $TOP-UNIP(X_k)$ is entirely the same). Let now $\Mu$ be an object of $UNIP(X_k)$. To check b) for it, we use again induction on the length of the underlying $W(k_1)$-module of the pull back of $\Mu$ through $z_1$. If $\Mu=(\Mo_{X_k},d)$, then $\Mu$ is the truncation mod $p$ of $(\Mo_{X},d)$. For the general case, we use a short exact sequence 
$$0\to (\Mo_{X_k},d)\hookrightarrow\Mu\twoheadrightarrow\Mu_1\to 0.$$ 
By induction, we consider an epimorphism $e_1:\Mu^1\twoheadrightarrow\Mu_1$, with $\Mu^1=(\Mf_1,d_1)$ an object of $p-UNIP(X_k)$. Let $e_2:\Mu^2\to\Mu$ and $\Mu^2\to\Mu^1$ be the morphisms of $a-p-UNIX(X_k)$ defining the fibre product of $\Mu\to\Mu_1$ and $e_1$. $p\Mu_2$ is an object of $p-UNIP(X_k)$ (this can be checked easily by ``evaluating" at $\tilde X^\wedge$'s as in 2.2.1.7 5)). Using the short exact sequence
$$0\to p\Mu_2\hookrightarrow\Mu_2\twoheadrightarrow\Mu_2/p\Mu_2\to 0,$$
and similar arguments involving pull backs, it is enough to show that $\Mu_2/p\Mu_2$ is the cokernel of an isogeny between objects of $p-UNIP(X_k)$. But this is or hypothesis. This ends the proof.
\medskip
{\bf Remark.} We assume $X$ is proper and geometrically connected. We can assume (cf. Grothendieck's algebraization theorem) that all objects of $p-UNIP(X_k)$ involve locally free $\Mo_X$-sheaves and connections on them. The hypothesis of b) is satisfied if different cohomological conditions hold (one such condition could be: $H^1(X,\Mf)$ is a torsion free $W(k)$-module and the restriction map $H^1(X,\Mf)\to H^1(X_k,\Mf/p\Mf)$ is surjective, for all locally free $\Mo_{X}$-sheaves which are obtained via short exact sequences from $\Mo_{X}$). 
\medskip
{\bf B. Ends.} Let $R$, $\Phi_R:R^\wedge\to R^\wedge$, $a$ and $b$ be as in 2.2.1 c). Let ${\got C}=(M,(F^i(M))_{i\in S(a,b)},(\vph_i)_{i\in S(a,b)})$ be an object of $\Mm\Mf_{[a,b]}(R)$. We can assume ${\rm Spec}(R/pR)$ has a finite number of connected components; so we view $M$ as an $R$-module (and not just as an $R^\wedge$-module). For $j\in SS(a,b)$, let 
$$F^j({\rm End}(M)):=\{m\in {\rm End}(M)|m(F^i(M))\subset F^{i+j}(M),\,\forall i\in S(a,b)\}.$$ 
As $F^i(M)$'s (as $R$-modules) are direct summands of $M$ (see [Fa1, 2.1 i)]), $F^j({\rm End}(M))$'s are as well direct summands (as $R$-modules) of ${\rm End}(M)$. The object of $\Mm\Mf_{[a-b,b-a]}(R)$
$$End({\got C}):=({\rm End}(M),(F^j({\rm End}(M)))_{j\in SS(a,b)},(\tilde\vph_j)_{j\in SS(a,b)})$$
is called the $End$ object of ${\got C}$. $\tilde\vph_j$'s are obtained from $\vph_i$'s in the logical way: if $m\in F^i({\rm End}(M))$ and $x\in F^j(M)$, then 
$$\tilde\vph_i(m)(\vph_j(x)):=\vph_{i+j}(m(x)).\leqno (ENDFR)$$ 
So, if $M$ is a projective $R/p^nR$-module for some $n\in\NN$, then $\tilde\vph_j$'s are obtained via the canonical identification ${\rm End}(M)=M\otimes_{R/p^nR} M^*$. The natural Lie bracket map ${\rm End}(M)\otimes_R {\rm End}(M)\to {\rm End}(M)$, defines a morphism 
$$End({\got C})\otimes_R End({\got C})\to End({\got C})$$
between objects of $\Mm\Mf_{[2a-2b,2b-2a]}(R)$.
\smallskip
Similarly, we define $End$'s of objects of $FC$ or of $p-FC$, for an arbitrary Fontaine category $FC$ of objects. In particular, if ${\got C}$ is a $p$-divisible object of $\Mm\Mf_{[0,1]}(W(k))$ (i.e. if ${\got C}$ is a filtered $\sg$-crystal), then any Lie $p$-divisible subobject of $End({\got C})$ is a filtered Lie $\sg$-crystal.
\smallskip
Similarly one defines $Hom$ objects (resp. $p$-divisible objects): they are direct summands of suitable $End$ objects (resp. $p$-divisible objects). Moreover, we have a natural, functorial (bilinear) identification $Hom({\got C},{\got C}_1)={\got C}^*\otimes {\got C}_1$.
\smallskip
If $X$ or $X^\wedge$ is not equipped with a Frobenius lift, then for the definition of $End$'s of objects of $\Mm\Mf_{[a,b]}^\nabla(X)$ we assume $b-a\le p-1$
and we think of $End({\got C})$ as the collection of $End$ objects indexed by pairs $(U,\Phi_U)$, where $U$ is an open subscheme of $X$ having a non-trivial special fibre and $\Phi_U$ is a Frobenius lift of $U^\wedge$, the $End$ object corresponding to such a pair being $End({\got C}_U)$, with ${\got C}_U$ as the pull back (restriction) of ${\got C}$ to $U$. The natural gluings (see [Fa1, p. 34-35]) between such restrictions result in natural gluings between $End({\got C}_U)$'s. See C and D below for a fancier presentation of such gluings.
\smallskip
Warning: taking $End$'s is not a functorial process. 
\medskip
{\bf Exercise.} {\bf a)} We assume $X_k$ has a finite number of connected components. Then for any 2 objects ${\got C}_1$ and ${\got C}_2$ of $\Mm\Mf(X)$, the $\ZZ_p$-module ${\rm Hom}({\got C}_1,{\got C}_2)$ is finite.
\smallskip
{\bf b)} Let $b\in\NN\cup\{0\}$. We assume $k=\bar k$. Let ${\got C}_1=(M_1,\vph_1)$ and ${\got C}_2=(M_2,\vph_2)$) be two objects of $p-\Mm_{[0,b]}(W(k))$. Then the $\ZZ_p$-module of morphisms between ${\got C}_1/p^n{\got C}_1$ and ${\got C}_2/p^n{\got C}_2$ which lift to morphisms between ${\got C}_1/p^{n+b}{\got C}_1$ and ${\got C}_2/p^{n+b}{\got C}_2$ is finite and by pull backs via the canonical morphism ${\rm Spec}(W(k_1))\to {\rm Spec}(W(k))$, with $k_1$ an algebraically closed field containing $k$, it remains the same. Here morphisms are interpreted in terms of $\Mo_{X^\wedge}$-modules endowed with Frobenius endomorphisms.
\medskip
{\bf Hints.} For a), using Teichm\"uller lifts we can assume $X={\rm Spec}(W(k))$; using devisage and the existence of $Hom$-objects the situation gets reduced to case when both ${\got C}_1$ and ${\got C}_2$ are annihilated by $p$. For b), as above we can speak about $Hom({\got C}_1,{\got C}_2)=({\rm Hom}(M_1,M_2),\vph)$ and about the maximal $W(k)$-submodule $M_0$ of ${\rm Hom}(M_1,M_2)$ taken by $\vph$ into ${\rm Hom}(M_1,M_2)$; each morphism we are interested in can be interpretated defines uniquely an element of $M_0/p^nM_0$ fixed by $\vph$. 
\medskip
{\bf C. Tensorial closures.} The problem we face (see the last part of B) is: if $X^\wedge$ is not equipped with a Frobenius lift, then just $\Mm\Mf_{[a,a+p-1]}^\nabla(X)$ is presently well defined; from a tannakian point of view this is slightly inconvenient (see E and F below; see the proof of J below for the use of the word slightly). So it is desirable to define in some way $\Mm\Mf_{[a,a+p+i]}^\nabla(X)$, with $i\in\NN\cup\{0\}$. There are many ways to proceed: accordingly, as below we deal with four such ways, the first (resp. the second, the third and the fourth) one will have the upper index $\nabla(tens)$ (resp. $\nabla(loc-tens)$, $\nabla(p+tens)$ and $\nabla(p+loc-tens)$) instead of $\nabla$. The constructions have to be performed for each connected component of $X^\wedge$; so we assume $X_k$ is connected (i.e. is integral). Let 
$$\UU(X)$$ 
be the set of pairs $(U,\Phi_U)$ as in B. As in what follows we work with two or more Frobenius lifts of the $p$-adic completion of some fixed open subscheme  of $X$, to emphasize which Frobenius lift we use to define Fontaine categories involving connections, we denote by $\Mm\Mf_{[a,b]}^\nabla(U,\Phi_U)$ and by $\Mm\Mf^\nabla(U,\Phi_U)$ the Fontaine categories defined by $(U,\Phi_U)\in\UU(X)$ and which previously were defined just by $\Mm\Mf_{[a,b]}^\nabla(U)$ and respectively by $\Mm\Mf^\nabla(U)$. The pull back functor
$$\Mm\Mf_{[a,a+p-1]}^{\nabla}(X)\to \Mm\Mf^\nabla(U,\Phi_U),$$ 
is faithful; in general it is not fully and this motivates what follows. If $\Phi_U^1$ is a second Frobenius lift of $U^\wedge$, then (see [Fa1, proof of 2.3]) we have a canonical isomorphism of categories
$$i(U,a,\Phi_U,\Phi_U^1):\Mm\Mf_{[a,a+p-1]}^\nabla(U,\Phi_U)\tilde\to\Mm\Mf_{[a,a+p-1]}^\nabla(U,\Phi_U^1)$$
and, if $\Phi_U^2$ is a third Frobenius lift of $U^\wedge$, we have the following cocycle condition
$$i(U,a,\Phi_U^1,\Phi_U^2)\circ i(U,a,\Phi_U,\Phi_U^1)=i(U,a,\Phi_U,\Phi_U^2).$$ 
\indent
{\bf Definition.} If $TC$ is a tensorial category which is small, then by the tensorial closure in $TC$ of a subcategory $SC$ of $TC$ we mean the smallest tensorial subcategory (i.e. the intersection of all tensorial subcategories) of $TC$ having $SC$ as a subcategory. 
\medskip
Let 
$$\Mm\Mf^{\nabla(tens)}(U,\Phi_U)$$
be the tensorial closure of the union of the categories $\Mm\Mf_{[c,c+p-1]}^\nabla(U,\Phi_U)$ taken inside $\Mm\Mf^\nabla(U,\Phi_U)$. We will not include here a big ado about tensorial closures of arbitrary abelian categories. We just point out that there is a standard way to check that $\Mm\Mf^{\nabla(tens)}(U,\Phi_U)$ does not depend (up to canonical isomorphism) on the choice of $\Phi_U$.
We have:
\medskip
{\bf Lemma.} {\it There is a uniquely determined isomorphism
$$i(U,\Phi_U,\Phi_U^1):\Mm\Mf^{\nabla(tens)}(U,\Phi_U)\tilde\to\Mm\Mf^{\nabla(tens)}(U,\Phi_U^1)$$
of tensorial categories which extends $i(U,c,\Phi_U,\Phi_U^1)$, $\forall c\in\ZZ$.} 
\medskip
{\bf Proof:} We consider the product category 
$$TS:=\Mm\Mf^\nabla(U,\Phi_U)\times \Mm\Mf^\nabla(U,\Phi_U^1)$$ and the smallest tensorial subcategory $SC$ of it containing the union of the categories $\Mm\Mf_{[c,c+p-1]}^\nabla(U,\Phi_U)$, $c\in\ZZ$, embedded into $TS$ via the natural inclusion 
$$i_{U,\Phi_U}^c:\Mm\Mf_{[c,c+p-1]}^\nabla(U,\Phi_U)\hookrightarrow \Mm\Mf^\nabla(U,\Phi_U)$$ and the composite of $i(U,c,\Phi_U,\Phi_U^1)$ with the similarly defined natural inclusion $i_{U,\Phi_U^1}^c$. We consider the natural (projection) functors $F:SC\to\Mm\Mf^{\nabla(tens)}(U,\Phi_U)$ and $F_1:SC\to\Mm\Mf^{\nabla(tens)}(U,\Phi_U^1)$. From very definitions we get: they are surjective (on morphisms). We also consider the functor $FF:\Mm\Mf^\nabla(U,\Phi_U)\to \Mm\Mf(W(k_1))$ (resp. $FF_1:\Mm\Mf^\nabla(U,\Phi_U^1)\to \Mm\Mf(W(k_1))$) defined by pull back through a $\Phi_U$-Teichm\"uller lift of $U^\wedge$ in a $k_1$-valued point of it, with $k_1$ an arbitrary perfect field. The fact that $F$ and $F_1$ are injective (on morphisms), results from the following two Facts.
\medskip
{\bf F1.} $\forall c\in\ZZ$, the image of $\Mm\Mf^\nabla_{[c,c+p-1]}(U,\Phi_U)$ in $\Mm\Mf(W(k_1))\times \Mm\Mf(W(k_1))$ defined naturally via $FF$ and $FF_1$ factors through $\Mm\Mf(W(k_1))$ diagonally embedded in $\Mm\Mf(W(k_1))\times \Mm\Mf(W(k_1))$.
\medskip
{\bf F2.} $\forall c\in\ZZ$, the image of $\Mm\Mf^\nabla_{[c,c+p-1]}(U,\Phi_U)$ in $COH(U)\times COH(U)$ defined naturally via the functor that associates to an object $Ob$ of $TS$ the underlying $\Mo_U$-sheaves of $FF(Ob)$ and of $FF_1(Ob)$ factors through $COH(U)$ diagonally embedded in $COH(U)$.
\medskip
Both these two Facts can be checked locally (so we can assume $U$ is affine) and so are a consequence of [Fa1, proof of 2.3]. We get: $F$ and $F_1$ are isomorphisms of categories and so we can take
$$i(U,\Phi_U,\Phi_U^1):=F_1\circ F^{-1}.$$
This ends the proof.
\medskip
{\bf Corollary.} {\it We have the following cocycle condition
$$i(U,\Phi_U^1,\Phi_U^2)\circ i(U,\Phi_U,\Phi_U^1)=i(U,\Phi_U,\Phi_U^2).$$}
\indent
We define $\Mm\Mf_{[a,a+p+i]}^{\nabla(tens)}(X)$ to be the category whose objects are obtained from the ones of the categories $\Mm\Mf_{[c,c+p-1]}^\nabla(X)$, $c\in\ZZ$, through the following 4 operations: 
\medskip
-- taking duals and tensor products as well as kernels and cokernels compatible with the standard gluings and whose underlying $\Mo_{X^\wedge}$-sheaves are having filtrations in the range $[a,a+p+i]$;
\medskip\noindent
its maps (morphisms) are defined in the logical way. Similarly, by not specifying the range of filtrations, we define $\Mm\Mf^{\nabla(tens)}(X)$.
\smallskip
To fully formalize the definition of $\Mm\Mf^{\nabla(tens)}(X)$, we recall (see [Fa1, p. 34-5]) that 
$\Mm\Mf_{[c,c+p-1]}(X)$ is the subcategory of
$$\FF(X):=\times_{(U,\Phi_U)\in\UU(X)} \Mm\Mf^\nabla(U,\Phi_U)$$
whose morphisms are families of morphisms $(m_{U,\Phi_U})_{(U,\Phi_U)\in\UU(X)}$, with each $m_{U,\Phi_U}$ as a morphism of $\Mm\Mf_{[c,c+p-1]}^\nabla(U,\Phi_U)$, such that for any functor $i(U,c,\Phi_U,\Phi_U^1)$ as above we have 
$$i(U,c,\Phi_U,\Phi_U^1)(m_{U,\Phi_U})=m_{U,\Phi_U^1}.$$
So we define $\Mm\Mf^{\nabla(tens)}(X)$ as the tensorial closure of the union of $\Mm\Mf_{[c,c+p-1]}(X)$'s ($c\in\ZZ$) in $\FF(X)$. 
\medskip
{\bf Remarks.} {\bf 1)} Each morphism of $\Mm\Mf^{\nabla(tens)}(X)$ defines a family of morphisms $(m_{U,\Phi_U})_{(U,\Phi_U)\in\UU(X)}$, with each $m_{U,\Phi_U}$ as a morphism of $\Mm\Mf^{\nabla(tens)}(U,\Phi_U)$, such that for any functor $i(U,\Phi_U,\Phi_U^1)$ as in the Lemma the following equality 
$$i(U,\Phi_U,\Phi_U^1)(m_{U,\Phi_U})=m_{U,\Phi_U^1}\leqno (GLUE)$$
makes sense and holds.
\smallskip
We denote by $\Mm\Mf^{\nabla(loc-tens)}(X)$
 the subcategory of $\FF(X)$ formed by all such families of morphisms. We do not stop to study when the natural inclusion of categories
$$\Mm\Mf^{\nabla(tens)}(X)\subset \Mm\Mf^{\nabla(loc-tens)}(X)$$
is an isomorphism.
\smallskip
{\bf 2)} If $X^\wedge$ has a Frobenius lift, then $\Mm\Mf^{\nabla(loc-tens)}(X)$ is naturally a subcategory of $\Mm\Mf^\nabla(X)$.
\smallskip
{\bf 3)} We consider an open cover $(U_i)_{i\in I}$ of $X$. In the definition of $\Mm\Mf^{\nabla(tens)}(X)$ (or of $\Mm\Mf^{\nabla(loc-tens)}(X)$) we can restrict to the subset of $\UU(X)$ whose elements are of the form $(U,\Phi_U)$, with $U$ an affine, open subscheme of $U_i$, for some $i\in I$. 
\smallskip
{\bf 4)} As the category $COH(X)$ is tensorial, following the pattern of F2 we get that for each object of $\Mm\Mf^{\nabla(tens)}(X)$ we can speak about its underlying $\Mo_X$-sheaf (it has the DC property) and about its filtration by direct summands. So the gluing performed by (GLUE) of 1) is just at the level of Frobenius endomorphisms.   
\medskip
{\bf D. Two other variants.} The problem with $\Mm\Mf^{\nabla(tens)}(X)$ is: the category $p-\Mm\Mf^\nabla(X)$ is already well defined (see 2.2.1.1 3)). But there is nothing to guarantee that under the operation of taking truncations mod $p^n$, $n\in\NN$, we get a functor from $p-\Mm\Mf^\nabla(X)$ to $\Mm\Mf^{\nabla(tens)}(X)$. So as a second variant, we define $\Mm\Mf^{\nabla(p+tens)}(X)$ to be the tensorial closure in $\FF(X)$ of the union of $\Mm\Mf^{\nabla(tens)}(X)$ and of the images of the truncation mod $p^m$ functors from $p-\Mm\Mf^\nabla(X)$ to $\FF(X)$, $m\in\NN$. 4) of C still applies to it. Similarly, we define $\Mm\Mf^{\nabla(p+loc-tens)}(X)$.
\medskip
{\bf Exercise.} a) We define $\Mm\Mf^{\nabla(tens)}_{[a,a+p-1]}(X)$ to be the full subcategory of $\Mm\Mf^{\nabla(tens)}(X)$ whose objections are involving filtrations in the range $[a,a+p-1]$. Show that we have a canonical identification  $\Mm\Mf^\nabla_{[a,a+p-1]}(X)=\Mm\Mf^{\nabla(tens)}_{[a,a+p-1]}(X)$. The same holds, with $\Mm\Mf^{\nabla(tens)}(X)$ being replaced by $\Mm\Mf^{\nabla(loc-tens)}(X)$.
\smallskip
b) The category of $p$-divisible objects of $\Mm\Mf^{\nabla(p+tens)}(X)$ (defined following the pattern of 2.2.1 c)) is naturally identified with $p-\Mm\Mf^\nabla(X)$.
\smallskip
c) Define $p-\Mm\Mf^{\nabla(tens)}(X)$. Show that it is a subcategory of $p-\Mm\Mf^\nabla(X)$. 
\smallskip
d) Perform c) as well for $\NN-pro-\Mm\Mf^{\nabla(tens)}(X)$ and for $isog-p-\Mm\Mf^{\nabla(tens)}(X)$.
\medskip 
{\bf Hint}: use the proof of the Lemma of C.
\medskip
{\bf E. Rigidity.} All Fontaine categories involving filtrations (of sheaves) in a finite range are not stable under tensor products and so are not rigid. However, if $X$ is as in 2.2.1 c) and $X^\wedge$ is equipped with a Frobenius lift (resp. is not equipped with a Frobenius lift) and if $FC$ is any one of the following two (resp. four) categories $\Mm\Mf(X)$ and $\Mm\Mf^\nabla(X)$ (resp. $\Mm\Mf^{\nabla(tens)}$, $\Mm\Mf^{\nabla(loc-tens)}$, $\Mm\Mf^{\nabla(p+tens)}$ and $\Mm\Mf^{\nabla(p+loc-tens)}$), then $FC[p^n]$, $p-FC$ and $isog-p-FC$ are rigid. 
\medskip
{\bf F. Fibre functors.} Let $O$ be DVR and let $u\in O$ be a uniformizer. Let $\Mt\Mc$ be a tensorial category over $O/u^nO$. [De5] and [Mi4] are stated over fields; however, the definition of a fibre functor as defined in [De5, 2.8] can be adapted for tensorial categories over (quotients of) $O$ as follows.
\medskip
{\bf Definitions. 1)}  Let $X$ be a flat $O/u^nO$-scheme. A tensorial subcategory $\Ms\Mc$ of $COH(\Mo_X)$ over $O/u^nO$ is said to be a good subcategory on $X$ if it has the property that each object of it, locally in the Zariski topology of $X$ has the DC property. 
\smallskip
{\bf 2)} By an almost fibre (resp. by a fibre) functor of $\Mt\Mc$ we mean an $O/p^nO$-linear functor which respects tensor products, identity elements and the standard constraints A, U, and C (see [De5, 2.7]) from $\Mt\Mc$ to $COH(X_1)$ (resp. to $LF(X_1)$), with $X_1$ an $O/u^nO$-scheme which is not necessarily flat, which is the composite of an $O/p^nO$-linear, exact, tensorial functor 
$$\Mt\Mc\to COH(X)$$ 
factoring through a good subcategory on $X$ (resp. $\Mt\Mc\to LF(X)$), with $X$ a flat $O/u^nO$-scheme, with the natural pull back functor $m^*:COH(X)\to COH(X_1)$ (resp. $m^*:LF(X)\to LF(X_1)$) defined by some $O/u^nO$-morphism $m:X_1\to X$. 
\smallskip
{\bf 3)} $\Mt\Mc$ is called tannakian if it has an almost fibre functor.
\medskip
{\bf Fact.} {\it Any almost fibre functor of $\Mt\Mc$ is exact.}
\medskip
{\bf Proof:} Let $\Ms\Mc$ be as in the Definition. We need to show that the restriction $m^*|\Ms\Mc$ to it of any pull back functor $m^*:COH(X)\to COH(X_1)$, with $m:X_1\to X$ an $O/u^nO$-morphism, is exact. We can assume $X$ is connected. We consider a short exact sequence 
$$0\to \Mo_1\hookrightarrow\Mo_2\twoheadrightarrow \Mo_3\to 0\leqno (SES)$$ 
in $\Ms\Mc$. To show that the complex 
$$0\to m^*(\Mo_1)\to m^*(\Mo_2)\to m^*(\Mo_3)\to 0\leqno (COM)$$ 
is exact, we proceed by induction on $q\in S(1,n)$ such that $u^q$ annihilates $\Mo_2$. It is enough to show that the complex $0\to m^*(\Mo_1)\to m^*(\Mo_2)$ is exact. If $q=1$ then $(SES)$ locally in the Zariski topology of $X$ (or of $X_k$) splits and so $(COM)$ is exact. The passage from $q-1$ to $q$ is argued as follows. We use a second induction on the rank ${\rm rk}(\Mo_2)$ of the locally free $\Mo_{X_{O/uO}}$-module 
$$\sum_{i=1}^q u^{i-1}\Mo_2/u^i\Mo_2.$$ 
We can assume $\Mo_3$ is non-zero. We have a standard short exact sequence 
$$0\to u^{q-1}m^*(\Mo_1)\hookrightarrow m^*(\Mo_1)\twoheadrightarrow m^*(\Mo_1/u^{q-1}\Mo_1)\to 0.$$ 
Denoting by $\Mo_{12}$ the image of $\Mo_1$ in $\Mo_2/u^{q-1}\Mo_2$ and by $\Mo_{12}^\prime$ the kernel of $\Mo_1\to\Mo_{12}$, we get a short exact sequence 
$$0\to\Mo_{12}^\prime\hookrightarrow\Mo_1\twoheadrightarrow\Mo_{12}\to 0.\leqno (SES1)$$ 
As $\Mo_3$ is non-zero, we have ${\rm rk}(\Mo_1)<{\rm rk}(\Mo_2)$. So by the induction statements, the pull back through $m^*$ of $(SES1)$ and of the exact complexes $0\to \Mo_{12}^\prime\hookrightarrow u^{q-1}\Mo_2$ and $0\hookrightarrow \Mo_{12}\to \Mo_1/u^{n-1}\Mo_1$, are exact. These imply that the complex $0\to m^*(\Mo_1)\hookrightarrow m^*(\Mo_2)$ is exact. This ends the proof. 
\medskip
For the language of groupoids we refer to [De5, p. 113-5 and \S 3] and [Mi4, Appendix A]; again we allow the base to be an arbitrary $\ZZ$-algebra $R_{\ZZ}$. For simplicity, as we deal with very explicit groupoids, a groupoid acting on an $R_{\ZZ}$-scheme $X$ is usually denoted just as an $R_{\ZZ}$-scheme. The definition of a representation of a $O/u^nO$-groupoid acting on an $O$-scheme $X$ is as in [De5, 3.3] (so an arbitrary quasi-coherent $\Mo_X$-sheaf $V$ is allowed to be the target --see [De5, p. 114]-- of such a representation). 
\smallskip
Let $X$ and $FC$ be as in E. In all that follows, for simplicity of presentation, we assume $X_k$ is connected. As the only global sections of $\Mo_{X_k}$ fixed by the Frobenius endomorphism of $X_k$ are the ones defined by elements of $\FF_p$, we have 
$${\rm End}_{FC}(X_{W_n(k)}(0))=\ZZ/p^n\ZZ;$$ 
here $X_{W_n(k)}(0)$ is the pull back of $W_n(k)(0)$ to $X$. Accordingly, we view as well $X_{W_n(k)}$ as a $\ZZ/p^n\ZZ$-scheme. 
\medskip
{\bf Fact.} {\it $FC[p^n]$ and all its rigid full subcategories (like $q-unip-FC[p^n]$, $solv-FC[p^n]$, etc.) are tannakian categories over $\ZZ/p^n\ZZ$.} 
\medskip
To see this we just need to point out (cf. the Fact) that there is a canonical almost fibre functor 
$$CAFF_n(FC):FC[p^n]\to COH(X_{W_n(k)});$$ 
it associates to any object of $FC[p^n]$ its underlying $\Mo_{X_{W_n(k)}}$-sheaf (the action on morphism being the logical one). We have the following compatibility property: identifying $COH(X_{W_n(k)})$ with a full subcategory of $COH(X_{W_{n+1}(k)})$, the restriction of $CAFF_{n+1}(FC)$ to $FC[p^n]$ is $CAFF_n(FC)$.
\smallskip
Similarly, if $X$ is affine (resp. is proper over $W(k)$), then $isog-p-FC$ is a tannakian category over $\QQ_p$: the canonical fibre functor 
$$CAFF_{\infty}(FC)$$ 
of it takes values in $LF(X^\wedge\otimes_{\ZZ_p} \QQ_p)$ (resp. takes values, cf. Grothendieck's algebraization theorem, in $LF(X_{B(k)})$). 
\medskip
{\bf G. Pull backs.}
All pull back operations of 2.2.1.3 are functorial, $\ZZ_p$-linear, respect tensor products and (when appropriate) identity objects. Moreover, they do extend to the context of $\Mm\Mf^{\nabla(tens)}(X)$ and of $\Mm\Mf^{\nabla(p+tens)}(X)$ (for instance, cf. 3) and 4) of C).
\medskip
{\bf H. Fontaine's comparison theory.} We assume that $p>2$ and that $X$ is an integral, regular, formally smooth $W(k)$-scheme having a connected special fibre. We also assume that we can speak about Fontaine's comparison theory for objects of $\Mm\Mf_{[c,c+p-2]}^\nabla(U)$, for any affine, open subscheme $U$ of $X$ which is faithfully flat over $W(k)$ and formally \'etale over $W(k)[x_1,...,x_m][{1\over {\prod_{i\in S(1,m)} x_i}}]$, with $m$ equal to the relative dimension of $X$ over $W(k)$ ([Fa1, 2.6] says that $X$ smooth over $W(k)$ would do; this can be easily extended to pro-\'etale covers of smooth $W(k)$-schemes). By this we mean that there is a natural fully faithful, contravariant functor $\DD_U$ from $\Mm\Mf_{[c,c+p-2]}^\nabla(U)$ to the category of representations of $\Gamma_U$ on $\ZZ_p$-modules of finite length which respects the type as defined in [Fa1, p. 37]; here $\Gamma_U$ is a suitable (like in [Fa1, p. 38-40]) quotient of the fundamental group $\pi_1((U^\wedge)_{B(k)},{\eta_U})$, with $\eta_U$ as the generic point of $U^\wedge$. We deduce the existence of a fully faithful, contravariant functor (still denoted by) $\DD_U$ from $p-\Mm\Mf_{[c,c+p-2]}^\nabla(U)$ to the category of representations of $\Gamma_U$ on free $\ZZ_p$-modules of finite rank. Let $\tilde\Gamma_X$ be a suitable fundamental group obtained as in [Fa1, rm. of p. 41] such that we can ``glue" the above functors $\DD_U$, to get a fully faithful, contravariant functor 
$\DD_X$ from $p-\Mm\Mf_{[c,c+p-2]}^\nabla(X)$ to the category of representations of $\tilde\Gamma_X$ on $\ZZ_p$-modules which are free of finite rank. For instance, if $X$ is proper over $W(k)$, then we can take (see [Fa1, 2,6*]) $\tilde\Gamma_X=\pi_1(X_{B(k)},\eta_X)$, with $\eta_X$ as the generic point of $X$.
\smallskip
We consider the subcategory $FCT-isog-p-\Mm\Mf^\nabla(X)$ of $isog-p-\Mm\Mf^\nabla(X)$ which is obtained through the four operations of C but performed starting from the categories $isog-p-\Mm\Mf_{[c,c+p-2]}^\nabla(X)$'s, $c\in\ZZ$. As in C, we can define it in a fancier way as the tensorial closure of the union of the categories $isog-p-\Mm\Mf_{[c,c+p-2]}^\nabla(X)$'s, $c\in\ZZ$, in $isog-p-\Mm\Mf^\nabla(X)$. We have a natural fibre functor 
$$\DD_X:FCT-isog-p-\Mm\Mf^\nabla(X)\to {\rm Rep}(\tilde\Gamma_X;\QQ_p)$$
into the category of $\QQ_p$-linear representations of $\tilde\Gamma_X$:
it is the tensorial extension of the restriction of $\DD$ to the union of the categories $isog-p-\Mm\Mf_{[c,c+p-2]}^\nabla(X)$'s, $c\in\ZZ$ (so it is obtained via the above operations --of taking kernels, cokernels, tensor products and duals-- from $\QQ_p$-representations of $\tilde\Gamma_X$ obtained by making $p$-invertible in the $\ZZ_p$-representations of it associated (via $\DD_X$) to objects of $p-\Mm\Mf_{[c,c+p-2]}^\nabla(X)$. Its existence can be checked as in (the proof of) C, using tensorial closures. 
\medskip
{\bf I. The quasi-solvable context.} The category $isog-unip-p-\Mm\Mf_{[0,0]}^\nabla(X)$ is a full, tensorial subcategory of $isog-q-unip-p-\Mm\Mf_{[0,0]}^\nabla(X)$, as well as this last category is a full, tensorial subcategory of $FCT-isog-p-\Mm\Mf^\nabla(X)$. So in H, we can speak about the restriction of $\DD_X$ to $isog-unip-p-\Mm\Mf_{[0,0]}^\nabla(X)$ or to $isog-q-unip-p-\Mm\Mf_{[0,0]}^\nabla(X)$ as being a fibre functor; we do not know if (or when) we can replace here $unip$ by $solv$. 
\smallskip
It is worth pointing out that:
\medskip
{\bf Fact.} {\it $isog-unip-p-\Mm\Mf^\nabla_{[0,0]}(X)$ can be naturally identified with the category of unipotent $F$-isocrystals on $X_k$.}
\medskip
{\bf Proof:} We just need to show that any unipotent $F$-isocrystal ${\got C}$ on $X_k$ is obtained from an object of $unip-p-\Mm\Mf^\nabla_{[0,0]}(X)$ by making $p$-invertible. Let $n$ be the rank of ${\got C}$. Let ${\got C}_1\subset {\got C}_2\subset ...\subset{\got C}_n={\got C}$ be a filtration of it by $F$-subisocrystals such that $\forall i\in S(1,n-1)$, ${\got C}_{i+1}/{\got C}_i$ is a trivial $F$-isocrystal of rank 1. Let ${\got C}^{\rm int}$ be an $F$-crystal on $X_k$ having ${\got C}$ as its $F$-isocrystal. We view it as a $p$-divisible object of $\Mm\Mf_{[0,0]}^\nabla(X)$ (this can be checked using Teichm\"uller lifts: the slopes of any unipotent $F$-isocrystal over a perfect field $k_1$ are all zero and so any $\sg_{k_1}$-crystal producing it can be viewed as a $p$-divisible object of $\Mm\Mf_{[0,0]}(W(k_1))$). But now, based on this and on Fact of 2.2.1.1 7), we get directly that the intersection of ${\got C}^{\rm int}$ with ${\got C}_i$ taken inside ${\got C}$ is a $p$-divisible subobject ${\got C}^{\rm int}_i$ of ${\got C}^{\rm int}$ and that the quotient ${\got C}_{i}^{\rm int}/{\got C}^{\rm int}_{i-1}$ is a trivial object of $p-\Mm\Mf^\nabla_{[0,0]}(X)$, $\forall i\in S(1,n)$. This ends the proof. 
\medskip
{\bf J. $\pi_1$-groupoids.} What follows is not used at all in \S2-4; however, we hope to come back soon to its important (unexploited) ideas. We refer to E. We assume $X_k$ is connected. We consider the $\ZZ/p^n\ZZ$-groupoid 
$$\pi_1(FC,n)$$ 
acting on $X_{W_n(k)}$ and defined (as in [De5, 1.11]) by automorphisms of the almost fibre functor $CAFF_n(FC)$. 
\smallskip
The case $n=1$ is simpler as $CAFF_1(FC)$ is in fact a fibre functor; so (cf. [De5, 1.12]) $\pi_1(FC,1)$ is faithfully flat over $X_k\times_{\FF_p} X_k$ and the tannakian category $FC[p]$ is isomorphic (in the tannakian sense) to the tannakian category of representations of $\pi_1(FC,1)$. Unfortunately loc. cit. deals only with the case of a base scheme $BS$ which is the spectrum of a field and with (see [De5, 2.8]) locally free sheaves of locally finite rank on some $BS$-scheme. So we do not stop to check here that similar properties hold for $n\ge 2$. As $FC[p^{n-1}]$ is a full subcategory of $FC[p^n]$ as well as a quotient category of $FC[p^n]$ (via the functor truncation mod $p^{n-1}$), we have a canonical identification 
$$\pi_1(FC,n)_{W_{n-1}(k)}=\pi_1(FC,n-1),$$ 
$\forall n\in\NN$, $n\ge 2$. So it makes sense to speak about the projective limit of $\pi_1(FC,n)$, $n\in\NN$; it is a $\ZZ_p$-groupoid 
$$\pi_1(FC,n\to\infty)$$ 
on the $p$-adic formal scheme $X^\wedge$.
\smallskip
Similarly, in the $isog-p$ context, we will not try to work with groupoids over analytic spaces: whenever we refer to a $\QQ_p$-context we assume $X$ is either an affine or a proper $W(k)$-scheme. Let now $X$ be affine (resp. proper $W(k)$-scheme). Starting from $isog-p-FC$, we similarly define 
$$\pi_1(FC,\infty)$$ 
to be the $\QQ_p$-groupoid acting on $X^\wedge[{1\over p}]$ (resp. on $X_{B(k)}$) defined by the fibre functor $CAFF_{\infty}(FC)$. Loc. cit. implies $\pi_1(FC,\infty)$ is faithfully flat over $X^\wedge[{1\over p}]\times_{B(k)} X^\wedge[{1\over p}]$ (resp. over $X_{B(k)}\times_{B(k)} X_{B(k)}$). 
\smallskip
We expect that in many cases $\pi_1(FC,\infty)$ is ``algebraizable" in a natural sense. To exemplify here what we mean by natural sense, we restrict t the case $X$ is affine. In such a case $\pi_1(FC,n\to\infty)$ is algebraizable, i.e. it is obtained naturally from a $\ZZ_p$-groupoid $\pi_1(FC,alg)$ on $X^\wedge$. Warning: $\pi_1(FC,alg)$ is not necessarily flat over $\ZZ_p$. 
\smallskip
However, denoting by $\pi_1(FC,alg)_{\QQ_p}$ the natural pull back of $\pi_1(FC,alg)$ to a $\QQ_p$-groupoid on $X^\wedge[{1\over p}]$, we have a natural homomorphism
$$m_{FC}:\pi_1(FC,alg)_{\QQ_p}\to \pi_1(FC,\infty);$$
it is defined by the natural ``inverting $p$" functor $a-p-FC\to isog-p-FC$ which passes to isogeny classes. 
See 7) below for a study of it. Here we just point out that if $X$ is a proper $W(k)$-scheme, one should similarly be able to construct a natural algebraization of $\pi_1(FC,n\to\infty)$ and a natural homomorphism $m_{FC}$. Warning: in this paragraph, in case $X^\wedge$ is not equipped with a Frobenius lift, we take $FC=\Mm\Mf^{\nabla(p+tens)}(X)$.
\smallskip
We have lots of variants and properties. 1) to 7) below are just samples of what one can do.
\smallskip
{\bf 1) The unip and solv variants.} Similarly we define $\pi_1^{\rm unip}(FC,n)$, $\pi_1^{\rm unip}(FC,\infty)$, $\pi_1^{\rm solv}(FC,n)$ and $\pi_1^{\rm solv}(FC,\infty)$, by working in the $unip$ or $solv$ context. For instance, if $X$ is proper over $W(k)$, $\pi_1^{\rm unip}(FC,\infty)$ is the groupoid acting on $X_{B(k)}$ defined by automorphisms of the restriction of $CAFF_{\infty}(FC)$ to $isog-unip-p-FC$. 
\smallskip
In general, we have natural homomorphisms of groupoids on $X_{W_n(k)}$
$$\pi_1(FC,n)\to \pi^{\rm solv}_1(FC,n)\to \pi^{\rm unip}_1(FC,n)$$
and in case $X$ is an affine or a proper $W(k)$-scheme 
$$\pi_1(FC,\infty)\to \pi^{\rm solv}_1(FC,\infty)\to \pi^{\rm unip}_1(FC,\infty).$$
\indent
{\bf 2) Cocharacters.} One of the nice features: whenever we consider a $W(k_1)$-valued point of $X$, with $k_1$ a perfect field containing $k$, [Wi] allows us to define cocharacters of the resulting pull back groupoids. For instance (see 2.2.14.1 below), we get a cocharacter of the $\ZZ/p^n\ZZ$-groupoid on $W_n(k_1)$ obtained by pulling back $\pi_1(FC,n)$ to a $\ZZ/p^n\ZZ$-groupoid on $W_n(k_1)$. Such pull backs are extremely useful in studying the flatness of different $\ZZ/p^n\ZZ$-groupoids we defined. 
\smallskip
{\bf 3) Fontaine's comparison theory variant.} Let now $X$ be as in H. Let $\pi_1(FCT(X))$ (resp. $\pi_1^{\rm unip}(FCT(X))$) be the group scheme over $\QQ_p$ defined as the automorphism of the fibre functor $\DD_X$ of H (resp. of its restriction to $isog-unip-p-\Mm\Mf_{[0,0]}(X)$). Similarly, let $\pi_1(\tilde\Gamma_X)$ be the group scheme over $\QQ_p$ which is defined by automorphisms of the tautological fibre functor of ${\rm Rep}(\tilde\Gamma_X;\QQ_p)$. We get natural homomorphisms 
$$\pi_1(\tilde\Gamma_X)\to\pi_1(FCT(X))\to\pi_1^{\rm unip}(FCT(X)).$$
\indent
{\bf 4) Pull backs.} If $m:X_1\to X$ is a morphism as in 2.2.1.3, then the above constructions of groupoids are compatible with pull backs. Just one example. We take $FC=\Mm\Mf^{\nabla(p+tens)}(X)$ and $FC_1=\Mm\Mf^{\nabla(p+tens)}(X_1)$. Then, corresponding to the pull back operation of 2.2.1.3 (see also G) we have a natural morphism 
$$\pi_1(FC_1,n)\to m^*(\pi_1(FC,n))$$
of $\ZZ/p^n\ZZ$-groupoids on $X_{1W_n(k)}$.
\smallskip
{\bf 5) The extra feature of the case $k=\FF_p$.} If $k=\FF_p$ and if $z\in X(\ZZ_p)$, then the pull back of $\pi_1(FC,n)$ via $z$, is a group scheme over $\ZZ/p^n\ZZ$ equipped naturally with a $\ZZ/p^n\ZZ$-valued point defined naturally by the Frobenius automorphisms $\vph_{\ZZ}$ of $\Mm\Mf(\ZZ_p)[p^n])$ (they are obtained as in 2.2.1 c) starting from direct sum decompositions produced by canonical split cocharacters; see 2.2.14.1 below for the mod $p^n$ version of such cocharacters).
\smallskip
{\bf 6) Some extra features of the $\Mm\Mf_{[0,0]}(*)$ case.} We just include two samples. The Frobenius endomorphism of $X_k$ gives birth, via the interpretation of Fact of I, to a ``Frobenius" endomorphism of $\pi_1^{\rm unip}(\Mm\Mf^{\nabla(p+tens)}(X),\infty)$. 
\smallskip
The previous paragraph is just a property of the $[0,0]$ range. If $FC$ is $\Mm\Mf(X)$ (resp. any one of the Fontaine categories of E involving connections) let $FC_0[p^n]:=\Mm\Mf_{[0,0]}(X)[p^n]$ (resp. let $FC_0[p^n]:=\Mm\Mf_{[0,0]}^\nabla(X)[p^n]$). $FC_0[p^n]$ is a full tensorial subcategory of $FC[p^n]$. So the restriction of $CAFF_n(FC)$ to it is an almost fibre functor: we denote by $\pi_1^0(FC,n)$ the $\ZZ/p^n\ZZ$-groupoid on $X_{W_n(k)}$ it defines. It is naturally equipped with a Frobenius lift. Moreover, we have a natural homomorphism of $\ZZ/p^n\ZZ$-groupoids
$$\pi_1(FC,n)\to\pi_1^0(FC,n).$$
\smallskip
{\bf 7) The $\ZZ_p$-context.} If $X^\wedge$ is not equipped with a Frobenius lift, then below we work with $FC=\Mm\Mf^{\nabla(p+tens)}(X)$. We have a variant $\pi_1(FC(ESS),n)$ of $\pi_1(FC,n)$, obtained by replacing $FC[p^n]$ with the tensorial closure of the image $FC(ESS)[p^n]$ of $p-FC$ in $FC[p^n]$ via the truncation mod $p^n$ functor. Here (ESS) refers to: essentialized. We have a similar natural identification (with $n\ge 2$)
$$\pi_1(FC(ESS),n)_{W_{n-1}(k)}=\pi_1(FC(ESS),n-1)$$
and a natural homomorphism (with $n\ge 1$)
$$\pi_1(FC,n)\to\pi_1(FC(ESS),n)$$
compatible with restrictions (from mod $p^n$ to mod $p^{n-1}$, with $n\ge 2$). So passing to limit, we get a homomorphism of $\ZZ_p$-groupoids
$$m_{FC,n\to\infty}:\pi_1(FC,n\to\infty)\to\pi_1(FC(ESS),n\to\infty).$$
If $X$ is affine, then $\pi_1(FC(ESS),n\to\infty)$ is as well algebraizable, i.e. is defined naturally by a $\ZZ_p$-groupoid $\pi_1(FC(ESS),alg)$ on $X^\wedge$. Moreover, $m_{FC,n\to\infty}$ is defined naturally by a homomorphism 
$$m_{FC,alg}:\pi_1(FC,alg)\to\pi_1(FC(ESS),alg).$$
It seems to us that the pull back of it to a $\QQ_p$-groupoid is an isomorphism. We have a natural factorization of $m_{FC}$ through ${m_{FC,alg}}_{\QQ_p}$: we denote it by $m_{FC(ESS)}$. We have the following philosophy:
\medskip
{\bf Ph.} {\it $\pi_1(FC(ESS),alg)$ is the ``standard $\ZZ_p$" part, while $\pi_1(FC,alg)$ ``captures" the ``exotic $p$-torsion" as well.}
\medskip
{\bf Proposition.} {\it We assume $X$ is proper and geometrically connected over $W(k)$. Then there is a uniquely determined $\ZZ_p$-groupoid $\pi_1(FC(ESS),alg)$ on $X$ having the properties:
\medskip
{\bf a)} Its $p$-adic completion is canonically identified with $\pi_1(FC(ESS),n\to\infty)$.
\smallskip
{\bf b)} Its pull back to a $\QQ_p$-groupoid on $X_{B(k)}$ is canonically identified with $\pi_1(FC,\infty)$.
\smallskip
{\bf c)} For any affine, open subscheme $U$ of $X$ having a non-empty special fibre, the pull back to $\Mu:=U^\wedge\times_{W(k)} B(k)$ of the identifications of a) and b) are compatible with the canonical homomorphism ${m_{FC(ESS)}}_{\Mu}:\pi_1(FC(ESS),n\to\infty)_{\Mu}\to\pi_1(FC,\infty)_{\Mu}.$}
\medskip
{\bf Proof:} In essence, this is obvious: so we go quickly through the main details. We use freely Grothendieck's algebraization theorem. We need the following obvious Fact. 
\medskip
{\bf Fact.} {\it Any finite subcategory of $p-isog-FC$ is obtained from a finite category of $p-FC$ by making $p$-invertible.}
\medskip
We recall that we adopted (see 2.1) the point of view that $FC$ is a small category. To avoid using the notion of universes (which we do not mention any way), the reader if desires, for what follows can work with a skeleton of $p-FC$. We consider the set $S(FC)$ formed by finite subcategories of $p-FC$. We order it under the relation of inclusion: it is a filtered set. For each $\al\in S(FC)$, we consider the functor
$$\Mf_{\al}:\al\to LF(X)$$
which associates to each object of $\al$ the algebraization of its underlying $\Mo_{X^\wedge}$-module and which acts on morphisms in the logical way (via algebraizations). Let $G_{\al}$ be the $\ZZ_p$-groupoid on $X$ parameterizing automorphisms of $\Mf_{\al}$. 
\smallskip
If $\al_1$, $\al_2\in S(FC)$ are such that $\al_1\subset\al_2$, then the restriction of $\Mf_{\al_2}$ to $\al_1$ is $\Mf_{\al_1}$. We deduce that we have a natural homomorphism of $\ZZ_p$-groupoids
$$f_{\al_2,\al_1}:G_{\al_2}\to G_{\al_1}.$$     
We take 
$$\pi_1(FC(ESS),alg):={\rm proj.}\, {\rm lim.}_{\al\in S(FC)} G_{\al}$$
via $f_{\al_2,\al_1}$'s. The Fact implies that its fibre over $\QQ_p$ is $\pi_1(FC,\infty)$. Similarly, the natural morphism
$$\pi_1(FC(ESS),n)\to \pi_1(FC(ESS),alg)_{W_n(k)}={\rm proj.}\, {\rm lim.}_{\al\in S(FC)} G_{\al W_n(k)}$$
is an isomorphism. The equality part is due to the fact that tensor products commute with inductive limits. The fact that the arrow is an isomorphism is a consequence of the fact that $FC(ESS)[p^n]$ is the tensorial closure of the image of $TR_n(FC)$ (of 2.2.1.6 7)): from the point of view of groupoids of automorphisms of almost fibre functors, the operation of taking the tensorial closure is irrelevant (cf. the exactitude of almost fibre functors and the list of 4 operations of C). This takes care of a) and b). c) can be checked locally; it boils down to: if $R$ is a flat $\ZZ_p$-algebra, then the only $\ZZ_p$-subalgebra $R_1$ of $R[{1\over p}]$ such that $R_1[{1\over p}]=R[{1\over p}]$ and $R_1^\wedge=R^\wedge$, is $R$ itself. This is a consequence of the fact that $R=R[{1\over p}]\cap R^\wedge$, the intersection being taken inside $R^\wedge[{1\over p}]$. This ends the proof.
\medskip
{\bf K. Some comparisons.} We will be brief. Whatever was done above in the context of Fontaine categories can be redone in the de Rham context (of $X$). One defines different tannakian categories of vector bundles on $X_{W_n(k)}$ (or on $X$ or $X_{B(k)}$, etc.) which are or are not endowed with some specific type of connections. So we get, as above, plenty of $\pi_1^{dR}$ groupoids. It is natural to try to ``compare" them. The unipotent case is by far the easiest, cf. the following general Fact.
\medskip
{\bf Fact.} {\it Let $O$ and $u$ be as in F. Let $K:=O[{1\over u}]$. Let $X$ be a proper, flat $O$-scheme whose ring of global sections is $O$. Then any pair $(\Mf,\nabla)$, with $\Mf$ a locally free $\Mo_{X_K}$-sheaf of finite rank $r$ and with $\nabla$ a connection on it, which is obtained from the pair $(\Mo_{X_K},d)$ (with $d$ as the connection annihilating $1$) via short exact sequences, is obtained, up to isomorphism, from a similar pair $(\Mf(X),\nabla)$ on $X$ by making $u$ invertible.}
\medskip
{\bf Proof:} We proceed by induction on $r$. The case $r=1$ is trivial. We now consider a short exact sequence 
$$0\to (\Mf_1,\nabla_1)\hookrightarrow (\Mf,\nabla)\twoheadrightarrow (\Mo_{X_K},d)\to 0,\leqno (SES)$$ 
such that $(\Mf_1,\nabla_1)=(\Mf_1(X),\nabla_1)_K$, with $(\Mf_1(X),\nabla_1)$ a pair obtained from $(\Mo_X,d)$ via short exact sequences. $M:={\rm Ext}^1(\Mo_X,\Mf_1(X))$ is a finitely generated $O$-module. Let $\gamma\in M[{1\over u}]$ be the class of (SES) viewed without connections. Considering the pull back of (SES) via the automorphism of $(\Mo_{X_K},d)$ defined by $u^m$, with $m\in\NN$ and big enough, we can assume that $\gamma\in M$ and that $\nabla$ takes the $\Mo_X$-sheaf $\Mf(X)$ (which is the extension of $\Mo_X$ through $\Mf_1(X)$) defined by $\gamma$ into $\Mf(X)\otimes_{\Mo_X} \Om_{X/O}$. This ends the proof.
\smallskip
The tannakian category $isog-unip-p-\Mm\Mf(X)$ used here is different from the one used in [Shi1-2]: loc. cit. uses unipotent isocrystals and not unipotent $F$-isocrystals, i.e (cf. the below Lemma) uses $isog-p-UNIP(X_k)$. 
\medskip
{\bf Lemma.} {\it Any object ${\got C}$ of $p-UNIP(X_k)$ is a unipotent crystal, i.e is obtained from short exact sequences of direct summand type from objects of $p-UNIP(X_k)$ which are crystals on $X_k$ in invertible sheaves.}
\medskip
{\bf Proof (slightly sketched):} Let ${\got D}_n$ be the maximal subobject of ${\got C}/p^n{\got C}$ which is obtained by pulling back a $W_n(k)$-module of finite type (and endowed with the trivial connection); its existence is implied by the fact that $UNIP(X_k)[p^n]$ is an abelian category having direct sums and of the fact that any quotient object of any such pull back is as well a pull back of a $W_n(k)$-module of finite type. By very definitions, the $\Mo_{X_{W_n(k)}}$-sheaf of the evaluation of ${\got D}_n$ at $X_{W_n(k)}$ is non-trivial. We consider the maximal subobject ${\got E}_n$ of ${\got D}_n$ with the property that the $\Mo_{X_{W_n(k)}}$-sheaf of its evaluation at $X_{W_n(k)}$ is free. For $m\in\NN$, we have a natural monomorphism $i_{n,m}:{\got E}_{n+m}[p^n]\hookrightarrow {\got E}_n$. Let ${\got F}_n$ be the intersection of the images of all $i_{n,m}$, $m\in\NN$. By reasons of ranks, we get a natural isomorphism $j_{n,1}:{\got F}_{n+1}[p^n]\hookrightarrow {\got F}_n$. The projective limit of $j_{n,1}$'s define a subobject ${\got F}$ of ${\got C}$. It is easy to see that ${\got F}$ is the pull back of a free $W(k)$-module of positive rank and that ${\got C}/{\got F}$ is an object of $p-UNIP(X_k)$. So the Lemma follows by induction on the rank of ${\got C}$.
\medskip
We can perform B to J in the context of $UNIP(X_k)$. As in E (resp. in F) we get that $UNIP(X_k)[p^n]$, $p-UNIP(X_k)$, $isog-p-UNIP(X_k)$ are rigid (resp. that $UNIP(X_k)[p^n]$ is tannakian over $W_n(k)$; its canonical almost fibre functor $CAFF_n(UNIP(X_k))$ is defined by evaluations at $X_{W_n(k)}$). As in G we speak about pull backs of objects of $UNIP(X_k)$, of $p-UNIP(X_k)$, etc., under morphisms $X_1\to X$. b) of Fact of A holds as well for $unip-\Mm\Mf^\nabla(X)$. So very often, for unipotent contexts we do not need to make distinction between some ``essentialized" or not (see 7) of J). However, to simplify the things below we work in ``essentialized" contexts. 
\smallskip
From now on we assume $X$ is proper. By replacing in the Proposition of J, $FC(ESS)$ by $unip-FC(ESS)$ and by $UNIP(X_k)(ESS)$ we speak respectively about $\pi_1^{\rm unip}(FC(ESS),alg)$ and about the unipotent, crystalline fundamental $W(k)$-groupoid $\pi_1^{{\rm unip},{\rm crys}}(X_k,X)$ on $X$ of $X_k$; so 
$$\pi_1^{{\rm unip},{\rm crys}}(X_k,X_{W_n(k)}):=\pi_1^{{\rm unip},{\rm crys}}(X_k,W(k))_{W_n(k)}$$ 
is the groupoid of automorphisms of the canonical almost fibre functor obtained by restricting $CAFF_n(UNIP(X_k))$ to $UNIP(X_k)(ESS)[p^n]$. The combination of the Fact and of the Lemma implies that its fibre over $B(k)$ is the same as the $B(k)$-groupoid on $X_{B(k)}$ we get using vector bundles endowed with a flat connection and which are unipotent in the natural sense.
\smallskip
We have a natural forgetful functor (which forgets the Frobenius endomorphisms)
$$unip-p-\Mm\Mf(X)\to p-UNIP(X_k)$$
and so, at the level of $W(k)$-groupoids we get a homomorphism
$$m_{\rm unip}:\pi_1^{{\rm unip},{\rm crys}}(X_k,X)\to\pi_1^{\rm unip}(\Mm\Mf^{\nabla(p+tens)}(X)(ESS),alg)_{W(k)}.$$
There are situations when $m_{\rm unip}$ is an isomorphism; however, in general we have to look at $\pi_1^{\rm unip}(\Mm\Mf^{\nabla(p+tens)}(X),alg)$ as a quotient (often trivial) of $\pi_1^{{\rm unip},{\rm crys}}(X_k,X)$. 
\smallskip
If $X^\wedge$ is not (resp. is) equipped with a Frobenius lift, then the fibre product 
$$\pi_1^F(X)$$ 
of $m_{\rm unip}$ with the natural homomorphism 
$$\pi_1(\Mm\Mf^{\nabla(p+tens)}(X)(ESS),alg)_{W(k)}\to \pi_1^{\rm unip}(\Mm\Mf^{\nabla(p+tens)}(X)(ESS),alg)_{W(k)}$$ 
(resp. with
$\pi_1(\Mm\Mf^{\nabla}(X)(ESS),alg)_{W(k)}\to \pi_1^{\rm unip}(\Mm\Mf^{\nabla(p+tens)}(X)(ESS),alg)_{W(k)}$) 
is called Fontaine's crystalline fundamental $W(k)$-groupoid of $X$; it is on $X$.
\medskip
{\bf 2.2.4.1. Complements.} {\bf 1)} $\pi_1^{{\rm unip},{\rm crys}}(X_k,X)$ is a $W(k)$-groupoid on $X$ and so it does depend on $X$. On the other hand $UNIP(X_k)$ depends only on $X_k$. To reconcile these two things, we can proceed in two ways. In the first way, we start with an arbitrary regular, formally smooth $k$-scheme $Y_k$. We define the unipotent crystalline fundamental crystal in $W(k)$-groupoids
$$\pi_1:=\pi_1^{{\rm unip},{\rm crys}}(Y_k)$$ 
as follows. We work with $CRIS(Y_k/W(k))$. If $(Z,i,\dl)$ is an arbitrary object of it (i.e. a thickening of some $Y_k$-scheme $S_k$; so $i:S_k\hookrightarrow Z$ is a closed embedding and $\dl$ is a PD-structure on the ideal sheaf defining $i$), to define $\pi_1(Z,i,\dl)$ we can assume that $Z$ is affine, that $p^n=0$ in $Z$, and that there is a $W(k)$-morphism $m_{Z,\tilde Y}$ from $Z$ to an affine, regular, formally smooth lift $\tilde Y$ of an open subscheme $\tilde Y_k$ of $Y_k$ whose restriction to $S_k$ is the $k$-morphism $S_k\to Y_k$. Then we take 
$$\pi_1(Z,i,\dl):=m_{Z,\tilde Y}^*(\pi_1^{{\rm unip},{\rm crys}}(\tilde Y_k,\tilde Y_{W_n(k)})).$$ 
As the objects of $UNIP(Y_k)$ are crystals in coherent sheaves, it is trivial to check that $\pi_1$ is indeed a crystal in $W(k)$-groupoids. For this definition of $\pi_1$, it is enough to have $X_k$ such that locally in the Zariski topology, $X_k$ has flat lifts to $W(k)$ which, when viewed naturally as thickenings of $X_k$, are final objects of $CRIS(X_k/W(k))$; in particular, the $\pi_1$ crystal in groupoids can be defined for arbitrary perfect $k$-schemes.
\smallskip
For the second way, we assume $Y_k$ is connected. Then the pull back of $\pi_1$ via any $y\in Y_k(k_1)$, with $k_1$ a perfect field containing $k$, is a crystal on ${\rm Spec}(k_1)$ in group schemes. So it can be identified with an affine group scheme $\pi_1(y)={\rm Spec}(R_y)$ over $W(k_1)$ (cf. also the Proposition of J; we view ${\rm Spec}(W(k_1))$ as a proper scheme over itself). We will not stop to check if (or when) $\pi_1(y)$ does not depend on $y\in X(W(k_1))$. It is well known that the generic fibre $\pi_1(y)_{B(k_1)}$ does not depend on $y$. Argument:  $\pi_1(y)_{B(k_1)}$ is a pro-unipotent group scheme and so $H^1(\Gamma_k,\pi_1(y)_{B(k_1)}(\overline{B(k)}))=0$; on the other hand [De5, 1.12] and [Mi1, Cor. A.9] tell us that $\pi_1(y_1)_{B(k_1)}$ is an inner form of $\pi_1(y)_{B(k_1)}$, $\forall y_1\in Y_k(k_1)$. However, this seems not to be enough to conclude that $\pi_1(y)$ does not depend on $y$. 
\smallskip
Denoting, $I_y:={\rm Ker}(R_y\to R_y\otimes_{W(k_1)} B(k_1))$), the group scheme naturally defined by ${\rm Spec}(R_y/I_y)$ is referred as the flat part $\pi_1^{\rm flat}(y)$ of $\pi_1(y)$.
We do not stop to check when $\pi_1^{\rm flat}(y)$ is of finite type (over $W(k_1)$). Even if its generic fibre is of finite type (over $B(k_1)$), we can not conclude from the proof of Proposition of J that it is of finite type; however, referring to the Fact of 2.2.4 K, we can always choose $(\Mf(X),\nabla)$ to be maximal in some sense and this is our guarantee that at least in many cases $\pi_1^{\rm flat}(y)$ is of finite type. However, $\pi_1$ is only one of the numerous possible variants: we can always restrict our attention to subcategories of $p-UNIP(X_k)$ of interest; for instance, referring to 2.2.4 J, we can work as well with some $G_{\al}$ having the property that its generic fibre is the group scheme (assumed to be of finite type) $\pi_1(y)_{B(k_1)}$.
\smallskip
{\bf 2)} We refer to 7) of  2.2.4 J. The study of $m_{FC,alg}$ is related to the study of objects of $FC$ which have lifts. If $0\to\Mo_1\hookrightarrow\Mo_2\twoheadrightarrow\Mo_3\to 0$ is a short exact sequence in $FC$, b) of 2.2.4 A points out that often $\Mo_2$ has a lift if $\Mo_1$ and $\Mo_3$ do have. In general, $\Mo_1$ has a lift if $\Mo_2$ does. Using this one gets that often any morphism between objects of $FC$ having a lift is the cokernel of a commutative (square) diagram in $p-FC$. Still this is not enough to get that $m_{FC,alg}$ has a section. The problem is with pairs of morphisms of $FC$ of the form
$$\Mo_1\rightarrow\Mo\leftarrow\Mo_2$$
where $\Mo_1$ and $\Mo_2$ have lifts while $\Mo$ does not. Usually we overcome this problem by considering only objects of $FC$ which are self dual. Briefly, we denote by $SD-FC$ the category:
\medskip
-- whose objects are pairs $(\Mo,f)$, with $\Mo$ an object of $FC$ and with $f:\Mo\tilde\to\Mo^t$ an isomorphism of $FC$;
\smallskip
-- whose morphisms are the logical ones (if $(\Mo_1,f_1)$ is another object of $(\Mo,f)$, then the morphisms from $(\Mo,f)$ to $(\Mo_1,f_1)$ are given by morphisms $h:\Mo\to\Mo_1$ of $FC$ such that we have $f_1\circ h=h^t\circ f$, with $h^t:\Mo_1\to\Mo$ as the dual of $h$).
\medskip
It is a tensorial, $\ZZ_p$-linear category. Similarly, we get $SD-FC(ESS)$. $m_{FC,alg}$ has an analogue in the context of $SD-FC(ESS)$; it can be checked that it has a section. This points out that we can speak about self dual (or orthogonal or symplectic) fundamental groupoids (cf. also 2.2.23 below).
\smallskip
{\bf 3)} All above can be redone over finite, ramified, DVR extensions of $W(k)$ of index of ramification at most $p-1$ (cf. also 2.2.1.4.3). To us, based on the Fact of I and of a logical ramified version of it, the comparisons in the unipotent $F$-context (of $isog-unip-p-\Mm\Mf(X)$) or of $FCT-isog-p-\Mm\Mf^\nabla(X)$ are all a consequence of Fontaine's comparison theory of [Fa2, th. 5*]. So to us the interesting cases of comparisons are the ones involving $\Mm\Mf^{\nabla(p+tens)}(X)$ or $p-\Mm\Mf^\nabla(X)$ or an index of ramification which is greater than $p-1$. We hope to come back to them. We end by mentioning that the whole of 2.2.4 can be adapted to the log-smooth context (see [Shi1-2] in connection to $isog-p-UNIP(X_k)$).
\smallskip
{\bf 4)} Warning: in the almost $p$-divisible context (for instance, of $a-p-FC$), though we have an identity object, we do not have duals; so, working with EDC instead of DC, we do not get useful extensions of 2.2.4 F. However, we can work with almost $p$-divisible objects in order to construct directly different algebraizations of ($W(k)$-) $\ZZ_p$-groupoids; for instance, we have a variant of the Proposition of 2.2.4 J in the context of $\pi_1(FC,n\to\infty)$: we just need to replace in its proof $p-FC$ by $a-p-FC$.
\smallskip
{\bf 5)} We refer to 2.2.4 C. We have a variant 
$$\Mm\Mf^{\nabla(big-tens)}(X)$$ 
(resp. $\Mm\Mf^{\nabla(big-p+tens)}(X)$) of $\Mm\Mf^{\nabla(tens)}(X)$ (resp. of $\Mm\Mf^{\nabla(p+tens)}(X)$) obtained by allowing besides the mentioned 4 operations (of 2.2.4 C) an extra one:
\medskip
-- taking subobjects compatible with all gluings (see (GLUE) of C) of objects obtained naturally by changes of Frobenius lifts of the $p$-adic completion of open, affine subschemes of $X$.
\medskip
Warning: we do not know when we can define pull backs for these variants. However, see 1.15.1 and 2.2.20.1 9) for their usage. The same apply in the context of $\nabla(loc-tens)$.
\smallskip
{\bf 6)} All of 2.2.4 can redone in the topologically (q-) solvable (or unipotent) context. Though it might look premature to say, it seems to us that such topological contexts are right contexts for making important steps towards the understanding of ``the whole" crystalline cohomology group of $X_k$'s as above (and not only of their unipotent counterparts). 
\medskip
{\bf 2.2.5. Definitions.} Let $X$ be a connected scheme.
\smallskip
{\bf 1)} A Shimura group pair over $X$ is a pair $(G,[\mu])$, with $G$ a reductive group over $X$ and with $[\mu]$ a $G(X^\prime)$-conjugacy class of cocharacters $\mu:\GG_m\to G_{X^\prime}$ defined over a connected \'etale cover $X^\prime$ of $X$, such that the induced (via inner conjugation) action of $\GG_m$ on ${\rm Lie}(G_{X^\prime})$ is defined (warning!) precisely
by the characters: the null, the identical and the inverse of the identical character of $\GG_m$. So $\mu$ is automatically a closed embedding. If moreover $G$ is an adjoint group, then we speak about a Shimura adjoint group pair.
\smallskip
{\bf 2)} Two such Shimura group pairs $(G,[\mu_1])$ and $(G,[\mu_2])$ over $X$ (with
$\mu_i:\GG_m\to G_{X_i}$ defined over a connected \'etale cover $X_i$ of $X$, $i=\overline{1,2}$) are said to be identical if there is a connected \'etale cover $X^\prime$ of $X$, dominating $X_1$ and $X_2$ and over which $\mu_1$ and $\mu_2$ define the same $G(X^\prime)$-conjugacy class of cocharacters of $G_{X^\prime}$. 
\smallskip
{\bf 3)} For any such Shimura group pair $(G,[\mu])$ over $X$, we denote by $[{\rm Lie}(\mu)]$ the $G^{\rm ad}(X^\prime)$-conjugacy class of the (rank one) direct summand $d\mu({\rm Lie}(\GG_m))$ of ${\rm Lie}(G_{X^\prime})$. We refer to the pair  $({\rm Lie}(G),[{\rm Lie}(\mu)])$ as the Shimura Lie pair (over $X$) attached to the Shimura group pair $(G,[\mu])$. 
\smallskip
{\bf 4)}
Two such Shimura Lie pairs $({\rm Lie}(G_1),[{\rm Lie}(\mu_1)])$ and $({\rm Lie}(G_2),[{\rm Lie}(\mu_2)])$ over $X$ (with $\mu_i$ as a cocharacter of $G_{X_i}$, for $X_i$ a connected \'etale cover of $X$), are said to be isomorphic if:
\medskip
a) there is an isomorphism of Lie algebras $\rho:{\rm Lie}(G_1)\tilde\to {\rm Lie}(G_2)$, taking the kernel of the natural Lie homomorphism (between locally free $\Mo_X$-sheaves endowed with a Lie structure) $q_1:{\rm Lie}(G_1)\to {\rm Lie}(G_1^{\rm ad})$ onto the kernel of the similarly defined natural Lie homomorphism $q_2:{\rm Lie}(G_2)\to {\rm Lie}(G_2^{\rm ad})$, and such that the resulting isomorphism ${\rm Im}(q_1)\tilde\to {\rm Im}(q_2)$ is obtained naturally from an isomorphism ${\rm Lie}(G_1^{\rm ad})\tilde\to {\rm Lie}(G_2^{\rm ad})$ induced by an isomorphism $G_1^{\rm ad}\tilde\to G_2^{\rm ad}$;
\smallskip
b) over a connected \'etale cover $X^\prime$ of $X$, dominating $X_1$ and $X_2$, $\rho$ takes $[{\rm Lie}({\mu_1}_{X^\prime})]$ onto $[{\rm Lie}({\mu_2}_{X^\prime})]$.
\medskip
{\bf 2.2.6. Remarks. 1)} Most common in 2.2.5 1), either:
\medskip
a) $G$ is split and $X$ is local henselian, or 
\smallskip
b) $X={\rm Spec}(\QQ_p)$ (resp. $X={\rm Spec}(\ZZ_p)$) and $G$ (resp. $G_{\QQ_p}$) is unramified over $\QQ_p$.
\medskip
In case a) we can take $X^\prime=X$. Argument: we consider a cocharacter of $G$ which over the spectrum of the separable closure $k_X^{\rm sep}$ of the residue field $k_X$ of $X$ is $G(k_X^{\rm sep})$-conjugate to $\mu_{k_X^{\rm sep}}$, and so we can assume $X$ is strictly henselian. We just need to show that any two cocharacters $\mu_1$ and $\mu_2$ of $G$ which over $k_X$ are identical, are in fact $G(X)$-conjugate. Based on [SGA3, Vol. III, 6.1 of p. 32] (applied to the centralizers of the images of $\mu_1$ and $\mu_2$ in $G$) and [SGA3, Vol. III, 1.5 of p. 329], we can assume $\mu_1$ and $\mu_2$ factor through the same torus of $G$. This implies $\mu_1=\mu_2$. 
\smallskip
For future references we also point out that 
$$G^{\rm ad}(X)=T_{\rm ad}(X)q(G(X)),\leqno (DER)$$ 
where $T_{\rm ad}$ is an arbitrary maximal torus of $G^{\rm ad}$ and $q:G\twoheadrightarrow G^{\rm ad}$ is the natural epimorphism (this formula can be obtained by simply using two short exact sequences $0\to Z(G^{\rm der})\hookrightarrow G^{\rm der}\twoheadrightarrow G^{\rm ad}\to 0$ and $0\to Z(G^{\rm der})\hookrightarrow T_{\rm der}\twoheadrightarrow T_{\rm ad}\to 0$, with $T_{\rm der}$ as the pull back of $T_{\rm ad}$ to $G^{\rm der}$, in the cohomological context of the \'etale or flat topology of $X$).
\smallskip
In case b) there is a smallest unramified cover $X^\prime$ of $X$ such that the Shimura group pair involved is definable by a cocharacter $\mu\colon \GG_m\to G_{X^\prime}$, cf. [Mi3, 4.6-7] (resp. cf. loc. cit., [Bo2, 15.14] and the fact --it can be deduced from [Ti2]-- that $G_{X^\prime}$ has a maximal split torus whose generic fibre is as well a maximal split torus).
\smallskip
{\bf 2)} We can work 2.2.5 in terms of the faithfully flat (or \'etale)
topology of $X$ but for this paper, the present form of 2.2.5 is all that we need.
\smallskip
{\bf 3)} Based on 2.2.5 1) and 2), we speak about morphisms
(maps) between two Shimura group pairs $(G_1,[\mu_1])$ and $(G_2,[\mu_2])$ over $X$: these are homomorphisms $f:G_1\to G_2$ such that for a simply connected pro-\'etale cover $X^\prime$ of $X$, $[{\mu_1}_{X^{\prime}}\circ f_{X^\prime}]=[{\mu_2}_{X^\prime}]$. We denote such a morphism by 
$$f:(G_1,[\mu_1])\to (G_2,[\mu_2]).$$ 
It is called an injective map (resp. an isomorphism), if as a group homomorphism is a closed embedding (resp. is an isomorphism). We get the category ${\rm ShGp}(X)$ of Shimura group pairs over $X$.
\smallskip
{\bf 4)} A study of Shimura group pairs over $\ZZ_p$ together with their symplectic embeddings (i.e. injective maps into $\bigl({\rm GSp}(W_{\ZZ_p},\psi),[\mu]\bigr)$; here $(W_{\ZZ_p},\psi)$ is a symplectic space over $\ZZ_p$, while $[\mu]$ is uniquely determined by requiring that $\mu:\GG_m\to GSp(W_{\ZZ_p},\psi)$ acts on $W_{\ZZ_p}$ via the null and the inverse of the identical character of $\GG_m$), is implicitly started in 3.10 and is fully developed in \S 7. See 4.6.7-8 below for samples.
\smallskip
{\bf 5)} Any Shimura variety ${\rm Sh}(G,X)$ of dimension $\ge 1$ gives birth to a Shimura group pair $(G,[\mu])$ over ${\rm Spec}(\QQ)$. Here $\mu:\GG_m\to G_{E^1(G,X)}$ is a cocharacter, with $E^1(G,X)$ a finite field extension of $E(G,X)$, which over $\CC$ under an $E(G,X)$-monomorphism $E^1(G,X)\hookrightarrow\CC$, is $G(\CC)$-conjugate to anyone of the cocharacters $\mu_x$, $x\in X$, def. in [Va2, 2.2]. The association: $(G,[\mu])$ to ${\rm Sh}(G,X)$, is functorial.  
\medskip
{\bf 2.2.7. Remark.} Let $\ell\in\NN$. If $X$ is local henselian and if $G$ is an $X$-simple, adjoint, split group, then the
Shimura group pairs $(G,[\mu])$ are implicitly classified by [De3, 1.3.9] (cf. also [Sa]); it is [SGA3, Vol. III, 6.1 of p. 32 and 1.5 of p. 329] which allows us to treat the situation as if $X$ is the spectrum of a field. We have
$\bigl[{\ell+1\over 2}\bigr]$ isomorphism classes if $G$ is of $A_\ell$ ($\ell\ge 1$) Lie type, 2 isomorphism classes
 if $G$ is of $D_\ell$ Lie type ($\ell\ge 5$), and one such class if $G$ is of $D_4$, $E_6$,
$E_7$, $B_\ell$ or $C_\ell$ Lie type ($\ell\ge 1$). There are no Shimura group pairs $(\tilde{G},[\tilde{\mu}])$ with $\tilde{G}$ an $X$-simple, adjoint group
of $E_8$, $F_4$ or $G_2$ Lie type. 
\smallskip
Moreover any Shimura group pair $(G_1,[\mu_1])$ over such an $X$, with $G_1$ an adjoint group, can be written uniquely (up to isomorphism) as a product of pairs $(G_2,[\mu_2])$ which are either having $\mu_2$ as the trivial cocharacter or are Shimura group pairs having $G_2$ as an $X$-simple, adjoint group. This is a consequence of the classification of adjoint groups over a field (see [Ti1]).
\medskip
{\bf 2.2.8. Basic definitions and notations. 1)} A Shimura filtered $\sg$-crystal over $k$ is a quadruple $(M,F^1,\vph,G)$,
with $(M,F^1,\vph)$ a filtered $\sg$-crystal such that $F^1$ is a proper direct summand of $M$ and with $G$ a quasi-split reductive subgroup of $GL(M)$, for which the following condition is satisfied:
\medskip
{\bf a)} There is a family of tensors $(t_\al)_{\al\in\Mj}$ of the $F^0$-filtration of $\Mt(M[{1\over p}])$ defined by $F^1$, such that $\vph(t_{\al})=t_{\al}$, $\forall\al\in\Mj$, and $G_{B(k)}$ is the subgroup of $GL\bigl(M[{1\over p}]\bigr)$ fixing $t_\al$, $\forall\al\in\Mj$.
\medskip
The following two properties (introducing extra notations) are implied by this condition a):
\medskip
{\bf b)} The quadruple $\bigl({\rm Lie}(G),\vph, F^0\bigl({\rm Lie}(G)\bigr),F^1\bigl({\rm Lie}(G)\bigr)\bigr)$ is a filtered Lie $\sg$-crystal; here
$$F^0\bigl({\rm Lie}(G)\bigr):={\rm Lie}(G)\cap F^0\bigl({\rm End}(M)\bigr)=\bigl\{x\in {\rm Lie}(G)\bigm|x(F^1)\subset F^1\bigr\}$$ 
and
$$F^1\bigl({\rm Lie}(G)\bigr):={\rm Lie}(G)\cap F^1\bigl({\rm End}(M)\bigr)=\bigl\{x\in {\rm Lie}(G)\bigm|x(F^1)=\{0\}\, {\rm and}\, x(M)\subset F^1\bigr\};$$
\medskip
{\bf c)} There is an injective cocharacter $\mu:\GG_m\hookrightarrow G$  producing a direct sum decomposition
$M=F^1\oplus F^0$, with $\be\in\GG_m(W(k))$ acting through $\mu$ on $F^i$ as the multiplication with $\be^{-i}$,
 $i=\overline{0,1}$.
\medskip
The argument for why a) implies b) and c) goes as follows. c) is a consequence of 2.2.1.2 and of the assumption that $F^1$ is a proper direct summand of $M$. From c) we get that $F^i\bigl({\rm Lie}(G)\bigr)\otimes_{W(k)} k$ is the intersection of $F^i\bigl({\rm End}(M)\bigr)\otimes_{W(k)} k$ with ${\rm Lie}(G)\otimes_{W(k)} k$. So $\vph\bigl(p{\rm Lie}(G)+F^0({\rm Lie}(G))+{1\over p}F^1({\rm Lie}(G)\bigr)$ is a direct summand of ${\rm End}(M)$; as it is included in ${\rm Lie}(G)[{1\over p}]$ (this is a consequence of a)), it is ${\rm Lie}(G)$ itself. c) implies $F^0({\rm Lie}(G))$ is a parabolic Lie subalgebra of ${\rm Lie}(G)$ having $F^1({\rm Lie}(G))$ as its unipotent radical and we have $[{\rm Lie}(G),F^1({\rm Lie}(G))]\subset F^0({\rm Lie}(G))$; so b) follows.
\medskip
{\bf 2)} A Shimura $\sg$-crystal over $k$ is a triple $(M,\vph,G)$ which can be extended to a
quadruple $(M,F^1,\vph,G)$ defining a Shimura filtered $\sg$-crystal over $k$. $F^1$ itself or this quadruple, is called a lift of $(M,\vph,G)$.
\smallskip
{\bf 3)} A generalized Shimura $p$-divisible object of $\Mm\Mf_{[a,b]}(W(k))$ is a quadruple 
$$(M,(F^i(M))_{i\in S(a,b)},\vph,G), 
$$
where the triple $(M,(F^i(M))_{i\in S(a,b)},\vph)$ is a $p$-divisible object of $\Mm\Mf_{[a,b]}(W(k))$ (so $a,b\in\ZZ$, with $a\le b$) and $G$ is a quasi-split reductive subgroup of $GL(M)$, for which the following three conditions hold:
\medskip
{\bf a')} The condition a) above holds under the modification that the $F^0$-filtration of $\Mt(M)$ is defined (see 2.1) by the filtration $(F^i(M))_{i\in S(a,b)}$ of $M$;
\medskip
{\bf b')} The condition b) holds but with $F^0\bigl({\rm Lie}(G)\bigr)$ and $F^1\bigl({\rm Lie}(G)\bigr)$ as direct summands of ${\rm Lie}(G)$ induced (see 2.1) from the tensor product filtration of $\Mt(M)$;
\medskip
{\bf c')} The cocharacter $\mu:\GG_m\to G$ obtained by applying 2.2.1.2 to the quadruple $(M,(F^i(M))_{i\in S(a,b)},\vph,(t_\al)_{\al\in\Mj})$ (cf. a')) and producing a direct sum decomposition
$M=\oplus_{i\in S(a,b)}\tilde{F}^i$, with $\be\in\GG_m(W(k))$ acting through $\mu$ on $\tilde F^i$ as the multiplication with $\be^{-i}$,
 $i=\overline{a,b}$, and with $F^i(M)=\oplus_{j\in S(i,b)}\tilde{F}^j$, is such that at least two of the $W(k)$-submodules $\tilde{F}^i$, $i=\overline{a,b}$, are non-zero, and so in particular it is injective. 
\medskip
{\bf 3a)} If in 3) we do not impose b') and we just assume that $G$ is reductive, then we speak about a $p$-divisible object with a reductive structure of $\Mm\Mf_{[a,b]}(W(k))$.
\smallskip
{\bf 4)} A generalized Shimura $p$-divisible object over $k$ is a triple $(M,\vph,G)$ which can be extended to a quadruple $(M,(F^i(M))_{i\in S(a,b)},\vph,G)$, called a lift of $(M,\vph,G)$, defining a generalized Shimura $p$-divisible object of $\Mm\Mf_{[a,b]}(W(k))$.
\smallskip
{\bf 4a)} Similarly to 4), starting from 3a) we speak about a $p$-divisible object with a reductive structure over $k$ (in the range $[a,b]$) and about its lifts.
\smallskip
{\bf 5)} Let $n\in\NN$. A Shimura $\sg^n$-isocrystal over $k$ is a triple $(M_{B(k)},\vph,G_{B(k)})$, where $(M_{B(k)},\vph)$ is a $\sg^n$-isocrystal over $k$ and $G_{B(k)}$ is a reductive subgroup of $GL(M_{B(k)})$, for which there is a family of tensors $(t_{\al})_{\al\in\Mj}$ of $\Mt(M_{B(k)})$ such that $\vph(t_{\al})=t_{\al}$, $\forall\al\in\Mj$, and $G_{B(k)}$ is the subgroup of $GL(M_{B(k)})$ fixing $t_{\al}$, $\forall\al\in\Mj$.
\smallskip
{\bf 6)} A Shimura filtered $\sg^n$-isocrystal over $k$ is a quadruple 
$$(M_{B(k)},(F^i(M_{B(k)}))_{i\in\ZZ},\vph,G_{B(k)}),$$ 
where $(M_{B(k)},\vph,G_{B(k)})$ is a Shimura $\sg^n$-isocrystal over $k$ and $(F^i(M_{B(k)}))_{i\in\ZZ}$ is a decreasing filtration of $M_{B(k)}$, such that, with the notations of 5), we can moreover choose the family of tensors $(t_{\al})_{\al\in\Mj}$ to be in the $F^0$-filtration of $\Mt(M_{B(k)})$ defined by the filtration $(F^i(M_{B(k)}))_{i\in\ZZ}$ of $M_{B(k)}$. When $n=1$ we drop it.
\smallskip
{\bf 7)} By an endomorphism of a Shimura (resp. Shimura filtered) $\sg$-crystal $(M,\vph,G)$ (resp. $(M,F^1,\vph,G)$) we mean an element $f\in {\rm Lie}(G)$ (resp. $f\in F^0({\rm Lie}(G))$) fixed by $\vph$. We denote it by $f:(M,\vph,G)\to (M,\vph,G)$ (resp. $f:(M,F^1,\vph,G)\to (M,F^1,\vph,G)$). The set of such endomorphisms has a natural structure of a Lie algebra over $\ZZ_p$. Similarly, we define endomorphisms for 3) to 6).
\medskip
{\bf 2.2.9. Remarks, facts and variants. 0)}  We often drop out mentioning ``over $k$".
\smallskip
{\bf 1)} We refer to 2.2.8 1). We saw that condition a) implies $\vph\bigl({\rm Lie}(G)[{1\over p}]\bigr)={\rm Lie}(G)[{1\over p}]$. It is ``almost" implied by b), cf. [Fa2, rm. iii) after th. 10]: in loc. cit. we do have to consider Tate-twists as well; moreover, according to us, the first (warning: not the second) part of loc. cit. involving the existence of the line $L_0$ is not correct as it stands (this is why we used condition a) to define Shimura filtered $\sg$-crystals). 
\smallskip
To argue this, let ${\got C}_0:=(M_0,F^1_0,\vph_0)$ be a filtered $\sg$-crystal. Let $\rho_0:\Gamma_k\to GL(N_0)(\ZZ_p)$ be the Galois representation defined by the dual of the Tate-module of any $p$-divisible group $D_0$ over $W(k)$ of whose associated filtered $\sg$-crystal is ${\got C}_0$. If $p=2$ we assume $D_0$ exists (it does exist after replacing $k$ by an abelian extension of degree dividing a power of $2$ which depends only on $(M_0,F^1_0)$, cf. the Corollary of the review 2.3.18.1 D below). So $N_0$ is a free $\ZZ_p$-module of the same rank as $M_0$. If $p\ge 3$, then $\rho_0$ is uniquely determined by ${\got C}_0$ (for instance, cf. [Fa2, th. 7]). Let $G_0$ be the flat, closed subgroup of $GL(N_0)$ which is the algebraic envelope of $\rho_0$ (see [Se1]). Let $G_1$ be an arbitrary integral, closed subgroup of $GL(N_0)$ normalized by $G_0$. Let $\tilde G_1$ be the integral, closed subgroup of $GL(M_0)$ corresponding to $G_1$ via Fontaine's comparison theory (for $p=2$ it is recalled in 2.3.18.1 E below). We assume it is smooth and its Lie algebra is the underlying module of a Lie $p$-divisible subobject of $End({\got C}_0)$ (if $p\ge 5$, the second part of this assumption is implied by the first one, cf. [Fa1, 2.6]). It is very easy to adjust the situation such that $\tilde G_1$ has a connected special fibre as well. But if a line $L_0$ exists as in [Fa2, rm. iii) after th. 10], then $G_0$ is contained in a subgroup of $GL(N_0)$ which is the extension of $\GG_m$ by $G_1$. This is not always so. So what is missing in loc. cit.? The mistake consists in: the $\ZZ_p$-submodules of $\Mt(N_0)$ of rank 1 and normalized by $G_1$, are not fixed by $G_0$ but just permuted (and the orbits under the resulting permutation are not necessarily finite). In other words, there are two ingredients missing in loc. cit.: condition c) as well as the passage to an algebraic field extension of $k$ (see the Fact below).
\smallskip
We come back to our situation of 2.2.8 1). As already pointed out c) results from a). On the other hand, we have:
\medskip
{\bf Fact.} {\it If $k=\bar k$ and the conditions b) and c) hold, then a) holds automatically.}
\medskip
{\bf Proof:} We just need to use a standard argument involving $\ZZ_p$-structures (it is recalled in 8) below) and the fact that $G_{B(k)}$ is reductive: with the independent notations of 8), $G_{\QQ_p}$ is the subgroup of $GL(M^a[{1\over p}])$ fixing some tensors of $\Mt(M^a)$) (cf. [De4, 3.1 c)]); viewing these tensors as tensors of $\Mt(M)$ (via the natural identification $M=M^a\otimes_{\ZZ_p} W(k)$) we do get that a) holds automatically. 
\medskip
So no doubt, a great part of the essence of def. 2.2.8 1) is expressed by b). 
All above applies as well to the context of 2.2.8 3).
\smallskip
{\bf 1')} 1)  suggests the introduction of pseudo Shimura (filtered) $\sg$-crystals: instead of a) of 2.2.8 1) we assume just b) of 2.2.8 1). Similarly, if b) and c) hold and a) holds after a passage to $W(k_1)$, with $k_1$ a finite field extension of $k$ (resp. with $k_1=\bar k$), then we speak about quasi (resp. potentially) Shimura (filtered) $\sg$-crystals. Warning: for these notions we still keep $G$ to be a quasi-split reductive group over $W(k)$. The greatest part of the terminology to be introduced below for Shimura $\sg$-crystals makes sense (and so it is used as well without extra comment) for these pseudo, quasi or potentially contexts. 
\smallskip
We have the following Galois interpretations. Not to introduce extra notations, we refer to the context of a pseudo Shimura filtered $\sg$-crystal $(M_0,F^1_0,\vph_0,\tilde G_1)$ as in 1)  (so $\tilde G_1$ is a quasi-split reductive group and $G_0$ normalizes $G_1$; for $p=2$ we assume $D_0$ exists). 
\medskip
{\bf Fact.} {\it $(M_0,F^1_0,\vph_0,\tilde G_1)$ is a Shimura (resp. quasi Shimura) filtered $\sg$-crystal iff $G_0$ (resp. the Zariski closure in $GL(N_0)$ of the connected component of the origin of $G_{0\QQ_p}$) is a subgroup of $G_1$.}
\medskip
{\bf Proof:} Below we use freely the passage (via Fontaine's comparison theory) from tensors in the $F^0$-filtration of $\Mt(M_0[{1\over p}])$ fixed by $\vph$ to tensors of $\Mt(N_0[{1\over p}])$ fixed by $\rho_0(\Gamma_k)$. For instance, it implies (via arguments involving torsors) $G_{1\QQ_p}$ is a reductive group. So loc. cit. applies again to the representation of $G_{1\QQ_p}$ on $N_0[{1\over p}]$. Based on this, the part ``without quasi" is obvious. We deal now with the part ``involving quasi". If $(M_0,F^1_0,\vph_0,\tilde G_1)$ is a quasi Shimura filtered $\sg$-crystal, then the image through $\rho_0$ of a compact, open subgroup of $\Gamma_k$ is contained in $G_1(\ZZ_p)$ (see 2.2.8 a)); so the connected component of the origin of $G_{0\QQ_p}$ is a subgroup of $G_{1\QQ_p}$. 
\smallskip
We assume now that the connected component of the origin of $G_{0\QQ_p}$ is a subgroup of $G_{1\QQ_p}$. So the canonical split cocharacter $\mu:\GG_m\to GL(M_0)$ of $(M_0,F^1_0,\vph_0)$, factors through $\tilde G_1$. But, as $\tilde G_1$ is reductive, using an argument involving $\QQ_p$-structures entirely similar to the one of the proof of the Fact of 1) but performed in the \'etale context, we get that the restriction of $\rho_0$ to $\Gamma_{\bar k}$ factors through $G_1(\ZZ_p)$. So, as the number of connected components of $G_{0\QQ_p}$ is finite, we deduce the existence of a finite Galois extension $k_1$ of $k$, such that the restriction of $\rho_0$ to $\Gamma_{k_1}$ factors through $G_1(\ZZ_p)$. So, based once more on a similar $\QQ_p$-argument as in the proof of 1) but performed in the \'etale context, we get that the extension of $(M_0,F^1_0,\vph_0)$ to $k_1$ is a Shimura filtered $\sg_{k_1}$-crystal. This ends the proof.
\medskip
As in 2.2.8 7), we define endomorphisms of pseudo or quasi or potentially (filtered) Shimura $\sg$-crystals.
\smallskip
{\bf 2)} If in 2.2.8 1) we have
$F^1({\rm Lie}(G))\ne 0$, then the pair $(G,[\mu])$ defines a Shimura group pair over ${\rm Spec}(W(k))$.
Similarly, the cocharacter $\mu^{\rm ad}:\GG_m\to G^{\rm ad}$ (induced naturally by the cocharacter $\mu$ of c') of 2.2.8 3)) is either trivial or gives birth (cf. 2.2.9 b')) to a Shimura adjoint group pair $(G^{\rm ad},[\mu^{\rm ad}])$ over ${\rm Spec}(W(k))$.
\smallskip
{\bf 3)} We refer to 2.2.8 1). The $G(W(k))$-conjugacy class $\Mc\Mc\Mp$ of the parabolic Lie subalgebra $F^0({\rm Lie}(G))$ of ${\rm Lie}(G)$, is uniquely determined by the Shimura $\sg$-crystal $(M,\vph,G)$. This is a consequence of the fact that the direct summand $F^1/pF^1$ of $M/pM$ is uniquely determined  by $\vph$ and of the following Fact:
\medskip
{\bf Fact 1.} {\it Two parabolic subgroups $P_1$ and $P_2$ of $G$ having the same special fibre, are $G(W(k))$-conjugate.}
\medskip
{\bf Proof:} As this well known result plays an important role in many parts of \S2-14, we include a proof of it. Let $T_k$ be a maximal torus of $P_{1k}=P_{2k}$ and so of $G_k$ (cf. [Bo2, 18.2 (i)]). From [SGA3, Vol. II, 3.6 of p. 48] we deduce  that $T_k$ lifts to a maximal torus $T_i$ of $P_i$, $i=\overline{1,2}$. To loc. cit. we need to add that:
\medskip
{\bf Subfact.} {\it Any closed subscheme of $G^\wedge$ which mod $p^m$ is a torus $\tilde T_m$, $\forall m\in\NN$, is obtained by completing $p$-adically a uniquely determined torus $\tilde T$ of $G$.}
\medskip
The argument goes as follows. We can assume $k=\bar k$ and $G=GL(M)$. We consider the direct sum decomposition $M/p^mM=\oplus_{\be\in CH_m} M_{\be}(m)$, with $CH_m$ as the set of characters of $\tilde T_m$, naturally defined by the action of $\tilde T_m$ on $M/p^nM$. As $\tilde T_{m+1}$ lifts $\tilde T_m$, $CH_m$ does not depend on $m$ and so we can denote it by $CH$; moreover, $\forall\be\in CH$, there is (as $M$ is $p$-adically complete) a direct summand $M(\be)$ of $M$ such that $M_{\be}(m)$ is the reduction mod $p^m$ of it, $\forall m\in\NN$. So we get an action of the (abstract) torus $\tilde T$ whose reductions modulo positive, integral powers of $p$ are $\tilde T_m$'s, on $\tilde M$: $\tilde T$ acts on $M(\be)$ via $\be$. It is a closed embedding (this is a matter of characters and so it follows from the fact that mod $p$ we have a closed embedding). The Subfact follows.
\smallskip
This addition to loc. cit. is not repeated in the subsequent places relying on it. 
\smallskip
Also from loc. cit. we deduce the existence of $g\in G(W(k))$ which mod $p$ is the identity, such that under its inner conjugation $T_1$ is taken onto $T_2$. So we can assume $T_1=T_2$. But $P_{1W(\bar k)}$ and $P_{2W(\bar k)}$ are determined by the set of roots (i.e. of characters of $T_{1W(\bar k)}$) of the action (via inner conjugation) of $T_{1W(\bar k)}$ on their Lie algebras. As these sets can be read out from this action mod $p$, we get $P_1=P_2$. This ends the proof.
\medskip
We call $\Mc\Mc\Mp$ as the filtration class  of the Shimura $\sg$-crystal $(M,\vph,G)$. About $\mu$ (or the $G(W(k))$-conjugacy class $[\mu]$ defined by $\mu$) we say it defines $\Mc\Mc\Mp$. 
\medskip
{\bf Fact 2.} {\it $[\mu]$ itself is uniquely determined by $(M,\vph,G)$.}
\medskip
{\bf Proof:}  
Let $\mu_1:\GG_m\to G$ be another cocharacter producing (as in c) of 2.2.8 1)) a direct sum decomposition $M=F_1^1\oplus F_1^0$ such that $(M,F_1^1,\vph,G)$ is a Shimura filtered $\sg$-crystal. Based on Fact 1 we can assume $F_1^0({\rm Lie}(G))=F^0({\rm Lie}(G))$. Let $P_0$ be the parabolic subgroup of $G$ having $F^0({\rm Lie}(G))$ as its Lie algebra; both $\mu$ and $\mu_1$ factor through it. Let $T$ be a maximal split torus of $P_0$: it is obtained by the usual lifting process (of loc. cit.), from a maximal split torus of $P_{0k}$. Applying [Bo2, 15.14] to $P_{0k}$, we can assume that the special fibres of $\mu$ and $\mu_1$ factor through $T_k$. Based again on [SGA3, Vol. II, 3.6 of p. 48] (applied to $P_0$ and its cocharacters $\mu$ and $\mu_1$) we can assume $\mu$ and $\mu_1$ themselves factor through $T$. In  particular, we get: their images are commuting to each other. As $F_1^1/pF^1_1=F^1/pF^1$, this implies $F_1^1=F^1$ and $F_0=F_0^1$; so $\mu=\mu_1$. This ends the proof.  
\medskip
{\bf Corollary.} {\it The $G^{\rm ad}(W(k))$-conjugacy class $[\mu^{\rm ad}]$ defined by the cocharacter $\mu^{\rm ad}:\GG_m\to G^{\rm ad}$ induced by $\mu$, is uniquely determined by $(M,\vph,G)$. Moreover, if $G$ is a torus, then $\mu$ itself is uniquely determined.}
\medskip
{\bf 4)} Generalized Shimura $p$-divisible objects over $k$ show up in this paper in: 3.1.7, the beginning paragraph of 3.4 (and so implicitly in many parts of 3.4-5), 3.6.1.5-6, 3.15.6, 4.5.6 8), 4.4.13, 4.9.18, 4.12.12.0, 4.14.2 and Appendix. They will be used much more frequently in \S5-10 (cf. also 4.5.4 below). Similarly, the ``with a reductive structure" context of 2.2.8 3a) and 4a) shows up concretely just in 3.13.7.8-9 and in Appendix. The whole of 1), 1') and of 3) as well as the whole of 5), 6), 8), 9) and 10) below can be entirely adapted to this last context: as no modifications are needed (besides always assuming in 1) that $\rho_0$ exists), we do not restate the things. Warning: when we speak about pseudo (or quasi or potentially) $p$-divisible objects with a reductive structure, we do not assume that the reductive group involved is quasi-split.
\smallskip
We use freely the terminology to be introduced for Shimura (filtered) $\sg$-crystals and which depends only on their attached (cf. 2.2.13 below) Shimura Lie $\sg$-crystals, in the context of generalized Shimura $p$-divisible objects over $k$ or of $\Mm\Mf(W(k))$; as a sample see 2.2.22 2) below. Any context involving generalized Shimura $p$-divisible objects is often referred as a generalized Shimura context.
\smallskip
{\bf 5)} When we want to emphasize a family of tensors $(t_\al)_{\al\in\Mj}$ as in a) of 2.2.8 1), we denote a Shimura
$\sg$-crystal by a quadruple $(M,\vph,G,(t_\al)_{\al\in\Mj})$. Similarly we denote a Shimura filtered $\sg$-crystal by a quintuple $(M,F^1,\vph,G,(t_\al)_{\al\in\Mj})$. Strictly speaking we should refer to $(M,\vph,G,(t_\al)_{\al\in\Mj})$ as a Shimura $\sg$-crystal with an emphasized family of tensors but we rarely do so; the same about the filtered case. 
\smallskip
If $m\in\NN\cup\{0\}$ is such that all homogeneous components of $t_{\al}$ are of degree at most $m$, $\forall\al\in\Mj$, then we say $(M,\vph,G)$ is of degree at most $m$; by the degree 
$${\rm deg}(M,\vph,G)\in\NN\cup\{0\}$$ 
we mean the smallest such possible value $m$, where the family $(t_{\al})_{\al\in\Mj}$ is allowed to vary subject to the constraint of a) of 2.2.8 1). If ${\rm deg}(M,\vph,G)=0$, then $G=GL(M)$. 
\smallskip
{\bf 6)} If $(M,\vph,G,(t_\al)_{\al\in\Mj})$ and $(M_1,\vph_1,G_1,(t_{1\al})_{\al\in\Mj_1})$
(resp. if $(M,F^1,\vph,G,(t_\al)_{\al\in\Mj})$ and $(M_1,F^1_1,\vph_1,G_1,(t_{1\al})_{\al\in\Mj_1}))$ are two Shimura (resp. two Shimura filtered) $\sg$-crystals, then an isomorphism between them is given by an
isomorphism $t:M\tilde\to M_1$ taking the triple $(\vph,G,(t_\al)_{\al\in\Mj})$ into $(\vph_1,G_1,(t_{1\al})_{\al\in\Mj_1})$ (resp. taking the quadruple $(F^1,\vph,G,(t_\al)_{\al\in\Mj})$ into $(F^1_1,\vph_1,G_1,(t_{1\al})_{\al\in\Mj_1})$). If $j:\Mj\tilde\to\Mj_1$ is a bijection and if the isomorphism $t$ takes $t_{\al}$ into $t_{1j(\al)}$, then we call $t$ a $j$-isomorphism. In particular, an $1_{\Mj}$-isomorphism between two Shimura $\sg$-crystals 
$(M,\vph,G,(t_\al)_{\al\in\Mj})$ and $(M,\vph_1,G,(t_\al)_{\al\in\Mj})$ is given by an element $g\in G(W(k))$ such that $\vph_1=g\vph g^{-1}$. 
An isomorphism of $(M,\vph,G)$ does not need to fix the family of tensors $(t_{\al})_{\al\in\Mj}$; easy examples can be constructed if $k=\FF_p$ and so, by natural extension, over any $k$.
\smallskip
Warning: the $1_{\Mj}$-automorphisms are not in general those endomorphisms which are invertible. 
\smallskip
{\bf 7)} The condition $G$ quasi-split, used in 2.2.8 1) to 4), is not really needed (cf. 3.8 below). It is inserted just to be able to be concrete in 3.2.3 below. Also Theorem 2.3.9 below makes it very convenient for applications to Shimura varieties of Hodge type. When we do not know that $G$ is quasi-split or when it is not, we speak about a Shimura (filtered) $\sg$-crystal of being not necessarily quasi-split. 
\smallskip
{\bf 8)} If in 2.2.8 1) we have $k=\bar k$, then writing $\vph=a\circ\mu({1\over p})$, $a$ becomes a $\sg$-linear automorphism of $M$ and so of $\Mt(M)$ (and in particular of ${\rm Lie}(G)$), and 
$$M^a:=\{x\in M|a(x)=x\}$$ 
is a free $\ZZ_p$-module of rank $\dim_{W(k)}(M)$. Moreover, 
$${\rm Lie}(G)^a:=\{x\in {\rm Lie}(G)|a(x)=x\}={\rm Lie}(G)\cap {\rm End}(M)^a\subset{\rm End}
(M^a)$$ 
is the Lie algebra of a reductive subgroup $G_{\ZZ_p}$ of $GL(M^a)$; the extension of $G_{\ZZ_p}$ to $W(k)$ is the reductive subgroup $G$ of $GL(M)$. As $\mu$ and $\vph$ fix $t_{\al}$, we have $t_{\al}\in\Mt(M^a)[{1\over p}]$, $\forall\al\in\Mj$.
\smallskip
Let $B(\FF_{p^d})$ be the field of definition of
the $G(B(k))$-conjugacy class of the cocharacter $\mu_{B(k)}:\GG_m\to G_{B(k)}$ (cf. [Mi3, 4.6-7] applied to $G_{\QQ_p}$; $G_{\QQ_p}$ splits over a finite unramified extension of $\QQ_p$). Here $d\in\NN$. Let $T_{\ZZ_p}$ be a maximal torus of $G_{\ZZ_p}$ which over $W(\FF_{p^d})$ contains a maximal split torus of $G_{W(\FF_{p^d})}$. Let $\mu_0:\GG_m\to G_{W(\FF_{p^d})}$ be a cocharacter which over $W(k)$ is $G(W(k))$-conjugate to $\mu$ and which factors through $T_{W(\FF_{p^d})}$. The existence of $\mu_0$ is guaranteed by loc. cit. and by the fact (see [Ti2]) that the generic fibre of a split torus of $G_{W(\FF_{p^d})}$ is a maximal split torus of $G_{B(\FF_{p^d})}$. We deduce the existence of an element $g\in G(W(k))$ such that the Shimura filtered $\sg$-crystal 
$\bigl(M,F^1,g\vph,G,{(t_\al)}_{\al\in\Mj}\bigr)$ is $1_{\Mj}$-isomorphic to the extension to
$k$ of a Shimura filtered $\sg$-crystal 
$$\bigl(M^a\otimes_{\ZZ_p} W(\FF_{p^d}),F^1_{\mu_0},a\circ \mu_0({1\over p}),G_{W(\FF_{p^d})},{(t_\al)}_{\al\in\Mj}\bigr)$$ 
over $\FF_{p^d}$, where $F^1_{\mu_0}\subset M^a\otimes_{\ZZ_p} W(\FF_{p^d})$ is the $F^1$-filtration associated as usual to $\mu_0$ (here we still denote by $a$ its restriction to $M^a\otimes_{\ZZ_p} W(\FF_{p^d})$). If $G$ is a torus, then $g=1_M$. 
\smallskip
So, from many points of view, for the study of Shimura $\sg$-crystals over a perfect field $k_0$ of characteristic $p$, we can assume $k_0$ is a finite field (if $k_0\ne\overline{k_0}$, then passing to a finite field extension $k_1$ of $k_0$, $M^a$, defined as above but working with $g(\vph\otimes 1)$ instead of $\vph$ for a suitable $g\in G(W(k_1))$, is still a free $\ZZ_p$-module of rank $\dim_{W(k)}(M)$).
\smallskip
{\bf 9)} Any filtered $\sg$-crystal $(M,F^1,\vph)$, with $F^1$ a proper direct summand of $M$, is identified with the Shimura filtered $\sg$-crystal $(M,F^1,\vph,GL(M))$, and any $\sg$-crystal $(M,\vph)$ which can be extended to a filtered $\sg$-crystal $(M,F^1,\vph)$, with $F^1$ a proper direct summand of $M$, is identified with the Shimura $\sg$-crystal $(M,\vph,GL(M))$. 
\smallskip
{\bf 10)} Any Shimura (filtered) $\sg$-crystal over $k$ gives birth to a Shimura (filtered) isocrystal over $k$, by making $p$ invertible. But not any Shimura (filtered) isocrystal over $k$ is obtained in this manner: for instance, the Shimura isocrystals over $k$ can have slopes which do not belong to the interval $[0,1]$ and the reductive groups involved are not necessarily quasi-split.
\smallskip
{\bf 11)} The condition (of 2.2.8 1)) $F^1$ is a proper direct summand of $M$ as well as the condition (of 2.2.8 3)) at least two of the submodules $\tilde{F}^i$'s are different from $\{0\}$, are inserted just to avoid the trivial cases. So the cocharacter $\mu$ of c) of 2.2.8 1) (or of c') of 2.2.8 3)) does not factor through the center of $GL(M)$. When we have $F^1=\{0\}$ or $F^1=M$ in 2.2.8 1) (resp. we have $M=\tilde F^i$ for some $i\in S(a,b)$ in 2.2.8 3)), we speak about a trivial Shimura (filtered) $\sg$-crystal (resp. about a trivial generalized Shimura $p$-divisible object of $\Mm\Mf(W(k))$). Without a special reference, we always assume we are in a non-trivial situation.
\smallskip
{\bf 12)} For the sake of flexibility it is important in 2.2.8 1) not to assume that $(M,\vph)$ has a principal quasi-polarization or that $\dim_{W(k)}(F^1)=\dim_{W(k)}(F^0)$ (as we have in situations emerging from principally polarized abelian varieties over $W(k)$). So we completely ignore the phenomenon of Tate-twists. However, when we have a principal quasi-polarization
$$p_M:M\otimes_{W(k)} M\to W(k)(1)$$ 
normalized by $G$, we speak about a principally quasi-polarized Shimura (filtered) $\sg$-crystal. As notation, we put $p_M$ as part of some $n$-tuple at the very end; example: $(M,F^1,\vph,G,(t_{\al})_{\al\in\Mj},p_M)$. Moreover we do not assume $G_{B(k)}$ is the subgroup of $GL(M[{1\over p}])$ fixing $t_{\al}$, $\forall\al\in\Mj$: we just assume it is the subgroup of ${\rm GSp}(M[{1\over p}],p_M)$ fixing these tensors. 
\smallskip
Let now $k=\bar k$. We consider a $\ZZ_p$-structure $G_{\ZZ_p}$ of $G$ obtained as in 8). As $\mu$ is injective, $G$ does not fix $p_M$. So $Z(G_{\ZZ_p})$ contains a split torus of dimension $1$ which normalizes $p_M$ without fixing it. The elements of $Z(G_{\ZZ_p})(\ZZ_p)$ are naturally identified with $1_{\Mj}$-automorphisms of $(M,F^1,\vph,G,(t_{\al})_{\al\in\Mj})$ but not of $(M,F^1,\vph,G,(t_{\al})_{\al\in\Mj},p_M)$.
\medskip
{\bf 2.2.10. Shimura filtered $F$-crystals.} Let $(M,F^1,\vph,G)$ be a Shimura filtered $\sg$-crystal. Let $H$ be a smooth subscheme of $G$ such that the origin of $G$ factors through $H$. Let $R$ be the $W(k)$-algebra of the completion of $H$ at this factorization $o$. We consider a $W(k)$-isomorphism 
$\tilde f:R\tilde\to W(k)[[x_1,\ldots,x_d]]$, with $d$ the relative dimension of $H$ over $W(k)$ in $o$; we view it as an identification. Let $\Phi_R$ be the Frobenius lift of $R$ taking (via this identification) $x_i$ to $x^p_i$,
$i=\overline{1,d}$. Let $\bar\Om_{R/W(k)}$ be the free $R$-submodule of $\Om_{R/W(k)}$ generated by $dx_1$,..., $dx_d$.
\smallskip
We consider the $p$-divisible object of $\Mm\Mf_{[0,1]}(R)$ defined by the following triple $(M\otimes_{W(k)} R, F^1\otimes_{W(k)} R,\Phi)$,
with 
$$\Phi:=h(\vph\otimes id)$$ 
defined by the universal element $h\in H(R)$ of $H$. [Fa2, th. 10] guarantees that there is a unique connection 
$$\nabla:M\otimes_{W(k)} R\to M\otimes_{W(k)} \bar\Om_{R/W(k)}$$ 
such that $\Phi$ is $\nabla$-parallel in the sense of [Fa2, ch. 7]; from [Fa2, th. 10] we also deduce it is integrable and nilpotent mod $p$. In 3.6.18.4 we reobtain in a completely new manner the existence and the uniqueness of $\nabla$, while in 3.6.18.4.1 we reobtain in a completely new manner the fact that $\nabla$ is automatically integrable and nilpotent mod $p$.  
\smallskip
The quadruple
$$(M\otimes_{W(k)} R, F^1\otimes_{W(k)} R,\Phi,\nabla)$$
is a filtered $F$-crystal over $R/pR$. We call it a Shimura filtered $F$-crystal over $R/pR$ (defined by $(M,F^1,\vph,G)$) and we denote it by a sextuple
$$(M,F^1,\vph,G,H,\tilde f).$$ 
When $H=G$, we denote it by a  quintuple $(M,F^1,\vph,G,\tilde f)$; when $d=0$ we keep denoting it by $(M,F^1,\vph,G)$.
Usually $H$ is a smooth subgroup of $G$. When we want to emphasize the family of tensors $(t_\al)_{\al\in\Mj}$, then we denote a Shimura filtered $F$-crystal also by a $7$-tuple
$\bigl(M,F^1,\vph,G,H,\tilde f,(t_\al)_{\al\in\Mj}\bigr)$ and then, by abuse of notation (cf. 2.1), we still write $(t_\al)_{\al\in\Mj}$ for this family of tensors, when viewed as being formed by tensors of $\Mt(M\otimes_{W(k)} R[{1\over p}])$. We have: 
\medskip
{\bf Fact 1.} {\it $\nabla(t_\al)=0$, $\Phi(t_\al)=t_\al$ and $t_\al$ belongs to the $F^0$-filtration of $\Mt(M\otimes_{W(k)} R[{1\over p}])$ defined by $F^1[{1\over p}]$, $\forall\al\in\Mj$.}
\medskip
The fact that $\nabla(t_\al)=0$ is implied by [Fa2, rm. ii) after th. 10], while the other parts of Fact 1 are obvious. So the $R$-linear endomorphism of $M\otimes_{W(k)} R$ taking $m\in M$ into $\nabla({\partial\over {\partial x_i}})(m)$ belongs to ${\rm Lie}(G)\otimes_{W(k)} R$, $\forall i\in S(1,d)$. Accordingly, we still denote by $\nabla$ the natural connection on $\Mt(M)\otimes_{W(k)} R$ or on $N\otimes_{W(k)} R$, with $N$ a direct summand of $\Mt(M)$ normalized by $G$, induced by the connection $\nabla$ on $M\otimes_{W(k)} R$; in particular, this applies to ${\rm Lie}(G)\otimes_{W(k)} R\subset (M\otimes_{W(k)} M^*)\otimes_{W(k)} R$. 
\smallskip
Fact 1 implies directly the following two Facts:
\medskip
{\bf Fact 2.} {\it The quadruple 
$$({\rm Lie}(G)\otimes_{W(k)} R,\Phi,F^0({\rm Lie}(G))\otimes_{W(k)} R,F^1({\rm Lie}(G))\otimes_{W(k)} R,\nabla)$$ 
is a filtered Lie $F$-crystal over $R/pR$}.
\medskip
{\bf Fact 3.} {\it For any point of ${\rm Spec}(R/pR)$ with values in a perfect field $k_1$, we get from $(M,F^1,\vph,G,H,\tilde f)$ by pulling back, a Shimura $\sg_{k_1}$-crystal over $k_1$.}
\medskip
Similarly to the end of 2.2.1 c) or to 2.2.9 12), we speak about principally quasi-polarized Shimura filtered $F$-crystals. Moreover, the conventions of 2.2.9 12) apply as well. Warning: if $(M,F^1,\vph,G,p_M)$ is a principally quasi-polarized Shimura filtered $\sg$-crystal, then whenever we want to keep $p_M$ as the cycle of a principal quasi-polarization of $(M,F^1,\vph,G,H,\tilde f)$, we need to take $H$ to be as well a subscheme of ${\rm Sp}(M,p_M)$.     
\medskip
{\bf 2.2.10.1. Comments.} 2.2.10 can be performed as well for the context of a $W(k)$-morphism $f_H:H\to G$ for which there is a $W(k)$-morphism $z_H:{\rm Spec}(W(k))\to H$ such that $z_H\circ f_H$ is the origin of $G$, with $H$ an arbitrary regular, formally smooth $W(k)$-scheme (so ${\rm Spec}(R)$ is now the completion of $H$ in $z_H$). The resulting quadruple $(M\otimes_{W(k)} R,F^1\otimes_{W(k)} R,\Phi,\nabla)$ is still referred as a Shimura filtered $F$-crystal over $R/pR$ but the notation is $(M,F^1,\vph,G,f_H\circ z_H,\tilde f)$. 
\smallskip
However, in such a generality we often prefer to speak about filtered $F$-crystals (with tensors) over the special fibre $S_H={\rm Spec}(R/pR)$ of the completion of $H$ in $z_H$: we feel inclined to reserve the terminology Shimura filtered $F$-crystals just for contexts where we do have some closed embedding into the completion of $G$ in its origin. Similarly, when we are in the context of a $p$-divisible object with tensors of $\Mm\Mf_{[0,1]}^\nabla(H)$ which is modeled on 2.2.10, then we often prefer to use the alternating terminology filtered $F$-crystals (with tensors) instead of the one of Shimura filtered $F$-crystals.
\smallskip
Any Shimura (filtered) $\sg$-crystal over $k$ is a Shimura (filtered) $F$-crystal over $k$; we use this second terminology only when $k$ is not specified. 
\medskip
{\bf 2.2.11. Definitions. 1)} A non-trivial Shimura filtered Lie $\sg$-crystal is a filtered Lie $\sg$-crystal $\bigl({\got g},\vph,F^0({\got g}), F^1({\got g})\bigr)$, with {\got g} the Lie algebra of a reductive group $G$ over $W(k)$, for which there is an injective cocharacter $\mu:\GG_m\hookrightarrow G^{\rm ad}$ such that:
\medskip
--  $(G^{\rm ad},[\mu])$ is a Shimura group pair over ${\rm Spec}(W(k))$;
\smallskip
-- it defines $F^0(\got g)$ and $F^1(\got g)$; in other words, we have 
${\got g}={\got g}^0\oplus{\got g}^1\oplus{\got g}^{-1}$, with ${\got g}^1\ne 0$, with $\be\in\GG_m(W(k))$ acting through $\mu$ on ${\got g}^i$ as the multiplication with $\be^{-i}$, $i=\overline{-1,1}$, and with $F^1({\got g})={\got g}^1$ and $F^0({\got g})={\got g}^0\oplus{\got g}^1$.
\medskip
{\bf 2)} A non-trivial Shimura Lie $\sg$-crystal is a Lie $\sg$-crystal $({\got g},\vph)$ which can be extended to a non-trivial Shimura filtered Lie $\sg$-crystal $\bigl({\got g},\vph,F^0({\got g}),F^1({\got g})\bigr)$. $\bigl({\got g},\vph,F^0({\got g}),F^1({\got g})\bigr)$ or just $F^0({\got g})$ itself is called a lift of $({\got g},\vph)$; as $F^1({\got g})$ is determined by $F^0({\got g})$, being its nilpotent radical, we never refer to $(F^0({\got g}),F^1({\got g}))$ as a lift of $({\got g},\vph)$. 
\smallskip
Writing $\vph=a\circ\mu({1\over p})$, $a$ is a $\sg$-linear automorphism of {\got g}, whose class 
$$[a]\in {\rm Aut}_{\sg-{\rm lin}}({\got g})/G^{\rm ad}(W(k))$$ 
(modulo inner automorphisms
of {\got g}) does not depend (cf. Corollary of 2.2.9 3)) on the choice of $\mu$. We call $[a]$ the automorphism class of $({\got g},\vph)$.
\smallskip
{\bf 3)} A Lie $\sg$-crystal $({\got g},\vph)$, with {\got g} the Lie algebra of a reductive group $G$ over $W(k)$ and with $\vph({\got g})={\got g}$, is called a trivial Shimura Lie $\sg$-crystal; in such a case the quadruple $({\got g},\vph,{\got g},\{0\})$ is referred as a trivial Shimura filtered Lie $\sg$-crystal.
When we do not bother about trivial or non-trivial, we speak about Shimura (filtered) Lie $\sg$-crystals.
\smallskip
{\bf 4)} If in 1) to 3), $G$ is an adjoint group then we speak about a Shimura adjoint (filtered) Lie $\sg$-crystal.
\medskip
{\bf 2.2.11.1. Remark.} We refer to 2.2.11 1) and 2). From the existence of $\mu$ we get: $F^0({\got g})$ is a parabolic Lie subalgebra of ${\got g}$, having $F^1({\got g})$ as its nilpotent radical. Facts 1 and 2 of 2.2.9 3) apply entirely, cf. the below general Fact and the fact that $F^0({\got g})/pF^0({\got g})$ is uniquely determined by $\vph$, being the kernel of the restriction to ${\got g}/p{\got g}$ of $p\vph$ mod $p$. So we use the same terminology: the $G^{\rm ad}(W(k))$-conjugacy class of $F^0({\got g})$ is called the filtration class of $({\got g},\vph)$ and we say that $[\mu]$ (or just $\mu$ or $F^0({\got g})$) defines it.  
\medskip
{\bf Fact.} {\it We consider two parabolic subgroups $P^1_H$ and $P^2_H$ of a reductive group $H$ over a field $L$. The following four statements are equivalent:
\medskip
a) $P^1_H$ is a subgroup of $P^2_H$;
\smallskip
b) ${\rm Lie}(P^1_H)\subset {\rm Lie}(P^2_H)$;
\smallskip
c) The unipotent radical $N^2_H$ of $P^2_H$ is a subgroup of the unipotent radical $N^1_H$ of $P^1_H$;
\smallskip
d) ${\rm Lie}(N^2_H)\subset {\rm Lie}(N^1_H)$.}
\medskip
{\bf Proof:} We can assume we are over an algebraically closed field. It is known that $P_H^1\cap P_H^2$ contains a maximal torus of $H$ (see [Bo2, 14.22 (ii)]). So the fact that a) is equivalent to b) follows from [Bo2, 14.17-18]. Similarly we get that c) is equivalent to d). From loc. cit. we also get that a) implies c). We now assume c) holds. From [Bo2, 14.22 (i)] we get that $P_H^1\cap P_H^2$ is a parabolic subgroup of $H$. Using [Bo2, 14.17-18] once more, we get that a) holds. This ends the proof.
\medskip
{\bf 2.2.12. Definitions.} Let $\bigl({\got g},\vph,F^0({\got g}),F^1({\got g})\bigr)$ be a Shimura filtered Lie $\sg$-crystal. We use the notations of 2.2.3 3) for $({\got g},\vph)$. $\bigl({\got g},\vph,F^0({\got g}),F^1({\got g})\bigr)$ is said to be:
\smallskip
{\bf a)} of parabolic type, if ${\got p}_{\ge 0}\subset F^0({\got g})$;
\smallskip
{\bf b)} of Borel type, if there is a Borel Lie subalgebra {\got b} of {\got g} contained in $F^0({\got g})$ and such that $\vph({\got b})\subset{\got b}$;
\smallskip
{\bf b$^\prime$)} of strong Borel type, if it is of Borel type and if the parabolic Lie subalgebra of {\got g} corresponding to non-negative slopes of $({\got g},\vph)$ is a Borel Lie subalgebra (so there is only one parabolic Lie subalgebra of {\got g} contained in $F^0({\got g})$ and taken by $\vph$ into itself);
\smallskip
{\bf c)} cyclic diagonalizable (if the permutation $\pi$ of ii) below is $1_A$, then we drop the word cyclic), if there is a $W(k)$-basis $\{e_i|i\in A\}$ of $\got g$, with $A:=S(1,\dim_{W(k)}(\got g))$ assumed to be non-empty, such that the following three conditions are satisfied
\medskip
\item{i)} for subsets $C\subset B\subset A$, $\{e_i|i\in C\}$ is a $W(k)$-basis of $F^1({\got g})$ and $\{e_i|i\in B\}$ is a $W(k)$-basis of $F^0({\got g})$;
\smallskip
\item{ii)} there is a permutation $\pi$ of $A$ for which we have $\vph(e_i)=p^{\vep_i}e_{\pi(i)}$, with $\vep_i\in\{-1,0,1\}$ equal to 1 if $i\in C$, equal to $0$ if $i\in B\bsl C$ and equal to $-1$ if $i\in A\bsl B$; 
\smallskip
\item{iii)} $\forall i,\,\ j\in A$, either $[e_i,e_j]\in W(k)e_{k(i,j)}$, for some $k(i,j)\in A$, or a power of $\vph$ acts trivially on $[e_i,e_j]$.
\medskip
{\bf d)} of toric type, if the Lie subalgebra of {\got g} corresponding to the slope $0$ of $({\got g},\vph)$ is contained in $F^0({\got g})$.
\medskip
{\bf 2.2.12.1. Remarks.} {\bf 1)} We refer to the context of 2.2.12 a). As ${\got p}_{\ge 0}\subset F^0({\rm Lie}(G))$, at the level of their unipotent radicals we have $F^1({\rm Lie}(G))\subset {\got p}_{>0}$. We get a filtered $\sg$-crystal 
$$({\got p}_{\ge 0},F^1({\rm Lie}(G)),\vph).$$ 
A $p$-divisible group $D$ (resp. $D_k$) over $W(k)$ (resp. over $k$) whose associated filtered $\sg$-crystal (resp. $\sg$-crystal) is $({\got p}_{\ge 0},F^1({\rm Lie}(G)),\vph)$ (resp. is $({\got p}_{\ge 0},\vph)$), is called a $p$-divisible group of $({\got g},\vph,F^0({\got g}),F^1({\got g})$) (resp. of $({\got g},\vph)$). $D_k$ is always uniquely determined (cf. [Fo, p. 152 and p. 160]). If $p\ge 3$, then $D$ itself is uniquely determined (cf. loc. cit. and [Me, ch. 4-5], or cf. 2.2.1.1 2); see also [Fo, ch. 4, \S 5]). For $p=2$, $D$ still exists: we can take it to be the direct sum of an \'etale $2$-divisible group with a $2$-divisible group having all slopes in the interval $(0,1)$ (cf. the fact that $F^1({\rm Lie}(G))\subset W_{\vep}({\got p}_{\ge 0},\vph)$ for $\vep>0$ small enough; see also 2.3.18.1 C below).
\smallskip
Similarly, in the context of 2.2.12 b), we get a filtered $\sg$-crystal $({\got b},F^1({\rm Lie}(G)),\vph)$. So the triple $({\got g}/{\got b},F^0({\got g})/{\got b},\vph)$ is a $p$-divisible object of $\Mm\Mf_{[-1,0]}(W(k))$. So ${\got p}_{>0}$ is contained in ${\got b}$ and so in $F^0({\got g})$.
\smallskip
{\bf 2)} From many points of view, the study of Shimura filtered Lie $\sg$-crystals of parabolic or Borel type, gets reduced to the study of suitable filtered $\sg$-crystals. In general, the simplest type of $p$-divisible objects of $\Mm\Mf_{[-1,1]}(W(k))$, are those which are the extension of a $p$-divisible object of $\Mm\Mf_{[-1,0]}(W(k))$ (resp. of $\Mm\Mf_{[0,1]}(W(k))$) by a $p$-divisible object of $\Mm\Mf_{[0,1]}(W(k))$ (resp. of $\Mm\Mf_{[-1,0]}(W(k))$).
\smallskip
{\bf 3)} A Shimura filtered Lie $\sg$-crystal of parabolic type is of toric type; but the converse is not true (for examples see 4.4.13 below).
\smallskip
{\bf 4)} If $k=\bar k$, then in 2.2.12 c) it is enough to assume $\vph(e_i)$ is $p^{\vep_i}e_{\pi(i)}$ times an invertible element of $W(k)$. The same applies to 2.2.1 d). 
\smallskip
{\bf 5)} Condition iii) of 2.2.12 c) is included in order to keep better track of the Lie structure (to be compared with 2.2.1 d)); we rarely make use of it. Warning: in iii) of 2.2.12 c) it is unreasonable to require $[e_i,e_j]\in\ZZ e_{k(i,j)}$ instead of $[e_i,e_j]\in W(k)e_{k(i,j)}$.  
\medskip
{\bf 2.2.13. Remarks.} Any Shimura $\sg$-crystal $(M,\vph,G)$ has attached to it (cf. b) of 2.2.8 1)) a Shimura Lie $\sg$-crystal
$({\rm Lie}(G),\vph)$. It also has attached to it a Shimura adjoint Lie $\sg$-crystal $\bigl({\rm Lie}(G^{\rm ad}),\vph\bigr)$ (we always identify ${\rm Lie}(G^{\rm ad})$ with a $W(k)$-Lie subalgebra of ${\rm Lie}(G^{\rm der})[{1\over p}]$ and so of ${\rm Lie}(G)[{1\over p}]$). Argument: if $F^1$ is a lift of $(M,\vph,G)$, then the quadruple 
$$\bigl({\rm Lie}(G^{\rm ad}),\vph,F^0({\rm Lie}(G^{\rm ad})),F^1({\rm Lie}(G^{\rm ad}))\bigr),$$ 
with $F^i({\rm Lie}(G^{\rm ad})):={\rm Lie}(G^{\rm ad})\cap F^i({\rm End}(M))[{1\over p}]$, $i=\overline{0,1}$, is a Shimura adjoint filtered Lie $\sg$-crystal (we can assume $k=\bar k$ and so we can refer to 2.2.9 8): ${\rm Lie}(G^{\rm ad}_{\ZZ_p})$ is included in ${\rm Lie}(G_{\QQ_p})$ and so is fixed by $a$; moreover $\mu({1\over p})$ normalizes ${\rm Lie}(G^{\rm ad})[{1\over p}]$). Similarly, for any Shimura filtered $\sg$-crystal we speak about the Shimura (adjoint) filtered Lie $\sg$-crystal attached to it. The same applies to the generalized Shimura context.
\smallskip
We refer to 2.2.2 4). We assume ${\got g}$ is the Lie algebra of an adjoint group $G$ over $R^\wedge$. Then we refer to $\Ml$ as a Shimura adjoint filtered Lie $F$-crystal over $R/pR$. The fact that this matches with 2.2.11 1) and 4) for $R=W(k)$ can be checked as follows. The composite homomorphism 
$$G\buildrel{\rm Ad}\over\to Aut({\got g})\hookrightarrow GL({\got g})$$ 
is a closed embedding (for the passage from $B(k)$ to $W(k)$, cf. [Va2, 3.1.2.1 c)]) and under it $G_{B(k)}$ is the connected component of the origin of $Aut({\got g})_{B(k)}$. So the canonical split cocharacter of $({\got g},F^0({\got g}),F^1({\got g}),\vph)$ factors through $Aut({\got g})$ (cf. the functorial aspect of it and the Lie structure of $\Ml$) and so through $G$. We also get: for $R=W(k)$, $G$ is uniquely determined.
\medskip
{\bf 2.2.13.1. An interpretation.} We consider a Shimura adjoint Lie $\sg$-crystal $({\got g},\vph)$. Let $G$ be as above. Then the triple $({\got g},\vph,G)$ is a quasi $p$-divisible object with a reductive structure over $k$ (starting from the above part referring to $Aut({\got g})_{B(k)}$, this can be checked as in Fact of 2.2.9 1')). Also, if $({\got g},\vph,F^0({\got g}),F^1({\got g}))$  is a lift of $({\got g},\vph)$, then the quintuple $({\got g},F^0({\got g}),F^1({\got g}),\vph,G)$ is a quasi $p$-divisible object with a reductive structure of $\Mm\Mf_{[-1,1]}(W(k))$. 
\medskip
{\bf 2.2.13.2. Terminology.} Whenever we want to emphasize $k$ in 2.2.11, we  speak about a Shimura (adjoint) (filtered) Lie $\sg$-crystal of being over $k$. If the reductive group $G$ is split, then we speak about split Shimura (adjoint) (filtered) (Lie) $\sg$-crystals. If just $G^{\rm ad}$ is split then we speak about split Shimura (adjoint) (filtered) Lie $\sg$-crystals. We also speak about Shimura (adjoint) (filtered) Lie isocrystals (over $k$), obtained from Shimura (adjoint) (filtered) Lie $\sg$-crystals, by making $p$ invertible. For instance, with the notations of 2.2.8 1), $({\rm Lie}(G)[{1\over p}],\vph)$ is a Shimura Lie isocrystal.
\medskip
{\bf 2.2.13.3. Remark.} Let $(M,F^1,\vph,G)$ be a Shimura $\sg$-crystal. Let $P$ be the parabolic subgroup of $G$ having $F^0({\rm Lie}(G))$ as its Lie algebra. Let $P^1$ be its image in $G^{\rm ad}$. From Fact 1 and Fact 2 of 2.2.9 3), as $P_{B(k)}$ is its own normalizer in $G_{B(k)}$, we get that the lifts of $(M,\vph,G)$ are in bijection with the $W(k)$-valued points of the completion of $G/P$ in its $W(k)$-valued point defined by the origin of $G$. Similarly, from 2.2.11.1 we get that the lifts of the Shimura (adjoint) Lie $\sg$-crystal attached to $(M,F^1,\vph,G)$ are in bijection with the $W(k)$-valued points of the completion of $G^{\rm ad}/P^1$ in its $W(k)$-valued point defined by the origin of $G^{\rm ad}$. As $P$ contains $Z(G)$, $G/P$ can be canonically identified with $G^{\rm ad}/P^1$. We get:
\medskip
{\bf Fact.} {\it Any lift of the Shimura (adjoint) Lie $\sg$-crystal attached to $(M,\vph,G)$ is defined uniquely (via the attachment process) by a lift of $(M,\vph,G)$.}
\medskip
{\bf 2.2.13.4. Variants.} Referring to 2.2.8 3a) and 4a), we similarly define the (adjoint) filtered $[a-b,b-a]$-Lie $\sg$-crystal attached to a $p$-divisible object with a reductive structure of $\Mm\Mf_{[a,b]}(W(k))$ and the (adjoint) $b-a$-Lie $\sg$-crystal attached to a $p$-divisible object with a reductive structure over $k$ in the range $[a,b]$. 
\medskip
{\bf 2.2.14. Truncations.} Let $n\in\NN$ and let $(M,F^1,\vph,G)$ be a Shimura filtered $\sg$-crystal. By its truncation mod $p^n$ we mean the quintuple 
$$(M/p^nM,F^1/p^nF^1,\vph,\vph_1,G_{W_n(k)}),$$ 
where $(M/p^nM,F^1/p^nF^1,\vph,\vph_1)$ is the truncation mod $p^n$ of $(M,F^1,\vph)$ (here as well as below we denote different reductions of $\vph$ still by $\vph$ and not by $\vph_0$). As in 2.2.9 6) we define the notion of isomorphism between truncations mod $p^n$ of two Shimura filtered $\sg$-crystals over $k$. 
\smallskip
For each $\al\in\Mj$ such that $t_{\al}\neq 0$, let $m(\al)\in\ZZ$ be the smallest integer such that $p^{m(\al)}t_{\al}\in\Mt(M)$; if $t_{\al}=0$ let $m(\al)=0$. Let $t_{\al}^n$ be the reduction mod $p^n$ of $p^{m(\al)}t_{\al}$: it is an element of $\Mt(M/p^nM)$. The quintuple 
$$(M/p^nM,F^1/p^nF^1,\vph,\vph_1,G_{W_n(k)}),(t_{\al}^n)_{\al\in\Mj})$$
 is called the truncation mod $p^n$ of $(M,F^1,\vph,G,(t_{\al})_{\al\in\Mj}).$  Similarly, if $V_M$ is the $\sg^{-1}$-linear endomorphism of $M$ such that $\vph\circ V_M=p1_M$, the quadruple (resp. quintuple)
$$(M/p^nM,\vph,V_M,G_{W_n(k)})$$ 
(resp. $(M/p^nM,\vph,V_M,G_{W_n(k)},(t_{\al}^n)_{\al\in\Mj})$) is called the truncation mod $p^n$ of $(M,\vph,G)$ (resp. of $(M,\vph,G,(t_{\al})_{\al\in\Mj})$). 
We refer to 
$(M/p^nM,F^1/p^nF^1,\vph,\vph_1,G_{W_n(k)})$ as a lift of $(M/p^nM,\vph,V_M,G_{W_n(k)})$; the same applies to the context with tensors.
 Above, we still write $V_M$ for its reductions mod $p^n$.
\smallskip
As in 2.2.8 7) and 2.2.9 6) we define endomorphisms and automorphisms and (in the context of tensors $1_{\Mj}$-automorphisms) of these truncations mod $p^n$.
Warning: in general, $G_{W_n(k)}$ is not the subgroup of $GL(M/p^nM)$ fixing $t_{\al}^n$, $\forall\sl\in\Mj$. So by an inner $1_{\Mj}$-automorphism of any of the mentioned truncations mod $p^n$ involving tensors, we mean an automorphism defined by an element of $G(W_n(k))$. 
\smallskip
We have variants of the above definitions, in the principally quasi-polarized context. Putting aside the part involving Verschiebung maps and tensors, we have a variant of truncations in the Shimura filtered Lie context. In particular, we speak (cf. end of 2.2.13) about inner automorphisms of Shimura adjoint (filtered) Lie $\sg$-crystals.
\medskip
{\bf 2.2.14.0. Exercise.} Let $\mu$ and $F^0$ be as in 2.2.8 c). Let $a$ be the $\sg$-linear automorphism of $M$ such that $\vph=a\circ \mu({1\over p})$. We still denote by $a$ its reduction mod $p$. If $g_1$, $g_2\in G(W(k))$ are such that the truncations mod $p$ of $(M,g_1\vph,G)$ and $(M,g_2\vph,G)$ are identical, show that $g_2^{-1}g_1$ mod $p$ belongs to the connected, smooth, unipotent subgroup $N_k$ of $G_k$ centralizing $a(F^0/pF^0)$ and $M/pM/a(F^0/pF^0)$ (so the Lie algebra of $N_k$ is the image of the $\sg$-linear endomorphism of ${\rm Lie}(G_k)$ defined by $p\vph$). Hint: we can assume $G=GL(M)$.
\medskip
{\bf 2.2.14.1. Remark.} Using the functorial aspect of the canonical split cocharacters of [Wi], 2.2.1.2 as well as such cocharacters make sense in the context of objects of $\Mm\Mf(W(k))$. In detail: any such object $Ob$ is the cokernel of an isogeny $m:Ob_1\hookrightarrow Ob_2$ between two $p$-divisible objects of $\Mm\Mf(W(k))$ (cf. Fact of 2.2.1.1 6)); using standard properties (as in the part of 2.2.1.0 referring to $A$ and $B$) of abelian categories (cf. 2.2.1.1 6)) the canonical cocharacter $\mu_{Ob}$ of the $GL$-group of the underlying $W(k)$-module of $Ob$ obtained from the canonical split cocharacters of the $GL$-groups of the underlying $W(k)$-modules of $Ob_1$ and $Ob_2$ by passing to quotients, is independent of the choice of $m$. 
\smallskip
In particular, working with $Ob:=(M/p^nM,F^1/p^nF^1,\vph,\vph_1)$ we get a canonical cocharacter 
$$\mu_{Ob}:\GG_m\to GL(M/p^nM)$$ 
factoring through $G_{W_n(k)}$ and producing a direct sum decomposition 
$$M/p^nM=F^1/p^nF^1\oplus F^0_n,$$ 
with $\be\in\GG_m(W_n(k))$ acting through $\mu_{Ob}$ trivially on $F^0_n$ and as the multiplication by $\be^{-1}$ on $F^1/p^nF^1$. $F^0_n$ is the reduction mod $p^n$ of the maximal direct summand of $M$ on which the canonical split cocharacter of $(M,F^1,\vph)$ acts trivially. Usually we use the canonical aspect of $\mu_{Ob}$ in the form: it does not depend on which Shimura filtered $\sg$-crystal $(M,F^1,\vph,G)$ we use to define $Ob$.  
\medskip
{\bf 2.2.14.2. Fontaine truncations.} The above approach to truncations of $(M,\vph,G)$ modulo powers of $p$ does not generalize to the context of $p$-divisible objects of $\Mm\Mf(W(k))$ or more generally of $\Mm\Mf^\nabla(X)$ or of $\Mm\Mf(X)$, with $X$ as in 2.2.1 c). So here we propose a second approach whose generalization is entirely trivial. Not to be long, we work directly with a $p$-divisible object $(M,(F^i(M))_{i\in S(a,b)},\vph,G)$ with a reductive structure of $\Mm\Mf_{[a,b]}(W(k))$. As for $s\in S(a,b)$, we have
$$(p^{-s}\vph)^{-1}(M)\cap M=\sum_{j=0}^{s-a} p^jF^{s-j}(M),$$
by decreasing induction on $s$ we get that the filtration $(F^i(M)/pF^i(M))_{i\in S(a,b)}$ of $M/pM$ is uniquely determined by $(M,\vph,G)$. Based on this, following entirely the proofs of 2.2.9 3), we also get: 
\medskip
{\bf Fact.} {\it Any other filtration $(F_1^i(M))_{i\in S(a,b)}$ of $M$ such that $(M,(F_1^i(M))_{i\in S(a,b)},\vph,G)$ is a $p$-divisible object with a reductive structure of $\Mm\Mf_{[a,b]}(W(k))$, is uniquely determined by $(M,\vph,G)$ up to conjugation by elements of $G(W(k))$ which take $F^i(M)$ into $\sum_{j=0}^{i-a} p^jF^{i-j}(M)$, $\forall i\in S(a,b)$.}
\medskip
Obviously, such elements of $G(W(k))$ normalize $(F^i(M)/pF^i(M))_{i\in S(a,b)}$. Warning: if $(M,\vph,G)$ is a generalized Shimura $p$-divisible object over $k$, then any element of $G(W(k))$ normalizing $(F^i(M)/pF^i(M))_{i\in S(a,b)}$ is such that, $\forall i\in S(a,b)$, it takes $F^i(M)$ into $\sum_{j=0}^{i-a} p^jF^{i-j}(M)$; but this does not hold in the general case. 
\smallskip
Let $\Mf$ be the set of all filtrations $(F_1^i(M))_{i\in S(a,b)}$ of $M$ as in the Fact. By Fontaine truncation mod $p^n$ of $(M,\vph,G)$ we mean the set 
$$F_n((M,\vph,G))$$ 
of isomorphism classes of truncations mod $p^n$ of lifts of $(M,\vph,G)$ defined by elements of $\Mf$. Similarly we define $F_n((M,\vph,G,(t_{\al})_{\al\in\Mj}))$ (resp. $F_n^{\rm inner}((M,\vph,G,(t_{\al})_{\al\in\Mj}))$) by speaking about isomorphism classes (resp. inner isomorphism classes) of truncations mod $p^n$ of (arbitrary) lifts of $(M,\vph,G,(t_{\al})_{\al\in\Mj})$. Here inner isomorphisms refers as: isomorphisms between two such lifts defined by elements of $G(W_n(k))$.
\smallskip
A last thing: the Fact of 2.2.13.3 can be entirely adapted to the context of lifts of $(M,\vph,G)$ and of its attached adjoint. 
\medskip
{\bf 2.2.14.3. Weak truncations.} If $(M,\vph,G,(t_{\al})_{\al\in\Mj})$ is a $p$-divisible object with a reductive structure of $\Mm\Mf_{[0,a]}(W(k))$, with $a\in\NN\cup\{0\}$, by its weak truncation mod $p^n$ we mean the quadruple $(M/p^nM,\vph,G_{W_n(k)},(t_{\al}^n)_{\al\in\Mj})$, with $t_{\al}^n$'s defined as in 2.2.14, and with $\vph$ still denoting its reduction mod $p^n$.
\medskip
{\bf 2.2.15. Comment.} Though G. Shimura never dealt with the type of crystals we have introduced, we still feel that the terminology Shimura crystals used in 2.2-14 is appropriate, as G. Shimura was the first one (see the report [Sh]) to consider abelian varieties endowed with (algebraic) cycles whose adequate integral crystalline realization are providing the very first examples of the type of $\sg$-crystals with extra structure defined in 2.2.8 2). Similarly we think the terminology Shimura group or Lie pair is justified. One could argue if defs. 2.2.8 5) and 6) are justified. We think it is a convenient terminology, very much related to the one introduced in 2.2.8 1) to 4); the study of such isocrystals, to our knowledge, was started (motivated by the context of Shimura varieties of Hodge type) implicitly in [Ko1]. 
\medskip  
{\bf 2.2.16. Lemma.} {\it We assume $k=\bar k$. Let $(M,F^1,\vph)$ be a filtered $\sg$-crystal over $k$, with $F^1$ a proper direct summand of $M$. It is cyclic diagonalizable iff there is a torus $T$ of $GL(M)$ such that $(M,F^1,\vph,T)$ is a Shimura filtered $\sg$-crystal.}
\medskip
{\bf Proof:}
If $(M,F^1,\vph)$ is cyclic diagonalizable, then we take $T$ to be the maximal torus of $GL(M)$ which normalizes the $W(k)$-submodule of $M$ generated by $e_i$, $\forall i\in S(1,{\rm dim}_{W(k)}(M))$, where $\{e_1,...,e_{{\dim}_{W(k)}(M)}\}$ is a $W(k)$-basis of $M$ as in 2.2.1 d). ${\rm Lie}(T)$ is included in $F^0({\rm End}(M))$ and is normalized by $\vph$; so it has a $W(k)$-basis formed by elements fixed by $\vph$. So the canonical split cocharacter $\mu:\GG_m\to GL(M)$ of $(M,F^1,\vph)$ factors through $T$. From the Fact of 2.2.9 1), we get: $(M,F^1,\vph,T)$ is a Shimura filtered $\sg$-crystal.  
\smallskip
If there is a torus $T_{W(k)}$ of $GL(M)$ such that $(M,F^1,\vph,T_{W(k)})$ is a Shimura filtered $\sg$-crystal, then the proof of cyclic diagonalizability of $(M,F^1,\vph)$ is entirely as of 4.1.2 below. In other words, if $T$ is the $\ZZ_p$-structure of $T_{W(k)}$ obtained as in 2.2.9 8), and if $W(k(v))$, with $k(v)$ a finite field, is the smallest Witt subring of $W(k)$ such that the cocharacter $\mu:\GG_m\to T_{W(k)}$ (obtained as in c) of 2.2.8 1)) is obtained from one of $T_{W(k(v))}$ by pull back, then we can move from ``the over $k$" context to ``the over $k(v)$" context: not to introduce extra notations we assume $k=k(v)$. So we can entirely copy (so this is not a forwarding but a convenient replacement) the notations of 4.1.1.1-4 which pertain just to the quintuple $(M,F^1,\vph,T,k(v))$ to get as in 4.1.2, that $(M,F^1,\vph)$ is a cyclic diagonalizable filtered $\sg$-crystal. The only difference: as we are not necessarily in a principal quasi-polarized context, the set $I^0_{\bar\vph}$ (see 4.1.1.4) can have 1 element as well. This proves the Lemma.
\medskip
{\bf 2.2.16.1. Remark.} In practice we take $T$ to be the smallest torus of $GL(M)$ such that ${\rm Lie}(T)$ is normalized by $\vph$ and contains $\mu({\rm Lie}(\GG_m))$, where $\mu:\GG_m\to GL(M)$ is the canonical split of $(M,F^1,\vph)$.
\medskip
For the case $p=2$ of 2.2.16.2-3 below see 2.3.18.2.
\medskip
{\bf 2.2.16.2. A Galois interpretation.} We assume $p\ge 3$ and $k=\bar k$. Let $(M,F^1,\vph)$ be a filtered $\sg$-crystal over $k$. Let 
$$\rho:\Gamma_k\to GL(N)(\ZZ_p)$$ 
be its attached Galois representation, obtained via Fontaine's comparison theory; so $N$ is a free $\ZZ_p$-module of rank equal to the rank of $M$ over $W(k)$. We have:
\medskip
{\bf Corollary.} {\it $(M,F^1,\vph)$ is cyclic diagonalizable iff $\rho$ factors through the group of $\ZZ_p$-valued points of a torus $T_1$ of $GL(N)$.}
\medskip
{\bf Proof:} We can assume $F^1$ is a proper direct summand of $M$. If $(M,F^1,\vph)$ is cyclic diagonalizable, then let $T$ be a maximal torus of $GL(M)$ such that $(M,F^1,\vph,T)$ is a Shimura filtered $\sg$-crystal (cf. the proof of 2.2.16). ${\rm Lie}(T)$ is generated by elements fixed by $\vph$. Any such element can be viewed as an endomorphism of $(M,F^1,\vph)$ and so (cf. the functorial aspect of Fontaine's comparison theory), as an endomorphism of $N$ fixed by ${\rm Im}(\rho)$. The Lie subalgebra $L$ of ${\rm End}(N)$ generated by such endomorphisms is the Lie algebra of a maximal torus $T_1$ of $GL(N)$: this can be seen moving from $GL(N)$ to $GL(N\otimes_{\ZZ_p} W(k))$ and using the fact that ${\rm Lie}(T)=L\otimes_{\ZZ_p} W(k)$ is generated by projectors of $M$ onto $W(k)$-submodules of $M$ of rank one. As $T_1$ is the centralizer of $L$ in $GL(N)$, $\rho$ factors through the $\ZZ_p$-valued points of $T_1$. 
\smallskip
The converse can be proved in entirely the same manner: we can assume $T_1$ is a maximal torus of $GL(N)$. This proves the Corollary.
\medskip
{\bf 2.2.16.2.1. Remark.} Referring to the proof of 2.2.16.2, we have a natural identification of $T_{1W(k)}$ with $T$: both can be recovered as the invertible elements of ${\rm Lie}(T)=L\otimes_{\ZZ_p} W(k)$ (viewed as endomorphisms of $M$ or of $N\otimes_{\ZZ_p} W(k)$).
\medskip
{\bf 2.2.16.3. Corollary.} {\it We assume $p\ge 3$ and $k=\bar k$. If we have a proper direct sum decomposition $(M,F^1,\vph)=(M_1,F_1^1,\vph_1)\oplus (M_2,F_2^1,\vph_2)$, then $(M,F^1,\vph)$ is cyclic diagonalizable iff $(M_1,F_1^1,\vph_1)$ and $(M_2,F_2^1,\vph_2)$ are cyclic diagonalizable.}
\medskip
{\bf Proof:} This is implied by the Corollary of 2.2.16.2.
\medskip
{\bf 2.2.16.4. Degree of definition.} Let $(M,F^1,\vph)$ be a cyclic diagonalizable $\sg$-crystal. Let $\mu:\GG_m\to GL(M)$ be the canonical split cocharacter of it. Let $M=F^1\oplus F^0$ be the direct sum decomposition it produces (so $\be\in\GG_m(W(k))$ acts through $\mu$ as the multiplication with $\be^{-i}$ on $F^i$, $i=\overline{0,1}$, cf. 2.1). Let $h$ be the projector of $M$ on $F^1$ associated to it. Let $d\in\NN$ be the smallest number such that $\vph^d(h)=h$. The existence of $d$ can be seen moving to $\bar k$ and applying the part of the proof of 2.2.16 referring to $k(v)$: from very definitions we have $d=[k(v):\FF_p]$. We refer to $d$ as the degree of definition of $(M,F^1,\vph)$. Passing to a perfect field containing $k$, it remains the same. Let $a$ be the $\sg$-linear automorphism of $M$ such that $\vph=a\circ \mu({1\over p})$. 
\medskip
{\bf Definition.} $(M,F^1,\vph)$ is called strongly cyclic diagonalizable, if the $\ZZ_p$-submodule $M^a$ of $M$ formed by elements fixed by $a$, has rank equal to $\dim_{W(k)}(M)$.  
\medskip
It is easy to see that the fact that $(M,F^1,\vph)$ is strongly cyclic diagonalizable implies $k$ contains $\FF_{p^d}$. Moreover, $d$ has the following (geometric) interpretation. 
\medskip
{\bf Fact.} {\it We assume $k=\bar k$. We have: 
\medskip
a) $(M,F^1,\vph)$ is the extension from $\FF_{p^d}$ to $k$ of a strongly cyclic diagonalizable filtered $\sg_{\FF_{p^d}}$-crystal;
\smallskip
b) if $(M,F^1,\vph)$ is the extension from $\FF_{p^m}$ to $k$ of a strongly cyclic diagonalizable filtered $\sg_{\FF_{p^m}}$-crystal $(\tilde M,F^1(\tilde M),\tilde\vph)$, then $d|m$.}
\medskip
{\bf Proof:} a) is a consequence of the part of 2.2.9 8) referring to the case of a torus and of the above interpretations of $d$ in terms of $k(v)$. For b), we first remark that $\tilde\vph^m$ is a $W(\FF_{p^m})$-linear endomorphism of $\tilde M$. Let $T$ be as in the proof of 2.2.16 and let $L$ be as in 2.2.16.2.1; so we identify $L$ with the set of elements of ${\rm Lie}(T)$ fixed by $\vph$. But, due to our hypothesis, any $e\in L$ belongs to ${\rm End}(\tilde M)$ and is fixed by the $B(\FF_{p^m})$-valued point of $GL(\tilde M)$ defined by $\tilde\vph^m$. So $\tilde\vph^m$ is a $B(\FF_{p^m})$-valued point of the torus of $GL(\tilde M)$ whose extension to $W(k)$ is $T$; so it fixes the image of the canonical split cocharacter of $(\tilde M,F^1(\tilde M),\tilde\vph)$. We get: $\vph^m(h)=\tilde\vph^m(h)=h$ and so $d|m$. This ends the proof.
\medskip
{\bf 2.2.16.5. About the (adjoint) Lie case.} All the above part of 2.2.16 can be adapted to the context of Shimura (adjoint) filtered Lie $\sg$-crystals; the only difference: the restriction $p\ge 3$ encountered in 2.2.16.2-3 has to be replaced by the restriction $p\ge 5$. We do not present here the details: we just mention some new important features, in the larger context of filtered Lie $\sg$-crystals.
\smallskip
We start with a filtered Lie $\sg$-crystal $\bigl({\got g},\vph,F^0({\got g}),F^1({\got g})\bigr)$. We assume $p\ge 5$ and ${\got g}\neq \{0\}$ but we do not require $k$ to be algebraically closed. Let $\rho:\Gamma_k\to GL(N)(\ZZ_p)$ be its attached Galois representation (cf. [Fa1, 2.6] applied modulo positive, integral powers of $p$). $N[{1\over p}]$ (resp. $N$) has --via Fontaine's comparison theory-- a natural Lie algebra structure if $p\ge 5$ (resp. if $p\ge 7$). Let $G_N$ be the Zariski closure in $GL(N)$ of the subgroup of $GL(N[{1\over p}])$ formed by Lie automorphisms of $N[{1\over p}]$. $\rho$ factors through the group of $\ZZ_p$-valued points of $G_N$. In general the structure of $G_N$ can be very complicated. So for simplicity, in all that follows we assume $\bigl({\got g},\vph,F^0({\got g}),F^1({\got g})\bigr)$ is a Shimura adjoint filtered Lie $\sg$-crystal. We have:
\medskip
{\bf Fact.} {\it In such a case ${\rm Lie}(G_N)[{1\over p}]$ is isomorphic to $N[{1\over p}]$ and so is a $\QQ_p$-form of ${\got g}\otimes_{W(k)} \overline{B(k)}$.}
\medskip
{\bf Proof:} The first (resp. second) part is a consequence of the fact that any semisimple Lie algebra over an algebraically closed field of characteristic $0$ has no differentiations which are not inner, i.e. which are not defined by its elements (resp. of Fontaine's comparison theory).
\medskip
It is worth pointing out: [Va2, 4.3.10 b)] can be used to get criteria when the Zariski closure in $G_N$ of the connected component of the origin of the generic fibre of $G_N$ is an adjoint group over $\ZZ_p$. However, from the point of view of cyclic diagonalizability this is not relevant: if $k=\bar k$ and moreover $\bigl({\got g},\vph,F^0({\got g}),F^1({\got g})\bigr)$ is cyclic diagonalizable, we consider the Zariski closure $T_0$ in $G_N$ of the connected component of the origin of the generic fibre of the intersection of $G_N$ with some maximal torus $T_1$ of $GL(N)$; so $T_0$ is (for instance, see [Va2, 4.3.9]) a subtorus of $T_1$ and so of $G_N$. Here $T_1$ is obtained as in the proof of 2.2.16.2 from a suitable maximal torus $T$ of $GL({\got g})$; so (cf. the proof of 2.2.16) $T$ is a maximal torus of $GL({\got g})$ whose Lie algebra is contained in $F^0({\rm End}({\got g}))$ and is mapped by $\vph$ onto itself. 
\smallskip
We now assume that $k$ is arbitrary (so $p\ge 2$). 
Let $G$ be the adjoint group over $W(k)$ of whose Lie algebra is ${\got g}$, cf. end of 2.2.13. If $\bigl({\got g},\vph,F^0({\got g}),F^1({\got g})\bigr)$ is cyclic diagonalizable, we can work as well the previous paragraph in the crystalline context: we consider the Zariski closure $T_G$ in $G$ of the connected component of the origin of the generic fibre of the intersection of $G$ with $T$. As above, $T_G$ is a torus of $G$; as in the proof of 2.2.16 we get that the canonical split of $\bigl({\got g},\vph,F^0({\got g}),F^1({\got g})\bigr)$ is a cocharacter of $GL({\got g})$ factoring through $T_G$. So, as in 2.2.16.4, if $({\got g},\vph)$ is non-trivial we define the degree of definition, to be referred as the $A$-degree of definition, of $\bigl({\got g},\vph,F^0({\got g}),F^1({\got g})\bigr)$.
\smallskip
Also, the above Fact still holds, provided we just state it rationally from the very beginning, i.e. we work with a Galois $\QQ_p$-representation $\rho_{\QQ_p}:\Gamma_k\to GL(N[{1\over p}])$ (we do assume that such a representation does exist: in most cases pertaining to $p|6$ this is obvious; [CF] and [La, th. 3.2] imply that in fact it always exists provided we define it using the semistable version of Fontaine's ring of $W(k)$).
\smallskip
\medskip
{\bf 2.2.17. Definitions.} Let $(M,\vph,G)$ be a Shimura $\sg$-crystal over $k$. By a lift of it of quasi CM type we mean a direct summand (filtration) $F^1$ of $M$ such that the quadruple $(M,F^1,\vph,G)$ is a Shimura filtered $\sg$-crystal with the property that the intersection of $F^0({\rm Lie}(G_{B(k)}))$ with the Lie subalgebra ${\got p}_{=0}$ of ${\rm Lie}(G_{B(k)})$ corresponding to the slope $0$ of $({\rm Lie}(G_{B(k)}),\vph)$, contains a Lie subalgebra of a maximal torus of $G_{B(k)}$ which is normalized by $\vph$. We also say: $(M,F^1,\vph,G)$ is a lift of quasi CM type of $(M,\vph,G)$. 
\smallskip
By a lift of an abelian variety $A_k$ over $k$ of quasi CM type we mean an abelian variety $A$ over $W(k)$ lifting $A_k$ and such that the filtered $\sg$-crystal of $A$ is a lift of quasi CM type of the $\sg$-crystal of $A_k$ (here we use the identifications of 2.2.9 9)). Above, the abbreviation CM refers to complex multiplication.
\medskip
{\bf 2.2.17.1. Remark.} A Shimura $\sg$-crystal can have more than one lift of quasi CM type. For instance, this is the case for the $\sg$-crystal of (a product of) supersingular elliptic curves over $\bar k$. Instead of slope ${1\over 2}$ we can work equally fast and well with any other slope ${a\over r}$, with $r\in\NN$, $r\ge 2$, $(a,r)=1$, and with any $\sg$-crystal of rank $r$ and pure slope ${a\over r}$ and which has a lift of quasi CM type: we get an infinite number of lifts of quasi CM type. These examples are very simple; but in some sense (see 4.4.13 below) they are   typical.  
\medskip
{\bf 2.2.18. Proposition.} {\it We assume $k=\bar k$. Let $F^1$ be a lift of quasi CM type of a Shimura $\sg$-crystal $(M,\vph,G)$. Then the filtered $\sg$-crystal $(M,F^1,\vph)$ is cyclic diagonalizable.}
\medskip
{\bf Proof:} 
Let $\mu:\GG_m\to GL(M)$ be the canonical split cocharacter of $(M,F^1,\vph)$. From its functoriality we deduce (cf. def. 2.2.8 1)) that $\mu$ factors through $G$. We consider (see 2.2.17) a $B(k)$-Lie subalgebra $LIE$ of $F^0({\rm Lie}(G_{B(k)}))$ normalized by $\vph$ and which is the Lie algebra of a maximal torus $T_{B(k)}$ of $G_{B(k)}$. All slopes of $(LIE,\vph)$ are $0$. So, as $k=\bar k$, we get: $LIE$ is generated by elements fixed by $\vph$. Again from the functoriality of $\mu$ we get: $\mu$ fixes such elements. So the generic fibre of $\mu$ factors through $T_{B(k)}$. Let $T^1_{B(k)}$ be the smallest subtorus of $T_{B(k)}$ such that the generic fibre of $\mu$ factors through it and its Lie algebra is normalized by $\vph$. The Zariski closure of $T^1_{B(k)}$ in $GL(M)$ is a torus $T^1$ of $G$: it is generated by $\GG_m$ subgroups of $GL(M)$, commuting among themselves and obtained from the image of $\mu$ via iterates of $\vph$. So $\mu$ factors through $T^1$. We get a Shimura filtered $\sg$-crystal $(M,F^1,\vph,T^1)$ (cf. Fact of 2.2.9 1)). Now everything results from 2.2.16. This ends the proof.
\medskip
{\bf 2.2.19. Remark.} There are plenty of examples of filtered $\sg$-crystals $(M,F^1,\vph)$, with $\dim_{W(k)}(M)=2\dim_{W(k)}(F^1)$, such that the $\sg$-crystal $(M,\vph)$ does not have a lift of quasi CM type. 
The idea of constructing such examples is: if $k=\bar k$, then a $\sg$-crystal $(M,\vph)$ which has a lift of quasi CM type is a direct sum of $\sg$-crystals having only one slope (cf. 2.2.1.1 4) and 2.2.18). Moreover, the different invariants (like Hodge numbers, etc.) of such summands are of particular type, and so very easy computable. So we just have to pass to an adequate isogeny. Here is an example involving slopes having $3$ as their denominators.
\medskip
{\bf Example 1.} Let $M$ be of rank 6 over $W(k)$. Let $\vph$ be such that it sends the elements of a $W(k)$-basis $\{e_1,...,e_6\}$ of $M$ into the $B(k)$-basis $\{pe_2,pe_3,e_1,pe_5,e_6,e_4\}$ of $M[{1\over p}]$, the ordering being preserved. Let $M^{\prime}$ be the $W(k)$-lattice of $M[{1\over p}]$ generated by $M$ and by ${{e_1+e_4}\over p}$. If $F^1:=<e_2,e_4,e_3+e_6>$, then the triple $(M^\prime,F^1,\vph)$ is a filtered $\sg$-crystal; its slopes are ${1\over 3}$ and ${2\over 3}$. From 2.2.1.1 4) and 2.2.18 we deduce ${\got C}:=(M^{\prime},\vph)$ has no lift of quasi CM type. Moreover, the $\sg$-crystal ${\got C}\oplus {\got C}^*(1)$ has a principal quasi-polarization but (the argument is the same; see also 2.2.16.3 and its $p=2$ version of 2.3.18.2 below) it has no lift of quasi CM type. 
\smallskip
Also, it is worth pointing out that the Lie subalgebra of ${\rm End}(M^\prime)$ corresponding to the slope $0$ (resp. to positive slopes) of $({\rm End}(M^\prime),\vph)$ is $L_0:={\rm End}(M_1)\oplus {\rm End}(M_2)$ (resp. is $L_{>0}:={\rm Hom}(M_2,M_1)$), where $M_1:=<e_1,e_2,e_3>$ and $M_2:=<e_4,e_5,e_6>$. But $L_0/pL_0$ and $L_{>0}/pL_{>0}$ have a non-zero intersection (inside ${\rm End}(M^\prime/pM^\prime)$); so $L_0$ is not the Lie algebra of a Levi subgroup of the parabolic subgroup of $GL(M^\prime)$ whose generic fibre has $L_0[{1\over p}]\oplus L_{>0}[{1\over p}]$ as its Lie algebra. In particular $L_0$ is not the Lie algebra of a reductive subgroup of $GL(M^\prime)$. 
\medskip
Also, we would like to recall (see Example 2 below) that there are filtered $\sg$-crystals $(M,F^1,\vph)$, with $\dim_{W(k)}(M)=2\dim_{W(k)}(F^1)$, such that there is a lattice $M^\prime$ of $M[{1\over p}]$ with the property that the pair $(M^\prime,\vph)$ is a $\sg$-crystal which can not be extended to a filtered $\sg$-crystal.  So above, we have to be careful in choosing ``adequate isogenies".
\medskip
{\bf Example 2.} Let $M$ be of rank 4 over $W(k)$. Let $\vph$ be such that it sends the elements $e_0$, $e_1$, $f_0$ and $f_1$ of a $W(k)$-basis of $M$, into $pe_0$, $e_1$, $pf_1$ and respectively $f_0$. Let $M^{\prime}$ be the lattice of $M[{1\over p}]$ generated by $e_0$, $e_1$, ${{e_0+f_0}\over p}$ and $f_1$. The $\sg$-crystal $(M^{\prime},\vph)$ can not be extended to a filtered $\sg$-crystal.
\medskip
{\bf 2.2.19.1. Definition.} A Shimura adjoint Lie $\sg$-crystal $({\got g},\vph)$ is said to be of reductive type, if $W(0)({\got g},\vph)$ is the Lie algebra of a reductive subgroup of the adjoint group whose Lie algebra is ${\got g}$.
\medskip
{\bf 2.2.19.2. Exercise.} Show that the Shimura adjoint Lie $\sg$-crystal attached to a Shimura $\sg$-crystal $(M,\vph,G)$ having a lift of CM type, is of reductive type. Hint: if $G=GL(M)$, then this follows from 2.2.18, 2.2.16 and 2.2.1.1 4); if $G$ is not $GL(M)$, use 2.2.3 3) and a relative argument of passing things from $GL(M)$ to $G$. 
\medskip
{\bf 2.2.20. Shimura $p$-divisible groups.}
Let $X$ be a regular, formally smooth, faithfully flat $\tilde W(k)$-scheme, where $\tilde W(k)$ is a DVR whose completion is $W(k)$. A Shimura $p$-divisible group over $X$ is a pair $(D,(t_{\al})_{\al\in\Mj})$ comprising from a $p$-divisible group $D$ over $X$ and a family of sections $(t_{\al})_{\al\in\Mj}$ of the essential tensor algebra of $\Mf\otimes_{\ZZ_{(p)}} \QQ$, where $\Mf$ is the first crystalline cohomology $\Mo_{X^\wedge}$-sheaf of $D\times D^t$, satisfying the following axiom: 
\medskip
{\bf AX.} There is an open, affine covering of $X$, such that for any member $U={\rm Spec}(R)$ of it whose special fibre is non-empty, denoting $M_R:=H^1_{\rm crys}(D_U/U^{\wedge})$, the mentioned sections produce sections (we keep the same notation) of $\Mt(M_R[{1\over p}])$ satisfying:
\medskip
i) they are crystalline sections, i.e. they are annihilated by the natural connection on $M_R$, are in the $F^0$-filtration of $\Mt(M_R[{1\over p}])$ defined by the natural $F^1$-filtration $F^1(M_R)$ of $M_R$, and are fixed by every $\Phi_R$-linear endomorphism $\Phi_{M_R}$ of $M_R$ defined (via $D_U$) by any Frobenius lift $\Phi_R$ of $R^\wedge$;
\smallskip
ii)  the Zariski closure of the subgroup of $GL(M_R[{1\over p}])$ fixing $t_{\al}$, $\forall\al\in\Mj$, in $GL(M_R)$ is a reductive group $\tilde G_R$ over $R^\wedge$;
\smallskip 
iii) $F^1(M_R)$ is a proper direct summand of $M_R$ in each geometric point of the special fibre of $U$. 
\medskip
{\bf 2.2.20.1. Remarks.}
{\bf 1)} Any $p$-divisible group whose associated filtered $F$-crystal has a structure of a Shimura filtered $\sg$-crystal (or of a Shimura filtered $F$-crystal) with an emphasized family of tensors, gets naturally the structure of a Shimura $p$-divisible group. 
\smallskip
{\bf 2)} The importance of this notion (of Shimura $p$-divisible groups) springs from 1.15.1. In \S1-4 we need it mostly just to simplify the presentation.
\smallskip
{\bf 3)} As in 2.2.9 6), we speak about isomorphisms and $1_{\Mj}$-isomorphisms between two Shimura $p$-divisible groups over the same base $X$.
\smallskip
{\bf 4)} Sometimes we need to speak about principally quasi-polarized Shimura $p$-divisible groups. The extra thing we need for this is: $D$ has a principal quasi-polarization $p_D$, which as a perfect form on $M_R$ is normalized by $\tilde G_R$, for any $U$ as above. Notation: $(D,(t_{\al})_{\al\in\Mj},p_D)$. Moreover we allow 2.2.20 ii) to be modified as in 2.2.9 12).
\smallskip
{\bf 5)} If $X_1$ is a regular, formally smooth, faithfully flat $\tilde W(k)$-scheme, then the pull back of a Shimura $p$-divisible group over $X$ through a $\tilde W(k)$-morphism $X_1\to X$, is also a Shimura $p$-divisible group (over $X_1$).
\smallskip
{\bf 6)} Instead of crystalline sections we can use, as a variant, \'etale sections (in the $\QQ_p$-context). But in the case when $X$ is not (the $p$-adic completion of) a pro-\'etale scheme over a smooth $\tilde W(k)$-scheme, it is more convenient to work with crystalline sections: see the way [Fa1, 7.1] is stated and proved and see the easiness of expressing an integral reductiveness condition in 2.2.20 ii). Moreover, in \S 5 we define Shimura $p$-divisible groups over larger classes of schemes (including ones of characteristic $p$) and so often the \'etale $\QQ_p$-context is not suited at all. 
\smallskip
{\bf 7)} It is convenient to allow the family of tensors $(t_{\al})_{\al\in\Mj}$ to vary under the constraint that the reductive groups $\tilde G_R$ of 2.2.20 ii) remain the same. In other words, by enlarging $\Mj$ or by passing to a subset of $\Mj$ subject to the mentioned constraint we consider that we have essentially the same Shimura $p$-divisible group. This leads to the notion of quasi-isomorphism: two Shimura $p$-divisible groups over $X$,
$E^1=(D^1,(t^1_{\al})_{\al\in\Mj^1})$ and
$E^2=(D^2,(t^2_{\al})_{\al\in\Mj^2})$, are said to be quasi-isomorphic if there is an isomorphism $t:D^1\tilde\to D^2$ such that for any $U$ as in 2.2.20, the reductive group $\tilde G^1_R$ over $R^\wedge$ we get via $E^1$ is mapped (via the isomorphism of $H^1_{\rm crys}$'s groups induced by $t_U$) isomorphically into the reductive group $\tilde G^2_R$ over $R^\wedge$ we get via $E^2$. Similarly for the principally quasi-polarized context.
\smallskip
{\bf 8)} We would like to point out that the terminologies of Shimura $p$-divisible groups and crystals do not quite match: in 2.2.20 we do not require $\tilde G_R$ to become quasi-split over geometric points of ${\rm Spec}(R/pR)$, while in 2.2.8 the accent does not fall on tensors. The reason we adopted these terminologies is of purely practical matters: in the case of Shimura $p$-divisible groups we are (besides problems involving classifications, reductions and Newton polygons) very much interested in understanding the tensors and the refined structures of the reductive groups (like $\tilde G_R)$ involved, while we mainly use Shimura $\sg$-crystals as a tool for accomplishing these goals (see 1.15 for a sample).
\smallskip
The condition iii) of 2.2.20 AX is inserted to match 2.2.9 11). The conventions of 2.2.9 11) extend naturally to the context of Shimura $p$-divisible groups.
\smallskip
{\bf 9)} To simplify the notations, we assume $\tilde W(k)=W(k)$. Let $\Md=(D,(t_{\al})_{\al\in\Mj})$ be a Shimura $p$-divisible group over $X$. We have:
\medskip
{\bf Fact.} {\it If $p\ge 3$ or if $p=2$ and $X^\wedge$ is equipped with a Frobenius lift (resp. if $p=2$ and $X^\wedge$ is not equipped with a Frobenius lift) then for any $n\in\NN$, $\Md_n:=\DD(D[p^n])$ (see 2.2.1.0) is an object of $\Mm\Mf_{[0,1]}^\nabla(X)$ with the property that the Lie object $End(\Md_n)$ of $\Mm\Mf_{[-1,1]}^\nabla(X)$ (resp. of $\Mm\Mf_{[-1,1]}^{\nabla(big-tens)}(X)$) has a natural Lie subobject ${\got g}(\Md_n)$: the underlying $R$-module of the pull back of ${\got g}(\Md_n)$ to $U$ is ${\rm Lie}(\tilde G_R)\otimes_R R/p^nR$.}
\medskip
{\bf Proof:}  We first show that the induced filtration of ${\rm Lie}(\tilde G_R)\otimes_R R/p^nR$ is defined by direct summands and ${\rm Lie}(\tilde G_R)\otimes_R R/p^nR$ is generated by $\Phi_{M_R}(F^0({\rm Lie}(\tilde G_R))\otimes_R R/p^nR)$, $\Phi_{M_R}({1\over p}F^1({\rm Lie}(\tilde G_R))\otimes_R R/p^nR)$ and $p\Phi_{M_R}({\rm Lie}(\tilde G_R))$. This is a local statement and so, by completing, we can assume $R$ is noetherian. So (by multiplying by some power of $p$) we can assume $t_{\al}\in \Mt(M_R)$, $\forall\in\Mj$. This implies that ${\rm Lie}(\tilde G_R)$ is the kernel of the evaluations of endomorphisms of $M_R$ at these sections. As ${\rm End}(M_R)$ is a $p$-adically complete and noetherian $R^\wedge$-module, we get that ${\rm Lie}(\tilde G_R)\otimes_R R/p^nR$ is the image in ${\rm End}(M_R/p^nM_R)$ of the kernel of the reduction mod $p^{n+m}$ of these evaluations, with $m\in\NN$ big enough. As these last evaluations can be viewed as morphisms between objects of $\Mm\Mf^\nabla(R)$ (and so, by forgetting the connections, of $\Mm\Mf(R)$), the statement follows from the strictness part of [Fa1, 2.1 iii)] (cf. also 2.2.1.1 6)). 
\smallskip
To end the proof, we just need to add: if $p=2$ and $X^\wedge$ is not equipped with a Frobenius lift, then  we are dealing with objects of $\Mm\Mf_{[-1,1]}^\nabla(U)$, with $U$ as in 2.2.20, which are obtained via a natural process from ones of $\Mm\Mf_{[0,1]}^\nabla(X)$, and so they canonically get glued together (cf. also 3) of 2.2.4 C).
\medskip
We refer to the pair 
$$(\Md_n,{\got g}(\Md_n))$$ 
(resp. to ${\got g}(\Md_n)$) as the crystalline counterpart (resp. as the Lie subobject) of the kernel of the multiplication by $p^n$ of $\Md$. Similarly, we use the same language in the principally quasi-polarized context. We also say: the pair $(\DD(D),(t_{\al})_{\al\in\Mj})$ is associated to $\Md$. 
\medskip
{\bf 2.2.21. A result of Faltings: universal deformations.} We consider a Shimura filtered $\sg$-crystal $(M,F^1,\vph,G,(t_{\alpha})_{\alpha\in\Mj})$ with an emphasized family of tensors. Let $P$ be the parabolic subgroup of $G$ normalizing $F^1$. Let $H$ be a smooth closed subscheme of $G$ such that we get naturally an open embedding $H\hookrightarrow G/P$ (for instance, we can take $H$ to be a nilpotent subgroup of $G$). We assume $H(W(k))$ contains the identity element of $G(W(k))$. Let $R$ and $\tilde f$ be as in 2.2.10.
\smallskip
Let $D$ be a $p$-divisible group over $W(k)$ such that its corresponding filtered $\sg$-crystal is $(M,F^1,\vph)$; for $p=2$ we assume it exists. We get a Shimura $p$-divisible group $\Md:=(D,(t_{\alpha})_{\alpha\in\Mj})$ over ${\rm Spec}(W(k))$. There is a uniquely determined Shimura $p$-divisible group $\Md_R=(D_R,(t_{\alpha})_{\alpha\in\Mj})$ over ${\rm Spec}(R)$ whose associated filtered $F$-crystal with tensors is $(M,F^1,\vph,G,H,\tilde f,(t_{\al})_{\al\in\Mj})$ and whose pull back through the closed embedding $z:{\rm Spec}(W(k))\to {\rm Spec}(R)$ defined by $x_i$'s being made to be $0$, is (identifiable with) $\Md$ itself. This is just a restatement of part of [Fa2, th. 10]. In [Fa2, p. 136-7] the uniqueness of $\Md_R$ is stated just in [Fa2, rm i) after th. 10]: however, going through the proof of [Fa2, th. 10] we do get the uniqueness part. Also it is worth pointing out that the uniqueness part follows as well from previous work: $D_{R/pR}$ is uniquely determined, cf. [BM, ch. 4] (or [dJ1, th. of intro.]). But the ideal of $p(x_1,...,x_d)$ of $R$ has nilpotent divided powers modulo any positive, integral power of $(x_1,...,x_d)$. As $R$ is complete w.r.t. the $(x_1,...,x_d)$-topology, starting from $D$, $D_{R/pR}$ and $(M,F^1,\vph,G,H,\tilde f)$ we construct uniquely $D_R$ (see [Me, ch. 4-5]). Moreover, [Fa2, rm. iii) after th. 10] can be restated as follows. We have (cf. also the comments of 2.2.9 1)) the following universal property.
\medskip
{\bf UP.} {\it For any Shimura $p$-divisible group $\Md_{R_1}$ over ${\rm Spec}(R_1)={\rm Spec}(W(k))[[y_1,...,y_m]]$, with $m\in\NN$, whose pull back through the closed embedding $z_1:{\rm Spec}(W(k))\to {\rm Spec}(R_1)$ defined by $y_j$'s being made to be $0$, is (identifiable with) $\Md$, there is a unique $W(k)$-morphism $z_R:{\rm Spec}(R_1)\to {\rm Spec}(R)$ such that $z_R\circ z_1=z$ and we have an $1_{\Mj}$-isomorphism $z_R^*(\Md_R)\tilde\to \Md_{R_1}$ which in $z_1$ achieves (i.e. is compatible with) the mentioned identifications.}
\medskip
In loc. cit., the above UP is (implicitly) stated in terms of filtered $F$-crystals with tensors. In what follows we refer to the UP above in both contexts: of Shimura $p$-divisible groups or of filtered $F$-crystals with tensors. To $\Md_R$ we refer as a (or as the) universal Shimura $p$-divisible group (defined by $\Md$); for $p\ge 3$ it is uniquely determined by the quadruple $(M,\vph,G,(t_{\alpha})_{\alpha\in\Mj})$ (see 2.3.18 why this is not always so for $p=2$). There is a natural variant of this UP in a principally quasi-polarized context; in particular, we speak about universal principally quasi-polarized Shimura $p$-divisible groups.
\medskip
{\bf 2.2.21.1. Remark.} 2.2.10 and 2.2.21 make sense for the context of not necessarily quasi-split Shimura filtered $\sg$-crystals.  
\medskip
{\bf 2.2.22. Complements on Shimura (Lie) $\sg$-crystals.} Below we deal with four such complements.
\medskip
{\bf 1) Extra terminology.} A Shimura (resp. Shimura filtered) $\sg$-crystal $(M,\vph,G)$ (resp. $(M,F^1,\vph,G)$) over $k$ is said to be: 
\medskip
-- quasi cyclic diagonalizable, if there is a torus $T$ of $G$ such that $(M,\vph,T)$ (resp. $(M,F^1,\vph,T)$) is a quasi Shimura (resp. quasi Shimura filtered) $\sg$-crystal;
\smallskip
-- cyclic diagonalizable, if there is a lift $F^1$ of it such that (resp. if) $(M,F^1,\vph)$ is cyclic diagonalizable in the sense of 2.2.1 d);
\smallskip
-- strongly cyclic diagonalizable, if there is a torus $T$ of $G$ such that $(M,\vph,T)$ (resp. $(M,F^1,\vph,T)$) is a cyclic diagonalizable Shimura (resp. Shimura filtered) $\sg$-crystal and $(M,F^1,\vph)$ is strongly cyclic diagonalizable (here, for the non-filtered context, $F^1$ is the unique direct summand of $M$ such that $(M,F^1,\vph,T)$ is a Shimura filtered $\sg$-crystal, cf. Corollary of 2.2.9 3));
\smallskip
-- potentially cyclic diagonalizable if over $\bar k$ it is cyclic diagonalizable.
\medskip
We have:
\medskip
{\bf Fact.} {\it For Shimura (filtered) $\sg$-crystals, the cyclic diagonalizability implies the quasi cyclic diagonalizability.}
\medskip
{\bf Proof:} Let $(M_0,F_0^1,\vph_0,\tilde G_1)$ be a cyclic diagonalizable Shimura filtered $\sg$-crystal. We use the notations of 2.2.9 1). We can assume $\rho_0$ is well defined; so $G_1$ contains $G_0$ (cf. Fact of 2.2.9 1')). In the first paragraph of the proof of 2.2.16, working with an arbitrary $k$, as the Weyl group of $G$ is finite, the quadruple $(M,F^1,\vph,T)$ is a quasi Shimura filtered $\sg$-crystal. So by passing to a finite field extension of $k$, we can assume we have a Shimura filtered $\sg$-crystal $(M_0,F_0^1,\vph_0,T_1)$ with $T_1$ a maximal torus of $GL(M_0)$. So the Fact follows from the Fact of 2.2.9 1') applied to the quadruple $(M_0,F^1_0,\vph_0,T_2)$, with $T_2$ as the Zariski closure of the connected component of the origin of the intersection $T_{1B(k)}\cap\tilde G_{1B(k)}$ in $\tilde G_1$ (as in 2.2.16.5 we argue $T_2$ is a torus).
\medskip
If $k=\bar k$, then the converse holds (cf. 2.2.16); this is not necessarily true if $k\neq\bar k$ (simple examples can be obtained by looking at ordinary elliptic curves over algebraic extensions of finite fields). 
\smallskip
If $G=GL(M)$, it is easy to see that the above notion of cyclic (resp. strongly cyclic) diagonalizable filtered $\sg$-crystal coincides (via 2.2.9 9)) to the previous one of 2.2.1 d) (resp. of Def. of 2.2.16.4). If $k=\bar k$, then the cyclic diagonalizability implies the strongly cyclic diagonalizability; this is not necessarily true if $k\neq\bar k$ (cf. the Fact of 2.2.16.4).
\smallskip
As in 2.2.16.4, we speak about the degree of definition of a potentially cyclic diagonalizable filtered $\sg$-crystal; in particular, we speak about the degree of definition of a potentially Shimura $\sg$-crystal $(M,\vph,T)$, with $T$ a torus. As any maximal torus of $G^{\rm ad}$ is its own centralizer in $G^{\rm ad}$, from the Fact of 2.2.9 1') we get: 
\medskip
{\bf Criterion.} {\it Referring to $(M,F^1,\vph,G)$, we assume there is a maximal torus of $G^{\rm ad}$ of whose Lie algebra is contained in $F^0({\rm Lie}(G^{\rm ad}))$ and is normalized by $\vph$ (resp. and is $W(k)$-generated by elements fixed by $\vph$). Then $(M,F^1,\vph,G)$ is potentially cyclic diagonalizable (resp. there is a maximal torus $T$ of $G$ such that $(M,F^1,\vph,T)$ is a Shimura filtered $\sg$-crystal).} 
\medskip
Following the above pattern, a Shimura filtered Lie $\sg$-crystal over $k$ is called quasi (resp. potentially) cyclic diagonalizable, if its extension to $k_1$ is cyclic diagonalizable, with $k_1$ a finite field extension of $k$ (resp. with $k_1=\bar k$). 
\smallskip
Let $n\in\NN$. The Shimura filtered $\sg$-crystal $(M,F^1,\vph,G)$ is said to be:
\medskip
-- cyclic (resp. quasi, potentially or strongly cyclic) diagonalizable of level $n$, if there is a cyclic (resp. quasi, potentially or strongly cyclic) diagonalizable Shimura filtered $\sg$-crystal $(M,F^1,g\vph,T)$, with $g\in G(W(k))$ and with $T$ of torus of $G$, such that the truncations mod $p^n$ of $(M,F^1,\vph,G)$ and of $(M,F^1,g\vph,G)$ are isomorphic.
\medskip
Similarly, we speak about Shimura filtered $\sg$-crystals $(M,F^1,\vph,G,(t_{\al})_{\al\in\Mj})$ with an emphasized family of tensors which are cyclic (resp. quasi, potentially or strongly cyclic) diagonalizable of level $n$: above we just need to replace isomorphic by $1_{\Mj}$-isomorphic (see 2.2.9 6)). 
\medskip
{\bf Exercise.} We assume $k=\bar k$. $(M,F^1,\vph,G)$ is cyclic diagonalizable of level $n$ iff all iterates under $\sg$ of the canonical split (cocharacter of $G_{W_n(k)}$) of its truncation mod $p^n$ (see 2.2.14.1) are commuting among themselves. Hint: just interpret 2.2.16 mod $p^n$ (it is [SGA3, Vol. II, p. 48] which allows us to lift tori of $G_{W_n(k)}$ and their cocharacters to tori of $G$ and respectively to their cocharacters).
\medskip
If $G=GL(M)$, it is easy to see that the truncation mod $p^n$ of $(M,F^1,\vph)$ is cyclic diagonalizable in the sense of 2.2.1 d) iff $(M,F^1,\vph,G)$ is cyclic diagonalizable of level $n$. Similarly, a Shimura filtered Lie $\sg$-crystal is said to be cyclic (resp. quasi or potentially cyclic) diagonalizable of level $n$ if its truncation mod $p^n$ is isomorphic (under an isomorphism respecting the Lie structures) to the truncation mod $p^n$ of a Shimura filtered Lie $\sg$-crystal which is cyclic (resp. quasi or potentially cyclic) diagonalizable. As an application of 2.2.16.5 we have:
\medskip
{\bf Corollary.} {\it We assume $k=\bar k$. Then $(M,F^1,\vph,G)$ is cyclic diagonalizable of level $n$ iff its attached Shimura adjoint filtered Lie $\sg$-crystal ${\got L}{\got C}$ is so.}
\medskip
{\bf Proof:} One implication is easy, as in general, the Shimura adjoint filtered Lie $\sg$-crystal attached to a strongly cyclic diagonalizable Shimura filtered $\sg$-crystal is cyclic diagonalizable (see the passage from 4.1.2 to 4.1.4 below). For the other implication, we just need to show (cf. the Exercise) that the canonical split cocharacter of the truncation mod $p^n$ of $(M,F^1,\vph,G)$ is such that its iterates under $\sg$ are commuting with each other. But to check this, it is enough to deal with their images in $G^{\rm ad}$ and so, cf. end of 2.2.13, with their images in $GL({\rm Lie}(G^{\rm ad}_{W_n(k)}))$. But the fact that these last images are commuting among themselves, is implied by the fact that ${\got L}{\got C}$ is cyclic diagonalizable of level $n$. This ends the proof.
\medskip
Also, from Exercise we get: if $k=\bar k$, $(M,F^1,\vph,G)$ is cyclic diagonalizable of level $n$ iff $(M,F^1,\vph,G,(t_{\al})_{\al\in\Mj})$ is.
\medskip
{\bf 2) Classes and isogeny classes.} Let $(M,F^1,\vph,G)$ be a Shimura filtered $\sg$-crystal. By the class of $(M,\vph,G)$ we mean the set $Cl(M,\vph,G)$ of isomorphism classes of Shimura $\sg$-crystals of the form $(M,g\vph,G)$, with $g\in G(W(k))$. By the isogeny class of $(M,\vph,G)$ we mean the set ${\rm Iso}(M,\vph,G)$ of isomorphism classes of Shimura $\sg$-crystals obtained from the set of Shimura filtered $\sg$-crystals of the form 
$$(g(M),g(F^1),\vph,G(g)),$$ 
with $g\in G(B(k))$ and with $G(g)$ as the Zariski closure of $G_{B(k)}$ in $GL(g(M))$, by forgetting the filtrations. Obviously $G(g)$ is a reductive group over $W(k)$ isomorphic to $G$ (see also the general result of [Ti2, end of 2.5]). For any representative $(g(M),\vph,G(g))$ of an isomorphism class of ${\rm Iso}(M,\vph,G)$, we have ${\rm Iso}(g(M),\vph,G)={\rm Iso}(M,\vph,G)$. If ${\rm Iso}(M,\vph,G)$ has only $1$ element, we say $(M,\vph,G)$ has constant isogeny class. As $(g(M),g(F^1),\vph,G(g))$ is isomorphic to $(M,F^1,g^{-1}{\vph}g,G)$, and as $g^{-1}{\vph}g$ can be put in the form $g_1\vph$, with $g_1\in G(W(k))$, $Cl(M,\vph,G)$ is a disjoint union of different isogeny classes. We say ${\rm Iso}(M,\vph,G)$ is an isogeny class of $Cl(M,\vph,G)$.
\smallskip
Similarly, we denote and define the class and the isogeny class of a Shimura (adjoint) Lie $\sg$-crystal and, even more generally, of a $p$-divisible object with a reductive structure over $k$.
\smallskip
If a property $\Mp$ pertaining to Shimura $\sg$-crystals defining the same class $Cl(M,\vph,G)$, is enjoyed by any representative of an isomorphism class of an isogeny class of $Cl(M,\vph,G)$ (resp. of an isomorphism class of $Cl(M,\vph,G)$) once it is enjoyed by a Shimura $\sg$-crystal $(M_1,\vph_1,G_1)$ defining this isogeny class (resp. defining $Cl(M,\vph,G)$) then we refer to it as an isogeny invariant (resp. as a class invariant) property. For instance, the automorphism class of $({\rm Lie}(G),\vph)$ is a class invariant. 
\smallskip
All above extend to the context of an emphasized family of tensors $(t_{\al})_{\al\in\Mj}$; so we speak about $Cl(M,\vph,G,(t_{\al})_{\al\in\Mj})$ (its elements are isomorphism classes up to $1_{\Mj}$-isomorphisms) and ${\rm Iso}(M,\vph,G,(t_{\al})_{\al\in\Mj})$, etc.
\medskip
{\bf Definition.} A class of Shimura adjoint Lie $\sg$-crystals is said to be torsionless, if any Shimura adjoint Lie $\sg$-crystal defining it has no automorphism of order $p$.
\medskip
For instance, if $p$ does not divide $t(G^{\rm ad})$ then the Shimura adjoint Lie $\sg$-crystal associated to $(M,\vph,G)$ has no inner automorphisms of order $p$. In \S13 we will see that there are examples (involving $A_n$ Lie types) with $k=\bar k$, when $p| t(G^{\rm ad})$ and still $Cl({\rm Lie}(G^{\rm ad}),\vph)$ is torsionless. 
\medskip
{\bf 3) Classification.} Any isocrystal over $\bar k$ can be put in a cyclic diagonalizable form, cf. Dieudonn\'e's classification of isocrystals over $\bar k$. So the simplest type of $\sg$-crystals over $k$ are those which can be put similarly in a cyclic diagonalizable form. Though def. 2.2.1 d) has a general form, for simplicity, in what follows we mostly focus our attention on filtered $\sg$-crystals (see the conventions of 2.1). However, see L and 3.15.7 A below. 
\medskip
{\bf A.} Let $n\in\NN$ and let $I_n:=S(1,n)$. Till J below we consider only $n$-tuples $(e_1,e_2,...,e_n)$ with $e_i\in\{0,1\}$, $\forall i\in I_n$. We are interested in them (cf. 2.1) only from the circular point of view: for us $(e_1,...,e_n)$ is ``equally good" as $(e_j,e_{j+1},...,e_n,e_1,...,e_{j-1})$, where $j\in S(2,n)$. By the dual of such an $n$-tuple $(e_1,e_2,...,e_n)$ we mean the $n$-tuple $(1-e_1,1-e_2,...,1-e_n)$.
\smallskip
An $n$-tuple $(e_1,...,e_n)$ is said to be in the standard form if 
$$
\sum_{i=1}^n e_i2^{n-i}\ge\sum_{i=1}^n d_i2^{n-i},\leqno (1)
$$ 
for any other $n$-tuple of the form $(d_1,...,d_n)=(e_j,...,e_n,e_1,...,e_{j-1})$, for some $j\in S(2,n)$. We denote by $(e_1^*,...,e_n^*)$ the $n$-tuple which is in the standard form, and which is obtained from the dual of $(e_1,...,e_n)$ via a circular rearranging. If $(e_1,...,e_n)=(e_1^*,...,e_n^*)$ then we say $(e_1,...,e_n)$ is self dual. If $(e_1,...,e_n)$ is self dual then $n$ is even and $\sum_{i\in I_n} e_i={n\over 2}$.
\medskip
{\bf B.} In what follows (cf. 2.1) the right lower indices of entries of an $n$-tuple are considered mod $n$. So $e_1=e_{n+1}$, etc. For any $n$-tuple $(e_1,...,e_n)$ which is in the standard form we can construct a circular diagonalizable filtered $\sg$-crystal $(M,F^1,\vph)$ as follows.
\smallskip
 {\bf Construction.} {$M$ is a free $W(k)$-module of rank $n$, $\vph(a_i)=p^{e_i}a_{i+1}$, where $\{a_i|i\in I_n\}$ is a $W(k)$-basis of $M$, while $F^1$ is the $W(k)$-submodule of $M$ generated by all $a_i$'s with $e_i=1$.} 
\smallskip
About $(M,F^1,\vph)$ we say: it is in the standard form associated to $(e_1,...,e_n)$ or that it has type $(e_1,...,e_n)$; about $(e_1,...,e_n)$ we say: it is a type of $(M,F^1,\vph)$. From the very def. 2.2.1 d) we get: any filtered $\sg$-crystal which is cyclic diagonalizable is a direct sum of circular diagonalizable ones; moreover, any filtered $\sg$-crystal which is circular diagonalizable has a type.  
\smallskip
We say $(e_1,...,e_n)$ is indecomposable if the extension of $(M,F^1,\vph)$ to $\bar k$ does not have a proper direct summand. Otherwise we say $(e_1,...,e_n)$ is decomposable. Warning: some of the types can be decomposed; example: the type $(1,0,1,0)$ ``corresponds" (over $\bar k$) to a ``direct sum" of two copies of the type $(1,0)$.
\medskip
{\bf C. Proposition.} {\it Any circular diagonalizable filtered $\sg$-crystals has a unique type.}
\medskip
{\bf Proof:} Let $(e_1,...,e_n)$ and $(d_1,...,d_n)$ be two $n$-tuples in the standard form. If $(e_1,...,e_n)=(0,...,0)$ then the statement is trivial. So we can assume $e_1=1$. Let $(M_1,F_1^1,\vph_1)$ and $(M_2,F_2^1,\vph_2)$ be respectively the circular diagonalizable filtered $\sg$-crystals having these types. We can assume that inequality $(1)$ above holds. We assume the existence of an isomorphism 
$$f:(M_1,F_1^1,\vph_1)\tilde\to (M_2,F_2^1,\vph_2).$$ 
We use a $W(k)$-basis $\{a_1,...,a_n\}$ of $M_1$ as above, as well as a similarly constructed $W(k)$-basis $\{b_1,...,b_n\}$ of $M_2$. We write $f(a_1)=\sum_{i=1}^n \al_ib_i$, with $\al_i\in W(k)$, $\forall i\in I_n$. Let $j\in I_n$ be such that $\al_j$ is non-zero mod $p$. In what follows we do not need that $f$ is an isomorphism: we just need that it is a morphism and that such a $j$ does exist. By induction on $l\in I_n$ we prove that $(e_1,...,e_l)=(d_j,...,d_{j+l-1})$; due to (1) and the standard forms of $(e_1,...,e_n)$ and $(d_1,...,d_n)$, this implies $(e_1,...,e_n)=(d_1,...,d_n)$. 
\smallskip
 As $\vph_1(a_1)=pa_2$ and $f$ is a morphism, we get $d_j=1$; so the case $l=1$ holds. Let now $l\ge 2$. We assume $(e_1,...,e_{l-1})=(d_j,...d_{j+l-2})$. This implies $f(a_{l})=\sum_{i=1}^n \al_{i,l-1}b_i$, with $\al_{j+l-1,l-1}$ non-zero mod $p$ (here $\al_{i,l-1}\in W(k)$, $\forall i\in I_n$). So if $e_l=1$, then $d_{j+l-1}$ must be also 1. If $e_l=0$ and $d_{j+l-1}=1$, then we reach a contradiction with the inequality $(1)$. So we have $e_l=d_{j+l-1}$. This ends the inductive argument and so the proof.
\medskip
From now on till the end of I below we assume $k=\bar k$; however, part of the arguments below (like the one of the below Corollary) just require that $k$ contains a ``sufficiently big" finite field.
\medskip
{\bf Corollary.} {\it $(M,F^1,\vph)$ has a principal quasi-polarization iff $(e_1,...,e_n)$ is self dual.} 
\medskip
{\bf Proof:} Let ${\got C}:=(M,F^1,\vph)$. Then ${\got C}^*(1)$ is in the standard form associated to $(e_1^*,...,e_n^*)$. So if we have an isomorphism ${\got C}\tilde\to{\got C}^*(1)$, then (cf. Proposition) $(e_1,...,e_n)$ is self dual. 
\smallskip
If $(e_1,...,e_n)$ is self dual, then ${\got C}$ and ${\got C}^*(1)$ are isomorphic. Let $\{a_1,...,a_n\}$ be a $W(k)$-basis of $M$ as above. Let $\{a_1^*,...,a_n^*\}$ be the dual $W(k)$-basis of $M^*$. The underlying $W(k)$-module of ${\got C}^*(1)$ is $M^*$ and the Frobenius endomorphism of it takes $a_i^*$ into $p^{1-e_i}a_{i+1}^*$. Let $h\in I_n$ be such that $(e_1^*,...,e_n^*)=(1-e_h,1-e_{h+1},...,1-e_n,1-e_1,...,1-e_{h-1})$. We take the smallest such $h$. We have $e_i+e_{i+h-1}=1$, $\forall i\in\NN$. So $e_i=e_{i+2h-2}$, $\forall i\in\NN$. Let $d=(2h-2,n)$. As $n$ is even, so $d$ is. We get $e_i=e_{i+d}$, $\forall i\in\NN$. As $h$ was taken to be the smallest possible value and as $d$ can not be $h-1$, we conclude $d=2h-2$. So $2h-2|n$. Let $m\in\NN$ be such that $n=m(2h-2)$.
\smallskip
If $m$ is odd, we consider an isomorphism $i_{{\got C}}:{\got C}\tilde\to{\got C}^*(1)$ that takes $a_1$ into $ua_{1+{n\over 2}}^*$, with $u\in\GG_m(W(k))$ such that $\sg^{n\over 2}(u)=-u$. It is easy to see that starting from $i_{{\got C}}$ we get a perfect alternating form $\psi:M\otimes_{W(k)} M\to W(k)(1)$.
\smallskip
If $m$ is even then $4|n$. Accordingly, for $j\in\{1,2\}$, let $e_1^j:=\sg^{n\over 2}(v_j)e_1+v_je_{1+{n\over 2}}$, with $v_j\in\GG_m(W(k))$ such that $\sg^n(v_j)=v_j$, and $v_1v_2^{-1}$ mod $p$ is not fixed by $\sg^{n\over 2}$. We get a direct sum decomposition $(M,F^1,\vph)=(M_1,F_1^1,\vph)\oplus (M_2,F_2^1,\vph)$, such that $e_1^j\in M_j$. The type of $(M_j,F_j^1,\vph)$ computed starting from $e_1^j$ is an ${n\over 2}$-tuple which is also self dual. So we can apply induction (the case $n=2$ is well known). This proves the Corollary. 
\medskip
{\bf D. The uniqueness property.} We assume now $p\ge 3$ (see 2.3.18.2 below for the case $p=2$). From very definitions, any cyclic diagonalizable filtered $\sg$-crystal is a direct sum of circular diagonalizable filtered $\sg$-crystals which have (cf. 2.2.16.3) an indecomposable type. Moreover we have:
\medskip
{\bf Proposition.} {\it Such a decomposition is unique up to isomorphisms.}
\medskip
{\bf Proof:} We start with a cyclic diagonalizable filtered $\sg$-crystal $(M,F^1,\vph)$. We can assume it has a unique slope $\alpha\in [0,1]$. Moreover, we can assume $\alpha\in (0,1)$: otherwise the uniqueness part is trivial. This implies $F^1$ is a proper summand of $M$. Let $T$ be a torus of $GL(M)$ obtained as in 2.2.16.1. Let $\rho$ and $N$ be as in 2.2.16.2. Let $T_0$ be the torus of $GL(N)$ corresponding to $T$ via Fontaine's comparison theory. $\rho$ factors through the $\ZZ_p$-valued points of $T_0$ and in fact, due to the smallness property expressed in 2.2.16.1, $T_0$ is the algebraic envelope of $\rho$. So to decompose $(M,F^1,\vph)$ into direct summands is the same as to decompose $N$ into direct summands normalized by $T_0$. Based on 2.2.16.2, the Proposition is equivalent to: $N$ can be decomposed uniquely (up to isomorphism) into a direct sum of indecomposable $T_0$-subrepresentations. But this is a consequence of the following two well known facts:
\medskip
-- the representation of $T_{0\FF_p}$ on $N/pN$ is semisimple;
\smallskip
-- any subrepresentation of the representation of $T_{0\FF_p}$ on $N/pN$ lifts to a subrepresentation of the representation of $T_0$ on $N$ defined by a direct summand of $N$. 
\medskip
This proves the Proposition.
\medskip
As in the above proof, the torus $T_0$ does not depend on which algebraically closed field $k$ we use to define it, we get: the notion of indecomposable type is well defined, without specifying such a field. 
\smallskip
Let $\Mt_n$ be the set of indecomposable $n$-tuples $(e_1,...,e_n)$ which are in the standard form. Let 
$$\Mt_{\infty}:=\cup_{n=1}^{\infty} \Mt_n.$$
The above Proposition says that to any cyclic diagonalizable filtered $\sg$-crystal $(M,F^1,\vph)$, with $M\neq\{0\}$, we can associate uniquely a formal sum
$$\sum_{\tau\in \Mt_{\infty}} m(\tau)\tau,$$
where $m(\tau)\in\NN\cup\{0\}$ is the number of distinct factors having type $\tau$ in a direct sum decomposition of $(M,F^1,\vph)$ in circular diagonalizable filtered $\sg$-crystals which have an indecomposable type. For all but finite $m(\tau)=0$; moreover, at least 1 number $m(\tau)$ is non-zero. From the Corollary of C we get:
\medskip
{\bf Corollary.} {\it $(M,F^1,\vph)$ has a principal quasi-polarization iff $m(\tau)=m(\tau^*)$, $\forall\tau\in\Mt_{\infty}$.}
\medskip
{\bf E. Examples.} We come back to $p\ge 2$.
\smallskip
a) The quadruples $(1,1,0,0)$ and $(1,0,1,0)$ are producing (under the construction of B) non-isomorphic filtered $\sg$-crystals, cf. C. Another way to see this: the Hodge numbers of their resulting $\sg^2$-crystals are not the same. Also the degree of definition of a cyclic diagonalizable filtered $\sg$-crystal having type $(1,1,0,0)$ (resp. $(1,0,1,0)$) is 4 (resp. 2), cf. def. 2.2.16.4.
\smallskip
b) There are precisely 10 isomorphism classes of cyclic diagonalizable filtered $\sg$-crystals $(M,F^1,\vph)$, with $M$ of rank $6$, whose isocrystals are ${1\over 2}$-symmetric. Four of them do involve the slopes $0$ and $1$, one involves the slopes ${1\over 3}$ and ${2\over 3}$, while the other five are involving just the slope ${1\over 2}$. Among these $5$, only $3$ of them have principal quasi-polarizations: the types $(1,1,0,1,0,0)$ and $(1,1,0,0,1,0)$ are dual to each other and so are not self dual; so the cyclic diagonalizable filtered $\sg$-crystals having any one of these two types do not have principal quasi-polarizations (cf. the Corollary of C or of D). So among these 10 classes, only 8 are principally quasi-polarizable.
\smallskip
c) An $n$-tuple $(e_1,...,e_n)$ such that $(n,\sum_{i\in I_n} e_i)=1$ is indecomposable. 
\medskip
{\bf F. Corollary.} {\it We assume $p\ge 3$. The degree of definition of a cyclic diagonalizable filtered $\sg$-crystal is equal to the least common multiple of the ranks of the circular filtered $\sg$-crystals which are direct summands of it and have an indecomposable type.}
\medskip
{\bf Proof:} We just have to show (cf. D) that if $(e_1,...,e_n)$ is an indecomposable type, then the cyclic diagonalizable filtered $\sg$-crystal $(M,F^1,\vph)$ constructed in B, has $n$ as its degree of definition. To see this, we can assume $k=\FF$. Obviously, the degree of definition $d$ of $(M,F^1,\vph)$ is a divisor of $n$. The fact that it is precisely $n$ can be read out from 4.1.1.4 below: as $(e_1,...,e_n)$ is indecomposable, the permutation $\bar\vph$ (constructed as in 4.1.1.4) is a cycle of some length $d_l$, and the $W(\FF_{p^{d_l}})$-module $_pF_{\bar f_l}$ (constructed as in 4.1.1.4) has rank 1. So $M$ has rank $d_l$. As $d_l|d$ (see 4.1.1.4), the conclusion follows. Warning: 4.1.1.4 can be read independently at any time, by just picking up the definition of $\bar h_i$'s in 4.1.1.1 below.
\medskip
{\bf G.} Let $p\ge 2$. We consider an $n$-tuple $(e_1,...,e_n)$ in the standard form. Let $(M,F^1,\vph)$ be the circular diagonalizable filtered $\sg$-crystal constructed in B. From the proof of 2.2.16 and from Corollary of 2.2.9 3), we deduce that the canonical split cocharacter $\mu:\GG_m\to GL(M)$ of it, is such that $\be\in\GG_m(W(k))$ acts through $\mu$ trivially on all those $a_i$'s for which $e_i=0$. So, from the very def. 2.2.16.4, we get: the degree of definition of $(M,F^1,\vph)$ is equal to the smallest number $h\in I_n$ such that $e_i=e_{i+h}$, $\forall i\in I_n$. We have:
\medskip
{\bf Corollary.} {\it  We assume $p\ge 3$. A type $(e_1,...,e_n)$ is indecomposable iff the circular diagonalizable filtered $\sg$-crystal having type $(e_1,...,e_n)$, has $n$ as its degree of definition.}
\medskip
{\bf Proof:} One implication follows from F. We assume now that the degree of definition is $n$. We can assume $e_1=1$. We assume we have a direct sum decomposition $(M,F^1,\vph)=(M_1,F_1^1,\vph_1)\oplus (M_2,F_2^1,\vph_2)$, with $M_1$ and $M_2$ as proper $W(k)$-submodules of $M$ and with $(M_1,F_1^1,\vph_1)$ an indecomposable cyclic diagonalizable filtered $\sg$-crystal. Let $(d_1,...,d_m)$ be its type. Its degree of definition is $m$ and divides $n$. Let $d_{s+ml}:=d_s$, where $l\in S(1,{n\over m}-1)$ and $s\in S(1,m)$. We get an $n$-tuple $(d_1,...,d_n)$. If inequality (1) of A holds, then following the first proof of C in the context of the projection of $(M,F^1,\vph)$ on $(M_1,F^1_1,\vph_1)$ we get that $(e_1,...,e_n)=(d_1,...,d_n)$. This contradicts the fact that the degree of definition of $(M,F^1,\vph)$ is $n$. If the inequality (1) does not hold, then we can assume $\sum_{i=1}^n d_i2^{n-i}>\sum_{i=1}^n e_i2^{n-i}$. Using the monomorphism $(M_1,F_1^1,\vph_1)\hookrightarrow (M,F^1,\vph)$, we similarly reach a contradiction. This ends the proof.  
\medskip
For the case $p=2$ of Corollaries of paragraphs F and G see 2.3.18.2.
\medskip
{\bf H. Combinatorics.} There are many problems of combinatorial nature involving cyclic diagonalizable filtered $\sg$-crystals. Usually such problems are either very elementary or too hard. In what follows, just to point out the flavor of the problems, we mention just two such problems of elementary nature. 
\smallskip
From $G$ above we get: the number $n\abs{\Mt_n}$ can be interpreted as the number of functions $\ZZ/n\ZZ\to\{0,1\}$ which are not periodic. So we have the following recursive formula computing the number of elements of $\Mt_n$:
$$\sum_{q|n, q\in\NN} q\abs{\Mt_q}=2^n.$$
In particular, $\abs{\Mt_1}=2$, $\abs{\Mt_4}=3$, $\abs{\Mt_6}=9$ and, for $q$ a prime, $\abs{\Mt_q}={{2^q-2}\over q}$.
\smallskip
The number of isomorphism classes of filtered $\sg$-crystals $(M,F^1,\vph)$, with $M$ of rank $2n$, which are cyclic diagonalizable and have a principal quasi-polarization, is $2^n$. We do not explain this fact here; we just mention that this number $2^n$ can be easily obtained by combining [Oo3] and 4.12.4 below (see \S 9-10 for details). This second problem is a particular case (corresponding to $G_M=GSp(M,\psi_M)$, for some principal quasi-polarization $\psi_M:M\otimes_{W(k)} M\to W(k)(1)$) of the first of the following two general problems (to be analyzed in \S 9-10).
\medskip
{\bf Problem 1.} Let $(M,F^1,\vph,G_M)$ be a Shimura filtered $\sg$-crystal. Find the number of isomorphism classes of cyclic diagonalizable filtered $\sg$-crystals of the form $(M,F^1,g\vph)$, with $g\in G_M(W(k))$.
\medskip
{\bf Problem 2.} Find the set $SDD(M,\vph,G_M)$ of degrees of definitions of the cyclic diagonalizable filtered $\sg$-crystals of the form $(M,F^1,g\vph)$, with $g\in G_M(W(k))$.
\medskip
{\bf Example.} If $G_M=GSp(M,\psi_M)$, then 
$SSD(M,F^1,G_M)$ is the set
$$\{l.c.m.[ba_1,...,ba_m]|m\in\NN,\, b\in\{1,2\}\, {\rm and}\, a_1,...,a_m\in\NN,\, {\rm with}\, 2\sum_{j=1}^m a_j=\dim_{W(k)}(M)\}.$$
This is a consequence of the Corollaries of D and F (for $p=2$ cf. also 2.3.18.2 below).   
\medskip
{\bf I. Classification and standard invariants.} Let $m\in\NN$ and let $q\in\NN\cup\infty$. Let ${\bf s}:=(d_1,...,d_m)$ be an arbitrary $m$-tuple, not necessarily in the standard form.  We assume $p\ge 3$ (see 2.3.18.2 for $p=2$). Let $(M,F^1,\vph)$ and $m(\tau)$'s be as in D. We refer to $m(\tau)$'s as the classification invariants. From them we can obtain many other invariants.
\smallskip
The number 
$$FL(M,F^1,\vph):=\sum_{\tau\in\Mt_{\infty}} m(\tau)$$ 
is referred as the factor length of $(M,F^1,\vph)$. The maximum 
$IL(M,F^1,\vph)$ 
of the numbers $m(\tau)$'s is referred as the isotype length of $(M,F^1,\vph)$. The number $SKL(M,F^1,\vph)$ of non-zero $m(\tau)$'s, is referred as the skeleton length of $(M,F^1,\vph)$. The greatest (resp. the number of) $n\in\NN$ such that there is $\tau\in\Mt_n$ with $m(\tau)\neq 0$, is referred as the cyclic length (resp. the slope denominator length) of $(M,F^1,\vph)$ and is denoted by $CL(M,F^1,\vph)$ (resp. by $SDL(M,F^1,\vph)$). These five types of length, together with the degree of definition of $(M,F^1,\vph)$ (i.e. the least common multiple of all such $n$'s), are the most useful simple types of invariants of $(M,F^1,\vph)$.
\smallskip
For any $n\in\NN$ and for every $\tau_0:=(e_1,...,e_n)\in\Mt_n$, let
${\bf s}(\tau_0)$ be the number of $i\in I_n$ such that $e_{i+j}=d_j$, for $j=\overline{1,m}$. By the ${\bf s}_q$-invariant of $(M,F^1,\vph)$ we mean the number
$${\bf s}_q(M,F^1,\vph):=\sum_{n=1}^q\sum_{\tau\in\Mt_n} m(\tau){\bf s}(\tau).$$ 
When $q=\infty$, we drop it as a right lower index of ${\bf s}$. 
\smallskip
Many of these invariants have geometric interpretation (see \S 9-10 for plenty of examples). Here we just point out that: if ${\bf s}=(1)$ then ${\bf s}(M,F^1,\vph)=\dim_{W(k)}(F^1)$, if ${\bf s}=(0)$ then ${\bf s}_1(M,F^1,\vph)$ is the $p$-rank of $(M,\vph)$, and if ${\bf s}=(1,0)$ or if ${\bf s}=(0,1)$ then ${\bf s}(M,F^1,\vph)$ computes the $a$-number of the $p$-divisible group over $k$ of whose associated $\sg$-crystal is $(M,\vph)$. 
\medskip
{\bf J. The adjoint Lie case.} Parts A to D above can be adapted to some extend to the case of Shimura adjoint filtered Lie $\sg$-crystals which are cyclic diagonalizable. There are two main new features. First, we have to consider $n$-tuples whose entries are elements of the set $\{-1,0,1\}$. Second, except for the trivial context (of 2.2.11 3)) we are never in a circular context (this can be checked starting from the part of 2.2.16.5 referring to $T_G$). So, starting from such an $n$-tuple, we can imitate the  construction of B, to get a $p$-divisible object of $\Mm\Mf_{[-1,1]}(W(k))$; but it is not at all clear how to put some Lie structure on it. So in what concerns B to D, it is more appropriate to concentrate on analyzing particular cases showing up in geometric situations, rather than trying to present a whole general theory. Another approach: we totally forget the adjoint Lie structure and proceed as in L below.
\medskip
{\bf K. Variant mod $p^m$.} Let $m\in\NN$. All of A to J above makes sense and remains true modulo $p^m$, i.e. for truncations mod $p^m$ of cyclic diagonalizable filtered $\sg$-crystals; no modifications of arguments are needed (besides using --see 2.2.1 c)-- $\vph_0$ and $\vph_1$ instead of $\vph$). 
\medskip
{\bf L. Extensions.} The arguments of A to I and of K above are such that the motivated reader can easily extend them to the context of cyclic diagonalizable $p$-divisible objects of $\Mm\Mf_{[0,m]}(W(k))$; the only modifications needed to be performed are as follows:
\medskip
i) the role of $2$ in A is just of a positive integer greater than $m$;
\smallskip
ii) the dual of an $n$-tuple $(e_1,...e_n)$ formed by elements of $S(0,m)$, is defined to be the $n$-tuple $(m-e_1,...,m-e_n)$;
\smallskip
iii) in the Corollaries of C and E above, we have to situate ourselves in the self dual context of 2.2.23 d) below;
\smallskip
iv) in paragraphs D, F and G we need $p\ge m+2$ (this can still be weaken, as their references to 2.3.18.2 for $p=2$ point out; see 2.3.18.3 below);
\smallskip
v) in H above, $2^n$ has to be replaced by $(m+1)^n$.
\medskip
Related to iv) we could add that 2.2.16-19 have analogues in the context of $p$-divisible objects $\Mm\Mf_{[0,m]}(W(k))$ (for $p\ge m+2$ no modifications of arguments are needed, just of language).
\medskip
{\bf M. The non-filtered context.} In the first proof of C, the fact that $f(F_1^1)=F_2^1$ played no role. So starting from this one can check that all invariants of $(M,F^1,\vph)$ defined in I above are in fact invariants of $(M,\vph)$ itself (in other words, if $k=\bar k$ and if $(M,F^1_1,\vph)$ and $(M,F^1,\vph)$ are both cyclic diagonalizable filtered $\sg$-crystals, then they are isomorphic). More generally, we have:
\medskip
{\bf Corollary.} {\it We assume $k=\bar k$. Any object $(N,\vph_N)$ of $p-\Mm(W(k))$ which is obtained from a cyclic diagonalizable object of $p-\Mm\Mf(W(k))$ by forgetting the filtration, is uniquely decomposed (up to isomorphism) as a direct sum of objects of $p-\Mm(W(k))$ which are obtained from circular diagonalizable objects of $p-\Mm\Mf(W(k))$ by forgetting the filtrations and which can not be decomposed further on in non-trivial direct summands having this property.}
\medskip
{\bf Proof:} We just need to show the uniqueness part. Let $$(N,\vph_N)=\oplus_{i=1}^{n^1} (N^1_i,\vph_{N^1_i})=\oplus_{i=1}^{n^2} (N^2_i,\vph_{N^2_i})$$ be two such maximal decompositions; so each $(N^j_i,\vph_{N^j_i})$ is obtained from a circular diagonalizable object ${\got C}_{i^j}$ of $p-\Mm\Mf(W(k))$ by forgetting the filtration and can not be decomposed into non-trivial direct summands which are obtainable in the same way, $j=\overline{1,2}$, $i\in S(1,n^j)$, with $n^1$, $n^2\in\NN$. Let $\tau_i^j=(e_i^j(1),...,e_i^j(a_i^j))$ be the type (defined as in A and B) of ${\got C}_{i^j}$, with $a_i^j\in\NN$. We can assume all entries of $\tau_i^j$'s are non-negative; let $a$ (resp. $b$) be the least common multiple (resp. the greatest) of $a_i^j$'s, $j=\overline{1,2}$, $i\in S(1,n^j)$. The essence of the argument is the proof of G. So let $\tau_i^j({\rm ext})$ be the $a$-tuple obtained by putting ${a\over {a_i^j}}$-copies of $\tau_i^j$ together. We are interested in them from the circular point of view and we order them w.r.t. base $b+1$ (in a way similar to inequality (1) of A involving base $2$); let $\tau_{i_0}^{j_0}=(e_1,...,e_a)$ be the greatest $a$-type among $\tau_i^j({\rm ext})$'s w.r.t. this ordering. We can assume $j_0=1$. We proceed by induction on $d:=\dim_{W(k)}(N)$; so we assume that the Corollary holds for $d$ smaller than some $r\in\NN$. We now prove it for $d=r$. 
\smallskip
Following the proof of G (and so the first proof of C) we get that $\tau_{i_0}^1({\rm ext})=\tau_{i_1}^2({\rm ext})$, for some $i_1\in S(1,n^2)$, and that $N_{i_0}^1$ projects surjectively onto $N_{i_1}^2$; here we do need to use the fact that $a_{i_0}^1$ is the smallest number $h\in\NN$ such that $e_i=e_{i+h}$, $\forall i\in S(1,a)$ (if it would not be, then using similar arguments as the ones of the last paragraph of C we reach a contradiction to the ``indecomposition" of $(N_{i_0}^1,\vph_{N_{i_0}^1})$). This implies: we can assume that $N_{i_0}^1=N_{i_1}^2$. By passing to quotients (defined by $N/N_{i_0}^1$) and using the inductive assumption, the Corollary follows. 
\medskip
Warning: the above proof does not apply to get an analogue of 2.2.16.3 in the context of $\Mm\Mf(W(k))$. 
\medskip
{\bf N. Exercise.} Show that the converse of 2.2.19.2 does not hold in general. Hint: first use the types $(1,0)$ and $(1,1,0,0)$ in the same way as in Example 1 of 2.2.19 and then use the non-filtered version of the proof of G, to argue that we do not get something which is cyclic diagonalizable.
\medskip
{\bf 4) The deformation dimension.} By the deformation dimension of the Shimura $\sg$-crystal $(M,\vph,G)$ we mean the number
$$dd((M,\vph,G)):=\dim_{W(k)}({\rm Lie}(G)/F^0({\rm Lie}(G)))=\dim_{W(k)}(F^1({\rm Lie}(G))).$$
Here $F^0({\rm Lie}(G))$ and $F^1({\rm Lie}(G))$ are defined as in 2.2.8 1) starting from any (it does not matter which one, cf. Fact 1 of 2.2.9 3)) lift $F^1$ of $(M,\vph,G)$. 
\medskip
{\bf 2.2.22.1. $a$-invariants.} Let $n\in\NN$. For any $n$-tuple $\tau=(e_1,...,e_n)$ of integers we denote by $a^+(\tau)$ (resp. by $a^-(\tau)$) the number of $j\in S(1,n)$ such that $(0,1)=(e_j,e_{j+1})$ (resp. such that $(1,0)=(e_j,e_{j+1})$). We refer to it as the $a^+$-invariant (resp. as the $a^-$-invariant) of $\tau$. The number
$$a(\tau):={{a^+(\tau)+a^-(\tau)}\over 2}\in {1\over 2}\ZZ$$
is referred as the $a$-invariant of $\tau$. 
\smallskip
We assume $k=\bar k$. For any cyclic diagonalizable $p$-divisible object ${\got C}$ of $\Mm\Mf(W(k))$ (and so implicitly for the object of $p-\Mm(W(k))$ obtained from it by forgetting its filtration) we define its $a^+$-invariant (resp. the $a^-$-invariant and the $a$-invariant) by the following two rules (cf. Corollary of M of 2.2.22 3)):
\medskip
\item{{\bf a)}} they are additive with respect to direct sum decompositions;
\smallskip
\item{{\bf b)}} if ${\got C}$ is circular and is not decomposable as a direct sum of two cyclic diagonalizable $p$-divisible objects of $\Mm\Mf(W(k))$, then it is the $a^+$-invariant (resp. the $a^-$-invariant and the $a$-invariant) of the type (it is defined as in A and B of 2.2.22 3)) of ${\got C}$.   
\medskip
{\bf Exercise.} Show that the $a$-invariant of a cyclic diagonalizable Shimura $\sg$-crystal $(M,\vph,GL(M))$ is the same as the $a$-number of the $p$-divisible group whose associated $\sg$-crystal is $(M,\vph)$. 
\medskip
If $k\neq\bar k$, by the $a^+$-invariant (or the $a^-$-invariant or the $a$-invariant) of any potentially cyclic diagonalizable Shimura (filtered) adjoint Lie $\sg$-crystal we mean the $a^+$-invariant (or the $a^-$-invariant or the $a$-invariant) of its extension to $\bar k$. It is easy to see that, under field extensions, all these invariants are the same. 
\medskip
{\bf Example.} Let $(M,F^1,\vph)$ be a cyclic diagonalizable filtered $\sg$-crystal. We use the notations of 2.2.1 d), with all $n(i)$'s belonging to $\{0,1\}$. Let $b_1$ (resp. $b_2$, $b_3$ and $b_4$) be the number of $i\in A$ such that $(n(i),n(\pi(i))$ is $(0,0)$ (resp. is $(0,1)$, $(1,0)$ and $(1,1)$). We have:
\medskip
i) $a_1:=b_2=b_3$ is the $a$-invariant of $(M,\vph)$;
\smallskip
ii) $b_1+b_4=\dim_{W(k)}(M)-2a_1$;
\smallskip
iii) the $a^+$-invariant (resp. the $a^-$-invariant) $a_2^+$ (resp. $a_2^-$) of $({\rm End}(M),\vph)$ is $b_1b_2+b_3b_4$ (resp. is $b_1b_3+b_2b_4$) and so $a_2^+=a_2^-$;
\smallskip
iv) the $a$-invariant $a_2$ of $({\rm End}(M),\vph)$ is $a_1\dim_{W(k)}(M)-2a_1^2$. 
\medskip
So fixing the rank of $M$, $a_2$ is determined by $a_1$. Moreover: knowing $a_2$, we have at most two choices for $a_1$. 
\medskip
{\bf 2.2.23. Complements on quasi-polarizations.}
In what follows we generalize the notion of principal quasi-polarization introduced in 2.2.1 c). Let $n\in\NN$.
Let $R$ be a regular, formally smooth $W(k)$-algebra. Let $\Phi_R$ be a Frobenius lift of $R^\wedge$. Let ${\got C}:=(M,(F^i(M))_{i\in S(a,b)},\vph)$ (resp. ${\got C}:=(M,(F^i(M))_{i\in S(a,b)},(\vph_i)_{i\in S(a,b)})$) be a $p$-divisible object (resp. an object) of $\Mm\Mf_{[a,b]}(R)$, whose underlying module $M$ is a projective $R^\wedge$-module (resp. projective $R/p^nR$-module). 
\medskip
{\bf A. Definitions.} a) ${\got C}$ is said to be a symmetric $p$-divisible object (resp. object) if there is a morphism ${\got C}\otimes_{R^\wedge} {\got C}\to R^\wedge(a+b)$ (resp. ${\got C}\otimes_{R/p^nR} {\got C}\to R/p^nR(a+b)$) defined (at the level of underlying modules) by a perfect symmetric form $p_M:M\otimes_{R^\wedge} M\to R^\wedge$ (resp. $p_M:M\otimes_{R/p^nR} M\to R/p^nR$). $p_M$ is referred as a symmetric principal quasi-polarization of ${\got C}$. 
\smallskip
b) ${\got C}$ is said to be an alternating $p$-divisible object (resp. object) if there is a morphism ${\got C}\otimes_{R^\wedge} {\got C}\to R^\wedge(a+b)$ (resp. ${\got C}\otimes_{R/p^nR} {\got C}\to R/p^nR(a+b)$) defined (at the level of underlying modules) by a perfect alternating form $p_M:M\otimes_{R^\wedge} M\to R^\wedge$ (resp. $p_M:M\otimes_{R/p^nR} M\to R/p^nR$). $p_M$ is referred as an alternating principal quasi-polarization or just as a principal quasi-polarization of ${\got C}$. 
\smallskip
c) If $M$ is a projective $R^\wedge$-module, then any morphism ${\got C}\otimes_{R^\wedge} {\got C}\to R^\wedge(a+b)$ defined by a symmetric form $p_M:M\otimes_{R^\wedge} M\to R^\wedge$ which becomes perfect after inverting $p$, is called a symmetric quasi-polarization of ${\got C}$. As in b) above, we have a variant where symmetric is replaced by alternating. 
\smallskip
d) ${\got C}$ is said to be a self dual $p$-divisible object (resp. object) if there is a morphism ${\got C}\otimes_{R^\wedge} {\got C}\to R^\wedge(a+b)$ (resp. ${\got C}\otimes_{R/p^nR} {\got C}\to R/p^nR(a+b)$) defined (at the level of underlying modules) by a perfect bilinear form $p_M:M\otimes_{R^\wedge} M\to R^\wedge$ (resp. $p_M:M\otimes_{R/p^nR} M\to R/p^nR$). $p_M$ is referred as a bilinear principal quasi-polarization of ${\got C}$ or as a perfect bilinear form of ${\got C}$. 
\smallskip
e) In a) to d) above we have a variant where $\Mm\Mf_{[a,b]}(R)$ is replaced by an arbitrary Fontaine category of objects.
\medskip
{\bf B. Examples.} We just restrict to situations emerging from the context of a Shimura filtered $\sg$-crystal $(M,F^1,\vph,G)$ over $k$. If there is a symmetric perfect form $p_M:M\otimes_{W(k)} M\to W(k)(1)$ such that $G$ normalizes $p_M$, then we refer to the quintuple $(M,F^1,\vph,G,p_M)$ as a symmetric principally quasi-polarized Shimura filtered $\sg$-crystal; this notion will play an important role in \S5-14. 
\smallskip
Let now ${\rm Tr}$ (resp. $KIL$) be the natural trace (resp. Killing) form on ${\rm Lie}(G)$ (resp. on ${\rm Lie}(G^{\rm ad})$). They are both fixed by $\vph$ and by inverting $p$ they become perfect. Let $g_{R^\wedge}\in G(R^\wedge)$. We consider the Lie $p$-divisible objects of $\Mm\Mf_{[-1,1]}(R)$
$${\got C}:=({\rm Lie}(G)\otimes_{W(k)} R^\wedge,g_{R^\wedge}(\vph\otimes 1),F^0({\rm Lie}(G))\otimes_{W(k)} R^\wedge,F^1({\rm Lie}(G))\otimes_{W(k)} R^\wedge)$$
and
$${\got C}^{\rm ad}:=({\rm Lie}(G^{\rm ad})\otimes_{W(k)} R^\wedge,g_{R^\wedge}(\vph\otimes 1),F^0({\rm Lie}(G^{\rm ad}))\otimes_{W(k)} R^\wedge,F^1({\rm Lie}(G^{\rm ad}))\otimes_{W(k)} R^\wedge).$$
We have symmetric quasi-polarizations
$${\rm Tr}:{\got C}\otimes_{R^\wedge} {\got C}\to R^\wedge(0)$$
and
$$KIL:{\got C}^{\rm ad}\otimes_{R^\wedge} {\got C}^{\rm ad}\to R^\wedge(0).$$
If ${\rm Tr}$ (resp. $KIL$) is a perfect form, then ${\got C}$ (resp. ${\got C}^{\rm ad}$) is a symmetric $p$-divisible object of $\Mm\Mf_{[-1,1]}(R)$ and ${\got C}/p^n{\got C}$ (resp. ${\got C}^{\rm ad}/p^n{\got C}^{\rm ad}$) is a symmetric object of $\Mm\Mf_{[-1,1]}(R)$.
\smallskip
Let now $({\got g},\vph,F^0({\got g}),F^1({\got g}))$ be a Shimura adjoint filtered Lie $\sg$-crystal over $k$. We assume all simple factors of {\got g} are either of classical Lie type or of $E_6$ or $E_7$ Lie type. If none of these simple factors are of $A_m$ or $C_m$ Lie type with $p$ dividing $2(m+1)$, or of $B_m$ Lie type with $p$ dividing $2(m-1)$, or of $D_m$ Lie type with $p$ dividing $2(2m-1)$, or of $E_6$ or $E_7$ Lie type with $p$ dividing $6$, then the Killing form on ${\got g}$ is perfect, cf. [Hu2, p. 48-9] for the $E_6$ and $E_7$ Lie types and cf. [Va2, 5.7.2 and 5.7.2.1] for the classical Lie types. Loc. cit. excludes the case $p=3$ just for the simple reason that [Va2] was focused on $p\ge 5$; but the argument for it is entirely the same (cf. also [Va2, 6.6.6]). 
\smallskip
If $p>2$ and all simple factors of {\got g} are of $B_m$ or $C_m$ Lie type with $m\ge 2$, or of $D_m$ Lie type with $m\ge 4$, then the classical trace form on {\got g} is perfect. To recall the definition of the classical trace form for the mentioned Lie types, we can assume all simple factors of {\got g} are of the same Lie type. Then the classical trace form on {\got g} is the multiple of the Killing form $KIL$ on {\got g} by some $r\in\GG_m(\QQ)$, such that its extension to $W(\bar k)$, when restricted to a simple factor $\Mf$ of ${\got g}\otimes_{W(k)} W(\bar k)$, is the trace form associated to the classical (it is uniquely determined up to isomorphism; for its complex version we refer to [He, \S 8 of ch. 3]) representation of $\Mf$ of dimension $d(m)$ (as $p>2$ {\got g} is the Lie algebra of any semisimple group of whose adjoint is the adjoint group over $W(k)$ --it is uniquely determined, cf. end of 2.2.13-- having {\got g} as its Lie algebra); here $d(m)$ is $2m+1$ or $2m$ depending on the fact that ${\got g}$ is or is not of $B_m$ Lie type. 
\smallskip
Similarly, if $p=2$ and all factors of {\got g} are of $A_{2m}$ Lie type, $m\in\NN$, then the classical trace form (it is defined similarly to the above paragraph) on {\got g} is perfect. So we get:
\medskip
{\bf Fact.} {\it If $p>3$ and {\got g} has no simple factor of $A_{pm-1}$ Lie type with $m\in\NN$, or if $p=3$ and {\got g} has no simple factor of $E_6$, $E_7$ or of $A_{pm-1}$ Lie type with $m\in\NN$, or if $p=2$ and all factors of {\got g} are of some $A_{2m}$ Lie type, $m\in\NN$, then $({\got g},F^0({\got g}),F^1({\got g}),\vph)$ is a symmetric $p$-divisible object of $\Mm\Mf_{[-1,1]}(W(k))$.}
\medskip
{\bf C. Exercise.} {\bf a)} We assume ${\rm Spec}(R/pR)$ has the $ALP$ property. Show that any bilinear principal quasi-polarization of an object ${\got C}$ of $\Mm\Mf_{[0,1]}^\nabla(R)[p]$ whose pull backs via any $W(k)$-morphism ${\rm Spec}(W(\bar k))\to {\rm Spec}(R)$ is symmetric, is in fact symmetric.
\smallskip
{\bf b)} We assume $p=2$. Show that any symmetric object ${\got C}$ of some category $\Mm\Mf_{[0,1]}(*)[2]$ is alternating.  
\medskip
{\bf Hints.} We can assume we are in the context of a) of 2.2.23 A, with $n=1$, $(a,b)=(0,1)$, $k=\bar k$ and $R=W(k)[[x_1,...,x_m]]$, for some $m\in\NN\cup\{0\}$; moreover $M$ is an $R/pR$-module. From very definitions, if $m^i=\vph_i(n^i)$ with $n^i\in F^i(M)$, we have $p_M(m^i,m^i)=0$, $i\in\{0,1\}$. But any $m\in M$ can be written as $m=a_0m^0+a_1m^1$, with $a_0$, $a_1\in R$. So, in case b) as $p_M$ is symmetric and $p=2$, we have $p_M(m,m)=0$. In case a) we can assume $\Phi_R$ takes $x_i$ into $x_i^p$. Using the fact that $p_M(m^1,m^0)=\Phi_R(p_M(n^1,n^0))$, by induction on $q\in\NN$ we get that $p_M$ is symmetric modulo the ideal $(x_1,...,x_m)^q$ of $R/2R$.   
\medskip
{\bf 2.2.24. Semisimple elements.} Let $(M,\vph,G)$ be a quasi $p$-divisible object with a reductive structure over a finite field $k$. Let $k_1$ be a finite field extension of $k$ such that the extension of $(M,\vph,G)$ to $k_1$ is a $p$-divisible object with a reductive structure over $k_1$. Let $m\in\NN$ be such that $k_1=\FF_{p^m}$. Then for any $n\in\NN$, $\vph^{mn}$ is a linear automorphism of $M[{1\over p}]$ and so (as it can be seen by extension to $B(k_1)$) it is identifiable with an element $h_{mn}\in G(B(k))$. Let $h_{mn}^s\in G(B(k))$ be its semisimple part. The slopes of $(M,\vph)$ and their multiplicities are computed in the standard way from the $p$-adic valuations and the multiplicities of the eigenvalues of $h_{mn}^s$ (see [Man]).
\medskip
{\bf 2.2.24.1.} If $k$ is arbitrary, then we can work as well with the semisimple element $h(\vph)$ of $GL(M)(B(k))$ which acts on $W(\al)(M[{1\over p}],\vph)$ as the multiplication with $p^{\al\dim_{W(k)}(M)!}$. As the slope decompositions of 2.2.3 3) are functorial w.r.t. morphisms between (latticed) isocrystals, we get (cf. def. 2.2.8 3a) and 4a)) that $h(\vph)\in G(B(k))$. The image of $h(\vph)$ in $G^{\rm ab}(B(k))$ is a class invariant (i.e. by replacing $\vph$ by $g\vph$, with $g\in G(W(k))$, it remains the same). We conclude:
\medskip
{\bf Fact.} {\it If the image of $h(\vph)$ in $G^{\rm ad}(B(k))$ is trivial, then $h(\vph)$ is uniquely determined by $Cl(M,\vph,G)$.}
\medskip
{\bf Proof:} Any element of $Z(G^{\rm der})(B(k))$ is defined by an element $z\in Z(G^{\rm der})(W(k))$ by making $p$ invertible; as $h(\vph)$ and $z$ commute, if $z$ is non-trivial then the eigenvalues of $h(\vph)z$ acting on $M[{1\over p}]$ are not all integral powers of $p$. From this the Fact follows.
\medskip\smallskip 
{\bf 2.3. The standard Hodge situation.} The whole of 2.3 is dedicated to the foundation of integral aspects of integral canonical models of Shimura varieties of Hodge type: a language, a context and many general important features are presented. Implicitly, parts of what follows represent the generalization of different parts of [Ko2] (especially of ch. 5 of loc. cit.). However, the things are far more technical then loc. cit., due to the passage from tensors of degree 2 (endomorphisms, polarizations, etc.) to tensors of arbitrary degree; here ``tensors" are thought as some realizations of Hodge cycles of abelian varieties. Later on (see  2.3.8 2)) it turns out that the context we introduce (see 2.3.4), though in the beginning might look unnatural or restricted, it is the most general one could think of in terms of Hodge embeddings in a reductive $\ZZ_p$-context. In Grothendieck's and Chevalley's spirit of foundations of algebraic geometry, we work in a context about which (at least for $p\ge 3$), in \S 6 we are able to show it is ``free" of any restriction (assumption) on tensors. We rely heavily on [Va2].  In 2.3.1-17 (resp. 2.3.18) we deal with the case $p\ge 3$ (resp. $p=2$).
\medskip
{\bf 2.3.1. The $\ZZ_{(p)}$ setting.} Let $(W,\psi)$ be a symplectic space over $\QQ$ and let 
$$f:{\rm Sh}(G,X)\hookrightarrow {\rm Sh}\bigl({\rm GSp}(W,\psi),S\bigr)
$$ 
be an injective map. Let $p\ge 3$ be a rational prime. We assume the existence of a $\ZZ_{(p)}$-lattice $L_{(p)}$ of $W$ such that $\psi$ induces a perfect form 
$$\psi:L_{(p)}\otimes_{\ZZ_{(p)}} L_{(p)}\to\ZZ_{(p)}$$ 
and the Zariski closure $G_{\ZZ_{(p)}}$ of $G$ in ${\rm GSp}(L_{(p)},\psi)$ is a reductive group over $\ZZ_{(p)}$ (if needed we can multiply $\psi$ by an element of $\GG_m(\QQ)$, cf. [Va2, 4.1.6]). We call such a $\ZZ_{(p)}$-lattice good w.r.t. the map $f$ (cf. [Va2, 5.8.3]). This assumption implies $G$ is unramified over $\QQ_p$. Let $v$ be a prime of $E(G,X)$ dividing $p$. It is unramified over $p$, cf. [Mi3, 4.7]. We always view $\CC$ as an $O_{(v)}$-algebra (cf. the definition of $E(G,X)$; for instance see [Va2, 2.6]). Let $K:={\rm GSp}(L_{(p)},\psi)(L_{(p)}\otimes_{\ZZ_{(p)}} \ZZ_p)$ and let $H:=K\cap G(\QQ_p)$. $K$ is a hyperspecial subgroup of ${\rm GSp}(W,\psi)(\QQ_p)$ (as $\psi:L_{(p)}\otimes_{\ZZ_{(p)}} L_{(p)}\to\ZZ_{(p)}$ is a perfect form) and $H$ is a hyperspecial subgroup of $G(\QQ_p)$ (as $G_{\ZZ_{(p)}}$ is a reductive group and as $H=G_{\ZZ_{(p)}}(\ZZ_p)$). Let 
$$
e:={\dim_\QQ(W)\over 2}\in\NN.
$$ 
\indent
Let $(v_\al)_{\al\in\Mj}$ be a family of homogeneous tensors of the essential tensor algebra $\Mt(W^*)$ of $W^*\oplus (W^*)^*=W^*\oplus W$, such that $G$
is the subgroup of ${\rm GSp}(W,\psi)$ fixing $v_{\al}$, $\forall\al\in\Mj$ (cf. [De4, 3.1 c)]). Sometimes it is more convenient to work with an extended family 
$$(v_{\al})_{\al\in\Mj^\prime}$$ 
(so $\Mj\subset\Mj^\prime$) of such tensors, such that $G$ is the subgroup of $GL(W)$ fixing $v_\al$, $\forall\al\in\Mj^\prime$. Also it is convenient to consider that all elements of ${\rm Lie}(Z(G))$ are part of this extended family (here we identify canonically ${\rm End}(W)$ with ${\rm End}(W^*)$). As $G(\QQ)$ contains the group of scalar automorphisms of $W$, $\forall\al\in\Mj^\prime$ we have 
$v_\al\in W^{\ast\otimes n_\al}\otimes_{\QQ} W^{\otimes n_\al}$, with $n_\al:={{{\rm deg}(v_{\al})}\over 2}$. 
\smallskip
Let $G^0_{\ZZ_{(p)}}$ be the maximal reductive subgroup of $G_{\ZZ_{(p)}}$ fixing $\psi$. So 
$$G_{\ZZ_{(p)}}^0={\rm Sp}(L_{(p)},\psi)\cap G_{\ZZ_{(p)}},$$ 
as the center of $G^0_{\QQ}$ has $-1_W$ as a $\QQ$--valued point: this is a consequence of the definition of $X$ as a $G(\RR)$-conjugacy class of homomorphisms ${\rm Res}_{\CC/\RR} \GG_m\to G_{\RR}$; so the maximal compact subgroup of ${\rm Res}_{\CC/\RR} \GG_m(\RR)$ is connected and, when viewed as a subgroup of $G(\RR)$, contains $-1_{W\otimes_{\QQ} \RR}$ and is contained in $G^0_{\ZZ_{(p)}}(\RR)$. 
\smallskip
We have a natural short exact sequence $0\to G^0_{\ZZ_{(p)}}\hookrightarrow G_{\ZZ_{(p)}}\twoheadrightarrow\GG_m\to 0$ and so a natural isogeny $\GG_m\times G^0_{\ZZ_{(p)}}\to G_{\ZZ_{(p)}}$ whose kernel is $\mu_2$.
\smallskip
To any $x\in X$ it corresponds an injective cocharacter $\mu_x:\GG_m\hookrightarrow G_{\CC}$ (see [Va2, 2.2]). We often prefer to view $G$ as a subgroup of $GL(W^*)$ or of ${\rm GSp}(W^*,\tilde\psi)$, where $\tilde\psi$ is the perfect form on $W^*$ induced naturally by $\psi$, and to emphasize this we denote $\mu_x$ by $\mu_x^\ast$; so $\mu_x^\ast$ defines, as in [Va2, 2.2], the Hodge structure on $W^*$ corresponding to $x\in X$. $\mu_x^\ast$ achieves a direct sum decomposition $W^*\otimes_{\QQ} \CC=F^1_x\oplus F^0_x$, with $\be\in\GG_m(\CC)$ acting through it on $F_x^i$ as the multiplication with $\be^{-i}$, $i=\overline{0,1}$.  
\medskip
{\bf 2.3.2. The moduli setting for the Siegel modular variety.} We get an injective map 
$$f:(G,X,H,v)\hookrightarrow\bigl({\rm GSp}(W,\psi),S,K,p\bigr)$$ 
of Shimura quadruples. Let $\Mm$ be the extension to $O_{(v)}$ of the integral canonical model of the Shimura quadruple $\bigl({\rm GSp}(W,\psi),$ $S,K,p\bigr)$. It is known (for instance see [Va2, 3.2.9 and 4.1]) that $\Mm$ parameterizes  isomorphism
classes of principally polarized abelian schemes of dimension $e$ over $O_{(v)}$-schemes having, in a compatible way, level-$N$ symplectic similitude structure for any $N\in\NN$
relatively prime to $p$; the choice of a $\ZZ$-lattice $L$ of $W$ such that $L\otimes_{\ZZ} \ZZ_{(p)}=L_{(p)}$ and  $\psi:L\otimes_{\ZZ} L\to\ZZ$ (makes sense and) is a perfect form is implicit here, cf. [Va2, 4.1]. Here, as well as everywhere else, the level structures and their compatibilities are as in [Va2, paragraph before 4.1.0].
\medskip
{\bf 2.3.3. The integral canonical model.} Let $\Mn$ be the integral canonical model of $(G,X,H,v)$
(for $p\ge 5$ cf. [Va2, 6.4.1] and the defs. of [Va2, 3.2.6]; for $p=3$ cf. \S 6). As $\Mm$ has the extension property (see [Va2, 2.3.3 3) and 6)]), we have (cf. also [Va2, 3.2.7 4)]) a natural $O_{(v)}$-morphism 
$$i_{\Mn}:\Mn\to\Mm.$$  
[Va2, 3.2.12] implies: $\Mn$ is the normalization of the Zariski closure of ${\rm Sh}_H(G,X)$ in $\Mm$. From the proof of [Va2, 3.4.1] we deduce: $i_{\Mn}$ is a finite morphism and not just pro-finite.  
\smallskip
The universal principally polarized abelian scheme over $\Mm$ gives birth, by pull back, to a principally polarized abelian scheme $(\Ma,\Mp_{\Ma})$ over $\Mn$. $\Ma$ is naturally endowed with a family $(w_{\al}^{\Ma})_{\al\in\Mj^\prime}$ of Hodge cycles. Argument: due to the existence of level structures, this can be read out from an \'etale $\QQ_l$-context, with $l$ a prime different from $p$; in other words, all Hodge cycles of the pull back of $\Ma$ to ${\rm Sh}(G,X)$ are already Hodge cycles of $\Ma$ and so [Va2, 4.1.1-2] applies. So the points of $\Mn$ with values in  fields $\Mk$ of characteristic $0$, give birth to principally polarized abelian varieties (obtained from $(\Ma,\Mp_{\Ma})$ by pull back), having a family $(w_\al)_{\al\in\Mj^\prime}$ of Hodge cycles (obtained from $(w_{\al}^{\Ma})_{\al\in\Mj^\prime}$ by pull back) and compatible level-$N$ symplectic similitude structures, for any $N\in\NN$ with $(N,p)=1$, such that some additional conditions are satisfied (cf. loc. cit.). 
\smallskip
We denote by $H_0$ an
arbitrary compact, open subgroup of $G(\AA^p_f)$ such that:
\medskip 
{\bf a)} the quotient morphism $\Mn\to\Mn/H_0$ is a pro-\'etale cover;
\smallskip
{\bf b)} there is a principally polarized abelian scheme $(\Ma_{H_0},\Mp_{\Ma_{H_0}})$ over $\Mn/H_0$, with $(\Ma,\Mp_{\Ma})$ obtained from it by pull back via the pro-\'etale cover $\Mn\to\Mn/H_0$.
\medskip
$\Mn/H_0$ is a quasi-projective, smooth $O_{(v)}$-scheme, cf. [Va2, 6.4.1]: this is a consequence of the fact that the morphism $i_{\Mn}$ is finite. Warning: we do not assume that the family of Hodge cycles $(w_{\al}^{\Ma})_{\al\in\Mj^\prime}$ is definable over $\Mn/H_0$ (i.e. is obtained by pulling back a family of Hodge cycles of $\Ma_{H_0}$); also we do not assume $H_0\times H$ is included in (some specific compact, open subgroup of) $GSp(\tilde L,\psi)(\tilde L\otimes_{\ZZ} \hat\ZZ)$ for some $\ZZ$-lattice $\tilde L$ of $W$ such that we have a perfect form $\psi:\tilde L\otimes_{\ZZ} \tilde L\to\ZZ$ (however, in many cases, the proof of [Va2, 6.7.2] can be used to show that we can assume that such an inclusion does exist). Also for the sake of future study (involving twists), we do not impose any restriction (besides b)) on $(\Ma_{H_0},\Mp_{\Ma_{H_0}})$; in particular, even if we have
$$H_0\times H\subset K(N):=\{g\in GSp(L,\psi)(L\otimes_{\ZZ} \hat\ZZ)|g\equiv 1_{L\otimes_{\ZZ} \hat\ZZ} \, {\rm mod}\, N\},\leqno (INCL)$$ for some $N\in\NN\setminus\{1,2\}$, with $(N,p)=1$, we do not assume it is the pull back of the universal principal polarized abelian scheme over $\Mm/K^p(N)$ having level-$N$ symplectic similitude structure. Here $K^p(N)\subset GSp(W,\psi)(\AA_f^p)$ is a subgroup such that we have a natural identification $K(N)=K^p(N)\times K$. However, by passing to a normal, open subgroup of $H_0$, $(\Ma_{H_0},\Mp_{\Ma_{H_0}})$ becomes isomorphic to such a pull back, under an isomorphism which over $\Mn$ it becomes the one of b). Warning: despite the notation $(\Ma_{H_0},\Mp_{\Ma_{H_0}})$, this pair does not depend only on $H_0$.
\smallskip
Let $\tilde H_0$ be an arbitrary compact, open subgroup of $G(\AA_f^p)$. We denote by $c_p(\tilde H_0)$ (resp. by $c_p(\tilde H_0,v)$) the number of connected components of the generic (resp. special) fibre of $\Mn/\tilde H_0$. From [Va2, 6.4.6 3)] (see also [Va2, 3.2.7 2)]) we get: this number is independent of which hyperspecial subgroup $H$ of $G(\QQ_p)$ we are dealing with and so the notations are justified.
\medskip
{\bf 2.3.3.1. The compact type case.} If $(G,X)$ is of compact type, then $\Mn/H_0$ is a projective, smooth $O_{(v)}$-scheme, cf. [Va2, 6.4.1.1 2)] (the proof of loc. cit. applies to the case $p=3$ as well).
\medskip
{\bf 2.3.3.2. The abelian and adjoint counterparts.} We denote by $\Mn^{\rm ad}$ (resp. by $\Mn^{\rm ab}$) the integral canonical model of the adjoint (resp. toric) Shimura quadruple $(G^{\rm ad},X^{\rm ad},H^{\rm ad},v^{\rm ad})$ (resp. $(G^{\rm ab},X^{\rm ab},H^{\rm ab},v^{\rm ab})$) of $(G,X,H,v)$. In connection to the existence of $\Mn^{\rm ad}$ for $p=3$ (resp. of $\Mn^{\rm ab}$ for $p\ge 3$) we refer to \S 6, cf. also [Va2, 6.2.2 b)] and AE.4.1 (resp. we refer to [Va2, 3.2.8]). Let $H_0^{\rm ad}$ (resp. $H_0^{\rm ab}$) be a compact, open subgroup of $G^{\rm ad}(\AA_f^p)$ (resp. of $G^{\rm ab}(\AA_f^p)$) such that under the natural epimorphism $G\twoheadrightarrow G^{\rm ad}$ (resp. $G\twoheadrightarrow G^{\rm ab}$) it contains the image of $H_0$. Warning: we do not assume $\Mn^{\rm ad}/H_0^{\rm ad}$ is smooth over $O_{(v^{\rm ad})}$. On the other hand, $\Mn^{\rm ab}/H_0^{\rm ab}$ is always an \'etale cover of ${\rm Spec}(O_{(v^{\rm ab})})$, cf. [Va2, 3.2.8]. In the same way we got the existence of $i_{\Mn}$, we get the existence of an $O_{(v)}$-morphism
$$q_{\Mn}:\Mn\to\Mn^{\rm ad}_{O_{(v)}}\times_{O_{(v)}} \Mn^{\rm ab}_{O_{(v)}}.$$ 
From [Va2, 6.4.5 b)] (cf. also [Va2, 6.2.2-3] and AE.4.1) we get that it is a pro-finite, pro-\'etale morphism (for $p=3$ cf. also \S6). If $\Mn^{\rm ad}$ is a pro-\'etale cover of the smooth $O_{(v^{\rm ad})}$-scheme $\Mn^{\rm ad}/H_0^{\rm ad}$, then the similarly defined $O_{(v)}$-morphism
$$q_{\Mn}(H_0):\Mn/H_0\to \Mn^{\rm ad}_{O_{(v)}}/H_0^{\rm ad}\times_{O_{(v)}} \Mn^{\rm ab}_{O_{(v)}}/H_0^{\rm ab}$$
is an \'etale cover of its image.
\medskip
{\bf 2.3.4. Definition.} The $\ZZ_{(p)}$-lattice $L_{(p)}$ of $W$ is said to be crystalline well positioned for the map $f$ w.r.t. $v$, if, besides the property mentioned in 
2.3.1 (i.e. it is good w.r.t. $f$), it also satisfies
the following property:
\medskip
{\bf $(\ast)$} {\it for any perfect field $k$ of characteristic $p$ and for every $z\in\Mn(W(k))$, the triple $(A,(w_\al)_{\al\in\Mj},p_A):=z^*(\Ma,(w_{\al}^{\Ma})_{\al\in\Mj},\Mp_{\Ma})$ is such that the Zariski closure in $GL(M)$
(with $M:=H^1_{\rm crys}\bigl(A/W(k)\bigr)=H^1_{dR}(A/W(k))$) of the subgroup of $GL(M[{1\over p}])$ fixing the perfect alternating form $p_A:M\otimes_{W(k)} M\to W(k)(1)$ (induced by the polarization $p_A$) and the crystalline (de Rham) component $t_{\al}$ of $w_\al$ ($t_\al$ is a homogeneous tensor of $\Mt(M[{1\over p}])$), $\forall\al\in\Mj$,
is a reductive group $\tilde G$ over $W(k)$.}
\medskip
If $L_{(p)}$ is crystalline well positioned for the map $f$ w.r.t. any prime $v$ of $E(G,X)$ dividing $p$, then we say $L_{(p)}$ is crystalline well positioned for the map $f$.
\medskip
{\bf 2.3.5. Basic definition.} A triple $(f,L_{(p)},v)$ as in 2.3.1, with $L_{(p)}$ a
$\ZZ_{(p)}$-lattice of $W$ which is crystalline well positioned for the map $f$ w.r.t. $v$, is called a standard
Hodge situation, abbreviated SHS. If moreover there is a $\ZZ_{(p)}$-subalgebra $\Mb$ of ${\rm End}(L_{(p)})$, which over $W(\FF)$ is a product of matrix $W(\FF)$-algebras, which is self dual w.r.t. $\psi$, and is such that the subgroup of ${\rm GSp}(W,\psi)$ fixing its elements has $G$ as its connected component of the origin, then the quadruple $(f,L_{(p)},v,\Mb)$ is called a standard PEL situation (we recall that PEL stands for polarization, endomorphisms and level structures). 
\medskip
{\bf 2.3.5.1. Exercise.} Let $(G_0,X_0)$ be an adjoint Shimura pair such that all its simple factors are of $A_n$ type for some $n\in\NN$. Let $p>2$ be a prime such that $G_0$ is unramified over $\QQ_p$. Show that there is an injective map $f\colon (G,X)\hookrightarrow ({\rm GSp}(W,\psi),S)$, with $(G^{\rm ad},X^{\rm ad})=(G_0,X_0)$, and a $\ZZ_{(p)}$-lattice $L_{(p)}$ of $W$ good w.r.t. $f$, such that for any prime $v$ of $E(G,X)$ dividing $p$ we get a standard PEL situation $(f,L_{(p)},v,\Mb)$, for a suitable $\ZZ_{(p)}$-subalgebra $\Mb$ of ${\rm End}(L_{(p)})$. Hint: [Va2, 6.5.1] allows us to assume $G_0$ is $\QQ$--simple and for such a case see [Va2, 6.5.1.1 and Case 2 of 6.6.5.1]; in other words, loc. cit. handles the case $p\ge 5$ and so, for $p=3$ we just have to literally copy the arguments for $p\ge 5$.
\medskip
{\bf 2.3.5.2. The factors.} Let $(f,L_{(p)},v)$ be a SHS. Let 
$$(G^{\rm ad},X^{\rm ad})=\prod_{j\in\Mf(G^{\rm ad})} (G_j,X_j)$$ 
be the direct factor decomposition in simple, adjoint Shimura pairs. Here $\Mf(G^{\rm ad})$ is a finite set of indices, indexing the $\QQ$--simple factors of $G^{\rm ad}$. We have:
\medskip
{\bf Fact 1.} {\it For any $j\in\Mf(G^{\rm ad})$ there is an injective map $f_j:(G^j,X^j)\hookrightarrow (G,X)$ such that:
\medskip
a) $G^{j{\rm der}}$ is the semisimple subgroup of $G$ naturally isogeneous to $G_j$ (so $(G^{j{\rm ad}},X^{j{\rm ad}})=(G_j^{\rm ad},X_j^{\rm ad})$);
\smallskip
b) $(f\circ f_j,L_{(p)},v_j)$ is a SHS, for any prime $v^j$ of $E(G^j,X^j)$ dividing $v$.}
\medskip
{\bf Proof:} [Va2, 3.2.7 11)] takes right away care of a) and to get b) we just need to use simple arguments. Following loc. cit., we can assume $G_j$ is generated by $G^{j{\rm der}}$ and by a maximal torus of $G$ whose Zariski closure in $G_{\ZZ_{(p)}}$ is a maximal torus of $G_{\ZZ_{(p)}}$; so (cf. [Va2, 3.1.6]) the Zariski closure of $G^j$ in $G_{\ZZ_{(p)}}$ is a reductive subgroup $G^j_{\ZZ_{(p)}}$. Let $H^j:=H\cap G^j(\QQ_p)$. Let $\Mn^j$ be the integral canonical model of the Shimura quadruple $(G^j,X^j,H^j,v^j)$. In the same way we got $i_{\Mn}$, we get naturally a finite $O_{(v^j)}$-morphism
$$i_{\Mn}^j:\Mn^j\to\Mn_{O_{(v^j)}}.$$
The fact that the triple $(f\circ f_j,L_{(p)},v_j)$ is a SHS, can be deduced immediately from [Va2, 3.1.6] via [Va2, 4.3.13] (to be compared with the proof of [Va2, 5.7.1]). To recall the details, we refer to 2.3.4, with $z$ factoring through $\Mn_j$. Defining $\tilde G_j$ in the similar way, [Va2, 4.3.13] ``takes care" of the maximal subtorus of $Z(G^j_{\ZZ_{(p)}})$ (i.e. it implies that the Zariski closure in $GL(M)$ of the connected component of the origin of $Z(\tilde G_{jB(k)})$ is a torus), while the fact that $(f,L_{(p)},v)$ is a SHS ``takes care" of the derived part of $G_{\ZZ_{(p)}}$ (i.e. it implies that the Zariski closure in $GL(M)$ of $\tilde G_{jB(k)}^{\rm der}$ is a semisimple group); it is [Va2, 3.1.6] which ``takes care" of combining the toric part with the derived part. This proves Fact 1.
\medskip
So from many points of view (for instance, see Fact 2 below), for the study of a SHS $(f,L_{(p)},v)$, we can assume $G^{\rm ad}$ is a simple $\QQ$--group; however, for the sake of generality, we never assume this a priori. 
\smallskip
Let $q_{G^{\rm ad}}$ be the order of the center of the simply connected group cover of $G^{\rm der}_{\CC}$. We have:
\medskip
{\bf Fact 2.} {\it A connected component $\Mc_0$ of $\Mn_{W(\FF)}$ is a quotient of the product $\times_{j\in\Mf(G^{\rm ad})} \Mc_0^j$, with ${\Mc}_0^j$ a connected component of $\Mn^j_{W(\FF)}$, through a pro-finite, $q_{G^{\rm ad}}^2$-torsion abelian group.}
\medskip
{\bf Proof:} First, it is irrelevant which connected components we choose, cf. [Va2, 3.3.2]. Second, we need to apply [Va2, 3.2.16] to the Shimura quadruple 
$$(\prod_{j\in\Mf(G^{\rm ad})} G_j,\prod_{j\in\Mf(G^{\rm ad})} X_j,\prod_{j\in\Mf(G^{\rm ad})} H^j,w),$$
where $w$ is the prime of the composite field of all reflex fields $E(G^j,X^j)$'s dividing all $v^j$'s. Third, we need to consider a cover $(\tilde G,\tilde X,\tilde H,\tilde v)\to (G,X,H,v)$, such that $\tilde G^{\rm der}=\prod_{j\in\Mf(G^{\rm ad})} G_j^{\rm der}$, cf. [Va2, 3.2.7 10)]. So, Fact 2 follows by combining [Va2, 6.2.3 and 6.2.3.1].
\medskip
The mentioned torsion abelian group is trivial if $G^{\rm der}=\times_{j\in\Mf(G^{\rm ad})} G^{j{\rm der}}$, cf. [Va2, 6.2.3]; in fact here the word if can be replaced by iff, as it can be easily deduced from [Va2, 3.3.1] and [De2, 2.1.7]. 
\medskip
{\bf 2.3.5.3. PEL-envelopes.} The PEL-envelope $f_1:(G_1,X_1)\hookrightarrow (GSp(W,\psi),S)$ of $f$ was introduced in [Va2, 4.3.12]: $G_1$ is the connected component of the origin of the subgroup of $GSp(W,\psi)$ centralizing the Lie algebra $\Mb_{\QQ}$ of endomorphisms of $W$ fixed by $G$, while $X_1$ is uniquely determined by requiring $f$ to factor through $f_1$. Let 
$$\Mb:=\Mb_{\QQ}\cap {\rm End}(L_{(p)}).$$ 
Let $v_1$ be the prime of $E(G_1,X_1)$ divided by $v$. If the quadruple $(f_1,L_{(p)},v_1,\Mb)$ is a standard PEL situation, then we refer to it as the PEL-envelope of $(f,L_{(p)},v)$. 
\medskip
{\bf 2.3.5.4. Automorphisms.} The group of automorphisms of the Shimura triple $(G,X,H)$ (or quadruple $(G,X,H,v)$) is a subgroup ${\rm Aut}((G,X,H))$ of $Aut(G_{\ZZ_{(p)}})(\ZZ_{(p)})$, cf. [Va2, 3.2.7 9)] and AE.0. We consider the natural homomorphism
$$a_{(G,X,H)}:{\rm Aut}((G,X,H))\to Aut(G^{\rm ad}_{\ZZ_{(p)}})(\ZZ_{(p)}).$$ 
If $G^{\rm ab}$ itself is the smallest subtorus of $G^{\rm ab}$ through whose extension to $\RR$ the homomorphism ${\rm Res}_{\CC/\RR} \GG_m\to G^{\rm ab}_{\RR}$ defining $X^{\rm ab}$ factors, then $a_{(G,X,H)}$ is injective and its image (identifiable with ${\rm Aut}((G,X,H))$) is a subgroup of finite index of $Aut(G^{\rm ad}_{\ZZ_{(p)}})(\ZZ_{(p)})$. So, (cf. also [De1, 1.6]) we have 
$${\rm Aut}((GSp(W,\psi),S,K))=Aut(GSp^{\rm ad}(L_{(p)},\psi))(\ZZ_{(p)}).$$ 
The subgroup of this last group normalizing $f$ (i.e., under this identification, normalizing $G$ and --over $\RR$-- taking $X$ onto itself) is called the automorphism group of $(f,L_{(p)},v)$ and is denoted by ${\rm Aut}(f,L_{(p)},v)$. 
\smallskip
Let now $g\in {\rm Aut}(f,L_{(p)},v)$. We denote by $a_g$ the $O_{(v)}$-automorphism of $\Mm$ it defines naturally (cf. 2.2.1.5.1 and natural extension of scalars). It leaves invariant the closed subscheme ${\rm Sh}_H(G,X)$ of $\Mm_{E(G,X)}$ and so it gives birth to an $O_{(v)}$-automorphism $b_g$ of $\Mn$; we have:
$$i_{\Mn}\circ b_g=a_g\circ i_{\Mn}.$$ 
\medskip
{\bf 2.3.5.5. Standard covers.} We consider a reductive subgroup $G_1$ of $GL(W)$ containing $G$ and such that $G^{\rm der}_1=G^{\rm der}$. Let $X_1$ be the $G_{1\RR}(\RR)$-conjugacy class of any element of $X$. We get a Shimura pair $(G_1,X_1)$ (see the axioms of [Va2, 2.3]). We have $E(G,X)=E(G_1,X_1)$. If $G_1\neq G$ (resp. if $G_1\neq G$ and $G_1$ is a subgroup of $GSp(W,\psi)$), then $(G_1,X_1)$ is an enlargement of $(G,X)$ (resp. of $(G,X)$ in $(GSp(W,\psi),S)$) in the sense of [Va3, 4.3.1 and 4.3.1.1].  
\smallskip
In practice $G_1$ is generated by $G$ and by a maximal torus of the centralizer of $G$ in $GL(W)$. In such a case $Z(G_1)$ is a torus: it is the group of invertible elements of the semisimple, commutative $\QQ$--subalgebra ${\rm Lie}(Z(G_1))$ of ${\rm End}(W)$. So we get:
\medskip
{\bf Fact 1.} {\it The adjoint map $a_1:(G_1,X_1)\to (G_1^{\rm ad},X_1^{\rm ad})$ is a cover (in the sense of [Va2, 2.4]). So $X_1=X^{\rm ad}$ (cf. [Mi1, 4.11]).}
\medskip
We refer to $a_1$ as a standard cover defined by $f$. We can assume $G_1$ is unramified over $\QQ_p$, cf. [Va2, 3.1.4]. If moreover the Zariski closure $G_{1\ZZ_{(p)}}$ of $G_1$ in $GL(L_{(p)})$ is a reductive group over $\ZZ_{(p)}$, then we refer to $a_1$ as a standard cover defined by $(f,L_{(p)})$; the Shimura quadruple $(G_1,X_1,H_1,v)$, with $H_1:=G_1(\QQ_p)\cap GL(L_{(p)})(\ZZ_{(p)})$, is referred as a standard extension of $(G,X,H,v)$ obtained via $(f,L_{(p)})$. 
\medskip
{\bf Fact 2.} {\it If the centralizer $C_{\ZZ_{(p)}}$ of $G_{\ZZ_{(p)}}^{\rm der}$ in $GL(L_{(p)})$ is a reductive group over $\ZZ_{(p)}$, then standard covers defined by $(f,L_{(p)})$ do exist.} 
\medskip
{\bf Proof:} We consider (see [Har, 5.5.3]) a maximal torus of $C_{\ZZ_{(p)}}$ and we take $G_1$ to be generated by its generic fibre and by $G^{\rm der}$. The fact that $G_{1\ZZ_{(p)}}$ is a reductive group, is implied by [Va2, 3.1.6].  This ends the proof.
\medskip
We assume now that $a_1$ is a standard cover defined by $(f,L_{(p)})$. Let $\Mn_1$ be the integral canonical model of $(G_1,X_1,H_1,v)$, cf. [Va2, 6.2.3]. Based on [Va2, 3.2.7 4) and 3.2.15], we have a natural open closed $O_{(v)}$-embedding $i_1:\Mn\hookrightarrow\Mn_1$. Moreover the natural $O_{(v)}$-morphism
$$q_{\Mn_1}:\Mn_1\to\Mn^{\rm ad}_{O_{(v)}}$$ 
is a pro-\'etale cover (cf. [Mi1, 4.13] and [Va2, 6.4.5.1 b)]; for $p=3$ cf. also \S 6). From the very definition of $X_1$, we have:
\medskip
{\bf Fact 3.} {\it The monomorphisms ${\rm Res}_{\CC/\RR} \GG_m\hookrightarrow G_{1\RR}$ defined by elements of $X_1$, factor through $G_{\RR}$ and so through $GSp(W\otimes_{\QQ} \RR,\psi)$.}      
\medskip
{\bf 2.3.5.6. Remark.} Fact 3 allows the interpretation of $\Mn_{1E(G,X)}$ as a moduli scheme of principally $\ZZ_{(p)}$-polarized abelian varieties over $E(G,X)$-schemes of dimension $e$, endowed with a family of Hodge cycles (still indexable by $\Mj^\prime$) and having compatibly level-$N$ structures (not necessarily symplectic similitude ones), $\forall N\in\NN$ with $(N,p)=1$, and satisfying some axioms. As this is very much the same as [Va2, 4.1 and 4.1.0] (we just need to allow principal $\ZZ_{(p)}$-polarizations as well as level-$N$ structures which are not necessarily symplectic similitude ones), we will not repeat the details here. 
\medskip
{\bf 2.3.5.6.1. Variants. A.} We refer to 2.3.5.5. We assume $G_1$ is unramified over $\QQ_p$ but the Zariski closure of $G_1$ in $GL(L_{(p)})$ is not a reductive group. From [Va2, 3.1.1 and 3.2.7.1] we deduce the existence of a hyperspecial subgroup $H_1$ of $G_1(\QQ_p)$ containing $H$ and which is the group of $\ZZ_p$-valued points of a reductive group $G_{1\ZZ_{(p)}}$ having $G_1$ as its generic fibre and having $G_{\ZZ_{(p)}}$ as a closed subgroup. So, we can still speak about $i_1$ and $q_{\Mn_1}$. Based on [Va2, 3.3.1-2], 2.3.5.6 still makes sense in such a situation. 
\smallskip
We consider a $\ZZ_{(p)}$-lattice ${L_1}_{(p)}$ of $W$ such that the Zariski closure of $G_1$ in $GL({L_1}_{(p)})$ is $G_{1\ZZ_{(p)}}$, cf. [Ja, 10.4 of Part I] and [Va2, 3.1.2.1 c)]. We refer to $a_1$ as a standard cover defined by $(f,{L_1}_{(p)})$. Even better, we consider the standard monomorphism
$${\rm sm}:SL({L_1}_{(p)})\hookrightarrow GSp({L_1}_{(p)}^2,\psi_1)$$ 
(as in [Va2, 6.6.5 d1)]; to be compared with [Va2, 6.5.1.1 v)]), with $\psi_1$ a suitable perfect alternating form on ${L_1}_{(p)}^2:={L_1}_{(p)}\oplus {L_1}_{(p)}$. Let $W_1:={L_1}_{(p)}^2\otimes_{\ZZ_{(p)}} \QQ$. We get naturally an injective map
$$\tilde f:(G,X)\hookrightarrow (GSp(W_1,\psi_1),S_1)$$ 
(the restriction of $\tilde f$ to $G^0$ is the composite of the monomorphism $G^0\hookrightarrow SL(W)$ with the generic fibre of ${\rm sm}$) such that $G_{1\ZZ_{(p)}}$ is naturally a subgroup of $GL({L_1}_{(p)}^2)$. 
\medskip
{\bf B.} We have a variant of 2.3.5.5-6 and A, where $G_1$ is generated by $G$ and by the center of the centralizer of $G$ in $GL(W)$. What we gain in this case: the torus $Z(G_1)$ is the group of invertible elements of a semisimple, commutative $\QQ$--subalgebra of ${\rm End}(W)$ left invariant by the involution of ${\rm End}(W)$ defined by $\psi$. So we have the following logical two extensions of the interpretation of 2.3.5.6.
\medskip
\item{(INTR1)} If moreover $G^{\rm der}$ is the derived subgroup of the connected component of the origin $G^\prime$ of the centralizer $C(G_1)$ of $Z(G_1)$ in $GSp(W,\psi)$, then a similar interpretation as in 2.3.5.6 holds for $\Mn_1$ itself: we just need to replace Hodge cycles by $\ZZ_{(p)}$-endomorphisms and to use Serre--Tate's deformation theory (to be compared with [Ko2, ch. 5]). Warning: if $C(G_1)\neq G^\prime$, then we still do not need to mention Hodge cycle which are not $\ZZ_{(p)}$-endomorphisms, as $C(G_1)$ is a subgroup of $G_1$. 
\medskip
\item{(INTR2)} If $G^{\rm der}$ is not the derived subgroup of $G^\prime$, then denoting by $X^\prime$ the $G^\prime(\RR)$-conjugacy class of homomorphisms ${\rm Res}_{\CC/\RR}\GG_m\to G^\prime_{\RR}$ generated by $X$, we get a PEL type embedding $f^\prime:(G^\prime,X^\prime)\hookrightarrow (GSp(W,\psi),S)$. Let $v^\prime$ be the prime of $E(G^\prime,X^\prime)$ divided by $v$. As $Z(G_1)$ is unramified over $\QQ_p$ (it splits over $\QQ_p^{\rm un}$ as $G_1$ does), $G^\prime$ is as well unramified over $\QQ_p$. If $G_1^\prime$ and respectively $\Mn_1^\prime$ have (by working with $v^\prime$) the same meaning for $G^\prime$ as $G_1$ and $\Mn_1$ have for $G$, we get naturally (see $i_{\Mn}$ of 2.3.3) a finite $O_{(v^\prime)}$-morphism $\Mn_1\to\Mn_1^\prime$ and so to any $\FF$-valued point of $\Mn_1$ it is attached naturally (cf. (INTR1)) a principally $\ZZ_{(p)}$-polarized abelian variety over $\FF$ of dimension $e$, endowed with a family of $\ZZ_{(p)}$-endomorphisms (indexed by the elements of ${\rm Lie}(Z(G_1))\cap {\rm End}(L_{(p)})$) and having compatibly level-$N$ structures (not necessarily symplectic similitude ones), $\forall N\in\NN$ with $(N,p)=1$.   
\medskip
{\bf 2.3.5.7. Remark.} We refer to the Existence Property of 1.10. From [Va2, 6.5.1.1 i) and the first part of 6.6.5.1] we deduce that we can moreover assume that the representation of $G^{0\rm der}_{W(\FF)}$ on $GL(L_{(p)}\otimes_{\ZZ_{(p)}} W(\FF))$ is a direct sum of irreducible representations such that:
\medskip
-- each simple (irreducible) direct factor of it is a representation associated to a minimal weight of a product factor of $G^{0\rm der}_{W(\FF)}$ having a simple adjoint, and so its special fibre is an irreducible subrepresentation of $G^{0\rm  der}_{\FF}$ (cf. also [Bo1, 6.4 and 7.3]).
\medskip
So in such a case $C_{\ZZ_{(p)}}$ is a reductive group, and so Fact 2 applies. Moreover, referring to 2.3.5.6.1 B, we have: as the Zariski closure of $Z(G_1)$ in $GL(L_{(p)})$ is a torus of $GL(L_{(p)})$ whose Lie algebra $\Mb$ is self dual w.r.t. $\psi$, the triple $(f^\prime,L_{(p)},v^\prime,\Mb)$ is a standard PEL situation.
\medskip
{\bf 2.3.5.8. An application: the proof of AE.4.2 a).} The case $p=3$ is very much the same (see \S6). So, to be simpler in references, we assume $p\ge 5$. Let $(G_1,X_1,H_1,v_1)$ be a Shimura quadruple of preabelian type, with $v_1$ dividing $p$. Let $\Mn_1$ (resp. $\Mn_1^{\rm ad}$) be its integral canonical model (resp. be the integral canonical model of its adjoint), cf. [Va2, 6.4.2] and AE.4.1. Let $H_{01}$ be a subgroup of $G_1(\AA_f^p)$ such that $H_{01}\times H_1$ is a subgroup of $G_1(\AA_f)$ smooth for $(G_1,X_1)$. If $p$ does not divide $t(G_1^{\rm ad})$, the fact that $\Mn_1$ is a pro-\'etale cover of $\Mn_1/H_{01}$ is implied by the Proposition of AE.4. We now assume $p| t(G_1^{\rm ad})$ and $H_{01}\times H_1$ is $p$-smooth for $(G_1,X_1)$. We need to show that $\Mn_1$ is a pro-\'etale cover of $\Mn_1/H_{01}$. As this is a problem of connected components of $\Mn_{1W(\FF)}$, we can assume (cf. Definition 1 of AE.4), that there is a prime $l$ different from $p$ and such that the image of $H_{01}$ in $G_1^{\rm ad}(\QQ_l)$ is contained in a compact, open subgroup $H_{01l}^{\rm ad}$ having no pro-$p$ subgroups. Let $H_{01}^{l{\rm ad}}$ be a compact, open subgroup of $G_1^{\rm ad}(\AA_f^{p,l})$ containing the image of $H_{01}$ in it; here $\AA_f^{p,l}$ denotes the ring of finite ad\`eles whose both $p$- and $l$-components are omitted.  
\smallskip
As the natural $O_{(v_1)}$-morphism $\Mn_1\to\Mn_{1O_{(v_1)}}^{\rm ad}$ is pro-\'etale (see [Va2, 6.4.5 b)]), and as the quotient $O_{(v_1)}$-morphism $\Mn_1\to\Mn_1/H_{01}$ factors through the natural morphism $\Mn_1\to\Mn_{1O_{(v_1)}}^{\rm ad}/\tilde H_{01}^{\rm ad}$, with 
$$\tilde H_{01}^{\rm ad}:=H_{01}^{l{\rm ad}}\times H_{01l}^{\rm ad},$$ 
we can assume $G_1$ is adjoint. Based on [Va2, 3.4.5.1], it is enough to show that the morphism $\Mn_1^{\rm ad}\to\Mn_1^{\rm ad}/H_{01}^{l\rm ad}$ is a pro-\'etale cover; here the role of $H_{01}^{l{\rm ad}}$ is that of an arbitrary compact subgroup of $G_1^{\rm ad}(\AA_f^{p,l})$. Based on this and on [Va2, 3.2.16], we can assume $G_1$ is $\QQ$--simple. So based on 2.3.5.7, not to introduce extra notations, (by moving back to a $G_1$ which is not adjoint but has a $\QQ$--simple adjoint) we can assume that:
\medskip
-- we have a SHS $(f,L_{(p)},v)$ and we are in the context of a standard cover $a_1:(G_1,X_1)\to (G_1^{\rm ad},X_1^{\rm ad})=(G^{\rm ad},X^{\rm ad})$ defined by $(f,L_{(p)})$;
\smallskip
-- with the notations of 2.3.5.5, we are dealing with a subgroup $\tilde H_{01}^{\rm ad}$ of $G_1^{\rm ad}(\AA_f^p)$ which has the above shape and properties (so in particular, $\tilde H_{01}^{\rm ad}\times G_{1\ZZ_{(p)}}^{\rm ad}(\ZZ_p)$ is $p$-smooth for $(G_1^{\rm ad},X_1^{\rm ad})$).
\medskip
A connected component of $\Mn_{1W(\FF)}^{\rm ad}/H_{01}^{l\rm ad}$ is a quotient of a connected component $\Mc^0$ of $\Mn_{1W(\FF)}$ through a group of automorphisms $GA$ of $\Mn_{1W(\FF)}$ leaving invariant $\Mc^0$ and defined by translations by a subgroup of $G_1(\AA_f^p)$ whose image in $G_1^{\rm ad}(\QQ_l)$ is trivial, cf. [Mi1, 4.13] and [Va2, 3.3.1]. If $h\in GA$ fixes $y\in\Mc^0(\FF)$, then as in AE.4.1 we get that $h$ acts trivially on $\Mc^0(\FF)$ (we need to work precisely with $l$). Warning: here we dot need to bother about Hodge cycles which are not defined by endomorphisms, i.e. in connection to the $\ZZ_{(p)}$-automorphism $a$ we get (as in AE.4.1) we are bothered just about $\ZZ_{(p)}$-polarizations, level structures (and if one desires, about $\ZZ_{(p)}$-endomorphisms), cf. 2.3.6.1 (INTR2). So $\Mc^0$ is a pro-\'etale cover of $\Mc^0/GA$. This ends the proof of AE.4.2 a).
\medskip
{\bf 2.3.5.8.1. Variant.} What follows is a natural extension of AE.4.1.1 and so provides a variant of the last paragraph of 2.3.5.8; so we can assume $G^{\rm ad}$ is $\QQ$--simple. If it is of some $A_n$ Lie type with $p|n+1$, then we can proceed as in AE.4.1.1 to get that we can assume that moreover $h\in G(\AA_f^p)$. If $G^{\rm ad}$ is not of $A_n$ Lie type with $p|n+1$, then $q_{G^{\rm ad}}$ (see 2.3.5.2) is relatively prime to $p$ and so for any $\tilde h\in G^{\rm ad}(\AA_f^p)$ there is $\tilde q\in\NN$, with $(\tilde q,p)=1$, such that $\tilde h^{\tilde q}$ belongs to the image of $G^{\rm der}(\AA_f^p)$ in $G^{\rm ad}(\AA_f^p)$. So, as in AE.4.1.1, we can assume that $h\in G(\AA_f^p)$; we conclude: 
\medskip
{\it Regardless of how $G^{\rm ad}$ is, in the last paragraph of 2.3.5.8, the use of ``$\ZZ_{(p)}$-" in front of polarizations (and isogenies) can be entirely avoided.}
\medskip
{\bf 2.3.6. Examples.} We refer to 2.3.1. If the family of tensors $(v_\al)_{\al\in\Mj}$ is $\ZZ_{(p)}$-well positioned w.r.t. $\psi$ for the group $G$, if it is enveloped by $L_{(p)}$ (see [Va2, 4.3.4]
for the meaning of these), and if it is of degree at most $2(p-2)$, then the triple $(f,L_{(p)},v)$ is a SHS
(cf. [Va2, 5.1 and 5.6.5]). More generally, if there is a family of tensors of $\Mt(L_{(p)}\otimes_{\ZZ_{(p)}} W(\FF))$ which is $W(\FF)$-well positioned for $G_{W(\FF)}$ and of degree at most $2(p-2)$, then the triple $(f,L_{(p)},v)$ is a SHS, cf. [Va2, 5.6.9].
\smallskip
In particular, if $f$ is a good embedding w.r.t. $p$ (see def. [Va2, 5.8.1]), then there is a $\ZZ_{(p)}$-lattice of $W$ which is crystalline well positioned for the map $f$. From this the Existence Property of 1.10 follows. 
\medskip
{\bf 2.3.6.1. The $d$-invariants.} Let $\tilde W$ be a finite dimensional vector space over a field $\tilde k$. We assume $\tilde k$ is either of characteristic $0$ or it is perfect of arbitrary positive characteristic. We consider a faithful representation $\tilde f:\tilde G\hookrightarrow GL(\tilde W)$, with $\tilde G$ a reductive group over $\tilde k$. We assume it is not an isomorphism. If $\tilde G$ contains the $1$ dimensional torus $\tilde T_1$ of scalar automorphisms of $\tilde W$, then we also assume the existence of a non-degenerate bilinear form $\tilde\psi$ on $\tilde W$ or on its dual normalized by $\tilde G$. The existence of the numbers introduced by the following definition is a consequence of [De4, 3.1 c)] and of the finite dimensionality of $\tilde W$.
\medskip
{\bf Definition.} Let $d(\tilde f)\in\NN\cup\{0\}$ (resp. $d_0(\tilde f)\in\NN\cup\{0\}$) be the smallest number such that there is a family of tensors of $\Mt(\tilde W)$ of degree at most $d(\tilde f)$ (resp. at most $d_0(\tilde f)$) such that the subgroup of $GL(\tilde W)$ fixing all its members and, in case $\tilde G$ contains $\tilde T_1$, normalizing $\tilde\psi$, is $\tilde G$ (resp. has $\tilde G$ as the connected component of its origin). We refer to $d(\tilde f)$ (resp. to $d_0(\tilde f)$) as the $d$-invariant (resp. as the $d_0$-invariant) of $\tilde f$. 
\medskip
We have $d(\tilde f)\ge d_0(\tilde f)$. By enlarging $\tilde k$, these numbers remain the same. So let $d_1(\tilde f)$ (resp. $d_2(\tilde f)$) be the $d$-invariant (resp. the $d_0$-invariant) of the natural representation of an arbitrary maximal torus of $\tilde G$ on $\tilde W$. Let $d_3(\tilde f)$ (resp. $d_4(\tilde f)$) be the $d$-invariant (resp. the $d_0$-invariant) of the natural representation of $\tilde G^{\rm der}$ on $\tilde W$. Let $d_5(\tilde f)$ (resp. $d_6(\tilde f)$) be the $d$-invariant (resp. the $d_0$-invariant) of the natural representation of the maximal subtorus of $Z(\tilde G)$ on $\tilde W$. 
\medskip
{\bf Example.} We refer to 2.2.9 5). We assume $G$ does not contain the $1$ dimensional split torus of scalar automorphisms of $M$. Then the $d$-invariant of the representation of $G_{B(k)}$ on $M[{1\over p}]$ is ${\rm deg}(M,\vph,G)$. The proof of this is very much the same as the proof of the Fact of 2.2.9 1').
\medskip
If $(f,L_{(p)},v)$ is a SHS, we refer to the numbers $d(f)$, $d_0(f)$, $d_1(f)$,..., $d_6(f)$, $d(f_p)$, $d_0(f_p)$, $d_3(f_p)$ and $d_4(f_p)$ as its $d$-invariants. Here $f$ is identified with a representation of $G$ on $W$ and 
$$f_p:G_{\FF_p}\hookrightarrow GL(L_{(p)}/pL_{(p)})$$ 
is the (alternating) faithful representation defined naturally by the closed embedding $G_{\ZZ_{(p)}}\hookrightarrow GL(L_{(p)})$. Warning: simple arguments involving representations of tori show that we have $d_i(f_p)=d_i(f)$ for $i\in\{1,2,5,6\}$; but we have no a priori reason to think that in general this automatically holds for $i\in\{0,3,4\}$ or that $d(f)=d_p(f)$. From the second paragraph of 2.3.1, we get: $d(f)$, $d_0(f)$, $d_1(f)$, $d_2(f)$, $d_5(f)$ and $d_6(f)$ are even numbers. Any one of the conditions $d(f)=0$ or $d_0(f)=0$ are equivalent to $f$ being an isomorphism. 
\smallskip
Using the language of [Va2, 4.3.1], we have: if $(G,X)$ is saturated in $(GSp(W,\psi),S)$, then $d_0(f)\le 4$. In general we do not know when $d_0(f)=4$. If we have a standard PEL situation $(f,L_{(p)},v,\Mb)$, then $d_0(f)\in\{0,2\}$ but (see the D case in [Ko2, \S 5]) it can happen that $d(f)\ge 4$. 
\medskip
{\bf 2.3.7. Convention.} Whenever we consider a SHS $(f,L_{(p)},v)$, we use freely the notations of 2.3.1-3 and of 2.3.5.2. 
\medskip
{\bf 2.3.8. Remarks.} 
{\bf 1)} We do not know when the converse of 2.3.6 is true. That is, we do not know for which SHS $(f,L_{(p)},v)$ there is a family $(v_\al)_{\al\in\Mj}$ of homogeneous tensors of the essential tensor algebra of $W^*\oplus W$, which is $\ZZ_{(p)}$-well positioned w.r.t. $\psi$ for $G$, is enveloped by $L_{(p)}$, and is of degree at most $2(p-2)$. It is not hard to see, that for $p=3$ we do need that at least one of the following two conditions is satisfied: 
\medskip
-- $d_0(f)$ is $0$ or $2$;
\smallskip
-- there is a saturated enlargement $(G_1,X_1)$ (see def. [Va2, 4.3.1]) of $(G,X)$ in $(GSp(W,\psi),S)$ such that the $d_0$-invariant of the representation of $G_1$ on $W$ is $2$.
\medskip
{\bf 2)} A major result of \S 6 states that any $\ZZ_{(p)}$-lattice $L_{(p)}$ of $W$ good w.r.t. $f$ is automatically crystalline well positioned for the map $f$. In particular, we will see in \S6 that in 2.3.5.6.1 A we get a SHS $(\tilde f,{L_1}_{{(p)}}^2,v)$.
\smallskip
{\bf 3)} Usually when we mention the polarizations we use the set of indices $\Mj$ (i.e. we would like to ``carry" as few Hodge cycles as possible); when we do not, then we use $\Mj^\prime$. If this is confusing to the reader, the reader can assume $\Mj=\Mj^\prime$.
\smallskip
{\bf 4)} We consider a quadruple $(f,L_{(p)},v,\Mb)$ which satisfies all requirements of a standard PEL situation, except that $(f,L_{(p)},v)$ is a SHS. Then such a quadruple is a standard PEL situation. For the case when $\Mb_{\QQ}:=\Mb\otimes_{\ZZ_{(p)}} \QQ$ is a simple $\QQ$--algebra, cf. [LR] and [Ko2, ch. 5] (cf. also [Va2, 5.6.3] and AE.1 for very detailed proofs in the spirit of arbitrary SHS's). Warning: in the way we defined a standard PEL situation, $\Mb_{\QQ}:=\Mb\otimes_{\ZZ_{(p)}} \QQ$ is not necessarily a $\QQ$--simple algebra and so $G^{\rm ad}$ does not necessarily have all its simple factors of the same Lie type. However, [Va2, 4.3.11] and its slight correction of AE.1, still apply (versus [Va2, 5.6.3]) to such a general $\ZZ_{(p)}$-algebra $\Mb$: just the third paragraph of [Va2, p. 469] has to be modified (we need to define separately a number $n\in\NN$ as in loc. cit. for each simple factor of $\Mb_{\QQ}$).  
\medskip
{\bf 2.3.8.1. Exercise.} Let $(G,X,H,v)$ be a Shimura quadruple of Hodge type, with $v$ dividing $p\ge 3$. Assuming 2.3.8 2), show that there is a SHS $(f,L_{(p)},v)$ for which we do have $H=G_{\ZZ_{(p)}}(L_{(p)}\otimes_{\ZZ_{(p)}} \ZZ_p)$. Hint: just restate the proof of [Va2, 6.7.2].
\medskip
{\bf 2.3.9. Theorem.} {\it Let $(f,L_{(p)},v)$ be a SHS. We consider the context and notations of 2.3.4. We have:
\smallskip
{\bf a)} $\tilde G$ is an inner form of $G_{W(k)}$; moreover $\tilde G^{\rm ab}$ is a split torus.
\smallskip
{\bf b)} The reductive group $\tilde G$ splits over a finite Galois cover ${\rm Spec}(W(k_1))$ of ${\rm Spec}(W(k))$ of order dividing a fixed power of $q_{G^{\rm ad}}$. Moreover, the Galois group ${\rm Gal}(k_1/k)$ is nilpotent of index of nilpotency at most two. If for each simple factor $(G_j,X_j)$ of $(G^{\rm ad},X^{\rm ad})$ of $D_n^{\RR}$ type with $n\in 2\NN$, $n\ge 6$, the semisimple subgroups of $G^{\rm der}_{\CC}$ naturally isogeneous to a simply factor of $G_{j\CC}$, are simply connected, then this Galois group is abelian.
\smallskip
{\bf c)} If all factors of $G^{\rm ad}$ are of $C_n$ or $A_n$ Lie type for some $n\in\NN$, and if for any simple factor $G_{\QQ}^\prime$ of $G^{\rm ad}_{\QQ}$ of $A_n$ Lie type, the semisimple subgroup of $G^{\rm der}_{\QQ}$ naturally isogeneous to it is a simply connected semisimple group, then $\tilde G$ is isomorphic to $G_{W(k)}$, regardless of the choice of the $W(k)$-valued point $z$ of $\Mn$.
\smallskip
{\bf d)} If $G^{\rm der}$ is a simply connected semisimple group, then the Galois group of b) is abelian of exponent 2.}
\medskip
{\bf Proof:} For the sake of future references, we index the main paragraphs by capital letters.
\smallskip
{\bf A.} Let $N\in\NN$ with $(N,p)=1$. As $(\Ma,\Mp_{\Ma})$ has level-N symplectic similitude structure, so does $(A_k,p_{A_k})$. So the $N$-th roots of unity are contained in $k$. We conclude: $k\supset\FF=\overline{k(v)}$. So also $B(k)\supset B(\FF)$. The Proposition is a consequence of this last inclusion and of Fontaine's comparison theory. We present the instructive details. We recall that:
\medskip
{\bf Theorem.} {\it Let $q$ be a prime. Any Hodge cycle of an abelian variety over the field of fractions $FR$ of a complete DVR of mixed characteristic $(0,q)$ having a perfect residue field, is a de Rham cycle, i.e. under the isomorphism of the de Rham conjecture, its de Rham component is mapped into the $q$-component of its \'etale component.}
\medskip
The case when the abelian variety is definable over a number field contained in $FR$, this result was known since long time (for instance, see [Bl, (0.3)]). This extra hypothesis was removed in [Va2, 5.2.16] (in the part of [Va2, 5.2] preceding [Va2, 5.2.16] an odd prime is used; however the proof of loc. cit. applies to all primes).  
\smallskip 
{\bf B.} As $G$ is unramified over $\QQ_p$, $G_{B(\FF)}$ is a split group. So $G_{W(k)}$ is a split group. Fontaine's comparison theory implies $\tilde G_{B(k)}$ is an inner form of $G_{B(k)}$ (cf. Theorem of A). But $\tilde G$ splits over a finite \'etale extension of $W(k)$. So over this extension it becomes isomorphic to $G_{W(k)}$, cf. the uniqueness of a split reductive group associated to a given root datum (see [SGA3, Vol. III, p. 313-4]). As $G^{\rm ad}_{W(k)}(W(k))$ is that maximal bounded subgroup of $G^{\rm ad}_{W(k)}(B(k))$ (see [Ti2, 3.2]) which normalizes ${\rm Lie}(G^{\rm ad}_{W(k)})$, we get: any automorphism of $G_{W(k)}$ which generically is inner, is inner. We deduce: $\tilde G$ is an inner form of $G_{W(k)}$. So $\tilde G^{\rm ab}$ is a split torus. This takes care of a).
\smallskip
{\bf C.} Now we prove simultaneously the other parts. Any element of $\Mt(L_{(p)}^*[{1\over p}])$ fixed by $G_{\QQ_p}$ gives birth (cf. the identification of [Va2, top of p. 473]) through Fontaine's comparison theory to an element of $\Mt(M[{1\over p}])$ fixed by $\tilde G$ (again cf. Theorem of A; as we are over $W(k)$ we can use as well Fontaine's comparison theory in a purely crystalline context). As $B(k)\supset B(\FF)$, we can use $B(\FF)$-linear combinations of such elements. We denote by 
$${\rm Fix}(G_{B(\FF)})$$ 
the set of tensors of $\Mt(L_{(p)}^*\otimes_{\ZZ_{(p)}} B(\FF))$ fixed by $G_{B(\FF)}$. So any element of ${\rm Fix}(G_{B(\FF)})$ has a counterpart in $\Mt(M[{1\over p}])$. As $G_{B(\FF)}$ is a split group and all simple Lie factors of ${\rm Lie}(G^{\rm ad}_{B(\FF)})$ are classical, split, simple Lie algebras over $B(\FF)$, we can often ``capture" the fact that $G_{B(\FF)}$ is split by just using the set ${\rm Fix}(G_{B(\FF)})$. The most useful elements of ${\rm Fix}(G_{B(\FF)})$ are its projectors: they single out different subrepresentations of the representation of $G_{B(\FF)}$ on $\Mt(L_{(p)}^*\otimes_{\ZZ_{(p)}} B(\FF))$ which can ``encode" useful information on $G^{\rm der}_{B(\FF)}$.
\smallskip
{\bf D.} We first assume that:
\medskip
\item{(*)} for any simple factor $AD$ of $G^{\rm ad}_{B(\FF)}$ ``coming" from a simple factor of $(G^{\rm ad},X^{\rm ad})$ of $A_n$ ($n\in\NN$, $n\ge 2$) or $D_n^{\RR}$ type ($n\in 2\NN$, $n\ge 6$), the semisimple subgroup $SSS$ of $G^{\rm der}_{B(\FF)}$ naturally isogeneous to it, is simply connected and the representation of $G_{B(\FF)}$ on $L_{(p)}^*\otimes_{\ZZ_{(p)}} B(\FF)$ does contain a subrepresentation $W_{SSS}$ (not necessarily irreducible) on which the connected component of the origin of $Z(G_{B(\FF)})$ acts as multiplication with scalars and on which $G^{\rm der}_{B(\FF)}$ acts through a quotient of it which is naturally isomorphic to $SSS$.
\medskip
The image of $G_{B(\FF)}$ in $GL(W_{SSS})$ is a reductive group $GSSS$ which has a quotient $QG$ which is either a $GL_{n+1}$-group or is a $GSO(2n)$-group. From [De4, 3.1 (a)] we deduce that the classical representation of $QG$ (defined as in 2.2.23 B) ``shows up" as a subrepresentation of $GSSS$ on $\Mt(W_{SSS})$ factoring through $QG$. From this, Weyl's complete reductibility theorem, the structure of the irreducible representations of a direct sum of absolutely simple, split, semisimple Lie algebras over $B(\FF)$, and the fact that the irreducible subrepresentations of the representation of a simple Lie factor of ${\rm Lie}(G_{B(\FF)}^{\rm ad})$ on $L_{(p)}^*\otimes_{\ZZ_{(p)}} B(\FF)$ are given by minimal weights (see [Sa] or [De2, p. 261]) we deduce: 
\medskip
{\bf Fact.} {\it For any simple factor ${\got g}_0$ of ${\rm Lie}(G^{\rm ad}_{B(\FF)})$, we can single out a subrepresentation of the representation of ${\rm Lie}(G^{\rm ad}_{B(\FF)})$ on $\Mt(L_{(p)}^*\otimes_{\ZZ_{(p)}} B(\FF))$ which is the tensor product $\Mt\Mp$ of the classical faithful representation (viewed as a monomorphism) ${\got g}_0\hookrightarrow  {\rm End}(W_0)$ (over $B(\FF)$; for its complex version we refer to [He, \S 8 of ch. 3]) with other representations of the same type of the other simple factors of ${\rm Lie}(G_{B(\FF)}^{\rm ad})$, just by using (a projector of) ${\rm Fix}(G_{B(\FF)})$; moreover, the bilinear form (if it is not zero, then it is non-degenerate and unique up to a non-zero scalar) $b_0$ on $W_0$ annihilated by ${\got g}_0$ (if ${\got g}_0$ is of $A_n$ Lie type, $n\in\NN\setminus\{1\}$, then this bilinear form is zero), gives birth to a bilinear form on $\Mt\Mp$ fixed by $G_{B(\FF)}^{\rm der}$ and normalized by $G_{B(\FF)}$. If ${\got g}_0={\rm Lie}(AD)$, then $\Mt\Mp=W_0$.}
\medskip
As $\tilde G$ is an inner form of $G_{W(k)}$, all simple factors of $\tilde G^{\rm ad}$ are absolutely simple factors. From the Fact and Fontaine's comparison theory, we get the existence of an absolutely simple factor $\tilde {\got g}_0$ of ${\rm Lie}(\tilde G_{B(k)}^{\rm ad})$, which is an inner form of ${\got g}_0$ and which has a faithful representation on $\tilde W_0:=W_0\otimes_{B(\FF)} B(k)$, annihilating a bilinear form $\tilde b_0$ on $\tilde W_0$. $\tilde b_0$ is a non-degenerate form iff the Lie type of ${\got g}_0$ is not $A_n$, with $n\in\NN$, $n\ge 2$. We recall that the dual of a tensor product of representations is the tensor product of the dual representations.
\smallskip
We now treat the case when all factors of $G^{\rm ad}$ are of some $A_n$ or $C_n$ Lie type ($n\in\NN$). By reasons of dimensions (resp. due to the existence of the alternating non-degenerate form $\tilde b_0$) in case ${\got g}_0$ is of $A_n$ (resp. of $C_n$) Lie type, we get: $\tilde {\got g}_0$ is isomorphic to ${\got g}_0\otimes_{B(\FF)} B(k)$. We deduce $\tilde {\got g}_0$ is a split Lie algebra over $B(k)$. So $\tilde G^{\rm der}$ is a split group. So (cf. a)) $\tilde G$ is a split group. As $\tilde G_{B(k)}$ and $G_{B(k)}$ are both split and inner forms of each other, they are isomorphic. We conclude (as in B): $\tilde G$ is isomorphic to $G_{W(k)}$. This ends the proof of c) under the assumption (*). 
\smallskip
{\bf E.} Coming back to the general situation (of the assumption (*)), it is enough to show: if $\tilde b_0$ is a non-degenerate symmetric form (i.e. if ${\got g}_0$ is of $B_n$ or $D_n$ Lie type), then $\tilde {\got g}_0$ and ${\got g}_0\otimes_{B(\FF)} B(k)$ become isomorphic to each other over a finite, unramified abelian extension $B(k_1)$ of $B(k)$ of order a divisor of $2^{N(W_0)}$, where $N(W_0)\in\NN$ depends only on $m:=\dim_{B(\FF)}(W_0)$. 
\smallskip
To see this, let $\tilde G_0$ be a semisimple group over $W(k)$ whose generic fibre is the connected component of the origin of $O(\tilde W_0,\tilde b_0)=Aut(\tilde W_0,\tilde b_0)$. Its existence is implied by the fact that $\tilde G$ is a reductive group. It is enough to show that $\tilde G_{0k}$ splits over a finite abelian extension $k_1$ of $k$ whose Galois group is a $2$-torsion group of whose order is bounded in terms of $m$. We consider a $W(k)$-lattice $\tilde L_0$ of $\tilde W_0$ such that $\tilde G_0$ is a closed subgroup of $GL(\tilde L_0)$, cf. [Ja, 10.4 of Part I] and [Va2, 3.1.2.1 c)]. The existence of $\tilde b_0$ implies (very easy exercise; see also the general result in [Ja, 10.9 of Part I]) the existence of a non-degenerate, symmetric bilinear form on $\tilde L_0/pL_0$ fixed by $\tilde G_{0k}$. With respect to some $k$-basis of $\tilde L_0/pL_0$, the quadratic form of this new bilinear form, can be written as a sum $\sum_{i=1}^m a_ix_i^2$, where $a_i$'s are non-zero elements of $k$. We get: $\tilde G_{0k}$ splits over the abelian extension $k_1$ of $k$ obtained by adjoining square roots of not more than $[{m\over 2}]$ elements of $k$. So $[k_1:k]|2^{[{m\over 2}]}$. This ends the proof of b) and d) under the assumption (*). 
\smallskip
{\bf F.} To prove b), c) and d) in general we use the standard techniques of [Va2, 6.2.3.1, 6.2.3, 6.4.2 and 6.4.5.1] (for the case $p=3$ cf. \S 6). In other words, cf. 2.3.6 and loc. cit., there is another SHS $(f^1,L_{(p)}^1,v^1)$ such that:
\medskip
-- the adjoint quadruples of Shimura quadruples $(G^1,X^1,H^1,v^1)$ and $(G,X,H,v)$ are isomorphic;
\smallskip
-- (*) is satisfied for it (for its last part pertaining to subrepresentations cf. also [Va2, 6.5.1.1 and 6.6.5]: they rely on [De2, proof of 2.3.10] and for loc. cit. it is obviously satisfied);
\smallskip
-- there is a natural isogeny $G^{1{\rm der}}\to G^{\rm der}$.
\medskip
For this new SHS we use the standard notations, except that we put an upper right index $1$ everywhere.
As $k$ is an $\FF$-algebra, from [Va2, 6.2.3 and 3.2.7 10)] we get: we can assume that each connected component of $\Mn^1_{W(k)}/H_0^1$ is an abelian cover of a connected component of $\Mn_{W(k)}/H_0$. Moreover, we can assume that the degree $d(H_0,H_0^1)$ of this cover divides a fixed power of $q_{G^{\rm ad}}$ (cf. [Va2, 6.2.3.1]) and that all Hodge cycles of $\Ma^1$ are as well Hodge cycles of $\Ma^1_{H_0^1}$.  
So the standard trick of using Segre's embeddings as in [Va2, Example 3 of 2.5] (it is fully recalled in 4.9 below; here we need just the rational version of 4.9.3-4 below and so from 4.9.2.0 only its C part is needed), allows us to shift our attention to the new SHS; by doing this we might need to replace $W(k)$ by $W(k_1)$, where $k_1$ is an abelian extension of $k$ of degree a divisor of $d(H_0,H_0^1)$.  But for the new situation we can apply the last paragraph of D and E. So we get b), c) and d), as the part of b) involving nilpotency is entirely trivial. This ends the proof of the Theorem.
\medskip
{\bf 2.3.9.1. Remark.} We do not know precisely for which SHS $(f,L_{(p)},v)$ the reductive group $\tilde G$ of 2.3.4 is automatically isomorphic to $G_{W(k)}$. The fixed power of $q_{G^{\rm ad}}$ referred in 2.3.9 b) is effectively computable, as the degree $d(H_0,H_0^1)$ is effectively computable (over $\CC$). 
\smallskip
If $G^{\rm der}$ is a simply connected semisimple group, then the Galois group of 2.3.9 b) is the product of at most $r(G^{\rm ad})$ copies of $\ZZ/2\ZZ$. Here $r(G^{\rm ad})$ is a sum indexed by the simple factors of $G^{\rm ad}_{\CC}$ of $B_n$ or $D_n$ Lie type (like the term of the sum corresponding to a $SO(2,m)^{\rm ad}_{\CC}$ factor is $[{{m+2}\over 2}]$, etc.).
\medskip
{\bf 2.3.9.2. Convention.} From now on, for the sake of simplifying the notations, we assume that whenever we consider a point ${\rm Spec}(k)\to\Mn$ or a point ${\rm Spec}(W(k))\to\Mn$, the reductive group $\tilde G$ we get as in 2.3.4 is (isomorphic to) the split group $G_{W(k)}$. This convention is motivated by 2.3.9 and 2.3.9.1. It does not apply to the similar type of points of a quotient of $\Mn$ (like $\Mn/H_0$).
\medskip
{\bf 2.3.10. Shimura crystals attached to points of $\Mn$.} Let $(f,L_{(p)},v)$ be a SHS. We use the notations of 2.3.4 (with $\tilde G=G_{W(k)}$). Let $\sg:=\sg_k$. Let $\vph$ be the Frobenius endomorphism of
$M=H^1_{\rm crys}(A/W(k))$. For future references, we state here explicitly as a Corollary, part of 2.3.9 A and C.  
\medskip
{\bf Corollary.} {\it We have $\vph(t_{\al})=t_{\al},\,\forall\al\in\Mj^\prime.$}
\medskip
Warning: once we know $\Mn$ is smooth, the Corollary follows as well from the previous forms of the Theorem of 2.3.9 A, as one can check using a density argument of $O_{(v)}^{\rm sh}$-valued points of $\Mn$.
\smallskip 
Moreover, $t_{\al}$ is a tensor of the $F^0$-filtration of $\Mt(M[{1\over p}])$ defined by the Hodge filtration $F^1$ of $M$ defined by $A$. So the triple $(M,\vph,G_{W(k)})$ (resp. the quadruple $(M,F^1,\vph,G_{W(k)})$) is (cf. 2.2.8 1) and 2)) a Shimura (resp. a Shimura filtered) $\sg$-crystal, called the Shimura (resp. the Shimura filtered) $\sg$-crystal attached to the point $y:{\rm Spec}(k)\to\Mn_{k(v)}$ defined by $z$ (resp. attached to $z$).
\smallskip
We deduce (see c) of 2.2.8 1)) the existence of an injective cocharacter $\tilde\mu:\GG_m\hookrightarrow G_{W(k)}$ such that we have a direct sum decomposition $M=F^1\oplus F^0$ with $\be\in\GG_m(W(k))$ acting through $\tilde\mu$ on $F^i$ as the multiplication with $\be^{-i}$, $i=\overline{0,1}$. Another proof of the existence of such a cocharacter $\tilde\mu$ can be obtained by entirely following the construction of a cocharacter in [Va2, 5.3.1], as 2.3.9.2 tells us that $G_{W(k)}$ is split.
\smallskip
Also, either from loc. cit. or from the proof of Fact 2 of 2.2.9 3), we deduce that if the cardinality of $W(k)$ as a set is not greater than the one of $\CC$, then for any $O_{(v)}$-monomorphism $W(k)\hookrightarrow\CC$, $\tilde\mu_{\CC}$ is $G(\CC)$-conjugate to the cocharacters $\mu_x^*$, $x\in X$, of 2.3.1.
\smallskip
When we want to emphasize the family of tensors $(t_{\al})_{\al\in\Mj^\prime}$ then we also write  down $(M,\vph,G_{W(k)},(t_\al)_{\al\in\Mj^\prime})$ and $(M,F^1,\vph,G_{W(k)},(t_{\al})_{\al\in\Mj^\prime})$. We also speak about the Shimura (adjoint) Lie $\sg$-crystal attached to $y$, about the Shimura  (adjoint) filtered Lie $\sg$-crystal attached to $z$, about the Shimura (adjoint) (Lie) isocrystal attached to $y$, and about the Shimura (adjoint) filtered (Lie) isocrystal attached to $z$ (cf. 2.2.13 and 2.2.13.2).
\smallskip
When we want to emphasize the principal polarization $p_A$, we speak about the principally quasi-polarized Shimura (resp. Shimura filtered) $\sg$-crystal attached to $y$ (resp. to $z$) and denote it by $(M,\vph,G_{W(k)},p_A)$ or by $(M,\vph,G_{W(k)},(t_{\al})_{\al\in\Mj^\prime},p_A)$ (resp. by $(M,F^1,\vph,G_{W(k)},p_A)$ or by $(M,\vph,G_{W(k)},(t_{\al})_{\al\in\Mj^\prime},p_A)$).
\smallskip
We consider now a morphism $z:{\rm Spec}(W(k))\to\Mn/H_0$, with $k$ a perfect field. Let $(A,p_A):=z^*(\Ma_{H_0},\Mp_{\Ma_{H_0}})$. Over a finite, \'etale extension $W(k_1)$ of $W(k)$, $A_{W(k_1)}$ gets naturally a family $(w_{\al})_{\al\in\Mj^\prime}$ of Hodge cycles (see 2.3.3). As the set of all Hodge cycles of $A_{W(k_1)}$ is invariant under the Galois group ${\rm Gal}(W(k_1)/W(k))={\rm Gal}(k_1/k)$, we get that the principally quasi-polarized filtered $\sg$-crystal $(M,F^1,\vph,p_A)$ attached to $(A,p_A)$ has a natural structure of a principally quasi-polarized quasi Shimura filtered $\sg$-crystal $(M,F^1,\vph,\tilde G_{W(k)},p_A)$ which is not necessarily quasi-split. Argument: we take $\tilde G_{W(k)}$ to be the Zariski closure in $GL(M)$ of the subgroup of $GL(M\otimes_{W(k)} B(k))$ whose extension to $B(k_1)$ is the subgroup of $GL(M\otimes_{W(k)} B(k_1))$ fixing the de Rham component of $w_{\al}$, $\forall\al\in \Mj^\prime$; b) and c) of 2.2.8 1) hold, as they hold after this extension, and moreover a) holds (cf. Corollary) after this extension. Warning: we have no reason to consider $\tilde G_{W(k)}$ to be an inner form of $G_{W(k)}$. However, if it is an inner form and if $k$ is finite, then we have $\tilde G_{W(k)}=G_{W(k)}$ (cf. Lang's theorem applied over $k$). 
\smallskip
This allows us to speak about the not necessarily quasi-split Shimura $\sg$-crystal or about the Shimura (adjoint) Lie $\sg$-crystal or isocrystal attached to a point of $\Mn/H_0$ with values in a perfect field, and about the not necessarily quasi-split Shimura filtered $\sg$-crystal or Shimura (adjoint) filtered Lie $\sg$-crystal or isocrystal attached to a point of $\Mn/H_0$ with values in the Witt ring of such a field. Warning: all of them depend on the choice of $\Ma_{H_0}$. 
\medskip
{\bf 2.3.10.1. Comment.} Shimura Lie $\sg$-crystals attached to points $y:{\rm Spec}(k)\to\Mn_{k(v)}$ play a central role in \S1-14. We can not refrain from commenting why they have not been used before. We think that mainly due to two reasons. The first reason is: in the case of special fibres of integral canonical models of Siegel modular varieties, all the information provided by the Shimura Lie $\sg$-crystal attached to such a point, can be read out directly from the principally quasi-polarized Shimura $\sg$-crystal ``producing it". The second reason is: the way of thinking of points of $\Mn_{k(v)}$ was mainly rational, i.e. was involving an analysis of the Shimura isocrystals defined by Shimura $\sg$-crystals attached to these points, even though for standard PEL situations (see [Ko2, ch. 5] and 2.3.5) was obvious (cf. [Va2, 4.3.11]; see also [LR] and [Ko2]) that ``we get" Shimura $\sg$-crystals. 
\medskip
{\bf 2.3.10.2. Example.} We consider a special point ${\rm Spec}(B(k))\to\Mn$ (see def. [Va2, 2.10]); as $\Mn$ has the extension property (see [Va2, 3.2.3 3) and 6)]), we get that it extends to a morphism $z_1:{\rm Spec}(W(k))\to\Mn$. As $\FF$ is a subfield of $k$ (see 2.3.9 A), and as we are dealing with special points, we can assume $k=\FF$. So, following the proof of 2.2.18 we get: the Shimura filtered $\sg$-crystal attached to $z_1$ is cyclic diagonalizable.
\medskip
{\bf 2.3.11. Some simple properties.} Let $(f,L_{(p)},v)$ be a SHS. All the local deformation theory presented in [Va2, 5.4] can be applied for it (for $\Mn$), cf. the reductiveness property expressed in 2.3.4. In particular, the fact that $\Mn$ is the normalization of the Zariski closure of ${\rm Sh}_H(G,X)$ in $\Mm$, results
also from [Va2, 5.4-5] (see [Va2, 5.4.8 and 5.5]). 
\smallskip
We have variants of [Va2, 5.4]: we work over an arbitrary perfect field $k$ (instead of $\FF$), or we work with $G^0_{W(k)}$ (instead of $G^{\rm der}_{W(k)}$), or we work with $\Mn/H_0$ (instead of $\Mn$, i.e. we work with a finite level symplectic similitude structure, instead of all level-$N$ symplectic similitude structures, $N\in\NN$, $(N,p)=1$). This is so due to the fact that [Fa2, th. 10] is true for any perfect field. Here we just state one such form of these variants. We consider a point $z\in \Mn_{W(k)}(W(k))$. Let $(M,F^1,\vph,G_{W(k)},(t_{\al})_{\al\in\Mj},p_M)$ be its attached principally quasi-polarized Shimura filtered $\sg$-crystal. We have:
\medskip
{\bf Theorem.} {\it There are formally smooth $W(k)$-morphisms from the completion $\tilde G$ of $G^{\rm der}_{W(k)}$ in its origin into $\Mn_{W(k)}$ such that:
\medskip
-- the origin of $\tilde G$ is mapped into $z$;
\smallskip
-- the principally quasi-polarized filtered $F$-crystal with tensors over $\tilde G_k$ defined naturally (via pull back) by $(\Ma,\Mp_{\Ma})$ and its natural family of (de Rham components of) Hodge cycles, is isomorphic to a principally quasi-polarized Shimura filtered $F$-crystal of the form 
$$(M,F^1,\vph,G_{W(k)},G^{\rm der}_{W(k)},\tilde f,(t_{\al})_{\al\in\Mj},p_M)$$under an isomorphism which in $z$ is the natural identification of $(M,F^1,\vph,G)$ with itself.}
\medskip
From [Va2, 5.4.8] we get (as in [Va2, 5.5]):
\medskip
{\bf Corollary.} {\it The finite $W(k)$-morphism $\Mn_{W(k)}\to\Mm_{W(k)}$ defined by $i_{\Mn}$ is a formal embedding in each $k$-valued point of $\Mn_{W(k)}$.} 
\medskip
We consider now the locally free $\Mo_{\Mn/H_0}$-sheaf 
$$\Ml_{H_0}:=R^1_{dR}(\Ma_{H_0}/\Mn/H_0)$$ 
of rank $2e$. $\Mp_{\Ma_{H_0}}$ makes it to be naturally equipped with a perfect alternating form $\psi_{\Ml_{H_0}}$. Let $\Ml$ be the pull back of $\Ml_{H_0}$ to $\Mn$. We have:
$$\Ml=R^1_{dR}(\Ma/\Mn).$$
The essential tensor algebra $\Mt(\Ml_{E(G,X)})$ of the direct sum of the pull back of $\Ml$ to $\Mn_{E(G,X)}$ and of the dual of it, comes equipped with a family of sections $(t_{\alpha}^\Ma)_{\alpha\in\Mj^\prime}$ defined by the de Rham components of the family of Hodge cycles 
$(w_{\alpha}^\Ma)_{\alpha\in\Mj^\prime}$ of $\Ma$. Let 
$$\Mg_{\Mn}$$ 
be the Zariski closure in $GL(\Ml)$ of the subgroup of $GL(\Ml_{E(G,X)})$ fixing these sections. Working in the faithfully flat topology, from Theorem we get:
\medskip
{\bf Fact 1.} {\it $\Mg_{\Mn}$ is a reductive subgroup of $GL(\Ml)$.} 
\medskip
Let now $H_1$ be a normal, open subgroup of $H_0$ such that $\Mg_{\Mn}$ is obtained from a subgroup of $GL(\Ml_{H_1})$ by pull back via the quotient morphism $\Mn\to\Mn/H_1$; here $\Ml_{H_1}$ is the pull back of $\Ml_{H_0}$ via the quotient morphism $i_{10}:\Mn/H_1\to\Mn/H_0$. We assume that all the sections of $\Mt(\Ml_{E(G,X)})$ we have considered are obtained by pull back from sections of the pull back of $\Mt(\Ml_{H_1})$ to the generic fibre of $\Mn/H_1$ (i.e. we assume $H_1$ is small enough). As the quotient morphism $\Mn/H_1\to\Mn/H_0$ is an \'etale cover and we can assume that the finite group $H_0/H_1$ permutes the (pull back to generic points of $\Mn/H_1$ of the) last family of sections, we get (see also 2.3.10):
\medskip
{\bf Fact 2.} {\it $\Mg_{\Mn}$ is obtained by pull back from a reductive subgroup $\Mg_{\Mn/H_0}$ of $GL(\Ml_{H_0})$.}
\medskip
We consider the principally quasi-polarized $p$-divisible object 
$$({\got L}_{H_0},p_{{\got L}_{H_0}})$$ 
of $\Mm\Mf_{[0,1]}^\nabla(\Mn_{W(k(v))}/H_0)$ corresponding to the principally quasi-polarized $p$-divisible group 
$$(\Md_{H_0},\Mp_{\Md_{H_0}})$$ of
$(\Ma_{H_0},\Mp_{\Ma_{H_0}})$. The underlying $\Mo_{(\Mn_{W(k(v))}/H_0)^\wedge}$-sheaf of ${\got L}_{H_0}$ is the $p$-adic completion of $\Ml_{H_0}$, while its connection is the $p$-adic completion of the Gauss--Manin connection on $\Ml_{H_0}$ defined by $\Ma_{H_0}$. From Fact 2 and 2.2.20.1 9) applied to $i_{10}^*({\got L}_{H_0})$ and from Corollary of 2.3.10 applied generically (in the context of $i_{10}^*(\Ml_{H_0})$, i.e. in the context of $\Mn/H_1$), we get: the $End$ $p$-divisible object $End({\got L}_{H_0})$ of 
$\Mm\Mf_{[-1,1]}^\nabla(\Mn_{W(k(v))}/H_0)$ has a natural Lie $p$-divisible subobject 
$${\got G}_{H_0},$$ 
whose underlying Lie $\Mo_{(\Mn_{W(k(v))}/H_0)^\wedge}$-sheaf is the $p$-adic completion of the Lie algebra of $\Mg_{\Mn/H_0}$. 
We conclude:
\medskip
{\bf Fact 3.} {\it For any morphism $U\to\Mn$ (resp. $U\to\Mn/H_1$), with $U$ a regular, formally smooth scheme over a faithfully flat $O_{(v)}$-algebra which is a DVR having a perfect residue field and index of ramification 1, we get naturally a principally quasi-polarized Shimura $p$-divisible group 
$$(D_U,(t_{\al}^D)_{\al\in\Mj},p_{D_U})$$ 
over $U$.}
\medskip
 Here $(D_U,p_{D_U})$ is the principally quasi-polarized $p$-divisible group over $U$ obtained by pulling back $(\Md_{H_0},\Mp_{\Md_{H_0}})$ via the composite morphism $U\to\Mn/H_0$, while the family of crystalline sections $(t_{\al}^D)_{\al\in\Mj}$ is naturally obtained from the mentioned sections of the generic fibre of $\Mt(\Ml_{H_1})$. The change of the original family of tensors $(v_{\al})_{\al\in\Mj}$ (for instance, if we work with $\Mj^\prime$ instead of $\Mj$), corresponds to a passage to quasi-isomorphic principally quasi-polarized Shimura $p$-divisible groups, cf. 2.2.20.1 7). We have a logical version of Fact 3 in terms of $F$-crystals: such a version encompasses the part of 2.3.10 referring to 2.3.4. 
\smallskip
From 2.2.21 and the above first two paragraphs (referring to [Va2, 5.4-5]), we get:
\medskip
{\bf Fact 4.} {\it If the morphism $U\to\Mn/H_1$ identifies $U$ with the completion of $\Mn_{W(k)}/H_1$ in a $W(k)$-valued point of it, then $(D_U,(t_{\al}^D)_{\al\in\Mj},p_{D_U})$ is a universal principally quasi-polarized Shimura $p$-divisible group.} 
\medskip
The pull back of $(\Ma,\Mp_{\Ma})$ via an automorphism of $\Mn$ defined by the right translation with an element of $G(\AA_f^p)$, gives birth to a principally polarized abelian scheme over $\Mn$ which is naturally $\ZZ_{(p)}$-isogeneous to $(\Ma,\Mp_{\Ma})$ (see [Va2, p. 454]); strictly speaking, cf. loc. cit., for this statement we have to consider principal polarizations up to $\GG_m(\ZZ_{(p)})$-multiples. So we have a natural continuous right action of $G(\AA_f^p)$ on $\Ml$, viewed as an affine scheme over $\Mn$. This implies that $\Ml$ does not depend on the choice of a $\ZZ$-lattice $L$ of $W$ as in 2.3.2. We can assume that under this action the considered global sections of $\Ml_{E(G,X)}$ are permuted. We get (cf. also AE.4.2 a)): 
\medskip
{\bf Fact 5.} {\it Let $\tilde H_0$ be a compact subgroup of $G(\AA_f^p)$ such that either $\tilde H_0\times H$ is $p$-smooth for $(G,X)$ or $p\not | t(G^{\rm ad})$ and $\tilde H_0\times H$ is smooth for $(G,X)$. Let $\Ml_{\tilde H_0}$ be the quotient of $\Ml$ under the action of $\tilde H_0$ on it. Then in Fact 2 we can replace $H_0$ by $\tilde H_0$.}
\medskip
{\bf Fact 6.} {\it Under right translation by elements of $G(\AA_f^p)$, the Shimura (resp. Shimura filtered) $\sg$-crystals attached to $k$-valued (resp. to $W(k)$-valued) points of $\Mn$ are the same. If $k=\bar k$ and if $\bigl({-1\over p}\bigr)=-1$, then the same holds  in the principally quasi-polarized context.}
\medskip
For the principally quasi-polarized context, we just need to add: if $\bigl({-1\over p}\bigr)=-1$, then (as $p\ge 3$) $\forall\be\in\GG_m(\ZZ_p)$, either $\be$ or $-\be$ is the square of an element of $\GG_m(\ZZ_p)$. So if $k=\bar k$ and if $\bigl({-1\over p}\bigr)=-1$, then $\forall\al_N\in\GG_m(\ZZ_p)$, a principally quasi-polarized Shimura $\sg$-crystal $(N,\vph_N,G_N,p_N)$ over $k$ is isomorphic to $(N,\vph_N,G_N,\al_Np_N)$ under a scalar automorphism of $N$ defined by an element of $\GG_m(\ZZ_p)$.
\medskip
Let now $O_{1}$ be a local, formally \'etale, faithfully flat $O_{(v)}$-algebra such that the representation of $G_{O_{1}}$ on $L_{(p)}^*\otimes_{\ZZ_{(p)}} O_{1}$ is a direct sum of two non-zero subrepresentations. This can be codified in terms of projectors which are $O_{1}$-linear combinations of Betti components of Hodge cycles which are defined by endomorphisms of $\Ma$. So, if these endomorphisms are obtained from endomorphisms of $\Ma_{H_0}$ via pull back, we get:
\medskip
{\bf Fact 7.} {\it The pull back of $\Ml_{H_0}$ to $\Mn_{O_1}/H_0$ gets naturally a direct sum decomposition into two non-zero summands.}
\medskip
Fact 7 just points out a first property of the general theory of automorphic vector bundles in mixed characteristic, whose elaboration will be started in \S5 (cf. end of 1.15).
\medskip
{\bf 2.3.12. Some extra standard notations.} We start with a closed point
$z:{\rm Spec}(W(k))\hookrightarrow\Mn_{W(k)}/H_0$. Let $y:{\rm Spec}(k)\hookrightarrow\Mn_{k}/H_0$ be defined by $z$. We use the notations of 2.3.10. So $F^1$ is the Hodge filtration of $M$ defined by $A$. For simplifying the presentation, we assume that the family of Hodge cycles $(w_{\al})_{\al\in\Mj^\prime}$ of $A_{W(k_1)}$ is definable over $W(k)$ and that the group $\tilde G_{W(k)}$ of 2.3.10 is in fact (isomorphic to) $G_{W(k)}$. We can always achieve  this by replacing $k$ by a finite field extension $k_1$ of it. Let $O_y$ be the local ring of $y$. So $O_y^{\rm h}$ and $\widehat {O_y}$ are its henselization and respectively its completion. Let $(A^{\rm h}_y,p_{A^{\rm h}_y})$ be the principally polarized abelian scheme over
${\rm Spec}(O^{\rm h}_y)$ obtained from $(\Ma_{H_0},\Mp_{\Ma_{H_0}})$ by pull back through the natural morphism ${\rm Spec}(O_y^{\rm h})\to\Mn/H_0$. Let 
$$
M^{\rm h}_y:=H^1_{dR}(A_y^{\rm h}/O^{\rm h}_y), 
$$
and let $F_y^{1{\rm h}}$ be its Hodge filtration defined by $A_y^{\rm h}$. 
\smallskip
As $\Mn$ is formally smooth over $\ZZ_{(p)}$, $O^{\rm h}_y$ is an integral domain. Let $K_y^{\rm h}$ be the field of fractions of $O^{\rm h}_y$.
Let $(w^{\rm h}_\al)_{\al\in\Mj^\prime}$ be the family of Hodge cycles with which  $A_{yK_y^{\rm h}}^{\rm h}$ is naturally endowed: as we assumed that the family $(w_{\al})_{\al\in\Mj^\prime}$ is defined over $W(k)$, the family $(w^{\rm h}_\al)_{\al\in\Mj^\prime}$ is defined over $O_y^{\rm h}$ and not only over its  strict henselization. This is a consequence of the following obvious fact:
\medskip
{\bf Fact.} {\it A Hodge cycle of a deformation which is $0$ in a point is $0$ (everywhere).} 
\medskip
Let 
$t^{\rm h}_\al$ be the de Rham component of $w_\al^{\rm h}$. It is a homogeneous tensor of $\Mt(M^{\rm h}_y\otimes_{O_y^{\rm h}} K_y^{\rm h})$.
\medskip
{\bf 2.3.13. Lemma.} {\it There is an isomorphism $f_y:M^{\rm h}_y\tilde\to M\otimes_{W(k)}O^{\rm h}_y$ taking $F^{1{\rm h}}_y$ onto $F^1\otimes_{W(k)} O^{\rm h}_y$, $p_{A_y^{\rm h}}$ into $p_A$,  and  $t^{\rm h}_\al$ into $t_\al$, $\forall\al\in\Mj^\prime$ (here we view $p_{A_y^{\rm h}}$ and $p_A$ as perfect forms on $M_y^{\rm h}$ and respectively on $M\otimes_{W(k)}O^{\rm h}_y$).}
\medskip
{\bf Proof:} Such an isomorphism exists over $\widehat {O_y}$, cf. Fact 4 of 2.3.11. The obstruction to the existence of $f_y$ is measured by an element $\ga\in H^1_{ff}(O^{\rm h}_y,P^0_{O^{\rm h}_y})$; here $P^0$ is the parabolic subgroup of $G^0_{W(k)}$ normalizing $F^1$, while the lower right index $ff$ refers to the faithfully flat topology. 
So $\ga$ becomes $0$ over $\widehat {O_y}$. As $O_y^{\rm h}$ is the henselization of the excellent ring $O_y$, using Artin's approximation theorem we get $\ga=0$. This proves the Lemma.
\medskip
{\bf 2.3.13.1. Corollary.} {\it There is an \'etale morphism $a:Y={\rm Spec}(R)\to\Mn_{W(k)}/H_0$ to which the point $y$ lifts, and  there is an isomorphism $M_R{\tilde\to} M\otimes_{W(k)}R$ taking $F^1_R$ onto $F^1\otimes_{W(k)}R$, $p_{A_Y}$ into $p_A$, and  $t^R_\al$ into $t_\al$, $\forall\al\in\Mj^\prime$ (here $M_R$, $F^1_R$, $(t^R_\al)_{\al\in\Mj^\prime}$, $p_{A_Y}$ have a meaning similar to the one in the case of $O^{\rm h}_y$; for instance, $M_R:=H^1_{dR}(A_Y/R$), where $A_Y:=\Ma_{H_0}\times_{\Mn/H_0} Y$,  the morphism $Y\to\Mn/H_0$ being naturally induced by $a$).}
\medskip
{\bf 2.3.13.2. Corollary.} {\it To check that 2.3.4 (*) holds, we just need to consider $z\in\Mn(W(\FF))$.} 
\medskip
{\bf 2.3.14. Warning: ${\rm Sp}$ or ${\rm GSp}$?}
Once for all we have to decide if we use ${\rm GSp}(W,\psi)$ or ${\rm Sp}(W,\psi)$ (and so $G_{W(k)}$ or $G_{W(k)}^0$). The use of ${\rm Sp}(W,\psi)$ corresponds to a permanent mentioning of polarizations (and so to a very precise presentation but more complicated notations). The use of ${\rm GSp}(W,\psi)$ corresponds to a loose treatment of polarizations (and so to a less precise presentation but simpler notations). We decided to go along with ${\rm Sp}(W,\psi)$: this explains the form of 2.3.15 below.
\smallskip 
The idea behind this comment is: we can work out [Va2, 5.4.4-5] using not  only (the $W(k)$-versions of) ${\rm Sp}(W,\psi)$ and $G^{\rm der}$ (or of $G^0$) but also (the $W(k)$-versions of) ${\rm GSp}$ and $G$. The use of ${\rm GSp}(W,\psi)$ means that the cycle of the polarization is not exactly $p_A$ but $p_A$ times some invertible element $u_p$. To exemplify what we mean by this, we consider a particular case: denoting by ${\rm Spec}(R_1)$ the completion of $GSp(M,p_A)$ in its origin, the natural universal element $g_{R_1}\in GSp(M,p_A)(R_1)$ takes $p_A$ into $\gamma_pp_A$, where $\gamma_p\in\GG_m(R_1)$ is congruent to the identity modulo the ideal $I_1$ of $R_1$ defining the origin of ${\rm Spec}(R_1)$. If $\Phi_{R_1}$ is a Frobenius lift of $R_1$ as in 2.2.10 (so it takes $I_1$ into $I_1^p$), then $u_p\in\GG_m(R_1)$ and $\gamma_p$ are related through the formula $\Phi_{R_1}(u_p)\gamma_p=u_p$. It is easy to see that there is a unique such $u_p$ which mod $I_1$ is the identity.
\smallskip
In general, this element $\gamma_p$ can be made to be the identity, by considering a square root $x$ of $\gamma_p$ (as $p>2$, $x\in R_1$); in our example involving $R_1$, the automorphism $g_{x^{-1}}$ of $M\otimes_{W(k)} R_1$ defined by scalar multiplication with $x^{-1}$, ``brings back" $\gamma_pp_A$ into $p_A$. So, to keep $p_A$ as the cycle of a principal quasi-polarization, we have to replace $g_{R_1}$ by $g_{x^{-1}}g_{R_1}\in Sp(M,p_A)(R_1)$.
\smallskip
In fact, due to the existence of the cocharacter $\tilde\mu$ of 2.3.10, in general, we can make $\gamma_p$ to be the identity element without using that $p$ is odd; for instance in the example involving $R_1$, instead of $g_{x^{-1}}$ we can use $\tilde\mu(\gamma_p)$ in order to bring back $\gamma_pp_A$ into $p_A$. 
\smallskip
So indeed the difference in the mentioned two choices is mostly of notations and of preciseness.      
\medskip
{\bf 2.3.15. Proposition.} {\it There is a smooth morphism $a:Y={\rm Spec}(R)\to\Mn_{W(k)}/H_0$ such that:
\medskip
a) The scheme $a^{-1}(z)$ contains a $W(k)$-valued point $z_1$;
\smallskip
b) We have an isomorphism $M_R\tilde\to M\otimes_{W(k)}R$ taking $F^1_R$ onto $F^1\otimes_{W(k)}R$, $p_{A_Y}$ into $p_A$, and $t^R_\al$ into $t_\al$, $\forall\al\in\Mj^\prime$ (here $M_R:=H^1_{dR}(A_Y/R)$, with $A_Y$ defined as in 2.3.13.1, etc.; the logical notations);
\smallskip
c) There is a Frobenius lift $\Phi_R$ of $R^\wedge$ of the form $z_i\to z^p_i$, $i\in I$, where the elements of the subset $\{z_i|i\in I\}$ of $R$ are forming a  regular system of parameters of $z_1$ in $Y$;
\smallskip
d) Based on b), the $p$-divisible object ${\got C}$ of $\Mm\Mf_{[0,1]}^\nabla(R)$ associated to the $p$-divisible group of 
$A_Y$, can be put in the form 
$\bigl(M\otimes_{W(k)} R^\wedge,F^1\otimes_{W(k)} R^\wedge, g_Y(\vph\otimes 1),\nabla_{R^\wedge}\bigr)$, with the perfect alternating form $p_A$ on $ M\otimes_{W(k)} R^\wedge$ corresponding to $p_{A_Y}$, with $\nabla_{R^\wedge}$ a connection on $M\otimes_{W(k)} R^\wedge$ annihilating $t_{\al}$, $\forall\al\in\Mj^\prime$, and with $g_Y\in G^0_{W(k)}(R^\wedge)$ defined by a morphism (still denoted by $g_Y$)
$$g_Y:Y^\wedge\to G^{0\wedge}_{W(k)}.$$ 
Moreover, we can assume the special fibre of $g_Y$ is \'etale (so $\abs{I}=\dim_{\QQ}(G^0_{\QQ})$) and the composite of (the $p$-adic completion of) $z_1$ with $g_Y$ defines the origin of $G^{0\wedge}_{W(k)}$.}
\medskip
{\bf Proof:} Let $\hat G^0={\rm Spec}(\hat R^0)$ be the completion of $G^0_{W(k)}$ in its origin. As in [Va2, 5.4.4] (see 2.3.11) we construct a formally smooth morphism 
$$\hat G^0\buildrel{\hat a}\over\to\Mn_{W(k)}/H_0$$ 
such that $z$ factors through the origin of $\hat G^0$ and the principally quasi-polarized filtered $F$-crystal with tensors $\hat {\got C}^0$ associated to the principally quasi-polarized $p$-divisible group over $\hat G^0$ obtained naturally from $(\Md_{H_0},\Mp_{\Md_{H_0}})$ (by pull back through the composite of $\hat a$ with the natural morphism $\Mn_{W(k)}/H_0\to\Mn/H_0$) and its natural family of de Rham sections (see Fact 3 of 2.3.11 and 2.3.12), is a principally quasi-polarized Shimura filtered $F$-crystal $(M,F^1,\vph,G_{W(k)},G_{W(k)}^0,\tilde f,(t_{\al})_{\al\in\Mj^\prime},p_A)$ (so $z^*(\hat {\got C}^0)$ is naturally identified with $(M,F^1,\vph,G,(t_{\al})_{\al\in\Mj^\prime})$). Here $\tilde f:W(k)[[x_1,...,x_{\abs{I}}]]\tilde\to\hat R^0$ is a $W(k)$-isomorphism. For the part referring to formal smoothness, cf. also [Va2, 5.4.7-8].
\smallskip
Let $n\ge 2$ be an integer and let $I$ be the ideal of $\hat R^0$ defining the origin of $\hat G^0$ (so $\hat R^0/I=W(k)$). Let (cf. Artin's approximation theorem) $Y={\rm Spec}(R)$ be a smooth $W(k)$-scheme such that:
\medskip
i) there is a smooth morphism $a:Y\to\Mn_{W(k)}/H_0$ such that the point $z$ lifts to a point $z_1:{\rm Spec}(W(k))\to Y$;
\smallskip
ii) if $I_1$ is the ideal of $R$ defining $z_1$, the natural morphism ${\rm Spec}(R/I_1^n)\to\Mn_{W(k)}/H_0$ defined by $a$ is obtained from the natural morphism ${\rm Spec}(\hat R^0/I^n)\to\Mn_{W(k)}/H_0$ defined by $\hat a$ via a $W(k)$-isomorphism $i_n:\hat R^0/I^n\tilde\to R/I_1^n$.
\medskip
Replacing if needed $Y$ with $Y^\prime={\rm Spec}(R^\prime)$, where $Y^\prime\to Y$ is an \'etale morphism such that $z_1$ lifts to a $W(k)$-valued point of $Y^{\prime}$, we can achieve (cf. 2.3.13.1) that the morphism $a$ endows the abelian scheme $A_Y$ with a family of Hodge cycles $(w_\al^R)_{\al\in\Mj^\prime}$ and that  b) and c) are satisfied (here $z_i\in R$ is such that $i_n(\tilde f(x_i))$ is congruent to $z_i$ modulo $I_1^n$). This takes care of a), b) and c). 
\smallskip
Now the principally quasi-polarized $p$-divisible object ${\got C}_R$ of $\Mm\Mf^\nabla_{[0,1]}(R)$ associated to the principally quasi-polarized $p$-divisible group of $(A_Y,p_{A_Y})$ can be put (cf. b)) in the form
$\bigl(M\otimes_{W(k)} R^\wedge,F^1\otimes_{W(k)} R^\wedge,g_Y(\vph\otimes 1),\nabla_{R^\wedge},p_A\bigr)$, with $g_Y\in G^0_{W(k)}(R^\wedge)$ (as $p_{A_Y}$ is fixed by $g_Y$ and as the tensors $t^R_\al$, $\al\in\Mj^\prime$, are fixed --this is a consequence of Corollary of 2.3.10 applied generically-- by the $\Phi_R$-linear endomorphism of $M\otimes_{W(k)} R^\wedge$ defined by ${\got C}_R$) and with $\nabla_{R^\wedge}$ a connection on $M\otimes_{W(k)} R^\wedge$ annihilating $t_\al$, $\forall\al\in\Mj^\prime$. 
\smallskip
Let $P^0$ be as in the proof of 2.3.13. Due the smoothness of $P^0$ over $W(k)$, any element of $P^0(R/I_1^n)$ congruent to the identity mod $I_1$, lifts, potentially after a replacement of $Y$ with an $Y^\prime$ as above, to an element of $P^0(R)$ and so to an element $h_{\rm corr}$ of $P^0(R^\wedge)$. From property ii), by replacing $g_Y(\vph\otimes 1)$ by $h_{\rm corr}g_Y(\vph\otimes 1)h_{\rm corr}^{-1}$ for a suitable such correction element $h_{\rm corr}$, we get that we can assume (simple argument at the level of tangent spaces) that $g_Y$ is defined by a $W(k)$-morphism $g_Y:Y^\wedge\to G^{0\wedge}_{W(k)}$, \'etale mod $p$ and having the property that the composite of the $p$-adic completion of $z_1$ with $g_Y$ defines the origin of $G_{W(k)}^{0\wedge}$. This takes care of d) and ends the proof.
\medskip
From the proof of 2.3.15, by removing extra regular parameters we get:
\medskip
{\bf 2.3.15.1. Corollary.} {\it There is an \'etale morphism $\tilde a:\tilde Y={\rm Spec}(\tilde R)\to\Mn_{W(k)}/H_0$ such that:
\medskip
a) The scheme $\tilde a^{-1}(z)$ contains a $W(k)$-valued point $\tilde z_1$;
\smallskip
b) We have an isomorphism $M_{\tilde R}\tilde\to M\otimes_{W(k)}\tilde R$ taking $F^1_{\tilde R}$ onto $F^1\otimes_{W(k)}\tilde R$, $p_{A_{\tilde Y}}$ into $p_A$,  and $t^{\tilde R}_\al$ into $t_\al$, $\forall\al\in\Mj^\prime$ (here $M_{\tilde R}:=H^1_{dR}(A_{\tilde Y}/\tilde R)$, with $A_{\tilde Y}$ defined as in 2.3.13.1, etc.; the logical notations);
\smallskip
c) There is a Frobenius lift $\Phi_{\tilde R}$ of $\tilde R^\wedge$ of the form $z_i\to z^p_i$, $i\in\tilde I$, where the elements of the subset $\{z_i|i\in\tilde I\}$ of $\tilde R$ are forming a system of regular parameters of $\tilde z_1$ in $\tilde Y$ (so $|\tilde I|=\dim_{\CC}(X)=\dim_{E(G,X)}(Sh(G,X))$);
\smallskip
d) Based on b), the $p$-divisible object of $\Mm\Mf_{[0,1]}^\nabla(\tilde R)$ associated to the $p$-divisible group of $A_{\tilde Y}$, can be put in the form 
$\bigl(M\otimes_{W(k)} \tilde R^\wedge,F^1\otimes_{W(k)} \tilde R^\wedge,g_{\tilde Y}(\vph\otimes 1),\nabla_{\tilde R^\wedge}\bigr)$, with the alternating form $p_A$ on $M\otimes_{W(k)} \tilde R^\wedge$ corresponds to $p_{A_{\tilde Y}}$, with the connection $\nabla_{\tilde Y^\wedge}$ on $M\otimes_{W(k)} \tilde R^\wedge$ annihilating $t_\al$, $\forall\al\in\Mj^\prime$, and with  $g_{\tilde Y}\in G^0_{W(k)}(\tilde R^\wedge)$ defined by a morphism (still denoted by $g_{\tilde Y}$) $g_{\tilde Y}:\tilde Y^\wedge\to G^{0\wedge}_{W(k)}$. Moreover, we can assume that:
\medskip
{\item{--}} the composite of (the $p$-adic completion of) $\tilde z_1$ with $g_{\tilde Y}$ defines the origin of $G^{0\wedge}_{W(k)}$;
\smallskip
{\item{--}} modulo the ideal $J:=p((z_i)_{i\in\tilde I})^2$ of the completion $\widehat{\tilde R}$ of $\tilde R$ w.r.t. its ideal $((z_i)_{i\in\tilde I})$, $g_{\tilde Y}$ defines an isomorphism $ISO$ from ${\rm Spec}(\widehat{\tilde R}/J)$ into the completion $\hat N$ of a commutative, unipotent subgroup $N$ of $G^0_{W(k)}$ in its origin, taken  modulo the ideal defined by $p$ times the second power of the ideal sheaf of $\hat N$ defining its origin; $N$ can be any such subgroup achieving an open embedding $N\hookrightarrow G_{W(k)}/P$ (so $g_{\tilde Y}$ mod $p$ is a formal closed embedding in the special fibre $\tilde y_1$ of $\tilde z_1$) and even more, if a Frobenius lift of $\hat N={\rm Spec}(W(k)[[\om_1,...,\om_{\dim_{\CC}(X)}]])$ taking $\om_l$ to $\om_l^p$, $l=\overline{1,\dim_{\CC}(X)}$, is chosen, we can also assume that such an isomorphism $ISO$ is compatible with the Frobenius lifts;
\medskip
e) The special fibre of $\tilde Y$ is connected.}
\medskip
The c) of 2.3.15 (resp. of 2.3.15.1) is automatically satisfied by passing to a suitable open, affine subscheme of $Y$ (resp. of $\tilde Y$); it is included in order to state its following d) accurately.
\medskip
{\bf 2.3.16. Corollary.} {\it There is an open, dense subscheme $U_k$ of $G_k$ such that: for any geometric point ${\rm Spec}(k_1)\buildrel{y_1}\over\to U_k$ 
(with $k_1=\overline{k_1}\supset k$),
there is a point $z_1\in \Mn/H_0(W(k_1))$, lifting a $k_1$-valued point of the same connected component of $\Mn_{\bar k}/H_0$ through which the point $y\circ i_k$ factors, with $i_k:{\rm Spec}(\bar k)\to {\rm Spec}(k)$ as the natural morphism, and having the property that the Shimura filtered
$\sg_{k_1}$-crystal attached to it is isomorphic to $\bigl(M\otimes_{W(k)} W(k_1),F^1\otimes_{W(k)} W(k_1),g_1(\vph_0\otimes 1),G_{W(k_1)}\bigr)$, where $g_1\in G_{W(k)}\bigl(W(k_1)\bigr)$ mod $p$ is (i.e. defines) the point $y_1$.}
\medskip
{\bf Proof:} Working with $G^0$ instead of $G$, this is a direct consequence of 2.3.15: we can take as $U^0_k$ the image of the special fibre of $Y^\wedge$ through the formally \'etale morphism $g_Y:Y^\wedge\to G^{0\wedge}_{W(k)}$ of 2.3.15; this can be read out from the existence of Teichm\"uller lifts ${\rm Spec}(W(\overline{k_1}))\to Y^\wedge$. Now, as $U_k$, we can take (cf. the Subfact of 2.3.11) an open subscheme of $G_k$ contained in the image of $U^0_k\times\GG_m$ through the natural isogeny (see 2.3.1) $G^0_k\times\GG_m\to G_k$. 
\smallskip
The argument for this (to be compared with 2.3.14) goes as follows: for any $\al\in\GG_m(W(k))$, there is $\be\in\GG_m(W(k))$ such that $\al=\sg(\be)\be^{-1}$; accordingly a Shimura $\sg$-crystal $(N,\vph_N,G_N)$ over $k$ is isomorphic to the Shimura $\sg$-crystal $(N,\al_N\vph_N,G_N)$, where $\al_N$ is the automorphism of $N$ defined by scalar multiplication with $\al$. This proves the Corollary.
\medskip
{\bf 2.3.16.1. Comment.} Warning: by working with $G_{W(k)}$ and not with $G^0_{W(k)}$, we have to give up mentioning the principal quasi-polarizations; so the isomorphism of 2.3.16 is just of Shimura filtered $\sg$-crystals. For applications to problems involving Newton polygons or truncations of Shimura Lie $F$-crystals this is equally good. If we want to keep track of the principal quasi-polarizations, we have to restate 2.3.16, with $U_k$ being replaced by $U_k^0$ of its proof.
\medskip
{\bf 2.3.17. Proposition.} {\it With the notations of 2.3.12, the Shimura filtered $\sg$-crystal
$(M,F^1_1,\vph,G_{W(k)})$ attached to another lift 
$z_1:{\rm Spec}(W(k))\hookrightarrow\Mn_{W(k)}/H_0$ of $y$, satisfies 
$$F^1_1=g(F^1),$$ 
with $g\in G_{W(k)}(W(k))$ such that mod $p$ it normalizing $F^1/pF^1$. Conversely, for any such $g$ there is a unique lift $z_g:{\rm Spec}(W(k))\hookrightarrow\Mn_{W(k)}$ of $y$, whose attached Shimura filtered $\sg$-crystal is $(M,g(F^1),\vph,G_{W(k)})$.} 
\smallskip
{\bf Proof:} As the parabolic subgroups of $G_{W(k)}$ normalizing $F^1$ and respectively $F^1_1$ are having the same special fibre, the first statement follows from Fact 2 of 2.2.9 3) applied in the context of the canonical split cocharacters (they factor through $G_{W(k)}$) of $(M,F^1,\vph)$ and of $(M,F^1_1,\vph)$. Also we would like to point out, that the approach of [Va2, 5.3.1] represents another way to see that these two parabolic subgroups are $G_{W(k)}(W(k))$-conjugate. 
\smallskip
The converse is a direct consequence of the smoothness of $G_{W(k)}$ and of $\Mn_{W(k)}$, of the Corollary of 2.3.11 (or of Fact 4 of 2.3.11) and of the local deformation theory of abelian varieties (over $k$) (see also the existence of the isomorphism $FIL_0$ of 2.4 below). This proves the Proposition.
\medskip
{\bf 2.3.17.1. Corollary.} {\it  Let $g\in G_{W(k)}(W(k))$ be such that mod $p$ it normalizes $F^1/pF^1\subset M/pM$. Then the principally quasi-polarized Shimura filtered $\sg$-crystal $\break (M,g(F^1),\vph,G_{W(k)},p_A)$ is induced from the principally quasi-polarized Shimura filtered F-crystal 
$$(M,F^1,\vph,G_{W(k)},G^0_{W(k)},\tilde f,p_A)$$ 
(of the proof of 2.3.15) through a $W(k)$-epimorphism $\hat R^0\twoheadrightarrow W(k)$.}
\medskip
The 2.2.21 UP suggests (to be compared with 3.12 below) variants of this Corollary, where we replace $G^0_{W(k)}$ by suitable smooth subschemes of $G_{W(k)}$ containing the origin.  
\medskip
{\bf 2.3.17.2. Exercise.} Show that 2.3.17.1 is a property of (principally quasi-polarized) Shimura $\sg$-crystals over a perfect field $k$ of characteristic $p\ge 3$. Hint one: just look at 2.4 below and use 2.2.21 UP. Hint two: if we feel more comfortable with the context of principally polarized abelian schemes, first assume $k=\bar k$ and follow the (independent) ideas of 4.12.12 and 4.12.12.0 or just 4.12.12.5 below, and then use Galois descent.
\medskip
{\bf 2.3.18. The $p=2$ SHS.} 
The case $p=2$ presents some particular features and so we felt it is appropriate to single it aside. 
\smallskip
{\bf A.} Working with $p=2$, we can keep 2.3.1-2. But for 2.3.3 we have to assume $\Mn$ exists: the proof of [Va2, 3.2.12] still applies (via [Va2, 3.4.1]) to give us that $\Mn$ is the Zariski closure of the normalization of ${\rm Sh}_H(G,X)$ in $\Mm$ (cf. 2.2.1.5.1). 2.3.3.1 needs no modification.
In connection to $\Mn^{\rm ad}$ (resp. to $\Mn^{\rm ab}$) of 2.3.3.2 we refer to 4.14.3.2.3 below (resp. to [Va2, 3.2.8]).
\smallskip
We keep 2.3.4-5 but we prefer to abbreviate the things differently: the resulting triple $(f,L_{(2)},v)$ (resp. quadruple $(f,L_{(2)},v,\Mb)$) is referred to as a $p=2$ SHS (resp. as a $p=2$ standard PEL situation). All of 2.3.6.1 makes sense for $p=2$. 2.3.5.1 remains true for $p=2$ as its hint does, while 2.3.5.3 makes sense for $p=2$ as well (cf. also B below). 
Related to Fact 1 of 2.3.5.2, we define $\Mn^j$ to be the normalization of the Zariski closure of ${\rm Sh}_{H^j}(G^j,X^j)$ in $\Mn$. The fact that $\Mn^j$ is smooth, $\forall j\in \Mf(G^{\rm ad},X^{\rm ad})$, is proved in the same manner as for the $A$ and $C$ cases of a $p=2$ standard PEL situation (this is reviewed in B below): $G^j_{\ZZ_{(2)}}$ is the centralizer of a torus $T$ of $G_{\ZZ_{(2)}}$ and so we have a relative PEL situation (in the sense of [Va2, 4.3.16]; see also [Va2, 4.3.14]); so one just needs to combine [Va2, 4.3.13] (applied in the context of the monomorphisms $T\hookrightarrow G_{\ZZ_{(2)}}\hookrightarrow GL(L_{(2)})$), with the Steps B4-10 of B below (i.e. with the fact that $(f,L_{(2)},v)$ is a $p=2$ SHS and with the deformation theory of endomorphisms of abelian schemes). 
Related to Fact 2 of 2.3.5.2 and to the part of 2.3.5.5 pertaining to $q_{\Mn_1}$, we refer to 4.14.3.2.3. 2.3.5.4 needs no modifications (as 2.2.1.5.1 handles the case $p=2$ as well). Except the mentioned part of 2.3.5.5, 2.3.5.5-6 and 2.3.5.6.1 A need no extra comment; the fact that 2.3.5.6.1 B holds for $p=2$ is implied by the Theorem of B below. However, we postpone to \S6 and [Va5] (resp. to 4.14.3.2) for a $p=2$ analogue of 2.3.5.7 (resp. of 2.3.5.8 and of 2.3.5.8.1). 
\smallskip
{\bf B.} The fact that 2.3.8 4) has a $p=2$ analogue
is just partially documented in the literature; for instance, in [Ko2] (see end of p. 391 of loc. cit.) the so called D case is avoided for $p=2$. The following Theorem fills out this gap in the literature and so it provides (via Lemma 2 of 4.6.4 below) the first new instances of $p=2$ SHS's.
\medskip
{\bf Theorem.} {\it We consider a quadruple $(f,L_{(2)},v,\Mb)$ which satisfies all properties of a $p=2$ standard PEL situation, except the fact that $(f,L_{(2)},v)$ is a $p=2$ SHS. Then the triple $(f,L_{(2)},v)$ is a $p=2$ SHS, and so the quadruple $(f,L_{(2)},v,\Mb)$ is a $p=2$ standard PEL situation.}
\medskip
{\bf Proof:} If $G$ is a torus, then the Theorem is easy: it is a consequence of [Va2, 4.3.13] (cf. also [Va2, 3.2.8 and 3.3.2]). From now on we assume $G$ is not a torus. We itemize the ideas by numbers attached to $B$.
\smallskip
{\bf B1.} We define $\Mn^\prime$ (resp. $\Mn$) to be the Zariski closure (resp. the normalization of the Zariski closure) of ${\rm Sh}_H(G,X)$ in $\Mm$. Let $(\Ma,\Mp_{\Ma})$ be defined as usual. In what follows, not to complicate the notations we still denote by $\Mb$ the set of $\ZZ_{(2)}$-endomorphisms with which $\Ma$ is naturally endowed and so, with which any abelian scheme obtained from $\Ma$ by pull back is endowed. Let $\Mb_1$ be the centralizer of $\Mb$ in ${\rm End}(L_{(2)})$. Let $*$ be the involution of $\Mb$ or of $\Mb_1$ defined by $\psi$. Till the end of 2.3.18, we use orderings of the form $(\Ma,\Mp_{\Ma},\Mb)$ instead of usual orderings of the form $(\Ma,\Mb,\Mp_{\Ma})$.
\smallskip
{\bf B2.} Let $n\in\NN$, $n\ge 4$. We consider the group scheme $\Md_n^1$ over $\ZZ_{(2)}$ fixing a quadratic form $QF:=x_1x_2+...+x_{2n-1}x_{2n}$ in $2n$ independent variables. It is known (see [Bo2, 23.6] for the picture over $\FF_2$) that the group $\Md_n$ over $\ZZ_{(2)}$ defined as the Zariski closure of the connected component of the origin of the generic fibre of $\Md_n^1$ in $\Md_n^1$, is a split, semisimple group of $D_n$ Lie type. So we get a $2n$ dimensional faithful representation 
$$\rho_n:\Md_n\hookrightarrow GL(\ZZ_{(2)}^{2n});$$
it is associated to the minimal weight $\om_1$ (see [De2] and [Bou2, planche IV]). Its special fibre is absolutely irreducible (simple argument at the level of dimensions) and (cf. [Bo2, 23.6]) it is alternating; so the perfect, symmetric bilinear form $BF$ on $L_n:=\ZZ_{(2)}^{2n}$ fixed by $\Md_n$ (and naturally defined by $QF$; it is unique up to multiplication with an invertible element of $\ZZ_{(p)}$), when taken mod $2$ is alternating. 
We have:
\medskip
{\bf Corollary.} {\it Let $R$ be a reduced, faithfully flat $\ZZ_{(2)}$-algebra. Let $M_R$ be a free $R$-submodule of $L_n\otimes_{\ZZ_{(2)}} R[{1\over 2}]$ of rank $2n$ and such that we get a perfect bilinear form $BF:M_R\otimes_R M_R\to R$ (so $M_R[{1\over 2}]=L_n\otimes_{\ZZ_{(2)}} R[{1\over 2}]$). If the Zariski closure $\Md_n(R)$ of ${\Md_n}_{R[{1\over 2}]}$ in $GL(M_R)$ is a semisimple group, then $BF$ restricted to $M_R/2M_R$ is alternating.}
\medskip
{\bf Proof:} Localizing in the \'etale topology, we can assume $\Md_n(R)$ is split and $R$ is local. We consider a Borel subgroup $B_R$ of ${\Md_n}_R$ and a maximal torus $T_R$ of $B_R$. Applying [SGA3, Vol. III, 1.5 of p. 329] over $R[{1\over 2}]$, we get that we can assume $T_{R[{1\over 2}]}$ and $B_{R[{1\over 2}]}$ extend to a torus and respectively to a Borel subgroup of $\Md_n(R)$. Based on [SGA3, p. 313-4 and 1.3 of p. 328] we deduce that we can assume $\Md_n(R)={\Md_n}_R$. But this implies that $2^{q(M_R)}M_R=L_n\otimes_{\ZZ_{(2)}} R$, for some $q(M_R)\in\ZZ$ (this can be checked immediately starting from the decompositions of $M_R$ and of $L_n\otimes_{\ZZ_{(2)}} R$ in irreducible $T_R$-modules and the fact that $\rho_n$ mod $2$ is irreducible). As $BF$ restricted to $M_R$ and to $L_n\otimes_{\ZZ_{(2)}} R$ are perfect, we get $q(M_R)=0$ and the Corollary follows. 
\medskip
We also get that the involution of $\Mb_1\otimes_{\ZZ_{(2)}} \FF$ has all its simple factors (i.e. the involutions restricted to simple factors of $\Mb_1\otimes_{\ZZ_{(2)}} \FF$ left invariant by $*$ or to products of two such simple factors permuted by $*$) are either of second type or of alternating first type. We refer to this property as the $ALT$ property. 
\smallskip
Parts a) and b) of the following Exercise complement parts of [Ko2] (like 7.2-3 of loc. cit., etc.). a) can be viewed as a converse of the Corollary (in our geometric context).
\medskip
{\bf B3. Exercise.} Let $R$ be a local $\ZZ_{(2)}$-algebra such that $\Mb\otimes_{\ZZ_{(2)}} R$ is a product of matrix $R$-algebras. Let $m_R$ be the maximal ideal of $R$; we assume $R/m_R$ has characteristic $2$. Let $\psi_1:L_{(2)}\otimes_{\ZZ_{(2)}} L_{(2)}\otimes_{\ZZ_{(2)}} R\to R$ be a perfect alternating form such that $\Mb\otimes_{\ZZ_{(2)}} R$ is self dual w.r.t. it and the resulting involution of $\Mb\otimes_{\ZZ_{(2)}} R$ is the same as the one defined by $\psi$. We assume the involution of $\Mb_1\otimes_{\ZZ_{(2)}} R/2R$ induced by $\psi_1$ has all its factors (in the same sense as above) either of second type or of alternating first type. We have:
\medskip
{\bf a)} We assume $R$ is a reduced, flat $\ZZ_{(2)}$-algebra. Then, the Zariski closure $\tilde G_R$ in $GL(L_{(2)}\otimes_{\ZZ_{(2)}} R)$ of the connected component of the origin of the subgroup of $GSp(L_{(2)}\otimes_{\ZZ_{(2)}} R[{1\over 2}],\psi_1)$ fixing the elements of $\Mb$, is reductive.
\smallskip
{\bf b)} We assume $R$ is a complete DVR of mixed characteristic (0,2) and with an algebraically closed residue field. We also assume the triples $(L_{(2)}\otimes_{\ZZ_{(2)}} R,\psi_1,\Mb\otimes_{\ZZ_{(2)}} R)$ and $(L_{(2)}\otimes_{\ZZ_{(2)}} R,\psi,\Mb\otimes_{\ZZ_{(2)}} R)$ become isomorphic by inverting $2$. Then there is an isomorphism of triples 
$$i_1:(L_{(2)}\otimes_{\ZZ_{(2)}} R,\psi_1,\Mb\otimes_{\ZZ_{(2)}} R)\tilde\to (L_{(2)}\otimes_{\ZZ_{(2)}} R,\psi,\Mb\otimes_{\ZZ_{(2)}} R).$$
\indent
{\bf c)} We assume $R$ is an artinian algebra of perfect residue field. Let $I$ be an ideal of $R$. Then any two representations of $\Mb\otimes_{\ZZ_{(2)}} R$ isomorphic modulo $I$, are isomorphic.
\smallskip
{\bf d)} Referring to c), if $\psi_1$ is congruent to $\psi$ modulo $I$, then there is a reductive subgroup of $GSp(L_{(2)}\otimes_{\ZZ_{(2)}} R,\psi_1)$ fixing each element of $\Mb\otimes_{\ZZ_{(2)}} R$ and of the same dimension as $G$.
\smallskip
{\bf e)} We assume all simple factors of $G^{\rm ad}$ are of some $A_n$ or $C_n$ Lie type ($n\in\NN$). Then $\Mb$ is a family of tensors strongly $\ZZ_{(2)}$-very well positioned w.r.t. $\psi$ for $G$ (in the sense of [Va2, 4.3.4]).
\medskip
{\bf Hints and proof of a) in the orthogonal case.} For the a) part use the splitness property in the context of the classification of the possible Cases (A, C or D) of [Ko2, top of p. 375 and p. 395]. The fact that the same classification holds over $R$, can be deduced from the splitness part by using geometric points of ${\rm Spec}(R[{1\over 2}])$. 
\smallskip
Related to the orthogonal case, we can assume we are dealing with an involution of a matrix algebra ${\rm End}(M_R)$ of rank $4n^2$ over $R$ which is defined (in the usual sense; see [KMRT, ch. 1]) by a perfect, symmetric bilinear form (still to be denoted $BF$) on the free $R$-module $M_R$ of rank $2n$. We can assume $R$ is local and $2$-adically complete. $BF$ restricted to $M_R/2M_R$, being alternating, has the same standard diagonal form: its blocks are formed by $2\times 2$ matrices whose entries are $0$ on the diagonal and $1$ otherwise. But if $u$, $v\in M_R/4M_R$ are such that $BF(u,v)=1$, $BF(u,u)=2b$ and $B(v,v)=2a$, with $a$, $b\in R/4R$, for any $x\in R/4R$ we have $BF(u+xv,u+xv)=2(b+x+x^2a)$ (here we still denote by $BF$ its reduction mod $4$). 
\smallskip
We get immediately that by passing to an \'etale extension of $R$ which mod $2$ is defined by at most $n$ equations of the form 
$$x+cx^2+d=0,$$ 
with $c$, $d\in R/2R$, and by completing $2$-adically, we can assume $BF$ restricted to $M_R/4M_R$ has the standard diagonal form. But, as $R$ is $2$-adically complete, by induction we get that we can assume that $BF$ itself has the standard diagonal form. So, as the operation of taking the Zariski closure is well behaved w.r.t. flat extensions, we can assume (by giving up the requirement that $R$ is $2$-adically complete) that $R$ is $\ZZ_{(2)}$ itself, that $M_R=L_n$ and that we are dealing with $BF$ of B2. But $QF(v):={{BF(v,v)}\over 2}$, $v\in L_n$, is exactly the quadratic form we considered in B2 and so the first paragraph of B2 applies.
\smallskip
b) follows from a): using an arbitrary maximal split torus of $\tilde G_R$ and the description of possible Cases in [Ko, p. 395], the situation gets reduced  to an abstract one, in a  context where $\Mb\otimes_{\ZZ_{(2)}} R$ is a matrix algebra and we are dealing with a direct sum of at most two copies of its standard representation. c) is trivial. d) follows from a), by using a natural lifting process. For e), just combine a) with the standard reduction steps of [Va2, 4.3.7 5')].
\medskip
{\bf B4.} We start with a morphism $y:{\rm Spec}(\FF)\to\Mn_{W(\FF)}$. Let
$$(A,p_A):=y^*\bigl((\Ma,\Mp_{\Ma})_{\Mn_{W(k)}}\bigr).$$ 
The problem we face is: the formal deformation space $DS$ of the triple $(A,p_A,\Mb)$ is formally smooth over $W(\FF)$ of relative dimension equal to $\dim_{\CC}(X)$ iff $G^{\rm ad}$ has all its simple factors of some $A_n$ or $C_n$ Lie type; if $G^{\rm ad}$ has factors of some $D_n$ Lie type ($n\ge 4$), then $DS_{\FF}$ has a tangent space of dimension greater then $\dim_{\CC}(X)$ (for instance, if $G^{\rm ad}$ is an absolutely simple $\QQ$--group of $D_n$ Lie type then the tangent space of $DS_{\FF}$ has dimension $n+\dim_{\CC}(X)$; this can be deduced from B2 and from the dimensions of hermitian symmetric domains denoted in [He, p. 518] as D III or as C I). To overcome this problem, let (as in [Va2, 5.1.2]) 
$$z:{\rm Spec}(V)\to\Mn_{W(\FF)}$$
 be a lift of $y$, with $V$ a finite, faithfully flat, DVR extension of $W(\FF)$; it is not a priori a closed embedding. Let 
$$(A_V,p_{A_V}):=z^*((\Ma,\Mp_{\Ma})_{W(\FF)}),$$ 
and let $(t_{\al})_{\al\in\Mj^\prime}$ be the family of de Rham components of the family of Hodge cycles with which $A_V$ is (as in 2.3.3) naturally endowed. Let $\tilde G_V$ be the subgroup of $GL(H^1_{dR}(A_V/V))$ obtained as the Zariski closure of the subgroup of $GL(H^1_{dR}(A_V/V)[{1\over 2}])$ fixing $t_{\al}$, $\forall\al\in\Mj^\prime$. Let $F^1$ be the Hodge filtration of $H^1_{dR}(A_V/V)$ defined by $A_V$.
\medskip
{\bf B5.} From Fontaine's comparison theory with $\QQ_2$ coefficients (used as in 2.3.9 and in [Va2, 5.2.17.2]), and the determinant condition of [Ko2, p. 389-90], we deduce that the extra condition of B3 b) (pertaining to making $2$ invertible) is satisfied for $R=V$ and for the context $CONT$ of the following two triples: $(L_{(2)}\otimes_{\ZZ_{(2)}} V,\psi,\Mb\otimes_{\ZZ_{(2)}} V)$ and $((H^1_{dR}(A_V/V))^*,p_{A_V},\Mb\otimes_{\ZZ_{(2)}} V)$; here we still denote by $p_{A_V}$ the perfect alternating form on $(H^1_{dR}(A_V/V))^*$ defined naturally by $p_{A_V}$ (via Fontaine's comparison theory or via de Rham cohomology). Strictly speaking, in order to appeal to B3, we need to choose an isomorphism $(L_{(2)}\otimes_{\ZZ_{(2)}} V,\Mb\otimes_{\ZZ_{(2)}} V)\tilde\to ((H^1_{dR}(A_V/V))^*,\Mb\otimes_{\ZZ_{(2)}} V)$; in what follows we will not mention it.
\smallskip
 However, in order to apply B3 b) we need to check that, modulo $2$, in the context of the second triple of CONT, involutions of orthogonal first type do not show up.
\medskip
{\bf B6. Exclusion of involutions of orthogonal first type.} We consider a simple factor $\Mb_2$ of $\Mb_1\otimes_{\ZZ_{(2)}} \ZZ_2$ left invariant by $*$ and such that the resulting involution of $\Mb_2[{1\over 2}]$ is of orthogonal first type. As a $\ZZ_2$-algebra, it is a matrix algebra $M_{2n}(W(\FF_{2^m}))$ of a $2n$ dimensional free $W(\FF_{2^m})$-module $V_{2n}$, with $m\in\NN$ and with $n$ as in B2 (the Brauer group of a finite field is trivial). To be consistent with the usual cohomological aspect of this paper, it is more convenient to view $\Mb$ and $\Mb_1$ (resp. $\Mb_2$) as being semisimple subalgebras of ${\rm End}(L_{(2)}^*)$ (resp. of ${\rm End}(L_{(2)}^*\otimes_{\ZZ_{(2)}} \ZZ_2)$). Let $\Gamma$ be the Galois group of $V[{1\over 2}]$ and let 
$$\rho:\Gamma\to GL_{2n}(W(\FF_{2^m}))$$ 
be the Galois representation obtained naturally from the Galois representation on the dual of the Tate-module of $A_{V[{1\over 2}]}$, by ``concentrating" only on the part corresponding to $\Mb_2$ (viewed in the dual context); the identifications of [Va2, top of p. 473] (performed in our present $p=2$ context) tell us that this makes sense.  
\smallskip
$\rho$ is a symmetric representation. Moreover, from $ALT$ we get that $\rho$ mod $2$ is alternating. We consider the $2$-divisible group $D_V$ of $A_V$. The $\ZZ_p$-algebra $M_{2n}(W(\FF_{2^m}))$ splits over $W(\FF_{2^m})$. So $D_V[2]\otimes_{\ZZ_2} W(\FF_{2^m})$ gets naturally decomposed (using a $W(\FF_{2^m})$-endomorphism which is a projector) as a direct sum $E_1\oplus E_2$ in such a way that $E_1[{1\over 2}]$ is naturally associated to $\rho_{W(\FF_{2^m})}$ mod $2$. So $E_1$ comes equipped naturally with a bilinear principal $W(\FF_{2^m})$-quasi-polarization. Here ``bilinear principal" is used as in d) of 2.2.23 A. Looking at $E_1[{1\over 2}]$, as $\rho$ mod $2$ is alternating, we get that the word ``bilinear" can be dropped. If one wants to avoid using tensor products of the form $D_V[2]\otimes_{\ZZ_2} W(\FF_{2^m})$, then one needs to restate everything in terms of symmetric quasi-polarizations (bilinear forms) on $W(\FF_{2^m})^{2n}$ viewed as a $\ZZ_2$-module, which are moreover $W(\FF_{2^m})$-bilinear.
\smallskip
Using the fact that the Dieudonn\'e's functor $\DD$ on $p-FF({\rm Spec}(V/2V))$ is faithfully flat (cf. the complete intersection situation of [BM, 4.3.2]), all these can be transferred into the crystalline cohomology (and so homology) context (for the convenience of the reader, a proof of this --which can be read at any time-- is included in B9 below). In particular, we get that the involution of each simple factor of $\Mb_2\otimes_{\FF_2} V/2V$ (viewed as a factor of the centralizer of $\Mb\otimes_{\ZZ_{(2)}} V/2V$ in ${\rm End}(H^1_{dR}(A_V/V)\otimes_V V/2V))$) is alternating. This ends the argument of the exclusion process.
\medskip
{\bf B7.} So the conditions of B3 a) are satisfied for $CONT$. We deduce: $\tilde G_V$ is a reductive group over $V$; as $V$ is complete having $\FF$ as its residue field, it is split and so we (can) identify it with $G_V$. So, (cf. also the $p=2$ analogue of 2.3.13.2) to show that $(f,L_{(2)},v))$ is a SHS, we just have to show that $\Mn$ is formally smooth over $O_{(v)}$. 
\smallskip
As in [Va2, 5.3.1] we deduce the existence of a cocharacter $\mu_V:\GG_m\to G_V$ giving birth to a direct sum decomposition $H^1_{dR}(A_V/V)=F^1\oplus F^0$, with $\beta\in\GG_m(V)$ acting through $\mu_V$ on $F^i$ as the multiplication with $\beta^{-i}$, $i=\overline{0,1}$. 
\smallskip
For future references, we point out that from B3 b) applied to the two triples of CONT we also get: 
\medskip
f) {\it They are in fact isomorphic.}
\medskip
{\bf B8.} To show that $\Mn$ is formally smooth over $O_{(v)}$ we follow very closely the approach of [Va2, 5.2-5]. We recall that some parts of loc. cit. were dealing just with odd primes. However, as we are dealing essentially just with Hodge cycles which are coming from endomorphisms, we can go around the limitations of loc. cit; here the use of ``essentially" has to do with the fact that the connected component of the origin of the subgroup $G_1$ of $GSp(W,\psi)$ fixing the elements of $\Mb$, is $G$. To explain this ``roundaboutness", we need quite a lot of preliminaries. 
\smallskip
Let $e(V):=[V:W(\FF)]$. Let $\tilde Re(V)$ (resp. $Re(V)$) be the subring of $B(\FF)[[T]]$ (with $T$ an independent variable) which is formed by formal power series $\sum_{n=0}^\infty a_nT^n$, with $a_n\in B(\FF)$ such that $b_n:=a_n[{n\over {e(V)}}]!\in W(\FF)$, $\forall n\in\NN\cup\{0\}$ (resp. such that the sequence $(b_n)_{n\in\NN}$ is formed by elements of $W(\FF)$ and converges to $0$). For $q\in\NN$, let $I(q)$ be the ideal of $\tilde Re(V)$ formed by formal power series with $a_0=a_1=...=a_{q-1}=0$. $Re(V)$ and $\tilde Re(V)$ have natural Frobenius lifts $\Phi_{Re(V)}$ and respectively $\Phi_{\tilde Re(V)}$: they take $T$ to $T^2$.
\smallskip
We follow closely [Va2, 5.2.1]. We choose a uniformizer $\pi_V$ of $V$. Let $f_{e(V)}\in W(\FF)[[T]]$ be the Eisenstein polynomial of degree $e(V)$ having $\pi_V$ as a root of it. So $V=W(\FF)[[T]]/(f_e(V))$. Let $S_{e(V)}$ be the subring of $B(\FF)[[T]]$ generated by $W(\FF)[[T]]$ and by ${{f_{e(V)}^n}\over {n!}}$, $n\in\NN$. As $f_{e(V)}$ is an Eisenstein polynomial, $S_{e(V)}$ can be defined also as the subring of $B(\FF)[[T]]$ generated by $W(\FF)[[T]]$ and by ${{T^{en}}\over {n!}}$, $n\in\NN$. Its $2$-adic completion is $Re(V)$. Warning: we do not have an interpretation for $\tilde Re(V)$, similar to the one of loc. cit. pertaining to odd primes. 
\smallskip
By mapping $T$ to $\pi_V$, as $V$ is $2$-adically complete, we get a $W(\FF)$-epimorphism 
$$q(\pi_V):Re(V)\twoheadrightarrow V.$$ 
We denote by $q_n(\pi_V)$ the $W_n(\FF)$-epimorphism obtained from $q(\pi_V)$ by tensoring with $W_n(\FF)$.  ${\rm Ker}(q_n(\pi_V))$ is an ideal of $Re(V)/2^nRe(V)$ endowed with a natural structure of divided powers (induced from the natural one of ${\rm Ker}(q(\pi_V))$); unfortunately we do not have nilpotent divided powers. Due to this and the fact that $q(\pi_V)$ does not factor through a $W(\FF)$-epimorphism from $\tilde Re(V)$ to $V$, [Va2, 5.2-5] needs substantial ``adjustments". We need the following obvious statement:
\medskip
{\bf Fact.} {\it $\forall n\in\NN$, $Re/2^nRe$ is the inductive limit of its local artinian $W_n(\FF)$-subalgebras taken by $\Phi_{Re(V)}$ mod $2^n$ into themselves (and ordered under the relation of inclusion).}
\medskip
{\bf B9. Preliminary assumption.} We first assume the triple $(A_V,p_{A_V},\Mb)$ lifts to a (similar) triple $(A_{Re(V)},p_{A_{Re(V)}},\Mb)$ over $Re(V)$. 
Let $M:=H^1_{dR}(A_{Re(V)}/Re(V))$, let $\psi_M$ be the perfect form on $M$ defined by $p_{A_{Re(V)}}$ and let $\nabla_M$ be the $2$-adic completion of the Gauss--Manin connection on $M$ defined by $A_{Re(V)}$. Let $\Phi_M$ be the $\Phi_{Re(V)}$-linear endomorphism of $M$ defined by $A_{Re(V)}$. $\Mb\otimes_{\ZZ_{(2)}} Re(V)$ acts naturally on $M$. Let $F^0({\rm End}(M))$ be the maximal direct summand of ${\rm End}(M)$ taking the Hodge filtration of $M$ defined by $A_{Re(V)}$ into itself. 
\smallskip
As ${\rm Ker}(q(\pi_V))$ has a natural divided structure, from the existence of the Dieudonn\'e's functor $\DD$ as a crystal (see [BM]), we get that the quintuple $(M,\Phi_M,\nabla_M,\psi_M,\Mb)$ is well defined without mentioning (i.e. assuming the existence of) $(A_{Re(V)},p_{A_{Re(V)}},\Mb)$. As the elements of $\Mb$, viewed as endomorphisms of $M$ are fixed by $\Phi_M$, we get that the involution of ${\rm End}(M)$ defined by $\psi_M$ leaves invariant $\Mb$, producing the same involution of it as $*$. In other words, it is easy to see that two elements of $F^0({\rm End}(M))$ fixed by $\Phi_M$ and which coincide mod ${\rm Ker}(q(\pi_V))$, are the same; one just needs to use the fact that 
$${\rm Ker}(q(\pi_V))\cap\cap_{s\in\NN} \Phi_{Re(V)}^s(Re(V))=\{0\}.$$ 
\indent
Moreover, the last paragraph of B6 extends automatically to the whole crystalline site ${\rm CRIS}({\rm Spec}(V/p^nV)/{\rm Spec}(\ZZ_2))$, $\forall n\in\NN$. We get that each simple factor of the involution of the reduction mod $2$ of the centralizer of $\Mb$ in ${\rm End}(M)$ induced by $\psi_M$ is either of second type or of alternating first type. We present an argument (in the spirit of this paper) for this well known fact (stated first in B6). We can assume we are in the context of a truncation mod $2$ $E$ of a $2$-divisible group over $V/2V$ which has a principal quasi-polarization $p_E:E\tilde\to E^t$. We have a natural involution 
$$INV:SP\tilde\to SP$$ 
on the space $SP$ of lifts of $E$ to (truncations mod $2$ of $2$-divisible groups) over some fixed artinian $W(\FF)$-algebra $AL_2$ lifting $AL_1:=V/2V$ in such a way that the kernel of the $W(\FF)$-epimorphism $AL_2\twoheadrightarrow AL_1$ is annihilated by the maximal ideal of $AL_2$: it takes such a lift $\tilde E$ into the lift $\tilde E^t$, viewed via $p_E$, as a lift of $E$. The elements of $SP$ fixed by $INV$ are in one-to-one correspondence to pairs $(\tilde E,p_{\tilde E})$, with $\tilde E$ as mentioned before and with $p_{\tilde E}:\tilde E\tilde\to \tilde E^t$ an isomorphism lifting $p_E$. But any $\FF$-linear automorphism of order $2$ of a non-trivial $\FF$-vector space, has non-zero fixed points. Based on this and on Grothendieck's deformation theorem of [Il, 4.4 c) and e)] (applied repeatedly to a projective system $(AL_n)_{n\in\NN}$) we get easily that the pair $(E,p_E)$ lifts to a pair $(E_1,p_{E_1})$ over $W(\FF)[[T]]$, with $E_1$ as the truncation mod $2$ of a $2$-divisible group over $W(\FF)[[T]]$ and with $p_{E_1}$ an isomorphism $E_1\tilde\to E_1^t$. So first a) and then b) of 2.2.23 C apply entirely to $\DD(E_1)$ and to its bilinear principal quasi-polarization defined naturally by $p_{E_1}$. Using the natural $W(k)$-monomorphism $W(\FF)[[T]]\hookrightarrow Re(V)$, we get that the bilinear form on the extension to $Re(V)$ of the underlying module of $\DD(E_1)$ is alternating.
\smallskip
So from B3 a) we deduce that we have a natural reductive subgroup of $GL(M)$; as in [Va2, 5.2.2.1 and 5.2.17] we get that it is isomorphic to $G_{Re(V)}$. 
\smallskip
From the Fact of B8 we get that the level-$N$ symplectic similitude structures of $(A_V,p_{A_V})$, with $(N,2)=1$, lift compatibly to level-$N$ symplectic similitude structures of $(A_{Re(V)},p_{A_{Re(V)}})$. 
\smallskip
Next we apply [Va2, 5.3.3]. As in loc. cit., to $(A_{Re(V)},p_{A_{Re(V)}})$ and its symplectic similitude structures, it corresponds a morphism
$$m_{Re(V)}:{\rm Spec}(Re(V))\to\Mm.$$
The affine transformation 
$$\tilde g_0:Re(V)\otimes_{W(\FF)} V\to Re(V)\otimes_{W(\FF)} V$$
which is $V$-linear and takes $T$ into $\pi_V^{e(V)-1}T+\pi_V$
is still well defined. Warning: this is not so if we replace $Re(V)$ by $\tilde Re(V)$. We consider an $O_{(v)}$-monomorphism $i_V: V\hookrightarrow\CC$. So, as in [Va2, 5.3.3.1] we construct a $W(\FF)$-monomorphism
$$\tilde g:Re(V)\hookrightarrow \CC[[T]],$$
which, when composed with the $\CC$-epimorphism $\CC[[T]]\twoheadrightarrow\CC$ defined by $T$ goes to $0$, becomes $i_V\circ q_{\pi_V}$. Warning: the naive way of defining (here or in loc. cit.) $\tilde g$ directly by $T$ goes to $T+i_V(\pi_V)$ leads into problems of convergence of formal power series with coefficients in $\CC$ and evaluated at $i_V(\pi_V)$.
\smallskip
As $G$ is the connected component of the origin of $G_1$ and as $Re(V)$ is an integral ring of characteristic $0$, $\nabla_M$ respects the $G_{Re(V)}$-action.
We explain what we mean by this. The triple $(M,\psi_M,\Mb)$ is obtained from a triple $(M_0,\psi_0,\Mb)$ over $W(\FF)$ by extensions of scalars. This can be proved in many ways: one way relies on the hint for b) of B3; another way can be deduced easily from the approach of B10 below. However, here we adopt a slightly different approach: denoting by $*|$ the restriction of the involution of ${\rm End}(M)$ to the centralizer of $\Mb$ in ${\rm End}(M)$, the triple $(M,*|,\Mb)$ is obtained from a (similar) triple $(M_0,*|,\Mb)$ over $W(\FF)$ by extensions of scalars. Argument: the only non-trivial (i.e the $D$) case is a consequence of the proof of a) of B3 ($Re(V)$ is strictly henselian and $2$-adically complete, cf. B8). So identifying $M=M_0\otimes_{W(\FF)} Re(V)$, (regardless of the fact that we use $(M_0,\psi_0,\Mb)$ or $(M_0,*|,\Mb)$) $\nabla_M$ is of the form $\dl+\be_{\rm end}$, with $\dl$ as the connection on $M$ annihilating $M_0$ and with (this is what we mean by $\nabla_M$ respects the $G_{Re(V)}$-action) $\be_{\rm end}\in {\rm Lie}(G_{Re(V)})\otimes dT$. 
\smallskip
So, as in [Va2, 5.3.3.1] we get that $m_{Re(V)}$ factors through $\Mn^\prime$. There is one detail which needs to be pointed out, as $G$ is not always the subgroup of $GSp(W,\psi)$ fixing the elements of $\Mb$; so in order to apply loc. cit. we need to add any one of the following two things:
\medskip
{\bf i)} [Va2, 5.2.2.1 and 5.2.15] make sense for $p=2$ as well (being just rational statements) and so we can use all de Rham components of the Hodge cycles with which $A$ is naturally endowed, in order to fully use [Va2, 4.1.5] exactly as mentioned in [Va2, 5.3.3.1];
\smallskip
{\bf ii)} ${\rm Lie}(G)$ is the Lie subalgebra of ${\rm Lie}(GSp(W,\psi))$ annihilating the elements of $\Mb$ and so we can still apply [Va2, 4.1.5] without mentioning the extra Hodge cycles of i), as loc. cit. relies (as [Fa2, end of rm. iii) after th. 10] does) just on facts involving Lie algebras. 
\medskip
The arguments below show that (the $p=2$ analogue of) [Va2, 5.3.4] can be skipped.
\smallskip
Let $(A_0,p_{A_0},\Mb)$ be the triple we get over $W(\FF)$ by pulling back $(A_{Re(V)},p_{A_{Re(V)}},\Mb)$ via the $W(\FF)$-monomorphism $z_0:{\rm Spec}(W(\FF))\hookrightarrow {\rm Spec}(Re(V))$, defined at the level of rings by: $T$ goes to $0$.     
\smallskip  
With its construction we are essentially done, as 
[Fa2, th. 10 and the remarks after] treat as well the case $p=2$. In other words, as in [Va2, 5.4.5] (see also 2.3.11) we construct a versal deformation $(A_R,p_{A_R},\Mb)$ of the triple $(A_0,p_{A_0},\Mb)$ over $R:=W(\FF)[[x_1,...,x_m]]$, with $m:=\dim_{\CC}(X)$. We consider the triple $(A_{\tilde Re(V)},p_{A_{\tilde Re(V)}},\Mb)$ over $\tilde Re(V)$ obtained from the triple $(A_{Re(V)},p_{A_{Re(V)}},\Mb)$ via the $W(\FF)$-morphism 
$$i_{e(V)}:{\rm Spec}(\tilde Re(V))\to {\rm Spec}(Re(V))$$ 
defined by the logical $W(\FF)$-monomorphism $Re(V)\hookrightarrow \tilde Re(V)$. It is easy to see that it is obtained from $(A_R,p_{A_R},\Mb)$ by pull back via a $W(\FF)$-morphism 
$$z_{\tilde Re(V))}:{\rm Spec}(\tilde Re(V))\to {\rm Spec}(R)$$ 
such that $z_{\tilde Re(V))}\circ \tilde z_0$ at the level of rings is defined naturally by: $x_i$'s go to $0$; here we denote by $\tilde z_0$ the natural factorization of $z_0$ through $i_{e(V)}$. To see this, we work in the crystalline cohomology context and mimic entirely [Fa2, th. 10 and rm. iii) after it]. Though $\tilde Re(V)$ is not isomorphic to $W(\FF)[[T]]$, it has all ring properties of $W(\FF)[[T]]$ needed in order to be able to mimic loc. cit.:
\medskip
{\bf iii)} $\tilde Re(V)$ is the projective limit of the quotients $\tilde Re(V)/I(q)$ (the transition $W(\FF)$-epimorphisms being the logical ones);
\smallskip
{\bf iv)} $\forall q\in\NN$, the quotient $I(q)/I(q+1)$, as a $W(\FF)$-module, is free of rank one;
\smallskip
{\bf v)} $\Phi_{\tilde Re(V)}$ takes $I(q)$ into $I(q+1)$, $\forall q\in\NN$. 
\medskip
It is the freeness part of iv) which guarantees that we have no obstructions from the fact that $G_1$ is not necessarily connected: as ${\rm Lie}(G_1)={\rm Lie}(G)$, the part of [Fa2, rm. iii) after th. 10] involving strictness for filtrations applies entirely in the same manner. 
\smallskip
To conclude that $\Mn_{W(\FF)}$ is formally smooth over $W(\FF)$ in $y$, we just have to remark that in the same way we applied [Va2, 4.1.5] in [Va2, 5.4.5], we can apply it in the present situation. In other words, the morphism $m_R:{\rm Spec}(R)\to\Mm$ corresponding to $(A_R,p_{A_R})$ and its natural symplectic similitude structures (lifting those of $(A_0,p_{A_0})$), factors through $\Mn_{W(\FF)}$ and (as in [Va2, 5.4.7-8]) the resulting $W(\FF)$-morphism is formally \'etale. We just need to add: the morphism ${\rm Spec}(W(\FF))\to\Mm$ defined by $(A_0,p_{A_0})$ and its natural symplectic similitude structures factors through $\Mn^\prime$ as $m_{Re(V)}$ does, and so $m_R$ factors through $\Mn^\prime$ and so through $\Mn_{W(\FF)}$. This proves the Theorem modulo the checking that the preliminary assumption of the beginning of this B9 always holds.
\medskip
{\bf B10. Argument for the preliminary assumption.} 
It goes in two main steps: first we work mod $2$ and then we use induction modulo (higher) powers of $2$. 
\medskip
{\bf The mod $2$ part.} First of all we need to show that the triple $\Mt\Mr\Mi_1$ obtained from $(A_V,p_{A_V},\Mb)$ by taken it mod $2$, lifts to a triple $(A_{Re(V)/2Re(V)},p_{A_{Re(V)/2Re(V)}},\Mb)$ over $Re(V)/2Re(V)$. $V/2V$ can be identified with $\FF[[T]]/(T^{e(V)})$ and so $q_1(\pi_V)$ factors through $\tilde Re(V)/2\tilde Re(V)$; we denote by 
$$z_{V/2V}:{\rm Spec}(V/2V)\to {\rm Spec}(\tilde Re(V)/2\tilde Re(V))$$
this factorization, viewed at the level of schemes. We recall from B9, that the quintuple $(M,\Phi_M,\nabla_M,\psi_M,\Mb)$ is well defined (without any assumption) and that we have naturally a reductive subgroup $G_{Re(V)}$ of $GL(M)$. We have:
\medskip
{\bf Lemma.} {\it The cocharacter $\mu_V$ lifts to a cocharacter $\mu_{Re(V)}$ of $G_{Re(V)}$.}
\medskip
The proof of this is very much the same as of [Va2, 5.3.2]. The only difference: we can not be so explicit in writing done some ideals and in fact we have to deal as well with $W(\FF)$-monomorphisms (between $W(\FF)$-algebras) and not just with $W(\FF)$-epimorphisms. So, based on the Fact of B8, we consider a sequence 
$$(A_n)_{n\in\NN},$$ 
with $A_n$ a $W_n(\FF)$-subalgebra of $Re(V)/2^nRe(V)$, such that $\forall n\in\NN$ we have:
\medskip
\item{{\bf P1}} $A_n$ is included in $A_{n+1}$ mod $2^n$;
\smallskip
\item{{\bf P2}} under $q_n(\pi_V)$, $A_n$ is mapped surjectively onto $V/2^nV$;
\smallskip
\item{{\bf P3}} $G_{Re(V)/2^nRe(V)}$ is obtained from a reductive group $G_{A_n}$ over $A_n$ by extension of scalars and the reduction of $G_{A_{n+1}}$ mod $2^n$ is naturally identified with the logical pull back (see P1) of $G_{A_n}$.
\medskip
By induction on $n\in\NN$ we get (as in loc. cit.) a cocharacter $\mu_n$ of $G_{A_n}$ lifting (see P2) the reduction of $\mu_V$ mod $2^n$ and such that $\mu_{n+1}$ mod $2^n$ is obtained by natural pull back (see P3) from $\mu_n$. So to get the Lemma we just need to take $\mu_{Re(V)}$ such that its reduction mod $2^n$ is obtained from $\mu_n$ by natural pull back (see P3), $\forall n\in\NN$.
\medskip
We consider the direct sum decomposition $M=F^1_M\oplus F^0_M$ produced by $\mu_{Re(V)}$ (so $\be\in \GG_m(Re(V))$ acts through $\mu_{Re(V)}$ on $F^i(M)$ by multiplication with $\be^{-i}$). $F^1_M$ lifts the direct summand $F^1$ of $H^1_{dR}(A_V/V)=M/{\rm Ker}(q(\pi_V))$. 
\smallskip
So, without having $(A_0,p_{A_0},\Mb)$ we can construct its filtered $\sg$-crystal endowed with a family of endomorphisms by just pulling back $(M,F^1_M,\Phi_M,\Mb)$ via $z_0$. So, as in [Va2, 5.4.5] (see also 2.2.21 UP) we can construct (without having $(A_R,p_{A_R},\Mb)$) a versal, principally quasi-polarized filtered $F$-crystal ${\got C}_R$ over $R/2R$ endowed with a family $\Mb$ of endomorphisms; here $R$ is as in B9. From [dJ1, th. of intro.] (the form of it we need is reproved in 3.14 B6 below) and from Serre--Tate's deformation theory we deduce the existence of a unique triple $(A_{R/2R},p_{A_{R/2R}},\Mb)$ lifting the special fibre $(A_V,p_{A_V},\Mb)_{\FF}$ of $(A_V,p_{A_V},\Mb)$ and whose $F$-crystal over $R$ is obtained from ${\got C}_R$ and its family of endomorphisms by forgetting the filtration of the underlying $R$-module of ${\got C}_R$. 
\smallskip
The extension of $(M,F^1,\Phi_M,\nabla_M,\Mb)$ to $\tilde Re(V)$ is obtained from $({\got C}_R,\Mb)$ by pull back via a $W(\FF)$-morphism $z_{\tilde Re(V)}:{\rm Spec}(\tilde Re(V))\to {\rm Spec}(R)$: this is nothing else but the corresponding part of B9, performed just in terms of principally quasi-polarized filtered $F$-crystals endowed with endomorphisms and not in terms of principally polarized abelian schemes with endomorphisms.  
\smallskip
As $V/2V$ is a complete intersection, again based on [BM, 4.3.2], we get that the pull back of $(A_{R/2R},p_{A_{R/2R}},\Mb)$ via the composite of $z_{V/2V}$ with $z_{\tilde Re(V)}$ mod $2$, is nothing else but $\Mt\Mr\Mi_1$. But from the Fact of B8, we get that the resulting $\FF$-morphism $R/2R\to V/2V$ lifts to an $\FF$-morphism $R/2R\to Re(V)/2Re(V)$. This ends the mod $2$ part.  
\medskip
{\bf The inductive part.} The second step is to lift, by induction on $n\in\NN$, the triple $\Mt\Mr\Mi(1)$ over $Re(V)/2Re(V)$ lifting $\Mt\Mr\Mi_1$ to a triple $\Mt\Mr\Mi(n+1)$ over $Re(V)/2^{n+1}Re(V)$ lifting naturally (i.e. in a compatible way) the reduction $\Mt\Mr\Mi_{n+1}$ mod $2^{n+1}$ of the triple $(A_V,p_{A_V},\Mb)$. We can assume $F^1_M$ mod $2$ is the filtration of $M/2M$ defined by the lift of $A_{V/2V}$ to $Re(V)/2Re(V)$ we got (we just have to redo the Lemma once more). Accordingly, we want $\Mt\Mr\Mi(n+1)$ in such a way that the filtration of $M/2^{n+1}M$ it defines is nothing else but $F^1_M$ mod $2^{n+1}$.
\smallskip
As the ideal $4Re(V)$ of $Re(V)$ has a natural divided power structure, nilpotent modulo any ideal $2^mRe(V)$ with $m\in\NN$, $m\ge 3$, we get that we can assume $n=1$. We consider two ideals of $Re(V)/4Re(V)$: $I_1$ is the ideal generated by $2$, while $I_2$ is ${\rm Ker}(q_2(\pi_V))$. Their intersection is nothing else but $2I_2$ and so its induced divided power structure is nilpotent. There is a unique triple 
$$(A_{(Re(V)/4Re(V))/2I_2},p_{A_{(Re(V)/4Re(V))/2I_2}},\Mb)$$ 
lifting $\Mt\Mr\Mi(1)$ and $\Mt\Mr\Mi_{2}$ at the same time. So, due to the just mentioned nilpotent part, there is a unique way of lifting $(A_{(Re(V)/4Re(V))/2I_2},p_{A_{(Re(V)/4Re(V))/2I_2}},\Mb)$ to $Re(V)/4Re(V)$ in the way prescribed by $F^1$ mod $4$. This takes care of $n=1$ and so ends the inductive part.  
\medskip
{\bf Conclusion.} Due to the mentioned compatibility the triple $(A_V,p_{A_V},\Mb)$ lifts to a triple over ${\rm Spf}(Re(V))$; as we are in a principally polarized context, we can replace (cf. Grothendieck's algebraization theorem) ``${\rm Spf}$" by ``${\rm Spec}$". So indeed the triple $(A_V,p_{A_V},\Mb)$ lifts to a triple $(A_{Re(V)},p_{A_{Re(V)}},\Mb)$ over ${\rm Spec}(Re(V))$. This ends the argument for the preliminary assumption and so the proof of the Theorem.
\medskip
{\bf B11. Remarks.} {\bf 1)} It is worth pointing out that the above proof does not use Fontaine's comparison theory (if in B9 we go ahead with ii) and not with i) and if we are not bothered to get f) of B7; warning: B5 can be worked out as well in the de Rham context). To be compared with the proof of [Va2, 5.1] which relies on such a theory.
\smallskip
{\bf 2)} B3 a) points out that the theory of [Va2, 4.3] of $O$-well positioned families of tensors (with $O$ a DVR) can be ``extended": 
\medskip
a) we use tensors which are not necessarily linear (like the function $QF(v)$ of the proof of it: we have $QF(\al v)=\al^2 v$, $\forall\al\in M_R$);
\smallskip
b) we impose restrictions on tensors modulo the maximal ideal of $O$.
\medskip
However, such theories are not easily adaptable to (i.e. usable in) the crystalline context of [Fa2, \S 4]. 
\smallskip
{\bf 3)} We refer to [Va2, 4.1] for notations. B9 ii) points out that [Va2, 4.1.3-5] have versions, where we replace the set $\Mj$ of loc. cit., with a subset $\Mj_0$ of it such that the subgroup of $GSp(W,\psi)$ fixing $s_{\al}$, $\forall\al\in\Mj_0$, has $G$ as its connected component of the origin. For these versions, we just need to replace in [Va2, 4.1.3] ``becomes a quadruple of $\Ma(G,X,W,\psi)$" by ``can be extended (by adding tensors indexed by $\Mj\setminus\Mj_0$) to a quadruple of $\Ma(G,X,W,\psi)$".
\smallskip
{\bf 4)}  $\tilde Re(V)$ is the simplest example of the class of $GD$-rings (here $GD$ stands for good deformation) to be axiomatized in \S 6 (starting from iii) to v) of B9).
\smallskip
{\bf 5)} Using [Va2, 5.6.4] we get immediately (cf. also B4) that $\Mn=\Mn^\prime$, provided $G^{\rm ad}$ does not have factors of some $D_n$ Lie type ($n\ge 4$). However, loc. cit. does not apply to the case of such factors.
\medskip
{\bf C.} 2.3.8 3) and conventions 2.3.7 and 2.3.9.2 apply but not 2.3.8 1). However, if we replace  $2(p-2)$ by ${\rm max}\{2,2(p-2)\}$ in 2.3.8 1) (in order to accommodate the case $p=2$), then the $p=2$ analogue of (the modified) 2.3.8 1) is the Theorem of B. 2.3.6 makes no sense for $p=2$. It is premature to state a morally supportive variant of 2.3.8.1; however, the proof of [Va2, 6.7.2] can be used with the same purpose for $p=2$ as well. Similarly, it is premature to deal with the $p=2$ analogue of 2.3.8 2).
\smallskip
{\bf D.} The parts a) and c) of 2.3.9 remain valid without any modification as (cf. A) 2.3.5.1 holds for $p=2$. Moreover, 2.3.9 d) applies as well provided (*) of 2.3.9 D holds: the part of the proof of 2.3.9 referring to [Ja, 10.4 of Part I] applies as well. Briefly, this goes as follows. The difference is that, with the notations of the proof of 2.3.9, $\tilde G_{0k}$ is the subgroup of $GL(L_0/pL_0)$ fixing a non-degenerate quadratic form $\tilde Q_0$ in $m$ variables (argument: this is so over $\bar k$, cf. [Bo2, 23.6] and its lift version to $W(\bar k)$; but such a quadratic form is unique modulo an invertible element and so, using Hilbert's Theorem 90, we can assume it is defined over $k$). Now after a finite (computable) number $n(m)$ of Galois extensions (the first one of $k$) of degree $2$ we can bring $\tilde Q_0$ in the standard form $x_0^1+x_1x_2+...+x_{2n-1}x_{2n}$ (resp. $x_1x_2+...+x_{2n-1}x_{2n}$) if $m=2n+1$ is odd (resp. if $m=2n$ is even). We postpone to refer to any $p=2$ equivalent of 2.3.9 b). 
\smallskip
{\bf E.} 2.3.10-16 remain valid: this is so due to the fact that [Fa2, th. 10] and the remarks following it are still true for $p=2$; so [Va2, 5.4] and its variants (see 2.3.11) remain true for $p=2$. 
But in connection to 2.3.17 we have to be more careful: for $p=2$ the second paragraph of the proof of 2.3.17 has to be modified slightly, as a $2$-divisible group over $W(k)$ is not determined by its associated filtered $\sg$-crystal; however it applies to 2.3.17 in the same way (see 2.3.18.1 B and C and 2.4.1 below), provided:
\medskip
-- for the uniqueness part of 2.3.17 we restrict to the case when $(M,\vph)$ does not have slope $0$;
\smallskip
-- for the existence part of the whole of 2.3.17 (inclusive 2.3.17.2) we restrict to the case when $k$ has no abelian extensions of degree $2$ or to the case when $(M,\vph)$ does not have slope $0$ or does not have slope $1$. 
\medskip
{\bf 2.3.18.1. A review.} Here we recall some known facts on $2$-divisible groups $D$ over a perfect field $k$ of characteristic $2$, for which we could not find a good reference. 
\smallskip
{\bf A.} We write $D$ as a product 
$$D=D_1\times D_2$$ 
such that $D_1$ is ordinary, while $D_2$ has all its slopes in the interval $(0,1)$. Let $(M,\vph)$ be its $\sg$-crystal and let $M=M_1\oplus M_2$ be the corresponding direct sum decomposition. For simplifying the presentation, we choose a lift of $D$ to a $2$-divisible group $D(2)$ over $W_2(k)$. Let $n\in\NN\cup\{0\}$ be such that $k^n$ is the $k$-linear space parameterizing lifts of $D$ to $W_2(k)$, with $D(2)$ corresponding to the origin. Let $T_{\rm def}$ be the vector group scheme it defines naturally. Let $k^n$ be the $k$-linear space which parameterizes lifts of the $F^1$-filtration of $M/2M$ (defined as the kernel of $\vph$ mod $2$) to $F^1$-filtrations of $M/4M$, with the $F^1$-filtration of $M/4M$ defined by $D(2)$ as the origin; let $T_{\rm fil}$ be the vector group scheme it defines naturally. 
\smallskip
In what follows we assume $n>0$: if $n=0$ all results below are trivial. There is a natural morphism (of schemes):
$$m_D:T_{\rm def}\to T_{\rm fil}.$$ 
Its existence can be seen easily using crystalline cohomology theory (i.e. using the $\DD$ functor in the context of the universal deformation space of $D$ over $W(k)[[x_1,...,x_n]]$; see 2.2.1.0 and [Il, 4.8]); see also 2.4 below for a second approach via PD-hulls and crystals.
We would like to point out that the role of $D(2)$ is irrelevant here: above, instead of $k$-vector spaces, we can consider as well affine spaces (in the classical sense of Euclidean Geometry); a change of $D(2)$ corresponds to a change of the origins of $T_{\rm def}$ and of $T_{\rm fil}$. We have:
\medskip
{\bf B. Fact.} {\it $m_D$ is a dominant, finite morphism. Generically, the field extension of the field of fractions $FF$ of $T_{\rm fil}$ we get, is the composite of at most $n$ extensions of degree $2$ of $FF$.}
\medskip
{\bf Proof:} We first mention the case when $D$ is an ordinary $2$-divisible group over $k$ and its lift $D(2)$ to $W_2(k)$ is modeled on (i.e. lifts to) its canonical lift. But in such a case this Fact can be checked without difficulty, starting from the theory of canonical crystalline coordinates. In other words, [Ka3, 4.3.1] and the standard formulas of pulling back $F$-crystals (see [De3, 1.1]; for mod $4$ versions of loc. cit. see also [Fa1, p. 36-37] and the reference to it in 2.2.1 c)) reduces the situation to the case when $D$ is the $2$-divisible group of an ordinary elliptic curve. For this last case we refer to [Og, 3.14] (from loc. cit. we get: $m_D$ can be identified with a morphism ${\rm Spec}(k[t])\to {\rm Spec}(k[t])$ which at the level of rings takes $t$ into $t^2-t$). We conclude: $m_D$ is a Galois cover of whose Galois group is $(\ZZ/2\ZZ)^n$.
\smallskip
To treat the general case without appealing to heavy computations, we use a specialization argument. So let $V(k)$ be the normalization of $k[[x]]$ in the algebraic closure of $k((x))$. Let $D_{V(k)}$ be a $2$-divisible group over $V(k)$ lifting $D$ and whose generic fibre is ordinary (for instance, cf. 3.1.8.1 below; it can be read at any time). Let $D_{V(k)}(2)$ be a lift of it to $W_2(V(k))$ which also lifts $D(2)$ (cf. [Il, 4.4]). Let $M_{V(k)}:=H^1_{\rm crys}(D_{V(k)}/W(V(k)))$. As in A, we denote by $T^{V(k)}_{\rm def}$ (resp. by $T^{V(k)}_{\rm fil}$) the $n$ dimensional affine $V(k)$-scheme whose $V(k)$-valued points parameterize lifts of $D_{V(k)}$ to $W_2(V(k))$ (resp. lifts. of the $F^1$-filtration of $M_{V(k)}/2M_{V(k)}$ defined naturally by $D_{V(k)}$ to an $F^1$-filtration of $M_{V(k)}/4M_{V(k)}$), with $D_{V(k)}(2)$ (resp. with the $F^1$-filtration of $M_{V(k)}/4M_{V(k)}$ defined by $D_{V(k)}(2)$) corresponding to its origin. We get (easy argument at the level of $V(k)$-valued points): $m_D$ is the special fibre of a natural morphism
$$m_{D_{V(k)}}:T^{V(k)}_{\rm def}\to T_{\rm fil}^{V(k)}$$
(whose generic fibre is constructed as $m_D$ is).  
\smallskip
It is easy to see that $m_D$ has finite fibres (see also E below). So $m_D$ is a dominant morphism. From Zariski's Main Theorem, we get that $T_{\rm def}$ is an open subscheme of the normalization of $T_{\rm fil}$ in the field of fractions of $T_{\rm def}$. But $\AA^n_k$ can not be an open subscheme of an integral, affine $k$-scheme $S$ such that $S\setminus\AA_k^n$ is non-empty: otherwise we would get that the polynomial ring $k[x_1,...,x_n]$ has invertible elements which do not belong to $k$. So $m_D$ is also a finite morphism. As the generic fibre of $m_{D_{V(k)}}$ is a Galois cover of whose Galois group is $(\ZZ/2\ZZ)^n$, the Fact follows.
\medskip
{\bf C.} From B we get that the natural map
$$m_D(\bar k):T_{\rm def}(\bar k)\to T_{\rm fil}(\bar k)$$
defined by $m_D$ is surjective. 
\medskip
In general, the fibres of $m_D(\bar k)$ are related to the passage to some particular isogenies (involving the $2$-torsion). So using (see E below) Fontaine's comparison theory with integral coefficients in the context of $2$-divisible groups over $W(k)$, such a fibre $m_D(\bar k)^{-1}(\alpha)$ can have more than $1$ element precisely in the case when $M_1$ has both slopes $0$ and $1$ with positive multiplicity. So, if $(M_1,\vph)$ has only the slope $0$ or only the slope $1$ or if $M_1=\{0\}$, then $m_D(\bar k)$ is a bijection. 
\smallskip
 {\bf D.} On the other hand, to lift a $2$-divisible group $D(m)$ over $W_m(k)$ lifting $D$, with $m\in\NN$, $m\ge 2$, to a $2$-divisible group over $W_{m+q}(k)$, with $q\in\NN$, we just have to lift the $F^1$-filtration of $M/2^mM$ defined by $D(m)$ to an $F^1$-filtration of $M/2^{n+m}M$. This is a consequence of [Me, ch. 4-5], as the ideal $2^mW_{m+q}(k)$ of $W_{m+q}(k)$ is equipped naturally with a nilpotent divided power structure (which is not the case for the ideal $2W_{m+q}(k)$ of $W_{m+q}(k)$). Combining this with the Fact of B, we get:
\medskip
{\bf Corollary.} {\it For any direct summand $F^1$ of $M$ such that the triple ${\got C}:=(M,F^1,\vph)$ is a filtered $\sg$-crystal, there is an abelian extension $k_1$ of $k$ of Galois group a subgroup of $\ZZ/2\ZZ^n$, such that the extension of ${\got C}$ to $k_1$ is associated to a $2$-divisible group over $k_1$.}
\medskip
{\bf E.} We now recall how Fontaine's comparison theory with integral coefficients in the context of $2$-divisible groups over $W(k)$ applies in C. We follow [Fa2, ch. 4 and ch. 6]. Warning: in [Fa2, ch. 1] it is stated that in [Fa2, ch. 2-6] an odd prime is used. However, in [Fa2, ch. 8] (especially see the ending paragraph of loc. cit.) it is explained why the things are fine for $p=2$ as well; so below we do refer to [Fa2, ch. 4 and ch. 6] even for $p=2$. Let $B^+(W(k))$ be the Fontaine's ring constructed as in [Fa2, ch. 4 and ch. 8], using the $2$-adic topology (for completing); the reason we do not use the PD topology (for completing) is that we could not convince ourselves that we still get  a natural Frobenius lift (see also [Fa2, p. 125] where the arguments refer only to odd primes). We recall: $B^+(W(k))$ is an integral, local $W(k)$-algebra, endowed with a separated, decreasing filtration $\bigl(F^i(B^+(W(k))\bigr)_{i\in\NN\cup\{0\}}$, with a Frobenius lift $FR$, and with a natural Galois action by $\Gamma_k$ (cf. also [FI, 1.3]). Moreover we have a natural $W(k)$-epimorphism compatible with the natural Galois actions by $\Gamma_k$
$$s_{W(k)}: B^+(W(k))\twoheadrightarrow \overline{W(k)}^\wedge,$$
where $\overline{W(k)}^\wedge$ is the $2$-adic completion of the normalization of $W(k)$ in $\overline{B(k)}$. Its kernel is $F^1(B^+(W(k)))$ and so we can identify 
$$gr^0:=F^0(B^+(W(k)))/F^1(B^+(W(k))=B^+(W(k))/F^1(B^+(W(k))$$ 
with $\overline{W(k)}^\wedge.$
So 
$$gr^1:=F^1(B^+(W(k))/F^2(B^+(W(k))$$ 
is naturally a $\overline{W(k)}^\wedge$-module: it is free of rank one; let $u$ be a generator of it. Let $\be_0\in F^1(B^+(W(k)))$ be defined as in [Fa2, p. 125] (in [FI, 1.2.3 and 1.3.3] it is denoted by $t$); its image in $gr^1$ is $2u$ times an invertible element $i_{\be_0}$ of $\overline{W(k)}^\wedge$ and, even more, 
$${{\be_0}\over 2}\in F^1(B^+(W(k)))$$
 and 
$$FR(\be_0)=2\be_0$$
 (cf. any one of the last two loc. cit.). 
\smallskip
Let $D_0$ be a $2$-divisible group over $W(k)$ lifting $D$. Let $H^1_{\acute et}(D_0)$ be the dual of its Tate-module. Let $F^1$ be the $F^1$-filtration of $M$ it defines naturally. Let $F^0$ be a direct supplement of $F^1$ in $M$. We have a $B^+(W(k))$-monomorphism
$$i_M(D_0):M\otimes_{W(k)} B^+(W(k))\hookrightarrow H^1_{\acute et}(D_0)\otimes_{\ZZ_2} B^+(W(k))$$
respecting the tensor product filtrations (the filtration of $H^1_{\acute et}(D_0)$ is defined logically by: $F^1(H^1_{\acute et}(D_0))=\{0\}$ and $F^0(H^1_{\acute et}(D_0))=H^1_{\acute et}(D_0)$), the Galois actions, and the Frobenius; all these Frobenius endomorphisms of modules are still denoted by $FR$. 
\smallskip
Its existence is just a particular case of [Fa2, th. 7] (cf. end of [Fa2, ch. 8]). From the way loc. cit is stated, clearly the prime $2$ is also treated. However, we need to point out that the ring $R_V$ used in loc. cit., as $p=2$, is not definable as in [Fa2, ch. 2] (for instance, if $k=\FF$, then the ring $R_V$ of loc. cit. is not the ring $\tilde Re(V)$ used in 2.3.18 B8 but the ring $Re(V)$ of the mentioned place). Moreover, our case is very simple: $W(k)$ has no ramification and so, as $B^+(W(k))$ is a $W(k)$-algebra, we can trivially define the crystalline cohomology group of $D_0$ relative to $B^+(W(k))$, using (cf. end of [Fa2, ch. 8]) $CRIS({\rm Spec}(k)/{\rm Spec}(W(k)))$. It can be checked that, working with the mentioned site, the homology and the cohomology of $\QQ_2/\ZZ_2$ can still be computed similarly to [Fa2, p. 132-3] (see [BM, ch. 2]).
\smallskip
The part of [Fa2, th. 7] expressing the fact that the cokernel of $i_M(D_0)$ is annihilated by the element $\beta_0$, in its simpler form, can be restated as (cf. also the strictness part of [Fa2, th. 5] and the above part pertaining to $i_{\be_0}$): 
\medskip
{\bf Key Property.} {\it The $\overline{W(k)}^\wedge$-linear map 
$$j_M(D_0):F^0\otimes_{W(k)} gr^1\oplus F^1\otimes_{W(k)} gr^0\to H^1_{\acute et}(D_0)\otimes_{\ZZ_2} gr^1,$$ induced naturally from $i_M(D_0)$ at the level of one gradings, is injective and its cokernel is annihilated by $2$.} 
\medskip
If now $D_0^1$ is another $2$-divisible group over $W(k)$ lifting $D$, producing the same $F^1$-filtration $F^1$ of $M$, then via the similarly constructed $B^+(W(k))$-monomorphism $i_M(D_0^1)$, $H^1_{\acute et}(D_1^0)$ is (cf. the Key Property) a $\ZZ_2$-lattice of 
$H^1_{\acute et}(D_0)[{1\over 2}]$ containing $2H^1_{\acute et}(D_0)$ and contained in ${1\over 2}H^1_{\acute et}(D_0)$. Using a well known theorem of Tate, we get that $D_0^1$ is determined by such a $\ZZ_2$-lattice, and this takes care (cf. also paragraph D) of the finiteness part of the fibres of $m_D$. If $H^1_{\acute et}(D_0)\neq H^1_{\acute et}(D_0^1)$, i.e. if $D_0\neq D_0^1$ as lifts of $D$, let 
$$H^1:=H^1_{\acute et}(D_0)\cap H^1_{\acute et}(D_0^1).$$ 
Let $H^1_0$ be a direct summand of $H^1$ contained in $2H^1_{\acute et}(E)$, with $E\in\{D_0,D_0^1\}$: we do have non-zero such summands; warning: we do not require $H^1_0$ to be $\Gamma_k$-invariant. So, via $j_M(E)$, $H^1_0\otimes_{\ZZ_p} gr^1$ is a direct summand of $F^0\otimes_{W(k)} gr^1\oplus F^1\otimes_{W(k)} gr^0$. We deduce the existence of an element 
$$a\in FF^1:=F^0\otimes_{W(k)} F^1(B^+(W(k))\oplus F^1\otimes_{W(k)} B^+(W(k))$$
such that the following three things hold:
\medskip
a) $FR(a)=2a$;
\smallskip
b) its image in $F^0\otimes_{W(k)} gr^1\oplus F^1\otimes_{W(k)} gr^0$ is non-zero modulo the maximal ideal $m_{\rm big}$ of $\overline{W(k)}^\wedge$;
\smallskip
c) its image in $H^1\otimes_{\ZZ_p} gr^1$ belongs to $H^1_0\otimes_{\ZZ_p} gr^1$ and moreover its image in $H^1_0\otimes_{\ZZ_p} gr^1/m_{\rm big}gr^1$ is non-zero.
\medskip
For instance, we can take $a:=i_M(E)({b\over 2}\be_0)$, where $b\in H_0^1$ generates a direct summand of it and moreover ${b\over 2}\in H^1_{\acute et}(E)$; as ${b\over 2}\be_0=b{{\be_0}\over 2}$, we get that a) to c) hold.
\smallskip 
We assume now that $(M,\vph)$ does not have integral slopes, i.e. we assume $M_1=\{0\}$. We want to reach a contradiction. We can assume $k=\bar k$. We choose a $W(k)$-basis $\Mb_1=\{e_1,...,e_{\dim_{W(k)}(M)}\}$ (resp. $\Mb_2=\{f_1,...,f_{\dim_{W(k)}(M)}\}$) of an arbitrary $W(k)$-submodule $N$ of $M$ (resp. of $M$) such that $\vph^l(e_i)=2^{n_i}e_{i}$, $\forall i\in S(1,\dim_{W(k)}(M))$, with $l\in\NN$ and with all $n_i$'s belonging to the set $S(1,l-1)$ (resp. such that it is formed by elements of $F^0$ and of $F^1$). Let $q\in\NN$ be such that $p^qM\subset N$. We write 
$$a=\sum_{i=1}^{\dim_{W(k)}(M)} \al_if_i,$$ 
with $\al_i$ belonging to $F^1(B^+(W(k))$ if $f_i\in F^0$ and to $B^+(W(k))$ if $f_i\in F^1$. Let $m:=l(q+1)$. We have $(l-1)(q+1)+q<m$. As $\{\vph^m(f_i)|i\in S(1,\dim_{W(k)}(M))\}$ is a $B(k)$-basis of $M[{1\over p}]$ formed by elements of $N\setminus p^{(q+1)(l-1)+q}M$, and as $FR^m(a)=2^ma\in 2^mM\otimes_{W(k)} B^+(W(k))$, we get that $\forall i\in S(1,\dim_{W(k)}(M))$ such that $f_i\in F^1$, $\al_i$ is not an invertible element. So the image of $a$ in $FF^1/m_{\rm big}FF^1$ belongs to the image of $F^0\otimes_{\ZZ_p} gr^1$ into this last $\overline{W(k)}^\wedge/m_{\rm big}$-module (i.e. $k$-vector space). 
\smallskip
From this, c) and the structure of $B^+(W(k))/F^2(B^+(W(k))$ (as an extension of $gr^0$ by $gr^1$), we deduce the existence of an element $b^0\in F^0$ such that $i_M(E)(b^0)\in H^1\otimes_{\ZZ_p} B^+(W(k))$ is non-zero modulo the maximal ideal of $B^+(W(k))$. This implies $(M,\vph)$ has slopes $0$. Contradiction. So $M_1$ is non-zero. 
\smallskip
We assume now that $(M,\vph)$ does not have slope 1. So $(M_1,\vph)$ has pure slope $0$. Let $H_0$ be the direct summand of $H^1_{\acute et}(D_0)$ and of $H^1_{\acute et}(D_0^1)$ such that $H_0[{1\over 2}]$ corresponds (via Fontaine's comparison theory) to the filtered isocrystal $(M_1[{1\over 2}],0,\vph_1)$. It is also a direct summand of $H^1$. So we can choose $H_0^1$ such that $H^1_0\oplus H_0$ is still a direct summand of $H^1$. Based on the same arguments as above, we can assume $b^0\in M_2$. So $(M,\vph)$ has the slope $0$ with multiplicity greater than the rank of $M_1$. Contradiction. So $(M,\vph)$ has the slope $1$ with positive multiplicity. Using the standard Cartier duality, we get that $(M,\vph)$ has as well the slope $0$ with positive multiplicity. This takes care of the bijectiveness part of C.
\medskip
{\bf 2.3.18.1.1. The fibres of $m_D$.} We assume now that $D_0$ is the product of the canonical lift of $D_1$ with an arbitrary lift $D^2$ of $D_2$. The same arguments of 2.3.18.1 E can be used to give that any lift $D_1^0$ of $D$ as in the mentioned place  is a product of a lift of $D_1$ with $D^2$ (we can also use the fact that the Galois representation associated to $D^2[2]$ does not fix any non-zero element and that $D^1_0$ is the extension of a $2$-divisible group by $D^2$).
In general we have:
\medskip
{\bf Fact.} {\it All fibres of $m_D(\bar k)$ (of 2.3.18.1 B) have exactly $2^{s(-1)}$ points, where $s(-1)$ is the multiplicity of the slope $-1$ for $({\rm End}(M),\vph)$.}
\medskip
{\bf Proof:} Let $D_0$ and $F^1$ be as in 2.3.18.1 E. $D_0$ is the extension of an \'etale $2$-divisible group $D_{01}$ by another $2$-divisible group $D_{02}$ over $W(k)$ whose special fibre does not have slope $0$. $D_{02}$ is uniquely determined by $F^1$, cf. 2.3.18.1 C. $D_{02}$ itself is the extension of a $2$-divisible group $D_{03}$ which does not have integral slopes by a $2$-divisible group $D_{04}$ of multiplicative type. We have a short exact sequence of $k$-vector spaces
$$0\to {\rm Ext}^1({D_{01}}_{W_2(k)},{D_{04}}_{W_2(k)})\hookrightarrow {\rm Ext}^1({D_{01}}_{W_2(k)},{D_{02}}_{W_2(k)})\twoheadrightarrow {\rm Ext}^1({D_{01}}_{W_2(k)},{D_{03}}_{W_2(k)})\to 0\leqno (EXT)$$
(the surjectivity part is checked by reasons of dimension). ${D_0}_{W_2(k)}$ is given by an element $\gamma\in {\rm Ext}^1({D_{01}}_{W_2(k)},{D_{02}}_{W_2(k)})$; its image in ${\rm Ext}^1({D_{01}}_{W_2(k)},{D_{02}}_{W_2(k)})$ is determined by $F^1/4F^1$ (again cf. 2.3.18.1 C). So (cf. also 2.3.18.1 D) the number of lifts $D_0$ of $D$ giving birth to $F^1$ are in one-to-one correspondence to the number of elements in ${\rm Ext}^1({D_{01}}_{W_2(k)},{D_{04}}_{W_2(k)})$ giving birth to the same $F^1$-filtration of $H^1_{\rm crys}(D_{04k}\times D_{01k}/W_2(k))$. So, as $s(-1)$ is the product of the ranks of $D_{01}$ and $D_{04}$, the Fact follows from the first paragraph of the proof of 2.3.18.1 B.
\medskip
{\bf 2.3.18.1.2. Exercise.} Using 2.3.18.1.1 show that the morphism $m_D$ of 2.3.18.1 B is a Galois cover. Hint: the situation gets reduced to the ordinary case.
\medskip
{\bf 2.3.18.2. The case $p=2$ for some parts of 2.2.} We assume $k$ is an algebraically closed field of characteristic $2$. First we refer to 2.2.16.2. Any cyclic diagonalizable filtered $\sg$-crystal $(M,F^1,\vph)$ is a direct sum of circular diagonalizable filtered $\sg$-crystals. From 2.3.18.1 C and E we deduce that if the Newton polygon of $(M,\vph)$ does not have slope $0$ or $1$, then we have attached to it a uniquely determined Galois representation with coefficients in $\ZZ_2$; so in such a context the Corollary of 2.2.16.2 still holds, as the part of its proof involving endomorphisms still holds. 
\smallskip
Accordingly, as we can always isolate aside the slopes 0 and 1 (as in 2.3.18.1 A), 2.2.16.3 remains true for $p=2$ as well. Similarly, we get that the Proposition of D of 2.2.22 3) remains true for $p=2$. This implies that the whole of D, F, G, H, I, K and M of 2.2.22 3) remains true for $p=2$.
\medskip
{\bf 2.3.18.3. Remark.} Let now $p$ be an arbitrary prime. 2.3.18.1 E can be adapted to $p$-divisible objects of $\Mm\Mf_{[a,p-1+a]}(W(k))$ as follows. Performing a Tate twist we can assume $a=0$. [Fa1, p. 35-36] explains how to associate Galois representations to truncations of such $p$-divisible objects; using projective limits we can associate a Galois representation $\rho$ to a $p$-divisible object ${\got C}$ of $\Mm\Mf_{[0,p-1]}(W(k))$. Of course, one needs to see when it involves a free $\ZZ_p$-module $N$ of the same rank as the underlying module of ${\got C}$. We will not stop to argue this (it is a consequence of [CF]); we will just mention that in the case ${\got C}$ is cyclic diagonalizable, it is easy to see based on loc. cit. that indeed $N$ is free and has the right rank. Moreover, following the pattern of 2.3.18.1 E it can be checked that, if ${\got C}$ either does not have direct summands which are $p$-divisible objects of $\Mm\Mf_{[0,0]}(W(k))$ or does not have direct summands which are $p$-divisible objects of $\Mm\Mf_{[p-1,p-1]}(W(k))$, then ${\got C}$ is uniquely determined by its associated Galois representation. So referring to iv) of L of 2.2.22 3), we can replace $p\ge m+2$ by $p\ge m+1$. 
\medskip\smallskip
{\bf 2.4. Some identifications.} What follows is thought as an appendix of \S 2: we detail [Va2, 5.5.1] (slightly reformulated) in the abstract context of Shimura $\sg$-crystals; moreover the case $p=2$ is included as well. Based on the end of 2.2.10 and of 2.2.21, there is no need to restate all things below in a principally quasi-polarized context. So, based on Fact 4 of 2.3.11 (and its $p=2$ version), everything below applies to the context of a ($p=2$) SHS. 
\smallskip
We start with a Shimura $p$-divisible group $(D,(t_{\al})_{\al\in\Mj})$ over $W(k)$, with $k$ an arbitrary perfect field. Let $(M,F^1,\vph,G,(t_{\al})_{\al\in\Mj})$ be its attached non-necessarily quasi-split Shimura filtered $\sg$-crystal. Let $R:=W(k)[[x_1,...,x_m]]$, with $m:=dd((M,\vph,G))$. Let $\Md_R:=(D_R,(t_{\al})_{\al\in\Mj})$ be a universal Shimura $p$-divisible group over ${\rm Spec}(R)$ which in the $W(k)$-valued point of ${\rm Spec}(R)$ defined by making all $x_i$'s to be $0$, is $(D,(t_{\al})_{\al\in\Mj})$. Let $P$ (resp. $P_M$) be the parabolic subgroup of $G$ (resp. of $GL(M)$) normalizing $F^1$. Let ${\rm Spec}(R_1)$ be the completion of $G/P$ in its $W(k)$-valued point defined by the origin of $G$. As a $W(k)$-algebra, $R_1$ is isomorphic to $R$ (cf. 2.2.22 4)).
\smallskip
Let $R^{PD}$ (resp. $R_1^{PD}$) be the $p$-adic completion of the PD-hull of the maximal ideal of $R$ (resp. of $R_1$). 
The pull back to ${\rm Spec}(R^{PD}/pR^{PD})$ of the $F$-crystal with tensors over $R/pR$ defined by $\Md_R$ is trivial, i.e. it can be identified with $(M\otimes_{W(k)} R^{PD},\vph\otimes 1,(t_{\al})_{\al\in\Mj})$ (see [BM]). So, corresponding to the $F^1$-filtration of $M\otimes_{W(k)} R^{PD}$ defined by the pull back of $\Md_R$ to ${\rm Spec}(R^{PD})$, we get naturally a $W(k)$-morphism
$$FIL:{\rm Spec}(R^{PD})\to GL(M)/P_M.$$
It factors through $G/P$ as it can be checked using $W(k)$-valued points of ${\rm Spec}(R^{PD})$ (see Fact 2 of 2.2.9 3)).
As $R^{PD}/p^nR^{PD}$ is an inductive limit of its artinian local $W_n(k)$-subalgebras (this is the same as the Fact of 2.3.18 B), $\forall n\in\NN$, we get that $FIL$ factors through ${\rm Spec}(R_1)$ and so through $R_1^{PD}$. The resulting $W(k)$-morphism 
$$FIL_0:{\rm Spec}(R^{PD})\to {\rm Spec}(R_1^{PD}),$$
is an isomorphism if $p\ge 3$ or if $p=2$ and each $W(k)$-valued point of ${\rm Spec}(R_1^{PD})$ lifts uniquely to a $W(k)$-valued point of ${\rm Spec}(R^{PD})$, as it can be checked immediately at the level of tangent spaces, starting from the fact that $\Md_R$ is universal. By the identification of [Va2, 5.5.1] we meant such $W(k)$-isomorphisms as $FIL_0$ (with $p\ge 3$). 
\smallskip
Any $W(k)$-automorphism $g_R:{\rm Spec}(R)\tilde\to {\rm Spec}(R)$ extends to a $W(k)$-automorphism $g_R^{PD}$ of ${\rm Spec}(R^{PD})$. We assume now that there is $g\in G^{\rm ad}$ which mod $p$ normalizes $F^0({\rm Lie}(G^{\rm ad}))/pF^0({\rm Lie}(G^{\rm ad}))$ and such that:
\medskip
-- the $F^0$-filtration $F^0_{g_R^*(z)}$ of ${\rm Lie}(G^{\rm ad})$ defined by the pull back through $g_R$ of any ${\rm Spec}(W(k))$-section $z$ of ${\rm Spec}(R)$ whose attached Shimura adjoint filtered Lie $\sg$-crystal is defined by an $F^0$-filtration $F^0_z$ of ${\rm Lie}(G^{\rm ad})$, is nothing else but $g^{-1}(F^0_z)$.
\medskip
Denoting still by $g$, the $W(k)$-automorphism of ${\rm Spec}(R_1)$ (resp. of ${\rm Spec}(R_1^{PD})$) defined by inner conjugation by $g$ (resp. defined by inner conjugation by $g$ and passage to the $p$-adic completion of $PD$-hulls), we get:
\medskip
{\bf Fact.} {\it $g\circ FIL_0=FIL_0\circ g_R^{PD}$.}
\medskip 
{\bf Proof:} It is enough to check this equality at the level of $W(k)$-valued points of ${\rm Spec}(R^{PD})$; but for such points, the equality holds, cf. the assumption.
\medskip
{\bf Corollary.} {\it We assume that either $p\ge 3$ or $p=2$ and each $W(k)$-valued point of ${\rm Spec}(R_1^{PD})$ lifts uniquely to a $W(k)$-valued point of ${\rm Spec}(R^{PD})$. Let $V$ be a faithfully flat $W(k)$-algebra which is a complete DVR of index of ramification at most $p-1$. If $g_R$ acts freely on ${\rm Spec}(R[{1\over p}])$, then $g$ acts freely on the set of $V$-valued points of ${\rm Spec}(R_1[{1\over p}])$.}
\medskip
{\bf Proof:} We just have to add that any $V$-valued point of ${\rm Spec}(R)$ (resp. of ${\rm Spec}(R_1)$) lifts uniquely to a $V$-valued point of ${\rm Spec}(R^{PD})$ (resp. of ${\rm Spec}(R_1^{PD})$).
\medskip
{\bf 2.4.1. The case $p=2$.} We assume now $p=2$. Let $T_{\rm def}^G$ (resp. $T_{\rm fil}^G$) be the tangent space of $R/pR$ (resp. of $R_1/pR_1$), viewed as a $dd((M,\vph,G))$ dimensional affine space over $k$. Its $k$-valued points are in one-to-one correspondence to lifts of the maximal point ${\rm Spec}(k)\to {\rm Spec}(R)$ to a $W_2(k)$-valued point of ${\rm Spec}(R)$ (resp. of lifts of $F^1/2F^1$ to a $F^1$-filtration $F^1_4$ of $M/4M$ defined as usual by a cocharacter of $G_{W_2(k)}$), with the lift provided by all $x_i$'s being $0$ (resp. by $F^1/4F^1$) corresponding to the origin. $FIL_0$ gives birth naturally to a $k$-morphism
$$m_D^G:T_{\rm def}^G\to T_{\rm fil}^G.$$
From 2.3.18.1 B we get that $m_D^G$ is dominant and finite, generically obtained by using a composite of extensions of degree $2$. 
\smallskip
Let ${\got C}_R$ be the Shimura filtered $F$-crystal with tensors over ${\rm Spec}(R/pR)$ associated to $\Md_R$. We have:
\medskip
{\bf Claim.} {\it If $k=\bar k$, then the induced map at the level of $k$-valued points $m_D^G(k):T_{\rm def}^G(k)\to T_{\rm fil}^G(k)$ is surjective and all its fibres have the same number of elements which is a non-negative, integral power of $2$.}
\medskip
{\bf Proof:} This is a consequence of 2.3.18.1.1.2. We now include a more useful second proof, based on 2.2.21 UP. It is more convenient to work with $W(k)$-valued points instead of $W_2(k)$-valued points (cf. 2.3.18.1 D). 
\smallskip
So let $z_i:{\rm Spec}(W(k))\to {\rm Spec}(R)$, $i=\overline{1,2}$, be two different $W(k)$-morphisms giving birth to the same filtration $F^1_0:=F^1_{z_1}=F^1_{z_2}$ of $M$ (if such $W(k)$-morphisms do not exist we have nothing to prove). From 2.2.21 UP we deduce the existence of a unique $W(k)$-isomorphism
$z_{12}:{\rm Spec}(R)\tilde\to {\rm Spec}(R)$ such that:
\medskip
i) $z_{12}\circ z_2=z_1$, and
\smallskip
ii) $z_{12}^*({\got C}_R)$ is isomorphic to ${\got C}_R$, under an isomorphism $i_{12}$ lifting (in $z_2$) the canonical identification $(M_0,F^1_{z_2},\vph,G,(t_{\al})_{\al\in\Mj})=(M_0,F^1_{z_1},\vph,G,(t_{\al})_{\al\in\Mj})$. 
\medskip
As ${\got C}_R$ has no automorphism whose pull back through $z_2$ is the identity, the isomorphism $i_{12}$ is canonical. So from 2.2.21 UP and from i) and ii) we get:
\medskip
{\bf a)} $z_{12}$ does not fix any $W(k)$-valued point of ${\rm Spec}(R)$;
\smallskip
{\bf b)} For any $z\in {\rm Spec}(R)(W(k))$, the $F^1$-filtrations of $M$ defined by $z$ and $z_{12}\circ z$ are the same;
\smallskip
{\bf c)} If we have a third point $z_3$ (not necessarily different from $z_1$ or $z_2$) such that $F^1_{z_3}=F^1_0$, then $z_{12}\circ z_{23}=z_{13}$ (the convention being that $z_{13}$ or $z_{23}$ is the identity automorphism of ${\rm Spec}(R)$ if respectively $z_1=z_3$ or $z_2=z_3$).
\medskip
In other words, if $\al\in T_{\rm fil}^G(k)$ is defined by $F^1_{z_1}/4F^1_{z_1}$, we have an injective map
(thought as a monodromy map)
$$m_{z_1}:m_D^G(k)^{-1}(\al)\hookrightarrow {\rm Aut}_{W(k)}({\rm Spec}(R)),$$
which takes the element of $m_D(k)^{-1}(\al)$ defining $z_3$ into $z_{13}$. 
The image of $m_{z_1}$ does not depend on $z_1$ (cf. a) and b)). So all fibres of $m_D^G(k)$ have the same number of elements. So the Claim follows from the sentence on extensions of degree 2 before it.
\medskip
{\bf 2.4.2. Corollary.} {\it $z_{12}$ is an automorphism of order $2$.}
\medskip
{\bf Proof:} The fact that $z_{12}$ is of finite order is implied by the finiteness of the fibres of $m_D^G(k)$ and by its uniqueness. From the fact that the special fibre of $D_R$ is a versal deformation, [dJ1, th. of intro.] (the form of it we need is reproved in 3.14 B6 below) implies $z_{12}$ mod $2$ is the trivial automorphism of ${\rm Spec}(R/2R)$. So $z_{12}^2$ mod $4$ is the trivial automorphism of ${\rm Spec}(R/4R)$ and so (being of finite order) is the trivial automorphism of ${\rm Spec}(R)$. 
\vfill\eject
\centerline{}
\bigskip
\bigskip
\centerline{\bigsll {\bf \S 3 The basic results}}
\bigskip\bigskip
\medskip
Let $k$ be a perfect field of characteristic a prime $p\ge 2$. Let $\sg$ be the Frobenius automorphism of $W(k)$. Let $(M,\vph,G)$ be a Shimura $\sg$-crystal
over $k$ (cf. 2.2.8 2)). From its very definition, we deduce (see 2.2.8 1) and 2)) the existence of a cocharacter $\mu:\GG_m\to G$
which produces a direct sum decomposition $M=F^1\oplus F^0$ (with $\be\in\GG_m(W(k))$ acting through $\mu$ on $F^i$ as the multiplication with $\be^{-i}$, $i=\overline{0,1}$) such that $(M,F^1,\vph,G)$ is a Shimura filtered $\sg$-crystal. We also deduce that $F^0$ and $F^1$ are proper direct summands of $M$. The $G(W(k))$-conjugacy class $[\mu]$ of the cocharacter $\mu$ is uniquely determined by $(M,\vph,G)$ (cf. Fact 2) of 2.2.9 3)). We consider a family
$(t_\al)_{\al\in\Mj}$ of tensors of $\Mt(M)[{1\over p}]$ such that $\vph(t_{\al})=t_{\al}$, $\forall\al\in\Mj$, and $G_{B(k)}$ is the subgroup of $GL\bigl(M[{1\over p}]\bigr)$ fixing $t_\al$, $\forall\al\in\Mj$. If convenient, one can assume that all tensors $t_{\al}$ are homogeneous. Let 
$${\got g}:={\rm Lie}(G).$$
\indent
Fixing the group $G$, the family of tensors $(t_\al)_{\al\in\Mj}$ and the $G(W(k))$-conjugacy class $[\mu]$, any other Shimura $\sg$-crystal
$(M,\vph_1,G)$, with $\vph_1$ fixing $t_\al$, $\forall\al\in\Mj$, and with $[\mu]$ defining its filtration class, up to isomorphism (cf. 2.2.9 6)), can be put in the form $\vph_1=g\vph$, with $g\in G(W(k))$. To see this, we consider two Shimura filtered $\sg$-crystals $(M,F^1,\vph,G,(t_{\al})_{\al\in\Mj})$ and $(M,F^1_1,\vph_1,G,(t_{\al})_{\al\in\Mj})$ having the same filtration class. Changing the second one by an isomorphism, we can assume $F^1=F^1_1$. But then $\vph_1=g\vph$ with $g\in GL(M)(W(k))$. As $\vph_1$ and $\vph$ are both fixing $t_\al$, $\forall\al\in\Mj$, we get that actually
$g\in G(W(k))$. Let
$$d_M:=\dim_{W(k)}(M).$$
\medskip
{\bf 3.0. The starting setting.} We consider Shimura $\sg$-crystals $(M,\vph_1,G)$ and Shimura filtered $\sg$-crystals $(M,F^1,\vph_1,G)$, with $\vph_1$ fixing $t_\al,\forall \al\in\Mj$, and with $[\mu]$ defining their filtration class. In other words, we consider Shimura $\sg$-crystals with an emphasized family of tensors whose class is $Cl(M,\vph,G,(t_{\al})_{\al\in\Mj})$ (cf. 2.2.22 2)) as well as their lifts. Occasionally, for the sake of accuracy, we repeat this. We use freely 2.2.3 3). 
\smallskip
In 3.1-13 we assume $p\ge 3$; see 3.14 for the modifications needed to be made for $p=2$. In 3.1-5, 3.7-9 and 3.11 we mostly deal with the generalization of the classical Serre--Tate (ordinary) theory to the context of $Cl(M,\vph,G,(t_{\al})_{\al\in\Mj})$; in particular, we also deal with problems involving Newton polygons. In 3.6 (resp. 3.12) we mostly deal with global (resp. local) deformations in the context of $Cl(M,\vph,G,(t_{\al})_{\al\in\Mj})$. In 3.10 we introduce a language and formulas pertaining to Shimura (adjoint) Lie $\sg$-crystals. In 3.13 we mostly deal with truncations mod $p$ in the context of $Cl(M,\vph,G,(t_{\al})_{\al\in\Mj})$. Often, in 3.1-14 we deal as well with more general contexts then the one of $Cl(M,\vph,G,(t_{\al})_{\al\in\Mj})$. In 3.15 we gather different refinements which can be obtained by combining different parts of 3.1-14.
\smallskip
Warning: from now on, without a special reference, we assume $p\ge 3$.
\medskip\smallskip
{\bf 3.1. The first group of basic results.}
\medskip
{\bf 3.1.0. Theorem. a)} {\it Among all  $\sg$-crystals $(M,g\vph)$ with $g\in G(W(k))$ (i.e. among all Shimura $\sg$-crystals $(M,\vph_1,G,(t_\al)_{\al\in\Mj})$ whose filtration class is defined by $[\mu]$), the ones which have the smallest Newton polygon $\Mp$ (in the sense that all others have a Newton polygon strictly above it) are precisely the ones $(M,g_0\vph)$ for which there is a Shimura
filtered $\sg$-crystal $(M,F^1_0,g_0\vph,G)$ such that ${\got p}_0:=W_0({\got g},g_0\vph)$ is contained in
the parabolic Lie subalgebra $F^0_0({\got g})$ of {\got g} formed by elements which take $F^1_0$ into itself.
\smallskip
{\bf b)} For any Shimura $\sg$-crystal $(M,g_0\vph,G)$ having $\Mp$ as its Newton polygon, there is a unique filtration $F^1_0$ of $M$, such that $(M,F^1_0,g_0\vph,G)$ is a Shimura filtered $\sg$-crystal satisfying (with the notations of a)) ${\got p}_0\subset F^0_0({\got g})$.
\smallskip
{\bf c)} Among all $\sg$-crystals $(M,g\vph)$ with $g\in G(W(k))$ (i.e. among all Shimura $\sg$-crystals $(M,\vph_1,G,(t_\al)_{\al\in\Mj})$ whose filtration class is defined by $[\mu]$), the ones which have the smallest Newton polygon 
$Lie_G(\Mp)$ (in the sense that all others have a Newton polygon strictly above it) of their attached
Shimura Lie $\sg$-crystals $({\got g},g\vph)$, are precisely the ones having $\Mp$ as their Newton polygon.
\smallskip
{\bf d)} There is a Shimura filtered $F$-crystal ${\got C}_F=(M,F^1,\vph,G,\tilde f)$ such that the Shimura $F$-crystal we get (cf. Fact 3 of 2.2.10) over a geometric point over the generic point of the scheme over which this ${\got C}_F$ is, has $\Mp$ as its Newton polygon.}
\medskip
{\bf 3.1.1. Definitions.} {\bf a)} Any Shimura $\sg$-crystal $\bigl(M,\vph_1,G,(t_\al)_{\al\in\Mj}\bigr)$ of whose filtration class is defined by $[\mu]$ and has $\Mp$ as its Newton polygon, is called a $G$-ordinary $\sg$-crystal (or a 
Shimura-ordinary $\sg$-crystal or, even better, a $(G,[\mu],[a])$-ordinary $\sg$-crystal, as $\Mp$ depends not only on the faithful representation $G\hookrightarrow GL(M)$ but also on the filtration class defined by $[\mu]$ and on the automorphism class [a] of $({\got g},\vph_1)$ as defined in 2.2.11 2)). 
\smallskip
{\bf b)} The formal isogeny type of a $G$-ordinary $\sg$-crystal $(M,g_0\vph,G)$ (resp. such a $G$-ordinary $\sg$-crystal), is called the $G$-ordinary type or the Shimura-ordinary type (resp. is called a $G$-ordinary $\sg$-crystal) produced by (or attached to) the Shimura
$\sg$-crystal $(M,\vph,G)$ we started with. 
\smallskip
{\bf c)} A Shimura filtered $\sg$-crystal
$(M,F^1_1,\vph_1,G)$ such that $W_0({\got g},\vph_1)$ is contained in the parabolic Lie subalgebra of {\got g} formed by elements taking $F^1_1$ into itself, is called the $G$-canonical lift (or the Shimura-canonical lift) of (the $(G,[\mu],[a])$-ordinary $\sg$-crystal) $(M,\vph_1,G)$. Occasionally, we refer to $F^1_1$ itself as the $G$-canonical lift of $(M,\vph,G)$. Pairs of the form $(\DD^{-1}(M,F^1_1,\vph_1),(t_{\al})_{\al\in\Mj})$ are called Shimura-canonical $p$-divisible groups over $W(k)$, while pairs of the form $(\DD^{-1}(M,F^1_1,\vph_1)_{k},(t_{\al})_{\al\in\Mj})$ are called Shimura-ordinary $p$-divisible groups (over $k$).
\medskip
{\bf 3.1.1.1. Exercise.} If $G=GL(M)$, then the $G$-ordinary type we get is an ordinary type (i.e. its associated Newton polygon has only the slopes 0 and 1).
Hint: this should be obvious; if it is not, then look at 3.1.8.1 or at 3.2.3 below.
\medskip
From 3.1.0 a) and b) and from defs. 2.2.8 7) and 2.2.9 6) we get:  
\medskip
{\bf 3.1.1.2. Corollary.} {\it Any endomorphism (resp. automorphism or $1_{\Mj}$-automorphism) of a Shimura-ordinary $\sg$-crystal over $k$, is as well an endomorphism (resp. automorphism or $1_{\Mj}$-automorphism) of its canonical lift.}
\medskip
{\bf 3.1.2. Remarks.} {\bf 1)} The terminology $G$-ordinary $\sg$-crystals and $G$-canonical lifts was suggested to us by F. Oort and R. Pink in June 1995. It is a convenient terminology for the abstract context; in geometric situations we will use gradually more and more the terminology Shimura-ordinary $\sg$-crystals and Shimura-canonical lifts. 
\smallskip
{\bf 2)} The class of Shimura-ordinary $\sg$-crystals is stable under perfect field extensions (cf. 3.1.0 b)). So, we speak about a Shimura-ordinary (resp. a Shimura-canonical) $p$-divisible group over $k$ even if by chance in 3.1.1 c) we do not have a quasi-split reductive group $G$.
\medskip
{\bf 3.1.2.1. Convention.} Whenever we pass to a perfect field $k_1$ containing $k$, for not overloading the notations, we still use the terminology $G$-ordinary, instead of $G_{W(k_1)}$-ordinary ($\sg_{k_1}$-crystals).
\medskip
{\bf 3.1.3. Reformulation of 3.1.0 in terms of sets.} 3.1.0 asserts the following five things:
\medskip
{\bf a)} Among all Newton polygons of $\sg$-crystals $(M,g\vph)$ with 
$g\in G(W(k))$, there is one $\Mp$ which is the smallest; let $A_G$ be the set of $g\in G(W(k))$ such that the Newton polygon of $(M,g\vph)$ is $\Mp$.
\smallskip
{\bf b)} Among all Newton polygons of Shimura Lie $\sg$-crystals $({\got g},g\vph)$ with $g\in G(W(k))$, there is one ${\rm Lie}_G(\Mp)$ which is the smallest; let $B_G$ be the set of $g\in G(W(k))$ such that the Newton polygon of $({\got g},g\vph)$ is ${\rm Lie}_G(\Mp)$.
\smallskip
{\bf c)} $A_G=B_G=C_G$, where $C_G$ is the set of $g\in G(W(k))$ such that $(M,g\vph)$ has the property that there is a Shimura filtered $\sg$-crystal $(M,F^1_g,g\vph,G)$ for which $W_0({\got g},g\vph)$ is contained in the parabolic Lie subalgebra of {\got g} formed by elements taking $F^1_g$ into itself (i.e. having the property that there is a Shimura filtered $\sg$-crystal $(M,F^1_g,g\vph,G)$ whose attached Shimura filtered Lie $\sg$-crystal $\bigl({\got g},g\vph,F^0_g({\got g}),F^1_g({\got g})\bigr)$ is of parabolic type in the sense of def. 2.2.12 a)).
\smallskip
{\bf d)} If $g\in C_G$, then there is a unique lift $F^1_g$ of $(M,g\vph,G)$ as in c).
\smallskip
{\bf e)} Any Shimura $\sg$-crystal can be deformed (using specific Shimura filtered $F$-crystals) so that it is the specialization of a $G$-ordinary $\sg_{k_1}$-crystal over a (suitable) perfect field $k_1$ containing $k$.
\medskip
{\bf 3.1.4. Proposition.} {\it For any $G$-ordinary $\sg$-crystal $(M,\tilde \vph,G)$, there is a unique cocharacter $\tilde\mu:\GG_m\to G$ producing a direct sum decomposition $M=\tilde F^1\oplus \tilde F^0$ (with $\be\in\GG_m(W(k))$ acting through $\tilde\mu$ on $\tilde F^i$ as the multiplication with $\be^{-i}$, $i=\overline{0,1}$) such that $(M,\tilde F^1,\tilde \vph,G)$ is a $G$-canonical lift and all elements of $W_0({\got g},\tilde\vph)$ (resp. of $W^0({\got g},\tilde\vph)$) take $\tilde F^1$ (resp. $\tilde F^0$) into itself.}
\medskip
For a first proof of this Proposition see 3.11.1-2 below. A second proof of it can be obtained by just copying the proof of a) of 4.4.1 3) below (where we work in the context of Shimura varieties of Hodge type) and so, for making \S 4 as little dependent on \S 3 as possible, it is reproduced only there. 
\medskip
{\bf 3.1.5. Exercise.} The cocharacter $\tilde\mu$ of 3.1.4, when viewed as a cocharacter of $GL(M)$, is the canonical split of $(M,\tilde F^1,\tilde\vph)$.
\smallskip
{\bf Solution.} We can assume $k=\bar k$. From 3.11.1 a) or c) below we deduce that we can assume $G$ is a torus (cf. also Fact of 2.2.22 1)). So Exercise follows from the functoriality of canonical splits (to be compared with 2.2.1.2) and from Corollary of 2.2.9 3). 
\medskip
A second solution can be obtained using the proof of b) of 4.4.1 3) below. Warning: 3.1.4-5 are not referred before 4.5, so no vicious circle is created (the terminology of 3.1.6 below is used in 3.11.2 just to introduce some extra terminology in 3.11.3-4 and so in 4.3.8.3). Also 3.1.5 can be deduced from 3.1.4, using just the functorial aspect of canonical split cocharacters (cf. 3.1.1.2 and the Claim of 2.2.3 3)). 
\smallskip
The uniqueness part of 3.1.4 implies:
\medskip
{\bf 3.1.5.1. Corollary.} {\it $\tilde\mu$ is fixed by any automorphism of $(M,\tilde\vph,G)$.} 
\medskip
{\bf 3.1.6. Definition.} The cocharacter $\tilde\mu$ of 3.1.4 is called the canonical split (cocharacter) of $(M,\tilde\vph,G)$ (or of $(M,\tilde F^1,\tilde\vph,G)$, cf. [Wi]).
\medskip
{\bf 3.1.7. Remark.} The whole of 3.1.0-6 (as well as the greatest part --i.e. the part which does not mention $p$-divisible groups or Shimura varieties-- of this chapter \S 3) remains valid in the context of generalized Shimura $p$-divisible objects over $k$. See 3.15.6 and Appendix for reformulation of some of the results in this extended context; one exception: 3.4-5 below are organized in such a way that this generalized context is handled simultaneously with the context of $(M,\vph,G)$. The extension of 3.2.1.1, 3.2.2, 3.2.3, 3.2.8, 3.3.1-3 below to generalized Shimura $p$-divisible objects over $k$ is trivial. But this is not quite so for the deformation theory of 3.6: occasionally the deformation theory in the generalized Shimura context requires some additional arguments (see 3.6.1.6, 3.6.18.7.1 c), 3.6.18.7.3 C and 3.15.6 below). 
\medskip
From 3.1.0 d) and 2.2.9 9) we get (cf. 3.1.1.1):
\medskip
{\bf 3.1.8. Corollary.} {\it Any $p$-divisible group $D$ over a perfect field $\tilde k$ of characteristic $p>2$ is the specialization of an ordinary $p$-divisible group.}
\medskip
{\bf 3.1.8.1. Comment.} 3.1.8, for the case when $D$ has a quasi-polarization over a field extension of $\tilde k$ was first obtained in [NO]. For the case when the $a$-number of $D$ is at most $1$ (resp. when the Newton polygon of $D$ has only one slope) see [Oo2] (resp. see [dJO, 5.15]). The last three loc. cit. handle the case $p=2$ as well. The general case, including the case $p=2$, though considered to be well known, it is untraceable in the literature. It can be easily deduced from [Fa2, th. 10] through a simple algebraization process (see 3.6.18.4.1 below for another form of it), based on an elementary analysis of Hasse--Witt invariants. This goes as follows. 
\smallskip
We can assume we are in a context involving Shimura $\sg_{\tilde k}$-crystals, cf. 2.2.9 9). So we use the previous notations of the beginning paragraph of \S 3, with $G=GL(M)$ and $\tilde k=k$, but with $p\ge 2$. If $p=2$ we need to choose $F^1$ such that $(M,F^1,\vph)$ is associated to a $2$-divisible group over $W(k)$. Writing $GL(M)={\rm Spec}(R_{GL})={\rm Spec}(W(k)[x_{11},...,x_{d_Md_M}][{1\over {\rm DET}}]$, with ${\rm DET}$ as the determinant of the matrix of whose entries are the variables $x_{ij}$'s, we choose a Frobenius lift of $R_{GL}^\wedge$ taking $x_{ij}$ into $x_{ij}^p$ if $i\neq j$, and taking $x_{ii}-1$ into $(x_{ii}-1)^p$, $\forall i,j\in S(1,d_M)$.
\smallskip
We work with the Frobenius endomorphism of $M\otimes_{W(k)} R_{GL}^\wedge$ defined by ${\rm gl}(\vph\otimes 1)$, with ${\rm gl}$ as the universal element of $GL(M)(R_{GL}^\wedge)$. Let $g\in GL(M)(W(k))$ be such that $g\vph(F^1)=pF^1$. So the multiplicity of the slope $0$ (resp. $1$) for $(M,g\vph)$ is $\dim_{W(k)}(M/F^1)$ (resp. is $\dim_{W(k)}(F^1)$). But ${\rm gl}$ mod $p$ specializes to $g$ mod $p$. So the Hasse--Witt invariant of ${\rm gl}(\vph\otimes 1)$ over the generic point of the completion of ${\rm Spec}(R_{GL})$ in its origin is as well $\dim_{W(k)}(M/F^1)$. So, based on the construction of 2.2.10 (involving a Shimura filtered $F$-crystal $(M,F^1,\vph,GL(M),\tilde f)$), the conclusion follows from 2.2.21 UP.
\smallskip
The same proof works to show that any principally quasi-polarized $p$-divisible group over $\tilde k$ is the specialization of a principally quasi-polarized, ordinary $p$-divisible group.
\medskip\smallskip
{\bf 3.2. The outline of the proof of 3.1.0.}
As the proof of 3.1.0 involves a couple of pages, we start outlying its main features.
\medskip
{\bf 3.2.1. Lemma.} {\it If $(N,\vph_i)$, $i=\overline{1,2}$, are two $\sg$-crystals over $k$ such that $\forall m\in\NN$, $\vph_1^m=h_m\vph^m_2$, with $h_m\in GL(N)(W(k))$, then they have the same Newton polygon. Moreover, $\forall\al\in [0,\infty)$ we have $W_{\al}(N,\vph_1)=W_{\al}(N,\vph_2)$.}
\medskip
{\bf Proof:} The condition $\vph^m_1=h_m\vph_2^m$, with $h_m\in GL(N)(W(k))$, implies that the Hodge
numbers of the $\sg^m$-crystals $(N,\vph^m_1)$ and $(N,\vph^m_2)$ are the same. But the
Newton polygon of a $\sg$-crystal is determined by the Hodge numbers of its iterates
[Ka2, 1.4.4]. This takes care of the first part.
\smallskip
To see the second part we can assume $k=\bar k$. Let $r\in\NN$ be such that $\vph_i^r$ acts diagonally w.r.t. a $B(k)$-basis $\{e_1^i,...,e_n^i\}$ of $N\otimes_{W(k)} B(k)$ formed by elements of $N\setminus pN$, $i\in\{1,2\}$; here $n:=\dim_{W(k)}(N)$. Let $i_N\in\NN\cup\{0\}$ be such that, $\forall i\in\{1,2\}$, $p^{i_N}N$ is contained in $<e_1^i,...,e_n^i>$. We assume that if $l,s\in S(1,n)$, with $l<s$, then $\vph_i^r(e_l^i)=p^{n_l^i}e_l^i$ and $\vph_i^r(e_s^i)=p^{n_s^i}e_s^i$, with $n_l^i$ and $n_s^i$ non-negative integers such that $n_l^i\ge n_s^i$. Let $\Ms\Ml$ be the finite set of slopes of $(N,\vph_1)$ or of $(N,\vph_2)$. Let $\al\in\Ms\Ml$. Let $W_i(\al):=W_{\al}(N,\vph_i)$. We use mathematical induction on the decreasing values of elements of $\Ms\Ml$ to show that $W_1(\al)=W_2(\al)$. 
\smallskip
We assume that $\forall\beta\in\Ms\Ml$, with $\beta>\al$, we have $W_1(\be)=W_2(\be)$. We need to show $W_1(\al)=W_2(\al)$. Let $t\in S(0,n)$ be such that $\{e_1^i,...,e_t^i\}$ is a $B(k)$-basis of $W_i(\be_0)$, with $\be_0$ as the smallest element of $\Ms\Ml\cap (\al,\infty)$ ($t=0$ iff this last intersection is the empty set). Let $s\in S(1,n)$ be such that $\{e_1^i,...,e_{s}^i\}$ is a $B(k)$-basis of $W_i(\al)$, $i=\overline{1,2}$; we have $s>t$. For $q\in\{t+1,...,s\}$ we write 
$$e_q^1=\sum_{l=1}^n c_{l,q}e_l^2,$$ 
with $c_{l,q}\in B(k)$. We apply the hypothesis with $m$ running through all positive multiples of $r$; so let $m=rm_0$, with $m_0\in\NN$. We get
$$
p^{m_0n_q^1}h_m^{-1}(e_q^1)=\sum_{l=1}^n \sg^{m}(c_{l,q})p^{m_0n_l^2}e_l^2.\leqno (HYP)
$$ 
So $p^{i_M}\sg(c_{l,q})p^{m_0(n_l^2-n_q^1)}\in W(k)$, $\forall l\in S(1,n)$. Taking $m_0\in\NN$ big enough (we need $m_0-i_M$ to be greater than the greatest $p$-adic valuation of all these coefficients $c_{l,q}$), we get: if $c_{l,q}\neq 0$, then $n_l^2\ge n_q^1$. But 
$${n_{u}^2\over r}<\al={n_{t+1}^1\over r}={n_{s}^1\over r},$$ 
if $u\in S(s+1,n)$. So $c_{l,q}=0$ for $l\in S(s+1,n)$. We deduce $e_q^1\in W_2(\al)$. We conclude $W_1(\al)\subset W_2(\al)$. By reasons of dimensions we get $W_1(\al)=W_1(\al)$. 
This proves the Lemma.
\medskip
{\bf 3.2.1.1. Corollary.} {\it If $(N,\vph_i,G_N)$, $i=\overline{1,2}$, are two Shimura $\sg$-crystals such that $\forall m\in\NN$, $\vph_1^m=h_m\vph_2^m$, with $h_m\in G_N(W(k))$, then, $\forall\al\in\QQ$, $W_{\al}({\rm Lie}(G_N),\vph_i)$, does not depend on $i\in\{1,2\}$; so the Newton polygon of $({\rm Lie}(G_N),\vph_i)$, does not depend on $i\in\{1,2\}$.}
\medskip
{\bf Proof:}
We just need to apply 3.2.1 to the $\sg$-crystals $({\rm Lie}(G_N),p\vph_i)$, $i\in\{1,2\}$.      
\medskip
{\bf 3.2.2. The trivial situation.} We come back to the Shimura filtered $\sg$-crystal $(M,F^1,\vph,G)$ considered in the beginning of \S 3. For the meaning of $F^0({\got g})$ and $F^1({\got g})$ see 2.2.8 1). The following assertions are (obviously) equivalent:
\medskip
\item{a)} $F^1({\got g})=\{0\}$;
\item{b)} $F^0({\got g})={\got g}$;
\item{c)} the cocharacter $\mu:\GG_m\to G$ factors through the center of $G$;
\item{d)} the $G(W(k))$-conjugacy $[\mu]$ of $\mu$ has a unique representative
$\mu:\GG_m\to G$.
\medskip
If these four assertions are true (for instance, this is so if $G$ is a torus), then for any $g\in G(W(k))$, $\vph\circ g=g_1\circ\vph$, with $g_1\in G(W(k))$. By induction, $\forall m\in\NN$ we have $(g\circ\vph)^m=g_m\circ\vph^m$, with $g_m\in G(W(k))$. So in this case, 3.1.0 is a direct consequence of 3.2.1 and 3.2.1.1.
\medskip
{\bf 3.2.2.1. Convention.} In all that follows till the end of \S 3, we assume the assertions a) to d) of 3.2.2 are not true; the only exceptions: in 3.11.1 and 3.11.2 A to C we do not assume this. So $(G,[\mu])$ is a Shimura group pair over ${\rm Spec}(W(k))$ (cf. def. 2.2.5). However, we point out that all results of 3.6, 3.11 and 3.14 below are valid without this assumption (for $G$ a torus, these results are either trivial or are not of any real use).
\medskip
{\bf 3.2.3. A special example.} Let $P$ be the parabolic subgroup of $G$ having
$F^0({\got g})$ as its Lie algebra. As $G$ is quasi-split, it has a Borel subgroup $B$. We can assume $B$ is a subgroup of $P$ (cf. [Bo2, 21.12] applied to $G_{B(k)}$ and the Fact of 2.2.3 3)). Let $T$ be a maximal torus of $B$. It contains a maximal split torus of $G$. As any two such maximal split tori are $G(W(k))$-conjugate (following the proof of Fact 1 of 2.2.9 3), this is a consequence of [Bo2, 15.14]), we can assume $\mu$ factors through $T$; so $\mu$ can be also viewed as a cocharacter $\mu:\GG_m\to T$. 
\smallskip
Let {\got t} (resp. {\got b}) be the Lie algebra of $T$ (resp. of $B$). We have inclusions ${\got t}\subset {\got b}\subset F^0({\got g})$, and $F^1({\got g})\subset{\got b}$. We get $\vph({\got t})\subset\vph({\got b})\subset{\got g}$. As ${\got t}/p{\got t}$ does not intersect $F^1({\got g})/pF^1({\got g})$, $\vph({\got t})$ is a direct summand of ${\got g}$. 
\smallskip
Let $T_0$ be a maximal torus of $G$ having $\vph({\got t})$ as its Lie algebra. The existence of $T_0$ as a maximal torus of $G$ can be checked as follows. A regular element of ${\got t}[{1\over p}]$ is mapped through $\vph$ into a regular element of ${\got g}[{1\over p}]$. We deduce the existence of a maximal torus of $G_{B(k)}$ having $\vph({\got t})[{1\over p}]$ as its Lie algebra; we define $T_0$ as its Zariski closure in $G$. If it is a torus of $G$, then its Lie algebra is automatically $\vph({\got t})$, as $\vph({\got t})$ is a direct summand of ${\got g}$. But to verify that indeed $T_0$ is a torus of $G$ we can assume $k=\bar k$. As $\mu$ factors through $T$, it is appropriate to use 2.2.9 8). Writing $\vph=a\circ\mu({1\over p})$, $a$ becomes a $\sg$-linear automorphism of $M$ and so of $\Mt(M)$. We have $\vph({\got t})=a({\got t)}$. As $a(t_{\al})=t_{\al}$, $\forall\al\in\Mj$, $a$ induces a $\sg$-linear Lie automorphism of ${\got g}$ as well. $a$ takes the Lie algebra of a torus of $GL(M)$ into the Lie algebra of a torus of $GL(M)$ (easy argument at the level of projectors of $M$ onto direct summands of it of rank $1$). So (cf. also the part of [Va2, 4.3.9] involving tori) $T_0$ is a torus of $G$. 
\smallskip
As $\vph$ is a $\sg$-linear automorphism of ${\got g}[{1\over p}]$, there is a Borel subgroup of $G_{B(k)}$ having $\vph({\got b})[{1\over p}]$ as its Lie algebra. Its Zariski closure in $G$ is (cf. the Fact of 2.2.3 3)) a Borel subgroup $B_0$ of $G$. ${\rm Lie}(B_0)$ contains $\vph({\got b})$. We recall that any two Borel subgroups of $G$ are $G(W(k))$-conjugate: this can be deduced from [Bo2, 15.14] (applied to $G_k$) and Fact 1 of 2.2.9 3). So there is an element $g_0\in G(W(k))$ such $g_0B_0g_0^{-1}=B$. Let $T_1:=g_0T_0g_0^{-1}\subset B$. Let $k_1$ be an algebraic field extension of $k$ such that $G_k$ splits over $k_1$; we usually take $k_1$ to be $\bar k$ or a finite field extension of $k$. So $T$ and $T_1$ split over $W(k_1)$. So similarly, we deduce the existence of $g_1\in B(W(k_1))$ such that 
$g_1{T_1}_{W(k_1)}g_1^{-1}=T_{W(k_1)}$. 
\smallskip
Let $\vph_0:=g_0\vph$ and let $\tau_0$ be the formal isogeny type of the  $\sg$-crystal $(M,\vph_0)$. We consider the following Shimura $\sg_{k_1}$-crystal $\bigl(M\otimes_{W(k)} W(k_1),\vph_1,G_{W(k_1)}\bigr)$ defined by
$$
\vph_1:=g_1(\vph_0\otimes 1).
$$
We have $\vph_0({\got b})\subset {\got b}$ and $\vph_0({\got b}[{1\over p}])=
{\got b}[{1\over p}]$. Similarly, we have 
$$\vph_1({\got b}\otimes_{W(k)} W(k_1))\subset
{\got b}\otimes_{W(k)} W(k_1)$$
 and 
$$\vph_1({\got t}\otimes_{W(k)} W(k_1))={\got t}\otimes_{W(k)} W(k_1).$$
\medskip
{\bf 3.2.4. Implications.} Theorem 3.1.0 is a direct consequence of the results 3.2.5-8 below and of the specialization theorem: 
\medskip
-- 3.2.6 and 3.2.7 b) imply 3.1.3 b) and that ${\rm Lie}_G(\tau)$ is the Newton polygon of $({\got g},\vph_0)$;
\smallskip
-- 3.2.6 implies 3.1.3 a) and e) and that $\Mp$ is the Newton polygon attached to $\tau_0$; 
\smallskip
-- 3.2.5 implies (via the above interpretation of ${\rm Lie}_G(\tau)$) $B_G\subset A_G$; 
\smallskip
-- 3.2.7 b) implies $A_G\subset C_G\cap B_G$;
\smallskip
-- 3.2.7 c) implies $C_G\subset A_G$ (so we get $A_G=B_G=C_G$, i.e. we get 3.1.3 c));
\smallskip
-- 3.2.8 implies 3.1.3 d). 
\medskip
In these implications, in order to benefit from the deformation aspect of 3.2.6, we used implicitly 3.2.5-7 for other perfect fields containing $k$. 
\medskip
{\bf 3.2.5. Claim.} {\it Among all Shimura Lie $\sg$-crystals $({\got g},g\vph)$ with
$g\in G(W(k))$, there is no one whose Newton polygon is strictly below the Newton polygon of $({\got g},\vph_0)$. Moreover, any $\sg$-crystal $(M,g\vph)$, with $g\in G(W(k))$, having the property that $({\got g},g\vph)$ and $({\got g},\vph_0)$ have the same Newton polygon, has the formal 
isogeny type $\tau_0$.} 
\medskip
{\bf 3.2.6. Lemma.} {\it For any $g\in G(W(k))$ there is a Shimura filtered $F$-crystal $(M,
F^1,g\vph,G,\tilde f)$ such that the resulting $F$-crystal over the algebraic closure of the field of fractions of the generic point of the special fibre of the completion of $G$ in its origin, has the formal isogeny type $\tau_0$.}
\medskip
{\bf 3.2.7. Proposition. a)} {\it ${\got p}_0:=W_0({\got g},\vph_0)$ is contained in} $F^0({\got g})$.
\smallskip
{\bf b)} {\it The Shimura $\sg$-crystals $(M,g_2\vph,G)$ (with $g_2\in G(W(k))$) having $\tau_0$ as the formal isogeny type are precisely the ones which, up to isomorphisms defined by elements of $G(W(k))$, can be extended to Shimura filtered $\sg$-crystals of the form $(M,F^1,\tilde p\vph_0,G)$, with $\tilde p$ a $W(k)$-valued point of the parabolic subgroup $P_0$ of $G$ satisfying
Lie$(P_0)={\got p}_0$. For any $\tilde p\in P_0(W(k))$, $W_0({\got g},\tilde p\vph_0)$ is ${\got p}_0$ (and so it is included in $F^0({\got g})$), and $({\got g},\vph_0)$ and 
$({\got g},\tilde p\vph_0)$ have the same Newton polygon.}
\smallskip
{\bf c)} {\it Any Shimura filtered $\sg$-crystal $(M,F^1_2,\vph_2,G)$ with the property that $W_0({\got g},\vph_2)$ is contained in $F^0_2({\got g})$, is such that $(M,\vph_2)$ has the formal isogeny type $\tau_0$.}
\medskip
{\bf 3.2.8. The uniqueness part.} The uniqueness of the filtration mentioned in 3.1.0, results from the fact that
any two such filtrations are conjugate by an element of $G(W(k))$ (cf. Fact 2 of 2.2.9 3)) and from the fact
that two $G(B(k))$-conjugate parabolic subgroups of $G_{B(k)}$, containing the same parabolic subgroup, are the same (cf. [Bo2, 14.22 (iii)]).
\medskip
{\bf 3.2.9. Warning.} We present two proofs of 3.2.5-7. The first one ends in 3.4.14 (via 3.5); from many points of view it is just the generalization of the last two paragraphs of 3.1.8.1. The second one ends in 3.7.1-5 for 3.2.6-7 and (cf. 3.3) in 3.4.13 for 3.2.5; it relies on 3.6.0-8.
\medskip
{\bf 3.2.10. An extra set.} Let $D_G$ be the set of $g\in G(W(k))$ such that 
$$g\vph_0=g_3p_0\vph_0g_3^{-1},\leqno (FORM)$$ 
with $p_0\in P_0(W(k))$ and with $g_3\in G(W(k))$ normalizing $F^1/pF^1$.
Let $D_G^1\subset G(W(k_1))$ be defined similarly but with $\vph_0$ replaced by $\vph_1$. 
\medskip\smallskip
{\bf 3.3. About the proof of 3.2.5.}
The statement of 3.2.5 is such that we can replace $k$ by a finite field extension $k_1$ of it. So by a suitable passage from $W(k)$ to $W(k_1)$, with $k_1$ as mentioned, we can assume $G$ is split. As $G$ is split, we can choose $T$ in 3.2.3 to be split, and so $T_0$ and $T_1$ are as well split. So we can take $g_1\in B(W(k))$ to define $\vph_1$ as in 3.2.3. We recall that any Borel subgroup of $G$ is its own normalizer in $G$. 
\smallskip
The claim 3.2.5 is a direct consequence of the following 4 items.
\medskip
{\bf 3.3.1.} There is $b\in B(W(k))$ such that $\vph_0=b\vph_1$: we can take $b=g_1^{-1}$. Even more, $\forall b_1\in B(W(k))$ we have $\vph_1 b_1=b_1^\prime\vph_1$, with $b_1^\prime\in B(W(k))$. Argument: as
$\vph_0({\got b}[{1\over p}])=\vph_1({\got b}[{1\over p}])=
{\got b}[{1\over p}]$, $b_1$ normalizes ${\got b}[{1\over p}]$ and so $b_1\in B(B(k))\cap G(W(k))=B(W(k))$. So by induction on $n\in\NN$, we get: $\vph_0^n=b(n)\vph_1^n$, for some $b(n)\in B(W(k))$.
\medskip
{\bf 3.3.2.} We use freely the Fact of 2.2.11.1 and of 2.2.3 3). $B$ is a subgroup of $P_1$, where $P_1$ is the parabolic subgroup of $G$ defined by the equality ${\rm Lie}(P_1)=W_0({\got g},\vph_1)$. Argument: as ${\got b}\subset F^0({\got g})$, we have $\vph_1({\got b})\subset{\got b}$; so $({\got b},\vph_1)$ is a $\sg$-crystal having only non-negative slopes and so ${\got b}\subset W_0({\got g},\vph_1)$. 
\smallskip
Moreover, from 3.2.1.1 and 3.3.1 we get ${\rm Lie}(P_1)={\got p}_0$. So $P_1=P_0$.
\medskip
{\bf 3.3.3.} 3.3.1 implies that $g_1\in D_G$ and that $D_G^1$ is the right translation of $D_G$ by $b_1^{-1}=g_1$. Let $g_2\in G(W(k))$. The Newton polygon of $({\got g},g_2\vph)$ is strictly below or equal to the one of $({\got g},\vph_0)$ iff $g_2\in D_G$. 
 Warning: this holds without the assumption $G$ is split. The first (resp. second) assertion is proved in 3.4.11 (resp. in 3.4.13).
\medskip
{\bf 3.3.4.} If $g\in D_G$, then $({\got g},g\vph_0)$ (resp. $(M,g\vph_0)$) has the same Newton polygon as $({\got g},\vph_0)$ (resp. as $(M,\vph_0)$). Argument: we can assume that
in 3.2.10, $g_3$ is the identity; then, as in 3.3.1, we have $(p_0\vph_0)^n=p_{n-1}\vph_1^n$, with $p_{n-1}\in P_0(W(k))=P_1(W(k))$, $\forall n\in\NN$, and so 3.2.1 and 3.2.1.1 apply. Here we can replace $\vph_0$ and $D_G$ respectively by $\vph_1$ and $D_G^1$, cf. 3.3.3.
\medskip\smallskip
{\bf 3.4. Split Shimura Lie $\sg$-crystals.} This section and the next one (3.5) are the very heart of the proof of 3.1.0. As a secondary goal of it, we introduce some notations and language to be consistently used in the rest of this paper and in the subsequent ones. Warning: not all cases to be considered below in 3.4-5 can show up in our context of Shimura $\sg$-crystals; we have in mind the so called special cases: cases e6) and e7) of 3.4.3.2, as well as the case d) of 3.4.3.2 with the set $A$ to be introduced in 3.5.1 having three elements or having two elements which via outer automorphisms of ${\got g}_1$ are not identifiable (cf. [Sa, p. 458-60] and its obvious passage from $\RR$ to an arbitrary field of characteristic $0$; see also [De2, 1.3.9-10] or [Se2, \S 3]). However, all such cases do show up in the generalized Shimura context. So, for the sake of completeness and of future references, in 3.4-5 we treat these cases as well, without any extra comment. The only two modifications needed to be made in the whole of 3.4-5 in order to accommodate this generalized Shimura context: 
\medskip
-- the element $\bar h$ of 3.4.0 has to be modified accordingly;
\smallskip
-- the $g_3$ elements of 3.2.10 (this is in connection to 3.4.11), 3.4.9 and 3.4.13 have to be modified accordingly (see the part of 2.2.14.2 referring to a generalized Shimura context). 
\medskip
So, the reader interested just in the context of Shimura $p$-divisible groups (i.e. in what we considered till now in 3.2-3), can entirely ignore these cases.
\smallskip
We refer to the situation of the beginning of 3.3: $G$ and $T$ are split and $g_1\in G(W(k))$. We have a direct sum decomposition
$${\got g}[{1\over p}]={\got z}[{1\over p}]\bigoplus\oplus_{i\in I}
{\got g}_i[{1\over p}]$$
of Lie algebras over $B(k)$. Here {\got z} is the Lie algebra of the maximal subtorus $Z(G)^0$ of $Z(G)$, while ${\got g}_i$ is the Lie algebra of a semisimple subgroup $G_i$ of $G$, such that 
$$
G^{\rm ad}=\prod_{i\in I}G_i^{\rm ad}
$$
(cf. [Ti1, p. 46] applied over $k$ and logically lifted to $W(k)$), with $G_i^{\rm ad}$ an absolutely simple, split, adjoint group over $W(k)$, $\forall i\in I$. It is defined by a Lie monomorphism
$$
i_G:{\got z}\oplus\bigoplus_{i\in I}{\got g}_i\hookrightarrow{\got g}
$$ 
which becomes an isomorphism by inverting $p$. $i_G$ is not always surjective; as we assumed $p\ge 3$, the failure of surjectivity can happen only when $p=3$ and there are factors $G_i$ of $E_6$ Lie type or when $p$ is arbitrary and there are factors $G_i$ of $A_\ell$ Lie type, with $p$ dividing $\ell+1$. 
\smallskip
Moreover, we have a natural isogeny $$
Z(G)^0\prod\prod_{i\in I} G_i\to G
$$
which at the level of Lie algebras defines $i_G$.
\medskip
{\bf 3.4.0. Notations.} $\vph_1$ normalizes {\got z} and permutes the factors  
${\got g}_i[{1\over p}]$, $i\in I$, producing a permutation $\ga$ of $I$. 
We write 
$$
\ga=\prod_{s\in\Mj_0}\ga_s
$$ 
as a product of disjoint cyclic permutations $\ga_s$. Warning: trivial permutations are allowed, i.e. each $i\in I$ fixed by $\gamma$ defines an element $s_i\in\Mj_0$ such that $\gamma_{s_i}$ is the trivial permutation of $I$ and so of $\{i\}$; for such an $i$, we still say $\gamma_{s_i}$ permutes cyclically the subset $\{i\}$ of $I$. As ${\got z}\subset F^0({\got g})$, the $\sg$-crystal $({\got z},\vph_1)$ has all slopes $0$.
We start working with just one cycle $\ga_0$ of $\ga$ (so we consider that $0\in\Mj_0$) permuting cyclically a subset $I_0$ of $I$. We can assume that $I_0=S(1,n)$, $n\in\NN$, and that $\ga_0$ is the
permutation: $i\to i+1$ (with $n+1=1$). Let 
$${\got g}_0:=\bigoplus_{i\in I_0}{\got g}_i.$$
We fix an element ${g}_2\in G(W(k))$; let
$$\vph_2:={g}_2\circ\vph_1.$$
\indent
Let $\bar h:M\to M$ be the endomorphism defined by: $\bar h(x)=0$ if $x\in F^0$ and 
$\bar h(x)=x$ if $x\in F^1$. We recall that the cocharacter $\mu:\GG_m\to G$ (of the beginning of \S 3) produces a direct sum decomposition
$M=F^1\oplus F^0$, with $\be\in\GG_m(W(k))$ acting through $\mu$ on $F^i$ as the multiplication with $\be^{-i}$,
$i=\overline{0,1}$; so 
$$\bar h\in d\mu ({\rm Lie}(\GG_m))\subset {\rm Lie}(T)=
{\got t}\subset {\rm End}(M).\leqno (HBAR)
$$ 
Let $\bar h_i$ be the component of $\bar h$ in 
${\got g}_i\bigl[{1\over p}\bigr]$ (corresponding to the direct sum decomposition ${\got g}[{1\over p}]={\got z}[{1\over p}]\bigoplus\oplus_{i\in I} {\got g}_i[{1\over p}]$). If $\bar h_i=0$ then ${\got g}_i\subset F^0({\got g})$. If $\bar h_i\ne 0$ then $\bar h_i$ is a semisimple element of ${\got g}_i\bigl[{1\over p}\bigr]$. When viewed as an endomorphism of ${\got g}_i[{1\over p}]$ (via $z_i\in {\got g}_i$ is mapped into $[\bar h_i,z_i]$), its eigenvalues are precisely $-1$, $0$ and $1$. 
\smallskip
We have $F^1({\got g})=\bigl\{x\in{\got g}\bigm|[\bar h,x]=x\bigr\}$; so 
$$
F^1({\got g}_i):=F^1({\got g})\cap{\got g}_i=\bigl\{x\in{\got g}_i\bigm|[\bar h_i,x]=x\bigr\}.
$$ 
Similarly, $F^0({\got g})=\bigl\{x\in{\got g}\bigm|[\bar h,x]=0\ {\rm or}\ x\bigr\}$ and so 
$$
F^0({\got g}_i):=F^0({\got g})\cap{\got g}_i=\bigl\{x\in{\got g}_i\bigm|[\bar h_i,x]=0\ {\rm or}\ x\bigr\}.
$$ 
Let $F^j({\got g}_0):=\bigoplus_{i\in I_0} F^j({\got g}_i)$,
$j=\overline{0,1}$. We get:
\medskip
{\bf 3.4.1. Lemma.} {\it The quadruple $\bigl({\got g}_0,\vph_j,F^0({\got g}_0),F^1({\got g}_0)\bigr)$ 
is a Shimura filtered Lie $\sg$-crystal, $j=\overline{1,2}$.}
\medskip
{\bf Proof:} As $\mu$ is a cocharacter over $W(k)$, we have $F^s({\got g}_i)/pF^s({\got g}_i)={\got g}_i/p{\got g}_i\cap F^s({\rm End}(M))/p F^s({\rm End}(M))$, $\forall s\in\{0,1\}$. So $\vph_j\bigl({1\over p}F^1({\got g}_i)+F^0({\got g}_i)+p{\got g}_i\bigr)$ is a direct summand of {\got g} of rank $\dim_{W(k)}({\got g}_i)=\dim_{W(k)}({\got g}_{i+1})$; it is included in ${\got g}_{i+1}\bigl[{1\over p}\bigr]$ and so it is ${\got g}_{i+1}$ ($\forall i\in I_0$ and $\forall j\in\{1,2\}$). So $\bigl({\got g}_0,\vph_j,F^0({\got g}_0),F^1({\got g}_0)\bigr)$ is a Lie $p$-divisible object of $\Mm\Mf_{[-1,1]}(W(k))$. It is a Shimura filtered Lie $\sg$-crystal,  as the cocharacter of $\prod_{i\in I_0} G_i^{\rm ad}$, obtained by composing $\mu$ with the natural epimorphism $G\twoheadrightarrow\prod_{i\in I_0} G_i^{\rm ad}$, shows that the conditions of 2.2.11 1) are satisfied. This ends the proof.
\medskip
{\bf 3.4.1.1. Remark.}
We can work equally well 3.4.1 (as well as all that follows) with
${\rm Lie}\bigl(\sum_{i\in I_0}G_i\bigr)$ (we have ${\got g}_0\subset {\rm Lie}(\sum_{i\in I_0}G_i)\subset{\got g}_0[{1\over p}]$) or with ${\rm Lie}\bigl((\sum_{i\in I_0}G_i)^{\rm ad}\bigr)=\oplus_{i\in I_0} {\rm Lie}(G_i^{\rm ad})$ instead of ${\got g}_0$ (cf.  2.2.13); here $\sum_{i\in I_0}G_i$ is the semisimple subgroup of $G$ generated by all $G_i$'s, $i\in I_0$.
\medskip
{\bf 3.4.1.2. The non-trivial context assumption.} 
We assume $F^1({\got g}_0)\ne 0$; otherwise there is nothing to be done (to be compared with 3.2.2): the Newton polygons of $({\got g}_0,\vph_1)$ and of $({\got g}_0,\vph_2)$ have only one slope $0$ and so are the same.
\medskip
{\bf 3.4.1.3. A convenient context.} Let $\tilde G_0$ be the subgroup of $G$ generated by the torus $T$ and by the semisimple subgroups $G_i$, $i\in I_0$. So ${\tilde G}^{\rm ad}_0=\prod_{i\in I_0} G_i^{\rm ad}$. In the same way we got $i_G$, we get that ${\got g}_0$ is (identifiable with) a Lie subalgebra of ${\rm Lie}(\tilde G_0)$. As the cocharacter $\mu$ of 3.2.3 can be viewed as well as a cocharacter of $\tilde G_0$ (still denoted by $\mu$), we get (cf. also 2.2.9 1) and 1')) that the triple  $(M,\vph_1,\tilde G_0)$ is a potentially Shimura $\sg$-crystal. Moreover, we get a Shimura group pair $(\tilde G_0,[\mu])$ over $W(k)$.  
\medskip
{\bf 3.4.2. Slopes (w.r.t. $I_0$).}
Let ${\got H}_2$ (resp. ${\got H}_1$) be the set of slopes of the isocrystal $\bigl({\got g}_0[{1\over p}],\vph_2\bigr)$ (resp. of the isocrystal
$\bigl({\got g}_0[{1\over p}],\vph_1\bigr)$). The sets ${\got H}_2$ and ${\got H}_1$ are subsets of the real interval $[-1,1]$. We have slope type direct sum decompositions 
$$
{\got g}_0[{1\over p}]=\bigoplus_{\al\in\Mh_j}{g}_j(\al),
$$  
with ${g}_j(\al):=W(\al)({\got g}_0[{1\over p}],\vph_j)$, $j=\overline{1,2}$ (cf. 2.2.3 3)). As $\vph^n_1$ and $\vph^n_2$ normalize any one of ${\got g}_i[{1\over p}]$, $i\in I_0$, we get:
$$g_j(\al)=\bigoplus_{i\in I_0}g_j(\al)\cap
{\got g}_i\bigl[{1\over p}\bigr],\leqno ({\bf 3.4.2.1})$$
$\forall j\in\{1,2\}$ and $\forall \al\in\Mh_j$. As $\vph_j$ maps $g_j(\al)\cap{\got g}_i[{1\over p}]$ bijectively onto $g_j(\al)\cap{\got g}_{i+1}[{1\over p}]$, we have:
\medskip
{\bf 3.4.2.2. Fact.} {\it $\forall j\in\{1,2\}$, the dimension of $g_j(\al)\cap{\got g}_i[{1\over p}]$ does not depend on $i\in I_0$.}
\medskip
{\bf 3.4.3. Roots.}
Let ${\got t}_i$ (resp. ${\got b}_i$) be the component of {\got t} (resp. of {\got b}) in
${\got g}_i$; so ${\got t}_i:={\got g}_i\cap{\got t}\bigl[{1\over p}\bigr]$
(resp. ${\got b}_i:={\got g}_i\cap{\got b}\bigl[{1\over p}\bigr]$). Let $\Dl_i$ be the base
of roots defined by the action (via inner conjugation) on ${\got g}_i$ of the image $T_{i\rm ad}$ in $G^{\rm ad}_i$ of the maximal torus $T_i$ of $G_i$ having ${\got t}_i$
as its Lie algebra, corresponding to the Borel Lie subalgebra ${\got b}_i$ of ${\got g}_i$. Let $\Phi_i$ be the set of roots of ${\got g}_i$ relative to the torus $T_{i\rm ad}$. So $\Phi_i$ is a finite set of characters of $T_{i\rm ad}$. Then $\Phi^+_i$ (the set of positive roots) and $\Phi^-_i$ (the set of negative roots) are defined w.r.t. $\Dl_i$. So $\Dl_i\subset\Phi^+_i$. Here $i\in I$.
\medskip
{\bf 3.4.3.0. Weyl's decompositions.}
We assume here $k=\bar k$. Let 
$${\got z}\oplus\bigoplus_{i\in I} {\got t}_i\oplus\bigoplus_{\scriptstyle i\in I\atop\scriptstyle \al\in\Phi_i}\tilde g_i(\al)\subset {\got g}$$
be the Weyl direct sum decomposition of ${\rm Lie}(Z(G)^0\times\times_{i\in I} G_i)$ associated to the natural action (via inner conjugation) of the image $\prod_{i\in I} T_{i\rm ad}$ of the maximal torus $T$ of $G$ in $G^{\rm ad}$. As $\vph_1({\got t})={\got t}$, $\vph_1$ permutes the $B(k)$-vector subspaces $\tilde g_i(\al)[{1\over p}]$ of ${\got g}[{1\over p}]$. We get a permutation $\ga_{\Phi}$ of the disjoint union $\cup_{i\in I} \Phi_i$. We have
$$
\vph_1\bigl(\tilde g_i(\al)\bigr)=p^{s_i(\al)}\tilde g_{\ga(i)}
\bigl(\ga_{\Phi}(\al)\bigr),\leqno (1) 
$$
where $s_i(\al)\in\{-1,0,1\}$ is such that the cocharacter $\mu:\GG_m\to T$ has the property that $\mu(\be)$ acts on elements of $\tilde g_i(\al)$ as the multiplication with $\be^{-s_i(\al)}$, $\forall\be\in\GG_m(W(k))$. As $\gamma_{\Phi}$ takes $\cup_{i\in I} \Phi_i^+$ onto $\cup_{i\in I} \Phi_i^+$ and $\cup_{i\in I} \Phi_i^-$ onto $\cup_{i\in I} \Phi_i^-$, and as $s_i(\al)\le 0$ (resp. $s_i(\al)\ge 0$) if $\al\in\Phi_i^-$ (resp. if $\al\in\Phi_i^+$) (these two facts are a consequence of end of 3.2.3), we deduce that Lie$(P_1)\subset F^0({\got g})$. So, forgetting the Lie structure, the pair $({\rm Lie}(P_1),\vph_1)$ is a $\sg$-crystal.
\medskip
{\bf 3.4.3.1. The number $e$ (w.r.t. $I_0$).}
$\vph^n_1$ normalizes ${\got g}_1[{1\over p}]$, ${\got b}_1$ and ${\got t}_1$ and so it permutes the $1$-dimensional $B(k)$-vector subspaces of ${\got b}_1[{1\over p}]$ normalized by ${\got t}_1[{1\over p}]$; let $\pi_1$ be the resulting permutation of $\Dl_1$. Let $e$ be the order of $\pi_1$. It is very convenient to work simultaneously with the Shimura Lie $\sg$-crystal 
$$
({\got g}_0^e,\vph_j^{(e)}), 
$$
where ${\got g}_0^e:={\got g}_0\oplus\cdots\oplus{\got g}_0$ ($e$ summands) and $\vph_j^{(e)}$ is the $\sg$-linear automorphism of ${\got g}_0^e[{1\over p}]$ acting like $\vph_j$ on any direct summand  ${\got g}_i$, $i\in I_0$,  of a 
${\got g}_0$ copy of ${\got g}_0^e$, but sending the ${\got g}_n$ direct summand  of the $s$-th  
${\got g}_0$ copy to the ${\got g}_1$ direct summand of the $(s+1)$-th ${\got g}_0$ copy,  $\forall s\in\ZZ/e\ZZ$. The proof of the fact that the pair $({\got g}_0^e,\vph_j^{(e)})$ is a Shimura Lie $\sg$-crystal is the same as of 3.4.1. The set of slopes of  $({\got g}_0^e,\vph_j^{(e)})$ is the set of slopes of $({\got g}_0,\vph_j)$, all multiplicities being $e$-times more, $j=\overline{1,2}$. 
\smallskip
The permutation of $\Dl_1$ defined by $(\vph_1^{(e)})^{ne}$ is the trivial one. \smallskip
$e$ can be 1, 2 or 3. It can be 3 only in the case when ${\got g}_1$ is of $D_4$ Lie type; it can be to 2 only when ${\got g}_1$ is of  $A_\ell$ ($\ell\ge 2$), $D_\ell$ ($\ell\ge 4$) or $E_6$ Lie type (and then $\pi_1$ is a usual involution of $\Dl_1$; outside of the $D_4$ Lie type we can replace here ``a" by the). All these are a consequence of the structure of the group of automorphisms of the Dynkin diagram of $\Phi_1$: it is trivial or isomorphic to $S_3$ or to $\ZZ/2\ZZ$. The value of $e$ does not depend on the choice of $i\in I_0$ used to define it (above we chose $1\in I_0$, i.e. we worked with $\Dl_1$).
\smallskip 
Let $\tilde I_0:=S(1,en)$. We like to think of ${\got g}_0^e$ as a direct sum $\oplus_{i\in\tilde I_0} {\got g}_i$. Here for $s\in S(1,e-1)$ and $i\in I_0$, ${\got g}_{i+sn}:={\got g}_i$; similarly ${\got t}_{i+sn}:={\got t}_i$, $\Dl_{i+sn}:=\Dl_i$, ${\got b}_{i+sn}:={\got b}_i$ and $\bar h_{i+sn}:=(\vph_1^{(e)})^{sn}(\bar h_i)\in {\got g}_{i+sn}[{1\over p}]$.
\medskip
{\bf 3.4.3.2. Simple roots (w.r.t. $I_0$).}
Let $\ell$ be
the rank of ${\got g}_1$. For $i\in\tilde I_0$, let
$\al_1(i),\ldots,\al_\ell(i)$ be the simple roots of $\Dl_i$, numbered in such a way that 
\medskip
-- $\vph_1^{(e)}$ sends the ($1$ dimensional) $B(k)$-vector subspace of ${\got g}_i[{1\over p}]$ corresponding to $\al_s(i)$ (i.e. on which $T_{i\rm ad}$ acts via its character defining $\al_s(i)$), $s\in S(1,\ell)$,
to the ($1$ dimensional) $B(k)$-vector subspace of ${\got g}_{i+1}[{1\over p}]$ corresponding to $\al_s(i+1)$ (with $en+1=1$),
and 
\smallskip
-- the numbering for $i=1$ is the one from [Bou2, planche I-VI]. 
\medskip
The maximal root of $\Dl_1$ (dropping the lateral right index (1) from now on) depends on the Lie type of ${\got g}_1$. It is (cf. loc. cit.):
\medskip
\item{{\bf a)}} $\al_1+\al_2+\cdots+\al_\ell$ if ${\got g}_1$ is of  $A_\ell$ Lie type ($\ell\ge 2$);
\smallskip
\item{{\bf b)}} $\al_1+2\al_2+\cdots+2\al_\ell$ if ${\got g}_1$ is of  $B_\ell$ Lie type ($\ell\ge 1$);
\smallskip
\item{{\bf c)}} $2\al_1+2\al_2+\cdots+2\al_{\ell-1}+\al_\ell$ if ${\got g}_1$ is of $C_\ell$ Lie type ($\ell\ge 2$);
\smallskip
\item{{\bf d)}} $\al_1+2\al_2+\cdots+2\al_{\ell-2}+\al_{\ell-1}+\al_\ell$ if ${\got g}_1$
is of  $D_\ell$ Lie type ($\ell\ge 4$);
\smallskip
\item{{\bf e6)}} $\al_1+2\al_2+2\al_3+3\al_4+2\al_5+\al_6$ if ${\got g}_1$ is of $E_6$ Lie type;
\smallskip
\item{{\bf e7)}} $2\al_1+2\al_2+3\al_3+4\al_4+3\al_5+2\al_6+\al_7$ if ${\got g}_1$ is of $E_7$ Lie type.
\medskip
${\got g}_1$ is one of these Lie types, cf. 2.2.7 and 3.4.1.2-3.
We usually prefer to consider $B_\ell(\ell\ge 1)$ and $A_\ell(\ell\ge 2)$,
i.e. we prefer to consider the Lie type $B_1$ instead of $A_1$; but we are not bothered if we are dealing with $B_2$ or with $C_2$. In what follows, depending on the Lie type of ${\got g}_1$, we refer to case a), to case b), etc.
\medskip
{\bf 3.4.3.3. Possibilities for $\bar h_1$.} Let 
$${\got g}_1={\got t}_1\oplus\bigoplus_{\be\in\Phi_1}{\got g}_\be$$ 
be the Weyl direct sum 
decomposition associated to the action of the maximal split torus $T_{1\rm ad}$ of $G_1^{\rm ad}$ on ${\got g}_1$ defined by inner conjugation. So $T_{1\rm ad}$ acts on ${\got g}_\be$ through the character defining the root $\be$. 
\smallskip
In cases b), c) and e7), $\bar h_1\in {\got t}_1[{1\over p}]$ is either $0$ or is defined as follows:
\medskip\noindent
-- in case b), $\bar h_1$ acts as zero on ${\got g}_{\al_s}$, $s\in S(2,\ell)$, and as identity on ${\got g}_{\al_1}$;
\medskip\noindent
-- in case c), $\bar h_1$ acts as zero on ${\got g}_{\al_s}$, $s\in S(1,\ell-1)$, and as identity on ${\got g}_{\al_\ell}$;
\medskip\noindent
-- in case e7), $\bar h_1$ acts as zero on ${\got g}_{\al_s}$, $s\in S(1,6)$, and as identity on ${\got g}_{\al_7}$.
\medskip
In case a), $\bar h_1\in {\got t}_1[{1\over p}]$ is either $0$ or there is $s_1\in S(1,\ell)$ such that
it acts as zero on ${\got g}_{\al_s}$, $s\in S(1,\ell)\setminus\{s_1\}$, and as identity on ${\got g}_{\al_{s_1}}$ 
(we have $\ell$ non-zero possibilities for $\bar h_1$).
\smallskip
In case d), $\bar h_1\in {\got t}_1[{1\over p}]$ is either $0$ or there is $s_1\in\{1,\ell-1,\ell\}$ such that
it acts as zero on ${\got g}_{\al_s}$, $s\in S(1,\ell)\setminus\{s_1\}$, and as identity on ${\got g}_{\al_{s_1}}$ (we have 3 non-zero possibilities for $\bar h_1$).
\smallskip
In case e6), $\bar h_1\in {\got t}_1[{1\over p}]$ is either $0$ or there is $s_1\in\{1,6\}$ such that it acts as zero on ${\got g}_{\al_s}$, $s\in S(1,6)\setminus\{s_1\}$, and as identity on ${\got g}_{\al_{s_1}}$ (we have 2 non-zero possibilities for $\bar h_1$).
\medskip
{\bf 3.4.3.4. The $\vep$'s (w.r.t. $I_0$).} For $i\in\tilde I_0$, we define: $\vep_i=0$ if $\bar h_i=0$ and $\vep_i=s$ if $\bar h_i$ acts as identity on the rank 1 direct summand of ${\got g}_i$ corresponding to the root
$\al_s(i)$; here $s\in S(1,\ell)$. So, fixing $i$, $\vep_i$ can take more than $1$ non-zero value only if we are in the cases a), d) or e6).
\smallskip
The assumption $F^1({\got g}_0)\ne 0$ implies that the set 
$$
I_1:=\{i\in I_0\bigm|\vep_i\ne 0\}
$$ 
is not empty. Let 
$$
\tilde I_1:=\{i\in\tilde I_0\bigm|\vep_i\ne 0\}.
$$ 
As $I_0\subset \tilde I_0$, we have $I_1\subset\tilde I_1$.
\medskip
{\bf 3.4.4. The smallest slope $\dl$ (w.r.t. $I_0$).} We have:
\medskip
{\bf Fact.} {\it The sets $\Mh_1$ and $\Mh_2$ are subsets of the interval
$\bigl[-{\abs{I_1}\over n},{\abs{I_1}\over n}\bigr]$; they are symmetric w.r.t. its mid-point $0$. Moreover, the smallest element of $\Mh_1$ is 
$$
\dl:=-{1\over n}\abs{I_1}<0.\leqno (MIN)
$$}
\smallskip
{\bf Proof:} {The first part is trivial, cf. the definition of $I_1$ and 2.2.3 1). For the second part, we just have to point out that ${\got g}_{-\be}$, with $\be\in\Phi_1^+$ as the maximal root w.r.t. the base $\Dl_1$, is included in $g_1(\dl)$.
\medskip
{\bf 3.4.5. Faltings--Shimura--Hasse--Witt shifts and maps (w.r.t. $I_0$).} Let $j\in\{1,2\}$. Let $\Psi^0_j:{\got g}_0\to{\got g}_0$ be the $\sg$-linear map defined by:
if $x\in{\got g}_i$, then
\medskip
\item{i)} $\Psi^0_j(x):=\vph_j(x)$ if $\vep_i=0$;
\smallskip
\item{ii)} $\Psi^0_j(x):=p\vph_j(x)$ if $\vep_i\ne 0$.
\medskip
Taking everything mod $p$, we get $\sg$-linear maps 
$$\bar\psi^0_j:{\got g}_0/p{\got g}_0\to
{\got g}_0/p{\got g}_0.$$ 
We have $\Psi^0_2=g_2\Psi^0_1$ and $\bar\psi_2^0=g_2\bar\psi_0^0$. 
\smallskip
Denoting $\tilde F^{l+1}({\got g}_i):=F^l({\got g}_i)$ if $\vep_i\neq 0$ and $\tilde F^{l+1}({\got g}_i):=\{0\}$ if $\vep_i=0$, $l\in\{1,2\}$, the quadruple 
$$\bigl({\got g}_0,\oplus_{i\in I_0} \tilde F^1({\got g}_i),\oplus_{i\in I_0} \tilde F^2({\got g}_i),\Psi_j^0\bigr)$$
is a $p$-divisible object of $\Mm\Mf_{[0,2]}(W(k))$; warning: $\Psi_j^0$ does not preserve the Lie bracket of ${\got g}_0$, cf. 3.4.1.2. We call it (resp. just the pair $({\got g}_0,\Psi_j^0)$) the Faltings--Shimura--Hasse--Witt shift of $\bigl({\got g}_0,\vph_j,F^0({\got g}_0),F^1({\got g}_0)\bigr)$ (resp. of $({\got g}_0,\vph_j)$). We also refer to $({\got g}_0,\Psi_j^0)$ as the Faltings--Shimura--Hasse--Witt shift of $(M,\vph_j,G)$ corresponding to $I_0$, and to $\Psi_j^0$ as a Faltings--Shimura--Hasse--Witt shift. Similarly, we refer to $\bar\psi^0_j$ as the Faltings--Shimura--Hasse--Witt map of $\bigl({\got g}_0,\vph_j,F^0({\got g}_0),F^1({\got g}_0)\bigr)$ or of $({\got g}_0,\vph_j)$, or as the Faltings--Shimura--Hasse--Witt map of $(M,\vph_j,G)$ w.r.t. $I_0$, or just as a Faltings--Shimura--Hasse--Witt map.
\smallskip
We can assume (in accordance to A of 2.2.22 3)) that $\vep_1\ne 0$. Let 
$$
N_2:=g_2(\dl)\cap{\got g}_1. 
$$
For any
$i\in\tilde I_1$, let ${\got p}_{\vep_i}$ be the parabolic Lie subalgebra of ${\got g}_1$
containing ${\got b}_1$ and such that its nilpotent radical ${\got s}_{\vep_i}$ is the $W(k)$-submodule of ${\got g}_1$ generated by all ${\got g}_\be$'s, with $\be\in\Phi_1^+$ having $+\al_{\vep_i}$ in its expression. Let
${\got q}_{\vep_i}$ be the parabolic Lie subalgebra of ${\got g}_1$, opposite to 
${\got p}_{\vep_i}$ (w.r.t. ${\got t}_i$ and $\Dl_i$) and let ${\got r}_{\vep_i}$ be its nilpotent  radical. The multiplicity 
$$m_{\dl}:=\dim_{B(k)} g_1(\dl)$$ 
of the slope $\dl$ for
$({\got g}_0,\vph_1)$ is (cf. the choice of numberings in 3.4.3.2) $n$ times the rank of the intersection 
$$\bigcap_{i\in \tilde I_1}{\got r}_{\vep_i}.$$
The multiplicity $\dim_{B(k)} g_2(\dl)$ of the slope $\dl$ for $({\got g}_0,\vph_2)$ is $n$ times the rank of $N_2$ (cf. 3.4.2.2). We denote by ${\got p}_{\vep_i}(i)$, 
${\got q}_{\vep_i}(i)$, ${\got r}_{\vep_i}(i)$ and respectively ${\got s}_{\vep_i}(i)$ the intersection of ${\got g}_i$ with $\vph_1^{i-1}\bigl({\got p}_{\vep_i}[{1\over p}]\bigr)$,
$\vph_1^{i-1}\bigl({\got q}_{\vep_i}[{1\over p}]\bigr)$, $\vph_1^{i-1}\bigl({\got r}_{\vep_i}[{1\over p}]\bigr)$ and respectively with $\vph_1^{i-1}\bigl({\got s}_{\vep_i}[{1\over p}]\bigr)$, $i\in\tilde I_1$. So ${\got p}_{\vep_i}(i)=F^0({\got g}_i)$, $\forall i\in I_1$; in particular, ${\got p}_{\vep_i}(1)={\got p}_{\vep_i}$. 
\medskip
{\bf 3.4.5.1. On multiplicities of $\dl$.} The pair $({\got g}_0,\Psi_2^0)$ is a $\sg$-crystal. We have $W^0({\got g}_0,\Psi_2^0)=g_2(\dl)\cap {\got g}_0$ and $W_{-2\dl}({\got g}_0,\Psi_2^0)=g_2(-\dl)\cap {\got g}_0$. We get:
\medskip
{\bf Fact.} {\it $\bigl(g_2(-\dl)\cap {\got g}_1\bigr)\otimes_{W(k)} k\subset {\got s}_{\vep_1}.$}
\medskip
{\bf Proof:} We can assume $k=\bar k$. Let $L:=\Psi_2^0$. Let $m\in\NN$ be a multiple of $en$ such that 
$$s:=-2m\dl\in\NN.$$
We consider the $\ZZ_p$-submodule $M_2(s)$ of $g_2(-\dl)\cap {\got g}_1$ generated by elements on which $L^m$ acts as scalar multiplication with $p^s$. Its rank is $\dim_{W(k)}(g_2(-\dl)\cap {\got g}_1)$. Let $N_2(s)$ be the $W(k)$-submodule of $g_2(-\dl)\cap {\got g}_1$ generated by $M_2(s)$. We can assume $L^m(g_2(-\dl)\cap {\got g}_1)\subset g_2(-\dl)\cap {\got g}_1$.
\smallskip
It is easy to see (based on 3.4.1 and the definition of $\vep_i$'s) that for any $x\in {\got g}_1\setminus p{\got g}_1$, $L^m(x)\notin p^{s+1}{\got g}_1$; moreover, if $L^m(x)\in p^{s}{\got g}_1$, then $x$ mod $p$ belongs to ${\got s}_{\vep_1}$ (if $m$ is big enough we also get that $x\in g_2(-\dl)\cap {\got g}_1$). In particular, all Hodge numbers of the $\sg^m$-linear endomorphism $L^m$ of $g_2(-\dl)\cap {\got g}_1$ are at most $s$. As the length of $g_2(-\dl)\cap {\got g}_1/L^m(g_2(-\dl)\cap {\got g}_1)$ and of $N_2(s)/L^m(N_2(s))$ as $W(k)$-modules are the same, we get that all these Hodge numbers must be $s$ and so the Fact follows.
\medskip       
{\bf A. Corollary.} {\it We have $N_2(s)=g_2(-\dl)\cap {\got g}_1$.}  
\medskip
{\bf B.} As $L^m(N_2)\subset N_2$ and the slopes of $(N_2,L^m)$ are all $0$, as in the above proof we get $L^m(N_2)=N_2$. So, if $k=\bar k$, $N_2$ is generated by elements fixed by $L^m$. On the other hand, it is well known that $N_2/pN_2$ is $\cap_{q\in\NN} \tilde L^q({\got g}_1/p{\got g}_1)$, where $\tilde L$ is the truncation mod $p$ of the restriction of $L^m$ to ${\got g}_1$. 
\medskip
{\bf 3.4.6. Theorem.} {\it If $\dl\in\Mh_2$, then the multiplicity of the slope $\dl$ of the Shimura Lie $\sg$-crystal $({\got g}_0,\vph_2)$ is less or equal to $m_{\dl}$. If it is $m_{\dl}$, then $\forall i\in I_1$, we have}
$$\bigl(\bigl(\bigoplus_{\al\in\Mh_2\cap[0,1]}g_2(\al)\bigr)\cap{\got g}_i\bigr)\otimes_{W(k)} k\subset {\got p}_{\vep_i}(i)\otimes_{W(k)} k.\leqno {\bf (SUB)}$$}
\indent
{\bf 3.4.7. Corollary.} {\it The fact that this multiplicity is $m_{\dl}$ or not depends only on the expression of $g_2$ mod $p$.}
\medskip
Before proving 3.4.6, we first mention some of its immediate implications. Not to complicate the notations by introducing extra indices, in 3.4.8-12 below we allow $g_2$ to vary as needed.
\medskip
{\bf 3.4.8. Proposition.} {\it The Newton polygon of $({\got g}_0,\vph_1)$
is not strictly above the Newton polygon of any other Shimura Lie $\sg$-crystal of the form
$({\got g}_0,g_2\vph_1)$, with $g_2\in\tilde G_0(W(k))$.}
\medskip
{\bf Proof:} We assume 3.4.8 is false; let $g_2\in\tilde G_0(W(k))$ for which 3.4.8 is not true. As $\dl\in\Mh_1$, from 3.4.4 and 3.4.6 we deduce that the multiplicity of the slope $\dl\in\Mh_2$ for $({\got g}_0,\vph_2)$ is also $m_{\dl}$. 
\smallskip
Let $P_2$ be the parabolic subgroup of $\tilde G_0$ whose Lie algebra is $W_0({\rm Lie}(\tilde G_0),\vph_2)$. $p\vph_2$ takes ${\rm Lie}(P_2)$ into itself. (SUB) implies that ${P_2}_k$ is included in $P^0_k$ (cf. also the Fact of 2.2.11.1); here $P^0$ is the parabolic subgroup of $\tilde G_0$ having $F^0({\got g})\cap {\rm Lie}(\tilde G_0)$ as its Lie algebra. We deduce the existence of an element $g_3\in{\rm Ker}\bigl(\tilde  G_0(W(k))\to\tilde G_0(k)\bigr)$ such that $P_3:=g_3^{-1}P_2g_3$ is a subgroup of $P^0$ (to be compared with Fact 1 of 2.2.9 3)). As 3.4.8 is a statement involving Newton polygons, by replacing $\vph_1$ by $\tilde a\vph_1\tilde a^{-1}$, with $\tilde a\in P^0(W(k))$, we can assume ${\got b}\cap {\rm Lie}(\tilde G_0)$ is included in ${\rm Lie}(P_3)$.
\smallskip
Let $\vph_3:=g_3^{-1}\vph_2g_3$. We can write $\vph_3=g_4\vph_1$, with
$g_4\in\tilde G_0(W(k))$. $({\got g}_0,\vph_3)$ and $({\got g}_0,\vph_2)$
have the same Newton polygon. Let $B_0$ be the Borel subgroup of $\tilde G_0$ whose Lie algebra is ${\got b}\cap {\rm Lie}(\tilde G_0)$. We have $B_0\subset P^0$, cf. 3.2.3. Repeating the argument of 3.2.3, we deduce the existence of an element $g_5\in P_3(W(k))$ such that $\vph_4=g_5\vph_3$ takes ${\got t}$
onto ${\got t}$ and ${\rm Lie}(B_0)$ into ${\rm Lie}(B_0)$. For $s\in\NN$ we get $\vph^s_4=p_{3,s}\vph^s_3$, with $p_{3,s}\in P_3(W(k))$ and $p_{3,1}=g_5$; also $\vph^s_4=b(s)\vph_1^s$, with $b(s)\in B_0(W(k))$ and $b(1)=g_5g_4$ (all these are as in 3.3.1 and 3.3.4; the fact that in 3.4.1.3 we are just in a potentially context does not represent an obstruction, as we can see moving to $\bar k$). This implies ${\rm Lie}(P_3)={\rm Lie}(\tilde G_0\cap P_1)$, cf. 3.2.1.1; so $P_3=P_{10}:=\tilde G_0\cap P_1$. This implies: $({\got g}_0,g_2\vph_1)$ has a lift of parabolic type.
\smallskip
From 3.2.1 and 3.2.1.1 we also get that $({\got g}_0,\vph_2)$ and $({\got g}_0,\vph_1)$  have the same Newton polygon. Contradiction. This ends the proof of the Proposition.
\medskip
{\bf 3.4.9. Corollary.} {\it If $g_2\in\tilde G_0(W(k))$ is such that the Shimura Lie $\sg$-crystals
$({\got g}_0,g_2\vph_1)$ and $({\got g}_0,\vph_1)$ have the same Newton polygon, then $(M,g_2\vph_1)$ and $(M,\vph_1)$ have the same Newton polygon and $g_2\vph_1$ is of the form $g_3p_1\vph_1 g_3^{-1}$, with $p_1\in P_{10}(W(k))$ ($P_{10}$ is as in the proof of 3.4.8) and
with $g_3\in\tilde G_0(W(k))$ such that mod $p$ it normalizes $F^1/pF^1$.}
\medskip
{\bf 3.4.10. Corollary.} {\it If $h_0\in\tilde G_0(W(k))$, then the
multiplicity of the slope $\dl$ for $({\got g}_0,h_0\vph_1)$ is less or equal to $m_{\dl}$. If we have equality, then
$({\got g}_0,h_0\vph_1)$ and $({\got g}_0,\vph_1)$ have the same Newton polygon.}
\medskip
3.4.9-10 are consequences of 3.4.6 and of the proof of 3.4.8. 3.4.11-12 below result from  3.4.6-10, by working with all cyclic permutations
$\ga_s$, $s\in\Mj_0$, at once (i.e. we work with the direct sum of $\sg$-linear maps $\oplus_{s\in\Mj_0} \Psi_j^s$, where $\Psi_j^s$ is defined as $\Psi_j^0$ of 3.4.5 but for an arbitrary $s\in\Mj_0$). 
\medskip
{\bf 3.4.11. Corollary.} {\it If $g_2\in G(W(k))$ is such that the Shimura Lie $\sg$-crystals
$({\got g},g_2\vph_1)$ and $({\got g},\vph_1)$ have the same Newton polygon, then the $\sg$-crystals $(M,g_2\vph_1)$ and $(M,\vph_1)$ have the same Newton polygon and $g_2\in D_G^1$ (see 3.2.10). Moreover the fact that $({\got g},g_2\vph_1)$ has the same Newton polygon as $({\got g},\vph_1)$ depends only on the expression of $g_2$ mod $p$.} 
\medskip
{\bf 3.4.12. Corollary.} {\it The Newton polygon of the Shimura Lie $\sg$-crystal $({\got g},\vph_1)$ is not strictly above the Newton polygon of any other Shimura Lie $\sg$-crystal of the form $({\got g},g_2\vph_1)$, with $g_2\in G(W(k))$.}
\medskip
{\bf 3.4.13. Remark.} 3.4.11-12 remain true if we replace $\vph_1$ and $D_G^1$ by $\vph_0$ and respectively by $D_G$ and we do not assume $G$ is split (assumption made at the beginning of 3.4). This is an easy consequence of 3.3.1, of the translation part of 3.3.3 and of the fact that two parabolic subgroups of $G$ which over $W(\bar k)$ are $G(W(\bar k))$-conjugate are in fact $G(W(k))$-conjugate (this fact results from [Bo2, 20.9] applied over $B(k)$ and Iwasawa's decomposition of [Ti2, 3.3.2]). So if $g_2\in G(W(k))$ is such that $g_2(\vph\otimes 1)=g_3p_0(\vph_0\otimes 1)g_3^{-1}$, with $p_0\in P_0(W(\bar k))$ and $g_3\in G(W(\bar k))$ normalizing $F^1\otimes_{W(k)} \bar k$, as every parabolic subgroup is its own normalizer, we get that we can assume $g_3\in G(W(k))$; this implies $p_0\in G(W(k))\cap P_0(W(\bar k))=P_0(W(k))$. This proves 3.3.3 and ends the proof of 3.2.5.
\medskip
{\bf 3.4.14. Exercise.} Using just the ideas and results of 3.1.8.1 and 3.4.6-13 prove the whole of 3.1.0 (i.e. end the first proof of 3.2.6-7). 
\medskip
This Exercise ends  (cf. 3.2.4) the first proof of 3.1.0. For a solution of it, which can be read out at any time, see 3.7.6.
\medskip\smallskip
{\bf 3.5. The proof of 3.4.6.} 
With the notations of 3.4.4-5 it is enough to show that
$$\dim_{W(k)} (N_2)\le\dim_{W(k)}(\bigcap_{i\in\tilde I_1}{\got r}_{\vep_i})={m_{\dl}\over n},\leqno (0)$$ 
and, if equality holds, then (SUB)
is true, with $i=1$.
\smallskip
For $i\in\tilde I_0$, let $\bar{\got g}_i:={\got g}_i\otimes_{W(k)} k$. 
For $i\in\tilde I_1$, let $\Pi_{\vep_i}$ be the projector of $\bar{\got g}_1$ of whose
kernel is $\bar{\got p}_{\vep_i}:={\got p}_{\vep_i}\otimes_{W(k)} k$ and of whose image is
$\bar{\got r}_{\vep_i}:={\got r}_{\vep_i}\otimes_{W(k)} k$; also we define $\bar{\got s}_{\vep_i}:=
{\got s}_{\vep_i}\otimes_{W(k)} k$ and $\bar{\got q}_{\vep_i}:={\got q}_{\vep_i}\otimes_{W(k)} k$.
\medskip
Following the image of $\bar{\got g}_1$ under the $s$-th iterate of $\bar\psi^0_2$, and
identifying $\bar{\got g}_i$ with $\bar{\got g}_i\otimes_k k$ (the homomorphism $k\to k$ being the Frobenius automorphism $\sg$) with $\bar{\got g}_1$ (this is allowed as we are dealing here just with images), we come across $k$-vector spaces of the form
$$^{h_s}\Pi_{\vep_{i_s}}\bigl(^{h_{s-1}}\Pi_{\vep_{i_{s-1}}}\bigl(\cdots
\bigl(^{h_2}\Pi_{\vep_{i_2}}(^{h_1}\bar{\got r}_{\vep_1})\cdots\bigr)\bigr)\bigr),\leqno (1)$$
where $h_1,\ldots,h_s\in G_1^{\rm ad}(k)$. Here and in what follows the left upper index $^{h_l}$, $l=\overline{1,s}$, means inner conjugation by $h_l$: it is used for conjugating groups as well as Lie subalgebras or elements of ${\got g}_1$; except for elements we do not use parentheses, i.e. we write $^{h_1}\bar{\got r}_{\vep_1}$ instead of $^{h_1}(\bar{\got r}_{\vep_1})$, etc. Also, in what follows, to simplify the presentation, we consider $k=\bar k$: the statements of 3.4.6 are such that this is allowed. So the homomorphism $G_1(k)\to G^{\rm ad}_1(k)$ is surjective; accordingly, we always assume $h_1,\ldots,h_s\in G_1(k)$.
\medskip
{\bf 3.5.0. Computation of $N_2\otimes_{W(k)} k$.} Taking $s\in\NN$ big enough and a multiple of $\abs{\tilde I_1}$, the $k$-vector space of (1) becomes nothing else but $N_2\otimes_{W(k)} k$ (cf. 3.4.5.1 B). We allow such an $s$ to vary. We  first prove (0), i.e.  that the dimension of the $k$-vector spaces of (1) computing $\dim_k (N_2\otimes_{W(k)} k)$ is less or equal to $\dim_{W(k)} (\bigcap_{i\in\tilde  I_1}
{\got r}_{\vep_i})$.
\medskip
{\bf 3.5.1. The set $A$.} Let 
$$A:=\{\vep_i|i\in\tilde I_1\}.$$ 
So $\cap_{i\in\tilde I_1} {\got r}_{\vep_i}=\cap_{a\in A} {\got r}_a$. The set $A$ can have more than $1$ element only when ${\got g}_1$ is in case a), d) or e6).
In case d) $A$ has at most 3 elements. In case e6) $A$ has at most 2 elements. In case a)
$A$ can have even more than 3 elements but in this case a), regardless of the value $\abs{A}$ we have
$$
\bigcap_{a\in A}{\got r}_a={\got r}_{a_0}\cap{\got r}_{a_1}, 
$$
with $a_0$ (resp. $a_1$) as the smallest (resp. as the greatest) element of $A$. This is so because any positive root containing 
$\al_{a_0}$ and $\al_{a_1}$ in its expression as sums of elements of $\Dl_1$, automatically contains
$\al_{a_2}$ in its expression, for any $a_2\in S(a_0,a_1)$ (see [Bou2, planche I]).
\medskip
{\bf 3.5.2. Notations.} Let $\bar G_1:={G_1}_k$, $\bar T_1:={T_1}_k$ and $\bar T_{1\rm ad}:={T_{1\rm ad}}_k$ (see 3.4.3). We call an
element $w\in\bar G_1(k)$ normalizing $\bar T_1$, a Weyl element. For any 
$\vep_i\in A$, let $\bar Q_{\vep_i}$ (resp. $\bar P_{\vep_i})$ be the parabolic subgroup of $\bar G_1$ having $\bar{\got q}_{\vep_i}$ (resp. $\bar{\got p}_{\vep_i})$ as its Lie algebra. It contains $\bar T_1$ as it can be checked immediately starting from the Fact of 2.2.11.1.
\medskip
{\bf 3.5.3. The first properties.} We study closer a $k$-linear endomorphism of $\bar{\got g}_1$ of the form
$\Pi_a\,^g\Pi_b$, with $a,b\in A$. We have the following main property:
$${\rm Im}(\Pi_a\,^g\Pi_b)={\rm Im}(\,^{\tilde g}\Pi_a\,^{\tilde w}\Pi_b),\ \ {\rm with}\ \ \tilde g\in\bar G_1(k)\ \ {\rm and}\,\ {\rm with}\ \ \tilde w\ \ {\rm a\, Weyl\,\, element}.\leqno (2)$$
\indent
We  also have the following inversion property:
$$\Pi_a\,^w\Pi_b=\,^w\Pi_b\,^w\Pi_a\,^{w^{-1}},\leqno (3)$$
for any Weyl element $w\in\bar G_1(k)$.
\smallskip
For any $g\in\bar G_1(k)$ we have the following elimination property
$${\rm Im}(\Pi_a\,^g\Pi_a)\subset {\rm Im}(\Pi_a).\leqno (4)$$
\indent
{\bf Proofs:} The elimination property is trivial. The inversion property can be checked easily, using Weyl's decomposition of $\bar{\got g}_1$ w.r.t. $\bar T_{1\rm ad}$. The proof of (2) is more involved: it needs some preliminaries.
\smallskip
Let $w_0\in\bar G_1(k)$ be a Weyl element that takes $\Phi_1^+$ into $\Phi_1^-$ (i.e. which takes the Borel subgroup of $\bar G_1$ having ${\got b}_1\otimes_{W(k)} k$ as its Lie algebra into its opposite $\bar B_1^{\rm opp}$ w.r.t. $\bar T_1$). Let $\Pi_{w_0^{-1}(a)}$ be the projector of $\bar{\got g}_1$ of whose kernel is $w_0^{-1}(\bar{\got p}_a)$ and of whose image is $w_0^{-1}(\bar{\got r}_a)$. So we have $$\Pi_a=\,^{w_0}\Pi_{w_0^{-1}(a)}\,^{w_0^{-1}}.$$
The Bruhat decomposition of $\bar G_1$ (see [Bo2, 14.12]) allows us to write
$w_0^{-1}g=b_awq_b$, with $q_b\in\bar Q_b(k)$, with $w\in\bar G_1(k)$ a Weyl element and with $b_a\in\bar B_1^{\rm opp}(k)$. ${\rm Lie}(\bar B_1^{\rm opp})$ is included in the kernel of $\Pi_{w_0^{-1}(a)}$. As $\bar B_1^{\rm opp}$ is the semidirect product of $\bar T_1$ and of its unipotent radical, using [Bo2, 14.5 (2)] applied to this unipotent radical (and $\bar T_1$), we can write 
$$b_a=b_a^1b_a^2$$
where $b_a^1\in\bar B_1^{\rm opp}(k)$ normalizes the image of $\Pi_{w_0^{-1}(a)}$ and $b_a^2$ is a $k$-valued point of the smooth, connected, unipotent subgroup of $\bar B_1^{\rm opp}$ generated by $\GG_a$ subgroups normalized by $\bar T_1$ and whose Lie algebra is the opposite of the image of $\Pi_{w_0^{-1}(a)}$ w.r.t. the action of $\bar T_{1\rm ad}$. Obviously
$$\Pi_{w_0^{-1}(a)}\,^{b_a^1}=\,^{b_a^1}\Pi_{w_0^{-1}(a)}.$$
Also we have  
$$\Pi_{w_0^{-1}(a)}\,^{b_a^2}=\Pi_{w_0^{-1}(a)}.$$ 
In other words, for any $x_a\in w_0^{-1}(\bar{\got r}_a)$, we have 
$$^{b_a^2}(x_a)=x_a+y_a,\leqno (5)$$ 
where $y_a\in w_0^{-1}(\bar{\got p}_a)$. To check this, we can assume (cf. loc. cit.) that $b_a^2$ is a $k$-valued point of a mentioned $\GG_a$ subgroup as well as $x_a\in {\got b}_1\otimes_{W(k)} k$ is normalized by $\bar T_1$. So formula (5) is an immediate consequence of the formulas of [BT, 4.2]. Besides loc. cit., we just need to remark that this formula (5) is trivial if $\bar G_1$ is of $A_1$ Lie type; the reason we have to remark this is: loc. cit. deals only with pairs of roots whose sum is non-zero.   
\smallskip
We conclude 
$$\Pi_{w_0^{-1}(a)}\,^{b_a}=\,^{b_a^1}\Pi_{w_0^{-1}(a)}.$$
\indent
To compute the image of $\Pi_a\,^g\Pi_b$, we can assume $q_b$ is the identity elements. We have:
$${\rm Im}(\Pi_a\,^g\Pi_b)={\rm Im}(\,^{w_0}\Pi_{w_0^{-1}(a)}\,^{w_0^{-1}g}\Pi_b)={\rm Im}(\,^{w_0}\Pi_{w_0^{-1}(a)}\,^{b_aw}\Pi_b)={\rm Im}(\,^{w_0b_a^1w_0^{-1}w_0}\Pi_{w_0^{-1}(a)}\,^{w}\Pi_b).$$
But this last ${\rm Im}$ can be rewritten as ${\rm Im}(\,^{w_0b_a^1w_0^{-1}}\Pi_a\,^{w_0w}\Pi_b).$
So we can take $\tilde w:=w_0w$ and $\tilde g:=w_0b_a^1w_0^{-1}$. This proves (2).
\medskip
{\bf 3.5.3.1. The proof of inequality (0).} We are now ready to prove, in four steps, the first statement of 3.4.6 (i.e. to prove 3.5 (0)). We can assume that any two consecutive $\vep_{i_j}$ and $\vep_{i_{j+1}}$ of 3.5 (1) are distinct (cf 3.5.3 (4)).
\medskip
{\bf 1)} If $\abs{A}=1$, we have nothing to prove (i.e. 3.5 (0) is a consequence of the Fact of 3.4.5.1 and of 3.5.0).
\smallskip
{\bf 2)} If $\abs{A}=2$, then in any sequence $\,^{h_s}\Pi_{\vep_{i_s}}...\,^{h_1}\Pi_{\vep_1}$ ``computing" $N_2\otimes_{W(k)} k$,
we get a subsequence $\Pi_a\,^g\Pi_b$, with $A=\{a,b\}$ and $g\in\bar G_1(k)$. So, we get that $\dim_{W(k)}(N_2)\le\dim_{W(k)}({\got r}_a\cap {\got r}_b)$, provided the inequality
$$\dim_k (\Pi_a(^g\bar{\got r}_b))\le\dim_k(\bar{\got r}_a\cap\bar{\got r}_b)
$$ 
is checked.
This last inequality can be proved in many ways. Using property (2), we are reduced to the case when $g$ is a Weyl element and so 
$$\dim_k(\Pi_a(^g\bar{\got r}_b))=\dim_k(\bar{\got r}_b)-\dim_k(^g\bar{\got r}_b\cap\bar{\got p}_a)=\dim_k(\bar{\got r}_b)-\dim_k(\bar{\got p}_a)+\dim_k(^g\bar{\got p}_b\cap\bar{\got p}_a)$$
 (the second equality due to the
fact that $g$ is a Weyl element). This proves the first part of the following Claim:
\medskip
{\bf Claim.} {\it The dimension of $\Pi_a(^g\bar{\got r}_b)$ is
maximal precisely when the dimension of the intersection $^g\bar{\got p}_b\cap\bar{\got p}_a$ of parabolic Lie subalgebras of $\bar{\got g}_1$ is maximal, i.e. when the intersection $^g\bar{\got p}_b\cap\bar{\got p}_a$ contains a Borel Lie subalgebra of $\bar{\got g}_1$.} 
\medskip
{\bf Proof:} The Fact of 2.2.11.2 points out that it is irrelevant if we work with inclusions between parabolic subgroups of $\bar G_1$ or between inclusions of their Lie algebras: this will be used freely in what follows. If the intersection $^g\bar{\got p}_b\cap\bar{\got p}_a$ contains a Borel Lie subalgebra of $\bar{\got g}_1$, then this intersection has the same dimension as the intersection $\bar{\got p}_b\cap\bar{\got p}_a$. Argument: we can assume ${\rm Lie}(\bar B_1)\subset ^g\bar{\got p}_b\cap\bar{\got p}_a$ and so the statement follows from [Bo2, 14.22 (iii)]. So if $\dim_k(^g\bar{\got p}_b\cap\bar{\got p}_a)< \dim_k(\bar{\got p}_b\cap\bar{\got p}_a)$, then $^g\bar{\got p}_b\cap\bar{\got p}_a$ does not contain the Lie algebra of a Borel subgroup of $\bar G_1$. Based on this, the part of the Claim pertaining to ``i.e." can be checked easily, case by case, starting from the fact that $g$ is a Weyl element and from [Bou2, planche I, IV and V].  
\smallskip
The case $A=\{1,\ell-1\}$, with ${\got g}_1$ of $D_{\ell}$ Lie type, is treated for the sake of convenience in the fourth paragraph of 4) below. The case e6) is left as an exercise: though in [Bou2, planche V] the structure of the Weyl group associated to the $E_6$ Lie type is not given, this exercise can be solved in the same manner as the two cases below, by working not in terms of vectors $\vep_i$'s (of loc. cit.) but directly in terms of roots; however see [Bou2, p. 219 and Exc. 2) of \S 4 of ch. VI] (also the case e6) can be checked immediately using computers). We are left with the following two cases. 
\medskip
{\bf Case 1: {\got g} is of $A_{\ell}$ Lie type.} So $A=\{i_0,i_1\}$ with $1\le i_0< i_1\le\ell+1$. We can assume $b=i_0$ and $a=i_1$, cf. 3.5.3 (3). For this paragraph we use the notations of [Bou2, planche I]: so $\vep_i$'s are not any more some integers in the set $S(0,\ell)$ but some vectors in $\RR^{\ell+1}$. The roots $\be\in\Phi_1^-$ such that ${\got g}_{\be}\subset {\got r}_c$, with $c\in A$, are of the form $-\vep_i+\vep_s$, with $i\le c< s\le\ell+1$. We consider a Weyl element $h\in \bar G_1(k)$. By abuse of notation, we identify it with an element of the Weyl group of loc. cit.; so we identify it with a permutation $h$ of the set $S(1,\ell+1)$ via the formula $h(e_i)=e_{h(i)}$. 
\smallskip
So the intersection $^h\bar{\got r}_{i_0}\cap\bar{\got r}_{i_1}$ has as many elements as the number of pairs $(i,j)$ such that $i\le i_0<j$ and $h(i)\le i_1\le h(j)$. Obviously the number of such pairs is of the form $\tilde i_0(\ell+1-\tilde i_1)$, where $\tilde i_0\in S(1,i_0)$ and $\tilde i_1\in S(i_1,\ell+1)$, and so it is at most $i_0(\ell+1-i_1)$. The equality holds precisely when $h$ maps the set $S(1,i_0)$ into the set $S(1,i_1)$ and maps a subset of the set $S(i_0+1,\ell+1)$ onto the set $S(i_1+1,\ell+1)$. But it is easy to see that any such permutation can be written as a product $h=h_1h_2$, where $h_2$ (resp. $h_1$) is the permutation of $S(1,\ell+1)$ leaving invariant each one of the sets $S(1,i_0)$ and $S(i_0+1,\ell+1)$ (resp. each one of the sets $S(1,i_1)$ and $S(i_1+1,\ell+1)$). As any Weyl element defining $h_2$ (resp. $h_1$) normalizes ${\got p}_{i_0}$ (resp. ${\got p}_{i_1}$), $^h\bar{\got p}_{i_0}\cap\bar{\got p}_{i_1}$ contains the Lie algebra of a Borel subgroup of $\bar G_1$ and we are through.
\medskip
{\bf Case 2: $A=\{\ell-1,\ell\}$, with ${\got g}_1$ of $D_{\ell}$ Lie type.} For this paragraph we use the notations of [Bou2, planche IV]. We have:
\medskip
--  {\it the roots of $\Phi_1^-$ corresponding to the action (via inner conjugation) of  $\bar T_{1\rm ad}$ on $\bar{\got r}_{\ell-1}$ are either of the form $-\vep_i+\vep_\ell$, with $1\le i<\ell$, or of the form $-(\vep_i+\vep_j)$, with $1\le i<j<\ell$;}
\smallskip
-- {\it the roots of $\Phi_1^-$ corresponding to the action (via inner conjugation) of $\bar T_{1\rm ad}$ on $\bar{\got r}_{\ell}$ are of the form $-(\vep_i+\vep_j)$, with $1\le i<j\le\ell$.}
\medskip
We assume $b=\ell$.
We consider a Weyl element $h\in \tilde G_1(k)$; we identify it with
an element of the Weyl group $W_{D_{\ell}}$ of loc. cit. If $h$ sends $\vep_{\ell}$ into $\vep_{\ell}$, then the intersection $^h\bar{\got r}_{b}\cap\bar{\got r}_a$ has at most ${{(\ell-2)(\ell-1)}\over 2}$ elements, with equality iff $h$ leaves invariant the set $\{\vep_1,...,\vep_{\ell-1}\}$; in such a case $h$ normalizes both $\bar{\got p}_{a}$ and $\bar{\got p}_{b}$. 
\smallskip
If $h$ sends $\vep_{\ell}$ into $-\vep_{\ell}$, then the intersection $^h\bar{\got r}_{b}\cap\bar{\got r}_a$ has $s+{{(s-1)s}\over 2}$ elements, where $s$ is the number of elements $i\in S(1,\ell-1)$ such that $h(\vep_i)\in\{\vep_1,...,\vep_{\ell-1}\}$; as $s\le\ell-2$ (cf. the structure of $W_{D_{\ell}}$ as a semidirect product mentioned in loc. cit.), this number is less or equal to ${{(\ell-2)(\ell-1)}\over 2}$, with equality iff $s=\ell-2$. If $s=\ell-2$, let $h_2\in W_{D_{\ell}}$ be the involution $(\vep_{\ell},\vep_{i_0})$, where $i_0\in S(1,\ell-1)$ is the unique element such that $h(\vep_{i_0})\in\{-\vep_1,...,-\vep_{\ell-1}\}$ (so $h(\vep_{i_0})=-\vep_{i_0}$). We have: $h_2$ (resp. $hh_2^{-1}$) normalizes $\bar{\got p}_{b}$ (resp. $\bar{\got p}_{a}$).
\smallskip
If $h$ sends $\vep_{\ell}$ into an element of the set $\{-\vep_1,...,-\vep_{\ell-1}\}$, then the intersection $^h\bar{\got r}_{b}\cap\bar{\got r}_a$ has $s+{{s(s-1)}\over 2}$ elements, where $s$ is the number of elements of the set $\{\vep_1,...,\vep_{\ell-1}\}$ mapped through $h$ into this set. As $s\le\ell-2$, this number is less or equal to ${{(\ell-2)(\ell-1)}\over 2}$. If $s=\ell-2$, then $h$ normalizes $\bar{\got p}_{a}$. If $h$ sends $\vep_{\ell}$ into $\vep_j$, with $j\in S(1,\ell-1)$, then the intersection $^h\bar{\got r}_{b}\cap\bar{\got r}_a$ has $s+{{(s-1)s}\over 2}$ elements, where $s$ is the number of $i\in S(1,\ell-1)$ such that $h(\vep_i)\in\{\vep_1,...,\vep_{j-1},\vep_{j+1},...,\vep_{\ell-1}\}$; as $s\le\ell-2$, this number is less or equal to ${{(\ell-2)(\ell-1)}\over 2}$, with equality iff $s=\ell-2$. If $s=\ell-2$, then $h$ normalizes $\bar{\got p}_{b}$. 
\smallskip
So in all cases $^h\bar{\got p}_{b}\cap\bar{\got p}_a$ contains the Lie algebra of a Borel subgroup of $\bar G_1$ and we are through.
\smallskip
This proves the Claim.
\medskip
From Claim and the first paragraph of its proof we get: $\dim_k(\,^g\bar{\got p}_b\cap\bar{\got p}_a)\le\dim_k(\bar{\got p}_b\cap\bar{\got p}_a)$.   
\smallskip
{\bf 3)} If $\abs{A}\ge 3$, then using repeatedly properties (2) and (3), we can bring
together a short sequence of the form $\Pi_a\,^g\Pi_b$, if ${\got g}_1$ is in case a), with $a$ the smallest element of $A$ and with $b$ the greatest element of $A$, or of the form
$\Pi_a\,^g\Pi_b\,^h\Pi_c$  if ${\got g}_1$ is in case d), with $A=\{a,b,c\}$. This proves 3.5 (0) if ${\got g}_1$ is in case a) (cf. 2)).
\smallskip
{\bf 4)} We assume we are in case d) with $\abs{A}=3$.
\smallskip
Let ${\got g}_1$ be of $D_\ell$ Lie type ($\ell\ge 4$). So $A=\{1,\ell-1,\ell\}$. Using property (2), we can assume $\Pi_{a}\,^g\Pi_{b}\,^h\Pi_{c}$ is such that $h$ is a Weyl element. Due to the circular form of the expressions of 3.5 (1), we can assume $c=1$ and $b=\ell-1$. But in this case the dimension of the image of $\Pi_{a}\,^g\Pi_{b}\,^h\Pi_{c}$ is less or equal to 
$$\ell-2=\dim_k(\bar r_a\cap\bar r_b\cap\bar r_c).$$
\indent
To check this, we use the notations of loc. cit. We have:
\medskip
--  {\it the roots of $\Phi_1^-$ corresponding to the action (via inner conjugation) of $\bar T_{1\rm ad}$ on $\bar{\got r}_1$ are of the form $-\vep_1\pm\vep_j$, with $2\le j\le\ell$.}
\medskip
So if the element $h$
of the Weyl group of loc. cit. sends $\vep_1$ into an element of the set $\{-\vep_1,-\vep_2,\ldots,-\vep_{\ell-1},\vep_\ell\}$ then we have
$\bar{\got r}_{\ell-1}\cap\,^h\bar{\got r}_1=\{0\}$ and we are done. If
$h$ sends $\vep_1$ into an element $\vep_{i_0}\in\{\vep_1,\vep_2,\ldots\vep_{\ell-1},-\vep_\ell\}$,
then $\bar{\got r}_{\ell-1}\cap\,^h\bar{\got r}_1$ has dimension exactly $\ell-1$. If $\dim_k(\bar{\got r}_{\ell-1}\cap\,^h\bar{\got r}_1)=\ell-1$, then $h=h_1h_2$ where $h_1$ is the transposition $(\vep_1,\vep_{i_0})$ and where $h_2:=h_1^{-1}h$ fixes $\vep_1$. As $h_2$ (resp. $h_1$) normalizes $\bar{\got r}_1$ and $\bar{\got p}_1$ (resp. $\bar{\got r}_{\ell-1}$ and $\bar{\got p}_{\ell-1}$) the intersection $\bar{\got q}_{\ell-1}\cap\,^h\bar{\got q}_1$ contains a Borel subalgebra of $\bar{\got g}_1$. 
\smallskip
The arguments of the proof of (2) show that we can assume $g$ is a Weyl element as well. So we can apply the same arguments as in the proof of Case 2 of the Claim of 2), in order to conclude that $\dim_k\bigl(\bar{\got r}_\ell\cap\,^g(\bar{\got r}_{\ell-1}\cap\,^h\bar{\got r}_1)\bigr)$ is at most $\ell-2$. Briefly, this goes as follows. We can assume $h$ is the identity element. If $g$ sends $\vep_1$ into $\{-\vep_1,...,-\vep_{\ell}\}$, then $\bar{\got r}_\ell\cap\,^g(\bar{\got r}_{\ell-1}\cap\bar{\got r}_1)=\{0\}$. If $g$ sends $\vep_1$ into $\{\vep_1,...,\vep_{\ell}\}$, then $\dim_k\bigl(\bar{\got r}_\ell\cap\,^g(\bar{\got r}_{\ell-1}\cap\bar{\got r}_1)\bigr)$ is equal to the number of elements of the set $\{-\vep_2,...,-\vep_{\ell-1},\vep_{\ell}\}$ sent by $g$ into the set $\{-\vep_1,...,-\vep_{\ell-1},-\vep_{\ell}\}$. Based on the structure of the Weyl group $W_{D_{\ell}}$ we get that the number of such elements is at most $\ell-2$. If it is $\ell-2$, we distinguish two situations. If $\vep_1$ is sent by $g$ into $\vep_{\ell}$, then $g$ normalizes $\bar{\got p}_{\ell}$. If $\vep_1$ is sent by $g$ into $\vep_{s}$, with $s\in S(1,\ell-1)$, we write $g$ as $g_1g_2$, with $g_2$ as the Weyl element which leaves invariant each $\vep_i$, $i\in S(1,\ell-1)\setminus\{s_0\}$ and which sends $\vep_{\ell}$ into the same element $-\vep_{s_0}$ as $g$ does. $g_1$ (resp. $g_2$) normalizes $\bar{\got p}_{\ell}$ (resp. $\bar{\got p}_{\ell-1}\cap \bar{\got p}_1$).  
\smallskip
This ends the proof of 3.5 (0) and so of the first statement of 3.4.6.
\medskip
{\bf 3.5.4. The proof of (SUB).} We assume now $\dim_{W(k)}(N_2)={m_{\dl}\over n}$. To prove (SUB) we first remark: 
\medskip
{\bf Fact.} {\it The Lie subalgebra $\bigl(g_2(-\dl)\cap{\got g}_1\bigr)\otimes_{W(k)} k$ of $\bar{\got g}_1$ is of the form
$^g\bigl(\cap_{a\in A}\bar{\got s}_a\bigr)$, with $g\in\bar P_{\vep_1}(k)$, and so its 
normalizer in $\bar G_1$ is $^g\bigl(\bigcap_{a\in A}\bar P_a\bigr)\subset\bar P_{\vep_1}$.}
\medskip
{\bf Proof:}
We start pointing out: if the Lie subalgebra $\bigl(g_2(-\dl)\cap{\got g}_1\bigr)\otimes_{W(k)} k$ of $\bar{\got g}_1$ is of the form $^g\bigl(\cap_{a\in A}\bar{\got s}_a\bigr)$, with $g\in \bar G_1(k)$, and if its normalizer in $\bar G_1$ is $^g\bigl(\bigcap_{a\in A}\bar P_a\bigr)$, then automatically $g\in \bar P_{\vep_1}(k)$. To see this, we first remark that $\bigl(g_2(-\dl)\cap{\got g}_1\bigr)\otimes_{W(k)} k$ is included in $\bar{\got s}_{\vep_1}=F^1(\bar{\got g}_1)$ (cf. the Fact of 3.4.5.1) and so the Lie algebra $F^1(\bar{\got g}_1)$ being abelian is included in the normalizer of $\bigl(g_2(-\dl)\cap{\got g}_1\bigr)\otimes_{W(k)} k$ in $\bar{\got g}_1$. Using Fact of 2.2.11.1 and [Bo2, 14.22 (i)] we get that $^g\bigl(\bigcap_{a\in A}\bar P_a\bigr)$ contains a Borel subgroup of $\bar P_{\vep_1}$; using $\bar P_{\vep_1}(k)$-conjugates, we can assume this Borel subgroup is also included in $\cap_{a\in A}\bar P_a$ and so [Bo2, 14.22 (iii)] implies $\cap_{a\in A}\bar P_a=^g\bigl(\cap_{a\in A}\bar P_a\bigr)$. As any parabolic subgroup of $\bar G_1$ is its own normalizer, we conclude: $g\in\bar P_{\vep_1}$. So we are left to check:
\medskip
{\bf Claim.} {\it The subspace $\bigl(g_2(-\dl)\cap{\got g}_1\bigr)\otimes_{W(k)} k$ of $\bar{\got g}_1$ is of the form
$^g\bigl(\cap_{a\in A}\bar{\got s}_a\bigr)$ and its normalizer in $\bar G_1$ is $^g\bigl(\bigcap_{a\in A}\bar P_a\bigr)$.}
\medskip
To prove this Claim we need some preliminaries. Let $\bar{\got g}_{\al}:={\got g}_{\al}\otimes_{W(k)} k$. We recall the following three well known properties:
\medskip
{\bf a)} {\it $[\bar{\got g}_{\al},\bar{\got g}_{-\al}]$ is a $1$ dimensional $k$-vector subspace of ${\rm Lie}(\bar T_1)={\got t}_1/p{\got t}_1$ (this is the first place where we need $p>2$);}
\medskip 
{\bf b)} {\it If $\abs{A}\ge 2$ and if $\al,\be,\al+\be\in\Phi_1$, then $\al$ and $\be$ generate a root system of rank 2 of $A_2$ Lie type (the inequality $\abs{A}\ge 2$ implies that all roots of $\Phi_1$ have the same length, cf. [Bou2, planche I, IV, and V]);}
\medskip
{\bf c)} {\it The differential map of the commutator map $c_a:\bar G_1\to\bar G_1$, defined by $c_a(g)=gag^{-1}a^{-1}$ is given by the formula $dc_a=1_{{\got g}_1}-Ad(a)$ (for instance, cf. [Bo2, 3.16]).}
\medskip
If $\abs{A}\ge 2$, from [BT, 4.2-3] and b) we get (via c)) the following formula
$$
[\bar{\got g}_\al,\bar{\got g}_\be]=\bar{\got g}_{\al+\be},\leqno (6) 
$$
where $\al,\be,\al+\be\in\Phi_1$. It can be as well checked directly, inside the $SL$-group of a 3 dimensional $k$-vector space; it is true even if $p=2$ (for an arbitrary simple, adjoint group over $k$). Also, regardless of the value of $\abs{A}$, from [BT, 4.2-3] we get (again via c)) the following inclusion
$$
[\bar{\got g}_\al,\bar{\got g}_\be]\subset \bar{\got g}_{\al+\be},\leqno (7) 
$$
where $\al,\be\in\Phi_1$ are such that $\al+\be\neq 0$; if ${\al+\be}\notin \Phi_1\cup\{0\}$ then $\bar{\got g}_{\al+\be}:=\{0\}$.
\smallskip
We are now ready to prove the Claim.
We first consider the case when $A=\{a\}$ has only 1 element. Let 
$$\Lambda_1:=\{\al\in\Phi_1^-|{\got g}_{\al}\subset {\got r}_a\}.$$ 
In this case, the first part of the Claim is a consequence of the Fact of 3.4.5.1, while for the second part we just need to remark: if $x_a\in\oplus_{\al\in\Lambda_1} \bar{\got g}_{\al}$ normalizes $\bar{\got s}_a$, then $\forall\al\in\Lambda_1$, from the inclusion $[\bar{\got g}_{-\al},x_a]\subset\bar{\got s}_a$, we get, cf. a) and (7), that the component of $x_a$ in $\bar{\got g}_{\al}$ is $0$; so $x_a=0$. 
\smallskip
If $\abs{A}\ge 2$, in order to benefit from the previous notations, we prove this Claim working with notations involving negative roots: the reductions to Weyl elements performed in 3.5.3 allows us to shift from positive to negative roots, as we like. We need to consider three cases.
\medskip
{\bf Case 1.} If $\abs{A}=2$, with $A=\{a,b\}$, we just have to show: 
\medskip
{\bf Subfact.} {\it If the image of the endomorphism 
$\Pi_a\,^w\Pi_b$, with $w\in\bar G_1(k)$ a Weyl element, has  dimension
$\dim_k(\bar{\got r}_a\cap\bar{\got r}_b)$, then $\bar{\got r}_b\cap\,^{w^{-1}}\bar{\got r}_a$
is $\bar G_1(k)$-conjugate to $\bar{\got r}_a\cap\bar{\got r}_b$, and the normalizer of $\bar{\got r}_a\cap\bar{\got r}_b$ is 
$\bar Q_a\cap\bar Q_b$.} 
\medskip
The first part of the Subfact is easy. Argument: $\bar Q_a$ and $^w\bar Q_b$ must contain a Borel subgroup of $\bar G_1$ (cf. Claim of 3.5.3.1 2) and Fact of 2.2.11.1); as any two Borel subgroups of $\bar G_1$ are $\bar G_1(k)$-conjugate, we can assume this Borel subgroup is $\bar B_1^{\rm opp}$ and so, based on [Bo2, 14.22 (iii)], we can assume $^w\bar Q_b=\bar Q_b$. For the second part of it, as any parabolic subgroup of a reductive group over $k$ is connected, it is enough to show that the Lie subalgebra of $\bar{\got g}_1$ normalizing $\bar{\got r}_a\cap \bar{\got r}_b$ is ${\rm Lie}(\bar Q_a\cap\bar Q_b)$. To check this, let 
$$\Lambda_2:=\{\al\in\Phi_1^+|{\got g}_{\al}\subset {\got s}_a\cup {\got s}_b\}.$$
Let $x_{ab}\in\oplus_{\al\in\Lambda_2} \bar{\got g}_{\al}$ normalizing $\bar{\got r}_a\cap \bar{\got r}_b$. Let $\al\in\Lambda_2$. Let $\be\in\Phi_1^-$ be such that:
\medskip
i) ${\got g}_{\be}\subset {\got r}_a\cap {\got r}_b$;
\smallskip
ii) either $\be+\al\in\Phi_1$ or $\be=-\al$.
\medskip
The existence of such a root can be read out from [Bou2, planche I, IV and V]. We include a discussion case by case, in a way needed to be used in 3.14 C (in connection to the $p=2$ case). 
\smallskip
If we are in case a) with $a<b$ and $\al=\sum_{s=m}^q \al_s$, where either $m\le a$ and $q\ge a$ or $m\le b$ and $q\ge b$, then we take respectively either $\be=\sum_{s=m}^{\max\{q,b\}} -\al_s$ or $\be=\sum_{s=\min\{m,a\}}^{q} -\al_s$. 
\smallskip
If we are in case e6) with $1=a<b=6$ and $\al=\sum_{s=1}^6 b_s\al_s$, where $b_s\in\{0,1,2,3\}$ and either $b_1$ or $b_6$ is $1$, then we can assume $b_1=1$ and so we take $\be=-\al_1-\al_6-(\sum_{s=2}^5 c_s\al_s)$ where:
\medskip
-- $c_3=c_4=c_5=1$ and $c_2=0$ if $b_2=b_6=0$;
\smallskip
-- $c_i=1$, $i=\overline{2,5}$, if $b_2=1$ and $b_6=0$;
\smallskip
-- $c_i=b_i$, $i=\overline{2,5}$, if $b_6=1$.
\medskip
For an exemplification in 3.14 C below, we also point out that we can also take:
\medskip
-- $c_i=1$, $i=\overline{2,5}$, if $b_2=0$ and $b_6=1$;
\smallskip
-- $c_2=1$ and $c_3=c_4=c_5=2$ if $b_2=b_6=1$ and $b_3$, $b_4$, $b_5\in\{1,2\}$ are such that either only one of them is $1$ or all three of them are $1$;
\smallskip
-- $c_2=c_3=c_5=2$ and $c_4=3$, if $b_6=b_2=1$ and either $b_4=3$ or $b_4=2$ and $b_3=b_5$;
\smallskip
-- $c_2=1$, $c_3=c_4=c_5=2$, if $b_2=2$ and $b_6=1$.
\medskip
If we are in case d), we need to consider two situations. First we assume $A=\{1,\ell-1\}$. If $\al$ does not contain $\al_{\ell}$ in its expression, we take $\be=-\sum_{s=1}^{\ell-1} \al_s$; if $\al=\sum_{s=1}^{\ell-1} \al_s$ then we can also take $\be=-\al-\al_{\ell}$. If $\al$ does not contain $\al_{\ell-1}$ in its expression but it does contain $\al_1$, then we take $\be=-\sum_{s=1}^\ell \al_s$. If $\al=\sum_{s=s_0}^{\ell} \al_s$, with $s_0\ge 2$, or if $\al$ is the maximal root, then we take $\be=-\sum_{s=1}^\ell \al_s$. If $\al$ contains $\al_{\ell}$, $\al_{\ell-1}$ and $\al_1$ in its expression but is not the maximal root, then we take $-\be$ to be this maximal root. If $\al$ contains $\al_{\ell}$, $\al_{\ell-1}$ and $\al_{s_0}$ in its expression, with $s_0\ge 2$, and does not contain $\al_{s_1}$, for any $s_1\in S(1,s_0-1)$, then we take $\be=-\al-\sum_{s=1}^{s_0-1} \al_s$. 
\smallskip
We consider now the case $A=\{\ell-1,\ell\}$. We can assume $\al$ contains $\al_{\ell-1}$ in its expression. If $\al$ does not contain $\al_{\ell}$ in its expression and it is different from $\al_{\ell-1}$, then we take $\be=-\al-\al_{\ell}$. If $\al=\al_{\ell-1}$, then we take $\be=-\sum_{s=\ell-2}^{\ell} \al_s$. If $\al$ contains $\al_{\ell}$ in its expression, then, not to introduce an $s_0$ as above, we can assume $\al$ contains as well $\al_1$ in its expression; so the choice of the previous paragraph of $\be$ in such a situation applies as well in the context of $A=\{\ell-1,\ell\}$.
\medskip
We now come back to $x_{ab}$. From the inclusion $[\bar{\got g}_{\be},x_{ab}]\subset\bar{\got r}_a\cap \bar{\got r}_b$, from formulas (6) and (7), from property a) and from i) and ii) we get that the component of $x_{ab}$ in $\bar{\got g}_{\al}$ is $0$. So $x_{ab}=0$. 
\medskip
{\bf Case 2.} If $\abs{A}\ge 3$ and ${\got g}_1$ is of $A_{\ell}$ Lie type, then  we have to work with $\Pi_a\,^w\Pi_b$,
with $a$ and $b$ the smallest and respectively the greatest element of $A$ and with $w\in\bar G_1(k)$ a Weyl element. So this case gets reduced to the case $\abs{A}=2$.
\medskip
{\bf Case 3.} If $\abs{A}=3$ and ${\got g}_1$ is of $D_{\ell}$ Lie type, then we have to work with $\Pi_a\,^w\Pi_b\,^{w_0}\Pi_c$, where $A=\{a,b,c\}$ and $w,w_0\in G_1(k)$ are both Weyl elements (cf. end of 3.5.3.1 4)). The same argument at the level of Lie algebras as in Case 1 applies: based on 3.5.3.1 4), we just need to show that the normalizer of $\bar{\got r}_a\cap \bar{\got r}_b\cap \bar{\got r}_c$ is $\bar Q_a\cap\bar Q_b\cap\bar Q_c$. To get this we just need to replace everywhere (in Case 1) ${\got r}_a\cap {\got r}_b$ by ${\got r}_a\cap {\got r}_b\cap {\got r}_c$, $\bar{\got r}_a\cap \bar{\got r}_b$ by $\bar{\got r}_a\cap \bar{\got r}_b\cap \bar{\got r}_c$, $\bar Q_a\cap\bar Q_b$ by $\bar Q_a\cap\bar Q_b\cap\bar Q_c$, $\bar s_a\cup\bar s_b$ by $\bar s_a\cup\bar s_b\cup\bar s_c$ and $x_{ab}$ by $x_{abc}$. The existence of a root $\be$ in this case is argued entirely as in the d) case of Case 1. To see this, let $\al\in\Phi_1^+$ be such that ${\got g}_{\al}\in {\got s}_1\cup {\got s}_{\ell-1}\cup {\got s}_{\ell}$. If $\al$ does not contain $\al_{\ell-1}$ and $\al_{\ell}$ in its expression, we take $\be=-\sum_{i=1}^{\ell} \al_i$. If $\al=\sum_{i=s}^{\ell-1}$, with $s\in S(1,\ell-1)$, then we take $\be=-\sum_{i=1}^{\ell} \al_i-\sum_{i=s}^{\ell-2}\al_i$. If $\al$ contains $\al_1$, $\al_{\ell-1}$ and $\al_{\ell}$ in its expression, then we take $\be$ as in Case 1 (for when we are in case d)). If $\al$ contains just $\al_{s}$'s in its expression, with $s\in S(s_0,\ell)$, with $s_0\in S(1,\ell-1)$, then we take $\be=-\al-\sum_{i=1}^{s_0-1}\al_i$. This ends the proof of the Claim and so of the Fact.
\medskip
{\bf 3.5.4.1. End of the proof of 3.4.6.} We are now ready to prove (SUB). 
The parabolic subgroup of $\bar G_1$ having $\bigl(\bigoplus_{\al\in\Mh_1\cap[0,1]}g_0(\al)\cap{\got g}_1\bigr)\otimes_{W(k)} k$
 as its Lie algebra normalizes $\bigl(g_2(-\dl)\cap{\got g}_1\bigr)\otimes_{W(k)} k$. So the Fact of 3.5.4 implies that it is included in  $\bar P_{\vep_1}$. 
This ends the proof of (SUB) and so of 3.4.6.
\medskip
{\bf 3.5.5. Variants.} 
It is easy to see that 3.4.1-8, 3.4.10, the logical variants of 3.4.9 and 3.4.11 (involving, cf. the beginning paragraph of 3.4.0, a generalized Shimura context), as well as the whole of 3.5 remain true if, instead of $g_2\in G(W(k))$ (resp. $h_0\in\tilde G_0(W(k))$ for 3.4.10) we work with any inner automorphism of $G_{W(k)}$ (resp. of ${\tilde G}_0$). Warning: in general this is not so if we use an outer automorphism, as such an automorphism can change the automorphism class of $({\got g},\vph)$, ``producing" a different $G$-ordinary type $\tau$. In particular, based on 3.4.11 and 3.4.13 we get: 
\medskip
{\bf Corollary.} {\it Let $g\in G^{\rm ad}(W(k))$.
\medskip
{\bf 1)} If $({\got g},\vph_0)$ and $({\got g},g\vph_0)$ have the same Newton polygon, then at the level of $\sg$-linear Lie automorphisms of ${\got g}$ we have $g\vph_0=g_3p_4\vph_0 g^{-1}_3$, where $g_3$ is a $W(k)$-valued point of $G^{\rm ad}$ which mod $p$ defines a $k$-valued point of the parabolic subgroup of $G^{\rm ad}$ having $F^0({\got g})[{1\over p}]\cap {\rm Lie}(G^{\rm ad})$ as its Lie algebra, and where $p_4$ is a $W(k)$-valued point of the parabolic subgroup of $G^{\rm ad}$ having ${\rm Lie}(P_0)[{1\over p}]\cap {\rm Lie}(G^{\rm ad})$ as its Lie algebra. Moreover, $({\got g},g\vph_0)$ has a lift of parabolic type.
\smallskip
{\bf 2)} The fact that $({\got g},\vph_0)$ and $({\got g},g\vph_0)$ have the same Newton polygon or not depends only on the expression of $g$ mod $p$.
\smallskip
{\bf 3)} 1) and 2) remain true if instead of ${\got g}$ we work with ${\rm Lie}(G^{\rm ad})$. Moreover, we have variants of these for when we work with ${\got g}_0$ instead of ${\got g}$ or with $\oplus_{i\in I_0} {\rm Lie}(G_i^{\rm ad})$ instead of ${\rm Lie}(G^{\rm ad})$.}
\medskip
In the above Corollary, the role of $({\got g},\vph_0)$ is precisely of an arbitrary Shimura Lie $\sg$-crystal having a lift of Borel type. 
\medskip
{\bf 3.5.6. An expectation.} We expect that the following statement is true: 
\medskip
{\bf Expectation.} {\it Let $H$ be a split, semisimple group over an arbitrary field $\tilde k$. Let $P_H(1)$,..., $P_H(n)$ be $n$ parabolic subgroups of it ($n\in\NN$). Then the dimension of the intersection of the Lie algebras of $n$ parabolic subgroups $\tilde P_1$,..., $\tilde P_n$ of $H$ which are respectively $H(\tilde k)$-conjugate to $P_H(1)$,..., $P_H(n)$, is maximal precisely in the case when there is a Borel subgroup of $H$ contained in the intersection $\cap_{i=1}^n\tilde P_i$.}
\medskip
The case $n=2$ is easy (to be compared with 3.5.3.1 2)). For $n>2$ we can prove it in many situations but we do not have a proof which works in the general case. 
\medskip\smallskip
{\bf 3.6. The second group of basic results: global deformations and some principles for different Fontaine categories.}
We start constructing global deformations: 3.6.0-11 form a sequence; they deal with global deformations in the context of Shimura $\sg$-crystals. However, along this sequence some digressions are as well included; they are 3.6.1.1.1 4) and 5), 3.6.1.1.2, 3.6.1.5-6, 3.6.6.2, 3.6.8.1 till 3.6.8.1.5 inclusive, 3.6.8.3 and 3.6.8.9; these digressions can be looked at just at the needed moments. As mentioned in 1.14.5, the readers should look at 3.6.8 from the general point of view from the very beginning. In 3.6.12-14 we deal with the geometric context of a SHS. 3.6.15-17 contain complements. 3.6.18-19 deal with different principles of Fontaine categories; warning: the global deformations in the generalized Shimura context are constructed in 3.15.6 (cf. also the relative moduli principle of 3.6.18.7.3 C). Some conclusions are gathered in 3.6.20.
\smallskip
We preserve the previous notations of 3.0-3. Let $G={\rm Spec}(R)$ and let $g_2\in G(W(k))$. Let $e_M:=\dim\bigl(GL(M)\bigr)=d_M^2$. Let $d$ be the relative dimension of $G$. We have $d\in S(4,e_M)$, cf. 3.1.2.1.
\smallskip 
Let $U_0={\rm Spec}(R_0)$ and
$U_2={\rm Spec}(R_2)$ be open, affine subschemes of $G$ through which the origin $a_0$ of $G$ factors (we still denote by $a_0$ this factorization) and respectively through which $g_2$ factors (for the sake of symmetry, this factorization is denoted by $a_2$), and having the following two properties:
\smallskip
{\bf i)} There are \'etale, affine $W(k)$-morphisms $U_0\buildrel{b_0}\over\to Y:=
{\rm Spec}\bigl(W(k)[z_1,\ldots,z_d]\bigr)$ and $U_2\buildrel{b_2}\over\to Y$, such that the $W(k)$-morphisms $b_0\circ a_0$ and $b_2\circ a_2$ are defined at the level of rings by: $z_i$ goes to $0$, $i=\overline{1,d}$;
\smallskip
{\bf ii)} There are \'etale, affine $W(k)$-morphisms $W_0\buildrel{c_0}\over\to Z:={\rm Spec}\bigl(W(k)[z_1,\ldots,
z_{e_M}]\bigr)$ and $W_2\buildrel{c_2}\over\to Z$, with $W_j={\rm Spec}(S_j)$ an open, affine subscheme of $GL(M)$ containing $U_j$, $j\in\{0,2\}$, such that $\tilde d\circ b_0=c_0\circ i_0$ and $\tilde d\circ b_2=c_2\circ i_2$, with $i_0:U_0\hookrightarrow W_0$ and $i_2:U_2\hookrightarrow W_2$ as the natural inclusions and with $\tilde d:Y\hookrightarrow Z$ as the $W(k)$-monomorphism which at the level of rings takes $z_i$ into $z_i$ (resp. into $0$) if $i\in S(1,d)$ (resp. if $i\in S(d+1,e_M)$). 
\medskip
{\bf 3.6.0. Some preliminaries.} From now on, till the end of 3.6.8, we assume $k$ is an infinite field. For the case of a finite field see 3.6.9 3).
\smallskip
$R/pR$ is an integral domain (cf. [Ti2, 3.8.1]). So the intersection $U_0\cap U_2$ is non-empty. Moreover, [Bo2, 18.3] guarantees the existence of a dense set of $k$-valued points of $(U_0\cap U_2)_k$. We choose a $W(k)$-monomorphism 
$$h:{\rm Spec}(W(k))\hookrightarrow U_0\cap U_2$$
such that the two $W(k)$-epimorphisms $W(k)[z_1,\ldots,z_d]\twoheadrightarrow W(k)$
associated to $b_0\circ i_{02}\circ h$ and $b_2\circ i_{20}\circ h$ (with $i_{02}$ and $i_{20}$ as the inclusions of $U_0\cap U_2$ into $U_0$ and respectively into $U_2$) send $z_i$, $i=\overline{1,d}$, into elements of
$\GG_m(W(k))$. Composing the $W(k)$-morphisms $c_0$ and $c_2$ with $W(k)$-automorphisms of $Z$ of the form
$z_i\to u_iz_i$, with $u_i\in\GG_m(W(k))$, $i=\overline{1,e_M}$, we can assume that both these two $W(k)$-epimorphisms
send $z_i$ to $1$, $i=\overline{1,d}$. 
\smallskip
Let $U^0_0$, $U^0_2$, $W^0_0$ and $W^0_2$ be the right translations of $U_0$, $U_2,W_0$ and respectively $W_2$ by $h^{-1}$. So the $W(k)$-monomorphism ${\rm Spec}(W(k))\hookrightarrow GL(M)$ defining $1_M$ factors through $U^0_0$, $U^0_2$, $W^0_0$ and $W^0_2$. Let $\Phi$ be the Frobenius lift of $Y$ and $Z$ taking
$z_i-1$ to $(z_i-1)^p$, $i=\overline{1,e_M}$. Let $\Phi_{R_0}$, $\Phi_{R_2}$, $\Phi_{S_0}$ and $\Phi_{S_2}$ be the Frobenius lifts of $U^{0\wedge}_0$, $U^{0\wedge}_2$, $W^{0\wedge}_0$ and respectively of $W^{0\wedge}_2$ such that the two diagrams below:
$$
\def\mapright#1{\smash{
\mathop{\longrightarrow}\limits^{#1}}}
\def\mapdown#1{\Big\downarrow
\rlap{$\vcenter{\hbox{$\scriptstyle#1$}}$}}
\matrix{U^{0\wedge}_j &\mapright{{r^\wedge_h\atop\sim}} &U^\wedge_j 
&\mapright{b^\wedge_j} &Y^{\wedge}\cr
\mapdown{\Phi_{R_j}} &&&&\mapdown{\Phi^\wedge}\cr
U^{0\wedge}_j &\mapright{{r^\wedge_h\atop\sim}} &U^\wedge_j 
&\mapright{b^\wedge_j} &Y^{\wedge},\cr}$$
and
\medskip
$$
\def\mapright#1{\smash{
\mathop{\longrightarrow}\limits^{#1}}}
\def\mapdown#1{\Big\downarrow
\rlap{$\vcenter{\hbox{$\scriptstyle#1$}}$}}
\matrix{W^{0\wedge}_j &\mapright{{r^\wedge_h\atop\sim}} &W^\wedge_j 
&\mapright{c^\wedge_j} &Z^\wedge\cr
\mapdown{\Phi_{S_j}} &&&&\mapdown{\Phi^\wedge}\cr
W^{0\wedge}_j &\mapright{{r^\wedge_h\atop\sim}} &W^\wedge_j 
&\mapright{c^\wedge_j} &Z^{\wedge},\cr}$$
$j\in\{0,2\}$, are cartesian; here we denote by $r_h$ any isomorphism defined by the right translation by $h$ via restrictions, while $b^\wedge_j$, $r^\wedge_h$, $c^\wedge_j$ and $\Phi^\wedge$ refer to morphisms defined by the $p$-adic completion of $b_j$, $r_h$, $c_j$ and respectively of $\Phi$. The isomorphisms $r_h$ allow us to identify $U^0_j={\rm Spec}(R_j)$ and $W^0_j={\rm Spec}(S_j)$, $j\in\{0,2\}$. Let $t_i(j)\in S_j$ be such that under the identification $W_j^0={\rm Spec}(S_j)$, the $W(k)$-morphism $c_j^\wedge\circ\tilde r_h^\wedge$, at the level of rings, takes $z_i-1$ into $t_i(j)$, $\forall j\in\{0,2\}$ and $\forall i\in S(1,e_M)$. So 
$$\Phi_{S_j}(t_i(j))=t_i(j)^p.$$ 
We also regard $t_i(j)$, $i\in S(1,d)$, as an element of $R_j$. 
\medskip
{\bf 3.6.1. Some initial $p$-divisible objects.} Let 
$$M_{R_j}:=\bigl(M\otimes_{W(k)}R_j^\wedge, F^1\otimes_{W(k)}R_j^\wedge,\Phi^r_j\bigr)$$ 
and 
$$M_{S_j}:=\bigl(M\otimes_{W(k)}S_j^\wedge,F^1\otimes_{W(k)}S_j^\wedge,\Phi_j^s\bigr),$$ 
with $\Phi^r_j:=r_j(h\vph_0\otimes 1)$ and $\Phi^s_j:=s_j(h\vph_0\otimes 1)$. Here $r_j\in G(R_j^\wedge)=G^\wedge(R_j^\wedge)$ (resp. $s_j\in GL(M)(S_j^\wedge)=GL(M)^\wedge(S_j^\wedge)$)
is the universal elements of $G^\wedge$ (resp. of $GL(M)^\wedge$), defined by the $p$-adic completion of the inclusion
$U^0_j\hookrightarrow G$ (resp. $W^0_j\hookrightarrow GL(M)$). $M_{R_j}$ (resp. $M_{S_j}$) is a $p$-divisible object of $\Mm\Mf_{[0,1]}(R_j)$ (resp. of $\Mm\Mf_{[0,1]}(S_j)$). So for any $n\in\NN$, $M_{R_j}/p^nM_{R_j}$ is an object of $\Mm\Mf_{[0,1]}(R_j)$ and $M_{S_j}/p^nM_{S_j}$ is an object of $\Mm\Mf_{[0,1]}(S_j)$. Above $j\in\{0,2\}$.  
\medskip
{\bf 3.6.1.1. Their pull backs.} For any formally \'etale, affine $W(k)$-morphism ${\rm Spec}(\tilde Q^\wedge_j)\buildrel{d_j}\over\to{\rm Spec}(R^\wedge_j)$, with $\tilde Q_j$ a $W(k)$-algebra, we denote by 
$$
M_{\tilde Q_j}:=\bigl(M\otimes_{W(k)}\tilde Q_j^\wedge,F^1\otimes_{W(k)}\tilde Q_j^\wedge, r_j\circ d_j(h\vph_0\otimes 1)\bigr)
$$ 
the pull back of $M_{R_j}$ through $d_j$, for $\Phi_{\tilde Q_j}$  the only Frobenius lift of $\tilde Q^\wedge_j$ satisfying 
$$d_j\circ\Phi_{\tilde Q_j}=\Phi_{R_j}\circ d_j.$$
 It is a $p$-divisible object of $\Mm\Mf_{[0,1]}(\tilde Q_j)$. We call it the extension of $M_{R_j}$ to ${\rm Spec}(\tilde Q_j)$ or to $\tilde Q_j$ (through $d_j$ or through the $W(k)$-homomorphism $R_j^\wedge\to \tilde Q_j^\wedge$ defining $d_j$). 
\smallskip
For any $n\in\NN$, $M_{\tilde Q_j}/p^nM_{\tilde Q_j}$ is an object of $\Mm\Mf_{[0,1]}(\tilde Q_j)$. 
\smallskip
In particular,
if ${\rm Spec}(\hat R^0_j)$ is the completion of $U^0_j$ in the origin, we denote by 
$$M_{\hat R^0_j}:=(M\otimes_{W(k)} \hat R^0_j,F^1\otimes_{W(k)} \hat R^0_j,\hat \Phi^r_j)$$ 
the extension of $M_{R_j}$ through the canonical $W(k)$-monomorphism $i_{R_j}:R_j\hookrightarrow\hat R^0_j$,
 for $\Phi_{\hat R^0_j}$  the natural Frobenius lift of $\hat R^0_j$ induced from $\Phi_{R_j}$ by
completion, $j\in\{0,2\}$. Here 
$$\hat\Phi^r_j:=\hat r_j(h\vph_0\otimes 1),$$ 
with $\hat r_j\in G(\hat R^0_j)$ defined by $r_j$ via $i_{R_j}$. $M_{\hat R_j^0}$ is a $p$-divisible object of $\Mm\Mf_{[0,1]}(\hat R_j^0)$.
\smallskip
Similarly, if $k_1$ is a perfect field containing $k$, we denote by $M_{{R_j}_{W(k_1)}}$ its pull back through the canonical affine $W(k)$-morphism ${\rm Spec}({R_j}_{W(k_1)})\to {\rm Spec}(R_j)$, the Frobenius lift of ${R_j}_{W(k_1)}$ being the logical one. We also call it the extension of $M_{R_j}$ to $W(k_1)$; we view ${R_j}_{W(k_1)}$ as a $W(k_1)$-algebra.
\smallskip
In a similar manner we define the pull back of any  object (or $p$-divisible object) of $\Mm\Mf_{[\tilde a,\tilde b]}(R_3)$, with $R_3$ a regular, formally smooth $W(k)$-algebra for which a Frobenius lift of $R_3^\wedge$ is chosen (fixed), through a formally \'etale, affine $W(k)$-morphism ${\rm Spec}(\tilde Q_3^\wedge)\to {\rm Spec}(R_3^\wedge)$, or through the natural affine $W(k)$-morphism ${\rm Spec}(R_{3W(k_1)})\to {\rm Spec}(R_3)$, with $k_1$ as above. 
\smallskip
All the above terminology conforms to 2.2.1.3.
\medskip
{\bf 3.6.1.1.1. Definitions and notations. 1)} Let $d_j:{\rm Spec}(\tilde Q_j^\wedge)\to {\rm Spec}(R_j^\wedge)$ be a formally \'etale, affine $W(k)$-morphism, with $\tilde Q_j$ a $W(k)$-algebra. We recall that (cf. beginning of \S 3) we have a direct sum decomposition $M=F^1\oplus F^0$. Let 
$$
\Phi(\tilde Q_j)^0:=r_j\circ d_j(h\vph_0\otimes 1)
$$ 
be the $\Phi_{\tilde Q_j}$-linear endomorphism of $M\otimes_{W(k)} \tilde Q_j^\wedge$ defined by $r_j\circ d_j$ and let  
$$
\Phi(\tilde Q_j)^1:(F^1+pM)\otimes_{W(k)}\tilde Q_j^\wedge\to M\otimes_{W(k)} \tilde Q_j^\wedge
$$ 
be the $\Phi_{\tilde Q_j}$-linear map defined by the formula 
$$p\Phi(\tilde Q_j)^1(m)=\Phi(\tilde Q_j)^0(m),$$ 
where $m\in (F^1+pM)\otimes_{W(k)} \tilde Q_j^\wedge$. For $n\in\NN$ we still denote by $\Phi(\tilde Q_j)^i$ its reduction mod $p^n$ ($i\in\{0,1\}$). 
\smallskip
{\bf 2)} We say $M_{\tilde Q_j}/p^nM_{\tilde Q_j}$ potentially can be viewed as an object of $\Mm\Mf_{[0,1]}^\nabla(\tilde Q_j)$ if there is a connection 
$$
\nabla:M\otimes_{W(k)} \tilde Q_j/p^n\tilde Q_j\to M\otimes_{W(k)} \Om_{(\tilde Q_j/p^n\tilde Q_j)/W_n(k)}
$$ 
such that we have:
$$\nabla\circ\Phi(\tilde Q_j)^0(m)=p\Phi(\tilde Q_j)^0\circ d\Phi_{\tilde Q_j*}/p\circ\nabla(m),\leqno (E_1)$$
if $m\in F^0\otimes_{W(k)} \tilde Q_j/p^n\tilde Q_j$, and
$$\nabla\circ\Phi(\tilde Q_j)^1(m)=\Phi(\tilde Q_j)^0\circ d\Phi_{\tilde Q_j*}/p\circ\nabla(m),\leqno (E_2)$$
if $m\in F^1\otimes_{W(k)} \tilde Q_j/p^n\tilde Q_j$. Here $d\Phi_{\tilde Q_j*}/p$ denotes the differential of the Frobenius lift $\Phi_{\tilde Q_j}$ of $\tilde Q_j^\wedge$ divided by $p$ and then taken mod $p^n$. About $\nabla$ we say it makes (or allows) $M_{\tilde Q_j}/p^nM_{\tilde Q_j}$ potentially to be viewed as an object of $\Mm\Mf_{[0,1]}^\nabla(\tilde Q_j)$; often (to be short) we also say $\nabla$ is a connection on $M_{\tilde Q_j}/p^nM_{\tilde Q_j}$. Warning: the system of equations obtained by putting $(E_1)$ and $(E_2)$ together, does not depend on the choice of the direct supplement $F^0$ of $F^1$ in $M$; moreover, as we chose $\Phi(\tilde Q_j)^0$ and $\Phi(\tilde Q_j)^1$ to be Frobenius maps (i.e. $\Phi_{\tilde Q_j}$-linear maps) and not $\tilde Q_j^\wedge$-linear maps, we got $m\in F^i\otimes_{W(k)} \tilde Q_j/p^n\tilde Q_j$ and not in $F^i\otimes_{W(k)} {}_{\sg}\tilde Q_j/p^n\tilde Q_j$.   
\smallskip
We say $\nabla$ respects the $G$-action, if $\forall i\in S(1,d)$, the $\tilde Q_j/p^m\tilde Q_j$-endomorphism of $M\otimes_{W(k)} \tilde Q_j/p^n\tilde Q_j$ that takes $m\in M$ into 
$$\nabla({{\partial}\over {\partial t_i(j)}})(m)\in M\otimes_{W(k)} \tilde Q_j/p^n\tilde Q_j,$$ 
is an element of ${\rm Lie}(G)\otimes_{W(k)} \tilde Q_j/p^n\tilde Q_j$.
\smallskip
We say $M_{\tilde Q_j}/p^nM_{\tilde Q_j}$ potentially is an object of $\Mm\Mf_{[0,1]}^\nabla(\tilde Q_j)$ if there is precisely one connection $\nabla$ on it; about this connection $\nabla$ we say it makes (or allows) $M_{\tilde Q_j}/p^nM_{\tilde Q_j}$ potentially to be an object of $\Mm\Mf_{[0,1]}^\nabla(\tilde Q_j)$. 
\smallskip
We use the same language for the case of $p$-divisible objects, i.e. when $n=\infty$. Warning: in such a case 
$$\nabla:M\otimes_{W(k)} \tilde Q_j^\wedge\to M\otimes_{W(k)} \Om_{\tilde Q_j/W(k)}^\wedge,$$
and we also say that $\Phi(\tilde Q_j)^0$ is $\nabla$-parallel.
\smallskip
If in any of the above cases the connection $\nabla$ is integrable, then we drop the word potentially. We prove later on (cf. 3.6.18.4.1), that any connection $\nabla$ which makes $M_{\tilde Q_j}/p^nM_{\tilde Q_j}$ potentially to be viewed as an object of $\Mm\Mf_{[0,1]}^\nabla(\tilde Q_j)$, is automatically integrable and so the word potentially is just for the time being.
\smallskip
{\bf 3)} Let $\Mz_k$ be the closed subscheme of $U_0^0$, $U_2^0$, $G$, $W_0^0$, $W_2^0$ and of $GL(M)$, defining the origin of the special fibre of $G_k$. So $\Mz_k={\rm Spec}(k)$. 
\smallskip
{\bf 4)} For an arbitrary regular, formally smooth $W(k)$-algebra $R_3$ as in the end of 3.6.1.1, and for any ($p$-divisible) object of $\Mm\Mf_{[\tilde a,\tilde b]}(R_3)$, with $\tilde a,\tilde b\in\ZZ$, $\tilde a\le\tilde b$, we use the same language as in 2) but always working with connections satisfying the Griffiths transversality condition. The equations needed to be satisfied are similar to the ones of $(E_1)$ and $(E_2)$ of 2) (for instance, see [Fa1, p. 33]). Warning: in such a generality, the word potentially is not for the time being (cf. 3.6.18.5.5 below). 
\smallskip
As a quick exemplification, we refer to the context of ${\got C}$ of 2.2.1 c), with $M$ a projective $R^\wedge$-module. We assume $F^{\tilde a}(M)=M$ and $F^{\tilde b+1}(M)=\{0\}$. The set of equations needed to be satisfied in order that a connection $\nabla:M\to M\otimes_{R^\wedge} \Om^\wedge_{R^\wedge/W(k)}$ makes ${\got C}$ potentially to be viewed as a $p$-divisible object of $\Mm\Mf_{[\tilde a,\tilde b]}^\nabla(R)$ are 
$$\nabla\circ\vph_i(m_i)=\vph_{i-1}\circ d\Phi_{R^\wedge *}/p\circ\nabla(m_i)  ,\leqno (E_i)$$
 with $m_i\in F^i(M)$, for $i\in S(\tilde a,\tilde b)$, and
$$\nabla(F^i(M))\subset F^{i-1}(M)\otimes_{R^\wedge} \Om^\wedge_{R^\wedge/W(k)}, \leqno (GT)$$
$i\in S(2+\tilde a,\tilde b)$. Taken them mod $p^n$, we get the set of equations needed to be satisfied in order that a connection $\nabla:M/p^nM\to M/p^nM\otimes_{R/p^nR} \Om_{R/p^nR/W_n(k)}$ makes ${\got C}/p^n{\got C}$ potentially to be viewed as an object of $\Mm\Mf_{[\tilde a,\tilde b]}^\nabla(R)$. Warning: we always look at (GT) as an equation (it can be expressed in terms of some morphisms of $R^\wedge$-modules being 0). 
\smallskip
{\bf 5)} Let $X$ and $X_1$ be regular, formally smooth $W(k)$-schemes or $p$-adic formal schemes over $W(k)$. We consider a formally \'etale morphism $m_X:X_1\to X$ and we assume $X$ (or $X^\wedge$ in case $X$ is a scheme) is equipped with a Frobenius lift $\Phi_X$. We denote by $\Phi_{X_1}$ the Frobenius lift of $X_1$ (or of $X_1^\wedge$, in case $X_1$ is a scheme) defined naturally by $\Phi_X$ via $m_X$. Let ${\got C}$ be an arbitrary ($p$-divisible) object of $\Mm\Mf(X)$. As in 4) we speak about connections on $m_X^*({\got C})$, i.e. about connections on the underlying module of $m_X^*({\got C})$ which make $m_X^*({\got C})$ potentially to be viewed as a ($p$-divisible) object of $\Mm\Mf^\nabla(X_1)$, or about $m_X^*({\got C})$ potentially to be a ($p$-divisible) object $\Mm\Mf^\nabla(X_1)$, etc. Warning: we always assume that these connections satisfy the Griffiths transversality condition.
\medskip
{\bf 3.6.1.1.2. Digression: the nilpotency mod $p$.} We refer to 3.6.8.1.1 4) with $n=1$. Let $\tilde M:=M/pM$ be endowed with the reduction mod $p$ of the filtration of $M$. We assume $R=W(k)[[x_1,...,x_m]]$ is a ring of formal power series. Let $N_i$ be the $R/pR$-submodule of $\tilde M$ generated by $\vph_i(F^i(\tilde M))$. Let $j\in S(1,m)$. We consider the $k$-endomorphism of $\tilde M$ defined by $D_j:=\nabla({\partial\over {\partial x_j}})$. We have
$$D_j(\vph_i(F^i(\tilde M))\subset N_{i-1},\leqno (INCL)$$cf. ($E_i$) and (GT). We get  
$$D_j^p(N_i)\subset\oplus_{s\in S(\tilde a,i-1)} N_s.\leqno (NIL1)$$ 
As $M=\oplus_{i\in S(\tilde a,\tilde b)} N_i$, we also get 
$$D_j^{p^{\tilde b-\tilde a+1}}=0.\leqno (NIL2)$$
\indent
In general we get: if $s\in\NN$, if $j_1$,..., $j_s\in S(1,m)$ and $i_1$,..., $i_s\in\NN\cup\{0\}$ are such that $i_1+i_2+...+i_s\ge p^{\tilde b-\tilde a+1}$, then (cf. (INCL))
$$D_{j_1}^{i_1}D_{j_2}^{i_2}...D_{j_s}^{i_s}=0.\leqno (NIL3)$$
We conclude:
\medskip
{\bf Corollary.} {\it All connections on some object of some Fontaine category of objects not involving connections, are nilpotent mod $p$, regardless of the fact that they are or are not integrable.}
\medskip
{\bf 3.6.1.1.3. The convenience assumption.} Sometimes it is convenient to make the assumption that the closed subscheme of ${\rm Spec}(R_j^0/pR_j^0)$ defined by $t_i(j)=0$, $i$ running through the elements of an arbitrary  subset of $S(1,d)$, is connected (so in particular, $\Mz_k$ is the closed subscheme of ${U^0_j}_k$ defined by $t_i(j)=0$, $\forall i\in S(1,d)$). 
\smallskip
This can always be achieved by removing a closed subscheme of $U_j^0$ which is not $U_j^0$ itself. This should be done before the choice of $h$ so that we do not lose track of the morphisms $a_0$ and $a_2$. It can be done as follows: we choose an open, affine subscheme $\Mu={\rm Spec}(R_{\Mu})$ through which the origin of $G$ factors and an \'etale morphism $b_{\Mu}:\Mu\to Y$, such that the closed subscheme ${\rm Spec}(R_{\Mu}/(p,(z_i)_{i\in I_d}))$ of $\Mu$, with $I_d$ an arbitrary subset of $S(1,d)$, is connected and the logical $W(k)$-morphism ${\rm Spec}(R_{\Mu}/((z_i)_{i\in I}))\to G$ is the origin. We can take now $U_0=\Mu$ and $U_2=g_2\Mu$; as the morphisms $b_0$ and $b_2$ we take $b_{\Mu}$ and respectively the composite of the left translation by $g_2^{-1}$ isomorphism $U_2\tilde\to U_0$ with $b_{\Mu}$.  
\smallskip
Warning: unless specifically stated we do not work under above assumption.
\medskip
{\bf 3.6.1.2. Lemma (the uniqueness principle in a convenient Shimura $\sg$-crystal context).} {\it For any formally \'etale, affine $W(k)$-morphism $l_j:{\rm Spec}(\tilde Q^\wedge_j)\to {\rm Spec}(R^\wedge_j)$, with $\tilde Q_j$ a $W(k)$-algebra such that all connected components of ${\rm Spec}(\tilde Q_j/p\tilde Q_j)$ have a non-empty intersection with $\ell_j^{-1}(\Mz_k)$, and for every $n\in\NN$, $M_{\tilde Q_j}/p^nM_{\tilde Q_j}$ potentially can be viewed in at most one way as an object of $\Mm\Mf_{[0,1]}^\nabla(\tilde Q_j)$; similarly, $M_{\tilde Q_j}$ potentially can be viewed in at most one way as a $p$-divisible object of $\Mm\Mf_{[0,1]}^\nabla(\tilde Q_j)$.}
\medskip
{\bf 3.6.1.3. Theorem (the $\nabla$, the moduli, and the surjectivity principle in a convenient Shimura $\sg$-crystal context). 1)} {\it For any $n\in\NN$ there is a formally  \'etale, affine $W(k)$-morphism 
$${\rm Spec}(Q^\wedge_{j,n})\buildrel{\ell_{j,n}}\over\to
{\rm Spec}(R^\wedge_j),$$ 
$j\in\{0,2\}$, with $Q_{j,n}$ a smooth $W(k)$-algebra, such that:
\medskip
\item{i)} $\ell_{j,n}^{-1}(\Mz_k)={\rm Spec}(k)$; 
\smallskip
\item{ii)} $M_{Q_{j,n}}/p^nM_{Q_{j,n}}$ is an object of $\Mm\Mf^\nabla_{[0,1]}(Q_{j,n})$, i.e. there is a unique 
connection
$$
\nabla^n_j:M\otimes_{W(k)}Q_{j,n}/p^nQ_{j,n}\to M\otimes_{W(k)} \Om_{(Q_{j,n}/p^nQ_{j,n})/W_n(k)}
$$
which makes $M_{Q_{j,n}}/p^nM_{Q_{j,n}}$ potentially to be viewed as an object of $\Mm\Mf_{[0,1]}^\nabla(Q_{j,n})$; $\nabla_j^n$ is integrable (so the word potentially can be disregarded), nilpotent mod $p$, and respects the $G$-action;
\smallskip
\item{iii)} for any other formally \'etale, affine $W(k)$-morphism ${\rm Spec}(\tilde Q^\wedge_j)\buildrel{\tilde\ell_j}\over\to {\rm Spec}(R^\wedge_j)$ (with $\tilde Q_j$ a $W(k)$-algebra) such that $M_{\tilde Q_j}/p^nM_{\tilde Q_j}$ potentially can be viewed as an object of $\Mm\Mf^\nabla_{[0,1]}(\tilde Q_j)$ and all connected components of ${\rm Spec}(\tilde Q_j/p\tilde Q_j)$ have a non-empty intersection with $\tilde\ell_j^{-1}(\Mz_k)$, there is a unique morphism ${\rm Spec}(\tilde Q^\wedge_j)\buildrel{\tilde\ell_j^n}\over\to {\rm Spec}(Q^\wedge_{j,n})$ such that $\tilde\ell_j=\ell_{j,n}\circ\tilde\ell_j^n$ and the connection (it is unique cf. 3.6.1.2) on $M\otimes_{W(k)} \tilde Q_j/p^n\tilde Q_j$ which makes $M_{\tilde Q_j}/p^nM_{\tilde Q_j}$ potentially to be viewed as an object of $\Mm\Mf_{[0,1]}^\nabla(\tilde Q_j)$ is the pull back of $\nabla_j^n$ through $\tilde\ell_j^n$ mod $p^n$;
\smallskip
\item{iv)} the special fibre of ${\rm Spec}(Q_{j,n})$ is a geometrically connected $k$-scheme.
\medskip
{\bf 2)} The unique $R_j^\wedge$-morphism (cf. ii) and iii)) 
$$\ell_j^n:{\rm Spec}(Q_{j,n+1}^\wedge)\to {\rm Spec}(Q_{j,n}^\wedge)$$ 
such that $\ell_{j,n+1}=\ell_{j,n}\circ\ell_j^n$ and the pull back of $\nabla^n_j$ through the morphism $\ell_j^n$ mod $p^n$ is $\nabla^{n+1}_j$ mod $p^n$, when taken mod $p$ (i.e. the resulting morphism at the level of special fibres) is: \'etale, quasi-finite, and generically an \'etale cover of degree $p^{d(n)}$, with $d(n)\in [0,d^2]$. Similarly, the morphism $\ell_{j,1}$ mod $p$ has all these properties.
\smallskip
{\bf 3)} The construction of ${\rm Spec}(Q_{j,n}^\wedge)$ commutes with (perfect) base field extensions, i.e. if $k_1$ is a perfect field containing $k$, then the $p$-adic completion of the extension of $\ell_{j,n}$ to $W(k_1)$ is the morphism of 1) obtained for the same $n\in\NN$  and  for the situation where $G$ is replaced by $G_{W(k_1)}$, $R_j$ is replaced by ${R_j}_{W(k_1)}$ and $M_{R_j}$ is replaced by its extension to $W(k_1)$.
\smallskip
{\bf 4)} Let ${\rm Spec}(Q_j)$ be the $p$-adic completion of the $\NN$-projective limit of ${\rm Spec}(Q_{j,n}^\wedge)$ under the transition morphisms $\ell_j^n$, $n\in\NN$, and let 
$$\ell_j:{\rm Spec}(Q_j)\to U_j^0={\rm Spec}(R_j)$$ 
be the resulting $W(k)$-morphism. There is a unique $W(k)$-morphism $a^j:{\rm Spec}(W(k))\to {\rm Spec}(Q_j)$ which composed with $l_j$ is $a_0$. $M_{Q_j}$ is a $p$-divisible object of $\Mm\Mf_{[0,1]}^\nabla(Q_j)$. The $k$-scheme ${\rm Spec}(Q_j/pQ_j)$ is geometrically connected. Moreover, $l_j$ has a universal property similar to the one described in iii) of 1): 
\medskip
{\bf UP.} For any formally \'etale, affine $W(k)$-morphism ${\rm Spec}(\tilde Q^\wedge_j)\buildrel{\tilde\ell_j}\over\to {\rm Spec}(R^\wedge_j)$ (with $\tilde Q_j$ a $W(k)$-algebra) such that $M_{\tilde Q_j}$ potentially can be viewed as a $p$-divisible object of $\Mm\Mf^\nabla_{[0,1]}(\tilde Q_j)$ and all connected components of ${\rm Spec}(\tilde Q_j/p\tilde Q_j)$ have a non-empty intersection with $\tilde\ell_j^{-1}(\Mz_k)$, there is a unique morphism ${\rm Spec}(\tilde Q^\wedge_j)\buildrel{\tilde\ell_j^{\infty}}\over\to {\rm Spec}(Q_j)$ such that $\tilde\ell_j=\ell_{j}\circ\tilde\ell_j^\infty$ and the connection $\tilde\nabla_j$ (it is unique cf. 3.6.1.2) on $M\otimes_{W(k)} \tilde Q_j^\wedge$ which makes $M_{\tilde Q_j}$ potentially to be viewed as a $p$-divisible object of $\Mm\Mf_{[0,1]}^\nabla(\tilde Q_j)$, when taken mod $p^n$, is the pull back of $\nabla_j^n$ through the composite of $\tilde\ell_j^\infty$ with the natural $W(k)$-morphism ${\rm Spec}(Q_j)\to {\rm Spec}(Q_{j,n})$, $\forall n\in\NN$.
\medskip
{\bf 5)} There is a reduced, closed subscheme $\Mb_j(k)$ of $G_k$, not containing the origin, depending on $h\vph_0$ and on $\Phi_{R_j}$ but not on how small or big $U_j^0$ is, and such that the morphism $\ell_j$ has non-empty fibres over (geometric) points of ${\rm Spec}(R_j/pR_j)$ which (when viewed as points of $G_k$) do not factor through $\Mb_j(k)$. For any perfect field $k_1$ containing $k$, we get similarly (cf. 3)) a reduced, closed subscheme $\Mb_j(k_1)$ of $G_{k_1}$ having similar properties as $\Mb_j(k)$; in fact we can take $\Mb_j(k_1):=\Mb_j(k)_{k_1}$.
\medskip
{\bf 6)} ${\rm Spec}(Q_j/pQ_j)$ is an $AG$ $k$-scheme.}
\medskip
The proofs of 3.6.1.2-3 are presented in 3.6.8, after we include some remarks (see 3.6.1.4-6), some direct applications (see 3.6.2-5), and some extra applications (see 3.6.7) which rely on some independent results of 3.6.6.
\medskip
{\bf 3.6.1.4. Remarks. 1)} Based on 3.6.18.4.2 below, it can be checked that in many situations we can take $Q_{j,n}=R_j$ (so $d(n)=0$), $\forall n\in\NN$. All these situations will be listed in \S 7, after we present (in \S 7) the classification of Shimura group pairs over $\ZZ_p$, and introduce (see 3.10) some extra terminology. See also 3.11.4 below.
\smallskip
{\bf 2)} The $W(k)$-morphism $\ell_{j,n}$ is an \'etale cover iff $Q_{j,n}^\wedge=R_j^\wedge$ (cf. i) of 3.6.1.3 1)). A simple example when $\ell_{j,1}$ is not an \'etale cover was suggested by A. Ogus. 
\smallskip
We take $d_M=2$ and $G=GL(M)$. We assume $(M,F^1,h\vph_0)$ is the filtered $\sg$-crystal of the canonical lift of an ordinary elliptic curve over $k$. We also assume the existence of a $W(k)$-basis $\{x,y\}$ of $M$ such that $h\vph_0$ takes $x$ and $y$ respectively into $x$ and $py$. We choose a Frobenius lift of $G^\wedge$ such that under the canonical identification $G={\rm Spec}(W(k)[x_{11},x_{12},x_{21},x_{22}][{1\over {x_{11}x_{22}-x_{21}x_{12}}}])$ defined by the above ordered basis of $M$, $x_{12}$ goes to $x_{12}^p$, $x_{21}$ goes to $x_{21}^p$, while $x_{ii}-1$ goes to $(x_{ii}-1)^p$, $i=\overline{1,2}$. We can take $U_0=U_2=G$. 
\smallskip
To see that $l_{j,1}$ is not an \'etale cover, it is enough to see that any natural pull back of $M_{R_j}/pM_{R_j}$ to a pro-\'etale cover of $\GG_a^\wedge$, with $\GG_a$ viewed as the subgroup of $G$ acting trivially on $x$ and on the quotient module $M/W(k)x$, has no connection on it. This is easy: we can assume $k=\bar k$; identifying $\GG_a={\rm Spec}(W(k)[t])$, with $t:=x_{12}$, the unique connection we get on $M\otimes_{W(k)} W(k)[[t]]$ making the pull back of $M_{R_j}$ through the natural $W(k)$-morphism ${\rm Spec}(W(k)[[t]])\to G$ to be a $p$-divisible object of $\Mm\Mf_{[0,1]}^\nabla({\rm Spec}(W(k)[[t]]))$, annihilates $x$ and takes $y$ into $\al(t)xdt$, where 
$$\al(t)=-\sum_{i=0}^\infty t^{p^i-1}.$$ 
The same holds mod $p$. But $\al(t)$ mod $p$ satisfies the equation (in $u$) 
$$t^{p-1}u^p=u+1.$$ 
This equation defines an \'etale scheme over ${\GG_a}_k={\rm Spec}(k[t])$ which is not an \'etale cover: it is obtained from the \'etale cover ${\rm Spec}(k[t,v]/(v^p-v-t))$ of ${\GG_a}_k$ by removing its closed subscheme (is a disjoint union of $p-1$ copies of ${\rm Spec}(k)$) defined by $t=0$, $v\neq 0$.  
\smallskip
In this example, coming back to the general notations,  the $k$-morphism 
$$\ell_j^\prime:\ell_j^{-1}({\rm Spec}(R_j/pR_j[{1\over\prod_{i=1}^d t_i(j)}])\to {\rm Spec}(R_j/pR_j[{1\over\prod_{i=1}^d t_i(j)}]):=\Mw_j$$
 is an $\NN$-pro-\'etale $k$-morphism of infinite degree. Here $d=4$, $t_i(j)=x_{ii}-1$ and $t_{i+2}(j)=x_{i,3-i}$, with $i\in\{1,2\}$. It defines an $\NN$-pro-\'etale cover above the open subscheme of $\Mw_j$ defined by $t_1(j)\neq -1$, as it can be easily checked starting from 3.6.18.4 B) below. 
\smallskip
{\bf 3)} The estimate $d(n)\le d^2$ is not the best possible. For instance if $G=GL(M)$ it can be checked (this is explained in 3.11.4 below) that we can replace $d^2$ by ${d^2\over 4}$. 
For better estimates of the real numbers $d(n)$, see 3.11.4. For the time being, we just remark that the bound $d^2$ is independent on $n$. We do not know (if or) when $d(n)=d(1)$, $\forall n\in\NN$.
\smallskip
{\bf 4)} We do not know if (or when) the $k$-morphism defined by the special fibre of $\ell_j$ is surjective. Similarly about the $k$-morphisms defined by the special fibres of $\ell_{j,n}$ and $\ell_j^n$, $n\in\NN$, $j\in\{0,2\}$. 3.6.18.5.2 below implicitly motivates why it is natural to expect (under some conditions) such $k$-morphisms to be surjective. 
\smallskip
{\bf 5)} $M_{Q_{j,m}}/p^nM_{Q_{j,m}}$ and $(M_{Q_{j,m}}/p^nM_{Q_{j,m}},\nabla_j^m\, {\rm mod}\, p^n)$, $m,n\in\NN$, $n\le m$, are rigid objects in the sense that they have few endomorphisms. For instance, their endomorphisms (called $G$-endomorphisms) which as elements of ${\rm End}(M/p^nM)\otimes_{W_n(k)} Q_{j,m}/p^nQ_{j,m}$ belong to ${\got g}/p^n{\got g}\otimes_{W_n(k)} Q_{j,m}/p^nQ_{j,m}$, are often such that these elements are the reduction mod $p^n$ of those elements of the Lie algebra of the maximal subtorus of $Z(G)$ which are fixed by $h\vph_0$. This remains true in many cases for the set of all endomorphisms, the group $Z(G)$ being replaced by the centralizer of $G$ in $GL(M)$. This important phenomenon is elaborated in \S 7, being an essential tool in the strategy of proving part of the split criterion of 1.15.4. There are variants of it when we pass to \'etale, affine ${\rm Spec}(Q_{j,m})$-schemes. In these variants extra endomorphisms can show up: like the ones springing from direct factors of ${\rm Lie}(G^{\rm ad}_{W(k_1)})$ contained in $W(0)({\rm Lie}(G^{\rm ad}_{W(k_1)}),h\vph_0\otimes 1)$; here $k_1$ is a finite field extension of $k$.
\smallskip
{\bf 6)} iii) and 3.6.1.3 UP can be formulated for arbitrary formally \'etale $R_j^\wedge$-schemes, not necessarily affine (cf. also 3.6.18.4.2 below).
\medskip
{\bf 3.6.1.5. $\sg$-$\Ms$-crystals.} Let $\tilde a,\tilde b\in\ZZ$, with $\tilde a\le\tilde b$. Let $(\tilde M,(F^i(\tilde M))_{i\in S(\tilde a+1,\tilde b)},\tilde\vph)$ be a $p$-divisible object of $\Mm\Mf_{[\tilde a,\tilde b]}(W(k))$. Let $\tilde G$ be a smooth, closed subgroup of $GL(\tilde M)$ having connected fibres and such that: 
\medskip
-- ${\rm Lie}(\tilde G)[{1\over p}]$ is stable under the action of $\tilde\vph$;  
\smallskip
-- ${\rm Lie}(\tilde G)$ together with the induced filtration (from the natural induced filtration of ${\rm End}(\tilde M)$; see 2.1) and Frobenius endomorphism $\tilde\vph$ of ${\rm Lie}(\tilde G)[{1\over p}]$, is a $p$-divisible subobject of the $p$-divisible object $({\rm End}(\tilde M),(F^i(\tilde M))_{i\in SS(\tilde a,\tilde b)},\tilde\vph)$ of $\Mm\Mf(W(k))$ (adding the Lie structure it is an $[\tilde a-\tilde b,\tilde b-\tilde a]$-filtered Lie $\sg$-crystal in the sense of 2.2.3 2)).
\medskip
The quadruple 
$$(\tilde M,(F^i(\tilde M))_{i\in S(\tilde a+1,\tilde b)},\tilde\vph,\tilde G)$$ 
is called a filtered $\sg$-$\Ms$-crystal (over $k$); here $\Ms$ stands for subobject and refers to the second of the above two conditions. As in 2.2.8 2) or 4) we get $\sg$-$\Ms$-crystals. As examples we mention the pseudo Shimura $\sg$-crystals (cf. 2.2.9 1')) and the $p$-divisible objects with a reductive structure over $k$ (cf. 2.2.8 4a)). We would like to point out that we did not assume the existence of a cocharacter $\mu:\GG_m\to\tilde G$, normalizing $F^i(\tilde M)$, $\forall i\in S(\tilde a,\tilde b)$, with $\beta\in\GG_m(W(k))$ acting through it on $F^i(\tilde M)/F^{i+1}(\tilde M)$ as the multiplication with $\beta^{-i}$ (to be compared with 2.2.8 3), where this is a consequence of the hypothesis involved). When such a cocharacter does exist, we speak about (filtered) $\sg$-$\Ms$-crystals for which the $W$-condition holds (the letter $W$ stands to honor [Wi]). It is easy to see that the Facts of 2.2.9 1) and 1') apply to the present context of (filtered) $\sg$-$\Ms$-crystals (in connection to 2.2.9 1') for $p\le 1+\tilde b-\tilde a$ we need to assume that a corresponding Galois representation does exist); keeping in mind what $\Ms$ stands for, we felt it is appropriate not to use the word ``pseudo" in front of the above terminology if the $W$-condition does not hold.
\smallskip
One could allow the fibres of $\tilde G$ not to be connected; but this can create problems. First, the Lie algebra of a smooth group $\tilde H$ over $W(k)$ is determined by its connected component of the origin. Second, if the special fibre is not connected, then we can get into trouble with problems involving Newton polygons, as for such problems it is important only how the $p$-adic completion of $\tilde H$ looks like. On the other hand: in [Fa2, rm. ii) after th. 10] the group involved does not necessarily have a connected special fibre.
\medskip
{\bf 3.6.1.6. Remark.} It is expected that everything in 3.6.1.3 (when properly formulated; see the paragraphs below) remains true in the context of $\sg$-$\Ms$-crystals. 
\smallskip
This is true for $(\tilde a,\tilde b)=(0,1)$ as well as for the context of generalized Shimura $p$-divisible objects of $\Mm\Mf(W(k))$. The proof in the case when $(\tilde a,\tilde b)=(0,1)$ is entirely the same (cf. also 3.6.8.4 1) and 3) below). 
\smallskip
In the other case we have to interpret 3.6.1.3 as a property of Shimura filtered Lie $\sg$-crystals. For instance, in many cases, we can construct a Shimura filtered $\sg$-crystal whose Shimura adjoint filtered Lie $\sg$-crystal is isomorphic to the one of our Shimura $p$-divisible object of $\Mm\Mf(W(k))$; for details see 4.6.7. The general case for this interpretation is explained in 3.6.18.7.3 C and 3.15.6. There is a second approach to this interpretation: the duality aspects of 2.2.23 B are the very root of it. In other words, the philosophy (to be elaborated in \S 7) is: 
\medskip
{\bf Ph.} {\it A great part of this 3.6 pertaining to $p$-divisible objects (or just objects) of some category $\Mm\Mf_{[0,1]}(*)$ can be extended to symmetric or alternating objects and to $p$-divisible objects which have symmetric or alternating quasi-polarizations, of some category $\Mm\Mf_{[-1,1]}(*)$.} 
\medskip
In the above two cases we can introduce as in 2.2.10 and 2.2.10.1 the notion of filtered $F$-$\Ms$-crystal and respectively of generalized Shimura (filtered) $F$-crystal or $p$-divisible object of some category $\Mm\Mf_{[\tilde a,\tilde b]}(*)$; here $*$ stands for regular, formally smooth $W(k)$-algebras or (when possible) ($p$-adic formal) $W(k)$-schemes.
\smallskip
What we meant above by ``properly formulated": following the discussion of 3.6.8.9 below, we might have to work with a different Frobenius endomorphism of the module underlying the $p$-divisible objects involved (i.e. we do not use a natural equivalent $e_j\in \tilde G(*^\wedge)$ of $r_j$ of 3.6.1 but some other element $\tilde e_j\in\tilde G(*^\wedge)$, very carefully chosen; warning: we do have in mind the requirement that $\tilde e_j$ mod $p$ is $e_j$ mod $p$).
\medskip  
{\bf 3.6.2. First consequences.} We now look at some immediate implications of 3.6.1.2-3. It is known (see 2.2.1.1 2); here --as well as in all subsequent places relying on it-- is the second place where we need $p\ge 3$) that the category  $\Mm\Mf^\nabla_{[0,1]}(Q_{j,n})$ is antiequivalent (via the $\DD^{-1}$ functor) to $p-FF({\rm Spec}(Q_{j,n}^\wedge))$,
$j\in\{0,2\}$. So corresponding to $M_{Q_{j,n}}/p^nM_{Q_{j,n}}$ and $\nabla_j^n$, we get a finite, flat, commutative group scheme 
$$D_{j,n}$$ 
over ${\rm Spec}(Q^\wedge_{j,n})$ of rank $p^{nd_M}$. Due to the uniqueness of $\nabla^n_j$, $(\ell_j^n)^*(D_{j,n})$ can be naturally identified with $D_{j,n+1}[p^n]$, $\forall n\in \NN$, $j\in\{0,2\}$. So the objects $M_{Q_j}/p^nM_{Q_j}$ together with the connection on $M\otimes_{W(k)} Q_j/p^nQ_j$ obtained from $\nabla_{j,n}$ by extension of scalars, $n\in \NN$, when put together (i.e. when we consider the limit of the $\NN$-projective system they define naturally) ``give birth" under the mentioned antiequivalence to a $p$-divisible group $D_j$ over ${\rm Spec}(Q_j)$: we have $D_j[p^n]={D_{j,n}}_{{\rm Spec}(Q_j)}$. Let $\nabla_j$ be the connection on $M\otimes_{W(k)} Q_j$ which makes $M_{Q_j}$ to be a $p$-divisible object of $\Mm\Mf_{[0,1]}^\nabla(Q_j)$. As the connections $\nabla_j^n$, $n\in\NN$, respect the $G$-action, $\nabla_j$ respects as well the $G$-action; so $\nabla_j(t_{\al})=0$, $\forall\in\Mj$. We conclude: 
\medskip
{\bf Corollary.} {\it The pair 
$$
D(j):=(D_j,(t_{\al})_{\al\in\Mj})
$$ 
is a Shimura $p$-divisible group over $Q_j$.}
\medskip
{\bf 3.6.2.0. A fact.} Let $\tilde\ell_j$ be as in 3.6.1.3 4). From the UP of 3.6.1.3 4) we get:
$$(M_{\tilde Q_j},\tilde\nabla_j)=\DD((\tilde\ell_j^{\infty})^{*}(D_j)).$$ 
\smallskip
{\bf 3.6.2.1. Some properties.} For simplifying the notations we assume now, till the end of 3.6.5, that $k=\bar k$. We also assume $h^{-1}$ (resp. $g_2h^{-1}$) mod $p$ does not define a point of $\Mb_0(k)$ (resp. of $\Mb_2(k)$). So we have the following properties (the first three, i.e. a) to c), are a consequence of the shape of the Frobenius lift $\Phi_{R_j}$ of $R_j^\wedge$ of 3.6.0).
\medskip
{\bf a)} $a_j^*(D(j))$ is the Shimura $p$-divisible group $(\DD^{-1}(M,F^1,h\vph_0),(t_{\al})_{\al\in\Mj}))$, $j\in\{0,2\}$.
\smallskip
{\bf b)} The pull back of $D(0)$ through any $W(k)$-morphism ${\rm Spec}(W(k))\to  {\rm Spec}(Q_0)$ which composed with $\ell_0$ is the $W(k)$-morphism  
${\rm Spec}(W(k))\to U^0_0$ defining $h^{-1}$, is the Shimura $p$-divisible group $(\DD^{-1}(M,F^1,\vph_0),(t_{\al})_{\al\in\Mj}))$.
\smallskip
{\bf c)} The pull back of $D(2)$ through any $W(k)$-morphism ${\rm Spec}(W(k))\to {\rm Spec}(Q_2)$ which composed with $\ell_2$ is the $W(k)$-morphism ${\rm Spec}(W(k))\to U^0_2$ defining $g_2h^{-1}$, is the Shimura $p$-divisible group $(\DD^{-1}(M,F^1,g_2\vph_0),(t_{\al})_{\al\in\Mj}))$.
\smallskip
{\bf d)} We denote by 
$$
O_0:=W(k)[[w_1,\ldots,w_d]]
$$ 
the $W(k)$-algebra of the completion of $G$ (or of
$U^0_0$ or $U^0_2$) in $a_0$. It is the same as the $W(k)$-algebra of the completion of ${\rm Spec}(Q_j)$ in $a^j$, $j\in\{0,2\}$; so we have $O_0=\hat R^0_0=\hat R^0_2$. The pull backs of $D(0)$ and $D(2)$ to 
${\rm Spec}(O_0)$ are denoted by $D^0_0$ and respectively by $D^0_2$. The $p$-divisible objects with tensors of $\Mm\Mf_{[0,1]}^\nabla(O_0)$ they define are both Shimura filtered $F$-crystals: $(M,F^1,h\vph_0,G,\tilde f_0,(t_{\al})_{\al\in\Mj})$ and respectively $(M,F^1,h\vph_0,G,\tilde f_2,(t_{\al})_{\al\in\Mj})$. Warning: the Frobenius lift of $O_0$ corresponding to $D^0_0$ (i.e. defined by $\tilde f_0$) is different from the Frobenius lift of $O_0$
corresponding to $D^0_2$ (i.e. defined by $\tilde f_2$). 
\smallskip
But these two Shimura filtered $F$-crystals are
induced one from each other, cf. [Fa2, i) and iii) after th. 10] (the mentioned two Frobenius lifts of $O_0$ are as needed for applying loc. cit., cf. the way we defined $\Phi_{R_j}$ in 3.6.0). This means that there are $W(k)$-homomorphisms
$$
O_0{{\scriptstyle{\om_0}\atop\rightrightarrows}\atop\scriptstyle{\om_2}}O_0
$$ 
such that the Shimura filtered
$F$-crystal with tensors associated to the pull back of $D^0_j$ through the morphism of schemes associated to $\om_j$, is $1_{\Mj}$-isomorphic (in the same sense as in 2.2.9 6)) to $(M,F^1,h\vph_0,G,\tilde f_{2-j},(t_{\al})_{\al\in\Mj})$, under an isomorphism lifting the identity automorphism of the Shimura $p$-divisible group of a), $\forall j\in\{0,2\}$.
\medskip
{\bf 3.6.2.2. Terminology.} We call this last fact, that $D_0^0$ and $D_2^0$ are induced one from each other, the very weak gluing principle. It is trivial to check, using 2.2.21, that we can actually choose $\om_j$, $j\in\{0,2\}$, to be isomorphisms of $W(k)$-algebras, and we call this fact the weak gluing principle. 
\medskip
{\bf 3.6.3. Pull backs to Witt rings.} We consider (cf. 3.6.2) the following Lie $p$-divisible object  of $\Mm\Mf_{[-1,1]}^\nabla(Q_j)$: 
$$End(M_{Q_j}):=\bigl({\rm End}(M)\otimes Q_j,r_j\circ\ell_j(h\vph_0\otimes 1),F^0({\rm End}(M))\otimes Q_j,F^1({\rm End}(M))\otimes Q_j,\nabla_j\bigr)$$
(all tensors products are over $W(k)$). Here, as in 2.2.10, we still denote by $\nabla_j$ the natural connection on ${\rm End}(M)\otimes_{W(k)}Q_j$. 
\smallskip
Similarly, we define ${\got g}_{Q_j}$ by working with ${\got g}$ instead of ${\rm End}(M)$. For any $n\in \NN$, $End(M_{Q_j}/p^nQ_j):=End(M_{Q_j})/p^nEnd(M_{Q_j})$ is a Lie object of
$\Mm\Mf^\nabla_{[-1,1]}(Q_j)$, having (as $\nabla_j$ respects the $G$-action and $r_j\ell_j\in G(Q_j)$) ${\got g}_{Q_j}/p^n{\got g}_{Q_j}$ as a Lie subobject, $j\in\{0,2\}$. As in Fact 3 of 2.2.10 (see also the general principle of 2.2.20 5)), we get: 
\medskip
{\bf Fact.} {\it For any $W(k)$-morphism $z:{\rm Spec}(W(k^1))\to {\rm Spec}(Q_j)$, with $k^1$ a perfect field containing $k$, the pull back of $(M_{Q_j},\nabla_j,(t_{\al})_{\al\in\Mj})$ through $z$ is a Shimura filtered $\sg_{k^1}$-crystal with an emphasized family of tensors
$$\bigl(M\otimes_{W(k)} W(k^1),F^1\otimes_{W(k)} W(k^1),\vph_j^1,G_{W(k^1)},(t_{\al})_{\al\in\Mj}\bigr),$$ 
where $\vph^1_j=g^1_j\vph_0$ for some element $g^1_j\in G(W(k^1))$.}
\medskip
So we can speak about the
Shimura filtered Lie $\sg_{k^1}$-crystal attached to it: its underlying $W(k^1)$-module is ${\rm Lie}(G_{W(k^1)})={\got g}\otimes_{W(k)} W(k^1)\subset {\rm End}(M\otimes_{W(k)} W(k^1))$.
\medskip
{\bf 3.6.4. Generic pull backs.} Considering an arbitrary point $z^j:{\rm Spec}(W(k_j))\to {\rm Spec}(Q_j)$, with $k_j$ a perfect field containing $k$, such that the induced point ${\rm Spec}(k_j)\to{\rm Spec}(Q_j/pQ_j)$ sits over the generic point of ${\rm Spec}(Q_j/pQ_j)$, we get (cf. 3.6.3) a Shimura filtered
$\sg_{k_j}$-crystal 
$${\got C}_j:=\bigl(M\otimes_{W(k)}W(k_j),F^1\otimes_{W(k)}W(k_j),\vph^j,G_{W(k_j)}\bigr), 
$$
with $\vph^j=g^j\vph_0$ for some element $g^j\in G(W(k_j))$. The  Shimura filtered Lie $\sg_{k_j}$-crystal
$$Lie({\got C}_j):=\bigl({\got g}\otimes_{W(k)} W(k_j),\vph^j,F^0({\got g})\otimes_{W(k)} W(k_j),F^1({\got g})\otimes_{W(k)} W(k_j)\bigr)$$ 
has a Newton polygon which does not depend on
$j\in\{0,2\}$ (cf. 3.6.2 d)) and is below the Newton polygon of $({\got g},\vph_0)$ (as ${\got C}_0$ specializes to $(M,F^1,\vph_0)$, cf. 3.6.2 b)). Moreover ${\got C}_2$ specializes to $(M,F^1,g_2\vph_0)$ and so $Lie({\got C}_2)$ specializes to $({\got g},g_2\vph_0)$ (cf. 3.6.2 c)).
\medskip
{\bf 3.6.5. Remarks. 1)} We could have worked to achieve the connection of the Shimura $\sg$-crystal $(M,g_2\vph_0,G)$
to a Shimura $\sg$-crystal expected to be $G$-ordinary working with just one deformation, cf. 3.6.7 below. We preferred to work with two global deformations (over $Q_0$ and $Q_2$) as we hope that this method will be useful in other contexts: we have in mind applications to the construction of integral models of Shimura varieties of special type by using moduli schemes of $p$-divisible groups when possible and when not of generalized Shimura $p$-divisible objects, as well as applications to $p$-adic uniformizations
of Shimura varieties, cf. also 3.6.19 below.
\smallskip
{\bf 2)} Theoretically it can happen that $h^{-1}$ mod $p$ belongs to $\Mb_0(k)$ or that $g_2h^{-1}$ mod $p$ belongs to $\Mb_2(k)$ (cf. the constructions of 3.6.8.2 below). Then we can  not go ahead with 3.6.2-4. So we have to proceed more carefully, cf. 3.6.7 below.
\smallskip
{\bf 3)} We denote the $p$-divisible group $D_j$ as well as the Shimura $p$-divisible group $D(j)$ over ${\rm Spec}(Q_j)$ introduced in 3.6.2 by 
$$
\Md(U_j^0,\Phi_{R_j},M,F^1,h\vph_0,G),
$$ 
in order to emphasize the data needed to construct them. 
\medskip
{\bf 3.6.6. Lemma.} {\it We consider the set $\Ml$ of $k$-valued points of $G$ obtained from the set $D_G$ of 3.2.10 by reduction mod $p$. We have: it is dense in $G_k$.}
\medskip
{\bf Proof:} Let $H$ be the Zariski closure of $\Ml$ in $G_k$. Writing $g=p_1n_1$, with $p_1$ normalizing $F^1$ and with $n_1\in {\rm Ker}(G(W(k))\to G(k))$, we get that $\Ml$ is the image at the level of $k$-valued points of a $k$-morphism 
$$m_H:\tilde H\to H$$ 
of varieties: we take $\tilde H:={P_0}_k\times_k P_k\times U_k$, with $P_k$ as the parabolic subgroup of $G_k$ normalizing $F^1/pF^1$ and with $U_k$ as the affine $k$-variety defined naturally by the $k$-vector space underlying ${\got g}/p{\got g}+F^0({\got g})$; the morphism $m_H$ is easily describable starting from 3.2.10 (FORM) by working mod $p$. As $\tilde H(k)$ is Zariski dense in $\tilde H$ (see [Bo2, 18.3]) we get that $m_H$ is a dominant $k$-morphism. 
\smallskip
For a Galois extension $k_1$ of $k$, we denote by $\Ml_1\subset G(k_1)$ and $H_1\hookrightarrow G_{k_1}$ the analogues of $\Ml$ and $H$ obtained by working over $W(k_1)$ instead of $W(k)$. $\Ml_1$ is stable under the natural action of the Galois group ${\rm Gal}(k_1/k)$ on $G(k_1)$ and moreover $\Ml$ is a subset of the subset of $\Ml_1$ formed by elements fixed by it. In fact the similarly constructed $k_1$-morphism $m_{H_1}$ is nothing else but the extension of $m_H$ to $k_1$ (see the end of the previous paragraph). So if $H_1=G_{k_1}$ then $H=G_k$. So we can assume $k=\bar k$ (in fact for most of below arguments it is enough to assume $G$ is split). This allows us to replace (in the definition of $\Ml$ and $H$) $\vph_0$ by $\vph_1$ (cf. 3.2.3, 3.3.1-2 and the translation part of 3.3.3).
\smallskip
We use the notations of 3.4.0-3. $H$ is $Z(G_k)$-invariant and $Z(G_k)\subset H$. It is enough to work with a subset $I_0$ of $I$ corresponding to a cycle $\gamma_0$ of $\gamma$; we just need to show that $H$ contains the subgroup of $G$ generated by $G_i$, $i\in I_0$. 
\smallskip
For $i\in I_0$, let $N_i^+$ be the integral, connected, smooth, unipotent, abelian subgroup of $G_i$ having $F^1({\got g}_i)$ as its Lie algebra and let $N_{i}^-$ be the integral, connected, smooth, unipotent, abelian subgroup of $G_i$ which is the opposite of $N_i^+$ w.r.t. the action of $T_i$ (via inner conjugation) on $G_i$. So we have
$${\rm Lie}(N_{i}^-)=\bigoplus_{\scriptstyle \al\in\Phi^+_i\atop\scriptstyle
{\tilde g}_i(\al)\subset {\rm Lie}(N_i^+)} {\tilde g}_i(-\al).$$
We consider the following two cases: 
\medskip
{\bf 1)} $\abs{\tilde I_0}=1$, and 
\smallskip
{\bf 2)} $\abs{\tilde I_0}\ge 2$. 
\medskip
We first treat the case 1). So $\tilde I_0=I_1=I_0=\{1\}=\{i\}$. Let ${\rm exp}$ be the exponential map which takes $x\in p{\rm Lie}(N_i^-)$ into the $W(k)$-valued point of $G$ defined by the automorphism $1_M+x$ of $M$. Let $P_0(i):=G_i\cap P_0$. Working with an arbitrary element $p_0\in P_0(i)(W(k))$ and with 
$$g\in \exp\bigl(p{\rm Lie}(N_i^-)\bigr),$$
we deduce (as $\vph_1(p{\rm Lie}(N_{i}^-))={\rm Lie}(N_i^-)$, cf. 3.4.3.2 and the fact that $\tilde I_0$ has only 1 element) that $H$ contains the product of ${P_0(i)}_k$ and of ${N_i^-}_k$; from [Bo2, 14.14] we get that $H$ contains ${G_i}_k$. 
\smallskip
The case 2) is very much the same. For $i\in \tilde I_0$, we define $G_i:=G_{i'}$, where $i'\in I_0$ is such that $n=\abs{I_0}$ divides $i-i'$; similarly we use the same approach mod $n$ for different other groups or Lie algebras or sets of roots indexed by $i\in \tilde I_0$ (so $\Phi_i^+=\Phi_{i'}^+$, etc.). Warning: sometimes, in what follows we specifically state that $i\in I_0$, i.e. sometimes the correct set of indices is $I_0$ and not $\tilde I_0$; when this is not stated explicitly, then $i\in\tilde I_0$.
\smallskip
Let 
$$M_i:=\sum_{\scriptstyle\al\in\Phi_i^-\atop
\scriptstyle{\tilde g}_i(\al)\not\subset {\rm Lie}(N_{i}^-)\oplus {\rm Lie}(P_0(i))} \GG_{a,\al(i)}.$$
Here $\GG_{a,\al(i)}$ is the $\GG_a$ subgroup of $G_i$ having $\tilde g_i(\al)$ as its Lie algebra, $P_0(i):=G_i\cap P_0$, while the sum sign $\sum$ refers to the fact that we take the subgroup of $G_i$ generated by these $\GG_a$ subgroups. $M_i$ is a connected, unipotent, smooth subgroup of $G_i$.
\smallskip
We choose $g=\prod_{i\in\tilde I_0} g_i$, with $g_i\in G_i(W(k))$ of the form $g_i=v_iw_i$, where 
$$v_i\in {\rm exp}\bigl(p{\rm Lie}(N_i^-)\bigr)$$ 
and 
$$w_i\in M_i(W(k)),$$ 
$\forall i\in\tilde I_0$. We also choose 
$$p_0=\prod_{i\in\tilde I_0} p_0^i,$$ 
with $p_0^i\in P_0(i)(W(k))$. For $\tilde g\in G(W(k))$ normalizing $F^1/pF^1$, we denote
$$\vph_1(\tilde g):=\vph_1\tilde g\vph_1^{-1}\in G(W(k)).$$
We have $e\in\{1,2\}$, cf. 3.4.3.1 and the beginning paragraph of 3.4.
To study $\tilde h$ defined by the equality $\tilde h\vph_1=gp_0\vph_1g^{-1}$, we consider two subcases. 
\smallskip
We first assume $\abs{I_0}=1$; so $\abs{\tilde I_0}=2$. We consider the $\ZZ_p$-structure $G_{\ZZ_p}$ of $G$ constructed as in 2.2.9 8) starting from $(M,F^1,\vph_1,G)$. Using it, it makes sense to speak about $\sg$ acting on $W(k)$-valued or $k$-valued points of $G$. We take $v_2$, $w_2$ and $p_0^2$ to be $1_M$. We get $\tilde h=v_1w_1p_0^1\vph_1(w_1^{-1})\vph_1(v_1^{-1})$. So $\tilde h$ mod $p$ is 
$$\bar w_1\bar p_0^1\sg(\bar w_1^{-1})\bar s_1,$$ 
where $\bar w_1$, $\bar p_0^1$ and $\bar s_1$ are the reduction mod $p$ of respectively $w_1$, $p_0^1$ and $\vph(v_1^{-1})$. The role of $\bar s_1$ is that of an arbitrary element of $\sg(N_{1}^-)(k)$. 
\smallskip
We write $\bar s_1=\bar m_1\bar s_{11}$, with $\bar m_1\in M_1(k)\subset \sg(N_1^-)(k)$ and with $\bar s_{11}\in\sg(N_1^-)(k)$. We have $\sg(\bar w_1^{-1})\bar m_1=\bar m_1\sg(\bar w_1^{-1})\tilde s_{11}$, with $\tilde s_{11}\in N_1^-\cap\sg(N_1^-)(k)$; this is an immediate consequence of the formula of [BT, 4.2 (1)]. Using [SGA3, Vol. III, 4.1.2 of p. 172], we can write 
$$\bar p_0^1\bar m_1=\bar m_{11}\bar p_0^{1,1},$$ 
with $\bar m_{11}\in M_1(k)$ and with $\bar p_0^{1,1}\in P_0(1)(k)$, provided we choose $\bar p_0^1$ to be a $k$-valued point of a small enough open subscheme of ${P_0(1)}_k$ containing its origin. We can always choose $\bar m_1$ such that $\bar w_1\bar m_{11}$ is the identity element. Denoting $\bar s_{111}:=\tilde s_{11}\bar s_{11}$, we get that $\tilde h$ mod $p$ can be written as
$$\bar p_0^{1,1}\sg(\bar w_1^{-1})\bar s_{111}.$$
$\bar p_0^{11}$, $\sg(\bar w_1^{-1})$ and $\bar s_{111}$ are running independently through the $k$-valued points of an open subscheme of ${\bar P_0(i)}_k$ containing the origin, of $\sg(M_1)_k$ and respectively of $\sg(N_1^-)_k$. So, from [Bo2, 14.14] we get that $\Ml$ contains the $k$-valued points of an open, dense subscheme of $\bar G_1$; so $H$ contains $\bar G_1$.  
\smallskip
From now on we assume $\abs{I_0}\ge 2$. Let 
$$
h_i:=v_iw_ip_0^i\vph_1(w_{i-1}^{-1})\vph_1(v_{i-1}^{-1}).
$$ 
\indent
It is enough to show that the Zariski closure $Z_H(I_0)$ in $\prod_{i\in I_0} {G^{\rm ad}_i}_k$ of the $k$-valued points of $\prod_{i\in I_0} {G^{\rm ad}_i}_k$ defined naturally by reductions mod $p$ of $\tilde h$'s, is $\prod_{i\in I_0} {G^{\rm ad}_i}_k$ itself. Let $\tilde h_i$ be the image of $h_i$ in ${G^{\rm ad}_i}_k$. If $e=1$, the image of $\tilde h$ in ${G^{\rm ad}_i}_k$ is $\tilde h_i$. If $e=2$ and if $p_0^1$,..., $p_0^n$ are all $1_M$, for $i\in I_0$, the image of $\tilde h$ in ${G^{\rm ad}_i}(k)$ is nothing else but the image of $v_iw_iv_{i+n}w_{i+n}p_0^{i+n}\vph_1(w_{i+n}^{-1})\vph_1(v_{i+n}^{-1})\vph_1(w_{i}^{-1})\vph_1(w_{i}^{-1})$ in ${G^{\rm ad}_i}(k)$ and so is of the form
$$\tilde h_i^1\tilde h_{i+n}\tilde h_i^2,\leqno (EXPR)$$ 
where $\tilde h_i^1$, $\tilde h_i^2\in G^{\rm ad}_i(k)$ depend only on $w_i$ and on $\vph_1(w_{i_1}^{-1})\vph_1(v_{i_1}^{-1})$ (i.e. only on ``ingredients" contributing to the expression of $\tilde h_i$).
\smallskip
Let $\tilde N_i$ be the integral, connected, smooth, unipotent subgroup of $G_i$ having 
$\vph_1(p{\rm Lie}(N_{i-1}^-))$ as its Lie algebra. The role of $\vph_1(v_{i-1}^{-1})$ is that of an arbitrary element $\tilde n_i$ of $\tilde N_i(W(k))$. 
\smallskip
Let $\tilde N_i^0$ be the integral, connected, smooth, unipotent subgroup of $G$ whose Lie algebra is $\vph_1({\rm Lie}(M_{i-1}))$. 
Let 
$$\tilde N_i^1:=M_i\cap\tilde N_i^0,$$
 and let 
$$\tilde N_i^2:=N_i^-\cap\tilde N_i^0.$$
\indent
The groups $\tilde N_i$, $\tilde N_i^0$, $\tilde N_i^1$ and $\tilde N_2$ are unipotent subgroups of $G_i$. $\tilde N_i^2$ is a normal subgroup of $\tilde N_i^0$. Moreover, $\tilde N_i^0$ is the semidirect product of $\tilde N_i^2$ and of $\tilde N_i^1$. The role of $\vph_1(w_{i-1}^{-1})$ is that of an arbitrary element of $\tilde N_i^0(W(k))$.  
So we can write 
$$\vph_1(w_{i-1}^{-1})=\tilde n_i^1\tilde n_i^2,$$
with $\tilde n_i^s\in\tilde N_i^s(W(k))$, $s\in\{1,2\}$. We get
$h_i=v_iw_ip_0^i\tilde n_i^1\tilde n_i^2\tilde n_i$.
\smallskip
Using again [SGA3, Vol. III, 4.1.2 of p. 172] we can write 
$$h_i=v_iw_i\bar n_i^1\tilde p_0^i\tilde n_i^2\tilde n_i,$$ 
with $\bar n_i^1\in\tilde N_i^1(W(k))$ and $\tilde p_0^i\in P_0(i)(W(k))$, provided we choose $p_0^i$ to be a $W(k)$-valued point of an open subscheme of $P_0(i)$ containing the origin and small enough. Loc. cit. can be applied as the product $P_0(i)\tilde N_i^1$ is an open subscheme of the parabolic subgroup of $G_i$ generated by $P_0(i)$ and $\tilde N_i^1$. Fixing some $i\in\tilde I_0$, $\bar n_i^1$, $\tilde p_0^i$, $\tilde n_i^2$ and $\tilde n_i$ are allowed to run independently through open, dense subschemes through which the origin of $G$ factors of subgroups of $G_i$ of whose Lie algebras are forming a direct sum decomposition of ${\got g}_i$; we refer to this property as IND. 
\smallskip 
We write $w_i=\tilde w_i^2\tilde w_i^1$, with $\tilde w_i^1\in\tilde N_i^1(W(k))$ and with $\tilde w_i^2$ a $W(k)$-valued point of the connected, smooth subgroup of $M_i$ of whose Lie algebra is the direct summand of ${\rm Lie}(\tilde N_i^1)$ in ${\rm Lie}(M_i)$ normalized by $T$. We refer to $\tilde w_i^1$ as the component of $w_i$ in $\tilde N_i^1(W(k))$. Similarly, below we refer to the component of any $W(k)$-valued point of $M_{i1}$ in $M_{i2}(W(k))$; here $M_{i1}$ is an arbitrary connected, smooth subgroup of $M_i$ which is the semidirect product of two connected subgroups of it normalized by $T$, one of them being $M_{i2}$. Warning: such components are always written down on the right side of products of two elements. Redenoting $\tilde p_i:=\tilde p_i^0$, we get 
$$
h_i=v_i\tilde w_i^2\tilde w_i^1\bar n_i^1\tilde p_i\tilde n_i^2\tilde n_i.
$$ 
\indent
{\bf Key point.} {\it The $1$ dimensional $B(k)$-vector spaces $\tilde g_i(\al)[{1\over p}]$, $i\in I_0$, which are not included in ${\rm Lie}(P_1)[{1\over p}]={\rm Lie}(P_0)[{1\over p}]$ (see 3.3.2 for this identification), are permuted by $\vph_1$ in such a way that in any resulting cycle there is such a $B(k)$-vector space contained in ${\rm Lie}(N_i^-)[{1\over p}]$, for some $i\in I_1$.}
\medskip
The key point is a consequence of 3.4.3.0 (1) and of the signs of $s_i(\al)$'s as explained in 3.4.3.0 (based on the end of 3.2.3). So we have:
\medskip
{\bf Claim.} {\it $\forall i\in\tilde I_0$, $\tilde h_i$ mod $p$ runs independently through the $k$-valued points of an open, dense subscheme of ${G_i}_k^{\rm ad}$.}
\medskip
To prove the Claim, we fix arbitrarily some $i\in\tilde I_0$; for simplifying the presentation, we assume it is the biggest element of $\tilde I_0$. As $v_i$ is congruent to the identity mod $p$, $\tilde h_i$ can be defined as well as the image of $v_i^{-1}h_i$ in $G_i^{\rm ad}(k)$. We choose $w_i$ to be an arbitrary point of an open, dense subscheme of $M_i$ containing the origin and small enough. We now choose $g_{i-1}$ and $p_0^i$ (i.e. we choose $\bar n_i^1$, $\tilde p_0^i$, $\tilde n_i^2$ and $\tilde n_i$, cf. IND) such that $v_i^{-1}h_i$ is an arbitrary element of an open, dense subscheme of $G_i$ containing the origin. So the Claim holds for $\tilde h_i$ mod $p$. We can choose $g_{i-1}$ (more precisely $\bar n_i^1$) to depend only on the component $\tilde w_i^1$ of $w_i$. Next, we choose inductively on $s\in S(2,\abs{\tilde I_0}-1\})$ elements $g_{i-s}$ and $p_0^{1+i-s}$ such that $v_{i+1-s}^{-1}h_{i+1-s}$ is an arbitrary element of an open, dense subscheme of $G_{i+1-s}$ through which the origin factors. Moreover, we can assume $g_{i-2}$ depends only on the component $\tilde w^{12}_i\in \tilde N_i^{12}(W(k))$ of $\tilde w_i^1$, where $\tilde N_i^{12}$ is the connected, smooth subgroup of $\tilde N_i^1$ whose Lie algebra is $\vph_1({\rm Lie}(\tilde N_{i-1}^1))\cap {\rm Lie}(\tilde N_i^1)$, etc. This takes care of the Claim except for $1$. The above key point implies: $v_{1}^{-1}h_{1}$ is as well an arbitrary element of an open, dense subscheme of $G_{1}$ containing the origin. In other words, a similar choice of $g_i$ (i.e. of $w_i$ and $v_i$) and of $p_0^{1}$ is possible, which does not depend on the initial choice of $w_i$ (i.e. the last component $\tilde w_i^{12...\abs{\tilde I_0}}$ --of $\tilde w_i^{12...\abs{\tilde I_0}-1}$-- we define inductively is $1_M$). So we can choose as the ``new" $w_i$, the initial element $w_i$ we started with. So the Claim holds as well for $1$; this proves the Claim.
\smallskip
From the Claim and (EXPR) we get: $Z_H(I_0)$ is $\prod_{i\in I_0} {G^{\rm ad}_i}_k$. This ends the proof of the Lemma.  
\medskip
{\bf 3.6.6.0. Exercise.} If $g_2\in D_G$ and $g_3\in G(W(k))$ mod $p$ is the identity element, then $(M,g_3g_2\vph_0)$ has $\tau_0$ as its formal isogeny type and $({\got g},g_3g_2\vph_0)$ has the same Newton polygon as $({\got g},\vph_0)$. Hint: use the ideas of the proof of 3.6.6 and of b) of 4.4.1 2) below; by induction on $m\in\NN$ show that we can assume $g_3$ mod $p^m$ is the identity element and so, choosing $m$ big enough, 3.3.4 and [Ka2, 1.4.4] apply. 
\smallskip
This Exercise is referred in what follows just to point out some alternatives or possible shortcuts in 3.6.6.1 2), 4.6.6 and Appendix, and so it is not used before Appendix. See 3.6.17 for a solution of it (deviating to some extent from its hint).
\medskip
{\bf 3.6.6.1. Remarks. 1)} Another proof of 3.6.6 (and actually of the stronger fact that $\Ml\supset U_1(k)$, with $U_1$ an open subscheme of $G_k$ containing its origin; warning: the above proof of 3.6.6 also obtains this stronger fact) can be obtained using 3.4.8, 3.4.11 and 3.6.10 below.
\smallskip
{\bf 2)} A fast proof of 3.1.3 a) and b) together with the inclusion $A_G\subset B_G$ can be obtained by putting together 3.6.6, 3.3.4, 3.4.11 and 3.6.10 below; also: the use of 3.4.11 can be substituted by 3.6.6.0.
\medskip
{\bf 3.6.6.2. Some extensions of 3.6.6.} 3.6.6 can be entirely adapted (cf. its proof) to the generalized Shimura context. One just needs to show that the case $e=3$ can be treated entirely similar. For $e=3$, the subcase $\abs{I_0}\ge 2$ of the Case 2) of 3.6.6 needs no modifications: one just needs to take $p_0^{n+1}$,..., $p_0^{2n}$ to be as well identity elements; the subcase $\abs{I_0}=1$ and $\abs{\tilde I_0}=3$ can be treated entirely similar to the first subcase of the Case 2) of 3.6.6. The exponential map as used in 3.6.6, still makes sense in the generalized Shimura context: it can be used with the same purpose (as we can see mod $p^2$); this is still part of the first place where we need $p\ge 3$. 
\medskip
{\bf 3.6.7. One global deformation.} Section 3.6 started with the desire of connecting the (arbitrary) Shimura filtered $\sg$-crystal $(M,F^1,g_2\vph_0,G)$ with the (extension to $\bar k$ of the) Shimura filtered $\sg$-crystal $(M,F^1,\vph_0,G)$ or with another  Shimura filtered $\bar\sg$-crystal which is --supposed to be-- a $G$-ordinary $\bar\sg$-crystal. This is achieved as follows.
\smallskip
Let $\Mu={\rm Spec}(R_1)$ be an open, affine subscheme of $G$, with $\Mu(W(k))$ containing $1_M$. We assume there is an \'etale $W(k)$-morphism $b_{\Mu}:\Mu\to Y$, with the origin of $G$ factoring through the closed subscheme of $\Mu$ obtained as the inverse image through $b_{\Mu}$ of the closed subscheme of $Y$ defined by $z_i=0$, $i=\overline{1,d}$. If needed, we replace $\Mu$ by a suitable open, affine subscheme of it, so that the condition ii) of the beginning of 3.6 is satisfied for $b_{\Mu}$. Let $\Phi_{R_1}$ be the Frobenius lift of $\Mu^\wedge$ obtained through $b_{\Mu}$ from the Frobenius lift of $Y$ which (at rings' level) takes $z_i$ to $z_i^p$, $i=\overline{1,d}$.
\smallskip
Let $D_1:=\Md(\Mu,\Phi_{R_1},M,F^1,g_2\vph_0,G)$ be the $p$-divisible group over the $p$-adic completion ${\rm Spec}(Q_1)$ of an $\NN$-pro-\'etale, affine scheme over $\Mu$, obtained as in 3.6.1.3, cf. 3.6.5 3). Sometimes, not to complicate the notation, we also refer to $D_1$ as a Shimura $p$-divisible group, the Shimura structure being the natural one.
From 3.6.1.3 5) and 3.6.6 we deduce (as in 3.6.3-4 and 2.3.16) that:
\medskip
{\bf 3.6.7.1. Corollary.} {\it The Shimura $\sg_{k_1}$-crystal ${\got C}_{1}$ associated to the pull back of $D_1$ to the spectrum of the algebraic closure $k_1$ of the field of fractions of $Q_1/pQ_1$, specializes to $(M,g_2\vph_0,G)$ and to a Shimura $\bar {\sg}$-crystal of the form $(M\otimes_{W(k)} W(\bar k),g_3(\vph_0\otimes 1),G_{W(\bar k)})$, with $g_3\in G(W(\bar k))$ which mod $p$ belongs to the set $\Ml(\bar k)$ defined as the set $\Ml$ of 3.6.6 but working over $W(\bar k)$ instead of over $W(k)$.} 
\medskip
{\bf 3.6.8. The proofs of 3.6.1.2-3.} From many points of view this 3.6.8 is the very heart of the whole of 3.6; so, though it uses very concrete situations, the reader should always keep in mind that it captures the very essence of most of what follows after it in 3.6.9-20. We work with just one index $j$: $0$. We have the following 14 steps, supported by the Lemma 3.6.8.2. 
\medskip
{\bf 1) The complete local case.} We recall $\Phi_{\hat R^0_0}$ is induced from $\Phi_{R_0}$.
 There is a unique connection $\nabla_{\hat R^0_0}$ on $M\otimes_{W(k)}\hat R^0_0$ which makes $M_{\hat R_0^0}$ potentially to be viewed as a $p$-divisible object of $\Mm\Mf_{[0,1]}^\nabla(\hat R_0^0)$; it is integrable and nilpotent mod $p$. For these facts cf. [Fa2, th. 10] (and its proof for the uniqueness part): the part of the paragraph before 3.6.1 involving Frobenius lifts allows us to apply loc. cit.
\medskip
{\bf 2) The $G$-action.} $\nabla_{\hat R^0_0}$ respects the $G$-action, i.e. it is of the form $\delta_0+\be$, with 
$$
\be\in\hat R^0_0\otimes_{R_0} ({\got g}\otimes_{W(k)}\Om_{R_0/W(k)})
$$ 
(cf. [Fa2, rm. ii) after th. 10]). Here and in what follows, for any formally \'etale $R_0$-algebra $\tilde R_0$, $\delta_0$ is the (truncation modulo some positive, integral power of $p$ of the) connection on $M\otimes_{W(k)} \tilde R_0$ annihilating $M$.
\medskip
{\bf 3) The general form of connections.} Let $S(M):=S(1,d_M)$. Let $S(G):=S(1,d)$. We choose a $W(k)$-basis $\{e_i|i\in S(M)\}$ of $M$. Let $\{e_{is}|i,s\in S(M)\}$ be the canonical $W(k)$-basis of ${\rm End}(M)$ w.r.t. this chosen $W(k)$-basis of $M$; so $e_{is}(e_s)=e_i$. For the sake of generality, for the time being, we make this choice arbitrarily. Later on we specialize to a proper choice (see 3.6.8.5) suitable for computations. $\Om_{R_0/W(k)}$ is a free $R_0$-module having $\{dt_l(0)|l\in S(G)\}$ as a basis.
\smallskip
Let $\tilde l_0:{\rm Spec}(\tilde Q_0^\wedge)\to {\rm Spec}(R_0)$ be a formally \'etale, affine $W(k)$-morphism, with $\tilde Q_0$ a $W(k)$-algebra. An arbitrary connection $\nabla$ on $M\otimes_{W(k)} \tilde Q_0^\wedge$
can be written in the form 
$$\nabla=\delta_0+\sum_{(i,s,l)\in S(M)\times S(M)\times S(G)} x_{isl}e_{is}\otimes dt_l(0),\leqno (0)$$ 
with $x_{isl}\in\tilde Q_0^\wedge$, $\forall (i,s,l)\in S(M)\times S(M)\times S(G)$. Similarly, if $n\in\NN$, and if $\nabla$ is a connection on $M\otimes_{W(k)} \tilde Q_0/p^n\tilde Q_0$, then we have the same expression for it but with $x_{isl}\in\tilde Q_0/p^n\tilde Q_0$.
\medskip
{\bf 4) The equations.} We start treating the case $n=1$. The condition that a connection $\nabla$ on $M\otimes_{W(k)} \tilde Q_0/p\tilde Q_0$ makes $M_{\tilde Q_0}/pM_{\tilde Q_0}$ potentially to be viewed as an object of $\Mm\Mf_{[0,1]}^\nabla(\tilde Q_0)$, is expressed (cf. 3.6.1.1.1 2)) by the following equations
$$\nabla\circ\Phi(\tilde Q_0)^0(m)=p\Phi(\tilde Q_0)^0\circ d\Phi_{\tilde Q_0*}/p\circ\nabla(m)=0,\leqno (1)$$
if $m\in F^0\otimes_{W(k)} \tilde Q_0/p\tilde Q_0$, and
$$\nabla\circ\Phi(\tilde Q_0)^1(m)=\Phi(\tilde Q_0)^0\circ d\Phi_{\tilde Q_0*}/p\circ\nabla(m),\leqno (2)$$
if $m\in F^1\otimes_{W(k)} \tilde Q_0/p\tilde Q_0$. 
\medskip
{\bf 5) A new way of looking at these equations.} The above equations are replaced, using the chosen form (0) of $\nabla$ (with $x_{isl}\in\tilde Q_0/p\tilde Q_0$), by equations involving $x_{isl}$, $(i,s,l)\in S(M)\times S(M)\times S(G)$, and the Frobenius transforms of some lifts of them to $\tilde Q_0$ taken mod $p$ (i.e. and their $p$-powers). These equations have coefficients in $R_0/pR_0$ and so they define an affine scheme $\Ms^1$ of finite type over $\Ms^0:={\rm Spec}(R_0/pR_0)$. $\Ms^1$ is a moduli scheme of connections. Let 
$$
\ell(0):\Ms^1\to\Ms^0
$$ 
be the $k$-morphism we get. 
\smallskip
 The nice thing is that these equations involving $x_{isl}$'s take the form
$$x_{isl}=L_{isl}(x_{111}^p,x_{112}^p,...,x_{d_Md_Md}^p)+a_{isl},\leqno (3)$$
where the form $L_{isl}$ is homogeneous and linear in each variable and has coefficients in $R_0/pR_0$ and where $a_{isl}\in R_0/pR_0$. We have such an equation for each triple $(i,s,l)\in S(M)\times S(M)\times S(G)$. This nice shape (3) of the equations involving $x_{isl}$'s is obtained as follows. Let $i\in S(M)$. Choosing a $W(k)$-basis $\{\tilde e_1^u,...,\tilde e_{\dim_{W(k)}(F^u)}^u\}$ of $F^u$, $u=\overline{0,1}$, we can write 
$$
e_i=\sum_{u=0}^1\sum_{q=1}^{\dim_{W(k)}(F^u)} a_q(u)\Phi(R_0)^u(\tilde e_q^u),\leqno (DIV)
$$ 
with all $a_q(u)$'s as elements of $R_0^\wedge$ (we recall: $M_{R_0}$ is a $p$-divisible object of $\Mm\Mf_{[0,1]}(R_0)$). Summing up as suggested by (DIV) the equations obtained by plugging $\tilde e_q^0$'s and $\tilde e_q^1$'s in (1) and respectively in (2) of 4), we get the mentioned shape for the equations involving $x_{isl}$, $(s,l)\in S(M)\times S(G)$, in the left hand side of (3). By doing this for any $i\in S(M)$, we obtain ``all" equations in the variables $x_{isl}$ (i.e. any other equation produced by (1) and (2) is a linear combination of the ones of (3)).  
\medskip
{\bf 6) The uniqueness argument.} The existence of another (i.e. different) connection on $M\otimes_{W(k)} \tilde Q_0/p\tilde Q_0$ making $M_{\tilde Q_0}/pM_{\tilde Q_0}$ potentially to be viewed as an object of $\Mm\Mf_{[0,1]}^\nabla(\tilde Q_0)$ corresponds to a non-trivial solution (with values in $\tilde Q_0/p\tilde Q_0$) of the system of equations 
$$x_{isl}=L_{isl}(x_{111}^p,x_{112}^p,...,x_{d_Md_Md}^p),\leqno (4)$$
$(i,s,l)\in S(M)\times S(M)\times S(G)$. Let $\Mi$ be the ideal of $R_0/pR_0$ generated by $t_l(0)$, $l\in S(G)$. The coefficients of the linear forms $L_{isl}$ are elements of $\Mi^{p-1}$. This is so due to the fact that $\Phi_{R_0}(t_l(0))=t_l(0)^p$. 
\smallskip
To prove 3.6.1.2 for $n=1$, we can assume $\tilde Q_0$ is a complete, local ring with the maximal ideal generated by $\Mi$ and residue field $\bar k$. By induction on $q\in\NN$, we get that any solution of the above system (4) of equations, belongs to $\Mi^{p-1+pq}\tilde Q_0/p\tilde Q_0$. So we get that any such solution must be the trivial one, given by $x_{isl}=0$, $\forall (i,s,l)\in S(M)\times S(M)\times S(G)$. This proves 3.6.1.2, for $n=1$.
\medskip
{\bf 7) The ``non-empty" \'etale part.} From the shape of the equations of (3) we get $\Om_{\Ms^1/\Ms^0}=\{0\}$. On the other hand, the criterion of formal smoothness is satisfied for $\ell(0)$. Argument: we need to show that for any $R_0/pR_0$-algebra $R(0)$ and for every ideal $I(0)$ of it such that $I(0)^2=\{0\}$, any solution $(y_{111}^0,...,y_{d_Md_Md}^0)$ of (3) in $R(0)/I(0)$ lifts to a solution of (3) in $R(0)$; but if $(y_{111},...,y_{d_Md_Md})$ is an arbitrary $d_M^2d$-tuple formed by elements of $R(0)$ and lifting $(y_{111}^0,...,y_{d_Md_Md}^0)$, then 
$$\bigl(L_{111}(y_{111}^p,...,y_{d_Md_Md}^p)+a_{111},...,L_{d_Md_Md}(y_{111}^p,...,y_{d_Md_Md}^p)+a_{d_Md_Md}\bigr)$$ 
is a solution of (3) in $R(0)$ lifting $(y_{111}^0,...,y_{d_Md_Md}^0)$. 
\smallskip
So $\Ms^1$ is \'etale over $\Ms^0$. Let $\Ms^1_0$ be the maximal open closed subscheme of $\Ms^1$ with the property that any connected component of it has a non-empty intersection with the closed subscheme $\ell(0)^{-1}(\Mz_k)$ of $\Ms^1$. From 6) we get that $\ell(0)^{-1}(\Mz_k)$ is either empty or ${\rm Spec}(k)$. So $\Ms_0^1$ is a connected, affine $k$-scheme. There are to ways to see that it is non-empty. For the first way we just need to point out: from 1) we deduce the existence of a natural formally \'etale $R_0/pR_0$-morphism 
$$m_0^0:{\rm Spec}(\hat  R_0^0/p\hat R_0^0)\to\Ms^1$$ 
factoring through $\Ms_0^1$. For the second way, avoiding the use of 1), see 8) below.
\smallskip
Lifting now $\Ms^1_0$ to a smooth, affine $W(k)$-scheme (this is standard), we get a formally \'etale, affine $W(k)$-morphism $\ell_{0,1}:{\rm Spec}(Q_{0,1}^\wedge)\to {\rm Spec}(R_0^\wedge)$, with $Q_{0,1}$ a smooth $W(k)$-algebra such that ${\rm Spec}(Q_{0,1}/pQ_{0,1})=\Ms^1_0$, which lifts the $k$-morphism $\Ms^1_0\to \Ms^0$ defined naturally by the $k$-morphism $\ell(0)$ of 5). From its construction and from 6) we deduce that ii) and  iii) of 3.6.1.3 1) are satisfied for $n=1$ (the fact that the connection $\nabla_0^1$ on $M\otimes_{W(k)} Q_{0,1}/pQ_{0,1}$ we get is integrable, nilpotent mod $p$ and respects the $G$-action results from the existence of $m_0^0$, cf. 1)). 
\smallskip
Next we present the second way to see that $\Ms^1_0$ is non-empty.
\medskip
{\bf 8) The ``non-empty" \'etale part: a new approach.} $M_{R_0}$ is the pull back of $M_{S_0}$ through the natural inclusion $i_0^0:U_0^0\hookrightarrow W^0_0$ respecting the Frobenius lifts. So to prove that $\Ms^1_0$ is non-empty, we can work with $M_{S_0}$ instead of $M_{R_0}$, i.e. we can assume
$G=GL(M)$ and $R_0=S_0$. So $d=d_M^2$. But in this case in 5) we have $d^2=d_M^4$ equations in $d^2$ variables. So the fact that $S_0^1$ is non-empty is handled by the below Lemma (applied --cf. the \'etaleness part of 7)-- to the reduction of (3) of 5) modulo $\Mi$); we return to the proofs of 3.6.1.2-3 in 3.6.8.2, after a short intermezzo exploiting it to a greater extend.
\medskip
{\bf 3.6.8.1. Fundamental Lemma.} {\it Let $m,l\in\NN$ and let $k_1$ be a field of characteristic $p>0$. We consider an affine $k_1$-scheme $\My_0$ defined via a system of $m$ equations
in $m+s$ variables $x_1,...,x_{m+s}$, $s\in\NN\cup\{0\}$, of the form
$$x_i=L_i(x_1,...,x_{m+s})^{p^l}+c(i),\leqno (5)$$  
$i=\overline{1,m}$, where $L_i$ is a homogeneous polynomial of degree $1$ with coefficients in $k_1$ and $c(i)\in k_1$, $\forall i\in S(1,m)$. Then  $\My_0$ is a non-empty, smooth $k_1$-scheme of pure dimension $s$. 
\smallskip
If $s=0$ then $\My_0$ is \'etale over $k_1$ and the number of elements of the set $\My_0(\overline{k_1})$ is finite and equal to $p^{lm_1}$, with $m_1\in S(0,m)$. Moreover, $m_1$ depends only on the coefficients of the homogeneous linear forms $L_i$, $i\in S(1,m)$, and not on the free coefficients $c(i)$, $i\in S(1,m)$. We also have: $m_1$ is less or equal to the rank $r_L$ of the $m\times m$ matrix $A_L$ obtained using the coefficients of $L_i$ as its rows, $i\in S(1,m)$.}
\medskip
{\bf Proof:} It is trivial to check (as in 7)) the criterion of formal smoothness for the $k_1$-morphism $\My_0\to {\rm Spec}(k_1)$. As $\My_0$ is of finite type over $k_1$ we deduce $\My_0$ is smooth over $k_1$. $\Om_{\My_0/k_1}$ is a free $\Mo_{\My_0}$-sheaf of rank $s$: $dx_{m+1}$,..., $dx_{m+s}$ is a basis of it; so $\My_0$ is of pure dimension $s$. For the other things we can assume $s=0$; so $\My_0$ is \'etale over $k_1$. For the part involving $\My_0(\overline{k_1})$, we can assume $k_1$ is algebraically closed. 
\smallskip
We use mathematical induction. The case  $m=1$ is obvious. Let now $m\ge 2$. We introduce a new variable $x_0$, so that the equations become homogeneous:
$$x_ix_0^{p^l-1}=L_i(x_1,...,x_m)^{p^l}+c(i)x_0^{p^l}.\leqno (6)$$ 
These equations are defining a closed subscheme $\My$ of $\PP_{k_1}^m$. By repeatedly using the projective dimension theorem (cf. [Ha, p. 48]), we get it is non-empty. If it has no point in the hypersurface of $\PP^m_{k_1}$ defined by $x_0=0$, then we are done: $\My_0=\My$ and $\abs{\My_0(k_1)}=p^{ml}$, cf. the intersection theory. If $\My$ has a point with $x_0=0$, then we get that the homogeneous polynomials $L_i$ of degree 1, $i=\overline{1,m}$, are linearly dependent, i.e. there are constants $a(i)\in k_1$, not all $0$, such that $\sum_{i=1}^m a(i)L_i=0$. So also 
$$\sum_{i=1}^m a(i)^{p^l}L_i^{p^l}=0.$$ 
We can assume $a(m)$ is different from $0$ and so we can assume it is 1. 
\smallskip
So the initial system of (non-homogeneous) equations is equivalent to the system of equations, where we preserve the first $m-1$ equations, while we replace the last equation (obtained for $i=m$) by the linear equation
$$x_m+\sum_{i=1}^{m-1} a(i)^{p^l}x_i=c(m)+\sum_{i=1}^{m-1} c(i)a(i)^{p^l}.$$
This last equation allows us to eliminate $x_m$. We come across a system of $m-1$ equations in $m-1$ variables of the same type as the original one. Moreover, the coefficients of the new homogeneous linear forms in the variables $x_1$, $x_2$,..., $x_{m-1}$ depend on $L_i$, $i\in S(1,m)$, but do not depend on $c(i)$, $i\in S(1,m)$. So we can proceed by induction. The last part about the estimate of $m_1$ is obvious: as above, we can eliminate $m-r_L$ variables at once. This proves the Lemma.
\medskip
{\bf 3.6.8.1.0. More on $m_1$.} $m_1$ can be smaller than $r_L$. In fact easy examples show that, if $r_L<m$, then $m_1$ can take any value of the set $S(0,r_L)$. Here, we just mention one example to illustrate this phenomenon. We take $m=2$ and we consider the system defined by the following two equations $x_1=x_1^{p^l}-x_2^{p^l}$ and $x_2=x_1^{p^l}-x_2^{p^l}$. So $r_L=1$. But $m_1=0$ as the only solution of this system of equations is the pair $(x_1,x_2)=(0,0)$.
\smallskip
On the other hand, if $r_L=m$ (i.e. if $A_L$ is an invertible matrix) then the system (5) of equations for $s=0$ is equivalent to a system of equations of the form
$$y_i^{p^l}=\tilde L_i(y_1,...,y_m),\leqno (7)$$
$i=\overline{1,m}$, with $\tilde L_i$ as non-homogeneous linear forms; so in this case $\My=\My_0$ and $m_1=m$.
\medskip
{\bf 3.6.8.1.1. Corollary.} {\it Let $m\in\NN$ and let $(E,\tau)$ be a pair comprising from a field $E$ and an automorphism $\tau$ of it. Let $L_i(x_1,...,x_m)$ be linear forms in $m$ variables with coefficients in $E$, $i=\overline{1,m}$. We consider the difference system $DS$ of equations
$$y_i=L_i(\tau(y_1),...,\tau(y_m)),\leqno (8)$$
$i=\overline{1,m}$. We have: $DS$ is non-empty, i.e. it has a solution in the difference closure of $(E,\tau)$.}
\medskip
{\bf Proof:}
We recall that the difference closure is a pair $(E_1,\tau_1)$, with $E_1$ a field containing $E$ and with $\tau_1$ an automorphism of it extending $\tau$, such that each difference variety over $(E,\tau)$ (see [Hr, \S 5] for a definition) has a point in $(E_1,\tau_1)$ and there is no subpair $(E_2,\tau_2)$ of it, with $E_1\neq E_2$, having all these properties. Its existence and uniqueness up to isomorphism can be proved in the same way as in the case of the algebraic closure of a field.
\smallskip
We can assume $(E,\tau)$ is equal to its difference closure. The emerging philosophy of E. Hrushovski and A. Macintyre (see [Hr]) allows us to assume (via ultraproducts of fields) that $E=\FF$ and that $\tau$ is an integral, non-zero power $q$ of its Frobenius automorphism having $\FF_p$ as its fixed field. As the subfield of $\FF$ generated by the coefficients of $L_i$'s is finite, we can assume $q>0$. So the Corollary follows from 3.6.8.1.
\medskip
{\bf 3.6.8.1.2. Corollary.} {\it Let $m,n\in\NN$. Let $R(0)$ be an $\FF_p$-algebra. Let $L_i(x_1,x_2,...,x_m)$ be homogeneous linear forms in $m$ variables with coefficients in $R(0)$, $i=\overline{1,m}$. We have:
\smallskip
{\bf a)} {\bf (the abstract surjectivity principle)} For any $m$-tuple $(c(1),...,c(m))$ of elements of $R(1)$, with $R(1)$ an $\FF_p$-algebra containing $R(0)$, the system of equations 
$$x_i=L_i(x_1^{p^n},x_2^{p^n},...,x_m^{p^n})+c(i),\leqno (9)$$
$i=\overline{1,m}$, defines an \'etale, affine ${\rm Spec}(R(1))$-scheme ${\rm Spec}(R(2))$. The resulting morphism $m(2):{\rm Spec}(R(2))\to {\rm Spec}(R(1))$ is surjective. 
\smallskip
Moreover, there is an open, dense subscheme $U(0)$ of ${\rm Spec}(R(0))$ such that, regardless of the choice of the $m$-tuple $(c(1),...,c(m))$, $m(2)$ defines an \'etale cover above $U(1):={\rm Spec}(R(1))\times_{{\rm Spec}(R(0))} U(0)$.
\smallskip
{\bf b)} {\bf (the abstract constructibility property)} Let $R(\infty)$ be the $\NN$-inductive limit of $\FF_p$-algebras $R(q)$, $q\in\NN$, $q\ge 2$, with $R(q)$ as the $R(q-1)$-algebra defined by a system of equations of the form of (9) but with $c(i)$ replaced by some $c(i)_{q}\in R(q-1)$, $i=\overline{1,m}$. Then the image of any connected (resp. irreducible) component of ${\rm Spec}(R(\infty))$ in ${\rm Spec}(R(1))$ is a constructible, dense subset of a connected (resp. irreducible) component of ${\rm Spec}(R(1))$. Moreover the image of any open closed subscheme of ${\rm Spec}(R(\infty))$ in ${\rm Spec}(R(1))$ is an open subscheme. 
\smallskip
{\bf b')} {\bf (the $AG$ property)} If $R(0)$ is of finite type over $k$, then ${\rm Spec}(R(\infty))$ is an $AG$ $k$-scheme and so has the $ALP$ property.
\smallskip
{\bf c)} If $R(0)$ is a strictly henselian local ring, then the additive map 
$$
\Ml(n): R(0)^m\to R(0)^m
$$ 
defined by 
$$\Ml(n)(x_1,...,x_m):=\bigl(x_1-L_1(x_1^{p^n},...,x_m^{p^n}),...,x_m-L_m(x_1^{p^n},...,x_m^{p^n})\bigr)
$$
is surjective. Moreover, all its fibres have the same number of elements equal to $p^{nm_1}$, for some $m_1\in S(0,m)$.}
\medskip
{\bf Proof:}
We first prove a). As the number of the coefficients of the forms $L_1$,..., $L_m$ is finite, we can assume $R(0)$ is a finitely generated $\FF_p$-algebra. In particular $R(0)$ is a noetherian ring. As in 7) we get that the system (9) of equations (of 3.6.8.1.2) defines an \'etale, affine morphism $m(2):{\rm Spec}(R(2))\to {\rm Spec}(R(1))$. 3.6.8.1 implies it is surjective. For the last part, we can assume $R(0)$ is an integral domain. So we just have to show that $U(2):={\rm Spec}(R(2))\times_{{\rm Spec}(R)} U(0)$ is finite over $U(1)$, with $U(0)$ a non-empty, open, affine subscheme of ${\rm Spec}(R(0))$ which depends only on the forms $L_i$, $i=\overline{1,m}$. We can also assume that the $p^n$-th roots of the coefficients of the linear forms $L_i$ are in $R(0)$. So we can write $L_i(x_1^{p_n},...,x_m^{p^n})=\tilde L_i(x_1,...,x_m)^{p^n}$, with $\tilde L_i$ as homogeneous linear forms. 
\smallskip
We follow the proof of 3.6.8.1. We proceed by induction on $m\in\NN$. The case $m=1$ is obvious. Let now $m\ge 2$. If the linear forms $\tilde L_i$ are linearly dependent over the field of fractions $K(0)$ of $R(0)$, by localizing $R(0)$ w.r.t. some non-zero element of it, we can eliminate one variable as in the proof of 3.6.8.1; so the induction applies. If the linear forms $\tilde L_i$ are linearly independent over $K(0)$, then they are linearly independent, over any point of a non-empty, open, affine subscheme $U(0)$ of ${\rm Spec}(R(0))$. So the systems of equations defining $U(2)$ over $U(1)$ can be put in a similar form to 3.6.8.1.0 (7) and so $U(2)$ is a finite $U(1)$-scheme; we get: $U(2)$ is an \'etale cover of $U(1)$. This proves the a) part.
\smallskip
b) and b') are a consequence of a) using standard noetherian induction on the topological space underlying the spectrum of the $\FF_p$-subalgebra of $R(0)$ generated by the coefficients of $L_1$,..., $L_m$. 
\smallskip
c) is a direct consequence of 3.6.8.1 and of the \'etale property expressed in a). This ends the proof of the Corollary.
\medskip
{\bf 3.6.8.1.3. Definitions.} We assume $R(0)$ is reduced. The maximal open subscheme $U(0)$ of ${\rm Spec}(R(0))$ having the property of 3.6.8.1.2 a) is called the Artin--Schreier open subscheme defined by the linear homogeneous forms $L_1$,..., $L_m$. Using the reductions of $L_1$,..., $L_m$ modulo the ideal of $R(0)$ defining the complement ${\rm Spec}(R(0))\setminus U(0)$ with its reduced scheme structure, we get similarly an Artin--Schreier open subscheme $U(1)$ of ${\rm Spec}(R(0))\setminus U(0)$: it is dense. In this way, by induction, we get a stratification of ${\rm Spec}(R(0))$ in reduced, locally closed subschemes; as the $\FF_p$-subalgebra of $R(0)$ generated by the coefficients of the forms $L_1$,..., $L_m$ is noetherian, this stratification has a finite number of strata $U(q)$, $q=\overline{1,n_L}$, where $n_L\in\NN$ depends on $L_i$'s. We call it the Artin--Schreier stratification of ${\rm Spec}(R(0))$ defined by the linear forms $L_1$,..., $L_m$. 
\smallskip
Let now $c(1)$,..., $c(m)\in R(0)$ be arbitrary elements. We take $R(1)$ to be the $R(0)$-algebra defined by the system $(9)$ of equations; it is reduced as $R(0)$ is. 
We similarly define the refined Artin--Schreier stratification of ${\rm Spec}(R(0))$ defined by the linear homogeneous forms $L_1$,..., $L_m$. The only difference: instead of $U(0)$, we take its open closed subscheme $U(0)^\prime$ such that the number of geometric points of ${\rm Spec}(R(1))$ above any geometric point $y$ of $U(0)^\prime$ is maximal (i.e. does not depend on $y$ and is strictly greater than the number of geometric points of ${\rm Spec}(R(1))$ above geometric points of ${\rm Spec}(R(0))\setminus U(0)^\prime$). So, as $R(1)$ is an \'etale $R(0)$-algebra, it is the stratification indexed by the number of geometric points of ${\rm Spec}(R(1))$ we get above geometric points of ${\rm Spec}(R(0))$ (see also [EGA IV, 15.5.9]). Its strata are similarly indexed: $U(q)^\prime$, $q=\overline{1,n_L^\prime}$, with $n_L^\prime\in\NN$. We have $n_{L^\prime}\ge n_L$. 
\smallskip
The nice think about above two stratifications of ${\rm Spec}(R(0))$ is: their strata are naturally indexed.
\medskip
{\bf 3.6.8.1.4. Lemma (the first form of the purity principle).} {\it All Artin--Schreier stratifications and all refined Artin--Schreier stratifications satisfy the purity condition.}
\medskip
{\bf Proof:} We can assume $R(0)$ is a finitely generated, reduced $\FF_p$-algebra. The following well known result allows us to assume $R(0)$ is as well integral and normal.
\medskip
{\bf Exercise.} Let $S$ be an affine, reduced scheme. An open subscheme $U$ of it is affine iff the normalization $U^n$ of $U_{\rm red}$ in a finite, normal extension of the ring of fractions of $S$ is. Hint: We need to show $U$ is affine if $U^n$ is. This is a consequence of Chevalley's theorem (see [Ha, Exc. 4.2 of p. 222]). Loc. cit. is stated with lots of unneeded hypotheses. Using the fact that the morphism $U^n\to U$ is surjective, we get that $U$ is quasi-compact. So using [EGA IV, \S 8] we can assume $S$ is finitely generated over a field; so loc. cit. applies (cf. also [Ma, 31.H]). 
\medskip
Using induction, we can assume $q=0$. Let $K(0)$ be the field of fractions of $R(0)$. Let $R(0)^n$ be the normalization of $R(0)$ in the ring of fractions $K(1)$ of $R(1)$. 
${\rm Spec}(R(1))$ is an open, affine subscheme of ${\rm Spec}(R(0)^n)$. So its complement $C(1)$ is either empty or of pure codimension 1 (easy consequence of [Ma, th. 38, p. 124]). The same holds for its image in ${\rm Spec}(R(0))$ which is ${\rm Spec}(R(0))\setminus U(0)$. As $R(0)$ is integral, $U(0)=U(0)^\prime$. This ends the proof.  
\medskip
{\bf 3.6.8.1.5. Remark.} 3.6.8.1 and 3.6.8.1.1-4 remain valid if $L_i$'s are not necessarily (homogeneous) linear forms (in $m$ variables) but sums of monomials which are some $p$-powers (not necessarily equal) of just one of the variables involved. The reason is: we can introduce new variables and equations of the form $y_i=x_i^p$, so that the new system of equations is of the same form as the one of 3.6.8.1.2 a), with $n=1$ and (potentially) with an increased number of variables. 
\medskip
We come back to the proof of 3.6.1.2-3.
\medskip
{\bf 3.6.8.2. Lemma (the explicit form of the surjectivity principle).} {\it There is a reduced, closed subscheme $\Mb_0(k)$ of $G_k$ not containing the origin, depending only on the linear forms $L_{isl}$ and not on the coefficients $a_{isl}$, $(i,s,l)\in S(M)\times S(M)\times S(G)$, and such that the fibres of $\ell_{0,1}$ above points of ${\rm Spec}(R_0/pR_0)$ not belonging to $\Mb_0(k)$, are non-empty.}
\medskip
{\bf Proof:} 
All schemes to be  introduced below are reduced.
There is a closed subscheme $\Mr(0)$ of ${\rm Spec}(R_0/pR_0)$, different from ${\rm Spec}(R_0/pR_0)$, and such that denoting by $\Mt(0)$ its open complement, the morphism $\ell_{0,1}^{-1}(\Mt(0))\to\Mt(0)$ is an \'etale cover. By elementary matrix operations we see that $\Mr(0)$ depends only on the linear forms $L_{isl}$ and not on the coefficients $a_{isl}$ (cf. 3.6.8.1.2 a) and its proof). We start constructing $\Mb_0(k)$.
\smallskip 
In what follows $\Mj(\tilde i)$ are finite set of indices, $\tilde i=\overline{0,d-1}$. All irreducible components of $\Mr(0)$ not passing through the origin of $G_k$ will be locally closed subschemes $S_u(0)$, $u\in\Mj(0)$, of $\Mb_0(k)$. For any irreducible component of $\Mr(0)$ passing through the origin of $G_k$, we have to repeat the argument. Let $\Mc^1$ be the union of all such irreducible components, viewed as a connected, closed subscheme of ${\rm Spec}(R_0/pR_0)$. There is (cf. 3.6.8.1.2 a)) an open, dense subscheme $\Mt(1)$ of $\Mc^1$, again depending only on the reduction of $L_{isl}$'s modulo the ideal of $R_0/pR_0$ defining $\Mc^1$ and not on the reduction modulo it of the coefficients $a_{isl}$, $(i,s,l)\in S(M)\times S(M)\times S(G)$, such that the $k$-morphism $\ell_{0,1}^{-1}(\Mt(1))\to\Mt(1)$ is an \'etale cover. Let $\Mr(1)$ be the subscheme of $\Mc^1$ defined by the complement of $\Mt(1)$ in $\Mc^1$. The irreducible components of $\Mr(1)$ not passing through the origin of $G_k$, will be subschemes $S_u(1)$, $u\in\Mj(1)$, of $\Mb_0(k)$. We keep going with those irreducible components of $\Mr(1)$ passing through the origin of $G_k$. 
\smallskip
By induction on $\tilde i\in S(0,d-1)$ we construct ($d$ is the dimension of $G_k$):
\medskip
{\bf a)} a connected, closed subscheme $\Mc^{\tilde i+1}$ of ${\rm Spec}(R_0/pR_0)$ of dimension at most $d-\tilde i-1$, and which is either empty or contains the origin (so $\Mc^d$ is either empty or is $\Mz_k$);
\smallskip
{\bf b)} a maximal open, dense subscheme $\Mt(\tilde i)$ of $\Mc^{\tilde i}$, with $\Mc^0:={\rm Spec}(R_0/pR_0)$, such that the $k$-morphism $\ell_{0,1}^{-1}(\Mt(\tilde i))\to\Mt(\tilde i)$ is an \'etale cover;
\smallskip
{\bf c)} closed subschemes $S_u(\tilde i)$ of ${\rm Spec}(R_0/pR_0)$, $u\in\Mj(\tilde i)$; they are the irreducible components not containing the origin of the closed subscheme of ${\rm Spec}(R_0/pR_0)$ defined as the complement $\Mr(\tilde i)$ of $\Mt(\tilde i)$ in $\Mc^{\tilde i}$. 
\medskip
The formula which allows us to apply the inductive argument is: $\Mc^{\tilde i+1}$ is the union of the irreducible components of $\Mr(\tilde i)$ containing the origin. So, from constructions, $\Mr(d-1)$ is either empty or has dimension $0$. 
\smallskip
$\Mb_0(k)$ is defined as the Zariski closure in $G_k$ of the (finite) union of these subschemes $S_u(\tilde i)$, $u\in\Mj(\tilde i)$, with $\tilde i$ running from 0 to $d-1$. It does not contain the origin of $G_k$.
\smallskip
From its construction we deduce that it does not depend on the coefficients $a_{isl}$, $(i,s,l)\in S(M)\times S(M)\times S(G)$, and that for any perfect field $k_1$ containing $k$ we have $\Mb_0(k_1)=\Mb_0(k)_{k_1}$. Moreover, we obviously have:
\medskip
{\bf Fact.} {\it The disjoint union $\cup_{\tilde i\in S(0,d)} \Mt(\tilde i)$ is a stratification of the complement of $\Mb_k$ in ${\rm Spec}(R_0/pR_0)$. Here $\Mt(d)$ is either empty or is $\Mz_k$, depending on the fact that there is or there is not an $\tilde i\in S(0,d-1)$ such that $\Mz_k$ is a closed subscheme of $\Mt(\tilde i)$. The Zariski closure of any non-empty $\Mt(\tilde i)$ contains $\Mz_k$.}
\medskip
From the Fact and b) the Lemma follows.
\medskip
{\bf 9) End of the proof for the case $n=1$.} From the fact that the coefficients of $L_{isl}$ are in $\Mi$ we deduce that the closed subscheme $\ell_{0,1}^{-1}(\Mz_k)$ of $\Ms_0^1$ is ${\rm Spec}(k)$. This implies iv) of 3.6.1.3 1). Related to the part of 3.6.1.3 2) referring to $\ell_{0,1}$, we just need to show that the fibres of $\ell_{0,1}$ over geometric points have at most $p^{d^2}$ points. This can be checked as follows. The equations of (3) of 5), are obtained working with the $W(k)$-basis $\{e_{is}|i,s\in S(M)\}$ of ${\rm End}(M)$. Using a second $W(k)$-basis of ${\rm End}(M)$ such that a subbasis of it is a $W(k)$-basis of $\got g$, as we are dealing with connections respecting the $G$-action, our variables $x_{ijl}$'s are as well solution of a system of $d_M^2d-d^2$-linear equations with coefficients in $R_0/pR_0$ whose rank is precisely $d_M^2d-d^2$ in each point of ${\rm Spec}(R_0/pR_0)$. As in the proof of 3.1.8.1, we can eliminate $d_M^2d-d^2$-variables at once: so $\Ms_0^1$ is a connected component of an ${\rm Spec}(R_0/pR_0)$-scheme obtained using a system of equations with (coefficients in $R_0/pR_0$ and with) the property that a subsystem of it is entirely similar to (3) of 5) but involving just $d^2$ variables; so 3.6.8.1 applies to give us the desired bound on the number of points of geometric fibres of $\ell_{0,1}$. This ends the proof of 3.6.1.3 1) and 2) involving $n=1$.
\medskip
{\bf 10) The inductive statement.} By induction on $n\in\NN$ we prove 3.6.1.2 and 3.6.1.3 1). We assume we managed to construct $Q_{0,m}$, for any $m\in\NN$ smaller than $n+1$. We now construct $Q_{0,n+1}$ and show that 3.6.1.2 holds for $n+1$.
\medskip
{\bf 11) The general form of a lift.} Let $\tilde\ell_0:{\rm Spec}(\tilde Q_0^\wedge)\to {\rm Spec}(Q_{0,n}^\wedge)$ be a formally \'etale, affine $W(k)$-morphism, with $\tilde Q_0$ a $W(k)$-algebra. The general form of a connection on $M\otimes_{W(k)} \tilde Q_0/p^{n+1}\tilde Q_0$ lifting the scalar extension of $\nabla_0^n$, is of the form
$$\nabla=\nabla_0^n+\sum_{(i,s,l)\in S(M)\times S(M)\times S(G)} p^nx_{isl}e_{is}\otimes dt_l(0),\leqno (10)$$ 
with all $x_{isl}$'s belonging to $\tilde Q_0/p\tilde Q_0$; here we still denote by $\nabla_0^n$ a fixed lift of the mentioned scalar extension to a connection on $M\otimes_{W(k)} \tilde Q_0/p^{n+1}\tilde Q_0$ respecting the $G$-action.
\medskip
{\bf 12) The key point.} The condition that this connection $\nabla$ is a connection on $M_{\tilde Q_0}/p^{n+1}M_{\tilde Q_0}$ is expressed by a system of  equations involving $x_{isl}$'s, which is of the same type as the system (3) of  5). The key point is:
\medskip
{\bf Fact.} {\it We get the same linear forms $L_{isl}$, $\forall (i,s,l)\in S(M)\times S(M)\times S(G)$!}
\medskip
However, the new feature is that the free coefficients $a_{isl}$ of 5) are replaced by new coefficients $a_{isl}(n)\in Q_{0,n}/pQ_{0,n}$. The fastest way to see this Fact and the new feature, is to use the same trick as in 8) (for $G=GL(M)$ they are obvious). 
\smallskip
As the coefficients $a_{isl}$ have played no role in the steps 5) to 9), we can repeat the arguments of 5) to 9). So we get a moduli scheme $\Ms^{n+1}$ of connections, which is an \'etale scheme over $\Ms_0^n:={\rm Spec}(Q_{0,n}/pQ_{0,n})$. We define $\Ms_0^{n+1}$ in the similar manner: it is the maximal open closed subscheme of $\Ms^{n+1}$ such that each connected component of it intersects the pull back of $\Mz_k$ through the natural $k$-morphism $\Ms^{n+1}\to {\rm Spec}(R_0/pR_0)$. 
\smallskip
This repetition takes care of 3.6.1.2 and of 3.6.1.3 1) for $n+1$ (again iv) is implied by i), while 1) implies --as in 7)-- that $\nabla_0^{n+1}$ is integrable, nilpotent mod $p$, and respects the $G$-action).
\medskip
{\bf 13) Some conclusions.} 3.6.1.3 4) results from 3.6.1.3 1). 3.6.1.3 3) is obvious. The fact that the morphism $\ell_0^n$ mod $p$ is: 
\medskip
-- \'etale, results from the same argument used in 7) (cf. also 3.6.8.1.2 a));
\smallskip
-- of finite type results from the construction of ${\rm Spec}(Q_{0,n+1}/pQ_{0,n+1})$ (it is a closed subscheme of a scheme of finite type over $\Ms_0^n$);
\smallskip
-- quasi-finite, with geometric fibres of cardinality at most $p^{d^2}$, results from the fact that $\nabla_0^{n+1}$ respects the  $G$-action (and so we actually have only $d^2$ variables) and from the intersection theory (cf. 9)).
\medskip
This takes care of 3.6.1.3 2).
\medskip
{\bf 14) End of the proof.} 3.6.1.3 5) and 6) are a consequence of the proof of 3.6.8.2 and of Fact of 12) (cf. 3.6.8.1.2 a) and b')). More precisely, for any $n\in\NN$ we construct an open subscheme $\Ms_1^n$ of $\Ms_0^n$, such that, under the $W(k)$-morphism $\ell_0^n$, $\Ms_1^{n+1}$ projects surjectively onto $\Ms_1^n$, as well as $\Ms^1_1$, under the $W(k)$-morphism $\ell_{0,1}$, projects surjectively onto the open subscheme of ${\rm Spec}(R_0/pR_0)$ defined by the complement of $\Mb_0(k)$. 
\smallskip
$\Ms^n_1$ is formed by removing from $\Ms_0^n$ the closed subscheme formed by the union of all irreducible components of $\ell_{0,n}^{-1}(\Mc^1_{\tilde i})$ not containing the closed subscheme $\ell_{0,n}^{-1}(\Mz_k)={\rm Spec}(k)$, for $\Mc^1_{\tilde i}$ running through all irreducible components of the scheme $\Mr(\tilde i)$ defined in 3.6.8.2, with $\tilde i$ as its turn running through the set $S(0,d-1)$.
\smallskip
We check now that $\Ms_1^{n+1}$ projects surjectively onto $\Ms_1^n$ (the argument for the fact that $\Ms_1^1$ projects surjectively onto the open subscheme of ${\rm Spec}(R_0/pR_0)$ defined by the complement of $\Mb_0(k)$ is entirely the same). First, any integral subscheme of $\Ms_1^{n+1}$ containing (resp. not containing) $\ell_{0,n+1}^{-1}(\Mz_k)={\rm Spec}(k)$ is mapped through $\ell_0^n$ into an integral subscheme of $\Ms_1^n$ containing (resp. not containing) $\ell_{0,n}^{-1}(\Mz_k)={\rm Spec}(k)$. So we do get a natural \'etale $k$-morphism 
$$\bar\ell_0^n:\Ms_1^{n+1}\to \Ms_1^n.$$ 
To check it is surjective, it is enough to show that any $\bar k$-valued point $y_n$ of $\Ms_1^n$ lifts to a $\bar k$-valued point of $\Ms_1^{n+1}$. 
\smallskip
Let $\tilde i_n\in S(0,d)$ be such that the image of $y_n$ in ${\rm Spec}(R_0/pR_0)$ is a point of $\Mt(\tilde i_n)$ (see the notations of 3.6.8.2). As $\bar\ell_0^n$ is an \'etale $k$-morphism (and so its image is open), as $(\ell_{0,m})^{-1}(\Mt(\tilde i_n))\to \Mt(\tilde i_n)$, with $m\in\{n,n+1\}$, is an \'etale cover (cf. 3.6.8.2 b) and 3.6.8.1.2 a)), and as $\Ms_1^{n+1}$ is a non-empty $k$-scheme (we have $\ell_{0,n+1}^{-1}(\Mz_k)={\rm Spec}(k)\subset \Ms_1^{n+1}$), we get that $y_n$ lifts to a $\bar k$-valued point of $\Ms_1^{n+1}$. 
\smallskip
This ends the proof of 3.6.1.2-3.
\medskip
{\bf 3.6.8.3. A digression.} What follows is a digression meant just for specialists; it is neither quoted nor used anywhere else, besides the isolated remark 3.6.8.4 4). 
\smallskip
A natural question arises: how to compute the connections $\nabla_j^n$ and so, how to construct the $W(k)$-morphisms $\ell_j^n$ (write down the equations defining $Q_{j,n}$)? This is not very difficult (cf. 3.6.8 5); see also 3.6.8.5 and 3.6.18.4 below). Here, as a digression,  we show how using Shimura varieties of $A_n$ type ($n\in\NN$) we can construct a formally \'etale, affine $W(k)$-morphism ${\rm Spec}(\tilde Q_{j,1}^\wedge)\to {\rm Spec}(R_j^\wedge)$ which factors through $\ell_{j,1}$. We can assume $k$ is the algebraic closure of a countably generated field over $\FF_p$; so in particular we have $k=\bar k$. We work again with only one index $j$: 0.
\medskip
{\bf 1)} $M_{R_0}$ is the pull back of $M_{S_0}$ through the inclusion $i_0^0:U_0^0\hookrightarrow
W^0_0$. So we can work with $M_{S_0}$ instead of $M_{R_0}$, i.e. we can assume from now on that
$G=GL(M)$ and that $R_0=S_0$. 
\medskip
{\bf 2)} Substituting $M$ with $\tilde M:=M\oplus W(k)$ and the direct sum decomposition $M=F^1\oplus F^0$ with  $\tilde M=F^1\oplus {\tilde F}^0$, where ${\tilde F}^0:=F^0\oplus W(k)$, we replace the Shimura $\sg$-crystal $(M,h\vph_0, GL(M))$ with the Shimura $\sg$-crystal $(\tilde M,\tilde\vph,GL(\tilde M))$ (with $\tilde\vph$ acting as $h\vph_0$ on $M$ and as $\sg$ on the $W(k)$ summand of $\tilde M$). This is achieved in the following 4 items.
\medskip
a) We choose arbitrarily a $W(k)$-basis $\{v_1,...,v_{d_M}\}$ of $M$; it gives birth to a $W(k)$-basis $\{v_1$,...,$v_{d_M+1}\}$ of $\tilde M$, where $v_{d_M+1}$ is 1 of the summand $W(k)$ of $\tilde M$.
\smallskip
b) We view $GL(M)$ as a subgroup of $GL(\tilde M)$ (with $b_M\in GL(M)(W(k))$ acting as $b_M$ on the summand $M$ of $\tilde M$ and as the identity on the other summand $W(k)$ of $\tilde M$).
\smallskip
c) We choose ${\tilde S}_0$ to be an open, affine subscheme of $GL(\tilde M)$ of the form $S_0\times\GG_a^{2d_M}\times\GG_m$ (the $\GG_a$ factors correspond to elements of $GL(\tilde M)(W(k))$, viewed as matrices, having 1 on the diagonal and 0 anywhere else, except at the intersection of an $i_1$-th row with an $i_2$-th column, where $i_1$ and $i_2$ are distinct elements of the set $\{1,2,...,d_M+1\}$ such that the intersection $\{i_1,i_2\}\cap\{d_M+1\}$ is non-empty; the $\GG_m$ factor is the subgroup of $GL(\tilde M)$ normalizing the summand $W(k)$ of $\tilde M$ and fixing $M$). In order that such a choice of $\tilde S_0$ to be possible, we might have to replace $S_0$ by a suitable open, affine subscheme of it (cf. [SGA3, Vol. III, 4.1.2 of p. 172]; so $S_0\times\GG_a^{2d_M}\times\GG_m$ is embedded in $GL(\tilde M)$ using products of elements, ordered adequately).
\smallskip
d) We choose a Frobenius lift of the $p$-adic completion of ${\tilde S}_0$ of the form $\Phi_{S_0}\times\Phi_1$ (corresponding to the above product decomposition of ${\tilde S}_0$), with $\Phi_1$ as the Frobenius lift of the $p$-adic completion of $\GG_a^{2d_M}\times\GG_m={\rm Spec}(W(k)[x_1,...,x_{2d_M+1}][{1\over x_{2d_M+1}}])$ which takes $x_s$ to $x_s^p$, $\forall s\in S(1,2d_M)$, and $x_{2d_M+1}-1$ to $(x_{2d_M+1}-1)^p$.
\medskip
{\bf 3)} The same argument which  allowed  us to assume $G$ is $GL(M)$, allows us to replace $GL(M)$ with $GL(\tilde M)$, if $p|d_M$. So we can assume $p$ does not divide $d_M$. 
\medskip
{\bf 4)} We remark that there is a SHS $(f,L_{(p)},v)$ defined by an injective map $f:(G_0,X)\hookrightarrow\bigl({\rm GSp}(W,\psi),S\bigr)$ (here we put exceptionally a lower right index
$0$, i.e. we write $G_0$, $2e_0=\dim_\QQ(W)$ and $\Mj_0^\prime$, not to create confusion with the meaning of $G=GL(M)$, $e_M=\dim_{W(k)}(G)=\dim_{W(k)}(GL(M))$ and $\Mj^\prime$) such that we have a natural monomorphism
$$G=GL(M)\hookrightarrow G_{0W(k)},$$ 
resulting into a natural identification  $G_{0W(k)}^{\rm der}=SL(M)$ (so the confusion between $G$ and $G_{0W(k)}$ would have not been too serious), and having the property that the Shimura Lie pairs attached to the Shimura group pairs $(G,[\mu])$ and $(G_1,[\tilde\mu^\ast_x])$,
are isomorphic (cf. def. 2.2.5 3) and 4)). Here $G_1$ is the subgroup of ${G_0}_{W(k)}$ generated by $G_{0W(k)}^{\rm der}$ and by the image of any cocharacter $\tilde\mu^\ast_x:\GG_m\to G_{0W(k)}$, whose extension, under an $O_{(v)}$-monomorphism $W(k)\hookrightarrow \CC$,
 is $G_0(\CC)$-conjugate to the cocharacters $\mu^\ast_x$, $x\in X$. 
\smallskip
So $G^{\rm der}_{0\RR}$ must be $SU(p_0,q_0)_{\RR}$, with $p_0:=\dim_{W(k)}(F^1)$ and $q_0:=\dim_{W(k)}(F^0)=\dim_{W(k)}(M/F^1)$. We can take $G^{\rm der}_0$ to be an arbitrary  $\QQ$--form of $G_{0\RR}^{\rm der}$ which
splits over $\QQ_p$ (to be compared with Example 4 of 4.6 below). [Va2, 6.4.2] guarantees the existence of a SHS as above if $p\ge 5$; the existence of such a SHS for $p=3$ is not a problem: we can assume we have a standard PEL situation $(f,L_{(p)},v,\Mb)$, cf. 2.3.5.1. For instance, if $p_0=1$ and $q_0=2$, we can work with a Picard surface ${\rm Sh}(G_0,X)$ (cf. 
[Go]) having a reflex field $E(G_0,X)$ in which $p$ is totally split. The general case is just an obvious extension of the case $p_0=1$, $q_0=2$ (we have $p_0q_0\ne 0$ due to 3.2.2.1). In what follows, as $p$ does not divide $d_M$, we assume $(f,L_{(p)},v)$ is obtained as in [Va2, 6.5.1.1]. 
\smallskip
We have a natural \'etale isogeny $G_1\to G=GL(M)$ (as $p$ does not divide $d_M$). It allows us to identify the completion of $G_1$ in its origin with the completion ${\rm Spec}(O_0)$ of $G$ in its origin. So ${\got g}_1:={\rm Lie}(G_1)={\got g}$.
\medskip
{\bf 5)} We now assume the existence of $z\in\Mn_{W(k)}/H_0(W(k))$ (cf. the standard notations 2.3.1-3 for the SHS introduced in 4)), such that its attached Shimura filtered
$\sg$-crystal $\bigl(M_z,F^1_z,\vph_z,{G_0}_{W(k)},(t^z_\al)_{\al\in\Mj^\prime_0}\bigr)$ (cf. 2.3.10 and the assumption $k=\bar k$) has the property that the Shimura filtered Lie $\sg$-crystal $\bigl({\got g}_1,\vph_z,F^0_z({\got g}_1),F^1_z({\got g}_1)\bigr)$, which is naturally a Lie $p$-divisible subobject of its attached Shimura filtered $\sg$-crystal (as ${\got g}_1$ is naturally a Lie subalgebra of ${\rm Lie}({G_0}_{W(k)})$), is isomorphic to the Shimura filtered Lie $\sg$-crystal
(${\rm End}(M),h\vph_0,F^0({\rm End}(M)$), $F^1({\rm End}(M)))$.
\medskip
{\bf 6)} From the assumed existence of $z$, we get that ${\got g}_{R_0}={\rm End}(M)_{S_0}$ is ($p$-adically) algebraizable. To explain what we mean by this, we start with a formally smooth $W(k)$-morphism $m_z:{\rm Spec}(O_0)\to\Mn_{W(k)}/H_0$ lifting $z$ and such that the pull back (through its composite with the natural morphism $m_{\rm nat}:\Mn_{W(k)}/H_0\to\Mn/H_0$) of the $p$-divisible group $\Md_{H_0}$ of $\Ma_{H_0}$ and of de Rham components of the family of Hodge cycles $\Ma_{H_0}$ is naturally endowed with), has associated to it a Shimura filtered $F$-crystal 
$$
{\got C}:=\bigl(M_z, F^1_z,\vph_z,G_{0W(k)},G_1,\bar f,(t^z_\al)_{\al\in\Mj^\prime_0}\bigr)
$$ 
over $O_0/pO_0$ (with $\bar f:O_0\tilde\to W(k)[[z_1,...,z_{e_M}]]$ the $W(k)$-isomorphism naturally induced by the $W(k)$-morphism $c_0$ of ii) before 3.6.0). The existence of $m_z$, is just [Va2, 5.4] (cf. the first two paragraphs of 2.3.11). Using Artin's approximation theorem (as in the proof of 2.3.15), we get:
\medskip
{\bf 7)} a smooth morphism
${\rm Spec}(S_1)\buildrel {\ell_1}\over\rightarrow\Mn_{W(k)}/H_0$ and a formally \'etale, affine morphism ${\rm Spec}(S^\wedge_1)\buildrel {\ell}\over\rightarrow W_1^0:={\rm Spec}(S_0)\times_G G_1$ such that the $p$-divisible object of
$\Mm\Mf^\nabla_{[0,1]}(S_1)$ associated to $(m_{\rm nat}\circ\ell_1)^{-1}(\Md_{H_0})$ can be put in the form
$$
{\got C}_{\ell_1}=\bigl(M_z\otimes_{W(k)} S^\wedge_1, F^1_z\otimes_{W(k)} S^\wedge_1, h_1(\vph_z\otimes 1),\nabla_z\bigr),
$$ 
with $h_1\in {G_0}_{W(k)}(S^\wedge_1)$ which mod $p$ is defined by  
$\ell$ (we view ${\rm Spec}(S_0)\times_G G_1$ as an open subscheme of $G_1$ and so as a locally closed subscheme of ${G_0}_{W(k)}$); if $pr_1:{\rm Spec}(S_1^\wedge)\to {\rm Spec}(S_0^\wedge)$ is the $p$-adic completion of the composite of $\ell$ with the projection $W_1^0\to {\rm Spec}(S_0)$, then we choose the Frobenius lift $\Phi_{S_1}$ of $S_1^\wedge$ such that the diagram
$$
\def\mapright#1{\smash{
\mathop{\longrightarrow}\limits^{#1}}}
\def\mapdown#1{\Big\downarrow
\rlap{$\vcenter{\hbox{$\scriptstyle#1$}}$}}
\matrix{{\rm Spec}(S^\wedge_1)&\mapright{pr_1} &{\rm Spec}(S_0^\wedge)\cr
\mapdown{\Phi_{S_1}}&&\mapdown{\Phi_{S_0}}\cr
{\rm Spec}(S^\wedge_1)&\mapright{pr_1}&{\rm Spec}(S_0^\wedge)\cr}$$
is cartesian. Warning: the fact that we chose $\Phi_{S_1}$ to be obtained via $\ell$ and $pr_1$, means that we know that $h_1$ only mod $p$ is defined by $\ell$;
\medskip
{\bf 8)} a point ${\rm Spec}(W(k))\buildrel{d_0}\over\hookrightarrow {\rm Spec}(S_1)$ such that $l_1\circ d_0=z$ and $\ell\circ d_0^\wedge$ is the origin of $W^0_1\subset G_1$;
\medskip
{\bf 9)} modulo a power $I^n_{d_0}$ ($n\in\NN\bsl\{1\})$ of the ideal $I_{d_0}$ of $S_1$ defining $d_0$, ${\got C}_{\ell_1}$ is isomorphic to ${\got C}$ modulo the $n$-th power of the ideal of $O_0$ defining the origin of ${\rm Spec}(O_0)$. 
\medskip
7) to 9) are just a variant of 2.3.15: the fact that we are not mentioning the polarizations, allows us (see also 2.3.14) to work such a variant in the context of ${G_0}_{W(k)}$ itself and not only of ${G_0^0}_{W(k)}$; by removing extra variables (in the same way we got 2.3.15.1 from 2.3.15), we can assume $h_1$ mod $p$ is defined by $\ell$.
\medskip
{\bf 10)} So the Lie $p$-divisible subobject of the Lie $p$-divisible object $End({\got C}_{\ell_1})$ of $\Mm\Mf^\nabla_{[-1,1]}(S_1)$, whose underlying $S^\wedge_1$-module is ${\got g}\otimes_{W(k)} S^\wedge_1$, is 
$${\got g}_{S_1}:=\bigl({\got g}\otimes_{W(k)} S^\wedge_1,h_1(\vph_z\otimes 1),F^0({\got g})\otimes_{W(k)} S^\wedge_1,F^1({\got g})\otimes_{W(k)} S^\wedge_1,\nabla_z\bigr).$$
\medskip
{\bf 11)} So the pull back of $End(M_{S_0})/pEnd(M_{S_0})={\got g}_{S_0}/p{\got g}_{S_0}$ through $pr_1$, can be viewed as a Lie object (identifiable with ${\got g}_{S_1}/p{\got g}_{S_1}$) of $\Mm\Mf^\nabla_{[-1,1]}(S_1)$. This  implies the existence of an \'etale $k$-morphism ${\rm Spec}(S_1/pS_1)\to\Ms_0^1$ which composed with $\ell_{0,1}$ (see 3.6.8 5) for notations) is the $k$-morphism ${\rm Spec}(S_1/pS_1)\to {\rm Spec}(S_0/pS_0)$ defined by taking $pr_1$ mod $p$. 
\smallskip
In other words, we can regain the connection of $M_{S_1}/pM_{S_1}$ from the one of $M_z\otimes_{W(k)} S_1^\wedge$ and from the one of $End(M_{S_1})/pEnd(M_{S_1})$. In fact, from the proof of [Va2, 6.5.1.1], as ${G_0}_{\ZZ_p}$ is a split group, we get that (cf. also 2.3.16.1 for the passage from the Lie algebras context of 5) to the level of $GL$-groups) we can assume there is a monomorphism 
$$i_z:(M,F^1,h\vph_0)\hookrightarrow (M_z,F_z^1,\vph_z)$$
between $p$-divisible objects of $\Mm\Mf_{[0,1]}(W(k))$, 
with $i_z(M)$ as a direct summand of $M_z$ normalized by ${G_0}_{W(k)}$ and such that the image of ${G_0}_{W(k)}$ in $GL(i_z(M))$ is $GL(i_z(M))$ itself.   
\medskip
{\bf 12)} If there is no $z\in\Mn_{W(k)}/H_0(W(k))$ as in 5), we have to choose
$$z_0\in\Mn_{W(k)}/H_0(W(k))$$ 
such that its attached Shimura filtered $\sg$-crystal
$\bigl(M_{z_0},F^1_{z_0},\vph_{z_0},{G_0}_{W(k)}\bigr)$ has the property that its corresponding Shimura Lie
$\sg$-crystal $\bigl({\got g},\vph_{z_0}, F^0_{z_0}({\got g}),F^1_{z_0}({\got g})\bigr)$
is isomorphic to a Shimura filtered Lie $\sg$-crystal 
$$
\bigl({\rm End}(M),g_3\vph_0, F^0(
{\rm End}(M)), F^1({\rm End}(M))\bigr), 
$$
with $g_3\in GL(M)(W(k))$ a $W(k)$-valued point of a ``small" open subscheme $U_3$ of $U_0\cap U_2$ containing $h$ (cf. 2.3.16). Let $U^0_3$ be the right translate of
$U_3$ by $h^{-1}$. 
\smallskip
Going once more through the arguments of 5) to 11) but performed just mod $p$, we can use this part involving a ``near by" point $z_0$, to get that the spectrum $\tilde\Mz$ of the normalization of the $S_0/pS_0$-subalgebra of $\hat R_0^0/p\hat R_0^0$ (we recall that $R_0=S_0$) generated by the coefficients of $\be$ of 3.6.8 1) taken mod $p$, is  generically \'etale over ${\rm Spec}(S_0/pS_0)$. 
\smallskip
As some trouble might result if $\Ms^1$ is not connected or if (theoretically) it is of dimension greater than $e_M$, this statement on $\tilde\Mz$ can be seen in two ways. Either by recalling that the system (4) of equations of 3.6.8 6) gives birth to an \'etale scheme over ${\rm Spec}(S_0/pS_0)$ (so somehow we go back to the general ideas) or by using the Kodaira--Spencer map (here is the place where we need $U_3$ to be ``small": we can assume that this map defined using the coefficients of the mentioned $\be$ still has a image of the right rank in $g_3$ mod $p$) and an argument of gluing (performed around a suitable point of ${U_3}_k$) entirely similar to the one of 3.6.2.1 d) but performed mod $p$. The need of working mod $p$ is implied by the fact that we do not know a priori that the spectrum of the $S_0$-algebra generated by the coefficients of $\be$, is smooth over $W(k)$ in a (suitable) $W(k)$-valued point of it lifting $g_3$; on the other hand we can assume that $g_3$ mod $p$ lifts to a smooth point of $\tilde\Mz$. So the lifting property of 3.6.18.5.2 below used in a relative context (see 3.6.18.7 below) implies that the gluing argument of 3.6.2.1 d) can be performed (i.e. taken) mod $p$ as well. 
\smallskip
We do not detail fully this second way here, as actually a point $z$ as in 5) always exists: as $k=\bar k$, 3.6.15 A below applies, as we can assume $\Mn$ has the completion property (cf. 4.12.12.0 below). We just point out that it is based on the fact that we can ``regain" ${\rm Spec}(S_0/pS_0)$, from the (coefficients w.r.t. some $k$-basis of $M/pM$ of the) Frobenius endomorphism of the $S_0/pS_0$-module underlying $M_{S_0}/pM_{S_0}$. The above mentioned gluing argument performed mod $p$, is the only way we can presently think of showing that $\tilde\Mz$ is an \'etale ${\rm Spec}(R_0/pR_0)$-scheme, without using the general theory (of equations) of 3.6.8 but just using the fact (here is the place where $z_0$ --and so $g_3$-- plays a role) that it has a suitable section (cf. 5) to 11) performed mod $p$) in the \'etale topology.
\medskip
{\bf 3.6.8.4. Remarks. 1)} The special form of the Frobenius endomorphism of $M\otimes_{W(k)} R_0$ played no role in the proofs of 3.6.1.2-3. This allows us to use the method employed in 3.6.8 in (all) other situations (for instance in 3.6.18.4 below).
\smallskip
{\bf 2)} In the proof of the surjectivity principle (for instance, of 3.6.8.2), the special shape of $U_0$ played no role. So it can be used in all other situations.
\smallskip
{\bf 3)} The Frobenius lift $\Phi_{R_0}$ of $R_0^\wedge$ did play a role in the above proofs of 3.6.1.2-3; it is not entirely essential: if we are dealing with an arbitrary Frobenius lift of $R_0^\wedge$, (DIV) of 3.6.8 5) still applies to give us systems of equations of the same shape as of (3) of 3.6.8 5). However, for the case of a general type of a Frobenius lift of $R_0^\wedge$, 3.6.1.3 i) and all the other parts of 3.6.1.3 depending on it have to be reformulated (so 3.6.1.2 has to be reformulated as well), see 3.6.18.8.3 below.
\smallskip
{\bf 4)} Let $(f,L_{(p)},v)$ be a SHS defined by an injective map $f:(G_0,X)\hookrightarrow ({\rm GSp}(W,\psi),S)$; in order to refer to 3.6.8.3, we use $G_0$ instead of $G$. If we already know that the Shimura filtered $\sg$-crystal $(M,F^1,h\vph_0,G)$ is isomorphic to the Shimura filtered $\sg$-crystal attached to a point $z:{\rm Spec}(W(k))\to\Mn_{W(k)}/H_0$, then 3.6.8.3 is significantly simplified: the replacement of $G$ by $GL(M)$ (see 3.6.8.3 1) to 3)) is not needed. This implies that (for 3.6.8.3) we do not need the property ii) of the beginning of 3.6, as well as in 3.6.0 we do not need the ``complication" with $W_0$, $W_2$, $W_0^0$, $W_2^0$, $\Phi_{S_0}$, $\Phi_{S_2}$, etc. 
\smallskip
{\bf 5)} We do not know if (or when) we have ${\rm Spec}(Q_{0,n}/pQ_{0,n})=\Ms^n$. We do expect this equality (at least in the majority of cases).
\medskip
{\bf 3.6.8.5. The closed subscheme $\Mb_0(k)$.} It is very much desirable to have a good description of the closed subscheme $\Mb_0(k)$ of $G_k$ introduced in 3.6.8.2. In what follows we explain a simpler approach for determining its  intersection with $U_0^0$.
\smallskip
It is enough to determine $\Mb_0(\bar k)$ (cf. 3.6.8.2); so we can assume $k=\bar k$. We choose the $W(k)$-basis $\{e_i|i\in S(M)\}$ of $M$ in such a way that: 
\medskip
{\it There is a permutation $\tilde\pi$ of the set $S(M)$ with the property that $h\vph_0(e_i)$ mod $p$ is either $0$ or is $e_{\tilde\pi(i)}$ mod $p$.}
\medskip
The existence of such a $W(k)$-basis is guaranteed by [Mu, cor. p. 143], applied to the $\sg$-linear map $M/pM\to M/pM$ defined by $h\vph_0$ mod $p$, and by the standard Jordan form of a nilpotent $\sg$-linear endomorphism (this is the same as in the $k$-linear context).
\medskip
{\bf 1)} Rewriting the equations of 3.6.8 5) skillfully, for the above chosen $W(k)$-basis of $M$, we come across equations of the form 
$$t_l(0)^{p-1}x_{isl}^p=\tilde L_{isl}(x_{11l},x_{12l},...,x_{d_Md_Ml}) +\tilde b_{isl},\leqno (11)$$
$\forall (i,s,l)\in S(M)\times S(M)\times S(G)$ such that $h\vph_0(e_s)$ is zero mod $p$ and $h\vph_0(e_i)$ is non-zero mod $p$, and of the form
$$\tilde L_{isl}(x_{11l},...,x_{d_Md_Ml})=\tilde b_{isl},\leqno (12)$$
$\forall (i,s,l)\in S(M)\times S(M)\times S(G)$ such that $h\vph_0(e_s)$ is non-zero mod $p$ or $h\vph_0(e_s)$ and $h\vph_0(e_i)$ are both $0$ mod $p$. Here $\tilde L_{isl}$ is a homogeneous linear form and $\tilde b_{isl}\in R_0/pR_0$, $\forall (i,s,l)\in S(M)\times S(M)\times S(G)$.
\smallskip
These equations are obtained by first plugging $e_s$ in the equations (1) and (2) of 3.6.8 4) and then identifying the coefficients of $r_0\circ\tilde l_0(e_i)dt_l(0)$, $(i,s,l)\in S(M)\times S(M)\times S(G)$ (we recall that $\Phi(\tilde Q_0)^0=r_0\circ\tilde l_0(h\vph_0\otimes 1)$, cf. 3.6.1.1.1).
\medskip
{\bf 2)} Denoting 
$$y_{isl}:=t_l(0)x_{isl},$$ 
$\forall (i,s,l)\in S(M)\times S(M)\times S(G)$, this system of equations (defined by (11) and (12)) when restricted to $R_0/pR_0[{1\over {\prod_{i=1}^d t_i(0)}}]$ is equivalent to a system of equations of the form
$$y_{isl}^p=\tilde L_{isl}(y_{11l},...,y_{d_Md_Ml})+t_l(0)\tilde b_{isl},\leqno (13)$$
if $(i,s)\in S(M)\times S(M)$ is such that $h\vph_0(e_s)$ mod $p$ is $0$ and $h\vph_0(e_i)$ is non-zero mod $p$, and of the form
$$\tilde L_{isl}(y_{11l},...,y_{d_Md_Ml})=t_l(0)\tilde b_{isl},\leqno (14)$$
if $(i,s)\in S(M)\times S(M)$ is such that $h\vph_0(e_s)$ mod $p$ is non-zero or $h\vph_0(e_s)$ and $h\vph_0(e_i)$ are both $0$ mod $p$; here $l$ runs through the elements of $S(G)$.
\smallskip
This new system of equations is simpler than the one in 3.6.8 5). It can be used to get a very good description of the intersection of $\Mr(0)$ (of the proof of 3.6.8.2) with ${\rm Spec}(R_0/pR_0[{1\over {\prod_{i=1}^d t_i(0)}}])$. 
\smallskip
{\bf 3)} Similarly, for any subset $S(t)$ of $S(G)$, reducing the system of equations of 2) to $R_0/(p,(t_i(0))_{i\in S(t)})[{1\over {\prod_{i\in S(G)\setminus S(t)} t_i(0)}}]$, we get a system of equations similar to the one in 2), involving $y_{isl}$, $(i,s,l)\in S(M)\times S(M)\times\bigl(S(G)\setminus S(t)\bigr)$. It is easier to determine the intersection of $\Mb_0(k)$ with  ${\rm Spec}(R_0/(p,(t_i(0))_{i\in S(t)})[{1\over {\prod_{i\in S(G)\setminus S(t)} t_i(0)}}])$ using it instead of another one, similar to the one of 3.6.8 5). This motivates why we get simplifications if we work under the convenience assumption of 3.6.1.1.3.
\smallskip 
{\bf 4)} Let now $\tilde G$ be a smooth subgroup of $G$ having connected fibres and such that the Lie algebra of its generic fibre is normalized by $h\vph_0$. The whole of 3.6.1.3 can be redone in this context (cf. 3.6.1.6), with the group $G$ replaced by $\tilde G$. Warning: not to complicate the story, below we use the same notations, with $G$ being replaced everywhere by $\tilde G$.
\medskip 
{\bf 3.6.8.6. The matrix form of the equations of 3.6.8.} For future references, and for small applications in 3.6.8.7 and 3.6.18.7 below, we include here as well the matrix form of the equations of 3.6.8 5). We consider a $W(k)$-basis $\{e_1,...,e_{d_M}\}$ of $M$ as in 3.6.8 2), such that $\{e_1,...,e_{\dim_{W(k)}}(F^1)\}$ is a $W(k)$-basis of $F^1$. 
\smallskip
For $l\in S(G)$, let $X_l$ be the square matrix whose entries are $x_{isl}$'s. Let $E_G$ be the square matrix with entries in $R_0/pR_0$ of the $R_0/pR_0$-endomorphism of $M\otimes_{W(k)} R_0/pR_0$ that takes $e_i$ mod $p$ into $\Phi(\tilde Q_0)^1(e_i)$ mod $p$ if $i\in S(1,\dim_{W(k)}(F^1))$ and into $\Phi(\tilde Q_0)^0(e_i)$ mod $p$ if $i\in S(1+\dim_{W(k)}(F^1),d_M)$, computed w.r.t. the $R_0/pR_0$-basis of $M\otimes_{W(k)} R_0/pR_0$ defined naturally by the chosen $W(k)$-basis $\{e_1,...,e_{d_M}\}$ of $M$. Let $M_{h\vph_0}$ be the square matrix of the $k$-endomorphism of $M/pM$ that takes $e_i$ mod $p$ into $h\vph_0({1\over p}e_i)$ mod $p$ if $i\in S(1,\dim_{W(k)}(F^1))$ and into $h\vph_0(e_i)$ mod $p$ if $i\in S(1+\dim_{W(k)}(F^1),d_M)$, computed similarly. We view $M_{h\vph_0}$ as having entries in $R_0/pR_0$. Let 
$$\pi(F^0):M_{d_M}(R_0/pR_0)\to M_{d_M}(R_0/pR_0)$$
be the $R_0/pR_0$-linear projector which takes a matrix $A\in M_{d_M}(R_0/pR_0)$ having entries $(a_{ij})_{i,j\in S(1,d_M)}$ into the matrix whose entries $(\tilde a_{ij})_{i,j\in S(1,d_M)}$ are defined by the rule: $\tilde a_{ij}$ is $a_{ij}$ if $i>\dim_{W(k)}(F^1)$ and $j\le \dim_{W(k)}(F^1)$ and is $0$ otherwise.
\smallskip
We write
$$d(E_G)=\sum_{l=1}^d C_ldt_l(0),$$
with $C_l$'s as square matrices with entries in $R_0/pR_0$ and with the $d$-operator acting on a matrix via its action on the entries. So now, the equations of 3.6.8 5) can be rewritten in the matrix form as follows
$$C_l+X_lE_G=E_G\pi(F^0)\bigl(t_l(0)^{p-1}X_l^{[p]}\bigr),\leqno (MF)$$
$\forall l\in S(G)$. Warning: for this simplified matrix form (MF), we do need that $\Phi_{S_0}(t_l(0))=t_l(0)^p$, $\forall l\in S(G)$. As $M_{R_0}$ is a $p$-divisible object, $E_G$ and $M_{h\vph_0}$ are invertible matrices. So (MF) can be put as well in the following convenient matrix form
$$C_lE_G^{-1}+X_l=E_G\pi(F^0)\bigl(t_l(0)^{p-1}X_l^{[p]}\bigr)E_G^{-1},\leqno (CMF)$$
$\forall l\in S(G)$.
\medskip
{\bf 3.6.8.6.1. Relative forms of the equations of 3.6.8.}  Let $D_G:=E_GM_{h\vph_0}^{-1}$. As the direct sum decomposition $M=F^1\oplus F^0$ is defined by the cocharacter $\mu$ of $G$, we have
$$D_G\in G(R_0/pR_0)$$
and 
$$\pi(F^0)({\got g}/p{\got g}\otimes_k R_0/pR_0)\subset {\got g}/p{\got g}\otimes_k R_0/pR_0.$$
For this last inclusion we use the standard identification of $M_{d_M}(R_0)$ with ${\rm End}(M\otimes_{W(k)} R_0)$ defined naturally by the chosen $W(k)$-basis of $M$.
So $C_lM_{h\vph_0}^{-1}=dD_G$ and $C_lE_G^{-1}=dD_GD_G^{-1}$ belong to ${\got g}/p{\got g}\otimes_k R_0/pR_0$. We have: 
\medskip
{\bf Fact.} {\it If $X_l\in {\got g}/p{\got g}\otimes_k R_0/pR_0$, then $E_G\pi(F^0)\bigl(t_l(0)^{p-1}X_l^{[p]}\bigr)E_G^{-1}\in {\got g}/p{\got g}\otimes_k R_0/pR_0$.} 
\medskip
{\bf Proof:} As we are dealing with $W(k)$-bases of $M$ and not with $R_0^\wedge$-bases of $M\otimes_{W(k)} R_0^\wedge$, the statement of the Fact has a linear aspect. So we can assume $k=\bar k$ and that $\{e_1,...,e_{d_M}\}$ is in fact a $\ZZ_p$-basis of a $\ZZ_p$-module $M_{\ZZ_p}$ of $M$ constructing as in 2.2.9 8) but for $(M,F^1,h\vph_0,G)$. In such a case $X_l^{[p]}\in {\got g}/p{\got g}\otimes_k R_0/pR_0$ and $M_{h\vph_0}\in G(W(k))$. The Fact follows. 
\medskip
So, in our context of connections respecting the $G$-action, we can view (CMF) as $\abs{S(G)}$ equalities between elements of ${\got g}/p{\got g}\otimes_k R_0/pR_0$ or as an equality between two elements of ${\got g}/p{\got g}\otimes_k \Om_{R_0/pR_0/k}$. 
\smallskip
Similarly, when passing from things mod $p^n$ to things mod $p^{n+1}$, we come across a system of equations $SE_n$ which can be put in a similar to (MF) or to (CMF) matrix form; the only difference: $C_l$ gets replaced by some other matrix $C_l(n)$, $l\in S(G)$. Moreover, in our context of connections respecting the $G$-action, we have $C_l(n)M_{h\vph_0}^{-1},C_l(n)E_G^{-1}\in {\got g}/p{\got g}\otimes_k R_0/pR_0$. So we can view as well $SE_n$ as $\abs{S(G)}$ equalities between elements of ${\got g}/p{\got g}\otimes_k R_0/pR_0$. 
\medskip
{\bf 3.6.8.6.2. Digression.} If we choose a different Frobenius lift of $R_0^\wedge$ which mod $p^2$ takes all ideals $(t_l(0))$'s into themselves, then following the pattern of 3.6.8 5), we get a system of equations which can be put similarly in the matrix form: the only difference we have is that the matrix which is $t_l(0)^{p-1}$ times the identity matrix $I_{d_M}$ of $M_{d_M}(R_0)$ is replaced by another scalar multiple of $I_{d_M}$. More generally, if we consider an arbitrary Frobenius lift of $R_0^\wedge$, then $\forall l\in S(G)$, the matrix $t_l(0)^{p-1}X_l^{[p]}$ has to be substituted by a linear combination 
$$\sum_{l^\prime\in S(G)} b_{ll^\prime}X_{l^\prime}^{[p]},\leqno (COMB)$$
where all $b_{ll^\prime}$'s are elements of $R_0/pR_0$. 
\medskip
{\bf 3.6.8.7. The abelian situation.}
We assume now that we are in the context of 3.6.8.5 4), with $\tilde G$ an abelian group. We speak about the abelian situation.
The nice thing in this case is:
\medskip
{\bf Fact.} {\it  $\forall (i,s,l)\in S(M)\times S(M)\times S(\tilde G)$, the coefficients of $\tilde L_{isl}$ are elements of $W(k)$.} 
\medskip
{\bf Proof:} This is a consequence of the fact that $\tilde G$ is abelian and that the connections $\nabla_j^n$ we get respect the $\tilde G$-action. In other words, referring to (MF) above with $E_G$ (resp. $D_G$) being denoted by $E_{\tilde G}$ (resp. by $D_{\tilde G}$), we have: $D_{\tilde G}$ commutes with $X_l$, $\forall l\in S(\tilde G)$; so (MF) can be rewritten as 
$$D_{\tilde G}^{-1}C_l+X_lM_{h\vph_0}=M_{h\vph_0}\pi(F^0)\bigl(t_l(0)^{p-1}X_l^{[p]}\bigr),$$
$\forall l\in S(\tilde G)$. The Fact follows.
\medskip
From Fact and 3.6.8.5 2) we get: in the abelian situation the open subscheme $\Mt(0)$ of ${\rm Spec}(R_0/pR_0)$ introduced in the proof of 3.6.8.2, contains ${\rm Spec}(R_0/pR_0[{1\over {\prod_{i\in S(\tilde G)} t_i(0)}}])$. The argument can be repeated for the closed subscheme ${\rm Spec}(R_0/(p,(t_i(0))_{i\in S(t)})$, with $S(t)$ an arbitrary subset of $S(\tilde G)$, to get the following result: 
\medskip
{\bf 3.6.8.8. Proposition.} {\it Under the convenience assumption, in the abelian situation, the $W(k)$-morphisms $\ell_{j,n}$ when taken mod $p$, $n\in\NN$, are defining \'etale covers of (i.e. above) the locally closed subschemes of ${\rm Spec}(R_j/pR_j)$ defined by making some of the variables $t_i(j)$ to be $0$ while inverting all others. In particular, we can take $\Mb_j(k)$ to be the empty scheme.}
\medskip
{\bf 3.6.8.9. The constancy property.}
We call the fact that when we pass in the proof of 3.6.1.2-3 from things mod $p$ to things mod $p^n$, we get the same homogeneous forms (see 3.6.8 (3) and (4), 3.6.8.5 (11), (12), (13) and (14), and 3.6.8 12)), as the constancy property. A system of equations with coefficients in an $\FF_p$-algebra and which can be put in the same shape as of (3) is said to be of first type. The system of equations defined by (11) and (12) is said to be of additive second type. The system of equations defined by (13) and (14) is said to be of adjusted additive second type. For extra terminology see 3.6.18.4.6.
\smallskip
As (implicitly) mentioned in 3.6.8.4 1) and 3), this constancy property remains true in the context of a (truncation of a) $p$-divisible object of $\Mm\Mf_{[0,1]}(R)$, with $R$ a regular, formally smooth $W(k)$-algebra whose $p$-adic completion is equipped with a Frobenius lift $\Phi_R$. We can replace everywhere $\Mm\Mf_{[0,1]}$ by $\Mm\Mf_{[c,c+1]}$, $c\in\ZZ$. So, the natural question arises:
\medskip
{\bf Q.}  {\it What about (truncations of) $p$-divisible objects of $\Mm\Mf_{[a,b]}(R)$, with $a,b\in\ZZ$, $b\ge a+2$, which are not $p$-divisible objects of $\Mm\Mf_{[c,c+1]}(R)$ for some $c\in S(a,b)$?}
\medskip
Let ${\got C}=(M_R,(F^i(M_R))_{i\in S(a,b)},\Phi_{M_R})$, with $M_R$ a projective $R^\wedge$-module, be such a $p$-divisible object. Due to the fact that the connections have to satisfy the Griffiths transversality condition, the counting argument (of equations and variables) performed before 3.6.8.1, gives us in this context (of ${\got C}$), more equations than variables. Assuming that $\Om_{R/pR/k}$ is a free $R/pR$-module of rank $r_R\in\NN$ and that $M_R$ is as well free, to get a connection on ${\got C}/p{\got C}$ is the same thing as to solve (in $R/pR$) a (suitable) system of equations of the following form:
$$x_i+c_i=L_i(x_1^p,...,x_{m}^p),\leqno (15)$$
for $i\in S(1,m)$, and 
$$c_i=L_i(x_1^p,...,x_{m}^p),\leqno (16)$$
for $i\in S(m+1,m+s)$; based on 3.6.1.1.1 4), the argument for this is entirely the same as the one of 3.6.8 5) relying on its (DIV). Here $c_i\in R/pR$ and $L_i$ is a homogeneous linear form with coefficients in $R/pR$, $\forall i\in S(1,m+s)$, while $m,s\in\NN\cup\{0\}$. We have:
\medskip
-- $m+s$ is $r_R$ times the rank of ${\rm End}(M_R)$;
\smallskip
-- $m$ is $r_R$ times the rank of the $F^{-1}$-filtration of ${\rm End}(M_R)$ defined by $(F^i(M_R))_{i\in S(a,b)}$. 
\medskip
Our assumption of the question (pertaining to $c\in\ZZ$) implies $s>0$. In general, if $s>0$ we say that a system of equations of the same type as (15) and (16) with coefficients in an $\FF_p$-algebra, is of third type (in $m$ variables and $s$ extra equations).  
\smallskip
In our situation (of ${\got C}$) $ms>0$ and so it easy to see that there are systems  of equations of the above type which do not have any solution or which do not define an \'etale scheme over ${\rm Spec}(R/pR)$. However the constancy property remains:
\medskip
{\bf The constancy property.} {\it Let $n\in\NN$. To lift a connection on ${\got C}/p^n{\got C}$ to a connection on ${\got C}/p^{n+1}{\got C}$, is the same thing as to solve (in $R/pR$) a similar system of equations as of the one of (15) and (16) above but where $c_i$ is replaced by $c_i(n)\in R/pR$, $\forall i\in S(1,m+s)$ (this holds even if ${\got C}$ is a $p$-divisible object of $\Mm\Mf_{[c,c+1]}(R)$ for some $c\in S(a,b)$: in  such a case $s=0$).} 
\medskip
This is still a very useful fact: see 3.6.18.5.5, 3.6.18.8 and 3.6.18.8.1 b) below. Also the following two approaches can produce lots of results. 
\medskip
{\bf A1.} (15) defines (cf. 3.6.8.1.2 a)) an \'etale, affine ${\rm Spec}(R/pR)$-scheme ${\rm Spec}(R_1)$. So to have (16) satisfied as well we just need to consider the closed subscheme ${\rm Spec}(R_2)$ of ${\rm Spec}(R_1)$ of solutions of (16); we are interested only in its maximal open closed subscheme which is \'etale over ${\rm Spec}(R/pR)$. We refer to ${\rm Spec}(R_1)$ as the \'etale part of the system of equations defined by (15) and (16).
\medskip
{\bf A2.} For many $W(k)$-algebras $R$ (like the local strictly henselian ones), the above system of equations defined by (15) and (16), when viewed as being in the $m+s$ variables $x_1$,..., $x_m$, $c_{m+1}$,..., $c_{m+s}$, has a solution in $R/pR$, regardless of the values of the coefficients $c_1$,..., $c_m$ or of the forms $L_i$'s. Using this, one could get nice criteria of when 3.6.8 above can be worked out for such a (truncation of a) $p$-divisible object ${\got C}$, after suitable modifications (corresponding to the new values of $c_{m+1}$,..., $c_{m+s}$) of $\Phi_{M_R}$. We hope to come back to these ideas in a future paper: the reason we do not persuade them here is related to the fact that we have not being able to fully understand how these ``suitable modifications" could be performed. Here we just mention two things: 
\medskip
a) these ``suitable modifications" should be performed only from lifting things mod $p^n$ to things mod $p^{n+1}$, $n\in\NN$, so that the shapes of the forms $L_i$'s do not get changed;
\smallskip
b) to perform them for lifting things in the context of a), everything comes down to solving some linear partial differential systems of equations of the form
$${{\partial f_0}\over {\partial w_j}}=a_{j}f_0+b_{j},\leqno (LPDSEQ)$$ 
$\forall j\in S(1,r_R)$, with $f_0$ and all $a_j$'s and $b_j$'s in ${\rm End}(M_R/pM_R)$ and with $w_1$,..., $w_{r_R}\in R/pR$ such that $\{dw_1,...,dw_{r_R}\}$ is an $R/pR$-basis of $\Om_{R/pR/k}$; here $b_j$'s are related to elements $c_{m+1}$,..., $c_{m+s}$ of $R$ obtained as above (i.e. such that they are part of a solution of suitable systems  of equations similar to the one defined by (15) and (16)). 
\medskip
Unfortunately, even if $R/pR=k[[w_1,...,w_{r_R}]]$ such a system (LPDSEQ) does not necessarily have a solution. 
\medskip
Coming back to systems of equations of third type, we have the following Corollary of 3.6.8.1.
\medskip
{\bf 3.6.8.9.0. Corollary.} {\it The number of solutions of a system of equations of third type in $m$ equations and $s$ extra variables with coefficients in an algebraically closed field $k_1$ of characteristic $p$, is either $0$ or of the form $p^{m_1}$, where $m_1\in S(0,m)$.}  
\medskip
{\bf Proof:} We use the notations of (15) and (16) above, with the role of $R$ being replaced by the one of $k_1$. If $L_{m+1}$,..., $L_{m+s}$ are all $0$, then the Corollary is obvious (based on 3.6.8.1). If one of these $s$ linear forms is non-zero and $m>2$, then taking the $p$-th root of it, as in 3.6.8.1 we can eliminate one variable and the sum $m+s$ drops by $1$ at least; warning: $s$ itself does not necessarily drop by $1$. Using induction on $m+s$, we are reduced to the case $m=1$, which is trivial, or to the case $s=0$, which is handled by 3.6.8.1. 
\medskip
Warning: under specialization, the number of solutions of a system of equations of third type can increase.
\medskip
{\bf 3.6.8.9.1. Comment.} It is worth commenting why for some time (see [Va13]) we thought that we can work out 3.6.1.3 using (mod $p$) just \'etale, affine morphisms instead of $\NN$-pro-\'etale, affine morphisms. Referring to 3.6.1.3, once we obtained $Q_{j,1}$, we thought that the free coefficients $a_{isl}(n)$ of 3.6.8 12) we obtain inductively for $n\in\NN$, do not depend on $n$. But most common this is not so. 
\smallskip
However, in many cases, proceeding as follows, we can avoid the passage (mod $p$) to $\NN$-pro-\'etale morphisms. Warning: what follows has to be interpreted up to a passage to an \'etale ${\rm Spec}(Q_{j,n})$-scheme, for some $n\in\NN$; so ${\rm Spec}(Q_{j,m})$, with $m\ge n$, are replaced accordingly by resulting pull backs. We have to make a decision: 
\medskip
{\bf Decision.} {\it Either we work with ``wild" connections (their behavior modulo high powers of the prime $p$ can not be controlled, i.e. can not be prescribed a priori) but with a very nice (for instance, like in 3.6.1) Frobenius endomorphism of $M\otimes_{W(k)} R_j$, and so later on we hope to get (sort of) a completion property (see 3.6.15 A) and (of) a slice principle (see \S 7), or we work with very nice connections and with a very ``wild" Frobenius endomorphism of $M\otimes_{W(k)} R_j$.}
\medskip
We explain this. We have the following obvious property:
\medskip
{\bf Fact.} {\it The homogeneous linear forms $L_{isl}$ of 3.6.8 (3) depend only on the expressions of $\Phi(R_j)^0$ and $\Phi(R_j)^1$ mod $p$ and of the expression of $\Phi_{R_j}$ mod $p^2$.}
\medskip
But the coefficients $a_{isl}(n)$ depend on the expressions of these three maps modulo higher powers of $p$. But we can try to modify $\Phi(Q_{j,1})^0$ and $\Phi(Q_{j,1})^1$ with something which is $0$ mod $p$ (so we do not change the forms $L_{isl}$, and so we do not affect the \'etale property expressed in 3.6.1 through $r_j$ mod $p$ being a universal element), in order to make $a_{isl}(n)=0$, $\forall (i,s,l)\in S(M)\times S(M)\times S(G)$ and $\forall n\in\NN$, $n\ge m$ for some given $m\in\NN$: first we can try to modify $\Phi(Q_{j,1})^0$ and $\Phi(Q_{j,1})^1$ modulo $p^2$ to get all $a_{isl}(1)$'s to be $0$ and then, looking at what results, we can try to modify it mod $p^3$ by something which is $0$ mod $p^2$ to get all $a_{isl}(2)$'s to be $0$, etc. Warning: we want all these modifications not to change $\Phi(Q_{j,1})^0$ and $\Phi(Q_{j,1})^1$ modulo the ideal defining the lift (it is unique) to ${\rm Spec}(Q_{j,1})$ of the origin $a_0$ of ${\rm Spec}(R_j)$.
\smallskip
Why should this be possible? Answer: there are four reasons. First, $\forall (l,n)\in S(G)\times\NN$, we can assume the endomorphism of $M\otimes_{W(k)} Q_{j,n}/pQ_{j,n}$ defined by the square matrix whose coefficients are $a_{isl}(n)$, $i,s\in S(M)$, is an element of ${\got g}\otimes_{W(k)} Q_{j,n}/pQ_{j,n}$ (see the part of 3.6.8.6.1 referring to 3.6.8.6 (CMF)). Second, any time we have to pass from things mod $p^n$ to things mod $p^{n+1}$, we lift the connection (we got mod $p^n$) (as in 3.6.8 11)), not blindly but in such a way that the resulting coefficients $a_{isl}(n)$ are linear combinations (with fixed coefficients of $R_j/pR_j$) of $p$-power elements of $Q_{j,n}/pQ_{j,n}$; after this careful lift we start to modify $\Phi(Q_{j,n})^i$, $i=\overline{0,1}$, by trying to achieve all $a_{ijl}(n)$'s to be $0$. Third (resp. fourth), we have Raynaud's (resp. Grothendieck's) theorem of [BBM, 3.1.1] (resp. of [Il, 4.4]). 
\smallskip
So the situation (after suitable modifications) should get stabilized: we expect (at least if $k$ is of some special type, i.e. ``almost" algebraically closed in some sense which we can not specify too much: just that it is related to the defs. of 3.6.18.4.6 F below) the existence of $m\in\NN$ such that (after suitable modifications) we can take $Q_{j,n}=Q_{j,m}$ $\forall n\in\NN$, $n\ge m$. From many points of view, it is preferable to have (if possible) $a_{isl}(n)=0$, $\forall n\in\NN$, $n\ge m$; so the problem of such a stabilization is a very important one. If such a stabilization is possible, then we do not have to pass (mod $p$) to $\NN$-pro-\'etale covers: this is so due to the fact that $Q_{j,n}/pQ_{j,n}$, $n\ge m$, is a connected component of the \'etale scheme over ${\rm Spec}(Q_{j,m}/pQ_{j,m})$ defined by the system (4) of equations of 3.6.8. This gives a very nice expression (fully controllable) of the resulting connections $\nabla_j^n$ but the control on $\Phi(Q_{j,n})^0$ and $\Phi(Q_{j,n})^1$ is much less; for instance, we can not regain the (sort of) completion property or (of) the slice principle.  
\smallskip
Once we gained the completion property and the slice principle, the second approach (trying to avoid the use of $\NN$-pro-\'etale morphisms), looks, from many points of view, more convenient. In all concrete cases we analyzed (for instance, this is automatically so if we are in a context involving polarized abelian varieties, or which can be reduced to such a context: see 4.12.12 and 4.12.12.0) the situation (after suitable modifications) gets stabilized. Warning: neither Grothendieck's nor Raynaud's theorems mentioned above are enough to imply such a stabilization, even for situations where we can assume that all tensors $(t_{\al})_{\al\in\Mj}$ are related to (i.e. are defined by) endomorphisms of $M$. We will come back to this in \S 10.   
\medskip
{\bf 3.6.9. Remarks. 1)} We could have worked 3.6.1 with $U_0$ and $U_2$ as open subschemes of $G^{\rm der}$. This is the most common case needed for applications. For example, a $\sg$-crystal
$(M,tg_2\vph_0)$, with $t\in Z(G)(W(k))$ and $g_2\in G(W(k))$, is very much the same as the $\sg$-crystal $(M,g_2\vph_0)$; so we often (for instance when $k=\bar k$ and $p$ does not divide the order of the center of $G^{\rm der}$) can consider
only elements $g_2\in G^{\rm der}(W(k))$. 
\smallskip
{\bf 2)} The assumption ii) of the beginning of 3.6 is not really needed: it has been inserted just to be able to be faster in 3.6.8 8) and 12) and so to give an alternative way in 3.6.8 7) which is not based on 3.6.8 1). Obviously, the property ii) of the beginning of 3.6, is implied by the property i) of the mentioned place, provided we replace $U_0$ and $U_2$ by suitable open, affine subschemes of them. This makes life easier (cf. 3.6.9.1 and 3.6.10 below); in particular 3.6.11 below can be expressed in a simpler form, without mentioning $U$ or $b$.
\smallskip
{\bf 3)} Looking back to 3.6.1-8 we see that the only two places we used the fact that $k$ is an infinite field are:
\medskip
{\it --} 3.6.0, where we needed the existence of an element $h\in\bigl(U_0\cap U_2\bigr)(k)$, which as a point of $U_0$ (resp. of $U_2$) is different from the point obtained by taking the special fibre of $a_0$ (resp. of $g_2$);
\smallskip
{\it --} 3.6.6, where we needed $\Ml$ to be infinite.
\medskip
So if we are not bothered about gluing things (``over $W(k)$") (cf. 3.6.5), we can work with only one deformation, cf. 3.6.7. So, in this context of just one deformation, the whole of 3.6.1 makes sense and remains true even in the case when $k$ is a finite field. 
\medskip
{\bf 3.6.9.1. Potential-deformation sheets.} Let $k$ be a perfect field of characteristic $p>2$. A locally closed, affine subscheme $U$ of a reductive group $\tilde G$ over $W(k)$ (resp. an affine  scheme $U$ over $W(k)$), with $U(W(k))$ containing the identity element of $\tilde G(W(k))$ (resp. such that a $z\in U(W(k))$ is fixed), is said to be a potential-deformation reductive sheet (resp. just sheet), if there is an \'etale morphism (resp. a formally \'etale morphism) 
$$
b_U:U\to {\rm Spec}(W(k)[x_1,...,x_{d(U)}]), 
$$
for some $d(U)\in\NN$, such that the identity element of $\tilde G(W(k))$ (resp. $z$) factors through the closed subscheme of $U$ defined by $x_i=0$, $i\in S(1,d(U))$. If moreover, we have:
\medskip
-- the closed subscheme of $U$ defined by the identity element of $\tilde G(W(k))$ (resp. defined by $z$) is the closed subscheme of $U$ defined by $x_i=0$, $i=\overline{1,d}$, and
\smallskip
-- the closed subscheme of $U_k$ defined by $x_i=0$, with $i$ running through an arbitrary subset of $S(1,d(U))$, is connected,
\medskip\noindent
then we speak about a convenient potential-deformation reductive sheet (resp. just sheet).
\smallskip
Very often we denote a potential-deformation reductive sheet (resp. just sheet) by a pair $(U,b_U)$ (resp.  by a triple $(U,b_U,z)$). For any potential-deformation sheet we denote by $\Phi_U$ the Frobenius lift of $U^\wedge$ which takes $x_i$ into $x_i^p$, $\forall i\in S(1,d(U))$.
\smallskip
This definition of a (convenient) potential-deformation (reductive) sheet is suggested by 3.6.9 2) and the convenience assumption 3.6.1.1.3.
\medskip
{\bf 3.6.10. Corollary.} {\it For any $g_3\in G(W(k))$ and for every open, affine subscheme $U_3={\rm Spec}(R_3)$ of $G$ which is a potential-deformation reductive sheet, there is a uniquely determined $\NN$-pro-\'etale, affine morphism 
$$
\ell_4:U_4={\rm Spec}(R_4)\to U_3
$$ 
such that:
\medskip
-- there is a closed subscheme of $G_k$ not containing the origin of $G_k$, with the property that the fibres of $\ell_4$ above points of the special fibre of $U_3$ not belonging to it, are non-empty; 
\smallskip
-- the special fibre of $U_4$ is a geometrically connected, $AG$ $k$-scheme;
\smallskip
-- $\ell_4^{-1}(\Mz_k)={\rm Spec}(k)$ and so the origin of $G$, when viewed as a $W(k)$-valued point of $U_3$, can be lifted uniquely to a point $a_4\in U_4(W(k))=U_4^{\wedge}(W(k))$; 
\smallskip
-- there is a $p$-divisible group $\Md(U_3,\Phi_{U_3},M,F^1,g_3\vph_0,G)$ over $U_4^\wedge$, whose associated $p$-divisible object of $\Mm\Mf_{[0,1]}^\nabla(U_4)$ (the Frobenius on $U_4^\wedge$ being the one induced from a fixed Frobenius lift $\Phi_{U_3}$ of $U_3^\wedge$ as in 3.6.9.1), is isomorphic to
$$
(M_{R_4},\nabla_4)=\bigl(M\otimes_{W(k)} R^\wedge_4, F^1\otimes_{W(k)} R^\wedge_4, h_4(g_3\vph_0\otimes 1),\nabla_4\bigr),
$$ 
with $h_4\in G(R_4^\wedge)$ as the universal
element of $G$ defined by $U_4^\wedge$ and with $\nabla_4$ a connection on $M\otimes_{W(k)} R_4^\wedge$ respecting the $G$-action;
\smallskip
\medskip
-- we have the following universal property:}
\medskip
\item{{\bf UP.}} {\it For any formally \'etale, affine morphism $\tilde\ell_4:\tilde U_4={\rm Spec}(\tilde R_4)\to U_3$, such that the closed subscheme $\tilde\ell_4^{-1}(\Mz_k)$ has a non-empty intersection with all connected components of the special fibre of $\tilde U_4$ and for which we get a $p$-divisible group $\tilde \Md_4$ over $\tilde U_4^\wedge$, defining (by forgetting the connection of $\DD(\tilde\Md_4)$) a $p$-divisible object $N_{\tilde R_4}$ of $\Mm\Mf_{[0,1]}(\tilde R_4)$ isomorphic to $M_{\tilde R_4}$, there is a unique morphism $\ell:{\rm Spec}(\tilde R_4^\wedge)\to {\rm Spec}(R_4^\wedge)$ such that $\ell_4^\wedge=\tilde\ell_4^\wedge\circ\ell$ and $\tilde\Md_4$ is isomorphic to $\ell^*(\Md(U_3,\Phi_{U_3},M,F^1,g_3\vph_0,G))$ under an isomorphism $IS$ such that, by forgetting the connections, $\DD(IS)$ is the mentioned isomorphism $N_{\tilde R_4}\tilde\to M_{\tilde R_4}$.}
\medskip
Here $M_{R_3}$ is defined as in 3.6.1.1: it allows us to get $M_{R_4}$ and $M_{\tilde R_4}$ as in 3.6.1.1. 
\medskip
{\bf 3.6.11. Corollary.} {\it Let $\tilde H$ be a smooth subgroup of $G$ having a connected special fibre. Let $e_{\tilde H}:=\dim_{W(k)}(\tilde H)$. Let
$U_{\tilde H}={\rm Spec}(R_5)$ be an open, affine subscheme of $\tilde H$ for which there is an open subscheme $U$ of $G$ containing $U_{\tilde H}$ such that there are \'etale morphisms $b$
and $c$ making the diagram
$$
\def\mapright#1{\smash{
\mathop{\longrightarrow}\limits^{#1}}}
\def\mapdown#1{\Big\downarrow
\rlap{$\vcenter{\hbox{$\scriptstyle#1$}}$}}
\matrix{U_{\tilde H} &\mapright{i_{\tilde H}} &U\cr
\mapdown{b} &&\mapdown{c}\cr
{\rm Spec}\bigl(W(k)[z_1,\ldots,z_{e_{\tilde H}}]\bigr) &\mapright{a}&{\rm Spec}\bigl(W(k)[z_1,\ldots,z_d]\bigr)\cr}$$
to be cartesian (with $i_{\tilde H}$ as the natural inclusion, with $a$ as the natural closed embedding defined by: $z_i$ goes to $z_i$, for $i\in S(1,e_{\tilde H})$, and to 0 otherwise) and the origin $a_0$ factors through the closed subscheme of $U_{\tilde H}$ defined by $z_i=0$, $i=\overline{1,e_{\tilde H}}$. Let the Frobenius lift $\Phi_{R_5}$ of $R_5^\wedge$ be defined by: $z_i$ (viewed via $c$ as an element of $R_5$) goes to $z_i^p$, $i=\overline{1,e_{\tilde H}}$.
\smallskip
Then for any $\tilde g\in G(W(k))$, there is an $\NN$-pro-\'etale, affine morphism 
$$
\ell_6:U_6={\rm Spec}(R_6)\to U_{\tilde H},
$$ 
with the special fibre of $U_6$ a geometrically connected, $AG$ $k$-scheme and with $\ell_6^{-1}(\Mz_k)={\rm Spec}(k)$ (so there is a unique $a_6\in U_6^\wedge$ lifting the factorization of $a_0$ through $U_{\tilde H}$), such that:
\medskip
-- there is a $p$-divisible group over $U_6^\wedge$ whose
associated $p$-divisible object of $\Mm\Mf_{[0,1]}^\nabla(R_6)$ is isomorphic to 
$$
(M_{R_6^\wedge},\nabla_6)=\bigl(M\otimes_{W(k)} R_6^\wedge, F^1\otimes_{W(k)} R_6^\wedge,h_6(\tilde g\vph_0\otimes 1),\nabla_6\bigr),
$$ 
with $h_6$ as the universal element of $\tilde H$ defined by the $p$-adic completion of $\ell_6$; $\nabla_6$ respects the $G$-action;
\smallskip
-- we have a universal property similar to the one described at the end of 3.6.10; 
\smallskip
-- there is a closed subscheme of $\tilde\Mb_k$ not containing the origin of $\tilde H_k$ and such that the fibres of $\ell_6$ above points of the special fibre of $U_{\tilde H}$ not belonging to it, are non-empty.}
\medskip
The results 3.6.10-11 are a direct consequence of 3.6.1.3 and its proof and of 3.6.9 2) and 3). We just need to mention two things (for instance in connection to 3.6.10). First, $\ell$ is uniquely determined by requiring $\DD(\tilde\Md_4)$ to be isomorphic to $\ell^*(M_4,\nabla_4)$. Second, $\tilde R_4/p\tilde R_4$ has a finite $p$-basis, as $R_3/pR_3$ does have, and so the passage from filtered $F$-crystals to $p$-divisible groups is achieved via the fully faithfulness property of [BM, 4.1.1] and via Grothendieck--Messing theory of [Me, ch. 4-5] (cf. also 3.6.2.0); the usage of this last theory is part of the second place where we need $p>2$.
\medskip
{\bf 3.6.11.1. Notations.} We denote by $\Md(U_{\tilde H},\Phi_{R_5},M,F^1,\tilde g\vph_0,G)$ the $p$-divisible group over $U_6^\wedge$ mentioned in 3.6.11. For $n\in\NN$, we denote by $U_{6,n}={\rm Spec}(Q_{6,n}^\wedge)$ the scheme constructed as in 3.6.1.3 1), starting from the data of the sextuple 
$$(U_{\tilde H},\Phi_{R_5},M,F^1,\tilde g\vph_0,G).$$ 
We use morphisms $\ell_{6,n}:{\rm Spec}(Q_{6,n}^\wedge)\to {\rm Spec}(R_6^\wedge)$, and $\ell_6^n:{\rm Spec}(Q_{6,n+1}^\wedge)\to {\rm Spec}(Q_{6,n}^\wedge)$, having the same significance as in 3.6.1.3. We can take as $U_6$ the limit of the $\NN$-projective system $\ell_6^n$, $n\in\NN$, and as $\ell_6$ we can take the resulting $W(k)$-morphism $U_6\to U_{\tilde H}$. Let $Q_6$ be the $W(k)$-algebra such that $U_6^\wedge={\rm Spec}(Q_6)$.
\medskip
{\bf 3.6.12. Remark.} We refer to 3.6.11. We can replace ``$G$-action" by $\tilde H^1$-action, where $\tilde H^1$ is the smallest smooth, connected subgroup of $G$ containing $\tilde H$ and the image of $\mu$ and whose Lie algebra is taken by $p\tilde g\vph_0$ into itself, cf. [Fa2, rm. ii) after th. 10] and Fact of 2.2.9 1). Moreover, it is enough to assume $U_{\tilde H}$ is a locally closed subscheme of $G$ through which $a_0$ factors, having a geometrically connected special fibre and smooth over $W(k)$, without assuming that it is the open subscheme of a smooth
subgroup $\tilde H$ of $G^{\rm der}$. However, in practice $\tilde H$ is either a subgroup containing $G^{\rm der}$ or is an abelian subgroup of $G$ isomorphic to $\GG^{e_{\tilde H}}_a$; in this last situation we can take $U_{\tilde H}=\tilde H$. 
\medskip
{\bf 3.6.13. The standard deformation setting in the context of a SHS.} As explained in 1.14.2, the reading of 3.6.13-14 should proceed only after the reading of 3.9 and 3.11-12. We start moving towards geometric contexts.
\smallskip
Let $(f,L_{(p)},v)$ be a SHS defined by an injective map $f:(\tilde G,X)\hookrightarrow \bigl({\rm GSp}(W,\psi),S)$. We use the standard notation of 2.3.1-3 but with $G$ replaced everywhere by $\tilde G$. Let $k$ be a perfect field of characteristic $p$ and let 
$$
z:{\rm Spec}(W(k))\to\Mn_{W(k)}/H_0. 
$$
Let $(M,F^1,g\vph_0,G^\prime,\tilde\psi_M)$ 
be its attached principally quasi-polarized (not necessarily quasi-split) quasi Shimura filtered $\sg$-crystal (cf. 2.3.10). Here $G^\prime$ is a form (not a priori inner) of $\tilde G_{W(k)}$. For simplifying the presentation we assume $G^\prime=\tilde G_{W(k)}$; so we do not have to mention any more ``not necessarily quasi-split". Let $G$ be the subgroup of $G^\prime$ such that under the last identification it is $\tilde G^0_{W(k)}$ (see 2.3.1 for the meaning of $\tilde G^0_{\ZZ_{(p)}}$). Here $g\in G(W(k))$; strictly speaking $\vph_0$ is obtained from $g\vph_0$, as in 3.2.3 (i.e. the Shimura adjoint filtered Lie $\sg$-crystal attached to $(M,F^1,\vph_0,\tilde G_{W(k)})$ is of Borel type). We use some of the previous notations of this 3.6, adjusted as need (the adjustments are mentioned at the right time). In particular we have a direct sum decomposition $M=F^1\oplus F^0$ and $a_0$ denotes the origin of $G$; in the case of this direct sum decomposition, the wording ``adjusted as needed" refers to the fact that we are not bothered that it is associated to a cocharacter $\mu$ of $\tilde G_{W(k)}$ and not of $G$ itself.
\smallskip 
Let $U_3={\rm Spec}(R_3)$ be an open, affine subscheme of $G=\tilde G^0_{W(k)}$ which is a potential-deformation reductive sheet. Let $(M_{R_4},\nabla_4)=\bigl(M\otimes_{W(k)} R_4^\wedge,F^1\otimes_{W(k)} R_4^\wedge,h_4(g\vph_0\otimes 1),\nabla_4\bigr)$ be the $p$-divisible object of $\Mm\Mf^\nabla_{[0,1]}(R_4)$ mentioned in 3.6.10 for $g_3=g\in G(W(k))$. We use the notations of 3.6.10; strictly speaking $(M,F^1,g\vph_0,G,\tilde\psi_M)$ is not a principally quasi-polarized Shimura filtered $\sg$-crystal but it is ``close enough" (cf. 2.2.9 1) and 1'): we are in a pseudo context) to justify this usage (cf. 3.6.1.6 or 3.6.11).
\smallskip
Let $A_z:={\rm Spec}(W(k))\,_z\times_{\Mn_{W(k)}/H_0}\Ma_{H_0}$.
For similar reasons, we also assume the subgroup $H_0$ of $G(\AA_f^p)$ is ``small enough (w.r.t. $z$)" so that the family of Hodge cycles ${(w_\al)}_{\al\in\Mj^\prime}$ with which the abelian variety $A_z\times_{W(k)} {\rm Spec}(W(\bar k))$ over $W(\bar k)$ is naturally endowed (through any lift ${\rm Spec}(W(\bar k))\to\Mn_{W(k)}$ of $z$), is defined over $W(k)$ and so it gives birth to a family ${(t_\al)}_{\al\in\Mj^\prime}$ of de Rham components of the $F^0$-filtration of $\Mt(M)$ and implicitly of the $F^0$-filtration of $\Mt(M\otimes_{W(k)} R_4)$. To be more concrete, we assume $H_0=H_1$, where $H_1$ is as in 2.3.11. This assumption is not really needed (cf. 3.6.14.2 below); however, it allows us to avoid using the word ``quasi". 
\smallskip
We have: $\nabla_4(t_{\al})=0$ (as $\nabla_4$ respects the $G$-action) and $h_4(g\vph_0\otimes 1)(t_{\al})=t_{\al}$, $\forall\al\in\Mj^\prime$. The fact that $G$ fixes $\tilde\psi_M$, implies that $\tilde\psi_M$ defines a principal quasi-polarization of $M_{R_4}$. Let 
$$
\tilde M_{R_4}:=(M_{R_4},\nabla_4,\tilde\psi_M).
$$
It is a principally quasi-polarized $p$-divisible object of $\Mm\Mf_{[0,1]}^\nabla(R_4)$. Sometimes we prefer to regard it as being with tensors (they are $t_{\al}$, $\al\in\Mj^\prime$). 
\smallskip
We denote by $m_{\rm nat}$ the natural morphism ${(\Mn_{W(k)}/H_0)}^\wedge\to\Mn/H_0$. 
\smallskip
{\bf SA1.} Let $\widehat{G}={\rm Spec}(\widehat{R_3})={\rm Spec}(\widehat{R_4})$ be the completion of $G$ in its origin (or the completion of $U_4$ in $a_4$). From [Va2, 5.4.4] (see 2.3.11 or the proof of 2.3.15) we deduce the existence of a $W(k)$-morphism $z_{\widehat{G}}:\widehat{G}\to\Mn_{W(k)}/H_0$ such that the completion of $a_4$ is mapped into $z$ and the two principally quasi-polarized Shimura $p$-divisible groups over $\widehat{G}$ obtained by pulling back the ones over $\Mn_{W(k)}/H_0$ (see Fact 3 of 2.3.11) and over ${\rm Spec}(R_4^\wedge)$ (see 3.6.10) are isomorphic, through an isomorphism lifting the identity in (the completion of) $a_4$. Warning: we assume $z_{\widehat{G}}$ mod $p$ factors through the henselization of the localization of ${\rm Spec}(R_4/pR_4)$ in the $k$-valued point defined by $a_4$.
\medskip
{\bf 3.6.14. Theorem (the deformation principle: the expanded form).} {\it There is an integral, affine scheme $V_G={\rm Spec}(S_G)$, which is the $p$-adic completion of an $\NN$-pro-\'etale scheme over a smooth $W(k)$-scheme, and there
are $W(k)$-morphisms $a:{\rm Spec}(W(k))\hookrightarrow V_G$, $b:V_G\to {(\Mn_{W(k)}/H_0)}^\wedge$ and $c:V_G\to U^\wedge_4$ such
that:
\medskip
\item{a)} $a_4=c\circ a$ and the $p$-adic completion of $z$ is $b\circ a$;
\smallskip
\item{b)} $b$ is obtained by taking the $p$-adic completion of the composite of an $\NN$-pro-\'etale morphism with a smooth morphism while $c$ is the $p$-adic completion of an $\NN$-pro-\'etale morphism;
\smallskip
\item{c)} The special fibre of $V_G$ is geometrically connected and the image of $c$ mod $p$ contains an open, dense subscheme of the special fibre of $U_4^\wedge$;
\smallskip
\item{d)} If $N_{S_G}=(N,F^1_{S_G},\Phi_N,\nabla_N,\tilde\psi_N)$ is the principally quasi-polarized $p$-divisible object of $\Mm\Mf^\nabla_{[0,1]}(S_G)$ defined by the principally quasi-polarized $p$-divisible group of the principally polarized abelian scheme $V_G\,_{b_0}\times_{\Mn_{W(k)}/H_0} (\Ma_{H_0},\Mp_{\Ma_{H_0}})$ over $V_G$, then there is an isomorphism of $N_{S_G}$ into $\tilde M_{S_G}$ (the pull back of $\tilde M_{R_4}$ through the $W(k)$-morphism $c$) taking the de Rham component $t_\al^{S_G}$ (of the Hodge cycle $w_\al^{S_G}$ of $V_G\,_{b_0}\times_{\Mn/H_0} \Ma_{H_0}$; it is viewed as a tensor of $N[{1\over p}]$ and is obtained similarly to 2.3.12-13) into $t_\al$, $\forall\al\in\Mj^\prime$ (here $b_0:=m_{\rm nat}\circ b$).}
\medskip
Let $\tilde H$ be the unipotent radical of the parabolic subgroup of $\tilde G^0_{W(k)}$ normalizing the direct summand $F^0$ of $M$. Let $U_{\tilde H}={\rm Spec}(R_5)$ be an open, affine subscheme of $\tilde H$ such that the conditions of 3.6.11 are satisfied. We  use the notations of 3.6.11 and 3.6.11.1, with $\tilde g=g$. So $e_{\tilde H}=\dim_{\CC}(X)$. Similarly to $\tilde M_{R_4}$ above, we define a principally quasi-polarized $p$-divisible object (sometimes viewed of being with tensors) of $\Mm\Mf_{[0,1]}^\nabla(R_6)$
$$
\tilde M_{R_6}:=(M_{R_6},\nabla_6,\tilde\psi_M).
$$ 
\indent
{\bf SA2.} As in SA1 we now assume that the similarly defined $W(k)$-morphism from the completion of $U_{\tilde H}$ in the origin into $\Mn_{W(k)}/H_0$, when taken mod $p$, factors through the  henselization of the localization of ${\rm Spec}(Q_6/pQ_6)$ in its $k$-valued point defined by $a_6$. 
\medskip
We refer to this assumption (as well as of the one of SA1) as the starting assumption.
\medskip
{\bf 3.6.14.1. Theorem (the deformation principle).} {\it There is an integral, affine scheme $V_{\tilde H}={\rm Spec}(S_{\tilde H})$, which is the $p$-adic completion of an $\NN$-pro-\'etale scheme over a smooth $W(k)$-scheme, and there
are $W(k)$-morphisms $a:{\rm Spec}(W(k))\hookrightarrow V_{\tilde H}$, $b:V_{\tilde H}\to {(\Mn_{W(k)}/H_0)}^\wedge$ and $c:V_{\tilde H}\to U^\wedge_6$ such
that:
\medskip
\item{a)} $a_6=c\circ a$ and the $p$-adic completion of $z$ is $b\circ a$;
\smallskip
\item{b)} $b$ and $c$ are $p$-adic completions of $\NN$-pro-\'etale morphisms;
\smallskip
\item{c)} The special fibre of $V_{\tilde H}$ is geometrically connected and the image of $c$ mod $p$ contains an open, dense subscheme of the special fibre of $U^\wedge_6$;
\smallskip
\item{d)} If $N_{S_{\tilde H}}=(N,F^1_{S_{\tilde  H}},\Phi_N,\nabla_N,\tilde\psi_N)$ is the principally quasi-polarized $p$-divisible object of
$\Mm\Mf^\nabla_{[0,1]}(S_{\tilde H})$ defined by the principally quasi-polarized $p$-divisible group of the principally polarized abelian
scheme $V_{\tilde H}\,_{b_0}\times_{\Mn_{W(k)}/H_0} (\Ma_{H_0},\Mp_{\Ma_{H_0}})$ over $V_{\tilde H}$, then there is an isomorphism of $N_{S_{\tilde H}}$ into $\tilde M_{S_{\tilde H}}$ (the pull back of $\tilde M_{R_6^\wedge}$ through the $W(k)$-morphism $c$) taking the de Rham component $t_\al^{S_{\tilde H}}$ (of the Hodge cycle $w_\al^{S_{\tilde H}}$
of $V_{\tilde H}\,_{b_0}\times_{\Mn/H_0} \Ma_{H_0}$; it is viewed as a tensor of $N[{1\over p}]$ and is obtained as in 2.3.12-13) into $t_\al$, $\forall\al\in\Mj^\prime$ (here $b_0:=m_{\rm nat}\circ b$).}
\medskip
{\bf Proof:} The proof of 3.6.14 is entirely similar to the proof of 3.6.14.1, just the presentation in the case of 3.6.14.1 is slightly more complicated. So we prove here only 3.6.14.1. We do not hesitate (though this makes the proof longer) to take the necessary time to explain in detail what the systems of equations of first type as defined in 3.6.8.9 (can) bring new (besides 3.6.1.3) to the deformation theory of Shimura $p$-divisible groups. Also, as mentioned in 1.14.2, we do not hesitate to refer to parts of 3.9 and 3.11-12. We can assume $\dim_{\CC}(X)>0$: if $\tilde G$ is a torus the above Theorems are trivial.
\smallskip
{\bf A.} Let $Y_1={\rm Spec}(S^1_{\tilde H})$ be a smooth, affine $W(k)$-scheme, with a geometrically connected special fibre, and such that there is a point
$a_1\in Y_1(W(k))=Y_1^\wedge(W(k))$ and there are $W(k)$-morphisms $b_1:Y_1\to\Mn_{W(k)}/H_0$ and $c_1:Y_1^\wedge\to U_{\tilde H}$ satisfying the properties:
\medskip
\item{A)} $a_6=c_1\circ a_1$ and $z=b_1\circ a_1$;
\smallskip
\item{B)} $c_1$ mod $p$ and $b_1$ are both \'etale morphisms;
\smallskip
\item{C)} Denoting $(A_{Y_1},p_{A_{Y_1}}):=b_1^*(\Ma_{H_0},\Mp_{\Ma_{H_0}})$, $A_{Y_1}$ has the family of Hodge cycles
$(w_\al^{S^1_{\tilde H}})_{\al\in\Mj^\prime}$ (with which it is naturally endowed) defined over $Y_1$ (and not only over an \'etale cover of $Y_1$) 
and there is an isomorphism 
$$d_1:N:=H^1_{dR}(A_{Y_1}/S^1_{\tilde H})\otimes_{S^1_{\tilde H}} {S^1_{\tilde H}}^\wedge\tilde\to M\otimes_{W(k)} {S^1_{\tilde H}}^\wedge$$ 
taking the Hodge filtration of $N$ defined by $A_{Y_1}$ into $F^1\otimes_{W(k)} {S^1_{\tilde H}}^\wedge$, $p_{A_{Y_1}}$ into $\tilde\psi_M$, and the de Rham
component $t_\al^{S^1_{\tilde H}}$ of $w_\al^{S^1_{\tilde H}}$ into $t_\al$, $\forall\al\in\Mj^\prime$;
\smallskip
\item{C')} There is a Frobenius lift of ${S^1_{\tilde H}}^\wedge$ taking $z_i$ into $z_i^p$, $i\in S(1,e_{\tilde H})$, where the elements of the set $\{z_1,...,z_{e_{\tilde H}}\}$ are forming a system of regular parameters of $a_1$ in $Y_1$, such that the resulting Frobenius endomorphism of $M\otimes_{W(k)} S^{1\wedge}_{\tilde H}$ (cf. the identification $d_1$ of $C$)) is $g_{Y_1}(g\vph_0\otimes 1)$, where $g_{Y_1}\in\tilde G^0_{W(k)}(S^{1\wedge}_{\tilde H})$ mod $p$ is defined by $c_1$;
\smallskip
\item{D)} Denoting by $I_z$ the ideal of ${S^1_{\tilde H}}^\wedge$ defining the point $a_1$, the principally quasi-polarized $p$-divisible group we get over ${\rm Spec}({S^1_{\tilde H}}^\wedge/I_z^2)$ together with its crystalline tensors is $1_{\Mj^\prime}$-isomorphic (in the same sense as of 2.2.9 6)) to the principally quasi-polarized $p$-divisible group over ${\rm Spec}(R_5/I_0^2)={\rm Spec}(Q_6/I_0^2Q_6)$ together with its tensors, where $I_0$ is the ideal of $U_{\tilde H}$ defining $a_6$. (Here $c_1$ allows us to identify ${\rm Spec}({S^1_{\tilde H}}^\wedge/I_z^2)={\rm Spec}(R_5/I_0^2)$.)
\medskip
The existence of such a scheme $Y_1$ and of $W(k)$-morphisms $a_1$, $b_1$ and $c_1$ with the above properties is just 2.3.15.1.
\smallskip
{\bf B.} Over 
$$
Y_{16}:=Y_1\times_{W(k)} {\rm Spec}(Q_6)
$$ 
we have two principally quasi-polarized $p$-divisible groups $(\Md_1,\Mp_{\Md_1})$ and $(\Md_2,\Mp_{\Md_2})$, obtained by pulling back through the natural projections from this product into its two factors, the principally quasi-polarized $p$-divisible group of $(A_{Y_1},p_{A_{Y_1}})$ and respectively the principally quasi-polarized $p$-divisible group 
$$
\bigl(\Md(U_{\tilde H},\Phi_{R_5},M,F^1,g\vph_0,G),\tilde\psi_M\bigr)
$$ 
over ${\rm Spec}(Q_6)$ (cf. 3.6.11.1); they can be viewed naturally as principally quasi-polarized Shimura $p$-divisible groups over $Y_{16}$. For any $n\in\NN$, let $Y_{16,n}$ be the moduli $Y_{16}$-scheme parameterizing isomorphisms between the two principally quasi-polarized, finite, flat, commutative group schemes over $Y_{16}$ defined by the kernels of the multiplication with $p^n$ on them. It is an affine scheme of finite type over $Y_{16}$. We get an $\NN$-projective system defined by the natural (``restriction") $W(k)$-morphisms 
$$
t^n:Y_{16,n+1}\to Y_{16,n}.
$$ 
Let $Y_{16,\infty}$ be its projective limit. 
Property D) implies: there is a natural $W(k)$-morphism 
$$
h_Z:Z:={\rm Spec}(W(k)[x_1,...,x_{e_{\tilde H}}]/(x_1,...,x_{e_{\tilde H}})^2)\to Y_{16,\infty}
$$ 
inducing $W(k)$-isomorphisms $Z\tilde\to {\rm Spec}({S^1_{\tilde H}}^\wedge/I_z^2)$ and $Z\tilde\to {\rm Spec}(Q_6/I_0^2Q_6)$; we can assume that under the first such $W(k)$-isomorphism $z_i$ mod $I_z^2$ is mapped into $x_i$ mod $(x_1,...,x_{e_{\tilde H}})^2$, $i=\overline{1,e_{\tilde H}}$. Here $x_i$'s are independent variables.
\smallskip
{\bf C.} Let $\tilde H^0$ be the completion of $\tilde H$ (or of $U_{\tilde H}$) in its origin. From [Fa2, th. 10 and the remarks after] (cf. D) and the proof of 2.3.15) we deduce the existence of a unique $U_{\tilde H}$-morphism 
$$
m_{\tilde H}:\tilde H^0={\rm Spec}(S^0_{\tilde H})\to Y_{16,\infty};
$$ 
$Y_{16,\infty}$ is naturally a ${\rm Spec}(Q_6)$-scheme and so a $U_{\tilde H}$-scheme. We have a natural factorization $f_Z$ of $h_Z$ through $m_{\tilde H}$. The key point is:
\medskip
\item{E)} {\it We can assume the morphism $m_{\tilde H}$ mod $p$ factors through the $Y_{16}$-scheme $\tilde Y_1\times_{Y_{16}} Y_{16,\infty}$, where the $Y_{16}$-scheme $\tilde Y_1$ is defined by $\tilde\ell_6$ of the following fibre product:
$$
\def\mapright#1{\smash{
\mathop{\longrightarrow}\limits^{#1}}}
\def\mapdown#1{\Big\downarrow
\rlap{$\vcenter{\hbox{$\scriptstyle#1$}}$}}
\matrix{\tilde Y_1 &\mapright{\tilde c_1} &{\rm Spec}(Q_6)\cr
\mapdown{\tilde\ell_6} &&\mapdown{\ell_6}\cr
Y_1^\wedge &\mapright{c_1}&U_{\tilde H}.\cr}$$}
\medskip
{\bf D.} E) is a consequence (in fact a restatement) of SA2. We stop to digress why the statement of the starting assumption is not automatically satisfied. We first remark that the composite of $f_Z$ with the morphism $\tilde H^0\to\tilde Y_1$ is a canonical closed embedding $a_Z:Z\hookrightarrow\tilde Y_1$ (cf. D) and B). The completion of $\tilde Y_1$ w.r.t. the ideal defining this embedding, can be identified with $\tilde H^0$ (cf. B)). Let $\tilde c_1^0$ and $\tilde\ell_6^0$ be the natural $W(k)$-morphisms (obtained through this identification) from $\tilde H^0$ to ${\rm Spec}(Q_6)$ and respectively to $Y_1^\wedge$. Let $m\in\NN$, $m\ge 2$. Let $\tilde I_0:=I_0S^0_{\tilde H}$. Referring to [Fa2, th. 10], when we have to pass from things modulo $\tilde I_0^m$ to things modulo $\tilde I_0^{m+1}$, we do not know a priori that the $k$-morphism $\tilde c_1^0$ mod $p$ is not changed. 
To argue this, let $P$ be the parabolic subgroup of $G$ normalizing $F^1$; so ${\rm Lie}(P)=F^0({\rm Lie}(G))$. The obstruction of lifting things modulo $\tilde I_0^m$ to things modulo $\tilde I_0^{m+1}$, is given by some element 
$$
\ga_m\in p{\rm Lie}(G)\otimes_{W(k)} \tilde I_0^m/\tilde I_0^{m+1}
$$ 
(cf. C') and D) for $m=2$). 
We can modify (as in loc. cit.) $\tilde c_1^0$ by something which is $0$ mod $p$ so that we have $\ga_m\in F^0({\rm Lie}(G))\otimes_{W(k)} \tilde I_0^m/\tilde I_0^{m+1}$. But this element corresponds (as in loc. cit.) to a replacement of a Frobenius endomorphism $\Phi_{\tilde H}^0$ of $M\otimes_{W(k)} {S_{\tilde H}^0}$ (defined for instance via $\Md_1$) by another one $(1-\ga_m)\Phi_{\tilde H}^0 (1-\ga_m)^{-1}$ isomorphic to it. The new obstruction we get is given by some new element
$$
\ga_{m+1}\in {\rm Lie}(G)\otimes_{W(k)} \tilde I_0^{m+1}/\tilde I_0^{m+2}.
$$
But $\gamma_{m+1}$ is not a priori divisible by $p$ (however: if $\gamma_m$ is divisible by $p^2$, then $\gamma_{m+1}$ is divisible by $p$).
 So the induction does not apply. That is why we need SA2. 
\smallskip
{\bf E.} Let now $n\in\NN$. The principally quasi-polarized, finite, flat, commutative group scheme we naturally have over ${\rm Spec}(Q_{6,n+1}^\wedge)$ (see 3.6.2) is defined over the spectrum of a $Q_{6,n+1}$-subalgebra $\tilde Q_{6,n+1}$ of $Q_{6,n+1}^\wedge$ of finite type. So $Y_{16,n+1}$ is obtained from a scheme $Y_{16,n+1,{\rm al}}$ of finite type over $Y_1\times_{W(k)} {\rm Spec}(\tilde Q_{6,n+1})$ by pull back through the natural $W(k)$-morphism ${\rm Spec}(Q_6)\to {\rm Spec}(\tilde Q_{6,n+1})$. From [BLR, th. 12 of p. 83] we deduce that the resulting $W(k)$-morphism $\tilde H^0\to Y_{16,n+1,{\rm al}}$ factors through a smooth, affine ${\rm Spec}(Q_{6,n+1})$-scheme $W_{16,n+1}$ (the fact that $W(k)$ is an excellent DVR and $U_3$ is a potential deformation reductive sheet implies that all hypotheses of loc. cit. are satisfied). 
\smallskip
There is an affine $Y_{16,n}$-scheme $Z_{16,n}$ which (via the natural projection ${\rm Spec}(Y_{16,n})\to {\rm Spec}(Q_6)$) is smooth over ${\rm Spec}(Q_6)$: for instance we can take $Z_{16,n}$ to be the pull back to ${\rm Spec}(Q_6)$ of the ${\rm Spec}(Q_{6,n+1}^\wedge)$-scheme
$$Z^1_{16,n}:={\rm Spec}(Q_{6,n+1}^\wedge)\times_{{\rm Spec}(Q_{6,n+1})} W_{16,n+1};$$ 
as $Z_{16,n}^1$ mod $p^{n+1}$ is the same as $W_{16,n+1}$ mod $p^{n+1}$, we have a natural $W(k)$-morphism $Z_{16,n}\to Y_{16,n}$ (cf. 2.2.1.1 2) and Fact 2 of 2.2.1.0). 
\smallskip
From this and D), using a slice, we get an affine $W(k)$-morphism 
$$
t(n):{\rm Spec}(T_n)\to Y_{16,n}
$$ 
such that, when composed with the projections of $Y_{16,n}$ on $Y_1$ and respectively on ${\rm Spec}(Q_6)$, it gives birth to a formally \'etale, affine morphism $q(n)$ and respectively to an \'etale, affine morphism $s(n)$ (the fact that $q(n)$ mod $p$ is formally \'etale can be also checked using arguments on Kodaira--Spencer maps). Moreover we can assume that the induced (by $h_Z$) $W(k)$-morphism $Z\to Y_{16,n}$ factors through $t(n)$, producing an affine $W(k)$-morphism 
$$
h_n:Z\to {\rm Spec}(T_n), 
$$
and that (cf. E)) $t(n)$ mod $p$ factors through $\tilde Y_1\times_{Y_{16}} Y_{16,n}$; so the special fibre of the $W(k)$-morphism $\tilde H^0\to Y_{16,n}$ defined by $m_{\tilde H}$ factors through the special fibre of $t(n)$. Also we can assume ${\rm Spec}(T_n/pT_n)$ is a geometrically connected $k$-scheme of finite type.
\smallskip
By induction on $n\in\NN$ we check that we can choose the things in such a way that  
there is an \'etale morphism 
$$
t_n:{\rm Spec}(T_{n+1})\to {\rm Spec}(T_n)
$$
for which we have:  
\medskip
{\item{F)} {\it $s(n+1)=s(n)\circ t_n$, $t_{n+1}\circ h_{n+1}=h_n$ and $t^n\circ t(n+1)\equiv t(n)\circ t_n$ mod $p^{n+1}$.}
\medskip
{\bf F. Proof of F): the ordinary case.}
The proof of F) is related to the following statement: for suitable $m\in\NN\cup\{0\}$, we can not have two formally \'etale morphisms from ${\rm Spec}(T_{n+m})$ into $Y^1$ which are different mod $p^n$ but identical mod $p$, whose images mod $p$ contain an open subset $U^1$ of the special fibre of $Y^1$ and are such that the pull backs of $A_{Y_1}[p^{n+m}]$ and of its principal quasi-polarization through them are the same (i.e. are isomorphic). We refer to this statement as the statement $(n,m)$. 
\smallskip 
The idea is: it is enough to check a statement $(n,m)$ at $\tilde G_{W(k)}$-ordinary points of $U^1$ (due to the fact that ${\rm Spec}(T_n/pT_n)$ is connected), by using (good) $\tilde G_{W(k)}$-multiplicative coordinates; here $\tilde G_{W(k)}$-ordinary points of $U^1$ are in the sense of 3.9.3.1 below. As the theory of such coordinates is developed just in 4.1-4 and 4.7 below, we have to work around. Such a theory of ordinary points and canonical coordinates is already available (cf. [De3] and [No1]) for the particular case in which the Shimura-ordinary type associated to $\Mn$ (see 4.1) is an ordinary type. So we first prove --based on [Ka3, ch. 2]-- statement $(n,n)$ and F) for the particular case in which we assume the existence of a dense set of ordinary points of $\Mn_{k(v)}$. 
\smallskip
We can assume $s(n+1)=s(n)\circ t_n$ and $t_{n+1}\circ h_{n+1}=h_n$. For the congruence mod $p^n$ of $t^n\circ t(n+1)$ and $t(n)\circ t_n$, it is enough to check that $q(n+1)$ is congruent to $q(n)\circ t_n$ mod $p^n$. So we start checking that $q(n)\circ t_{n+1}$ is congruent to $q(n+1)$ modulo $p^{m(n)}$ for a sufficiently large $m(n)\in\NN$. From [Ka3, p. 151] we deduce that, at the ring levels, $q(n)\circ t_{n}$ is congruent to $q(n+1)$ modulo the $[{{n+1}\over 2}]$-power of any maximal ideal of $T_n$ corresponding to an ordinary point ${\rm Spec}(T_n)(\bar k)$ (with the notations of loc. cit., if $x\in A(k)[p^n]$, then $p^{2n}\tilde x=0$). But the ordinary points of ${\rm Spec}(T_n/pT_n)$ are dense (cf. our density assumption and the fact  that ordinary points can be recognized working mod $p$). So we can take $m(n)=[{{n+1}\over 2}]$. This takes care of the statement $(n,n)$. After a relabeling (like $T_{2^n}$ becomes $T_n$), we get F). This ends the proof of F) in the mentioned particular (ordinary) case. 
\medskip
{\bf G.} We come back to the general case and assume (for the time being) that property F) holds. Let ${\rm Spec}(T_{\infty})$ be the limit of the $\NN$-projective system defined by $t_n$, $n\in\NN$. From F), we deduce the existence of a $W(k)$-morphism 
$$
{\rm Spec}(T_{\infty}^\wedge)\to Y_{16,\infty}
$$
such that $h_Z$ factors through it. Moreover, ${\rm Spec}(T_{\infty}^\wedge)$ together with the $W(k)$-morphism ${\rm Spec}(W(k))\to {\rm Spec}(T_{\infty}^\wedge)$ defined via $h_Z$ and the logical $W(k)$-morphisms ${\rm Spec}(T_{\infty}^\wedge)\to U_6^\wedge$ and ${\rm Spec}(T_{\infty}^\wedge)\to (\Mn_{W(k)}/H_0)^\wedge$, satisfies all requirements of 3.6.14.1, except of c). So to get 3.6.14.1 c), we have to work in a more refined way: we have to show that ${\rm Spec}(T_{\infty}/pT_{\infty})$ can be obtained from ${\rm Spec}(T_1/pT_1)$ using systems of equations ``closed in spirit" to those of first type, in the same manner as in the proof of the last part of 3.6.8.1.2 a); on the way of achieving this we also prove F) in the general case. The needed refinement is present in H to K below.
\medskip
{\bf H. Proof of F): the general case.} 
We work in the general case.
We reached the situation: over ${\rm Spec}(T_1)$ there are two principally quasi-polarized Shimura $p$-divisible groups such that their kernels under the multiplication by $p$ are isomorphic; as ${\rm Spec}(T_1)$ is a regular, formally smooth $W(k)$-scheme it makes sense to speak about an isomorphism preserving the tensors: by this we mean a filtered isomorphism between filtered, projective $T_1^\wedge$-modules underlying the filtered $F$-crystals of these two $p$-divisible groups, preserving the quasi-polarizations and (in the sense of $1_{\Mj^\prime}$-isomorphisms) the tensors, and which mod $p$ results in an isomorphism $\Mi_1$ between truncations mod $p$ of principally quasi-polarized $p$-divisible objects of $\Mm\Mf_{[0,1]}^\nabla(T_1)$.
We have:
\medskip
\item{G)} {\it There is a formally \'etale, $p$-adically complete, affine ${\rm Spec}(T_1)$-scheme $\tilde Z_{\infty}$ through which the $W(k)$-morphism $h_1$ factors, such that over the special fibre $\tilde Y_{\infty}$ of $\tilde Z_{\infty}$, the two principally quasi-polarized $p$-divisible groups we got naturally over ${\rm Spec}(T_1/pT_1)$, are isomorphic, under an isomorphism which, when viewed (by forgetting the filtrations) as an isomorphism between $p$-divisible objects of $\Mm_{[0,1]}^\nabla(\tilde Z_{\infty})$, lifts the extension of $\Mi_1$ to $\tilde Z_{\infty}$, preserves the tensors and in the $Z$-valued point defined by the factorization of $h_1$, is the one defined naturally via $h_Z$.}
\medskip
This is a consequence of F). As we proved property F) only in the ordinary case, we need to point out that G) is also a consequence of SA2. To see this we just need to remark that SA2 implies directly:
\medskip
{\bf Fact.} {\it Let $n\in\NN$. The integral closed subscheme of the ${\rm Spec}(T_1/pT_1)$-scheme parameterizing isomorphisms between the kernels of the multiplication by $p^n$ of the two principally quasi-polarized $p$-divisible groups we got over ${\rm Spec}(T_1/pT_1)$, through which the composite morphism of $\tilde m_H$ mod $p$ and of the natural $k$-morphism ${Y_{16,\infty}}_k\to {Y_{16,1}}_k$ factors naturally, is generically finite over ${\rm Spec}(T_1/pT_1)$ and so \'etale over ${\rm Spec}(T_1/pT_1)$ in its $k$-valued point obtained from the one of $\tilde H^0$, via this factorization.} 
\medskip
It is worth pointing out that the generically finiteness part of this Fact results also from the following obvious Subfact (which is a particular case of b) of 2.2.4 B): 
\medskip
{\bf Subfact.} {\it Let $k_1\subset k_2$ be two algebraically closed fields containing $k$. We consider isomorphisms over $W(k_2)$ between truncations mod $p^n$ of two Shimura $p$-divisible groups $\tilde D(1)=(D(1),(t(1)_{\al})_{\al\in\Mj})$ and $\tilde D(2)=(D(2),(t(2)_{\al})_{\al\in\Mj})$ over $W(k_1)$, preserving the tensors, in the sense that they lift to $1_{\Mj^\prime}$-isomorphisms (as in 2.2.9 6)) between $\tilde D(1)$ and $\tilde D(2)$. Then the number of such isomorphisms is finite and all of them are defined over $W(k_1)$.}
\medskip
We come back to the proof of F). We now show that G) implies F). From the Fact we deduce that we can assume $\tilde Z_{\infty}$ is the $p$-adic completion of an affine, $\NN$-pro-\'etale ${\rm Spec}(T_1)$-scheme. We write $\tilde Z_\infty={\rm Spec}(\tilde T_\infty)$. From G) we deduce the existence of a natural morphism $m^{\rm fil}_k:{\tilde Z_{\infty k}}\to Y_{16,\infty}$. We consider the underlying $\tilde T_{\infty}$-module $\tilde M_{\infty}$ of the $F$-crystal over ${\tilde Z_{\infty k}}$ associated to the $p$-divisible group over ${\tilde Z_{\infty k}}$ (defined by the pull back of any of the two $p$-divisible groups over ${\rm Spec}(T_1/pT_1)$). The two $F^1$-filtrations $F^1_1(\tilde M_{\infty})$ and $F^1_2(\tilde M_\infty)$ of $\tilde M_{\infty}$ defined respectively via $\Md_1$ and $\Md_2$, might not be the same. We consider the principally polarized abelian scheme $PPA$ over $\tilde Z_{\infty}$ we get from $(A_{Y_1},p_{A_{Y_1}})$ via natural pull backs. The ``replacement" of $F^1_1(\tilde M_\infty)$ by $F^1_2(\tilde M_\infty)$, corresponds to a ``replacement" of $PPA$ by another principally polarized abelian scheme $PPA^\prime$ over $\tilde Z_{\infty}$, which still lifts the special fibre of $PPA$, cf. Serre--Tate's deformation theory (this is part of the second place where we need $p>2$). The pull backs (via $h_1$) of $PPA$ and $PPA^\prime$ to ${\rm Spec}(W(k))$ are the same. 
\medskip
{\bf Claim.} {\it Passing to an \'etale cover of ${\tilde Z_{\infty}}$ through which $h_1$ still factors, $PPA^\prime$ becomes isomorphic (under an isomorphism which in $h_1$ is the identity) to the pull back of $(\Ma_{H_0},\Mp_{\Ma_{H_0}})$ through a uniquely determined $O_{(v)}$-morphism
$m^{\rm fil}_1:\tilde Z_{\infty}\to\Mn/H_0$ lifting the $k(v)$-morphism ${\tilde Z_{\infty k}}\to \Mn_{k(v)}/H_0$ naturally defined by $m^{\rm fil}_k$.} 
\medskip
The existence of $m^{\rm fil}_1$, with $\tilde Z_{\infty}$ being replaced by its completion in its $k$-valued point naturally defined by the factorization of $h_1$, is a consequence of 2.2.21 UP (cf. Fact 4 of 2.3.11). By passing to an \'etale cover of $\tilde Z_{\infty}$ defined by replacing $H_0$ by a compact, open, normal subgroup of it such that $(\Ma_{H_0},\Mp_{\Ma_{H_0}})$ is the natural pull back of the universal principally polarized abelian scheme over $\Mm/K^p(N)$ (for notations cf. 2.3.3), the existence of $m^{\rm fil}_1$ becomes a statement in the faithfully flat topology and so it is a consequence of variants of the deformation theory of [Va2, 5.4.4-6] mentioned in 2.3.11. By combining the last two sentences, the Claim follows.   
\medskip
We conclude: $m^{\rm fil}_1$ defines naturally a $W(k)$-morphism 
$$m^{\rm fil}:\tilde Z_{\infty}\to Y_{16,\infty}$$ 
through which $h_Z$ factors and whose special fibre is $m^{\rm fil}_k$. As $\tilde Z_{\infty}$ is the $p$-adic completion of an affine, $\NN$-pro-\'etale ${\rm Spec}(T_1)$-scheme, we have a version of $m^{\rm fil}$ mod $p^{n+1}$ in which $Y_{16,\infty}$ is replaced by $Y_{16,n}$ and $\tilde Z_{\infty}$ is replaced by an \'etale ${\rm Spec}(T_1)$-scheme (again cf. 2.2.1.1 2) and Fact 2 of 2.2.1.0). This ends the proof of F) in the general case.
\smallskip  
{\bf I.} To explain in detail this property G), we need some notations. 
Let $\Phi_2$ be the Frobenius endomorphism of $M\otimes_{W(k)} T_1^\wedge$ induced from the Frobenius endomorphism of $M\otimes_{W(k)} Q_6$ and let $\Phi_1$ be the Frobenius endomorphism of $M\otimes_{W(k)} T_1^\wedge$ defined via C) and the pull back of the natural $p$-divisible group on $Y_1$ through $q(1)$. For computing $\Phi_1$ and $\Phi_2$, we choose (for instance) the Frobenius lift $\Phi_{T_1}$ of $T_1^\wedge$ to be the one induced from the one of $Q_6$. Let ${\rm Spec}(T_1^{\rm h})$ be the henselization of the localization of ${\rm Spec}(T_1)$ w.r.t. the $k$-valued point defined naturally by $h_1$. 
\smallskip
We can assume the isomorphism $\Mi_1$ (see H) is defined by $1_{M\otimes_{W(k)} T_1/pT_1}$. Property G) is equivalent to: there is $\tilde h_{\infty}\in G(T_1^{\rm h\wedge})$ such that mod $p$ is the identity and over $T_1^{\rm h\wedge}$ we have: 
$$\tilde h_{\infty}\circ (\Phi_1\otimes 1)=(\Phi_2\otimes 1)\circ\tilde h_{\infty}.\leqno(17)$$ 
Working as in 3.6.1.3, modulo powers of $p$, starting from equation (17) we construct an $\NN$-pro-\'etale scheme over ${\rm Spec}(T_1/pT_1)$ as follows.
We start working mod $p^2$. 
So we are trying to find an element 
$\tilde h_1\in G(T_1^\wedge)$ which mod $p$ is the identity and for which we have:
$$\Phi_2\circ\tilde h_1(m/p)\equiv\tilde h_1\circ\Phi_1(m/p)\,\,\,  {\rm mod}\,\,\, p^2,\leqno(18)$$ 
for $m\in F^1$, and
$$\Phi_2\circ\tilde h_1(m)\equiv\tilde h_1\circ\Phi_1(m)\,\,\,   {\rm mod}\,\,\, p^2,\leqno(19)$$ 
if $m\in F^0$. The above congruencies are between elements of $M\otimes_{W(k)} T_1^\wedge$.
\smallskip
For solving the system of equations defined by (18) and (19), only the expression of $\tilde h_1$ mod $p^3$ is important. Moreover as we already know the existence of $\Mi_1$, we can assume $\tilde h_1$ mod $p^2$ defines an element of $P(T_1/pT_1)$. So let $\ga_{\tilde h_1}\in {\rm Lie}(P)\otimes_{W(k)} T_1/pT_1$ be such that $\tilde h_1$ mod $p^2$ is $1_{M\otimes_{W(k)} T_1/p^2T_1}+p\ga_{\tilde h_1}$ (warning: here $p\ga_{\tilde h_1}\in p{\rm Lie}(P)\otimes_{W(k)} T_1/pT_1={\rm Lie}(P)\otimes_{W(k)} pT_1/p^2T_1$). We choose arbitrarily a $W(k)$-basis $\{\tilde e_1,...,\tilde e_{\dim_{W(k)}(P)}\}$ of ${\rm Lie}(P)$ and a $W(k)$-basis $\{n_1,...,n_{e_{\tilde H}}\}$ of ${\rm Lie}(\tilde H)$. We can assume $\{\tilde e_1,...,\tilde e_{e_{\tilde H}}\}$ is a $W(k)$-basis of $F^1({\rm Lie}(G))$. We write 
$$
\ga_{\tilde h_1}=\sum_{i=1}^{\dim_{W(k)}(P)} x_i\tilde e_i,
$$
with $x_i$ as variables (thought as elements of $T_1/pT_1$). Moreover we can write $\tilde h_1$ mod $p^3$ as follows
$$1_{M\otimes_{W(k)} T_1/p^3T_1}+p\tilde\gamma_{\tilde h_1}+p^2\sum_{j=1}^{e_{\tilde H}} y_jn_j,$$
with $y_j$ as extra variables (thought as elements of $T_1/pT_1$) and with $1_{M\otimes_{W(k)} T_1/p^3T_1}+p\tilde\gamma_{\tilde h_1}\in P(T_1/p^3T_1)$  which mod $p^2$ is $1_{M\otimes_{W(k)} T_1/p^2T_1}+p\gamma_{\tilde h_1}$.
The above equations (18) and (19) get translated in a system of equations involving $x_i$'s and $y_j$'s. It defines an affine ${\rm Spec}(T_1/pT_1)$-scheme ${\rm Spec}(T_1^\prime)$. Usually it is not reduced; below we are interested just in the reduced scheme it defines. 
\smallskip
{\bf J.} 3.9.3.1 allows us to speak about $\tilde G_{W(k)}$-ordinary points of ${\rm Spec}(T_1/pT_1)$; their definition does not depend on which of the two $p$-divisible groups over ${\rm Spec}(T_1)$ (i.e. on which of the morphisms $q(1)$ and $s(1)$) we use to define them (cf. the existence of $t(1)$). From 3.12.1 (applied either in the context of 2.3.11 or in the abstract context of the pull back of $(\Md(U_{\tilde H},\Phi_{R_5},M,F^1,g\vph_0,G),\tilde\psi_M)$ to $\tilde H^0$) we get that they are dense. We have:
\medskip
\item{H)}
{\it ${\rm Spec}(T_1^\prime)$ has finite fibres over $\tilde G_{W(k)}$-ordinary points of ${\rm Spec}(T_1/pT_1)$ with values in $\bar k$, having the same number of points which is a non-negative, integral power of $p$.} 
\medskip
We get this by looking at the system of equations $SE_{\rm init}$ in $x_i$'s and $y_j$'s we got in the following standard manner:
\medskip
-- the right (resp. left) hand sides of (18) and (19) when combined, can be viewed as a product $\tilde h_1g_1$ (resp. $\tilde h_2g_1$) mod $p^2$, where $g_1$ is the $\Phi_{T_1}$-linear endomorphism of $M\otimes_{W(k)} T_1^\wedge$ defined by the formula $\Phi_1=g_1\mu({1\over p})$, with $\mu$ is as in 3.6.13, and where $\tilde h_2\in G(T_1^\wedge)$ is congruent to the identity mod $p$;
\smallskip
-- considering a sum $1_{M\otimes_{W(k)} T_1^\wedge}+p\bigl(\sum_{i=1}^{\dim_{W(k)}(P)} x_i(2)\tilde e_i\bigr)+p\bigl(\sum_{j=1}^{e_{\tilde H}} y_j(2)n_j\bigr)$ which mod $p^2$ is $\tilde h_2$.
\medskip
So, by computing $x_i(2)$'s and $y_j(2)$'s (in terms of Frobenius transforms of $x_i$'s and $y_j$'s) and by identifying the coefficients of $\tilde e_i$'s and $n_j$'s, $SE_{\rm init}$ is obtained by combining the equations of (20) and (21) below
$$x_i=L_i(x_{e_{\tilde H}+1}^p,...,x_{\dim_{W(k)}(P)}^p,y_1^p,...,y_{e_{\tilde H}}^p),\leqno(20)$$
$i\in S(1,\dim_{W(k)}(P))$, and 
$$L_i(x_{e_{\tilde H}+1}^p,...,x_{\dim_{W(k)}(P)}^p,y_1^p,...,y_{e_{\tilde H}}^p)=0,\leqno(21)$$
$i\in S(\dim_{W(k)}(P)+1,\dim_{W(k)}(G))$; warning: it is not of third type. Here $L_i$ is a non-homogeneous linear form with coefficients in $T_1/pT_1$, $i\in S(1,\dim_{W(k)}(G))$. A system of equations with coefficients in an $\FF_p$-algebra having the same form as (20) and (21), is called of fourth type.
\smallskip
The first $e_{\tilde H}$ equations of (20) allow us to eliminate the variables $x_1$,..., $x_{e_{\tilde H}}$. So, over the perfection ${\rm Spec}(T_1/pT_1)^{\rm perf}$ (of ${\rm Spec}(T_1/pT_1)$), the last $\dim_{W(k)}(P)-e_{\tilde H}$ equations of (20) together with the $p$-th roots of the equations of (21), define a system of equations SE in $\dim_{W(k)}(G)-e_{\tilde H}$ variables, which is very close in spirit to one of first type. 
\smallskip
In fact, in any $\tilde G_{W(\bar k)}$-ordinary point of ${\rm Spec}(T_1/pT_1)$ with values in $\bar k$ (it defines uniquely a $\bar k$-valued point of ${\rm Spec}(T_1/pT_1)^{\rm perf}$), SE becomes a system of equations with coefficients in $\bar k$ which is reduced-equivalent to one of first type. Here reduced-equivalent means that the reduced schemes of their affine $\bar k$-schemes of solutions are the same (i.e. are isomorphic). To see this let $(M\otimes_{W(k)} W(\bar k),g(y_G)(\vph_0\otimes 1),\tilde G_{W(\bar k)})$ be the Shimura $\bar\sg$-crystal of the pull back of $\Md_2$ through a $\tilde G_{W(\bar k)}$-ordinary point $y_G:{\rm Spec}(\bar k)\to {\rm Spec}(T_1/pT_1)$. We consider the $\Phi_{T_1}$-Teichm\"uller lift $z_G:{\rm Spec}(W(\bar k))\to {\rm Spec}(T_1^\wedge)$ of $y_G$. The proof of 3.11.1 c) below shows that there is $\tilde h_3\in G(W(k))$ congruent to the identity mod $p$ and such that $\tilde h_3z_G^*(\Phi_1\otimes 1)=z_G^*(\Phi_2\otimes 1)\tilde h_3$. As the replacement of $z_G^*(\tilde h_1)$ by $z_G^*(\tilde h_1)\tilde h_3$ corresponds to a change of the extension to $W(\bar k)$ of the chosen $W(k)$-basis of ${\rm Lie}(P)$ and of ${\rm Lie}(\tilde H)$, we can assume we are dealing with $W(\bar k)$-bases of $P_{W(\bar k)}$ and of $\tilde H_{W(\bar k)}$ such that the $W(\bar k)$-submodules they generate are permuted by $pg(y_G)(\vph_0\otimes 1)$ as in 3.4.3.0. So we get that $y_G^*(SE)$ is equivalent to a system of equations (in $u_i$'s) with coefficients in $\bar k$ of the form
$$a_iu_i=u_{\pi(i)}^p,$$
$i=\overline{e_{\tilde H}+1,\dim_{W(k)}(P)+e_{\tilde H}}$; we have: $a_i=0$ iff $i\ge 1+\dim_{W(k)}(P)$. Here $\pi$ is a suitable permutation of the set $S(1,\dim_{W(k)}(P)+e_{\tilde H})$. Obviously, it is reduced-equivalent to one of first type. So the part of H) pertaining to a power of $p$ (resp. to the same number of points) follows from 3.1.8.1 (resp. from 3.11.1 c)).
\medskip
{\bf K.} G) guarantees that ${\rm Spec}(T_1^\prime)$ is not an empty scheme and that $h_1$ mod $p$ factors through the normalization ${\rm Spec}(T_1^{0\prime})$ of the integral, closed subscheme of ${\rm Spec}(T_1^\prime)$ through which the special fibre of $\tilde Z_{\infty}$ factors naturally. So ${\rm Spec}(T_1^{0\prime})$ is a geometrically connected ${\rm Spec}(T_1/pT_1)$-scheme. Moreover, it is an \'etale cover over an open, dense subscheme $\Mu_1$ of the $\tilde G_{W(k)}$-ordinary locus $\Ml_1$ of ${\rm Spec}(T_1/pT_1)$.  
\smallskip
The same argument as the one used in the proof 3.6.8.1.2 a) applies to give us that we can choose $\Mu_1$ to depend only on the coefficients of the homogeneous parts of the forms $L_i$, $i=\overline{1,\dim_{W(k)}(G)}$. We state this explicitly as a Lemma; its notations are independent of the previous ones and its proof (just an application of 3.6.8.9.0 and 3.6.8.1.2 a)) is left as an exercise. 
\medskip
{\bf Lemma.} {\it Let $m,n\in\NN$, with $m\le n$. We consider a system of equations of fourth type 
$$u_i=L_i(u_1^p,...,u_n^p),$$
$i=\overline{1,m}$ and 
$$0=L_i(u_1^p,...,u_n^p),$$
$i=\overline{m+1,n}$, with coefficients in an integral $\FF_p$-algebra $R$. We assume that its affine ${\rm Spec}(R)$-scheme of solutions has finite non-empty fibres. 
\smallskip
Then there is an open, dense subscheme $U$ of ${\rm Spec}(R)$ such that for any affine, dominant, integral ${\rm Spec}(R)$-scheme ${\rm Spec}(R_1)$ and for every $n$-tuple $(c_1,...,c_n)$ of elements of $R_1$ such that the affine ${\rm Spec}(R_1)$-scheme ${\rm Spec}(R_2)$ of solutions of the system  of equations obtained by replacing above the linear form $L_i$ by $L_i+c_i$, $\forall i\in S(1,n)$, has generically a finite number of points, ${\rm Spec}(R_2)\times_{{\rm Spec}(R)} U$ is a finite scheme over ${\rm Spec}(R_1)\times_{{\rm Spec}(R)} U$.}
\medskip
So we can repeat the arguments, for things modulo $p^n$ ($n\in\NN$, $n>1$). We obtain new systems of equations, involving similar linear forms $L_i(n)$; moreover the homogeneous part of $L_i(n)$ is the same as of $L_i$, $\forall i\in S(1,\dim_{W(k)}(G))$ and $\forall n\in\NN\setminus\{1\}$. 
\smallskip
So, entirely similar to 3.6.8.1.2 a), we deduce that we can assume $\tilde Y_{\infty}$ is an affine ${\rm Spec}(T_1/pT_1)$-scheme which is a non-trivial $\NN$-pro-\'etale cover above $\Mu_1$. In other words, coming back to the paragraph G, we can assume the morphism ${\rm Spec}(T_{\infty}/pT_{\infty})\to {\rm Spec}(Q_6/pQ_6)$ has an image which contains an open subscheme of ${\rm Spec}(Q_6/pQ_6)$. This takes care of 3.6.14.1 c), ending the proof of 3.6.14.1.
\medskip
{\bf 3.6.14.2. Exercise.} Show that 3.6.14 and 3.6.14.1 remain true, as mentioned in 3.6.13, without the assumptions made in 3.6.13. Hint: the only difference is that we have to work with different de Rham tensors, which become linear combination of de Rham components of Hodge cycles, after passage to an \'etale morphism.
\medskip 
{\bf 3.6.14.3. Exercise.} Restate 3.6.14.1 and 3.6.14.1 in terms of a $k$-valued point $y$ of $\Mn_{W(k)}/H_0$, instead of a $W(k)$-valued point $z$ of $\Mn_{W(k)}/H_0$.
\medskip
{\bf 3.6.14.4. Variants.} We do not know when the starting assumptions SA1-2 are satisfied. Instead of open subschemes of $G=G^0_{W(k)}$ or of $\tilde H$, we can work with \'etale schemes over them; so the chances of being able to verify such starting assumptions, in new contexts which can be still of great use, are much higher. When there are no hopes to get any useful starting condition, in 3.6.14 d) and 3.6.14.1 d) we have to work with truncations mod $p^n$ of $p$-divisible objects involved ($n\in\NN$); in such a situation we do not need to pass to $\NN$-pro-\`etale morphisms. There are three variants of working with truncations mod $p^n$: below we detail them in the context of 3.6.14.1 (so that we can use the notations of the proof of 3.6.14.1). Warning: for below variants, no starting assumption is needed. 
\smallskip
In the first variant we entirely ignore the extra Shimura structure and keep track just of principally quasi-polarized, finite, flat, commutative group schemes; it stops at the existence of the $W(k)$-morphism $t(n)$ of 3.6.14.1 E. 
\smallskip
In the second variant we essentially ignore the tensors but we keep track of the extra Shimura (Lie) structures. Let $\tilde Y_{16,n}={\rm Spec}(\tilde R_{16,n})$ be the integral, closed, affine subscheme of $Y_{16,n}$ through which the composite of $m_{\tilde H}$ with the natural
$W(k)$-morphism $Y_{16,\infty}\to Y_{16,n}$ factors. $\tilde Y_{16,n}$ might not be formally smooth over $W(k)$. However, it is faithfully flat over $\ZZ_p$ and so we can still speak about $H^1_j:=H^1_{\rm crys}({\Md_j}_{\tilde Y_{16,n}}/\tilde Y_{16,n}^\wedge)$ and so about the $\tilde R_{16,n}$-submodule ${\got g}_j$ of ${\rm End}(H^1_j/p^nH^1_j)$ obtained naturally by pulling back the Lie submodule of the kernel of the multiplication by $p^n$ of $\Md_1$, $j=\overline{1,2}$. We have a natural identification of $H^1_1$ mod $p^n$ with $H^1_2$ mod $p^n$ and so of their ${\rm End}$'s; these identifications preserve the filtrations. We consider the maximal closed subscheme $\tilde Z_{16,n}$ of $\tilde Y_{{16,n}W_n(k)}$ with the property that under this identification of ${\rm End}$'s, ${\got g}_1$ mod $p^n$ is ${\got g}_2$ mod $p^n$; as ${\got g}_2$ and ${\rm End}(H^1_2)/p^n{\rm End}(H^1_2)$ are free $\tilde R_{16,n}/p^n\tilde R_{16,n}$-modules, its construction is obvious. 
\smallskip
Using $n$ times [BLR, th. 12 of p. 83] we deduce, cf. the existence of $m_{\tilde H}$, that we can assume $t(n)$ factors through $\tilde Y_{16,n}$ in such a way that mod $p^n$ it factors through $\tilde Z_{16,n}$: first time, we achieve this factorization mod $p$, then mod $p^2$, etc. (as in 3.6.14.1 E we can ``put" everything in a finite type context).
We conclude (cf. 2.2.20.1 9)): 
\medskip
{\bf Fact.} {\it We can assume that over ${\rm Spec}(T_n)$, the crystalline counterparts of the kernel of the multiplication by $p^n$ of $(\Md_1,\Mp_{\Md_1})$ and of $(\Md_2,\Mp_{\Md_2})$ become isomorphic.}
\medskip
In the third variant we fully keep track of the tensors: it is easy to state the above second variant in a way that keeps full track of the tensors, provided we either assume $t_{\al}\in\Mt(M)$, $\forall\al\in\Mj$, or we treat the tensors as in 2.2.14 (we just have to define $\tilde Z_{16,n}$ in a way that keeps track of these tensors as well). 
\medskip
{\bf 3.6.14.4.1.} It seems to us possible to prove directly, using (at least in the case $(G,X)$ is of compact type) just the techniques of 3.6.14, the following expected result.
\medskip
{\bf Expectation.} {\it We assume $k=\bar k$. Let $n\in\NN$. There is a $W(k)$-valued point of $\Mn_{W(k)}/H_0$ whose special fibre factors through the same connected component of $\Mn_k$ as $y$, such that the crystalline counterpart (see 2.2.20.1 9)) of the kernel of the multiplication by $p^n$ of the Shimura $p$-divisible group we get over $W(k)$ through it, is isomorphic to the crystalline counterpart of the same type of kernel of the Shimura $p$-divisible group we get through an arbitrary $W(k)$-valued point of ${\rm Spec}(Q_6)$.}
\medskip
{\bf 3.6.15. The completion property.} For the sake of presenting most of the twelve fundamental principles (of 1.2.3) in a reasonable compact form, we now include the following expected result. We continue to refer to the context of 3.6.13 and to assume $k=\bar k$.
Let $y:{\rm Spec}(k)\to\Mn_{W(k)}/H_0$ be defined by $z$. 
\medskip
{\bf A. Expectation (the completion property).} {\it Giving $g_1\in\tilde G^0_{W(k)}(W(k))$, there is a point
$z_1:{\rm Spec}(W(k))\to\Mn_{W(k)}/H_0$ lifting a $k$-valued point of the same connected component $\Mc^0$ of $\Mn_k/H_0$ through which $y$ factors and such that its attached
 principally quasi-polarized Shimura filtered $\sg$-crystal $\bigl(M_1,F^1_1,\vph_1,\tilde G_{W(k)},(t_{1\al})_{\al\in\Mj^\prime},\tilde\psi_{M_1}\bigr)$ is $1_{\Mj^\prime}$-isomorphic to
$\bigl(M,F^1,g_1\vph_0,\tilde G_{W(k)},(t_{\al})_{\al\in\Mj^\prime},\tilde\psi_M\bigr)$.}
\medskip
This Expectation is proved in \S 11; see also 4.12 below. The initial motivation for this expectation is based on 3.6.14.4.1 and on the following Fundamental Lemma.
\medskip
{\bf B. Fundamental Lemma.} {\it Let $(M,\vph,G)$ (resp. $(M,\vph,G,p_M)$) be an arbitrary $\sg$-$\Ms$-crystal over $k$ (resp. an arbitrary $\sg$-$\Ms$-crystal over $k$ endowed with a perfect bilinear form $p_M:(M,\vph)\otimes_{W(k)} (M,\vph)\to W(k)(n(\vph))$ normalized by $G$ and for which the $W$-condition holds; here $n(\vph)\in\ZZ$). Then there is $n\in\NN$, such that for any $h_n\in G(W(k))$ congruent to the identity mod $p^n$ (resp. congruent to the identity mod $p^n$ and centralizing $p_M$), $(M,h_n\vph,G)$ (resp. $(M,h_n\vph,G,p_M)$) is isomorphic to $(M,\vph,G)$ (resp. to $(M,h_n\vph,G,p_M)$) under an isomorphism defined by an element of $G(W(k))$.}
\medskip
{\bf Proof:}
We can assume $v:=\dim_{B(k)}({\rm Lie}(G_{B(k)}))$ is at least 1. We first deal with the context without $p_M$. Let $s_L(\vph)\in\NN\cup\{0\}$ be the smallest number such that $p^{s_L(\vph)}\vph$ takes ${\rm Lie}(G)$ onto itself. Using Dieudonn\'e's classification of isocrystals over $k$ for the isocrystal $({\rm Lie}(G_{B(k)}),\vph)$, we consider a $B(k)$-basis 
$$\Mb=\{e_1,...,e_v\}$$ 
of ${\rm Lie}(G_{B(k)})$ formed by elements of ${\rm Lie}(G)$ such that there is a permutation $\pi_L$ of $S(1,v)$ with the property that it can be decomposed into a product of disjoint cycles $\Mc_1$,..., $\Mc_s$, with $s\in\NN$, each such cycle $\Mc_{i_0}=(a_1,...,a_q)$, with $i_0\in S(1,s)$, having the property that $\vph(e_{a_j})=e_{a_{j+1}}$ if $j\in S(1,q-1)$ and $\vph(e_{a_q})=p^{n((a_1,...,a_q))}{e_{a_1}}$, with 
$$n((a_1,...,a_q))\in\ZZ$$ 
relatively prime to $q\in\NN$. Let $m\in\NN$ be such that $p^m{\rm Lie}(G)$ is included in the $W(k)$-span $E$ of $\Mb$. We take $n\in\NN$ such that 
$$n\ge 2m+1+s_L(\vph).\leqno (EST)$$ 
\indent  
We want to show that there is $h\in G(W(k))$ congruent to the identity mod $p^{m+1+s_L(\vph)}$ and such that $hh_n\vph=\vph h$. By induction on $t\in\NN$ we show that there is $h\in G(W(k))$ congruent to the identity mod $p^{n+t-1-m}$ and such that $hh_n\vph h^{-1}\vph^{-1}\in G(W(k))$ is congruent to the identity mod $p^{n+t}$. It is enough to deal with the case $t=1$. We use what we call the stairs method (modeled on $<e_1,...,e_v>$).
\smallskip
We take $h=\tilde h^1\tilde h^2$, where 
$$\tilde h^1:=\prod_{l\in S(1,v)} (1_M+x_lp^{u_l}e_l),\leqno (PROD)$$ 
with $u_l\in\NN$ depending only on the cycle $\Mc$ of $\pi_L$ to which $l$ belongs and with $x_l\in W(k)$, and where $\tilde h^2\in GL(M)(W(k))$ is congruent to the identity mod $p^{n+1+s_L(\vph)}$. We always take $u_l$ as the maximal possible value allowed by $h_n$ in the sense that $h_n$, up to elements of $GL(M)(W(k))$ congruent to the identity mod $p^{n+1+s_L(\vph)}$, can be similarly written down as such a product $h_n=h_n^1h_n^2$. In particular, we always have $u_l\ge n-m\ge m+1+s_L(\vph)$ and so the order in which products of the above form (PROD) are taken is irrelevant mod $p^{2(n-m)}$ and so mod $p^{n+1+s_L(\vph)}$; this also guarantees that $h_n$, up to elements congruent to the identity mod $p^{n+1+s_L(\vph)}$, can be similarly written down as such a product. Warning: $\tilde h^1$ is not in general an element of $G(W(k))$; however, mod $p^{n+1+s_L(\vph)}$ it is and so, as $G$ is smooth, we can always choose $\tilde h^2$ such that $h\in G(W(k))$. The expression of $hh_n\vph h^{-1}\vph^{-1}\in G(W(k))$ mod $p^{n+1}$ does not depend on $\tilde h^2$.
\smallskip
Moreover, if $h^\prime=\prod_{l\in S(1,v)} (1+x_l^{\prime}p^{u_l}e_l)$, with all $x_l^\prime$'s in $W(k)$, then 
$hh^\prime$ is congruent mod $p^{n+1+s_L(\vph)}$ to $\prod_{l\in S(1,v)} [1+(x_l^\prime+x_l)p^{u_l}e_l]$. We call this property as $ADD$.
\smallskip
We need to show that we can choose $x_l$'s such that $hh_n\vph h^{-1}\vph^{-1}\in G(W(k))$ is congruent to the identity mod $p^{n+1}$. This boils down to: if $h_n$ is not congruent to the identity mod $p^{n+1}$, then we can choose $x_l$'s such that writing $hh_n\vph h^{-1}\vph^{-1}\in G(W(k))$ as a product as in (PROD), we can replace $u_l$ by $u_l+1$, $\forall l\in S(1,v)$. Concentrating on a fixed cycle $\Mc:=\Mc_{i_0}$ of $\pi_L$, we need to show the existence of solutions of a system of equations over $k$ of the following circular form
$$b_j\tilde x_{a_j}=d_j\tilde x_{a_{j-1}}^p+c_j,\leqno (22)$$
 $j\in S(1,q)$ (here $q+1=1$), where $c_j$, $d_j\in k$, and with:
\medskip
-- $b_j=1$ for all $j\in S(1,q)$, if $n(\Mc)\ge 0$, and with 
\smallskip
-- $b_j=1$ for all $j\in S(1,q-1)$ and with $b_q=0$ and all $d_j$'s equal to $1$, if $n(\Mc)<0$. 
\medskip
The above system  of equations is viewed over $k$ and its shape can be obtained as usual, by working mod $p$ in the context of the $W(k)$-submodule of $E$ generated by all these $p^{u_l}e_l$'s. In other words, based on $ADD$, the variable $\tilde x_{a_j}$ is defined as:
\medskip
-- $x_{a_j}$ mod $p$, if $n(\Mc)\ge 0$;
\smallskip
-- $x_{a_j}$ mod $p$ for $j\neq q$ and as $y_{a_q}$ mod $p$ for $j=q$, where $y_{a_q}\in W(k)$ is defined by the identity $x_{a_q}=p^{-n(\Mc)}y_{a_q}$ (i.e. we take $x_{a_q}$ to be divisible by $p^{-n(\Mc)}$), if $n(\Mc)<0$.
\medskip
So, we just need to show that such a system  of equations does have a solution: in case $n(\Mc)\ge 0$, 3.6.8.1 (with $s=0$ and $l=1$) applies directly (as $k=\bar k$), while in case $n(\Mc)<0$ we can just apply the fact that $k$ is perfect; if $n(\Mc)\ge 0$ it is easy to see, without appealing to 3.6.8.1, that (22) does have solutions in $k$. This ends the proof in the context without $p_M$.
\smallskip
We now deal with the context involving $p_M$. We can assume $G$ does not fix $p_M$. Let $T$ be a 1 dimensional torus of $Z(G)$ which does not fix $p_M$ and which has a split $\ZZ_p$-structure $T_{\ZZ_p}$ such that the elements of $T_{\ZZ_p}(\ZZ_p)=\GG_m(\ZZ_p)$ are naturally isomorphisms of $(M,\vph,G)$; this is nothing else but end of 2.2.9 12) adapted to $(M,\vph,G,p_M)$. 
Let $r_f$ be the rank of the maximal finite, flat group subscheme of $T$ fixing $p_M$. Let 
$$m(r_f)\in\NN$$ 
be the smallest number such that any element of $\GG_m(\ZZ_p)$ congruent to $1$ mod $p^{m(r_f)}$ is the $r_f$-th power of an element of $\GG_m(\ZZ_p)$. We get:
\medskip
{\bf Fact.} {\it $\forall\al\in\GG_m(\ZZ_p)$ congruent to the identity modulo $p^{m(r_f)}$, $(M,\vph,G,p_M)$ and $(M,\vph,G,{\al}p_M)$ are isomorphic under an isomorphism defined by an element of $T(\ZZ_p)$.} 
\medskip
So the Lemma in the context involving $p_M$ follows from this Fact and from its part in the context without $p_M$: we just need to take $n\in\NN$ such that
$$n\ge\max\{2m+1+s_L(\vph),m(r_f)\}.\leqno (ESTP)$$
\medskip
{\bf 3.6.16. Definitions and remarks. 1)} We say $\Mn$ has the completion property, if 3.6.15 A holds. 
\smallskip
{\bf 2)} The smallest possible value of the number $n$ of 3.6.15 B, is called the isomorphism deviation of $(M,\vph,G)$ and is denoted by $isom-d(M,\vph,G)$.
\smallskip
{\bf 3)} To our knowledge, 3.6.15 B was not known previously even for the classical situation when $G=GSp(M,p_M)$ and $(M,\vph)$ is an object of $p-\Mm_{[0,1]}(W(k))$, with $M$ of rank big enough (like greater than $7$). 
\medskip
{\bf 3.6.16.1. Terminology.} All connections used in 3.6.1-16 as well as in 2.2-3 are called (in the honor of their source [Fa2, th. 10 and the remarks after]) Faltings connections. 
\medskip
{\bf 3.6.17. Solution of 3.6.6.0.} We proceed by induction. We can assume $g_2\in P_0(W(k))$ and $k=\bar k$. So we can replace $\vph_0$ by $\vph_1$ (we recall that $P_0=P_1$, cf. 3.3.2). But $(M,F^1,g_2\vph_1)$ is itself cyclic diagonalizable: based on 3.3.4 and 2.2.3 3), the elements of ${\got g}[{1\over p}]$ fixed by $g_2\vph_1$ take $F^1[{1\over p}]$ into itself and form a reductive $\QQ_p$-Lie subalgebra of ${\got g}[{1\over p}]$ of the same rank as ${\got g}[{1\over p}]$. So 2.2.18 applies and so (cf. its proof or cf. the Fact of 2.2.22 1)) there is a torus of $P_0$ whose Lie algebra is normalized by $g_2\vph_0$. Using $\ZZ_p$-structures as in 2.2.9 8), we get that we can assume it is a maximal torus. Arguments as in the proof of Fact 1 of 2.2.9 3), show that, up to isomorphisms defined by elements of $P_0(W(k))$, we can assume it is $T$ itself. So $g_2\in P_0(W(k))$ normalizes $T$ and so $g_2\vph_1$ has all properties of $\vph_1$ we needed to get 3.4.3.0. Not to introduce extra notations, we can assume $g_2=1_M$.
\smallskip
Let $g_3\in G(W(k))$ be congruent to the identity mod $p^m$, $m\in\NN$. So $g_3$ mod $p^{m+1}$ can be written down as $1_M+p^{m}x$, with $x=x_1+x_2$, where $x_1\in {\got p}_0$ and $x_2$ belongs to the Lie algebra ${\got p}_{<0}$ of ${\got g}$ corresponding to negative slopes of $({\got g},\vph_1)$. We take $h\in G(W(k))$ such that mod $p^{m+2}$ is 
$$1_M+p^{m}e_0+p^{m}\sum_{s=1}^{s_1} y_se_s+p^{m+1}\sum_{s=s_1+1}^{\dim_{W(k)}({\got p}_{<0})} y_se_s,$$ 
where $e_0\in {\got p}_0$, $y_s$'s are elements of $W(k)$, $s_1:=\sum_{i\in I} \dim_{W(k)}(M_i)$ and $\{e_1,...,e_{s_1}\}$ (resp. $\{e_{s_1+1},...,e_{\dim_{W(k)}({\got p}_{<0})}\}$) is a $W(k)$-basis of $\prod_{i\in I} M_i$ (resp. of $\prod_{i\in I} N_i^-$) formed by elements normalized by $T$ (cf. the proof of 3.6.6 for notations). 
\smallskip
Let $\Mb_0$ be an arbitrary $W(k)$-basis of ${\got p}_0$. Let $\Mb_1$ be the $W(k)$-basis of ${\rm Lie}(G)$ obtained by putting $\Mb_0$ and $\{e_1,...,e_{\dim_{W(k)}({\got p}_{<0})}\}$ together. As in the proof of 3.6.15 B, the congruence $hg_3\vph_1h^{-1}-1_M\equiv 0$ mod $p^{m+1}$ gets translated in a system of equations $SE$ with coefficients in $k$ in the $d$ variables obtained by taking the reduction mod $p$ of all $y_s$'s and of the coefficients $z_s$'s, $s\in S(\dim_{W(k)}({\got p}_{<0})+1,d)$, of $e_0$ w.r.t. $\Mb_0$; the equations of $SE$ are obtained by identifying with $0$ the coefficients of $hg_3\vph_1h^{-1}-1_M$ mod $p^{m+1}$ w.r.t. reduction mod $p^{m+1}$ of $\Mb_1$. $SE$ is a direct sum of two subsystems $SE_1$ (involving just the variables $z_s$'s) and $SE_2$ (involving just the variables $y_s$'s). $SE_1$ is a system of equations of first type in $\dim_{W(k)}({\got p}_0)$ variables. Moreover, as $k=\bar k$, $SE_2$ is a disjoint union of systems of equations of exactly the same type as the ones obtained in the proof of 3.6.15 B, corresponding to $n(\Mc)<0$ (cf. the Weyl decomposition of 3.4.3.0). So as in the mentioned proof, $SE_1$ and $SE_2$ have solutions in $k$. So using induction, up to inner isomorphisms, we can assume $g_3=1_M$. This solves 3.6.6.0.
\medskip
{\bf 3.6.18. General principles for connections in Fontaine categories.} Let $m\in\NN\cup\{0\}$. The notations to be used in 3.6.18-20 are entirely independent from the previous notations of \S 3, except that $k$ and $\sg$ will have the same meaning. We keep assuming $p\ge 3$, though all arguments below, except very few exceptions, work for $p=2$ as well; these exceptions are listed in 3.14. In this section
$$R:=W(k)[[z_1,...,z_m]]$$ 
is the ring of formal power series in $m$ variables with coefficients in $W(k)$. Let $\Phi_R$ be an arbitrary Frobenius lift of $R$. 
\smallskip
Any object of a category $\Mm\Mf_{[a,b]}(R)$ lifts to a $p$-divisible object (see Fact of 2.2.1.1 6) and 2.2.1 f)). So we prefer to look at such objects as being obtained from truncations (mod $p^n$, for some $n\in\NN\cup\{0\}$) of $p$-divisible objects, by a repeated finite process of extensions (via short exact sequences). 
\smallskip
We start with a $p$-divisible object 
$$
{\got C}:=(M\otimes_{W(k)} R,F^1\otimes_{W(k)} R,\Phi),
$$ 
of $\Mm\Mf_{[0,1]}(R)$. Here $M$ is a free $W(k)$-module of finite rank, $F^1$ is a direct summand of $M$, while $\Phi$ is a $\Phi_R$-linear endomorphism of $M\otimes_{W(k)} R$ producing an isomorphism (still denoted by $\Phi$) of $R$-modules 
$$\Phi:(M+{1\over p}F^1)\otimes_R {}_{\Phi_R}R\tilde\to M\otimes_{W(k)} R.$$ 
A natural question arises:
\medskip
{\bf Q} Is there a $p$-divisible group over $R$, whose corresponding $p$-divisible object of $\Mm\Mf_{[0,1]}^\nabla(R)$, when viewed just as a $p$-divisible object of $\Mm\Mf_{[0,1]}(R)$, is $\got C$?
\medskip
If its answer is yes, we say ${\got C}$ is induced from a $p$-divisible group. Similarly we speak about $\got C$ being induced (cf. also 2.2.1.3) from a $p$-divisible object of $\Mm\Mf_{[0,1]}^\nabla(R_1)$ via a $W(k)$-homomorphism $R_1\to R$, where $R_1$ is a regular, formally smooth $W(k)$-algebra.
\medskip
{\bf 3.6.18.0. Types of Frobenius lifts of $R$.}  [Fa2, th. 10] can be restated: the answer to {\bf Q} is yes provided $\Phi_R$ is of the form $\Phi_R(z_i)=z_i^p$, $i=\overline{1,m}$. This assumption (on the special shape of $\Phi_R$) has been made in [Fa2, ch. 7] (quoting) ``for technical reasons".
\smallskip
Let $S(m):=S(1,m)$ and let $S(M):=S(1,\dim_{W(k)}(M))$. Any Frobenius lift $\Phi_R$ of $R$ is uniquely determined by equations
$$\Phi_R(z_i)=z_i^p+pP_i(z_1,...,z_m)+pL_i(z_1,...,z_m)+p\ga_i,\leqno (23)$$ 
where $P_i=\sum_{l=2}^\infty Q_{il}$, with $Q_{il}$ a homogeneous polynomial of degree $l\ge 2$ in $m$ variables, where $L_i$ is a linear homogeneous polynomial and $\ga_i\in W(k)$, $\forall i\in S(m)$. It is easy to see that by replacing $z_i$ by $z_i+p\dl_i$, with $\dl_i\in W(k)$ conveniently chosen, $\forall i\in S(m)$, we can assume $\ga_i=0$, $\forall i\in S(m)$: by induction on $l\in\NN$ we can check that we can assume $p^l$ divides $p\ga_i$; this is the general argument showing the existence of (the well known) Teichm\"uller lifts. We assume this from now on. 
\smallskip
We can write the Frobenius endomorphism $\Phi$ of $M\otimes_{W(k)} R$ in the form 
$$
\Phi=g_R(\vph\otimes 1),$$ 
where $g_R\in GL(M)(R)$ modulo the ideal $(z_1,...,z_m)$ of $R$ is the identity element $1_M$ of $G(W(k))$ and where $\vph$ is the $\sg$-linear endomorphism of $M$ obtained from $\Phi$ by taking $z_i=0$, $\forall i\in S(m)$. So what we are trying to do (for answering {\bf Q}), is nothing else but to deform over $R$ the $p$-divisible group $D_{W(k)}:=\DD^{-1}(M,F^1,\vph)$ over $W(k)$, in the way ``prescribed" by $\got C$. The use of $\DD^{-1}$ is part of the second place (see 3.6.2) where we need $p>2$.  
\smallskip
We write 
$$L_i(z_1,...,z_m)=\sum_{j=1}^m a_{ij}z_j,$$ 
with $a_{ij}\in W(k)$, $\forall i,j\in S(m)$. We obtain a matrix $A(\Phi_R)\in M_m(W(k))$, having $a_{ij}$ as its entries. Let 
$$\bar A(\Phi_R)\in M_m(k)$$ 
be the matrix obtained by taking $A(\Phi_R)$ mod $p$. We call $A(\Phi_R)$ (resp. $\bar A(\Phi_R)$) the linear part (resp. the $k$-linear part) of $\Phi_R$ w.r.t. $z_1,...,z_m$. To study how $A(\Phi_R)$ changes under a $W(k)$-automorphism of $R$ which modulo its ideal $(z_1,...,z_m)$ is the identity, we can restrict to linear $W(k)$-automorphisms, i.e. to $W(k)$-automorphisms of $R$ taking $z_i$ to $\sum_{j=1}^m b_{ij}z_j$, with $b_{ij}\in W(k)$ defining an invertible matrix $B\in GL_m(W(k))$ (argument: a $W(k)$-automorphism of $R$ which is the identity modulo $(z_1,...,z_m)^2$, does not changes $A(\Phi_R)$). Under such a linear $W(k)$-automorphism the matrix $A(\Phi_R)$ is replaced by $\sg(B)A(\Phi_R)B^{-1}$; so $\bar A(\Phi_R)$ is replaced by 
$$\bar B^{[p]}\bar A(\Phi_R)\bar B^{-1},$$ 
with $\bar B\in GL_m(k)$ equal to $B$ mod $p$ (see end of 2.1 for the meaning of $\bar B^{[p]}$). 
\smallskip
In what follows it is important only how $\bar A(\Phi_R)$ looks like.
\medskip
{\bf 3.6.18.0.1. Terminology.} If $\Phi_R(z_i)=z_i^p$ (resp. if $\Phi_R(z_i)=(z_i+1)^p-1$), $\forall i\in S(m)$, we speak about an additive (resp. multiplicative) type Frobenius lift of $R$. If $A(\Phi_R)$ (resp. $\bar A(\Phi_R)$) is an invertible (resp. a nilpotent) matrix, we speak about an essentially multiplicative (resp. additive) type Frobenius lift of $R$. We also say that $\Phi_R$ is of (essentially) multiplicative (or additive) type. 
\medskip
{\bf 3.6.18.0.1.1. The global context.} Here, as a digression, we globalize the above terminology. We consider now an arbitrary regular, formally smooth $W(\tilde k)$-scheme $\tilde X$ equipped with a Frobenius lift $\Phi_{\tilde X}$, where $\tilde k$ is an arbitrary perfect field. We consider a closed embedding $\tilde y:{\rm Spec}(\tilde k_1)\hookrightarrow\tilde X$, with $\tilde k_1$ an algebraic field extension of $\tilde k$. We say $\Phi_{\tilde X}$ is a multiplicative (or additive, or essentially multiplicative or additive) type Frobenius lift of $\tilde X$ in $\tilde y$, if the Frobenius lift of (the $W(\tilde k)$-algebra of global sections of) the completion of $\tilde X$ in $\tilde y$ induced naturally by $\Phi_{\tilde X}$, is of the corresponding type.
\medskip
{\bf 3.6.18.1. Notations.} Let $f(0):=\dim_{W(k)}(M/F^1)$ and $f(1):=\dim_{W(k)}(F^1)$. We denote by  $s(i)$ the multiplicity  of the slope $i$ for $(M,\vph)$, $i\in\{0,1\}$. So $s(0)s(1)$ is the multiplicity of the slope $-1$ of the Lie $\sg$-crystal $({\rm End}(M),\vph)$. Let $r$ be the number of non-zero eigenvalues (counted with their multiplicities) in $\bar k$ of $\bar A(\Phi_R)$; we refer to it as the rank of $\Phi_R$. 
\smallskip
{\bf 3.6.18.1.1. Remarks.} It is easy to see that $r$ is well defined regardless of all choices (cf. also 3.6.18.4 B) below). In particular the notion of essentially additive type Frobenius lift of $R$ (case $r=0$) is well defined. We can have $r=0$ without $\bar A(\Phi_R)$ being $0$. So, if $m\ge 2$, we have plenty of essentially additive type Frobenius lifts of $R$ which are not of additive type. 
\smallskip
On the other, if $k=\bar k$ it can be easily checked that, up to a suitable change of variables, the notions of multiplicative type and of essentially multiplicative type coincide. We sketch the argument for this well known fact. If $\Phi(z_i)=(z_i+1)^p-1$, $i=\overline{1,m}$, and if $\Phi_1$ is an essentially multiplicative type Frobenius lift of $R$ taking $(z_1,...,z_m)$ into itself, then by induction on $q\in\NN$ we show  that (up to a change of variables) we can assume $\Phi(z_i)-\Phi_1(z_i)\in p(z_1,...,z_m)^{q+1}$, $\forall i\in S(1,m)$. For $q=1$ this is a consequence of the fact that for any $A\in GL_m(W(k))$, there is $B\in GL_m(W(k))$ such that $\sg(B)B^{-1}=A$. Assuming that the statement is true for $q-1$, to prove it for $q$, we show by induction on $l\in\NN$ that we can assume $\Phi(z_i)-\Phi_1(z_i)\in p^{l+1}(z_1,...,z_m)^{q}+p(z_1,...,z_m)^{q+1}$, $\forall i\in S(1,m)$: if $\Phi(z_i)-\Phi_1(z_i)\equiv p^lP_i(z_1,...,z_m)$ modulo $p^{l+1}(z_1,...,z_m)^q+p(z_1,...,z_m)^{q+1}$, with $P_i$'s as homogeneous forms of degree $q$, then replacing $z_i$ by $z_i-p^{l-1}P_i(z_1,...,z_m)$ we get that we can assume that $P_i$'s are all $0$. As $R_m$ is $(p,z_1,...,z_m)$-adically complete, the conclusion follows.
\smallskip
But this is not so if $k\neq\bar k$; example: if $k=\FF_p$ then $\Phi_R(z_i)=(z_i+1)^p-1+p^2z_i$ is not of multiplicative type, as its linear part w.r.t. $z_1,...,z_m$ is not conjugate over $\ZZ_p$ to the identity matrix of $M_m(\ZZ_p)$.
\medskip
{\bf 3.6.18.1.2. Exercise.} We assume $r=0$. Show that there is a sequence of ideals $(I_n)_{n\in\NN\cup\{0\}}$ of $R$ such that $I_0=R$, $\cap_{n\in\NN} I_n=pR$, $I_{n-1}/I_{n}$ is isomorphic to $k$ as an $R$-module and $d_{\Phi_R*}/p(I_n)\subset I_{n+1}\Om_{R/W(k)}^\wedge$, $\forall n\in\NN$. Hint: we can assume $\bar A(\Phi_R)$ is nilpotent and lower triangular; so $d_{\Phi_R*}/p(z_i)\in \bigl(p,z_{i+1},...,z_m)+(z_{js})_{j,s\in S(1,i)}\bigr)\Om_{R/W(k)}^\wedge$, $\forall i\in S(1,m)$ (so we can take $I_1=(p,z_1^2,z_2,...,z_m)$, $I_2=(p,z_1^2,z_2^2,z_1z_2,z_3,z_4,....,z_m)$,..., $I_m=(p)+(z_1,...,z_m)^2$, $I_{m+1}=(p)+(z_2,...,z_m)^2+(z_1^3)+z_1(z_2,...,z_m)$, etc.).
\medskip
{\bf 3.6.18.2. The generic situation.} This is obtained by requiring:
\medskip
a) $(M,\vph)$ to have all slopes $0$ and $1$ and $(M,F^1,\vph)$ to be diagonalizable;
\smallskip
b) $\Phi_R(z_i)=(z_i+1)^p-1$, $\forall i\in S(m)$; 
\smallskip
c) $k=\bar k$.
\medskip
We have:
\medskip
{\bf Fact.} {\it For any $n\in\NN$, there are precisely 
$$p^{nmf(0)f(1)}$$ connections on ${\got C}/p^n{\got C}$, regardless of the shape of $g_R$. They are all integrable and nilpotent mod $p$.}
\medskip
{\bf Proof:} The first part is just a very particular case of 3.6.18.4 B) below (which is proved independently of this proof). However, we felt it is appropriate to include a separate proof of it, to make the connection to previous work (in particular, see [De3]) and to prepare the reader for the general case. We consider a $W(k)$-basis $\Mb=\{e_1,...,e_{\dim_{W(k)}(M)}\}$ of $M$ such that $\vph(e_i)=p^{\vep_i}e_i$, with $\vep_i\in\{0,1\}$, $\forall i\in S(1,\dim_{W(k)}(M))$. Let $I^0$ (resp. $I^1$) be the set of those $i$'s such that $\vep_i=0$ (resp. $\vep_i=1$); so $f(0)=\abs{I^0}$ and $f(1)=\abs{I^1}$. We can assume that $g_R$ fixes all $e_i$'s with $i\in I^0$ and that, $\forall i\in I^1$ we have 
$$g_R(e_i)-e_i=\sum_{j\in I^0} a_{ij}e_j,$$ 
with all $a_{ij}$'s in $R$. The argument for this is entirely the same as the inductive argument of 3.6.18.1.1.
\smallskip
But any connection $\nabla_n$ on $M\otimes_{W(k)} R/p^nR$ which makes ${\got C}/p^n{\got C}$ potentially to be viewed as an object of $\Mm\Mf_{[0,1]}(R)$ annihilates all $e_i$'s with $i\in I^0$ and, $\forall i\in I^1$ it takes $e_i$ into 
$$\sum_{l=1}^m \sum_{j\in I^0} a_{ijl}e_jdz_l,\leqno (NABLA)$$ 
with $a_{ijl}\in R/p^nR$ such that mod $p^n$ we have
$$a_{ijl}=\Phi_R(a_{ijl})(z_l+1)^{p-1}+b_{ijl},\leqno (24)$$
$\forall (j,l)\in I_0\times S(1,m)$. 
Here $b_{ijl}$'s are elements of $R$ such that we have equalities 
$$\sum_{l=1}^mb_{ijl}dz_l=da_{ij},\leqno (25)$$
$\forall (i,j)\in I_1\times I_0$. Starting from the shape of $\nabla$ we get that, $\forall l\in S(1,m)$ and $\forall f_r\in R$, we have:
$$\nabla({\partial\over {\partial z_l}})^{p}(f_re_i)\equiv 0 \, {\rm mod}\,\, p\leqno (N0)$$
if $i\in I_0$ and
$$\nabla({\partial\over {\partial z_l}})^{2p-1}(f_re_i)\equiv 0 \, {\rm mod}\,\, p\leqno (N1)$$
if $i\in I_1$. So $\nabla_n$ is nilpotent mod $p$ (we could just quote 3.6.1.1 2) but (N1) is slightly more precise than loc. cit.). To check that $\nabla_n$ is integrable, we just need to show that the differential form 
$$w:=\sum_{(i,j,l)\in I_1\times I_0\times S(1,m)} a_{ijl}e_{ij}dz_l\leqno (DF)$$
is closed, where $\{e_{ij}|i,j\in S(1,\dim_{W(k)}(M))\}$ is the $W(k)$-basis of ${\rm End}(M)$ defined naturally (as in 3.6.8 3)) by $\Mb$. But from (24) we get that 
$$dw=\sum_{(i,j,l)\in I_1\times I_0\times S(1,m)} \Phi_R(da_{ijl})e_{ij}(z_l+1)^{p-1}dz_l=\Phi_{R*}/p(dw).\leqno (EQ)$$
So the differential $2$-form $dw$ is zero mod $p$. Plugging this in the right hand side of (EQ) we get that $dw$ is zero mod $p^2$, etc. Using mathematical induction we get $dw=0$.
\smallskip  
As $\Phi_R(a_{ijl})$ mod $p$ is $a_{ijl}^p$ mod $p$ and as $z_l+1$ is invertible in $R$, the number of such connections $\nabla_1$ is $p^{ms(0)s(1)}$ (cf. (24)). But using the constancy property (see 3.6.8.9), 3.6.8.1.2 a) and c) and induction on $n\in\NN$, the Fact gets reduced to the case $n=1$. This proves the Fact.
\medskip
{\bf 3.6.18.2.1. Two philosophies.} The reason we call this the generic situation is related to [De3, 1.4.7]. There are two philosophies at play: 
\medskip
-- any $p$-divisible group over $k$ is the specialization of an ordinary $p$-divisible group (cf. 3.1.8 and 3.1.8.1);   
\smallskip
-- whatever we can prove for the standard situation and does not involve the precise non-negative, integral power of $p$ of the number of solutions we get for different systems of equations (in the other cases we get smaller powers of $p$), can be proved in the general case (in particular we always get solutions).
\medskip
As samples for this second philosophy see 3.6.18.4 below. This second philosophy works also in the more general context where the category $\Mm\Mf_{[0,1]}(R)$ is replaced by another one $\Mm\Mf_{[a,b]}(R)$, with $a,b\in\ZZ$, $a<b$ (so we keep 3.6.18.2 b) and c), while 3.6.18.2 a) is modified accordingly; see 3.6.18.5.5 below).  
\medskip
{\bf 3.6.18.2.2. Remark.} Let ${\got C}_0$ be an object of $\Mm\Mf_{[0,1]}(R)$ which is the extension of an object of $\Mm\Mf_{[1,1]}(R)$ by an object of $\Mm\Mf_{[0,0]}(R)$. Then the same arguments of the proof of 3.6.18.2 apply to get that any connection on ${\got C}_0$ is integrable.
\medskip
{\bf 3.6.18.3. Example.} It is easy to see that the question {\bf Q} does not always have a positive answer. Here is a simple example for which the answer to {\bf Q} is no. Let $R=\ZZ_p[[t]]$ with $\Phi_R(t-1)=(t-1)^p$. We take $M$ to be a free $\ZZ_p$-module of dimension 2. Let $\{x,y\}$ be a $\ZZ_p$-basis of $M$ such that $\{y\}$ is a $\ZZ_p$-basis of $F^1$. Let $\Phi(x)=x$ and $\Phi(y)=py+ptx$. There is no connection $\nabla$ on $M\otimes_{\ZZ_p} R$ for which $\Phi$ is $\nabla$-parallel. Argument: such a connection must annihilate $x$ and takes $y$ into $axdt$ (to be compared with 3.6.1.4 2)), where $a\in R$, when taken mod $p$, must satisfy the equation $u^p(t-1)^{p-1}=u+1$ in $u$ with coefficients in $R/pR$; as the equation $v^p=v+1$ has no solution in $\FF_p$, no such $a$ exists. On the other hand, it can be checked by directly computing $a$ (for instance, cf. the proof of 3.6.18.2; see also 3.6.18.4 below) that, if we move from $\ZZ_p$ to $W(\FF)$ (i.e. if we replace $\FF_p$ by its algebraic closure), $(M\otimes_{\ZZ_p} W(\FF)[[t]],F^1\otimes_{\ZZ_p} W(\FF)[[t]],\Phi\otimes 1)$ is induced from a $p$-divisible group over $W(\FF)[[t]]$. So the ``obstruction" in the existence of $\nabla$ is: the field $\FF_p$ is not ``big enough". 
\medskip
{\bf 3.6.18.3.1. Lemma.} {\it If there is an integrable  connection $\nabla$ on $M\otimes_{W(k)} R$ which makes $\got C$ to be viewed as a $p$-divisible object of $\Mm\Mf_{[0,1]}^\nabla(R)$, then ${\got C}$ is induced from a $p$-divisible group.}
\medskip
{\bf Proof:} The connection $\nabla$ is also nilpotent mod $p$, cf. [Fa1, p. 34] or 3.6.1.1.2. So it allows us to change the Frobenius lift of $R$ to one of additive type, cf. [De3, (1.1.2.1)]. Now everything results from 2.2.21 UP.
\medskip
{\bf 3.6.18.4. Theorem (the universal local $\nabla$ principle for $\Mm\Mf_{[0,1]}$).} {\it We assume $k=\bar k$. We have:
\smallskip
{\bf A)} There are connections on $\got C$.
\smallskip
{\bf B)} Let $n\in\NN$. The number of connections on ${\got C}/p^n{\got C}$ is $p^{nrs(0)s(1)}$; this number is smaller or equal to $p^{nmf(0)f(1)}$. Each such connection lifts to precisely $p^{rs(0)s(1)}$ connections on ${\got C}/p^{n+1}{\got C}$.}
\medskip
{\bf Proof:} A) follows from B); in fact its proof is contained in the proof of 3.6.1.3, cf. 3.6.8.4 1) to 3) and the constancy property of 3.6.8.9. For future references the proof of B) is divided into parts, indexed by numbers attached to the right of the letter $P$. 
\medskip
{\bf P1.} The inequality part of B) is obvious. From 3.6.8.1.2 c) and the constancy property of 3.6.8.9 we deduce that for the counting part of B) it is enough to show that the number of connections on ${\got C}/p{\got C}$ is precisely $p^{rs(0)s(1)}$. As $k=\bar k$, from 3.6.18.0 we deduce that we can assume the matrix $\bar A(\Phi_R)$ is upper triangular. To construct (get) a connection on ${\got C}/p{\got C}$ is the same thing (see 3.6.8 (3)) as getting a solution in $R/pR$ of a system of equations of the form
$$x_{ijl}=L_{ijl}(x^p_{111},...,x^p_{\dim_{W(k)}(M)\dim_{W(k)}(M)m})+c_{ijl},$$  
$(i,j,l)\in S(M)\times S(M)\times S(m)$, where $L_{ijl}$ are linear homogeneous forms in $m\dim_{W(k)}(M)^2$ variables and where $c_{ijl}\in R/pR$, $\forall (i,j,l)\in S(M)\times S(M)\times S(m)$. 
\medskip
{\bf P2.} This system of equations defines an \'etale, affine scheme over ${\rm Spec}(R/pR)$, cf. 3.6.8.1.2 a). We consider the ideal $I:=(z_1,...,z_m)$ of $R/pR$. As $R/pR$ is a strictly henselian ring, it is enough to look at this system of equations modulo $I$. Due to the upper triangular form of $\bar A(\Phi_R)$, this system of equations modulo $I$ can be ``separated":
\medskip
{\bf Fact.} {\it In the expression of $L_{ijl}$ modulo $I$ only those variables $x_{ijl^\prime}$ show up for which $l^\prime\le l$.}
\medskip
So this allows us to assume $m=1$: first we solve the subsystem of equations we get in variables $x_{ij1}$'s; we use any solution of it to get a new system  of equations in variables $x_{ij2}$'s (in other words, $x_{ij1}$'s become constants for this new system of equations involving just $x_{ij2}$'s); etc. (the induction applies; warning: all these are modulo $I$). From now on, as $m=1$, we drop the index $l$.  
So modulo $I$ we get a system $\Ms$ of $\dim_{W(k)}(M)^2$ equations in $\dim_{W(k)}(M)^2$ variables with coefficients in $k$.
\medskip
{\bf P3.} If $\bar A(\Phi_R)\in k$ is $0$, then the system $\Ms$ is of the form $x_{ij}=d_{ij}$, $(i,j)\in S(M)\times S(M)$, with $d_{ij}\in k$. Obviously it has a unique solution. So we can assume $\be:=\bar A(\Phi_R)\in k$ is non-zero (in fact we can assume $\be=1$). 
\smallskip
Modulo $I$, $g_R$ contributes only to the coefficients $d_{ij}$ (obtained by taking $c_{ij}$ modulo $I$). So 3.6.8.1.2 c) implies we can assume $g_R$ is the identity element of $GL(M)(R)$. 
\medskip
{\bf P4.} We present first the particular case when $(M,F^1,\vph)$ is cyclic diagonalizable. We consider a $W(k)$-basis $\{e_i|i\in S(M)\}$ of $M$. We can assume there is a subset $S(F^1)$ of $S(M)$ such that $\{e_i|i\in S(F^1)\}$ is a $W(k)$-basis of $F^1$ and there is a permutation $\pi$ of $S(M)$ such that $\vph(e_i)=pe_{\pi(i)}$, if $i\in S(F^1)$, and $\vph(e_i)=e_{\pi(i)}$, if $i\notin S(F^1)$ (cf. 2.2.1 d)); let $d(\pi)\in\NN$ be the order of $\pi$. We always take $\pi$ such that: if $\vph^v(e_i)=p^ue_i$ for any $i$ belonging to a $v$-cycle of $\pi$, with $v\in\NN$ and $u\in\{0,v\}$, then this cycle has just 1 element (i.e. $v=1$). The system $\Ms$ is of the form 
$$x_{\pi(i)\pi(j)}+d_{ij}=\be\al_{ij}x_{ij}^p,\leqno (26)$$
$i,j\in S(M)$, where $\al_{ij}$ is 1 or 0 depending on the fact that $(i,j)$ is or is not in $(S(M)\setminus S(F^1))\times S(F^1)$ and where $d_{ij}\in k$ (cf. $(E_1)$ and $(E_2)$ of 3.6.1.1.1 2)). So if $\pi(i)\neq i$ or if $\pi(j)\neq j$, or if $\pi(i)=i$ and $\pi(j)=j$ but $(i,j)\notin (S(M)\setminus S(F^1))\times S(F^1)$, then in the sequence of elements $\al_{ij}$, $\al_{\pi(i)\pi(j)}$,..., $\al_{\pi^{d(\pi)-1}(i)\pi^{d(\pi)-1}(j)}$, at least one is $0$; so such an $x_{ij}$ is uniquely determined. For any pair $(i,j)$, with $\pi(i)=i$ and $\pi(j)=j$ and with $i\notin S(F^1)$ and $j\in S(F^1)$, we have precisely $p$ (independent) choices for $x_{ij}$. The number of such pairs is $s(0)s(1)$. This handles the particular case.
\medskip
{\bf P5.} We come back to the general case. We write 
$$(M,\vph)=(M_1,\vph)\oplus (M_0,\vph)\oplus (M_{01},\vph)$$
such that all slopes of $(M_1,\vph)$, $(M_0,\vph)$ and $(M_{01},\vph)$  are 1, 0 and respectively belong to the interval $(0,1)$ (cf. the existence of the filtration $F^1$ of $M$). Let $B(\vph)\in M_{\dim_{W(k)}(M)}(k)$ be the matrix formed naturally (i.e. it acts naturally on $M/pM$) by a solution of the system $\Ms$. From the shape of $(E_1)$ and $(E_2)$ of 3.6.1.1.1 2) we get that $B(\vph)$ annihilates the image of $\vph(M)$ in $M/pM$. This implies that $B(\vph)$ annihilates $M_0/pM_0\oplus M_{(0,1)}/pM_{(0,1)}$: if $y_1\in M_{(0,1)}$ is non-zero mod $p$, then in the sequence of elements $y_u\in M_{(0,1)}\setminus pM_{(0,1)}$, $u\in\NN$, $u\neq 1$, defined inductively by the rule $\vph(y_{u})=p^{\vep_u}y_{u-1}$, with $\vep_u\in\{0,1\}$, there is a first element $y_{u_0}$ such that $\vep_{u_0}=0$; as in the particular case we (inductively) get that $B(\vph)$ annihilates the reductions mod $p$ of $y_{u_0-1}$,..., $y_1$. 
\smallskip
So we just need to look at the action of $B(\vph)$ on the image in $M_1/pM_1$ of elements $e\in M_1$ satisfying $\vph(e)=pe$. Writing $\nabla({\partial\over {\partial z_1}})(e)$ mod $I$ as $e_0+e_1+e_{01}$, with $e_0\in M_0/pM_0$, $e_1\in M_1/pM_1$ and $e_{01}\in M_{01}/pM_{01}$, we get (cf. ($E_2$) of 3.6.1.1.1 2)) that $e_0+e_1+e_{01}=\be\vph(e_0+e_1+e_{01})$. So $e_1$ and $e_{01}$ are both $0$. So $B(\vph)(e)$ belongs to $M_0/pM_0$ and so $B(\vph)$ belongs to a $s(0)s(1)$ dimensional $k$-vector subspace of ${\rm End}(M/pM)$. The same argument of 3.6.8 9) (based on 3.6.8.1 and pertaining to number of points) shows that we have at most $p^{s(0)s(1)}$ possibilities for $B(\vph)$. From this we deduce (via 3.6.18.2) that there are precisely $p^{s(0)s(1)}$ possibilities for $B(\vph)$ and so solutions of $\Ms$.  
\medskip
{\bf P6.} Another way (shorter and faster) to get this last sentence goes as follows. We can assume $\Phi_R(z_1)=z_1^p+pz_1$ (so $\be=1$); so working modulo $(z_1)$ instead of modulo $(p,z_1)$, we get similarly a matrix $\tilde B(\vph)\in M_{\dim_{W(k)}(M)}(W(k))$ lifting $B(\vph)$ and such that, when viewed as an endomorphism of $M$, we have (cf. $(E_1)$ and $(E_2)$ of 3.6.1.1.1 2) and the fact that $g_R$ is the identity element of $G(R)$)
$$p\vph(\tilde B(\vph))=\tilde B(\vph);\leqno (27)$$  
it is 3.6.8.1.2 c) and the constancy part of 3.6.8.9 which allow us to consider such a lift $\tilde B(\vph)$. (27) says that $\tilde B(\vph)$ is an element of ${\rm End}(M)$ fixed by $p\vph$. But the $\ZZ_p$-submodule $SL_{-1}$ (resp. the $W(k)$-submodule) of ${\rm End}(M)$ generated by such elements is of rank $s(0)s(1)$ (resp. is a direct summand of rank $s(0)s(1)$). 
\medskip
{\bf P7.} One can combine the above two ways: without using any lifting argument (of $B(\vph)$ to $\tilde B(\vph)$), we get (cf. (ENDFR) of 2.2.4 B): $B(\vph)$ is fixed by the reduction mod $p$ of the $\sg$-linear endomorphism $p\vph$ of ${\rm End}(M)$ and so it belongs to $SL_{-1}$ mod $p$.
\smallskip 
This ends the proof of B) and so of the Theorem.  
\medskip
{\bf P8. Exercise.} Starting from 3.6.18.1.2, give a second proof of B).    
\medskip
{\bf 3.6.18.4.1. Theorem (the integrability principle for $\Mm\Mf_{[0,1]}(R)$).} {\it All connections $\nabla$ proved to exist in 3.6.18.4 are integrable and nilpotent mod $p$.}
\medskip
{\bf Proof:} Let $R^{\rm al}:=W(k)[z_1,...,z_m]$. Here the upper right index ``al" stands for algebraization. The fact that a connection is integrable is expressed through some algebraic equations. In our case it is enough to check (regardless of which $n\in\NN\cup\{\infty\}$ we work with) that these equations are satisfied modulo the ideal $J_r:=(p,z_1,...,z_m)^r$ of $R$, for any $r\in\NN$. But these equations modulo $J_r$ depend only on the expressions of $\Phi$ and $\Phi_R$ modulo $J_{r+2}$.  This allows us to algebraize the things. From the Chinese Remainder Theorem we deduce the existence:
\medskip
-- of a Frobenius lift $\Phi_{R^{\rm al}}$ of $R^{\rm al}$ which modulo its ideal $J_{r+2}^{\rm al}:=(p,z_1,...,z_m)^{r+2}$ is the same as $\Phi_R$ modulo $J_{r+2}$ and which modulo its ideal $I_1^{r+2}$, with $I_1:=(p,z_1-1,...,z_m-1)$, is isomorphic to the one coming from a generic situation as in 3.6.18.2;
\smallskip
-- of an element $g_{R^{\rm al}}\in GL(M)(R^{\rm al})$ such that modulo $J_{r+2}^{\rm al}$ it is $g_R$ modulo $J_{r+2}$, and modulo $I_1^{r+2}$ is an element $g_1\in GL(M)(W_n(k))\subset GL(M)(R^{\rm al}/I_1^{r+2})$ for which the Frobenius endomorphism of $M\otimes_{W(k)} R^{\rm al}/I_1^{r+2}$ defined by $g_1(\vph\otimes 1)$ is isomorphic to one coming from a generic situation.
\medskip
We assume $n\in\NN$. Let $\nabla_n$ be a connection on $M\otimes_{W(k)} R/p^nR$ which makes ${\got C}/p^n{\got C}$ potentially to be viewed as an object of $\Mm\Mf_{[0,1]}(R)$. As in 3.6.1.3 (and its proof) we construct for the $p$-divisible object 
$$
{\got C}^{\rm al}:=(M\otimes_{W(k)} R^{\rm al},F^1\otimes_{W(k)} R^{\rm al},g_{R^{\rm al}}(\vph\otimes 1))
$$ 
of $\Mm\Mf_{[0,1]}(R^{\rm al})$ a moduli affine ${\rm Spec}({R^{\rm al}}^\wedge)$-scheme 
$$M_n({\got C}^{\rm al})={\rm Spec}(R_n)$$ 
of connections which make the extension of ${\got C}^{\rm al}/p^n{\got C}^{\rm al}$ to $*$ potentially to be viewed as an object of $\Mm\Mf_{[0,1]}^{\nabla}(*)$; the Frobenius lift of $R^{\rm al}$ is $\Phi_{R^{\rm al}}$ and $*$ stands for a $p$-adically complete, formally \'etale $R^{\rm al}$-algebra. Mod $p$ it defines an \'etale, affine scheme $\Mr(n):={\rm Spec}(R_n/pR_n)$ over $\Mr(0):={\rm Spec}(R^{\rm al}/pR^{\rm al})$. Moreover the fibres of $\Mr(n)$ over $k$-valued points of $\Mr(0)$ have a number of points less or equal to the number $n_1$ of $\bar k$-valued points of the fibre $\Mf_1$ above the $k$-valued point of $\Mr(0)$ defined by $z_i=1$, $i\in S(m)$, cf. 3.6.18.4 B). So the \'etale, affine $k$-morphism $r_n:\Mr(n)\to\Mr(0)$ is an \'etale cover above the maximal open subscheme of $\Mr(0)$ with the property that the fibres of $r_n$ above geometric points of it have precisely $n_1$ points. We conclude:
\medskip
{\bf Fact 1.} {\it Any connected component of $\Mr(n)$ contains points of $\Mf_1$.}
\medskip  
From this and from 3.6.18.2 we get:
\medskip
{\bf Fact 2.} {\it $M_n({\got C}^{\rm al})$ is a moduli scheme of integrable connections, nilpotent mod $p$.}
\medskip
The last thing we need is:
\medskip
{\bf Fact 3.} {\it $\nabla_n$ modulo the ideal $(p,z_1,...,z_m)^r$ of $R/p^nR$, is equal to a connection obtained from the universal one on $M\otimes_{W(k)} R_n/p^nR_n$, through a formally \'etale, affine $W(k)$-morphism ${\rm Spec}(R/p^nR)\to M_n({\got C}^{\rm al})$.}
\medskip
The easiest way to argue Fact 3 is to remark that ${\rm Spec}(R_n/(J_r^{\rm al}+p^nR^{\rm al})R_n)$ is the moduli scheme of maps (thought as ``pseudo-connections")
$$\nabla_{r,n}:M\otimes_{W(k)} R^{\rm al}/J_r^{\rm al}+p^nR^{\rm al}\to M\otimes_{W(k)} FREE_{r,n},\leqno (PSCONN)$$
with $FREE_{r,n}$ as the free $R^{\rm al}/J_r^{\rm al}+p^nR^{\rm al}$-module having $dz_1$,..., $dz_m$ as a basis, satisfying the equations obtained from the equations satisfied by the universal connection on $M\otimes_{W(k)} R_n/p^nR_n$ by taking them modulo the ideal $J_r^{\rm al}+p^nR^{\rm al}$ of $R^{\rm al}$. This is a consequence of the fact that property (DIV) of 3.6.8 5) still holds modulo the ideal $J_{r+2}^{\rm al}+p^nR^{\rm al}$ of $R^{\rm al}$.  
\smallskip
From Facts 2 and 3 we get (we recall that $r\in\NN$ was arbitrary): $\nabla_n$ is integrable and nilpotent mod $p$. As $n\in\NN$ is arbitrary the same holds for the case of $p$-divisible objects (i.e. for $n=\infty$). This ends the proof of the Theorem. 
\medskip
{\bf 3.6.18.4.1.1. Corollary.} {\it Let $n\in\NN$. We consider a connection $\nabla$ on ${\got C}$. For $l\in S(1,m)$, we consider the $W(k)$-linear endomorphism $D_l:=\nabla({\partial\over {\partial z_l}})$ of $M\otimes_{W(k)} R$. We have:
\medskip
{\bf a)} $D_lD_s=D_sD_l$, $\forall s,l\in S(1,m)$;
\smallskip
{\bf b)} $D_l^{2p-1}$ mod $p$ is $0$, $\forall l\in S(1,m)$.}
\medskip
{\bf Proof:} a) is just a restatement of 3.6.18.4.1. We now argue b). From the proof of 3.6.18.4.1 we get that we can assume we are in the context of a generic situation. So everything results from (N0) and (N1) of the proof of 3.6.18.4.2. 
\medskip
b) represents a slight improvement to 3.6.1.1.2 (NIL2).
\medskip
{\bf 3.6.18.4.2. Corollary (the moduli principle).} {\it For any regular, formally smooth $W(k)$-scheme or $p$-adic formal scheme $S$ over ${\rm Spec}(W(k))$ equipped with a Frobenius lift and having a connected special fibre, every object (resp. every $p$-divisible object) $\got C$ of $\Mm\Mf_{[0,1]}(S)$ defines a moduli $p$-adic formal scheme $M(\got C)$ over $S$ of connections (on the underlying sheaves of modules of the pull back of ${\got C}$ to formally \'etale, $p$-adic formal schemes over $S$) which make the extension of $\got C$ to such a $p$-adic formal scheme $S_1$ to be viewed as an object (resp. as a $p$-divisible object) of $\Mm\Mf^{\nabla}_{[0,1]}(S_1)$. We have:
\medskip
{\bf a)} $\forall m\in\NN$, the $S_{W_m(k)}$-scheme $M({\got C})_{W_m(k)}$ is \'etale and affine (resp. is $\NN$-pro-\'etale and affine);
\smallskip
{\bf b)} Locally in the Zariski topology of $S_k$, the special fibre of $M({\got C})_{W_m(k)}$ is obtained through a finite sequence (resp. infinite sequence indexed by elements of $\NN$) of systems of equations of first type and defined using linear homogeneous forms having coefficients in (the ring of global sections of) $S_k$;
\smallskip
{\bf c)} The fibres of $M({\got C})_k$ over points with values in fields of $S_k$ are non-empty. In the case of an object of $\Mm\Mf_{[0,1]}(S)$, these fibres are \'etale schemes over fields defined by algebras of dimension (as vector spaces over these fields) a non-negative, integral power of $p$.}
\medskip
{\bf Proof:}
c) follows from b), cf. 3.6.8.1 (applied with $s=0$ and $l=1$). Due to the universal property (implicit in the reference of $M({\got C})$ as a moduli $p$-adic formal scheme), we can work locally in the Zariski topology of $S_k$. So we can assume that $S={\rm Spec}(\tilde R)$ is an affine scheme, that $\tilde R=\tilde R^\wedge$, that $\Om_{S_k/k}$ is a free $\Mo_{S_k}$-sheaf of finite rank, and that the underlying $\tilde R$-module $\tilde M$ of ${\got C}$ has the DC property, i.e. it is isomorphic to a direct sum 
$$\sum_{i=1}^q (\tilde R/p^i\tilde R)^{m_i},$$ 
with $q\in\NN$ and $m_i\in\NN\cup\{0\}$, $i=\overline{1,q}$, cf. [Fa1, 2.1 ii)] (resp. and that $\tilde M$ is free). It is the connectedness of $S_k$ which allows us to consider such a $q\in\NN$ if ${\got C}$ is an object.
\smallskip
 Working with ``potentially to be viewed" instead of ``to be viewed", a) and b) are a consequence of 3.6.8.4 1) to 3) and of 3.6.8.1.2 a), cf. also the below extra features for the context of objects. But (the proof of) 3.6.18.4.1 allows us to drop the word ``potentially".
\smallskip
When we pass from $p$-divisible objects to objects, there are three extra features (differences) showing up, which are worth being pointed out in detail. In what follows we assume ${\got C}=(\tilde M,F^1(\tilde M),\vph_0,\vph_1)$ is an object. It can happen that there are two distinct elements $i,j\in S(1,q)$ such that $m_i$ and $m_j$ are both non-zero. In such a situation, ${\rm End}(\tilde M)$ is not a free $\tilde R/p^l\tilde R$-module for some $l\in S(0,q)$ and so it is much more convenient (in some sense compulsory) to work inductively on $l\in S(1,q)$ with ${\rm End}(\tilde M)/p^l{\rm End}(\tilde M)$ instead of $\tilde M/p^l\tilde M$. This is the first significant difference in comparison with the proof of 3.6.1.3 and so the first extra feature.
\smallskip
We can assume $F^1(\tilde M)$ has also the DC property (see proof of 2.2.1.1 6)). We consider a set $\Mb:=\{e_1,...,e_r\}$, with $r:=\sum_{i=1}^q m_i$, of elements of $\tilde M$ such that $\tilde M$ is the direct sum of its cyclic $\tilde R$-submodules generated by them (the proof of 2.2.1.1 6) implies that the annihilator of any element of $\Mb$ is generated by a power of $p$). We can assume $F^1(\tilde M)$ is generated by the elements of a subset of $\Mb$. We refer to such a set $\Mb$ as an $\tilde R$-basis of ${\got C}$. Let $\Me:=\{f_1,...,f_u\}$ be the set of elements of ${\rm End}({\tilde M})$ such that:
\medskip
-- ${\rm End}({\tilde M})$ is the direct sum of its cyclic $\tilde R$-submodules generated by them;
\smallskip
-- it is naturally defined by elements of $\Mb$, i.e. for each $s\in S(1,u)$ there are $s_1$, $s_2\in S(1,r)$ such that $f_u(e_i)=0$ if $i\neq s_1$ and $f_u(e_{s_1})=p^{n(u)}e_{s_2}$, with $n(u)\in S(0,q)$.
\medskip
Also we can assume there are elements $z_1$,..., $z_v$ of $\tilde R$ such that $\{dz_1,..., dz_v\}$ is an $\tilde R$-basis of the $p$-adic completion of $\Om_{\tilde R/W(k)}$. The existence of a connection $\nabla$ on $\tilde M$ making ${\got C}$ to be viewed as an object $\Mm\Ms_{[0,1]}^\nabla(S)$ can be codified (i.e. restated) in finding solutions of $v$ equations involving elements of ${\rm End}(\tilde M)$. So writing
$$\nabla=\dl_0+\sum_{j^\prime=1}^v\sum_{j=1}^r x_{jj^\prime}f_jdz_{j^\prime},\leqno (28)$$
with $\dl_0$ an arbitrary but fixed connection on $\tilde M$, we need to find solutions of a system of equations $SE$ in the $rv$ variables $x_{jj^\prime}$'s. When viewed mod $p$ it is of first type (this, as in 3.6.8 5), is a consequence of the fact that $\tilde M$ is $\tilde R$-generated by $\vph(\tilde M)$ and by $\vph_1(F^1(\tilde M))$), and so 3.6.8.1.2 a) applies as in the proof of 3.6.1.3 (for instance, see 3.6.8 7)). Working by induction on $l\in S(1,q)$ modulo $p^l$ and following 3.6.8 10) to 12), we get a) and the part of b) which does not refer to linear homogeneous forms. For the integrability part, we need to remark that based on 3.6.18.2.2 and the above part of the proof we can perform entirely the algebraization process of the proof of 3.6.18.4.1 in the context of arbitrary objects of $\Mm\Mf_{[0,1]}(R)$ (see also 3.6.18.5.2 below for a second possibility of arguing this, via 3.6.18.4.1 and the Fact of 2.2.1.1 6)).
\smallskip
We now point out the second significant difference in comparison with the proof of 3.6.1.3 and so the second extra feature. By passing from things mod $p^l$ to things mod $p^{l+1}$, we get a system of equations in at most the same number of variables (the precise number of variables is $v$ times the rank of $p^l{\rm End}(\tilde M)/p^{l+1}{\rm End}(\tilde M)$ as an $\tilde R/p\tilde R$-module), and so potentially of a completely different form of what we got working mod $p^l$, $l\in S(1,q-1)$. But the situation is the same: we still get systems of equations of first type. 
\smallskip
Also it is worth pointing out that the $\tilde R/p\tilde R$-module $p^l{\rm End}(\tilde M)/p^{l+1}{\rm End}(\tilde M)$ is naturally a submodule of ${\rm End}(\tilde M)/p{\rm End}(\tilde M)$, and so the new systems of equations of first type are related to specific subobjects of $End({\got C})/pEnd({\got C})=End({\got C}/p{\got C})$; so the ``completely different forms" are not ``too far" from the initial form of $SE$. This is the third new feature. It can be reformulated as: the linear homogeneous forms of these systems of equations depend only on $End({\got C}/p{\got C})$ and on its (warning: non-necessarily Lie) subobjects. 
\smallskip
Conclusion: (we can assume) these systems of equations of first type have coefficients in $\tilde R/p\tilde R$. This proves a) and b) and ends the proof. 
\medskip
{\bf Examples.} As examples of such moduli $p$-adic formal schemes (resp. of special fibres of them) for connections, we have the $p$-adic completions (viewed as formal schemes) of the $W(k)$-schemes ${\rm Spec}(Q_{j,n}^\wedge)$ (resp. the $k$-schemes $\Ms^n$), $n\in\NN$, of 3.6.1.3 (resp. of 3.6.8). 
\medskip
As in Fact 2 of the proof of 3.6.18.4.1 we get:
\medskip
{\bf 3.6.18.4.3. Corollary (the touching property).} {\it We assume we are in the context of 3.6.18.4.2, with ${\got C}$ an object of $\Mm\Mf_{[0,1]}(S)$. Let $p^m$, with $m\in\NN\cup\{0\}$, be the maximal number of points ${\rm Spec}(k_1)\to M({\got C})_k$ lifting the same $k_1$-valued point of $S_k$, with $k_1$ an algebraically closed field containing $k$. Then any connected component of $M({\got C})_k$ intersects all fibres of $M({\got C})_k$ over geometric points of $S_k$ which are defined by \'etale algebras of dimension, as vector spaces over residue fields of $S_k$, equal to $p^m$.}
\medskip
{\bf 3.6.18.4.4. Remarks.} {\bf 1)} 3.6.18.4 remains true for the case when $R/pR$ is a local henselian ring and $R$ is a regular, formally smooth $W(k)$-algebra, having $k=\bar k$ as its residue field.
\smallskip
{\bf 2)} From now on, in connection to the categories $\Mm\Mf_{[0,1]}(*)$, we do not use (cf. the proof of 3.6.18.4.2) anymore the word potentially (see 3.6.1.1.1).
\smallskip
{\bf 3)} The moduli schemes of 3.6.18.4.2 do depend on the choice of $\Phi_S$ (for instance, cf. 3.6.18.4 B)).
\medskip
From 3.6.18.4.2 b) and 3.6.8.1.2 b), as we are in a regular context, we get:
\medskip
{\bf 3.6.18.4.5. Corollary (the surjectivity principle).} {\it Under the assumptions of 3.6.18.4.2, any connected component of $M({\got C})_k$, with $\got C$ a $p$-divisible object, maps into an open, dense subscheme of $S_k$.}
\medskip
{\bf 3.6.18.4.6. Extra terminology.} 
What follows in this section is an abstract algebraic digression supplementing the first paragraphs of 3.6.8.9. In that paragraphs, we defined abstractly systems of equations of first and third type, as well as, in connection to 3.6.8.5, we defined systems of equations of (adjusted) additive second type, while working with an additive type Frobenius lift of a regular, formally smooth $W(k)$-algebra $R$. 
\smallskip
Similarly, if the Frobenius lift involved is of multiplicative (resp. essentially multiplicative or additive) type (in some maximal point of ${\rm Spec}(R)$, see 3.6.18.0.1.1), we define systems of equations of multiplicative (resp. essentially multiplicative or additive) second type; in 3.6.18.4 P4 (see (26) of it) we came across such systems of equations but we preferred to view them as of first type. The main difference from what we got in 3.6.8.5 1) is: the coefficients of $x_{isl}^p$'s we get are not any more some $p-1$ powers of non-invertible elements (see 3.6.8.5 (11)) but they are of a more general nature; in particular they can be invertible. Also in 3.6.14.1 J we introduced systems of equations of fourth type; in what follows they are not used. 
\smallskip
We index different parts using capital letters. 
\smallskip
{\bf A.} Let $\bar R$ be an arbitrary $\FF_p$-algebra. Let $n\in\NN$. 
A system of equations of the form
$$x_i=L_i(x_1^p,x_2^p,...,x_n^p)+c_i,\leqno (*)
$$ 
$i=\overline{1,n}$, where $L_i$ is a linear homogeneous form with coefficients  in $\bar R$ and $c_i\in\bar R$, $\forall i\in S(1,n)$, is said to be of first type with coefficients in $\bar R$ (cf. the terminology of 3.6.8.9); we also call it a quasi Artin--Schreier system of equations in $n$ variables with coefficients in $\bar R$. It defines naturally an \'etale, affine ${\rm Spec}(\bar R)$-scheme (cf. 3.6.8.1.2 a)). If this system of equations defines a finite, \'etale morphism (at the level of schemes), then we omit the word quasi. The $n\times n$ matrix obtained naturally (as in 3.6.8.1) from the coefficients of the linear forms $L_i$, $i=\overline{1,n}$, is called the matrix of the system $(*)$ of equations. 
\smallskip
We came across such quasi Artin--Schreier systems of equations in the proofs of 3.6.1.3 and 3.6.18.4.2. For future references, we define as well (abstractly) systems of equations of second type: a system of equations of the form
$$a_ix_i^p=L_i(x_1,x_2,...,x_n)+c_i,\leqno (29)$$
$i=\overline{1,m}$, and 
$$0=L_i(x_1,x_2,...,x_n)+c_i,\leqno (30)$$
$i=\overline{m+1,n}$, where $m\in S(1,n)$, and where $L_i$'s and $c_i$'s are as above, is said to be of second type with coefficients in $\bar R$, if at least one of the elements $a_1$,..., $a_m$ of $\bar R$ is non-zero. However, as we do not use them again in the rest of this paper or in \S5-14, we do not stop to define abstractly systems of equations of (adjusted) additive (or multiplicative, etc.) second type.
\smallskip
{\bf B. Artin--Schreier fundamental groups.} An \'etale cover $q_1:{\rm Spec}(\bar R_1)\to {\rm Spec}(\bar R)$ is called crystalline elementary $n$-admissible if it is defined by an open closed subscheme of an \'etale cover of ${\rm Spec}(\bar R)$ obtained using an Artin--Schreier system of equations in at most $n$ variables; here the ``at most" part is motivated by the second extra feature of the proof of 3.6.18.4.2. We say $q_1$ is crystalline $n$-admissible if it is obtained through (i.e. it is a composite of) a finite sequence of consecutive crystalline elementary $n$-admissible \'etale covers. 
\smallskip
We assume ${\rm Spec}(\bar R)$ is connected. For all that follows, we fix a geometric point of ${\rm Spec}(\bar R)$ in the usual way related to fundamental groups; not to complicate the notations and as below it is convenient to look at fundamental groups as groups of automorphism of Galois pro-\'etale covers, we never mention it explicitly. The automorphism group of the Galois pro-\'etale cover of ${\rm Spec}(\bar R)$ generated by connected components of crystalline $n$-admissible \'etale covers of ${\rm Spec}(\bar R)$, is called the level $n$ Artin--Schreier fundamental group of $\bar R$ or of ${\rm Spec}(\bar R)$ and it is denoted by $\Pi_n^{AS}(\bar R)$. We have natural epimorphisms 
$$
\Pi_{n+1}^{AS}(\bar R)\twoheadrightarrow\Pi_n^{AS}(\bar R).  
$$
Variant: above we consider only finite sequences of length at most $m$, where $m\in\NN$ is fixed; we obtain pro-finite Galois groups $\Pi_{n,m}^{AS}(\bar R)$ and natural epimorphisms 
$$
\Pi_{n_1,m_1}^{AS}(\bar R)\twoheadrightarrow\Pi_{n,m}^{AS}(\bar R),
$$ 
with $n_1,m_1\in\NN$, $n_1\ge n$ and $m_1\ge m$. We denote by $\Pi_{\infty}^{AS}(\bar R)$ the projective limit of $\Pi_n^{AS}(\bar R)$. It is a quotient of the first fundamental group of ${\rm Spec}(\bar R)$. We call it the Artin--Schreier fundamental group of $\bar R$ or of ${\rm Spec}(\bar R)$. We have:
\medskip
{\bf C. Lemma (the fundamental lemma of the $\pi_1$-theory in positive characteristic).} {\it If $k_1$ is a field (not necessarily perfect) containing $\FF_p$, then $\Pi_{\infty}^{AS}(k_1)$ is exactly the Galois group of $k_1$.}
\medskip 
{\bf Proof:}
Let $k_2$ be a finite Galois extension of $k_1$. Let $n_0:=[k_2:k_1]$ and let $E=\{e_1,...,e_{n_0}\}$ be an arbitrary $k_1$-basis of $k_2$. An arbitrary element $x\in k_2$ can be written uniquely as a sum $x=\sum_{i=1}^{n_0} x_ie_i$, with $x_i\in k_1$, $i=\overline{1,n_0}$. We consider the equation 
$$x^p=x.\leqno (31)$$ 
It gets translated into an Artin--Schreier system of equations of the form 
$$x_i=L_i(x_1^p,...,x_{n_0}^p),\leqno (**)$$ 
$i=\overline{1,n_0}$, with $L_i$ as linear homogeneous forms. The square matrix $A_E$ formed by the coefficients of these linear forms is invertible (this can be seen by extension to $k_1^{\rm perf}$). So 3.6.8.1.2 a) applies. So $k_2$ splits over the Galois extension $k_3$ of $k_1$ generated by fields of whose spectra are connected components of the \'etale $k_1$-scheme defined naturally by this Artin--Schreier system of equations (argument: the number of solutions of the equation $x^p=x$ in $k_2\otimes_{k_1} k_3$ is --cf. 3.4.8.1.0-- on one hand $p^n$ and on the other hand it is $p$ to the power the number of connected components of ${\rm Spec}(k_2\otimes_{k_1} k_3)$). This ends the proof. 
\medskip
{\bf D. The global context.} In the above Lemma the condition $k_1$ is a field can be weaken: we just need that any projective $k_1$-module of finite rank is free. In particular it holds for local rings, for PID's and for polynomial rings over $k$. Warning: we do not need to assume we are in some normal context; if $s$ is a non-zero solution of (31) in any $\FF_p$-algebra, then $s^{p-1}$ is a projector and so it still ``isolates" some connected components. So to understand the first fundamental groups of arbitrary connected $\FF_p$-schemes we just need to understand how to ``glue" Artin--Schreier systems of equations. This goes as follows.
\smallskip
The proof of C above suggests looking only at Artin--Schreier systems of equations which are defined by homogeneous linear forms whose coefficients are producing an invertible matrix. Referring to the proof of C, the change of the $k_1$-basis $E$, results in a replacement of $A_E$ by $BA_E(B^{[p]})^{-1}$, with $B\in GL_{n_0}(k_1)$ (warning: this is different from 3.6.18.0). In general, two invertible matrices $A_1$ and $A_2$ of rank $n$ with coefficients in $\bar R$, are called $F$-equivalent (over $\bar R$), if there is an invertible matrix $B$ with coefficients in $\bar R$ such that $A_2=BA_1(B^{[p]})^{-1}$. So the determinants of $A_1$ and $A_2$ differ by a $p-1$ power of an invertible element of $\bar R$.
\medskip
{\bf Definition.} Let $X_p$ be arbitrary separated, connected $\FF_p$-scheme. A Galois--Artin--Schreier (or an Artin--Schreier homogeneous invertible) system of equations in $n$ variables over $X_p$ is defined by: 
\medskip
-- an open cover $\Mo\Mc$ of $X_p$ by affine schemes $X_p^i={\rm Spec}(R^i)$, $i\in\Mi$,\smallskip
-- an invertible $n\times n$ matrix $A^i$ with coefficients in $R^i$, $\forall i\in\Mi$,
\smallskip
-- an invertible $n\times n$ matrix $B_{ij}$ with coefficients in $R^{ij}$, $\forall i,j\in\Mi$, where $X_p^i\cap X_p^j={\rm Spec}(R^{ij})$,
\medskip\noindent
such that the following two things hold
\medskip
a) $\forall i,j\in\Mi$, over $X_p^i\cap X_j^p$ we have $B_{ij}A_j(B_{ij}^{[p]})^{-1}=A_i$, and
\smallskip
b) $\forall i,j,l\in\Mi$, over $X_p^i\cap X_p^j\cap X_p^l$ we have $B_{ij}B_{jl}=B_{il}$. 
\medskip
The equivalence of two such Galois--Artin--Schreier systems of equations in $n$ variables over $X_p$ is defined in the standard manner (of passing to refinements of $\Mo\Mc$, etc.). Following the pattern of (**) above, we can use them to define an \'etale cover of $X_p$ of degree $p^n$: above $X_p^i$ it is the \'etale cover defined by the Artin--Schreier system of equations in $n$ variables $x_1^i$,.., $x_n^i$ whose linear forms are homogeneous and whose matrix is $A_i$ itself, $\forall i\in\Mi$; the identities a) and b) allow us to glue logically such \'etale covers (i.e. $B_{ij}$ are the matrices of the linear combination transitions from $x_s^i$'s to $x_s^j$'s, $s=\overline{1,n}$, $\forall i,j\in\Mi$: if $x^i$ is the column vector having $x_s^i$'s as its entries, $s=\overline{1,n}$, then we have $x^i=B_{ij}x^j$). 
\smallskip
Following the pattern of B), we use such \'etale covers to define a pro-finite Galois group $\Pi^{GAS}_{n,m}(X_p)$, with $m\in\NN$ having the same role as above, as well as their $\NN$-projective limit $\Pi^{GAS}_{n}(X_p)$. We refer to $\Pi^{GAS}_{n}(X_p)$ as the level $n$ Galois--Artin--Schreier fundamental group of $X_p$. Moreover, by identifying an $n\times n$ matrix with entries $(a_{ij})_{i,j\in S(1,n)}$ with a $n+1\times n+1$ matrix whose entries $(b_{ij})_{i,j\in S(1,n+1)}$ are defined by $b_{ij}=a_{ij}$ and $b_{n+1,i}=b_{j,n+1}=0$ if $1\le i,j\le n$, and by $b_{n+1,n+1}=1$, we get a natural epimorphism $\Pi^{GAS}_{n+1}(X_p)\twoheadrightarrow\Pi^{GAS}_{n}(X_p)$, and so we define the $\NN$-projective limit 
$$\Pi^{GAS}_{\infty}(X_p)$$
 of these groups $\Pi^{GAS}_{n}(X_p)$, $n\in\NN$. 
Following the proof of C, we get:
\medskip
{\bf Lemma (the fundamental lemma of the $\pi_1$-theory in positive characteristic: the global form).} {\it The fundamental group of $X_p$ is $\Pi^{GAS}_{\infty}(X_p)$.}
\medskip
If $X_p$ is moreover integral, regular and noetherian, the above Lemma represents (in our opinion) a practical way of computing the fundamental group of $X_p$.
As a particular case of this Lemma we get the following criterion:
\medskip 
{\bf Corollary (a criterion of simply connectivity).} {\it The fundamental group of $X_p$ is trivial iff $\forall n\in\NN$, any Galois--Artin--Schreier system of equations in $n$ variables over $X_p$ is equivalent to the trivial one, defined by $n\times n$ matrices which are all identity.}
\medskip
{\bf E. Remarks.} {\bf 1)} The proof of C above implies ${\rm Gal}(k_2/k_1)$ is naturally a quotient of $\Pi^{AS}_{n_0,1}(k_1)$. It is very much desirable to replace here the lower right index ${n_0,1}$ with another one ${n,1}$, with $n$ as small as possible. We do not know how to get in general better upper bounds for such an $n$. If $n_0\ge 2$, by using 3.6.8.1.5 and the fact that $1$, $x$, $x^p$,..., $x^{p^{n_0-1}}$ are linearly dependent over $k_1$, for any $x\in k_2$, we can take $n=n_0-1$; in this way we can obtain a second proof of C (but not of D).
\smallskip
{\bf 2)} Let $l\in\NN$. Using $x_i^{p^l}$'s instead of $x_i^p$'s in (*), we define similarly level $(n,m,l)$ Artin--Schreier fundamental groups $\Pi_{n,m,l}^{AS}(\bar R)$; if $l=1$, we always drop mentioning it as an index. Obviously (see 3.6.8.1.5) we have a natural epimorphism:
$$\Pi^{AS}_{nl,m}(\bar R)\twoheadrightarrow\Pi^{AS}_{n,m,l}(\bar R).$$
Similarly, in the context of D above, we define the level $(n,m,l)$ Galois--Artin--Schreier fundamental group $\Pi_{n,m,l}^{GAS}(X_p)$ of $X_p$ (so we have to speak about the $F^l$-equivalence relation between invertible matrices with coefficients in a given $\FF_p$-algebra, etc.).
\smallskip
{\bf 3)} There are $n(n+1)$ coefficients involved in the definition of a quasi Artin--Schreier system of equations in $n$ variables; they are entirely independent.  One can learn a lot by looking at quasi Artin--Schreier systems of equations in a small number of variables, like 2 or 3. For $n=2$ we have 6 coefficients, while for $n=3$ we have 12.
\smallskip
{\bf 4)} If a quasi Artin--Schreier system of equations with coefficients in $\bar R$ defines a Galois cover of an open, dense subscheme of ${\rm Spec}(\bar R)$, then the \'etale scheme over ${\rm Spec}(\bar R)$ naturally associated to it is an open subscheme of a Galois cover of ${\rm Spec}(\bar R)$. This is an immediate consequence of 3.6.8.1.2 a). 
\smallskip
{\bf 5)} The Lemma of D (as well as the ``$\rho$" part of 3.6.18.7.2 below) can be reformulated in terms of truncations mod $p$ of unit $F$-crystals. However, B (and so implicitly F below) goes beyond the reach of the language of such truncations: there are plenty of Artin--Schreier systems of equations whose matrices are not invertible. Here are two quick examples (with $k$ as usual). 
\medskip
{\bf Example 1.} Let $R=k[t_1,t_2]/(t_1t_2)^n$, with $n\in\NN\setminus\{1\}$; we take the system of equations (in $z_1$ and $z_2$) defined by $z_1+c_1=t_1t_2z_1^p+t_1z_2^p$ and $z_2+c_2=t_2z_1^p+z_2^p$. 
\medskip
{\bf Example 2.} Let $U:=k[x_1,x_2,x_3,x_4]/(x_1x_3-x_2x_4)$; we take the system of equations $z_1=x_1z_1^p+x_2z_2^p$ and $z_2=x_3z_1^p+x_4z_2^p$ and we consider the maximal affine, open subscheme ${\rm Spec}(\bar R)$ of ${\rm Spec}(U)$ above which this system of equations defines an \'etale cover $q:{\rm Spec}(\bar R_1)\to {\rm Spec}(\bar R)$, cf. 3.6.8.1.3-4. We have: 
\medskip
{\bf 1)} The complement of ${\rm Spec}(\bar R)$ in ${\rm Spec}(U)$ is irreducible and is not given by one single equation.
\medskip
{\bf 2)} The \'etale cover $q$ is of degree $p$ and is not definable by an Artin--Schreier system of equations in 1 variable.
\medskip
1) is left as an exercise (cf. [Ha, Exc. 6.5 of p. 147]). To argue 2), we assume $q$ is given by the equation $y=ay^p+b$, with $a$, $b\in\bar R$. $a$ is invertible in $\bar R$ and so, based on 1) and on the fact that $U$ is normal, it is an invertible element of $U$. As any such element is an element of $k$, by passing to $\bar k$, we can assume $a=1$. So $q$ is a Galois extension; but it is easy to see that this is not so. Contradiction. So 2) follows.
\medskip
Example 2 motivates why in all that follows we do not mention at all such truncations. However, the terminology of crystalline (elementary) $n$-admissible \'etale covers used in B is motivated by the way we came across 3.6.8.1 and by such a possible reformulation.
\medskip
{\bf F. Definitions.} Let $\Pi_0^{AS}(\bar R)$ (resp. $\Pi_{\infty+1}^{AS}(\bar R)$) be the trivial group (resp. the fundamental group of ${\rm Spec}(\bar R)$). The Artin--Schreier dimension of $\bar R$, denoted $d^{AS}(\bar R)$, is the smallest number $q\in\{0\}\cup\NN\cup\{\infty\}\cup\{\infty+1\}$ such that $\Pi_q^{AS}(\bar R)$ is the whole fundamental group of ${\rm Spec}(\bar R)$. 
\smallskip
Let $q\in\NN\cup\{\infty\}$. If the group $\Pi_q^{AS}(\bar R)$ is trivial we say ${\rm Spec}(\bar R)$ is weakly $q$-simply connected. If moreover $\Pi_q^{GAS}({\rm Spec}(R))$ itself is trivial, then we drop the word weakly. Similarly, referring to D with $X_p$ non-affine, if $\Pi_q^{GAS}(X_k)$ is trivial we say $X_p$ is $q$-simply connected. 
\medskip
{\bf G. Homework.} Compute the level $n$ Artin--Schreier fundamental groups and the Artin--Schreier dimensions of algebraic field extensions of $\FF_p$.
\medskip 
{\bf H. Artin--Schreier systems of equations in mixed characteristic.} We do not stop to present here a full theory: we are just interested to point out the main idea. Let $O$ be a DVR of mixed characteristic $(0,p)$ and let $k_1$ be its residue field. Let $R_O$ be an integral, regular, faithfully flat, formally smooth $O$-algebra. We assume that the smallest open subscheme of ${\rm Spec}(R_O)$ containing ${\rm Spec}(R_O\otimes_O k_1)$ is ${\rm Spec}(R_O)$ itself. We situate ourselves in the context of A above, with $c_i\in R_O$ and with $L_i$ having coefficients in $R_O$, $\forall i\in S(1,n)$. We assume the matrix of the coefficients of $L_i$'s is invertible. Then, as usual, we get an \'etale cover of ${\rm Spec}(R_O)$ (this is a consequence of the assumption on ${\rm Spec}(R_O)$, via the classical purity theorem of [SGA1, p. 275]). Imitating B we define the level $n$ Artin--Schreier fundamental group $\Pi^{AS}_n(R_O)$ of $R_O$, as well as their projective limit, the Artin--Schreier fundamental group $\Pi^{AS}(R_O)$ of $R_O$. Variant: we can define the level $(n,m,l)$ Artin--Schreier fundamental group $\Pi_{n,m,l}^{AS}(R_O)$ of $R_O$. We do believe that, in this way, one can learn a lot about Galois groups of number fields; this is our first  (see [Va10] for the second) approach towards the understanding of finite quotients of such Galois groups.
\medskip
{\bf I. Artin--Schreier systems of equations in unequal positive characteristic.} The last sentence of H above suggests treating the situation where we deal with a system (*) of equations as in A and whose matrix is invertible but over affine schemes which are of characteristic a prime $l\neq p$. The first thing to be done is to get criteria (using the classical purity theorem) when such systems of equations (defining finite, flat schemes) do define \'etale covers. 
\medskip
 We hope to come back in a future paper with the main (expected or known to be true) properties of these Artin--Schreier fundamental groups and dimensions. We end by mentioning that we do believe that there is a deep connection between I, [Bu2] and the theory of ordinary reductions of a proper, smooth scheme over a number field.  
\medskip
We come back to the starting context of 3.6.18.
\medskip
{\bf 3.6.18.5. Theorem (the inducing property).} {\it The answer to the question {\bf Q} is yes if any one of the following conditions is satisfied:
\medskip
\item{a)} $k=\bar k$; 
\smallskip
\item{b)} $\Phi_R$ is of essentially additive type; 
\smallskip
\item{c)} the Newton polygon of $(M,\vph)$ does not have either the slope 0 or the slope 1.}
\medskip
{\bf Proof:} If a) holds, this is a consequence of 3.6.18.3.1 and 3.6.18.4.1. If b) or c) holds, then based on 3.6.18.4 B), there is a unique connection on the extension of ${\got C}$ to $W(\bar k)$. Using Galois descent, due to its uniqueness it is obtained from a connection $\nabla$ on ${\got C}$ by extension of scalars. So the Theorem follows from 3.6.18.3.1 and 3.6.18.4.1.  
\medskip
{\bf 3.6.18.5.1. Corollary.} {\it We assume $\Phi_R$ is of essentially additive type. Then we have an equivalence of categories between $p-\Mm\Mf_{[0,1]}(R)$ and $p-\Mm\Mf_{[0,1]}^\nabla(R)$, and we have an antiequivalence (via the $\DD^{-1}$ functor) of categories between these two categories and $p-DG({\rm Spec}(R))$.} 
\medskip
{\bf Proof:} For the equivalence part, from 3.6.18.4 (we have $r=0$) we deduce that there is a unique way to associate to a $p$-divisible object ${\got C}$ of $\Mm\Mf_{[0,1]}(R)$ a $p$-divisible object $({\got C},\nabla)$ of $\Mm\Mf_{[0,1]}^\nabla(R)$. So we just need to check that this association is well behaved w.r.t. morphisms between $p$-divisible objects. 
\smallskip
Everything boils down to (keeping the previous notations): if $f\in F^0({\rm End}(M))\otimes_{W(k)} R$ is fixed by $\Phi$, then it is annihilated by the (unique) connection $\nabla$ on $M\otimes_{W(k)} R$ for which $\Phi$ is $\nabla$-parallel. Let $\nabla(f)=f_1\in {\rm End}(M)\otimes_{W(k)} \bar\Om_{R/W(k)}$. It satisfies the equation
$$f_1=p\Phi\circ d\Phi_{R*}/p(f_1)=\Phi\circ d\Phi_{R*}(f_1).\leqno (32)$$
We claim: $f_1=0$. It is enough to show that $f_1$ mod $p$ is $0$. Choosing an arbitrary $W(k)$-basis $\Mb$ of ${\rm End}(M)$, by identifying the $R/pR$-coefficients of both sides of (32) w.r.t. it and the $R$-basis $dz_1$,..., $dz_m$ of $\bar\Om_{R/W(k)}$, we get a system of equations $SE$ of first type in $m\dim_{W(k)}(M)^2$ variables and with coefficients in $R/pR$. As in the proofs of 3.6.1.2 and of 3.6.18.4 B), due to the fact that $\Phi_R$ is of essentially additive type, we get: $SE$ has a unique solution; so $f_1$ mod $p$ is $0$.
\smallskip
The antiequivalence part is implied by the equivalence part and by 3.6.18.5 b): we can change the Frobenius lift $\Phi_R$ (cf. [De3, (1.1.2.1)]) to one of additive type $\Phi_R(z_i)=z_i^p$; so Grothendieck--Messing--Berthelot's theory (see [Me, ch. 4-5] and [BM, ch. 4]) and [Fa2, th. 10] apply. In other words:
\medskip
a) each $p$-divisible object of $\Mm\Mf_{[0,1]}^\nabla(R)$ is associated to a (uniquely) $p$-divisible group over $R$ (cf. 2.2.21 UP);
\smallskip
b) this association is functorial and fully faithful (cf. [BM, ch. 4] and [Me, ch. 4-5]).
\medskip
This ends the proof.  
\medskip
{\bf 3.6.18.5.2. Corollary (the liftability property).} {\it We assume that $k=\bar k$ or that $\Phi_R$ is of essentially additive type. Let $n\in\NN$. Let $\got C$ be an object of $\Mm\Mf_{[0,1]}(R)$ such that we have an epimorphism $e:{\got C}_1\twoheadrightarrow {\got C}$, with ${\got C}_1$ as the truncation mod $p^n$ of some $p$-divisible object of $\Mm\Mf_{[0,1]}(R)$ (resp. let $\got C$ be a $p$-divisible object of $\Mm\Mf_{[0,1]}(R)$). Then any connection $\nabla$ on ${\got C}$ (resp. on ${\got C}/p^n{\got C}$) can be lifted to a connection on ${\got C}_1$ (resp. on ${\got C}$). 
\smallskip
If $\Phi_R$ is of essentially additive type, then such a lift is unique. If $k=\bar k$ and we are in the context of $e$, then the number of such lifts is a non-negative, integral power of $p$.}   
\medskip 
{\bf Proof:} We first consider the case of a $p$-divisible object ${\got C}$. The subcase $k=\bar k$ is a consequence of (the proof of) 3.6.18.4.2 b) and of 3.6.8.1.2 a). The subcase when $\Phi_R$ is of essentially additive type is a consequence of 3.6.18.4 B). 
\smallskip
We consider now the case of an object ${\got C}$. We follow closely the proof of 3.6.18.4.2. Let ${\got C}=(M,F^1,\vph_0,\vph_1)$ and let ${\got C}_1=(M_1,F^1_1,\vph_{01},\vph_{11})$. Let $q:M_1\twoheadrightarrow M$ be the $R$-epimorphism defining $e$. We view $\nabla$ as a connection $\nabla_1:M_1\to M\otimes_R \bar\Om_{R/W(k)}$ for which $e$ is $\nabla$-parallel in the sense that we have 
$$\nabla_1\circ \vph_{11}(m)=\vph_0\circ d\Phi_{R*}/p\circ\nabla_1(m)\leqno (33)$$
 if $m\in F^1_1$, and 
$$\nabla_1\circ \vph_{01}(m)=p\vph_0\circ d\Phi_{R*}/p\circ\nabla_1(m)\leqno (34)$$
if $m\in M_1$.
 So we need to lift $\nabla_1$ to a connection $\tilde\nabla:M_1\to M_1\otimes_{R} \bar\Om_{R/W(k)}$ making ${\got C}_1$ to be viewed as an object of $\Mm\Mf_{[0,1]}^\nabla(R)$.
\medskip
We consider an $R$-basis $\Mb_1=\{e_1,...,e_s\}$ of ${\got C}_1$ as in the proof of 3.6.18.4.2. From [Fa1, 2.1 iii)] we get $F^1=q(F_1^1)$. Using this we deduce (via standard induction on the length of $M_1/(z_1,...,z_m)M_1$ as a $W(k)$-module) that we can assume that the non-zero elements of the form $q(e_i)$ are forming an $R$-basis of ${\got C}$. As $M_1$ is a free $R/p^nR$-module, the natural $R$-linear map ${\rm End}(M_1)\to {\rm Hom}(M_1,M)$ is an epimorphism. We consider its kernel $N$ and the $R/p^nR$-basis $\Mf_1=\{f_1,...,f_{s^2}\}$ of ${\rm End}(M_1)$ defined, as in the proof of 3.6.18.4.2, by $\Mb_1$. Let $\nabla^1$ be an arbitrary lift of $\nabla_1$ to a connection on $M_1$.
\smallskip
$N$ is the underlying module of the object ${\got N}$ of $\Mm\Mf_{[-1,1]}(R)$ which is the kernel of the natural epimorphism 
$$End({\got C}_1)\twoheadrightarrow Hom({\got C}_1,{\got C}).$$ 
So $N$ has an $R$-basis formed by non-zero elements of the form $p^{n(i)}f_i$, where all $n(i)$'s belong to $S(0,n)$ and are uniquely determined by $q$. Let $I_0$ be the set of those $i$'s such that $n(i)\le n-1$.
An arbitrary connection $\tilde\nabla$ on $M_1$ lifting $\nabla_1$ can be written in the form
$$\tilde\nabla=\nabla^1+\sum_{i\in I_0}\sum_{j=1}^m x_{ij}p^{n(i)}f_idz_j,\leqno (35)$$
with all $x_{ij}$'s in $R/p^nR$.
The fact that $\tilde\nabla$ makes ${\got C}_1$ to be viewed as an object of $\Mm\Mf_{[0,1]}^\nabla(R)$ gets translated in $m$ equalities between elements of $N$, i.e. in a system of equations $SE$ involving the variables $x_{ij}$'s and their Frobenius transforms.
\smallskip
So we can follow entirely the proof of 3.6.18.4.2: by induction on $l\in S(1,n)$ we show that we have a solution of the system of equations obtained from $SE$ by working mod $p^l$ (i.e. of the system of equations we get by considering $m$ equalities between elements of $N/p^lN$). But, as $M_1$ is $R$-generated by $\vph_{01}(M_1)$ and by $\vph_{11}(F^1_1)$, as in 3.6.8 5) we get that $SE$ mod $p$ is a system of equations of first type. So if $k=\bar k$, 3.6.8.1 applies. If $\Phi_R$ is of essentially additive type, then we can proceed as in 3.6.18.4 P2-3 to separate the variables and to show that $SE$ mod $p$ has a unique solution. This takes care of the case $l=1$. The same things apply when we move from things mod $p^l$ to things mod $p^{l+1}$, $l\in S(1,n-1)$, cf. the proof of 3.6.18.4.2. As $f_i$'s and $z_j$'s are defined over $R$ we can assume that all linear forms have coefficients in $R/pR$; this takes care (cf. 3.1.8.1.2 c)) of the ``power of $p$" part. This proves the Corollary.   
\medskip
{\bf 3.6.18.5.3. Corollary.} {\it We assume the Frobenius lift $\Phi_R$ is of essentially additive type. Then we have an equivalence of categories between $\Mm\Mf_{[0,1]}(R)$ and $\Mm\Mf_{[0,1]}^\nabla(R)$ and we have an antiequivalence (via the $\DD^{-1}$ functor) of categories between these two categories and $p-FF({\rm Spec}(R))$.} 
\medskip
{\bf Proof:} Any object ${\got C}$ of $\Mm\Mf_{[0,1]}(R)$ is the cokernel of an isogeny $m_{12}:{\got C}_1\hookrightarrow {\got C}_2$ between two $p$-divisible objects of $\Mm\Mf_{[0,1]}(R)$ (cf. Fact of 2.2.1.1 6)). From 3.6.18.5.1 (applied to $m_{12}$) we get that there is a connection $\nabla_{\got C}$ on ${\got C}$. 3.6.18.5.2 (applied twice to the natural epimorphisms ${\got C}_2\twoheadrightarrow {\got C}_2/p^n{\got C}_2\twoheadrightarrow {\got C}$, with $n\in\NN$ big enough) implies (via 3.6.18.5.1) $\nabla_{\got C}$ is unique. So the equivalence part is obtained as in the proof of 3.6.18.5.1. 
\smallskip
For the antiequivalence part, as we have 3.6.18.5.1, we just have to remark, besides the reference to the Fact of 2.2.1.1 6), one extra thing. Any finite, flat, commutative group scheme $G_R$ over $R$ of rank a power of $p$, is a closed subgroup of a $p$-divisible group over $R$ (cf. M. Raynaud's theorem of [BBM, 3.1.1]). So, based on standard techniques pertaining to abelian categories (see 2.2.1.0 and 2.2.1.1 6)), the Corollary follows.
\medskip
{\bf 3.6.18.5.4. Remarks.} {\bf 1)} It is easy to see that the results 3.6.18.5.1-3 extend to a more general context, inspired from 3.6.18.4.2 (cf. also 3.6.18.4.4). We leave this to the reader (to be compared with 3.6.18.8 below). In particular,  one can  state a moduli variant of 3.6.18.5.2 following the pattern of 3.6.18.4.2: it is obtained by just putting together 3.6.18.4.2 and 3.6.18.5.2. We could think of it as an explicit form of particular cases of Raynaud's theorem mentioned above. 
\smallskip
{\bf 2)} 3.6.18.5.1 and 3.6.18.5.3 are the third place where we need $p>2$.
\smallskip
{\bf 3)} The condition $k=\bar k$ of 3.6.18.4 (resp. of --part of-- 3.6.18.5.2) can be weaken as follows: we just need that ${\rm Spec}(k)$ is $\dim_{W(k)}(M)^2$-simply (resp. is $\abs{I_0}$-simply) connected (see defs. 3.6.18.6 D). The same applies in other situations (for instance, in the proof of 3.6.15 B, assuming that the base $\Mb$ exists, we just need $k$ to be $C$-simply connected, with $C$ as the maximum of the orders of all cycles $\Mc$ of $\pi_L$ with $n(\Mc)\ge 0$).
\medskip
{\bf 3.6.18.5.5. The case of categories $\Mm\Mf_{[a,b]}(R)$.} A natural question arises. What about the categories $\Mm\Mf_{[a,b]}(R)$, with $a,b\in\ZZ$, $b>a+1$? In what follows, if the reader desires can assume $a=0$.
\medskip
{\bf A.} First of all, as in 3.2.3 and 3.1.1.1 we can speak about generic situations. The second philosophy of 3.6.18.2.1 says that it is enough to look at the case of such generic situations. The conclusion is: due to reasons explained in 3.6.8.9, none of the results of 3.6.18.4-5 remain true (in general), except one. We start explaining this exception. 
\smallskip
Let ${\got C}=(M_R,(F^i(M_R))_{i\in S(a+1,b)},\Phi_{M_R})$ be a $p$-divisible  object of $\Mm\Mf_{[a,b]}(R)$. As in 3.6.18.1, for $i\in S(a,b)$, we define numbers $f(i)\in\NN\cup\{0\}$ by the formula
$$f(i):=\dim_R(F^i(M_R)/F^{i+1}(M_R))).$$ 
\indent
We assume $k=\bar k$. By a generic situation modeled on ${\got C}$ we mean (to be compared with 3.6.18.2) that we are dealing with a $p$-divisible object 
$${\got C}^1=(M_R,(F^i(M_R))_{i\in S(a+1,b)},\Phi_{M_R}^1)$$
of $\Mm\Mf_{[a,b]}(R)$ such that:
\medskip
a) defining a $p$-divisible object ${\got C}_{W(k)}^1:=(M,(F^i(M))_{i\in S(a+1,b)},\vph^1)$ of $\Mm\Mf_{[a,b]}(W(k))$ as in 3.6.18.0 (i.e. via the pull back of ${\got C}^1$ through the Teichm\"uller lift ${\rm Spec}(W(k))\hookrightarrow {\rm Spec}(R)$), there is a $W(k)$-basis $\{e_1,...,e_{\dim_{W(k)}}(M)\}$ of $M$ for which we have $\vph(e_i)=p^{\vep_i}e_i$, with $\vep_i\in S(a,b)$, and $\forall j\in S(a+1,b)$ the set of those $e_i$'s such that $\vep_i\ge j$ is a $W(k)$-basis of $F^j(M)$;
\smallskip
b) $\Phi_R$ is of multiplicative type. 
\medskip
We have:
\medskip
{\bf Fact.} {\it The number of connections on ${\got C}/p^n{\got C}$ is at most
$$
p^{nm\sum_{i=a}^{b-1} f(i)f(i+1)}.
$$
Moreover this number is attained for some generic situations modeled on ${\got C}$.} 
\medskip
{\bf Proof:} Using the fact that we are dealing only with connections satisfying the Griffiths transversality condition, the first part of the Fact is a consequence of the constancy property of 3.6.8.9 and of the estimate of $m_1$ in 3.6.8.9.0. 
\smallskip
Following the pattern of 3.6.18.2, we can write (via short exact sequences) ${\got C}^1$ as extensions of $p$-divisible objects of $\Mm\Mf_{[i,i]}(R)$, where $i\in S(a,b)$ is such that $f(i)\neq 0$. So, using induction on $b-a$, the first part of the Fact for ${\got C}^1$ is as well a consequence of 3.6.18.2. Here are the details. We can assume $a=0$ and accordingly we write ${\got C}^1$ as the extension of a $p$-divisible object ${\got C}_1^1$ of $\Mm\Mf_{[1,b]}(R)$ by a $p$-divisible object of $\Mm\Mf_{[0,0]}(R)$. The number of connections on ${\got C}^1/p^n{\got C}^1$ is at most equal to the sum of the number of connections on ${\got C}_1^1/p^n{\got C}_1^1$ and of the number of connections on ${\got C}^1_{< 2}/p^n{\got C}^1_{< 2}$, where ${\got C}^1_{< 2}$ is the maximal $p$-divisible subobject of ${\got C}^1$ which is a $p$-divisible object of $\Mm\Mf_{[0,1]}(R)$; moreover, the equality holds if $\Phi_M^1$ takes $F^1(M_R)$ into itself. 
\smallskip
So if ${\got C}^1$ is a direct sum of $p$-divisible objects of $\Mm\Mf_{[i,i]}(R)$, for various $i\in S(0,b)$, then the number of connections on ${\got C}^1/p^n{\got C}^1$ is (by induction) precisely $p^{nm\sum_{i=a}^{b-1} f(i)f(i+1)}$. This ends the proof.
\medskip
But if $m\ge 2$ and if the sum
$$\sum_{i=a}^{b-2} f(i)f(i+1)f(i+2)$$ 
is not $0$, not all of them are integrable; this creates problems as well as beauty. Obviously, if $f(i)f(i+1)f(i+2)=0$, $\forall i\in S(a,b-2)$, then all these connections are integrable. 
\smallskip
We now assume $\Phi_R(z_l+1)=(z_l+1)^p$, $\forall l\in S(1,m)$. Let $\tilde F^i(M)$, $i\in S(a,b)$, be direct summands of $M$ as in 2.2.1 c). We can assume we have an identity $M_R=M\otimes_{W(k)} R$ such that $F^i(M_R)$ is nothing else but the direct sum $\oplus_{j\ge i} \tilde F^j(M)\otimes_{W(k)} R$. We also assume $\vep_i\ge\vep_{i^\prime}$, if $i\ge i^\prime$, $i$, $i^\prime\in S(1,\dim_{W(k)}(M))$. Let $g_R\in {\rm Ker}(GL(M_R)\to GL(M))$ be defined as in 3.6.18.0 but for ${\got C}^1$. Let 
$$N(n,m,f(a),f(a+1),...,f(b),g_R)$$ 
be the number of integrable connections on ${\got C}^1/p^n{\got C}^1$. 
\smallskip
We consider the maximal subgroup $N_R^{-1}$ (resp. $N_R^{-2}$) of $GL(M_R)$ with the property that its $R$-valued points are those $R$-linear automorphisms $n_R$ of $M_R$ having the property: $1_{M_R}-n_R$ takes $\tilde F^i(M)$ into $\tilde F^{i-1}(M)\otimes_{W(k)} R$ (resp. into $\tilde F^{i-2}(M)$), $\forall i\in S(a,b)$. It is a smooth, unipotent, connected subgroup and its relative dimension is 
$$\sum_{i=a}^{b-1} f(i)f(i+1)$$ 
(resp. is $\sum_{i=a}^{b-2} f(i)f(i+2)$). $N_R^1$ is the extension from $W(k)$ to $R$ of a subgroup $N_{W(k)}^{-1}$ of $GL(M)$. From the proof of the above Fact we get: we can assume 
$$g_R\in {\rm Ker}(N_W(k)^{-1}(R)\to N_{W(k)}^{-1}(W(k))).$$  
\smallskip
We consider connections $\nabla$ on $M_R$ of the following particular form
$$\nabla=\dl+\sum_{l\in S(m)} x_ldz_l,$$
with $x_l\in {\rm Lie}(N_R^{-1})$ and with $\dl$ the connection on $M_R$ that annihilates $M$. We choose arbitrarily an $R$-basis $\Mb_{-2}$ of 
$$LIE_{-2}:={\rm Lie}(N_R^{-1})\oplus {\rm Lie}(N_R^{-2}).$$
Working as usual modulo powers of $p$, 
the systems of equations needed to be satisfied by:
\medskip
-- (the coefficients of) a connection on $M_R/p^{n}M_R$ lifting a fixed connection on $M_R/p^{n-1}M_R$ which makes ${\got C}^1/p^{n-1}{\got C}^1$ potentially to be viewed as an object of $\Mm\Mf_{[a,b]}^\nabla(R)$ (here $n\in\NN$), in order to be on ${\got C}^1/p^n{\got C}^1$,
\medskip\noindent
express the equality between two elements of $LIE_{-2}\otimes_R \Om_{R/W(k)}^\wedge$; identifying the coefficients of these two elements w.r.t. the $R$-basis of $LIE_{-2}\otimes_R \Om_{R/W(k)}^\wedge$ defined naturally by $\Mb_2$ and by the $R$-basis of $\{dz_1,...,dz_m\}$ of $\Om_{R/W(k)}^\wedge$, we get systems of equations of third type (with coefficients in $R/pR$). Warning: we do need to consider $LIE_{-2}$ and not just ${\rm Lie}(N_R^{-1})$ and this is why we do not get in general systems of equations of first type. 
\smallskip
However, we can still study to some extend such numbers $N(n,m,f(a),f(a+1),...,f(b),g_R)$ following the pattern of 3.6.18.4.1 of algebraizing things (via Chinese Reminder Theorem and work performed modulo $r$-powers of maximal ideals, $r\in\NN$), cf. D below. But, based on (24) of 3.6.18.2, it is easy to see that such numbers $N(n,m,f(a),f(a+1),...,f(b),g_R)$ do depend on $g_R\in {\rm Ker}(N_{W(k)}^{-1}(R)\to N_{W(k)}^{-1}(W(k)))$. For instance, if $(n,m,a,b)=(1,2,0,2)$ then 
$$N(1,2,1,1,1,1_{M_R})=(p-1)^3+(2p-1)^2$$ 
is the number of quadruples $(q_1,...,q_4)$ formed by elements of $\ZZ/p\ZZ$ and such that $q_1q_2=q_3q_4$ (this can be easily read out from (24) of 3.6.18.2); but there are many values of $g_R$ for which we have $N(1,2,1,1,1,g_R)=0$. This suggests the introduction of types of generic situations (modeled on ${\got C}$). Not to be too long, we propose here only one such type.
\medskip
{\bf Definition.} The generic situation ${\got C}^1$ is said to be authentic if  $N(n,m,f(a),f(a+1),...,f(b),g_R)>0$. 
\medskip
Despite the fact that none of the other results of 3.6.18.4-5 remain true in general for $\got C$, many of them have weaker forms, which are also very useful. We explain them one by one. 
\medskip
{\bf B.} We start with a general definition. We come back to the general case, i.e. we are not any more in the context of a generic situation. Let ${\got C}_{W(k)}$ be defined as in a) but starting from $\got C$.
\medskip
{\bf Definitions.} We assume $b\ge a+1$ and $k=\bar k$. Let $i\in S(a,b)$. The dimension of the $\FF_p$-vector subspace 
$$V_i({\got C}_{W(k)}):=\{x\in F^i(M)/pF^i(M)|\vph_i(x)=x\}$$ 
of $F^i(M)/pF^i(M)$, 
with $\vph_i:F^i(M)\to M$ as the $\sg$-linear map obtained from $\vph$ as in 2.2.1 c) (we still denote by $\vph_i$ its reduction mod $p$), is called the pseudo-multiplicity of the slope $i$ of ${\got C}_{W(k)}$ (or of ${\got C}$ in the maximal point of ${\rm Spec}(R)$). Let $S_i({\got C}_{W(k)})$ be the $\ZZ_p$-submodule of $M$ formed by elements fixed by $p^{-i}\vph$. Let 
$$W_i({\got C}_{W(k)}):=V_i({\got C}_{W(k)})\cap S_i({\got C}_{W(k)})/pS_i({\got C}_{W(k)}),$$ 
the intersection being taken inside $M/pM$. $\dim_{\FF_p}(W_i({\got C}_{W(k)}))$ is called the virtual pseudo-multiplicity of the slope $i$ of ${\got C}_{W(k)}$ (or of ${\got C}$ in the maximal point of ${\rm Spec}(R)$). Similarly, if $X$ is as in 2.2.1 c) and if ${\got C}_X$ is an object (resp. a $p$-divisible object) of $\Mm\Mf_{[a,b]}(X)$ we define the pseudo-multiplicity (resp. the pseudo-multiplicity and the virtual pseudo-multiplicity) of the slope $i$ of ${\got C}_X$ in a geometric point of ${\rm Spec}(X_k)$. 
\medskip
Warning: in general the (virtual) pseudo-multiplicities and the multiplicities of the slope $i\in S(a,b)$ do not coincide. However, they do coincide if $i\in\{a,b\}$ (for instance, if $b=a+1$ they always coincide).
\smallskip
We come back to the initial assumption $b>a+1$, with $k$ arbitrary.
The role of $s(0)s(1)$ in 3.6.18.4 B) has to be replaced by the pseudo-multiplicity $s_{pm}(-1)$ of the slope $-1$ of the extension of $End({\got C}_{W(k)})$ to $\bar k$. As in the proof of 3.6.18.4 (see P7 of it), we get the following improvement of the first part of the Fact of A:
\medskip
{\bf Corollary.}  {\it The number ${\rm CONN}({\got C},n)$ of connections on ${\got C}/p^n{\got C}$ is at most
$$p^{nrs_{pm}(-1)}.$$}
\indent
This Corollary is a particular case of 3.6.18.7.1 b) below (and so, as it is not used below, we refer to the mentioned place).
Even if $k=\bar k$ it can happen that actually we have less than predicted such connections. To construct such examples with $g_R$ as the trivial element, we just need to consider cases when $V_i(End({\got C}_{W(k)}))$ has dimension 1 over $\FF_p$ and $W_i(End({\got C}_{W(k)}))=\{0\}$: the proof of 3.6.18.7.1 b) below shows that if in the Corollary we have equality $\forall n\in\NN$, then each connection counted by ${\rm CONN}({\got C},n)$ lifts to a connection on ${\got C}$; so we can follow 3.6.18.4 P6-7 to get that, under the assumptions on the ranks of $V_i(End({\got C}_{W(k)}))$ and of $W_i(End({\got C}_{W(k)}))$, ${\rm CONN}({\got C},1)$ has at most 1 element. 
\medskip
{\bf Example.} We consider a $W(k)$-basis $\{x,y,z\}$ of a free $W(k)$-module $\tilde M$ of rank $3$ and define: $\tilde\vph(px)=x+py$, $\tilde\vph(p^3y)=z$ and $\tilde\vph(z)=p^3y$. We take $F^0(\tilde M)=F^3(\tilde M)=<z>$, $F^{-1}(\tilde M)=F^{-2}(\tilde M)=<x,z>$, and $F^{-3}(\tilde M)=\tilde M$. Let $\tilde{\got C}_{W(k)}:=(\tilde M,(F^i(\tilde M)_{i\in S(-3,3)}),\tilde\vph)$. Its slopes are $0$, $0$ and $-1$, its pseudo-multiplicity of the slope $-1$ is $1$ but there are no $a_1$, $a_2\in W(k)$ such that $p\tilde\vph(x+a_1y+a_2z)=x+a_1y+a_2z$. So $W_i(\tilde{\got C}_{W(k)})=\{0\}$. We take ${\got C}$ to be the pull back to ${\rm Spec}(R)$ of the direct sum $\tilde {\got C}_{W(k)}\oplus W(k)(0)$.
$End({\got C})$ has the pull back to ${\rm Spec}(R)$ of $\tilde{\got C}_{W(k)}$ as a direct summand. We get (as in 3.6.18.4 P6-7) that there is only one connection on $M/pM\otimes_k R/pR$, with $M:=\tilde M\oplus W(k)$, which lifts to a connection on $M\otimes_{W(k)} R$ making ${\got C}$ potentially to be viewed as a $p$-divisible object of $\Mm\Mf_{[-3,3]}^\nabla(R)$, even if $r\neq 0$: it is the connection annihilating $M/pM$.
\medskip
We do not know when there is a constant $f({\got C})\in\NN\cup\{0\}$ such that ${\rm CONN}({\got C},n)=p^{nf({\got C})}$, $\forall n\in\NN$. Warning: in general there is no such constant $f({\got C})$.
\medskip
{\bf C.} Similarly we define the number ${\rm INT}({\got C},n)$ counting the connections of B which are moreover integrable. It seems to us that we always have ${\rm INT}({\got C},n)=N(n,m,f(a),f(a+1),...,f(b),g_R)$ for some $g_R\in {\rm Ker}(N_{W(k)}^{-1}(R)\to N_{W(k)}^{-1}(W(k)))$; warning: we do not know how to prove this using the ideas of A and of the proof of 3.6.18.4.1. However, the mentioned ideas give us (it is an easy exercise, cf. also D below): 
\medskip
{\bf Corollary.} {\it Let $q\in\NN$. If ${\rm CONN}({\got C},n)=p^{nm\sum_{i=a}^{b-1} f(i)f(i+1)}$, then ${\rm INT}({\got C},n)\ge N(n,m,f(a),f(a+1),...,f(b),g_R)$, for some $g_R\in N_{W(k)}^{-1}(R)$ which mod $(z_1,...,z_m)^q$ is the identity element.} 
\medskip
It seems to us (based on the Corollary of B), that in the above Corollary we are actually in the context of a generic situation. 
\medskip
{\bf D.} 3.6.18.4.2 has a weaker version. We use its notations but replacing the lower right index $[0,1]$ everywhere by $[a,b]$. For $m\ge 2$ we have two moduli $p$-adic formal schemes: $M({\got C})$ of integrable connections and $M_1({\got C})$ of connections; this is a consequence of 3.6.8.9 A1. Always we have an open closed natural embedding
$$M({\got C})\hookrightarrow M_1({\got C})$$
of formally \'etale, $p$-adic formal schemes over $S^\wedge$. 
But we can not say anything about their fibres. In particular they can be empty. Also, in case of an object, there is no a priori reason (besides 3.6.8.9.0) to always have all non-empty fibres having a number of elements which is a power of $p$. 
\medskip
{\bf E.} In accordance to D, if $M_1({\got C})$ or $M({\got C})$ is not empty and the special fibre of $S$ is noetherian and connected, then the surjectivity principle of 3.6.18.4.5 is preserved for it. In the same context, the purity result of 3.6.8.1.4, the touching property of 3.6.18.4.3, and the surjectivity principle of 3.6.18.6 b) below are preserved. All these can be read out from their proofs and from 3.6.8.9 A1. 
\medskip
{\bf F.} The inducing property is a much subtler thing, as it can happen to have ${\rm INT}({\got C},n)=0$ for some $n\in\NN$ (see end of A). We do not know any simple and sufficiently general criterion which could guarantee that ${\rm INT}({\got C},n)\ge 1$, $\forall n\in\NN$.
\medskip
{\bf 3.6.18.5.6. Problem.} With the notations of 3.6.18.0-2, if $r\ge 1$, study the behavior of connections (obtained as in 3.6.18.4 and 3.6.18.4.2) w.r.t. (short) exact sequences of objects of $\Mm\Mf_{[0,1]}(R)$. Hint: look at the proofs of 3.6.18.5.2-3. 
\smallskip
The similar problem in the context of 3.6.18.5.5 is much harder.
\medskip
{\bf 3.6.18.5.7. A variant.} Using 3.6.18.5 c) instead of 3.6.18.5 b), we obtain logical variants of 3.6.18.5.1 and 3.6.18.5.3. We state these variants, in a language inspired from 3.6.18.5.5 B and so suitable for future generalizations. Regardless of how $\Phi_R$ and $k$ are, we have:
\medskip
{\bf Corollary.} {\it The full subcategory of $p-FF({\rm Spec}(R))$ (resp. of $p-DG({\rm Spec}(R))$) whose objects are such that either their \'etale part or their multiplicative type part is trivial, is antiequivalent (via the $\DD$ functor) to the full subcategory $s_{pm}(-1)=0-\Mm\Mf_{[0,1]}(R)$ (resp. $s_{pm}(-1)=0-p-\Mm\Mf_{[0,1]}(R)$) of the category $\Mm\Mf_{[0,1]}(R)$ (resp. of the category $p-\Mm\Mf_{[0,1]}(R)$) whose objects are having the pseudo-multiplicity of the slope $-1$ of the $End$ objects of their truncations mod $p$ equal to $0$. Moreover, $s_{pm}(-1)=0-\Mm\Mf_{[0,1]}(R)$ (resp. $s_{pm}(-1)=0-p-\Mm\Mf_{[0,1]}(R)$) is equivalent to the similarly defined full subcategory $s_{pm}(-1)=0-\Mm\Mf_{[0,1]}^\nabla(R)$ (resp. $s_{pm}(-1)=0-p-\Mm\Mf_{[0,1]}^\nabla(R)$) of $\Mm\Mf_{[0,1]}^\nabla(R)$ (resp. of $p-\Mm\Mf_{[0,1]}^\nabla(R)$).}
\medskip
{\bf Proof:} The proof of 2.2.1.1 6) shows that any object of $\Mm\Mf_{[0,1]}(R)$ has a strict lift in the sense of 2.2.1 e); so any object of $s_{pm}(-1)=0-\Mm\Mf_{[0,1]}(R)$ is the epimorphism of an isogeny between two objects of $s_{pm}(-1)=0-p-\Mm\Mf_{[0,1]}(R)$. So, copying the proofs of 3.6.18.5.1 and 3.6.18.5.3, we just have to deal with the part involving morphisms, in the context of $p$-divisible objects and groups. But over ${\rm Spec}(R)$, there is no non-trivial homomorphism from a finite, flat, commutative group scheme $G_1$ of $p$-power order whose pull back to ${\rm Spec}(R^{\rm sh})$ is connected into a finite, \'etale, commutative group scheme $G_2$ of $p$-power order; moreover, as $p\ge 3$, there is no non-trivial homomorphism $G_2\to G_1$. So, all morphisms between ($p$-divisible) objects to be considered are such that (the proof of) 3.6.18.4 c) applies to get that all elements $f_1$ as in the proof of 3.6.18.5.1 are $0$.
\medskip
3.6.18.5.7 is part of the third place where we need $p\ge 3$.
\medskip
{\bf 3.6.18.5.8. An equivalence of categories.} Let $S$ be as in 3.6.18.4.2. As in 3.6.18.5.7 we define the categories $s_{pm}(-1)=0-\Mm\Mf_{[0,1]}(S)$, $s_{pm}(-1)=0-p-\Mm\Mf_{[0,1]}(S)$, $s_{pm}(-1)=0-\Mm\Mf_{[0,1]}^\nabla(S)$, $s_{pm}(-1)=0-p-\Mm\Mf_{[0,1]}^\nabla(S)$. By just combining 3.6.18.4.2, 3.6.18.5 c) and 3.6.18.5.7 we get:
\medskip
{\bf Corollary.} {\it The category $s_{pm}(-1)=0-\Mm\Mf_{[0,1]}(S)$ (resp. $s_{pm}(-1)=0-p-\Mm\Mf_{[0,1]}(S)$) is naturally equivalent to the category $s_{pm}(-1)=0-\Mm\Mf_{[0,1]}^\nabla(S)$ (resp. to $s_{pm}(-1)=0-p-\Mm\Mf_{[0,1]}^\nabla(S)$).}
\medskip
{\bf 3.6.18.5.9. A formula.} We refer to the morphism $e$ of 3.6.18.5.2, in the case $k=\bar k$. Following the proofs of 3.6.18.5.2 and 3.6.18.4 B), we get  that the number of lifts of $\nabla$ to connections making ${\got C}_1$ to be viewed as an object of $\Mm\Mf_{[0,1]}^\nabla(R)$ is precisely 
$$p^{r\sum_{i=0}^{n-1}a_i},$$ 
where $a_i$ is the pseudo-multiplicity of the slope $-1$ of $p^i{\got N}/p^{i+1}{\got N}$ (in the maximal point of ${\rm Spec}(R)$). Here ${\got N}$ is as in the proof of 3.6.18.5.2. As $p^{i+1}{\got N}/p^{i+2}{\got N}$ is naturally a subobject of $p^i{\got N}/p^{i+1}{\got N}$, we have:
$$a_0\ge a_1\ge...\ge a_{n-1}.$$ 
\smallskip
For future references, we state here a combined version of 3.6.1.3 and 3.6.18.4.2.
\medskip
{\bf 3.6.18.6. Theorem.} 
{\it Let ${\rm Spec}(R_1)$ be a regular, formally smooth, affine $W(k)$-scheme. We assume that its special fibre is connected and that it is part of a triple $({\rm Spec}(R_1),\Phi_{R_1},z)$ defining a potential-deformation sheet. Let ${\got C}_1$ be a $p$-divisible object of $\Mm\Mf_{[0,1]}(R_1)$. 
\smallskip
{\bf a)} {\bf (the universal global $\nabla$ principle for $\Mm\Mf_{[0,1]}(*)$ for potential-deformation sheet contexts)}  The whole of 3.6.1.3 remains true (except the fact that we get connections respecting the $G$-action which does not make any sense here) in this context, provided we state the things in terms of $p$-adic formal schemes. In particular,  there is a moduli $p$-adic formal scheme $M_0({\got C}_1)$ having a geometrically connected special fibre and there is a formally \'etale, affine morphism $d_1:M_0({\got C}_1)\to {\rm Spec}(R_1^\wedge)$ such that $d_1^{*}({\got C}_1)$ is a $p$-divisible object of $\Mm\Mf_{[0,1]}^\nabla(M_0({\got C}_1))$ and which are universal in the following sense:
\medskip
{\bf UP} For any formally \'etale, $p$-adic formal scheme $X$ over ${\rm Spec}(R_1^\wedge)$ such that all connected components of its special fibre contain points mapping into the special fibre $z_k$ of $z$ and the pull back ${\got C}_1^X$ of ${\got C}_1$ to it is a $p$-divisible object of $\Mm\Mf_{[0,1]}^\nabla(X)$, there is a unique ${\rm Spec}(R_1^\wedge)$-morphism $\ell_X:X\to M_0({\got C}_1)$ (so ${\got C}_1^X$ is $(d_1\circ\ell_X)^*{\got C}_1$). 
\medskip
The special fibre of $d_1$ is an $\NN$-pro-\'etale, affine morphism. Moreover, the point $z$ lifts uniquely to a point ${\rm Spf}(W(k))\to M_0({\got C}_1)$.
If the underlying module of ${\got C}_1$ is free, then $d_1$ is defined naturally by the $p$-adic completion of an $\NN$-pro-\'etale morphism $\ell_1:{\rm Spec}(Q_1)\to {\rm Spec}(R_1)$.
\smallskip
{\bf b)} {\bf (the surjectivity principle)} The fibres of $d_1$ over an open subscheme of the special fibre of ${\rm Spec}(R_1^\wedge)$ containing $z_k$, are not empty. 
\smallskip
The same remains true for objects: we just need to replace the words $\NN$-pro-\'etale by the word \'etale (b) becomes trivial). If $R_1/pR_1$ is of finite type over $k$, then $M_0({\got C}_1)_k$ is an $AG$ $k$-scheme.}
\medskip
The proof of this is entirely the same as the proof of 3.6.1.3 (cf. 3.6.18.4.1 for the integrability part),  and so it is omitted. We just mention two things. The localization process mentioned in the proof of 3.6.18.4.2, allows us to assume that we are in the context described in the last sentence of a). $M_0({\got C}_1)$ is the (only) connected component of the moduli $p$-adic formal scheme $M({\got C}_1)$ of 3.6.18.4.2 to which $y$ lifts.
\medskip
{\bf 3.6.18.7. The relative situation.} 3.6.18.4-6 works in the relative context, i.e. when we start with a Shimura filtered $\sg$-crystal or, even more generally, with a suitable filtered $\sg$-$\Ms$-crystal over $k$. So here, with the notations of the opening paragraphs of 3.6.18, we assume that there is a smooth  subgroup $\tilde G$ of $GL(M)$ having connected fibres and the property that there is  a family of tensors $(t_{\al})_{\al\in\Mj}$ of the $F^0$-filtration of $\Mt(M[{1\over p}])$ defined by $F^1$ such that $\vph(t_{\al})=t_{\al}$, $\forall\al\in\Mj$, and $\tilde G_{B(k)}$ is the subgroup of $GL(M[{1\over p}])$ fixing $t_{\al}$, $\forall\al\in\Mj$. Till the end of 3.6.18.7.1 we assume $k=\bar k$. We have:
\medskip
{\bf Theorem.} {\it If $g_R\in\tilde G(R)$, then there is an integrable connection $\nabla$ on $M\otimes_{W(k)} R$, which respects the $\tilde G$-action (in the similar sense as of 3.6.1.1.1) and makes $\got C$ to be viewed as a $p$-divisible object of $\Mm\Mf_{[0,1]}^\nabla(R)$. Moreover, this $p$-divisible object is induced from a filtered $F$-$\Ms$-crystal of the form $(M,F^1,\vph,\tilde G,\tilde f)$ (it is defined as in 2.2.10).}
\medskip
{\bf Proof:} This is a consequence of 3.6.18.4.1 and [Fa2, rm. iii) after th. 10] (to be compared with 3.6.18.3.1), once we show that there is a connection on $M\otimes_{W(k)} R$ which makes $\got C$ to be viewed as a $p$-divisible object of $\Mm\Mf_{[0,1]}^\nabla(R)$ and which respects the $\tilde G$-action. To see this, following the ideas of the steps 4) to 7), 11) and 12) of 3.6.8, we get that it is enough to show that there is a connection on $M\otimes_{W(k)} R/pR$ which makes ${\got C}/p{\got C}$ to be viewed as an object of $\Mm\Mf_{[0,1]}^\nabla(R)$ and which respects the $\tilde G$-action (i.e. the liftability property of 3.6.18.5.2 works also in this relative context pertaining to filtered $\sg$-$\Ms$-crystals: this can be deduced from the way the paragraph below is presented). 
\smallskip
We start with a connection $\nabla=\delta_0+\be$, where $\delta_0$ is the connection on $M\otimes_{W(k)} R/pR$ annihilating $M$ and where
$$
\be\in LIE:={\rm Lie}(\tilde G)\otimes_{W(k)} \Om_{(R/pR)/k}.
$$ 
The equations $(E_1)$ and $(E_2)$ of 3.6.1.1.1 2), in this case (as $\Phi(t_\al)=t_\al$, $\forall\al\in\Mj$), express for each $l\in S(M)$ the equality between two elements of $LIE$: this can be read out from 3.6.8.6.1-2 (cf. also 2.2.1.2). So we get a system of $m\dim_{W(k)}(({\rm Lie}(\tilde G))$ equations in the same number of variables; the resulting system of equations is of first type and so 3.6.8.1.2 a) applies. This ends the proof of the Theorem.
\smallskip
Warning: not all connections on ${\got C}/p{\got C}$ respect the $\tilde G$-action; simple examples can be constructed starting from the generic situation described in 3.6.18.2. However, for $n\in\NN$, we have:
\medskip
{\bf 3.6.18.7.0. Formulas.} {\it The number of connections on $M\otimes_{W(k)} R/p^nR$ that respect the $\tilde G$-action and make ${\got C}/p^n{\got C}$ to be viewed as an object of $\Mm\Mf_{[0,1]}^\nabla(R)$ is $p$ to the power $nr$ times the number $s(-1)$ of slopes $-1$ of the Lie $\sg$-crystal $({\rm Lie}(\tilde G),\vph)$. Each such connection lifts in precisely $p^{rs(-1)}$ ways to a connection on $M\otimes_{W(k)} R/p^{n+1}R$ respecting the $\tilde G$-action and making ${\got C}/p^{n+1}{\got C}$ to be viewed as an object of $\Mm\Mf_{[0,1]}^\nabla(R)$.}
\medskip
The proof of this is entirely the same as the proof of 3.6.18.4 B).
\medskip
{\bf 3.6.18.7.1. Formulas and estimates (the general form).} We refer to the filtered $\sg$-$\Ms$-crystal of 3.6.1.5. Choosing an element $g_R\in\tilde G(R)$, let 
$${\got C}:=(\tilde M\otimes_{W(k)} R,(F^i(\tilde M)\otimes_{W(k)} R)_{i\in S(\tilde a+1,\tilde b)},g_R(\tilde\vph\otimes 1)).$$ 
It is a $p$-divisible object of $\Mm\Mf_{[\tilde a,\tilde b]}(R)$. We have:
\medskip
{\bf Theorem. a)} {\it If $\tilde a=0$ and $\tilde b=1$, then 3.6.18.7.0 still holds in this general context.
\smallskip
{\bf b)} In general, the number of connections on $M\otimes_{W(k)} R/p^nR$ which respect the $\tilde G$-action and make ${\got C}/p^n{\got C}$ potentially to be viewed as an object of $\Mm\Mf_{[\tilde a,\tilde b]}^\nabla(R)$ is at most equal to $p^{nrs_{pm}(-1)}$, where $s_{pm}(-1)$ is the pseudo-multiplicity of the slope $-1$ of $({\rm Lie}(\tilde G),(F^i({\rm Lie}(\tilde G)))_{i\in SS(\tilde a,\tilde b)},\tilde\vph)$.
\smallskip
{\bf c)} We assume the filtration of ${\rm Lie}(\tilde G)$ is in the range $[-1,1]$. Then in b) we have equality.}
\medskip
{\bf Proof:} Let $\tilde d:=m\dim_{W(k)}(F^{-1}({\rm Lie}(\tilde G)))$. All systems and subsystems of equations to be considered below are in $\tilde d$ variables and have coefficients in $R/pR$. For b) we just need to show: for any $i\in S(0,n-1)$, every connection $\nabla_0$ on $\tilde M\otimes_{W(k)} R/p^iR$ respecting the $\tilde G$-action and making ${\got C}/p^i{\got C}$ potentially to be viewed as an object of $\Mm\Mf_{[\tilde a,\tilde b]}^\nabla(R)$, can be lifted in at most $p^{rs_{pm}(-1)}$ ways to a connection on $\tilde M\otimes_{W(k)} R/p^{i+1}R$ having the properties:
\medskip
i) it respects the $\tilde G$-action, and
\smallskip
ii) it makes ${\got C}/p^{i+1}{\got C}$ potentially to be viewed as an object of $\Mm\Mf_{[\tilde a,\tilde b]}^\nabla(R)$. 
\medskip
We first remark that, as in 3.6.8.9, we can ``capture" the fact that a connection $\nabla$ on $\tilde M\otimes_{W(k)} R/p^{i+1}R$ lifting $\nabla_0$ satisfies i) and ii), by using a system of equations $SE$ of third type. We consider their subsystems which are of first type (see 3.6.8.9 A1). To be more precise, for instance, we can choose an $R$-basis $\Mb$ of $\tilde M$ such that suitable subbases of it are $W(k)$-bases of its direct summands $F^i(\tilde M)$, $i\in S(a,b)$. Using it we put $SE$ in a convenient matrix form as in 3.6.8.6 (CMF) (see also 3.6.8.6.1) and we consider just the subsystem defined by entries which correspond to the direct summand $F^{-1}({\rm Lie}(\tilde G))$ of ${\rm End}(\tilde M)$. So the choice of $\Mb$ gives birth naturally to a direct supplement $DS$ of $F^{-1}({\rm Lie}(\tilde G))$ in ${\rm End}(\tilde M)$ and we ``view" our systems of equations ``modulo it" (strictly speaking modulo $DS\otimes_{W(k)} \Om_{R/W(k)}^\wedge$). However, in what follows, it is important to allow such a direct supplement to be arbitrary (i.e. to vary) and in particular, not necessarily to be defined by a $W(k)$-basis of $\tilde M$; let $SDS$ be the set of such $DS$'s. Let $\pi_{DS}$ be the reduction mod $p$ of the projector of ${\rm End}(\tilde M)$ on $F^{-1}({\rm Lie}(\tilde G))$ having $DS$ as its kernel. So for any $i\in S(0,n-1)$ and for each such $DS\in SDS$ we obtain (via the projection of ${\rm End}(\tilde M)\otimes_{W(k)} \Om_{R/W(k)}^\wedge$ on $F^{-1}({\rm Lie}(G))\otimes_{W(k)} \Om_{R/W(k)}^\wedge$ ``along" $DS\otimes_{W(k)} \Om_{R/W(k)}^\wedge$) a subsystem $SS_{DS}(i)$ of $SE$ formed by precisely $\tilde d$ equations; it is of first type and it is uniquely determined up to equivalence.
\smallskip
From 3.6.8.9 A1 and 3.6.8.1.2 c), as $R/pR$ is strictly henselian, we get that the number $SOL_{DS}$ of solutions of $SS_{DS}(i)$ is independent on $i$ and $n$; so in what follows we just write $SS_{DS}$ instead of $SS_{DS}(i)$. Obviously, the number of lifts of $\nabla_0$ as above (i.e. which have properties i) and ii)) is at most $p^{SOL_{DS}}$. 
\smallskip
So, as these subsystems of equations are of first type, we can follow entirely the proof of 3.6.18.4 B) (see P1 to P3 and P7 of it): $SOL_{DS}$ is equal to $p$ to the power $r$ times the pseudo-multiplicity $s_{pm}^{DS}(-1)$ of the slope $-1$ of $({\rm Lie}(\tilde G),(F^i({\rm Lie}(\tilde G)))_{i\in SS(\tilde a,\tilde b)},\tilde\vph)$ ``along DS". In other words,
$s_{pm}^{DS}(-1)$ is the dimension of the $\FF_p$-vector space
$$V_{-1}^{DS}({\got C}):=\{x\in F^{-1}({\rm Lie}(\tilde G))/pF^{-1}({\rm Lie}(\tilde G))|\pi_{DS}\circ\tilde \vph_{-1}(x)=x\},$$
where $\tilde\vph_{-1}$ is (the reduction mod $p$ of) the $\sg$-linear map $F^{-1}({\rm End}(\tilde M))\to {\rm End}(\tilde M)$ defined as in 2.2.1 c) by $p\tilde\vph$ (see (ENDFR) of 2.2.4 B). As we are allowed to vary $DS$, we jut need to remark:
$$s_{pm}(-1)={\rm min} \{s_{pm}^{DS}(-1)|DS\in SDS\}.\leqno (36)$$
To prove this formula, let $DS_0\in SSD$ be such that $s_{pm}^{DS_0}(-1)$ attains the above minimum. Let $\{e_1,...,e_{\tilde q}\}$ be an $\FF_p$-basis of $V_{-1}^{DS_0}({\got C})$, and let $\{e_1,...,e_{\dim_{W(k)}(F^{-1}({\rm Lie}(\tilde G)))}\}$ be a $k$-basis of $F^{-1}({\rm Lie}(\tilde G))/pF^{-1}({\rm Lie}(\tilde G))$ w.r.t. which $\pi_{DS_0}\circ\tilde \vph_{-1}$ has an upper triangular form. Obviously, $s_{pm}^{DS_0}(-1)=\tilde q\ge s_{pm}(-1)$. If $\tilde q>s_{pm}(-1)$, then there is $j\in S(1,\tilde q)$ such that $\tilde \vph_{-1}(e_j)=e_j+f_j$, with $f_j\in DS_0$ a non-zero element. We can assume $j=1$. Let $\{f_1,...,f_{{\dim_{W(k)}(\tilde M)}^2-\dim_{W(k)}(F^{-1}({\rm Lie}(\tilde G)))}\}$ be a $k$-basis of $DS_0$. Replacing $DS_0$ by another 
$$DS_1\in SDS$$ 
such that $\pi_{DS_1}$ is the sum of $\pi_{DS_0}$ with the $k$-linear endomorphism of ${\rm End}(\tilde M/p\tilde M)$ which annihilates all $e_u$'s ($u\in S(1,\dim_{W(k)}(F^{-1}({\rm Lie}(\tilde G))))$) and all $f_s$'s ($s\in S(2,\dim_{W(k)}(\tilde M))$) and which takes $f_1$ into $-e_1$, we get $s_{pm}^{DS_1}(-1)\le -1+s_{pm}^{DS_0}(-1)$. Contradiction. This proves (36), and so b). 
\medskip
For a) and c), let $LIE$ be as in 3.6.18.7. To show that equality holds if $\tilde a=0$ and $\tilde b=1$, we just have to show that in fact the systems of equations of third type we come across are in fact of first type. 
\smallskip
Let $N_0$ be the Zariski closure of the connected component of the origin of the subgroup of $GL(\tilde M[{1\over p}])$ normalizing ${\rm Lie}(\tilde G)[{1\over p}]$ in $GL(\tilde M)$. So $N_0$ normalizes $\tilde G$. The canonical split cocharacter $\tilde\mu:\GG_m\to GL(\tilde M)$ of $(\tilde M,F^1(\tilde M),\tilde\vph)$ factors through $N_0$, cf. our assumptions of 3.6.1.5. So we can write $\tilde\vph=\tilde a\tilde\mu\bigl({1\over p}\bigr)$, with $\tilde a$ as a $\sg$-linear automorphism of $\tilde M$, normalizing ${\rm Lie}(\tilde G)$. But 3.6.8.6.1-2 can be entirely adapted to this ``normalizing" context, so that 3.6.8.6.1-2 apply as in 3.6.18.7: the systems of equations we get are expressing the equality between two elements of $LIE$ and so, as ${\got C}$ is a $p$-divisible object of $\Mm\Mf_{[0,1]}(W(k))$, they are of first type. This proves a).
\smallskip
The proof of c) is entirely the same as the proof of a). We just have to remark that (GT) of 3.6.1.1 4) represents no restrictions on $\be$ (i.e. writing $\nabla=\nabla^0+\be$, with $\nabla^0$ a fixed connection on $\tilde M\otimes_{W(k)} R/p^{n+1}R$ lifting $\nabla_0$, $\be\in p^m{\rm Lie}(\tilde G)\otimes_{W(k)} \Om_{R/W(k)}^\wedge/p^{m+1}{\rm Lie}(\tilde G)\otimes_{W(k)} \Om_{R/W(k)}^\wedge\tilde\to LIE$ involves as many coefficients --viewed as elements of $R/pR$-- as $\dim_{R/pR}(LIE)$); so the fact that $\nabla$ is a connection on ${\got C}/p^{n+1}{\got C}$ gets translated in: these coefficients are solutions of a quasi Artin--Schreier system of equations in as many variables as $\dim_{R/pR}(LIE)$. This ends the proof.
\medskip
We now combine 3.6.18.7.1 with a variant of 3.6.18.5.9. Let $a$, $b$ and $R$ have the same significance as in 3.6.18.5.5. Let ${\got C}$ be an arbitrary object of $\Mm\Mf_{[a,b]}(R)[p^n]$. Let 
$$s_{pm}(-1,{\got C})$$
be the pseudo-multiplicity of the slope $-1$ of
$$\oplus_{i=0}^{n-1} p^iEnd({\got C})/p^{i+1}End({\got C}).$$
\indent
{\bf 3.6.18.7.1.1. Corollary.} {\it a) We assume $(a,b)=(0,1)$. The number of connections on ${\got C}$ is precisely $p^{rs_{pm}(-1,{\got C})}$.
\smallskip
b) We do not make any assumption on $a$ and $b$. The number of connections on ${\got C}$ is at most $p^{rs_{pm}(-1,{\got C})}$.} 
\medskip
a) follows by combining the part of the proof of 3.6.18.4.2 referring to extra features with the proof of 3.6.18.7.1. b) follows from the proof of 3.6.18.7.1 b). The only difference from this proof: when we start lifting connections ``pertaining to" $End({\got C})/p^iEnd({\got C})$ to connections ``pertaining to" $End({\got C})/p^{i+1}End({\got C})$ (see the proof of 3.6.18.4.2), we come across systems of equations in a potentially less number of variables and so we are dealing with the pseudo-multiplicity of the slope $-1$ of $p^iEnd({\got C})/p^{i+1}End({\got C})$ (and not just of $End({\got C})/pEnd({\got C})$).
\medskip
{\bf 3.6.18.7.2. Relative forms of 3.6.18.4.6.} We come back to the starting setting of 3.6.18.7. It is worth detailing what 3.6.18.4.6 becomes in the relative context. There are two variants: one modeled on the faithful representation $\rho$ of $\tilde G_k$ on $M/pM$ and another one modeled on $(M,F^1,\vph,\tilde G)$. We first detail the first one and then we hint the idea of the second one.
\smallskip
Let $\tilde m:=\dim_{W(k)}(M)$. Let $\bar R$ be a $k$-algebra such that ${\rm Spec}(\bar R)$ is connected. In the proof of 3.6.18.4.6 C we came across Artin--Schreier systems of equations in $n$ variables ($n\in\NN$) with coefficients in $\bar R$, whose $n\times n$ matrix $A_L$ is invertible (i.e. $A_L\in GL_n(\bar R)$). If now $\tilde m=n$ and $A_L\in\tilde G_k(\bar R)$ (under the representation $\rho$), then we say that the Artin--Schreier system (*) of equations of 3.6.18.4.6 A is modeled on $\rho$, and any open closed subscheme of the ${\rm Spec}(\bar R)$-scheme defined by it is said to be a crystalline elementary $\rho$-admissible \'etale cover of ${\rm Spec}(\bar R)$. As in 3.6.18.4.6 we define crystalline $\rho$-admissible \'etale covers of ${\rm Spec}(\bar R)$, as well as Galois groups
$$\Pi^{\rho-AS}_{\tilde m,m}$$
and 
$$\Pi^{\rho-AS}_{\tilde m},$$
by using only Artin--Schreier systems of equations which are modeled on $\rho$. Here $m\in\NN$ has the same role as in 3.6.18.4.6 B, i.e. in the definition of $\Pi^{\rho-AS}_{\tilde m,m}$ we consider sequences of crystalline $\rho$-admissible \'etale covers of length at most $m$. Warning: all Artin--Schreier systems of equations modeled on $\rho$ are in precisely $\tilde m$ variables. We refer to $\Pi^{\rho-AS}_{\tilde m}$ (resp. to $\Pi^{\rho-AS}_{\tilde m,m}$) as the $\rho$ (resp. as the level $m$ $\rho$) Artin--Schreier fundamental group of ${\rm  Spec}(\bar R)$ (the lower right index $\tilde m$ just has the role of keeping track of the dimension of the faithful representation $\rho$). Similarly, following the pattern of 3.6.18.4.6, we can define (using in 3.6.18.4.6 D just matrices $A_i$'s and $B_{ij}$'s which are defined by points of $\tilde G_k$) for any separated, connected $k$-scheme $X_p$:
\medskip
-- a $\rho$ Galois--Artin--Schreier fundamental group $\Pi^{\rho-GAS}_{\tilde m}(X_p)$; 
\smallskip
--  a level $m$ $\rho$ Galois--Artin--Schreier fundamental group $\Pi^{\rho-GAS}_{\tilde m,m}(X_p)$. 
\medskip
The most interesting cases are when $\rho$ is an irreducible representation.   For instance, we would like to understand the case when $\rho$ is the spinorial (resp. semi-spinorial) representation of a split semisimple group over $k$ of $B_l$ (resp. of $D_l$) Lie type, with $l\in\NN$, $l\ge 3$  (resp. $l\ge 4$). As in 3.6.18.6 F, we define the notion of (weakly) $\rho$-simply connectivity. 
\smallskip
The second variant imitates the first variant as well as the case described in the equations (3) and (4) of 3.6.8; in other words, we require that some specific entries of $A_L$ are $0$, without requiring that $A_L$ is invertible or that it defines an element of $\tilde G_k$. We intend to come back in a future paper to detail this second variant, as we think it is equally important as the first one. 
\medskip
{\bf 3.6.18.7.3. The relative moduli principle.} Let $X$ be a regular, formally smooth, faithfully flat $W(k)$-scheme equipped with a Frobenius lift $\Phi_X$ and having a connected special fibre. We consider a $p$-divisible object $\got C$ of $\Mm\Mf_{[a,b]}(X)$ and a family of tensors $(t_{\al})_{\al\in\Mj}$ of the essential tensor algebra of the underlying $\Mo_{X^\wedge}$-sheaf $\Mf$ of $\got C$ such that there is an open, affine cover of $X$ with the property that for any member $U={\rm Spec}(R)$ of it whose special fibre is non-empty, denoting by $\Phi_U$ the Frobenius lift of $U^\wedge$ naturally defined by $\Phi_X$ and denoting by $(M_R,(F^i(M_R))_{i\in S(a,b)},\Phi_{M_R})$ the pull back of ${\got C}$ to $U^\wedge$ (so $M_R$ is a projective $R^\wedge$-module), the mentioned sections become sections (we keep the same notation) of $\Mt(M_R)$ satisfying:
\medskip
{\bf i)} they are in the $F^0$-filtration of $\Mt(M_R)$ defined by $(F^i(M_R))_{i\in S(a,b)}$ and are fixed by $\Phi_{M_R}$;
\smallskip
{\bf ii)} the Zariski closure of the subgroup of $GL(M_R[{1\over p}])$ fixing $t_{\al}$, $\forall\al\in\Mj$, in $GL(M_R)$ is a reductive group $\tilde G_R$ over $R^\wedge$.
\medskip
The intersection $F^i({\rm Lie}(G_R)):=F^i({\rm End}(M_R))\cap {\rm Lie}(G_R)$ is a direct summand of ${\rm Lie}(G_R)$, $i=\overline{0,1}$. The argument for this is entirely the same as the one of the first paragraph of the proof of 2.2.20.1 9).
\smallskip
Using Teichm\"uller lifts ${\rm Spec}(W(\tilde k))\to {\rm Spec}(R^\wedge)$, with $\tilde k$ a perfect field containing $k$, and applying 2.2.1.2 to the context of pull backs of $({\got C},(t_{\al})_{\al\in\Mj}$ via them, we get that (the pull back) of $F^0({\rm Lie}(G_R))$ through eacg geometric point of ${\rm Spec}(R/pR)$ or dominating ${\rm Spec}(R^\wedge[{1\over p}])$ is the Lie subalgebra of a parabolic subgroup of (the corresponding pull backs of $G_R$). Even better, we have:
\medskip
{\bf A. Claim.} {\it $F^0({\rm Lie}(G_R))$ is the Lie algebra of a uniquely determined parabolic subgroup $P_R$ of $G_R$. $F^1({\rm Lie}(G_R))$ is the Lie algebra of the unipotent radical of $P_R$.} 
\medskip
{\bf Proof:} The second part follows from the first part, via pull backs through Teichm\"uller lifts as above. The uniqueness of $P_R$ follows from the Fact of 2.2.11.1. 
\smallskip
To see that $P_R$ exists, we start remarking that we have a natural identification
$${\rm Lie}(G^{\rm ad})=\{x\in {\rm End}({\rm Lie}(G_R))|p^nx\in {\rm ad}({\rm Lie}(G_R)), {\rm for}\, {\rm some}\, n\in\NN\}.$$ 
So ${\rm Lie}(G^{\rm ad})$ can be interpreted as well as the underlying $R^\wedge$-module of the kernel (see 2.2.1.1 6)) of a morphism between $p$-divisible objects of $\Mm\Mf(R)$. So we similarly get direct summands $F^0({\rm Lie}(G^{\rm ad}_R))$ and $F^1({\rm Lie}(G^{\rm ad}_R))$ of ${\rm Lie}(G^{\rm ad}_R)$. As above, in each geometric point of ${\rm Spec}(R)$ factoring trough a Teichm\"uller lift, $F^0({\rm Lie}(G^{\rm ad}_R))$ becomes the Lie algebra of a parabolic subgroup of (the corresponding pull back of $G^{\rm ad}_R$) whose unipotent radical has as its Lie algebra the corresponding pull back of $F^1({\rm Lie}(G^{\rm ad}_R))$; in particular, this applies to all geometric points of ${\rm Spec}(R/pR)$. 
\smallskip
As $F^0({\rm Lie}(G^{\rm ad}_R))$ is a direct summand of ${\rm Lie}(G_R^{\rm ad})$, we can speak about the subgroup $N$ of $G_R^{\rm ad}$ normalizing it. We now check that $N_{R/pR}$ is a parabolic subgroup of $G^{\rm ad}_{R/pR}$ having $F^0({\rm Lie}(G^{\rm ad}_R))/pF^0({\rm Lie}(G^{\rm ad}_R))$ as its Lie algebra. 
\smallskip
We first assume $R=W(k)$. We consider a maximal torus $T_{W(k)}$ of the parabolic subgroup $PS$ of $G^{\rm ad}_{W(k)}$ whose Lie algebra is $F^0({\rm Lie}(G^{\rm ad}_{W(k)}))$. It acts on ${\rm Lie}(G^{\rm ad}_{W(k)}$ via inner conjugation. Let $F^{-1}$ be the opposite of $F^1({\rm Lie}(G^{\rm ad}_{W(k)}))$ w.r.t. this action. As we are in an adjoint context, $T_{W(k)}$ acts on $F^{-1}$ through characters which are not multiples by elements of $\NN\setminus\{1\}$ of characters of $T_{W(k)}$; so no element of $F^{-1}/pF^{-1}$ can take ${\rm Lie}(T_k)$ into ${\rm Lie}(PS_k)$. So ${\rm Lie}(N_k)={\rm Lie}(PS_k)$. As $N_k$ contains $PS_k$, we have (see Fact of 2.2.11.1) $N_k=PS_k$.
\smallskip
We come back to the general $R$. Using Teichm\"uller lifts we get that all fibres of $N_{R/pR}$ are parabolic subgroups and so are smooth and connected; they also have the same dimension. We get: $N_{R/pR}$ is smooth (using translations, it is enough to check that $N_{R/pR}$ is smooth in the origins of its fibres; but the argument for this is the same as the one of the second paragraph after [Va2, Claim 1 of p. 463]) and so a parabolic subgroup of $G_{R/pR}$. Based on this on the fact that each point of ${\rm Spec}(R)$ specializes to a point of ${\rm Spec}(R/pR)$, we get that all fibres of $N$ are smooth and have the same dimension. This implies: $N$ is smooth in all points of the connected component of the origin of any one of its fibres. So, as each point of ${\rm Spec}(R)$ specializes to a point of ${\rm Spec}(R/pR)$, we get that all such fibres are connected and so $N$ is smooth. Using Teichm\"uller lifts, we get that its Lie algebra is $F^0({\rm Lie}(G_R^{\rm ad}))$. 
\smallskip
As $N$ is a smooth subgroup of $G^{\rm ad}_R$, as $N_{R/pR}$ is a parabolic subgroup of $G_{R/pR}$, and as each point of ${\rm Spec}(R)$ specializes to a point of ${\rm Spec}(R/pR)$, $N$ is a parabolic subgroup of $G^{\rm ad}_R$. Argument: we can assume that $N$ contains a maximal torus $T_R$ of $G_{R}^{\rm ad}$ which is split, cf. [SGA3, Vol. III, 6.1 of p. 32]; as the set of characters of the action of $T_R$ on ${\rm Lie}(N)$ can be read out mod $p$, all fibres of $N_{R[{1\over p}]}$ are --cf. [Bo2, 7.1, 14.17-18 and 11.16]-- parabolic subgroups of corresponding fibres of $G_{R[{1\over p}]}^{\rm ad}$. So the inverse image of $N$ under the natural epimorphism $G_R\twoheadrightarrow G_R^{\rm ad}$ is a parabolic subgroup $P_R$ of $G_R$. Using Teichm\"uller lifts, we get: its Lie algebra is $F^0({\rm Lie}(G_R))$. This ends the proof.
\medskip
Let now $\Mg$ be the reductive group over $X^\wedge$ obtained by gluing naturally these reductive groups $G_R$. Let $\Mp$ be the parabolic subgroup of $\Mg$ obtained by gluing naturally these parabolic subgroups $P_R$. Here reductive group (resp. parabolic subgroup) is understood in the following sense: $\Mg$ (resp. $\Mp$) is a $p$-adic formal group scheme (resp. formal subgroup scheme of $\Mg$) over $X^\wedge$ such that $\Mg_{W_m(k)}$ (resp. $\Mp_{W_m(k)}$) is a reductive group scheme over $X_{W_m(k)}$ (resp. is a parabolic subgroup of $\Mg_{W_m(k)}$), $\forall m\in\NN$. 
\medskip 
{\bf Lemma.} {\it Locally in the \'etale topology of ${\rm Spec}(R)$ the filtration $(F^i(M_R))_{i\in S(a,b)}$ of $M_R$ is defined (as in 2.2.8 c')) by a cocharacter $\mu_R$ of $G_R$.}
\medskip
{\bf Proof:} We can assume $P_R$ is split, cf. [SGA3, Vol. III, 6.1 of p. 32]. Let $m_1:{\rm Spec}(W(k_1))\to X$ be a Teichm\"uller lift whose special fibre is a dominant geometric point of ${\rm Spec}(R/pR)$; so $m_1$ is as well a dominant morphisms. The filtration of $M_R\otimes_R W(k_1)$ is defined (as in 2.2.8 c')) by a cocharacter of $\tilde P_{W(k_1)}$. As we can replace such a cocharacter by any $P(W(k_1))$-conjugate of it, we can assume that it is obtained from a cocharacter $\mu_R$ of $P_R$ by pull back. If $M_R=\oplus_{i\in S(a,b)} \tilde F_R^i$ is the direct sum decomposition defined (as in 2.2.8 c')) by $\mu_R$, then we have $F_R^i=\oplus_{j=i}^b\tilde F^j(M_R)$, as this holds after tensorization with $W(k_1)$. This ends the proof.
\medskip
{\bf Proposition.} {\it There is an $\NN$-pro-\'etale morphism $f_1:{\rm Spec}(R_1)\to {\rm Spec}(R)$, surjective mod $p$, and such that there is a connection $\dl_{R_1}$ on $M_R\otimes_{R^\wedge} R_1^\wedge$ annihilating $t_{\al}$, $\forall\al\in\Mj$.} 
\medskip
{\bf Proof:} We can assume $\mu_R$ exists. We write $\Phi_{M_R}=g_R\mu_R({1\over p})$, with $g_R$ a $\Phi_U$-linear automorphism of $M_R$ fixing $t_{\al}$, $\forall\al\in\Mj$. We can assume $M_R$ is a free $R^\wedge$-module. The equation
$$\Phi_{M_R}(x)=x,\leqno (EQ)$$
with $x\in M_R$, mod $p$ defines a system of equations of first type with coefficients in $R/pR$ whose matrix is invertible. By passing (as usual) from things mod $p^m$ to things mod $p^{m+1}$, $m\in\NN$, we come across again on a system of equations of first type whose matrix is not changed. So from end of 3.6.8.1.0 and from 3.6.8.1.2 a) (to be compared with the proof of 3.6.1.3) we deduce the existence of $f_1$ of the type mentioned such that the $\ZZ_p$-module $M_{\ZZ_p}$ of solutions of (EQ) in $M_R\otimes_{R^\wedge} R_1^\wedge$ is such that the natural $R_1^\wedge$-linear map $M_{\ZZ_p}\otimes_{\ZZ_p} R_1^\wedge\to M_R\otimes_{R^\wedge} R_1^\wedge$ is an isomorphism. So $M_{\ZZ_p}$ defines a $\ZZ_p$-structure of $M_R\otimes_{R^\wedge} R_1^\wedge$; with respect to it we have $t_{\al}\in\Mt(M_{\ZZ_p})$, $\forall\in\Mj$. So we can take as $\dl_{R_1}$ the connection on $M_R\otimes_{R^\wedge} R_1^\wedge$ annihilating $M_{\ZZ_p}$. This ends the proof.
\medskip
{\bf B. Definitions.} We say a connection $\nabla$ on $\Mf$ respects the $\Mg$-action, if $\nabla(t_{\al})=0$, $\forall\al\in\Mj$. We say a connection on $\Mf/p^n\Mf$ making ${\got C}/p^n{\got C}$ potentially to be viewed as an object of $\Mm\Mf_{[a,b]}^\nabla(X)$ respects the $\Mg$-action, if locally in the $\NN$-pro-\'etale topology of $X$, it lifts to a connection on $\Mf$ which makes ${\got C}$ potentially to be viewed as a $p$-divisible object of $\Mm\Mf_{[a,b]}^\nabla(X)$ and which respects the $\Mg$-action. Similarly, we define connections respecting the $\Mg$-action, for pull backs of ${\got C}$ or of ${\got C}/p^n{\got C}$.
\medskip
{\bf C. Theorem (the relative form of the moduli principle).} {\it We assume the filtration of ${\rm Lie}(G_R)$ is in the range $[-1,1]$. Then there is a moduli $p$-adic formal scheme $M_{\Mg}({\got C}/p^n{\got C})$ (resp. $M_{\Mg}(\got C)$) of connections on pull backs of $\Mf/p^n\Mf$ (resp. of $\Mf$) through formally \'etale, $p$-adic formal schemes over $X^\wedge$, which make the extension of ${\got C}/p^n{\got C}$ (resp. of $\got C$) to such a $p$-adic formal scheme $X_1$ over $X^\wedge$ potentially to be viewed as an object (resp. as a $p$-divisible object) of $\Mm\Mf^{\nabla}_{[0,1]}(X_1)$ and which respect the $\Mg$-action. 
Moreover, we have:
\medskip
{\bf a)} $\forall m\in\NN$, the $X_{W_m(k)}$-scheme $M_{\Mg}({\got C}/p^n{\got C})_{W_m(k)}$ (resp. $M_{\Mg}({\got C})_{W_m(k)}$) is \'etale and affine (resp. is $\NN$-pro-\'etale and affine);
\smallskip
{\bf b)} Locally in the \'etale topology of $X_k$, the special fibre of $M_{\Mg}({\got C}/p^n{\got C})$ (resp. of $M_{\Mg}({\got C})$) is obtained through a finite sequence (resp. an infinite sequence indexed by $\NN$) of systems of equations of first type defined using linear homogeneous forms with coefficients in (the ring of global sections of the original --i.e. before localization--) $X_k$;
\smallskip
{\bf c)} The fibres of ${M_{\Mg}({\got C}/p^n{\got C})}_k$ (resp. of ${M_{\Mg}({\got C})}_k$) over points with values in fields of $X_k$ are non-empty. The fibres of ${M_{\Mg}({\got C}/p^n{\got C})}_k$ are \'etale schemes over fields defined by algebras of dimension (as vector spaces over these fields) a non-negative, integral power of $p$.}
\medskip
{\bf Proof:} Let $m:X_1\to X$ be a formally \'etale morphism, with $X_1$ a $p$-adic formal scheme. We first assume that $(a,b)=(0,1)$. The fact that a connection on $\Mm:=(m^\wedge)*(\Mf)$ respects the $\Mg$-action, is expressed through some algebraic equations. So as $M_{\Mg}(\got C)$ (resp. as $M_{\Mg}({\got C}/p^n{\got C})$) we can take the maximal open closed, formal subscheme of the moduli $p$-adic formal scheme $M(\got C)$ (resp. $M({\got C}/p^n{\got C})$) constructed in 3.6.18.4.2, over which we get connections respecting the $\Mg$-action (resp. which is the image of $M_{\Mg}(\got C)$ in $M({\got C}/p^n{\got C})$ under the natural $X^\wedge$-morphism $M({\got C})\to M({\got C}/p^n{\got C})$). This takes care of the existence part, as well as of a). 
\smallskip
b) follows from 3.6.18.7: in the mentioned place we used, cf. 3.6.8.6.1 (and end of 3.6.18.7.1), the fact that the direct sum decomposition $M=F^1\oplus F^0$ is defined (as in 2.2.8 c)) by a cocharacter of $\tilde G$ and that $\dl_0$ exists (i.e. and that we can introduce $\be$). But working in the $\NN$-pro-\'etale topology of $X^\wedge$, the Lemma and the Proposition of A take care of the existence of a corresponding cocharacter $\mu_R$ of $G_R$ and respectively of a connection $\dl_{R_1}$ substituting $\dl_0$. 
\smallskip 
c) results from b) (cf. 3.6.8.1 applied in the context of geometric points of $X_k$). This handles the case $(a,b)=(0,1)$. 
\smallskip
The general case is entirely the same. The only difference comes from the fact that we can not appeal to 3.6.18.4.2 and so we have to argue slightly differently the existence of $M_{\Mg}({\got C}/p^n{\got C})$ (resp. of $M_{\Mg}({\got C})$). Due to the universal property implicit in the moduli aspect of $M_{\Mg}({\got C}/p^n{\got C})$ (resp. of $M_{\Mg}({\got C})$) and due to the affineness part of a), using descent we can work locally in the \'etale (resp. $\NN$-pro-\'etale) topology of $X$. So we can assume that $X={\rm Spec}(R)$ is affine, that the filtration of $\Mf$ is defined by cocharacters of $\Mg$, that $\Om_{X_k/k}$ is a free $\Mo_{X_k}$-sheaf, that $\Mf$ is a free $\Mo_{X^\wedge}$-sheaf and that we have --to be compared with the last proof of A-- a $\ZZ_p/p^n\ZZ_p$-structure (resp. $\ZZ_p$-structure) of $\Mf/p^n\Mf$ (resp. of $\Mf$) formed by elements fixed by $\Phi_{M_R}$. Once we have this, the proof of 3.6.18.7.1 c) implies that we can construct $M_{\Mg}({\got C}/p^n{\got C})$ (resp. of $M_{\Mg}({\got C})$) following the pattern of 3.6.18.4.2:
\medskip
-- defining $M_{\Mg}({\got C}/{\got C})$ to be $X$, $M_{\Mg}({\got C}/p^n{\got C})_k$ is obtained from $M_{\Mg}({\got C}/p^{n-1}{\got C})_k$ using a quasi Artin--Schreier system of equations;  
\smallskip
-- $M_{\Mg}({\got C}/p^n{\got C})$ is the $p$-adic completion of the affine, \'etale $M_{\Mg}({\got C}/p^{n-1}{\got C})$-scheme whose special fibre is the $M_{\Mg}({\got C}/p^{n-1}{\got C})_k$-scheme $M_{\Mg}({\got C}/p^n{\got C})_k$;
\smallskip
--  $M_{\Mg}({\got C})$ is the $p$-adic completion of the $\NN$-projective limit of $M_{\Mg}({\got C}/p^n{\got C})$'s.
\medskip
This ends the proof. 
\medskip
{\bf 3.6.18.7.4. Comments.} {\bf 1)} We do not know how to compute the number of connected components of $M_{\Mg}({\got C}/p^n{\got C})$.
\smallskip
{\bf 2)} We do not know when in 3.6.18.7.3 b) we can work just in the Zariski topology of $X_k$.
\smallskip
{\bf 3)} We do not know when $M_{\Mg}({\got C}/p^n{\got C})$ is $X$-algebraizable, i.e. is obtained from an \'etale $X$-scheme $X_1$ by taking the $p$-adic completion. If such an $X_1$ exists, then we have a universal property for the morphism $X_1\to X$ involving formally \'etale $X$-schemes, similar to the one of iii) of 3.6.1.3 1). 
\smallskip
{\bf 4)} We have variants of 3.6.18.7.3 (and of 1) to 3)): for instance, we allow $X$ to be a $p$-adic formal scheme or we work in the context of 3.6.18.7.0 or of 3.6.18.7.1 a) (the existence of $\Mp$ in the last two contexts can be deduced as well from [Fa2, rm. iii) after th. 10]). 
\medskip
Till end of 3.6.18.10 we work with a regular, formally smooth, affine $W(k)$-scheme ${\rm Spec}(R_1)$ such that $Y_1:={\rm Spec}(R_1/pR_1)$ is a connected scheme. We fix a Frobenius lift $\Phi_{R_1}$ of $R_1^\wedge$. Let $a$, $b\in\ZZ$, with $a\le b$. 3.6.18.8-10 below can be stated as well in a non a priori affine context; however, as their essence ``boils down" to an affine context, they are stated just in this last context. The following two theorems are consequences of 3.6.18.7.1.1 b). However, for the sake of future references and due to their importance, they are as well stated separately, with arguments.
\medskip
{\bf 3.6.18.8. Theorem (the universal uniqueness principle).} {\it We assume there is a maximal ideal $m_{R_1}$ of $R_1$ such that $R_1/m_{R_1}$ is a finite field and the resulting Frobenius lift of the completion $\widehat{R_1}$ of $R_1$ w.r.t. is of essentially additive type. Then for any ($p$-divisible) object ${\got C}_1$ of $\Mm\Mf_{[a,b]}(R_1)$, there is at most one connection on it. So $a-p-\Mm\Mf_{[a,b]}^\nabla(R_1)$ is a full subcategory of $a-p-\Mm\Mf_{[a,b]}(R_1)$.}
\medskip
The proof of the first part of the Theorem is the same as the proof of 3.6.1.2 (cf. also the extra features of the proof of 3.6.18.4.2): we work inductively modulo $p^n$, $n\in\NN$. The difference from 3.6.8 6) is that we need to apply the same upper triangular trick used in 3.6.18.4 P2 (for the resulting Frobenius lift of $\widehat{R_1}$) for ``separating" the variables. The case of $p$-divisible objects is also a consequence of Corollary of 3.6.18.5.5 B. As in the proof of 3.6.18.5.1 we get that $a-p-\Mm\Mf_{[a,b]}^\nabla(R_1)$ is a full subcategory of $a-p-\Mm\Mf_{[a,b]}(R_1)$.
\medskip
{\bf 3.6.18.8.1. Theorem (the universal slope type uniqueness principle).} {\it Let $m_1:{\rm Spec}(W(k_1))\to {\rm Spec}(R_1)$ be a $\Phi_{R_1}$-Teichm\"uller lift, with $k_1$ an algebraic field extension of $k$. 
\smallskip
{\bf a)} Let ${\got C}_1$ be a $p$-divisible object of $\Mm\Mf_{[0,1]}(R_1)$ such that $m_1^*({\got C}_1)$ does not have slope $0$ or slope $1$. Then for any object (or $p$-divisible object) ${\got C}_1^\prime$ of $\Mm\Mf_{[0,1]}(R_1)$ whose truncation mod $p$ is isomorphic to ${\got C}_1/p{\got C}_1$, there is at most one connection on it.
\smallskip
{\bf b)} Let ${\got C}_1$ be a $p$-divisible object of $\Mm\Mf_{[a,b]}(R_1)$ such that $End(m_1^*({\got C}_1))$ has the pseudo-multiplicity of the slope $-1$ equal to $0$. Then for any object (or $p$-divisible object) ${\got C}_1^\prime$ of $\Mm\Mf_{[a,b]}(R_1)$ whose truncation mod $p$ is isomorphic to ${\got C}_1/p{\got C}_1$, there is at most one connection on it.}
\medskip
{\bf Proof:} The part of a) involving $p$-divisible objects is a consequence of the formula of 3.6.18.4 B). It is convenient to view the part of a) referring to objects as a particular case of b). Now b) is a consequence of 3.6.8.9 (cf. 3.6.18.5.5) and of the proofs of 3.6.18.4 B) and of 3.6.18.4.2. In other words, based on 3.6.18.5.5 D, the proofs of 3.6.18.4 B) and of 3.6.18.4.2 can be entirely adapted to the context of ${\got C}_1$. We just need to remark: the fact that $End(m_1^*({\got C}_1))$ has the pseudo-multiplicity of the slope $-1$ equal to $0$, can be read out from its truncation $\Mt$ mod $p$; so what we know for $\Mt$ remains true for all its subobjects. This proves the Theorem.  
\medskip
{\bf 3.6.18.8.2. Examples.}
We refer to 3.6.18.6. As usual we can speak about the $p$-rank of ${\got C}_1/p{\got C}_1$ at any maximal (or geometric) point of ${\rm Spec}(R_1/pR_1)$. The $p$-rank is $0$ in such a point $y$ iff the $F$-crystal we get (by pulling back ${\got C}_1$ via a Teichm\"uller lift of ${\rm Spec}(R_1)$ whose special fibre factors through $y$) does not have slope $0$. So if ${\got C}_1$ has $p$-rank $0$ w.r.t. any maximal point of ${\rm Spec}(R_1/pR_1)$, then in the construction of $M({\got C}_1)_k$ we do not need to pass to $\NN$-pro-\'etale covers (or in the case of an object to \'etale covers): we have $M_0({\got C}_1)=M({\got C}_1)={\rm Spf}(R_1^\wedge)$, i.e. $d_1$ mod $p$ is an isomorphism (cf. 3.6.18.5.8). 
\smallskip
Similarly, we can ``detect" the multiplicity of the slope $1$ of ${\got C}$ at a geometric point of ${\rm Spec}(R_1/pR_1)$ by just looking at ${\got C}_1/p{\got C}_1$: it is equal to the $p$-rank of $({\got C}/p{\got C})^*(1)$ at the point. So, similarly, if ${\got C}_1$ has slope $1$ with multiplicity $0$ at
 any maximal point of ${\rm Spec}(R_1/pR_1)$, then $M({\got C}_1)={\rm Spf}(R_1^\wedge)$. 
\medskip
{\bf 3.6.18.8.3. Some variants.} 3.6.1.2-3 remain true for any Frobenius lift $\Phi_{R_j}$ of $R_j^\wedge$ such that $\Phi_{\hat R_j^0}$ is of essentially additive type (cf. 3.6.18.4 B); see also 3.6.18.5.1 and 3.4.18.7.0). 
This gives us more freedom. To exemplify it, we assume we are in a context in which 3.6.8.4 5) applies, i.e. we automatically have ${\rm Spec}(Q_{0,n}/pQ_{0,n})=\Ms^n$, regardless of the choices involved; for instance, this is the case if the $G$-ordinary type defined by $(M,\vph,G)$ has no slope $0$ or no slope $1$, as in this case $\Ms^n={\rm Spec}(R_0/pR_0)$, $\forall n\in\NN$ (cf. 3.6.18.8.2). In such a context, we can always adjust the things (in a way similar to the algebraization process of the proof of 3.6.18.4.1) so that
in 3.6.2 we ``end up" over $k$, i.e. the $W(k)$-morphisms of 3.6.2.1 b) and c) do exist without the three assumptions of the beginning paragraph of 3.6.2.1 (that $k=\bar k$, etc.). 
\smallskip
For an arbitrary Frobenius lift $\Phi_{R_j}$ of $R_j^\wedge$, we lose the uniqueness part, so 3.6.1.2-3 have to be modified to a significant extend (for instance not all connections we get respect the $G$-action, $\ell_{j,n}^{-1}(\Mz_k)$ is not necessarily ${\rm Spec}(k)$, etc.). We do not stop to rewrite what remains true out of 3.6.1.2-3, in this general context. We just point out that we still have a universal property w.r.t. connections respecting the $G$-action (cf. 3.6.18.7.3 C), but the corresponding $p$-adically complete moduli $R_j^\wedge$-schemes, can (theoretically) have non-connected special fibres.
\smallskip
Moreover, 3.6.18.6 has a similar version, in which we work with an arbitrary Frobenius lift $\Phi_{R_1}$ and we are not in the context of a potential-deformation sheet. 
\medskip
{\bf 3.6.18.9. Global estimates.} We assume that for any reduced, closed subscheme of $Y_1$, its $\bar k$-valued points are dense in it. Let ${\got C}$ be an object of $\Mm\Mf_{[0,1]}(R_1)$. We denote by $\nabla({\got C})$ the number of connections on ${\got C}$. 
\smallskip
Let $y\in Y_1(\bar k)$. Let $r_y\in\NN\cup\{0\}$ be the rank (see 3.6.18.1) of the Frobenius lift of the completion of $R_1$ w.r.t. its maximal ideal defining the point of $Y_1$ through which $y$ factors, defined naturally by $\Phi_{R_1}$. 
\medskip
{\bf Claim.} {\it The values of $r_y$'s achieve a stratification $\Ms$ of $Y_k$ in reduced, locally closed subschemes in such a way that under specialization of strata, their values decrease. $\Ms$ satisfies the purity property.}
\medskip
{\bf Proof:} We consider a $p$-divisible group over $R_1$ which is a direct sum of $\QQ_p/\ZZ_p$ and of its dual. We apply 3.6.18.4.2 in the context of the $p$-divisible object of $\Mm\Mf_{[0,1]}^\nabla(R_1)$ associated to it. So the Claim follows from 3.6.18.4 B) and 3.6.8.1.3-4. 
\medskip
As in 3.6.18.1 we associate to ${\got C}/p{\got C}$ (and $y$) numbers $s_y(0),s_y(1)\in\NN\cup\{0\}$. Let $n\in\NN\cup\{0\}$ be the smallest number such that $p^n$ annihilates ${\got C}$. Let 
$$d_y({\got C}):=nr_ys_y(0)s_y(1).$$ 
From 3.6.18.4 B) and the extra features of the proof of 3.6.18.4.2 we get:
$$\nabla({\got C})\le min\{d_y({\got C})|y\in Y_1(\bar k)\}.\leqno (EST1)$$ 
We have:  
\medskip
{\bf Fact.} {\it If the number $r_ys_y(0)s_y(1)$ does not depend on $y\in Y_1(\bar k)$ and if the underlying module of ${\got C}$ is a projective $R_1/p^nR_1$-module, then the special fibre of the moduli formal scheme $M({\got C})$ (of 3.6.18.4.2) is an \'etale cover of $Y_1$. Moreover, if $r_ys_y(0)s_y(1)\neq 0$, then each one of the numbers $r_y$, $s_y(0)$ and respectively $s_y(1)$ independently do not depend on $y\in Y_1(\bar k)$.}
\medskip
{\bf Proof:} The first part follows from 3.6.18.4 B), 3.6.18.4.2, 3.6.8.1.2 a) and standard properties on the number of points of geometric fibres of an \'etale, affine morphism (see [EGA IV, 15.5.9]). The second part follows using specialization and generic arguments (cf. the above Claim for $r_y$'s). This ends the proof.
\medskip
Let now ${\got C}$ be an object of $\Mm\Mf_{[a,b]}(R_1)$. Let $n$ be as above. Let $s_{pm}^y(-1,{\got C})$ be the pseudo-multiplicity of the slope $-1$ of the following object of $\Mm\Mf_{[a,b]}(R_1)$ 
$$\sum_{i=0}^{n-1} p^iEnd({\got C})/p^{i+1}End({\got C}).$$
From 3.6.18.7.1.1 b) we get 
$${\rm log}_p(\nabla({\got C}))\le min\{r_ys_{pm}^y(-1,{\got C})|y\in Y_1(\bar k)\}.\leqno (EST2)$$
From 3.6.18.7.1.1 a) and loc. cit. we also get:
\medskip
{\bf Corollary.} {\it We assume $b=a+1=1$. The products $r_ys_{pm}^y(-1,{\got C})$ do not depend on $y\in Y_1(\bar k)$ iff the special fibre of the moduli formal scheme $M({\got C})$ (of 3.6.18.4.2) is an \'etale cover of $Y_1$.}
\medskip
{\bf 3.6.18.10. $p$-divisible objects of bounded $\nabla$-deviations.}
We assume $Y_1$ has the $ALP$ property. Let ${\got C}$ be a $p$-divisible object of $\Mm\Mf_{[0,1]}(R_1)$. We say ${\got C}$ is of bounded (resp. of totally bounded) $\nabla$-deviation if at least one connected component (resp. if all connected components) of the special fibre of $M({\got C})$ is (resp. are) of finite type over $Y_1$. The interesting cases are when $Y_1$ is of finite type over $k$ or when $R_1=W(k)[[z_1,...,z_m]]$, with $k\neq\bar k$. Plenty of examples of $p$-divisible objects of  bounded $\nabla$-deviation can be obtained following the strategy of modifications of 3.6.8.9.1. Similarly, using 3.6.18.8-9, we get plenty of examples of $p$-divisible objects of totally bounded $\nabla$-deviation; for instance, ${\got C}$ is of totally bounded deviation if ${\got C}/p{\got C}$ has $p$-rank $0$ at all $\bar k$-valued point of $Y_1$. Two problems are dear to us.
\medskip
{\bf P1.} In the interesting cases mentioned above, determine all $p$-divisible objects of bounded $\nabla$-deviation.
\smallskip
{\bf P2.} Determine conditions under which, in the case when $R_1$ is a smooth $W(k)$-algebra, a $p$-divisible object of $\Mm\Mf_{[0,1]}^\nabla(R_1)$ is such that:
\smallskip
{\item {a)}} it is associated to a $p$-divisible group over $R_1$ (and not only --see 2.2.1.1 2)-- to one over $R_1^\wedge$), or
\smallskip
{\item {b)}} when viewed as an object of $\Mm\Mf_{[0,1]}(R_1)$, it is of totally bounded $\nabla$-deviation. 
\medskip
{\bf 3.6.19. The gluing principle.} 3.6.14 and 3.6.14.1 are particular cases of a more general principle pertaining to Shimura $p$-divisible groups, which we call the gluing principle. However, we singled them out, as they are in the context of abelian varieties, are the most important cases needed in applications, and moreover we had to give a preliminary motivation for the Expectation of 3.6.15 A. In what follows we present different forms of the gluing principle.
\smallskip
{\bf A.} Let ${\rm Spec}(R^{i\wedge})$ be the $p$-adic completion of a pro-\'etale, affine scheme ${\rm Spec}(R^i)$ over a smooth, affine $W(k)$-scheme, $i\in\overline{1,2}$. Let $z^i:{\rm Spec}(W(k))\to {\rm Spec}(R^i)$, and let 
$$
\Md^i=(D^i,(t^i_{\al})_{\al\in\Mj})
$$ 
be a Shimura $p$-divisible group over ${\rm Spec}(R^{i\wedge})$. We still denote by $z^i$ the $p$-adic completion of $z^i$. We assume that:
\medskip
i) the two Shimura $p$-divisible groups $z_1^*(\Md^1)$ and $z_2^*(\Md^2)$ over $W(k)$ are isomorphic;
\smallskip
ii) the relative dimension $d_1$ of ${\rm Spec}(R^1)$ in $z^1$ is greater or equal to the relative dimension $d_2$ of ${\rm Spec}(R^2)$ in $z^2$; moreover $\Md^i$ is a versal deformation in $z^i$, $i=\overline{1,2}$;
\smallskip
iii) ${\rm Spec}(R^i/pR^i)$ is a connected $k$-scheme, $i=\overline{1,2}$.
\medskip
 Here the versality condition of ii) is defined as follows.
\medskip
{\bf B. Definitions.} 
Let $X$ be a regular, formally smooth $W(k)$-scheme. A Shimura $p$-divisible group $\Md_X$ over $X$ is said to be a versal (resp. uni plus versal) deformation in a closed point $y_X:{\rm Spec}(k_1)\hookrightarrow X$, with $k_1$ an algebraic field extension of $k$, if the dimension of the image of the Kodaira--Spencer map attached to $\Md_X$ (and so computed relative to $W(k)$) in this point (resp. if this map is injective and the mentioned dimension) is equal to the deformation dimension of the (non-necessarily quasi-split) Shimura $\sg_{k_1}$-crystal of $y_X^*(\Md_X)$. 
\smallskip
We assume now that $X_k$ has the $ALP$ property. We say $\Md_X$ is a versal (resp. uni plus versal) deformation if it is a versal (resp. uni plus versal) deformation in all closed points $y_X$ as above. 
\smallskip
An object or a $p$-divisible object of $\Mm\Mf^\nabla(X)$ is said to be a uni plus quasi-versal deformation in a closed point $y_X$ as above (resp., assuming that $X_k$ has the $ALP$ property, is said to be a uni plus quasi-versal deformation), if the Kodaira--Spencer map in $y_X$ (resp. in each point $y_X$ as above) is injective. 
\smallskip
Similarly, assuming that $X_k$ has the $ALP$ property, an object or a $p$-divisible object of $\Mm\Mf^\nabla(X)$ is said to be a quasi-versal deformation if the Kodaira--Spencer map in each point $y_X$ as above has an image whose dimension depends only on the connected component of $X_k$ to which $y_X$ belongs. Similarly, we define Shimura $p$-divisible groups over $X$ which are (uni plus) quasi-versal deformations. 
\medskip
{\bf C.} If $\Md_i$ has a principal quasi-polarization $\Mp_{\Md_i}$, $i=\overline{1,2}$, then we say that we are in a principally quasi-polarized context; if we are in such a context then we assume that the isomorphism of i) respects the principal quasi-polarizations involved.
\smallskip
{\bf D.} Let $(M,F^1,\vph,G)$ be the Shimura filtered $\sg$-crystal associated to the two Shimura $p$-divisible groups over $W(k)$ obtained in i). Let ${\rm Spec}(\widehat{R^1})$ (resp. ${\rm Spec}(R^{1h})$) be the completion of ${\rm Spec}(R^1)$ (resp. be the henselization of ${\rm Spec}(R^1)$) in $z^1$. We still denote by $z^i$ the resulting $W(k)$-valued point of ${\rm Spec}(\widehat{R^i})$ or of ${\rm Spec}(R^{ih})$. We assume the existence of a formally smooth $W(k)$-morphism 
$$f_{12}:{\rm Spec}(\widehat{R^1})\to {\rm Spec}(R^{2\wedge})$$ 
such that:
\medskip
\item{R1)} $f_{12}\circ z^1=z^2$;
\smallskip
\item{R2)} the image of $R^2/pR^2$ in $\widehat{R^1}/p\widehat{R^1}$ is contained in $R^{1h}/pR^{1h}$;
\smallskip
\item{R3)} $f_{12}^*(\Md^2)$ is isomorphic to the pull back of $\Md^1$ through the natural morphism ${\rm Spec}(\widehat{R^1})\to {\rm Spec}(R^{1\wedge})$, through an isomorphism lifting the one of i). 
\medskip
We also assume that:
\medskip
\item{EXTRA)} For any formally smooth $W(k)$-morphism $\tilde f_{12}:{\rm Spec}(\tilde R^{1\wedge})\to {\rm Spec}(R^{2\wedge})$, with ${\rm Spec}(\tilde R^1)$ a pro-\'etale, affine scheme over ${\rm Spec}(R_1)$ such that $z_1$ has a lift $\tilde z_1$ with the properties that $\tilde f_{12}\circ\tilde z_1=z_2$ and the objects with tensors of $p-\Mm^\nabla_{[0,1]}(\widehat{R^1})$ defined by $\tilde f_{12}^*(\Md^2)$ and $\Md^1_{{\rm Spec}(\tilde R^{1\wedge})}$ are isomorphic under an isomorphism lifting (in $\tilde z_1$) the one of i) but viewed in the non-filtered context, after replacing $\tilde R_1$ by an \'etale, affine scheme over it to which $\tilde z_1$ lifts we can modify $\tilde f_{12}$ by something which is $0$ mod $p$ so that $\tilde f_{12}^*(\Md^2))$ and $\Md^1_{{\rm Spec}(\tilde R^{1\wedge})}$ themselves are isomorphic under an isomorphism lifting (in $\tilde z_1$) the one of i).
\medskip
Following the proofs of 2.3.15 and 3.6.14.1 we can state (cf. also 3.6.2.2):
\medskip
{\bf Theorem (the gluing principle: the global form).} {\it There is an affine $W(k)$-scheme ${\rm Spec}(R^3)$, with $R^{3\wedge}=R^3$, and there is a formally \'etale $W(k)$-morphism $q_1:{\rm Spec}(R^3)\to {\rm Spec}(R^{1\wedge})$ and a formally smooth $W(k)$-morphism $q_2:{\rm Spec}(R^3)\to {\rm Spec}(R^{2\wedge})$, such that for a suitable $W(k)$-morphism $a^3:{\rm Spec}(W(k))\to {\rm Spec}(R^3)$ we have:
\medskip  
a) $q_i\circ a^3=z^i$, $i=\overline{1,2}$;
\smallskip
b) $q_1^*(\Md^1)$ is isomorphic to $q_2^*(\Md^2)$ through an isomorphism lifting (cf. a)) the one of i);
\smallskip
c) ${\rm Spec}(R^3/pR^3)$ is connected and a pro-\'etale scheme over a smooth, affine $k$-scheme.
\medskip
If $R^i$ is a smooth $W(k)$-algebra, $i=\overline{1,2}$,
then we can also assume that:
\medskip
d) the image of ${\rm Spec}(R^3/pR^3)$ in ${\rm Spec}(R^{3-i}/pR^{3-i})$ contains an open, dense subscheme of it.}
\medskip
{\bf Proof:} The proof of this Theorem is contained in the proof 3.6.14.1. We need to add just two extra things. First, [BLR, th. 12 of p. 83] applies entirely (as in 3.6.14.1 E), due to standard arguments on projective systems as in [EGA IV, \S 8]; one can quote as well b) of the following general Fact (for $X={\rm Spec}(R^i)$).
\medskip
{\bf Fact.} {\it {\bf a)} Let $Y$ be a quasi-compact, quasi-separated, regular, formally smooth $W(k)$-scheme. If $X$ is a projective limit of \'etale, affine $Y$-schemes $X_{\al}$, $\al\in\Mj(X)$, then the category $\Mm\Mf^{\nabla(p+tens)}(X)$ is the projective limit of the categories $\Mm\Mf^{\nabla(p+tens)}(X_{\al})$ (see 2.2.4 C and D) for the definition of these categories). If all these schemes are endowed with Frobenius lifts compatible with the projective system, then the similar thing can be said about the categories $\Mm\Mf(X)$ and $\Mm\Mf^\nabla(X)$.
\medskip
{\bf b)} If in a), $Y$ is moreover a smooth $W(k)$-scheme, then $p-FF(X^\wedge)$ (resp. $p-DG(X^\wedge)$) is antiequivalent (via the $\DD$ functor) to $\Mm\Mf_{[0,1]}^\nabla(X)$ (resp. to $p-\Mm\Mf_{[0,1]}^\nabla(X)$).}  
\medskip
The proof of a) is standard (see [EGA IV, 8.5.2]), based on [Fa1, 2.1 (ii)] and on the fact that any finitely generated projective module over a ring is finitely presented; for the $\nabla$ context we need to add: as $\Om_{{X_{\al}}_k/k}$ is a locally free $\Mo_{X_{\al}}$-sheaf of finite rank, $\forall\al\in\Mj(X)$, compatible with pull backs under special fibres of transition morphisms of the projective system, the connections involved can be viewed, locally in the Zariski topology of $Y$, as a finite sequence of linear maps between locally free sheaves of modules of finite ranks. 
\smallskip
b) follows from a) and from 2.2.1.1 2) (for the faithfulness part of $\DD$ over $X^\wedge$, cf. [BM, 4.1.1] applied over $X_k$; $X_k$ has locally in the Zariski topology a finite $p$-basis).
This Fact is the fourth place where we need $p\ge 3$.
\smallskip
The second thing we need to add is related to the fact that $d_i$ can be greater than the deformation dimension $d$ of the Shimura $p$-divisible groups of i). But the refinement part of the proof of 3.6.14.1 applies entirely in this situation (cf. EXTRA)). Using these two additions, the Theorem follows.  
\medskip
{\bf E. Variants.} There are many other variants of the gluing principle:
\medskip
iv) we work in a context involving Shimura filtered $F$-crystals instead of Shimura $p$-divisible groups;
\smallskip
v) we weaken the part of ii) involving the versality condition;
\smallskip
vi) we work in a principally quasi-polarized context;
\smallskip
vii) we work modulo a power of $p$.
\medskip
The variants iv) and vi) are straightforward, while in connection to vii) we refer to 3.6.14.4. We detail now what we mean by this v). Let $Q^i$ be the completion of $R^i$ in $z^i$. Let $\tilde\Md^i:=\Md^i_{{\rm Spec}(Q^i)}$. Let $Q^0$ be the regular $W(k)$-algebra of the universal deformation space (see 2.2.21 and 2.2.21.1) of the Shimura $p$-divisible $\Md(W(k))$ showing up in i). Corresponding to $\tilde\Md^i$ we naturally get $W(k)$-homomorphisms $q_i:Q^0\to Q^i$, $i=\overline{1,2}$. What we want to say in v): it is enough to assume that
\medskip
a) ${\rm ker}(q_1)={\rm ker}(q_2)$, that
\smallskip
b) over dominant geometric points ${\rm Spec}(\tilde k)\to {\rm Spec}(Q^i/pQ^i)$ we get naturally from $\tilde\Md^i$ (see Fact 3 of 2.2.10) Shimura-ordinary $\tilde\sg$-crystals, and that
\smallskip
c) the resulting $W(k)$-monomorphism $Q^0/{\rm ker}(q_i)\hookrightarrow Q^i$ ($i\in\{1,2\}$) is formally smooth. 
\medskip
We express the combination of a) to c) by: $\Md^1$ and $\Md^2$ unfold (in $z^1$ and respectively $z^2$) the same quasi-versal deformation of (the Shimura filtered $\sg$-crystal underlying) $\Md(W(k))$, which generically (in the special fibre) is Shimura-ordinary.    
\medskip
{\bf F.} We do not know how to check in practice that EXTRA) holds. In a context of polarized abelian varieties, EXTRA) is in general a consequence of Serre--Tate's deformation theory (even in the context of v)) and of [Va2, 4.1.5] (see 3.6.14.1 E). 
\medskip
{\bf 3.6.20. Final remarks for 3.6. 1)} We might wonder why we preferred to work in 3.6.1-14 with open subschemes of the group $G$ and not with open subschemes of the (abelian) unipotent subgroup of $G$ acting trivially on $F^0$ and on $M/F^1$. The reason is: the full generality is useful in other situations, gives us more flexibility and, as it can be easily checked, the particular case (of an abelian group) can be often obtained from the general case, by taking slices (cf. also \S 7).
\smallskip
Also one might wonder: Why global deformations? No doubt, from many points of view, local deformations offer equally useful information as the global ones. However, from the deeper point of view of the classification of $p$-divisible groups (or objects when appropriate) over arbitrary $\ZZ_{(p)}$-schemes, the global deformations represent a major tool. Also, from many simple points of view, they offer advantages. Just one sample: referring to 3.6.18.6, the Newton polygon stratification of $M_0({\got C}_1)_k$ defined by $d_1^*({\got C}_1)$, is always the pull back of a stratification of ${\rm Spec}(R_1/pR_1)$; this last stratification can be obtained significantly easier in the case when ${\rm Spec}(R_1/pR_1)$ is of finite type over $k$ (often we just have to use simple reductive group properties; see \S 9-10 for details and computations).  
\smallskip
{\bf 2)} We needed (3.1.8 and 3.4.14 below are afterthoughts) the $\nabla$ principle for Shimura $\sg$-crystals for proving 3.1.0 and Milne's conjecture (see 1.15) in their full generality (and not only in the context of a SHS $(f,L_{(p)},v)$, where it is not needed: 2.3.16 is enough, as the scheme $\Mn$ itself provides the needed global deformation; see the proof of 4.2 below). The other  principles (of 3.6.18.4-8) are just a natural extension of it. The deformation principle is a logical continuation of the $\nabla$-principle. Its importance stems from the fact that, when applies, we often can take slices (cf. \S 7; see also 3.6.8 8) and 3.6.8.3), and its versions of 3.6.14.4 are very convenient for different global and local computations, cf. \S 10. It also led us to 3.6.15 B, which is the very essence of the boundedness principle of 3.15.7 below. Also the rigidity property (see 3.6.1.4 5) and \S 7) is a useful tool. We will use it in \S 8, in connection to 1.15.4. The gluing principle and its natural extensions (see [Va6]) to the generalized Shimura context and to Fontaine categories $\Mm\Mf_{[-1,1]}(*)$, is our straightforward approach to prove the existence of integral canonical models of Shimura varieties of special type (cf. also [Va2, 3.2.7 8)], 3.6.1.6 and 4.14.2 below). No doubt we were influenced (and inspired) by [Fa1, 7.1] and [Fa2, th. 10].
\smallskip
{\bf 3)} In some parts of 3.6 we used [Fa2, th. 10 and the remarks after] but not in 3.6.18.4. It is easy to see that using 2.3.17.2, 3.6.18.4.1, 3.6.18.7.0 and [Fa1, 7.1] we obtain another (completely different) proof of a slightly stronger form of [Fa2, th. 10 and the remarks after] for $p\ge 3$. To detail this, we use the notations of 3.6.18.5.1 with $p\ge 3$; so we are dealing with a Frobenius lift of $R$ which is essentially of additive type. Based on the proof of the equivalence part of 3.6.18.5.1, any $p$-divisible object of $\Mm\Mf_{[0,1]}(R)$ is a $p$-divisible object of $\Mm\Mf_{[0,1]}^\nabla(R)$. So, as in the proof of 3.6.18.3.1 we can assume $\Phi_R(z_i)=z_i^p$, $i=\overline{1,m}$.
\medskip
{\bf The algebraization process.} {\it Let $q\in\NN$. We consider the ideal $I_R(q):=(z_1^q,...,z_m^q)$ of $R$. For any $p$-divisible object ${\got C}_0=({\got C}^0,\nabla_0)$ of $\Mm\Mf_{[0,1]}^\nabla(R)$, we can algebraize the things similarly to the proof of 3.6.18.4.1. First we consider a $p$-divisible object ${\got C}^{0\rm al}$ of $\Mm\Mf_{[0,1]}(R^{\rm al})$ which mod $I_R(q+2)$ is ${\got C}^0$ mod $I_R(q+2)$; here $R^{\rm al}$ is as in the proof of 3.6.18.4.1. Second, as in the proof of 3.6.18.4.1, we consider the moduli formal scheme $M_0({\got C}^{0\rm al})$ defined by ${\got C}^{0\rm al}$ and obtained via 3.6.18.4.6; it is algebraizable (cf. end of 3.6.18.4.6 a)), and so not to introduce extra notations, we view it as an affine $R^{\rm al}$-scheme (and not as a formal scheme) which is $p$-adically complete. 
\smallskip
From b) of the Fact of 3.6.19 (i.e. --implicitly-- from [Fa1, 7.1]) and from the $\NN$-pro-\'etale part of 3.6.18.4.6 a), we deduce that there is a $p$-divisible group over $M_0({\got C}^{0\rm al})$ whose associated $p$-divisible object of $\Mm\Mf_{[0,1]}^\nabla(M_0({\got C}^{0\rm al}))$ is the pull back of ${\got C}^{0\rm al}$ to $M_0({\got C}^{0\rm al})$ together with the unique connection (cf. 3.6.18.8) on it. 
\smallskip
But $\nabla_0$ mod $I_R(q)$ is uniquely determined by ${\got C}^0$ mod $I_R(q+2)$ (the argument is the same as the one of the proof of 3.6.1.2; see 3.6.8 6) for things mod $p$). As we have a uniquely determined $R^{\rm al}$-morphism ${\rm Spec}(R)\to M_0({\got C}^{0\rm al})$ (see 3.6.18.4.6 a)), we get that ${\got C}_0$ modulo $I_R(q)$ is associated to a $p$-divisible group $D_{q}$ over $R/I_R(q)$. $D_{q}$ is uniquely determined (for instance, cf. [BM, 4.3.2 (i)] and [Me, ch. 4-5]; this can be also checked starting from generalities on deformation theory as presented in [Fa2, \S 7] or from [dJ1, th. of intro.]). So $D_{q+1}$ mod $I_R(q)$ is $D_q$ and so (as $R$ is complete w.r.t. the topology defined by $I(q)$'s) we deduce the existence of a uniquely determined $p$-divisible group $D$ over $R$ such that $\DD(D)$ is ${\got C}_0$.}
\medskip
It is 2.3.17.2 (strictly speaking a variant of it adapted to the context of [Fa2, rm. iii) after th. 10]: see 2.4) which takes care of the mentioned context involving tensors. But the proofs in [Fa2, \S 7] (work for $p=2$ as well and) are of unsurpassable beauty. 
\smallskip
Also [Fa1, 7.1] can be deduced from [Fa2, th. 10], cf. 3.15.3 4) below. Moreover, the antiequivalence part of 3.6.18.5.1 (as well as of 3.14 B6 below) can be deduced from [dJ1, th. of intro.]; in fact it is easy to see (based on [Me, ch. 4-5]) that this mentioned result for the smooth context (of $R/pR$) and [Fa2, th. 10] are equivalent to each other, even for $p=2$ (see the part of 2.2.21 referring to some uniqueness). Warning: despite this equivalence, it is worth stating explicitly that the methods of [dJ1] can not be adapted to the generalized Shimura context, while the proof of [Fa2, th. 10] does (see 3.15.6 below). See \S 6 and [Va5-8] for general Dieudonn\'e theories in the context of Shimura $p$-divisible groups over very general base schemes (not necessarily regular).
\smallskip
As [Fa1, 7.1] is stated for $p\ge 3$, in 3.14 below we do not come back to it. 
\smallskip
{\bf 4)} There are variants of many of the above principles. For instance, we get plenty of variants by working with:
\medskip
-- almost $p$-divisible objects (see 2.2.1.7 1)), or with
\smallskip
-- finitely generated modules, which are not necessarily projective, or with 
\smallskip
-- (some) $W(k)$-algebras which are not necessarily regular, formally smooth. 
\medskip
We do not stop to state results in this generality, as anyone can state them easily, when needed.  We just point out that in the case of a $W(k)$-algebra $R$, with $R=R^\wedge$, we always have to:
\medskip
-- either work with connections which are w.r.t. an a priori specified locally free $R$-submodule of the $p$-adic completion of $\Om_{R/W(k)}$, taken by a fixed Frobenius lift of $R$ into itself, and then we have to deal with $\Phi_R$-linear endomorphisms of $R$-modules, which ``keep" us inside it, or to
\smallskip
-- work with ``pseudo-connections" as in the part of the proof of 3.6.18.4.1 referring to its Fact 3. 
\medskip
{\bf 5)} b) of the Fact of 3.6.19 applies if $X$ is the henselization of the localization of a polynomial $W(k)$-algebra $W(k)[x_1,...,x_m]$ w.r.t. its ideal $(x_1,...,x_m)$. If moreover $k=\bar k$, when we combine this with 3.6.18.5.4 1) and 3.6.18.7, we get a very elementary way to construct Shimura $p$-divisible groups over $X^\wedge$. See [Va11] for details and applications.
\smallskip
{\bf 6)} It seems to us (see also 3.6.18.10) that there is a great similarity between 3.6.1.3, the theory of $p$-jets of A. Buium (see [Bu1]) and the theory of almost \'etale extensions of G. Faltings (see [Fa1] and [Fa3]). We think it is possible to obtain significant progress in the related areas by ``fitting" these three theories together. 
\smallskip
{\bf 7)} We view the part $s=0$ of 3.6.8.1 and 3.6.8.1.2 c) as the purely algebraic analogue of Lang's theorem.
\smallskip
{\bf 8)} There are many other ways (besides 3.6.19 B) to express the versality (resp. uni plus versality) condition of 3.6.19 ii). One way is to use the proof of 3.12.1 below. Another way (more fashionable) is inspired from  2.3.17.2: any Shimura $p$-divisible group over $W(k)$, whose associated Shimura $\sg$-crystal over $k$ is (identifiable with) the one obtained from $\Md^i$ through the $k$-valued point $y^i$ of ${\rm Spec}(R^{i\wedge})$ defined by $z^i$, is induced from $\Md^i$ through a (resp. through a uniquely determined) $W(k)$-valued point of ${\rm Spec}(R^{i\wedge})$ lifting $y^i$. The fact that all these definitions are equivalent is proved in the same manner as for the case of a SHS (see 2.3.17.2 and [Va2, 5.4-5]; see also 2.4). This is the fifth place where we need $p\ge 3$. 
\smallskip
{\bf 9)} In [Va6] we will use the gluing principle (see 3.6.2.2 and 3.6.19 D) to construct (see also \S 6 for the first steps; see $\Ml_{r,d-r}$ of 1.12 as a first sample) the theory of Shimura envelopes. This theory generalizes the theory of (quotients of) integral canonical models of Shimura varieties.
\smallskip
{\bf 10)} In positive characteristic, regardless of the fact that we are interested in some fundamental groups (see 3.6.18.4.6 C), or on some connections (see the whole of 3.6.8 and of 3.6.18), or on some stratifications (see 3.6.8.1.4, the Claim of 3.6.18.9 and 4.5.18 below), or on some deformations (see 3.6.14), or on some homomorphisms between two $p$-divisible groups (see the proofs of 3.6.14.1 and of 3.6.17), the quasi Artin--Schreier systems of equations are at the very root of the subject.
\smallskip
{\bf 11)} In 3.6.18 we fixed $\Phi$ and $\Phi_R$ and varied the connections, subject to the equations $(E_1)$ and $(E_2)$ of 3.6.1.1.1 2). It would be interesting (though not so useful) to study what happens if we fix an integrable, nilpotent mod $p$ connection and one of the two $\Phi_R$ and $\Phi$, while allowing the other one to vary. The problem is: the connections involve two few variables (to be compared with 3.6.18.2), to truly determine (up to \'etale extensions) what has been allowed to vary.
\vfill\eject
{\bf 3.7. The proofs of 3.2.6-7.}
\medskip
{\bf 3.7.1. Proof of 3.2.6.} Let $g_2\in G(W(k))$. We use the Shimura $\sg_{k_1}$-crystal
${\got C}_1$ considered in 3.6.7.1. We recall that ${\got C}_1$ specializes to $(M\otimes_{W(k)} W(\bar k),g_3(\vph_0\otimes 1),G_{W(\bar k)})$, with $g_3\in G_{W(\bar k)}(W(\bar k))$ such that mod $p$ belongs to the subset $\Ml(\bar k)$ of $G(\bar k)$ defined as the set $\Ml$ of 3.6.6 but for $\bar k$. We deduce from 3.4.11-13 that the Shimura Lie $\sg_{k_1}$-crystal $Lie({\got C}_1)$ 
attached to ${\got C}_1$ has the same Newton polygon as
$({\got g},\vph_0)$ and so as $({\got g},\vph_1)$ (cf. also 3.3.1 and 3.3.4). From the fact that ${\got C}_1$
specializes to $(M,g_2\vph_0)$ (cf. 3.6.7.1), we deduce
that the Newton polygon of $({\got g},g_2\vph_0)$ is above the Newton polygon of 
$({\got g},\vph_0)$. Similarly,
working just with ${\got g}_0$ instead of ${\got g}$ itself, we get that the Newton polygon of $({\got g}_0,g_2\vph_0)$ is above the Newton polygon of 
$({\got g}_0,\vph_0)$. 
\smallskip
As $Lie({\got C}_1)$ and 
$({\got g},\vph_1)$ have the same Newton polygon, we get (cf. 3.4.11 applied over $k_1$) that the formal isogeny type of ${\got C}_1$ is $\tau_0$. As ${\got C}_1$ specializes to $(M,g_2\vph_0,G)$ and
$g_2\in G(W(k))$ was arbitrary, we get 3.2.6.
\medskip
{\bf 3.7.2. Proof of 3.2.7 a).} 3.2.7 a) results directly from 3.3.2 and 3.4.3.0. So (cf. also 3.2.2), in the end of 3.5.5, the expression ``of Borel type" can be substituted by ``of parabolic type".
\medskip
{\bf 3.7.3. Proof of 3.2.7 c).} Let $(M,F^1_2,\vph_2,G)$ be a Shimura filtered $\sg$-crystal with the property that $W_0({\got g},\vph_2)$ is contained in $F^0_2({\got g})$. As in 3.2.3, let ${}_2B$ be a Borel subgroup of $G$ containing the image of a cocharacter $\mu_2:\GG_m\to G$ of the form $g\mu g^{-1}$ (with $g\in G(W(k))$ such that $g(F^1)=F^1_2$) and such that we have $F^1_2({\got g})\subset {\rm Lie}({}_2B)\subset {}_2{\got p}\subset F^0_2({\got g})$. Also, as in 3.2.3, let $g_3\in {}_2P(W(k))$, with ${}_2P$ as the parabolic subgroup of $G$ having ${}_2{\got p}$ as its Lie algebra, be such that
$g_3\vph_2\bigl({\rm Lie}({}_2B)\bigr)\subset {\rm Lie}({}_2B)$. But now
$(M,F^1_2,g_3\vph_2,G)$ has the same properties as $(M,F^1,\vph_0,G)$ which allowed us
to prove 3.2.5-6. We get that $(M,g_3\vph_2)$ has the formal isogeny type $\tau_0$. As
$(g_3\vph_2)^s=g_{3,s}\vph_2^s$, with $g_{3,s}\in {}_2P(W(k))$, $\forall s\in\NN$, we deduce from 3.2.1, that $(M,\vph_2)$ has as well the formal isogeny type $\tau_0$. This proves 3.2.7 c).
\medskip
{\bf 3.7.4. Proof of 3.2.7 b).} Let $(M,F^1_2,\vph_2,G)$ be a Shimura filtered $\sg$-crystal having $\tau_0$ as its formal isogeny type. 3.2.6 implies that the Lie $\sg$-crystal $\bigl({\rm End}(M),\vph_2\bigr)$ has the smallest Newton polygon 
${\rm End}(\Mp)$ among all  Lie $\sg$-crystals $\bigl({\rm End}(M),g\vph_2\bigr)$, with $g\in G(W(k))$. But for any $g\in G(W(k))$, the Newton polygon of
$\bigl({\rm End}(M),g\vph_2\bigr)$ is obtained (in the logical way) from the Newton polygon of $({\got g},g\vph_2)$ and from the Newton polygon of the quotient $\bigl({\rm End}(M)/{\got g},g\vph_2\bigr)$ (this quotient can be extended to a $p$-divisible object of $\Mm\Mf_{[-1,1]}(W(k))$). From this and the specialization part of the last statement of 3.7.1 we get that $({\got g},\vph_2)$ has the same Newton polygon as $({\got g},\vph_0)$. From 3.3.3 we get that (up to isomorphism) we can assume $\vph_2=p_0\vph_0$, with $p_0\in P_0(W(k))$. So, applying 3.2.1.1 in the context of 3.3.4, 3.2.7 b) follows. This (cf. 3.7.2-3) ends the proof of 3.2.7. 
\medskip
{\bf 3.7.5. End of the second proof of 3.1.0.} 3.2.5, 3.2.6 and 3.2.7 are proved respectively in 3.4.13, 3.7.1 and 3.7.2-4. Based on 3.2.4 (and 3.2.8), 3.1.0 follows. When we combine 2.3.4 with 3.2.7 b) we get:
$$A_G=B_G=C_G=D_Gg_0.\leqno (ALL)$$
\medskip
{\bf 3.7.6. The solution of 3.4.14.} Among the above four sections 3.7.1-4, only 3.7.1 needs to be slightly modified. For this modification, we follow the pattern of 3.1.8.1. Let ${\rm Spec}(R)$ be the completion of $G$ in its origin. We consider a Shimura filtered $F$-crystal ${\got C}$ over $R/pR$ of the form $(M,F^1,g_2\vph_0,G,\tilde f)$. Let $\Mk$ be the algebraic closure of the field of fractions of $R/pR$. We apply Fact 3 of 2.2.10 to the natural $k$-morphism ${\rm Spec}(\Mk)\to {\rm Spec}(R/pR)$: we get a Shimura $\sg_{\Mk}$-crystal ${\got C}_{\Mk}=(M\otimes_{W(k)} W(\Mk),g_{\Mk}(\vph_0\otimes 1),G_{W(\Mk)})$ specializing to $(M,F^1,g_2\vph_0)$. 
\smallskip
We use a specialization argument in the following way. We can assume $g_{\Mk}$ mod $p$ is defined naturally by the resulting $k$-morphism ${\rm Spec}(\Mk)\to G$; so $g_{\Mk}$ mod $p$ specializes to $g_0$ mod $p$. We can assume $k=\bar k$; so $G_R$ is split. Let $\tilde G_R^{\rm der}$ be the product of the semisimple subgroups of $G_R^{\rm der}$ having simple adjoints. So we can define a $\Phi_R$-linear map 
$$\bar\psi_{R/pR}^0:{\rm Lie}(\tilde G_{R/pR}^{\rm der})\to {\rm Lie}(\tilde G_{R/pR}^{\rm der})$$ 
as in 3.4.5 but in the context of ${\got C}$. Let $\bar\psi^0_{\Mk}$ be its extension via the $k$-homomorphism $R\to\Mk$. From the mentioned specialization, we get that $\bar\psi^0_{\Mk}$ specializes to $\bar\psi^0_0$ and so, $\forall s\in\NN$ we have
$$\dim_{\Mk}{\rm Im}(\bar\psi^0_{\Mk})^s\ge\dim_k{\rm Im}(\bar\psi^0_0)^s.\leqno (SPEC)$$ 
For $s$ big enough, the right hand side of this inequality is $m_{\dl}$. So, from 3.4.6 applied over $\Mk$, we get that the multiplicity of the slope $\dl$ of $({\got g}_0\otimes_{W(k)} W(\Mk),g_{\Mk}(\vph_0\otimes 1))$ is $m_{\dl}$. 
\smallskip
The same remains true if we work with another cycle of the permutation $\gamma$ of 3.4.0. From 3.4.8-12 we conclude: the Newton polygons of $({\got g}\otimes_{W(k)} W(\Mk),g_{\Mk}(\vph_0\otimes 1))$ and of $({\got g},\vph_0)$ are the same. So the Newton polygon of $({\got g},g_2\vph_0)$ is above the Newton polygon of $({\got g},\vph_0)$. The last paragraph of 3.7.1 does not need to be modified: we just need to replace ${\got C}_1$ by ${\got C}_{\Mk}$ and $k_1$ by $\Mk$. These replacements have to be performed in connection to 3.7.4 as well. This ends the solution of 3.4.14.   
\medskip\smallskip
{\bf 3.8. The non quasi-split context.} The basic results of 3.1 remain true, if we 
do not assume $G$ is quasi-split (assumption made in the def. 2.2.8). To see this, we
first remark that if $G$ is not quasi-split, then the field $k$ must be infinite (any reductive group over a finite field is quasi-split, cf. [Bo2, 16.6]). 
Let $k_1$ be a finite Galois extension of $k$ such that $G_{k_1}$ is quasi-split. Let $\sg_1$ be the Frobenius automorphism of $W(k_1)$. We know 3.1.0 is true for $G_{W(k_1)}$ (i.e. for when  we work with Shimura $\sg_1$-crystals $(M\otimes_{W(k)} W(k_1),g(\vph\otimes 1),G_{W(k_1)})$, with $g\in G(W(k_1))$). So to get 3.1.0 for $G$ (i.e. for when we work with Shimura $\sg$-crystals $(M,g\vph,G)$, with $g\in G(W(k))$) we just have to remark two simple Facts: 
\medskip
{\bf 3.8.1. Fact.} There is $g\in G(W(k))$ such that $(M\otimes_{W(k)} W(k_1),g\vph\otimes 1,G_{W(k_1)})$ is a $G_{W(k_1)}$-ordinary $\sg_1$-crystal.
\medskip
{\bf 3.8.2. Fact.} If $(M,g_0\vph)$, with $g_0\in G(W(k))$, is such that its extension to $k_1$ is a $G_{W(k_1)}$-ordinary $\sg_1$-crystal, due to the uniqueness assertion of 3.1.0 b), the filtration $F_0^1$ of $M\otimes_{W(k)} W(k_1)$ describing the $G_{W(k_1)}$-canonical lift of the $G_{W(k_1)}$-ordinary $\sg_1$-crystal $\bigl(M\otimes_{W(k)} W(k_1),g_0\vph\otimes 1,G_{W(k_1)}\bigr)$,  is definable over $W(k)$ (i.e. $F^1_0=\tilde F^1_0\otimes_{W(k)} W(k_1)$, with $\tilde F^1_0\subset M$).
\medskip
3.8.2 is trivial. 3.8.1 results from:
\medskip
-- the fact that $G(k)$ is dense in $G_k$ (cf. [Bo2, 18.3]) and so its image in $G(k_1)$ is dense in $G_{k_1}$, and from 
\smallskip
-- the fact that there is an open subscheme $U$ of $G_{k_1}$ such that for any $g\in G(W(k_1))$ lifting a $k_1$-valued point of $U$, $\bigl(M\otimes_{W(k)} W(k_1),g(\vph_0\otimes 1)\bigr)$ is
a $G_{W(k_1)}$-ordinary $\sg_1$-crystal. 
\medskip
This last fact is a consequence of the stronger form of 3.6.6 mentioned in 3.6.6.1 1), of 3.3.4 and of the part of 3.4.11 referring to things mod $p$; as pointed out in 3.6.6.1 1), it is also a consequence of 3.6.10, 3.4.8 and 3.4.11 (via the specialization theorem).
\medskip\smallskip
{\bf 3.9. Supplements to 3.1 and applications.}
\medskip
{\bf 3.9.1. The Lie stable $p$-rank.} Let $(M,\vph_2,G)$ be a Shimura $\sg$-crystal. For not overloading the notation, we first assume $G$ is split. We use the previous notations of 3.1-3, 3.4.0 and 3.4.5. One difference: we take $\vph_2:=g_2\vph_0$, with $g_2\in G(W(k))$. 
\medskip
Let $\Psi:\bigoplus_{i\in I}{\got g}_i\to\bigoplus_{i\in I}{\got g}_i$ be the $\sg$-linear map defined by the rule:
$$\Psi(x):=\cases{p\vph_2(x) &if $x\in{\got g}_i\ \ {\rm and}\ \ F^1({\got g}_i)\ne \{0\}$\cr
\vph_2(x) &if $x\in{\got g}_i\ \ {\rm and}\ \ F^1({\got g}_i)=\{0\}\;.$\cr}$$
Let $\bar\Psi:\bigoplus_{i\in I}{\got g}_i\otimes_{W(k)} k\to\bigoplus_{i\in I}{\got g}_i\otimes_{W(k)} k$ be obtained by taking $\bar\Psi$ mod $p$. Let 
$$N_{\vph_2}:=\bigcap_{m\in \NN} {\rm Im}(\bar\Psi^m).$$ 
Based on 3.4.5, we refer to $\Psi$ (resp. $\bar\Psi)$ as the Faltings--Shimura--Hasse--Witt shift (resp. map) of $(M,\vph_2,G)$ or of $({\got g},\vph_2)$.
\smallskip
We define 
$$R_{\vph_2}:=\dim_k(N_{\vph_2}).$$ 
We call $N_{\vph_2}$ (resp. $R_{\vph_2}$) the Lie stable $k$-vector space (resp. the Lie stable $p$-rank or the Faltings--Shimura--Hasse--Witt invariant) of $(M,\vph_2,G)$.
\smallskip
Similarly we define the Lie stable $k$-vector space and the Lie stable $p$-rank of any split Shimura Lie $\sg$-crystal over $k$.
\medskip
{\bf 3.9.1.0. Cyclic ranks.} As $\bar\Psi$ is a $\sg$-linear endomorphism of a finite dimensional $k$-vector space, we deduce the existence of $m\in\NN$ such that 
$$N_{\vph_2}={\rm Im}(\bar\Psi^m);$$
for instance we can take $m$ to be the following value which does not depend on $g_2$
$$l.c.m.\{o(\gamma_s)r_s|s\in\Mj_0\},$$ 
where $o(\gamma_s)$ is the order of the permutation $\gamma_s$ of $I$ and where $r_s$ is the relative dimension of $G_{i_s}$, with $i_s\in I$ such that $\gamma_s$ permutes cyclically a subset of $I$ containing $i_s$. The smallest such value of $m$ is called the cyclic rank of $(M,\vph_2,G)$ or of $({\got g},\vph_2)$. Similarly, we define the cyclic rank of any split Shimura Lie $\sg$-crystal over $k$.
\medskip
{\bf 3.9.1.1. Remark.} The Lie stable $p$-ranks, the Lie stable $k_1$-vector spaces and the cyclic ranks are well defined for any Shimura $\sg_{k_1}$-crystal $(M_1,\vph^1,G_1)$ over a perfect field $k_1$. This is so due to the fact that the map $\Psi_1$ defined as in 3.9.1 but starting from the extension of $(M_1,\vph^1,G_1)$ to $\overline{k_1}$ is in fact definable over $W(k_1)$, cf. the existence of a cocharacter $\mu_1:\GG_m\to G_1$ defining the filtration class of $(M_1,\vph^1,G_1)$ and cf. the classification (see [Ti1]) of $k_1$-simple (and so implicitly of $W(k_1)$-simple) adjoint groups. In other words we have: 
\medskip
-- ${\rm Lie}(G_1^{\rm ad})$ is a direct sum of $W(k_1)$-simple Lie factors, each one $\Mf$ of them having the property that all simple Lie factors of $\Mf\otimes_{W(k_1)} W(\overline{k_1})$ are producing (as in 3.4.3.4) $\vep_i$'s numbers which are either all positive or are all $0$;
\smallskip
-- the group cover of $G^{\rm der}_{1W(\overline{k_1})}$ defined naturally using the product of all semisimple subgroups of $G^{\rm der}_{1W(\overline{k_1})}$ having simple adjoints, is obtained from a group cover of $G^{\rm der}_1$ by extensions of scalars.
\medskip
The same applies to the Shimura Lie context. So from now on we do not assume $G$ is split. 
\medskip
{\bf 3.9.2. Theorem.} {\it A Shimura $\sg$-crystal $(M,\vph_2,G)$ is a $G$-ordinary $\sg$-crystal iff its Lie stable $p$-rank is maximal (i.e. $R_{{\vph}_2}\ge R_{{\vph}_3}$, for any $\vph_3=g_3\vph_2$, with $g_3\in G(W(k))$) and so equal to $R_{{\vph}_0}$.}
\medskip
{\bf Proof:} 
We first remark that $({\got g},\vph_2)$ has the smallest Newton polygon among all Shimura Lie $\sg$-crystals $({\got g},g\vph_2)$ iff for any cycle $\ga_0$ of the permutation $\ga$ of $I$, $({\got g}_0,\vph_2)$ has the smallest Newton polygon among all Shimura Lie $\sg$-crystals $({\got g}_0,g\vph_2)$  (here $g\in G(W(k))$, $\ga_0$ is a cyclic, transitive permutation of a subset $I_0$ of $I$, and ${\got g}_0=\oplus_{i\in I_0} {\got g}_i$). One implication is trivial, while the other one results from 3.4.8 and 3.1.0 c) applied to all classes of the form $Cl(M,\vph_0,\tilde G_0)$ (so defined by elements $s\in\Mj_0$). We have:
$$R_{\vph_2}:=\sum_{s\in\Mj_0} m_{\dl_s},$$ 
with each $m_{\dl_s}$ as a suitable multiplicity of some slope 
$$\dl_s\in[-1,0].$$ 
For $s\in\Mj_0$, $\dl_s$ and $m_{\dl_s}$ are defined in the same manner as the slope $\dl$ and its multiplicity $m_{\dl}$ of 3.4.4-5; for instance, $\dl_0=\dl$ and $m_{\dl_0}$ is (cf. 3.4.5.1 B) the multiplicity of the slope $\dl$ for $({\got g}_0,\vph_2)$. Warning: the formula of 3.4.4 remains true even if $I_1$ is the empty set. So the Theorem follows from 3.4.10.
\medskip
{\bf 3.9.3. Corollary.} {\it The fact that a Shimura $\sg$-crystal $(M,g\vph,G)$, with $g\in G(W(k))$, is $G$-ordinary depends only on the value of $g$ mod $p$.}
\medskip
{\bf 3.9.3.1. An application.} We refer to 3.6.1.3. We consider a point $y\in {\rm Spec}(R_j)(\bar k)$; let $z\in {\rm Spec}(R_j)(W(\bar k))$ be its Teichm\"uller lift. We consider the Shimura filtered $\bar\sg$-crystal obtained from $M_{R_j}$ by pull back via $z$ (cf. 3.6.3); if by forgetting its filtration,  we get a $G$-ordinary $\bar\sg$-crystal, then we refer to $y$ as a $G$-ordinary (or Shimura-ordinary) point. 3.9.3 guarantees that in fact 
we can decide if $y$ is or is not a $G$-ordinary point, without considering 
$z$. 
\smallskip
From 3.1.0 c) we get that the multiplicity of the slope $-1$ of the Shimura Lie $\sg$-crystal of $(M,\vph_0,G)$ is greater or equal to the multiplicity of the slope $-1$ of any Shimura Lie $\bar\sg$-crystal of the form $(M\otimes_{W(k)} W(\bar k),g(\vph_0\otimes 1),G_{W(\bar k)})$, with $g\in G(W(\bar k))$. So from 3.6.18.7.0, from 3.6.18.4.3 and from the construction of the special fibre of ${\rm Spec}(Q_j)$ (see 3.6.8 1) to 14); see also 3.6.8.1.3 and 3.6.8.2) we get:
\medskip
{\bf Corollary.} {\it If $y$ is a $G$-ordinary point and if $r_y=\dim_{W(k)}(G)$ (see 3.6.18.9 for the definition of $r_y$), then $y$ lifts to a point $y_{\infty}\in {\rm Spec}(Q_j)(\bar k)$.} 
\medskip
{\bf 3.9.4. The refined Lie stable $p$-rank.} We use the notations of the beginning of 3.4, with $G$ split. We choose a bijection $f_{\Mj_0}$ between $\Mj_0$ and the set $S(1,\abs{\Mj_0})$. Then $\forall s\in\Mj_0$, we can define the Lie stable $p$-rank $R_{\vph_2}(s)$ of $(M,\vph_2,G)$ w.r.t. $s$ or w.r.t. the subset of $I$ of which $\gamma_s$ is a cyclic permutation. For instance, if we work
with the cycle $\ga_0$, $R_{\vph_2}(0)$ is nothing else but the Lie stable $p$-rank of
$({\got g}_0,\vph_2)$. We denote by $R^r_{\vph_2}$ the $\abs{\Mj_0}$-tuple $(a_1,\ldots,a_{\abs{\Mj_0}})$ of non-negative integers defined by the formula:
$$a_l:={R_{\vph_2}\bigl(f^{-1}_{\Mj_0}(l)\bigr)
\over{\rm length}\bigl(\ga_{f^{-1}_{\Mj_0}(l)}\bigr)},$$
$\forall l\in S(1,\abs{\Mj_0})$. 3.4.2 guarantees that $a_l$'s are non-negative integers; for instance, if $l=f_{\Mj_0}(0)$, then $R_{\vph_2}(0)$ is the multiplicity of the slope $\dl$
of $({\got g}_0,\vph_2)$ (cf. 3.4.5.1 B) and so is a multiple of $\abs{I_0}={\rm length}(\gamma_0)$. If $G$ is not split, then we define the refined Lie stable $p$-rank of $(M,\vph,G)$ to be the refined Lie stable $p$-rank of its extension to $\bar k$.
\smallskip
On the set of $\abs{\Mj_0}$-tuples of non-negative integers we
introduce a partial order by the rule  
$$(b_1,b_2,\ldots,b_{\abs{\Mj_0}})\le
(c_1,c_2,\ldots,c_{\abs{\Mj_0}})$$ 
iff 
$$b_i\le c_i, ,\forall i\in S(1,\abs{\Mj_0}).$$
From the proof of 3.9.2 we get: 
\medskip
{\bf 3.9.5. Corollary.} {\it A Shimura $\sg$-crystal $(M,\vph_2,G)$ is a $G$-ordinary $\sg$-crystal iff its refined Lie stable $p$-rank is maximal, i.e. is equal to $R^r_{\vph_0}$.}
\medskip
{\bf 3.9.6. The adjoint context.} Similarly to 3.9.4, we define the refined Lie stable $p$-rank of any split Shimura Lie $\sg$-crystal. The refined Lie stable $p$-rank of a split Shimura (Lie) $\sg$-crystal over $k$ is the same, if instead of working with Shimura
Lie $\sg$-crystals of the above form $({\got g}_0,\vph_2)$, we work with their adjoint Lie
$\sg$-crystals $({\got g}_0^{\rm ad},\vph_2)$, where ${\got g}_0^{\rm ad}:={\rm Lie}(\tilde G_0^{\rm ad})\subset {\rm Lie}(\tilde G_0)[{1\over p}]$ (see 2.2.13; see 3.4.1.3 for notations). 
\smallskip
To see this, we use the notations of 3.4. If we are dealing with a non-empty subset of $I$ whose elements are permuted transitively by $\gamma$ and is such that the Lie algebras of the semisimple subgroups of $G$ indexed by it, are contained in $F^0({\got g})$, the statement is obvious. So we can assume we are dealing with the subset $I_0$ of $I$ of 3.4.0 (cf. the assumption of 3.4.1.2). The kernel of the natural map ${\got g}_0\otimes_{W(k)} k\to {\got g}_0^{\rm ad}\otimes_{W(k)} k$, has trivial intersection with the Lie algebra of any connected, smooth, unipotent subgroup of $G_k$. So, $\forall m\in\NN$, $\bar\Psi^m({\got g}_1)$ can be naturally identified with the $k$-vector subspace of ${\got g}_0^{\rm ad}$ defined similarly, using ${\rm Lie}(G_1^{\rm ad})\otimes_{W(k)} k$ and the adjoint variant of $\bar\Psi$. 
\smallskip
In what follows, we refer to this adjoint variant of $\Psi$ (resp. of $\bar\Psi$) as the Faltings--Shimura--Hasse--Witt adjoint shift (resp. adjoint map) of $(M,\vph_2,G)$ or of $({\got g},\vph_2)$. Similarly, we speak about the Faltings--Shimura--Hasse--Witt adjoint shift (resp. adjoint map) of $(M,\vph_2,G)$ w.r.t. $I_0$ (or to an element $s\in\Mj_0$ or to the permutation $\gamma_s$), or more generally w.r.t. a non-empty subset of $\Mj_0$ (i.e. w.r.t. a non-empty subset of $I$ left invariant by $\gamma$).
\medskip
{\bf 3.9.7. The zero, the positive, and the non-negative rank and type.} By the non-negative rank (resp. type) of a Shimura $\sg$-crystal $(M,\vph,G)$ over $k$ (or of its attached Shimura Lie $\sg$-crystal) we mean the rank as a $W(k)$-module (resp. we mean the $G(W(k))$-conjugacy class) of ${\got p}_{\ge 0}:=W_0({\rm Lie}(G),\vph)$. Similarly we define the positive and the zero ranks and types, by using positive slopes and respectively the slope $0$ (the wording ``parabolic" has to be removed or replaced accordingly). 
\smallskip
Warning: to know the non-negative rank is the same as to know the zero or the positive rank; but this is not so if we replace rank by type. These definitions extend to all Shimura (adjoint) Lie $\sg$-crystals or isocrystals or to a filtered context. 
\medskip
{\bf 3.9.7.1. The refined non-negative type.} We consider a lift $F^1$ of $(M,\vph,G)$. Let $P$ be the parabolic subgroup of $G$ normalizing $F^1$. 
As in 3.2.3, we fix a maximal torus $T$ of $P$ such that the canonical split cocharacter of $(M,F^1,\vph)$ factors through $T$. 
For 
$$\al\in SS:=[0,1]\cap {1\over {\dim_{W(k)}(M)!}}\ZZ,$$ 
let $P_{\ge \al}$ (resp. $P_{=\al}$) be the integral, closed subgroup of $G$ such that the Lie algebra of its generic fibre is $W_{\al}({\rm Lie}(G)[{1\over p}],\vph)$ (resp. is $W(\al)({\rm Lie}(G)[{1\over p}],\vph)$) (see 2.2.3 3)). $P_{\ge 0}$ is a parabolic subgroup of $G$ and $P_{=\al}$ and $P_{\ge \al}$ are subgroups of $P_{\ge\be}$ if $\al\in [\be,1]$, with $\al$, $\be\in SS$. Let $CONG$ be the subgroup of $G(W(k))$ normalizing $F^1/pF^1$. 
\medskip
{\bf Definition.} By the refined non-negative type of $(M,\vph,G)$ we mean the $CONG$-conjugacy class of the (indexed) family of subgroups of $G$ formed by $(P_{\ge\al})_{\al\in SS}$ and by $(P_{=\al})_{\al\in SS}$. It depends only on the inner isomorphism class of $(M,\vph,G)$.
\medskip
The intersection $P_k\cap {P_{\ge 0}}_k$ contains a maximal torus of $P_k$. As in the proof of Fact 1 of 2.2.9 3), we deduce that, up to inner isomorphism of $(M,\vph,G)$ defined by an element $el$ of $G(W(k))$ which mod $p$ is the identity, we can assume that the intersection $P\cap P_{\ge 0}$ contains a maximal torus $T_0$ of $G$. Till end of 3.9.7.2 we assume $k=\bar k$. We can assume that, up to inner isomorphisms defined by elements of $P(W(k))$, $T_0=T$ (to be compared with 3.2.3). As $T$ was fixed, we have a finite number of possibilities for $P_{\ge 0}$: the number of parabolic subgroups of $G$ containing $T$ is finite. 
\medskip
{\bf Questions.} Which parabolic subgroups of $G$ containing $T$ can be the non-negative type of $(M,g\vph,G)$, for some $g\in G(W(k))$? What general principles are governing the refined non-negative types of $(M,g\vph,G)$, for $g\in G(W(k))$ allowed to vary? When $(M,\vph,G)$ is uniquely determined (up to inner isomorphism) by its (refined) non-negative type?
\medskip
{\bf 3.9.7.2. The inductive property.} 
Let $R_0$ be the quotient of $P_{\ge 0}$ by its unipotent radical $P_{>0}$. We view it as well as a subgroup of $G$ via: it is naturally isomorphic to the only Levi subgroup of $P_{\ge 0}$ containing $T$. Let 
$$g\in P_{\ge 0}(W(k))$$ 
be such that $g\vph({\rm Lie}(T))={\rm Lie}(T))$ (this is the same as 3.2.3). Let $h$ be the image in $R_0(W(k))$ of $g^{-1}$. The quadruple $(M,g\vph,T)$ is a Shimura $\sg$-crystal (cf. Fact 1 of 2.2.9 1)). So $(M,hg\vph,R_0)$ is as well a Shimura $\sg$-crystal. Moreover, $$\vph=n_+(1)(hg\vph)$$
with $n_+(1)\in P_{>0}(W(k))$. By induction on $m\in\NN$ we get that 
$$\vph^m=n_+(m)(hg\vph)^m,$$
with $n_+(m):=\vph^m(hg\vph)^{-m}\in P_{>0}(B(k))$. So the semisimple elements $b_1$ and $b_2$ of $G(B(k))$ defined by $\vph^m$ and respectively by $(hg\vph)^m$ (see 2.2.24 and 2.2.24.1) are $P_{>0}(B(k))$-conjugate. Argument: $b_1$ and $b_2$ have the same image in $R_0(B(k))$ and so can be viewed as $B(k)$-valued points of maximal tori $T_1$ and respectively $T_2$ of $P_{\ge 0 B(k)}$ having the same image in $R_{0B(k)}$; based on the conjugation part of [Bo2, 20.5] we can assume $T_1=T_2$ and so $b_1=b_2$. We conclude:
\medskip
{\bf Corollary.} {\it The Newton polygons of $(M,\vph)$ and of $(M,hg\vph)$ are the same.}
\medskip
{\bf 3.9.7.2.1. Remarks.} {\bf 1)} 3.9.7.1-2 extend automatically to the generalized Shimura context of 2.2.8 3) and 4) but (warning) not to the context of 2.2.8 3a) and 4a) (the reason is: in the context of 2.2.8 3a) and 4a), not any element $el$ as in 3.9.7.1 can be used as an inner isomorphism).
\smallskip
{\bf 2)} The duality of language mentioned at the end of 2.2.3 3) allows us to redo all of 3.9.7.1-2 in the context of non-positive slopes. From the point of view of 3.9.7.2, we presently can not see any advantage of using the (refined) non-positive types instead of (refined) non-negative types.  
\medskip
{\bf 3.9.7.3. An example.} Let $\al\in W(k)$. We take $\dim_{W(k)}(M)=6$ and $G=GL(M)$. Let $\vph_{\al}$ take the elements of a $6$-tuple $(e_1,...,e_6)$ formed by elements of a $W(k)$-basis of $M$ into $(e_2,pe_3,pe_1,e_5,e_6+\al e_1,pe_4)$. The slopes of $(M,\vph_{\al})$ are $1\over 3$ and $2\over 3$. The non-negative type of $(M,\vph_{\al},G)$ is defined by the parabolic subgroup $P_{\ge 0}$ of $G$ normalizing the $W(k)$-span $<e_1,e_2,e_3>$. Let $T$ be the maximal torus of $G$ normalizing $e_i$, $i=\overline{1,6}$. Let $P_{>0}$, $P_{\ge 0}$ and $R_0$ be as abve. Let $N_{16}$ be the unipotent subgroup of $G$ that fixes $e_i$, $i=\overline{1,5}$, and takes $e_6$ into the $W(k)$-span $<e_1,e_6>$. Its dimension is $1$. We have:
\medskip
{\bf Claim.} {\it If $\al_1$, $\al_2\in W(k)$ are such that $\al_1-\al_2\not\in pW(k)$, then $(M,\vph_{\al_1})$ and $(M,\vph_{\al_2})$ are not isomorphic (even if we consider pull backs to other algebraically closed fields containing $k$).}
\medskip
{\bf Proof:} This can be easily checked. First, any isomorphism between 2 members of the family $\bigl((M,\vph_{\al},G)\bigr)_{\al\in W(k)}$ is defined by an element $h_0\in P_{\ge 0}(W(k))$. Second, we can write $h_0=n_+r_0$, with $n_+\in P_{\ge 0}(W(k))$ and $r_0\in R_0(W(k))$. Third, $r_0$ is an automorphism of $(M,\vph_0)$. Fourth, a simple computation shows that $r_0$ mod $p$ (under inner conjugation) centralizes ${N_{16}}_k$. Fifth, $P_{>0}$ is abelian and so all its $W(k)$-valued points (under inner conjugation) act trivially on $N_{16}$. So, $h_0\vph_{\al_1}h_0^{-1}=n_{12}\vph_{\al_2}$, with $n_{12}\in P_{>0}(W(k))$ such that its component in $N_{16}$ is non-zero mod $p$; here the components of $n_{12}$ are with respect to $\GG_a$ subgroups of $P_{>0}$ normalized by $T$. This ends the proof.
\medskip
Conclusion: the number of isomorphism classes of $p$-divisible groups of rank $6$ and dimension $3$ over $k$ is not finite. Based on 3.6.15 B we also get: there is $n\in\NN$ such that the number of isomorphism classes of tuncations mod $p^n$ of $p$-divisible groups of rank $6$ and dimension $3$ over $k$ is not finite (it seems to us that we can take $n=2$ but we will not stop to check this).
\smallskip
The above example can be put in a family: we just need to apply 3.6.11 in the context of $N_{16}$ and of $(M,\vph_{0},G)$.   
\medskip
{\bf 3.9.8. Truncations in the adjoint Lie $\sg$-crystal context.} Not to be too long, we refer just to the adjoint context; the notations are independent of the previous ones. Let $({\got g},\vph,F^0({\got g}),F^1({\got g}))$ be an arbitrary Shimura adjoint Lie $\sg$-crystal. Let $({\got g},\tilde F^1({\got g}),\tilde F^2({\got g}),\Psi)$ be its Faltings--Shimura--Hasse--Witt shift (cf. 3.9.1.1 and 3.4.5). In 2.2.14 we defined the truncation mod $p^n$, where $n\in\NN$, of $({\got g},\vph,F^0({\got g}),F^1({\got g})$) as well as of $({\got g},\tilde F^1({\got g}),\tilde F^2({\got g}),\Psi)$. The pair 
$$({\got g}/p^n{\got g},\Psi)$$
 is referred as the truncation mod $p^n$ of $({\got g},\vph)$ (see 3.13.7.4 2) below for a justification; here we still denote by $\Psi$ its reduction mod $p^n$). As in 2.2.14.2, as the adjoint group over $W(k)$ having ${\got g}$ as its Lie algebra is uniquely determined (see 2.2.13), we speak about inner automorphisms of such truncations.
\medskip
{\bf 3.9.9. Deformations in the non-smooth context.} We consider the following situation. Let $(M,\vph_0,G,(t_{\al})_{\al\in\Mj})$ be a Shimura-ordinary $\sg$-crystal with an emphasized family of tensors. Let $g\in G(W(k))$. We consider a set $\Mj_g$ containing $\Mj$ and an extra family of tensors $(t_{\al})_{\al\in\Mj_g\setminus \Mj}$ of $\Mt(M[{1\over p}])$ fixed by $g\vph_0$. Let $G_1$ be the Zariski closure in $GL(M)$ of the subgroup of $GL(M[{1\over p}])$ fixing $t_{\al}$, $\forall\in\Mj_g$. We assume it is integral. The question is:
\medskip
{\bf Q.} {\it Under what conditions, the quadruple $(M,g\vph_0,G,(t_{\al})_{\al\in\Mj_g})$ is the specialization of a similar quadruple over $k[[x]]$, with $x$ an independent variable, such that by forgetting the extra tensors (indexed by $\Mj_g\setminus\Mj$) and by pull back under a geometric, dominant point of ${\rm Spec}(k[[x]])$, we get a Shimura-ordinary $F$-crystal (isomorphic to the corresponding pull back of $(M,\vph_0,G,(t_{\al})_{\al\in\Mj})$)?}
\medskip
For instance, if we have a quasi-polarization $p_M:(M,g\vph_0)\otimes_{W(k)} (M,g\vph_0)\to W(k)(1)$ and if $G$ is the maximal integral subgroup of $GL(M)$ normalizing $p_M$, then the first main result of [NO] (in its abstract formulation pertaining to $p$-divisible groups) can be interpreted as: no extra conditions are required. We have a similar question, where $k[[x]]$ is replaced by the henselization of the localization of $k[x]$ w.r.t. $(x)$. 
\smallskip
To handle such questions, we propose an approach in two steps. First step would be: use B6 of 3.14 below to deform $(M,\vph_0,G,(t_{\al})_{\al\in\Mj_g})$ to a similar quadruple $(M\otimes_{W(k)} W(k_1),g_1(\vph_0\otimes 1),G_{W(k_1)},(t_{\al})_{\al\in\Mj_g}))$ over the algebraic closure $k_1$ of $k((x))$ which lifts to $W(k_1)$, in the sense that there is a lift $F^1_1$ of $(M\otimes_{W(k)} W(k_1),g_1(\vph_0\otimes 1),G_{W(k_1)})$ such that $t_{\al}$ is in the $F^0_1$-filtration of $\Mt(M[{1\over p}])$ defined by $F^1_1$, $\forall\al\in\Mj_g\setminus\Mj$. Second step would be:   show that the quintuple $(M\otimes_{W(k)} W(k_1),g_1(\vph_0\otimes 1),G_{W(k_1)})$ can be deformed over the $p$-adic completion of an $\NN$-pro-\'etale scheme over a smooth $W(k_1)$-scheme, which has a connected special fibre and has $k_1$-valued points with the property that the pull back through them of the corresponding deformation are, when the extra tensors are ignored, Shimura-ordinary $\sg_{k_1}$-crystals. 
\smallskip
In what follows we deal only with the second step; it gives a lot of extra information (warning: the below proof is part of the third place where we need $p>2$). So we assume there is a lift $F^1$ of $(M,g\vph_0,G)$ such that $t_{\al}\in F^0(\Mt(M[{1\over p}]))$, $\forall\al\in\Mj_g\setminus\Mj$. We also assume there is a smooth, affine $W(k)$-scheme $G_2$ and a $W(k)$-morphism $m_1:G_2\to G_1$ such that the following condition holds:
\medskip
{\bf COND.} $g\in G_1(W(k))$ and there are $W(k)$-valued points $z_1$ and $z_2$ of $G_2$ lifting respectively $1_M$ and $g^{-1}$ and whose special fibres factor through the same connected component of $G_{2k}$.
\medskip
{\bf A. Theorem.} {\it There is an $\NN$-pro-\'etale scheme $G_3$ of $G_2$ to which $z_1$ lifts and whose special fibre is a connected, $AG$ $k$-scheme, and there is a deformation (defined by an object of $p-\Mm\Mf_{[0,1]}^\nabla(G_3)$) of $(M,F^1,g\vph_0,G,(t_{\al})_{\al\in\Mj_g})$ over $G_3^\wedge$ which generically in this special fibre is (in the sense of forgetting all tensors indexed by elements of $\Mj_g\setminus\Mj$) Shimura-ordinary.}
\medskip
{\bf Proof:} We write $G_2={\rm Spec}(R_2)$. Localizing, we can assume we have a potential-deformation sheet $(G_2,b_2,z_1)$ and that $G_{2k}$ is connected; warning: it does not matter if $z_2$ does not factor any more through $G_2$ under this localization. Let $\Phi_{R_2}$ be the Frobenius of $R_2^\wedge$ we get as in 3.6.9.1 from $b_2$. We apply 3.6.18.6 to the following $p$-divisible object of $\Mm\Mf_{[0,1]}(R_2)$
$${\got C}_2:=(M\otimes_{W(k)} R_2^\wedge,F^1\otimes_{W(k)} R_2^\wedge,g_2^{\rm univ}(g\vph_0\otimes 1)),$$
with $g_2^{\rm univ}\in G(R_2^\wedge)\subset GL(M)(R_2)$ the universal element logically defined by $m_1$. Let 
$$\ell_2:G_3\to G_2$$ 
be as $\ell_1$ in the end of 3.6.18.6 a); so $G_3^\wedge:=M_0({\got C}_2)$ as a $G_2$-scheme. 3.9.2 implies that the condition of Shimura-ordinariness is an open condition. So from COND we deduce that there is an open subscheme $U_2$ of $G_{2k}$ such that for any $W(\bar k)$-valued Teichm\"uller lift $z$ of $G_2^\wedge$ whose special fibre factors through $U_2$, $(z^*({\got C}_2),G_{W(\bar k)})$ is a lift of a Shimura-ordinary $\bar\sg$-crystal. From this and 3.1.8.1.2 b), the Theorem follows (to be compared with 3.9.3.1). 
\medskip
{\bf B. Example.} We assume that we are in the context of a quasi-polarized filtered $\sg$-crystal $(M,F^1,g\vph_0,p_M)$. So $G=GL(M)$ and $G$ is the maximal integral subgroup of $G$ normalizing $p_M$. We assume $g\in G_1(W(k))$ (if $k=\bar k$, due to the assumption on $F^1$ this represents no restriction as we can see through simple arguments involving $\ZZ_p$-structures as in 2.2.9 8)). Let $e\in\NN$ be such that $\dim_{W(k)}(M)=2e$. It is an elementary fact that there is a $W(k)$-basis $\{e_1,e_2,...,e_{2e}\}$ of $M$ such that $p_M(e_i,e_j)=0$ unless $\abs{j-i}=e$. We can assume $p_M(e_i,e_{i+e})=p^{n_i}$, where
$$1=n_1=n_2=...=n_{i_1}<n_{i_1+1}=...=n_{i_2}<...<n_{i_{m-1}+1}=...=n_{i_m},$$ with $m\in\NN$ and $i_m=e$. Let $i_0:=0$. There is a connected subgroup $S_1$ of $GL(M/pM)$ such that the special fibres of all $W(\bar k)$-valued points of $G_1$, do factor through $S_1$ and the resulting set of $\bar k$-valued points of $S_1$ is dense in $S_1$. This is an immediate consequence of the below two properties (P1) and (P2) and of the fact that the standard representation of any symplectic group over a field is absolutely irreducible (this last thinks takes care of the connectedness part). We have:
$$S_1^{\rm der}=\prod_{j=1}^m Sp(V_j,p^{-n_j}p_M)_k,\leqno (P1)$$
where 
$$V_j:=<e_{i_{j-1}+1},...,e_{i_j},e_{i_{j-1}+1+e},...,e_{i_j+e}>.$$
Moreover, 
$$s(V_j/pV_j)\subset \oplus_{l=j}^m V_j/pV_j\otimes_k \bar k,\leqno (P2)$$  
$\forall s\in S_1(\bar k)$. So using a sequence of dilatation (see [BLR, 3.2-4]; the first one is the dilatation of $S_1$ in $G_1$), we deduce the existence of a morphism $\ell_2$ as above for which COND holds; $\ell_2$ is in fact a homomorphism (cf. [BLR, (d) of p. 64]). We conclude: 
\medskip
{\bf Corollary.} {\it Any quasi-polarized $p$-divisible group over $k$ which lifts to a quasi-polarized $p$-divisible group over $W(k)$, can be deformed over the $p$-adic completion of an $\NN$-pro-\'etale scheme $X$ over a smooth $W(k)$-scheme, with $X_k$ as a geometrically connected, $AG$ $k$-scheme, to a quasi-polarized, ordinary $p$-divisible group over $\bar k$.}
\medskip
There are plenty of examples (for instance, see [Nor, ch. 4]) of quasi-polarized $\bar\sg$-crystals which are not obtained from quasi-polarized filtered $\bar\sg$-crystals by forgetting the filtration. Here is a very simple numerical version of loc. cit. We take $e=2$ and $k=\FF_p$. We assume $p_M(e_1,e_3)=1$ and $p_M(e_2,e_4)=p$. We denote $g\vph_0$ by $\vph$ and we assume $\vph$ takes $(e_1,e_2,e_3,e_4)$ into $(e_3,e_4,pe_1,pe_2)$. So all slopes of $(M,\vph)$ are ${1\over 2}$. It is trivial to see that there is no direct summand $F^1$ of $M\otimes_{\ZZ_p} W(\FF)$ lifting the $\FF$-submodule of $M\otimes_{W(k)} \FF$ generated by $e_3$ and $e_4$ mod $p$ and such that $p_M(a,b)=0$, $\forall a,b\in F^1$.
\medskip
{\bf C. Remarks.} {\bf 1)} We come back to the general context of a lift $F^1$ of $(M,g\vph_0,G)$ such that $t_{\al}\in F^0(\Mt(M[{1\over p}]))$, $\forall\al\in\Mj_g\setminus\Mj$. As in C we define $S_1$. Warning: in general $S_1$ is not connected. However, $S_1$ is a smooth $k$-group. We consider a homomorphism $m_1^\prime:G_2^\prime\to G_1$, with $G_2^\prime$ smooth over $W(k)$, which is obtained as in C, using dilatations. Using the universal property of dilatations (see [BLR, p. 63]), we deduce that we can assume $m_1$ factors through $m_1^\prime$. If COND holds for $m_1$ then it holds for $m^\prime$ and so we can assume $m_1=m_1^\prime$; however, for different computational problems it is still useful not to assume $m_1=m_1^\prime$. The Kodaira--Spencer map of the deformation we get over a $G_3^{\prime\wedge}$ as in the above proof (but working with $m_1^\prime$), in each $\bar k$-valued point $y$ of $G_3^\prime$ has an image whose dimension is exactly the dimension $d$ of the image of ${\rm Im}({\rm Lie}(G_{2k}^\prime)\to {\rm Lie}(S_{1k})$) in ${\rm End}(M)/F^0({\rm End}(M))\otimes_{W(k)} k$. So forgetting the tensors $t_{\al}$, $\al\in\Mj_g\setminus\Mj$, by taking suitable slices we get (over resulting formally smooth subschemes of $G_3^{\prime\wedge}$ of relative dimension $d$) Shimura $p$-divisible groups which are uni plus quasi-versal. 
\smallskip
{\bf 2)} We do not assume anymore the existence of $F^1$ as in 1). Still $S_1$ can be defined as in C. Let $S_2$ be the reduced group subscheme of $S_1$ which normalizes the kernel $\bar F_1$ of $\vph_0$ mod $p$. If $S_2$ has no characters producing as usual a direct sum decomposition $M/pM=\bar F^1\oplus\bar F^0$, then an $F^1$ as in 1) does not exist.
\medskip\smallskip
{\bf 3.10. Terminology and formulas.}
\medskip
All terminologies and formulas to be introduced below in the context of Shimura Lie $\sg$-crystals, will be used as well in the rational context of Shimura (adjoint) (filtered) Lie isocrystals attached to Shimura (adjoint) (filtered) Lie $\sg$-crystals. 
\medskip
{\bf 3.10.0. Cyclic Shimura Lie $\sg$-crystals.} The notations that follows are similar to the previous ones but, without special reference, they are independent. We start considering an arbitrary perfect field $k$. A Shimura adjoint Lie $\sg$-crystal $({\got g},\vph)$ over $k$ is said to be cyclic, if {\got g} is the Lie algebra of an adjoint group $G$ over $W(k)$ such that
considering the product decomposition $G_{W(\bar k)}=\prod_{i\in\tilde I}G_i^\prime$ in simple groups over $W(\bar k)$, $\vph\otimes 1$ permuting cyclically the direct summands ${\rm Lie}(G_i^\prime)[{1\over p}]$'s of
${\got g}\otimes_{W(k)} W(\bar k)[{1\over p}]$. 
\smallskip
 $({\got g},\vph)$ is said to be trivial if
$\vph({\got g})={\got g}$; otherwise it is said to be non-trivial. $({\got g},\vph)$ is said to be of $A_\ell$, $B_\ell$, $C_\ell$, $D_\ell$, $E_6$ or $E_7$ Lie type if the factors $G_i^\prime$ are of this Lie type; in what follows, we ignore the other Lie types, as they force $({\got g},\vph)$ to be trivial, cf. 2.2.7.
\smallskip
As in 3.2.3, let $g\in G(W(\bar k))$ be such that
$g(\vph\otimes 1)$ takes the Lie algebras {\got b} and {\got t} of a Borel subgroup $B_{W(\bar k)}$ of $G_{W(\bar k)}$ and respectively of a maximal torus $T_{W(\bar k)}$ of $B_{W(\bar k)}$ into themselves. The permutation of the
basis of roots of ${\rm Lie}(G_i^\prime)$ w.r.t. the maximal torus $T_i^\prime:=T_{W(\bar k)}\cap G_i^{\prime}$ of $G_i^\prime$ and corresponding to the Borel subgroup $B_i^\prime:=B_{W(\bar k)}\cap G_i^\prime$ of $G_i^\prime$, induced (as in 3.4.3.1) by the restriction of $(g(\vph\otimes 1))^{\abs{\tilde I}}$ to ${\rm Lie}(B_i^\prime)$, can be non-trivial. This permutation can be (cf. 3.4.3.1):
\medskip
-- an involution (this can happen only when {\got g} is of $A_\ell$, $D_\ell$ or
$E_6$ Lie type);
\smallskip
-- a cycle of order 3 (this can happen only when {\got g} is of $D_4$ Lie type);
\smallskip
--  trivial. 
\medskip
We say respectively that $({\got g}\otimes_{W(k)} W(\bar k),g(\vph\otimes 1))$
(or $({\got g},\vph)$ itself) is with involution, or with a 3-cycle, or without a cycle. 
\smallskip
All above definitions are used as well in a filtered context or in a non-necessarily adjoint context (i.e. when we work with the Lie algebra of a semisimple group over $W(k)$ which is not necessarily adjoint). In particular we also speak about cyclic Shimura Lie $\sg$-crystals which are not necessarily adjoint.
\medskip
{\bf 3.10.0.1. Remark.} 
If we just require the simple factors of ${\rm Lie}(G)$ to be permuted cyclically by $\vph$, then we say $({\got g},\vph)$ is weakly cyclic. These notions of weakly cyclic and cyclic coincide over finite fields but not in general. To exemplify this, we restrict to the trivial context, as the adjustments for the non-trivial context are easily made. Let $k$ be such that there is a Galois extension of it $k_1$ of whose Galois group is not cyclic. Let $\gamma$ be an automorphism of $W(k_1)$ whose restriction to $W(k)$ is $\sigma$. Let $H_1$ be a split, simple adjoint group over $\ZZ_p$. Let $H:={\rm Res}_{W(k_1)/W(k)} H_{1W(k_1)}$; it is a $W(k)$-simple, adjoint group. Its Lie algebra (over $W(k)$) is ${\rm Lie}(H_1)\otimes_{\ZZ_p} W(k_1)$. We consider the $\sg$-linear automorphism $a$ of it which acts identically on ${\rm Lie}(H_1)$ and as $\gamma$ on $W(k_1)$. $H_{W(k_1)}$ is a disjoint union of $[k_1:k]$ copies of $H_{1W(k_1)}$; but, as ${\rm Gal}(k_1/k)$ is not cyclic, the Lie algebras of these copies of $H_{1W(k_1)}$ are not permuted transitively by $a\otimes 1$.  
\smallskip
In what follows, without special reference, we deal only with the cyclic notion.
\medskip
{\bf 3.10.1. Definitions: cyclic factors.} Any split Shimura adjoint (filtered) Lie $\sg$-crystal ${\got C}^{\rm ad}$ is isomorphic to a direct sum of cyclic split Shimura adjoint (filtered) Lie $\sg$-crystals, which are called the cyclic  factors of ${\got C}^{\rm ad}$. The direct sum of the non-trivial cyclic factors is called the non-trivial part of ${\got C}^{\rm ad}$, while the direct sum of the trivial ones is called the trivial part of ${\got C}^{\rm ad}$. 
Similarly, for any Shimura filtered $\sg$-crystal $(M,F^1,\vph,G_M)$ or Shimura $\sg$-crystal $(M,\vph,G_M)$ (resp. for any Shimura filtered Lie $\sg$-crystal $({\rm Lie}(G_M),\vph,F^0({\rm Lie}(G_M)),F^1({\rm Lie}(G_M)))$ or Shimura Lie $\sg$-crystal $({\rm Lie}(G_M),\vph)$), with $G^{\rm ad}_M$ a split group, we define  its cyclic Lie factors (resp. its cyclic factors), as well as its cyclic adjoint factors.  
\smallskip
For instance, with the notations of the part of 3.4 ending with 3.4.1, the Shimura filtered Lie $\sg$-crystal $({\got g}_0,\vph_2,F^0({\got g}_0),F^1({\got g}_0))$ of 3.4.1 is cyclic and is called a cyclic Lie factor of $(M,F^1,\vph_2,G)$ (or a cyclic factor of $({\got g},\vph_2,F^0({\got g}),F^1({\got g}))$); its adjoint (defined naturally using the Lie algebra of $\prod_{i\in I_0} G_i^{\rm ad}$) is called a cyclic adjoint factor of $(M,F^1,\vph_2,G)$ or of $({\got g},\vph_2,F^0({\got g}),F^1({\got g}))$. 
\smallskip
Similarly we speak about the non-trivial part (or the trivial part) of a split Shimura (filtered) Lie $\sg$-crystal $\Ml$. Warning: they are not necessarily defined by direct summands of the underlying module of $\Ml$.
\smallskip
Similarly, for any Shimura (adjoint) (filtered) Lie $\sg$-crystal (resp. any Shimura (filtered) $\sg$-crystal), we speak about its weakly cyclic (adjoint) factors (resp. about its weakly cyclic Lie factors and its weakly cyclic adjoint factors) and about these factors being trivial or non-trivial. The only difference is: the adjoint groups we get are not always products of absolutely simple, adjoint groups (cf. 3.10.0.1). Whenever possible, we drop the word weakly.
\medskip
{\bf 3.10.2. A context.} From now on we work with a cyclic Shimura adjoint filtered Lie $\sg$-crystal
$\bigl({\got g},\vph,F^0({\got g}),F^1({\got g})\bigr)$ which is non-trivial and has the property that there are Lie algebras {\got b} and {\got t} of a Borel subgroup $B$  of $G$ and respectively of a maximal torus $T$ of $B$, such that $F^1({\got g})\subset{\got b}\subset F^0({\got g})$,  $\vph({\got t})={\got t}$ and $\vph({\got b})
\subset{\got b}$ (to be compared with 3.2.3). Here $G$ is as in 3.10.0.
\medskip
All the terminology to be introduced below for the case $G$ split is also used:
\medskip
a) in the context of the above paragraph (via the passage to $W(\bar k)$), and 
\smallskip
b) for cyclic factors of  split Shimura filtered Lie $\sg$-crystals defined by Shimura-canonical lifts (i.e. for factors like
$\bigl({\got g}_0,\vph_1,F^0({\got g}_0),F^1({\got g}_0)\bigr)$ of 3.4.1). 
\medskip
{\bf 3.10.3. Notations.} From now on we also assume $G$ is a split group. For not overloading the notations, we consider $G=\times_{i\in I_0} G_i$, ${\got g}={\got g}_0,F^0({\got g})=F^0({\got g}_0),F^1({\got g})=
F^1({\got g}_0)$ (``with ${\got g}_0$ of 3.4.0") and  $\vph=\vph_1$ (``of 3.2.3"). So we use fully the notations and conventions introduced in 3.4.0-3 and 3.4.5. So
$\tilde I=I_0$ and $G_i^{\prime}=G_i$, $\forall i\in I_0$. Warning: here $G_i$ are adjoint groups; so only the notations are as in the mentioned places (this motivates the above use of quotation marks). 
\medskip
{\bf 3.10.4. The $\vep$-type.} The $m$-tuple $(\vep_1,\ldots,\vep_m)$ is called the $\vep$-type of
$({\got g}_0,\vph_1)$. Here 
$$m:=\abs{\tilde I_0}=e\abs{I_0}\in\{\abs{I}_0,2\abs{I_0},3\abs{I_0}\}.$$ 
If $({\got g}_0,\vph_1)$ is with involution or with a 3-cycle, then we replace $({\got g}_0,\vph_1$) by 
$$({\got g}_0^e,\vph_1^{(e)})$$
 (cf.  3.4.3.1) ($e$ being 2 or 3, depending on $\pi_1$ being an involution or a 3-cycle). Here
$$\vep_i\in S(0,\ell).$$
\indent
{\bf 3.10.4.1. Warning.}
We replace $({\got g}_0,\vph_1)$ by $({\got g}_0^e,\vph_1^{(e)})$ only for defining the $\vep$-type of $({\got g}_0,\vph_1)$. For formulas that follows (see 3.10.6) this replacement is not allowed (performed).
\medskip
{\bf 3.10.4.2. On the value of $m$.} If $({\got g}_0,\vph_1)$ is without a cycle, then the $m$-tuple $(\vep_1,\ldots,\vep_m)$ is uniquely determined up to a cyclic permutation and up to a renumbering of the roots of $\Dl_1$ (cf. 3.4.3.2). We have $m=n=\abs{I_0}$. 
\smallskip
If $({\got g}_0,\vph_1)$ is with involution then the $m$-tuple
$(\vep_1,\ldots,\vep_m)$ is uniquely determined up to a cyclic permutation and, in the case when ${\got g}_0$ is of $D_4$ Lie type,  up to a renumbering of the roots of $\Dl_1$. We have $m=2n$.
\smallskip
If $({\got g}_0,\vph_1)$ is with a 3-cycle then the $m$-tuple $(\vep_1,\ldots,\vep_m)$ is uniquely determined up to cyclic permutation and up to a renumbering of the roots of $\Dl_1$. We have $m=3n$. 
\medskip
{\bf 3.10.4.3. Convention.} An $\vep$-type $(\vep_1,...,\vep_m)$ is said to be in the standard form, if $\sum_{i=1}^m \vep_i(\ell+1)^{m-i}$ is the greatest number among all numbers obtained in this manner using an $\vep$-type obtained from $(\vep_1,...,\vep_m)$ through the rules of 3.10.4.2 (i.e. up to renumberings and cyclic permutations). Here the role of $\ell+1$ is irrelevant: we just need an integer greater than all $\vep_i$'s (to be compared with 2.2.22 3)). We always prefer to use $\vep$-types which are in the standard form (cf. also 2.2.22 3)). For instance, we prefer to consider the $\vep$-type $(1,1,0,0)$ rather than the $\vep$-type $(0,0,1,1)$ and to consider the $\vep$-type $(\ell,\ell-1,0,1,0)$ rather than $(\ell-1,0,1,0,\ell)$ or $(1,0,\ell,\ell-1,0)$ (here $\ell\ge 4$).
\medskip
{\bf 3.10.4.4. Spreadings.} We define $A:=\{\vep_i|i\in I_1\}$ and $B:=\{\vep_i|i\in\tilde I_1\}$. We have $\abs{A}\le\abs{B}\le e\abs{A}$. The numbers $\abs{A}$ and $\abs{B}$ are called the apparent spreading and respectively the spreading of $({\got g}_0,\vph_1)$. 
\medskip
{\bf 3.10.5. Definitions.}
If $\vep_i\ne 0$, $\forall i\in S(1,m)$, then $({\got g}_0,\vph_1)$ is said to be of totally non-compact type. If there is $i\in I_0$ such that $\vep_i=0$, then we say $({\got g}_0,\vph_1)$ has compact factors or is of compact type. 
\smallskip
If all $\vep_i$, $i\in S(1,m)$, are equal and different from $0$ (i.e. when $\abs{B}=1$ and $I_0=I_1$), we say $({\got g}_0,\vph_1)$ is totally non-compact of constant type, or totally non-compact of type:
\medskip
i) 1 if ${\got g}_0$ is of  $B_\ell$, $C_\ell$ ($\ell\ge 1$), $D_4$, $E_6$ or  $E_7$ Lie type, or of $D_\ell$ Lie type ($\ell\ge 5$) with $\vep_1=1$.
\smallskip
ii) 2 if ${\got g}_0$ is of  $D_\ell$ Lie type, with $\ell\ge 5$, and $\vep_1\in\{\ell-1,\ell\}$.
\smallskip
iii) $\min\{\vep_1,1+\ell-\vep_1\}$ if ${\got g}_0$ is of $A_\ell$ Lie type.
\medskip
In case i) for ${\got g}_0$ of $D_\ell$ Lie type, with $\ell\ge 4$, we also say $({\got g}_0,\vph_1)$ is totally non-compact of $D^\RR_\ell$ type. In case ii) we also say $({\got g}_0,\vph_1)$ is totally non-compact of simple $D^{\HH}_\ell$ type. Similarly, if all non-zero $\vep_i$'s are equal ($i\in S(1,m)$) but there is $i_0\in I_0$ such that $\vep_{i_0}=0$, we say $({\got g}_0,\vph_1)$ is of constant type, or of type 1 (or of $D_{\ell}^{\RR}$ type), or of type 2 (or of simple $D_{\ell}^{\HH}$ type), or of type $\min\{\vep_i,1+\ell-\vep_i\}$ (with $i\in I_0$ such that $\vep_i\neq 0$) depending on which of the situations i) to iii) we are, but dropping (sometimes replacing by the expression compact) the expression of totally non-compact type.
\medskip
If $\abs{B}\ge 2$ we say $({\got g}_0,\vph_1)$ is of non-constant type. This can happen only when ${\got g}_0$ is of $A_\ell$ ($\ell\ge 2$), $D_\ell$ ($\ell\ge 4$) or $E_6$ Lie type. 
\smallskip
If ${\got g}_0$ is of $D_\ell$ Lie type, with $B=\{\ell-1,\ell\}$, or if ${\got g}_0$ is of $D_4$ Lie type, with $\abs{B}=2$, we say $({\got g}_0,\vph_1)$ is of non-simple $D_\ell^\HH$ type. The non-simple and the simple $D_{\ell}^{\HH}$ types are also referred as $D_{\ell}^{\HH}$ types. In the cases when ${\got g}_0$ is of $D_4$ Lie type with $\abs{B}=3$, or of $D_\ell$ Lie type, with $\ell\ge 5$ and  $\abs{B}\ge 2$ but $B\ne\{\ell-1,\ell\}$, we say $({\got g}_0,\vph_1)$ is of mixed $D_\ell$ type. If ${\got g}_0$ is of $D_\ell$ Lie type with $\abs{B}=3$, that is with $B=\{1,\ell-1,\ell\}$, we say $({\got g}_0,\vph_1)$ is of strongly mixed $D_\ell$ type.
\smallskip
If ${\got g}_0$ is of $A_\ell$ Lie type with $\abs{B}\ge 2$, we say $({\got g}_0,\vph_1)$
is of $\{a,b\}$ $A_\ell$ type, where $1\le a<b\le\ell$ and $a\le\ell+1-b$ and, if potentially changing $\ell-i$ with $i+1$, $i=\overline{0,\ell-1}$, $a$ is the smallest element of $B$ and $b$ is the greatest element of $B$. 
\smallskip
$q_1:=\abs{I_1}$ is called the non-compact length of $({\got g}_0,\vph_1)$ and $q_0:=n=\abs{I_0}$ is called the cyclic length of ${\got g}_0$. As we assumed $({\got g}_0,\vph_1)$ is non-trivial (see 3.10.2), we have $q_1\ge 1$. By the rank of $({\got g}_0,\vph_1)$ we mean the rank $\ell\in\NN$ of any simple Lie factor of ${\got g}_0$.  
\smallskip
Most of the above terminology is inspired from and conforms to [De2].
\medskip
{\bf 3.10.5.1. The concentrated $\vep$-type.} The $eq_0$-tuple $\tau^{\rm c}$ obtained from the $\vep$-type $(\vep_1,...,\vep_m)$ by removing all zeros is referred as the concentrated $\vep$-type of $({\got g}_0,\vph_1)$. $\tau^{\rm c}$ is called oscillating if $eq_0\ge 2$ and all its two consecutive elements (in the circular sense) are distinct. If $\abs{B}=2$ and $\tau^{\rm c}$ is oscillating, we say $\tau^{\rm c}$ is alternating. If $\tau^{\rm c}$ is alternating and if $B$ is stable under the group of outer automorphisms of $\Phi_1$ leaving invariant $\Phi_1^+$, we say $\tau^{\rm c}$ is involutive. For the $A_2$ and $E_6$ Lie type any oscillating concentrated $\vep$-type is automatically involutive (but this is not so for the $A_{\ell}$ and $D_{\ell}$ Lie types, with $\ell\ge 3$).    
\medskip 
{\bf 3.10.6. Formulas.} 
The formulas below are based:
\medskip
-- on the list of hermitian symmetric domains and of their real dimension (see [He, p. 518]), and on 
\smallskip
-- the standard tables of roots of the simple, split, semisimple Lie algebras over a field of characteristic $0$ (see [Bou2, planche I-VI]). 
\medskip
There is nothing deep about the computations involved and so we mostly just list the results. All below multiplicities are obtained by counting the roots of $\Phi_1^+$ having in their expression (as sums with positive, integral coefficients of roots of $\Dl_1$) a specific subset of $\Dl_1$.
\medskip
{\bf i) Let ${\got g}_0$ be of $B_\ell$, $C_\ell$ or $E_7$ Lie type.}
\medskip
The slopes of $({\got g}_0,\vph_1)$ are precisely $-{q_1\over q_0}$, $0$ and ${q_1\over q_0}$.
The multiplicity of ${q_1\over q_0}$ is $q_0\dim_{\CC}(X)$, where $X$ is an irreducible hermitian symmetric domain whose group of automorphisms is of the same Lie type as ${\got g}_0$.
\smallskip
If ${\got g}_0$ is of $C_\ell$ Lie type, then $\dim_{\CC}(X)={\ell(\ell+1)\over 2}=\dim_{W(k)}({\got g}_1){{\ell+1}\over {4\ell+2}}$.
\smallskip
If ${\got g}_0$  is of $B_\ell$ Lie type, then $\dim_{\CC}(X)=2\ell-1=\dim_{W(k)}({\got g}_1){{2\ell-1}\over {2\ell^2+\ell}}$.
\smallskip
If ${\got g}_0$  is of $E_7$ Lie type, then $\dim_{\CC}(X)=27={27\over 133}\dim_{W(k)}({\got g}_1)$.
\medskip
{\bf ii) Let ${\got g}_0$ be of $E_6$ Lie type.}
\medskip
{\bf Case 1.} $({\got g}_0,\vph_1)$ is of type 1; then everything is as in i), with
$\dim_{\CC}(X)=16={8\over 39}\dim_{W(k)}({\got g}_1)$.
\medskip
{\bf Case 2.} $({\got g}_0,\vph_1)$ is of non-constant type, without a cycle and with
$\abs{I_2}\ne\abs{I_1\bsl I_2}$, where $I_2:=\{i\in I_1|\vep_i=1\}$. We can assume
$q_2:=\abs{I_2}>\abs{I_1\bsl I_2}$. The slopes of 
$({\got g}_0,\vph_1)$ are precisely $-{q_1\over q_0}$, $-{q_2\over q_0}$, ${-q_1+q_2\over q_0}$, $0$, ${q_1-q_2\over q_0}$, ${q_2\over q_0}$ and ${q_1\over q_0}$.
\smallskip
The multiplicities of the slopes ${q_1\over q_0}$, ${q_2\over q_0}$ and ${q_1-q_2\over q_0}$ are all equal to $8q_0$.
\medskip
{\bf Case 3.} $({\got g}_0,\vph_1)$ is of non-constant type with involution or is without a cycle but with $\abs{I_2}=\abs{I_1\bsl I_2}$ (with $I_2$ as in Case 2). The slopes of 
$({\got g}_0,\vph_1)$ are precisely ${-q_1\over q_0}$, ${-q_1\over 2q_0}, 0, {q_1\over 2q_0}$ and ${q_1\over q_0}$.
\medskip
The multiplicity of ${q_1\over 2q_0}$ is $16q_0$. The multiplicity of ${q_1\over q_0}$ is $8q_0$.
\medskip
In all these three cases we have at most $3$ positive slopes.
\medskip
{\bf iii) Let ${\got g}_0$ be of $D_\ell$ Lie type.}
\medskip
{\bf Case 1.} $({\got g}_0,\vph_1)$ is of type 1 (i.e. of $D_\ell^\RR$ type); then the situation is as in i), with $\dim_{\CC}(X)=2(\ell-1)=\dim_{W(k)}({\got g}_1){{2\ell-2}\over {2\ell^2-\ell}}$.
\medskip
{\bf Case 2.} $({\got g}_0,\vph_1)$ is of type $2$ (i.e. it is of simple $D_l^{\HH}$ type); then the
situation is as in i), with $\dim_{\CC}(X)={\ell(\ell-1)\over 2}=\dim_{W(k)}({\got g}_1){{\ell-1}\over {4\ell-2}}$.
\medskip
{\bf Case 3.} $({\got g}_0,\vph_1)$ is of non-simple $D_\ell^\HH$ type. We can assume $B=\{\ell-1,\ell\}$ even for $\ell=4$. Let $I_2:=\{i\in I_1|\vep_i=\ell-1\}$. We can assume $q_2:=\abs{I_2}\ge\abs{I_1\bsl I_2}$. 
\smallskip
{\bf Subcase 1.} $({\got g}_0,\vph_1)$ is without a cycle and $q_2\ne q_1-q_2$. The slopes are like in Case 2 of ii). 
\smallskip
The multiplicity of the slope ${q_1\over q_0}$ is
${(\ell-1)(\ell-2)\over 2}q_0$. 
\smallskip
The multiplicities of the slopes ${q_2\over q_0}$ and ${{-q_2+q_1}\over q_0}$ are equal to $(\ell-1)q_0$.
\smallskip
{\bf Subcase 2.} $({\got g}_0,\vph_1)$ is without a cycle  and $\abs{I_2}=\abs{I_1\bsl I_2}$ or is with involution. The slopes are like in Case 3 of ii). \smallskip
The multiplicity of
${q_1\over q_0}$ is ${(\ell-1)(\ell-2)\over 2}q_0$. 
\smallskip
The multiplicity of ${q_1\over 2q_0}$ is $2(\ell-1)q_0$.
\medskip
{\bf Case 4.} $({\got g}_0,\vph_1)$ is of mixed $D_{\ell}$ type.
\smallskip
{\bf Subcase 1.} $({\got g}_0,\vph_1)$ is of strongly mixed type.
Let $I_2:=\{i\in I_1|\vep_i=1\}$, $I_3:=\{i\in I_1|\vep_i=\ell-1\}$ and $I_4:=\{i\in I_1|\vep_i=\ell\}$. Let $q_j:=\abs{I_j}$, $j=\overline{2,4}$. We have $q_1=q_2+q_3+q_4$.
We can assume $q_3\ge q_4$.
\smallskip
{\bf Subcase 1.1.}  $({\got g}_0,\vph_1)$ is without a cycle. So $q_j\ge 1$, $\forall j\in\{2,3,4\}$. The set of slopes
is formed from $0$ and $\pm a$, where $a$ is of the form ${\sum_{j\in\Mj_0}q_j\over q_0}$, with $\Mj_0$ a non-empty subset of $\{2,3,4\}$.
\smallskip
The multiplicity of 
$$a={\sum_{j\in\Mj_0}q_j\over q_0}$$
 is equal to 
$$\bigl(\sum_{\Mj\in\Mj_1}r_\Mj\bigr)q_0,$$ 
where $\Mj_1$ is the set of non-empty subsets of $\{2,3,4\}$, with the property that $\sum_{j\in\Mj_1}q_j=aq_0$. Here $r_{\{2,3,4\}}=\ell-2$, $r_{\{2,3\}}=r_{\{2,4\}}=1$, $r_{\{3,4\}}={(\ell-3)(\ell-2)\over 2}$ and $r_{\{2\}}=r_{\{3\}}=r_{\{4\}}=\ell-2$. In other words, $r_{\{2,3,4\}}$ is the number of roots of $\Phi_1^+$ having $\alpha_1$, $\alpha_{\ell}$ and $\alpha_{\ell-1}$ in their expression (as a sum of elements of $\Delta_1$), $r_{\{2,3\}}$ is the number of roots of $\Phi_1^+$ having $\alpha_1$ and $\alpha_{\ell-1}$ while not having $\alpha_{\ell}$ in their expression, etc. (we recall that we are using the notations of 3.4, cf. 3.10.3). 
\medskip
{\bf Subcase 1.2.} $({\got g}_0,\vph_1)$ is with a 3-cycle. This implies $\ell=4$. The slopes are $-{q_1\over q_0}, {-2q_1\over 3q_0}, {-q_1\over 3q_0},0,{q_1\over 3q_0}, {2q_1\over 3q_0}$ and ${q_1\over q_0}$. 
\smallskip
The multiplicity of the slope ${q_1\over q_0}$ is $2q_0$. 
\smallskip
The multiplicity of the slope ${2q_1\over 3q_0}$ is $3q_0$.
\smallskip
The multiplicities of the slopes ${q_1\over 3q_0}$ and $0$ are $6q_0$.
\smallskip
{\bf Subcase 1.3.} $({\got g}_0,\vph_1)$ is with involution. Let ${\rm Inv}$ be the involution of the basis of roots of ${\got g}_1$ (w.r.t. ${\got g}_1\cap {\got b}$) which takes $\al_{\ell-1}$ into $\al_{\ell}$ (this fits the assumption $q_3\ge q_4$ for  $\ell=4$). This implies $q_2\ge 1$. The slopes are $0,\pm a$, with $a$ of the form $\bigl(\sum_{j\in\Mj_0}q_j\bigr){1\over 2q_0}$, for 
$$\Mj_0\in SS:=\{\{2,3,4,5,6,7\},\{2,3,4,5\},\{3,4,6,7\},\{2,5\},\{3,4\}\}.$$ 
Here $q_{3+i}=q_i$ for $i=\overline{2,4}$, while the elements of $SS$ are keeping track of the orbits of the natural action of ${\rm Inv}$ on the non-empty subsets of the set $\{\al_1,\al_{\ell-1},\al_{\ell}\}$; so $\{3,4\}$ corresponds to the orbit formed by the two sets $\{\al_{\ell-1}\}$ and $\{\al_{\ell}\}$, $\{2,5\}$ corresponds to the orbit formed by the set $\{\al_1\}$, $\{3,4,6,7\}$ corresponds to the orbit formed by the set $\{\al_{\ell-1},\al_{\ell}\}$, $\{2,3,4,5\}$ corresponds to the orbit formed by the two sets $\{\al_1,\al_{\ell-1}\}$ and $\{\al_1,\al_{\ell}\}$, while $\{2,3,4,5,6,7\}$ corresponds to the orbit formed by the set $\{\al_1,\al_{\ell-1},\al_{\ell}\}$.
\smallskip
The multiplicity of 
$$a={1\over 2q_0}\bigl(\sum_{j\in\Mj_0}q_j\bigr)$$ 
is 
$$\bigl(\sum_{\Mj_1\in\Mj_2}r_{\Mj_1}\bigr)q_0,$$ 
where $\Mj_2:=\{\Mj\in SS|\sum_{j\in\Mj}q_j=2q_0a=\sum_{j\in\Mj_0}q_j\}$. Here $r_{\Mj_1}$ is $\ell-2$ if $\Mj_1=\{2,3,4,5,6,7\}$ or $\Mj_1=\{2,5\}$, is 2 if $\Mj_1=\{2,3,4,5\}$, is ${(\ell-3)(\ell-2)\over 2}$ if $\Mj_1=\{3,4,6,7\}$ and is $2(\ell-2)$ if $\Mj_1=\{3,4\}$. These numbers $r_{\Mj_1}$ have an entirely similar significance as the similar numbers of Subcase 1.1.
\smallskip
{\bf Subcase 2.} $({\got g}_0,\vph_1)$ is not of strongly mixed type. So $({\got g}_0,\vph_1)$ is without an involution and $\ell\ge 5$. This subcase is entirely analogous to Subcase 1.1, for the value $q_4=0$.
\medskip
In all these four cases the number of positive slopes is at most $7$.
\medskip
{\bf iv) Let ${\got g}_0$ be of $A_\ell$ Lie type $(\ell\ge 2)$.}
\medskip
{\bf Case 1. $\abs{B}=1$.} This happens when $\abs{A}=1$ and $({\got g}_0,\vph_1)$ is
without a cycle or when $({\got g},\vph_1)$ is with involution and $A=\{{\ell+1\over 2}\}$ (in this second situation $\ell$ must be odd).
\medskip
The slopes are $-{q_1\over q_0}$, $0$ and ${q_1\over q_0}$.
\medskip
The multiplicity of the slope ${q_1\over q_0}$ is $i_0(\ell+1-i_0)q_0$, where $B=\{i_0\}$. It is less or equal to ${{q_0(\ell+1)^2}\over 4}=\dim_{W(k)}({\got g}_0){{\ell^2+2\ell+1}\over {4(\ell^2+2\ell)}}$.
\medskip
{\bf Case 2. $\abs{B}=2$.} We distinguish three subcases.
\medskip
{\bf Subcase 1.} $A=\{i_0\}$, with $i_0\ne{\ell+1\over 2}$, and $({\got g}_0,\vph_1)$
is with involution. We can assume $i_0<{\ell+1\over 2}$. The slopes are
$-{q_1\over q_0},{-q_1\over 2q_0},0,{q_1\over 2q_0}$ and ${q_1\over q_0}$.
\smallskip
The multiplicity of ${q_1\over q_0}$ is $i^2_0q_0$.
\smallskip
The multiplicity of ${q_1\over 2q_0}$ is $2i_0(\ell+1-2i_0)q_0$.
\smallskip
{\bf Subcase 2.} $A=\{i_0,i_1\}$, with $1\le i_0\le{\ell+1\over 2}$, $i_0<i_1$, and $({\got g}_0,\vph_1)$ is without a cycle. Let $I_2:=\{j\in I_1|\vep_j=i_0\}$ and let $q_2:=\abs{I_2}$.
\medskip
{\bf Subcase 2.1.} $2q_2\neq q_1$. The slopes are $\pm{q_1\over q_0}$, $\pm{q_2\over q_0}$, $\pm{q_1-q_2\over q_0}$ and $0$.
\smallskip
The multiplicity of ${q_1\over q_0}$ is $i_0(\ell+1-i_1)q_0$.
\smallskip
The multiplicity of ${q_2\over q_0}$ is $i_0(i_1-i_0)q_0$.
\smallskip
The multiplicity of ${q_1-q_2\over q_0}$ is $(\ell+1-i_1)(-i_0+i_1)q_0$.
\smallskip
{\bf Subcase 2.2.} $2q_2=q_1$. The slopes are $\pm{q_1\over q_0},\pm{q_1\over2q_0}$ and $0$.
\smallskip
The multiplicity of ${q_1\over q_0}$ is $i_0(\ell+1-i_1)q_0$.
\smallskip
The multiplicity of ${q_1\over 2q_0}$ is $(i_1-i_0)(\ell+1-i_1+i_0)q_0$.
\smallskip
{\bf Subcase 3.} $A=\{i_0,\ell+1-i_0\}$ and $({\got g}_0,\vph_1)$ is with involution. This
is entirely similar to Subcase 1. 
\medskip
{\bf Case 3. $\abs{B}\ge 3$.} So $\abs{A}\ge 2$. We do not treat this case here. We just mention two things. The multiplicity of the greatest positive slope ${q_1\over q_0}$ is $i_0(\ell+1-i_1)$, where $i_0$ (resp. $i_1$) is the smallest (resp. the greatest) element of $B$. The number of positive slopes is at most $${{\abs{B}^2+\abs{B}}\over 2}.$$
We get this by counting the equivalence classes of the following equivalence relation $r_B$ on $\Phi_1^+$. $\al$, $\be\in\Phi_1^+$ are in relation $r_B$ iff in their expression as a sum of elements of $\Delta_1$ $\al_i$ shows up with the same coefficient, $\forall i\in B$: their number is exactly the number of positive slopes of the $A_{\abs{B}}$ Lie type. Moreover, if $({\got g}_0,\vph_1)$ is with involution then in fact the number of positive slopes is at most
$${{\abs{B}^2+\abs{B}+2[{{\abs{B}+1}\over 2}]}\over 4},$$
i.e. it is at most equal to the number of orbits of the natural action of the involution of $\Phi_1$ leaving $\Phi^+$ invariant on the set of such equivalence classes. If $q_1$ is big enough w.r.t. $\abs{B}$ one can check immediately that both these upper bounds can be attained.
\medskip
{\bf 3.10.7. Some facts.} If the $\vep$-type of $({\got g}_0,\vph_1)$ is  $(\vep_1,\ldots,\vep_m)$, then ${\got p}_0:=W_0({\got g}_0,\vph_1)$ (resp. $W^0({\got g}_0,\vph_1)$) is a direct sum of its intersection with the factors ${\got g}_i$ of ${\got g}_0$, $i\in I_0$, while its intersection with ${\got g}_i$ is the $W(k)$-module generated by ${\got t}_i$ and by the rank 1 $W(k)$-submodules corresponding --in the sense of 3.4.3.2-- to all roots of $\Phi^+_i$ (resp. of $\Phi^-_i$) and to those negative (resp. positive) roots not having $-\al_{\vep_j}(i)$ (resp. $\al_{\vep_j}(i)$) in their expression, for any $j\in S(1,m)$ such that $\vep_j\ne 0$; this is a consequence of 3.4.3.0 via the numbering in 3.4.3.2. 
\smallskip
Using this we get directly the following Facts:
\medskip
{\bf F0.} {\it We fix the Lie type of ${\got g}_0$. Then ${1\over q_0}$ times the number of slopes $0$ of $({\got g}_0,\vph_1)$ depends only on the set of non-zero elements $\vep_j$'s ($j\in S(1,m)$).}
\medskip
{\bf F1.} {\it ${\got p}_0\cap{\got g}_i$ is a maximal parabolic Lie subalgebra of ${\got g}_i$ iff $({\got g}_0,\vph_1)$ is of constant type.} 
\medskip
{\bf F2.} {\it ${\got p}_0=F^0({\got g}_0)$
iff $({\got g}_0,\vph_1)$ is totally non-compact of constant type (i.e. iff $({\got g}_0,\vph_1)$ is an ordinary Shimura Lie $\sg$-crystal).}
\medskip
{\bf F3.} {\it ${\got p}_0$ is a Borel Lie subalgebra of $\got g$ iff ${\got g}_0$ is of $A_\ell$ Lie type with $B=S(1,\ell)$. Here $\ell\ge 1$.}
\medskip
{\bf F4.} {\it The Lie stable $p$-rank of $({\got g}_0,\vph_1)$ is always positive and divisible by $q_0$. It is $q_0$ (i.e. the refined Lie stable $p$-rank $({\got g}_0,\vph_1)$ is the $1$-tuple $(1)$) iff ${\got g}_0$ is of $A_\ell$ Lie type with $\{1,\ell\}\subset B$.}
\medskip
When the situation of $F3$ takes place, we say that $\bigl({\got g}_0,\vph_1,F^0({\got g}_0),F^1({\got g}_0)\bigr)$ is of strong Borel
type (cf. 2.2.12 c)) or that it has maximal $A$-spreading. Also, for future references we state explicitly the following obvious fact:
\medskip
{\bf F5.} {\it $({\got g}_0,\vph_1)$ is totally non-compact iff its multiplicity of slope $-1$ (or $1$) is non-zero.}
\medskip
{\bf 3.10.8. An interpretation in terms of $p$-divisible groups.}
We consider a $p$-divisible group $D$ over $W(k)$ having $({\got b},F^1({\got g}_0),\vph_1)$ as its associated filtered $\sg$-crystal, cf. 2.2.12.1 1). The positive slopes and their multiplicities as listed in 3.10.6, are precisely the positive slopes and their multiplicities of the special fibre $D_k$ of $D$. 
\medskip
{\bf 3.10.8.1. The $a$-invariant.} Let $a({\got g}_0,\vph_1):=a(D_k)$ (cf. 2.1). We refer to it as the $a$-invariant of $({\got g}_0,\vph_1)$. As the $a$-number is not changed by field extensions, we similarly define the $a$-invariant of any Shimura Lie $\sg$-crystal which has a lift of parabolic type.
\medskip
{\bf Exercise.} Show that the $a$-invariant of $({\got g}_0,\vph_1)$ is the same as the $a$-invariant of the extension to $\bar k$ of its canonical lift as defined in 2.3.18.3 2). Hint: use Exercise of 2.2.22.1.
\medskip\smallskip
{\bf 3.11. The case $k=\bar k$ and applications.} We use the notations of 3.1-2. We assume $(M,\vph,G,(t_{\al})_{\al\in\Mj})$ is a $G$-ordinary $\sg$-crystal.
\medskip
{\bf 3.11.1. Theorem.} {\it Let ${\got C}=(M,F^1,\vph,G,{(t_\al)}_{\al\in\Mj})$ be the $G$-canonical lift of $(M,\vph,G,(t_{\al})_{\al\in\Mj})$. We assume $k=\bar k$. We have:
\medskip
{\bf a)} ${\got C}$ is cyclic diagonalizable;
\smallskip
{\bf b)} its Shimura filtered Lie $\sg$-crystal is cyclic diagonalizable and of Borel and parabolic type;
\smallskip
{\bf c)} it is $1_{\Mj}$-isomorphic to $\bigl(M,F^1,\vph_1,G,{(t_\al)}_{\al\in\Mj}\bigr)$ mentioned in 3.2.3.}
\medskip
{\bf Proof:} There are many ways in proving c). One way is entirely the same as the proof of b) of 4.4.1 2) (cf. also the extra details needed for the proof of b) of 4.4.1 3)) below, made in the context of a SHS. As the referred proof is self contained (based on 3.4.11), this way is not repeated here (cf. the policy of 1.14.2). A second way is to combine 3.6.17 with 3.4.11. 
\smallskip
For the convenience of the reader we present a third way, based on 3.6.17, 3.4.11 and 3.1.0 c). We can assume $\vph=g\vph_1$, with $g\in P_1(W(k)))$ (see 3.3.2 for the meaning of $P_1$). We need to find $h\in P_1(W(k))$ such that 
$$hg\vph_1=\vph_1h.\leqno (1)$$
The case when $g$ is congruent to the identity mod $p$ is entirely the same as the part of the proof of 3.6.17 referring to the $W(k)$-base $\Mb_0$ and to induction. 
\smallskip
Let $R_1$ be the Zariski closure of the reductive subgroup of $P_{1B(k)}$ whose Lie algebra is $W(0)({\got g},\vph_1)$. From Fact 1 of 2.2.9 1) we get that $(M,F^1,\vph_1,T)$ is cyclic diagonalizable. So from 2.2.19.2 we get that $R_1$ is a Levi subgroup of $P_1$ (more precisely: $R_1$ is the only Levi subgroup of $P_1$ containing $T$). So $P_1$ is the semidirect product of its unipotent radical $N_1$ and of $R_1$.
To treat the general case, we write
$$g=g_1n_1$$
with $g_1\in R_1(W(k))$ and $n_1\in N_1(W(k))$.  As $k=\bar k$, $R_1$ has a canonical $\ZZ_p$-structure: the part of 3.4.3.0 (1) referring to $s_i(\al)$'s implies that $\mu$ centralizes $R_1$ and so this $\ZZ_p$-structure is obtained as in 2.2.9 8); so, as $R_1$ is canonically identified with $P_1/N_1$, [Bo2, 16.4] implies that we can assume $g_1$ is congruent to the identity mod $p$. $N_1$ is a normal subgroup of $P_1$ normalized by $T$. Moreover it is unipotent; so there is a finite sequence of characteristic subgroups 
$$N_m\subset N_{m-1}\subset ...\subset N_1$$ 
of $N_1$ such that $N_i/N_{i+1}=\GG_a^{m_i}$, $\forall i\in S(1,m)$, with $m$, $m_i\in\NN\cup\{0\}$ and with $N_{m+1}$ as the trivial subgroup. Let $U_i$ be a subgroup of $N_i$ normalized by $T$ and such that the resulting homomorphism $U_i\to N_i/N_{i+1}$ is an isomorphism. We have $\vph({\rm Lie}(U_i))\subset {\rm Lie}(U_i)$ (cf. 3.4.3.0 (1) and the fact that each $N_i$ is a characteristic subgroup of $N_1$). Moreover, $\forall i\in S(1,m)$, we have 
$$\vph^{l_i}({\rm Lie}(U_i))\subset p{\rm Lie}(U_i),$$
for some $l_i\in\NN\cup\{0\}$ (cf. 3.4.3.0 (1); this is just the version of the Key Fact of the proof of 3.6.6 which involves positive --versus negative-- slopes). Using this, by repeatedly replacing $g\vph$ with $n_2g\vph n_2^{-1}$, with $n_2\in U_i(W(k))$, where $i\in S(1,m)$ is the biggest number such that $n_1$ mod $p$ is a $k$-valued point of $N_i$, we get that we can assume $n_1$ is congruent to the identity mod $p$. So $g$ is congruent to the identity mod $p$. This proves c).
\smallskip
To prove a) and b), we recall (cf. 3.2.3) that $\vph_1({\got t}\otimes_{W(k)} W(\bar k))={\got t}\otimes_{W(k)} W(\bar k)$. So, as $k=\bar k$, {\got t} is generated by elements fixed by $\vph_1$. Considering the family of tensors of $\Mt(M[{1\over p}])$ formed by such elements and by the members of the family $(t_{\al})_{\al\in\Mj}$, we get a concrete way to see that the quadruple $(M,F^1,\vph_1,T)$ is a Shimura filtered $\sg$-crystal. 
So a) results from this and from c) (cf. 2.2.16). b) results from c) and the fact that $k=\bar k$ (cf. 3.4.3.0 (1) and 2.2.12.1 4)). This proves the Theorem.  
\medskip
{\bf 3.11.2. First applications. A.} From 3.11.1 c) we get 3.1.4. This goes as follows. The uniqueness part is implied by 3.2.8; it allows us to assume that in 3.1.4 we have $k=\bar k$. But referring to $(M,F^1,\vph_1,G)$, we can take $\tilde\mu=\mu$, cf. 3.4.3.0. 
\medskip
{\bf B.} Let now $\mu_y:\GG_m\to G$ be the canonical split of $(M,\vph,G)$ (cf. 3.1.6). We write $\vph=a\circ \mu_y({1\over p})$, with $a$ as a $\sg$-linear automorphism of $M$. We assume $k=\bar k$. So $M^a:=\{x\in M|a(x)=x\}$ defines a $\ZZ_p$-structure on $M$. Let $\bar h_1\in {\rm End}(M)$ be defined by: $\bar h_1(x)=x$ if $x\in F^1$ and $\bar h_1(x)=0$ if $\GG_m$ through $\mu_y$ acts trivially on $x$. We have $\bar h_1\in {\rm End}\bigl(M^a\otimes_{\ZZ_p} W(\FF)\bigr)$ and  there is $d\in\NN$ such that $\vph^d(\bar h_1)=\bar h_1$ but $\forall i\in \NN$, $i<d$, we have $\vph^i(\bar h_1)\ne\bar h_1$; all these are a consequence of 3.11.1 a) (cf. also 2.2.16.4).
So $\mu_y$ can be viewed as an injective cocharacter 
$$\mu_y:\GG_m\hookrightarrow G_{W(\FF_{p^d})}\subset GL\bigl(M^a\otimes_{\ZZ_p} W(\FF_{p^d})\bigr).$$ 
Let $F^1_d:=F^1\cap M^a\otimes_{\ZZ_p} W(\FF_{p^d})$; we have $F^1=F^1_d\otimes_{W(\FF_{p^d})} W(k)$. Let $\vph^i(\bar h_1)=\bar h_{i+1}$, $\forall i\in \NN$. So $h_{i+1}=a^i(\bar h_1)$ and $\bar h_{d+1}=\bar h_1$. The $\bar h_i$'s are commuting among themselves (as they belong to ${\rm Lie}(T)$). As in 2.2.9 8), let $G_{\ZZ_p}$ be the subgroup of $GL(M^a)$ whose extension to $W(k)$ is the subgroup $G$ of $GL(M^a\otimes_{\ZZ_p} W(k))=GL(M)$. 
\smallskip
We get (cf. 2.2.16.4 a)): $(M,F^1,\vph,G,(t_{\al})_{\al\in\Mj})$ is isomorphic to the extension to $W(k)$ of the following strongly cyclic diagonalizable Shimura filtered $\sg_{\FF_{p^d}}$-crystal
$$\Mm\Mm:=\bigl(M^a\otimes_{\ZZ_p} W(\FF_{p^d}),F^1_d,\sg_{\FF_{p^d}}\circ\mu_y({1\over p}),G_{W(\FF_{p^d})},(t_{\al})_{\al\in\Mj}\bigr);$$ 
here $\sg_{\FF_{p^d}}$ acts trivially on $M^a$ and as Frobenius on $W(\FF_{p^d})$. We refer to $\Mm\Mm$ as the minimal model of $(M,F^1,\vph,G,(t_{\al})_{\al\in\Mj})$, cf. also Fact of 2.2.16.4; from its construction we get it is uniquely determined. 
\smallskip
$\bigl({\sg_{\FF_{p^d}}\circ\mu_y({1\over p})}\bigr)^d$ is a $B(\FF_{p^d})$-linear automorphism of $M^a\otimes_{\ZZ_p} B(\FF_{p^d})$; so it can be identified with an element of $G_{W(\FF_{p^d})}(B(\FF_{p^d}))$ and so of $G(K_0)$. As $\bar h_i$'s are commuting among themselves, we get immediately that in fact $\bigl({\sg_{\FF_{p^d}}\circ\mu_y({1\over p})}\bigr)^d$ is a $\QQ_p$-linear automorphism of $M^a[{1\over p}]$ defined by a $\QQ_p$-valued point of the smallest torus of $G_{\QQ_p}$ through whose extension to $B(\FF_{p^d})$ $\mu_y$ factors; in particular, it is semisimple. 
\medskip
{\bf C.} Let $P_{=0}$ be the subgroup of $G_{\ZZ_p}$ whose extension to $W(\FF_{p^d})$ is the centralizer of $\bar h_1$,..., $\bar h_d$; it is connected as its extension to $W(\FF_{p^d})$ is the centralizer of the torus of $G_{W(\FF_{p^d})}$ generated by the image of $\mu$ and by the conjugates under $\sg_{\FF_{p^d}}$ of this image (cf. [Bo2, 11.12] applied over $B(\FF_{p^d})$ and over $\FF_{p^d}$). ${\rm Lie}(P_{=0})\subset {\got g}$ is the Lie algebra over $\ZZ_p$ of endomorphisms of $(M,\vph,G)$ or of $(M,F^1,\vph,G)$ (one inclusion is obvious from constructions; the other one results from 3.1.1.2 and 3.1.4-5); even more, its extension to $W(\FF_{p^d})$ is $W(0)({\rm Lie}(G_{W(\FF_{p^d})}),\sg_{\FF_{p^d}}\circ\mu_y({1\over p}))$: all these can be read out as well from 3.4.3.0 (via 3.11.1 c)). From 3.4.3.0 (1) we also get that we have the following slope sign direct sum decomposition
$$({\rm Lie}(G_{W(\FF_{p^d})}),F^0({\rm Lie}(G_{W(\FF_{p^d})})),F^1({\rm Lie}(G_{W(\FF_{p^d})})),\sg_{\FF_{p^d}}\circ\mu_y({1\over p}))=\leqno (SSDSD)$$
$$=({\rm Lie}(P_{=0W(\FF_{p^d})}),{\rm Lie}(P_{=0W(\FF_{p^d})}),0,\sg_{\FF_{p^d}}\circ\mu_y({1\over p}))\oplus$$ 
$$\oplus({\rm Lie}(N_{>0}),{\rm Lie}(N_{>0}),F^1({\rm Lie}(G_{W(\FF_{p^d})}))\cap {\rm Lie}(N_{>0}),\sg_{\FF_{p^d}}\circ\mu_y({1\over p}))\oplus$$
$$\oplus ({\rm Lie}(N_{<0}),F^0({\rm Lie}(G_{W(\FF_{p^d})}))\cap {\rm Lie}(N_{<0}),0,\sg_{\FF_{p^d}}\circ\mu_y({1\over p}))$$
as $p$-divisible objects of $\Mm\Mf_{[-1,1]}(W(\FF_{p^d}))$: the first (resp. the second and the third) direct sum factor is in fact a $p$-divisible object of $\Mm\Mf_{[0,0]}(W(\FF_{p^d}))$ (resp. is a $p$-divisible object of $\Mm\Mf_{[0,1]}(W(\FF_{p^d}))$ and of $\Mm\Mf_{[-1,0]}(W(\FF_{p^d}))$); here $N_{<0}$ (resp. $N_{>0}$) is the integral, connected, smooth, unipotent subgroup of $G_{W(\FF_{p^d})}$ whose Lie algebra is the Lie subalgebra of ${\rm Lie}(G_{W(\FF_{p^d})})$ corresponding to negative (resp. positive) slopes of $({\rm Lie}(G_{W(\FF_{p^d})}),\sg_{\FF_{p^d}}\circ\mu_y({1\over p}))$. $N_{>0}$ and $N_{<0}$ have natural $\ZZ_p$-structures: they are obtained from subgroups of $G_{\ZZ_p}$ by extension of scalars (cf. end of B above).
\smallskip
Even better, 3.4.3.0 (1) implies: we actually get slope type direct sum decompositions 
$$({\rm Lie}(G_{W(\FF_{p^d})}),F^0({\rm Lie}(G_{W(\FF_{p^d})})),F^1({\rm Lie}(G_{W(\FF_{p^d})})),\sg_{\FF_{p^d}}\circ\mu_y({1\over p}))=\oplus_{\al\in S} {\got C}_{\alpha},\leqno (STDSD)$$ 
where $S$ is the set of slopes of $({\rm Lie}(G_{W(\FF_{p^d})}),\sg_{\FF_{p^d}}\circ\mu_y({1\over p}))$ and ${\got C}_{\alpha}$ is a $p$-divisible object of $\Mm\Mf_{[0,1]}(W(\FF_{p^d}))$ if $\al>0$, of $\Mm\Mf_{[0,0]}(W(\FF_{p^d}))$ if $\al=0$, and of $\Mm\Mf_{[-1,0]}(W(\FF_{p^d}))$ if $\al<0$; moreover the Newton polygon of ${\got C}_{\alpha}$ has only one slope: $\alpha$. So, STDSD ``isolates" entirely all slopes of $S$ and not only their signs (as the decomposition (SSDSD) does). Declaring ${\got C}_s$ to be the zero $p$-divisible object if $s\in [-2,2]\setminus S$, w.r.t. the Lie structure of ${\rm Lie}(G_{W(\FF_{p^d})})$ we have natural inclusions
$$[{\got C}_{\al},{\got C}_{\be}]\subset {\got C}_{\al+\be}$$
of $p$-divisible objects of $\Mm\Mf_{[-2,2]}(W(k))$; here $\al$, $\be\in S$.
\medskip
{\bf D.} Let $B$ be a Borel subgroup of $G$ normalizing $F^1$ and such that $\vph({\rm Lie}(B))\subset {\rm Lie}(B)$ (cf. the Borel type property of 3.11.1 b)). There are many ways to see that $\mu_y$ factors through $B$ (for instance, using the end of 2.2.12.1 1) and the fact that $B$ is its own normalizer in $G$, or using the functoriality of canonical splits, etc.). So $a({\rm Lie}(B))={\rm Lie}(B)$. Let $B_{\ZZ_p}$ be the Borel subgroup of $G_{\ZZ_p}$ which is the $\ZZ_p$-structure of $B$ defined naturally by $a$. So $B_{\ZZ_p}$ normalizes $\vph^i(F^1)$, $\forall i\in\NN$, and so we have: 
\medskip
{\bf Fact.}  {\it $\bar h_i$ is $G(B(k))$-conjugate to $\bar h_1$, iff $\bar h_i=\bar h_1$, i.e. iff $d|i-1$.}    
\medskip
Let $_{A}\bar h_1$ (resp. $_{T}\bar h_1$) be the image of $\bar h_1$ in ${\rm Lie}(G^{\rm ad})$ (resp. in ${\rm Lie}(G^{\rm ab})$),
and let $d_1\in \NN$ (resp. $d_2\in\NN$) be the smallest positive integer such that $_A\bar h_1$ (resp. $_T\bar h_1$) is the image of $\bar h_{d_1+1}$ in
${\rm Lie}(G^{\rm ad})$ (resp. the image of $\bar h_{d_2+1}$ in ${\rm Lie}(G^{\rm ab})$). Then $d_1$ and $d_2$ divide $d$, and in fact $d=l.c.m.[d_1,d_2]$.  
\medskip
{\bf 3.11.3. Degrees of definition.} Referring to 3.11.2 D, the number $d$ is called the degree of definition of $(M,\vph,G)$, $d_1$ is called the degree of definition of $\bigl({\rm Lie}(G^{\rm ad}),\vph\bigr)$ or the $A$-degree of definition of $(M,\vph,G)$, and $d_2$ is called the degree of definition of $\bigl({\rm Lie}(G^{\rm ab}),\vph\bigr)$ or the $T$-degree of definition of $(M,\vph,G)$. Here $A$- and $T$- stand for adjoint and respectively toric. $d$ is the same as the degree of definition of $(M,F^1,\vph)$ as def. in 2.2.16.4. Moreover,  as $N_{>0}$ is a non-trivial group (see 3.2.2.1), the $A$-degree of definition of $(M,\vph,G)$ is equal to the $A$-degree of definition of the (non-trivial) Shimura adjoint filtered Lie $\sg$-crystal attached to $(M,F^1,\vph,G)$ as defined in the end of 2.2.16.5. 
\medskip
{\bf 3.11.3.1. Corollary.} {\it Let $g\in G(W(k))$ and $T_g$ a maximal torus of $G$ be such that $(M,g\vph,T_g)$ is a quasi Shimura $\sg$-crystal. Then the degree of definition $d_g$ of $(M,g\vph,T_g)$ is a multiple of $d$.}
\medskip
{\bf Proof:} This Corollary is very much related to 2.2.16.4; however, it is not implied directly by the mentioned place. We can assume $k=\bar k$; so we can drop the word quasi and moreover we can assume $T_g=T$. Let $\mu_T:\GG_m\to T$ be the unique cocharacter producing a direct sum decomposition $M=F^1_T\oplus F^0_T$ (with $\be\in\GG_m(W(k))$ acting through it on $F^i_T$ as the multiplication with $\be^{-i}$, $i=\overline{0,1}$) such that $(M,F^1_T,g\vph,T)$ is a Shimura filtered $\sg$-crystal. Writing $g\vph=a_T\circ \mu_T({1\over p})$, $a_T$ is a $\sg$-linear automorphism of $M$. As in 3.11.2 B, let $T_{\ZZ_p}$ (resp. $\tilde G_{\ZZ_p}$) be the $\ZZ_p$-structure of $T$ (resp. of $G$) it defines naturally. $\tilde G_{\ZZ_p}$ is an inner form of $G_{\ZZ_p}$ and so from Lang's theorem, we get that in fact we have a natural identification $G_{\ZZ_p}=\tilde G_{\ZZ_p}$. As $T_{\ZZ_p}$ is a torus of $\tilde G_{\ZZ_p}$, from 2.2.16.4 a) we get (cf. also Fact 2 of 2.2.9 3)) that the $G(W(k))$-conjugacy class of the cocharacter $\mu_y$ of 3.11.2 is definable over $W(\FF_{p^{d_g}})$. Let $d_1:=(d_g,d)$. From [Mi3, 4.6-7] we deduce that this class is definable over $W(\FF_{p^{d_1}})$. So, as $\bar h_i=\sg_{\FF_{p^d}}^{i-1}(\bar h_1)$ (see 3.11.2 B), from the Fact of 3.11.2 D we get $d_1=d$. This proves the Corollary.  
\medskip
{\bf 3.11.3.2. Remark.} We consider a potentially cyclic diagonalizable Shimura filtered $\sg_{k_1}$-crystal over a finite field $k_1$. Warning: the degree of definition of its extension to $\overline{k_1}=\FF$ is not necessarily a divisor of the degree $[k_1:\FF_p]$; easy examples can be obtained by considering an elliptic curve $E$ over $\QQ$ which has supersingular good reduction w.r.t. $p$ and is such that $E_{\overline{\QQ}}$ has complex multiplication.
\medskip
{\bf 3.11.4. An application: more on the numbers $d(n)$ of 3.6.1.3 2).} We come back to the notations of 3.1-4. 3.11.1 allows us to define the number 
$$r((M,\vph,G))\in\NN\cup\{0\}$$ 
as the number of slopes $-1$ of the Lie $\sg$-crystal of any $G$-canonical lift produced by $(M,\vph,G)$. Directly from 3.6.18.7.0 we get: $\forall n\in\NN$ we have the estimate 
$$d(n)\le\dim_{W(k)}(G)r((M,\vph,G)).\leqno (EST1)$$ 
The numbers $r((M,\vph,G))$ can be easily computed using 3.10.6. As a very gross estimate we have $3r((M,\vph,G))\le \dim_{W(k)}(G)$ (cf. the formulas of 3.10.6: see i), ii), Cases 1 and 2 of iii), and Case 1 of iv); see 3.10.7 for why the situation gets reduced to the constant type context of 3.10.5). So we get 
$$d(n)\le {\dim_{W(k)}(G)^2\over 3};\leqno (EST2)$$ 
if $G=GL(M)$, then in the above inequality we can replace the denominator $3$ by $4$ (cf. the last part of Case 1 of 3.10.6 iv)). In many situations, we can get much better estimates than the above gross estimate: for instance, see Case 2) of 3.10.6 ii).
\smallskip
In particular, if $r((M,\vph,G))=0$ then $d(n)=0$, $\forall n\in\NN$ (his result matches 3.6.18.8.1 a)). So in such a case we can take $Q_{j,n}=R_j$, $\forall n\in\NN$. There are plenty of examples with $r((M,\vph,G))=0$; we mention here two such situations of general nature:
\medskip
{\bf 1)} the case when none of the adjoint cyclic factors of the $G$-canonical lift $(M,F^1,\vph_0,G)$ of $(M,\vph_0,G)$ are totally non-compact;
\smallskip
{\bf 2)} the situations pertaining to a SHS to be referred in 4.6.1 1) below.
\medskip
{\bf 3.11.5. Remark.} 3.11.1 c) and 3.1.3 a) imply: $(M,\vph,G)$ has constant isogeny class (cf. def. of 2.2.22 2)).
\medskip
{\bf 3.11.6. The Lie variant of 3.11.1.} Referring to the abstract context of 2.2.11 (i.e. we move away from the context of Shimura $\sg$-crystals, even if we use very similar notations) we have the following Lie variant of 3.11.1:
\medskip
{\bf Corollary.} {\it {\bf A)} We assume $k=\bar k$. If $({\got g},\vph,F^0({\got g}),F^1({\got g}))$ is a Shimura filtered Lie $\sg$-crystal of Borel type, then for any $g\in G^{\rm ad}(W(k))$ such that $({\got g},g\vph)$ has the same Newton polygon as $({\got g},\vph)$, the Shimura Lie $\sg$-crystals $({\got g},g\vph)$ and $({\got g},\vph)$ are both isomorphic under an isomorphism defined by an element of $G^{\rm ad}(W(k))$. 
\smallskip 
{\bf B)} A Shimura filtered Lie $\sg$-crystal $({\got g},\vph,F^0({\got g}),F^1({\got g}))$ of Borel type is of parabolic type. If $k=\bar k$, then the converse holds. Moreover, for $k$ arbitrary and $({\got g},\vph)$ non-trivial, there is a unique injective cocharacter $\mu:\GG_m\hookrightarrow G^{\rm ad}$ as in 2.2.11 1) and such that the parabolic Lie subalgebra of ${\got g}$ corresponding to non-positive (resp. to non-negative) slopes of $({\got g},\vph)$ is contained in ${\got g}^0\oplus {\got g}^{-1}$ (resp. in ${\got g}^0\oplus {\got g}^1$); it is the canonical split cocharacter of $({\got g},\vph,F^0({\got g}),F^1({\got g}))$.
\smallskip
{\bf C)} If $k=\bar k$, then any Shimura filtered Lie $\sg$-crystal of Borel type is cyclic diagonalizable.}
\medskip
{\bf Proof:} A) is a consequence of 3.5.5: we just need to use the same argument as in the proof of 3.11.1 c). Using A) and just interpreting (i.e. restating) 3.2.3, 3.2.8, 3.4.3.0 and 3.11.2 A in terms of Lie algebras, B) and C) follow.
\medskip
{\bf 3.11.6.1. Definitions.} A non-trivial Shimura Lie $\sg$-crystal which has a lift of parabolic type, is called a Shimura-ordinary Lie $\sg$-crystal, while this lift (as in 3.2.8, based on the part of 2.2.11.1 referring to 2.2.9 3), we get that it is unique) is called its Shimura-canonical lift. A cocharacter $\mu$ as in B) above is called the canonical split of $({\got g},\vph)$. Similarly, following the pattern of 3.1.1, we speak about generalized Shimura $p$-divisible objects over $k$ which are Shimura-ordinary, as well as about their Shimura-canonical lifts and about their canonical split cocharacters.
\medskip
{\bf 3.11.6.2. Remarks.} {\bf 1)} Referring to the Corollary of 3.11.6, it is worth pointing out another way to get the first sentence of B). Let ${\got p}_{\ge 0}$, ${\got p}_{>0}$ and ${\got p}_{<0}$ be the Lie subalgebras of ${\got g}$ corresponding respectively to non-negative, positive and negative slopes of $({\got g},\vph)$. We can assume $k=\bar k$. The Borel type part implies $F^1({\got g})\subset {\got p}_{>0}$ (cf. 2.2.12.1 1)); so ${\got p}_1\subset {\got p}_2$ (cf. Fact of 2.2.11.1). We conclude (based on 2.2.3 3)) ${\got p}_{\ge 0}\subset F^0({\got g}_0)$.
\smallskip
{\bf 2)} We refer to 2.2.12.1 1). The $p$-divisible group $D^+$ over $W(k)$ whose associated filtered $\sg$-crystal is $({\got p}_{>0},F^1({\got g}),\vph)$ is called the positive $p$-divisible group of $({\got g},\vph)$.
\smallskip
Similarly, as 
$$F^1({\got g}_0)\subset {\got p}_{>0}\subset {\got p}_{\ge 0}\subset F^0({\got g}_0),$$ 
the triple
$$({\got p}_{<0},F^0({\got g})\cap {\got p}_{<0},p\vph)$$ 
is a filtered $\sg$-crystal; warning: here $F^0({\got g})\cap {\got p}_{<0}$ is viewed as an $F^1$-filtration of ${\got p}_{<0}$. The $p$-divisible group $D^-$ over $W(k)$ having it as its associated filtered $\sg$-crystal, is called the negative $p$-divisible group of $({\got g},\vph)$. From the description of the permutation $\gamma_{\Phi}$ of 3.4.3.0, we get: 
\medskip
{\bf Fact.} {\it If $k=\bar k$, then $D^-$ is the dual of $D^+$.} 
\medskip
{\bf 3)} 3.11.6 forces us to come back to 3.10.8. We assume $k=\bar k$ and we use the notations of 1) and 2). 3.11.6 C) implies that the filtered $\sg$-crystal $({\got p}_{>0},F^1({\got g}),\vph)$ is cyclic diagonalizable. It is easy to compute its classification invariants starting from 3.10.6. For instance, if $({\got g},\vph)$ is cyclic of $B_{\ell}$, $C_{\ell}$ or $E_7$ Lie type, from 3.10.6 i) we get: $({\got p}_{>0},F^1({\got g}),\vph)$ is a direct sum of $\dim_{\CC}(X)$-copies of the circular diagonalizable filtered $\sg$-crystal having its type equal to the $\vep$-type of $({\got g},\vph)$, and so we just have to decompose the $\vep$-type into indecomposable types (as in 2.2.22 3)); here $X$ is the hermitian symmetric domain associated to the Lie type of ${\got g}$, while $\ell\in\NN$. 
\smallskip
Similar direct sum decompositions can be obtained for the context of 3.10.6 ii) to iv). However, in the Case 3 of 3.10.6 iv), they can be very complicated; in particular, the $\vep$-type is not enough to describe (directly, i.e. without using [Bou2, planche I]) the isomorphism classes of the direct summands (see D of 2.2.22 3)) of $({\got p}_{>0},F^1({\got g}),\vph)$.
\medskip
{\bf 3.11.7. Groups of automorphisms.} Let $(M,\vph,G,(t_{\al})_{\al\in\Mj})$ be a $G$-ordinary $\sg$-crystal and let $F^1$ be its canonical lift. We first assume $k=\bar k$. Let $\bar h_i$'s, $P_{=0}$ and $\Mm\Mm$ be as in 3.11.2 B and C. Let $\Mm\Mm_{\FF}$ be the extension of $\Mm\Mm$ to $\FF$. The group ${\rm Aut}_0$ (resp. ${\rm Aut}_1$) of automorphisms (resp. of $1_{\Mj}$-automorphisms) of $(M,F^1,\vph,G)$ (resp. of $(M,F^1,\vph,G,(t_{\al})_{\al\in\Mj})$) is the same as the group of automorphisms (resp. of $1_{\Mj}$-automorphisms) of $(M^a\otimes_{\ZZ_p} W(\FF),F^1_d\otimes_{W(\FF_{p^d})} W(\FF),\sg_{\FF_{p^d}}\circ\mu_y({1\over p})\otimes 1,G_{W(\FF)}\bigr)$ (resp. of $\Mm\Mm_{\FF}$). From 3.11.2 C we get that ${\rm Aut}_1$ is the group of $\ZZ_p$-valued points of $P_{=0}$. Similarly, let $P_{=0}^1$ be the maximal closed subgroup of the normalizer of $G_{\ZZ_p}$ in $GL(M^a)$, whose extension to $W(\FF_{p^d})$ centralizes $\bar h_i$, $i=\overline{1,d}$. From 3.1.5.1 we get: 
\medskip
{\bf Corollary.} {\it ${\rm Aut}_0$ is the group of $\ZZ_p$-valued points of $P_{=0}^1$.}
\medskip
Let ${\rm Aut}_2$ be the group of automorphisms of the positive (or negative) $p$-divisible group of $({\got g},F^0({\got g}),F^1({\got g}),\vph)$ or of the positive (or negative) $p$-divisible group of the Shimura filtered $\sg$-crystal attached to $\Mm\Mm_{\FF}$. It can be as well interpreted as the group of  $\ZZ_p$-valued points of an integral, affine group scheme over $\ZZ_p$ (cf. the $\ZZ_p$-structures of $N_{>0}$ and of $N_{<0}$ mentioned in 3.11.2 C). The group ${\rm Aut}_s$, with $s\in\{0,1,2\}$, is referred as the $s$-th group of automorphisms of $(M,F^1,G,(t_{\al})_{\al\in\Mj})$.   
\medskip
{\bf 3.11.8. The case of a perfect field.} With $k$ just perfect, similarly to 3.11.2-3, for $(M,\vph,G,(t_{\al})_{\al\in\Mj})$ we define its degree of definition $d\in\NN$ as well as its $A$-degree and $T$-degree of definition. We have:
\medskip
{\bf Lemma.} {\it $k$ contains $\FF_{p^d}$.}
\medskip
{\bf Proof:} If $k$ is finite, this can be checked using entirely the same ideas as in the proof of 3.11.3.1. To see the general case, it is enough to show that the smallest Lie subalgebra of ${\rm Lie}(G_{B(k)})$ taken by $\vph$ into itself and containing $d\mu$, with $\mu:\GG_m\to G$ as the canonical split of $(M,\vph,G)$, is generated by elements fixed by $\vph$. But this Lie algebra is a Lie subalgebra of the Lie algebra of the maximal torus of the center of $P_{\ge 0}\cap P_{\le 0}$, with $P_{\ge 0}$ (resp. $P_{\le 0}$) as the parabolic subgroup of $G_{B(k)}$ whose Lie algebra is $W_0({\rm Lie}(G_{B(k)}),\vph)$ (resp. is $W^0({\rm Lie}(G_{B(k)}),\vph)$) (this can be seen over $\bar k$; so 3.11.2 C applies). As $P_{\ge 0}$ and $P_{\le 0}$ are parabolic subgroups of $G_{B(k)}$, the Lemma follows by transferring everything (via Fontaine's comparison theory) into the \'etale context with $\QQ_p$-coefficients (cf. also 2.2.12.1 1) and 3.11.6.2 2); see also (SSDSD) of 3.11.2 C). 
\medskip
Let $\Mm\Mm:=\bigl(M^a\otimes_{\ZZ_p} W(\FF_{p^d}),F^1_d,\sg_{\FF_{p^d}}\circ\mu_y({1\over p}),G_{W(\FF_{p^d})},(t_{\al})_{\al\in\Mj}\bigr)$ be the minimal model of the extension of $(M,\vph,G,(t_{\al})_{\al\in\Mj})$ to $\bar k$. Let ${\rm Aut}_s$ be the $s$-th group of automorphisms of this extension.
\smallskip
The Lemma implies that we have a natural action of ${\rm Gal}(k)$ on ${\rm Aut}_i$, $i\in\{0,1,2\}$: it is the trivial action (cf. the existence of $\ZZ_p$-structures in 3.11.7). How ``far" is $(M,F^1,\vph,G)$ (resp. $(M,\vph,G,(t_{\al})_{\al\in\Mj})$) from being isomorphic (resp. $1_{\Mj}$-isomorphic) to the extension of $\Mm\Mm$ to $k$, is measured by a class 
$$\gamma\in H^1({\rm Gal}(k),{\rm Aut}_0)$$
(resp. $\gamma_1\in H^1({\rm Gal}(k),{\rm Aut}_1)$).
\smallskip
As in 3.11.6.2 we define the positive $D^+$ and the negative $D^{-}$ $p$-divisible group of $({\got g},F^0({\got g}),F^1({\got g}),\vph)$. As above, how ``far" is $D^+$ from being isomorphic to the extension to $k$ of the positive $p$-divisible group of the Shimura filtered $\sg$-crystal attached to $\Mm\Mm$, is measured by a class
$$\gamma_2\in H^1({\rm Gal}(k),{\rm Aut}_2).$$
We say $(M,\vph,G)$ or $(M,F^1,\vph,G)$ is balanced if $\gamma_2$ is the $0$ class. It is easy to see that there are examples when $(M,\vph,G)$ is balanced but $\gamma_1$ is not the $0$ class; for instance, this is the case if $k=\FF_p$, $M$ is of rank 2 over $\ZZ_p$, $G=GL(M)$ and we replace $(M,F^1,\vph)$ by $(M^a,F^1_d,\be_M\mu_y({1\over p}))$, where $\be_M$ is the scalar automorphism of $M$ defined by an invertible element of $\ZZ_p$ which is not a $p-1$-th power. 
\smallskip
It can be checked using quite involved matrix computations that the positive and the negative $p$-divisible group of $({\got g},F^0({\got g}),F^1({\got g}),\vph)$ are dual to each other. However, here we take a faster approach, though slightly restricted. We have the following obvious result:
\medskip
{\bf Fact.} {\it We assume $({\got g},F^0({\got g}),F^1({\got g}),\vph)$ or its adjoint is a symmetric $p$-divisible object of $\Mm\Mf_{[-1,1]}(W(k))$. Then the positive and the negative $p$-divisible group of $({\got g},\vph,F^0({\got g}),F^1({\got g}))$ are dual to each other.} 
\medskip
See 2.2.23 and the beginning paragraph of 3.4 for how practical this Fact is. For instance, it applies if the Killing form on ${\rm Lie}(G^{\rm ad})$ is perfect.
\medskip
{\bf 3.11.8.1. Direct sum decompositions.} We have slope sign (resp. slope type) direct sum decompositions similar to (SSDSD) (resp. of (STDSD)) of 3.11.2 C, for any Shimura filtered Lie $\sg$-crystal $({\got g},\vph,F^0({\got g}),F^1({\got g}))$ of parabolic type over $k$. This is a consequence of 3.11.6 B, by natural passage from $k$ to $\bar k$ (cf. also Lemma of 2.2.3 3)). Such (SSDSD)'s represent an improvement to 2.2.12.1 1). Based on (STDSD), we can uniquely associate a Galois representation to $({\got g},F^0({\got g}),F^1({\got g}),\vph)$: it is a direct sum of Galois representations indexed by the slopes of $({\got g},\vph)$. For instance, referring to (STDSD) of 3.11.2 C, it is the direct sum
$$\oplus_{\al\in S} \rho_{\al},$$ where $\rho_{\al}$ is the Galois representation associated (via Fontaine's comparison theory) to ${\got C}_{\al}$. 
\smallskip
We assume now that $({\got g},\vph,F^0({\got g}),F^1({\got g}))$ is the Shimura filtered Lie $\sg$-crystal attached to a $G$-canonical lift $(M,F^1,\vph,G,(t_{\al})_{\al\in\Mj})$. Let $\rho:\Gamma_k\to GL(N)(\ZZ_p)$ be the Galois representation associated to $(M,F^1,\vph)$ via Fontaine's comparison theory; this is the sixth place where we need $p\ge 3$. Let $v_{\al}\in\Mt(N[{1\over p}])$ correspond to $t_{\al}$ via Fontaine's comparison theory and let $G_N$ be the Zariski closure in $GL(N)$ of the subgroup of $GL(N[{1\over p}])$ fixing $v_{\al}$, $\forall\al\in\Mj$. We denote by $\rho_G$ the factorization of $\rho$ through $G_N(\ZZ_p)$. Corresponding to (STDSD) of $({\got g},F^0({\got g}),F^1({\got g}),\vph)$, we get via Fontaine's comparison theory a direct sum decomposition ${\rm Lie}(G_N)[{1\over p}]=\oplus_{\al\in S} \tilde L_{\al}$ in $\QQ_p$-vector spaces left invariant (under ${\rm Ad}(\ZZ_p)\circ\rho_G$) by $\Gamma_k$. 
\smallskip
To show that it extends to an \'etale slope type direct sum decomposition
$${\rm Lie}(G_N)=\oplus_{\al\in S} L_{\al}\leqno (ESTDSD)$$
of $\Gamma_k$-submodules, with $L_{\al}:=\tilde L_{\al}\cap {\rm Lie}(G_N)$, we can assume $k=\bar k$. Let $P_{=0W(k)}$ be the subgroup of $G$ defined as the extension to $W(k)$ of a reductive group $P_{=0}$ over $\ZZ_p$ constructed as in 3.11.2 C (as a subgroup of a natural $\ZZ_p$-model of $G$), and let $T$ be a maximal torus of it whose Lie algebra is generated by elements fixed by $\vph$. As in the end of the proof of 3.11.1 we get that $(M,F^1,\vph,T)$ is a Shimura filtered $\sg$-crystal. Moreover, as in 2.2.16.5 (cf. also [Va2, 4.3.13]), based on [Va2, 4.3.9] and the functorial aspect of Fontaine's comparison theory, we get that via this theory, to $T$ it corresponds a maximal torus $T_0$ of $G_N$. As in the proof of 2.2.16.2 we have a natural identification
$${\rm Lie}(T)={\rm Lie}(T_{0W(k)}).\leqno (TOR)$$
Argument: both sides of (TOR) can be identified with the extension to $W(k)$ of the Lie algebra of endomorphisms of $(M,F^1,\vph,T)$.
So we can assume $\rho_G$ factors through the 
group of $\ZZ_p$-valued points of $T_0$. Moreover, $L_{\al}$ is normalized by $T_0$ (this can be checked by moving to the crystalline context; see (LBP) of 2.2.3 3)). As the set $S(\al)$ of characters of $T_{0W(\FF)}$ through which $T_{0W(\FF)}$ acts on $L_{\al}\otimes_{\ZZ_p} W(\FF)$ is uniquely determined by $\al$ and we have $S(\al)\cap S(\al_1)=\varnothing$, if $\al\neq\al_1$, we get that we have a direct sum decomposition ${\rm Lie}(G_N)\otimes_{\ZZ_p} W(\FF)=\oplus_{\al\in S} L_{\al}\otimes_{\ZZ_p} W(\FF)$; from this (ESTDSD) follows.
\medskip
Declaring $L_s=\{0\}$ if $s\in [-2,2]\setminus S$, we have a natural inclusion
$$[L_{\al},L_{\be}]\subset L_{\al+\be}$$
of $\Gamma_k$-modules, where $\al$ and $\be$ belong to $S$.
If $k=\bar k$, then $L_0$ is the maximal $\ZZ_p$-submodule of ${\rm Lie}(G_N)$ on which $\Gamma_k$ acts trivially. From the functorial aspect of Fontaine's comparison theory we get that we have a natural identification
$$L_0={\rm Lie}(P_{=0})\leqno (ZERO)$$ 
\indent
We assume now that $G_N$ is a reductive group (1.15.1 points out implicitly that this is always so). Using the same arguments, we get that we have an \'etale slope type direct sum decomposition of ${\rm Lie}(G_N^{\rm ad})$. Moreover, we have a logical variant of (TOR) and (ZERO) in the adjoint context. To argue this last thing, we consider the Zariski closure $P_{=0}^\prime$ in $G_N$ of the reductive subgroup of $G_{N\QQ_p}$ having $L_0[{1\over p}]$ as its Lie algebra. It is a reductive subgroup of $G_N$, as it contains the maximal torus $T_0$ of $G_N$ and as its generic fibre is. From Fontaine's comparison theory we get that we have a natural identification of the generic fibres of $P_{=0}^\prime$ and of $P_{=0}$; under this identification, $\ZZ_p$-valued points of $P_{=0}^\prime$ are mapped into $\ZZ_p$-valued points of $P_{=0}$. From [Ti2, 3.4.1] we get that in fact we have an identification of $P_{=0}^\prime$ with $P_{=0}$. This implies that we have a natural identification of the extension to $W(k)$ of the image $\tilde P_{=0}$ of $P_{=0}^\prime$ in $G^{\rm ad}_N$ and the image of $P_{0W(k)}$ in $G$. 
From this 
the adjoint variants of (TOR) and (ZERO) follow. 
\medskip
{\bf 3.11.8.1.1. Remark.} If $p\ge 5$, then 3.11.8.1 holds in the general context of Shimura-canonical lifts of Shimura-ordinary (adjoint) Lie $\sg$-crystals, cf. 2.2.16.5.
\medskip
{\bf 3.11.9. The cyclic diagonalizable context.} Most of the decompositions of 3.11.2 C and 3.11.8 make sense, under proper formulation, in a cyclic diagonalizable context. The same applies in connection to 3.11.2 B and to ${\rm Aut}_i$ of 3.11.7, with $i\in\{0,1\}$. To list the two differences showing up, let ${\got C}=(M,F^1,\vph,G,(t_{\al})_{\al\in\Mj})$ be a cyclic diagonalizable Shimura filtered $\sg$-crystal. We do have a slope sign direct sum decomposition for ${\got C}$ as in 3.11.2 C but its summands are $p$-divisible objects of $\Mm\Mf_{[-1,1]}(W(k))$. 3.11.8.1 (ZERO) holds iff the Shimura filtered Lie $\sg$-crystal attached to ${\got C}$ is of toric type.  
\medskip\smallskip
{\bf 3.12. Shimura $F$-crystals and generic points.}
Let $\bigl(M,F^1,\vph,G,{(t_\al)}_{\al\in\Mj}\bigr)$ be a Shimura filtered $\sg$-crystal over $k$.
Let ${\got C}_H:=(M,F^1,\vph,G,H,\tilde f)$ be a Shimura filtered $F$-crystal with $H$ a smooth, locally closed subscheme of $G$ through which the origin factors and such that the $W(k)$-morphism $H\to G/P$ is smooth in a open neighborhood of this factorization; here $P$ is the
parabolic subgroup of $G$ normalizing $F^1$. With the notations of 2.2.10, let $\eta:{\rm Spec}(k_\eta)\to{\rm Spec}(R)$ be a geometric point sitting over the generic point of the special fibre of ${\rm Spec}(R)$. Let 
$$
{\got C}_\eta:=\bigl(M\otimes_{W(k)} W(k_\eta),\vph_\eta,G_{W(k_\eta)}\bigr)
$$ 
be the Shimura $\sg_{k_\eta}$-crystal obtained from
$(M,F^1,\vph,G,H,\tilde f)$ by pull back through $\eta$ (cf. Fact 3 of 2.2.10).
\medskip
{\bf 3.12.1. Theorem.} {\it ${\got C}_\eta$ is a $G_{W(k_\eta)}$-ordinary $\sg_{k_\eta}$-crystal.}
\medskip
{\bf Proof:} There is a Shimura filtered $F$-crystal ${\got C}_G:=(M,F^1,\vph,G,\tilde f_1)$ with the property that Shimura $F$-crystals obtained from it by pull back through geometric points over the generic point of the special
fibre of the scheme over which ${\got C}_G$ is defined are $G$-ordinary (cf. 3.1.0 d)). 
\smallskip
Let $\ell:=dd((M,\vph,G))$ (see 2.2.22 4)).
We consider a Shimura filtered $F$-crystal ${\got C}_1$ over the special fibre of ${\rm Spec}(R_1)$, where $R_1:=W(k)[[x_1,...,x_{\ell}]]$, which is a universal deformation of $(M,F^1,\vph,G,(t_{\al})_{\al\in\Mj})$ (see 2.2.21 UP). ${\got C}_H$ and ${\got C}_G$ are both induced from ${\got C}_1$ through (uniquely determined; see 2.2.21 UP) $W(k)$-homomorphisms $R_1\to R$ and respectively $R_1\to R_G$ compatible with the $W(k)$-epimorphisms $R_1\twoheadrightarrow W(k)$, etc., defining the origins; as the natural $W(k)$-morphism $H\to G/P$ is smooth around the origin and as the natural $W(k)$-morphism $G\to G/P$ is smooth, we get (to be compared with [Va2, 5.4.4-8]) that these two $W(k)$-homomorphisms have left inverses and so are formally smooth. From this the Theorem follows.
\medskip\smallskip
{\bf 3.13. Deviations of Shimura (filtered) (Lie) $\sg$-crystals.} In this section we express our hopes for obtaining a reasonable theory of good ``$G$-pseudo-canonical lifts" of Shimura $\sg$-crystals $(M,\vph,G)$ which are not potentially cyclic diagonalizable. We mostly restrict to present some definitions, to restate some of the previous results and problems in terms of them, to formulate some expectations and to list some of the most interesting problems arising.
However, see section 3.13.7: it is dedicated to a study of truncations modulo positive integral powers of $p$ of objects of $\Mm\Mf(W(\bar k))$; though we mostly deal with (the classification of) truncations mod $p$ in a context related to Shimura (Lie) $\sg$-crystals, the general context of such truncations mod $p$ is as well mentioned (see 3.13.7.8) and moreover some paragraphs refer to things modulo positive integral powers of $p$.
\smallskip
Let ${\got C}_0:=(M,\vph,G)$ be a Shimura $\sg$-crystal over $k$. Let $\Mf$ be the set of direct summands $F^1$ of $M$ such that the quadruple $(M,F^1,\vph,G)$ is a Shimura filtered $\sg$-crystal. $\Mf$ is the $G(W(k))$-conjugacy class of an arbitrary $F^1\in\Mf$ (cf. Fact 2 of 2.2.9 3)). To fix the notations we choose (arbitrarily) an element $F^1\in\Mf$. We get a Shimura filtered $\sg$-crystal ${\got C}:=(M,F^1,\vph,G)$. As usual, let $F^0({\rm Lie}(G))$ be the parabolic Lie subalgebra of ${\rm Lie}(G)$ formed by elements taking $F^1$ into itself. 
\smallskip
We assume $\Mf$ has more than 1 element; so $G$ is not a torus (cf. 3.2.2). Let $\got p$ be the parabolic Lie subalgebra of ${\got s}{\got l}(M)$ corresponding to non-negative slopes of $({\got s}{\got l}(M),\vph)$. Let 
$$
{\got p}_G:={\rm Lie}(G^{\rm der})\cap\got p. 
$$
Let $\pi$ be a projector of ${\got g}{\got l}(M[{1\over p}])$ on ${\rm Lie}(G_{B(k)})$ fixed by $\vph$ and $G$. For instance we can take $\pi$ such that its kernel is the $B(k)$-vector subspace of ${\got g}{\got l}(M[{1\over p}])$ perpendicular on ${\rm Lie}(G_{B(k)})$ w.r.t. the trace form on ${\got g}{\got l}(M[{1\over p}])$ (see [Va2, 4.2] and the reference to it of AE.0). In what follows, except (perhaps) the undetailed second part of 3.13.3 3) below, it is irrelevant which such projector $\pi$ we choose. Let 
$${\got p}_G^0:={\got p}\cap {\rm ker}(\pi).$$ 
The $B(k)$-vector spaces ${\got p}_G^0[{1\over p}]$ and ${\got p}_G[{1\over p}]$ are normalized by $\vph$.
\medskip
{\bf 3.13.1. Definitions and notations. 1)} By the deviation of ${\got C}$ we mean the rational number $$d({\got C}):={{\dim\bigl({\got p}_G\cap F^0({\rm Lie}(G))\bigr)}\over {\dim_{W(k)}({\got p}_G)}}.$$
\smallskip
{\bf 2)} By the complementary deviation of ${\got C}$ we mean the rational number $$d^\perp({\got C)}:={{\dim\bigl({\got p}_G^0\cap F^0({\got s}{\got l}(M))\bigr)}\over {\dim_{W(k)}({\got p}_G^0)}}.$$
\smallskip
{\bf 3)} By the global deviation of ${\got C}$ we mean the rational number $$gd({\got C}):={{\dim\bigl({\got p}\cap F^0({\got s}{\got l}(M))\bigr)}\over {\dim_{W(k)}({\got p})}}.$$
\smallskip
{\bf 4)} By the deviation of ${\got C}_0$ we mean the rational number $d({\got C}_0)$ defined as the maximum of all numbers $d({\got C})$, with the $F^1$-filtration $F^1$ of $M$ running through the set $\Mf$.  
\smallskip
{\bf 5)} By the complementary (resp. global) deviation of ${\got C}_0$ we mean the rational number $d^\perp({\got C}_0)$ (resp. $gd({\got C}_0)$) defined as the maximum of all rational numbers $d^\perp({\got C})$ (resp. $gd({\got C})$), with the $F^1$-filtration $F^1$ of $M$ running as in 4).
\smallskip
{\bf 6)} By a $G$-pseudo-canonical lift of ${\got C}_0$ we mean a Shimura filtered $\sg$-crystal $\tilde {\got C}=(M,F^1_0,\vph,G)$ such that $d({\got C}_0)=d(\tilde {\got C})$.
\smallskip
{\bf 7)} By the $CM$-deviation of ${\got C}_0$, we mean the number $CM-d({\got C}_0)\in\NN\cup\infty$, defined as follows: if there is a finite, flat DVR extension $V$ of $W(\bar k)$  such that there is a maximal $\QQ_p$-torus of $G_{B(\bar k)}$ whose Lie algebra is formed by elements of ${\rm Lie}(G_{B(\bar k)})$ fixed by $\vph\otimes 1$ and whose extension to $V[{1\over p}]$ is a torus of the generic fibre of a parabolic subgroup of $G_V$ lifting the parabolic subgroup of $G_{\bar k}$ normalizing $F^1\otimes_{W(k)} \bar k$, then $CM-d({\got C}_0)$ is the smallest possible value of $[V:W(\bar k)]$ for such a $V$; if such a $V$ does not exist, then $CM-d({\got C}_0)=\infty$.  
\medskip
{\bf 3.13.2. Remarks. 1)} ${\got C}_0$ is a $G$-ordinary $\sg$-crystal iff $d({\got C}_0)=1$, cf. 3.1.0 a). ${\got C}$ is a $G$-canonical lift iff $d({\got C})=1$, cf. 3.1.0 a) and b). This motivates our terminology: deviations of Shimura (filtered) $\sg$-crystals. So we view the numbers introduced in 3.13.1 1) to 5) as a measure of how far or close  a Shimura $\sg$-crystal $(M,\vph,G)$ (resp. a Shimura filtered $\sg$-crystal $(M,F^1,\vph,G)$) is from being a $G$-ordinary $\sg$-crystal (resp. a $G$-canonical lift).
\smallskip  
{\bf 2)} A $G$-pseudo-canonical lift is not necessarily unique. This can be seen looking at $\sg$-crystals associated to supersingular abelian varieties over $k$.
\smallskip
{\bf 3)} The $CM$-deviations are measuring how far a Shimura $\sg$-crystal is from being potentially cyclic diagonalizable: the $CM$-deviation of a Shimura $\sg_{\bar k}$-crystal is 1 iff it has a lift of quasi CM type (cf. 2.2.18).
\medskip
{\bf 3.13.3. Variants. 1)} 3.13.1 has been stated using non-negative slopes. The same thing can be stated using positive, negative or non-positive slopes. In such a case, as well as in others (to be defined below), we  use a left upper index of the form $+$, $-$, or $0-$ for the deviations introduced, to emphasize that they have been defined using respectively positive, negative or non-positive slopes. For instance we get rational numbers ${}^+d({\got C})$, $p^n-{}^{0-}d({\got C})$, etc. 
\smallskip
{\bf 2)} In some cases it is useful to define deviations of higher order. By this we mean that in 3.13.1 1) to 5) we replace $F^0({\rm Lie}(G^{\rm der}))$ or $F^0({\got s}{\got l}(M))$ with their iterates under the action of $\vph$; for instance, we replace $F^0({\got s}{\got l}(M))$ by 
$${\got s}{\got l}(M)\cap\vph^i(F^0({\got s}{\got l}(M[{1\over p}]))),$$ 
$i\in\NN\cup\{0\}$. So we obtain sequences of deviation numbers $d_i({\got C})$, $d_i^\perp({\got C})$, $gd_i({\got C})$, $d_i({\got C}_0)$, $d_i^\perp({\got C}_0)$, $gd_i({\got C}_0)$, $i\in\NN\cup\{0\}$, which for $i=0$ become the numbers introduced in 3.13.1.
\smallskip
{\bf 2')} Replacing in 2) $\vph^i(F^0({\got s}{\got l}(M[{1\over p}])))$ by 
$$\cap_{n=0}^i \vph^n(F^0({\got s}{\got l}(M[{1\over p}]))),$$ 
we obtain $\cap$ (intersection) type deviations of higher order. The notations to be used: $d_{\cap i}({\got C})$, $d_{\cap i}^\perp({\got C}_0)$, etc.  
\smallskip
In connection to 2) and 2'), we always drop the lower right index $i=0$.
\smallskip
{\bf 3)} Definitions 3.13.1 1) to 6) make sense for Shimura (filtered) isocrystals over $k$. Of course, for a Shimura isocrystal $(M^1,\vph^1,G^1)$ we need to choose a set $\Mf^1$ of ``admissible" filtrations of $M^1$; often $\Mf^1$ is a $G^1(B(k))$-conjugacy class of a fixed filtration. 
\smallskip
Sometimes, for the sake of refinement, we also work with $p^n$-deviations of Shimura (filtered) $\sg$-crystals, $n\in\NN$. Not to make the story too long, here we just define the $p^n$-deviation of $\got C$, $n\in\NN$. It is the rational number defined by:
$$p^n-d({\got C}):={{{\rm length}\bigl({\got p}_G\otimes_{W(k)} W_n(k)\cap F^0({\rm Lie}(G)\otimes_{W(k)} W_n(k))\bigr)}\over {{\rm length}\bigl({\got p}_G\otimes_{W(k)} W_n(k)\bigr)}}.$$
Here ${\rm length}(*)$ is the usual length of a finitely generated module $*$ over the artinian local ring $W_n(k)$. 
\smallskip
{\bf 4)} All deviations introduced above which are not of global or complementary type, make sense for Shimura (filtered) Lie $\sg$-crystals.
\smallskip
{\bf 5)} We can define isogeny deviations. They are denoted by putting $isog-$ in front of the notation of the involved deviation. For instance, by the isogeny deviation of ${\got C}_0$, denoted $isog-d({\got C}_0)$, we mean the greatest number $d({\got C}_0^i)$, with ${\got C}_0^i$ running through the representatives of $Isog({\got C}_0)$. For $isog-CM-({\got C}_0)$, the word ``greatest" has to be replaced by the word ``smallest".
\smallskip
{\bf 6)} We can define  strong $CM$-deviations, denoted by replacing $CM$ by $SCM$, by replacing (in 3.13.1 7)) the expression ``there is a maximal" by: ``for any maximal".
\medskip
{\bf 3.13.3.1. Remarks 1).}
${\got C}_0$ is a $G$-ordinary $\sg$-crystal iff $p-d({\got C}_0)=1$. ${\got C}$ is a $G$-canonical lift iff  $p^n-d({\got C})=1$, $\forall n\in\NN$. These two facts are a consequence of the proof of 3.4.8 (cf. also 3.1.0 a)).
\smallskip
{\bf 2)} We can introduce also the numbers 
$$p^n-d_f({\got C}):={{{\rm length}\bigl({\got p}_G\otimes_{W(k)} W_n(k)\cap F^0({\rm Lie}(G^{\rm der})\otimes_{W(k)} W_n(k))\bigr)}\over {{\rm length}\bigl(F^0({\rm Lie}(G^{\rm der}))\otimes_{W(k)} W_n(k)\bigr)}}.$$
But these numbers can be greater or smaller than the ones attached to a $G$-canonical lift. Using this, we can easily construct examples for which these numbers for $G$-canonical lifts do not characterize them.  
\smallskip
{\bf 3)} In 1995 we conjectured: 
\medskip
{\bf Conjecture.} {\it All $isog-SCM$-deviations of Shimura (Lie) $\sg$-crystals are finite numbers.}
\medskip
 This conjecture is very much related to Langlands--Rapoport's conjecture referred to in 1.15.9 and so it will play a central role in \S 14.
\smallskip
{\bf 4)} We do believe that the above deviations allow a relatively uniform treatment of all Shimura (adjoint) (filtered) (Lie) $\sg$-crystals, regardless of the fact they are or are not potentially cyclic diagonalizable.    
\medskip
{\bf 3.13.4. Problems.} Many problems (questions) are arising. We  mention just few of them.
\smallskip
{\bf 1)} When $G$-pseudo-canonical lifts can be defined as well in terms of $d^\perp(*)$ or of $gd(*)$ numbers?
\smallskip
{\bf 2)} Fixing ${\got C}_0$, can we determine all possible values of the numbers $d({\got C})$, $d^\perp({\got C})$, $gd({\got C})$, where the filtration $F^1$ of $M$ runs as in 3.13.1 4)?
\smallskip
{\bf 3)} Determine all Shimura $\sg$-crystals ${\got C}_0$ for which we can define uniquely up to isomorphisms a $G$-pseudo-canonical lift of them, by using the numbers introduced in 3.13.3 2) and 3).
\smallskip
{\bf 4)} How are  $p^n$-deviations ($n\in\NN$) related to deviations? 
\smallskip
{\bf 5)} How deviations ($p^n$-deviations, $CM$-deviations, etc.) vary when we replace $\vph$ by $g\vph$ with $g\in G(W(k))$? What are their possible values? 
\smallskip
{\bf 6)} What is the connection between $CM$-deviations and deviations?
\smallskip
{\bf 7)} Let $A$ be an abelian variety over a number field $E$. For any unramified prime $v$ of $E$ of good reduction for $A$, we can define a $CM$-deviation $CMD_A(v)$ and a $SCM$-deviation $SCMD_A(v)$ (of the $\sg_{k(v)}$-crystal of the reduction of $A$ w.r.t. $v$). How are they varying in terms of $v$? Is it true that all these $CM$-deviations are 1, except for a finite number? See 4.6.2.3 2) for a similar question in the relative context. 
\smallskip
{\bf 8)} How are the $CM$-deviations and the $SCM$-deviations related?
\medskip
{\bf 3.13.5. Wintenberger groups and deviations.} We consider the canonical split cocharacter $\mu(F^1):\GG_m\to GL(M)$ of the filtered $\sg$-crystal $(M,F^1,\vph)$. Due to its functoriality, it factors through $G$. The Zariski closure in $G$ of the smallest subgroup of $G_{B(k)}$ whose Lie algebra is normalized by $\vph$ and contains the Lie algebra of $\mu(F^1)(\GG_m)$, is called (in honor of [Wi, 4.2.3]) the Wintenberger group of $(M,F^1,\vph)$; it does not depend on $G$ and so is denoted by $W(M,F^1,\vph)$. From the proof of 2.2.18 we get: 
\medskip
{\bf Fact.} {\it $W(M,F^1,\vph)$ is a torus of $G$ iff $(M,F^1,\vph,G)$ is potentially cyclic diagonalizable.} 
\medskip
Let $\Mw_1$ be the subgroup of $G$ generated by $\mu(F^1)$, with $F^1$ running as in 3.13.1 4). Let $\Mw_2$ be the Zariski closure in $G$ of the  smallest (so connected) subgroup of $G_{B(k)}$ with the property that its Lie algebra is stable under $\vph$ and contains ${\rm Lie}(\Mw_1)$. Let $\Mw_{i+2}$ be the Zariski closure in $GL(M)$ of the centralizer of the generic fibre of $\Mw_i$ in $GL(M[{1\over p}])$, $i=\overline{1,2}$. We call $\Mw_i$ as the $i$-th Wintenberger group of ${\got C}_0$, $i=\overline{1,4}$. 
\medskip
{\bf Exercise.} If $(M,\vph)$ is an ordinary $\sg$-crystal, show that $\Mw_1$ is the extension of $\GG_m$ by a smooth, commutative, unipotent group of relative dimension equal to $dd({\got C}_0)$. Hint: use 3.6.18.2.
\medskip
{\bf 3.13.5.1. Versal endomorphisms of ${\got C}_0$.} Let  $(t_{\al})_{\al\in\Mj}$ be a family of tensors of $\Mt(M[{1\over p}])$ such that $({\got C},(t_{\al})_{\al\in\Mj})$ is a Shimura filtered $\sg$-crystal with an emphasized family of tensors. Let $\Md$ be a Shimura $p$-divisible group over $W(k)$ whose associated $p$-divisible object with tensors of $\Mm\Mf_{[0,1]}(W(k))$ is $({\got C},(t_{\al})_{\al\in\Mj})$. We consider an endomorphism $h:M\to M$ fixed by $\vph$. We say $h$ is a versal endomorphism of ${\got C}_0$ if any one of the following equivalent conditions is satisfied:
\medskip
{\bf i)} $h$ normalizes any element of $\Mf$;
\smallskip
{\bf ii)} $h$ is fixed by $\Mw_1$;
\smallskip 
{\bf iii)} $h$ is fixed by $\Mw_2$;
\smallskip
{\bf iv)} $h$ defines naturally an endomorphism of the $p$-divisible group $DG$ underlying the universal Shimura $p$-divisible defined (see 2.2.21) by $\Md$ (over the formal power series ring in $dd({\got C}_0)$ variables with coefficients in $W(k)$).
\medskip
The use of iv) is part of the third place where we need $p\ge 3$. Similarly we speak about versal automorphisms $h:(M,\vph)\tilde\to (M,\vph)$ of ${\got C}_0$; in this case the above equivalent conditions are also equivalent to:
\medskip
{\bf v)} $h\in\Mw_3(W(k))$;
\smallskip
{\bf vi)} $h\in\Mw_4(W(k))$.
\medskip
Warning: a versal automorphism is not necessarily an automorphism in the sense of 2.2.9 6).
\medskip
{\bf 3.13.5.2. Comments.} We do not know when $\Mw_i$ is a smooth group over $W(k)$, or when the special fibre of $\Mw_i$ is a connected group, or when the generic fibre of $\Mw_3$ or of $\Mw_4$ is a connected group. We do not know how to compute the relative dimension $w_i$ of $\Mw_i$ or how to determine the structure of the generic fibre of $\Mw_i$. We do not know when $w_i$ is an isogeny invariant.
\medskip
{\bf 3.13.5.3. $\Mw$-deviations.} 
Let $g\in G(W(k))$ be such that $(M,g\vph,G)$ is a $G$-ordinary $\sg$-crystal. We define $w_i(G-ord)\in\NN$ as the relative dimension of the $i$-th Wintenberger group of the extension of $(M,g\vph,G)$ to $\bar k$. The number 
$$wd_i({\got C}_0):={w_i\over {w_i(G-ord)}}\in (0,\infty)\cap\QQ$$
is called the $i$-th $\Mw$-deviation or the $i$-th Wintenberger deviation of ${\got C}_0$; it is well defined, cf. 3.11.1 c). Warning: $wd_i({\got C}_0)$ is not always $1$, as it can be checked using supersingular $\sg$-crystals. Some of the questions and problems of 3.13.2-4 can be adapted to the context of $\Mw$-deviations. 
\medskip
{\bf 3.13.5.4. Relative rigidity.} Let $G_1$ be a reductive subgroup of $GL(M)$ containing $G$ and such that the triple $(M,\vph,G_1)$ is a Shimura $\sg$-crystal. We say ${\got C}_0$ is rigid w.r.t. $G_1$, if for any versal automorphism of ${\got C}_0$ defined by an element $h\in G_1(W(k))$ we have: $h$ centralizes ${\rm Lie}(G^{\rm der})$. When $G_1=GL(M)$ we drop the words ``w.r.t. $G_1$". 
\medskip
{\bf Exercise.} We assume that $\dim_{W(k)}(M)=8$, that $G$ contains the group of scalar automorphisms of $M$ and that $G$ is a split $GL_2$-group whose representation on $M$ is a direct sum of 4 copies of the standard $2$ dimensional irreducible representation of $G$. If $dd({\got C}_0)=1$ (i.e. if it is non-zero) and if all slopes of $(M,\vph)$ are $0$ or $1$, then ${\got C}_0$ is rigid.
Hint: use the fact that any endomorphism of $M$ fixed by a Borel subgroup of $G$ is fixed by $G$ itself.
\medskip
{\bf 3.13.6. Three other types of deviations.} By the isomorphism deviation $isom-d({\got C}_0)$ of ${\got C}_0$ we mean the isomorphism deviation of its extension to $\bar k$ (cf. 3.6.16 2)). By the Newton polygon deviation (resp. Newton polygon Lie deviation) of ${\got C}_0$, we mean the smallest number $NP-d({\got C}_0)\in\NN$ (resp. $NPL-d({\got C}_0)\in\NN$) such that $(M,g\vph)$ and $(M,\vph)$ (resp. such that $({\rm Lie}(G),g\vph)$ and $({\rm Lie}(G),\vph)$) have the same Newton polygon, provided $g\in G(W(k))$ is congruent to the identity mod $p^{NP-d({\got C}_0)}$ (resp. mod $p^{NPL-d({\got C}_0)}$). We have $isom-d({\got C}_0)\ge NP-d({\got C}_0)$ and $isom-d({\got C}_0)\ge NPL-d({\got C}_0)$.
\smallskip
If $NP-d({\got C}_0)=1$ (resp. if $NPL-d({\got C}_0)=1$) we say ${\got C}_0$ is an $NP$-canonical (resp. $NP$-Lie-canonical) Shimura $\sg$-crystal. If $isom-d({\got C}_0)=1$, then we say ${\got C}_0$ is an $isom$-canonical Shimura $\sg$-crystal. Similarly we define $isom$-deviations of Shimura adjoint Lie $\sg$-crystals and $isom$-canonical Shimura adjoint Lie $\sg$-crystals.
\medskip
{\bf Examples.} Any Shimura-ordinary $\sg$-crystal is $NP$-canonical and $NP$-Lie-canonical, cf. 3.9.3 (and 3.1.0 a) and c)); even more: it is also an $isom$-canonical Shimura $\sg$-crystal, cf. 3.11.1 c). If $\dim_{W(k)}(M)=2$, then ${\got C}_0$ is automatically $NP$-canonical, $NP$-Lie-canonical and $isom$-canonical Shimura $\sg$-crystal.
\medskip
{\bf Problem.} Study how the isomorphism and Newton polygon (Lie) deviations are related to the other deviations introduced in 3.13.1-5. 
\medskip
{\bf 3.13.6.1. NP variation functions.} We assume $n_1:=NP-d({\got C}_0)>1$ (resp. $n_2:=NPL-d({\got C}_0)>1$). For $n\in S(1,n_1-1)$ (resp. $n\in S(1,n_2-1)$), we denote by $f_1(n)$ (resp. $f_2(n)$) the number of (distinct) Newton polygons of the Shimura (resp. Shimura Lie) $\sg$-crystals of the form $(M,g\vph)$ (resp. of the form $({\rm Lie}(G),g\vph)$), with $g\in G(W(k))$ congruent to the identity mod $p^n$. We refer to $f_1$ (resp. $f_2$) as the first (resp. the second) NP variation function attached to (or of) ${\got C}_0$.
\medskip
{\bf 3.13.7. CM levels.} Let $n\in\NN$. Till end of 3.13.7-9 we assume $k=\bar k$. We now assume ${\got C}_0$ is not cyclic diagonalizable. We say ${\got C}_0$ is of CM level at least $n$ if there is a lift of it $F^1\in\Mf$ such that $\got C$ is cyclic diagonalizable of level $n$ (in the sense of 2.2.22 1)).  
 We denote by 
$$n({\got C}_0)\in\NN\cup\{0,\infty\}$$ 
the smallest supremum of the set of those $n$ such that ${\got C}_0$ is of CM level at least $n$. We refer to $n({\got C}_0)$ as the CM level of ${\got C}_0$. End of 2.2.22 1) points out that using emphasized families of tensors we get the same CM levels.
\smallskip
If $k\neq\bar k$ and ${\got C}_0$ is not potentially cyclic diagonalizable, we define the CM level $n({\got C}_0)$ of ${\got C}_0$ to be the CM level of its extension to $\bar k$. 
\medskip
{\bf Comments.} {\bf 1)} It seems to us possible to formulate extremely challenging variants of 3.13.4 7), involving CM levels instead of $CM$-deviations. 
\smallskip
{\bf 2)} We assume that $k=\bar k$ and that ${\got C}_0$ is not cyclic diagonalizable. As the Weyl group of $G$ is finite, we deduce the existence of a finite number of isomorphism classes of cyclic diagonalizable Shimura filtered $\sg$-crystals of the form $(M,F^1,g\vph,T)$, with $g\in G(W(k))$ and with $T$ a maximal torus of $G$. So, from this and from 3.6.16 2) we get: 
\medskip
{\bf Corollary.} {\it There is $M\in\NN\cup\{0\}$, such that $n((M,g\vph,G))\le M$, for any $g\in G(W(k))$ with the property that $(M,g\vph,G)$ is not cyclic diagonalizable.}
\smallskip
The smallest such number $M$ is called the maximal CM level of $Cl({\got C}_0)$.
\smallskip
{\bf 3)} From the Corollary of 2.2.22 1) and from 2.2.13.3 we get that we can speak as well about CM levels of Shimura adjoint Lie $\sg$-crystals. 
\medskip
{\bf 3.13.7.1. Expectation (the CM level one property).} {\it Always $n({\got C}_0)\in\NN$.}
\medskip
{\bf An approach towards its proof.} One needs to show $n({\got C}_0)\neq 0$. Let $T$, $B$ and $\vph_1$ have the same significance as in 3.2.3. So we can write $\vph=g\vph_1$, with $g\in G(W(k))$ (with the notations of 3.2.3, $g:=(g_1g_0)^{-1}$). Let 
$$\bar\psi_1:{\rm Lie}(G^{\rm ad}_k)\to {\rm Lie}(G^{\rm ad}_k)$$ 
be the Faltings--Shimura--Hasse--Witt adjoint map of $(M,\vph_1,G)$. 
Let $P^+$ be the parabolic subgroup of $G$ having $F^0({\got g})$ as its Lie algebra and let $P^-$ be its opposite w.r.t. $T$. Let $N^+$ (resp. $N^-$) be the unipotent radical of $P^+$ (resp. of $P^-$). $P^{0}:=P^+\cap P^-$ is a Levi subgroup of both $P^+$ and $P^-$.
\smallskip
We consider the $\ZZ_p$-structure $G_{\ZZ_p}$ of $G$ defined as in 3.11.2 B for $(M,F^1,\vph_1,G)$. So for a smooth subscheme $H$ of $G$ (resp. of $G_k$), $\sg(H)$ is the smooth subscheme of $G$ (resp. of $G_k$) obtained from $H$ by applying $\sg$ to elements of $H(W(k))$ (resp. of $H(k)$). We have $\sg(T)=T$. Similarly we use $\sg$ in the adjoint context. As in 3.5.2, elements of $G(W(k))$ (resp. of $G(k)$) normalizing $T$ (resp. $T_k$) are called Weyl elements. $B$ is invariant under $\sg$, cf. 3.11.1 D.  
\smallskip
In connection to the Expectation, only the value $\bar g\in G(k)$ of $g$ mod $p$ is important: the reduction mod $p$ of $(M,F^1,g\vph_1,G)$ depends only on $\bar g$ and on ``the choice" of $F^1$ mod $p^2$. 
Changing $\vph:M\to M$ by another $\sg$-linear endomorphism $h{\vph}h^{-1}:M\to M$ isomorphic to it, with $h\in G(W(k))$ which mod $p$ normalizes $F^1/pF^1$, $g$ gets replaced by $hg\vph_1h^{-1}\vph_1^{-1}$. As we are interested only in $\bar g$, we can assume 
$$h=h_1h_3h_2,$$ 
with $h_1\in N^+(W(k))$, with $h_2\in P^{0}(W(k))$, and with $h_3\in N^-(W(k))$ congruent to the identity mod $p$. As $\vph_1=\sg\mu({1\over p})$, with $\mu:\GG_m\to G$ as the canonical split of $(M,\vph_1,G)$ (cf. 3.11.2 B), $\bar g$ gets replaced by the reduction mod $p$ of $hg\vph_1h^{-1}\vph_1^{-1}$ which is:
$$m_1m_2\bar g\sg(m_2^{-1})\sg(m_3^{-1});\leqno (*)$$
here $m_1\in N^+(k)$ and $m_2\in P^0(k)$ are the reduction mod $p$ of $h_1$ and respectively of $h_2$, while $m_3\in N^-(k)$ is uniquely determined by $h_3$; when $h_3$ vary, $m_3$ varies through all $k$-valued points of $N^-$. 
\smallskip
We consider the $k$-variety
$$H:=N^+_k\times_k P^0_k\times_k N^-_k,$$
and endow it with a new group structure by the following rule: if $(m_1,m_2,m_3)$, $(r_1,r_2,r_3)\in H(k)$ then
$$(m_1,m_2,m_3)(r_1,r_2,r_3):=(m_1m_2r_1m_2^{-1},m_2r_2,m_3m_2r_3m_2^{-1})\in H(k).$$
This new group structure is different from the product one; however, the group structures induced on $N^+_k\times_k N^-_k$ and on $P^0_k$ (embedded naturally in $H$) are the same as the initial ones. Moreover, we have a short exact sequence
$$0\to N^+_k\times_k N^-_k\hookrightarrow H\twoheadrightarrow P^0_k\to 0.\leqno (SES)$$ 
\indent
(*) gives us an action 
$$\TT:H\times G_k\to G_k$$
of $H$ on $G_k$ (viewed just as a variety) by the rule 
$$\TT_{(m_1,m_2,m_3)}(\bar g):=m_1m_2\bar g\sg(m_2^{-1})\sg(m_3^{-1}).$$ 
\indent
Let $P^{00}$ (resp. $P^{+0}$) be the image of $P^0$ (resp. of $P^+$) in $G^{\rm ad}$. Let 
$$\Mx:=N^+_k\setminus G_k^{\rm ad}/\sg({N^-_k}).$$
We consider the action
$$\TT^0:{P^{00}_k}\times\Mx\to\Mx$$
 of ${P^{00}_k}$ on $\Mx$ defined on points as follows: if $x\in\Mx$ is defined by $\bar g\in G^{\rm ad}(k)$, then for $h\in {P^{00}_k}(k)$, $\TT^0_h(x)$ is defined by the element $h\bar g\sg(h^{-1})$. We have another action $\TT^1$ of ${P^{00}_k}$ on $\Mx$: $\TT^1_h(x)$ is defined by the element $h\bar g$. We have:
\medskip
{\bf 3.13.7.1.1. Corollary.} {\it a) The orbits of $\TT$ are in one-to-one correspondence to inner $1_{\Mj}$-isomorphism classes of reduction mod $p$ of Shimura $\sg$-crystals of the form $(M,g\vph,G,(t_{\al})_{\al\in\Mj})$, with $g\in G(W(k))$.
\smallskip
b) The orbits of $\TT^0$ are in one-to-one correspondence to inner isomorphism classes of Faltings--Shimura--Hasse--Witt adjoint maps of the form $h\bar\psi_1$, with $h\in G_k^{\rm ad}(k)$ (here ``inner" refers to the fact that the isomorphisms are as well defined by elements of $G_k^{\rm ad}(k)$).}
\medskip 
{\bf Proof:} a) is a consequence of 2.2.14.0. For b), as the homomorphism $G(k)\to G^{\rm ad}(k)$ is surjective and due to (*), we just need to remark that:
\medskip
i) for $m_3\in G(k)$ we have $m_3\bar\psi_1=\bar\psi_1$ iff $m_3\in\sg(N^-_k)(k)$: one implication is trivial, as $\sg(N^-)$ is a commutative, unipotent group, while  the other one is the content of a) of Step 2 (it can be read out at any time) of 3.13.7.3 below;
\smallskip
ii) $P^{+0}_k$ is the subgroup of $G^{\rm ad}_k$ normalizing the kernel of $h\bar\psi_1$ (as any parabolic subgroup of $G_k$ is its own normalizer); so, as $g\bar\psi_1$ and $h\bar\psi_1$ have the same kernel ${\rm Lie}(P^{+0}_k)$, any inner isomorphism between them is defined by an element of $P^{+0}(k)$ (i.e. of $P^+(k)$);
\smallskip
iii) for any $\bar m\in P^{+0}(k)$ and every $h\in G^{\rm ad}(k)$, $\bar mh\bar\psi_1\bar m^{-1}$ is the Faltings--Shimura--Hasse--Witt adjoint map of the Shimura $\sg$-crystal $(M,\tilde m\tilde h\vph_1\tilde m^{-1},G)$, where $\tilde m$ (resp. $\tilde h$) is an element of $G(W(k))$ whose reduction mod $p$ lifts $\bar m$ (resp. $h$). 
 \medskip
{\bf 3.13.7.1.1.1. Exercise.} Let $m\in\NN$. Show that the truncations mod $p^m$ of two Shimura $\sg$-crystals with emphasized families of tensors $(N,g_N^i\vph_N,G_N,(t_{\al})_{\al\in\Mj_N})$ ($i=\overline{1,2}$) over a perfect field of characteristic $p\ge 3$ are inner isomorphic iff their Fontaine truncations mod $p^m$ are so. Hint: one implication (the if part) is obvious (cf. Fact 1 of 2.2.1.0); for the second one use induction on $m\in\NN$ (the above proof handles the case $m=1$ while the inductive step relies on 2.2.1.0 (VPHIONE) and on 2.2.14.0). 
\medskip
{\bf 3.13.7.1.2. Corollary (the mod $p$ property).} {\it $(M,g\vph_1,G,(t_{\al})_{\al\in\Mj})$ is $G$-ordinary iff its truncation mod $p$ is inner $1_{\Mj}$-isomorphic to the truncation mod $p$ of $(M,\vph_1,G,(t_{\al})_{\al\in\Mj})$.}
\medskip
{\bf Proof:} One implication follows from 3.11.1 c). The second implication follows from 3.9.2 and 3.13.7.1.1 a) and b). This ends the proof.
\medskip
The Expectation 3.13.7.1 can be reformulated as the first part of: 
\medskip
{\bf 3.13.7.1.3. Expectation (Bruhat decomposition in the Shimura $\sg$-crystal context).} {\it For any $\bar g\in G(k)$, there is an element $\bar g_1$ in the orbit $o(\bar g)$ of $\bar g$ under $\TT$ such that $\bar g_1(T_k)\bar g_1^{-1}=T_k$. Moreover, the distinct orbits of $\TT$ are in one-to-one correspondence to the elements of the quotient set $W_{P^+}\setminus W_G$, where $W_G$ is the Weyl group of $G$ w.r.t. $T$, while $W_{P^+}$ is its subgroup leaving $P^+$ invariant.}
\medskip
We have the following adjoint variant of 3.13.7.1.3:
\medskip
{\bf 3.13.7.1.4. Expectation (Bruhat decomposition in the Shimura adjoint Lie $\sg$-crystal context).} {\it For any $\bar g\in \Mx(k)$, there is an element in the orbit $o_0(\bar g)$ of $\bar g$ under $\TT_0$ defined by an element $\bar g_1\in G^{\rm ad}(k)$ such that $\bar g_1(T_k)\bar g_1^{-1}=T_k$. Moreover, the distinct orbits of $\TT_0$ are in one-to-one correspondence to the elements of $W_{P^+}\setminus W_G$.}
\medskip
{\bf 3.13.7.2. The connection between the orbits of $\TT$ and $\TT^0$.} The orbits of $\TT$ are in one-to-one correspondence to the ones of $\TT^0$ iff $\forall\bar g\in G(k)$ and $\forall z\in Z(G_k)$, we have $z\bar g\in o(\bar g)$. This is automatically so if the homomorphism $h(Z(G_k)):Z(G_k)\to Z(G_k)$ that takes $z\in Z(G_k)(k)$ into $z\sg(z^{-1})$ is an epimorphism; for instance, this is so if $Z(G_k)$ is a torus. 3.13.7.1 (or 3.13.7.1.3 or 3.13.7.1.4) once checked, implies that this always holds, regardless of how $h(Z(G_k))$ is: this is a consequence of 3.11.1 c) applied to Shimura filtered $\sg$-crystals of the form $(M,F^1,g\vph,T)$, with $g\in G(W(k))$. More precisely, from 3.11.1 c) we get:
\medskip
{\bf Fact.} {\it The pull back under the natural epimorphism $G_k\twoheadrightarrow \Mx$ of the orbit of $\TT^0$ defined by a Weyl element $w\in G(k)$ is the  orbit of $\TT$ defined by $w$.} 
\medskip
{\bf 3.13.7.3. Groups of automorphisms of Faltings--Shimura--Hasse--Witt adjoint maps.} We use the notations of 3.13.7.1, though below we take a slightly different approach. What follows is a digression in three Steps on such groups of automorphisms.
\medskip
{\bf Step 1. A reduction to Weyl's elements.} Let $Q:=\sg(P^-_k)$. Let $w_0\in G(k)$ be a Weyl element such that $w_0Qw_0^{-1}\cap P_k^+$ contains $B_k$. Using Bruhat decomposition (see [Bo2, 14.12]), we get that $\bar gw_0^{-1}=\bar pw_1w_0\bar qw_0^{-1}$, with $\bar p\in P^+(k)$, with $w_1\in G(k)$ a Weyl element and with $\bar q\in Q(k)$. So $\bar g={\bar p}w{\bar q}$, with $w:=w_1w_0$ a Weyl element. Let $\tilde p\in P^+(W(k))$ lifting $\bar p$. Replacing $(M,g\vph_1,G)$ by the Shimura $\sg$-crystal $(M,\tilde p^{-1}g\vph_1\tilde p,G)$ isomorphic to it, we can assume that the Faltings--Shimura--Hasse--Witt adjoint map of $(M,\vph,G)$, has the same image as $w\bar\psi_1$. This represents a quicker (but less useful) proof of 3.5.3 (2).
\medskip
{\bf Step 2. A lemma.} The notations of the following general Lemma, are entirely independent of all previous notations.
\medskip
{\bf Lemma.} {\it Let $k_1$ be an arbitrary field. Let $H$ be
an absolutely simple, adjoint group over $k_1$. Let $P$ be a parabolic subgroup of $H$. We consider the subgroup $C$ (resp. $P_C$) of $H$ centralizing (resp. normalizing) the Lie algebra ${\got n}$ of the unipotent radical $N$ of $P$. We have:
\medskip
a) If the characteristic of $k_1$ is not $2$, then $P_C=P$. If moreover $N$ is a commutative group, then $C=N$; 
\smallskip
b) If the characteristic of $k_1$ is $2$, then ${P_C}_{\rm red}=P$. If moreover $N$ is a commutative group, then $C_{\rm red}=C\cap P=N$;
\smallskip
c) The kernel of the adjoint map ${\rm Ad}:H\to GL({\rm Lie}(H))$ is trivial (i.e. the center of ${\rm Lie}(H)$ is trivial).}
\medskip
{\bf Proof:} We can assume the characteristic of $k_1$ is positive; we denote it by $p$. We can assume $k_1=\overline{k_1}$. The part c) is a well known result; for instance, it is a consequence of [Va2, 3.1.2.1 c)], applied to a lift of $H$ to an adjoint group $\tilde H$ over $W(k)$ (the characteristic $0$ case of c) is trivial). The existence of such a lift is a consequence of the uniqueness theorem of [SGA3, Vol. III, p. 313-4] and of the existence --for instance, see [Hu1, 25.5] or [SGA3, Vol. III, Exp. XXV]-- of Chevalley group schemes over $\ZZ$. 
\smallskip
$P_C$ contains $P$ and so ${P_C}_{\rm red}$ is a parabolic subgroup of $H$. Let $T$ be a maximal torus of $P$. Let $N^-$ be the unipotent radical of the opposite of $P$ w.r.t. $T$. As $T$ normalizes ${\got n}$ we have: $T$ normalizes $P_C$ and so it normalizes ${\rm Lie}(P_C)$. We now assume $p\ge 3$. To show that $P=P_C$, we just need to show that ${\rm Lie}(P_C)={\rm Lie}(P)$. If this is not so, then we deduce the existence of a non-zero element $x\in {\rm Lie}(P_C)\cap {\rm Lie}(\GG_a(x))$, with $\GG_a(x)$ as a $\GG_a$ subgroup of $N^-$ normalized by $T$; so $[x,{\got n}]\subset {\got n}$. Applying 3.5.4 a) to $[x,y]$, with $y$ a non-zero element of the Lie algebra of the $\GG_a$ subgroup of $N$ which is opposite (w.r.t. $T$) to $\GG_a(x)$, we reach a contradiction. 
\smallskip
If $p=2$ and ${P_C}_{\rm red}$ is not $P$, then from the structure of parabolic subgroups of $H$ containing $P$ (see [Bo2, 14.17-18]), we deduce that there is a $PSL_2$ or an $SL_2$-subgroup of ${P_C}_{\rm red}$ normalizing the Lie algebra of the unipotent radical of a Borel subgroup of it. Contradiction (there are no non-trivial 1 dimensional representations of such a group). 
\smallskip
We assume now $N$ is a commutative group; so $C$ contains $N$. It is enough to show that the reductive group $P/N$ acts faithfully on ${\got n}$. As the characteristic $0$ case is trivial, using again [Va2, 3.1.2.1 c)] in the context of the natural (inner) representation of $\tilde P/\tilde N$ on ${\rm Lie}(\tilde N)$, with $\tilde P$ as a lift of $P$ to a parabolic subgroup of $\tilde H$ and with $\tilde N$ as its unipotent radical, the Lemma follows.
\medskip
{\bf Exercise.} Using c) and standard noetherian arguments show that at the end of 2.2.13 the assumption $R=W(k)$ is not needed. 
\medskip
{\bf Step 3. Some computations.} We have:
\medskip
{\bf Fact.} {\it {\bf a)} $\forall h\in G^{\rm ad}(k)$, $hw\bar\psi_1$ has the same image as $w\bar\psi_1$ iff $h\in wQw^{-1}(k)$.
\smallskip
{\bf b)} A Weyl element $w_1\in G(k)$ is such that $w_1\bar\psi_1$ is inner isomorphic to $hw\bar\psi_1$ for some $h\in wQw^{-1}(k)$ iff the intersection $P^+(k)w_1\cap wQ(k)$ is non-empty, i.e. iff $w_1\in P^+(k)wQ(k)$
\smallskip
{\bf c)} The actions $\TT$ and $\TT^0$ split up in terms of double cosets of $G(k)$ of the form $P^+(k)wQ(k)$.
\smallskip
{\bf d)} Up to inner isomorphism, the image of $g\bar\psi_1$ depends only on the double coset $P^+(k)gQ(k)$.}
\medskip
{\bf Proof:} a) is a consequence of a) of Step 2 (the use of a) of Step 2 is the seventh place where we need $p>2$). Based on this, b), c) and d) follow from 3.13.7.1.1 ii).
\medskip
Based on a) and on Step 1) we can assume the Faltings--Shimura--Hasse--Witt adjoint map of $(M,\vph,G)$ is $hw\bar\psi_1$, with $h\in wQw^{-1}(k)$. c) implies: for the study of $\TT$ and of $\TT_0$, we always need to fix some element of $DC:=P^+(k)\setminus G(k)/Q(k)$; we presently work with the one defined by $w\in G(k)$. 
\smallskip
 Let $Q^0$ be the image of $Q$ in $G^{\rm ad}_k$. In what follows we identify $N^+_k$ with the unipotent radical of $P^{+0}_k$; similarly we identify $N^-_k$ with a subgroup of $G^{\rm ad}_k$. Let $N^0:=\sg(N^-_k)$ is the unipotent radical of $Q^0$. The quotient variety 
$$S:=wQ^0w^{-1}/wN^0w^{-1},$$
is the moduli variety of maps $hw\bar\psi_1$, with $h\in wQw^{-1}(k)$ (cf. a) and c) of Step 2). 
\smallskip
We consider the subgroup 
$$M_0:=P_k^{+0}\cap wQ^0w^{-1}$$
 of $G^{\rm ad}_k$. It is known (see [SGA3, Vol. III, 4.1.1-2 of Exp. XXVI]), that $M_0$ is a smooth, connected $k$-group; it contains the image ${T_{\rm ad}}_k$ of $T$ in $G^{\rm ad}_k$. For any $m\in M_0(k)$, $mw\bar\psi_1m^{-1}$ is of the form $h_mw\bar\psi_1$, with $h_m\in wQ^0w^{-1}(k)$, and so it defines a unique element of $S(k)$. Let $Z$ be the subgroup (warning: it is not necessarily reduced) of $M_0$ acting trivially on $w\bar\psi_1$ via inner conjugation. We get a well defined orbit morphism
$$m_{\rm orb}:M_0/Z\to S,$$ 
which is an embedding. We still denote by $m_{\rm orb}$ the morphism $M_0\to S$ defining it. Let $s_0\in S(k)$ be defined by $h$ being identity.
\medskip
In what follows, when we want to disregard the Lie algebra structure of the tangent space of the origin of a smooth group $\tilde H$ over $k$, we denote it by $T_0(\tilde H)$. Similarly, we denote by $T_{s_0}(S_1)$ the tangent space at $s_0$ of any subscheme $S_1$ of $S$ containing $s_0$.
\smallskip
Let $M_1$ (resp. $M_2$) be the intersection of $wQ^0w^{-1}$ with $N^+_k$ (resp. with $L^0:=P^{00}_k$). $M_1$ is a normal subgroup of $M_0$ and $M_0$ is the semidirect product of $M_1$ and $M_2$ (this is a consequence of the fact that $w$ is a Weyl element and of the fact that $M_0$ is connected). So $M_1$ and $M_2$ are as well smooth, connected. Any $m_0\in M_0(k)$ can be written uniquely as a product $m_0=m_1m_2$, with $m_i\in M_i(k)$, $i=\overline{1,2}$. From the part of 3.5.3 referring to $b_a$ (or from 3.13.7.1 (*)), we get that 
$$m_0w\bar\psi_1m_0^{-1}=m_1m_2w\sg(m_1^{-1})w^{-1}w\bar\psi_1.$$
For $\tilde m\in G^{\rm ad}(k)$, let 
$$s(\tilde m):=w\sg(\tilde m)w^{-1}.$$
So the map associating $h_{m_0}$ to $m_0=m_1m_2$ is given by: the pair $(m_1,m_2)$ is mapped into
$$m_1m_2s(m_2)^{-1}.\leqno (FOR)$$ 
Moreover, $s(L^0)$ is a Levi subgroup of $wQ^0w^{-1}$; here and in what follows we apply the same conventions pertaining to $\sg$, to the homomorphism $s:G_k\to G_k$. 
\smallskip
We can identify $T_{s_0}(S)$ with $T_0(s(L^0))$. As the tangent map in the origin of the Frobenius map $\sg$ is trivial, under this identification we have 
$$T_0(Z)=T_0(P^{+0}_k\cap s(N_k^-)).\leqno (1)$$ 
Let 
$$Z_1:=N^+_k\cap s(N_k^-).$$ 
It is a connected, smooth normal subgroup of $M_0$, contained in $M_1$ and in $Z$. $M_1/Z_1$ can be identified with $M_3:=N^+_k\cap s(L^0)$; $M_3$ is a smooth, connected group. The connectedness part of $Z_1$ and $M_3$ is implied by the fact that $M_1$ is the semidirect product of $M_3$ and of $Z_1$. Similarly we get that $M_4:=L^0\cap s(N_k^-)$ is connected and smooth and can be identified with a normal subgroup of $M_0/Z_1$. If $y\in M_4(k)$ belongs as well to $Z_{\rm red}/Z_1(k)$, then $ys(y)^{-1}$ belongs to $s(N_k^-)(k)$; so $s(y)\in s(N_k^-)(k)$. So $y$ is trivial. We get that we have a natural homomorphism 
$$q_Z:Z_{\rm red}/Z_1\to M_5:=M_0/Z_1M_4$$
whose kernel is finite and has no $k$-valued points besides the identity.
$q_Z$ is purely inseparable, cf. (1). The study of $Z$ can be pushed forward, using formally closed morphisms ${\rm Spec}(k[[T]])\to Z_{\rm red}/Z_1$ (to be compared with [MFK, ch. 2]). This will not be done here (see \S 9-10 for details). Here we just mention the case when the image of $q_Z$ has dimension $0$. In such a case $\dim_k(Z)=\dim_k(Z_1)$ and so $M_0/Z$ has dimension equal to the sum
$$\dim_k(M_4)+\dim_k(P^{+0}_k\cap s(L^0)).$$
But both $\dim_k(M_4)$ and $\dim_k(s(L^0)\cap N^-_k)$ are equal to
${{\dim_k(L^0)-\dim_k(s(L^0)\cap L^0)}\over 2}$. We conclude: $M_0/Z$ has the same dimension as $L^0$ and so the same dimension as $S=wQ^0w^{-1}/wN^0w^{-1}$. So $m_{\rm orb}$ is an open embedding.
\medskip
{\bf Example.} We assume $L^0\cap s(N_k^-)$ is trivial. Then $Z$ is smooth and $Z_1$ is the connected component of the origin of $Z$. So the image of $q_Z$ has dimension $0$ and so $m_{\rm orb}$ is an open embedding.
\medskip
{\bf 3.13.7.3.1. Standard inequalities.} From the existence of $m_{\rm orb}$ we get
$$\dim_k(Q^0)-\dim_k(N^0)\ge\dim_k(M_0)-\dim_k(Z).$$ 
As $\dim_k(L^0\cap s(P^+_k))=\dim_k(L^0\cap s(P^-_k))$, we get  
$$\dim_k(S)=\dim_k(L^0)=\dim_k(M_2)+\dim_k(M_4).\leqno (2)$$ 
Also from Step 3 of 3.13.7.3 (see def. of $Z_1$ and (1) of it) we get
$$dd((M,\vph,G))=\dim_k(N_k^-)\ge\dim_k(P^{+0}_k\cap s(N_k^-))\ge\dim_k(Z)\ge\dim_k(N^+_k\cap s(N_k^-)).\leqno (3)$$
\medskip
{\bf 3.13.7.4. Complements. A.} We come back to 3.13.7.1. It is easy to see that all stabilizers of points of $\Mx$ w.r.t. $\TT^1$ are smooth. Moreover, the tangent spaces of the orbit maps of $\TT^1$ and $\TT^0$ are the same (as the differential of the Frobenius map $\sg$ is the zero map). We conclude:
\medskip
{\bf Fact.} {\it $\forall x\in\Mx$, its orbit under $\TT^0$ has dimension greater or equal to the dimension of its orbit under $\TT^1$.}
\medskip
{\bf Lemma.} {\it Let $X$ be a $k$-variety (not necessarily connected) which has two stratifications $\Ms_1$ and $\Ms_0$ into locally closed, regular, connected subvarieties. We assume that the number $n(\Ms_1)$ of strata of $\Ms_1$ is finite and that $\forall x\in X(k)$ the dimension of the stratum of $\Ms_0$ to which $x$ belongs is greater or equal to the similarly defined stratum of $\Ms_1$. Then the number of strata of $\Ms_0$ is as well finite.}
\medskip
We leave it to the reader to prove this Lemma using noetherian induction on $n(\Ms_1)$. But the orbits of $\TT^1$ are in general not finite. So the above Fact and Lemma are not enough to conclude that $\TT^0$ has a finite number of orbits. However, 3.3.4, 3.6.6.1 1) and 3.11.1 c) imply:
\medskip
{\bf Corollary.} {\it $\TT$ (and so also $\TT^0$) has a dense open orbit.}
\medskip
{\bf B. The elimination and the insertion processes.} Let $G_i$, $i\in I$, be the simply factors of $G^{\rm ad}_k$. As usual (see 3.4.3.1 and the proof of 3.6.6; see also c) of Step 3) we can assume they are permuted transitively by $\sg$. So we can assume that $I_0:=I=S(1,n)$ and that $\sg$ takes $G_i$ into $G_{i+1}$, with $n+1=1$. Let $I_1$ be as in 3.4.3.4 (i.e. $I_1$ is the subset of $I_0$ corresponding to simple factors of $G_k^{\rm ad}$ not included in $P^{+0}_k$). We can assume $I^1$ is non-empty (otherwise 3.13.7.1.4 follows from Lang's theorem and so 3.13.7.1.3 follows from 3.13.7.2). We also assume $I_0\neq I_1$ (otherwise what follows is void). We order the elements of $I_1$ as 
$$a_1<a_2<...<a_{\abs{I_1}}.$$ 
We can assume $a_1=1$. For each $i\in I_1$, let 
$$n_i:=a_{i+1}-a_i,$$ 
with $a_{\abs{I_1}+1}:=n+a_1$, and let $N^{-i}_k:=N^-_k\cap G_i$. Let 
$$\Mx_{\rm nc}:=N_k^+\setminus \prod_{i\in I_1} G_i/\prod_{i\in I_1} \sg^{n_i}(N^{-i}_k).$$
Let $$P^{00}_k(I_1):=P^{00}_k\cap (\prod_{i\in I_1} G_i).$$
We have an action $\TT^0_{\rm nc}$ of $P^{00}_k(I_1)$ on $\Mx_{\rm nc}$ by the rule: for $h\in P^{00}_k(I_1)\cap G_i(k)$ and $x\in \Mx_{\rm nc}$ defined by an element $g_i\in G_i(k)$, 
${\TT^0_{\rm nc}}_h(x)$ is defined by the element
$$hg_i\sg^{n_i}(h^{-1}).$$
Similarly to $\TT^1$ we define 
$$\TT^1_{\rm nc}: P^{00}_k(I_1)\times\Mx_{\rm nc}\to \Mx_{\rm nc}$$ 
using left translations. It is trivial to see that:
\medskip
{\bf Fact.} {\it The orbits of $\TT^0$ (resp. of $\TT^1$) are in one-to-one correspondence with the orbits of $\TT^0_{\rm nc}$ (resp. of $\TT^1_{\rm nc}$). In particular, each orbit of $\TT^0$ is defined by a $k$-valued point of $\Mx$ defined by an element of $\prod_{i\in I_1} G_i(k)$.}
\medskip
 In other words the $G_i$'s factors defined by $i\in I_0\setminus I_1$ can be ``eliminated".
\medskip
We have a reversed process of ``inserting" non-compact factors. For each $i\in I_0\setminus I_1$, let $j\in I_1$ be the greatest element which is smaller than $i$. We define the following subgroups of $\prod_{i\in I_0} G_i=\prod_{i\in I_1} G_i\times\prod_{i\in I_0\setminus I_1} G_i$:
$$N_k^{+i}:=\sg^{j-i}(N_k^+\cap G_j),$$
$$N_k^{-i}:=\sg^{j-i}(N^{-j}_k),$$
$$P_k^{00}(i):=\sg^{j-i}(P_k^{00}(I_1)\cap G_j),$$
$$N_k^+(I_0):=N_k^+\times \prod_{i\in I_0\setminus I_1} N_k^{+i},$$
$$P_k^{00}(I_0):=P_k^{00}(I_1)\times \prod_{i\in I_0\setminus I_1} P_k^{00}(i).$$
Let 
$$\Mx_{I_0}:=N_k^+(I_0)\setminus \prod_{i\in I_0} G_i/\sg(N_k^-(I_0)),$$ 
where 
$$N_k^-(I_0):=N_k^-\times \prod_{i\in I_0\setminus I_1} N^{-i}_k$$ 
is the opposite of $N_k^+(I_0)$ w.r.t. the action of the image $T_k(I_0)$ of $T_k$ in $\prod_{i\in I_0} G_i$ via inner conjugation. We have a natural action $\TT^0(I_0)$ of $P_k^{00}(I_0)$ on $\Mx_{I_0}$: it is entirely defined in the same manner as $\TT^0$. We have a natural embedding
$$i_{\Mx}:\Mx_{\rm nc}\hookrightarrow \Mx_{I_0}$$
defined by the natural inclusion of $\prod_{i\in I_1} G_i$ into $\prod_{i\in I_0} G_i$ and a monomorphism
$$i_{I_0}:P_k^{00}(I_1)\hookrightarrow P_k^{00}(I_0)$$
which takes $h_i\in (P_k^{00}(I_1)\cap G_i)(k)$ into 
$$\prod_{l=0}^{n_i-1} \sg^l(h_i).$$
The pair $(i_{I_0},i_{\Mx})$ is compatible with the two actions $\TT^0_{\rm nc}$ and $\TT^0(I_0)$, i.e. we have
$$i_{\Mx}({\TT^0_{\rm nc}}_h(x))=\TT^0(I_0)_{i_{I_0}(h)}(i_{\Mx}(x)),$$
with $x\in\Mx_{\rm nc}(k)$ and $h\in P_k^{00}(I_1)(k)$.
\medskip
{\bf Exercise.} {\bf a)} Two $k$-valued points of $\Mx_{\rm nc}$ belong to the same orbit under $\TT^0(I_0)$ iff they belong to the same orbit under $\TT^0_{\rm nc}$ (and so the number of orbits of $\TT^0_{\rm nc}$ is finite if the number of orbits of $\TT^0(I_0)$ is finite). 
\smallskip
{\bf b)} If a $k$-valued point of $\Mx_{\rm nc}$ belongs to the orbit under $\TT^0(I_0)$ of a $k$-valued point of $\Mx_{I_0}$ defined by some Weyl element $w\in\prod_{i\in I_0} G_i(k)$ normalizing $T_k(I_0)$, then it belongs to the orbit under $\TT^0_{\rm nc}$ of some Weyl element $\tilde w\in\prod_{i\in I_1} G_i(k)$ normalizing $T_k(I_0)$.
\smallskip
{\bf Hints:} a) (resp. b)) boils down to the fact that there is no element of $P_k^{00}(i)(k)$ belonging to $N_k^{+i}(k)N_k^{-i}(k)$ (resp. the only Weyl elements normalizing $T_k(I_0)$ and belonging to $N_k^{+i}(k)P_k^{00}(i)(k)N_k^{-i}(k)$ are those belonging to $P_k^{00}(i)(k)$). So for b) one needs to use the classical Bruhat decomposition in the form: $G_i(k)$ is a disjoint union of sets of the form $N_k^{+i}(k)P_k^{00}(i)(k)N_k^{-i}(k)w$, with $w$ a Weyl element of $G_i(k)$ normalizing $T_k(I_0)$. Here $i\in I_0\setminus I_1$.
\medskip
The moral of this Exercise is (cf. also the Fact): if 3.6.13.7.1.4 holds in the context of $\TT^0(I_0)$ then it also holds in the context of $\TT^0$. So to check that 3.6.13.7.1.3 holds in the abstract context of Shimura adjoint Lie $\sg$-crystals, we can essentially restrict to the context of totally non-compact (see 3.10.5) cyclic factors. Here ``essentially" refers to the fact that we still need to check that the maximal number of Weyl elements of $\prod_{i\in I_1} G_i(k)$ normalizing $T_k(I_0)$ and defining elements of $\Mx_{\rm nc}$ belonging to distinct orbits of $\TT^0_{\rm nc}$ is as expected (in practice this is very easy).
\medskip
{\bf C.} We come back to 3.13.7.1. We recall that we are interested only on $g$ mod $p$. Using the classical Bruhat's decomposition (over $k$; see [Bo2, 14.12]), we can assume $g=b_1wb_2$, with $b_1$, $b_2\in B(W(k))$ and with $w$ a Weyl element. Changing $g\vph_1$ by $b_1^{-1}g\vph_1b_1$ and using the fact (see 3.3.1) that $\vph_1b_1=b_1^\prime\vph_1$, with $b_1^\prime\in B(W(k))$, we can assume $g=wb$, with $b\in B(W(k))$. We can assume $b$ is a product of elements of distinct $\GG_a$ subgroups of $B$ normalized by $T$. It is known that there is an arbitrary choice of the order in which such a product is taken (cf. [SGA3, Vol. III, 4.1.2 of p. 172]). So one can try to apply (increasing) induction on the number of factors of such a product which are not identity: one can easily check that many of such factors can be assumed to be trivial. This approach is very practical, if one has to deal with a specific (concrete) situation.   
\medskip
{\bf D. Remarks.} {\bf 1)} From 3.13.7.1.1 a) and the classification of principally quasi-polarized finite, flat, commutative group schemes over $k$ annihilated by $p$ and liftable to $W(k)$ (for instance, see [EO] and [Oo3]) we get that 3.13.7.1 holds if $G=GSp(M,\psi_M)$, with $\psi_M:M\otimes_{W(k)} M\to W(k)(1)$ a principal quasi-polarization of $(M,\vph)$. So 3.13.7.1.3-4 hold provided $G^{\rm ad}$ is a product of simple groups left invariant by $\sg$ and of some $C_n$ Lie type. 
\smallskip
{\bf 2)} One can construct group actions similar to the ones $\TT$ and $\TT^0$ of 3.13.7.1.1, by working modulo $p^m$, with $m\in\NN$ (i.e. by considering truncations mod $p^m$ of Shimura $\sg$-crystals and of Faltings--Shimura--Hasse--Witt shifts). Moreover, 3.13.7.2 has analogues for $m\ge 2$ (and this motivates 3.9.8). See \S 9-10 for details. One can use the language of stacks as well but we think the language of orbit maps (like $\TT^0$'s) is much more natural and practical. 
\smallskip
{\bf 3)} The Bruhat type decompositions hinted at in 3.13.7.1.3-4 can be formulated for many other general contexts. We do expect that all of them can be proved formally, by using a $\sg$-language in the context (see [Bou2, ch. IV, \S 2]) of Tits' systems. We now state a more general form of 3.13.7.1.3 in the context of reductive groups. 
\smallskip
Let $H$ be an an arbitrary reductive group over an algebraically closed field $k$ of characteristic $p$. Let $F:H\to H$ be an epimorphism fixing only a finite number of elements (for instance, cf. [Hu2, \S 8]). Let $P_H$ be an arbitrary parabolic subgroup of $H$ and let $N_H^+$ be its unipotent radical. We assume the existence of a maximal torus $T_H$ of $P_H$ taken by $F$ into itself. Let $P^0_H$ be the Levi subgroup of $P_H$ containing $T_H$. Let $N_H^-$ be the opposite of $N_H^+$ w.r.t. the action of $T_H$ (via inner conjugation) on ${\rm Lie}(H)$. Warning: we not need to assume the existence of a Borel subgroup of $P_H$ taken by $F$ into itself.
\medskip
{\bf Expectation (the Bruhat decomposition in the $F$-context for reductive groups).} {\it Any element $g\in H(k)$ can be written down as a product
$$g=n_+p_0wF(p_0^{-1})F(n_-),$$
with $n_+\in N_H^+(k)$, $n_-\in N_H^-(k)$, $p_0\in P^0_H(k)$ and with $w\in H(k)$ normalizing $T_H$. Moreover the minimal number of such $w$'s needed is precisely the number of elements of the quotient set $W_{P_H}\setminus W_H$, where $W_H$ is the Weyl group of $H$ normalizing $T_H$ and $W_{P_H}$ is its subgroup formed by elements leaving $P_H$ invariant.}  
\medskip
Behind this Expectation and its adjoint variant, there are group actions $\TT_H$ and $\TT^0_H$ as in 3.13.7.1.1. Not to repeat the things, here we just mention that $\TT^0_H$ is the action of the image $P^{00}_H$ of $P^0_H$ in $H^{\rm ad}$ on
$\Mx_H:=N^+_H\setminus H^{\rm ad}/F(N^-_H)$ which is obtained by passing to quotients from the action of $P^{00}_H$ on $H^{\rm ad}$ defined by: $p_{00}\in P^{00}_H(k)$ takes $g\in H^{\rm ad}(k)$ into $p_{00}gF^{\rm ad}(p_{00}^{-1})$. Here $F^{\rm ad}$ is the epimorphism of $H^{\rm ad}$ naturally defined by $F$. See 3.13.7.8 for $\TT_H$ in the $\sg$-context.
\medskip
{\bf 3.13.7.5. Example.} We assume that all simple factors of $G^{\rm ad}$ are of $A_1$ Lie type. We refer to 3.13.7.4 B, with $I=I_0$. As we are dealing with the $A_1$ Lie type, $P_k^{00}(I_1)$ is just a maximal torus of 
$$G(I_1):=\prod_{i\in I_1} G_i.$$ 
So the subgroup $B^+$ (resp. $B^-$) of $G(I_1)$ generated by $N_k^+$ (resp. by $N_k^-$) and by $P_k^{00}(I_1)$ is a Borel subgroup of $G(I_1)$. The classical Bruhat decomposition (applied as in the beginning of Step 1 of 3.13.7.3) tells us that $G(I_1)(k)$ is a disjoint union of double cosets $B^+(k)wB^-(k)=B^+(k)wN_k^-(k)$, with $w$ running through (representatives of) elements of the Weyl group of $G(I_1)$ normalizing $P_k^{00}(I_1)$. As the $A_1$ Lie type does not allow involutions, we have $N_k^-=\prod_{i\in I_1} \sg^{n_i}(N_k^-\cap G_i)$. 
\smallskip
So to check 3.13.7.1.4 in the present context we just need to show (cf. the elimination process of 3.13.7.4 B) that each element of the form $tw$, with $t\in P_k^{00}(I_1)(k)$ and with $w\in G(I_1)(k)$ normalizing $P_k^{00}(I_1)$, can be written as a product $\bigl(\prod_{i\in I_1} t_i\bigr)w\prod_{i\in I_1} \sg^{n_i}(t_i^{-1})$, where $t_i\in (G_i\cap P_k^{00}(I_1))(k)$. Restatement: $t$ can be written as 
$$\Bigl(\prod_{i\in I_1} t_i\Bigr)w\Bigl(\prod_{i\in I_1} \sg^{n_i}(t_i^{-1})\Bigr)w^{-1}.\leqno (FT)$$ But this is a consequence of [Hu2, th. of 18.3], applied to the epimorphism of $P_k^{00}(I_1)$ that takes $\tilde t_i\in (G_i\cap P_k^{00}(I_1))(k)$ into $w\sg^{n_i}(\tilde t_i)w^{-1}$. We can also use systems of equations, close in spirit to the ones of first type, to get this restatement: identifying $P^{00}_k(I_0)$ with a finite product of $\GG_m$'s groups, the actual type of systems of equations we get over $k$ is
$$x_s=c_sx_{\pi(s)}^{p^{n_s}},$$
where $s$ belongs to some set having as many element as the dimension of $P_k^{00}(I_1)$, with $\pi$ a permutation of it, with all $n_s$'s belonging to $\ZZ\setminus\{0\}$ and with all $c_s$'s as non-zero elements of $k$; obviously it has a solution formed by non-zero elements of $k$. Moreover, the map that takes $t$ into the element of (FT) is separable.
\smallskip
From this last paragraph and 3.13.7.2 we conclude:
\medskip
{\bf Corollary 1.} {\it The Expectations of 3.13.7.1 and 3.13.7.1.3-4 hold for the $A_1$ Lie type. Moreover, all stabilizers of $\TT^0$ are smooth.}   
\medskip
Similarly,  we get (based on loc. cit. and not on arguments involving systems of equations) that the Expectation of 3) of 3.13.7.4 D holds for the $A_1$ Lie type. More generally, the same arguments give us:
\medskip
{\bf Corollary 2 (the essentially Borel subgroup case).} {\it The Expectation of 3) of 3.13.7.4 D holds for the case when the image of $P_H$ in each simple factor $SF$ of $H^{\rm ad}$ is either a Borel subgroup of $SF$ or is $SF$ itself. Moreover, all stabilizers of $\TT^0_H$ are smooth.}
\medskip
Moreover, in the context of $\TT^0$ and $\TT^1$ for the $A_1$ Lie type, the Lemma of 3.13.7.4 A does apply (i.e. we do get based on it that $\TT^0$ has a finite number of orbits). 
\medskip
{\bf 3.13.7.6. A variant of the elimination process.} Let $I$ be as in 3.13.7.4 B. We do not assume any more $I=I_0$. Let $J$ be the subset of $I$ formed by elements $i$ such that $G_i$ is not included in the image of $P^+$ in $G^{\rm ad}_k$. Let $\pi$ be the permutation of $I$ defined naturally by $\sg$; so $\sg(G_i)=G_{\pi(i)}$. Based on the Fact and a) of 3.13.7.4 B, we can assume all components of $g$ in $G_i$, with $I\setminus J$, are identity elements. Instead of $\bar\psi_1$ and $g\bar\psi_1$ we can work equally well with their non-compact variants $\bar\psi_1^{\rm nc}$ and $g\bar\psi_1^{\rm nc}$. Here
$$\bar\psi_1^{\rm nc}:\oplus_{i\in J} {\rm Lie}(G_i)\to \oplus_{i\in J} {\rm Lie}(G_i)$$
takes an element $x_i\in {\rm Lie}(G_i)$, with $i\in J$, into $\bar\psi_1^{n_i}(x_i)$, where $n_i\in\NN$ is the smallest element such that $\pi^{n_i}(i)\in J$. We refer to $g\bar\psi_1^{\rm nc}$ as the non-compact Faltings--Shimura--Hasse--Witt adjoint map of $(M,g\vph_1,G)$ or of its attached  Shimura Lie $\sg$-crystal; accordingly the other terminology of 3.9 pertaining to types of Faltings--Shimura--Hasse--Witt maps extend to such a non-compact context. By an inner automorphism of $g\bar\psi^{\rm nc}$, we mean an automorphism of it defined by an element of $\prod_{i\in J} G_i(k)$. Let
$$i_J(p):\prod_{i\in J} G_i\hookrightarrow\prod_{i\in I} G_i$$
be the monomorphism defined as follows: for $g_i\in G_i(k)$, with $i\in J$, 
$$i_j(g_i):=\prod_{l=0}^{n_i-1} \sg^l(g_i)$$ 
(to be compared with $i_{I_0}$ of 3.13.7.4 B). Based on c) of Step 2 of 3.13.7.3 we get:
\medskip
{\bf Fact.} {\it The image under $i_J(p)$ of the group of inner automorphisms of $g\bar\psi_1^{\rm nc}$ is the group of inner automorphisms of $g\bar\psi_1$. In particular, these two groups of inner automorphisms have the same dimension.}
\medskip
So we can perform the elimination process not only in connection to finding the orbits of $\TT^0$ but also in connection to finding groups of inner automorphisms of Faltings--Shimura--Hasse--Witt adjoint maps (over $k$). From this Fact and Corollary 1 of 3.13.7.5 we get:
\medskip
{\bf 3.13.7.6.0. Corollary.} {\it We assume $G_i$ is of $A_1$ Lie type, $\forall i\in I$. Then the number of orbits of $\TT^0$ is $2^{\abs{J}}$: each such orbit $o(w)$ is defined by a Weyl element $w\in\prod_{i\in J} G_i(k)$ normalizing the image $T(J)_k$ of $T_k$ in $\prod_{i\in J} G_i$; moreover, if $w_1\in \prod_{i\in J} G_i(k)$ is another Weyl element, we have $o(w)=o(w_1)$ iff $w_1=tw_2$ for some $t\in T(J)_k(k)$. So each such orbit $o(w)$ defines a function 
$$f_w:J\to\{0,1\}$$
by the rule: $f(i)=0$ iff the component of $w$ in $G_i$ does belong to $T(J)_k(k)$. 
\smallskip
If $B^+$ (resp. $B^-$) is the image of $P^+$ (resp. of $P^-$) in $\prod_{i\in J} G_i$, then under the projection of $G^{\rm ad}_k$ on $\prod_{i\in J} G_i$, $o(w)$ is mapped onto 
$$B^+(k)wB^-(k).$$ 
The group of automorphisms of $w\bar\psi_1^{\rm nc}$ (and so of $w\bar\psi_1$) has dimension equal to 
$$\sum_{i\in J} f_w(i).$$ 
Moreover $o(w)$ specializes to $o(w_1)$ iff 
$$f_{w_1}(i)\ge f_w(i),\, \forall i\in J.\leqno (SPEC)$$}
\medskip
{\bf Proof:} The only thing left to be remarked is that the dimension $d_w$ of the group of inner automorphisms of $w\bar\psi_1^{\rm nc}$ is nothing else but 
$$\dim_k(B^+)-\dim_k(B^+wB^-/B^-).\leqno (DIM)$$
The argument goes as follows. In order to have $g\bar\psi_1^{\rm nc}g^{-1}$ of the form $h\bar\psi_1^{\rm nc}$ (here $g$, $h\in\prod_{i\in J} G_i(k)$), we need $g\in B^+(k)$, cf. the normalizing part of a) of Step 2 of 3.13.7.3 (see also Step 3 of 3.13.7.3). Moreover, $h\bar\psi_1^{\rm nc}=\bar\psi_1^{\rm nc}$ iff $h$ is a $k$-valued point of the unipotent radical of $B^-(k)$, cf. the centralizing part of a) of Step 2 of 3.13.7.3; so the variety parameterizing maps of the form $h\bar\psi_1^{\rm nc}$, with (see c) of Step 3 of 3.13.7.3) $h\in B^+(k)wB^-(k)$, has the same dimension as (in fact it is isomorphic to) $B^+wB^-/B^-$. So $\dim_k(B^+)-d_w=\dim_k(B^+wB^-/B^-)$. This ends the proof.
\medskip
{\bf 3.13.7.6.1. The non-compact context.} The whole of 3.13.7.3 can be redone in the context of $\bar\psi_1^{\rm nc}$. We just need to:
\medskip
--  replace the role of $\sg$ with the one of the endomorphism $F:\prod_{i\in J} G_i\to \prod_{i\in J} G_i$ that takes $g_i\in G_i(k)$ into $\sg^{n_i}(g_i)$, with $n_i\in\NN$ as in 3.13.7.6;
\smallskip
-- consider a Weyl element $w\in\prod_{i\in J} G_i(k)$ normalizing the image $T_k(J)$ of $T$ in $\prod_{i\in J} G_i$;
\smallskip
-- replace $s$ by the endomorphism $s_J$ of $\prod_{i\in J} G_i$ that takes $g_i\in G_i(k)$ into $wF(g_i)w^{-1}$. 
\medskip
Let $P^{+}(J)$ (resp. $Q(J)$) be the image of $P_k^{+}$ (resp. of $F(P^-_k)$) in $\prod_{i\in J} G_i$. Let $P^-(J)$ be the opposite of $P^+(J)$ w.r.t. the action (via inner conjugation) of $T_k(J)$; so $Q(J):=F(P^-(J))$. Let $\TT^0_{\rm nc}$ and $\Mx_{\rm nc}$ be as in 3.13.7.4 but obtained working with $J$ instead $I_1$ (i.e. we do not assume that we are in a cyclic context). Let $P^0(J)$ be the image of $P^0$ in $\prod_{i\in J} G_i$. Let $N(J)$ be the unipotent radical of $Q(J)$. As in c) of the Fact of Step 3 of 3.13.7.3, the action $\TT^0_{\rm nc}$ splits in terms of the image in $\Mx_{\rm nc}$ of double cosets of the form $P^{+}(J)wQ(J)$. So for its study we need to fix an element of
$$DC_{\rm nc}:=P^{+}(J)(k)\setminus \prod_{i\in J} G_i(k)/Q(J)(k).$$
For future references, we point out:
\medskip
{\bf A. Fact.} {\it Let $g\in\prod_{i\in J} G_i(k)$. If we are in the cyclic context (i.e. we have $I=I_0$ and $J=I_1$), then the Lie stable $p$-rank of $g\bar\psi_1^{\rm nc}$ is ${{\abs{I_1}}\over {\abs{I_0}}}$ times the Lie stable $p$-rank of $g\bar\psi_1$.}
\medskip
Let 
$$M_0^{\rm nc}:=P^{+}(J)\cap wQ(J)w^{-1}=P^{+}(J)\cap s_J(P^-(J))$$ 
and let $Z^{\rm nc}$ be its subgroup fixing $w\bar\psi_1^{\rm nc}$. Let 
$$Z_1^{\rm nc}:=N_k^+\cap s_J(N_k^-).$$ 
It is a smooth, connected subgroup of $Z^{\rm nc}$. As in 3.13.7.3.1 we get
$$dd((M,\vph,G))\ge\dim_k(P^{+}(J)\cap s_J(N_k^-))\ge\dim_k(Z^{\rm nc})\ge\dim_k(N_k^+\cap s_J(N_k^-)),\leqno (4)$$
and
$$\dim_k(Q(J))-\dim_k(N(J))\ge\dim_k(M_0^{\rm nc})-\dim_k(Z^{\rm nc}).\leqno (5)$$
\indent
Let $w_0\in\prod_{i\in J} G_i(k)$ be a Weyl element with the property that the intersection
$$P^+(J)\cap w_0Q(J)w_0^{-1}$$
contains the image $B_k(J)$ of $B$ in $\prod_{i\in J} G_i$; the image of $w_0$ in $DC_{\rm nc}$ is uniquely determined. We have:
\medskip
{\bf B. Claim.} {\it If $dd((M,\vph,G))=\dim_k(Z^{\rm nc})$, then $w\in P^{+}(J)(k)w_0Q(J)(k)$.}
\medskip
{\bf Proof:} From (4) we get $s_J(N_k^-)$ is a subgroup of $P^{+}(J)$. Based on [Bo2, 14.22 (i)], this is equivalent to: the intersection $P^+(J)\cap s_J(P^-(J))=P^+(J)\cap wQ(J)w^{-1}$ contains a Borel subgroup of $\prod_{i\in J} G_i$. As we are dealing with double cosets, we can assume that this Borel subgroup is $B_k(J)$. The Claim follows. 
\medskip
{\bf C. Theorem.} {\it The following two statements are equivalent:
\medskip
a) $N_k^+=s_J(N_k^-)$;
\smallskip
b) $w\in P^+(J)w_0Q(J)$, the image in $\Mx_{\rm nc}$ of the double coset $P^+(J)w_0Q(J)$ is one orbit of the action $\TT^0_{\rm nc}$ and we have $dd((M,\vph,G))=\dim_k(Z^{\rm nc})$.}
\medskip
{\bf Proof:} If $N_k^+=s_J(N_k^-)$, then $P^+(J)=s_J(P^-(J))=wQ(J)w^{-1}$; so $P^0(J)=s_J(P^0(J))$ and $w\in P^+(J)w_0Q(J)$. We get: $P^{+}(J)(k)w_0Q(J)(k)=P^+(J)(k)w_0$. So the image $IM$ of $P^{+}(J)w_0Q(J)$ in $\Mx_{\rm nc}$ can be naturally identified with $P^0(J)$ (cf. the Hint of a) of 3.13.7.4 B). Under this identification, the automorphism of $IM$ defined naturally by $h\in P^0(J)(k)$ via the action $\TT^0_{\rm nc}$, is nothing else but the automorphism of $P^0(J)$ that takes $h_1\in P^0(J)(k)$ into 
$$hh_1s_J(h)^{-1}$$ 
(this is the same as (FOR) of Step 3 of 3.13.7.3). So from Lang's theorem in the form of [Hu2, th. 8.3] (applied to the homomorphism $F=s_J$) we get that $IM$is one orbit of $\TT^0_{\rm nc}$. The last part of b) is a consequence of (4).
\smallskip
We assume now that b) holds. To check a) we can assume $w=w_0$. Due to the ``one orbit part", in (5) we have equality. As $\dim_k(N(J))=dd((M,\vph,G))$, we get $\dim_k(M_0^{\rm nc})=\dim_k(Q(J))=\dim_k(P^+(J))$. So $M_0^{\rm nc}=P^+(J)$ and so by reasons of dimensions $P^+(J)=s_J(P^-(J))$. So a) holds. This ends the proof.
\medskip
{\bf D.} The double coset $P^+(J)w_0Q(J)$ is referred as the final coset and any Shimura $\sg$-crystal whose non-compact Faltings--Shimura--Hasse--Witt map is inner isomorphic to $g_{\rm fin}\bar\psi_1^{\rm nc}$, with $g_{\rm fin}\in P^+(J)w_0Q(J)$ is referred as a quasi-final Shimura $\sg$-crystal. Similarly, we speak about Shimura adjoint Lie $\sg$-crystals which are quasi-final. 
\medskip
{\bf 3.13.7.6.2. The list.} The equality $N_k^+=s_J(N_k^-)$ is equivalent to the fact that each cyclic adjoint factor $\Mc$ of $(M,\vph_1,G)$ is of one of the following types:
\medskip
{\bf a)} $B_{\ell}$, $C_{\ell}$ or $D_{\ell}^{\RR}$ type;
\smallskip
{\bf b)} of $A_{\ell}$ type with ${\ell}$ odd and which is also (of constant) type ${{{\ell}+1}\over 2}$;
\smallskip
{\bf c)} of $A_{\ell}$ type whose concentrated $\vep$-type is involutive;
\smallskip
{\bf d)} of $D_{\ell}^{\HH}$ type with ${\ell}$ odd and greater than $4$, whose concentrated $\vep$-type is involutive;
\smallskip
{\bf e)} of simple $D_{\ell}^{\HH}$ type with ${\ell}$ even and greater than $3$.
\medskip
This list conforms with the expression of $w_0$'s in [Bou2, planche I to IV]. 
\medskip
{\bf 3.13.7.6.3. The cyclic diagonalizability property.} We assume that the two statements of 3.13.7.6.1 C hold. We consider a Weyl element $\om$ of $G(W(k))$ whose image in $G^{\rm ad}_k$ is the same as the image of $w_0$. We get a cyclic diagonalizable  
Shimura filtered $\sg$-crystal $(M,F^1,\om\vph_1,T)$. We have:
\medskip
{\bf Theorem.} {\it $(M,\om\vph_1,G)$ is an $isom$-canonical Shimura $\sg$-crystal.}
\medskip
{\bf Corollary (the cyclic diagonalizability property).} {\it Any Shimura $\sg$-crystal $(M,g\vph_1,G)$, with $g\in G(W(k))$, whose Faltings--Shimura--Hasse--Witt map is inner isomorphic to the one of $(M,\om\vph_1,G)$, is isomorphic to $(M,\om\vph_1,G)$ and so is cyclic diagonalizable.}
\medskip
{\bf Proofs:} To prove the Theorem, we need to show: if $g\in {\rm Ker} (G(W(k))\to G(k))$, then $(M,\om\vph_1,G)$ and $(M,g\om\vph_1,G)$ are isomorphic under an isomorphism defined by an element of $G(W(k))$. We follow the pattern of the proof of 3.6.15 B. By induction on $m\in\NN$ we show that we can assume $g\in {\rm Ker} (G(W(k))\to G(W_m(k)))$. It is enough to show that: if $g\in {\rm Ker} (G(W(k))\to G(W_m(k)))$, then there is $h\in G(W(k))$ congruent to the identity mod $p^{m-1}$ and such that the element $\tilde g\in G(W(k))$ defined by the equality $\tilde g\om\vph_1=hg\om\vph_1h^{-1}$ is congruent to the identity mod $p^{m+1}$. Let $\tilde\vph_1:=\om\vph_1$.
\smallskip 
We use the notations of 3.13.7.6 to index the simple factors of $G^{\rm ad}$. So $G^{\rm ad}=\prod_{i\in I} G_{iW(k)}$, with the special fibre of the simple group $G_{iW(k)}$ being $G_i$. 
Let $P^0_{W(k)}(J)$ be the reductive subgroup of $G$ normalized by $T$ and whose special fibre is the pull back to $G_k$ of $P^0(J)$. Let $N^+$ and $N^-$ be as in 3.13.7.1.1. Let $P^0(I)$ be the smallest reductive subgroup of $G$ containing $P_{W(k)}^0(J)$ and stable under the endomorphism $i_{\om}\circ\sg$ of $G$ which takes $u\in G(W(k))$ into $\om\sg(u)\om^{-1}$. The images of $P^0(I)$ and of $P_{W(k)}^0(J)$ in $\prod_{i\in J} G_{iW(k)}$ are the same (cf. a) of 3.13.7.6.1 C). We choose $W(k)$-basis $\Mb_0$, $\Mb_+$ and $\Mb_-$ of the Lie algebras of $P^0(I)$, $N^+$ and respectively of $N^-$. Conjugating $g\tilde\vph_1$ by elements of $G(W(k))$ congruent to $1_M$ mod $p^m$ and whose images in $G^{\rm ad}(W(k))$ belong to $\prod_{i\in I\setminus J} G_{iW(k)}(W(k))$, we can assume
$$g=g_cg_-g_+g_0$$ 
where $g_+\in N^+(W(k))$, $g_-\in N^-(W(k))$ and $g_0\in P^0(I)(W(k))$ are congruent to the identity mod $p^m$, and where $g_c\in G(W(k))$ is congruent to the identity mod $p^{m+1}$. The order in which such a product decomposition of $g$ is taken is entirely irrelevant mod $p^{2m}$ and so mod $p^{m+1}$. 
\smallskip
We take $h$ of the form $h=h_3h_2h_1$, where $h_3\in G(W(k))$ and $h_1\in P^0(I)(W(k))$ are congruent to the identity mod $p^m$ and where 
$$h_2=i_J(h_+)i_J(h_-);$$
here 
$$i_J:\prod_{i\in J} G_{iW(k)}\hookrightarrow G^{\rm ad}$$ is a monomorphism defined entirely as $i_J(p)$ of 3.13.7.6 and $h_+\in N^+(W(k))$ (resp. $h_-\in N^-(W(k))$) is (warning!) congruent to the identity mod $p^{m-1}$ (resp. mod $p^m$). $i_J(h_+)$ and $i_J(h_-)$ still make sense as the unipotent subgroups of $G$ are in natural one-to-one correspondence to unipotent subgroups of $G^{\rm ad}$. $h_1$, $g_-$ and $g_+$ are commuting mod $p^{m+1}$. So we get
$$\tilde g=h_3h_2h_1g\tilde\vph_1h_1^{-1}\tilde\vph_1^{-1}\tilde\vph_1h_2^{-1}h_3^{-1}\tilde\vph_1^{-1}=\tilde g_ch_3h_2g_-g_+\tilde g_0\tilde\vph_1h_2^{-1}h_3^{-1}\tilde\vph_1^{-1},$$
where $\tilde g_c\in G(W(k))$ is congruent to the identity mod $p^{m+1}$ and where 
$$\tilde g_0:=h_1g_0\tilde\vph_1h_1^{-1}\tilde\vph_1^{-1}.$$
 As $P^0(I)$ is stable under the endomorphism $i_{\om}\circ\sg$ of $G$, ${\rm Lie}(P^0(I))$ is stable under $\tilde\vph_1$. So writing $1_{M}-h_1$ mod $p^{m+1}$ as $p^m$ times a linear combination of elements of $\Mb_0$ with coefficients in $W(k)$, the condition that $\tilde g_0$ is congruent to the identity mod $p^{m+1}$ gets translated in choosing the reduction mod $p$ of these coefficients to form a solution of a quasi Artin--Schreier system of equations in $\dim_k(P^0(I)_k)$ variables with coefficients in $k$. So based on 3.6.8.1 (with $l=1$ and $s=0$) we can assume $\tilde g_0$ is congruent to the identity mod $p^{m+1}$. So $\tilde g$ is of the form $\tilde g_ch_3\tilde g_{+-}\tilde\vph_1h_3^{-1}\tilde\vph_1^{-1}$, where 
$$\tilde g_{+-}:=h_2g_-g_+\tilde\vph_1h_2^{-1}\tilde\vph_1^{-1}.$$ 
We consider two cases.
\medskip
{\bf Case 1: $m\ge 2$.} Every element of $G(W(k))$ which mod $p^m$ is $1_M$ commutes mod $p^{2m-1}$ with any element of $G(W(k))$ which mod $p^{m-1}$ is $1_M$. So $\tilde g_{+-}$ mod $p^{2m-1}$ (and so mod $p^{m+1}$) can be rewritten as $\tilde g_+\tilde g_-$, where 
$$\tilde g_+:=i_J(h_+)g_+\tilde\vph_1i_J(h_-^{-1})\tilde\vph_1^{-1}$$
 and 
$$\tilde g_-:=i_J(h_-)g_-\tilde\vph_1i_J(h_+^{-1})\tilde\vph_1^{-1}.$$ 
Based on a) of 3.13.7.6.1 C), we get that $\tilde g_+\in N^+(W(k))$ and $\tilde g_-\in N^-(W(k))$. We write $1_{M}-h_-$ mod $p^{m+1}$ as $p^m$ times a linear combination $C_-$ of elements of $\Mb_-$ with coefficients in $W(k)$. We write $h_+=h_+^1h_+^2$, with $h_+^1\in N^+(W(k))$ congruent to the identity mod $p^m$ and with $h_+^2\in N^+(W(k))$ congruent to the identity mod $p^{m-1}$. We write $1_M-h_+^1$ (resp. $1_M-h_+^2$) as $p^m$ (resp. as $p^{m-1}$) times a linear combination $C_+^1$ (resp. $C_+^2$) of elements of $\Mb_+$ with coefficients in $W(k)$. The conditions that both $\tilde g_+$ and $\tilde g_-$ are congruent to the identity mod $p^{m+1}$ get translated in choosing the reduction mod $p^2$ of all these coefficients to form a solution of a system $SE$ of equations with coefficients in $W_2(k)$. We split the study of its solutions in two steps. First we consider its subsystem $SUBSE$ of equations with coefficients in $k$ in the variables which are the reduction mod $p$ of all coefficients of $C_-$ and of $C_+^2$ and obtained by requiring $\tilde g_-$ to be congruent to the identity mod $p^{m+1}$ and $\tilde g_+$ to be congruent to the identity mod $p^m$: it is a quasi Artin--Schreier system of equations. We can not easily give a formula for the number of variables; this is so due to the fact that the definition of $s_J$ involves iterates of $\sg$ and so we need to appeal to 3.6.8.1.5; if we allow monomials as in the mentioned place, then we are dealing with systems of equations (not a priori of first type) which involve precisely $2\dim_k(N^+_k)$ variables. Second we pick up a solution of $SUBSE$ and we plug it (via a fixed lift $h_+^1$ of its truncation mod $p^m$) in $SE$. We get a new system of equations with coefficients in $k$ in the variables which are the reduction mod $p$ of all coefficients of $C_+^1$ and obtained by requiring $\tilde g_+$ to be congruent to the identity mod $p^{m+1}$: it is a quasi Artin--Schreier system of equations. The situation with the number of variables is as above, with $2\dim_k(N^+_k)$ replaced by $\dim_k(N^+_k)$. 
\smallskip
So we can assume $\tilde g_{+-}$ is congruent to the identity mod $p^{m+1}$. Choosing $h_3$ to be $1_M$, we get $\tilde g$ is congruent to the identity mod $p^{m+1}$.
\medskip
{\bf Case 2: $m=1$.} This case is a little bit complicated as we can not assume that we are dealing with elements which commute mod $p^2$ (see the choice of $h_+$). However, based on 3.5.3 (5) 
$$g_+\bigl(\tilde\vph_1i_J(h_-^{-1})\tilde\vph_1^{-1}\bigr)\bigl(\tilde\vph_1i_J(h_+^{-1})\tilde\vph_1^{-1}\bigr)$$ 
mod $p^2$ can be rewritten as 
$$g_+\bigl(\tilde\vph_1i_J(h_+^{-1})\tilde\vph_1^{-1}\bigr)\bigl(\tilde\vph_1i_J(h_-^{-1})\tilde\vph_1^{-1}\bigr)p_{+-0}$$
and so as
$$\tilde\vph_1i_J(h_+^{-1})\tilde\vph_1^{-1}g_+\tilde\vph_1i_J(h_-^{-1})\tilde\vph_1^{-1}p_{+-0},$$ 
with $p_{+-0}\in P^0(I)(W(k))$ congruent to the identity mod $p$. So $\tilde g_{+-}$ mod $p^2$ is congruent to 
$$i_J(h_+)EL_-g_+\tilde\vph_1i_J(h_-^{-1})\tilde\vph_1^{-1}p_{+-0},$$ 
where 
$$EL_-:=i_J(h_-)g_-\tilde\vph_1i_J(h_+^{-1})\tilde\vph_1^{-1}.$$
\indent
So, writing $1_{M}-h_-$ (resp. $1_M-h_+$) mod $p^{2}$ (resp. mod $p$) as $p$ times (resp. as) a linear combination of elements of $\Mb_-$ (resp. of $\Mb_+$) with coefficients in $W(k)$, the conditions that $EL_-$ is congruent to the identity mod $p^2$ and $i_J(h_+)g_+\tilde\vph_1i_J(h_-^{-1})\tilde\vph_1^{-1}$ is congruent to the identity mod $p$ get translated in choosing the reduction mod $p$ of these coefficients to form a solution of a quasi Artin--Schreier system of equations with coefficients in $k$; the situation with the number of variables is as in above part involving $SUBSE$. So we can assume $\tilde g_{+-}$ mod $p^2$ is of the form $g_+^\prime p_{+-0}$, where $g_+^\prime\in N^+(W(k))$ is congruent to  the identity mod $p$. Choosing $h_3$ to be the product of an element of $P^0(I)(W(k))$ congruent to the identity mod $p$ with $(g_+^\prime)^{-1}$, the arguments of the paragraph before Case 1 referring to $\tilde g_0$  allow us to assume $\tilde g_{+-}$ (and so $\tilde g$) is congruent to the identity mod $p^2$.  
\medskip
So regardless of how $m$ is, we can assume $\tilde g$ is congruent to the identity mod $p^{m+1}$. This proves the Theorem. We now refer to the Corollary. We can assume the truncations mod $p$ of $(M,g\vph_1,G)$ and $(M,w\vph_1,G)$ are inner isomorphic, cf. Fact of 3.13.7.2. So from the construction of $\TT$ we get: we can assume $g$ is congruent to $w$ mod $p$. So the Corollary follows from the Theorem. This ends the proofs.
\medskip
{\bf 3.13.7.6.3.1. An approach.} We include now an approach to be used to present some examples in 3.15.7 H and K below.
\smallskip
We follow the pattern of the proof of 3.13.7.6.3 and we use the previous notations. Let $w$ be a Weyl element of the final double coset. We use Weyl's decomposition ${\rm Lie}(G)={\rm Lie}(T)\oplus_{\al\in\Phi} {{\got g}_{\al}}$ w.r.t. $T$. So $\Phi$ is the set of cocharacters of the action of $T$ on ${\rm Lie}(G)$ and ${\got g}_{\al}$ is the direct summand of ${\rm Lie}(G)$ on which $T$ acts via $\al$. Let $\tilde\vph_1:=w\vph_1$. $\tilde\vph_1$ permutes ${\got g}_{\al}[{1\over p}]$'s; so wet a permutation $\pi$ of $\Phi$. We concentrate on one cycle $(\al_1,...\al_s)$ of it; whatever we do below for it is implicitly done for all other cycles of $\pi$. For $i\in S(1,s)$, let $u_i\in S(-1,1)$ be such that $\be\in\GG_m$ acts through $\mu$ on ${\got g}_{\al}$ as the multiplication with $p^{-u_i}$.
\smallskip
As $s_J(N_k^-)$ is a subgroup of $P^+(J)$, $s_J(N_k^+)$ is a subgroup of $P^-(J)$. This implies:
\medskip
{\bf Key property.} {\it There is no $i\in S(1,s)$ such that $u_i$ and $u_{i+1}$ are both $1$ or both $-1$ (here $u_{s+1}=u_1$).}
\medskip
We consider the following modification 
$${\rm Lie-mod}(G):={\rm Lie}(T)\oplus p{\rm Lie}(N_-)\oplus {\rm Lie}(N^+)\oplus pL^{0,-}\oplus L^{0,0}\oplus L^{0,+}$$
of ${\rm Lie}(G)$. Here $L^{0,+}$, $L^{0,0}$ and $L^{0,-}$ are direct sums of ${\got g}_{\al}$'s subject to the following rules:
\medskip
-- their direct sum is ${\rm Lie}(P^0)$;
\smallskip
-- ${\got g}_{\al_i}\subset L^{0+}$ (resp. ${\got g}_{\al_i}\subset L^{00}$) iff $n_i=0$ and the first non-zero term of the sequence $(n_j)_{j\ge i}$ is $1$ (resp. iff $(n_1,...,n_s)=(0,...,0)$). 
\medskip
We choose a generator $x_i$ of ${\got g}_{\al_i}\cap {\rm Lie-mod}(G)$, $i\in S(1,s)$. We obviously have:
\medskip
{\bf Fact.} {\it We assume that there is no cycle $(\al_1,...,\al_s)$ of $\pi$ such that $\sum_{i=1}^s n_i\ge 0$ and $(-1,0,0,...,0,-1)=(n_j,n_{j+1},...,n_l)$ for some $j\in S(1,n)$ and $l\in S(j,j+s-1)$. Then for $i\in S(1,s)$, writing $\tilde\vph_1(x_i)=p^{v_i}\gamma_ix_{i+1}$, with $\gamma_i\in\GG_m(W(k))$ and $v_i\in S(-1,1)$, all $v_i$, $i\in S(1,s)$, have the same sign, i.e. either all of them are non-negative or all of them are non-positive.}
\medskip
{\bf 3.13.7.6.3.2. Remark.} 3.13.7.6.3 and 3.13.7.6.1 remain true in a principally quasi-polarized context (cf. its proof and 3.6.15 B and its proof).
\medskip
{\bf 3.13.7.7. Remark.} 3.13.7.1-4, 3.13.7.6 and 3.13.7.6.1 make sense as well in a purely Lie context involving Shimura adjoint Lie $\sg$-crystals: no modifications of arguments are required (besides what mentioned in the beginning paragraph of 3.4). So to the list of 3.13.7.6.2 we can add (again cf. the expression of $w_0$ in [Bou, planche V and VI]) the following two extra types:
\medskip
{\bf f)} of $E_7$ Lie type;
\smallskip
{\bf g)} of $E_6$ Lie type whose concentrated $\vep$-type is involutive.
\medskip
Moreover, 3.13.7.6.3 can be entirely adapted to Shimura-ordinary $p$-divisible objects over $k$ whose attached Shimura adjoint Lie $\sg$-crystals are such that f) or g) holds for each one of their cyclic factors.
\medskip
{\bf 3.13.7.8. The case of Fontaine categories $p-\Mm\Mf_{[a,b]}(W(k))$.} We refer to 3) of 3.13.7.4 D. Its Expectation has interpretations similar to the ones of 3.13.7.1 and of 3.13.7.1.1. In what follows, we exemplify this phenomenon just in the general context pertaining to Fontaine categories. So we assume $H$ is the extension to $k$ of a reductive group $H_{\FF_p}$ over $\FF_p$ in such a way that $F$ is nothing else but the endomorphism $\sg$ of $H=H_k$. 
\smallskip
Let $H_{\ZZ_p}$ be a reductive group over $\ZZ_p$ whose special fibre is $H_{\FF_p}$. We consider a faithful representation 
$$\rho:H_{\ZZ_p}\hookrightarrow GL(M_{\ZZ_p}),$$ 
with $M_{\ZZ_p}$ a free $\ZZ_p$-module of finite rank. We consider a cocharacter $\mu:\GG_m\to H_{W(k)}$ such that $P_H$ lifts to a parabolic subgroup $P_{W(k)}$ of $H_{W(k)}$ in such a way that ${\rm Lie}(P_{W(k)})$ is the maximal Lie subalgebra of ${\rm Lie}(H_{W(k)})$ on which $\mu$ acts (via inner conjugation) through non-negative, integral powers of the identical character of $\GG_m$. $[\mu]$ is not uniquely determined. Its existence can be deduced immediately from the structure of parabolic subgroups of $H$ (see [Bo2, 14.17]): as in the proof of Fact 1 of 2.2.9 3) we can work just over $k$ and we can assume $\mu_k$ factors through a maximal torus $T$ of $P_H$; we need the composite of $\mu_k$ with a specified set $CHAR$ of independent (over $\ZZ$) characters of $T$ to be specified a priori. Here $CHAR$ is a set of characters of $T$ forming a bases of the root system of the action of $T$ on ${\rm Lie}(H)$, such that any Lie subalgebra of ${\rm Lie}(H)$ corresponding to an element of $CHAR$ is included in ${\rm Lie}(P_H)$.
\smallskip
For $g\in H_{W(k)}(W(k))$ we construct a $p$-divisible object ${\got C}=(M,(F^s(M))_{s\in S(a,b)},\vph)$ of $\Mm\Mf_{[a,b]}(W(k))$, with $a$, $b\in\ZZ$, $b\ge a$, as follows. Let $M:=M_{\ZZ_p}\otimes_{\ZZ_p} W(k)=\oplus_{s=a}^b \tilde F^s$ be the direct sum decomposition defined by $\mu$ as in 2.2.1.2. Let $F^s(M):=\oplus_{j=s}^b \tilde F^j$. We take $\vph$ to be $g\circ\sg\circ\mu({1\over p})$, where $\sg$ is identified with the $\sg$-linear automorphism of $M$ fixing $M_{\ZZ_p}$ and acting as $\sg$ on $W(k)$. 
The quadruple 
$$(M,(F^s(M))_{s\in S(a,b)},\vph,H_{W(k)})$$ 
is a $p$-divisible object with a reductive structure of $\Mm\Mf_{[a,b]}(W(k))$ (see 2.2.8 3a)). We recall (see 2.2.8 4a) and 2.2.1 d)) that we speak about lifts of $(M,\vph,H_{W(k)})$ and about ${\got C}$ or some truncations of it being cyclic diagonalizable. As in 3.13.7, provided $(M,\vph,H_{W(k)})$ has no lift which is cyclic diagonalizable, we speak about the CM level of $(M,\vph,H_{W(k)})$. 
\smallskip
We consider the direct sum decomposition of 
$${\rm Lie}(H_{W(k)})=\oplus_{j\in S(\mu)} {\got n}(j)$$
produced by $\mu$ via inner conjugation; so $\be\in\GG_m(W(k))$ acts through $\mu$ on ${\got n}_j$ as the multiplication with $\be^{-j}$. Here $S(\mu)\subset\ZZ$ is a finite set. Let $T_{W(k)}$ be a maximal torus of $H_{W(k)}$ lifting $T$ and through which $\mu$ factors. Let $N_0$ be the centralizer of $\mu$ in $H_{W(k)}$. It is a reductive subgroup of $H_{W(k)}$ whose Lie subalgebra is ${\got n}_0$. For $j\in S(\mu)\setminus\{0\}$, let $N_j$ be the integral, connected, smooth subscheme of $H_{W(k)}$ which is the product (taken in some order) of the distinct $\GG_a$ subgroups of $H_{W(k)}$ normalized by $T_{W(k)}$ and whose Lie algebras are included in ${\got n}(j)$. For what follows it is irrelevant the order in which such products are taken, cf. [SGA3, Vol. III, 4.1.2 of p. 172]. The tangent space of $N_j$ in the origin of $H_{W(k)}$ is ${\got n}_j$. Let $N^+$ (resp. $N^-$) be the integral, connected, smooth, unipotent subgroup of $H_{W(k)}$ generated by all $N_j$'s with $j>0$ (resp. with $j<0$). 
\smallskip
Let $h\in H_{W(k)}(W(k))$ be such that $(M,(F^s(M))_{s\in S(a,b)},h\vph h^{-1},H_{W(k)})$ is still a $p$-divisible object with a reductive structure of $\Mm\Mf_{[a,b]}(W(k))$. From Fact of 2.2.14.2 we get that $h(F^s(M))\subset \sum_{v=0}^{s-a} p^vF^{s-v}(M)$, $\forall s\in S(a,b)$. It is an easy exercise to see that this implies that we can write 
$$h=h_+h_-h_0,$$
where 
$$h_+=\prod_{j\in S(\mu)\cap (0,\infty)} h_j$$
and 
$$h_-=\prod_{j\in S(\mu)\cap (-\infty,0)} h_j=\prod_{j\in S(\mu)\cap (-\infty,0)} (1_M+p^{-j}u_j),$$
with $h_j\in N_j(W(k))$ congruent to the identity mod $p^{{\rm max}\{0,-j\}}$, $\forall j\in S(\mu)$.
It is loc. cit. (resp. loc. cit. and the fact that, $\forall j_1$, $j_2\in S(\mu)\cap (-\infty,0)$, the commutator $n_{j_1}n_{j_2}n_{j_1}^{-1}n_{j_2}^{-1}$, with $n_{j_i}\in N_{j_i}(W(k))$, $i=\overline{1,2}$, is a $W(k)$-valued point of the subgroup of $H_{W(k)}$ generated by $N_j$'s with $j\in S(\mu)\cap (-\infty,j_1+j_2)$) which tells us that it is irrelevant (for what follows) the order in which the product of $h_+$ (resp. of $h_-$) is taken. 
\smallskip
Based on the previous paragraph, as in 3.13.7.1 we get: $hg\vph h^{-1}\vph^{-1}\in H_{W(k)}(W(k))$ mod $p$ is of the same form $m_1m_2\bar g\sg(m_2^{-1})\sg(m_3^{-1})$ as (*) of 3.13.7.1, $\forall g\in H_{W(k)}(W(k))$. $m_1$, $\bar g$ and $m_2$ are the reduction mod $p$ of $h_+$, $g$ and respectively $h_0$, while $m_3\in N^-(k)$ is uniquely determined by $h_-$ (it is the reduction mod $p$ of $\prod_{j\in S(\mu)\cap (-\infty,0)} (1_M+u_j)$). 
So we come across entirely the same type of action $\TT_H$ of $N^+_k\times_k N_{0k}\times_k N^-_k$ on $H$ as $\TT$ of 3.13.7.1. 
So the Expectation of 3) of 3.13.7.4 D gets interpreted as:
\medskip
{\bf Expectation (the general form of the CM level one property in the context of Fontaine categories).} {\it We assume $(M,\vph,H_{W(k)})$ has no lift which is cyclic diagonalizable. Then its CM level is at least one.} 
\medskip
We have a variant of this in the adjoint context similar to the one of 3.13.7.1.4: we just have to restate everything in terms of the adjoint representations of $H_{\ZZ_p}$ and of $H_{W(k)}$. We denote by $\TT^0_H$ the group action of $N_{0k}$ on $N^+_k\setminus H_k/\sg(N_k^-)$ defined by $\TT_H$ by passing to quotient. 
\smallskip
Corollary 2 of 3.13.7.5 provides us with many examples when this Expectation holds; in particular it holds if $P_H$ is a Borel subgroup of $H$.
\medskip
{\bf 3.13.7.8.1. Example.} Let $n\in\NN$. We assume that $M$ has rank $n+1$, that $(a,b)=(0,n)$, that $H_{W(k)}=GL(M)$ and that $\dim_{W(k)}(F^s(M))=n+1-s$; so the subgroup $P_{W(k)}$ of $GL(M)$ normalizing $F^s(M)$, $i\in S(1,n)$, is a Borel subgroup. We say ${\got C}$ (resp. any one of its truncations) is a flag $p$-divisible object (resp. a flag object) of rank $n+1$. As $H_{W(k)}=GL(M)$, in what follows we do not mention it. We get (cf. end of 3.13.7.8): 
\medskip
{\bf Corollary.} {\it There is a lift of $(M,\vph)$ such that its truncation mod $p$ is cyclic diagonalizable. The number of isomorphism classes of cyclic diagonalizable flag objects annihilated by $p$ and of rank $n+1$ is equal to the number of elements of the Weyl group of $GL(M)$ and so it is $(n+1)!$.}
\medskip
{\bf 3.13.7.8.2. Remark.} If $P_H$ is as in Corollary 2 of 3.13.7.5, then the whole of 3.13.7.6.0 can be redone in such a context. Denoting by $H_i$, $i\in J$, those simple factors of $H^{\rm ad}$ with the property that the image of $P_H$ in them is a Borel subgroup, the only modifications needed to be made are:
\medskip
{\bf a)} the distinct orbits of $\TT_H$ (or of $\TT^0_H$) are in one-to-one correspondence to the elements of the Weyl group $W_J:=\prod_{i\in J} W_i$ of $\prod_{i\in J} H_i$;
\smallskip
{\bf b)} denoting by $B^+(i)$ the image of $P_H$ in $H_i$ and by $B^-(i)$ its opposite w.r.t. the image $T_J$ of $T$ in $\prod_{i\in J} H_i$, $\forall i\in J$, the orbit corresponding to $w_1\in W_J$ specializes to the orbit corresponding to another $w_2\in W_J$ if and only if $\forall i\in J$, the double coset $B^+(i)w_1(i)B^-(i)$ specializes to the double coset $B^+(i)w_2(i)B^-(i)$ (here we identify $W_J$ with a set of elements of $\prod_{i\in J} H_i(k)$ normalizing $T_J$, while $w_s(i)$ is the component of $w_s$ in $W_i$, $\forall s\in\{1,2\}$ and $\forall i\in J$);
\smallskip
{\bf c)} instead of group of automorphisms (of some Faltings--Shimura--Hasse--Witt maps) we need to speak about stabilizer subgroups of $\TT^0_H$.
\medskip
{\bf 3.13.7.8.3. Example.} We assume $H$ is simple, adjoint. If $H$ is of $G_2$ (or $F_4$ or $E_8$ Lie type), then the number of orbits of $\TT_H$ (or of $\TT^0_H$) is 12 (resp. is $2^73^2$ or $2^{14}3^55^27$), cf. [Bou2, planche VII to IX]. 
\medskip
{\bf 3.13.7.8.4. Remarks.} {\bf 1)} 3.13.7.8.3 is included just to ``compensate" the ``exclusion" of these Lie types of 2.2.7. It nourishes our hope that using Fontaine categories one can handle (at least in the context of unramified extensions of $\ZZ_p$; in fact 2.2.1.4.3 points out that even in the context of ramified extensions of $\ZZ_p$ there is a lot one can hope to do) uniformly all (generic fibres of) reductive groups, from many points of view (like inverse --local-- Galois problems, local Langlands's correspondences, etc.). The Shimura envelopes (hinted at in 3.6.20 9)) are (as well) the very first step in this direction.
\smallskip
{\bf 2)} The whole of 3.13.7.6.1-2 can be redone in the context of 3) of 3.13.7.4: as this needs quite a lot of reformulations, it will not be done here; here we will just mention that in the cyclic context (i.e. when $F^{\rm ad}$ permutes transitively the simple factors of $H^{\rm ad}$), in order to expect to be in a case similar to the ones listed in 3.13.7.6.2 and 3.13.7.7, we need to assume that all images of $P_H$ in a simple factor of $H^{\rm ad}$ which are not reductive, up to isomorphisms of such factors, are isomorphic. 
\smallskip
{\bf 2')} Warning: 3.13.7.6.3 can be adapted to the context of 3) of 3.13.7.4 or of 3.13.7.8 but often slightly modified; the need of modifications is implied by the fact that in Case 2 of the proof of 3.13.7.6.3 we did use (cf. its reference to 3.5.3 (5)) that we are in a generalized Shimura context. So either we try to adapt this Case 2 or we state results in the form: some $isom$-deviations (in the mentioned contexts) are not bigger than some natural numbers.
\smallskip
{\bf 3)} The Expectation of 3.13.7.8 and 3.13.7.8.1 can be restated in terms of Fontaine truncations mod $p$ (see 2.2.14.2) in a similar way to 3.13.7.1.1.1. This supports 1). 
\smallskip
{\bf 4)} Once the Expectation of 3.13.7.8 is checked, one can define the $a$-invariant of $(M,\vph,H_{W(k)})$ by: it is the $a$-invariant of any cyclic diagonalizable $p$-divisible object with a reductive structure whose Fontaine truncation mod $p$ is isomorphic to the one of $(M,\vph,H_{W(k)})$. One checks that this is well defined by performing the proof of M of 2.2.22 3) mod $p$.
\medskip
{\bf 3.13.7.9. The general form of Faltings--Shimura--Hasse--Witt shifts, invariants and maps.} With some extra work, we can be more precise in 3.13.7.8.2 c) (or in connection to its variants provided by cases where the Expectation of 3.13.7.8 is checked). Not to be too long, we situate ourselves in the filtered adjoint context (i.e. we leave to the reader to restate everything as in 3.9.1 or 3.9.4 in the general reductive context, filtered or non-filtered, by looking at derived subgroups). So we use the notations of 3.13.7.8 and we assume that $H$ is an adjoint group and that $\rho$ is the adjoint representation. So $M$ is a Lie algebra over $W(k)$. We still speak about $F^s(M)$, $s\in S(a,b)$ (though in fact we can take $a+b=0$). Let $H_i$, $i\in I$, be the simple factors of $H_{W(k)}$. For $i\in I$, let 
$$m_i\in\ZZ\cap (-\infty,0]$$ be the smallest number such that $H_i\cap N_{m_i}$ is a non-trivial group. Let
$$\Psi:{\rm Lie}(H_{W(k)})\to {\rm Lie}(H_{W(k)})$$
be the $\sg$-linear map that takes $x_i\in {\rm Lie}(H_i)$ into $p^{-m_i}\vph(x_i)$. We refer to it as the Faltings--Shimura--Hasse--Witt quasi-shift of $(M,\vph)$. A great part of 3.13.7.1-8 and of 3.9 can be redone for it. 
\smallskip
First, let $\bar\psi$ be its truncation mod $p$. Second, as in 3.13.7.6 we introduce $\Psi^{\rm nc}$ and $\bar\psi^{\rm nc}$. Third, as in 3.9.1 (resp. 3.9.4) we define the Faltings--Shimura--Hasse--Witt invariant of $(M,\vph)$ (or of $\bar\psi$): it is $\cap_{n\in\NN} {\rm Im}(\bar\psi^n)$. Warning: in the present generality, these invariants are not so useful; for instance, referring to 3.13.7.8.1, there are $(n-1)!$ isomorphism classes for which the Lie stable $p$-rank is maximal (i.e. equal to $1$). Accordingly, the subgroup of $H$ whose $k$-valued points are those $g$ such that $g\bar\psi=\bar\psi$, is in general ``too big" to be used as in 3.13.7.3. So we do not think of $\bar\psi$ as a generalized Faltings--Shimura--Hasse--Witt map. 
\smallskip
To us, the right generalization of the Faltings--Shimura--Hasse--Witt shifts and maps of 3.9 is as follows. For $s\in\ZZ$, let $\vph_s:=p^{-s}\vph$. Let $i\in I$. We denote by 
$$\Psi_{i,m_i}:{\rm Lie}(H_i)\to {\rm Lie}(H)$$ 
the restriction of $\Psi$ to ${\rm Lie}(H_i)$. If $m_i<0$, then for any $j\in S(m_i-1,-1)$ such that $H_i\cap N_j$ is not (resp. is) the origin of $H_i$, let 
$$\Psi_{i,j}:F^j({\rm Lie}(H_{i}))\to {\rm Lie}(\sg(H_{i}))$$
be the restriction of $\vph_j$ to $F^j({\rm Lie}(H_{i}))$ (resp. be the trivial map). The sequence of maps
$$(\Psi_{i,j})_{i\in I,\, j\in S(m_i,{\rm max}\{m_i,-1\})}\leqno (SEQ)$$
(resp. the sequence $(\bar\psi_{i,j})_{i\in I,\, j\in S(m_i,{\rm max}\{m_i,-1\})}$ formed by their truncation mod $p$) is referred as the Faltings--Shimura--Hasse--Witt shift (resp. map) of $(M,(F^s(M))_{s\in S(a,b)},\vph)$ (resp. of $(M,\vph)$). If $(M,\vph)$ is a Shimura adjoint Lie $\sg$-crystal, then for each $i\in I$ only one map $\Psi_{i,j_i}$ is defined; so we can view (SEQ) as a sum 
$$\sum_{i\in I} \Psi_{i,j_i}:{\rm Lie}(H)\to {\rm Lie}(H)$$ 
(resp. as the reduction mod $p$ of this sum): it is the usual Faltings--Shimura--Hasse--Witt shift (resp. map) as defined in 3.4.5 and 3.9. 
\smallskip
In general, 
$$FSHW(M,\vph):=(\bar\psi_{i,j})_{i\in I,\, j\in S(m_i,{\rm max}\{m_i,-1\})}$$ 
does not depend on the choice of $(F^s(M))_{s\in S(a,b)}$. If $h\in H_{W(k)}(W(k))$, then by an inner isomorphism between $FSHW(M,\vph)$ and $FSHW(M,h\vph)$, we mean an isomorphism defined by an element of $N^+N^0(k)$ (cf. also c) and the normalizing part of a) of Step 2 of 3.13.7.3). Following the pattern of 3.13.7.6, we similarly define the non-compact variant 
$$FSHW^{\rm nc}(M,\vph):=(\bar\psi_{i,j}^{\rm nc})_{i\in J,\, j\in S(m_i,{\rm max}\{m_i,1\})},$$
where $J$ is the subset of $I$ formed by elements $i$ such that $H_i$ is not a subgroup of $N_0$. It is easy to see that all of 3.13.7.1.1 b), 3.13.7.1.1.1, 3.13.7.3, 3.13.7.4 B, 3.13.7.6, 3.13.7.8.2 and 3.13.7.8.4 3) can be restated in terms of (isomorphism classes) of $FSHW(M,h\vph)$ or of $FSHW^{\rm nc}(M,h\vph)$. We just add two thinks:
\medskip
-- in 3.13.7.8.2 c), we can replace stabilizers by group of inner automorphisms of Faltings--Shimura--Hasse--Witt (non-compact) maps;
\smallskip
-- for the parts of a) and b) of Step 2 of 3.13.7.3 referring to centralizers, we need to define $C$ as the centralizer of the factors of the lower central series of ${\got n}$.
\medskip
A final thought: there are a couple of ways to define the Lie stable $p$-rank of $(M,\vph)$ (or of $FSHW(M,\vph)$) which can be used to get the analogue of 3.9.4. As we were not able to decide on time which is the most practical such way, we postpone to future work the definition of these invariants.
\medskip
{\bf 3.13.8. Reductive deviations.} We come back to the beginning of 3.13. So $k$ is arbitrary. We saw in Example 1 of 2.2.19 that the Shimura adjoint Lie $\sg$-crystal attached to ${\got C}_0$ is not necessarily of reductive type. The number $1\over n$, with $n\in\NN\cup\{0,\infty\}$ as the smallest number such that the Shimura adjoint Lie $\sg$-crystal attached to $(M,g\vph,G)$ is of reductive type, for $g\in G(W(k))$ congruent to the identity mod $p^n$, is referred as the first reductive deviation of ${\got C}_0$. Here ${1\over\infty}=0$ and ${1\over 0}=\infty$; if $n=\infty$, then $g=1_M$.
\smallskip
We now assume $k=\bar k$. Let ${\got p}_{=0}:=W(0)({\rm Lie}(G),\vph)$ and let ${\got p}_{=00}$ be the $\ZZ_p$-Lie subalgebra of ${\got p}_{=0}$ formed by elements fixed by $\vph$. The length of ${\got p}_{=0}/{\got p}_{=00}\otimes_{\ZZ_p} W(k)$ as a torsion $W(k)$-module is called the second reductive deviation of ${\got C}_0$. 
\medskip
{\bf Example.} The first (resp. second) reductive deviation of the $\sg$-crystal associated to a supersingular elliptic curve over $k=\bar k$ is $0$ (resp. is $1$). 
\medskip
{\bf 3.13.9. Toric deviations.}  We assume $k=\bar k$. Let ${\got p}_{=0}$ and ${\got p}_{=00}$ be as in 3.13.8. By the first (resp. second) toric deviation of ${\got C}$, we mean the rational number 
$${{\dim_{W(k)}(F^0({\rm Lie}(G))\cap {\got p}_{=0})}\over {\dim_{W(k)}({\got p}_{=0})}}$$
 (resp. ${{\dim_{\ZZ_p}(F^0({\rm Lie}(G))\cap {\got p}_{=00})}\over {\dim_{\ZZ_p}({\got p}_{=00})}}$). The first (resp. second) toric deviation of ${\got C}_0$ is defined as the greatest first (resp. greatest second) toric deviation of ${\got C}$, with the $F^1$-filtration $F^1$ of $M$ defining ${\got C}$ running through all elements of $\Mf$ (of the beginning of 3.13).
\smallskip
3.13.8-9 extend naturally to the context of Shimura (filtered) Lie $\sg$-crystals or isocrystals.
\medskip\smallskip
{\bf 3.14. The case $p=2$.} Referring to 3.1-13, we essentially used the fact that $p>2$ in precisely seven places: twice in 3.5.4 (see also 3.6.6.2), in 3.6.2 (and the subsequent places depending on it, including 3.6.10-11, 3.6.18.0 and 3.9.9 A), in 3.6.18.5.1 and 3.6.18.5.3 (and so in 3.6.18.5.7 and 3.13.5.1), in the Fact of 3.6.19, in 3.6.20 8), in 3.11.8.1, and in Step 2 of 3.13.7.3. We now deal, in a convenient order, with the $p=2$ analogues of these seven places.
\medskip
{\bf A.} We start with some small remarks. We repeat, cf. 2.2.1 c), that [Fa1, 2.3] works for $p=2$ as well (see also the $p=2$ analogue of 3.6.18.4.1.1). In 3.6.18.0 (resp. 3.10.8), the $p$-divisible group $D_{W(k)}$ (resp. $D$) does not necessarily exist, and if exists, it is not always uniquely determined for $p=2$ (resp. it is not always uniquely determined for $p=2$: its existence is argued in 2.2.12.1 1)). Accordingly, for the existence part in 3.6.18.0, for $p=2$, we usually need to assume that $k$ is $1$-simply connected (i.e. it has no abelian extension of degree $2$) or that $(M,\vph)$ either does not have slope $0$ or does not have slope $1$, cf. 2.3.18.1 B and C. So 3.6.18.5 c) does not need to be modified. In the paragraph before 3.4.0, for $p=2$ the list of cases when the Lie monomorphism $i_G$ is not an isomorphism is much longer; in fact, in most situations, $i_G$ is not an isomorphism. For $p=2$, a) of 3.6.18.10 P1 has to be formulated as follows: ``is associated to a $2$-divisible group over $R_1$ or over $R_1^\wedge$". The exponential map (in the generalized Shimura context) referred to in 3.6.6.2 for $p=2$ does not have the same simplified form as in the proof of 3.6.6 (i.e. $1_M+x$ mod $4$ has to be replaced by a sum $1_M+\sum_{s=1}^{n} {x^s\over {s!}}$, where $n\in\NN$ is such that $x^n=0$, $\forall x\in 2{\rm Lie}(N_{-i})$; however, its usage does not need to be modified).
\medskip
{\bf B.} The equivalence parts of 3.6.18.5.1 and 3.6.18.5.3 remain true for $p=2$ but not their antiequivalence parts. However, due to the fact that [Fa2, th. 10] is true as well for $p=2$, the mentioned antiequivalence parts can be reformulated in weaker forms as follows.
\medskip
{\bf B0. Exercise.} Any filtered $\sg$-crystal over $k$ is associated to a $2$-divisible group iff $k$ is $1$-simply connected. Hint: one implication is a consequence of 2.3.18.1 B and D; for the converse use the part of the proof of 2.3.18.1 B referring to [Og].
\medskip
{\bf B1. Proposition.} {\it We refer to 3.6.18.5.1 and 3.6.18.5.3 with $p=2$. We assume $k$ is $1$-simply connected. We have: 
\medskip
{\bf A)} the quotient of $2-DG({\rm Spec}(R))$ ``under the functor $\DD$" (of 2.2.1.0) is antiequivalent to the category $2-\Mm\Mf_{[0,1]}(R)$ and so to the category $2-\Mm\Mf_{[0,1]}^\nabla(R)$;
\smallskip
{\bf B)} any object of $\Mm\Mf_{[0,1]}^\nabla(R)$ is associated (via $\DD$) to a finite, flat, commutative group scheme of $2$-power order over $R$.} 
\medskip
{\bf Proof:} As any object of $\Mm\Mf_{[0,1]}(R)$ is the cokernel of an isogeny between $2$-divisible objects of $\Mm\Mf_{[0,1]}(R)$ (see Fact of 2.2.1.1 6)) and as the equivalence parts of 3.6.18.5.1 and 3.6.18.5.3 hold for $p=2$, B) follows from A). To see A) it is enough to show: 
\medskip
{\bf B2. Lemma.} {\it Any morphism $n_{12}:{\got C}_1\to {\got C}_2$ between two $2$-divisible objects of $\Mm\Mf_{[0,1]}^\nabla(R)$ is associated (via $\DD$) to a morphism $m_{12}:D_2\to D_1$ between two $2$-divisible groups over $R$.}
\medskip
We first consider the case $R=W(k)$. Let (cf. our assumption on $k$ and B0) $D_i$ be a $2$-divisible group over $W(k)$ such that $\DD(D_i)={\got C}_i=(M_i,F_i^1,\vph_i)$, $i=\overline{1,2}$. We write $D_i$ as the extension of an \'etale $2$-divisible group $E_i$ by a $2$-divisible group $F_i$ not having slope $0$. $E_i$ and $F_i$ are uniquely determined by ${\got C}_i$ (cf. 2.3.18.1 C and the fact that the complex 
$$0\to\DD(E_i)\to\DD(D_i)\to\DD(F_i)\to 0$$ 
is a short exact sequence). So $n_{12}$ determines uniquely morphisms $e_{12}:E_2\to E_1$ and $f_{12}:F_2\to F_1$. We need to show: we can choose $D_1$ and $D_2$ such that $e_{12}$ and $f_{12}$ are determined naturally by a morphism $m_{12}:D_2\to D_1$ which satisfies $\DD(m_{12})=n_{12}$. Let $D_3:=D_1\times_{E_1} E_2$; it is a $2$-divisible group over $k$. Let $(M_3,F_3^1,\vph_3)$ be its associated filtered $\sg$-crystal. 
\smallskip
Based on 2.3.18.1 D, we can work over $W_2(k)$. What follows is just a functorial version of 2.3.18.1 B (cf. also 2.3.18.1.1). We consider the affine scheme  $FP_{\rm def}$ defined naturally by the $k$-vector space defined by the fibre product of the natural $k$-linear maps 
$$
{\rm Ext}^1(E_{1W_2(k)},F_{1W_2(k)})\to {\rm Ext}^1(E_{2W_2(k)},F_{1W_2(k)})
$$ 
and 
$$
{\rm Ext}^1(E_{2W_2(k)},F_{2W_2(k)})\to {\rm Ext}^1(E_{2W_2(k)},F_{1W_2(k)}).
$$ 
Similarly we consider the affine scheme $FP_{\rm fil}$ defined naturally by the $k$-vector space of the similarly constructed fibre product obtained working in the crystalline context (i.e. using lifts mod $4$ of the filtrations $F^1_1$, $F^1_2$ and $F^1_3$ mod $2$, etc.). The natural morphism
$$FP_{\rm def}\to FP_{\rm fil}$$
is a Galois cover (cf. 2.3.18.1.2). So the resulting map 
$$FP_{\rm def}(k)\to FP_{\rm fil}(k)$$ 
is surjective. $n_{12}$ defines naturally an element 
$$\gamma_{n_{12}}\in FP_{\rm fil}(k).$$ 
We just need to choose $D_1$ and $D_2$ to be defined naturally by some element 
$$\gamma_{m_{12}}\in FP_{\rm def}(k)\subset {\rm Ext}^1(E_{1W_2(k)},F_{1W_2(k)})\times {\rm Ext}^1(E_{2W_2(k)},F_{2W_2(k)})$$ 
mapping into $\gamma_{n_{12}}$.   
\smallskip
We come back to the general case: so $m\in\NN\cup\{0\}$ (of 3.6.18) is arbitrary. We consider the Teichm\"uller lift $z:{\rm Spec}(W(k))\hookrightarrow {\rm Spec}(R)$; let $I_R$ be the ideal of $R$ defining $z$. We consider a morphism $m_{12}^z:D_2^z\to D_1^z$ between two $2$-divisible groups over $W(k)$ such that $\DD(m_{12}^z)$ is identifiable with $z^*(n_{12})$. From [Fa2, th. 10] we deduce the existence of a uniquely determined $2$-divisible group $D_i$ over $R$ such that $z^*(D_i)=D_i^z$ and $\DD(D_i)$ is identifiable with ${\got C}_i$ in a way compatible with the identification $\DD(D_i^z)=z^*({\got C}_i)$, $i=\overline{1,2}$. There is a uniquely determined morphism $m_{12}:D_2\to D_1$ such that $\DD(m_{12})$ is identifiable with $n_{12}$ in a way compatible with all previous identifications (this is just the endomorphism aspect of the uniqueness argument of 2.2.21 presented in connection to [BM, ch. 4] and [Me, ch. 4-5]). This proves the Lemma and so the Proposition.   
\medskip
{\bf B3. Remark.} We still assume that we are in the context of B1 above. The above proof can be entirely adapted to show that any commutative diagram in the category $2-\Mm\Mf_{[0,1]}^\nabla(R)$ is associated via the $\DD$ functor to a commutative digram in $2-DG({\rm Spec}(R))$. But any morphism between two objects of $\Mm\Mf_{[0,1]}^\nabla(R)$ can be written as a natural morphism (defined by $f_3$ and $f_4$) ${\rm Coker}(f_1)\to {\rm Coker}(f_2)$, where $f_1$, $f_2$, $f_3$ and $f_4$ are morphisms between objects of $2-\Mm\Mf_{[0,1]}^\nabla(R)$ such that $f_3\circ f_1=f_2\circ f_4$ (cf. 2.2.1.1 6) and the first sentence of the above proof). So one gets the following improvement of B) of B1:
\medskip
{\bf Corollary.} {\it The quotient of $2-FF({\rm Spec}(R))$ ``under $\DD$" is $\Mm\Mf_{[0,1]}^\nabla(R)$.}
\medskip
{\bf B4. Warning.} In A) of B1 above (resp. in the Corollary of B3) some $2$-divisible groups are ``forced" to be isomorphic (resp. some finite, flat, commutative group schemes of $2$-power order are forced to be isomorphic, while some morphisms between them are ``forced" to be zero).
\medskip
{\bf B5.} We come back to an arbitrary $k$. The variant 3.6.18.5.7 can be adapted for $p=2$, cf. 2.3.18.1 C: we just have to work either just with trivial \'etale parts or just with trivial multiplicative type parts (i.e. we get two equivalences and two antiequivalences, by working with the corresponding full subcategories of the categories of 3.6.18.5.7). One sample: The full subcategory of $p-FF({\rm Spec}(R))$ whose objects have trivial \'etale parts is antiequivalent to the full subcategory of $p-\Mm\Mf_{[0,1]}(R)$ whose objects are such that the Newton polygons of their pull backs to $\sg$-crystals over $k$ have no zero slopes.
\medskip
{\bf B6. Corollary.} {\it We refer to 3.6.18.5.1 with $p\ge 2$. The category $p-DG({\rm Spec}(R/pR))$ is antiequivalent (via $\DD)$ to the category $p-\Mm_{[0,1]}(R)$ as defined in 2.2.1.6.}
\medskip
{\bf Proof:} See [BM, 4.2.4] for the fully faithfulness part. We consider an object $(M_R,\Phi_{M_R},\nabla_{M_R})$ of $p-\Mm_{[0,1]}(R)$. We consider a direct summand $F^1(M_R)$ of $M_R$ such that ${\got C}:=(M_R,F^1(M_R),\Phi_{M_R},\nabla_{M_R})$ is an object of $p-\Mm_{[0,1]}(R)$. If $p=2$, we need to assume that $F^1(M_R)$ is such that the pull back of $(M_R,F^1(M_R),\Phi_{M_R})$ via the $W(k)$-valued Teichm\"uller lift of ${\rm Spec}(R)$, is associated to a $2$-divisible group over $W(k)$. From 2.2.21 UP (resp. from 3.6.18.5.3) for $p\ge 2$ (resp. for $p\ge 3$), we get the essential surjectivity part of $\DD$ on objects. This ends the proof.
\medskip
We refer to 3.6.18.5.1 with $p\ge 2$. We consider an ideal $I_R$ of $R$ taken by $\Phi_R$ into itself and such that $S:=R/I_R$ is flat over $W(k)$. Let $\Phi_S$ be the resulting Frobenius lift of $S$. 
We recall (cf. [BM]) that to any $p$-divisible group $D$ over ${\rm Spec}(S/pS)$ it is functorially associated a crystal $\DD(D)$ on the crystalline site $CRIS(S/pS/{\rm Spec}(\ZZ_p))$. So as the ideal $pS$ of $S$ is naturally equipped with a divided power structure and due to the existence of $\Phi_S$, evaluating $\DD(D)$ at ${\rm Spec}(S)$ we get a triple $(M_D,\Phi_D,\nabla_D)$, with $M_D$ a free $S$-module of finite rank, with $\nabla_D$ an integrable, nilpotent mod $p$ connection on $M_D$ and with $\Phi_D$ a $\Phi_S$-linear endomorphism of $M_D$ which is $\nabla_D$-parallel. As $D$ lifts to a $p$-divisible group over $S$, this triple is an object of the category $p-\Mm^\nabla(S)$ defined in 2.2.1.7 4). So we get a contravariant functor 
$$\DD_{S/pS}:p-DG({\rm Spec}(S/pS))\to p-\Mm^\nabla(S).$$ 
\medskip
{\bf B7. Corollary.} {\it  $\DD_{S/pS}$ is onto on objects.}
\medskip
{\bf Proof:} Let $({\got C}_{S},\nabla)$ be an object of $p-\Mm_{[0,1]}^\nabla(S)$. As in Fact of 2.2.1.1 6), ${\got C}_{S}$ is the reduction mod $I_R$ of an object ${\got C}$ of $p-\Mm_{[0,1]}(R)$. $\nabla$ is the unique connection on the underlying $S$-module $M_S$ of ${\got C}_S$ such that the pair $({\got C}_{S},\nabla)$ is an object of $p-\Mm^\nabla(S)$. To see this, let $\nabla_1$ be a second such connection. To show that $\nabla-\nabla_1\in {\rm End}(M_S)\otimes_S\Om_{S/W(k)}^\wedge$ is zero, it is enough to show that $\nabla$ and $\nabla_1$ coincide mod $p$. Using the notations of 3.6.18.1.2, we consider the sequence of ideals $(J_n)_{n\in\NN\cup\{0\}}$, with $J_n:=I_n+I_R$ of $R$; by induction on $n\in\NN$ we get easily that $\nabla-\nabla_1\in J_n{\rm End}(M_S)\otimes_S\Om_{S/W(k)}^\wedge$. As $\cap_{n\in\NN} J_n=I_R$, we get that $\nabla$ and $\nabla_1$ coincide mod $p$. So $\nabla=\nabla_1$.
\smallskip
So $({\got C}_{S},\nabla)=\DD_{S/pS}(D_{{\rm Spec}(S/pS)})$, where (cf. B6) $D$ is the $p$-divisible group over $R/pR$ corresponding to ${\got C}$. This ends the proof.
\medskip
{\bf B8. Remark.} B6-7 were first obtained (using a slightly different language) in [dJ1, th. of intro.]. B7 can be combined with [BM] to reobtain the full form of [dJ1, th. of intro.]. Moreover, we have variants of B7 where $S$ is not necessarily flat.
\medskip
{\bf C.} It is easy to see that in 3.5.4 it was not at all essential that $p>2$. To see this we use the notations of 3.5.4, with $p=2$. The formula (6) of 3.5.4 remains true for $p=2$ (it can be checked inside the $SL$-group of a 3 dimensional vector space over $k$) but formula 3.5.4 a) does not. In fact it is easy to see that it fails precisely in the following cases:
\medskip
\item{(*)} {\it $\bar G_1$ is an adjoint group of $B_{\ell}$ Lie type, with $\ell\ge 1$, and $\al$ is a short root.}
\medskip
As (*) is obvious if $\bar G_1$ is of $B_1$ Lie type, to argue (*) we can assume the rank $\ell$ of $\bar G_1$ is at least 2. Using the fact that all roots of $\Phi_1$ of the same length are conjugate under the Weyl group attached to it (for instance, see [Hu1, p. 53]), for the $A_{\ell}$ ($\ell\ge 2$), $D_{\ell}$ ($\ell\ge 4$), $E_6$ and $E_7$ Lie types, the situation gets reduced to the $A_2$ Lie type, for which 3.5.4 a) obviously holds. Similarly, for the $B_{\ell}$ (resp. $C_{\ell}$) Lie type case, with $\ell\ge 2$ (resp. $\ell\ge 3$), the situation gets reduced to the $B_2$ Lie type (resp. to the simply connected $C_2$ Lie type): we just need to remark the following two things (valid for $\ell\ge 2$).
\medskip
{\bf a)} The standard monomorphism $SO(2\ell+1)\hookrightarrow SO(2\ell+3)$ (over any algebraically closed field) lifts to a monomorphism ${\rm Spin}(2\ell+1)\hookrightarrow {\rm Spin}(2\ell+3)$ between simply connected semisimple groups. Here, for $i\in\{1,3\}$, $SO(2\ell+i)$ denotes a simple adjoint group of $B_{\ell+j}$ Lie type, where $j=0$ if $i=1$ and $j=1$ if $i=3$. 
\smallskip
{\bf b)} We have a similar standard monomorphism from the simply connected semisimple group over $k$ of $C_{\ell-1}$ Lie type into the adjoint semisimple group over $k$ of $C_{\ell}$ Lie type.
\medskip
The remaining cases, can be studied (the argument of Step 2 of 3.13.7.3 referring to [Va2, 3.1.2.1 c)], allows us to shift from positive characteristic to characteristic $0$) in characteristic $0$: the fact that 3.5.4 a) holds or not for $\bar G_1$ can be restated in terms of some semisimple group of $B_1$ Lie type are or are not all simply connected). But the situation pertaining to the $B_2$ Lie type (adjoint or not, can be read out from [Bou3, (II) of p. 201]). This ends the argument for (*).
\medskip  
On the other hand [BT, 4.2-3] and [Bo2, 3.16] treat the case $p=2$ as well and so formula (7) still holds for $p=2$. So 3.5.4 holds for $p=2$ provided we are in the context of Shimura $\sg$-crystals or in a generalized Shimura context not involving adjoint groups of $B_{\ell}$ Lie type. To see that 3.5.4 holds as well for the generalized Shimura context involving groups whose derived subgroups have normal, simple subgroups which are adjoint and of $B_{\ell}$ Lie type, we can proceed in many ways, like:
\medskip
i) we shift from an adjoint context to a simply connected context;
\smallskip
ii) we replace the argument based on 3.5.4 a) and made at the level of Lie algebras in characteristic $p$, either by a refined one or by an argument at the level of reductive groups in characteristic $0$;
\smallskip
iii) we refer to J below (for $\ell\ge 2$). 
\medskip
For the convenience of the reader we include here one way to argue things via ii). In the adjoint $B_{\ell}$ Lie type context, the Fact of 3.5.4 has to be restated as follows:
\medskip
{\bf Fact.} {\it The Lie subalgebra $\bigl(g_2(-\dl)\cap{\got g}_1\bigr)\otimes_{W(k)} k$ of $\bar{\got g}_1$ is 
$\bar{\got s}_a$ and  $\bigl(\bigl(\bigoplus_{\al\in\Mh_2\cap[0,1]}g_2(\al)\bigr)\cap{\got g}_1\bigr)\otimes_{W(k)} k$ is $\bar P_{\vep_1}$.}
\medskip
{\bf Proof:} We keep using the notations of 3.5.3, with $p=2$. As we are dealing with the $B_{\ell}$ Lie type, $\abs{A}=1$. For the rest of the argument we just use this fact. By reasons of dimensions we have (cf. the Fact of 3.4.5.1) 
$$\bigl(g_2(-\dl)\cap{\got g}_1\bigr)\otimes_{W(k)} k=\bar{\got s}_{\vep_1}=F^1({\got g}_1)\otimes_{W(k)} k.\leqno (1)$$ 
We consider the parabolic subgroup $\bar P_{\ge 0}^1$ of $\bar G_1$ whose Lie algebra is $\bigl(\bigl(\bigoplus_{\al\in\Mh_2\cap[0,1]}g_2(\al)\bigr)\cap{\got g}_1\bigr)\otimes_{W(k)} k$. Using (1) and Fact of 2.2.11.1 we get that $\bar P_{\ge 0}^1$ is a subgroup of $\bar P_{\vep_1}$. So, as in the proof of 3.4.8, replacing $\vph_2$ by $\tilde g_2\vph_2\tilde g_2$, with $\tilde g_2\in P_{\vep_1}(W(k))$, we can assume 
$$\bigl(\bigl(\bigoplus_{\al\in\Mh_2\cap[0,1]}g_2(\al)\bigr)\cap{\got g}_1\bigr)\subset {\got p}_{\vep_1}.\leqno (INCL)$$ 
Performing the same thing for every $i\in I_1$, we get that we can assume $\bigoplus_{\al\in\Mh_2\cap[0,1]}g_2(\al)\subset F^0({\got g}_0)$. This together with (1) and its analogues with $i\in I_1$, implies that 
$$(\bigoplus_{\al\in\Mh_2\cap[0,1]}g_2(\al)/g_2(-\dl),\bigoplus_{\al\in\Mh_2\cap[0,1]}g_2(\al)/g_2(-\dl),\vph_2)$$
is a $p$-divisible object of $\Mm\Mf_{[0,0]}(W(k))$. 
\smallskip
So $({\got g}_0,\vph_2)$ has no slopes in the interval $(0,-\dl)$ and so $\bigl(g_2(-\dl)\cap{\got g}_1\bigr)\otimes_{W(k)} k$ is the Lie algebra of the unipotent radical of $\bar P_{\ge 0}^1$. So from (1) and the Fact of 2.2.11.1 we get that in (INCL) we have in fact equality. From this the Fact follows
\medskip
{\bf D.} In 3.6.2 it was crucial that $p>2$. So everything in 3.6.1-18, which is not based on 3.6.2 or on the antiequivalence parts of 3.6.18.5.1 and 3.6.18.5.3, remains true for the case $p=2$. In other words, everything in 3.6.1-18 except 3.6.2, the antiequivalence parts of 3.6.18.5.1 and 3.6.18.5.3, 3.6.10-11 and the whole of 3.6.14, remains true for the case $p=2$; in connection to 3.6.18.5.7 see B5. In 3.13.5.1 iv), for $p=2$ we need to replace $DG$ by its filtered $F$-crystal.
\medskip
{\bf E.} Despite D above, it is  possible to modify the statement and the proofs of 3.6.14.1-4 accordingly (i.e. without relying on 3.6.2). We have a couple of variants. 
\smallskip
First, if the Shimura-ordinary type associated to the $p=2$ SHS $(f,L_{(2)},v)$ has no integral slopes, then we do not need any modification to 3.6.2 and so to the whole of 3.6.14 (cf. D and 2.3.18.1 C and D; in such a context the part of the proof of 3.6.14.1 depending on 2.2.1.1 2) still holds for $p=2$ and so the part of the proof of 3.6.14.1 referring to [BLR] still applies). Warning: this has limited applications; however $p=2$ variants of 4.6.1 1) below (for instance involving the $A_1$ Lie type, cf. 2.3.18 A) point out that occasionally it does apply. 
\smallskip
Second, we replace assumptions SA1-2 by the restatement of property G) of the proof of 3.6.14.1, where we do not mention anything about $2$-divisible groups; so we work entirely in terms of (non-filtered) $F$-crystals with tensors. So 3.6.14 and 3.6.14.1, under the logical (abstract) restatement involving just $F$-crystals, are true as well for a $p=2$ SHS (in terms of this restatement, the mentioned part referring to [BLR] can be skipped entirely). 
\smallskip
Third, the proof of Theorem 2 of 3.15.1 below points out that we can still work in the context of $2$-divisible groups over $k$-schemes (see also 3.15.2 below). In other words, in the second variant, we can replace $F$-crystals by $p$-divisible groups (over $k$-schemes).
\smallskip
Fourth, we have a variant in which we already assume that the things ``are fine" mod $4$, i.e. (besides S1 or S2) we assume that a scheme ${\rm Spec}(T_2)$ as in the proof of 3.6.14.1 does exist a priori. 
\smallskip
Similarly, 3.6.10-11 have to be stated just in terms of filtered $F$-crystals and not of $2$-divisible groups. The fact that 3.9.9 A, B and C still holds for $p=2$ is implied by the proof of Theorem 2 of 3.15.1 below: we just need to be in a context of smooth group schemes over $W(k)$ so that the part of the mentioned proof pertaining to lifts of cocharacters still applies (cf. also 3.15.3 5) below).
\medskip
{\bf F.} In 3.6.20 8), for the versal (resp. uni plus versal) context we need to assume $k$ is $1$-simply connected (resp. that the $2$-divisible groups over $k$ defined by special fibres of the ones considered in 3.6.19 i), do not have either slope $0$ or slope $1$), cf. 2.3.18.1 B and C.
\medskip
{\bf G.} We refer to 3.6.19. a) of its Fact remains true. Similarly, its Theorem remains true: it is stated already in terms of $p$-divisible groups and so (cf. also the extra assumption (EXTRA) of 3.6.19) for $p=2$ we do not get into trouble with D or E above. Similarly, its variants iv) to vi) mentioned in 3.6.19 E still remain true. In connection to its variant vii) we refer to 4.14.3 K and to [Va5]. 
\medskip
{\bf H.} In connection to 3.11.8.1 for $p=2$, it is 2.13.8.2 which points out that we can still define $\rho_{\al}$'s. Related to $\rho$ we need to assume that either $(M,\vph)$ does not have slope $0$ or $1$ (cf. 2.3.18.1 C) or $\rho$ is defined using a $2$-divisible group over $W(k)$ as in I below. 
\medskip
{\bf I.} Let $(M,F^1,\vph,G,(t_{\al})_{\al\in\Mj})$ be a Shimura-canonical lift of a (non-necessarily quasi-split) Shimura $\sg$-crystal with an emphasized family of tensors over a perfect field $k$ of characteristic $2$. Let $D_k$ be the $2$-divisible group over $k$ having $(M,\vph)$ as its associated $\sg$-crystal. As pointed out in 2.3.18.1, it can happen that there is more than one $2$-divisible group over $W(k)$ lifting $D_k$ and such that its associated filtered $\sg$-crystal is $(M,F^1,\vph)$. However, from 3.11.1 a) (and 2.3.18.1 C) we get that there is a unique such lift $D$ which is a direct sum of: an \'etale $2$-divisible group, of a multiplicative type $2$-divisible group and of a $2$-divisible group whose slopes are rational numbers of the interval $(0,1)$. We refer to the pair $(D,(t_{\al})_{\al\in\Mj})$ (resp. $(D_k,(t_{\al})_{\al\in\Mj})$) as a Shimura-canonical (resp. Shimura-ordinary) $2$-divisible group over $W(k)$ (resp. over $k$).  
\medskip
{\bf J.} To show that 3.13.7.3 (and so implicitly all of 3.13.7.1-2 and of 3.13.7.4-8) still holds for $p=2$, we just need to point out that b) of Step 2 of 3.13.7.3 applies in Step 3 of 3.13.7.3 in the same manner, as we are interested in $k$-valued points of $G^{\rm ad}_k$ and so we can work equally well with ${P_C}_{\rm red}$ instead of $P_{C}$. However, referring to b) of Step 2 of 3.13.7.3, we have the following improvement of it (for the Shimura adjoint Lie $\sg$-context):
\medskip
{\bf Supplement.} {\it We assume $p=2$ and that there is a Shimura pair $(H,[\mu])$ such that $P$ is the parabolic subgroup of $H$ whose Lie algebra is the maximal Lie subalgebra of ${\rm Lie}(H)$ on which $\mu$ acts via the identity and the trivial cocharacter of $\GG_m$. If moreover $H$ is not of $A_1$ Lie type, then ${P_C}_{\rm red}=P_C$ and so $C=N$.}
\medskip
{\bf Proof:} Let $T$ be a maximal torus of $P$. Let $\Phi$ be the set of roots of the action $AC$ of $T$ (via inner conjugation) on ${\rm Lie}(H)$. Let $\Phi_0$ (resp. $\Phi_1$) be the subset of $\Phi$ formed by the roots of the action of $T$ on ${\rm Lie}(P)$ (resp. on ${\got n}$). As $H$ is adjoint, $AC$ is faithful. Based on the proof of b) of Step 2 of 3.13.7.3, we just need to show ${P_C}_{\rm red}=P_C$. So, following the pattern of the part of the Claim of 3.5.4 referring to $\abs{A}=1$ and involving normalization, we just need to show that:
\medskip
{\bf Fact.}  {\it For any $\al\in\Phi\setminus\Phi_0$, there is $\be\in\Phi_1$, such that $\al+\be\in \Phi_0$ and the pair $(\al,\be)$ is not irregular in the sense of [BT, 4.3] (so the $p=2$ analogue of 3.5.4 (6) holds for it).}
\medskip
This Fact can be easily checked (as in Case 1 of 3.5.4) starting from [Bou2, planche I to VI] and [BT, 4.3 (i)]. Based on the list of irregular cases of loc. cit., the only situations which need extra details are the ones when $H$ is of $B_l$ or $C_l$ Lie type, $l\in\NN$, with $l\ge 2$. Choosing a Borel subgroup of $P$ containing $T$, we get an ordering of $\Phi$ such that the set $\Phi^+$ of positive roots is included in $\Phi_0$. Based on loc. cit., we can assume $\al$ is short. For future references we assume just that $\al\in\Phi\setminus\Phi_1$. We first deal with the $B_l$ Lie type case. With the standard notations of [Bou2, planche II], if $\al=-\sum_{i=m}^l \al_i$, with $m\ge 1$, then we can take $\be=(\sum_{i=1}^{m-1}\al_i)+2(\sum_{i=m}^l \al_i)$ if $m\ge 2$ and the maximal positive root if $m=1$. If $\al=\sum_{i=m}^l \al_i$, with $m\ge 2$, then we can take $\be=\sum_{i=1}^{m-1}\al_i$. In both situations, $\be$ is long, and we are not in an irregular situation. 
\smallskip
We now assume $H$ is of $C_l$ Lie type. If $\al=-\sum_{i=m}^j \al_i$, with $1\le m\le j<l-1$, then we can take $\be=\al_l+2\sum_{i=m}^{l-1} \al_i$. We now assume $\al=-\sum_{i=m}^l \al_i-\sum_{i=n}^{l-1} \al_i$, with $1\le m\le n\le l-1$. If $n\neq 1$, then we can take $\be$ to be the maximal root. If $n=1$, then we take $\be=\sum_{i=1}^l\al_i+\sum_{i=j}^{l-1} \al_i$, where $j\in S(1,l)$ is such that $\al$ and $\be$ are not orthogonal. In this last case, $\be$ is short but is not perpendicular on $\al$; so we are through (cf. [BT, 4.3 (i)]). This proves the Fact and ends the proof of the Supplement.
\medskip\smallskip
{\bf 3.15. Some conclusions.} This section is formed by different conclusions to 3.1-14; they are thought as refined versions of parts of 3.1-14 which can be obtained by thoroughly combining different parts of it. In 3.15.1-6 and 3.15.9 we deal with Dieudonn\'e theories, versal deformations and applications of them. In 3.15.7-8 and 3.15.10 we exploit to a much greater extend 3.6.15 B. In particular, we get:
\medskip
-- a new principle (the boudedness one, see 3.15.7);
\smallskip
-- a new proof (see 3.15.8) of the specialization theorem;
\smallskip
-- a new proof in a slightly more general context (see 3.15.10.1) of [dJO, 4.1].
\medskip
{\bf 3.15.1. Versal global deformations.} 3.6.1.3 allows us (cf. 3.6.9 3)) to strengthen 3.1.8 and its $p=2$ analogue of 3.1.8.1. Let $p\ge 2$ be an arbitrary prime and let $\tilde k$ be an arbitrary perfect field of characteristic $p$. We have:
\medskip
{\bf Theorem 1.}  {\it Any $p$-divisible group over $\tilde k$ admits a versal (global) deformation over the $p$-adic completion of an $\NN$-pro-\'etale, affine scheme ${\rm Spec}(R_1)$ over a smooth, affine $W(\tilde k)$-scheme ${\rm Spec}(R)$ (and so also over an $\NN$-pro-\'etale, affine scheme ${\rm Spec}(\bar R_1)$ over a smooth, affine $\tilde k$-scheme ${\rm Spec}(\bar R)$), of whose Kodaira--Spencer map is an isomorphism. We can choose $R_1$ such that ${\rm Spec}(R_1/pR_1)$ is a geometrically connected, $AG$ $k$-scheme and the $\tilde k$-morphism ${\rm Spec}(R_1/pR_1)\to {\rm Spec}(R/pR)$ is surjective. Moreover the $p$-divisible group over the generic point of ${\rm Spec}(R_1/pR_1)$ (or of $\bar R$) is ordinary.} 
\medskip
Similarly, in the relative context, we have:
\medskip
{\bf Theorem 2.} {\it Let $\Md(W(\tilde k))$ be a Shimura $p$-divisible group over $W(\tilde k)$. Then there is a smooth $W(\tilde k)$-algebra $R$, an $\NN$-pro-\'etale, affine $R$-scheme ${\rm Spec}(R_1)$, and a Shimura $p$-divisible group $\Md$ over ${\rm Spec}(R_1^\wedge)$ such that:
\smallskip
a) in a $W(\tilde k)$-valued point $z$ of ${\rm Spec}(R_1^\wedge)$, $\Md$ becomes (isomorphic to) $\Md(W(\tilde k))$;
\smallskip
b) $\Md$ is a uni plus versal (global) deformation (see 3.6.19 B); 
\smallskip
c) ${\rm Spec}(R_1/pR_1)$ is a geometrically connected, $AG$ $k$-scheme and the morphism ${\rm Spec}(R_1/pR_1)\to {\rm Spec}(R/pR)$ is surjective.
\smallskip
Moreover, the resulting Shimura $\sg_{\tilde k_g}$-crystal ${\got C}$ over the algebraic closure $\tilde k_g$ of the field of fractions of $R_1/pR_1$ is Shimura-ordinary.}
\medskip
{\bf Proof:} The parts involving (Shimura-) ordinariness of the above Facts are a consequence of 3.12.1 and its proof (cf. also 3.14 for $p=2$). The case $p\ge 3$ is covered, by taking slices, by 3.6.1.3, 3.6.2 and 3.12.1: we just have to perform 3.6.11 in the context when $\tilde H$ is as in the first paragraph after 3.6.14. The fact that the above Theorems hold as well for $p=2$, can be checked as follows. 
\smallskip
We can refer just to Theorem 2. Copying 3.6 i), parts of 3.6.0-1, and the case $n=1$ of 3.6.1.3, we need to start with a smooth $W(\tilde k)$-algebra $R$ such that:
\medskip
{\bf i)} it has a connected special fibre;
\smallskip
{\bf ii)} its spectrum has a $W(\tilde k)$-valued point $z$ and there is a triple $({\rm Spec}(R),b_R,z)$ which is a potential-deformation sheet in the sense of 3.6.9.1 (so the $R/pR$-module $\Om_{R/pR/\tilde k}$ is free);
\smallskip
{\bf iii)} its $2$-adic completion is equipped with the Frobenius lift $\Phi_R$ obtained via $b_R$ in the same way $\Phi_U$ of 3.6.9.1 is obtained (via $b_U$); so (in $z$) we have a property similar to the one of 2.3.15 c): in particular, $z$ is a Teichm\"uller lift;
\smallskip
{\bf iv)} there is a $2$-divisible object ${\got C}$ with cycles of $\Mm\Mf_{[0,1]}(R)$, which in $z$ becomes (through pull back) the $2$-divisible object ${\got C}_z$ with cycles of $\Mm\Mf_{[0,1]}(W(\tilde k))$ defined by $\Md(W(\tilde k))$, which is modeled (as in 2.2.10 or 3.6.0-1) on ${\got C}_z$ and whose truncation mod $2$ is equipped with a uni plus versal connection, respecting the extra ``Shimura structure" (so the relative dimension of $R$ over $W(\tilde k)$ is $dd({\got C})$);
\smallskip
{\bf v)} the underlying $R^\wedge$-module of ${\got C}$ together with its cycles, is obtained by extension of scalars from a pair $(\tilde M,(\tilde t_{\al})_{\al\in\Mj})$ over $W(\tilde k)$; so it is free, and the referred ``Shimura structure" is defined by a reductive subgroup $\tilde G$ of $GL(\tilde M)$ (so the connection of iv) respects the $\tilde G_R$-action in the sense of 3.6.1.1.1 2)).  
\medskip
In connection to iv) we add: working with $\tilde H$ as mentioned, we always can assume that the part on uni plus versality holds (perhaps after passage to an open, affine subscheme of ${\rm Spec}(R)$ through which $z$ still factors), cf. c) of 3.6.18.7.3 C. 
\smallskip 
So we can ``perform 3.6.1.3" for ${\got C}$ (for $n\ge 2$), cf. 3.14 and 3.6.18.4.2. We get: there is an $\NN$-pro-\'etale morphism ${\rm Spec}(R_1)\to {\rm Spec}(R)$ to which $z$ lifts uniquely (we still denote by $z$ this unique lift), such that $R_1/2R_1$ is an integral ring and (cf. also 3.6.18.8.1 and iii) above) the pull back ${\got C}_{R_1}$ of ${\got C}$ to ${\rm Spec}(R_1)$ is a uni plus versal $2$-divisible object with tensors of $\Mm\Mf_{[0,1]}^\nabla(R_1)$. Let ${\got C}_1$ be the $F$-crystal over $R_1/2R_1$ obtained from ${\got C}_{R_1}$ and its natural connection by forgetting its filtration.
\smallskip
Let $\Mk_1$ be the field of fractions of $R_1/2R_1$. It is a field which is an $\NN$-ind-\'etale extension of a finitely generated field over $\tilde k$ and so (see [BM, 1.1.2 (ii)]) it has a finite $2$-basis; $R_1/2R_1$ itself has a finite basis (cf. loc. cit. and ii)) but this is irrelevant for what follows. So the first main result of [dJ2] applies: there is a unique $2$-divisible group $\Mg_{\Mk_1}$ over $\Mk_1$ whose $F$-crystal is obtained from ${\got C}_1$ by natural pull back. It is worth pointing out that the existence of $\Mg_{\Mk_1}$ can be deduced as well from [Fa2, th. 10] using standard descent. We sketch the argument. The $\overline{\tilde k}$-valued points of ${\rm Spec}(R_1/2R_1)$ are dense: for instance, cf. 3.6.8.1.2 a) and 3.6.18.4.2 b); to be compared with 3.6.1.3 5). Even better, we can assume that each point of ${\rm Spec}(R_1/pR_1)$ of codimension 1 specializes to a $\overline{\tilde k}$-valued point $\tilde y$, cf. 3.6.8.1.2 b'). Moreover, from loc. cit. we get (see also 3.6.18.3.1) that the $F$-crystal over the completion $\widehat{R_{\tilde y}}$ of the local ring $R_{\tilde y}$ of (any such point) $\tilde y$ obtained naturally by pulling back ${\got C}_1$, is the $F$-crystal of a $2$-divisible group; it is uniquely determined (for instance, cf. [BM, 4.2.4]). Using this and standard Galois descent in an affine context (based on the fact that the natural morphism $m_{\tilde y}:{\rm Spec}(\widehat{R_{\tilde y}})\to {\rm Spec}(R_{\tilde y})$ is regular, see [Ma]), we deduce that the $2$-divisible group over the perfection of $\Mk_1$ obtained naturally via ${\got C}_1$ and classical Dieudonn\'e theory, is in fact naturally definable over $\Mk_1$. 
\smallskip
But now, the fact that the $2$-divisible group we got over $\Mk_1$ extends uniquely to a $2$-divisible group $\Mg_{R_1/2R_1}$ over $R_1/2R_1$ in the way prescribed by ${\got C}_1$, is standard. One very fast way to prove this goes as follows. We consider first a local ring $V$ of $R_1/2R_1$ which is a DVR. Let $\pi_V$ be a uniformizer of $V$ and let $V_1$ be a faithfully flat, finite $V$-algebra which is a DVR having $\pi_V$ as a uniformizer, a residue field which is the perfection of the residue field of $V$, and whose field of fractions is an algebraic extension of $V[{1\over {\pi_V}}]$. Using the above part referring to $m_{\tilde y}$, we know that there is a local, faithfully flat, $V$-algebra $V_2$, which is a DVR such that:
\medskip
-- the pull back $\Mg_{V_2[{1\over {\pi_V}}]}$ of $\Mg_{\Mk_1}$ to ${\rm Spec}(V_2[{1\over {\pi_V}}])$ extends to a $2$-divisible group over ${\rm Spec}(V_2)$;
\smallskip
-- the morphism $m_2:{\rm Spec}(V_2)\to {\rm Spec}(V)$ is regular. 
\medskip
From 3.14 B6 we get that there is a unique $2$-divisible group $\Mg_{\widehat{V_1}}$ over the completion $\widehat{V_1}$ of $V_1$ whose $F$-crystal over $\widehat{V_1}$ is obtained from ${\got C}_1$ by natural pull back. As the field of fractions of $\widehat{V_1}$ has a finite $2$-basis, using again [dJ2] (or using descent as above or [BM, 4.2.4]) we get that the notations are fitting each other, i.e. the generic fibre of $\Mg_{\widehat{V_1}}$ is indeed obtained from $\Mg_{\Mk_1}$ by natural pull back. So (simple argument at the level of $V_1$-lattices) $\Mg_{V_1[{1\over {\pi_V}}]}$ extends to a $2$-divisible group $\Mg_{V_1}$ over $V_1$. Though we do not need this, it is worth pointing out that there is a unique such extension, cf. [dJ3, 1.2]; the reason we do not need this: the $F$-crystal of $\Mg_{V_1}$ is automatically (cf. the construction of $\Mg_{V_1}$) the logic pull back of ${\got C}_1$. Standard descent involving the $V$-algebras $V_1$ and $V_2$ (we can assume $V_2$ is complete and so, as $m_2$ is regular, $V_2\otimes_V V_1$ is a DVR), shows that $\Mg_{\Mk_1}$ extends to a $2$-divisible group over $V$. Using [dJ3, 1.1] ($V$ has a $2$-basis, cf. [BM, 1.1.2 (ii) and (v)] and the fact that ${\rm Spec}(R_1/2R_1)$ is an  $\NN$-pro-\'etale ${\rm Spec}(R/2R)$-scheme) or using descent as above, we get that this extension is as predicted by ${\got C}_1$. 
But now, the passage from points of codimension 1 to the whole of ${\rm Spec}(R_1/2R_1)$ is a local statement for the faithfully flat topology of ${\rm Spec}(R_1/2R_1)$ and so we can apply again the part of the previous paragraph referring to [Fa2, th. 10] (or to [dJ2], as the field of fraction of $W(\tilde k_1)[[x_1,...,x_{dd({\got C}_z))}]]$, with $\tilde k_1$ an algebraic field extension of $\tilde k$, has a finite $2$-basis).   
\smallskip
The trouble with $p=2$ is in lifting things modulo $2$ to things modulo $4$ (cf. 2.3.18.1 A to C), i.e. we have trouble with filtrations. We can assume $\Mg_{R_1/2R_1}[4]$ is obtained by pulling back a finite, flat, commutative group scheme $\Mg(4)_{R/2R}$ over $R/2R$. Let $\hat R$ be the completion of $R$ (or of $R_1$) w.r.t. its maximal ideal defining $z$. We can now go around this ``filtration trouble", using [Fa2, th. 10 and rm. iii) after it]: from loc. cit. we deduce (cf. also 2.2.9 1) and 2.2.1.2) the existence of a cocharacter $\mu_{\hat R}:\GG_m\to\tilde G_{\hat R}$ producing as usual a direct sum decomposition $\tilde M\otimes_{W(\tilde k)} {\hat R}=F^1_{\hat R}\oplus F^0_{\hat R}$ such that the pull back of ${\got C}_{1}$ and its tensors through the logical $W(\tilde k)$-morphism ${\rm Spec}(\hat R)\to {\rm Spec}(R_1)$, when its underlying $\hat R$-module is endowed with the filtration defined by $F^1_{\hat R}$, is a $2$-divisible object with tensors of $\Mm\Mf_{[0,1]}^\nabla(\hat R)$ associated to a Shimura $2$-divisible group over ${\rm Spec}(\hat R)$ lifting $\Md(W(\tilde k))$ and lifting (cf. 3.14 B6) $\Mg_{{\hat R}/2{\hat R}}$. From this and [BLR, th. 12 of p. 83] (applied very much the same as in 3.6.14.4 but only once) we get that by passing to an \'etale $W(\tilde k)$-morphism ${\rm Spec}(R')\to {\rm Spec}(R)$, we can assume that: 
\medskip
{\bf vi)} $\Mg(4)_{R/2R}$ lifts to a finite, flat, commutative group scheme $\Mg(4)$ over $R$, such that the object ${\got C}(4):=\DD(\Mg(4)_{R/2R})$ of $\Mm\Mf_{[0,1]}^\nabla(R)$ has the property that the filtration $F^1_4$ of the underlying $R/4R$-module $\tilde M\otimes_{W(\tilde k)} R/4R$ of ${\got C}(4)$ respects the extra Shimura structure in the logical way (i.e. it is defined as usual by a cocharacter $\mu_4:\GG_m\to \tilde G_{R/4R}$);
\smallskip
{\bf vii)} $\Mg(4)$ in $z$ is the kernel of the multiplication by $4$ of (the $2$-divisible group of) $\Md(W(\tilde k))$.
\medskip
Based on properties of lifting cocharacters (see [SGA3, Vol. II, 3.6 of p. 48]; to be compared with [Va2, 5.3.2]) $F^1_4\otimes_{R/4R} R_1/4R_1$ lifts to a filtration $F^1$ of the $R_1^\wedge$-module $\tilde M\otimes_{W(\tilde k)} R_1^\wedge$ underlying ${\got C}_1$ in such a way that we still get a $2$-divisible object with tensors ${\got C}_2$ of $\Mm\Mf_{[0,1]}^\nabla(R_1)$. Moreover, we can assume that ${\got C}_2$ in the factorization $\tilde z$ of $z$ through ${\rm Spec}(R_1)$ is still ${\got C}_z$. By replacing $R/2R$ as well by a suitable localization of it, from 3.6.18.4.5 we get: we can assume that c) of Theorem 2 holds. So the above two Facts follow from the following Claim: 
\medskip
{\bf Claim.} {\it $\Mg(4)_{R_1^\wedge}$ lifts, as prescribed by ${\got C}_{2}$, to a $2$-divisible group $\Mg_{R_1^\wedge}$ over ${\rm Spec}(R_1^\wedge)$ in such a way that, when endowed naturally with tensors, in $\tilde z$ it is $\Md(W(\tilde k))$.} 
\medskip
{\bf Proof:} 
[Me, ch. 5, 1.6 and 2.3.4-5] implies there is a precisely one way to lift $\Mg_{R_1/2R_1}$ to a $2$-divisible group $\Mg_{R_1/4R_1}$ over ${\rm Spec}(R_1/4R_1)$ such that its $4$-torsion is the pull back of $\Mg(4)$ to ${\rm Spec}(R_1/4R_1)$. From [Me, ch. 5, 1.6]  we deduce, via a natural limit process, that there is a unique way to lift $\Mg_{R_1/4R_1}$ to a $2$-divisible group $\Mg_{R_1}$ over ${\rm Spec}(R_1^\wedge)$ such that its associated $2$-divisible object of $\Mm\Mf_{[0,1]}^\nabla(R_1)$ is ${\got C}_2$. The second part of the Claim (involving $\Md(W(\tilde k))$) is obvious (cf. 2.3.18.1 D and the construction of $\Mg_{R_1}$). This proves the Claim.
\medskip
{\bf 3.15.2. The variant of the Fact of 3.6.19 modulo $p$.} Let $X$ be a $W(\tilde k)$-scheme which locally in the Zariski topology is a pro-\'etale scheme over a smooth $W(\tilde k)$-scheme. We consider the category $$\Mm_{[0,1]}^\nabla(X)$$ 
defined as follows. If $X={\rm Spec}(R)$ is affine and $\Phi_R$ is a Frobenius lift of $R^\wedge$, then:
\medskip
-- its objects are formed by quadruples $(M,\Phi_M,\Phi_M^1,\nabla_M)$, where $M$ is an $R$-module endowed with a $\Phi_R$-linear endomorphism $\Phi_M$ and with a $\Phi_R$-linear map $\Phi_M^1:\tilde F^1(M)\to M$, with $\tilde F^1(M)$ as the pull back of the $R/pR$-submodule of $M/pM$ which is the kernel of $\Phi_M$ mod $p$ via the natural $R$-epimorphism $M\twoheadrightarrow M/pM$, and where $\nabla_M$ is a connection on $M$, with the property that locally in the Zariski topology of $X$ it can be extended to a quintuple $(M,F^1,\Phi_M,\Phi_M^1,\nabla_M)$, with $F^1\subset \tilde F^1(M)$, such that denoting by $\Phi_M^1|F^1$ the restriction of $\Phi_M^1$ to $F^1$, the quadruple $(M,F^1,\Phi_M,\Phi_M^1|F^1,\nabla_M)$ is an object of $\Mm\Mf_{[0,1]}^\nabla(R)$;
\smallskip
-- its morphisms from such a quadruple into a similar one $(M_1,\Phi_{M_1},\Phi_{M_1}^1,\nabla_{M_1})$ are given by parallel (in the natural sense) $R$-linear maps $f:M\to M_1$ such that $\Phi_{M_1}\circ f=f\circ \Phi_M$ (so $f(\tilde F^1(M))\subset \tilde F^1(M_1)$) and $\Phi_{M_1}^1\circ f=f\circ \Phi_M^1$. 
\medskip
If $X$ is not affine, then we consider an arbitrary cover of $X$ by open, affine subschemes, and we glue triples as above using the standard gluing arguments of [Fa1, 2.3] (cf. also 2.2.1 c) for $p=2$); due to existence of $F^1$'s, these arguments apply entirely to the non-filtered context of $\Mm_{[0,1]}^\nabla(X)$.  
\medskip
We also consider the category 
$$\Mm\Mv_{[0,1]}^\nabla(X)$$ 
defined as follows. If $X={\rm Spec}(R)$ is affine and $\Phi_R$ is as above, then the objects of $\Mm\Mv_{[0,1]}^\nabla(X)$ are quadruples 
$$(M,\Phi_M,V_M,\nabla_M),$$ 
with $M$, $\Phi_M$ and $\nabla_M$ as above, with $V_M:M\to M\otimes_R\, _{\Phi_R} R$ an $R$-linear map such that by identifying $\Phi_M$ with an $R$-linear map $\Phi_M:M\otimes_R\, _{\Phi_R} R\to M$ we have 
$$\Phi_M\circ V_M=p1_M\,\, {\rm and}\,\, V_M\circ\Phi_M=p1_{M\otimes_R\, _{\Phi_R} R}$$ 
(i.e. $V_M$ is a Verschiebung map), which locally in the \'etale topology of $X$, are obtained from an object of $\Mm_{[0,1]}^\nabla(X)$ in the logical way (here by ``logical" we mean nothing else but: if $F^1$ and $\Phi_M^1$ are as above, then we have $V_M\circ\Phi_M^1(x)=x$, $\forall x\in F^1$, cf. 2.2.1.0 (VPHIONE)). If $X$ is not affine, then we consider an arbitrary cover of $X$ by open, affine subschemes, and we glue quadruples as above using standard gluing arguments (of crystals on $X_k$ in coherent sheaves). The morphisms of $\Mm\Mv_{[0,1]}^\nabla(X)$ are (as above) the logical ones. 
\smallskip
Fact 1 of 2.2.1.0 implies: we have a natural identification of the category $\Mm_{[0,1]}^\nabla(X)$ with a subcategory (not necessarily full) of $\Mm\Mv_{[0,1]}^\nabla(X)$. We have:
\medskip
{\bf Theorem.} {\it The category of finite, flat, commutative group schemes of $p$-power rank over $X_{\tilde k}$ which locally in the \'etale topology of $X$ lift to finite, flat, commutative group schemes over $X$, is antiequivalent to the category $\Mm\Mv_{[0,1]}^\nabla(X)$.}   
\medskip
{\bf Proof:}
The functor achieving this antiequivalence of categories is the crystalline Dieudonn\'e functor $\DD$ (see [BBM, ch. 3]) evaluated at (the thickening of $X_k$ defined by) $X^\wedge$. $\DD$ is fully faithful, cf. [BM, 4.2.6]. Based on this and standard descent, to check that it is essentially surjective on objects, we can work locally in the ($\NN$-pro-) \'etale topology; in particular we can assume $X={\rm Spec}(R)$ is an affine scheme. We fix a Frobenius lift $\Phi_R$ of $R^\wedge$. So, the same argument based on [EGA IV, 8.5.2] and used to prove the Fact of 3.6.19, allows us to assume $R$ is a smooth $W(\tilde k)$-algebra; warning: not to complicate notations, even when we use below the $\NN$-pro-\'etale topology, we still refer to $X$. The case $p\ge 3$ is a consequence of 2.2.1.1 2), cf. the above definition of $\Mm_{[0,1]}^\nabla(X)$ involving ``locally in the Zariski topology". 
\smallskip
The case $p=2$ can be deduced from the proof of 3.15.1 and the moduli principle in the form of 3.14 D, 3.6.18.4.2 and 3.6.18.5.2. This goes as follows. We can assume we are dealing with an object $\Mo$ of $\Mm_{[0,1]}^\nabla(X)$, which is obtained from an object $\Mo_f$ of $\Mm\Mf_{[0,1]}^\nabla(X)$ by forgetting the filtration. We can also assume (cf. Fact of 2.2.1.1 6)) that $\Mo_f$, when viewed without connection, lifts to a $2$-divisible object ${\got C}$ of $\Mm\Mf_{[0,1]}(X)$, i.e. it is the cokernel of an isogeny $m_f:{\got C}_s\hookrightarrow {\got C}$ between two $2$-divisible objects of $\Mm\Mf_{[0,1]}(X)$. But in the $\NN$-pro-\'etale topology of $X$, ${\got C}$ can be viewed as a $2$-divisible object of $\Mm\Mf_{[0,1]}^\nabla(X)$ lifting $\Mo_f$ (i.e., using connections, $m_f$ can be viewed locally in the $\NN$-pro-\'etale topology of $X$ as an isogeny $m_f^\nabla$ between $2$-divisible objects of $\Mm\Mf_{[0,1]}^\nabla(X)$ having $\Mo_f$ as its cokernel; argument: this is just an abstract extension of the proof of 3.6.18.5.2, cf. 3.6.18.5.4 1)). 
\smallskip
So, from the proof of Theorem 2 of 3.15.1 involving the existence of $\Mg_{R_1/2R_1}$ and from the fully faithfulness of $\DD$ (in the context of $2$-divisible groups over $R_1/2R_1$; see [BM, 4.2.6]), we get that in the $\NN$-pro-\'etale topology of $X$, the object of $\Mm\Mv_{[0,1]}^\nabla(X)$ naturally defined by $\Mo$ is associated (via $\DD$) to a finite, flat, commutative group scheme $G$ of $2$-power order over $X_{\tilde k}$: $G$ is the kernel of an isogeny $i_{12}:D_1\to D_2$ between two $2$-divisible groups over $X_{\tilde k}$. $i_{12}$ is defined via: $\DD(i_{12})$ is the morphism $m_f^\nabla$ but viewed (by forgetting filtrations) as a morphism of $p-\Mm_{[0,1]}^\nabla(X)$. We can assume $G$ lifts to a finite, flat, commutative group scheme over $X$: using Artin's approximation theorem this is a consequence of A) and B6 of 3.14 (applied in the context of the pull back of $m_f^\nabla$ over completions of henselizations of localizations of $X$ w.r.t. its maximal points); but as any finite, flat group scheme is of finite presentation, we can replace the $\NN$-pro-\'etale topology of $X$ by the \'etale topology of $X$. This proves the Theorem. 
\medskip
{\bf 3.15.3. Comments and variants.} {\bf 1)} The expectation of 3.6.8.9.1 gets translated: in the two Theorems of 3.15.1, (at least if $\tilde k$ is ``almost" algebraically closed in some sense) it should be possible to avoid the ``complication" with $\NN$-pro-\'etale morphisms (and so to avoid introducing $R_1$). For the case $\tilde k=\bar{\tilde k}$, in connection to Theorem 1, this is achieved (using extra tools) in (a stronger form in) 4.12.12 and 4.12.12.2. Here we just add: the (ineffective) existence of $\Mg(4)$ can be deduced as well from [Il, 4.8] and Artin's approximation theorem.
\smallskip
{\bf 2)} The argument of 3.15.2 involving the essential surjectivity part (on objects) used for $p=2$ works as well for $p\ge 3$.
\smallskip
{\bf 3)} We have a logical variant of 3.15.2 in the context of the spectrum of a regular $W(\tilde k)$-algebra of formal power series.
It can be proved in the same way: instead of referring to 3.15.1, we can just quote directly 3.6.18.5.1 (for $p\ge 3$) and A) of 3.14 B1 (for $p=2$). Even better, we have: 
\medskip
{\bf Corollary 1.} {\it 3.15.2 remains true, if instead of $X$ we use its completion $\hat X_Y$ along a closed subscheme $Y$ which is formally smooth over $W(\tilde k)$.}
\medskip
{\bf Proof:} $\Mm\Mv_{[0,1]}^\nabla(\hat X_Y)$ is defined as follows. We can restrict (via gluing arguments as in 3.15.2) to local charts. So we can assume $\hat X_Y={\rm Spf}(R)$, with $R$ a $W(\tilde k)$-algebra. We take $\Mm\Mv_{[0,1]}^\nabla(\hat X_Y)$ to be $\Mm\Mv_{[0,1]}^\nabla({\rm Spec}(R))$ defined as in 3.15.1. It is well known that, locally in the Zariski topology of $Y_{\tilde k}$, the completion of $X_{\tilde k}$ along $Y_{\tilde k}$ still has a finite $p$-basis (for instance, see [dJ2, 1.1.3]); so the fully faithfulness part of the Corollary follows from [BM, 4.1.1]. Using an entirely similar algebraization process as in 3.6.20 3), we get (from Theorem of 3.15.2) the essential surjectivity part of the Corollary (the part involving lifts of finite, flat, commutative groups schemes over ${\hat X}_{Y_k}$ to ones over $\hat X_Y$ in the \'etale topology is as in the last paragraph of 3.15.2). 
\smallskip
{\bf 4)} Using the proofs of 3.15.1-2, we get a second proof of b) of the Fact of 3.6.19. It has the advantage that it entirely avoids the use of the language of Fontaine's comparison theory with $W_m(\FF_p)$-coefficients, $m\in\NN$ (we recall that the proof of [Fa1, 7.1] relies on such a theory). Warning: in the way we have presented the things, by avoiding Fontaine's comparison theory, we implicitly use [Fa2, th. 10] (or alternatively [dJ2-3]). 
\smallskip
{\bf 5)} Putting aside the Shimura-ordinariness part of Theorem 2 of 3.15, a) to c) of it can be entirely adapted to the context of 3.6.18.7.0 and even more generally to the context of 3.6.18.7.1 a).
\smallskip
{\bf 6)} The spirit of this paper is ``centered" on perfect fields. However, it can be easily checked that 3.15.2 and Corollary 1 of 3) can be entirely adapted to the case of fields having a finite $p$-basis: we just have to use (besides [BM, 4.1.1 or 4.2.6]) an elementary variant of b) of the Fact of 3.6.19 (which allows us to pass 3.15.2 from perfect fields to fields having a finite $p$-basis), in the context of the categories $\Mm\Mv_{[0,1]}^\nabla(*)$. Warning: if we are dealing with a field $k_0$ having a finite $p$-basis, the sheaves of relative differentials are relative to $\cap_{m\in\NN} k_0^{p^m}$ and not to $k_0$ and so, even when we are dealing with $\Mm\Mv_{[0,1]}^\nabla(W(k_0))$, non-trivial connections do show up. 
\smallskip
For the convenience of the reader we include here a weaker form of the referred elementary variant, which is enough for applications in connection to 3.15.2 and to Corollary 1 of 3). We start with an \'etale, affine ${\rm Spec}(W(k_0)[x_1,...,x_m])$-scheme $X$, with $m\in\NN$; we choose the Frobenius lift of $X^\wedge$ to be the natural one defined by: $x_i$ goes to $x_i^p$, $i=\overline{1,m}$. We write $k_0$ as an inductive limit of smooth, integral, finitely generated $\FF_p$-algebras $R_0^{\al}$, $\al\in\Mj$, with flat but not-necessarily \'etale transition homomorphisms. Let $k_0^{\rm al}$ be the field of fractions of $R_0^{\al}$. We can assume that each $R_0^{\al}$ lifts to a smooth, integral, finitely generated $\ZZ_p$-algebra $R^{\al}$. 
\smallskip
Let $n\in\NN$. We write $X_{W_n(k_0)}$ as the projective limit of $X^{\al}_n$, $\al\in\Mj$, where each $X^{\al}_n$ is an \'etale, affine ${\rm Spec}(W_n(k_0^{\al}))[x_1,...,x_m])$-scheme. Localizing $X$, we can assume that each $X_n^{\al}$ lifts to an \'etale, affine ${\rm Spec}(W(k_0^{\al}))[x_1,...,x_m])$-scheme $X^{\al}$; warning: we do not assume the existence of any transition morphisms between $X^{\al}$'s. The category $\Mm\Mv_{[0,1]}^\nabla(X)$ is defined as in 3.15.2. Let ${\got C}$ be an object of it annihilated by $p^n$. We can assume that the connection on the underlying module of ${\got C}$ is involving $\Om_{X_{W_n(k_0)}/W_n(k_0^{p^n})}$. As $\Om_{X_{W_n(k_0)}/W_n(k_0^{p^n})}$ is a finitely generated $\Mo_{X_{W_n(k_0)}}$-module which is the projective limit of the finitely presented $\Mo_{{X^{\al}_n}}$-modules $\Om_{{X^{\al}_n}/W_n({k_0^{\al}}^{p^n})}$, as in the proof of the Fact of 3.6.19, we get that ${\got C}$ is obtained by pulling back an object ${\got C}^{\al_0}$ of $\Mm\Mv_{[0,1]}^\nabla(X^{\al_0})$, for some $\al_0\in\Mj$. Localizing $R_0^{\al_0}$, we can assume $X_k^{\al_0}$ is obtained from an \'etale, affine scheme over ${\rm Spec}(R_0^{\al_0}[x_1,...,x_n])$ by natural pull back and that this \'etale scheme lifts to an \'etale, affine scheme $Y^{\al_0}$ over ${\rm Spec}(R^{\al_0}[x_1,...,x_n])$. Localizing $R_0^{\al_0}$ further on, we can assume ${\got C}^{\al_0}$ lifts to an object ${\got C}^{\al_0}_{\rm lift}$ of $\Mm\Mv_{[0,1]}^\nabla(Y^{\al_0})$, via the natural homomorphisms 
$$R^{\al_0}/p^nR^{\al_0}\tilde\to W_n(R_0^{\al_0})\rightarrow W_n(k_0^{\al_0})$$ (the first one, cf. the smoothness of $R^{\al_0}/p^nR^{\al_0}$, just lifts the natural identification $R^{\al_0}/pR^{\al_0}=W_1(R_0^{\al_0})$; it is trivial to check that, as denoted, it is an isomorphism). Applying Theorem of 3.15.2 to ${\got C}^{\al_0}_{\rm lift}$, we get that ${\got C}$ is associated to a finite, flat, group scheme annihilated by $p^n$ over $X_k$ which locally in the \'etale topology lifts to a group scheme over $X^\wedge$ (and so --via [BLR, th. 12 of p. 83]-- over $X$).  
\smallskip
The passage from the smooth context to the formally smooth context is entirely the same as the algebraization process of 3.6.20 3): we end up not with $Y^{\al_0}$ but with the infinitesimal neighborhoods (in $Y^{\al_0}$) of a closed subscheme of it which is smooth over $\ZZ_p$. 
\smallskip
As a conclusion, using the above mentioned adaptation, we reobtain in a completely new manner [dJ2, first main result]: even for $p=2$ we can avoid entirely the use of [dJ2-3] in the proofs of the two Theorems of 3.15.1 (cf. their references to descent). Moreover, 3) points out a weakening of the conditions (restrictions) in [dJ2, first main result]. We get:
\medskip
{\bf Corollary 2.} {\it In loc. cit. one can replace ``is locally of finite type over a field with a finite $p$-basis" by: is locally a pro-\'etale (formal) scheme over a (formal) scheme of finite type over a field with a finite $p$-basis.}
\medskip
Warning: Corollary 2 is by no chance a trivial improvement. For instance, there are pro-\'etale schemes over a scheme of finite type over a field $k_0$ with a finite $p$-basis, which have closed subschemes of codimension 1 having no $\overline{k_0}$-valued points. 
\medskip
{\bf 3.15.4. More on b) of the Fact of 3.6.19.} We assume $p\ge 3$. Let $k$ be a perfect field of characteristic $p$. Let $X$ be a regular, formally smooth $W(k)$-scheme which locally in the Zariski topology is a pro-\'etale scheme over a smooth $W(k)$-scheme. Let $\hat X_Y$ be the completion of $X$ along a closed subscheme $Y$ which is formally smooth over $W(k)$. Let $Z$ be a pro-\'etale, formal scheme over $\hat X_Y$; locally $Z$ is the completion of a pro-\'etale scheme over a smooth $W(k)$-scheme along a closed subscheme of it which is formally smooth over $W(k)$. $Z^\wedge$ is the formal scheme obtained by completing $Z$ $p$-adically (so, if $Z=\hat X_Y$, then $Z^\wedge$ is the completion of $X$ along $Y_k$). The category $\Mm\Mf_{[0,1]}^\nabla(Z)$ (and so implicitly $p-\Mm\Mf_{[0,1]}^\nabla(Z)$) is defined as follows. We can restrict (via standard gluing arguments) to local charts, i.e. we can assume $Z$ is the formal scheme defined by completing a pro-\'etale, affine scheme ${\rm Spec}(A)$ over a smooth, affine $W(k)$-scheme along a closed subscheme of it defined by an ideal $I_A$ of $A$ such that ${\rm Spec}(A/I_A)$ is formally smooth over $W(k)$. Denoting by $\hat A$ the completion of $A$ in the $I_A$-topology, we take 
$$\Mm\Mf^\nabla_{[0,1]}({\rm Spf}(\hat A)):=\Mm\Mf_{[0,1]}^\nabla(\hat A).$$ 
So $p-\Mm\Mf^\nabla_{[0,1]}({\rm Spf}(\hat A))=p-\Mm\Mf_{[0,1]}^\nabla(\hat A)$. We have:
\medskip
{\bf Corollary.} {\it The category $p-DG(Z^\wedge)$ (resp. $p-FF(Z^\wedge)$) is antiequivalent (via the $\DD$ functor) to $p-\Mm\Mf_{[0,1]}^\nabla(Z)$ (resp. to $\Mm\Mf_{[0,1]}^\nabla(Z)$).} 
\medskip 
{\bf Proof:}
Starting from Corollary 2 of 3.15.3 6) and Grothendieck--Messing theory of lifting $p$-divisible groups (see [Me, ch. 4-5]) the part referring to $p$-divisible groups follows. The part referring to finite, flat, commutative group schemes of $p$-power rank, results from its part referring to $p$-divisible groups, once we remark three things:
\medskip
a) $\DD$ is faithful (working in the faithfully flat topology, this is a consequence of 3.6.18.5.3);
\smallskip
b) any object (resp. morphism) of $\Mm\Mf_{[0,1]}^\nabla(Z)$, locally in the $\NN$-pro-\'etale topology of $Z$ lifts to a $p$-divisible object (resp. to a morphism between two $p$-divisible objects) of $\Mm\Mf_{[0,1]}^\nabla(Z)$: using Fact of 2.2.1.1 6) and the $p\ge 3$ analogue of 3.14 B3, this is a consequence of 3.6.18.4.2 (resp. of 3.6.18.5.4 1));
\smallskip
c) Galois descent allows us to replace ``in the $\NN$-pro-\'etale topology" by ``in the Zariski topology" in b) (cf. a)).
\medskip
In fact, using the equivalent of Corollary 1 of 3.15.3 3) over $k$, 3.6.18.5.3 and [BLR, th. 12 of p. 83], in b) we can work locally in the \'etale topology, without ``appealing" to a $p$-divisible context.
\medskip
{\bf 3.15.4.1. Variants.} 3.15.4 and 3.6.18.5.8 can be combined. Similarly, if $p=2$, the variant of 3.6.18.5.7 pointed out in 3.14 B5, can be combined with Corollary 2 of 3.15.3 6), to get a variant of 3.15.4 which disregards either all multiplicative type parts or all \'etale parts. Also, 3.15.4 has a version in which $k$ is replaced by any field having a finite $p$-basis: [dJ2, 2.2.2-3] is the main ingredient needed to show that the category $\Mm\Mf_{[0,1]}^\nabla(\hat A)$ is well defined and does not depend on the choice of a Frobenius lift of $\hat A^\wedge$; see also [dJ2, 2.3.4 and 2.4.8] for interpretations of $p-\Mm\Mf_{[0,1]}^\nabla(\hat A)$ in terms of suitable filtered $F$-crystals. Warning: in the case of a field having a $p$-basis, 3.15.4 b) has to be replaced by the paragraph following 3.15.4 c).
\medskip
{\bf 3.15.4.1.1. Some extensions to $p$-divisible contexts.} For future references we also point out the ``$p$-divisible" form of 3.15.2 and of Corollary 1 of 3). We have:
\medskip
{\bf Corollary.} {\it We refer to 3.15.2 (resp. to Corollary 1 of 3)). The category $p-DG(X_{\tilde k})$ (resp. $p-DG({\hat X}_{Y_{\tilde k}})$) is antiequivalent (via $\DD$) to $p-\Mm\Mv_{[0,1]}^\nabla(X)$ (resp. to $p-\Mm\Mv_{[0,1]}^\nabla(X_Y^\wedge)$).}
\medskip
{\bf 3.15.5. Open questions.}  Let $O$ be as in 2.2.1.4 and let $(Y,U)$ be as in 2.2.1.4.1. Let $n\ge 2$, $n\in\NN$. Let $f:G_1\to G_2$ be a morphism between truncations mod $p^n$ of two $p$-divisible groups over $Y_{O/pO}$. In connection to 3.6.14.4.1, 2.2.1.4 and 3.15.2, we have the following question.
\medskip
{\bf Q1.} {\it If $f$ restricted to $U$ is an isomorphism (resp. is a closed embedding), what extra conditions assure us that $f$ is an isomorphism (resp. is a closed embedding)?}
\medskip
We now present an approach which we hope will lead to an answer to Q1. Following the proof of 2.2.1.4 we can assume $Y={\rm Spec}(W(k)[[T]])$, with $k=\bar k$, and its Frobenius lift is defined (at the level of $W(k)$-algebras) via: $T$ goes to $T^p$. Using [Il, 4.8] we get that $G_1$ and $G_2$ lift to truncations mod $p^n$ of $p$-divisible groups over $R:=W(k)[[T]]$. So we can apply 3.15.2 to $f$. We can write $\DD(f)$ as a morphism $m_{12}$ between truncations mod $p^n$ of two objects of $p-\Mm_{[0,1]}(Y)$; here we identify $p-\Mm_{[0,1]}(Y)$ with $p-\Mm_{[0,1]}^\nabla(Y)$ (cf. the logical non-filtered version of 3.6.18.5.1 and its $p=2$ analogue). As $n\ge 2$, often the truncation mod $p$ of $m_{12}$ can be viewed as well as a morphism between two objects of $\Mm\Mf_{[0,1]}(Y)$. If this is so, then as in the proof of 2.2.1.4 we conclude that the induced morphism $f:G_1[p]\to G_2[p]$ is an isomorphism (resp. is a closed embedding) and so using Nakayama's lemma, we get that $f$ is an isomorphism (resp. is a closed embedding). 
\medskip
We now assume $Y^\wedge$ is equipped with a Frobenius lift and that $O=W(k)$. For future references we also point out the following two general forms of Q1 which involve ``arbitrary" $p$-divisible objects and for which the previous paragraph can be naturally adapted. Let $y_0:{\rm Spec}(k_0)\to Y$ be the maximal point of $Y$. Let ${\got C}_1$ and ${\got C}_2$ be two objects of $p-\Mm_{[0,a]}(Y)$, with $a\in\NN$. Let $m:{\got C}_1/p^{n_1}{\got C}_1\to {\got C}_2/p^{n_2}{\got C}_2$ be a morphism between two of their truncations; we view it as an $\Mo_Y$-linear map between $\Mo_Y$-sheaves endowed with Frobenius endomorphisms. Let $i_U$ be the inclusion of $U$ in $Y$. We assume $n_2\ge n_1\ge a+1$. 
\medskip
{\bf Q2.} {\it If $n_1=n_2$ and if $i_U^*(m)$ is an isomorphism (resp. an epimorphism), what extra conditions assure us that $m$ is an isomorphism (resp. an epimorphism)?}
\medskip
{\bf Q3.} {\it We assume the kernel of $i_U^*(m)$ is included in $p^h$ times the underlying $\Mo_{U}$-sheaf of $i_U^*({\got C}_1/p^{n_1}{\got C})$, with $h\in S(a+2,n_1)$. Under what conditions the underlying $W(k_0)$-module $M_0$ of the kernel of $y_0^*(m)$ is included in $p^{h-a-1}$ times the underlying $W(k_0)$-module $M_{01}$ of $y_0^*({\got C}_1/p^{n_1}{\got C}_1)$?}
\medskip
{\bf 3.15.6. Complements on the deformation theory in the generalized Shimura context.} In what follows we point out that 2.2.21 UP, 2.4, 3.2.4, a great part of 3.6 and 3.15.1 hold in the generalized Shimura context, provided we state everything in terms of (filtered) $F$-crystals. To avoid repetitions we refer to parts of 3.6 and 3.11 even if $p=2$, without mentioning each time that this is allowed by (different parts of) 3.14. 
\smallskip
We start with a generalized Shimura $p$-divisible object ${\got C}=(M,(F^i(M))_{i\in S(a,b)},\vph,G)$ of $\Mm\Mf_{[a,b]}(W(k))$; let $(t_{\al})_{\al\in\Mj}$ be a family of tensors of $\Mt(M[{1\over p}])$ as in 2.2.9 3a). Let $\mu:\GG_m\to G$ be the canonical split cocharacter of ${\got C}$ (see 2.2.1.2). Let $N$ be the integral, closed subgroup of $G$ whose Lie algebra is the Lie algebra of ${\rm Lie}(G)$ on which $\GG_m$ acts via $\mu$ (and inner conjugation) through the identical character. It is commutative and unipotent. Let $R$ be the completion of $N$ in the origin. We choose an isomorphism $\tilde f:R\tilde\to W(k)[[x_1,...,x_d]]$, with $d:=\dim_{W(k)}(N)$. Let $\Phi_R$ be the Frobenius of $R$ as in 2.2.10. We consider the triple
$${\got C}_R:=(M\otimes_{W(k)} R,(F^i(M)\otimes_{W(k)} R)_{i\in S(a,b)},n(\vph\otimes 1)),$$
with $n\in N(R)$ as the universal element. There is a unique connection $\nabla$ on ${\got C}_R$ (cf. 3.6.18.7.1 b) and c); the rank of $\Phi_R$ is $0$). 
Let $P_{\le 0}$ (resp. $P_{\ge 0}$) be the parabolic subgroup of $G$ having $W^0({\rm Lie}(G),\vph)$ (resp. $W_0({\rm Lie}(G),\vph)$) as its Lie algebra. Let $P_{<0}$ be the unipotent radical of $P_{\le 0}$. 
\medskip
{\bf A.} We first assume ${\got C}$ is a Shimura-canonical lift in the sense of 3.11.6.1: to motivate what follows in D and E below and to first deal with shorter proofs, we treat till end of B the simplest case of such a ${\got C}$. 
\smallskip
From 3.11.6 B) and the Fact of 2.2.11.1 we get that $N$ is a subgroup of $P_{<0}$. Moreover, $\nabla$ is of the form $\dl+\be$, with $\dl$ as the connection on $M\otimes_{W(k)} R$ annihilating $M$ and with 
$$\be\in {\rm Lie}(P_{<0})\otimes_{W(k)} \Om_{R/W(k)}^\wedge.\leqno (BETA)$$
(BETA) is very much the same as [Fa2, rm. ii) after th. 10], being a property of $\sg$-$\Ms$-crystals; it is a consequence of 3.6.18.7.1 b) and c). 
\smallskip
We consider the triple
$${\got C}_R^{\le 0}:=({\rm Lie}(P_{\le 0})\otimes_{W(k)} R,F^{0}({\rm Lie}(P_{\le 0}))\otimes_{W(k)} R,n(p\vph\otimes 1)).$$
It is an object of $p-\Mm\Mf_{[0,1]}(R)$. The generic fibre of the action (via inner conjugation) of $P_{<0}$ on ${\rm Lie}(P_{\le 0})$ is faithful (this can be seen easily by looking at the restriction of ${\rm ad}(x)$, with $x\in {\rm Lie}(P_{<0})$, to the Lie algebra of a maximal torus of $P_{\le 0}$). So $\nabla$ is integrable if the connection on ${\got C}_R^{\le 0}$ it induces naturally (and with which it can be identified) is integrable. So from 3.6.18.4.1 we get: 
\medskip
{\bf Corollary.} {\it $\nabla$ is integrable and so $({\got C}_R,\nabla)$ is a $p$-divisible object of $\Mm\Mf_{[a,b]}^\nabla(R)$.}
\medskip
{\bf B.} The Kodaira--Spencer map of $\nabla$ is injective (this, as in [Va2, 5.4.8], is a consequence of the shape of $N$). Let $R_1:=W(k)[[z_1,...,z_m]]$, with $m\in\NN$, and let $\Phi_{R_1}$ be its Frobenius lift taking $z_i$ into $z_i^p$. We consider a $p$-divisible object $({\got C}_1,(t_{\al}^\prime)_{\al\in\Mj})$ with tensors of $\Mm\Mf_{[a,b]}(R_1)$ which under pull back via the Teichm\"uller lift $z_1:{\rm Spec}(W(k))\to {\rm Spec}(R_1)$ is $1_{\Mj}$-isomorphic (in the sense of 2.2.9 6)) to ${\got C}$; we view this $1_{\Mj}$-isomorphism as an identification. We have:
\medskip
{\bf Theorem.} {\it There is a unique connection $\nabla_1$ on ${\got C}_1$. It is integrable, nilpotent mod $p$ and respects the $G$-action. Moreover, there is a uniquely determined $W(k)$-morphism $z_R:{\rm Spec}(R_1)\to {\rm Spec}(R)$ such that $z_R\circ z_1$ is the Teichm\"uller lift ${\rm Spec}(W(k))\to {\rm Spec}(R)$ and $({\got C}_1,\nabla_1,(t_{\al}^\prime)_{\al\in\Mj})$ is isomorphic to $z_R^*({\got C}_R,\nabla,(t_{\al})_{\al\in\Mj})$ through an isomorphism which in $z_1$ is the mentioned identification.}
\medskip
{\bf Proof:} 3.6.18.7.1 b) and c) imply the existence and the uniqueness of $\nabla_1$ as well as that $\nabla_1$ respects the $G$-action. The rest is entirely the same as the arguments of [Fa2, th. 10 and rm. iii) after it], as the Kodaira--Spencer map of $\nabla$ is injective and $d=\dim_{W(k)}({\rm Lie}(G)/F^0({\rm Lie}(G)))$ (the essence of loc. cit. is captured by these two properties and by Corollary of A). This ends the proof.
\medskip
{\bf C.} We do not assume any more that $(M,\vph,G)$ is Shimura-ordinary. As in 3.2.3 we deduce the existence of $g\in G(W(k))$ such that $({\rm Lie}(G),g\vph,F^0({\rm Lie}(G)),F^1({\rm Lie}(G)))$ is of Borel type. From 3.11.6 B) we get that $(M,(F^i(M))_{i\in S(a,b)},g\vph,G)$ is a Shimura-canonical lift of $(M,\vph,G)$. As above we construct ${\got C}_R$. As 3.4-5 treat the generalized Shimura context, from 3.4.6 we get that the slope $-1$ of $({\rm Lie}(G),\vph)$ is at most equal to the slope $-1$ of $({\rm Lie}(G),g\vph)$. From 3.6.18.5.5 B we get that the pseudo-multiplicity of $({\rm Lie}(G),h\vph)$ is the same as the multiplicity of the slope $-1$ of $({\rm Lie}(G),h\vph)$, $\forall h\in G(W(k))$ (cf. also 3.4.5.1 A and B). 
\medskip
{\bf D. Theorem.} {\it In 3.6.18.7.3 C the word potentially can be dropped (i.e. $M(\got C)$ and $M({\got C}/p^n{\got C})$ are moduli schemes of integrable connections).}
\medskip
{\bf Proof:} We refer to the notations of 3.6.18.7.3. We can assume that $R=W(k)[[x_1,...,x_m]]$, with $m\in\NN$ and $k=\bar k$, and that $X={\rm Spec}(R)$. So in Proposition of 3.6.18.7.3 A we have $R_1=R=R^\wedge$. From 3.6.18.7.1 b) and c) we get: any connection $\nabla$ on ${\got C}$ which respects the $G_R$-action is of the form $\dl_R+\be$, with $\be\in {\rm Lie}(G^{\rm der})\otimes_R \bar\Om_{R/W(k)}$. We need to show: $\nabla$ is integrable. 
\smallskip
Based on C, the proof of 3.6.18.4.1 can be adapted to a great extend to the context of 3.6.18.7.3 C. Warning: there is one significant difference which does not allow us just to copy the mentioned proof. The obstacle to be overcome is: 
\medskip
{\bf ob)} In 3.6.18.2 we had to deal only with integral slopes which made its proof very simple. The slopes of $(M,g\vph)$ (of C) are often ``far" from being integral. 
\medskip
The overcome ob), we first pass to an adjoint context as follows.
Due to the fact that $\dl_R$ is defined via the $\ZZ_p$-structure $M_{\ZZ_p}$ of the last proof of 3.6.18.7.3 A, to show that $\nabla$ is integrable it is enough to show that the connection on ${\rm Lie}(G_R^{\rm ad})$ induced (as in 2.2.10) by $\nabla$ is integrable. So we can assume that $G_R$ is adjoint and that the representation of $G_R$ on $M_R$ is the adjoint representation; so $M_R={\rm Lie}(G_R)$ and $(a,b)=(-1,1)$.  
\smallskip
$G_R$ is the pull back of an adjoint subgroup $G_{\ZZ_p}$ of $GL(M_{\ZZ_p})$. So ${\rm Lie}(G_R)={\rm Lie}(G_{\ZZ_p})\otimes_{\ZZ_p} R$. Based on the Lemma of 3.6.18.7.3 A, we can assume there is a cocharacter $\mu_R:\GG_m\to G_R$ defining (as in 2.2.11 1)) the filtration $(F^i({\rm Lie}(G_R))_{i\in S(-1,1)}$ of $M_R$; even more we can assume that $\mu_R$ is obtained from a cocharacter $\mu$ of $G_{W(k)}$ by natural pull back. So we can write $\Phi_{M_R}=g_R(\vph\otimes 1)$, where $\vph:=\sg\mu({1\over p})$, with $\sg$ viewed as a $\sg$-linear automorphism of $M:={\rm Lie}(G_{W(k)})$ fixing ${\rm Lie}(G_{\ZZ_p})$.
\smallskip
Based on this expression of $\Phi_{M_R}$ and on C we can algebraize the things as in the proof of 3.4.18.6.1: so we ``connect" ${\got C}$ with a generic situation defined by a triple 
$${\got C}_{R}:=(M\otimes_{W(k)} R,(F^i(M))_{i\in S(-1,1)}\otimes_{W(k)} R,g_{R}(g\vph\otimes 1))$$
which is an object of $p-\Mm\Mf_{[-1,1]}(R)$. 
Here $R$ is endowed with a Frobenius lift $\Phi_R^1$ of multiplicative type (possible different from $\Phi_U$  of 3.6.18.7.3), $g_R\in G_R(R)$ mod $(x_1,...,x_m)$ is $1_M$, $F^i(M)$ is the direct summand of $M$ such that $F^i(M_R)=F^i(M)\otimes_{W(k)} R$, $i=\overline{0,1}$, while $g\in G_{W(k)}(W(k))$ is (as in C) such that the Shimura adjoint Lie $\sg$-crystal $(M,g\vph)$ is Shimura-ordinary. So the Theorem follows from the following Lemma.
\medskip
{\bf Lemma.} {\it Any connection $\nabla_R$ on ${\got C}_{R}$ respecting the $G_R$-action (i.e. of the form $\dl+\be$, with $\dl$ as in A and with $\be\in{\rm Lie}(G_R)\otimes_R \bar\Om_{R/W(k)}$) is integrable.}
\medskip
{\bf Proof:} We can assume $(M,g\vph)$ is cyclic (equivalently, that $G_{\ZZ_p}$ is $\ZZ_p$-simple). As 3.4 handles the case of $Cl(M,\vph,G_{W(k)})$ as well, from 3.4.8 we get that the Newton polygons of pull backs of ${\got C}_R$ through geometric points of ${\rm Spec}(R/pR)$ all have the same Newton polygon and so are Shimura adjoint Lie $F$-crystals which are Shimura-ordinary; here the pull backs are obtained via Teichm\"uller lifts. So based on this we can imitate the arguments of B to overcome the obstacle ob). 
\smallskip 
For this we consider the $\Phi_R^1$-linear map $\Psi:M_R\to M_R$ obtained as $\Psi$ of 3.4.5; we call it the Faltings--Shimura--Hasse--Witt shift of ${\got C}_{R}$. Let $N$ be the $\ZZ_p$-submodule of $M_R$ formed by elements fixed by $\Psi$. $N\otimes_{\ZZ_p} R$ is a direct summand of $M\otimes_{W(k)} R$, cf. [De3, 1.2.4] and the above part on (constant) Newton polygons. 
\smallskip
If $G_{W(k)}$ is not (resp. is) of $B_n$ Lie type, let $\tilde P_{\le 0}$ be the subgroup of $G_R$ (resp. of the simply connected semisimple group cover $G^{\rm sc}_R$ of $G_R$) normalizing $N\otimes_{\ZZ_p} R$; we can view $N\otimes_{\ZZ_p} R$ as well as a Lie subalgebra of ${\rm Lie}(G^{\rm sc}_R)$. Based on the Claim of 3.5.4 (for $p=2$, cf. also 3.14 C and J), by reasons of dimensions we get that for any Teichm\"uller lift $TL:{\rm Spec}(W(k_1))\to {\rm Spec}(R)$, with $k_1$ a perfect field containing $k$, the image of $\tilde P_{{\le 0}W(k_1)}$ in $G_{W(k_1)}$ is the parabolic subgroup having $W^0(M_R\otimes_{R} W(k_1),TL^*(g_R)(g\vph\otimes 1))$ as its Lie algebra (it normalizes $N\otimes_{\ZZ_p} W(k_1)$, cf. 2.2.3 3)); this justifies our notation with ${\le 0}$ as a lower right index. So as in the proof of 3.6.18.7.3 A involving the adjoint context we get that $\tilde P_{\le 0}$ is a parabolic subgroup of $G_R$ (resp. of $G^{\rm sc}_R$). Let $P_{\le 0}$ be the image of $\tilde P_{\le 0}$ in $G_R$. $pg_{R}(g\vph\otimes 1)$ takes ${\rm Lie}(P_{\le 0})$ onto itself. 
\smallskip
We consider the intersection $H_{R/pR}$ of $P_{\le 0R/pR}$ with $P_{R/pR}$. Its fibres are smooth, connected and of same dimension: using the same type of pull backs as above, this follows easily from 3.11.6 A and B (more precisely, the smoothness part follows from the fact that the intersection $P_{\le 0W(k_1)}\cap P_{W(k_1)}$ contains a maximal torus of $G_{W(k_1)}$; see [SGA3, Vol. III, 4.1.1-2 and 4.5 of Exp. XXVI]). So (as in the part of 3.6.18.7.3 A refering to $N_{R/pR}$) the group scheme $H_{R/pR}$ is smooth. As $H_k$ contains a maximal torus of $G_k$, from [SGA3, Vol. II, 3.6 of p. 48] we get $H_{R/pR}$ contains a maximal torus of $G_{R/pR}$. Following the proof of Fact 2 of 2.2.9 3), we can assume (after suitable $G(R)$-conjugation) that $H:=P_{\le 0}\cap P$ contains a maximal torus of $G_R$ and that $\mu_R$ factors through it. We get that the triple
$$({\rm Lie}(P_{\le 0}),F^0({\rm Lie}(P_{\le 0})),g_{R}(g\tilde\vph\otimes 1)),$$
with $\tilde\vph:=p\vph$, is an object of $p-\Mm\Mf_{[0,1]}(R)$. As in A we get that $\be\in {\rm Lie}(P_{\le 0})\otimes_R\bar\Om_{R/W(k)}$ and that $\nabla_R$ is integrable. This ends the proof of the Lemma and so of the Theorem.
\medskip
{\bf Warning.} The way we got this Corollary does not apply to get that all connections of 3.6.18.7.1 c) are integrable. But we do expect that all such connections are integrable. 
\medskip
{\bf E.} D achieves the deformation theory in the generalized Shimura context. So many parts of 3.1-4 and 3.6.7-14 hold (occasionally under minor restatements) for the generalized Shimura context. In particular, we get:
\medskip
{\bf Corollary.} {\it B holds without assuming that $(M,g\vph,G)$ is Shimura-ordinary.}
\medskip
 For the sake of convenience, the part involving Newton polygons (i.e. pertaining to 3.1.0, 3.2.4, etc.) is deferred to Appendix. Here we just point out that:
\medskip
-- 3.1.0 d), 3.7.6 and 3.12 remain valid in the generalized Shimura context (no change of arguments are needed);
\smallskip
--  3.6.1.2-3, 3.6.18.2-4 and 3.6.18.8 for the generalized Shimura context are already ``encompassed" by 3.6.18.7.1 c), 3.6.18.7.3 C and D; 
\smallskip
-- one can state forms of 3.6.18.5.1, 3.6.18.5.3 and of 3.6.18.5.7-8 for the generalized Shimura context; as this context is not well suited from the point of view of morphisms, these forms are not stated here: we just mention that they are a consequence of 3.6.18.7.1 c), 3.6.18.7.3 C and D;
\smallskip
-- the inducing property of 3.6.18.5 still holds for the generalized Shimura context: we just need to restate its c) part in terms of $-1$ slopes of attached Shimura adjoint Lie $\sg$-crystals;
\smallskip
-- Corollary is the natural extension of 2.2.21 UP;
\smallskip
-- all of 2.4 can be adapted for $p\ge 5$ to the generalized Shimura context (we just need to work --cf. 2.2.13.4-- with filtrations in the adjoint context and so involving the range $[-1,1]$); moreover, we have weaker versions of these adaptations for $p\in\{2,3\}$ (to be compared with 2.2.16.5 and 2.3.18.3). 
\medskip
{\bf 3.15.7. The boundedness principle.} There is a widely spread opinion that $p$-divisible groups (and so by extrapolation $p$-divisible objects) involve an infinite process. The Fundamental Lemma of 3.6.15 B points out that this opinion is not quite accurate. Our philosophy (already hinted at in 3.6.18.10) is: 
\medskip
{\bf Ph.} {\it $p$-divisible objects involve a ``bounded infinite" process.}
\medskip
Till the end of 3.15.7 we formalize what we call the boundedness principle and so ``implicitly" explain what we mean by a ``bounded infinite" process. Different forms of it will be indexed by numbers attached to $BP$. 
\smallskip
We start with definitions meant to ease the statement of results. Till end of \S 3 we work with $p\ge 2$ and with $k$ a perfect field of characteristic $p$. In all that follows $M$ is a free $W(k)$-module of finite rank $d_M\in\NN$. 
\medskip
{\bf A. Definitions.} {\bf a)} By an elementary Dieudonn\'e (resp. elementary Dieudonn\'e--Fontaine) object of $p-\Mm(W(k))$ we mean an object $(M,\vph)$ for which there is a $W(k)$-basis $\{e_1,...,e_{d_M}\}$ of $M$ such that $\vph(e_i)=p^{n_i}e_{i+1}$, $i=\overline{1,d_M}$, with $e_{d_M+1}=e_1$, $n_i=0$ if $i\in S(1,d_M-1)$, and $n_{d_M}\in\ZZ$ is relatively prime to $d_M$ (resp. with $e_{d_M+1}=e_1$ and $n_i$'s being either all non-negative or all non-positive, and which can not be written as a direct sum of non-trivial $p$-divisible objects whose underlying $W(k)$-modules have such $W(k)$-basis). Such a $W(k)$-basis $\{e_1,...,e_{d_M}\}$ is called a standard $W(k)$-basis of $(M,\vph)$. An object of $p-\Mm(W(k))$ is called a Dieudonn\'e (resp. Dieudonn\'e--Fontaine) object if it is a direct sum of elementary Dieudonn\'e (resp. of elementary Dieudonn\'e--Fontaine) objects of $p-\Mm(W(k))$.  
\smallskip
{\bf b)} By the Dieudonn\'e volume (resp. torsion) of a latticed isocrystal ${\got C}=(M,\vph)$ we mean the smallest number $DV({\got C})\in\NN\cup\{0\}$  (resp. $DT({\got C})\in\NN\cup\{0\}$) such that there is an isogeny $m:{\got C}_1\hookrightarrow {\got C}_{W(\bar k)}$ (of latticed isocrystals), with ${\got C}_1$ a Dieudonn\'e object of $p-\Mm(W(\bar k))$, defined at the level of $W(\bar k)$-modules by a monomorphism $m:M_1\hookrightarrow M$ such that $M/m(M_1)$ has length $DV({\got C})$ (resp. such that $p^{DT({\got C})}M\subset m(M_1)$).
\smallskip
{\bf c)} Similarly to b) we define Dieudonn\'e--Fontaine volume $DFV({\got C})$ (resp. torsion $DFT({\got C})$) of ${\got C}$.
\smallskip
{\bf d)} We consider an $n$-tuple $\tau:=(a_1,...,a_n)$ formed by integers. If $\sum_{i=1}^n a_i$ is non-negative (resp. non-positive) then by the non-negative (resp. non-positive) sign deviation of $\tau$ we mean the maximum value of $-1$ (resp. of $1$) times sums which are of the form $\sum_{i=m}^{s} a_i$, where $m\in S(1,n)$ and $s\in S(m,n+m-1)$, and which have the property that all its subsums of the form $\sum_{i=m^1}^{s} a_i$, with $m_1\in S(m,s)$, are (warning!) non-positive (resp. are non-negative). If $\sum_{i=1}^n a_i$ is positive (resp. is negative), then by the sign deviation $SD(\tau)$ of $\tau$ we mean its non-negative (resp. its non-positive) sign deviation. If $\sum_{i=1}^n a_i=0$, then by the sign $SD(\tau)$ of $\tau$ we mean the minimum between its non-negative and non-positive sign deviations. 
\smallskip
{\bf d')} As in d) we define the (non-negative or the non-positive) value deviation $VD(\tau)$ of $\tau$, by considering sums of different entries (not necessarily consecutive) of $\tau$ which have the same sign. 
\smallskip
{\bf e)} The entry $a_q$ of $\tau$ is called quasi-special, if the sum $\sum_{i=q+s}^q a_i$ has a sign which does not depend on $s\in S(-n+1,0)$. If all these sums are negative (so also $a_q<0$), we refer to $a_q$ as a special negative entry of $\tau$. Similarly we define special positive entries of $\tau$. By the negative (resp. positive) $sp$-invariant of $\tau$ we mean the number 
$$sp_-(\tau)\in\NN\cup\{0\}$$ 
(resp. the number $sp_+(\tau)\in\NN\cup\{0\}$) of special negative (resp. of special positive) entries of $\tau$.   
\medskip 
The existence of $DV({\got C})$ and so of $DT({\got C})$ is equivalent to Dieudonn\'e's classification of isocrystals over $\bar k$ and this explains our terminology. Obviously 
$$DFV({\got C})\le DV({\got C})$$
and 
$$DFT({\got C})\le DT({\got C}).$$
\indent
{\bf Examples.} $SD(-1,1,-1,-1,1,1,0,-1)=1+1=2$, $SD(1,1,-2,1,3)=2$ and $SD(-1,1,-1)=1$. $sp_-(-1,1,-1,-1,1,1,-1,0,-1)=2$ (the special negative entries being the last two $-1$).
\medskip
3.6.15 B can be restated as:
\medskip
{\bf BP0.} {\it The isomorphism deviation of any $\sg$-$\Ms$-crystal $(M,\vph,G)$ over $k$ is a finite number which is not bigger than $2DT({\rm Lie}(G),\vph)+1+s_L(\vph)$, where $s_L(\vph)\in\NN\cup\{0\}$ is the smallest number such that $p^{s_L(\vph)}\vph({\rm Lie}(G))\subset {\rm Lie}(G)$.}
\medskip
{\bf B. Remarks.} {\bf 0)} It is an easy exercise to check that in the proof of 3.14.5 B we can in fact use Dieudonn\'e--Fontaine objects. So in BP0 we can replace $DT({\rm Lie}(G),\vph)$ by $DFT({\rm Lie}(G),\vph)$. This represents considerable improvements in practical calculations. For instance, referring to a) of A, if $(n_1,n_1,...,n_{d_M})=(1,1,...,1,-1)$, with $d_M\ge 2$, then $DFT({\got C})=1$ while $DT({\got C})={\rm max}\{1,d_M-2\}$. More generally, we have the following sequence of 3 inequalities
$$DFT({\got C})\le SD(n_1,...,n_{d_M})\le VD(n_1,...,n_{d_M})\le\sum_{i=1}^{d_M}\abs{n_i}.\leqno (INEQ)$$
The second and the third inequality follow from very definitions.
\smallskip
We now check the first inequality in the case $\sum_{i=1}^n n_i>0$ (the other cases are entirely the same) and not all $n_i$'s are non-negative. We consider the biggest $u\in\NN$ such that there is $m\in S(1,n)$ with the property that $a_{m-v,m}:=\sum_{i=m-v}^m n_i\le 0$, $\forall v\in S(0,u)$; so $n_{m+1}$ and $n_{m-u-1}$ are positive and $n_m\le 0$. For $v\in\ S(0,u)$, we replace $n_{m-v}$ by 
$$p^{-a_{m-v,m}}n_{m-v}.$$ 
As we are in a circular context, we can assume $m-u=1$. If $d_M=u+1$ or if all $n_i$'s, with $i\in S(u+2,...,d_M)$, are non-negative we are done, as by very definitions all $-a_{m-v,m}$'s belong to $S(0,SD(n_1,...,n_{d_M}))$. If this is not the case, we next deal with the remaining segment, formed by the integers $n_{u+2}$,..., $n_{d_M}$. We repeat the operation. We choose the biggest $u_1\in\NN$ such that there is $m_1\in S(u+2,d_M)$ with the property that $a_{m_1-v_1,m_1}:=\sum_{i=m_1-v_1}^{m_1} n_i\le 0$, $\forall v_1\in S(0,u_1)$; so $n_{m_1+1}$ and $n_{m_1-u_1-1}$ are positive and $n_{m_1}\le 0$. Due to the choice of $u$ (of being the biggest), we have $m_1-u_1\ge u+2$. For $v_1\in\ S(0,u_1)$, we replace $n_{m_1-v_1}$ by $p^{-a_{m_1-v_1,m_1}}n_{m_1-v_1}$ and we repeat the operation. So by induction on the number of remaining entries (they do not need to be indexed by a set of consecutive numbers of $S(1,d_M)$), we get that the first inequality holds.
\smallskip
In particular, if $n_i\in S(-1,1)$, $\forall i\in S(1,n)$, then (as a substitute of (INEQ)) it is convenient to work with the inequality
$$DFT({\got C})\le\min\{n^-,n^+\},\leqno (CONVINEQ)$$
where $n^-$ (resp. $n^+$) is the number of $i$'s such that $n_i=-1$ (resp. such that $n_i=1$).
However, for being quicker (as the goal below is not to be sharp), for the general estimates of C below we still work with $DT$'s instead of $DFT$'s.   
\smallskip
{\bf 1)} The union of the categories $p-\Mm_{[0,a]}(W(k))$, $a\in\NN$, is a full subcategory of the category of $\sg$-crystals over $k$ but is not equal to it. However, 3.6.15 B and BP0 can be entirely adapted to the context where we work with triples $(M,\vph,G)$, where $(M,\vph)$ is a latticed isocrystal and $G$ is an arbitrary integral subgroup of $GL(M)$ such that ${\rm Lie}(G_{B(k)})$ is normalized by $\vph$. Argument: using dilatations as in 3.9.9 B and C, we replace the closed subgroup $G$ of $GL(M)$ by a homomorphism $G_1\to GL(M)$ factoring through $G$, with $G_1$ a smooth group scheme over $W(k)$ having the same generic fibre as $G$; as there is $s\in\NN$ such that $p^s({\rm Lie}(G_{B(k)})\cap {\rm End}(M))\subset {\rm Lie}(G_1)$, we can proceed as in the proof of 3.6.15 B to show (by working with the triple $(M,\vph,G_1)$) that the isomorphism deviation of $(M,\vph,G)$ is still a finite number. 
In particular, this applies to latticed isocrystals over $k$. 
\smallskip
Also one can similarly check that the extra condition of 3.16.6 B pertaining to the $W$-condition is not needed.
\smallskip
{\bf 2)} We refer to the proof of 3.6.15 B. In practice it is more convenient to work with the following estimate
$$n-m\ge \max\{m+1+s_L(\vph),n(\Mc_1),...,n(\Mc_s)\}\leqno (CONVEST)$$
instead of (EST) of the mentioned place.
The reason is: for such an $n$, $\vph(\tilde h^1)\in {\rm End}(M)$, provided we have all $u_l$'s greater than $n-m$. This has a version in the context of Dieudonn\'e--Fontaine objects (cf. 0)).
\smallskip
{\bf 3)} We refer to 3.6.16 2). We think it is an interesting problem to compute $isom-d(M,\vph,G)$ (or at least to get sharp estimates of it) perhaps in terms of Hodge numbers of $(M,\vph^s)$ and (or) of $({\rm Lie}(G),\vph^s)$, $s\in\NN$. We do not think (EST) of 3.6.15 B is a very sharp estimate (for instance, for the case of the $\sg$-crystal of a supersingular elliptic curve we get $n\ge 4$ and not $n\ge 1$); however, we do not think it can be significantly improved (in general). The same applies to the context of principal bilinear quasi-polarizations.
\smallskip
{\bf 4)} As in 2.2.22.1, we define additively the negative and the positive $sp$-invariant of any cyclic diagonalizable object of $p-\Mm\Mf(W(k))$ and of the object of $p-\Mm(W(k))$ defined by it by forgetting the filtration.     
\medskip
{\bf C. Estimates.} Let $(M,\vph)$ be a latticed isocrystal over $k$. We want to estimate $DT(M,\vph)$. We can assume $k=\bar k$. Let $SL$ be the set of slopes of $(M[{1\over p}],\vph)$. For $\al\in SL$, we write 
$$\al={a_{\al}\over b_{\al}},$$ 
with $a_{\al},b_{\al}\in\ZZ$, $b_{\al}>0$ and $(a_{\al},b_{\al})=1$; let $m_{\al}\in\NN$ be such that the multiplicity of the slope $\al$ is $m_{\al}b_{\al}$. Let 
$$s\in\NN\cup\{0\}$$ 
be the smallest number such that $\vph(p^sM)\subset M$. 
We refer to it as the $s$-number of $(M,\vph)$. Let $d:=l.c.m.\{b_{\al}|\al\in SL\}$. Let $\al_{\rm max}\in\NN\cup\{0\}$ (resp. $e\in\NN$) be the maximum of all $sb_{\al}+a_{\al}$'s (resp. of all $m_{\al}b_{\al}$'s). 
Using 2.2.3 3)
we get a natural monomorphism
$$q_{SL}:\oplus_{\al\in SL} W(\al)(M,\vph)\hookrightarrow M.$$
Let 
$$h\in\NN\cup\{0\}$$ 
(resp. $h_{\al}\in\NN\cup\{0\}$) be the greatest Hodge number of $(M,p^s\vph)$ (resp. of $(W(\al)(M,\vph),p^s\vph)$). Let 
$$\tilde h:=\max\{h_{\al}|\al\in SL\}.$$ 
We refer to $\tilde h$ (resp. to $h$) as the small $h$-number (resp. as the $h$-number) of $(M,\vph)$. We have
$$\tilde h\le h$$
and the equality holds if $q_{SL}$ is an isomorphism; this
motivates the use of the word small. 
\medskip
{\bf BP1.} {\it For any triple $(a,b,c)\in\NN\times\NN\cup\{0\}\times \NN\cup\{0\}$, there is a smallest number 
$$n(a,b,c)\in\NN\cup\{0\}$$ 
such that $DT({\got C})\le n(a,b,c)$, for any latticed isocrystal ${\got C}$ over $k$ of rank $a$, $s$-number $b$ and $h$-number $c$. In particular, $DT(M,g\vph)\le n(d_M,s,h)$, $\forall g\in GL(M)$.}
\medskip
{\bf Proof:} We use induction on $a$: the case $a=1$ is trivial. We work in the context of $(M,\vph)$. We have
$$DT(M,\vph)\le es+DT(M,p^s\vph).\leqno (EST1)$$
Argument: if $N$ is a $W(k)$-submodule of $M$ such that the pair $(N,p^s\vph)$ is an elementary Dieudonn\'e object of $p-\MM(W(k))$ and if $\{b_1,...,b_n\}$ is a standard $W(k)$-basis of $(N,p^s\vph)$, then the pair $(<p^{ns}b_1,p^{ns-s}b_2,...,p^sb_{n-1},b_n>,\vph)$ is an elementary Dieudonn\'e object of $p-\Mm(W(k))$; so we just have to remark that $n\le e$. 
\smallskip
So from now on (till the Fact below inclusive) we assume $s=0$. 
\smallskip
We first consider the case when $\abs{SL}=1$ and $m_{\al}=1$. We consider a $B(k)$-basis $\{e_1,...,e_{b_{\al}}\}$ of $M[{1\over p}]$ formed by elements of $M$ and such that $\vph(e_i)=e_{i+1}$, $i=\overline{1,b_{\al-1}}$, and $\vph(e_{b_{\al}})=p^{a_{\al}}e_1$. We also assume $e_1\in M\setminus pM$. Let $t\in S(0,b_{\al}-1)$. 
Let $e_0:=\vph^{-1}(e_1)\in M[{1\over p}]$. 
\medskip
{\bf Claim.} {\it There is an increasing sequence $(c_t)_{t\in S(0,b_{\al}-1)}$ of elements of $\NN$ which depends only on $\al$ and $h_{\al}$ and such that
$$\vph(e_t)\notin \sum_{i=1}^t W(k)e_i+p^{c_t}M,\leqno (EST2)$$
$\forall t\in S(0,b_{\al}-1)$.}
\medskip
To check the Claim, we use a second induction on $t\in S(0,b_{\al}-1)$. Taking $c_0:=1$, (EST2) holds for $t=0$. We assume now that (EST2) holds for $t\le r-1\in S(0,b_{\al}-2)$. As $c_0\le c_1\le...\le c_{r-1}$ and as $\vph(e_{r-1})=e_r$, we have
$$p^{c_{r-1}-1}\bigl(M\cap (<e_1,...,e_r>[{1\over p}])\bigr)\subset <e_1,...,e_r>.\leqno (IND)$$ 
We now prove that (EST2) holds for $t=r$. We write  
$$\vph(e_r)=\sum_{i=1}^r a_ie_i+p^{n_r}x,\leqno (EQ)$$
with $x\in M\setminus pM$ and all $a_i$'s in $W(k)$. By our initial induction (on ranks), we can speak about $n(t,0,h_{\al})$. So let
$$
c_r:=c_{r-1}+b_{\al}(h_{\al}+c_{r-1})+2rb_{\al}!(h_{\al}+c_{r-1})^2\bigl(1+2{\rm max}\{n(r,0,j)|j\in S(0,h_{\al}+c_{r-1}-1)\}\bigr). \leqno (RECEQ)
$$
We assume $n_r\ge c_r+1$. As $\vph(e_r-\sum_{i=1}^{r-1}\sg^{-1}(a_{i+1})e_i)=p^{n_r}x+a_1e_1$, from (IND) we get that $a_1$ is not divisible by $p^{h_{\al}+c_{r-1}}$. We consider the Frobenius endomorphism $\vph_r$ of $M_r$ that takes $e_i$ into $e_{i+1}$ if $i\le r-1$ and which takes $e_r$ into $\sum_{i=1}^r a_ie_i$. As $a_1\neq 0$, by inverting $p$, $\vph_r$ becomes an isomorphism. So $(M_r,\vph_r)$ is a $\sg$-crystal whose $h$-number is the $p$-adic valuation of $a_1$ and so it is at most $h_{\al}+c_{r-1}$. Let 
$$q:=1+2r(b_{\al}-1)!(h_{\al}+c_{r-1}){\rm max}\{n(t,0,j)|j\in S(0,h_{\al}+c_{r-1}-1)\}.$$ 
\indent
By the very definition of $n(r,0,j)$, we know that 
$$DT(M_r,\vph_r)\le{{q-1}\over {2r(b_{\al}-1)!(h_{\al}+c_{r-1})}}.\leqno (IND_0)$$ 
\smallskip
Let $N_r$ be a $W(k)$-submodule of $M_r$ such that $M_r/N_r$ is annihilated by $p^{DT(M_r,\vph_r)}$ and $(N_r,\vph_r)$ is a Dieudonn\'e object of $p-\Mm(W(k))$. Writing $(N_r,\vph_r)$ as a direct sum of elementary Dieudonn\'e objects of $p-\Mm(W(k))$, we choose a simple factor $(N_r^1,\vph_r)$ of it such that the component $e_1^1$ of $p^{DT(M_r,\vph_r)}e_1$ in $N_r^1$ (w.r.t. to the resulting direct sum decomposition of $N_r$) is not divisible inside $N_r^1$ by $p^{1+DT(M_r,\vph_r)}$. 
\smallskip
From (EQ) we get that 
$p^{qa_{\al}}e_1=\vph^{b_{\al}q}(e_1)$ taken mod $p^{n_r}$ is congruent to $\vph_r^{qb_{\al}}(e_1)$. So (cf. also (IND)) $\vph_r^{qb_{\al}}(e_1^1)-p^{qa_{\al}}e_1^1\in p^{n_r-c_{r-1}-DT(M_r,\vph_r)}N_r^1$. 
\smallskip
On the other hand, as $(N_r^1,\vph_r)$ is an elementary Dieudonn\'e object of $p-\Mm(W(k))$, $\vph_r^{qb_{\al}}(e_1^1)$ mod $p^{n_r}$ is divisible inside $N_r^1$ by $p^{m_r}$, where 
$$m_r\in (\NN\cup\{0\})\cap [qb_{\al}u-DT(M_r,\vph_r)-r(h_{\al}+c_{r-1}-1),qb_{\al}u+DT(M_r,\vph_r)+r(h_{\al}+c_{r-1}-1)],\leqno (0)$$ 
with $u\in [0,h_{\al}+c_{r-1}-1]$ as the slope of $(N_r^1,\vph_r)$; so $vu\in\NN\cup\{0\}$ for some $v\in S(1,r!)$. 
\smallskip
As $n_r\ge c_r+1$, it is easy to see that $n_r-c_{r-1}-DT(M_r,\vph_r)$ is greater than $qb_{\al}u+DT(M_r,\vph_r)+r(h_{\al}+c_{r-1}-1)$ as well as than $qh_{\al}$ (and so then $qa_{\al}$). From this, $(IND_0)$ and the last two paragraphs we get: 
$$\abs{qa_{\al}-qb_{\al}u}\le {{q-1}\over {2r(b_{\al}-1)!(h_{\al}+c_{r-1})}}+r(h_{\al}+c_{r-1}-1).$$ 
A simple computation shows that $\abs{\al-u}<{1\over {b_{\al}!}}$. On the other hand, as $r\le b_{\al}-1$, from the existence of $v$ we get $b_{\al}!\abs{\al-u}\in\NN\cup\{0\}$. We conclude: $\al=u$. So $b_{\al}$ divides $a_{\al}$ and so $b_{\al}=1$. For $b_{\al}=1$ we have $r=0$ and this contradicts the inequality $r\ge 1$. So $n_r\le c_r$. As obviously $c_r$ depends only on $\al$, $h_{\al}$ and $c_{r-1}$, this ends the argument for the Claim. We conclude:
\medskip
{\bf Corollary.} {\it There is $d(\al,h_{\al})\in\NN$ depending only on $\al$ and $h_{\al}$ and such that $DT(M,\vph)\le d(\al,h_{\al}).$}
\medskip
Next we consider the case when $\abs{SL}=1$ and $m_{\al}$ is arbitrary. By induction on $m_{\al}\in\NN$ we show that
$$DT(M,\vph)\le (2m_{\al}-1)d(\al,h_{\al}).\leqno (EST3)$$
Using a short exact sequence
$0\to (M_1,\vph_1)\hookrightarrow (M,\vph)\twoheadrightarrow (M_2,\vph_2)\to 0$, with $\dim_{W(k)}(M_2)=b_{\al}$, we just need to remark the following 2 things:
\medskip
{\bf a)} The $h$-number of $(M_i,\vph_i)$ is at most the same as the $h$-number of $(M,\vph)$, $i\in S(1,2)$ (for $i=1$ this is obvious while for $i=2$ this follows by passing to duals of latticed isocrystals). 
\smallskip
{\bf b)} Any isogeny from an elementary Dieudonn\'e object $(M_0,\vph_0)$ of $p-\Mm(W(k))$ to $(M_2,\vph_2)$ factors, after multiplication with $p^{2({m_{\al}-1})d(\al,h_{\al})}$, through $(M,\vph)$.
\medskip
To check b), let $m_0\in M_0$ be such that $\vph_0^{b_{\al}}(m_0)=p^{a_{\al}}(m_0)$ and $M_0$ is $W(k)$-generated by $m_0$, $\vph_0(m_0)$,..., $\vph^{b_{\al}-1}_0(m_0)$. Let $m\in M$ be such that via the $W(k)$-epimorphism $M\twoheadrightarrow M_2$ is mapped into $m_0$; here we view $M_0$ as a $W(k)$-submodule of $M_2$. We have $\vph^{b_{\al}}(m)=p^{\al}m+\tilde m$, with $\tilde m\in M_1$. We need to show that, eventually after multiplication with $p^{2({m_{\al}-1})d(\al,h_{\al})}$, we can choose $m$ such that $\tilde m=0$. Using induction on $m_{\al}$, this boils down to the following obvious thing:
\medskip
{\bf Fact.} {\it For any $x\in p^{a_{\al}}M_0$, the equation $\vph_0^{b_{\al}}(y)=p^{a_{\al}}y+x$ has a solution $y\in M_0$}. 
\medskip
We come back to an arbitrary latticed isocrystal $(M,\vph)$; so we do not assume anymore that $s=0$. Let 
$$u:={\rm max}\{(2m_{s+\al}-1)d(s+\al,h_{\al})|\al\in SL\}.$$
\indent
{\bf Exercise.} Show that the cokernel of $q_{SL}$ is annihilated by $p^{2d_M\al_{\rm max}du+u}$. 
\medskip
{\bf Hint.} Let $N$ be a $W(k)$-submodule of $M$ such that $(N,p^s\vph)$ is a Dieudonn\'e object of $p-\Mm(W(k))$ and $M\cap N_0[{1\over p}]/N_0$ is annihilated by $p^u$, for any elementary Dieudonn\'e object of $p-\Mm(W(k))$ which is a simple factor of $(N,p^s\vph)$. If the statement of the Exercise is not true then there is $y\in M$ such that $p^{2d_M\al_{\rm max}du+u+v}y\in N\setminus pN$ for some $v\in\NN$. Using standard $W(k)$-bases of all such $(N_0,p^s\vph)$'s, procede by induction on the number $NR$ of non-zero coefficients of $p^{2d_M\al_{\rm max}du+u+v}y$ w.r.t. the $W(k)$-basis of $N$ obtained by the union of all these standard $W(k)$-bases (for this induction replace $2d_M$ by $2NR$).  
\medskip
From Exercise, (EST3), (EST1) and Corollary we get:
$$DT(M,\vph)\le 2d_M\al_{\rm max}du+se+2u.\leqno (TOTALEST)$$
So BP1 is a consequence of (TOTALEST); for the part pertaining to $g\in GL(M)(W(k))$, we need to point out that such $g$'s do not change the $s$-numbers or the $h$-numbers, while $m_{\al}$'s, $d$ and $e$ have upper bounds which depend only on $d_M$ and $h$. This ends the proof of BP1. 
\medskip
{\bf BP2.} {\it We assume $k=\bar k$. Let $a,b\in\ZZ$, $b\ge a$. Let $r\in\NN$. Then there is a smallest number
$$f(r,a,b)\in\NN\cup\{0\}$$ 
 such that the isomorphism class (and so also the Newton polygon) of any $\sg$-$\Ms$-crystal $(M,\vph,G)$ over $k$ whose underlying module is of rank at most $r$ and such that $(M,\vph)$ is an object of $p-\Mm_{[a,b]}(W(k))$, depends only on its Fontaine truncation mod $p^{f(r,a,b)}$. Similarly, if $a\ge 0$ then there is a smallest number
$$f^\prime(r,a,b)\in\NN\cup\{0\}$$ such that the isomorphism class of $(M,\vph,G)$ depends only on its weak truncation mod $p^{f^\prime(r,a,b)}$.}
\medskip
This is a consequence of BP0 and BP1; as a gross estimate we have: 
$$f^\prime(r,a,b)\le 1+b-a+2g(r^2,b-a,2(b-a)),$$
where for $(a_1,b_1,c_1)\in\NN\times\NN\cup\{0\}\times\NN\cup\{0\}$, 
$$g(a_1,b_1,c_1):={\rm max}\{n(j,s,h)|(j,s,h)\in S(1,a_1)\times S(0,b_1)\times S(0,c_1)\}.$$ 
The connection between these $f$'s and $f^\prime$'s is given (with $a\ge 0$) by the inequalities
$$f^\prime(r,a,b)\le f(r,a,b)+b+1\le f^\prime(r,a,b)+b+1;\leqno (WFTRUNC)$$
the factor $b+1$ is added here as the Fontaine truncation mod $p^{f(r,a,b)}$ of $(M,\vph,G)$ involves division of elements of the form $\vph(x)$ by $p^c$, with $c\in S(a,b)$; here $x\in M$. 
\medskip
{\bf BP2'.} {\it A $p$-divisible group over $\bar k$ of rank $r$ is uniquely determined by its truncation mod $p^{f^\prime(r,0,1)}$.}
\medskip
For $r\in\NN\cup\{0\}$, let
$$d(r)\in\NN\cup\{0\}$$
be the smallest number such that $DFT({\rm End}(M_1),\vph_1)\le d(r)$, for any $\sg$-crystal $(M_1,\vph_1)$ of rank $r$ which can be extended to a filtered $\sg$-crystal. We have (cf. BP0):
$$f^\prime(r,0,1)\le f(r,0,1)+1\le 3+2d(r).$$
\medskip
{\bf D. Global deformations.} By a global deformation of an object ${\got C}$ of $p-\Mm_{[a,b]}(W(k))$ (resp. of $p-\Mm\Mf_{[a,b]}(W(k))$) we mean a pair $(X,{\got C}_X)$, with $X$ a $W(k)$-scheme which locally in the Zariski topology is a pro-\'etale scheme over a smooth $W(k)$-scheme and which has a special fibre which is a geometrically connected, $AG$ $k$-scheme, and with ${\got C}_X$ an object of $p-\Mm^\nabla_{[a,b]}(X)$ (resp. of $p-\Mm\Mf_{[a,b]}^\nabla(X)$), such that the pull back of ${\got C}_X$ via a $k$-valued (resp. $W(k)$-valued) point of $X^\wedge$ is ${\got C}$. For our conventions on groupoids we refer to 2.2.4 F.
\medskip
{\bf BP3.} {\it We consider a global deformation $(X,{\got C}_X)$ of an object of $p-\Mm_{[a,b]}(W(k))$. Let $r$ be the local rank of its underlying $\Mo_{X^\wedge}$-module. 
\smallskip
{\bf a)} If the number of isomorphism classes of Fontaine truncations mod $p^{f(r,a,b)}$ of pull backs of ${\got C}_X$ through points of $X_k$ with values in an arbitrary algebraically closed field $k_1$ containing $k$ is finite and does not depend on the choice of $k_1$, then the number of isomorphism classes of pull backs of ${\got C}_X$ through $k_1$-valued points of $X_k$ is as well finite and does not depend on $k_1$. 
\smallskip
{\bf b)} In general, there is a $k$-groupoid $Y_k$ on $X_k$ such that the morphism $m_k:Y_k\to X_k\times X_k$ which is part of its definition is affine and of finite presentation and the pull backs of ${\got C}_X$ through two $k_1$-valued points $y_1$ and $y_2$ of $X_k$ are isomorphic iff the $k_1$-valued of $X_k\times_k X_k$ defined by the pair $(y_1,y_2)$ lifts to $m_k$.}
\medskip
{\bf Proof:} a) follows from BP2. To see b), we can assume $a\ge 0$. We consider the pull backs of ${\got C}$ to $X\times_{W(k)} X$ via the two projections on $X$ and the reduced $X_k\times_k X_k$-scheme $Y_{f^\prime(r,a,b)}$ parameterizing isomorphisms between their weak truncations mod $p^{f^\prime(r,a,b)}$. 
\smallskip
The construction of $Y_{f^\prime(r,a,b)}$ is standard. As we work mod $p^{f^\prime(r,a,b)}$, localizing and using pull back arguments we can assume $X$ is affine and of finite type over $W(k)$. Even more, we can assume $X_{W_{f^\prime(r,a,b)}(k)}$ is $W_{f^\prime(r,a,b)}(R/pR)$, with $R$ a smooth $W(k)$-algebra (this is as in 3.15.3 6)); so we can use arguments at the level of (suitable conjugacy classes) of matrices of $M_r(W_{f^\prime(r,a,b)}(R/pR))$ to get the existence of $Y_{f^\prime(r,a,b)}$. Here suitable refers to the fact that we need to conjugate with invertible matrices which have  a shape similar to the shape of $h$ of 3.13.7.8. BP2 implies: we can take $Y_k=Y_{f^\prime(r,a,b)}$. This ends the proof.  
\medskip
{\bf E. Remark.} We have a variant of BP3 in the relative context of 3.6.18.7.1 b): we just need to replace isomorphisms by inner isomorphisms. This is so as the part of the proof of D pertaining to matrices can be performed in the context of (suitable conjugacy classes) of matrices defined by suitable $W_{f^\prime(r,\tilde a,\tilde b)}(k)$-valued points of a (fixed) smooth subgroup of a $GL$-group (over an open subscheme of $X_{W_{f^\prime(r,\tilde a,\tilde b)}(k)}$).  
\medskip
{\bf F. The case of Shimura $p$-divisible groups.} We consider a Shimura $p$-divisible group $(D_X,(t_{\al})_{\al\in\Mj})$ over $X^\wedge$, with $X$ a regular, formally smooth $W(k)$-scheme having a connected special fibre. Let $r\in\NN\cup\{0\}$ be the rank of $D_X$. Let $\EE(D_X)$ be the object of $p-\Mm_{[0,1]}(X)$ obtained by forgetting the filtration of $\DD(D_X)$. As $X_k$ is connected, we speak (as in 2.2.9 6)) about $1_{\Mj}$-isomorphism classes of pull backs of $(*(D_X),(t_{\al})_{\al\in\Mj})$, with $*\in\{\DD,\EE\}$, via geometric points of $X_k$ (cf. also 3.6.18.7.3 A and 3.15.6 B). So we have the following particular combination of D and E:
\medskip
{\bf BP4.} {\it There is a $k$-groupoid $Y_k$ on $X_k$ such that the $X_k\times_k X_k$-scheme $Y_k$ is affine and of finite presentation and the pull backs of $(\EE(D_X),(t_{\al})_{\al\in\Mj})$ through two geometric points $y_1$ and $y_2$ of $X_k$ with values in the same algebraically closed field $k_1$ are isomorphic iff the point of $X_k\times_k X_k$ defined by the pair $(y_1,y_2)$ lifts to a $k_1$-valued point of $Y$.} 
\medskip
This is entirely the same as the proof of D: we just need to remark that we can define $Y_k$ using as well schemes of isomorphisms of truncations mod $p^{f^\prime(r,0,1)}$ of the pull backs of $D_{X_k}$ to $X_k\times_k X_k$ via the natural two projections of $X_k\times_k X_k$ on $X_k$.
\medskip
{\bf G. Some invariants.} Let $q\in\NN$. Let $(M,(F^i(M))_{i\in S(a,b)},\vph,G)$ be a filtered $\sg_{\FF_{p^q}}$-$\Ms$-crystal over $\FF_{p^q}$. So $(M,\vph)$ is an object of $p-\Mm_{[a,b]}(W(\FF_{p^q}))$. We assume the existence of a family of tensors $(t_{\al})_{\al\in\Mj}$ in $F^0(\Mt(M[{1\over p}])\otimes_{B(\FF_{p^q})} B(\FF))$ such that $\vph\otimes 1(t_{\al})=t_{\al}$, $\forall\al\in\Mj$, and $G_{B(\FF)}$ is the subgroup of $GL(M\otimes_{W(\FF_{p^q})} B(\FF))$ fixing $t_{\al}$, $\forall\al\in\Mj$. See 2.2.22 2) for the meaning of $Cl(M,\vph,G)$. The number of isomorphism classes of Fontaine truncations mod $p^{n}$ of $\sg_{\FF_{p^q}}$-$\Ms$-crystals over $\FF_{p^q}$ of the form $(M,g\vph,G)$, with $g\in W(\FF_{p^q})$, is at most $\abs{G(W_n(\FF_{p^q}))}$ and so it is finite (here $n\in\NN$). We have:
\medskip
{\bf BP5.} {\it We assume $k=\bar k$. The subset 
$$Cl_{q}(M\otimes_{W(\FF_{p^q})} W(k),\vph\otimes 1,G_{W(k)},(t_{\al})_{\al\in\Mj})$$
 of $Cl(M\otimes_{W(\FF_{p^q})} W(k),\vph\otimes 1,G_{W(k)},(t_{\al})_{\al\in\Mj})$ formed by inner isomorphism classes defined by $\sg_k$-$\Ms$-crystals with an emphasized family of tensors which are definable over $\FF_{p^q}$ is finite and its number of elements does not depend on $k$.}
\medskip
{\bf Proof:} Let $g\in G(W(k))$. As in 2.2.9 8) we construct a triple $(M_{\ZZ_p},G_{\ZZ_p},(t_{\al})_{\al\in\Mj})$ starting from $(M\otimes_{W(\FF_{p^q})} W(k),g(\vph\otimes 1),G_{W(k)},(t_{\al})_{\al\in\Mj})$. As $G_{\FF_p}$ is connected, from Lang's theorem we get that this triple does not depend on $g$; moreover, its extension to $W(\FF)$ is $1_{\Mj}$-isomorphic to $(M\otimes_{W(\FF_{p^q})} W(k),G_{W(k)},(t_{\al})_{\al\in\Mj})$. So BP5 follows from BP2 and the finiteness of the set $G_{\ZZ_p}(W_{f^\prime(d_M,a,b)}(\FF_{p^q}))$ (the ``indepence of $k$ is obvious; see also b) of 2.2.4 B).
\medskip
From BP5, the First Main Corollary of 1.6.4 follows. We refer to the numbers $D(r,d,q)$ of 1.6.4 as the classical Dieudonn\'e numbers. Similarly, we refer to 
$Cl_{q}(M\otimes_{W(\FF_{p^q})} W(k),\vph\otimes 1,G_{W(k)},(t_{\al})_{\al\in\Mj})$
as the level-$q$ Dieudonn\'e set of $Cl(M\otimes_{W(\FF_{p^q})} W(k),\vph\otimes 1,G_{W(k)},(t_{\al})_{\al\in\Mj})$. Its number of elements is referred as the level-$q$ Dieudonn\'e number of $Cl(M\otimes_{W(\FF_{p^q})} W(k),\vph\otimes 1,G_{W(k)},(t_{\al})_{\al\in\Mj})$ and is denoted by 
$$D_q(Cl(M\otimes_{W(\FF_{p^q})} W(k),\vph\otimes 1,G_{W(k)},(t_{\al})_{\al\in\Mj})).$$
\indent
As a very gross (though morally useful) estimate we have
$$D_q(Cl(M\otimes_{W(\FF_{p^q})} W(k),\vph\otimes 1,G_{W(k)},(t_{\al})_{\al\in\Mj}))\le p^{q\dim_{W(k)}(G)(f^\prime(d_M,a,b)-1)}\abs{G_{\ZZ_p}(\FF_{p^q})};\leqno (GREST)$$
the right hand side of (GREST) is $\abs{G_{\ZZ_p}(W_{f^\prime(d_M,a,b)}(\FF_{p^q}))}$. In particular, in the context of $p$-divisible groups we get
$$D(r,d,q)\le p^{qr^2(f^\prime(r,1,0)-1)}\abs{GL_r(\FF_{p^q})}\le p^{qr^2f^\prime(r,1,0)}.\leqno (GRESTP-DIV)$$
\medskip
{\bf H. $isom$-constant deformations.} The deformation (resp. the Shimura $p$-divisible group) of BP3 (resp. of BP4) is called $isom$-constant, if the $k$-morphism $Y_k\to X_k\times X_k$ is surjective. 
\medskip
{\bf Problem.} Given a Shimura $p$-divisible group over $W(k)$, construct explicitly (following the pattern of 3.6.11) global deformations of it which are $isom$-constant and uni plus quasi-versal.
\medskip
{\bf Example.} We consider a cyclic diagonalizable Shimura filtered $\sg$-crystal $(M,F^1,\tilde\vph_1,G)$. We use the notations of 3.13.7.6.3.1 but we do not assume $(M,\tilde\vph_1,G)$ is quasi-final. So $k=\bar k$. We assume there is a cycle $(\al_1,...,\al_s)$ of $\pi$ such that $u_s=-1$ is a negative special entry of $(u_1,...,u_s)$. Let $U_i$ be the $\GG_a$ subgroup of $G$ having ${\got g}_{\al_i}$ as its Lie algebra, $i\in S(1,s)$. Let $NU$ be the integral, closed subgroup of $G$ generated by $G_i$'s; it is included in the unipotent radical of the parabolic subgroup of $G$ having $W^0({\rm Lie}(G),\tilde\vph_1)$ as its Lie algebra. We write $U_s={\rm Spec}(W(k)[t])$ and we work with the Frobenius lift of $W(k)[[t]]$ which maps $t$ into $t^p$. We consider the $p$-divisible object 
$$M_{U_s}:=(M\otimes_{W(k)} W(k)[t],F^1\otimes_{W(k)} W(k)[t],u(\tilde\vph_1\otimes 1),(t_{\al})_{\al\in\Mj})$$ 
of $\Mm\Mf_{[0,1]}(U_s)$, with $u$ as the universal element of $G(U_s)=G(W(k)[t])$. We apply 3.6.18.6 to it. We deduce the existence of an $\NN$-pro-\'etale morphism $\tilde U\to U_s$ such that:
\medskip
{\bf 1)} the special fibre of $\tilde U$ is connected, has an image in ${U_s}_k$ which is open and has a unique $k$-valued point mapping into the origin of ${U_s}_k$;
\smallskip
{\bf 2)} the pull back $M_{\tilde U}$ of $M_{U_s}$ to $\tilde U^\wedge$ is a $p$-divisible object of $\Mm\Mf_{[0,1]}^\nabla(\tilde U^\wedge)$.
\medskip
As $u_s=-1$, we can assume the Kodaira--Spencer map naturally associated to $M_{\tilde U}$ is injective in each $k$-valued point of $\tilde U_k$. The pull back of $M_{\tilde U}$ through such a point, when viewed without filtration, is a Shimura $\sg_k$-crystal of the form (cf. the existence of Teichm\"uller lifts) $(M,h_1\tilde\vph_1,G)$, with $h_1\in U_s(W(k))$. It is easy to see that it is isomorphic to $(M,\tilde\vph_1,G)$, under an isomorphism defined by an element of $NU(W(k))$: as $u_s$ is a negative special entry we can apply entirely the stairs method of the proof of 3.6.15 B. More precisely, the case $s=1$ is trivial and if $s\ge 2$ then by induction on $j\in\NN$ we replace $h_j\tilde\vph_1$ by its (inner) conjugate $h_{j+1}\tilde\vph_1$ under an element (it is uniquely determined) of $U_{s-j}(W(k))$, where $h_{j+1}$ is a $W(k)$-valued point of $U_{s-j}$; here the indices of $U_i$'s are taken mod $s$.  We get: $\forall m\in\NN$, $h_{1+sm}$ is congruent mod $p^m$ to $1_M$. So (cf. BP4):
\medskip
{\bf Corollary.} {\it $(M_{\tilde U},(t_{\al})_{\al\in\Mj})$ is an $isom$-constant deformation of $(M,\tilde\vph_1,G)$ whose Kodaira--Spencer maps in $k$-valued points (resp. in the $k$-valued point of 1) are injective (resp. can be identified naturally with ${\rm Lie}(U_s)/p{\rm Lie}(U_s)$).}
\medskip
{\bf Remarks.} {\bf 1)} This example is the generalization of the well known classical case when $s=1$ and $u_s=-1$. The simplest case of general nature covered by this example is: the case when all $u_i$'s are non-positive and $u_s=-1$.
\smallskip
{\bf 2)} The condition that $(M,\tilde\vph_1,G)$ is cyclic diagonalizable can be significantly weaken (we just need to get a ``circular property" on the Lie algebras of $\GG_a$ subgroups $U_i$'s of $G$ as above); in particular, this always applies if the Lie stable $p$-rank of the Shimura adjoint Lie $\sg$-crystal attached to $(M,\tilde\vph_1,G)$ is positive (cf. 1) and 3.4.5.1 B).
\medskip
{\bf I. Minimal degrees of definition.} We take $k$ arbitrary. By a particular object (with tensors) of $p-\Mm(W(k))$ we mean a non-zero object ${\got C}$ (with tensors) of $p-\Mm(W(k))$ whose extension ${\got C}_{\bar k}$ to $\bar k$ is definable over a finite field. For instance:
\medskip
-- any cyclic diagonalizable Shimura $\sg$-crystal is particular (cf. 2.2.16.4 a));
\smallskip
-- if $k\subset\FF$ then all objects (with tensors) of $p-\Mm(W(k))$ are particular (this follows from the proof of BP5).  
\medskip
By the minimal degree of definition of a particular object ${\got C}$ (with tensors) of $p-\Mm(W(k))$ we mean the smallest $q\in\NN$ such that ${\got C}_{\bar k}$ is definable over $\FF_{p^q}$. Lemma of 3.11.8 implies that the minimal degree of definition of any Shimura-ordinary $\sg$-crystal is the same as its degree of definition. Similarly, by adapting the proof of the mentioned Lemma in the context of 2.2.16.5 we get: the minimal degree of definition of any Shimura-ordinary adjoint Lie $\sg$-crystal is the same as its $A$-degree of definition.
\medskip
{\bf J. Remarks.} {\bf 1)} In the terminology to be introduced in [Va11], BP4 can be restated as: there is a reduced almost stack of $p$-divisible groups of fixed rank over $\FF_p$-schemes. It is easy to see that it is of finite type (for instance, cf. BP1 or cf. Theorem 13 of 14). Here by ``reduced almost" we want to emphasize that we are dealing only with geometric points (as BP4 is stated just in terms of such points).
\smallskip
{\bf 2)} Simple examples (like 3.9.7.3) point out that in general, with $k=\bar k$, we have 
$${\rm lim}_{q\to\infty} D_q(Cl(M\otimes_{W(\FF_{p^q})} W(k),\vph\otimes 1,G_{W(k)},(t_{\al})_{\al\in\Mj}))=\infty.$$ 
\indent
{\bf K.} For future references we state explicitly the following combined versions of BP0 and of the first two inequalities of (INEQ) of rm. 0) of B. We assume $k=\bar k$.
\medskip
{\bf Corollary.} {\it {\bf a)} The Dieudonn\'e--Fontaine torsion of any cyclic diagonalizable Shimura Lie $\sg$-crystal is at most the maximum of the deformation dimension of its cyclic factors. 
\smallskip
{\bf b)} The $isom$-deviation of any object of $p-\Mm_{[a,b]}(W(k))$ which has a lift which is cyclic diagonalizable and whose underlying $W(k)$-module has rank $r$ is at most 
$$1+(b-a)(2r^2+1).$$} 
\indent
{\bf Examples.} We refer to 3.13.7.6.3.1. Then under the assumption of its Fact we have $DFT({\rm Lie}(G),\tilde\vph_1)=1$ (cf. ineq. two of (INEQ) of rm. 0) of B). So (cf. BP0) the $isom$-deviation of $(M,\tilde\vph_1,G)$ is at most $4$. For instance, this applies if $(M,\tilde\vph_1)$ is a direct sum of a finite number of copies of the cyclic diagonalizable $\sg$-crystal attached to the type $(1,0,0,...,0)$ (or to its dual); however, it is easy to see that in this case the $isom$-deviation of $(M,\tilde\vph_1,G)$ is in fact $1$. 
\smallskip
If moreover ${\rm Lie-mod}(G)$ is a Lie-algebra, then we can replace $4$ by $3$. This is so as the Case 2 of 3.13.7.6.3 can be entirely adapted to the context of ${\rm Lie-mod}(G)$. Warning: in many cases b) can be significantly improved (to be compared with a)).
\medskip
{\bf 3.15.8. New proof of the specialization theorem.}
The estimates of 3.15.7 are different (though close in nature) to the ones of [Ka2, 1.4-5]. No doubt, the estimates of loc. cit. are much sharper. However, the estimates of 3.15.7 give some extra information (insight) and in general they can be considerably refined. These refinements will not be presented here: to us, it seems far more useful to present a proof of the specialization theorem which uses just very simple (i.e. the most basic) estimates. We recall (see [Gr]; see also [Dem, p. 91] and [Ka2]) that Grothendieck proved that under specialization the Newton polygons of $F$-crystals ``go up", while Katz added (see [Ka2, 2.3.1-2]) that in fact for any $F$-crystal ${\got C}$ over an integral $\FF_p$-scheme $S$ there is an open, dense subscheme $U$ of $S$ such that the Newton polygons of pull backs of ${\got C}$ through geometric points of $U$ are all the same. Below we concentrate on reproving the existence of $U$.
\smallskip
We start by presenting a homomorphism form of 3.6.15 B. We assume $k=\bar k$. Let $(M_1,\vph_1)$ and $(M_2,\vph_2)$ be two $\sg$-crystals over $k$. Let $d_i:=DT({\rm End}(M_i),\vph_i)$ and let $s_i$ be the $s$-number of $({\rm End}(M_i),\vph_i)$, $i=\overline{1,2}$. Let $s$ be the $s$-number of $({\rm End}(M_1\oplus M_2),\vph_1\oplus\vph_2)$. Let $d:=DT({\rm End}(M_1\oplus M_2),\vph_1\oplus\vph_2)$. 
\smallskip
Let $q\in\NN\cup\{0\}$. From the proof of 3.6.15 B we get: that the images of the following two reduction group homomorphisms
$${\rm Aut}(M_1/p^{2d_1+1+s_1+q}M_1,\vph_1)\to {\rm Aut}(M_1/p^{d_1+1+s_1+q}M_1,\vph_1)$$
and
$${\rm Aut}(M_1,\vph_1)\to {\rm Aut}(M_1/p^{d_1+1+s_1+q}M_1,\vph_1)$$
are the same. Passing from automorphisms to homomorphisms in the standard way we get 1) of the following Corollary.
\medskip
{\bf Corollary (the homomorphism property).} {\it {\bf 1)} Let $\tilde s:=s+q+d+1$. The following $\ZZ_p$-linear reduction maps
$${\rm Hom}((M_1/p^{d+\tilde s}M_1,\vph_1),(M_2/p^{d+\tilde s}M_2,\vph_2))\to {\rm Hom}((M_1/p^{\tilde s}M_1,\vph_1),(M_2/p^{\tilde s}M_2,\vph_2))$$
and
$${\rm Hom}((M_1,\vph_1),(M_2,\vph_2))\to {\rm Hom}((M_1/p^{\tilde s}M_1,\vph_1),(M_2/p^{\tilde s}M_2,\vph_2))$$
have the same image.
\smallskip
{\bf 2)} There are numbers $n$, $e\in\NN$ which depend only on the rank $r_i$ of $M_i$ and of the maximal Hodge number $h_i$ of $(M_i,\vph_i)$, $i=\overline{1,2}$, such that for $\tilde n:=n+q$, the images of the two natural restriction maps
$${\rm Hom}((M_1,g_1\vph_1),(M_2,g_2\vph_2))\to {\rm Hom}((M_1/p^{\tilde n}M_1,g_1\vph_1),(M_2/p^{\tilde n}M_2,g_2\vph_2))$$ 
and
$${\rm Hom}((M_1/p^{\tilde n+e}M_1,g_1\vph_1),(M_2/p^{\tilde n+e}M_2,g_2\vph_2))\to {\rm Hom}((M_1/p^{\tilde n}M_1,g_1\vph_1),(M_2/p^{\tilde n}M_2,g_2\vph_2))$$
are the same, $\forall g_i\in GL(M_i)(W(k))$, $i=\overline{1,2}$.
\smallskip
{\bf 3)} Let $D_1$ and $D_1$ be two $p$-divisible groups over $k$. There is $n\in\NN$ which depends only on the rank $r_i$ of $D_i$, $i=\overline{1,2}$, such that a homomorphism $D_1[p^{n+q}]\to D_2[p^{n+q}]$ lifts to a homomorphism $D_1\to D_2$ iff it lifts to a homomorphism $D_1[p^{n+q+d(r_1+r_2)}]\to D_2[p^{n+q+d(r_1+r_2)}]$.}
\medskip
{\bf Proof:} 2) and 3) are a consequence of 1) and of BP1 (see end of C for the meaning of $d(r_1+r_2)$).
\medskip
We are now ready to proof the existence of $U$. Let $K$ be the field of fractions of $S$ and let ${\got C}_{\bar K}=(M,\vph)$ be the pull back of ${\got C}$ to $\bar K$. We consider an isogeny 
$$i:(M_1\otimes_{W(\FF_p)} W(\bar K),\vph_1\otimes 1)\to (M,\vph),$$ 
with $(M_1,\vph_1)$ a Dieudonn\'e object of $p-\Mm(W(\FF_p))$. Let $h_1$ be the $h$-number of $(M_1,\vph_1)$. We take $q\in\NN$ such that we have $p^qM\subset i(M_1)$. Localizing, we can assume that the $h$-numbers of pull backs of ${\got C}$ through geometric points of $S$ are all bounded above by the same number $h_{\got C}$. Let $m:=\dim_{W(\bar K)}(M)$. Let 
$$f:=1+\max\{h_1,h_{\got C}\}+2g(4m^2,\max\{h_1,h_{\got C}\},2\max\{h_1,h_{\got C}\})$$ 
(cf. the paragraph after BP2). The $s$-number (resp. the $h$-number) of the $End$ of the direct sum of $(M_1\otimes_{W(\FF_p)} W(k_1),\vph_1\otimes 1)$ with the pull back of ${\got C}$ through any $k_1$-valued point of $S$ is at most $\max\{h_1,h_{\got C}\}$ (resp. at most $2\max\{h_1,h_{\got C}\}$); here $k_1$ is an arbitrary algebraically closed field containing $\FF_p$. Let $\tilde q\in\NN$ be greater or equal to $e$ obtained in 2) for $h_2:=h_{\got C}$ and $r_2=r_1=m$.
\smallskip
We consider the reduction $i(q+\tilde q+f)$ of $i$ mod $q+\tilde q+f$; we view it as a morphism between (non-filtered) $W(\bar K)$-modules endowed with connections and Frobenius endomorphisms subject to axioms similar to the ones of $\Mm\Mf^\nabla(W(\bar K))$. The connections are of course initially (i.e. over $W(\bar K)$) trivial; though in what follows, they do not play any relevant role, we still feel appropriate to keep track of them. $i(q+\tilde q+f)$ is defined over a finite field extension $K_1$ of $K$ (even if $K$ does not have a finite $p$-basis). Replacing $S$ by its normalization $S_1$ in $K_1$, as the continuous map underlying the morphism $S_1\to S$ is proper, we can assume $K_1=K$. Localizing, we deduce the existence of an affine, open, non-empty subscheme $U={\rm Spec}(\bar R)$ of $S$ such that we have a morphism
$$i(q+\tilde q+f)_U:(M_1/p^{q+\tilde q+f}M_1,\vph_1)_U\to {\got C}_U/p^{q+\tilde q+f}{\got C}_U,$$
with the lower index $U$ meaning the pull back to $U$; we view it (cf. also [Be, 1.6.4 of p. 246 and 1.2.1 of p. 91]) as being a morphism between truncations mod $p^{q+\tilde q+f}$ of objects of $p-\Mm\Mo\Md^\nabla(W(\bar R))$. This last category is defined similarly to $p-\Mm^\nabla(W(\bar R))$ of 2.2.1.7 4) but this time we do not assume:
\medskip
-- that the $\Mo_U$-sheaf $\Om_{\bar R/k}$ is locally free of locally finite rank; 
\smallskip
-- or that we can ``put" some adequate filtrations on the underlying sheaves of its objects (this conforms to 1) of 3.15.7 B). 
\medskip
We can assume the cokernel of $i(q+\tilde q+f)_U$ is annihilated by $p^q$. From 2) we deduce that for any geometric point $y_1:{\rm Spec}(k_1)\to U$, the truncation mod $p^{q+f}$ of $y_1^*(i(q+\tilde q+f)_U)$ lifts to a morphism 
$$i_{y_1}:(M_1\otimes_{W(\FF_p)} W(k_1),\vph_1\otimes 1)\to y_1^*({\got C})$$ 
of $\sg$-crystals. As the cokernel of its reduction mod $p^{q+\tilde q+f}$ is annihilated by $p^q$, we get $i_{y_1}$ is injective and so, by reasons of ranks, it is an isogeny. So the Newton polygon of $y_1^*({\got C})$ does not depend on $y_1$. This ends the proof of the existence of $U$.
\medskip
{\bf 3.15.9. Three main specialization approaches.} Here, in connection to the specialization theorem, we include three main specialization approaches; not to make this paper too long, we just list their principles without performing any explicit computations (see \S 9-10 for such computations). We assume $k=\bar k$; however, everything else except the application of 3.15.7 BP2 to be presented below holds for an arbitrary $k$. We use the notations of 3.13.7.1. So $(M,F^1,\vph,G,(t_{\al})_{\al\in\Mj})$ is a Shimura filtered $\sg$-crystal with an emphasized family of tensors over $k$, $\mu$ is a cocharacter of $G$ defining its filtration class, $\vph_1:=g^{-1}\vph$ where $g\in G(W(k))$ is such that $(M,F^1,\vph_1,G)$ is a Shimura-canonical lift, $G_{\ZZ_p}$ is a $\ZZ_p$-structure of $G$ defined naturally by $(M,F^1,\vph_1,G)$, while $N^+$, $P^0$ and $P^-$ are subgroups of $G$ defined naturally by $\mu$. $\sg$ acts on $W(k)$-valued points of $G$ via this $\ZZ_p$-structure. Let $n:=1+f(d_M,1,0)$. Let 
$$
H_{W(k)}:=N^+\times_{W(k)} P^0\times_{W(k)} N^-.
$$
Its special fibre is $H$ as defined in 3.13.7.1. Let $\TT$ have the same significance as in 3.13.7.1. For $\tilde g\in G(W(k))$ we denote by $\tilde g(n)$ its reduction mod $p^n$. Let $g_1\in G(W(k))$ be such that $g(1)$ factors through the Zariski closure of the orbit $o_1$ of $g_1(1)$ under $\TT$ in $G_k$. 
\smallskip
We consider the morphism 
$$m(n):H_{W_n(k)}\to G_{W_n(k)}$$ 
which at the level of $W(k)$-valued points takes a triple 
$$(h_1,h_2,h_3)\in N^+(W_n(k))\times P^0(W_n(k))\times N^-(W_n(k))$$ 
into the $W_n(k)$-valued point of $G$ defined by 
$$h_1h_2h_3g_1(n)\sg(h_3^{-1})\sg(h_2^{-1})\sg(ph_1^{-1})g^{-1}(n);$$
$N^+$ is a direct sum of copies of $\GG_a$ and so it makes sense to multiply its elements by $p$; this defines $\sg(ph_1^{-1})$. The image of $m(n)_k$ is $t_1:=o_1g^{-1}(1)$. We consider a $W_n(k_1)$-valued point $z_n^0$ of $H_{W_n(k)}$, with $k_1$ a field (not necessarily perfect), such that:
\medskip
a) the $k_1$-valued point of $G$ defined by $m(n)\circ z_n^0$ mod $p$ sits over the generic point $s$ of $t_1$;
\smallskip
b) the $k_1$-valued point of $H$ we get is a closed point $s_h$ such that the resulting field extension $k_0\hookrightarrow k_1$ at the level of residue fields of $s$ and $s_H$ is totally inseparable and of finite degree; let $m$ be its degree.
\medskip
For the totally inseparable part of b) we just have to recall that $\TT$ is a group action (if the stabilizer subgroup of $g_1(1)$ under $\TT$ is smooth, then we can take $m=1$). Let $z_n$ be the composite of $\sg_{k_1}^{mn}$, viewed as an endomorphism of ${\rm Spec}(W(k_1))$, with $z_n^0$. $m(n)\circ z_n$ factors through ${\rm Spec}(W_n(k_0))$ (cf. b)); let $w_n:{\rm Spec}(W_n(k_0))\to G_{W_n(k)}$ be this factorization. Let $\tilde k:=\overline{k_0}$.
\smallskip
$s$ specializes to the origin of $G_k$. So we consider an affine, open subscheme $U$ of $G$ through which the origin of $G$ factors and such that we have an \'etale morphism from it into $\AA_{W(k)}^{\dim_{W(k)}(G)}$. Identifying $W_n(U_k)$ with $U_{W_n(k)}$ (this is as in 3.15.3 6)), often it is possible to show the existence of a Frobenius lift $\Phi_U$ of $U^\wedge$ such that its reduction $\Phi_U(n)$ mod $p^n$ is compatible with the factorization (still denoted by $w_n$) of $w_n$ through $U$, i.e. we have 
$$w_n\circ \sg_{k_0}(n)=\Phi_U(n)\circ w_n,\leqno (COMP)$$ 
where $\sg_{k_0}(n)$ is the Frobenius endomorphism of ${\rm Spec}(W_n(k_0))$. In what follows we assume such an $\Phi_U$ exists. 
\smallskip
We apply a variant of 3.6.1.3 in the context of $(M,F^1,\vph,G,(t_{\al})_{\al\in\Mj})$ and of $(U,\Phi_U)$; the difference from 3.6.1.3, 3.6.10 or 3.6.18.6 is: $\Phi_U$ is not necessarily of essentially additive type in the origin of $U_k$. We deduce the existence of a Shimura $p$-divisible group $(\Md,(t_{\al})_{\al\in\Mj})$ over the $p$-adic completion of an $\NN$-pro-\'etale scheme $U_1$ over $U$ which is universal in some specific sense. We do not need this universality. What we need: $U_{1k}$ is an $AG$ $k$-scheme and so each $\tilde k$-valued point of $U_k$ lifts to a $\tilde k$-valued of $U_{1k}$; here $\tilde k$ is an algebraically closed field containing $k$. We consider an arbitrary $k$-valued point $y_1$ of $U_{1k}$ lifting the origin of $U_k$. We get:
\medskip
{\bf 1)} The Shimura $\sg$-crystal associated to $y_1^*(\Md,(t_{\al})_{\al\in\Mj})$ is inner isomorphic to the Shimura $\sg$-crystal $(M,g_2\vph,G,(t_{\al})_{\al\in\Mj})$, with $g_2\in G(W(k))$ congruent to the identity mod $p$.
\medskip
{\bf 2)} The Shimura $\sg_{\tilde k}$-crystal associated to the pull back of $(\Md,(t_{\al})_{\al\in\Mj})$ through a $W(\tilde k)$-valued point $y_2$ of $U_1$ such that its special fibre specializes to $y_1$ and the $W_n(\tilde k)$-valued point of $G$ it defines naturally factors through $w_n$, is isomorphic to $(M\otimes_{W(k)} W(\tilde k),g_3(g_1\vph\otimes 1),G_{W(\tilde k)},(t_{\al})_{\al\in\Mj})$, with $g_3\in G(W(\tilde k))$ congruent to the identity mod $n-1=f(d_M,1,0)$.
\medskip
2) is a consequence of (COMP) and of the construction of $m(n)$ (we recall that $s$ specializes to the origin of $U_k$). Using 3.4.15 BP2' we get:
\medskip
{\bf 3)} In 2) we can assume $g_3$ is $1_M$.
\medskip
{\bf Remarks. 1)} This forms the going down form of the inductive specialization approach. At least if $m=1$ there are variants mod $p^n$ of this approach. As we do not stop to analyze when $\Phi_U$ indeed exists, there is no point to include these variants here. 
\smallskip
{\bf 2)} There is as well a going up form of the inductive specialization approach; in essence it is nothing else but the combination of 3.6.1.3 with the touching property of 3.6.18.4.3 and with the possibility of concretely constructing (in many situations; see H) $isom$-constant deformations.
\smallskip
{\bf 3)} There is as well a parallel approach: this is entirely similar to the algebraization process of the proof of 3.6.18.4.1. So using either 3.6.18.4.3 or the Fact of 3.6.18.9 as well as the Chinese Reminder Theorem (as in the mentioned) proof, one tries to connect two Shimura $\sg$-crystals through a global deformation which has some specific properties like (it is $isom$-constant, or has a constant Newton polygon, etc.).      
\medskip
{\bf 3.15.10. The second (Newton polygon) form of the purity principle.} 3.6.15 B was obtained independently of [dJO, \S 2-4]. However loc. cit. and the estimates of 3.15.7 share one common think: they both use Dieudonn\'e's objects of $p-\Mm(W(k)$ (i.e. Dieudonn\'e's classification); but we would like to point out that we think it is more appropriate to use in general (in connect to 3.15.7-10 and loc. cit.) Dieudonn\'e--Fontaine's objects instead of Dieudonn\'e's objects (of $p-\Mm(W(k)$). Moreover, the ideas of 3.6.15 B (in the form of the Corollary of 3.15.8) can be combined with [Ka2, 2.7.4] to get slightly better results (with much simpler proofs) than what [dJO, 4.1] gives. The goal of this section is to prove the following Theorem.
\medskip
{\bf Theorem.} {\it Let $S$ be an integral $\FF_p$-scheme. Let ${\got C}$ be an $F$-crystal on $S$. We assume there is $b_0\in\NN$ such that the $h$-number of all pull backs of ${\got C}$ through geometric points of $S$ are all bounded above by $b_0$. Let $U$ be the maximal open, dense subscheme of $S$ such that the Newton polygons of pull backs of ${\got C}$ through geometric points of $U$ are all the same (see 3.15.8). We also assume that, locally in the Zariski topology of $S$, the normalization $U^n$ of $U$ in the field of fractions of $S$ is regular in codimension 1 and any dominant, normal, affine $U^n$-scheme which is finite above spectra of local rings of $U^n$ which are DVR, is the normalization of $U^n$ in a finite field extension of the field of fractions $K$ of $S$. Then $U$ is an affine $S$-scheme.}
\medskip
{\bf Proof:} Let $r$ be the rank of ${\got C}$. We can assume $S={\rm Spec}(R)$ is affine and (cf. Exercise of 3.6.8.1.4) normal. So $U=U^n$. To prove that $U$ is affine, we are allowed (as the conditions on $U$ allow us) to replace $S$ by its normalization in any finite field extension of $K$. Let ${\got C}_0$ be a Dieudonn\'e object of $p-\Mm(\FF_p)$ such that we have an isogeny ${\got C}_{0\bar K}\hookrightarrow {\got C}_{\bar K}$. Let $b$ be the maximum of $b_0$ and of the maximal slope of ${\got C}_{0}$. Let 
$$n:=1+b+2g(4r^2,b,2b)$$ 
be as in the paragraph after 3.15.7 BP2. Its first role is: the Newton polygons of pull backs of ${\got C}$ through geometric points of $S$ depend only on their truncations mod $p^n$ (cf. BP0-1; here the role of $b$ is the same as of $\max\{h_1,h_{\got C}\}$ in 3.15.8). 
In what follows we take $m\in\NN$, $m\ge n$ (to be specified at the right time). 
\smallskip
Warning: till the end of the proof, in order to fully accommodate the context of non-perfect $\FF_p$-schemes, we use the terminology $F$-crystals in the sense of (truncations) of objects of $p-\Mm\Mo\Md^\nabla(W(*))$, with $*$ standing for $\FF_p$-algebras or affine $\FF_p$-schemes (this new terminology, is not an impediment from the point of view of Newton polygons as they can be computed via pull backs through geometric points); $p-\Mm\Mo\Md^\nabla(W(*))$ is interpreted as in 3.15.8 (so its objects are locally free of locally finite ranks $\Mo_{{\rm Spec}(W(*))}$-sheaves endowed with integrable, nilpotent mod $p$ connections and with Frobenius endomorphisms which satisfy the logical axiom connecting the Frobenius endomorphisms with the connections). We use right lower indices by schemes for pull backs of $F$-crystals.
\smallskip
Localizing and using Witt rings, we can assume there is a finitely generated $\FF_p$-subalgebra $R_0$ of $R$ such that ${\got C}/p^m{\got C}$ is obtained as the pull back to $S$ of an $F$-crystal ${\got C}[m]_0$ on $S_0:={\rm Spec}(R_0)$ in coherent sheaves annihilated by $p^m$. Warning: here we already use the convention of the previous paragraph.
\smallskip 
Using [Ka2, 2.7.4] we deduce that for any DVR, local ring $V$ of $U$, there is a perfect, faithfully flat $V$-algebra $V_1$ such that we have an isogeny $i_{V_1}:{\got C}_{0V_1}\to {\got C}_{V_1}$ of $F$-crystals. Let $n_V\in\NN$ be the smallest number such that the inclusion $p^{n_V}{\got C}_{V_1}\hookrightarrow {\got C}_{V_1}$ factors through $i_{V_1}$. 
\medskip
{\bf Claim 1.} {\it There is $q\in\NN$, depending only on ${\got C}_0$ and not on $V$, such that we can assume $n_V\le q$.}
\medskip
The proofs of [Ka2, 2.6.1-2 and 2.7.1] (they are forming the ingredient needed to get [Ka2, 2.7.4]) imply that there is $\tilde q\in\NN$ which does not depend on $V$ ($\tilde q$ depends only on the Newton polygon of ${\got C}_0$) such that we have an isogeny ${\got C}^V_{0V_1}\to {\got C}_{V_1}$ whose cokernel is annihilated by $p^{\tilde q}$, with ${\got C}_0^V$ an $F$-crystal on $\FF_p$. But the $h$-number (resp. the $s$-number) of ${\got C}_0^V$ (resp. of $End({\got C}_0^V)$) is at most $b+\tilde q$ (resp. at most $b+2\tilde q$). So Claim 1 follows from 3.15.7 BP1. 
\smallskip
We can assume $q$ is greater than $2b$ and $2(c+1)$, where $c\in\NN$ is such that ${\got C}_0$ is an object of $p-\Mm_{[0,c]}(W(k))$. We take (for instance)
$$m:=7q+2n+1.$$ 
Let 
$$i_{V_1}(m):{\got C}_{0V_1}/p^m{\got C}_{0V_1}\to {\got C}_{V_1}/p^{m}{\got C}_{V_1}$$ 
be the truncation mod $p^m$ of $i_{V_1}$. Using standard arguments of algebraic geometry we deduce the existence of a finite field extension $K_V$ of $K$ such that we can assume we have a similar morphism 
$i_{U_1(V)}(m)$, with ${\rm Spec}(U_1(V))$ as an open subscheme of the normalization of $U$ in $K_V$ through which ${\rm Spec}(V_1)\to U$ factors, whose cokernel is annihilated by $p^q$ and which is the truncation mod $p^m$ of a morphism ${\got C}_{0U_1(V)}/p^{m+2q}{\got C}_{0U_1(V)}\to {\got C}_{U_1(V)}/p^{m+2q}{\got C}_{U_1(m)}$. 
But the number of morphisms
${\got C}_{0\bar K}/p^m{\got C}_{0\bar K}\to {\got C}_{\bar K}/p^{m}{\got C}_{\bar K}$ which are the truncation mod $p^m$ of morphisms ${\got C}_{0\bar K}/p^{m+n}{\got C}_{0\bar K}\to {\got C}_{\bar K}/p^{m+n}{\got C}_{\bar K}$ is finite (this follows from the fact that any such morphism lifts to a morphism ${\got C}_{0\bar K}\to {\got C}_{\bar K}$, cf. 3.15.8 1) and the definition of $n$). So we can assume $K_V$ does not depend on $V$ and so we can assume $K_V=K$. 
\smallskip
We reached the following situation. There is a finite number of open subschemes $U_i$, $i\in I$, of $U$ such that:
\medskip
1) each point of $U$ of codimension $1$ belongs to some $U_i$;
\smallskip
2) for each $i\in I$, we have a morphism $i_{U_i}(m):{\got C}_{0U_i}/p^m{\got C}_{0U_i}\to {\got C}_{U_i}/p^{m}{\got C}_{U_i}$, whose cokernel and kernel are annihilated by $p^q$.
\medskip
What follows next is a method of ``gluing" $i_{U_i}(m)$'s, $i\in I$ (which we hope to be useful in other situations). It is enough to ``glue" these morphisms generically. So let
$$l_i(m):{\got C}_{0\bar K}/p^m{\got C}_{0\bar K}\to {\got C}_{\bar K}/p^m{\got C}_{\bar K},$$
be the morphism naturally defined by $i_{U_i}(m)$, $i\in I$. From 1) of the Corollary of 3.15.8 we get that there is $l_i:{\got C}_{0\bar K}\to {\got C}_{\bar K}$ such that its reduction mod $p^{m-n}$ lifts the truncation mod $p^{m-n}$ of $l_i(M)$. $l_i$'s are isogenies and in fact their cokernels are still annihilated by $p^q$. Let $l:{\got C}_{0\bar K}\to {\got C}_{\bar K}$ be an isogeny such that:
\medskip
3) its image lies inside the intersection of the images of all $l_i$'s ($i\in I$);
\smallskip
4) its cokernel is annihilated by $p^{2q}$.
\medskip
For 4) we need to point out that we can take as $l$ any $p^ql_i(m)$, with $i\in I$. Let $s_i$ be the factorization of $l$ through $l_i$. As in above part referring to $K_V$'s, we can assume that the truncations mod $p^{m-q}$ of $l$ and of all $s_i$'s are all definable over $K$ and not only over $\bar K$. As $s_i$'a are definable over finite fields, we can assume as well that $s_i$ extends (uniquely) to a morphism 
$s_{U_i}:{\got C}_{0U_i}/p^{m-q}{\got C}_{0U_i}\to {\got C}_{U_i}/p^{m-q}{\got C}_{U_i}$. So we can glue the composite of $s_{U_i}$ with the truncation mod $p^{m-q}$ of $i_{U_i}(m)$, $i\in I$, to get a morphism
$$i_{U_0}(m-q):{\got C}_{0U_0}/p^{m-q}{\got C}_{0U_0}\to {\got C}_{U_0}/p^{m-q}{\got C}_{U_0},$$
with $U_0$ an open subscheme of $U$ such that the codimension of $U\setminus U_0$ in $U$ is at least 2, whose kernel and cokernel are annihilated by $p^{2q}$.
\smallskip
We consider the $S$-scheme $S_1$ parameterizing morphisms from ${\got C}_{0S}/p^{m-2q}{\got C}_{0S}$ to ${\got C}/p^{m-2q}{\got C}$. It is an affine $S$-scheme, locally of finite presentation. Let $U^\prime$ be the normalization of the Zariski closure of $U_0$ in $S_1$ defined naturally by $i_{U_0}(m-q)$ and let $i_{U^\prime}(m-q)$ be the natural extension of $i_{U_0}(m-q)$ to it. From our hypotheses we get: $U$ is naturally an open subscheme of $U^\prime$. So the Theorem follows from the following Claim.
\medskip
{\bf Claim 2.} {\it $U$ is an affine subscheme of $U^\prime$.}
\medskip
{\bf Proof:} We can assume $U^\prime$ is an open subscheme of $S$ (otherwise we replace $S$ by $U^\prime$). We can assume $U^\prime$ and $i_{U^\prime}(m-q)$ are obtained by pull backs from an open subscheme $U^\prime_0$ of $S_0$ and respectively from a morphism 
$$i_{U_0^\prime}(m-q):{\got C}_{0U_0^\prime}/p^{m-q}{\got C}_{0U_0^\prime}\to {\got C}[m]_0/p^{m-q}{\got C}[m]_0.$$ 
\smallskip
We can assume we have a similar morphism
$$\tilde i_{U_0}(m-q):{\got C}^t_{0U_0}/p^{m-q}{\got C}^t_{0U_0}\to {\got C}_{U_0}^t/p^{m-q}{\got C}^t_{U_0},$$
whose cokernel is annihilated by $p^{2q}$ (by shrinking, we can assume we have the same $U_0$). Here the right upper index $t$ means (taking) the dual. Strictly speaking we have to perform first Tate twists: we need to tensor all these duals by the pull back of $\FF_p(q)$ to $S$; not to introduce extra notations we do not mention it. So, starting from it, we similarly define $\tilde U^\prime$, $\tilde U^\prime_0$. Moreover, we similarly get that $\tilde U^\prime$ is affine and that $U$ is naturally an open subscheme of $\tilde U^\prime$. We can assume $\tilde U^\prime$ is as well an open subscheme of $S$.
\smallskip
As $S$ is separated, $U^\prime\cap\tilde U^\prime$ is an affine scheme. So it is enough to show that $U=U^\prime\cap\tilde U^\prime$. We assume this is not so. So there is a geometric point $y:{\rm Spec}(k)\to U^\prime\cap\tilde U^\prime\setminus U$. 
We consider a morphism $m_y:{\rm Spec}(k[[T]])\to U^\prime_0\cap\tilde U_0^\prime$, whose special fibre factors through $y$ and whose generic fibre factors through the generic point of $U^\prime_0\cap\tilde U_0^\prime$. We consider the pull back morphism 
$$l_y(m-2q):{\got C}_{0k}/p^{m-2q}{\got C}_{0k}\to y^*({\got C}/p^{m-2q}{\got C}).$$   
\indent
We consider the composite $co$ of $i_{U_0}(m-q)$ with the dual of $\tilde i_{U_0}(m-q)$. The cokernel of $y^*(co)$ is annihilated by $p^{2q+2q+q}=p^{5q}$: $q\ge 2(c+1)>c$ and as $co$ is a morphism between constant $F$-crystals, this is checked (cf. b) of 2.2.4 B) by pulling back $co$ via geometric, dominant points of $U_0$; an extra $p^q$ is added here due to Tate twists. So the kernel of $y^*(co)$ is as well annihilated by $p^{5q}$. We conclude: the kernel and so also the cokernel of $l_y(m-2q)$ are annihilated by $p^{5q}$.
\smallskip
From 1) of the Corollary of 3.15.8 we deduce the existence of a morphism 
$l_y^{\infty}:{\got C}_{0k}\to y^*({\got C})$ lifting the truncation mod $p^{m-n-2q}$ of $l_y(m-2q)$. As the cokernel of $l_y(m-2q)$ is annihilated by $p^{5q}$ and as $m-2q-n\ge 5q+1$, $l_y^\infty$ is injective and so an isogeny. So $y$ factors through $U$. Contradiction. This ends the proof of the Claim and so of the Theorem.
\medskip
{\bf 3.15.10.1. Remarks.} {\bf 0)} The assumption on the existence of $b_0\in\NN$ such that the $h$-number of all pull backs of ${\got C}$ through geometric points of $S$ are all bounded above by $b_0$ is automatically satisfied locally, as one can see by evaluating the $r$-th exterior power of ${\got C}$ at (the thickening of $S$ defined naturally by) $W(S)$. So the Theorem holds without it; however, without assuming it, the estimates of its proof on $n$, $m$, etc., hold only locally.
\smallskip
{\bf 1)} The second assumption on $U$ (of the Theorem) is satisfied, for instance, if $U$ is a Krull scheme (i.e. if it has an open cover by affine, open subschemes which are spectra of Krull $\ZZ$-algebras); see [Sam, \S3 and 4.5 of ch. 1]. In particular, the Theorem applies if $S$ is (the normalization of) a locally noetherian, integral scheme $NS$ in a finite field extension of the field of fractions of $NS$ (cf. [Ma, rm. 2 at the end of ch. 7]). We have variants of the Theorem, by working (i.e. by putting conditions on $U$) in the faithfully flat topology (of $S$).
\smallskip
{\bf 2)} The two conditions on $U$ are not needed at least if we are in a context in which liftings of truncations (in a context involving connections) can be performed. To exemplify what we mean by this, we check that this condition is not needed for the context of an object ${\got C}$ of $p-\Mm_{[0,1]}(S)$ (and so, in particular, for the context of a $p$-divisible group over $S$). 
\smallskip
We refer to the above proof. For any point $\tilde y$ of $S_0$ of codimension 1, there is an open, affine subscheme $\tilde U(\tilde y)$ of $S_0$ containing it and an $\NN$-pro-\'etale morphism $\tilde U(\tilde y;pro)\to \tilde U(\tilde y)$ such that ${\got C}[m]_{0\tilde U(\tilde y;pro)}$ is the truncation mod $p^m$ of an object of $p-\Mm_{[0,1]}^\nabla(\tilde U(\tilde y;pro))$ (cf. Fact of 2.2.1.1 6) and 3.6.18.6 and its $p=2$ version; see the proof of 3.15.1). The Newton polygon stratification of $\tilde U(\tilde y;pro)$ defined by such an object is obtained from a stratification of $\tilde U(\tilde y)$ by pull back (cf. 3.6.20 1) and the fact that $m\ge n$). Let $U_{00}^\prime$ be the open, dense subscheme of $S_0$ over which, using such lifts, we get the same Newton polygon as of ${\got C}_0$. We assume now $\tilde y$ is a point of $U_{00}^\prime$. Denoting by $\tilde V$ the local ring of $\tilde y$, based on the mentioned lift, we can still construct morphisms $i_{\tilde U_1(\tilde V)}:{\got C}_{0\tilde U_1(\tilde V)}/p^m{\got C}_{0\tilde U_1(\tilde V)}\to {\got C}_{\tilde U_1(\tilde V)}/p^m{\got C}_{\tilde U_1(\tilde V)}$, with $\tilde U_1(\tilde V)$ as an open subscheme of the normalization of $U_{00}^\prime$ in a finite field extension of the field of fractions $K_0$ of $R_0$, whose kernels and cokernels are annihilated by $p^q$. Once we have them, we can proceed as in the above proof to get that $U^\prime_{00}$ is affine. So, using pull backs we get $U$ is affine.
\smallskip
Similarly, using 3.15.6 D, we get that these two conditions on $U$ are not needed if we are in a generalized Shimura context. Warning: for the context of $p$-divisible groups, we can use as well standard arguments of algebraic geometry and [Il, 4.8] to reduce the situation to the context of a local, complete, noetherian $S$; but this is not possible for the generalized Shimura context.
\smallskip
{\bf 3)} The proof of 3.15.10 is one of the main extra ingredients needed to get global variants of 3.6.15 B. See \S 9-10 for such variants. 
\medskip
\vfill\eject
\centerline{}
\vfill\eject
\centerline{}
\medskip
\bigskip
\centerline{\bigsll {\bf \S 4 Applications of the basic results to integral canonical models of}}\par
\centerline{\bigsll {\bf Shimura varieties of preabelian type, to $p$-divisible groups,}}
\par
\centerline{\bigsll {\bf and to abelian varieties}}
\bigskip\bigskip
\medskip
4.1-8 and 4.11 form a sequence continuing 2.3. They represent a foundation of integral aspects of Shimura varieties of Hodge type in cases of good reduction. In 4.9-10 and 4.12-13 we extend this foundation from the Hodge type to the preabelian type. In 4.1-13, with very few exceptions, we assume $p\ge 3$. Some overviewing remarks are included in 4.14. In particular, 4.14.3 handles to a great extend the case $p=2$. So in 4.1-14 we deal with many generalizations of classical theories (of Serre--Tate, Manin, Dwork, etc.).
\smallskip
\S 3 contains much more than what is needed for applications to integral canonical models of Shimura varieties of preabelian type. In 4.1-9 and 4.11-14 we try to make as little use of \S 3 as possible. In particular, we avoid the use of 3.6, except of its independent Lemma 3.6.6; even the use of 3.6.6 is in fact entirely avoidable, cf. 3.4.14. By making as little use as possible of \S3, inevitably parts of \S 3 are getting in \S4 a second or third proof. 
\smallskip
We start with a SHS $(f,L_{(p)},v)$. We recall that conventions 2.3.7 and 2.3.9.2 apply. In all that follows $k$ is a perfect field of characteristic $p$.
\medskip\smallskip
{\bf 4.1. The Shimura-ordinary type.} 
Let $L_p^\ast:=L^\ast_{(p)}\otimes_{\ZZ_{(p)}} \ZZ_p$. We also view $G_{\ZZ_p}$ as a subgroup of $GL(L_p^\ast)$. Let $T$ be a maximal torus of $G_{\ZZ_p}$ such that:
\medskip
1) there is a cocharacter 
$\mu:\GG_m\hookrightarrow T_{W(k(v))}$ whose extension to $\CC$ under an $O_{(v)}$-monomorphism $W(k(v))\hookrightarrow\CC$, is $G(\CC)$-conjugate to the cocharacters $\mu^\ast_x:\GG_m\hookrightarrow G_\CC$, $x\in X$;
\smallskip
2) there is a Borel
subgroup $B$ of $G_{\ZZ_p}$, whose Lie algebra is such that its elements take  the $F^1$-filtration of
$L^\ast_p\otimes_{\ZZ_p} W(k(v))$ defined by $\mu$ into itself, i.e. we have ${\rm Lie}(B)\subset F^0\bigl({\rm Lie}(G_{W(k(v))})\bigr)$. 
\medskip
We recall that $\mu$ gives birth to a direct sum decomposition 
$$L^\ast_p\otimes_{\ZZ_p} W(k(v))=F^1\oplus F^0,$$ with $\be\in\GG_m\bigl(W(k(v))\bigr)$ acting through $\mu$ on $F^i$ as the multiplication with $\be^{-i}$, $i=\overline {0,1}$.
\smallskip
For instance, we can take $B$ to be an arbitrary Borel subgroup of $G_{\ZZ_p}$ and  $T$ to be a maximal torus of $B$. So $T$ (resp. $T_{W(k(v))}$) contains a maximal split $\ZZ_p$-torus (resp. $W(k(v))$-torus) of $G_{\ZZ_p}$ (resp. of $G_{W(k(v))}$) and $T_{B(k(v))}$ contains a maximal split torus of $G_{B(k(v))}$ (for instance, see [Ti2]). From the definition of the reflex field (for instance, see [Va2, 2.6]) we get (see [Mi3, 4.6-7]): $B(k(v))$ is the field of definition of the $G(\CC)$-conjugacy class of cocharacters $\mu_x^\ast$, $x\in X$, of the extension of $G_{\QQ_p}$ to $\CC$ under the composite of the inclusion $\QQ_p\hookrightarrow B(k(v))$ with an  $O_{(v)}$-monomorphism $B(k(v))\hookrightarrow\CC$. So from loc. cit. and the fact that any two maximal split tori of $G_{B(k(v))}$ are $G(B(k(v)))$-conjugate (see [Bo2, 20.9 (ii)]) we get that we can choose a cocharacter $\mu$ as in 1); as any two Borel subgroups of $G_{W(k(v))}$ containing a fixed maximal split torus of $G_{W(k(v))}$ are conjugate under an element of $G_{W(k(v))}(W(k(v))$ normalizing this torus (based on [Bo2, 20.9 (i) and (ii)] this can be proved following the pattern of the proof of Fact 1 of 2.2.9 3)), we can assume that 2) holds as well. 
\smallskip
In what follows, $T$ and $B$ are fixed subject to 1) and 2): we do not assume a priori that they are chosen as described in the previous paragraph. So $B$ not necessarily contains $T$; however, it must contain the smallest subtorus $T_{\mu}$ of $T$ with the property that the cocharacter $\mu$ of $T_{W(k(v))}$ factors through ${T_{\mu}}_{W(k(v))}$. So $T$ can be any maximal torus of the centralizer of $T_{\mu}$ in $G_{\ZZ_p}$; in particular, it is not necessarily unique up to $G_{\ZZ_p}(\ZZ_p)$-conjugation. However, as any two maximal tori of $B$ are $B(\ZZ_p)$-conjugate and as any two Borel subgroups of $G_{\ZZ_p}$ are $G_{\ZZ_p}(\ZZ_p)$-conjugate (these statements can be checked by just following the pattern of the proof of Fact 1 of 2.2.9 3), starting from [Bo2, 19.2 and 20.9 (i)] applied over $\FF_p$), we get that $T_{\mu}$ and $\mu$ are uniquely determined up to $G_{\ZZ_p}(\ZZ_p)$-conjugation.
\smallskip
We refer to $(G_{\ZZ_p},[\mu])$ as the Shimura group pair defined by $(G,X,H,v)$. We similarly define the Shimura group pair defined by any Shimura quadruple $(G_1,X_1,H_1,v_1)$.
\medskip
{\bf 4.1.1. The definition.} The Shimura-ordinary type $\tau$ associated to (or defined by) the SHS $(f,L_{(p)},v)$ is the formal isogeny type associated to the $\sg$-crystal over $k(v)$ 
$$
{\got C}_{(f,v)}:=(L^\ast_p\otimes_{\ZZ_p} W(k(v)),\vph)
$$ 
defined by $\vph:=\sg\mu\bigl({1\over p}\bigr)$. Here the Frobenius automorphism 
$$\sg:=\sg_{k(v)}$$ 
of $W(k(v))$ is as well identified with the $\sg$-linear automorphism of $L^\ast_p\otimes_{\ZZ_p} W(k(v))$ fixing $L^\ast_p$. As the notation suggests, ${\got C}_{(f,v)}$ does not depend on the choice of $L_{(p)}$ (producing a SHS $(f,L_{(p)},v)$) or on the choices of $T$ or $B$ (subject to 1) and 2) of 4.1): it depends only on $T_{\mu}$, on $\mu$ and on the representation of $T_{\mu\QQ_p}$ on $L^*[{1\over p}]$, i.e. (cf. 4.1) it depends only on $f$ and $v$. 
\smallskip
Let $r(f,v)\in\NN\cup\{0\}$ be the multiplicity of its slope $1$. Let $FL(f,v)$, $IL(f,v)$, $SKL(f,v)$, $CL(f,v)$ and $SDL(f,v)$ be respectively the factor length, the isotype length, the skeleton length, the circular length and the slope denominator length of the extension of ${\got C}_{(f,v)}$ to $\FF$ (cf. I of 2.2.22 3)).
\smallskip
In 4.1 1) we can replace $k(v)$ by any other finite field: such a field must contain $k(v)$ and in fact $\mu$ is automatically definable over $W(k(v))$, cf. [Mi3, 4.6-7]. Moreover, if in 4.1 2) we replace $B$ by a parabolic subgroup $PAR$ of $G_{\ZZ_p}$, then we regain the same Shimura-ordinary type: this is a consequence of the fact that $PAR$ contains a Borel subgroup of $G_{\ZZ_p}$ (cf. Fact of 2.2.3 3) and [Bo2, 21.12] applied to $G_{\QQ_p}$). 
\medskip
{\bf 4.1.1.1. Simple properties.} Let $M:=L^\ast_p\otimes_{\ZZ_p} W(k(v))$ and let ${\got t}:={\rm Lie}\bigl(T_{W(k(v))}\bigr)$. Let 
$$\bar h_1:M\to M$$ 
be the endomorphism defined by:
$\bar h_1(x)=0$ if $x\in F^0$ and $\bar h_1(x)=x$ if $x\in F^1$. $\bar h_1$ is an element of
{\got t}. Let $d(v)\in\NN$ be such that $\abs{k(v)}=p^{d(v)}$; so $k(v)=\FF_{p^{d(v)}}$ and $\FF=\overline{k(v)}$. In what follows, as $d(v)$ is used very often, we denote it in a simpler way by $d$. For $i\in S(1,d-1)$, let $\bar h_{i+1}:=\vph^i(\bar h_1)$. We have:
\medskip
{\bf Fact 1.} {\it $\vph^d(\bar h_1)=\bar h_1$ and $\bar h_{i+1}\ne\bar h_1$,
$\forall i\in S(1,d-1)$.}
\medskip
{\bf Proof:} We have $\vph^i(\bar h_1)=\sg^i(\bar h_1)$. So, Fact 1 follows from [Mi3, 4.6-7]: $[\mu]$ can not be defined over $W(k_1)$, with $k_1$ a subfield of $k(v)$ different from $k(v)$.
\medskip
 Let 
$$h^0:={1\over d}\sum^d_{i=1}\bar h_i\in {\got t}[{1\over p}]$$ 
We have:
\medskip
{\bf Fact 2.} {\it $\vph(h^0)=\sg(h^0)=h^0$. Moreover, the slopes of the $\sg$-crystal
$(M,\vph)$ are the eigenvalues of $h^0$ (acting on $M[{1\over p}]$) and the multiplicities are the same.}
\medskip
{\bf Proof:} The first part is obvious. The second part results from the fact that $\vph^d$ is a $B(k(v))$-valued of $T_{W(k(v))}$ fixed by the action of $\sg$ on $T_{W(k(v))}$ and the eigenspaces of its action on $M[{1\over p}]$ are the same as the ones of $dh^0$, the eigenvalues (counted with multiplicities) being $p$ to those powers which are eigenvalues (counted with multiplicities) of $dh^0$. This is a consequence of the fact that ${\bar h}_i\in {\got t}$, $\forall i\in S(1,d)$, and so of the fact that the endomorphisms ${\bar h}_i$'s are commuting among themselves.
\smallskip
More precisely, if $\tilde x\in M_i:=\bigl\{x\in M|h^0(x)={i\over d}x\bigr\}$, with $i\in S(0,d)$, and if $\tilde x\not\in pM$, then $\vph^d(\tilde x)=p^i\tilde y$, with $\tilde y\in M_i\setminus pM$. Moreover, we have 
$M=\bigoplus_{0\le i\le d}M_i$. To check these last two statements, we first remark that each element of ${\rm Lie}(T)$ is fixed by $\vph$; adding these elements to the family $(v_{\al})_{\al\in\Mj^\prime}$, we get that $(M,F^1,\vph,T_{W(k(v))})$ is a Shimura filtered $\sg$-crystal. From Corollary of 2.2.9 3), we get that $\mu$ is the canonical split of $(M,F^1,\vph)$. So the statements follow from 2.2.1.1 4) and the first paragraph of G of 2.2.22 3).
\medskip
{\bf 4.1.1.2. The Lie counterpart.} Let ${\rm Lie}_G(\tau)$ be the formal isogeny type associated to the Shimura Lie $\sg$-crystal 
$$
\bigl({\rm Lie}(G_{W(k(v))}),\vph\bigr) 
$$
attached to the Shimura $\sg$-crystal $(M,\vph,G_{W(k(v))},(v_{\al})_{\al\in\Mj^\prime})$. As above, the slopes of it are the eigenvalues of $h^0$ (acting on ${\rm Lie}(G_{B(k(v))})$ via the rule: $x$ is mapped into $[h^0,x]$) and the multiplicities are the same. We call it the Shimura-ordinary Lie type associated to (or defined by) the SHS $(f,L_{(p)},v)$. Let $r(v^{\rm ad})\in\NN\cup\{0\}$ be the multiplicity of its slope $-1$. Let $SS(v^{\rm ad})$ be the set of its slopes; it is a subset of $[-1,1]\cap\QQ$, symmetric w.r.t. $0$ (cf. 2.2.3 1)); moreover, always $0\in SS(v^{\rm ad})$.
\medskip
{\bf 4.1.1.3. Remark.} Similarly to 4.1.1, ${\rm Lie}_G(\tau)$ depends only on the adjoint Shimura pair $(G^{\rm ad},X^{\rm ad})$, on the rank of $G^{\rm ab}$ and on $v^{\rm ad}$. So, as the notations suggest, $r(v^{\rm ad})$ and $SS(v^{\rm ad})$ depend only on $(G^{\rm ad},X^{\rm ad})$ and on $v^{\rm ad}$.
\medskip
{\bf 4.1.1.4. Some decompositions.} For any $i\in S(1,d)$ and $j\in\{0,1\}$ we denote by $^iF^j$ the free
$W(k(v))$-submodule of $M$ on which $\bar h_i$ acts trivially if $j=0$
and identically if $j=1$. For any function $\bar f:S(1,d)\to\{0,1\}$ we denote by 
$$F_{\bar f}:=\bigcap_{i\in S(1,d)}\,^iF^{\bar f(i)}.$$ 
Let $\Ml$ be the set of functions $\bar f$ as above for which $F_{\bar f}\ne 0$. We get a direct sum decomposition 
$$M=\oplus_{\bar f\in\Ml}F_{\bar f}$$ 
due to the fact that $\bar h_i$'s are commuting among themselves. Let 
$$
\bar\vph:\Ml\to\Ml
$$
be the bijection defined by the rule: $\bar\vph(\bar f)(i)=\bar f(i-1)$; here $\bar f(0):=\bar f(d)$. From the very definition of $F_{\bar f}$'s we get: 
\medskip
a) $\vph(F_{\bar f})=F_{\bar\vph(\bar f)}$ if $\bar f(1)=0$;
\smallskip
b) $\vph(F_{\bar f})=pF_{\bar\vph(\bar f)}$ if $\bar f(1)=1$. 
\medskip
Let $\bar\vph=\prod_{l\in I_{\bar\vph}}\bar\vph_l$ be written as a product of disjoint cyclic permutations. We allow trivial cyclic permutations. So we have a disjoint union 
$$I_{\bar\vph}=I^1_{\bar\vph}\cup I^0_{\bar\vph};$$ 
$l\in I_{\bar\vph}$ belongs to $I^1_{\bar\vph}$ iff $\bar {\vph}_l$ is a non-trivial permutation (and then any $\bar f\in\Ml$ such that $\bar {\vph}_l(\bar f)\ne\bar f$ is said to be associated to $\bar {\vph}_l$). So $I^0_{\bar\vph}:=\{\bar f\in\Ml|\bar\vph (\bar f)=\bar f\}$; any $\bar f\in I^0_{\bar\vph}$ is said to be associated to $\bar {\vph}_{\bar f}$.  $I^0_{\bar\vph}$ has either $0$ or $2$ elements depending on the fact that $r(f,v)$ is $0$ or not. As $\bar\vph^d=1_\Ml$, we get that the order $d_l$ of the cyclic permutation $\bar\vph_l$ divides $d$. For any $l\in I_{\bar\vph}$, we choose arbitrarily an element $\bar f_l$ of $\Ml$ which is associated to $\bar\vph_l$. We have
$$\vph^{d_l}(F_{\bar f_l})=p^{{d_l\over d}\ell(\bar f_l)}F_{\bar f_l},$$ 
with 
$$\ell(\bar f_l):=\sum^d_{i=1}\bar f_l(i),$$ 
and $\vph^{d_l}$ acts on $F_{\bar f_l}$ as
$p^{{d_l\over d}\ell(\bar f_l)}\sg^{d_l}$. Let 
$$
_pF_{\bar f_l}:=\bigl\{x\in F_{\bar f_l}|\sg^{d_l}(x)=x\bigr\}.
$$ 
It is a well known fact that $_pF_{\bar f_l}$ is a $W(\FF_{p^{d_l}})$-free submodule of $F_{\bar f_l}$ such that $F_{\bar f_l}={}_pF_{\bar f_l}\otimes_{W(\FF_{p^{d_l}})} W(k(v))$ (for instance, this follows easily from [Bo2, Prop. of p. 30] applied over $B(\FF_{p^{d_l}})$ to $F_{\bar f_l}[{1\over p}]$).
\medskip
{\bf 4.1.2. Proposition.} {\it The filtered $\sg$-crystal $(M,F^1,\vph)$ is cyclic diagonalizable.}
\medskip
{\bf Proof:} We choose arbitrarily a $W(\FF_{p^{d_l}})$-basis $\{\ell_s|s\in\Mj_l\}$ of ${}_pF_{\bar f_l},\forall l\in I_{\bar\vph}$. We view it as a $W(k(v))$-basis of $F_{\bar f_l}$. For any cyclic permutation $\bar\vph_l$ of length $\ge 2$ (i.e. for when we deal with an $l\in I^1_{\bar\vph}$) and for every element $\bar f\in\Ml$ associated to $\bar\vph_l$ but different from $\bar f_l$,
let $j\in\NN$ be the smallest number such that $\bar f=\bar\vph^j(\bar f_l)$. We get a $W(k(v))$-basis of $F_{\bar f}$ by taking
${1\over p^{e_{jl}}}$ of the images of $\vph^j(\ell_s)$, $s\in\Mj_l$, with 
$$e_{jl}:=\sum^j_{i=1}\bar f_l(i).$$ 
This expression of $e_{jl}$'s is a consequence of the following iterating formula:
$$\vph^j=\bigl(\prod_{m=1}^j \sg^m(\mu({1\over p}))\bigr)\sg^j,\leqno (IT)$$
where $\sg^m(\mu):=\sg^m\circ\mu({1\over p})\circ\sg^{-m}$, $\forall m\in\ZZ$. 
\smallskip
With respect to the $W(k(v))$-basis of $M$ obtained by putting together the chosen $W(k(v))$-bases of $F_{\bar f}$, $\bar f\in\Ml$, $(M,F^1,\vph)$ gets the desired cyclic diagonalizable form. This ends the proof of the Proposition.
\medskip
{\bf 4.1.2.1. Remark.} By very definition (see 2.2.22 1)) the Shimura filtered $\sg$-crystal $(M,F^1,\vph,T_{W(k(v))})$ is strongly cyclic diagonalizable.
\medskip
{\bf 4.1.2.2. Remark.} For what follows we refer to D and I of 2.2.22 3). All classification and standard invariants of (the extension to $\FF$) of ${\got C}_{(f,v)}$ are encoded in the permutation $\bar\vph$ of $\Ml$. For instance, the cyclic length is the maximum of the $d_l$'s numbers, while the skeleton length is equal to $\abs{I_{\bar\vph}}$. In particular, the associated formal sum to this extension is of the form
$$\sum_{q\in \NN;\, q|d} \sum_{\tau_q\in\Mt_q} m_{f,v}(\tau_q)\tau_q.$$ 
\indent
{\bf 4.1.3. Example.} $\Ml$ has 2 elements iff $\tau$ is an ordinary type or a supersingular type, i.e. iff $d=1$ or if $d=2$ and $(M,\vph)$ does not have integral slopes.
\medskip
{\bf 4.1.4. Proposition.} {\it The Shimura filtered $\sg$-crystal 
$$\bigl({\rm Lie}(G_{\ZZ_p})\otimes_{\ZZ_p} W(k(v)),\vph,F^0({\rm Lie}(G_{W(k(v))})),F^1({\rm Lie}(G_{W(k(v))}))\bigr)$$
 and its adjoint are cyclic diagonalizable and of Borel type.}
\medskip
{\bf Proof:} The fact that the mentioned two Shimura filtered Lie $\sg$-crystals are of Borel type results from the fact that $\vph({\rm Lie}(B_{W(k(v))}))\subset {\rm Lie}(B_{W(k(v))})$ (we recall that $B$ is a Borel subgroup of $G_{\ZZ_p}$). The
proof of their cyclic diagonalizability is entirely analogous to the proof of 4.1.2 (we just have to deal with functions from $S(1,d)$ to $S(-1,1)$), and so it is omitted. 
\medskip
Above we wrote ${\rm Lie}(G_{\ZZ_p})\otimes_{\ZZ_p} W(k(v))$ instead of ${\rm Lie}(G_{W(k(v))})$ just to emphasize how $\vph=\sg\mu\bigl({1\over p}\bigr)$ acts on it: $\sg$ fixes ${\rm Lie}(G_{\ZZ_p})$ and acts as the Frobenius automorphism on $W(k(v))$, while $\mu\bigl({1\over p}\bigr)$ acts via inner conjugation.
\medskip
{\bf 4.1.4.1. Remark.} We assume $B$ contains $T$; let $B^{\rm opp}$ be the opposite of $B$ w.r.t. $T$. As $\vph$ takes ${\rm Lie}(T_{W(k(v))})$ and ${\rm Lie}(B_{W(k(v))})$ into themselves, $p\vph$ takes ${\rm Lie}(B_{W(k(v))}^{\rm opp})$ into itself.
\medskip
{\bf 4.1.5. $\om$'s types.}  Let 
$N$ be the normalizer of $T$ in $G_{\ZZ_p}$. Let 
$$
W_G:=(N/T)\bigl(W(\overline{k(v)})\bigr)
$$ be the Weyl group of 
$G_{W(\FF)}$ w.r.t. $T_{W(\FF)}$. It depends only on $G^{\rm ad}_{\CC}$ and not on $v$. 
For $\om\in W_G$, let $g_{\om}\in N(W(\FF))$ be such that its image in $W_G$ is $\om$. Let
$${\got C}_{\om}:=(M\otimes_{W(k(v))} W(\FF),F^1\otimes_{W(k(v))} W(\FF),g_{\om}(\sg\mu({1\over p})\otimes 1),G_{W(\FF)})$$ be the Shimura filtered $\bar\sg$-crystal defined by $g_{\om}$. We denote by $\Mp_{\om}$ (resp. by ${\rm Lie}_G(\Mp_{\om})$) the Newton polygon of ${\got C}_{\om}$ (resp. of the Shimura Lie $\sg$-crystal attached to ${\got C}_{\om}$). Let $\tau_{\om}$ be the formal isogeny type of ${\got C}_{\om}$; we refer to it as the $\om$'s type associated to (or defined by) $(f,L_{(p)},v)$; if $\om$ is the identity element of $W_G$, we regain the Shimura-ordinary type of 4.1.1. We have the following extension of 4.1.2 and 4.1.4.
\medskip
{\bf 4.1.5.1. Proposition.} {\it ${\got C}_{\om}$ is cyclic diagonalizable, its attached Shimura (adjoint) filtered Lie $\sg$-crystal is cyclic diagonalizable; moreover it (and so also $\Mp_{\om}$ and ${\rm Lie}_G(\Mp_{\om})$) depends only on $\om$ and not on the choice of $g_{\om}$.}
\medskip 
{\bf Proof:} The proof of the cyclic diagonalizability of ${\got C}_\om$ or of its attached Shimura filtered (adjoint) Lie $\sg$-crystal is entirely analogous to the proof of 4.1.2, and so it is omitted. The cyclic diagonalizability results as well from 2.2.16, once we remark (cf. Fact of 2.2.9 1)) that $(M\otimes_{W(k(v))} W(\FF),F^1\otimes_{W(k(v))} W(\FF),g_{\om}(\sg\mu({1\over p})\otimes 1),T_{W(\FF)})$ is a Shimura filtered $\sg$-crystal; so the last part of the Proposition is a very particular case of 3.11.1 c). This ends the proof.
\medskip
{\bf 4.1.5.2. Remark.} In 4.1.5-6, instead of $\FF$, we can work as well with the smallest finite field extension $k(v^{\rm sp}_G)$ of $k(v)$ over which $G_{k(v)}$ splits; warning: in general, in connection to 4.1.5.1 we have to replace accordingly the wording cyclic diagonalizable by ``potentially cyclic diagonalizable". 
\medskip
{\bf 4.1.5.3. Degrees of definitions.} Let $d_{\om}$ be the degree of definition of ${\got C}_{\om}$. The set $SDD(f,v)$ of all $d_{\om}$'s is precisely the set $SDD(M\otimes_{W(k(v))} W(\FF),\sg\mu({1\over p})\otimes 1,G_{W(\FF)})$ defined in Problem 2 of 2.2.22 H. If $G_{\ZZ_p}$ is a split group, then $SDD(f,v)$ is a subset of the set of orders of elements of $W_G$ (cf. 4.1.5.2). From 3.11.3.1 we get that $d|d_{\om}$.  
\medskip
{\bf 4.1.5.4. Simple properties.} Let $\sg_{\om}:=g_{\om}\bar\sg$. It is a $\bar\sg$-linear automorphism of $M\otimes_{W(k(v))} W(\FF)$. We have $g_{\om}(\sg\mu({1\over p})\otimes 1)=\sg_{\om}\mu_{W(\FF)}({1\over p})$. The canonical split of ${\got C}_{\om}$ is $\mu_{W(\FF)}$ itself (cf. the proof of 4.1.5.1 and the Corollary of 2.2.9 3)). The elements of $M\otimes_{W(k(v))} W(\FF)$ fixed by $\sg_{\om}$ give birth to a $\ZZ_p$-structure of $(M\otimes_{W(k(v))} W(\FF),G_{W(\FF)},(v_{\al})_{\al\in\Mj^\prime})$: from Lang's theorem we get it is precisely $(L_p^*,G_{\ZZ_p},(v_{\al})_{\al\in\Mj^\prime})$. We deduce (cf. 2.2.9 8) applied to $(M\otimes_{W(k(v))} W(\FF),F^1\otimes_{W(k(v))} W(\FF),g_{\om}(\sg\mu({1\over p})\otimes 1),T_{W(\FF)})$) that $({\got C}_{\om},(v_{\al})_{\al\in\Mj^\prime})$ is $1_{\Mj^\prime}$-isomorphic to the extension to $\FF$ of a strongly cyclic diagonalizable Shimura filtered $F$-crystal over $\FF_{p^{d_{\om}}}$
$$(M\otimes W(\FF_{p^{d_{\om}}}),F^1_{\om},\sg_{\FF_{p^{d_{\om}}}}\mu_{\om}({1\over p}),G_{W(\FF_{p^{d_{\om}}})},(v_{\al})_{\al\in\Mj^\prime})),$$
where $\mu_{\om}$ is a cocharacter of $G_{W(\FF_{p^{d_{\om}}})}$ such that $[\mu_{\om}]=[\mu_{W(\FF_{p^{d_{\om}}})}]$ and where $F^1_{\om}$ is the maximal $W(\FF_{p^{d_{\om}}})$-submodule of $M\otimes W(\FF_{p^{d_{\om}}})$ on which $\mu_{\om}$ acts as the inverse of the identity character of $\GG_m$; its canonical split cocharacter is $\mu_{\om}$. The main property of $\mu_{\om}$ is (cf. the cylic diagonalizability; see also Exercise of 2.2.22 1)):
\medskip
{\bf Fact.} {\it It commutes with $\sg^s_{\FF_{p^{d_{\om}}}}\mu_{\om}\sg^{-s}_{\FF_{p^{d_{\om}}}}$, $\forall s\in\NN$.}
\medskip
Let $T_{\mu_{\om}}$ be the smallest torus of $G_{\ZZ_p}$ through which $\mu_{\om}$ factors. Let $T^{\om}$ be a maximal torus of $G_{\ZZ_p}$ containing it.
\medskip
{\bf Exercise.} Let $\om_1$, $\om_2\in W_G$. Show that $({\got C}_{\om_1},(v_{\al})_{\al\in\Mj^\prime})$ and $({\got C}_{\om_2},(v_{\al})_{\al\in\Mj^\prime})$ are $1_{\Mj^\prime}$-isomorphic, iff the extensions to $W(\FF)$ of $\mu_{\om_1}$ and $\mu_{\om_2}$ are $G_{\ZZ_p}(\ZZ_p)$-conjugate (and not only $G_{\ZZ_p}(W(\FF))$-conjugate), i.e. iff the extensions to $\FF$ of $\mu_{\om_1}$ and $\mu_{\om_2}$ are $G_{\ZZ_p}(\FF_p)$-conjugate. Hint: an $1_{\Mj^\prime}$-isomorphism between $({\got C}_{\om_1},(v_{\al})_{\al\in\Mj^\prime})$ and $({\got C}_{\om_2},(v_{\al})_{\al\in\Mj^\prime})$ must preserve their canonical lifts, i.e. it is given by an element of the centralizer of the image of $\mu_{W(\FF)}$ in $G_{W(\FF)}$.
\medskip
{\bf 4.1.5.5. The Faltings--Shimura--Hasse--Witt relation.} Let $R_G(v)$ be the equivalence relation on $W_G$ defined by: $\om_1$, $\om_2\in W_G$ are in relation $R_G(v)$, iff the Faltings--Shimura--Hasse--Witt adjoint maps of ${\got C}_{\om_1}$ and ${\got C}_{\om_2}$ are isomorphic, under an isomorphism defined by an element of $G^{\rm ad}_{\ZZ_p}(\FF)$.
\medskip
{\bf 4.1.6. The minimal model.} We view the alternating form $\tilde\psi$ of 2.3.1  also as a perfect alternating form $\tilde\psi: M\otimes_{W(k(v))} M\to W(k(v))(1)$. So
$$\Mm\Mm:=(M,F^1,\vph,G_{W(k(v))},(v_{\al})_{\al\in\Mj^\prime},\tilde\psi)$$ 
is a principally quasi-polarized Shimura filtered $\sg$-crystal. From 4.1.4, as in 3.4.3.0, we get (cf. also 3.11.6 B)) that its attached Shimura filtered Lie $\sg$-crystal is of parabolic type. From 3.1.0 a) we get that it is a Shimura-canonical lift. The cocharacter $\mu$ is the canonical split of $(M,\vph,G_{W(k(v))})$ (cf. the proof of Fact 2 of 4.1.1.1 and 3.1.5). From the very definitions of 3.11.2 B we get that $\Mm\Mm$ is the minimal model of its extension to $\FF$. 
\medskip\smallskip
{\bf 4.2. $G$-ordinary points.} 
\medskip
{\bf 4.2.1. Theorem.} {\it There is a $G(\AA_f^p)$-invariant, open, dense subscheme $\Mu$ of $\Mn_{k(v)}$ such that for any point $y:{\rm Spec}(k)\to\Mn_{k(v)}$ we have:
\medskip
a) $\tau$ is the formal isogeny type of the Shimura $\sg_k$-crystal attached to it (cf. 2.3.10) iff it factors through $\Mu$;
\smallskip
b) ${\rm Lie}_G(\tau)$ is the formal isogeny type of the Shimura Lie $\sg_k$-crystal attached to it (cf. 2.3.10) iff it factors through $\Mu$;
\smallskip
c) the Shimura $\sg_k$-crystal attached to it is Shimura-ordinary iff $y$ factors through $\Mu$.}
\medskip
{\bf Proof:} Let $V_0:=W(\FF)$. Let
$\Mn_{V_0}:=\Mn\times_{O(v)}V_0$ and let $\Mn_\FF$ be its special fibre. Let $y_0:{\rm Spec}(\FF)\hookrightarrow\Mn_\FF$ be a closed point and let $z_0:{\rm Spec}(V_0)\hookrightarrow\Mn_{V_0}$ be a lift of it. Let $\Mn^0_{\FF}$ be the connected component of $\Mn_\FF$ containing $y_0$. 
\smallskip
Let $(A_0,p_{A_0}):=z_0^*(\Ma,\Mp_{\Ma})$. Let $\bigl(M_0,F^1_0,\vph_0,G_{V_0},(t_\al)_{\al\in\Mj^\prime},p_{A_0}\bigr)$ be the principally quasi-polarized Shimura filtered $\bar\sg$-crystal attached to $z_0$. So (cf. 2.3.10) $M_0:=H^1_{dR}(A_0/V_0)$,
$F^1_0$ is its Hodge filtration defined by $A_0$ and $(t_\al)_{\al\in\Mj^\prime}$ is the family of de Rham components of the family $(w_\al)_{\al\in\Mj^\prime}$ of Hodge cycles with which $A_0$ is (cf. 2.3.3) naturally endowed. From 2.3.13.1 we deduce that for any point $z_1:{\rm Spec}(V_0)\hookrightarrow\Mn_{V_0}$ lifting an $\FF$-valued point $y_1$ of $\Mn^0_{\FF}$, the Shimura filtered $\bar\sg$-crystal attached to it is of the form $\bigl(M_0,F^1_0,\vph_1,G_{V_0},(t_\al)_{\al\in\Mj^\prime}\bigr)$.
We have $\vph_1=g_1\vph_0$, with $g_1\in G_{V_0}(V_0)$ (cf. the beginning of \S 3). From 2.3.16 we get:
\medskip
{\bf (4.2.1.1)} {\it For any $\FF$-valued point $\bar g_1$ of a non-empty open subscheme of $G_\FF$, there is $g_1\in G_{V_0}(V_0)$ lifting it and there is $z_1\in\Mn_{V_0}(V_0)$ lifting a $y_1\in\Mn^0_{\FF}(\FF)$ for which we have $\vph_1=g_1\vph_0$.}
\medskip
Let $\Mu_0$ (resp. $\Mu_1$) be the open subscheme of $\Mn^0_{\FF}$ of whose geometric points have the smallest Newton polygon of their attached Shimura $F$-crystals (resp. Shimura Lie $F$-crystals). $\Mu_0$ and $\Mu_1$ exist, cf. the specialization theorem.
\smallskip
As all $\bar\sg$-crystals $(M_0,g_1\vph_0)$, with $g_1\in G_{V_0}(V_0)$, give birth to 1/2-symmetric isocrystals we get:
\medskip 
{\bf (4.2.2)} {\it The fact that the  $\bar\sg$-crystal $(M_0,\vph_0)$ attached to the point $y_0\in\Mn^0_{\FF}(\FF)$ has the smallest Newton polygon among all $\bar\sg$-crystals attached to $\FF$-valued points of $\Mn^0_{\FF}$, is equivalent to the fact that the Lie $\bar\sg$-crystal $\bigl({\rm End}(M_0),\vph_0\bigr)$ has the smallest Newton polygon among all Lie $\bar\sg$-crystals $\bigl({\rm End}(M_0),\vph_1\bigr)$ defined naturally by $\bar\sg$-crystals attached to $\FF$-valued points of $\Mn^0_{\FF}$.}
\medskip
So, as the Newton polygon of $\bigl({\rm End}(M_0),\vph_1\bigr)$ is obtained from the Newton polygon of $\bigl({\rm Lie}(G_{V_0}),\vph_1\bigr)$ and from the Newton polygon of
$\bigl({\rm End}(M_0)/{\rm Lie}(G_{V_0}),\vph_1\bigr)$, $\Mu_0$ is an open subscheme of $\Mu_1$. $\bigl({\rm End}(M_0)/{\rm Lie}(G_{V_0}),\vph_1\bigr)$ is not a Lie $\bar\sg$-crystal; but the quadruple 
$$\bigl({\rm End}(M_0)/{\rm Lie}(G_{V_0}),F_0^0({\rm End}(M_0))/F_0^0({\rm Lie}(G_{V_0})),F_0^1({\rm End}(M_0))/F_0^1({\rm Lie}(G_{V_0})),\vph_1\bigr)$$  
is a $p$-divisible object of $\Mm\Mf_{[-1,1]}(V_0)$.
\medskip
{\bf 4.2.3. The use of Fontaine's comparison theory.} In what follows we refer to [Va2, \S 5]. Loc. cit. was worked out under the condition (*) of [Va2, 5.1]. Such a condition is not necessarily satisfied in the context of our SHS $(f,L_{(p)},v)$ (cf. 2.3.8 1)). However, its role was just to reach to a reductive context as of 2.3.4 and 2.3.11 (see [Va2, 5.2.12]). So below we refer to parts of [Va2, \S 5] without any extra comment. Let $K_0:=V_0[{1\over p}]$. 
\smallskip
It is known that there is a $K_0$-isomorphism 
$$
H^1_{\acute et}(A_{0\overline{K_0}},\QQ_p)\otimes_{\QQ_p}K_0\tilde\to M_0\otimes_{V_0} K_0
$$ 
taking the
$p$-component of the \'etale component $u_\al$ of $w_\al$ into $t_\al$, $\forall\al\in\Mj^\prime$, and taking $\tilde\psi$ into $p_{A_0}$ (cf. [Va2, 5.1.4 and 5.2.17.2]). But $H^1_{\acute et}(A_{0\overline{K_0}},\QQ_p)$ and the family of tensors $(u_\al)_{\al\in\Mj^\prime}$ of $\Mt\bigl(H^1_{\acute et}(A_{0\overline{K_0}},\QQ_p)\bigr)$ can be respectively identified (non-canonically) with $L_p^*\otimes_{\ZZ_p} \QQ_p$ and the family of tensors $(v_\al)_{\al\in\Mj^\prime}$ of $\Mt(L_p^*\otimes_{\ZZ_p} \QQ_p)$ (cf. [Va2, top of p. 473]). So
we get a $K_0$-isomorphism 
$$\rho_0:M\otimes_{W(k(v))} K_0\tilde\to M_0\otimes_{V_0} K_0$$
 taking $v_\al$ into $t_\al$, $\forall\al\in\Mj^\prime$, and taking $\tilde\psi$ into $p_{A_0}$.
\smallskip
The isomorphism $\rho_0$ allows us to identify the extension from $B(k(v)$) to $K_0$ of the $\sg$-linear automorphism $\vph$ of $M[{1\over p}]$ defined in 4.1, with a $\bar\sg$-linear automorphism of $M_0\otimes_{V_0} K_0$. 
So we have $\vph\otimes 1=\tilde g\vph_0$, with $\tilde g\in G^0(K_0)$.
There is a cocharacter $\mu_0:\GG_m\to G_{V_0}$ which produces a direct sum decomposition $M_0=F^1_0\oplus F^0_0$, with $\be\in\GG_m(V_0)$ acting through $\mu$ on $F^i_0$ as the multiplication with $\be^{-i}$ (cf. 2.3.10); under extension via an $O_{(v)}$-monomorphism $V_0\hookrightarrow\CC$, $\mu_0$ becomes $G(\CC)$-conjugate to the cocharacters $\mu^\ast_x$, $x\in X$. As $G_{K_0}$ is a split group, $\mu_0$ and $\mu_{K_0}$ are in fact $G_{V_0}(K_0)$-conjugate. So we can assume 
$$\rho_0(F^1\otimes_{W(k(v))} K_0)=F_0^1\otimes_{V_0} K_0.$$
\medskip
{\bf 4.2.3.1. The case of a torus.} We first consider the case when $G$ is a torus. $\Mn_{V_0}$ as a scheme is a disjoint union of copies of ${\rm Spec}(V_0)$, cf. [Va2, 3.2.8]. These copies are permuted transitively by $G(\AA_f^p)$, cf. [Va2, 3.3.1-2]. So taking $\Mu=\Mn_{k(v)}$, we just have to show that the formal isogeny type of the $\bar\sg$-crystal attached to $y_0$ is $\tau$. This is well known (for instance, this is implicitly contained in [Ko1, \S 2]). However, we recall one way to prove this, using the language of this paper. As $(M_0,F_0^1,\vph_0,G_{V_0})$ is a cyclic diagonalizable Shimura filtered $\bar\sg$-crystal (see 2.2.16), as in 4.1.1.1 we get that its formal isogeny type can be computed from the way $\vph_0^s$'s, $s\in\NN$, act on $\mu_0$ (as $G$ is a torus, $\mu_0$ is uniquely determined, cf. Corollary of 2.2.9 3)). But from the identifications of 4.2.3 we get: this way is the same as the way $\bar\sg^s$'s, $s\in\NN$, act on $\mu_{K_0}$. This proves 4.2.1 in the case when $G$ is a torus. 
\smallskip
From now on we assume $G$ is not a torus.
\medskip
{\bf 4.2.4. Some simplifications.} We got 2 particular hyperspecial subgroups of $G(K_0):G_{V_0}(V_0)$ (defined via the $V_0$-lattice $M_0$ of $M_0\otimes_{V_0} K_0$) and $G_{\ZZ_{(p)}}\bigl(M\otimes_{W(k(v))}V_0\bigr)$ (defined via the $V_0$-lattice $M\otimes_{W(k(v))}V_0$ of $M\otimes_{W(k(v))} K_0$). But any two hyperspecial subgroups of $G(K_0)$ are conjugate by an element $g_0\in G^{\rm ad}(K_0)$ (cf. [Ti2, p. 47]). Changing the Shimura adjoint filtered Lie $\bar\sg$-crystal attached to $z_0$ by the inner automorphism defined by $g_0$, it can be put in the form
$\bigl({\rm Lie}(G^{\rm ad}_{V_0}),\vph^\prime_0,F^0_0\bigl({\rm Lie}(G^{\rm ad}_{V_0})\bigr),F^1_0\bigl({\rm Lie}(G^{\rm ad}_{V_0})\bigr)\bigr)$, where $\vph^\prime_0=g(\vph\otimes 1)$, with $g\in G_{V_0}^{\rm ad}(K_0)$. Warning: here, in accordance with the original notation $G_{\ZZ_{(p)}}$ (of 2.3.1), $G_{V_0}^{\rm ad}$ is w.r.t. the second hyperspecial subgroup; we had no reason to write ${\rm Lie}(G_{\ZZ_p})\otimes_{\ZZ_p} V_0$ instead of ${\rm Lie}(G_{V_0}^{\rm ad})$.
\smallskip
As $\mu$ and $\mu_0$ are $G_{V_0}(K_0)$-conjugate, changing everything by an isomorphism (defined by an element of $G_{V_0}^{\rm ad}(V_0)$, cf. Iwasawa's decomposition of [Ti2, 3.3.2]), we can assume
$F^0_0\bigl({\rm Lie}(G^{\rm ad}_{V_0})\bigr)=F^0\bigl({\rm Lie}(G^{\rm ad}_{V_0})\bigr)$ and $F^1_0\bigl({\rm Lie}(G^{\rm ad}_{V_0})\bigr)=F^1\bigl({\rm Lie}(G^{\rm ad}_{V_0})\bigr)$. As $g$ normalizes ${\rm Lie}(G_{V_0}^{\rm ad})$ and as $G^{\rm ad}_{V_0}(V_0)$ is that maximal bounded subgroup of $G^{\rm ad}_{V_0}(K_0)$ (being a hyperspecial subgroup of it, see [Ti2, 3.2]) normalizing ${\rm Lie}(G_{V_0})$, $g\in G_{V_0}^{\rm ad}(V_0)$. So $({\rm Lie}(G_{V_0}),g(\vph\otimes 1))$ is as well the Shimura Lie $\bar\sg$-crystal attached to $y_0$.
\smallskip
This last fact allows us, via the first sentence of 3.5.5, to appeal to the results of 3.1-3, 3.6.6 and 3.4.
\medskip
{\bf 4.2.5. The use of (4.2.1.1).} (4.2.1.1) implies: there is a non-empty open subscheme $U$ of $G^{\rm ad}_{\FF}$ with the property that for any element $\bar g_0\in U(\FF)$, there is a
$V_0$-valued point of $\Mn_{V_0}$ lifting an $\FF$-valued point of $\Mn^0_{\FF}$ and such that its attached  Shimura adjoint filtered Lie $\bar\sg$-crystal is isomorphic to
$${\got l}_{g_1}:=\bigl({\rm Lie}(G^{\rm ad}_{V_0}),g_1g(\vph\otimes 1),F^0\bigl({\rm Lie}(G^{\rm ad}_{V_0})\bigr),F^1\bigl({\rm Lie}(G^{\rm ad}_{V_0})\bigr)\bigr),$$ with $g_1$ mod $p$ equal to $\bar g_0$.
\medskip
{\bf 4.2.6. The definition of $\Mu_{\FF}\cap\Mn^0_{\FF}$.} Now 4.2.1 for $\Mn^0_{\FF}$ is a direct consequence of 4.1.1.4, 4.2.5, 3.2.5, 3.4.3.0, 3.4.11, 3.4.13, 3.5.5, 3.6.6, 3.7.2 and 4.2.7 below. We present the details. From \S 3 and 4.1 we need:
\medskip
{\bf a)} The fact (see 4.1.4) that the Shimura (adjoint) filtered Lie $\sg$-crystal attached to $\bigl(M,F^1,\vph,G_{W(k(v))}\bigr)$ is of Borel type.
\smallskip
{\bf b)} 3.4.11 and 3.6.6 applied to the extension to $\FF$ of $\bigl(M,F^1,\vph,G_{W(k(v))}\bigr)$. We get (via 3.5.5 1) and 2)) that there is a Zariski dense set $\Ml$ of $\FF$-valued points of $G_{V_0}^{\rm ad}$ such that ${\got l}_{g_1}$ is (see 3.11.6.1) Shimura-ordinary and has the same Newton polygon as $({\rm Lie}(G^{\rm ad}_{V_0}),\vph\otimes 1)$, $\forall g_1\in G_{V_0}^{\rm ad}(V_0)$ lifting an element of $\Ml$. 
\smallskip
{\bf c)} As the intersection $U(\FF)\cap\Ml$ is non-empty, from b) and 4.2.5 we get that there are $\FF$-valued points of $\Mn^0_{\FF}$ whose attached Shimura Lie $\bar\sg$-crystals are Shimura-ordinary and have ${\rm Lie}_G(\tau)$ as their formal isogeny type; we recall that ${\rm Lie}\bigl(Z(G_{V_0})\bigr)$ contributes only with slopes $0$ to the formal isogeny type of $\bigl({\rm Lie}(G_{V_0}),g_1\vph_0\bigr)$, $\forall g_1\in G_{V_0}^{\rm ad}(V_0)$.
\smallskip
{\bf d)} The first sentence of 3.2.5 (applied to the class $Cl(M_0,\vph_0,G_{V_0})$) which together with c) implies that the $\FF$-valued points of $\Mu_1$ have ${\rm Lie}_G(\tau)$ as the formal isogeny type of their attached Shimura Lie $\bar\sg$-crystals. The role of $({\got g},\vph_0)$ in 3.2.5 is that of a Shimura-ordinary Lie $\sg$-crystal (attached to a Shimura-ordinary $\sg$-crystal), cf. 3.4.3.0, end of 3.5.5 and 3.7.2 (see also 3.11.6 A) and B)).
\smallskip
{\bf e)} The second sentence of 3.2.5 (applied to $Cl(M_0,\vph_0,G_{V_0})$) which (together with the fact that $\Mu_0$ is an open subscheme of $\Mu_1$) implies $\Mu_0=\Mu_1$.
\medskip
It is worth pointing out that in fact the use of 3.6.6 can be avoided, cf. 3.4.14. To define $\Mu$ it is enough to define $\Mu_{\FF}$.  We take $\Mu_{\FF}\cap\Mn^0_{\FF}$ to be $\Mu_1=\Mu_0$. 
From 3.1.0 a) and c) (applied to $Cl(M_0,\vph_0,G_{V_0})$) and c) we also get that the Shimura $\bar\sg$-crystal attached to an $\FF$-valued point of $\Mn_{\FF}^0$ is Shimura-ordinary iff it factors through $\Mu_1$.
\medskip
{\bf 4.2.7. A claim.} The fact that the $\FF$-valued points of $\Mu_1$ give birth to Shimura-ordinary $\bar\sg$-crystals having precisely $\tau$ as their formal isogeny type results from the following Claim:
\medskip
{\bf Claim.} {\it Under the isomorphism $\rho_0$ of 4.2.3, the slopes of $(M,\vph)$ together with their multiplicities,
and the slopes of the Shimura $\bar\sg$-crystal attached to an element of $\Mu_1(\FF)$ together with their multiplicities, are computed using $p$-adic valuations of the eigenvalues (and their multiplicities) of
two semisimple elements $s_1$ and respectively $s_2$ of $G_{K_0}(K_0)$ which are $G^{\rm ad}_{V_0}(K_0)$-conjugate, and so when viewed as semisimple elements of $G_{V_0}(\overline{K_0})$ are $G_{V_0}(\overline{K_0})$-conjugate. Here eigenvalues are computed as endomorphisms of $M_0\otimes_{V_0} K_0$.}
\medskip
{\bf Proof:}  As the homomorphism $G_{V_0}(\overline{K_0})\to G^{\rm ad}_{V_0}(\overline{K_0})$ is surjective, we just have to deal with the $G_{V_0}^{\rm ad}(K_0)$-conjugacy part. For being able to use the previous notations, we assume $y_0$ factors through $\Mu_1$. A sufficiently high power $n$ of $\vph$ (resp. of $\vph_0$) acts diagonally w.r.t. a $W(\FF)$-basis of $M\otimes_{W(k(v))} W(\FF)$ (resp.  of $M_0$), producing (for instance, via 3.11.2 B or 2.2.24.1) a semisimple element $s_1$ (resp. $s_2$) of $G_{K_0}(K_0)$ (we use the identification of 4.2.3 achieved by $\rho_0$). Choosing $n$ big enough and a positive, integral multiple of $d[k_{A_0}:\FF_p]$, where $k_{A_0}$ is the finite field such that $A_{0}$ and its family of Hodge cycles $(w_{\al})_{\al\in\Mj^\prime}$ are defined over $W(k_{A_0})$, we can define such semisimple elements without referring to 3.11.2 B or to 2.2.24.1, cf. 2.2.24. 
\smallskip
In what follows it is irrelevant if we use 3.11.2 B (equivalently 2.2.24.1) or 2.2.24 in order to define them provided, in case we define them via 2.2.24, we allow such semisimple elements to be replaced by a product of them with elements of $G_{V_0}(V_0)$ which do not change the $p$-adic valuations of the eigenvalues involved. However, to be short we assume $s_1$ and $s_2$ are defined via 3.11.2 B. 
\smallskip
As the cocharacter of $G^{\rm ab}_{K_0}$ obtained by composing $\mu_{K_0}$ with the natural epimorphism $m_{qu}:G_{K_0}\twoheadrightarrow G^{\rm ab}_{K_0}$, is the same as the one obtained by composing $\mu_0$ with $m_{qu}$ (as $\mu_{K_0}$ and $\mu_0$ are $G_{K_0}(K_0)$-conjugate), the images of $s_1$ and $s_2$ in $G^{\rm ab}_{K_0}(K_0)$ are the same (cf. the constructions in 3.11.2 B and 4.2.3). 
\smallskip
Moreover, we can assume that the images of $s_1$ and $s_2$ in $G^{\rm ad}_{K_0}(K_0)$ are $G^{\rm ad}_{K_0}(K_0)$-conjugate: this is a consequence of 4.2.6 b) and c) and of 3.11.6 A) (cf. also 4.2.5). So $s_1^{q}$ and $s_2^{q}$ are $G_{K_0}^{\rm ad}(K_0)$-conjugate. Here $q$ is the order of the isogeny $G^{\rm der}\to G^{\rm ad}$. This proves the Claim.
\medskip
{\bf 4.2.8. End of the proof.} The formal isogeny types $\tau$ and ${\rm Lie}_G(\tau)$ depend only on $f$ and $v$. So what we got for $\Mn^0_{\FF}$ remains true for all other connected components of $\Mn_\FF$. So (cf. also 4.2.6-7) the subscheme $\Mu$ of $\Mn_{k(v)}$ defined by 4.2.1 a) exists and is open, is also defined by 4.2.1 b), and it is dense in $\Mn_{k(v)}$. Moreover, from the end of 4.2.6, we also get that 4.2.1 c) holds. From Fact 6 of 2.3.11 we get that $\Mu$ is $G(\AA^p_f)$-invariant. This ends the proof of 4.2.1. 
\smallskip
Another way of passing from one connected component of $\Mn_{\FF}$ to all others is to use the following fact:
\medskip
{\bf 4.2.8.1.* Fact.}
{\it The connected components of $\Mn_\FF$ are permuted transitively by $G(\AA^p_f)$ (cf. [Va2, 3.3.2] and the existence --see [Va3]-- of smooth toroidal compactifications of $\Mn$).}
\medskip
Warning: If $(G,X)$ is of compact type, then this Fact is a consequence of 2.3.3.1 and of [Va2, 3.3.2].
\medskip
{\bf 4.2.9. Definition.} The points of $\Mu_{k(v)}$ (resp. of $\Mu_{k(v)}/\tilde H_0$) with values in fields are called $G$-ordinary (or $G_{W(k(v))}$-ordinary or Shimura-ordinary) points of $\Mn_{k(v)}$ (resp. of $\Mn_{k(v)}/\tilde H_0$). Here $\tilde H_0$ is an arbitrary compact (not necessarily open) subgroup of $G(\AA_f^p)$.
\medskip
{\bf 4.2.10. Corollary.} {\it The Lie algebra underlying the Shimura adjoint Lie $\sg_k$-crystal $LIE$ attached to a $k$-valued point of $\Mn_{k(v)}$ is naturally identified, up to inner isomorphisms (defined by elements of $G^{\rm ad}_{W(k)}(W(k))$, with ${\rm Lie}(G^{\rm ad}_{W(k)})$. Under this identification, we have:
\medskip
-- the automorphism class (see 2.2.11) of $LIE$ is given by the $\sg_k$-linear Lie automorphism of ${\rm Lie}(G^{\rm ad}_{W(k)})={\rm Lie}(G^{\rm ad}_{\ZZ_p})\otimes_{\ZZ_p} W(k)$ acting trivially on ${\rm Lie}(G^{\rm ad}_{\ZZ_p})$;
\smallskip
-- the filtration class (see 2.2.11.1) of $LIE$ is defined by $F^0({\rm Lie}(G^{\rm ad}_{\ZZ_p})\otimes_{\ZZ_p} W(k))$ (i.e. is defined by the composite of $\mu_{W(k)}$ with the epimorphism $G_{W(k)}\to G_{W(k)}^{\rm ad}$).}
\medskip
{\bf Proof:}
This is nothing else but 4.2.5 (or 4.2.4). The fact that we can assume $k$ is just perfect is implied by the convention 2.3.9.2 and so by the fact that $G_{W(k)}$ is split; so the isomorphism $\rho_0$ of 4.2.3 exists as well in this case (as it can be checked using standard properties of the Galois cohomology).
\medskip\smallskip
{\bf 4.3. Cyclic factors and (refined) Lie stable $p$-ranks.}
\medskip
{\bf 4.3.1. Cyclic factors.} In what follows we apply 3.9-10 to the context of 4.2. Let 
$$G^{\rm ad}_{\ZZ_p}=\times_{i\in\Mh}G_i$$ 
be the product decomposition of $G^{\rm ad}_{\ZZ_p}$ in $\ZZ_p$-simple factors. So any ${G_i}_{\QQ_p}$ is a $\QQ_p$-simple, adjoint group. We consider a disjoint decomposition
$$\Mh=\Mh^{\rm c}\cup\Mh^{\rm nc}$$
defined by: $i\in\Mh$ belongs to $\Mh^{\rm c}$ iff the composite $\mu(i)$ of $\mu$ with the natural epimorphism $G_{W(k(v))}\twoheadrightarrow G_{iW(k(v))}$ is trivial. $\Mh$ is the empty set iff $\dim_{\CC}(X)=0$, i.e. iff $G$ is a torus (cf. axiom [Va2, 2.3 (SV2)]).
\smallskip
We fix a point 
$$y:{\rm Spec}(k)\to\Mn_{k(v)}.$$ 
The  Shimura adjoint Lie $\sg_k$-crystal $({\rm Lie}(G^{\rm ad}_{W(k)}),\vph_y)$ attached to it is a product 
$$\prod_{i\in\Mh} \bigl({\rm Lie}({G_i}_{W(k)}),\vph_y\bigr),$$
 each simple factor of it being a cyclic Shimura adjoint Lie $\sg_k$-crystal. In other words, $\forall i\in\Mh$ $\vph_y$ permutes (cyclically) transitively the Lie algebras of the $B(k)$-simple factors of ${G_i}_{B(k)}$; this is a consequence of 4.2.10 and of the fact that $G_i$ is the Weil restriction of an absolutely simple, adjoint group over an unramified, finite extension of $\ZZ_p$. We call these factors the cyclic adjoint factors attached to the point $y$; so in what follows we refer to the $i$-th cyclic adjoint factor attached to $y$ as well as (cf. 3.9.1 and 3.9.6) to the $i$-th Faltings--Shimura--Hasse--Witt (adjoint) shift or map attached to $y$. Similarly, let
$$\bar\psi(i):{\rm Lie}({G_i}_k)\to {\rm Lie}({G_i}_k)$$
be the $i$-th Faltings--Shimura--Hasse--Witt adjoint map of the extension of $(M,\vph,G_{W(k(v))})$ to $k$. Based on 4.2.10, we always write the $i$-th Faltings--Shimura--Hasse--Witt adjoint map of $y$ in the form $g_y\bar\psi(i)$, with $g_y\in G_i(k)$ acting on ${\rm Lie}({G_i}_k)$ by inner conjugation; accordingly, by an inner isomorphism between two such $i$-th Faltings--Shimura--Hasse--Witt adjoint maps attached to $k$-valued points of $\Mn_{k(v)}$ we always mean an isomorphism defined by an element of $G_i(k)$.
\smallskip
We also speak about the cyclic Lie factors attached to $y$ (cf. def. 3.10.1). We naturally extend this terminology to points $z: {\rm Spec}(W(k))\to\Mn$ as well as to the context of the quotient $\Mn/H_0$. 
\smallskip
Warning: in the context of points of $\Mn/H_0$ with values in perfect fields, the adjoint groups whose Lie algebras are defining such cyclic factors are not necessarily products of absolutely simple ones and so (at least theoretically) we have to work with the notion of weakly cyclic factors (cf. 3.10.0.1); in other words, referring to 2.3.10, we do not know when $\tilde G_{W(k)}^{\rm ad}$ is a product of adjoint groups, each one of them with the property that the Lie algebras of simple factors of its extension to $W(\bar k)$ are permuted transitively by $\vph$. 
\medskip
{\bf 4.3.1.1. Notations and simple facts.}
We write
$${G_i}_{W(k)}=\times_{j_i\in\Mh_i}G_{j_i},$$ 
with $G_{j_i}$ an absolutely simple, adjoint $W(k)$-group and with $\Mh_i$ a set of indices. So
$$G^{\rm ad}_{W(k)}=\times_{\scriptstyle i\in\Mh} \times_{\scriptstyle j_i\in\Mh_i} G_j=\times_{\ell\in I_p(G^{\rm ad})}G_\ell,$$ 
where $I_p(G^{\rm ad})$ is the (disjoint) union of $\Mh_i$, $i\in\Mh$. $\abs{I_p(G^{\rm ad})}$ does not depend on $p$ (or on $k$ or $y$); however, it is appropriate to use the right lower index $p$ in the definition of $I_p(G^{\rm ad})$.
\smallskip
$\vph_y$ permutes the factors ${\rm Lie}(G_\ell)$, $\ell\in I_p(G^{\rm ad})$, achieving a permutation $\ga_p(G^{\rm ad})$ of $I_p(G^{\rm ad})$; it does not depend on $y$, cf. 4.2.10. $\ga_p(G^{\rm ad})$ is a product of $\abs{\Mh}$-disjoint cyclic permutations $\ga_i$ ($i\in\Mh)$, where $\ga_i$ is the cyclic permutation of $\Mh_i$ defined logically by the action of $\vph_y$ on ${\rm Lie}({G_i}_{B(k)})$. For $i\in\Mh$, let
$$d_i:=\abs{\Mh_i}\in\NN.$$
From the structure of simple, adjoint groups over fields (see [Ti1]) we get:
$$G_i={\rm Res}_{W(\FF_{p^{d_i}})/\ZZ_p} G^i,$$
with $G^i$ an absolutely simple, adjoint $W(\FF_{p^{d_i}})$-group.
\smallskip
For what follows we refer to 3.11.2-3. 
\medskip
{\bf Fact.} {\it The $T$-degree (resp., in case $G$ is not a torus, the $A$-degree) of definition of any Shimura-ordinary $\sg_{\bar k}$-crystal attached to a $\bar k$-valued $G$-ordinary point of $\Mn_{k(v)}$ is $d_T:=[k(v^{\rm ab}):\FF_p]$ (resp. is $d_A:=[k(v^{\rm ad}):\FF_p]$).}
\medskip
{\bf Proof:} We have a logical variant of 4.2.3 for any morphism $z:{\rm Spec}(W(k))\to\Mn$, cf. the proof of 4.2.10. So, not to introduce extra notations, we can assume $k=\FF$. The part involving $T$-degrees is a consequence of their definition of 3.11.3 and of the $G_{V_0}(K_0)$-conjugacy of the cocharacters $\mu_{K_0}$ and $\mu_0$ of 4.2.3 (cf. also the definition of $E(G^{\rm ab},X^{\rm ab})$). 
\smallskip
The $A$-degree of definition of $(M,\vph,G)$ is $[k(v^{\rm ad}):\FF_p]$, as it can be checked easily starting from very definitions (see Fact of 3.11.2 D and its logical adjoint version; see also [Mi3, 4.6-7]). So the part involving $A$-degrees results from 3.11.6 A) and 4.2.5-6 (see b) of 4.4.1 2) below for a second approach). This ends the proof.
\medskip
For $i\in\Mh^{\rm nc}$, let $k(v_i)$ be the subfield of $k(v)$ such that $B(k(v_i))$ is the field of definition of the $G_i(B(k(v)))$-conjugacy class of $\mu(i)_{B(k(v))}$ (the initial group being $G_{i\QQ_p}$).
\smallskip
For $i\in\Mh^{\rm nc}$, we define as in 2.2.16.5 and 3.11.3, an $A$-degree of definition 
$$d_A(i):=[k(v_i):\FF_p].$$
As in the above proof, it is the $A$-degree of definition of the $i$-th cyclic adjoint factor attached to any $\bar k$-valued $G$-ordinary point of $\Mn_{k(v)}$.
\smallskip
If $G$ is a torus, it is convenient to define $d_A:=1$. So always
$d_A$ is the least common multiple of $d_A(i)$, $i\in\Mh^{\rm nc}$, while $d=[d_A,d_T]$. We refer to $d_A(i)$ as the $i$-th $A$-degree of definition of (the Shimura-ordinary type produced by) the SHS $(f,L_{(p)},v)$. We also denote $d_A$ (resp. $d_T$) by $d_A(v)$ (resp. by $d_T(v)$).
\medskip
{\bf 4.3.1.2. Different invariants.} We denote by $R_y$ (resp. by $R^r_y$) the Lie (resp. the refined Lie) stable $p$-rank of $({\rm Lie}(G_{W(k)}^{\rm ad},\vph_y)$ (cf. 3.9.1, 3.9.1.1, 3.9.4 and 3.9.6). We refer to it as the Lie (resp. the refined Lie) stable $p$-rank of $y$.
\smallskip
We denote by $R(\tau)$ (resp. $R^r(\tau)$) the Lie (resp. the refined Lie) stable $p$-rank of the Shimura adjoint Lie $\sg$-crystal attached to $\bigl(M,\vph,G_{W(k(v))}\bigr)$ of 4.1.1.2.
\smallskip
For $i\in\Mh^{\rm nc}$, we denote by $\vep_i(v^{\rm ad})$ (resp. by $a_i(v^{\rm ad})$) the $\vep$-type (resp. the $a$-invariant) of the $i$-th cyclic adjoint factor of the extension of ${\got C}_{(f,v)}$ to $\FF$ (cf. 3.10.4 and 3.10.8.1). The concentrated $\vep$-type of $\vep_i(v^{\rm ad})$ is denoted by $\vep_i^{\rm c}(v^{\rm ad})$ (see 3.10.5.1). Let
$$a(v^{\rm ad}):=\sum_{i\in\Mh^{\rm nc}} a_i(v^{\rm ad}).$$
We refer to it as the $a$-invariant of $(f,L_{(p)},v)$. 
\smallskip
For $i\in\Mh^{\rm nc}$, let
$$D_i^+(v^{\rm ad})$$ 
(resp. $D_i^-(v^{\rm ad})$) be the positive (resp. negative) $p$-divisible group of the $i$-th cyclic adjoint factor of the extension of ${\got C}_{(f,v)}$ to $\FF$ (cf. 3.11.6.2 2)).  
\medskip
{\bf 4.3.2. The refined $a$-invariant.} $R_y$ is a non-negative integer. Fixing a bijection 
$${f_\Mh}:\Mh\tilde\to S(1,\abs{\Mh}),$$ 
$R^r_y$ is a sequence of non-negative integers, indexed by elements of $S(1,\abs{\Mh})$. These elements ``count" the $\QQ_p$-simple factors of $G^{\rm ad}_{\QQ_p}$. The $n$-th number of the sequence is the Lie stable $p$-rank of the cyclic Shimura adjoint Lie $\sg_k$-crystal $\bigl({\rm Lie}(G_{f^{-1}_\Mh(n)})\otimes_{\ZZ_p} W(k),\vph_y\bigr)$ divided by $d_{f_{\Mh}^{-1}(n)}$ (here $n\in S(1,\abs{\Mh})$). For $i\in\Mh^{\rm c}$, let $a_i(v^{\rm ad}):=0$. Let
$$a^r(v^{\rm ad}):=\bigl(a_{f^{-1}_{\Mh}(1)}(v^{\rm ad}),...,a_{f^{-1}_{\Mh}(\abs{\Mh})}(v^{\rm ad})\bigr);$$
we refer to it as the refined $a$-invariant of $(f,L_{(p)},v)$.
\medskip
{\bf 4.3.3. Remark.} The set $SLSp-rk(\Mn_{(k(v))})$ of Lie stable $p$-ranks of geometric points of $\Mn_{k(v)}$ is not necessarily a set of consecutive non-negative integers. For instance, if $\Mn=\Mm$ the set of values of such Lie stable $p$-ranks is $\{1,3,6,\ldots,{e(e+1)\over 2}\}$ (we recall that $\dim_\QQ(W)=2e$).
\medskip
{\bf 4.3.4. Theorem.} {\it The following statements are equivalent:}
\medskip
\item{a)} $y$ {\it factors through\/} $\Mu$;
\smallskip
\item{b)} $R_y=R(\tau)$;
\smallskip
\item{c)} $R^r_y=R^r(\tau)$.
\medskip
{\bf Proof:} This is a direct consequence of 3.9.2,  3.9.4, 4.2.1 and 4.2.4, applied to any connected component of $\Mn_\FF$ (cf. also 3.9.6).
\medskip
{\bf 4.3.5. Corollary.} $R(\tau)$ (resp. $R^r(\tau)$) {\it is the greatest number (resp. the greatest sequence) among all numbers $R_{\tilde y}$ (resp. among all sequences of numbers $R^r_{\tilde y}$) associated to points $\tilde y\in\Mn_{k(v)}(\tilde k))$ (with $\tilde k$ a perfect field).}
\medskip
{\bf Proof:} This is a consequence of 4.2.1 c), 4.2.4 and 3.9.5.
\medskip
{\bf 4.3.6. Corollary.} {\it $\Mn_{k(v)}\bsl \Mu$ is of pure codimension 1 in $\Mn_{k(v)}$.}
\medskip
{\bf Proof:} The proof of this is entirely analogous to the one in the case of the ordinary locus of $\Mm_{k(v)}$. If $\dim_{\CC}(X)=0$, then the Corollary is obvious: we have $\Mu=\Mn_{k(v)}$ (cf. 4.2.3.1).
We assume now $\dim_{\CC}(X)\ge 1$. 4.3.4-5 assure us (to be compared also with 4.3.7 below) that $\Mn_{k(v)}\bsl \Mu$ is of pure codimension 1, if it is non void. 
\smallskip
To see that $\Mn_{k(v)}\bsl \Mu$ is non void, we first treat the case when $\Mn$ has the completion property, just to emphasize how 3.6.15 A can be applied in practice.
So, assuming that $k=\bar k$ and that $\Mn$ has the completion property, the fact that $\Mn_{k(v)}\bsl \Mu$ is non void can be checked, by just remarking: If 
$(M_y,F^1_z,\vph_y,G_{W(k)})$ is the Shimura-canonical lift of the Shimura $\sg_k$-crystal $(M_y,\vph_y,G_{W(k)})$ attached to a $G$-ordinary point $y\in\Mu(k)$, then $F^1_z({\rm Lie}(G_{W(k)})\ne \{0\}$ (corresponding to the fact $\dim_{\CC}(X)>0$) and there is $g\in G^0_{W(k)}(W(k))$ such that $(M_y,g\vph_y,G_{W(k)})$ is not a $G_{W(k)}$-ordinary $\sg_k$-crystal.  
\smallskip
To argue this let $\ell_0\in I_p(G^{\rm ad})$ be such that  
$F^1_z\bigl({\rm Lie}(G_{\ell_0})\bigr)\ne 0$. Let $T$, $B$ and $B^{\rm opp}$ be as in 4.1.4.1. Let $\om\in G^0_{W(k)}(W(k))$ normalizing $T_{W(k)}$ and such that $wB_{W(k)}w^{-1}=B^{\rm opp}_{W(k)}$.
\smallskip
We can assume that $B_{W(k)}$ normalizes $F^1_z$ and that $\vph_y$ takes ${\rm Lie}(T_{W(k)})$ and ${\rm Lie}(B)\otimes_{\ZZ_p} W(k)$ into themselves (cf. 3.2.3 and either 3.6.15 A or 3.11.1 c)). The Lie stable $p$-rank of the  cyclic Shimura adjoint Lie $\sg_k$-crystal 
$\bigl({\rm Lie}({G_i}_{W(k)}),\om\vph_y\bigr)$ (with $i\in\Mh$ such that $G_{\ell_0}\subset G_{iW(k)}$) is $0$. So $(M_y,\om\vph_y,G_{W(k)})$ is not a $G_{W(k)}$-ordinary $\sg_k$-crystal, cf. F4 of 3.10.7. 
\smallskip
To treat the general case, we use special quadruples (cf. [Va2, 3.2.10]). Let $\tilde T$ be a maximal $\QQ$--torus of $G$ such that:
\medskip
{\bf i)} the Zariski closure of $\tilde T_{\QQ_p}$ in $G_{\ZZ_p}$ is a maximal torus of $G_{\ZZ_p}$ with the property that, $\forall i\in\Mh$, its image in $G_i$, when extended over $W(\FF_{p^{d_i}})$ has a quotient (via an isogeny) which is a product $T_i^1\times T_i^2$, where $T_i^1$ (resp. where $T_i^2$) is itself a product of tori which are $1$ dimensional and non-split (resp. which are such that their extensions to $W(\FF_{p^{2d_i}})$ have $0$ $W(\FF_{p^{2d_i}})$-ranks);
\smallskip
{\bf ii)} over $\RR$ it is the extension of a compact torus by a $1$ dimensional split torus. 
\medskip
The existence of such a torus is obtained by applying [Har, 5.5.3] to $G$: we just need to show the existence of a maximal torus of $G^i$ (of 4.3.1.1) whose Weil restriction from $W(\FF_{p^{d_i}})$ to $\ZZ_p$ is a torus of $G_i$ to each i) applies. But this is a consequence of the following two simple Facts: 
\medskip
{\bf Fact 1.} {\it $\forall i\in\Mh$, $G^i$ is either split or splits over $W(\FF_{p^{2d_i}})$.}
\medskip
{\bf Fact 2.} {\it We consider an absolutely simple, adjoint group $\tilde S$ of classical Lie type over a finite field $\tilde k$ of odd characteristic which splits over the quadratic extension $\tilde k_2$ of $\tilde k$. Then, if it is not a split group of $A_n$ or $D_{2n+1}$ Lie type or a non-split group of $D_{2n}$ Lie type, with $n\ge 2$ (resp. if it so), $\tilde S$ has a maximal torus $\tilde T$ which is isogeneous to a direct product of $1$ dimensional non-split tori (resp. to a product of $1$ dimensional non-split tori and of tori whose ranks over $\tilde k_2$ are $0$).} 
\medskip
{\bf Proofs:} Fact 1 is trivial if $G^i$ is not of $D_4$ Lie type; but this special case is a consequence of the fact that no Shimura variety of $D_4^{\rm mixed}$ type is of Hodge type (see [De2]; see also [Se2, cor. 2 of p. 182]). 
\smallskip
Fact 2 can be proved easily, case by case. If $\tilde S$ is of $C_n$ Lie type, then it is split. So to check the existence of $\tilde T$, we can assume $n=1$; but this case is trivial. If $\tilde S$ is of $B_n$ Lie type, then it is split, and again we can assume $n=1$, which is already covered by the $C_1$ Lie type case. If $\tilde S$ is of $D_n$ Lie type, with $n\ge 4$, then $\tilde S$ is the adjoint group of a special orthogonal group $SO$ of a non-degenerate quadratic form on a $2n$ dimensional $\tilde k$-vector space. If $\tilde S$ is split and $n$ is even or if $\tilde S$ is non-split and $n$ is odd we can take as $\tilde T$ the image in $\tilde S$ of a product of tori of $SO$, each one of them leaving invariant a quadratic form in two variables $x_1$ and $x_2$ which is of the form $x_1^2-a_{22}x_2^2$, with $a_{22}\in \tilde k$ not a square (argument: the quadratic form $y_1y_2+y_3y_4$ is equivalent to the quadratic form $x_1^2-a_{22}x_2^2-x_3^2+a_{22}x_4^2$). If $\tilde S$ is split and $n$ is odd or if $\tilde S$ is non-split and $n$ is even, then using a product of $n-3$ tori of $SO$ of the same type and considering the image in $\tilde S$ of its centralizer in $SO$, the situation gets reduces to the split $D_3$ Lie type case, i.e. to the split $A_3$ Lie type case.  
\smallskip
If $\tilde S$ is split of $A_n$ Lie type, $n\ge 2$, then we can argue at the level of representations: $T_n:={\rm Res}_{\tilde k_{n+1}/\tilde k} \GG_m/\GG_m$, with $\tilde k_{n+1}$ the finite field extension of $\tilde k$ such that $[\tilde k_{n+1}:\tilde k]=n+1$, is (isomorphic to) a maximal torus of $\tilde S$; if $n$ is even then $T_n$ has rank $0$ over $\tilde k_2$ and if $n$ is odd then $T_n$ is isogeneous to a product of a $1$ dimensional, non-split torus with a torus whose rank over $\tilde k_2$ is $0$. If  $\tilde S$ is non-split of $A_n$ Lie type, $n\ge 2$, then we can argue at the level of non-degenerate hermitian forms in a way similar to the $C_n$ Lie type case; see [Ro, th. of p. 243]. This proves the two Facts. 
\medskip
As in [Va2, 3.2.11], we deduce the existence of a special quadruple $(\tilde T,\{\tilde h\},\tilde H_{\tilde T},\tilde v)\hookrightarrow (G,X,H,v)$. We fix some $i\in\Mh^{\rm nc}$. We consider the cocharacter $\mu_i$ of $T_{iW(\FF)}^1\times T_{iW(\FF)}^2$ defined naturally by $\tilde h$ via extensions and logical composites. The $A$-degree of definition $d_A(i)$ of $\bigl({\rm Lie}({G_i}_{W(k)}),\vph\otimes 1\bigr)$ is smaller or equal to the $A$-degree of definition $d_{\tilde z}(i)$ of the $i$-th cyclic adjoint factor $\Mc_i(\tilde z)$ attached to an $W(\FF)$-valued point $\tilde z$ of $\Mn$ factoring through the integral canonical model of $(\tilde T,\{\tilde h\},\tilde H_{\tilde T},\tilde v)$ (this is the same as 3.11.3.1). From the proof of 2.2.18 we get that the Shimura filtered $\bar\sg$-crystal $(\tilde M,F^1(\tilde M),\tilde\vph,G_{W(\FF)})$ attached to $\tilde z$ is cyclic diagonalizable. We first assume that the mentioned two $A$-degrees of definition are not equal. So, if the special fibre of $\tilde z$ factors through $\Mu$, $(\tilde M,\tilde\vph,G_{W(\FF)})$ has two distinct lifts which are cyclic diagonalizable (cf. 2.3.17 and 3.11.1 a) for the second one). But it is an easy exercise to see that this is not possible (it is solved independently in 4.4.13.2-3 below). Contradiction. 
\smallskip
If $\mu_i$ does not factor through $T_{iW(\FF)}^1$, then $d_{\tilde z}(i)>2d_i$; on the other hand, Fact 1 implies that $d_A(i)\le 2d_i$, and so we get $d_A(i)\neq d_{\tilde z}(i)$. So from now on we assume $\mu_i$ factors through $T_{iW(\FF)}^1$. 
This implies (cf. ii); we need to use a formula similar to the one of $\vph^j$ in 4.1.2): $\Mc_i(\tilde z)$ has all slopes equal to $0$. So (cf. formulas of 3.10.6 and 4.2.1 a)) the special fibre of $\tilde z$ does not factor through $\Mu$. So $\Mn_{k(v)}\setminus\Mu$ is non-empty. This ends the proof of the Corollary.
\medskip
{\bf 4.3.7. Examples.}
We concentrate just on one $\ZZ_p$-simple factor 
$G_i$ of $G^{\rm ad}_{\ZZ_p}$. We can assume $\Mh_i=S(1,d_i)$. Let $z:{\rm Spec}(W(k))\to\Mn$ be a lift of $y$ and let
$$
_i{\got C}_z:=\Bigl({\rm Lie}\bigl({G_i}_{W(k)}\bigr),\vph_y,F^0\bigl({\rm Lie}({G_i}_{W(k)})\bigr),F^1\bigl({\rm Lie}({G_i}_{W(k)})\bigr)\Bigr)
$$ 
be the $i$-th cyclic adjoint factor attached to $z$. If $F^1\bigl({\rm Lie}({G_i}_{W(k)}\bigr)=\{0\}$ --i.e. if $i\in\Mh^{\rm c}$-- there is nothing to be done. From now on we assume $i\in\Mh^{\rm nc}$ and that $F^1\bigl({\rm Lie}(G_1)\bigr)\ne 0$. Let 
$$\bar\psi_i:{\rm Lie}({G_i}_k)\to {\rm Lie}({G_i}_k)$$ 
be the $i$-th Faltings--Shimura--Hasse--Witt adjoint map attached to $y$ (so $\bar\psi_i$ is the Faltings--Shimura--Hasse--Witt map of $({\rm Lie}({G_i}_{W(k)}),\vph_y)$). 
\medskip
{\bf Case 1. $_i{\got C}_z$ is without involution.} The $\sg_k^{d_i}$-linear map 
$$_1\overline{\psi_i}^{d_i}:\bigl({\rm Lie}(G_1)/F^0({\rm Lie}(G_1))\bigr)\otimes_{W(k)} k\to\bigl({\rm Lie}(G_1)/F^0({\rm Lie}(G_1))\bigr)\otimes_{W(k)} k$$ 
obtained from $\bar\psi_i^{d_i}$ by passage to quotients, still computes (in the sense of 3.5) ${1\over d_i}$ times the Lie stable $p$-rank of $_i{\got C}_z$.
\smallskip
If moreover $_i{\got C}_z$ is
of $B_\ell$, $C_\ell$ or $D_\ell^\RR$ type, then
this Lie stable $p$-rank is $d_i$ times 
the $f_{\Mh}(i)$-entry of $R^r(\tau)$ iff the $\sg_k^{d_i}$-linear map $_1\overline{\psi_i}^{d_i}$ has a non-zero determinant (cf. the formulas of 3.10.6 i)). So if all simple factors of $(G^{\rm ad},X^{\rm ad})$ are of $B_{\ell}$, $C_{\ell}$ or $D_{\ell}^{\RR}$ type, then the fact that for a point $y\in\Mn_{k(v)}(k)$ we have $R_y=R(\tau)$, can be expressed in terms of some determinants being non-zero. This is very close to the intuition provided by the (usual) Hasse--Witt matrices.
\medskip
{\bf Case 2. $_i{\got C}_z$ is with involution.} So $_i{\got C}_z$ is of $A_{\ell}$ or of $D_{\ell}^{\HH}$ type.  The $\sg_k^{2d_i}$-linear map 
$$_1\overline{\psi_i}^{2d_i}:\bigl({\rm Lie}(G_1)\bigr)/F^0\bigl({\rm Lie}(G_1)\bigr)\otimes_{W(k)} k\to\bigl({\rm Lie}(G_1)\bigr)/F^0\bigl({\rm Lie}(G_1)\bigr)\otimes_{W(k)} k$$ 
obtained from $\bar\psi_i^{2d_i}$ by passage to quotients, still computes (in the sense of 3.5) ${1\over d_i}$ times the Lie stable $p$-rank of $_i{\got C}_z$. 
\smallskip
So in the cases when $_i{\got C}_z$ is of $A_{\ell}$ or of $D_{\ell}^{\HH}$ type, with or without involution, the fact that the Lie stable $p$-rank of this $\sg_k^{m_i}$-linear map $_1\overline{\psi_i}^{m_i}$, with $m_i\in\{d_i,2d_i\}$, is (resp. is not) $d_i$ times the $f_{\Mh}(i)$-entry $n_i$ of $R^r(\tau)$, can be expressed (cf. 4.3.5) in terms of some matrices of rank at most 1 being non-zero (resp. being zero); such matrices (with entries in $k$) are obtained by taking the $n_i$-th exterior power of (cf. 3.9.1.0) $(_1\overline{\psi_i}^{m_i})^{\dim_{\QQ}(G)}$.
As a conclusion: 
\medskip
{\bf Property.} {\it We can always express the fact that for a point $y\in\Mn_{k(v)}(k)$ we have $R_y^r=R^r(\tau)$, in terms of some matrices (indexed by $i\in\Mh^{\rm nc}$) of rank at most 1 being non-zero.}   
\medskip
{\bf 4.3.8. Remarks. 1)} We denote by $p-{\rm Lie}_G(y)$ or by $r^y$ the Lie $p$-rank of $y$, defined as the multiplicity of the slope $-1$ of $({\rm Lie}(G^{\rm ad}_{W(k)}),\vph_y)$. If $p-{\rm Lie}_G(y)\ne 0$, then the Lie stable $p$-rank of $y$ is different from $0$ and the $p$-rank of $A_y:=y^*(\Ma)$ is different from $0$. But the converses of these two implications are not true,
i.e. there are examples (in the case of the first implication, if $y$ is a Shimura-ordinary point, then $G^{\rm ad}_{\RR}$ must have compact factors, cf. the formulas of 3.10.6) when the Lie stable $p$-rank of $y$ or the $p$-rank of $A_y$ is different from $0$ but the Lie $p$-rank of $y$ is $0$.
\smallskip
{\bf 2)} If $f$ is an isomorphism, i.e. if $\Mn=\Mm$, to know $r^y$ is equivalent to know the $p$-rank $r_y$ of $A_y$: we have 
$$r^y={r_y(r_{y+1})\over 2}.$$ 
But this is not true in the general context of Shimura varieties of Hodge type.
\smallskip
{\bf 3)} In the computations of 4.3.7 involving ranks and determinants it is irrelevant which factor $G_1$ of $G_{iW(k)}$ we choose, subject to $F^1({\rm Lie}(G_1))\neq 0$.
\medskip
{\bf 4.3.9. Definitions.} The $\sg_k^{m_i}$-linear map $_1\overline{\psi_i}^{m_i}$, $m_i\in\{d_i,2d_i\}$, defined in 4.3.7, is called the $i$-th Faltings--Shimura--Hasse--Witt reduced map attached to $y$, and the matrix defined by it (via the choice of some $k$-basis) is referred as the $i$-th Faltings--Shimura--Hasse--Witt reduced matrix attached to $y$ and to $1\in\Mh_i$. Exterior powers of such maps and matrices are called modified $i$-th Faltings--Shimura--Hasse--Witt reduced maps or matrices attached to $y$. 
\smallskip
Similarly, for any cyclic Shimura Lie $F$-crystal we define (cf. 3.9.1.1) such maps and matrices. The only difference: $m_i$ can be also $3d_i$ (cf. 3.10.0). It is 3.9.6 which points out that these Faltings--Shimura--Hasse--Witt reduced maps are the same as the ones defined in a semisimple (not necessarily adjoint) context. 
\medskip\smallskip
{\bf 4.4. $G$-canonical lifts.}
\medskip
{\bf 4.4.0. Definition.} A point $y:{\rm Spec}(k)\to\Mn_{k(v)}$ is called a parabolic (resp. a Borel or reductive) point, if there is a lift $z:{\rm Spec}(W(k))\to\Mn$ of it such that the Shimura filtered Lie $\sg_k$-crystal attached to $z$ is of parabolic (resp. of Borel or reductive) type.
\medskip
{\bf 4.4.1. Theorem. 1)} {\it A point $y:{\rm Spec}(k)\to\Mn_{k(v)}$ factors through $\Mu$ iff $y$ is a parabolic point.
\medskip
{\bf 2)} We assume now that $y$ factors through $\Mu$. We have:
\medskip
a) There is a unique lift $z:{\rm Spec}(W(k))\to\Mn$ of $y$ with the property that the Shimura filtered Lie $\sg_k$-crystal attached to $z$ is of parabolic type;
\smallskip
b) If $k=\bar k$, then the Shimura adjoint filtered Lie $\sg_k$-crystal attached to this unique lift $z$ is isomorphic to 
$$\Bigl({\rm Lie}(G^{\rm ad}_{\ZZ_p})\otimes_{\ZZ_p} W(k),\sg\mu({1\over p})\otimes 1, F^0\bigl({\rm Lie}(G^{\rm ad}_{\ZZ_p})\otimes_{\ZZ_p} W(k)\bigr), F^1\bigl({\rm Lie}(G^{\rm ad}_{\ZZ_p})\otimes_{\ZZ_p} W(k)\bigr)\Bigr),$$ 
and so it is cyclic diagonalizable and of Borel type ($\mu$ being as in 4.1);
\smallskip
c) If $k=\bar k$ then the Shimura filtered  Lie $\sg_k$-crystal attached to $z$ is isomorphic to 
$$\Bigl({\rm Lie}(G_{\ZZ_p})\otimes_{\ZZ_p} W(k),\sg\mu({1\over p})\otimes 1,F^0\bigl({\rm Lie}(G_{\ZZ_p})\otimes_{\ZZ_p} W(k)\bigr), F^1\bigl({\rm Lie}(G_{\ZZ_p})\otimes_{\ZZ_p} W(k)\bigr)\Bigr),$$
and so it is cyclic diagonalizable and of Borel type.
\medskip
{\bf 3)} Let $\bigl(M_y,F^1_z,\vph_y,G_{W(k)},(t_\al)_{\al\in\Mj^\prime},\psi_y\bigr)$ be the principally quasi-polarized Shimura filtered $\sg_k$-crystal attached to $z$ above. We have:
\medskip
a) There is a unique cocharacter $\mu_y:\GG_m\hookrightarrow G_{W(k)}$ which
produces a direct sum decomposition $M_y=F^1_z\oplus F^0_z$ such that the elements of the parabolic Lie subalgebra of
${\rm Lie}(G_{W(k)})$ corresponding to non-negative (resp. to non-positive) slopes of $\bigl({\rm Lie}(G_{W(k)}),\vph_y\bigr)$ take $F^1_z$ (resp. $F^0_z$) into itself. We have $[\mu_y]=[\mu_{W(k)}]$;
\smallskip
b) If $k=\bar k$, then the filtered $\sg_k$-crystal $(M_y,F^1_z,\vph_y)$ is cyclic diagonalizable and its isomorphism class does not depend on the choice of $y$. Moreover, the principally quasi-polarized Shimura filtered $\sg_k$-crystal $\bigl(M_y,F^1_z,\vph_y,G_{W(k)},(t_\al)_{\al\in\Mj^\prime},\psi_y\bigr)$ is isomorphic to any sextuple
$\bigl(M_{y_1},F^1_{z_1},\vph_{y_1},G_{W(k)},(t^1_\al)_{\al\in\Mj^\prime},\psi_{y_1}\bigr)$ defined similarly but starting from another $G$-ordinary point $y_1:{\rm Spec}(k)\to \Mu$ factoring through the same connected component of $\Mu_k$ through which $y$ factors;
\smallskip
c) If $k=\bar k$, then there is $n\in\NN$ such that the principally quasi-polarized Shimura filtered $\sg_k^n$-isocrystal $\bigl(M_y[{1\over p}],F^1_z[{1\over p}],\vph_y^n,G_{B(k)},(t_\al)_{\al\in\Mj^\prime},\psi_y\bigr)$ is $1_{\Mj^\prime}$-isomorphic to the principally quasi-polarized Shimura filtered $\sg_k^n$-isocrystal $\bigl(M\otimes_{W(k(v))} B(k),F^1\otimes_{W(k(v))} B(k),\break
(\vph\otimes 1)^n,G_{B(k)},(v_\al)_{\al\in\Mj^\prime},\tilde\psi\bigr)$;
\smallskip
d)* $\bigl(M_y,F^1_z,\vph_y,G_{W(k)},(t_\al)_{\al\in\Mj^\prime},\psi_y\bigr)$ is $1_{\Mj^\prime}$-isomorphic (in the sense of 2.2.9 6)) to 
$$\bigl(M\otimes_{W(k(v))} W(k),F^1\otimes_{W(k(v))} W(k),\vph\otimes 1,G_{W(k)},(v_\al)_{\al\in\Mj^\prime},\tilde\psi\bigr).$$ 
So in c) we can take $n=1$.}
\medskip
{\bf Proofs:} 1) is a direct consequence of 3.1.0 b) and 2.3.17, cf. 4.2.1 c). 
\smallskip
4.2.3-4 were stated for $\FF$. As these results remain obviously true by working with an algebraically closed field containing $\FF$, in what follows we refer to them (with no extra comment) in such an extended context.
\smallskip
a) (resp. b)) of 2) follows from 3.2.8 and 2.3.17 (resp. from 3.11.6 A) and 4.2.5-6). We include here a second approach to b) of 2). We can assume (cf. 4.2.4, 3.4.11 and 3.5.5) that the  Shimura adjoint filtered Lie $\sg_k$-crystal attached to $z$ is 
$$
\Bigl({\rm Lie}(G^{\rm ad}_{\ZZ_p})\otimes_{\ZZ_p} W(k), g(\sg\mu({1\over p})\otimes 1),F^0\bigl({\rm Lie}(G_{\ZZ_p}^{\rm ad})\otimes_{\ZZ_p} W(k)\bigr),F^1\bigl({\rm Lie}(G^{\rm ad}_{\ZZ_p})\otimes_{\ZZ_p} W(k)\bigr)\Bigr), 
$$
with $g\in\tilde P_0(W(k))$, where $\tilde P_0$ is the parabolic subgroup of $G^{\rm ad}_{W(k)}$ such that $\tilde{\got p}_0:={\rm Lie}(\tilde P^0)$ is $W_0\bigl({\rm Lie}\bigl(G^{\rm ad}_{W(k)}\bigr),\sg\mu({1\over p})\otimes 1\bigr)$. Let $\vph_y:=g(\sg\mu({1\over p})\otimes 1)$. $\tilde{\got p}_0$ is also $W_0\bigl({\rm Lie}(G^{\rm ad}_{W(k)}),\vph_y\bigr)$ (this is the same as 3.3.4). Let ${\got n}_0$ be the nilpotent radical of 
$\tilde{\got p}_0$. The quadruples $({\got n}_0,F^1({\got n}_0),\vph_y)$ and $\bigl({\got n}_0,F^1({\got n}_0),\sg\mu({1\over p})\otimes 1\bigr)$ are filtered Lie $\sg_k$-crystals as well as $p$-divisible
objects of $\Mm\Mf_{[0,1]}(W(k))$. We have $F^1({\got n}_0)=F^1({\rm Lie}(G^{\rm ad}_{W(k)}))$, cf. 2.2.12.1 1).
\smallskip
Let ${\got p}^0_0:=W(0)(\tilde{\got p}_0,\sg\mu({1\over p})\otimes 1)$ (resp. ${\got p}_0^1:=W(0)(\tilde{\got p}_0,\vph_y))$. As ${\got p}_0^1\subset F^0({\got g})$, we have $\vph_y({\got p}_0^1)\subset {\got p}_0^1$ and so 
$$\cap_{n\in\NN} \vph_y^n(\tilde {\got p}_0)\subset{\got p}_0^1.\leqno (1)$$ 
We have:
\medskip
{\bf Fact.} {\it ${\got p}_0^0$ and ${\got p}^1_0$ are both Lie algebras of two Levi subgroups $\tilde P_0^0$ and respectively $\tilde P_0^1$ of $\tilde P_0$ (so ${\got p}_0^1$, as a $W(k)$-module, is a direct summand of $\tilde{\got p}_0$). Moreover, in (1) we have equality.} 
\medskip
{\bf Proof:} The second part of the Fact is a consequence of the fact that $({\got p}_0^1,\vph_y)$ is a $\sg_k$-crystal having only slopes $0$. We now consider the Zariski closure $\tilde P_0^1$ in $\tilde P_0$ of the Levi subgroup of $\tilde P_{0B(k)}$ having ${\got p}_0^1[{1\over p}]$ as its Lie algebra; it is a group scheme over $W(k)$. Let $\tilde N_0$ be the unipotent radical of $\tilde P_0$.
\smallskip
For $(\tilde{\got p}_0,\sg\mu({1\over p})\otimes 1)$ the first part of the Fact is obvious, cf. its definition (cf. also 3.11.2 C and 4.1.6). For $(\tilde{\got p}_0,\vph_y)$ this first part follows from the fact that for any $n\in\NN$, we have $\vph_y^n=a_n\vph^n$, with $a_n\in\tilde P_0(W(k))$. In other words, for any $m\in\NN$, taking $n$ big enough, $\vph_y^n(\tilde{\got p}_0)$ mod $p^m$ is on one side included in ${\got p}_0^1$ mod $p^m$ (cf. (1)) and on the other hand it is the inner conjugate of ${\got p}_0^0$ mod $p^m$ by $a_n$ mod $p^m$.  
So ${\got p}_0^1$ mod $p^m$ is the Lie algebra of a Levi subgroup of $\tilde P_{0W_m(k)}$ and so is the Lie algebra of a uniquely determined reductive subgroup $R_m$ of $\tilde P_{0W_m(k)}$. The uniqueness of $R_m$ can be proved in many ways. For instance, we have:
\medskip
{\bf Subfact.} {\it $R_m$ is the subgroup $N_m$ of $\tilde P_{0W_m(k)}$ normalizing ${\got p}_0^1$ mod $p^m$.}
\medskip
The argument for this Subfact goes as follows. As we are in an adjoint context, the arguments of the proof of the Claim of 3.6.18.7.3 referring to $F^{-1}$ and to $T_{W(k)}$ gives us that ${\rm Lie}(N_m)$ is ${\got p}_0^1$ mod $p^m$. So, as $N_m$ contains $R_m$, we get that $N_m$ is smooth, having $R_m$ as its connected component of the origin. So to check that $N_m=R_m$, we just need to show that $N_m$ is connected. We can assume $m=1$. But $N_1$ is normalized by any maximal torus of $R_1$ and so, as $\tilde N_{0k}$ can be identified with the Lie algebra of ${\got n}_0/p{\got n}_0$ and as the group of connected components of $N_1$ can be identified with a finite, \'etale subgroup of $\tilde N_{0k}$, we get $N_1=R_1$.
\smallskip
The Subfact implies $R_m={R_{m+1}}_{W_m(k)}$. So as $R_m$ contains a maximal torus of $G^{\rm ad}_{W_m(k)}$ and as $m\in\NN$ is arbitrary, the Subfact of 2.2.9 3) can be entirely adapted to the sequence $(R_m)_{m\in\NN}$. We get the existence of a Levi subgroup $R_0$ of $G^{\rm ad}_{W(k)}$ whose reduction mod $p^m$ is $R_m$, $\forall m\in\NN$. As ${\rm Lie}(R_0)={\got p}_0^1$, the generic fibres of $R_0$ and $\tilde P_0^1$ are the same. So $R_0=\tilde P_0^1$. This proves the Fact. 
\medskip
Any two Levi subgroups of $\tilde P_0$ are conjugate by an element $g_0\in\tilde P_0(W(k))$: this can be proved in the same manner as the Fact 1 of 2.2.9 3), starting from [Bo2, 14.19] applied over $k$. Changing $\vph_y$ (viewed as a $\sg_k$-linear automorphism of ${\rm Lie}(G_{W(k)}^{\rm ad})$) with $g_0^{-1}\vph_yg_0$ (i.e. --as $\tilde{\got p}_0\subset F^0({\rm Lie}(G_{W(k)}^{\rm ad}))$-- changing it by one which is inner isomorphic to it), we can assume ${\got p}^0_0={\got p}^1_0$. But this implies $g\in\tilde P^0_0(W(k))$ (cf. the uniqueness part of loc. cit. applied over $\overline{B(k)}$). This means that $g$ fixes (under conjugation) the image of $\bar h_s$ in 
${\rm Lie}\bigl(G^{\rm ad}_{W(k)}\bigr)$, $\forall s\in S(1,d)$. We consider the $\sg_k$-linear automorphism $\sg_1:=g(\sg\otimes 1)$ of ${\rm Lie}(G_{W(k)}^{\rm ad})$. Let
$\tilde G^{\rm ad}_{\ZZ_p}$ be the adjoint group over $\ZZ_p$ having  
$$
{\rm Lie}(G^{\rm ad}_{W(k)})^{\sg_1}:=\bigl\{x\in {\rm Lie}(G^{\rm ad}_{W(k)})|\sg_1(x)=x\bigr\} 
$$
as its Lie algebra (as $k=\bar k$, a $\ZZ_p$-basis of
${\rm Lie}(G^{\rm ad}_{W(k)})^{\sg_1}$ is a $W(k)$-basis of ${\rm Lie}(G^{\rm ad}_{W(k)})$).
\medskip
{\bf Claim.} {\it $\tilde G^{\rm ad}_{\ZZ_p}$ is naturally isomorphic to $G^{\rm ad}_{\ZZ_p}$.}
\medskip
{\bf Proof:} Let $\tilde\mu$ be the composite of $\mu$ with the epimorphism $G_{W(k(v))}\twoheadrightarrow G^{\rm ad}_{W(k(v))}$. We recall that $V_0=W(\FF)$. $\tilde P^0_0$ is defined over $\ZZ_p$ (for instance, cf. Fact 2 of 4.1.1.1).
As $g\in\tilde P_0^0(W(k))$, we have $d\tilde\mu\bigl({\rm Lie}(\GG_m)\bigr)\subset {\rm Lie}(\tilde G^{\rm ad}_{V_0})$; so we can assume $k=\FF$. As the map $\tilde P^0_0(\FF)\to\tilde P^0_0(\FF)$ defined by $\tilde g\to\tilde g\bar\sg(\tilde g)^{-1}$ is surjective (cf. [Bo2, 16.4]), we can assume that $g$ mod $p$ is the identity element of $\tilde P_0^0(\FF)$. Argument: as $\tilde P^0_0$ commutes with the image of $\tilde\mu_{W(\FF)}$, when we replace $g(\sg\tilde\mu({1\over p})\otimes 1)$ with $g_1g(\sg\tilde\mu({1\over p})\otimes 1)g_1^{-1}$, with $g_1\in\tilde P^0_0(V_0)$, we get $g_1g\bar\sg(g_1)^{-1}(\sg\tilde\mu({1\over p})\otimes 1)$. 
\smallskip
$g\equiv 1$ mod $p$ implies $\tilde G^{\rm ad}_{\FF_p}=G^{\rm ad}_{\FF_p}$. But two adjoint semisimple groups over $\ZZ_p$ having isomorphic special fibres, are isomorphic. Argument: as the root datum of an adjoint group over $W(\FF)$ can be read out from its special fibre, from [SGA3, Vol. III, 5.3 of p. 314] we get that they are isomorphic over $W(\FF)$; and so we can apply the fact that a torsor of a smooth group scheme over $\ZZ_p$ having an $\FF_p$-valued point is trivial. 
\smallskip
So $\tilde G^{\rm ad}_{\ZZ_p}$ is isomorphic to $G^{\rm ad}_{\ZZ_p}$, under an isomorphism lifting the identity $\tilde G^{\rm ad}_{\FF_p}=G^{\rm ad}_{\FF_p}$ and taking the image $^{\rm ad}T_{\mu}$ of $T_{\mu}$ (of 4.1) in $G^{\rm ad}_{\ZZ_p}$ into the similarly defined (via $\tilde\mu$) torus $^{\rm ad}\tilde T_{\tilde\mu}$ of $\tilde G_{\ZZ_p}^{\rm ad}$; cf. [SGA3, Vol. II, 3.6 of p. 48] for the last part. This proves the Claim.
\medskip
But now $\bigl({\rm Lie}(G^{\rm ad}_{V_0}),\sg_1\tilde\mu({1\over p})\otimes 1\bigr)$ w.r.t. the new $\ZZ_p$-structure defined by ${\rm Lie}(\tilde G^{\rm ad}_{\ZZ_p})$, is exactly
$\bigl({\rm Lie}(G^{\rm ad}_{V_0}),\sg\tilde\mu({1\over p})\otimes 1\bigr)$ with the old $\ZZ_p$-structure defined by ${\rm Lie}(G^{\rm ad}_{\ZZ_p})$. This ends the second proof of b) of 2), as
$\bigl({\rm Lie}(G^{\rm ad}_{\ZZ_p})\otimes_{\ZZ_p} V_0,\sg\tilde\mu({1\over p})\otimes 1,F^0({\rm Lie}(G_{\ZZ_p}^{\rm ad})\otimes_{\ZZ_p} V_0),F^1({\rm Lie}(G_{\ZZ_p}^{\rm ad})\otimes_{\ZZ_p} V_0)\bigr)$ is of Borel type and cyclic diagonalizable (cf. 4.1.4). 
\smallskip
c) of 2) can be obtained entirely as b) of 2): we just have to replace everywhere ${\rm Lie}(G_{\ZZ_p}^{\rm ad})$ by ${\rm Lie}(G_{\ZZ_p})$ and to consider parabolic and Levi subgroups of $G_{W(k)}$ instead of $G_{W(k)}^{\rm ad}$ (Levi subgroups of $G_{W(k)}$ are in one-to-one correspondence to the Levi subgroups of $G_{W(k)}^{\rm ad}$ and so the part of the proof of the above Fact referring to the adjoint context of 3.6.18.7.3 still applies).
\smallskip
For the proof of a) of 3) we first remark that $\mu_y$, if  exists, is unique
($F^0_z$ is uniquely determined, cf. the argument of 3.2.8 for the uniqueness of $F^1_z$) and so, due to the uniqueness assertion, we can assume $k=\bar k$. From Fact 2 of 2.2.9 3) we get: if $\mu_y$ exists then $[\mu_y]=[\mu_{W(k)}]$. So the existence of $\mu_y$ results directly from b) of 2) (i.e. a cocharacter $\mu_y$ of $G_{W(k)}$, subject to $[\mu_y]=[\mu_{W(k)}]$, is uniquely determined by the cocharacter $\tilde\mu_y:\GG_m\to G^{\rm ad}_{W(k)}$ obtained by composing it with --see 4.2.7-- ${q_{\rm qu}}_{W(k)}$, cf. (DER) of 2.2.6 1)).
\smallskip
Using 3.4.11 and 2.3.13.1, for the part of b) of 3) not involving cyclic diagonalizability, we just have to show that for any $g\in P_0^0(W(k))$, there is $h\in G_{W(k)}^0(W(k))$ normalizing $F^1/pF^1$ and such that $g\vph_y=h\vph_yh^{-1}$; here $P_0^0$ is the parabolic subgroup of $G^0_{W(k)}$ whose image in $G^{\rm ad}_{W(k)}$ is $\tilde P_0$. But the existence of $h$ can be obtained following entirely the proof of b) of 2) (i.e. following the above Fact and Claim). In other words, the same arguments as above involving Levi subgroups allow us to assume that the image of $g$ in $G^{\rm ad}_{W(k)}(W(k))$ belongs to $\tilde P^0_0(W(k))$. Using this, we define $\tilde G_{\ZZ_p}$ in the same way we defined $\tilde G^{\rm ad}_{\ZZ_p}$ above; as a reductive group over $\ZZ_p$ is uniquely determined by its special fibre and as a torsor of a reductive group over $\ZZ_p$ which is trivial over $\FF_p$, is trivial, the proof of the Claim applies entirely in the context of the reductive group $\tilde G_{\ZZ_p}$. So a natural variant of the first paragraph after the proof of the Claim, applies. 
\smallskip
It is an easy exercise to see that the Shimura filtered $\sg_k$-crystal $(M_y,F^1_z,\vph_y,G_{W(k)})$ is cyclic diagonalizable. Hint: $\mu_y:\GG_m\to G_{W(k)}$ defines an element $\bar h_1^z\in {\rm End}(M_y)$ by $\bar h^z_1(x)=x$ if $x\in F^1_z$ and $\bar h^z_1(x)=0$ if $x\in F^0_z$; denoting  $\bar h^z_s:=\vph_y^{s-1}(\bar h_1^z)$, $s\in \NN$, we again get that
$\bar h^z_{d+1}=\bar h^z_1$, i.e. $\vph_y^d(\bar h^z_1)=\bar h^z_1$, and the whole argument of the proof of 4.1.2 can be repeated. b) of 2) implies that the Shimura filtered Lie $\sg_k$-crystal attached to $(M_y,F^1_z,\vph_y,G_{W(k)})$ is of parabolic type and so this exercise is (via 3.1.0 b)) as well a particular case of 3.11.1 a).  
\smallskip
To see c) of 3) we write $\vph_y=a_y\mu_y({1\over p})$, with $a_y$ as a $\sg_k$-linear automorphism of $M_y$. $M_y$ and $G_{W(k)}$ get a natural $\ZZ_p$-structure $M_y(\ZZ_p)$ and respectively $\tilde G_{\ZZ_p}$, by considering elements fixed by $a_y$ (see 2.2.9 8)). The degree of definition of $(M_y,F_z^1,\vph_y)$ is $d$ (cf. end of 3.11.2 D and the Fact of 4.3.1.1); so $\mu_y$ can be viewed as a cocharacter of $\tilde G_{W(k(v))}$ and so we can view $\bar h_s^z$'s as belonging to ${\rm Lie}(\tilde G_{W(k(v))})$. Moreover, $t_{\alpha}\in\Mt(M_y(\ZZ_p)[{1\over p}])$, $\forall\alpha\in\Mj^\prime$. $\tilde G_{\QQ_p}$ is (see 4.2.3) an inner form of $G_{\QQ_p}$, defined by a class 
$$\gamma\in H^1({\rm Gal}(\QQ_p),G_{\QQ_p}(\overline{\QQ_p})).$$ 
Using standard arguments (to be compared with [Va2, 5.2.17.2] and the proof of 2.3.13), we get that the image of $\gamma$ in $H^1({\rm Gal}(\QQ_p^{\rm un}),G_{\QQ_p}(\overline{\QQ_p}))$ is 0 (i.e. is the trivial class), where $\QQ_p^{\rm un}$ is the maximal unramified algebraic extension of $\QQ_p$. From the proof of b) of 2) we get that the image $\gamma^{\rm ad}$ of $\gamma$ in $H^1({\rm Gal}(\QQ_p),G_{\QQ_p}^{\rm ad}(\overline{\QQ_p}))$ is $0$; this implies: $\tilde G_{\QQ_p}$ is naturally isomorphic to $G_{\QQ_p}$ (so $G_{\ZZ_p}$ is isomorphic to $\tilde G_{\ZZ_p}$, cf. [Ti2]). 
\smallskip
We now check the following statement: there is $n\in\NN$ a multiple of $d$ such that the extension to $W(\FF_{p^n})[{1\over p}]$ of the representation of $G_{\QQ_p}$ on $L_p^*$ becomes isomorphic to the extension to $W(\FF_{p^n})[{1\over p}]$ of the representation of $\tilde G_{\QQ_p}$ on $M_y(\ZZ_p)[{1\over p}]$ and, under such an isomorphism $\rho(y,n)$, $v_{\al}$ goes into $t_{\al}$, $\forall\al\in\Mj^\prime$, $\tilde\psi$ goes into $\psi_y$, and $\bar h_s$ goes into $\bar h_s^z$, $\forall s\in S(1,d)$. 
\smallskip
First we need (this takes care of the part of the statement not mentioning $\bar h_s$'s):
\medskip
-- $n$ to be a multiple of $d$ such that the image of $\gamma$ in $H^1({\rm Gal}(W(\FF_{p^n})[{1\over p}]),G_{\QQ_p}(\overline{\QQ_p}))$ is $0$ and $G_{B(\FF_{p^n})}$ is split. 
\medskip
It is the proof of b) of 2), which takes care of the part involving
$\bar h_s$'s and $\bar h_s^z$'s; the two obstacles of getting this, i.e. the fact that the map $G(B(\FF_{p^n}))\to G^{\rm ad}(B(\FF_{p^n}))$ is not necessarily surjective and the fact that we do not want to use the fact --it is argued in 3.11.2 C-- that the centralizer in $G_{B(k(v))}$ of all $\bar h_i$'s is connected, are overcome by the fact that the centralizer $C$ of $\bar h_i$'s in the parabolic subgroup of $G_{B(k(v))}$ whose Lie algebra is $W_0({\rm Lie}(G_{B(k(v))}),\vph)$ is reductive (cf. the above part pertaining to Levi subgroups) and splits (as $G_{\QQ_p}$ does) over $B(\FF_{p^n})$. So the above considerations pertaining to $\gamma$ (and its images) apply to a similarly defined class $\gamma_c\in H^1({\rm Gal}(B(k(v))),C(\overline{\QQ_p}))$. So second we need:
\medskip
-- $\gamma_c$ has $0$ image in $H^1({\rm Gal}(B(\FF_{p^n})),C(\overline{\QQ_p}))$.
\medskip
The fact that the extension of $\rho(y,n)$ to $B(k)$ takes $(\vph\otimes 1)^n$ into $\vph_y^n$ is obvious, as $\vph^n$ and $\vph_y^n$ are expressible in terms of $\bar h_i$'s and respectively of $\bar h_i^z$'s. This proves c) of 3).
\smallskip
From 4.1.2 we get that the isomorphism class of the extension of $(M,F^1,\vph)$ to $k$ is determined by the elements $\bar h_i$'s and can be computed working just with rational coefficients (like over $B(k(v))$). Using this and the above proof of c) of 3), we get that the isomorphism class of the cyclic diagonalizable filtered $\sg_k$-crystal $(M_y,F^1_z,\vph_y)$ is the same as of the extension of $(M,F^1,\vph)$ to $k$ and so it is independent of $y$. This ends the proof of b) of 3). 
\smallskip
The proof of part d) of 3) will be  presented in \S 5 as it is related to Milne's
conjecture mentioned in 1.15. 
Here we just remark that we do not need to assume $k=\bar k$ (cf. 2.3.9.2).
\medskip
{\bf 4.4.2. Definition.} The unique lift $z:{\rm Spec}(W(k))\to\Mn$ (resp. $z:{\rm Spec}(W(k))\to\Mn/H_0$) of a $G$-ordinary point $y:{\rm Spec}(k)\to \Mn$ (resp. of a $G$-ordinary point $y:{\rm Spec}(k)\to\Mn/H_0$) with the property that its attached Shimura
filtered Lie $\sg_k$-crystal is of parabolic type is called the $G$-canonical or the Shimura-canonical lift of $y$.
\medskip
{\bf 4.4.2.1. Remark.}  The right translation by an element of $G(\AA_f^p)$ of a $G$-canonical lift, is again a $G$-canonical lift (cf. a) of 4.4.1 2) and Fact 6 of 2.3.11). This allows us to replace $H_0$ in the above definition by any compact subgroup $\tilde H_0$ of $G(\AA_f^p)$.
\medskip
{\bf 4.4.3. Definition.} Let $y:{\rm Spec}(k)\to\Mn_{k(v)}/H_0$ (resp. $z:{\rm Spec}(W(k))\to\Mn/H_0$) be an arbitrary $k$-valued (resp. $W(k)$-valued) point. Let $A_y:=y^*(\Ma_{H_0})$ (resp. $A_z:=z^*(\Ma_{H_0})$). An endomorphism $f_y$ (resp. $f_z$) of $A_y$ (resp. of $A_z$) is said to be a $G$-endomorphism if, as an element of ${\rm End}\bigl(H^1_{\rm crys}(A_y/W(k))\bigr)$, is an endomorphism of the quasi Shimura $\sg_k$-crystal $(M_y,\vph_y,\tilde G_{W(k)})$ attached to $y$ (see end of 2.2.9 1') and 2.3.10). We denote the set of $G$-endomorphisms of $A_y$ (resp. of $A_z$) by ${\rm End}_G(A_y)$ (resp. by ${\rm End}_G(A_z)$). ${\rm End}_G(A_y)$ and ${\rm End}_G(A_z)$ have natural structures of Lie algebras over $\ZZ$.
\smallskip
An automorphism $a_y$ (resp. $a_z$) of $A_y$ (resp. of $A_z$) is said to be a $G$-automorphism if, as an automorphism of $H^1_{\rm crys}(A_y/W(k))$, is an element of $\tilde G_{W(k)}(W(k))$. Notation ${\rm Aut}_G(A_y)$ (resp. ${\rm Aut}_G(A_z)$). ${\rm Aut}_G(A_y)$ and ${\rm Aut}_G(A_z)$ have a group structure. Warnings: 
\medskip
-- in general, ${\rm Aut}_G(*)$, as a set, is not a subset of ${\rm End}_G(*)$; here $*\in \{A_y,A_z\}$;
\smallskip
-- ${\rm End}_G(A_y)$, ${\rm End}_G(A_z)$, ${\rm Aut}_G(A_y)$ and  ${\rm Aut}_G(A_z)$ depend on the choice of $\Ma_{H_0}$.
\medskip
Similarly we speak about $G$-$R_{\ZZ}$--endomorphisms (i.e. about elements of ${\rm End}_G(A_y)\otimes_{\ZZ} R_{\ZZ}$), with $R_{\ZZ}$ an arbitrary $\ZZ$-algebra, or about $G$-$\QQ$--automorphisms.
\medskip
{\bf 4.4.4. Corollary.} {\it Let $y:{\rm Spec}(k)\to \Mu/H_0$ be a $G$-ordinary point and let $z:{\rm Spec}(W(k))\to\Mn/H_0$ be its $G$-canonical lift. Then 
$${\rm End}_G(A_y)={\rm End}_G(A_z)$$ 
and
$${\rm Aut}_G(A_y)={\rm Aut}_G(A_z).$$
In other words, any $G$-endomorphism (resp. any $G$-automorphism) of $A_y$ lifts to an endomorphism (resp. to an automorphism) of $A_z$.}
\medskip
{\bf Proof:} This is a direct consequence of a) of 4.4.1 2). See also 3.1.1.2.
\medskip
{\bf 4.4.5. Corollary.} {\it If $(f,L_{(p)},v,\Mb)$ is a standard PEL situation, then any $G$-canonical lift $z:{\rm Spec}\bigl(W(\FF)\bigr)\to\Mn/H_0$ of a $G$-ordinary point $y:{\rm Spec}(\FF)\to \Mu$ gives birth to a special point 
$z_0:{\rm Spec}(B(\FF))\to {\rm Sh}_{H_0\times H} (G,X)$ (i.e. $A_z$ has complex multiplication).}
\medskip
{\bf Proof:} We view ${\rm End}(A_y)$ as a $\ZZ$-subalgebra of ${\rm End}(M_y)$, with $M_y$ as in 4.4.3. So ${\rm End}_G(A_y)$ is the intersection of the centralizer of $\Mb$ in ${\rm End}(A_y)$ with the normalizer of the $B(\FF)$-vector space of the principal polarization $p_{A_y}$ of $A_y$, obtained from $\Mp_{\Ma}$ by pull back through $y$ and viewed as a cycle in the crystalline context. So, as any abelian variety over a finite field has complex multiplication, the $\QQ$--group defined by invertible elements of ${\rm End}_G(A_y)\otimes_{\ZZ} \QQ$ contains a maximal torus of rank equal to the rank of $G$. So, from the identity ${\rm End}_G(A_y)={\rm End}_G(A_z)$ of 4.4.4, we get that the Mumford--Tate group of (the generic fibre) of $A_z$ is a torus. This proves the Corollary.
\medskip
{\bf 4.4.6. Remark.} In \S 6 and \S 12 we will see that 4.4.5 remains true without assuming that $(f,L_{(p)},v)$ comes from a standard PEL situation $(f,L_{(p)},v,\Mb)$. This will easily imply that $\Mu$ is a locally closed subscheme of $\Mm_{k(v)}$ (cf. \S 13).
\medskip
{\bf 4.4.7. Remark.} We assume $k=\bar k$. With the notations of the proof of c) of 4.4.1 3),
$M_y(\ZZ_p):=\{x\in M_y|a_y(x)=x\}$ is a $\ZZ_p$-structure on $M$ w.r.t. which $\mu_y$ can be viewed as an injective cocharacter $\mu_y:\GG_m\hookrightarrow GL\bigl(M_y(\ZZ_p)\otimes_{\ZZ_p} W(k(v))\bigr)$; we have $t_{\al}\in\Mt(M_y(\ZZ_p))$, $\forall\al\in\Mj^\prime$. So $\bigl(M_y,F^1_z,\vph_y,{(t_\al)}_{\al\in\Mj^\prime}\bigr)$ is isomorphic to
$$(M_y(\ZZ_p)\otimes_{\ZZ_p}W(k),{^\prime F^1_z}\otimes_{W(k(v))} W(k),\sg\mu_y({1\over p})\otimes 1,(t_{\al})_{\al\in\Mj^\prime});$$ 
here we still denote by $\sg$ the $\sg$-linear automorphism of $M_y(\ZZ_p)\otimes_{\ZZ_p} W(k(v))$ fixing $M_y(\ZZ_p)$, while $^\prime F^1_z:=F^1_z\cap M_y(\ZZ_p)\otimes_{\ZZ_p} W(k(v))$. 
\smallskip
Moreover, the Zariski closure in $GL(M_y(\ZZ_p))$ of the subgroup of $GL(M_y(\ZZ_p)[{1\over p}]$ fixing $t_{\al}$, $\forall\al\in\Mj^\prime$, is $G_{\QQ_p}$, and the degree of definition of
$\bigl(M_y,F^1_z,\vph_y,G_{W(k)}\bigr)$ is $d$ (we recall that $k(v)=\FF_{p^d}$). For all these cf. the proof of c) of 4.4.1 3).
\medskip
{\bf 4.4.8. Exercises. 1)} We assume 4.4.6. If $y:{\rm Spec}(\FF)\to \Mu$, prove that there is a unique lift (the $G$-canonical lift of $y$, cf. 4.4.6) $z:{\rm Spec}\bigl(W(\FF)\bigr)\to\Mn$ such that
$A_z$ has complex multiplication. Hint: if $k(v)=\FF_p$ this is well known, cf. 4.6 P1 below and cf. the theory of lifts of ordinary abelian varieties over $\FF$ to abelian varieties over $W(\FF)$ having complex multiplication; if $k(v)\ne\FF_p$, first copy the proof of 2.2.18 and then use 3.4.3.0.
\smallskip
{\bf 2)} Show that the property expressed in 1) is a property of Shimura-ordinary $F$-crystals over perfect fields (we do need to work in the context of lifts of quasi CM type over an algebraically closed field of arbitrary positive characteristic). Hint: use 2.2.18 and the proof of b) of 4.4.1 2). 
\smallskip
{\bf 3)} This is just a particular case of 1). We assume $k=\bar k$. Let $z:{\rm Spec}(W(k))\to \Mn$ be the $G$-canonical lift
of a $G$-ordinary point $y:{\rm Spec}(k)\to \Mu$. Let $\bigl({\rm Lie}(G_{W(k)}^{\rm ad}),\vph_y\bigr)$ be as in 4.3.1. If $W_0({\rm Lie}(G_{W(k)}^{\rm ad}),\vph_y)$ is a Borel Lie subalgebra of ${\rm Lie}(G_{W(k)}^{\rm ad})$, then $z$ is the unique lift of $y$ to $W(k)$, whose attached Shimura filtered $\sg_k$-crystal is a lift of quasi CM type (cf. def. 2.3.17) of the Shimura $\sg_k$-crystal attached to $y$.
Hint: use 3.10.7.
\smallskip
For a theory behind (and a complete solution of) 4.4.8 see 4.4.13.2.
\medskip
{\bf 4.4.9. Remark.} Coming back to 4.4.4, we would like to mention that not always
we have ${\rm End}(A_y)={\rm End}(A_z)$. This is always the case when $k(v)=\FF_p$, cf. properties $P1$ and 4.6 P2 below. But when $k(v)\ne\FF_p$, we can have ${\rm End}(A_y)\mathrel{\mathop\supset\limits_{\ne}} {\rm End}(A_z)$ for some $G$-ordinary points $y$. For instance this is so when $G$ is a torus and the Shimura-ordinary type we get is a 
supersingular type (concrete example: when ${\rm Sh}\bigl({\rm GSp}(W,\psi),S\bigr)$ is the elliptic modular curve, $G$ is a torus and the prime $v$ is such that the Shimura-ordinary type is (1,1)). There are many examples with $G$ a torus such that ${\rm End}(A_y)\ne {\rm End}(A_z)$, and starting from this we can obtain examples with $G$ a reductive group which is not a torus and when not always ${\rm End}(A_y)$ is equal to ${\rm End}(A_z)$.
\smallskip
However, in \S 14, we will see that 
$${\rm End}(A_y)\otimes_{\ZZ} \QQ={\rm End}_G(A_y)\otimes_{\ZZ} \QQ \oplus {\rm End}(A_y)^{\perp_G}\otimes_{\ZZ} \QQ$$ 
and similarly 
$${\rm End}(A_z)\otimes_{\ZZ} \QQ={\rm End}_G(A_z)\otimes_{\ZZ} \QQ \oplus {\rm End}(A_z)^{\perp_G}\otimes_{\ZZ} \QQ,
$$ 
where ${\rm End}(A_y)^{\perp_G}$ (resp. ${\rm End}(A_z)^{\perp_G}$) is the set of endomorphisms of $A_y$ (resp. of $A_z$) with the property that inside ${\rm End}\bigl(H^1_{\rm crys}(A_y/W(k))\bigl[{1\over p}\bigr]\bigr)$ they are perpendicular to ${\rm Lie}(G_{B(k)})$ w.r.t. the trace form on ${\rm End}\bigl(H^1_{\rm crys}(A_y/W(k))\bigl[{1\over p}\bigr]\bigr)$.
\medskip
{\bf 4.4.10. Problem.} Referring again to 4.4.4, study the structure of the Lie algebras ${\rm End}_G(A_y)\otimes_{\ZZ} \QQ$ and ${\rm End}_G(A_y)$ for various $G$-ordinary points $y:{\rm Spec}(k)\to\Mu/H_0$. The most interesting case is when $k$ is a subfield of $\FF$.
\medskip
{\bf 4.4.11. Remark.} Not only $d({\got C})$ but also the numbers $d^\perp({\got C})$ and $gd({\got C)}$ (defined in 3.13.1), with $\got C$ the Shimura filtered $\sg_k$-crystal attached to a $G$-canonical lift $z:{\rm Spec}(W(k))\to\Mn$, do not depend on the choice of $z$ (cf. c) of 4.4.1 3)).
\medskip
{\bf 4.4.12. Remark.} We assume that d) of 4.4.1 3) holds and that $k=\bar k$. Let $V(k)$ be a complete DVR which is a finite, faithfully flat extension of $W(k)$. Let $m_{V(k)}\colon {\rm Spec}(V(k))\to\Mn$ be a morphism. Let $(A_{V(k)},p_{A_{V(k)}}):=m_{V(k)}^*(\Ma,\Mp_{\Ma})$ and let $(t_{\alpha}^{V(k)})_{\alpha\in\Mj^\prime}$ be the family of Hodge cycles with which $A_{V(k)}$ is naturally endowed. From 2.3.13.1 and the density part of 4.2.1, we deduce the existence of a $V(k)$-isomorphism 
$$L_p^*\otimes_{\ZZ_p} V(k)=M\otimes_{W(k(v))} V(k)\tilde\to H^1_{dR}(A_{V(k)}/V(k))$$
taking $\tilde\psi$ into $p_{A_{V(k)}}$ and taking $v_{\alpha}$ into $t_{\alpha}^{V(k)}$, $\forall\al\in\Mj^\prime$.
\medskip
{\bf 4.4.13. The $U$-ordinariness.} 2.2.17.1 raises the question: when a lift of quasi CM type of a Shimura $F$-crystal over a perfect field is unique? To answer this question we start with some definitions. In what follows $k$ is an arbitrary perfect field.
\medskip
{\bf 4.4.13.1. Definitions.} {\bf a)} A potentially cyclic diagonalizable Shimura $\sg_k$-crystal $(\tilde M,\tilde\vph,\tilde G)$ over $k$ is called a $U$-ordinary $\sg_k$-crystal, if its extension to $\bar k$ has a unique lift of quasi CM type. 
\smallskip
{\bf b)} Due to this uniqueness, this lift of quasi CM type is defined (via natural extension) by a direct summand $F^1(\tilde M)$ of $\tilde M$. $F^1(\tilde M)$ itself or the Shimura filtered $\sg_k$-crystal $(\tilde M,F^1(\tilde M),\tilde\vph,\tilde G)$ is called the $U$-canonical lift of $(\tilde M,\tilde\vph,\tilde G)$.
\smallskip
{\bf c)} A $U$-ordinary Shimura $\sg_k$-crystal $(\tilde M,\tilde\vph,\tilde G)$ is called $T$-ordinary if the Lie algebra $\Mt:=W(0)({\rm Lie}(\tilde G),\tilde\vph)$ is abelian (due to the potentially cyclic diagonalizable part implicit in the definition of $U$-ordinariness, this is equivalent --cf. 2.2.19.2-- to: $\Mt$ is the Lie algebra of a maximal torus of $\tilde G$).
\smallskip
{\bf d)} In the same way we defined $G$-ordinary points of $\Mn_{k(v)}$, we define $U$-ordinary and $T$-ordinary points of it, as well as $U$-canonical lifts of $U$-ordinary points of $\Mn_{k(v)}$ (or of $\Mn_{k(v)}/\tilde H_0$) with values in perfect fields.
\smallskip
{\bf e)} a) to c) extend naturally to the context of Shimura (filtered) Lie $\sg_k$-crystals or to the context of generalized Shimura $p$-divisible objects. In particular we speak about $U$-ordinary Shimura Lie $\sg_k$-crystals and about their $U$-canonical lifts (to be compared with 3.11.6.1).
\medskip
We have the following characterization of $U$-ordinary $\sg_k$-crystals:
\medskip
{\bf 4.4.13.2. Proposition.} {\it A Shimura $\sg_k$-crystal $(\tilde M,\tilde\vph,\tilde G)$ over $k$ is $U$-ordinary iff it can be extended to a Shimura filtered $\sg_k$-crystal $(\tilde M,F^1(\tilde M),\tilde\vph,\tilde G)$ whose attached Shimura filtered Lie $\sg_k$-crystal is (see def. 2.2.12 d)) of toric type.}
\medskip
{\bf Proof:} We first assume the existence of an $F^1$-filtration $F^1(\tilde M)$ of $\tilde M$ such that the Shimura filtered Lie $\sg_k$-crystal $\bigl({\rm Lie}(\tilde G),\tilde\vph,F^0({\rm Lie}(\tilde G)),F^1({\rm Lie}(\tilde G))\bigr)$ is of toric type. Let $Z$ be the integral, closed subgroup of $\tilde G$ such that the Lie algebra of its generic fibre is $W(0)({\rm Lie}(\tilde G)[{1\over p}],\tilde\vph)$. We can assume $k=\bar k$. So, the generic fibre of $Z$ has a natural $\QQ_p$-structure $Z_{\QQ_p}$: the Lie algebra of $Z_{\QQ_p}$ is formed by elements of ${\rm Lie}(Z)[{1\over p}]$ fixed by $\tilde\vph$. $Z_{\QQ_p}$ is a reductive group having a maximal torus of the same dimension as the rank of $G_{\CC}$ (cf. 2.2.3 3) and the $\QQ_p$-version of 2.2.9 8)). 
\smallskip
Let $\tilde\mu:\GG_m\to\tilde G$ be the canonical split of $(\tilde M,F^1(\tilde M),\tilde\vph,\tilde G)$. As in the proof of 2.2.18 we get that $F^1(\tilde M)$ is a lift of quasi CM type and that $\tilde\mu_{B(k)}$ factors through the center of $Z$. Let now $F_1^1(\tilde M)$ be another lift of quasi CM type of $(\tilde M,\tilde\vph,\tilde G)$. Let $\tilde\mu_1:\GG_m\to\tilde G$ be the canonical split of $(\tilde M,F^1_1(\tilde M),\tilde\vph,\tilde G)$. Its generic fibre factors through a maximal torus of $Z_{B(k)}$ and so it commutes with $\tilde\mu$. This implies $\tilde\mu=\tilde\mu_1$ (to be compared with Fact 2 of 2.2.9 3)) and so we have $F_1^1(\tilde M)=F^1(\tilde M)$.
\smallskip
To prove the other implication, we can assume $k=\bar k$. We just have to show that for any cyclic diagonalizable Shimura filtered $\sg_k$-crystal $(\tilde M,F^1(\tilde M),\tilde\vph,\tilde G)$ whose Shimura filtered Lie $\sg_k$-crystal is not of toric type, there is a lift $F_1^1(\tilde M)$ of quasi CM type of $(\tilde M,\tilde\vph,\tilde G)$ which is different from $F^1(\tilde M)$. Let $Z_{\QQ_p}$ be as above. It is easy to see (using just topological arguments) that there is $g\in Z_{\QQ_p}(\QQ_p)$ which takes $\tilde M$ onto itself and which does not normalize $F^1(\tilde M)$; as any such $g$ is fixed by $\vph$, it moreover normalizes $F^1(\tilde M)/pF^1(\tilde M)$. But now, from the very def. 2.2.17, we get that $g(F^1(\tilde M))$ is another lift of quasi CM type of $(\tilde M,\tilde\vph,\tilde G)$. This ends the proof.
\medskip
{\bf 4.4.13.2.1. Corollary.} {\it A Shimura $\sg_k$-crystal is $U$-ordinary iff its Shimura adjoint filtered Lie $\sg$-crystal is.}
\medskip
This is a consequence of 4.4.13.2 and of the Fact of 2.2.13.3. 
\medskip
{\bf 4.4.13.3. Examples.} From 4.4.13.2 and 3.1.0 a) we get: any Shimura-ordinary $\sg_k$-crystal is $U$-ordinary and any Shimura-canonical lift is a $U$-canonical lift. This solves all 3 exercises of 4.4.8 (for 4.4.8 1) cf. also 2.3.17). But there are plenty of $U$-ordinary $\sg_k$-crystals which are not Shimura-ordinary; below we present the simplest type of such examples, directly in the adjoint context.
\medskip
{\bf 4.4.13.3.1. New examples.} Let $\tilde G_0$ be a split, simple, adjoint group over $\ZZ_p$ of some $C_n$, $B_n$, $E_7$ or (with $n\ge 4$) $D_n$ Lie type. Let $\tilde G:={\rm Res}_{W(\FF_{p^3})/\ZZ_p} \tilde G_{0W(\FF_{p^3})}$; so $\tilde G_0$ is naturally a subgroup of $\tilde G$. $\tilde G_{W(\FF_{p^3})}$ is a product of three copies $\tilde G_1$, $\tilde G_2$ and $\tilde G_3$ of $\tilde G_{0W(\FF_{p^3})}$. Let $\tilde\mu:\GG_m\hookrightarrow\tilde G$ be an injective  cocharacter factoring through $\tilde G_0$ and such that the pair $(\tilde G,[\tilde\mu])$ is a Shimura group pair; if $\tilde G_0$ is of $D_n$ Lie type, we assume $(\tilde G,[\tilde\mu])$ is of $D_n^{\RR}$ type (i.e. the derived subgroup of the centralizer of $\tilde\mu$ in $\tilde G_0$ is of $D_{n-1}$ Lie type). We consider a maximal torus $\tilde T$ of $\tilde G$ through which $\tilde\mu$ factors and a Borel subgroup of $\tilde G$ containing $\tilde T$ and such that ${\rm Lie}(\tilde B)$ is included in the $F^0$-filtration of ${\rm Lie}(\tilde G)$ produced as usual --see 2.2.8 3)-- by $\tilde\mu$.  Let $w\in \tilde G(W(\FF_{p^3}))$ be a Weyl element (w.r.t. $\tilde T$) such that its component in $\tilde G_3$ is trivial, while its component in $\tilde G_i$ $(i\in\{1,2\}$) takes $\tilde B_{W(\FF_{p^3})}\cap\tilde G_i$ into its opposite (w.r.t. the maximal torus $\tilde T_{W(\FF_{p^3})}\cap\tilde G_i$ of $\tilde G_i$). 
\smallskip
As in 4.1.5 we construct a Shimura adjoint filtered Lie $\sg_{\FF_{p^3}}$-crystal 
$$\bigl({\rm Lie}(\tilde G_{W(\FF_{p^3})}),w(\sigma\tilde\mu({1\over p})\otimes 1),F^0({\rm Lie}(\tilde G_{W(\FF_{p^3})})),F^1({\rm Lie}(\tilde G_{W(\FF_{p^3})}))\bigr)$$ 
which is (to be compared with 4.1.5.2) potentially cyclic diagonalizable, is of toric type but, as its Lie stable $p$-rank is $0$, it is not of parabolic type. So any generalized Shimura filtered $\sg_{\FF_{p^3}}$-crystal ${\got C}$ having it as its attached adjoint is a $U$-canonical lift, without being a Shimura-canonical lift; if $\tilde G$ is of $A_1$ Lie type, then the $U$-ordinary $\sg_{\FF_{p^3}}$-crystal underlying ${\got C}$ is in fact a $T$-ordinary $\sg_{\FF_{p^3}}$-crystal.  
\smallskip
Similar examples can be constructed modeled on some $A_n$, $D_n^{\HH}$, $E_6$ or $D_n^{\rm mixed}$ type. There are plenty of examples involving as well adjoint factors which are included in the corresponding $F^0$-filtration, as well as examples which involve other Weyl elements. See [Va8] for a minute analysis and classification of such examples. 
\medskip
{\bf 4.4.13.4. Remarks.} {\bf 1)} 4.4.13.2 implies that $U$-ordinary points which are not Shimura-ordinary do not show up in the classical setting of Siegel modular varieties. But they do show up in the context of some of the classical standard PEL situations of [Ko2, ch. 5] (for instance, cf. the $C_n$ Lie type part of 4.4.13.3.1 and 4.12.12.6 below; see also 4.12.12.6.6 3) below). To our knowledge, not even one example of a $U$-ordinary point which is not Shimura-ordinary has been previously identified as such. 
\smallskip
{\bf 2)} 4.4.3-4, 4.4.5-6 and 4.4.9-10 extend to the context of $U$-ordinary points; in particular, 4.4.5-6 give us completely new ways to construct abelian varieties having complex multiplication (in connection to 4.4.6 and 4.4.9-10 see \S6 and respectively \S14). 
\smallskip
{\bf 3)} The letter $T$ in 4.4.13.1 stands for toric, while the letter $U$ stands for uniqueness as well as for the way (serpent like) the Weyl elements (see 4.4.13.3.1 for a sample) have to act on different Lie algebras of (some) adjoint groups in order to produce $U$-ordinary (Lie) $\sg_k$-crystals.
\smallskip
{\bf 4)} We do not know when the $U$-ordinariness is an isogeny invariant or when 3.11.5 holds for $U$-ordinary $\sg_{\bar k}$-crystals. Also we do not know when we can recognize a $U$-ordinary (adjoint Lie) $\sg_k$-crystal from its Newton polygon. It seems to us, that the answer would be: only when we are in a Shimura-ordinary context.
\smallskip
{\bf 5)} As in 3.11.3 we define the degree, the $T$-degree and (when appropriate) the $A$-degree of definition of a $U$-ordinary $\sg_k$-crystal $(\tilde M,\tilde\vph,\tilde G)$. Following 3.1.5, the canonical split cocharacter $\tilde\mu$ of its $U$-canonical lift is also called the canonical split cocharacter of $(\tilde M,\tilde\vph,\tilde G)$. It has interpretations similar to the ones of 3.1.4-6: 
\medskip
{\bf INT.} {\it $\tilde\mu$ is the unique cocharacter of $\tilde G$ defining the filtration class of $(\tilde M,\tilde\vph,\tilde G)$ and such that its image fixes the elements of $W(0)({\rm Lie}(\tilde G),\tilde\vph)$ (cf. the proof of 4.4.13.2). It is fixed by any automorphism of $(\tilde M,\tilde\vph,\tilde G)$. Any endomorphism of $(\tilde M,\tilde\vph,\tilde G)$ is as well an endomorphism of its $U$-canonical lift.}
\smallskip
Similarly we define the canonical split cocharacter of any $U$-ordinary Shimura Lie $\sg_k$-crystal ${\got L}$. It has an entirely similar interpretation as of INT.
\medskip
{\bf 4.4.13.5. Remark.} We have the following variant of the Proposition of AE.4:
\medskip
{\bf Corollary.} {\it Let $g\in {\rm Aut}(G,X,H)$. If $a_g$ fixes an $\FF$-valued $U$-ordinary point of $\Mn$, then it fixes its $U$-canonical lift $z:{\rm Spec}(W(\FF))\to\Mn$.}
\medskip
{\bf Proof:} In such a case $g_0$ (defined as in [Va2, p. 495] starting from the non-trivial part of the Shimura adjoint filtered Lie $\sg_{\FF}$-crystal ${\got C}_z^{\rm ad}$ attached to $z$), is an automorphism of ${\got C}_z^{\rm ad}$ itself (see 4.4.13.4 5)) and so (cf. the uniqueness part of 2.3.17) $a_g$ fixes $z$.
\medskip\smallskip
{\bf 4.5. Different stratifications.} 
\medskip
{\bf 4.5.1. Definition.} The stratification of $\Mn_{k(v)}$ in $G(\AA^p_f)$-invariant, reduced, locally closed subschemes indexed and defined by Newton polygons of Shimura Lie $\sg_k$-crystals attached to points ${\rm Spec}(k)\to\Mn_{k(v)}$, is called the canonical Lie stratification.
\medskip
{\bf 4.5.2. Definition.} The stratification of $\Mn_{k(v)}$ in $G(\AA^p_f)$-invariant, reduced, locally closed subschemes, indexed and defined by sequences of length $\abs{\Mh}$ formed by Newton polygons of the cyclic factors (their number is $\abs{\Mh}$) of
Shimura adjoint Lie $\sg_k$-crystals attached to points ${\rm Spec}(k)\to\Mn_{k(v)}$, is called the
refined canonical Lie stratification.
\medskip
{\bf 4.5.3. Remarks. 1)} $\Mu$ is the open, dense stratum of both stratifications of $\Mn_{k(v)}$ defined above (cf. 4.2.1 and 4.3.4). In the refined context, in all that follows we can as well replace $\Mh$ by $\Mh^{\rm nc}$.
\smallskip
{\bf 2)} Let $i\in\Mh$. Let $\pi_i$ be the projector of ${\rm End}(L_p^*[{1\over p}])$ on ${\rm Lie}(G_i)$ having as its kernel the elements of ${\rm End}(L_p^*[{1\over p}])$ perpendicular on ${\rm Lie}(G_i)$ with respect to the trace form on ${\rm End}(L_p^*[{1\over p}])$ (to be compared with [Va2, 4.2] and AE.0). It is fixed by $G_{\QQ_p}$ and so it can be viewed (cf. [Va2, 4.1] and the identifications of [Va2, top of p. 473]) as a $\QQ_p$-linear combination of Betti realizations of Hodge cycles of abelian varieties over fields of characteristic zero obtained by pulling back $\Ma_{H_0}$. We can assume as well that $\pi_i$ is a $\QQ_p$-linear combination of $v_{\al}$'s of 2.3.1. So the $G(\AA^p_f)$-invariant part of 4.5.1-2 follows from this and the paragraph after Fact 4 of 2.3.11.
\medskip
{\bf 4.5.4. $\rho$ $p$-divisible objects.} We assume $k=\bar k$. Any representation 
$$\rho:G_{\ZZ_p}\to GL(N),$$ 
with $N$ a non-zero, free $\ZZ_p$-module of finite rank, allows us to attach to any point $z:{\rm Spec}(W(k))\to\Mn$, a $p$-divisible object
$_\rho{\got C}_z$ 
of $\Mm\Mf_{[a_\rho,b_\rho]}(W(k))$; here $a_\rho,b_\rho\in\ZZ$ are defined by the fact that 
$$\mu_\rho:=\rho_{W(k(v))}\circ\mu:\GG_m\to GL\bigl(N\otimes_{\ZZ_p} W(k(v))\bigr)$$ 
produces a direct sum decomposition 
$$N\otimes_{\ZZ_p} W(k(v))=\oplus_{i\in\Mi_{\rho}}F^i_\rho,$$
with all $F^i_\rho\ne 0$, with $a_\rho,b_\rho\in\Mi_{\rho}$ such that $\Mi_{\rho}\subset S(a_\rho,b_\rho)$, and where $\be\in\GG_m(W(k(v)))$ acts through $\mu_{\rho}$ on $F^i_{\rho}$ as the multiplication with $\be^{-i}$. This goes as follows.
\smallskip
Let $y$ be the special fibre of $z$. We consider a $1_{\Mj^\prime}$-isomorphism (in the sense of 2.2.9 6)) of the quadruple $\bigl(M\otimes_{W(k(v))} W(k),\mu_{W(k)},(v_\al)_{\al\in\Mj^\prime},\tilde\psi\bigr)$ with the quadruple
$\bigl(M_y,\mu_y,(t_\al)_{\al\in\Mj^\prime},p_{A_y}\bigr)$, cf. d) of 4.4.1 3) and 4.4.12. If we do not want to use d) of 4.4.1 3), then we can use $\ZZ_p$-structures as in 2.2.9 8): at the end we get the same thing. 
\smallskip 
Here $(A_z,p_{A_z}):=z^*(\Ma,\Mp_{\Ma})$, $A_y$ is the special fibre of $A_z$ and $M_y:=H^1_{\rm crys}\bigl(A_y/W(k)\bigr)$; $(t_\al)_{\al\in\Mj^\prime}$ is the family of de Rham (crystalline) components of the family of Hodge cycles with which  $A_z$ is naturally endowed. $\mu_y:\GG_m\to GL(M_y)$ is a cocharacter  fixing $t_\al$, ${\forall\al\in\Mj^\prime}$; it produces a direct sum decomposition $M_y=F^1_z\oplus F^0_z$, with $F^1_z$ as the Hodge filtration of $M_y$ defined by $A_z$ and with $\be\in\GG_m(W(k))$ acting through $\mu_y$ on $F^i_z$ as the multiplication with $\be^{-i}$.
\smallskip
We can write the Frobenius
$\vph_y$ of $M_y$, under this isomorphism, as $g(\sg\mu({1\over p})\otimes 1)$ with
$g\in G^0_{\ZZ_p}(W(k))\subset GL(M)(W(k))$. Now we define
$$_\rho{\got C}_z:=\bigl(N\otimes_{\ZZ_p} W(k),(F^i_\rho\otimes_{W(k(v))} W(k))_{i\in\Mi_{\rho}},\vph_{y,\rho}\bigr),$$ 
with 
$$\vph_{y,\rho}:=\rho_{W(k)}(g)(\sg\mu_\rho({1\over p})\otimes 1).$$ 
Obviously, $_\rho{\got C}_z$ is uniquely determined up to inner isomorphism (regardless of all choices). In \S 5 we will see that there is a $p$-divisible object 
$$_{\rho}{\got C}\Mn$$ 
of $\Mm\Mf^\nabla_{[a_\rho,b_\rho]}(\Mn)$ attached to $\rho$ in such a manner that in $z$ it becomes $_\rho{\got C}_z$. $_\rho{\got C}_z$ is called the $\rho$ $p$-divisible
object of $\Mm\Mf_{[a_\rho,b_\rho]}(W(k))$ attached to $z$. 
\smallskip
The $p$-divisible object $_\rho{\got C}\Mn$ defines a  Newton polygon stratification of
$\Mn_{k(v)}$ in $G(\AA_f^p)$-invariant, reduced, locally closed subschemes (cf. \S 5). We call it the $\rho$-stratification of $\Mn_{k(v)}$. As here we do not show the existence of $_{\rho}{\got C}\Mn$, we briefly indicate how we can define the $\rho$-stratification of $\Mn_{k(v)}$ directly, by just sketching the construction of $_{\rho}{\got C}\Mn$. We can work in the \'etale topology of $\Mn_{W(k)}$, around the point through which $y$ factors. In order to apply 2.3.15, we still denote by $z$ the resulting $W(k)$-morphism ${\rm Spec}(W(k))\to \Mn_{W(k)}/H_0$. Let $a:Y={\rm Spec}(R)\to\Mn_{W(k)}/H_0$ be as in 2.3.15. We assume $a$ factors through an affine subscheme of $\Mn_{W(k)}/H_0$ whose $p$-adic completion has a Frobenius lift and it is compatible with the two Frobenius lifts. Let $(M_y\otimes_{W(k)} R^\wedge, F^1_z\otimes_{W(k)} R^\wedge, g_Y(\vph_y\otimes 1),\nabla_{R^\wedge})$ have the analogue meaning of 2.3.15 d). Following the same pattern of the construction of $_\rho{\got C}_z$, we consider the following $p$-divisible object $_\rho{\got C}_R$ of $\Mm\Mf_{[a_\rho,b_\rho]}^\nabla(R)$ 
$$(N\otimes_{\ZZ_p} R^\wedge,(F^i_\rho\otimes_{W(k(v))} R^\wedge)_{i\in\Mi_{\rho}},\rho_R(g_Y)(\vph_{y,\rho}\otimes 1),\nabla_{R^\wedge}\bigr).$$   
Here, as $\nabla_{R^\wedge}$ respects the $G_R$-action, we still denote by $\nabla_{R^\wedge}$ the connection on $N\otimes_{\ZZ_p} R$ it defines naturally (via the Lie homomorphism ${\rm Lie}(G_{\ZZ_p})\to {\rm End}(N)$ associated to $\rho$).
Using Teichm\"uller lifts ${\rm Spec}(W(\FF))\to Y$, we get that the Newton polygon stratification of $R/pR$ defined by $_\rho{\got C}_R$, is the pull back via the special fibre of $a$ of a Newton polygon stratification of the image $IM$ of the special fibre of $a$ in $\Mn_k/H_0$; it is nothing else but the Newton polygon stratification defined by the Newton polygons of $_\rho{\got C}_z$ and so it is the pull back to $IM$ of the $\rho$-stratification of $\Mn_{\FF}/H_0$. Standard arguments (based on 4.5.6 1) below and the part of 2.3.9 A referring to $\FF$) show that in the last sentence we can replace $\Mn_{\FF}/H_0$ by $\Mn_{k(v)}/H_0$.
\medskip
{\bf 4.5.5. Definition.} The intersection of all $\rho$-stratifications of $\Mn_{k(v)}$ is called the absolute stratification of $\Mn_{k(v)}$.
\medskip
{\bf 4.5.6. Exercises. 0)} The $\rho$-stratification of $\Mn_{k(v)}$ is trivial (i.e. formed just by one stratum) iff $\mu_{\rho}$ factors through the center of the image of $\rho_{W(k(v))}$ in $GL(N\otimes_{\ZZ_p} W(k(v)))$. Hint: use 3.2.2 and the proof of 4.3.6.
\smallskip
{\bf 1)} Any $\rho$-stratification is $G(\AA^p_f)$-invariant. Hint: adapt 4.5.3 2).
\smallskip
{\bf 2)} The canonical Lie stratification is the ${\rm Ad}$-stratification, where ${\rm Ad}:G_{\ZZ_p}\to GL({\rm Lie}(G^{\rm ad}_{\ZZ_p}))$ is the adjoint representation.
\smallskip
{\bf 3)} If ${\rm Ad}_i:G_{\ZZ_p}\to GL({\rm Lie}(G_i))$, $i\in\Mh$, is the (irreducible) subrepresentation of ${\rm Ad}$ induced on the $\ZZ_p$-direct summand  ${\rm Lie}(G_i)$ of ${\rm Lie}(G^{\rm ad}_{\ZZ_p})$ (see 4.3.1), then the refined  canonical Lie  stratification is the intersection of all ${\rm Ad}_i$-stratifications, $i\in\Mh$.
\smallskip
{\bf 4)} The absolute stratification is well-defined, i.e. its strata are locally closed subschemes and their number is finite. In fact, always there is a $\rho$ such that the absolute stratification and the $\rho$-stratification coincide. Hints: for the first part use that the set of $G(B(\FF))$-conjugacy classes of semisimple elements obtained as in 2.2.24.1 but working with a $Cl(M_0,\vph_0,G_{W(\FF)})$ as in 4.2.6 is finite; for the second part, starting from a finite number $\rho_1$, $\rho_2$,..., $\rho_m$ of representations of $G_{\ZZ_p}$ such that the absolute stratification is the intersection of the $\rho_i$'s stratifications, $i\in S(1,m)$, use a sum $\sum_{s=1}^m \rho_s(n_s)$, with all $n_s$'s belonging to $\NN$ and such that $n_1<<n_2<<...<<n_m$; here $\rho_s(n_s)$ is the tensor product representation of $\rho_s$ with the one dimensional representation achieving the Tate twist by $\ZZ_p(n_s)$ (variant, if all images of $\rho_s$'s are semisimple, we can equally well work with a  sum $\sum_{s=1}^m n_s\rho_s$). 
\smallskip
{\bf 5)} If $\Mn=\Mm$, the usual Newton polygon and absolute stratifications coincide. Hint: see 4.5.6.1 below.
\smallskip
{\bf 6)} If $\rho_0:G_{\ZZ_p}\hookrightarrow GL(L^\ast_p)$ is the faithful representation introduced in 4.1, then the $\rho_0$-stratification is the stratification of $\Mn_{k(v)}$ defined by
Newton polygons of Shimura $F$-crystals attached to points of $\Mn_{k(v)}$ with values in perfect fields.
\smallskip
{\bf 7)} Show by examples that 5) is not true for an arbitrary SHS $(f,L_{(p)},v)$. Hint: either use products of Shimura varieties of Hodge type or look at the case when $\abs{\Mh^{\rm nc}}\ge 2$.
\smallskip
{\bf 8)} $\Mu$ is an open subscheme of the open, dense stratum of any $\rho$-stratification and so it is the open, dense stratum  of the absolute stratification. Hint: use b) of 4.4.1 2).
\smallskip
{\bf 9)} Show that the $\rho$-stratification depends only on the restriction of $\rho$ to $G^{\rm der}_{\QQ_p}$. Hint: use the part of 4.2.7 involving abelianizations of groups.
\smallskip
All the above Exercises are very easy.
\medskip
{\bf 4.5.6.1. Proposition.} {\it We assume $G^{\rm ad}$ has no simple factor of:
\medskip
-- $A_n$ or $D_n$ Lie type, $n\in\NN$, $n\ge 2$, or of 
\smallskip
-- $B_n$, $n\in\NN$, $n\ge 4$, which has a simple factor over $\QQ_p$ corresponding to an $i\in\Mh^{\rm nc}$ such that $\vep_i^{\rm c}(v^{\rm ad})$ is not an $1$-tuple. 
\medskip
Then the absolute stratification coincides with the refined Lie canonical stratification.}
\medskip
{\bf Proof:} Let $y_l:{\rm Spec}(k)\to\Mn_{k(v)}$, $l=\overline{1,2}$, be two points of the same stratum of the refined canonical Lie stratification of $\Mn_{k(v)}$. We need to show that they belong to the same stratum of any $\rho$-stratification of $\Mn_{k(v)}$. Let $(M_l,\vph_l,G_{W(k)})$ be the Shimura $\sg_k$-crystal attached to $y_l$. Let $m\in\NN$ be such that $\vph_l^m$ acts diagonally on elements of a $B(k)$-basis of $M_l[{1\over p}]$. We can assume that $k=\FF$ and so that the Shimura $\sg_k$-crystals involved are defined over finite fields; so we can always choose $m$ so that 2.2.24 applies: we get a semisimple element $b_l\in G(B(\FF))$, whose eigenvalues, as an endomorphism of $M_l$ are $p$ to powers $m$ times the slopes of $(M,\vph_l)$, multiplicities corresponding. For future references, from now on we still use the general notation $k$ instead of $\FF$ (cf. 2.2.24.1).
\smallskip
 We choose $m$ to be a multiple of 6 times the product of the lengths of the cyclic adjoint factors of the Shimura Lie $\sg_k$-crystal attached to $y_1$ (or $y_2$; 4.2.10 or 4.4.1 2) implies this product does not depend on $y_1$). The key fact is:
\medskip
{\bf Claim.} {\it Replacing $m$ by a suitable multiple of itself, we can assume $b_1$ and $b_2$ are $G(B(k))$-conjugate.}
\medskip
To prove this claim, let $\tilde G$ be the product of the simply connected semisimple group cover of $G^{\rm der}$ and of the connected component of the origin of $Z(G)$. Replacing $m$ by a multiple of itself, $b_l$ is replaced by a positive power of itself; so we can assume $b_l$ lifts to an element $\tilde b_l\in\tilde G(B(k))$. We consider a representation  
$$\tilde G_{B(k)}\to GL(V_{B(k)})$$ 
factoring through $\tilde G_{B(k)}^{\rm der}$ and such that when restricted to $\tilde G^{\rm der}_{B(k)}$ it is a direct sum, indexed by the simple factors of $G^{\rm ad}_{B(k)}$, of the standard  embeddings (see [He, \S 8 of ch. 3]) of the classical groups defined by suitable group covers of these simple factors (these group covers are quotients of $\tilde G^{\rm der}_{B(k)}$). For instance, if we have a factor ${\rm Sp}(\tilde V_{B(k)},\tilde\psi_{B(k)})$ of $\tilde G_{B(k)}^{\rm der}$, with $(\tilde V_{B(k)},\tilde\psi_{B(k)})$ a symplectic space over $B(k)$, then we get the natural representation of ${\rm Sp}(\tilde V_{B(k)},\tilde\psi_{B(k)})$ on $\tilde V_{B(k)}$. Similarly, if we have a simple factor of $G^{\rm ad}_{B(k)}$ of $B_{\ell}$ Lie type, we get the classical orthogonal faithful representation of dimension $2\ell+1$ of the split $SO(2\ell+1)$-group (over $B(k)$). 
\smallskip
$\tilde b_l$ acts diagonally on $V_{B(k)}$; we can choose $\tilde b_l$ so that its eigenvalues, as an automorphism of $V_{B(k)}$, are integral powers of $p$. $\tilde b_l$ also acts diagonally (via inner conjugation) on each simple factor $\Mf$ of the Lie algebra of the image of ${\rm Lie}(\tilde G^{\rm der}_{B(k)})$ in ${\rm End}(\tilde  V_{B(k)})$. With its eigenvalues, as an automorphism of $\Mf$, we construct a Newton polygon $NP_l(\Mf)$. 3.4.1 and 3.4.2.2 imply that the images of $b_l$ in a set of $GL$-groups of Lie algebras of simple factors of $G^{\rm ad}_{B(k)}$ which are permuted transitively by $\vph_l$, up to isomorphism of these groups, are the same. So, the fact that $y_1$ and $y_2$ are points of the same stratum of the refined canonical Lie stratification of $\Mn_{k(v)}$, gets restated in: for each such simple factor $\Mf$ we have
$$NP_1(\Mf)=NP_2(\Mf).\leqno (EQ)$$ 
\indent
{\bf Fact.} {\it Let $C$ be a classical semisimple, split group over $B(k)$ which has a simple adjoint and is not of $A_l$ Lie type with $l\ge 2$ (so $C$ is $SO(n)$ or ${\rm Sp}(n)$, for some adequate $n\in\NN$). Let $b^1$ and $b^2$ be semisimple elements of $C(B(k))$ such that their eigenvalues (under the classical representation of $C$) are integral powers of $p$. We assume that the Newton polygons defined by the eigenvalues of their actions (via inner conjugation) on ${\rm Lie}(C)$ are the same. We have:
\medskip
{\bf a)} if $C$ is of $C_n$ or $B_3$ Lie type, then $b^1$ and $b^2$ are $C(B(k))$-conjugate;
\smallskip
{\bf b)} if $C=SO(2n)$ (i.e. if $C$ is of $D_n$ Lie type) and if the greatest eigenvalue of $b^1$ and the greatest eigenvalue of $b^2$ both have multiplicity greater than $1$, then $b^1$ and $b^2$ are conjugate by a $B(k)$-valued point of the orthogonal group $O(2n)$;
\smallskip
{\bf c)} if $C=SO(2n+1)$ (i.e. if $C$ is of $B_n$ Lie type) and if the greatest eigenvalue of $b^1$ and the greatest eigenvalue of $b^2$ both have multiplicity greater than $1$, then $b^1$ and $b^2$ are $C(B(k))$-conjugate.}
\medskip  
The simple proof of this Lemma is left as an exercise.
\smallskip
Let $m_{\rm iso}=2$ be the exponent of the kernel of the natural isogeny from $\tilde G_{B(k)}$ into the product of $\tilde G^{\rm ab}_{B(k)}$ with the image of $\tilde G^{\rm der}_{B(k)}$ in $GL(V_{B(k)})$. As in 4.2.7, we can assume that the images of $\tilde b_1$ and $\tilde b_2$ in $\tilde G^{\rm ab}_{B(k)}$ are the same. Let $i\in\Mh^{\rm nc}$ be such that $G_i$ is of $B_n$ Lie type. We take $\Mf$ to be included in ${\rm Lie}(G_{iB(k)})$. If the $i$-th cyclic adjoint factor of $y_1$ or of $y_2$ is Shimura-ordinary, then from (EQ) and from 4.4.1 2) and 4.2.1 b) we get that we can assume the images of $\tilde b_1$ and $\tilde b_2$ in the simple factor of $\tilde G^{\rm ad}_{B(k)}$ corresponding to $\Mf$ are the same. If the $i$-th cyclic adjoint factors of $y_1$ and of $y_2$ are not Shimura-ordinary and if $\vep_i^{\rm c}(v^{\rm ad})$ is an $1$-tuple, then (easy exercise; it is implicitly solved by 4.12.12.6.4.1 a) below) we are in a situation where c) applies. So, from (EQ) and a) and c) we conclude: $\tilde b_1^{m^2_{\rm iso}}$ and $\tilde b_2^{m^2_{\rm iso}}$ are $\tilde G(B(k))$-conjugate. This proves the Claim.
\smallskip
The Claim implies: $y_1$ and $y_2$ are points of the same stratum of any $\rho$-stratification of $\Mn_{k(v)}$ (cf. the constructions in 4.5.4). This ends the proof of the Proposition. 
\medskip
{\bf 4.5.6.2. Examples.} {\bf A.} It seems to us that the restriction on $B_n$ Lie types in 4.5.6.1 is not needed. If $i\in\Mh^{\rm nc}$ is such that $G_i$ is of $B_n$ Lie type, then $G^i$ is split and so it is the special orthogonal group of a perfect, symmetric bilinear form on a free $W(\FF_{p^{d_i}})$-module $V(i)$ of rank $2n+1$. So, viewing $V(i)$ as a $\ZZ_p$-module, we consider the representation 
$$\tilde\rho_i^{\rm ad}:G_{\ZZ_p}\to GL(V(i))$$
which factors through $G_i$ inducing the tautological representation of $G_i$ on $V(i)$. From the proof of 4.5.6.1 we get:
\medskip
{\bf Corollary.} {\it If $G_i$ is of some $B_n$ Lie type, $\forall i\in\Mh^{\rm nc}$, then the absolute stratification of $\Mn_{k(v)}$ is the intersection of the $\tilde\rho_i$-stratifications, $i\in\Mh^{\rm nc}$.} 
\medskip
{\bf B.} On the other hand, 4.5.6.1 does not necessarily hold for the case when $G^{\rm ad}$ has simple factors of $A_n$ Lie type, with $n\ge 2$, or of $D_n$ Lie type, with $n\ge 4$. Here is a very simple example, in the abstract context of Shimura $\sg_k$-crystals involving the $A_3$ Lie type: the proof of 4.12.12 (cf. also 4.12.12.6) below, tells us that this abstract context can be easily adapted to the context of a SHS. Let $(\tilde M,F^1(\tilde M),\tilde\vph,GL(\tilde M))$ be a Shimura filtered $\sg_k$-crystal over $k=\bar k$, with $\tilde M$ (resp. $F^1(\tilde M)$) of rank $4$ (resp. of rank $2$). Let $g_1,g_2\in GL(\tilde M)$ be such that:
\medskip
-- $(\tilde M,g_1\tilde\vph)$ has slope $1\over 3$ with multiplicity $3$ and slope $1$ with multiplicity $1$;
\smallskip
-- $(\tilde M,g_2\tilde\vph)$ has slope $2\over 3$ with multiplicity $3$ and slope $0$ with multiplicity $1$.
\medskip
The Lie isocrystals $({\rm End}(\tilde M)[{1\over p}],g_1\tilde\vph)$ and $({\rm End}(\tilde M)[{1\over p}],g_2\tilde\vph)$ are isomorphic. So the Newton polygon of 
$$({\rm Lie}(PSL(M)),g_i\tilde\vph)$$ 
does not depend on $i\in\{1,2\}$. Similar examples can be constructed for any other $A_{2+m}$ (resp. $D_{3+m}$) Lie type, with $m\in 1+2\NN$ (resp. $m\in\NN$). For examples involving the $A_{2m}$ Lie types ($m\in\NN$), we usually need groups over $\ZZ_p$ whose adjoints are simple but not absolutely simple.
\smallskip
Warning: in many situations 4.5.6.1 remains true even if $G^{\rm ad}$ has simple factors of some $A_n$ or $D_n$ Lie types; for instance, this is so if $G^{\rm ad}_{\ZZ_p}$ is a split group and $G^{\rm ad}_{\RR}=SU(1,2)^{\rm ad}_{\RR}$.
\medskip
{\bf 4.5.7. Some sets of Newton polygons.} Let $y:{\rm Spec}(\FF)\to\Mn_{k(v)}$ be an arbitrary point and let $\bigl(M_0,\vph_0,G_{W(\FF)}\bigr)$ be its attached Shimura $\bar\sg$-crystal. Let $NP(\Mn_{k(v)})$ be the set of Newton polygons defined by $\bar\sg$-crystals of the form $(M_0,g\vph_0)$, with
$g\in G^0_{W(\FF)}\bigl(W(\FF)\bigr)$. Let $LNP(\Mn_{k(v)})$ be the set of Newton polygons defined by Shimura adjoint Lie $\bar\sg$-crystals of the form $\bigl({\rm Lie}(G^{\rm ad}_{W(\FF)}),g\vph_0\bigr)$, with $g\in G^0_{W(\FF)}\bigl(W(\FF)\bigr)$. Similarly, we define $RLNP(\Mn_{k(v)})$ as the set of sequences of length $\abs{\Mh}$ of Newton polygons attached
to cyclic factors (their number is $\abs{\Mh}$) of the above Shimura adjoint Lie $\bar\sg$-crystals.
Let
$$N(G^{\rm ad},X^{\rm ad},v^{\rm ad}):=\abs{LNP(\Mn_{k(v)})},$$
and let
$$N_1(G^{\rm ad},X^{\rm ad},v^{\rm ad}):=\abs{RLNP(\Mn_{k(v)})}.$$
We also write $N(\Mn_{k(v)})=N(G^{\rm ad},X^{\rm ad},v^{\rm ad})$ and $RLNP(\Mn_{k(v)})=RLNP(G^{\rm ad},X^{\rm ad},v^{\rm ad})$. 
\medskip
{\bf 4.5.8. Remarks.} {\it {\bf 1)}} The sets $LNP(\Mn_{k(v)})$ and $RLNP(\Mn_{k(v)})$ and so also the natural numbers $N(G^{\rm ad},X^{\rm ad},v^{\rm ad})$ and $N_1(G^{\rm ad},X^{\rm ad},v^{\rm ad})$, depend only on the isomorphism class of the Shimura group pair $(G^{\rm ad}_{\ZZ_p},\mu^{\rm ad})$, where $\mu^{\rm ad}:\GG_m\to G^{\rm ad}_{W(k(v^{\rm ad}))}$ has the similar meaning as $\mu$ of 4.1 but associated to the Shimura variety ${\rm Sh}(G^{\rm ad},X^{\rm ad})$ (cf. b) of 4.4.1 2)). This justifies our notations and shows that the set $RLNP(G^{\rm ad},X^{\rm ad},v^{\rm ad})$ can be defined directly without any reference to some SHS. 
\smallskip
{\it {\bf 2)}} In \S 10 we will use different $\rho$-stratifications of $\Mn_{k(v)}$ to compute the numbers $N(G^{\rm ad},X^{\rm ad},v^{\rm ad})$ and $N_1(G^{\rm ad},X^{\rm ad},v^{\rm ad})$.
\medskip
{\bf 4.5.9. The quasi-affineness property.} The canonical stratification of the special fibre $\Mm_{k(v)}$ of $\Mm$ (cf. [EO] and [Oo3, 7.2]; in [Oo3] is also referred as the EO stratification) gives birth (by pull back) to an $f$-canonical stratification of $\Mn_{k(v)}$, in $G(\AA^p_f)$-invariant, reduced, locally closed, quasi-affine subschemes (as $\Mn\to\Mm$ is a finite morphism, $\Mn_{k(v)}\to\Mm_{k(v)}$ is a finite morphism). We do not deal here with the interplay between the $f$-canonical and the canonical Lie stratifications 
of $\Mn_{k(v)}$. We just mention that $\Mu$ is an open subscheme of the open, dense stratum of the
$f$-canonical stratification of $\Mn_{k(v)}$, cf. b) or d) of 4.4.1 3), the definition of the canonical stratification of $\Mm_{k(v)}$, and (in case we do not want to use d) of 4.4.1.3) the following Exercise.
\medskip
{\bf Exercise.} Let $(A_s,p_{A_s})$, $s=\overline{1,2}$, be 2 principally polarized abelian varieties over an algebraically closed field of positive characteristic. If $A_1[p]$ is isomorphic to $A_2[p]$ then $(A_1[p],p_{A_1})$ is isomorphic to $(A_2[p],p_{A_2})$. Hint: use the classification of [EO] and [Oo3, 7.2] and Corollary of D of 2.2.22 3).
\medskip
3.9.3 suggests that this open, dense stratum should be $\Mu$ itself. As the quasi-affineness is well behaved under the operations of passing to finite morphisms or to open subschemes and taking quotients (via finite group actions), from [Oo3, 7.2] we get:
\medskip
{\bf Corollary.} {\it For any compact subgroup $\tilde H_0$ of $G(\AA_f^p)$, the $G$-ordinary locus of $\Mn_{k(v)}/\tilde H_0$ is a quasi-affine $k(v)$-scheme.}
\medskip
The interpretation of the $f$-canonical stratification of $\Mn_{k(v)}$ in terms of finite, flat group schemes over $\FF$ annihilated by $p$ (and liftable to $W(\FF)$) and endowed with some endomorphisms and principal polarizations should lead (via the use of PEL-envelopes as defined in 2.3.5.3) to the definition of a refined $f$-canonical stratification of $\Mn_{k(v)}$. 4.5.6.2 B implicitly points out that in general, the refined $f$-canonical stratification of $\Mn_{k(v)}$ once defined, can be different from the $f$-canonical stratification of $\Mn_{k(v)}$. 
\medskip
{\bf 4.5.10. Remark.} We do not stop here to prove that if $y_1,y_2\in\Mn_{k(v)}(k)$ are two distinct points giving birth to the same point $y\in\Mm_{k(v)}(k)$, then $y_1$ and $y_2$ are always in the same stratum, for any $\rho$-stratification of $\Mn_{k(v)}$, as a proof of Langlands--Rapoport's conjecture for $\Mn_{k(v)}$ in the context provided by $(f,L_{(p)},v)$ (see [Mi2] and [Va2, 1.7]) implies that such a situation $y_1\to y\gets y_2$ does not occur (cf. \S 14; to be compared also with [Va2, 5.6.4] and 4.4.6).
\medskip
{\bf 4.5.11. Toric points.} Here $k$ is again an arbitrary perfect field. The isomorphism of b) of 4.4.1 2) or just 4.2.10 leads us to define:
\medskip
a) a map $f_G:W_G\to LNP(\Mn_{k(v)})$ associating to $\om\in W_G$ the Newton polygon ${\rm Lie}_G(\Mp_\om)$ (of 4.1.5);
\smallskip
b) the set of $G(\om)$-ordinary points $y:{\rm Spec}(k)\to\Mn_{k(v)}$ by the rule: the Shimura adjoint Lie $\sg_k$-crystal attached to them, over $\bar k$, are inner isomorphic to the extension to $\bar k$ of the Shimura adjoint Lie $\bar\sg$-crystal of ${\got C}_\om$ of 4.1.5; here inner refers to isomorphisms defined by elements of $G^{\rm ad}_{W(k)}$, cf. 4.2.10. 
\medskip
{\bf 4.5.11.1. Three open problems.} If $\Mn$ has the completion property, then the set of $G(\om)$-ordinary points is non-empty. We assume that the set of $G(\om)$-ordinary points is non-empty (in 4.12.12.6 below we prove that this is always so). The $G(\om)$-ordinary points are points of a uniquely determined stratum ${\got s}_\om$ of the refined canonical Lie stratification of $\Mn_{k(v)}$; using the last part of 4.1.5.1 and the same argument as in the proof of the first sentence of b) of 4.4.1 3), we get that the $G(\om)$-ordinary points are points of a uniquely determined stratum ${\got s}^a_\om$ of the absolute stratification of $\Mn_{k(v)}$. 
\smallskip
We do not know (if or) when:
\medskip\noindent
$ALL$ all points of ${\got s}_\om$ (or of ${\got s}^a_\om$) are $G(\om)$-ordinary points, or when
\smallskip
$D(\om)$ $G(\om)$-ordinary points are points of an open, dense subscheme of ${\got s}_\om$ (or of ${\got s}^a_\om$), or when
\smallskip
\item{$ECC$} (at the opposite pole) each connected component of ${\got s}_\om$ (or of ${\got s}^a_\om$) has $G(\om)$-ordinary points.
\medskip
{\bf 4.5.11.2. Definition.} The points of $\Mn_{k(v)}$ with values in perfect fields which are $G(\om)$-ordinary, for some $\om\in W_G$, are called toric points of $\Mn_{k(v)}$.
\medskip 
{\bf 4.5.11.2.1. Exercise.} Show that toric points can be defined in terms of Galois representations: a point $y:{\rm Spec}(\bar k)\to\Mn_{k(v)}$ is a toric point, iff there is a $W(\bar k)$-valued point $z$ of $\Mn$ lifting $y$ such that the natural $p$-adic Galois representation 
$$\rho_z:\Gamma_{\bar k}\to GL(H^1_{\acute et}(z^*(\Ma),\overline{B(k)}))(\ZZ_p)$$ 
factors through the group of $\ZZ_p$-valued points of a maximal torus of $G_{\ZZ_p}$;
here we view $G_{\ZZ_p}$ as a subgroup of $GL(L_p^*)$ and identify (cf. [Va2, top of p. 473]) $L_p^*$ with $H^1_{\acute et}(z^*(\Ma),\overline{B(k)})$. Hint: use 2.2.13.3, 2.2.16.2, the Criterion of 2.2.22 1), 2.3.17 and the natural variant of 4.2.4 with $\FF$ replaced by $\bar k$. 
\medskip
{\bf 4.5.11.3. Remark.} We will see in \S10 that toric points are ``pillars" of the (refined) canonical Lie stratification of $\Mn_{k(v)}$ as well as of all other stratifications of 4.5.15-16 below. The most important ones are the $U$-ordinary points introduced in 4.4.13.
\medskip
{\bf 4.5.12. Problems. 1)} The map $f_G$ is usually not injective. Compute the number of elements of its fibres.
\smallskip
{\bf 2)} Compute the refined Lie stable $p$-ranks of toric points.  
\medskip
{\bf 4.5.13. CM levels of non-toric points.} Let $n\in\NN$. We consider a point $y$ of $\Mn_{k(v)}$ with values in a perfect field $k$ which is not a toric point. We say $y$ is of CM level at least $n$ if there is a morphism $z:{\rm Spec}(W(\bar k))\to\Mn$ whose special fibre factors through $y$ and whose attached Shimura filtered $\sg_{\bar k}$-crystal is (see defs. of 2.2.22 1)) cyclic diagonalizable of level $n$. Based on the equivalent of 4.2.4 over $\bar k$ instead of $\FF$ and on Corollary of 2.2.22 1), the second condition on $z$ can be restated in a more practical form as: 
\medskip
{\it The truncation mod $p^n$ of the Shimura adjoint filtered Lie $\sg_{\bar k}$-crystal ${\got C}_z$ attached to $z$ is isomorphic to the truncation mod $p^n$ of the extension to $\bar k$ of the Shimura adjoint filtered Lie $\bar\sg$-crystal attached to any one of ${\got C}_{\om}$'s of 4.1.5.}
\medskip
 We denote by 
$$n(y)\in\NN\cup\{0,1+2\dim(X)\}$$ 
the smallest supremum of the set of those $n$ such that $y$ is of CM level at least $n$; it is a) of 3.15.7 K and 3.15.7 BP0 which motivates the use of $1+2\dim(X)$ here. We refer to $n(y)$ as the CM level of $y$. 
Based on 2.3.17, we get that $n(y)$ is equal to the CM level of its attached Shimura $\sg_k$-crystal (see 3.13.7). The Expectation of 3.13.7.1 gets restated as: we expect that $n(y)\ge 1$. 
\medskip
{\bf 4.5.14. Problem.} Compute the number of strata and the dimensions of the connected components of the strata of the $f$-canonical stratification of $\Mn_{k(v)}$ introduced in 4.5.9. The first part of this problem is a particular case of Problem 1 of H of 2.2.22 3).
\medskip
{\bf 4.5.15. Other Lie type stratifications.} Using the (refined) Lie stable $p$-ranks attached to points of $\Mn_{k(v)}$ with values in fields, we obtain a (refined) Lie stable stratification of $\Mn_{k(v)}$ in $G(\AA_f^p)$-invariant, reduced, locally closed subschemes (cf. 3.9.1.1 and the part of 2.3.11 referring to $\Mg$ and to ${\got G}_{H_0}$). We can use as well Lie $p$-ranks (see 4.3.8 1)) to define another Lie type stratification of $\Mn_{k(v)}$; but the refined Lie stable stratification is in general more refined than it. 
\medskip
{\bf Example.} The fact F5 of 3.10.7 points out that there are many situations when the stratification of $\Mn_{k(v)}$ by Lie $p$-ranks has only one stratum, while 4.3.5-6 point out that, if $G$ is not a torus, this is not so for the refined Lie stable stratification of $\Mn_{k(v)}$. 
\medskip
We have many variants: using the maps $_1\overline{\psi_i}^{d_i}$ of 4.3.7 and their variants, or their exterior or symmetric powers, or Faltings--Shimura--Hasse--Witt (adjoint) maps, we can define a whole bunch of Lie type $G(\AA_f^p)$-invariant stratifications of $\Mn_{k(v)}$. Before pointing out one such stratification which is very useful, we explain what we mean here by ``their variants". In 4.3.7 we have selected an absolutely simple factor of $G_{iW(k)}$ with the property that the $F^1$-filtration of its Lie algebra defined --via 4.2.10-- by any $W(k)$-morphism $z:{\rm Spec}(W(k))\to\Mn$, is non-zero; but we have no reason to choose one particular such factor and so the choice of another one produces similar such maps (to be referred as variants of the maps of 4.3.7). 
\smallskip
For instance, using the ranks of images of the maps of 4.3.7 and of their variants and not just the ranks of their stable images, we get the Lie non-stable stratification: it is the intersection of as many stratifications of $\Mn_{k(v)}$ as non-compact factors of $G^{\rm ad}_{\RR}$ we have. In general it is different from (in some sense unrelated to) the refined Lie stable stratification of $\Mn_{k(v)}$. 
\smallskip
Here is an example in the abstract form. We refer to 4.4.13.3.1. We choose another Weyl element $w_1\in \tilde G(W(\FF_{p^3}))$ (resp. $w_2\in \tilde G(W(\FF_{p^3}))$) whose component in $\tilde G_1$ is the same as the component of $w$ and whose components in $\tilde G_2$ and $\tilde G_3$ are trivial (resp. whose components in $\tilde G_1$ and $\tilde G_2$ are the same as the ones of $w$ and whose component in $\tilde G_3$ is defined similarly to its other two components). Then the Shimura adjoint Lie $\sg_{\FF_{p^3}}$-crystals obtained as in 4.4.13.3.1 but using these two Weyl elements have the same Lie stable rank (equal to $0$) but their set of three maps defined as $_1\overline{\psi_i}^{d_i}$ of 4.3.7 are different (in the case of $w_2$ they are all zero maps, while in the case of $w_1$ only 2 of them are zero maps). 
\smallskip
There are many questions arising and which involve these Lie type stratifications (for instance: how are they related?, what is the number of their strata?, when are all strata smooth or quasi-affine?, etc.). To us, the most interesting ones (and so worth stating them separately) are:
\medskip
{\bf Q$_1$)} When the connected components of a given stratum of the Lie stable stratification of $\Mn_{k(v)}$ (or of $\Mn_{\FF}$) are permuted transitively by $G(\AA_f^p)$?
\smallskip
{\bf Q$_2$)} Using all stratifications hinted at in this 4.5.15 or introduced in 4.5.1-5, we can define their intersection. We call it the pseudo-ultra stratification of $\Mn_{k(v)}$. It is easy to see that in general it is a significant refinement of the absolute stratification. Examples: one example is already provided in the paragraph above referring to 4.4.13.3.1; a similar example can be obtained without changing ``the type" of the Weyl element $w$ of 4.4.13.3.1: referring to 4.4.13.3.1, if we choose a similar type of Weyl element but whose component in $\tilde G_1$ (or in $\tilde G_2$) is trivial, then this change is fully ``recorded" by the pseudo-ultra stratification. Can we in some way keep accurate track of this refinement?        
\smallskip
{\bf Q$_3$)} The above stratifications are very much related to the incipient theory of deviations of 3.13: the $\sg_k$-linear maps of 3.9.1 and other ones constructed naturally from them --like the ones of 4.3.7-- give birth (see above and below for samples) to a big bunch of stable, or of non-stable, or of anti-stable deviations of Shimura (filtered) (Lie) $\sg_k$-crystals. Which other type of ``points" of $\Mn_{k(v)}$ (besides the $G$-ordinary ones) and of their lifts (to $\Mn$), can be singled out using the $\sg_k$-linear maps of 3.9.1 (and eventually some deviations similar to the ones introduced in 3.13)? Can we in this way single out all toric points (which are not $U$-ordinary)?
\medskip
{\bf 4.5.15.0. Anti-stability.} To exemplify the ideas behind these new types of deviations and what we mean by anti-stable, we refer to 3.4.5 (so we work in a context modeled on 3.9.1 and not on its simplified, i.e. reduced, form of 4.3.7). Let $n,m\in\NN$ and let $I_0^s$ be a subset of $I_0$. As we do need to use explicitly iterates of functions, the $\sg_k$-linear map $\Psi_j^0$ of 3.4.5 is denoted here just by $\Psi_j$. We consider the image 
$$\Psi_j^n(\oplus_{i\in I_0^s} {\got g}_i)\subset {\got g}_0$$
and we take it mod $p^m$: we get a $W_m(k)$-submodule ${\rm Im}_{(n,m,I_0^s)}$ of ${\got g}_0/p^m{\got g}_0$ of finite length. Its length $l_{(n,m,I_0^s)}$ is called the $(n,m,I_0^s)$-deviation of the Shimura $\sg_k$-crystal $(M,\vph_j,G)$ or of its cyclic factor $({\got g}_0,\vph_j)$. Similarly, for $I^s$ an arbitrary subset of the set $I$ of the paragraph before 3.4.0, we can define the $(n,m,I^s)$-deviation of $(M,\vph_j,G)$ or of the Shimura Lie $\sg_k$-crystal $({\got g},\vph_j)$. 
\smallskip
We can ``complicate" the situation even more, by working in a filtered context $(M,F^1,\vph_j,G)$ and considering the intersection
$${\rm Im}_{(n,m,I_0^s)}\cap F^0({\got g}_0)/p^mF^0({\got g}_0);$$
we denote its length by $l_{(n,m,I_0^s,F^0({\got g}_0))}$. If $l_{(n,m,I_0^s)}\neq 0$, let
$$d(n,m,I_0^s,F^0({\got g}_0)):={{l_{(n,m,I_0^s,F^0({\got g}_0))}}\over {l_{(n,m,I_0^s)}}}.$$
If $l_{(n,m,I_0^s)}=0$, then $d(n,m,I_0^s,F^0({\got g}_0)):=0$. $d(n,m,I_0^s,F^0({\got g}_0))$ is called the $(n,m,I^s_0)$-deviation of $(M,F^1,\vph_j,G)$ or of $({\got g}_0,\vph_j,F^0({\got g}_0),F^1({\got g}_0))$. We have logical variants when $I_0^s$ is replaced as above by $I^s$. 
\smallskip
Working in the adjoint context, i.e. with the Shimura adjoint (filtered) Lie $\sg_k$-crystal attached to $(M,F^1,\vph_j,G)$, we similarly define such deviations: they will be denoted by adding ${\rm ad}$ as an upper right index; like $d^{\rm ad}(n,m,I_0^s,F^0({\got g}_0))$, $l^{\rm ad}_{(n,m,I_0^s)}$, etc. 3.9.6 points out: for any $n\in\NN$ and for every subset $I^s$ of $I$ we have
$$l^{\rm ad}_{(n,1,I^s)}=l_{(n,1,I^s)}.$$ 
\smallskip
If $I_0^s\neq I_0$ or if $I^s$ is not stable under $\ga$, then these deviations are referred as anti-stable deviations. All the above type of deviations can be defined in the abstract context of Shimura (adjoint) (filtered) Lie $F$-crystals over perfect fields. 
\medskip
{\bf 4.5.15.1. A concrete description of the pseudo-ultra stratification of $\Mn_{k(v)}$ and some variants.} For future references, we now include a canonical description of the pseudo-ultra stratification introduced in 4.5.15 $Q_2)$. We consider two points $y_1$ and $y_2$ of $\Mn_{k(v)}$ with values in the same algebraically closed field $k$. Let (see 4.3.1)
$$\bar\psi_j(i)=g_{j,i}\bar\psi(i)$$ 
be the $i$-th Faltings--Shimura--Hasse--Witt adjoint map attached to $y_j$, $j=\overline{1,2}$, $i\in\Mh$. Here $g_{j,i}\in G_i(k)$. For $n,m\in\NN$ and for a subset $I^s$ of the set $I_p(G^{\rm ad})$ of 4.3.1.1, we can define a number $l_{n,m,I^s,y_j}\in\NN\cup\{0\}$ as in 4.5.15.0. Using these numbers we define the pseudo-ultra stratification of $\Mn_{k(v)}$ as follows. $y_1$ and $y_2$ belong to the same stratum of this stratification iff the following two conditions hold:
\medskip
{\bf a)} {\it They belong to the same stratum of the absolute stratification;}
\smallskip
{\bf b)} {\it $l_{n,1,I^s,y_1}=l_{n,1,I^s,y_2}$, for any $n\in\NN$ and for every subset $I^s$ of $I_p(G^{\rm ad})$.}
\medskip
In order to define a refinement of it, we need a basic assumption (approach): 
\medskip
{\it Either we assume that the below Expectation holds or we allow stratifications of $\Mn_{k(v)}$ in potentially infinite number of strata (see 2.1; we recall that our convention in the last case is to describe the reduced, locally closed subschemes of $\Mn_{k(v)}$ or of $\Mn_{\FF}$ defined by such a stratification $\Ms$ of $\Mn_{k(v)}$, as the passage to an arbitrary $\bar k$ --i.e. the description of the reduced, locally closed subschemes of $\Mn_{\bar k}$ which are implicitly part of the definition of $\Ms$-- is automatic)}.
\medskip
{\bf Expectation (the finiteness property).} {\it For $i\in\Mh$, the set of inner isomorphism classes of Faltings--Shimura--Hasse--Witt adjoint maps of the form $({\rm Lie}({G_i}_k),g_i\bar\psi_1(i))$, with $g_i\in G_i(k)$, is finite and does not depend on $k$.}
\medskip
Based on 3.13.7.1, we do believe that this finiteness property always holds and that it can be proved (for instance) by just refining 3.5.3 slightly. Under this basic assumption we can define the ultra stratification in the same way we defined the pseudo-ultra stratification, by just replacing b) by the following requirement:
\medskip
{\bf c)} {\it For any $i\in\Mh$, $\bar\psi_1(i)$ and $\bar\psi_2(i)$ are inner isomorphic.}
\medskip
If we work only with c) (i.e. we do not impose a) as well), we speak about the quasi-ultra stratification or the Faltings--Shimura--Hasse--Witt stratification of $\Mn_{k(v)}$; it is a stratification which depends only on adjoint Lie $F$-crystals attached to points of $\Mn_{k(v)}$ with values in perfect fields. The argument of why c) always defines locally closed subschemes of $\Mn_{k(v)}$ is presented in 4.5.15.2. The Expectation should be compared with [EO] (or [Oo3]), where a similar finiteness property in the context of the $p$-torsion of principally quasi-polarized $p$-divisible groups over $k$ is used; in particular, loc. cit. implies (via 3.13.7.2 and 1) of 3.13.7.4 D) that the above finiteness property holds if all $G_i$'s of 4.3.1 are split of some $C_n$ Lie type.
\smallskip
We do not know what is the right connection between the quasi-ultra stratification and the refined canonical Lie stratification. Also, we do not know when the quasi-ultra (or the pseudo-ultra) and the ultra stratifications coincide. 
\smallskip
Warning: in what follows, whenever we refer to the (quasi-) ultra stratification, we implicitly assume that we work under the above basic assumption (approach).
\medskip
{\bf 4.5.15.2. Formulas pertaining to the quasi-ultra stratification.} From their very definition, the pseudo-ultra and the quasi-ultra stratifications are (once checked to exist) $G(\AA_f^p)$-invariant (see Fact 6 of 2.3.11). So from the part of 2.3.9 A pertaining to $\FF$ we get: to show that the quasi-ultra stratification of $\Mn_{k(v)}$ is well defined, we can work with the quotient $\Mn_{\FF}/H_0$ (in particular, in what follows we speak about the quasi-ultra stratification of $\Mn_{k_1}/\tilde H_0$ for any compact subgroup $\tilde H$ of $G(\AA_f^p)$ and for every perfect field $k_1$ which is either algebraically closed or is an algebraic extension of $k(v)$). For the sake of convenience we keep working with $k=\bar k$. 
\smallskip
The below Formula tells us that the quasi-ultra stratification is a very elementary concept (i.e. it can be easily studied, using few mathematical tools). Its proof is a trivial application of 2.3.15.1 (or of 3.6.14.4). We use the notations of 4.5.15.1. Let $\bar P_i(y_1)$ be the parabolic subgroup of ${G_i}_k$ having as its Lie algebra the $F^0$-filtration of ${\rm Lie}({G_i}_k)$ defined naturally via any lift of $y_1$ to a $W(k)$-valued point of $\Mn$ (i.e. defined by the reduction mod $p$ of the $F^0$-filtration of ${\rm Lie}(G_{W(k)}^{\rm ad})$ defined by any such lift; see 2.3.10 and 2.2.13). Let $U_i$ be an arbitrary smooth, unipotent subgroup of ${G_i}_k$ such that we have a natural open embedding 
$$U_i\hookrightarrow {G_i}_k/\bar P_i(y_1).$$
\indent 
We denote by $\TT^0$ the action we get as in 3.13.7.1 starting from the Shimura $\sg$-crystal associated to an arbitrary Shimura-ordinary $k$-valued point of $\Mn_{k(v)}$. It is 2.3.13.1, 4.2.1 and b) of 4.4.1 2) which allow us to make this choice arbitrarily: we can assume that $\forall i\in\Mh$, the $i$-th Faltings--Shimura--Hasse--Witt adjoint map associated to any such point is $\bar\psi(i)$ (of 4.3.1). $\TT^0$ gets decomposed into $\abs{\Mh}$-actions $\TT^0(i)$, $i\in\Mh$: $\TT^0(i)$ is the action of the centralizer $C_i(\mu)_k$ in ${G_i}_k$ of the image of $\mu_{k}$ in ${G_i}_k$ on 
$$\Mx(i):=N_i^+(\mu)_k\setminus {G_i}_k/\sg(N_i^-(\mu))_k\leqno (SP(i))$$
defined by restricting $\TT^0$ to $C_i(\mu)_k$. Here $N_i^-$ (resp. $N_i^+$) is the maximal unipotent subgroup of ${G_i}_{W(k(v))}$ with the property that $\mu$ acts on its Lie algebra (via inner conjugation) via the identical (resp. the inverse of the identical) character of $\GG_m$. 
\smallskip
We consider the locally closed subscheme $V_i(1)$ of $U_i$ whose $k$-valued points are those $u_i\in U_i(k)$ such that the pairs $({\rm Lie}({G_i}_k),\bar\psi_1(i))$ and $({\rm Lie}({G_i}_k),u_i\bar\psi_1(i))$ are inner isomorphic. Its existence is implied by the existence of $\TT^0$, as the condition on $u_i$ can be reformulated as: $u_ig_{1,i}\bar\psi(i)$ and $g_{1,i}\bar\psi(i)$ are inner isomorphic. Let $$d_i(y_1)$$ 
be the dimension of the connected component of $V_i(1)$ containing the origin of $U_i$. For $i\in\Mh^{\rm c}$, we have $d_i(y_1)=0$. Let $z\in\Mn_{k(v)}(W(k))$ be an arbitrary lift of $y_1$. Let $\tilde a$ and $g_{\tilde Y}$ have the same significance as in 2.3.15.1 (applied to $z$, viewed as a $W(k)$-valued point of $\Mn_{W(k)}/H_0$ lifting the $k$-valued point --still denoted by $y_1$-- of $\Mn_{k}/H_0$ defined by $y_1$; so $y$ of 2.3.12.1 is denoted here by $y_1$), with $N$ such that its special fibre is (naturally identifiable with) 
$$\prod_{i\in\Mh} U_i.$$ 
As the special fibre $\tilde a_k$ of $\tilde a$ is an \'etale morphism, we get that the image under $\tilde a_k$ of 
$$g_{\tilde Y}^{-1}(\prod_{i\in\Mh} V_i(1)),$$ 
is a locally closed subscheme of $\Mn_k/H_0$: it is the intersection of the stratum $s(y_1)$ of the quasi-ultra stratification of $\Mn_k/H_0$ to which $y_1$ belongs with the image of the special fibre of $\tilde a$. So the quasi-ultra stratification of $\Mn_k/H_0$ is well defined and we have:
\medskip
{\bf Formula.} {\it The dimension of the stratum of the quasi-ultra stratification of $\Mn_{k(v)}$ to which $y_1$ belongs is equal to 
$$d_{FSHW}(y_1):=\sum_{i\in\Mh} d_i(y_1)=\sum_{i\in\Mh^{\rm nc}} d_i(y_1).\leqno (1)$$} 
\medskip
{\bf 4.5.15.2.1. Variant.} Using 2.3.15 instead of 2.3.15.1, we can express even better the numbers $d_i(y_1)$'s. To explain this we need some more notations. Let 
$$d(i):=\dim_k({G_i}_k/{\bar P_i(y_1)});$$ 
for $i\in\Mh^{\rm c}$ it is $0$. Let $(M_1,\vph_1,G_{W(k)})$ is the Shimura $\sg_k$-crystal attached to $y_1$. From very definitions we get the second equality of
$$\dim_{\CC}(X)=dd((M_1,\vph_1,G_{W(k)}))=\sum_{i\in\Mh^{\rm nc}} d(i);\leqno (2)$$ 
the first equality is argued as in [Va2, end of 5.4.7]. We consider the group
$${\rm Inn}({\rm Lie}({G_i}_k),\bar\psi_1(i))$$
of inner automorphisms of the $i$-th Faltings--Shimura--Hasse--Witt adjoint map attached to $y_1$. It is an algebraic group over $k$. Let $d^i(y_1)$ be its dimension. We have:
\medskip
{\bf Proposition.} {\it $d_i(y_1)+d^i(y_1)=d(i)$, $\forall i\in\Mh$.}
\medskip
{\bf Proof:} We use the notations of 2.3.15 (as above, the resulting $W(k)$-valued point of $\Mn/H_0$ defined by $z$ is still denoted by $z$, the lift $z_1$ mentioned in 2.3.15 a) is as well still denoted by $z$ while its special fibre is denoted by $y_1$ instead of $y$). The morphism $a$ is of relative dimension $\dim_{\QQ}(G^0_{\QQ})-\dim_{\CC}(X)$. So the pull back of $s(y_1)$ to $Y_k$ is of dimension $$\dim_{\QQ}(G^0_{\QQ})-\dim_{\CC}(X)+d_{FSHW}(y_1).\leqno (3)$$ 
The key fact is that the special fibre of $g_Y$ is \'etale. We get 
$$\dim_{\QQ}(G^0_{\QQ})-\dim_{\CC}(X)+d_{FSHW}(y_1)=\dim_k(S^0(y_1)),\leqno (4)$$
where $S^0(y_1)$ is the connected component of the origin of the maximal reduced, locally closed subscheme $S(y_1)$ of $G^0_k$ whose $k$-valued points are those $g\in G^0(k)$ such that $\bar\psi_1(i)$ and $g\bar\psi_1(i)$ are inner isomorphic, $\forall i\in\Mh$. Defining similarly a reduced, locally closed subscheme $S_i(y_1)$ of ${G_i}_k$ (it is nothing else but the image of $S(y_1)$ in ${G_i}_k$), the connected component of it containing the origin is the image of $S^0_i(y_1)$ in ${G_i}_k$. We get that 
$$\dim_k(S^0(y_1))=\dim_k(Z(G^0_k))+\sum_{i\in\Mh} \dim_k(S^0_i(y_1)).\leqno (5)$$
Combining the above formulas (1) to (5) and splitting up the discussion in terms of $i\in\Mh$, we get
$$\dim_k({G_i}_k)-d(i)+d_i(y_1)=\dim_k(S^0_i(y_1)).\leqno (6)$$
\indent
Let $Z_i(y_1)$ be the subscheme of $S_i(y_1)$ centralizing (via left translations) $\bar\psi_1(i)$; we view it as a (non-necessarily reduced) group scheme over $k$. $S_i(y_1)$ is invariant under right translations by elements of $Z_i(y_1)(k)$. The $k$-valued points of the quotient variety $S_i(y_1)/Z_i(y_1)$ are in one-to-one correspondence to the set of distinct maps $h\bar\psi_1(i)$, with $h\in G_i(k)$, which are inner isomorphic to $\bar\psi_1(i)$. As $h\bar\psi_1(i)$ and $\bar\psi_1(i)$ have the same kernel ${\rm Lie}(\bar P_i(y_1))$, $\forall h\in G_1(k)$, an inner isomorphism between them is defined by an element of $\bar P_i(y_1)(k)$ (cf. also a) of Step 2 of 3.13.7.3). So $S_i(y_1)/Z_i(y_1)$ is naturally isomorphic to $\bar P_i(y_1)/{\rm Inn}({\rm Lie}({G_i}_k),\bar\psi_1(i))$. So all connected components of $S_i(y_1)$ have the same dimension and are smooth and moreover we have 
$$\dim_k(S^0_i(y_1))=\dim_k(\bar P_i(y_1))-d^i(y_1)+\dim_k(Z_i(y_1)).\leqno (7)$$
As $\dim_{W(k)}(G_i)=\dim_k(\bar P_i(y_1))+d(i)$, from (6) and (7) we get
$$d(i)-d_i(y_1)-d^i(y_1)=d(i)-\dim_k(Z_i(y_1)).\leqno (8)$$
\indent
We are just left to show that $d(i)=\dim_k(Z_i(y_1))$. From its definition, $Z_i(y_1)(k)$ is formed by those elements of $G_i(k)$ centralizing the image of $\bar\psi_1(i)$. This image is isomorphic (under the inner isomorphism defined by $g_{1,i}^{-1}$) to the image of $\bar\psi(i)$. So the equality $d(i)=\dim_k(Z_i(y_1))$ follows from a) and c) of Step 2 of 3.13.7.3. This ends the proof.
\medskip
{\bf Remark.} The above proof shows: we can read out which stratum of the quasi-ultra stratification of $\Mn_{k}$ specializes to which strata by just looking at which orbits of $\TT^0$ specialize to which orbits of it.
\medskip
{\bf 4.5.15.2.2. Remark.} $Z_i(y_1)$ is connected, $\forall i\in\Mh$, cf. a) of Step 2 of 3.13.7.3. So $S^0_i(y_1)=S_i(y_1)$, $\forall i\in\Mh$; so if $Z(G_k)$ is connected, $S(y_1)$ itself is connected. Moreover, as $S_i(y_1)$ is smooth, $\forall i\in\Mh$, we get that $S(y_1)$ is smooth. So, as $a$ (resp. $g_Y$ mod $p$) is smooth (resp. is \'etale), we get that $s(y_1)$ is smooth over $k$. We conclude:
\medskip
{\bf Corollary.} {\it The strata of the quasi-ultra stratification of $\Mn_{k(v)}$ are regular and formally smooth over $k(v)$.}
\medskip
{\bf 4.5.15.2.3. Inequalities.} Let $P_i^-(\mu)$ (resp. $P_i^+(\mu)$) be the subgroup of $G_{iW(k(v))}$ on which $\mu$ acts via inner conjugation through the trivial and the identity (resp. the trivial and the inverse of the identity) cocharacter of $\GG_m$. So $N_i^-(\mu)$ (resp. $N_i^+(\mu)$) as defined in 4.5.15.2, is the unipotent radical of $P_i^-(\mu)$ (resp. of $P_i^+(\mu)$). Let $\om$ be as in 4.1.5. We denote by 
$$d^i(\om)$$ 
the dimension of the automorphism group of the $i$-th Faltings--Shimura--Hasse--Witt adjoint map $w\bar\psi(i)$ attached to ${\got C}_{\om}$. 
\smallskip
We assume now that there is a $k$-valued point of $\Mn$ such that the truncation mod $p$ of the Shimura adjoint filtered Lie $\sg$-crystal associated to a lift $z\in\Mn(W(k))$ of $y$, is inner isomorphic to the truncation mod $p$ of the extension of the Shimura adjoint filtered Lie $\bar\sg$-crystal associated to ${\got C}_{\om}$ (of 4.1.5) to $k$. From the first sentence of 3.13.7.3.1, splitting (in our context) its statement in terms of $i\in\Mh$, we get
$$\dim_k(P_i^-(\mu))-\dim_k(N_i^-(\mu))\ge\dim_k(P_i^+(\mu)\cap w\sg({\rm Lie}(P_i^-(\mu)))w^{-1})-d^i(y_1).\leqno (9)$$
Combining (9) with the Proposition of 4.5.15.2.1 we get the second (the first one follows from very definitions) inequality of:
$$\dim_{W(k)}(w\sg(N_i^-(\mu))w^{-1}\cap N_i^-(\mu))\le d_i(y_1)\le\dim_{W(k)}(w\sg(P_i^-(\mu))w^{-1}\cap N_i^-(\mu)).\leqno (10)$$
\medskip
{\bf 4.5.15.2.4. Remarks.} {\bf 1)} The quasi-ultra stratification is a refinement of the Lie stable stratification of $\Mn_{k(v)}$. From 4.3.4 c) we get: they both have the same open, dense stratum (it is $\Mu$ of 4.2.1). 
\smallskip
{\bf 2)} It is crucial to define quasi-ultra stratifications in terms of inner automorphisms and not just in terms of automorphisms: otherwise, we obtain (in general) less refined stratifications, as 4.5.6.2 B points out. So, 4.5.15.1 c)
 implies 4.5.15.1 b) but, in general, the converse does not hold. For instance, if $G^{\rm ad}_{\ZZ_p}$ is split and $G^{\rm ad}_{\RR}=SU(1,2)^{\rm ad}_{\RR}$, 4.5.15.1 b) defines only $2$ strata of $\Mn_{k(v)}$, while 4.5.15.1 c) defines $3$.
\smallskip
{\bf 3)} It is worth pointing out that in the above study of the quasi-ultra stratification of $\Mn_{k(v)}$, the word Verschiebung does not show up at all (cf. also 3.13.7). So we view Faltings--Shimura--Hasse--Witt adjoint maps as the adjoint Lie analogue (in any relative context pertaining to $p$-divisible objects with a reductive structure) of truncations mod $p$ of $p$-divisible groups over $\FF_p$-schemes, cf. 3.13.7.9. 
\smallskip
A great part of the theory (see 3.13.7 and above) of Faltings--Shimura--Hasse--Witt adjoint maps can be redone (particularly see 3.13.7.1.2) entirely in terms of truncations mod $p$ of Shimura $\sg_k$-crystals as defined in 2.2.14; but, in the context of a SHS $(f,L_{(p)},v)$, till the proof of d) of 4.4.1 3) or at least of 4.2.8.1 is not written down, we can get into trouble with the possibility that a connected component of $\Mn_{W(\FF)}$ might have a special fibre which is not connected, and so the situation would not be satisfactorily enough. Moreover, the language of Faltings--Shimura--Hasse--Witt adjoint maps is universal (cf. 3.13.7.8-9; for instance it applies immediately to other classes of varieties or to the $E_6$, $E_7$ and $D_n^{\rm mixed}$ types). Also we do believe, that the computations are much easier (transparent) in the adjoint context (not involving any Verschiebung maps).  
\smallskip
{\bf 4)} We get significantly better estimates than (9) or (10) above, if we use the non-compact Faltings--Shimura--Hasse--Witt adjoint maps attached to points of $\Mn_{k(v)}$ (they are definable as the usual ones, starting from 3.13.7.6 and 3.9.6).
\smallskip
{\bf 5)} In general, the inner isomorphism classes of Fontaine truncations mod $p$ of (pull backs via $k$-valued points) of a global deformation (over a regular, formally smooth $W(k)$-scheme $X$) of a $p$-divisible object with a reductive structure of $\Mm\Mf_{[a,b]}(W(k))$, do define locally closed subschemes of $X_k$ regardless of how $a$ and $b$ are: the proof is the same as in 4.5.15.2 (cf. 3.13.7.8-9). However, due to limitations explained in 3.6.8.9, presently we can not say when Proposition of 4.5.15.2.1 generalizes outside (cf. 3.15.6) of the generalized Shimura context. 
\medskip
{\bf 4.5.15.2.5. A supplement to 4.5.9.} Let $(M_y,\vph_y,G_{W(k)},p_{M_y})$ be the principally quasi-polarized Shimura $\sg_k$-crystal associated to a point $y\in\Mn(k)$. Let $C$ be the center of the centralizer of $Z(G^0_{W(k)})$ in $GL(M_y)$. We first assume there is a torus $C_0$ of $C$ containing $Z(G^0_{W(k)})$ and contained in $Sp(M_y,p_{M_y})$. So for any $g\in Z(G^0_{W(k)})(W(k))$ the truncations mod $p$ of $(M_y,\vph_y,G_{W(k)},p_{M_y})$ and of $(M_y,g\vph_y,G_{W(k)},p_{M_y})$ are isomorphic under an isomorphism defined by a $k$-valued point of $C_0$ (cf. 3.13.7.2). As the map $G^0_k(k)\to G^{\rm ad}_k(k)$ is surjective, we get that any $k$-valued point of the same connected component of $\Mn_k$ through which $y$ factors and whose Faltings--Shimura--Hasse--Witt adjoint map is inner isomorphic to the one of $y$, maps (via the special fibre of $i_{\Mn}$ of 2.3.2) into the same stratum of the canonical stratification of $\Mm_{k(v)}$ as defined in [EO] (see also [Oo3]; if $y$ is a toric point, then the use of $C_0$ is entirely avoidable here, cf. Fact of 3.13.7.2).
\smallskip
From b) of 4.4.1 3) and from the Exercise of 4.5.9, we get that this remains true even if such a torus $C_0$ does not exist. So, as in 4.5.9, we get that each stratum of the quasi-ultra stratification of $\Mn_{k(v)}$ is a quasi-affine scheme. 
\medskip
{\bf 4.5.15.3. Automorphism invariants.} We do not assume anymore $k$ algebraically closed. Let $z_1\in\Mn/H_0(k_1)$ be such that the group $\tilde G_{W(k)}$ obtained as in 2.3.10 is $G_{W(k)}$ itself. Denoting by $y_1$ its special fibre, let $\bar\psi_1(i)$ and ${\rm Inn}({\rm Lie}({G_i}_k),\bar\psi_1(i))$
be as in 4.5.15.1 and 4.5.15.2.1 (cf. 3.9.1.1). If $k$ is infinite (resp. finite), let $H_i(y_1)$ be the connected component of the origin (resp. be the trivial subgroup) of ${\rm Inn}({\rm Lie}({G_i}_k),\bar\psi_1(i))$. The group 
$${\rm aut}_i(y_1):={\rm Inn}({\rm Lie}({G_i}_k),\bar\psi_1(i))(k)/H_i(y_1)(k)$$
is called the $i$-th automorphism invariant of $y_1$. 
\smallskip
Such an automorphism invariant can be defined for any (Faltings--Shimura--Hasse--Witt adjoint map of a) cyclic Shimura (adjoint) Lie $\sg_k$-crystal. In particular, for any $\om\in W_G$, we denote by ${\rm aut}_i(\om)$ the automorphism invariant of the $i$-th cyclic adjoint factor of the Shimura filtered $\bar\sg$-crystal ${\got C}_{\om}$-crystal of 4.1.5. 
\smallskip
In case $k=\bar k$, it seems to us that these invariants are unrelated to the $d_i(y_1)$'s dimensions introduced in 4.5.15.2. Often, ${\rm aut}_i(y_1)$ can be naturally interpreted as a subgroup of the group of $\FF_p$-valued points of a connected group scheme $A_i(y_1)$ over $\FF_p$; it seems to us that even in such cases the dimension of $A_i(y_1)$ is unrelated to $d_i(y_1)$.
\medskip
{\bf Example.} We situate ourselves in the abstract context of 3.4.5; so the $E_6$ and $E_7$ Lie types are also allowed. We assume $I_0=I_1=\{1\}$ and $k=\bar k$. With the notations of 3.5, we consider the reductive subgroup ${P_{00}}_k$ of ${G_1}_k$ of whose Lie algebra is $\bar{\got q}_{\vep_1}\cap \bar{\got p}_{\vep_1}$. Using the $\ZZ_p$-structures of 3.11.2 C, we get that it has a natural $\FF_p$-structure $P_{00\FF_p}$ and that it makes sense to speak about elements of $\bar{\got q}_{\vep_1}$ fixed by $\sg$. The group ${\rm Inn}({\got g}_0/p{\got g}_0,\bar\psi_0^0)(k)$ is naturally identifiable with the group of elements of ${P_{00}}(k)$ normalizing the $\FF_p$-Lie subalgebra of $\bar{\got q}_{\vep_1}$ of elements fixed by $\sg$. Using the Lemma of Step 2 of 3.13.7.3 we get: ${\rm Inn}({\got g}_0/p{\got g}_0,\bar\psi_0^0)(k)$ is finite and is a subgroup of the group of $\FF_p$-valued points of the image of $P_{00\FF_p}$ in the $GL$-group of this abelian $\FF_p$-Lie algebra.
\medskip
{\bf 4.5.16. Faltings--Shimura--Dieudonn\'e (adjoint, principal or standard) stratification.}
Again, in what follows we allow stratifications of $\Mn_{k(v)}$ in potentially an infinite number of strata. Let $n\in\NN$. We consider the equivalence relation $E_n$ on the set of $\FF$-valued points of $\Mn_{k(v)}$ defined by the rule (again 4.2.10 allows us to use the word inner):
\medskip
\item{{\bf $E_n$}} {\it $y_1$, $y_2\in\Mn(\FF)$ are in relation $E_n$ iff the truncations mod $p^n$ (see 3.9.8) of their Shimura adjoint Lie $\bar\sg$-crystals attached to them are inner isomorphic.}
\medskip
The case $n=1$, corresponds to the equivalence relation on $\Mn(\FF)$ defined naturally by the quasi-ultra stratification of $\Mn_{k(v)}$. It is easy to see that similarly, the equivalence classes of $E_n$ define naturally constructible subsets of $\Mn_{k(v)}$. We call $E_n$ as the level-$n$ Faltings--Shimura--Dieudonn\'e equivalence relation on $\Mn(\FF)$; it is (cf. Fact 6 of 2.3.11) $G(\AA_f^p)$-invariant and so we still denote by  $E_n$ the induced equivalence relation on $\Mn/H_0(\FF)$. We need an extra equivalence relation $E_{\infty}$ on $\Mn(\FF)$:
\medskip
\item{{\bf $E_{\infty}$}} {\it $y_1$, $y_2\in\Mn(\FF)$ are in relation $E_{\infty}$ iff their Shimura adjoint Lie $\bar\sg$-crystals attached to them are inner isomorphic.} 
\medskip
Similarly we define the equivalence relations $E_n(l)$ and $E_{\infty}(l)$ on $\Mn()$ or on $\Mn/H_0(l)$, with $l$ an algebraically closed field containing $k(v)$. Let $E_{\infty}(k(v)):=E_{\infty}$. We have:  
\medskip
{\bf Theorem.} {\it There is a uniquely determined stratification $ST$ of $\Mn_{k(v)}$ in $G(\AA_f^p)$-invariant, reduced, locally closed subschemes having the following properties:
\medskip
{\bf a)} it has a stratum which is an open, dense subscheme of $\Mn_{k(v)}$;
\smallskip
{\bf b)} its strata are regular and quasi-affine; 
\smallskip
{\bf c)} all connected components of a given stratum have the same dimension;
\smallskip
{\bf d)} If $l$ is $k(v)$ or an algebraically closed field containing $k(v)$, two $\bar l$-valued points of $\Mn_l$ belong to a locally closed subscheme of $\Mn_l$ which is a stratum of $ST$ iff they are in relation under $E_{\infty}(l)$.}
\medskip
{\bf Proof:} From 3.15.7 BP2 we get that there is $n\in\NN$ effectively computable and such that $E_n=E_{\infty}$. We consider the adjoint variant ${\got G}^{\rm ad}_{H_0}$ of the Lie $p$-divisible object ${\got G}_{H_0}$ of $\Mm\Mf_{[-1,1]}(\Mn/H_0)$ introduced in 2.3.11 (cf. also the first proof of 3.6.18.7.3 A). From b) of 3.15.7 D and 3.15.7 E applied to it, we get: there is an affine $k(v)$-morphism 
$$m:Y_{k(v)}\to \Mn_{k(v)}/H_0\times \Mn_{k(v)}/H_0$$
of finite type such that $y_1$, $y_2\in\Mn_{k(v)}/H_0(\FF)$ are in relation $E_{n}$ iff the $\FF$-valued point $(y_1,y_2)$ of $\Mn_{k(v)}/H_0\times \Mn_{k(v)}/H_0$ lifts to an $\FF$-valued point of $Y_{k(v)}$. Let $p_1$ and $p_2$ be the $k(v)$-morphisms $Y_{k(v)}\to \Mn_{k(v)}/H_0$ naturally defined by the projections of $\Mn_{k(v)}/H_0\times \Mn_{k(v)}/H_0$ on $\Mn_{k(v)}/H_0$, $s=\overline{1,2}$. So $y_2$ is in relation $E_{\infty}$ with $y_1$ iff $y_2$ is in the image $I_{y_1}$ of the restriction of $p_2$ to $p_1^{-1}(y_1)$. $I_{y_1}$ is a constructible set. As $E_n=E_{\infty}$, from the Theorem and Fact 4 of 2.3.11 we get that the local geometry of $I_{y_1}$ in each $\FF$-valued point of it is the same. So all connected components of $I_{y_1}$ are regular and have the same dimension; so $I_{y_1}$ has a natural structure of a reduced, locally closed subscheme of $\Mn_{k(v)}/H_0$ (and not just of a constructible set).  From Fact 6 of 2.3.11 we get that the pull back of $I_{y_1}$ to $\Mn_{k(v)}$ is $G(\AA_f^p)$-invariant. 
\smallskip
The same remains true if we work with an arbitrary algebraically closed field $l$ containing $k(v)$. So $ST$ exists, is uniquely determined by d), is $G(\AA_f^p)$-invariant, its strata are regular, and c) holds for it. The second part of a) is implied by b) of 4.4.1 2). The quasi-affineness part of b) is implied by 4.5.15.2.5.  
This ends the proof.
\medskip
{\bf 4.5.16.0. Remarks.} {\bf 1)} We refer to $ST$ as the Faltings--Shimura--Dieudonn\'e adjoint stratification of $\Mn_{k(v)}$. Warning: its number of strata is in general infinite (cf. 3.9.7.3); however, this is not a serious handicap, as often it is easy to see that locally in the Zariski topology, its strata which are locally closed subschemes of $\Mn_{k}$ are obtained as the fibres of a suitable finite number of fibrations. 
\medskip
{\bf 2)} We have two main variants of the Theorem of 4.5.16. In the first (resp. second) one we replace Shimura adjoint Lie $\bar\sg$-crystals by Shimura (resp. by principally quasi-polarized Shimura) $\bar\sg$-crystals; the first one is stated in Theorem 6 of 1.6.3. Their proofs need only one modification: in the principally quasi-polarized context, if $({{-1}\over p})=1$ we can not say that we have $G(\AA_f^p)$-invariant equivalence relations and so we need to choose $H_0$ small enough (for instance, we need $H_0$ such that 2.3.3 (INCL) holds). We get the Faltings--Shimura--Dieudonn\'e standard (resp. principal) stratification of $\Mn_{k(v)}$. Often we drop the word standard.
\smallskip
If the homomorphism $Z(G_{\ZZ_p})(W(\FF))\to Z(G_{\ZZ_p})(W(\FF))$ that takes an element $a\in Z(G_{\ZZ_p})(W(\FF))$ into $a\bar\sg(a^{-1})$ is surjective and if 4.2.8.1 holds, then the Faltings--Shimura--Dieudonn\'e adjoint and standard stratifications of $\Mn_{k(v)}$ coincide. Moreover, if $({{-1}\over p})=-1$ they also coincide with the Faltings--Shimura--Dieudonn\'e principal stratification of $\Mn_{k(v)}$ (to be compared with the proof of Fact 6 of 2.3.11). 
\medskip
{\bf 4.5.16.1. The third form (i.e. the isomorphism form for a SHS) of the purity principle.} We have:
\medskip
{\bf Corollary.} {\it The Faltings--Shimura--Dieudonn\'e (adjoint) stratification satisfies the purity property.}
\medskip
{\bf Proof:} It is enough to deal with the adjoint context: the arguments for the non-adjoint context are entirely the same. Let $SN$ be the normalization of the Zariski closure $\bar s(y_1)$ in $\Mn_{\FF}/H_0$ of a stratum $s(y_1)$ of the Faltings--Shimura--Dieudonn\'e adjoint stratification defined by some $y_1\in\Mn(\FF)$. Let $n\in\NN$ be such that $E_n=E_{\infty}$ and let $m\in\{n,n+3\}$. Let ${\got L}(m)$ be the reduction mod $p^{m}$ of the Shimura adjoint Lie $\bar\sg$-crystal attached to $y_1$; we view it as an object of $\Mm_{[0,2]}(W(\FF))$ whose underlying $W_m(\FF)$-module has a natural Lie structure. We consider the scheme $ISO_m(SN)$ parameterizing inner isomorphisms between:
\medskip
i) the pull back of ${\got L}(m)$ through the natural morphism $SN\to {\rm Spec}(\FF)$,
\smallskip
ii) and the pull back of ${\got G}^{\rm ad}_{H_0}(1)/p^m{\got G}^{\rm ad}_{H_0}(1)$, viewed without filtration, through the natural morphism $SN\to\Mn_{k(v)}$. 
\medskip
Here we view the pull backs of i) and ii) as $F$-crystals on $SN$ in coherent sheaves endowed with Lie structures. $ISO_m(SN)$ is an affine $SN$-schemes (even if $SN$ is not regular): locally in the Zariski topology of $SN$ it is constructed via evaluations at $W_m(SN)$; to be compared with the proof of 3.15.7 D). 
\smallskip
We have a natural reduction morphism
$$REDUC:ISO_{n+3}(SN)\to ISO_n(SN).$$
We consider the maximal reduced, closed subscheme of $ISO_n(SN)_{\rm red}$ through which $REDUC_{\rm red}$ factors. It is a quasi-finite, affine $SN$-scheme $SN_{n,n+3}$; the quasi-finiteness follows (via specialization) from b) of 2.2.4 B. 
\smallskip
From the part of the above proof referring to local geometry, we get that all connected components of $SN_{n,n+3}$ have the same dimension, are irreducible and normal. Moreover, the number of points of each fibre of $SN_{n,n+3}$ above an $\FF$-valued point of $SN$ mapping into $s(y_1)$ (resp. not mapping into $s(y_1)$) depends only on $y_1$ (resp. is $0$). Let $SN^n$ be the normalization of $SN$ in the ring of fractions $RF$ of $SN_{n,n+3}$. From Zariski's Main Theorem we get that $SN_{n,n+3}$ is an open subscheme of $SN^n$. As $SN_{n,n+3}^n$ and $SN^n$ are affine $SN$-schemes, the complement $CO$ of $SN_{n,n+3}^n$ in $SN^n$ is either empty or of pure codimension 1. The complement of $s(y_1)$ in $\bar s(y_1)$ is the image of $CO$ in $\bar s(y_1)$; so it is either empty or of pure codimension 1. This (cf. also Exercise of 3.6.8.1.4) ends the proof.  
\medskip
{\bf 4.5.16.2. Remarks.} {\bf 1)} It seems to us that the types of 3.9.7 can be used as well to define a stratification of $\Mn_{k(v)}$ (automatically in a finite number of strata).  
\smallskip
{\bf 2)} The proof of 4.5.16.1 applies as well to give us that the quasi-ultra stratification of $\Mn_{k(v)}$ satisfies the purity property. 
\medskip
{\bf 4.5.17. Non Lie type stratifications.} Forgetting the extra (Shimura) structure of $F$-crystals attached to points of $\Mn_{k(v)}$ with values in perfect fields, we can define other (not necessarily of Lie type!) stratifications of $\Mn_{k(v)}$, closer in spirit to the stratifications of $\Mm_{k(v)}$ (like the one using $p$-ranks, see [NO], etc.). The slightly unpleasant feature of such stratifications of $\Mn_{k(v)}$: very often they are not canonical (i.e. they depend on the SHS producing them; on the contrary the Lie type stratifications are canonical, see 4.9 below) or they are not refined enough, as it can be seen through examples in which we have $\abs{\Mh^{\rm nc}}\ge 2$.  
\medskip
{\bf 4.5.18. Artin--Schreier stratifications.} By the sum of two systems of equations in $n_1$ and respectively in $n_2$ variables, we mean the system of equations in $n_1+n_2$ variables obtained by ``putting" them together. Let $X_p$ be an arbitrary reduced $\FF_p$-scheme. A stratification $\Ms$ of it in reduced, locally closed subschemes is called an Artin--Schreier stratification (to be abbreviated as: an AS stratification), if there is $m\in\NN$, such that locally in the Zariski topology of $X_p$, $\Ms$ is an AS stratification obtained as in 3.6.8.1.3, using quasi Artin--Schreier systems of equations in $m$ variables (see 3.6.18.4.6 A for the definition of these systems of equations). Similarly we define a refined Artin--Schreier stratification (to be abbreviated as a RAS stratification) of $X_p$. The smallest such $m$ is called the minimal degree of definition of $\Ms$ and is denoted by $md(\Ms)$. From the estimate of $m_1$ in 3.6.8.1 (with $l=1$) and from the end of 3.6.8.1.3 we get:
\medskip
{\bf Fact.} {\it $md(\Ms)-1$ is greater or equal to the minimum number of strata we get by restricting $\Ms$ to open, affine subschemes of $X_p$.}
\medskip
{\bf Example 1.} We assume $X_p$ is a smooth, connected $k$-scheme and there is a quasi-polarized $p$-divisible group $(D_p,p_{D_p})$ over $X_p$. We consider the stratification $\Ms(D_p)$ of $X_p$ defined by $p$-ranks of $D_p$ over geometric points of $X_p$. From 3.6.18.4 B) and the moduli principle of 3.6.18.4.2 we deduce  that it is a RAS stratification (we can take $m$ to be $\dim_k(X_p){\rm rk}(D_p)^2$). Using a suitable direct sum of $D_p$ with an \'etale $p$-divisible group of sufficiently high rank, based again on 3.6.18.4 B) and the moduli principle of 3.6.18.4.2, we get that $\Ms(D_p)$ is a RAS stratification, without assuming that $D_p$ is quasi-polarized. Using Hasse--Witt invariants, we can reobtain this fact (we can take $m={\rm rk}(D_p)$), again without using that $D_p$ is quasi-polarized. 
\smallskip
If $X_p$ is a connected component of $\Mm_{k(v)}$ (of 2.3.2) and we are dealing with the standard principally quasi-polarized $p$-divisible group over it, from [NO, 4.1] we deduce that this RAS stratification is as well an AS stratification. 
\medskip
{\bf Example 2.} The Lie $p$-rank stratification of $\Mn_{k(v)}$ is a RAS stratification. As in Example 1, this follows from 3.6.18.7.0 and Corollary of 3.6.18.7.3 (via 3.6.18.4.2). 
\medskip
Warning: not always a RAS stratification is an AS stratification.
\medskip
{\bf Example 3.} Let $R:=k[x_1,x_2]$. We consider the quasi Artin--Schreier system of equations in 2 variables with coefficients in $R$ defined by the equations $z_1=x_1x_2z_1^p+x_1z_2^p$ and $z_2=x_1x_2z_1^p+2x_1z_2^p$. The RAS stratification of ${\rm Spec}(R)$ defined by it is not an AS stratification: the complement of the open, dense stratum has 2 strata of the same dimension.
\medskip
{\bf Example 4.} From 3.9.1.0, following the proof of 3.6.18.4.6 C, we get that the Lie stable stratification of $\Mn_{k(v)}$ is a RAS stratification. This is a particular case of the following principle. Let $q\in\NN$ and let $F_{X_p}$ be the Frobenius endomorphism of $X_p$. We consider a locally free $\Mo_{X_p}$-sheaf $\Mf$ together with an $\Mo_{X_p}$-linear map $l:{F_{X_p}^q}^*(\Mf)\to\Mf$. Working similarly to 3.9.1 we can define the stable $p$-ranks of $l$ w.r.t. geometric points of $X_p$ and use them to define a stratification $\Ms(X_p)$ of $X_p$ in reduced, locally closed subschemes. As in the particular case, we get that it is a RAS stratification. Moreover, 3.6.8.1.4 applies to it.
Examples 1 and 2 are as well particular samples of this principle: using Hasse--Witt maps and their versions at the level of ${\rm End}$'s, in reduced forms similar to the ones of 4.3.7, we get that always a Lie $p$-rank (resp. a $p$-rank) stratification in a context involving Shimura $p$-divisible groups (resp. $p$-divisible groups) is a RAS stratification.
\smallskip
In general, we can not perform the same construction, if we restrict $l$ to ${F_{X_p}^q}^*(\Mf_1)$, with $\Mf_1$ an $\Mo_{X_p}$-subsheaf of $\Mf$ which locally in the Zariski topology is a non-trivial direct summand: to compute similar stable ranks, we come across systems of equations of third type (see 3.6.8.9) and so (see 3.6.8.9.0) in general they can not be used to define (following the pattern of 3.6.8.1.3) stratifications of $X_p$.
\medskip
{\bf Example 5.} We assume $X_p$ is regular. Using systems of equations of the form
$$x_i=f_ix_i^p,$$
$i=\overline{1,m}$, with $f_i$'s as elements of some $\FF_p$-algebras of global sections of open, affine subschemes of $X_p$, we get that the natural stratification of $X_p$ defined by any divisor of it with normal crossings, is a RAS stratification. Warning: Example 3 shows that not all RAS stratifications of an $\FF_p$-scheme can be obtained in such a way, i.e. starting from a divisor with normal crossings.
\medskip
{\bf Example 6.} We consider two quasi Artin--Schreier systems of equations $S_1$ and $S_2$ with coefficients in a reduced $\FF_p$-algebra $R$; let $\Ms_1$ and respectively $\Ms_2$ be the RAS stratifications of ${\rm Spec}(R)$ they define. Using the sum of $n_1$ copies of $S_1$ with $n_2$ copies of $S_2$, for suitable $n_1$, $n_2\in\NN$ (which depend only on the number of variables of $\Ms_1$ and $\Ms_2$), we get that the intersection $\Ms_1\cap\Ms_2$ is a RAS stratification. We conclude: always the intersection of a finite number of RAS stratifications is a RAS stratification.  
\medskip
{\bf 4.5.18.1. Conjecture.} {\it Always a Newton polygon stratification associated to a $p$-divisible group over $X_p$ (or to a $p$-divisible object of some Fontaine category) is a RAS stratification.}
\medskip
The motivation for this Conjecture is based on:
\medskip
-- examples;
\smallskip
-- on the proof of 3.15.10; 
\smallskip
-- on 3.6.20 4); 
\smallskip
-- and on the fact that taking different exterior powers of objects of $\Mm\Mf_{[a,b]}(*)$, we expect formulas similar to the ones of 3.6.18.4 B) and 3.6.18.7 and which involve (besides different pseudo-multiplicities) as well slopes which are not necessarily integral values; these different exterior powers should be able (at least in many cases; often we might have to use direct sums or products as in Example 1) to ``capture the shape" of Frobenius endomorphisms of (the ``involved") vector bundles modulo higher powers of $p$ (we recall that in getting the formula of 3.6.18.4 B) only the value of $\Phi$ mod $p^2$ was important). 
\medskip
We have a variant of this conjecture in the context (see 3.15.7 D and E) of isomorphism classes of $p$-divisible objects with a reductive structure. So, using 3.6.8.1.4 we should be able to reobtain the forms of the purity principle presented in 3.15.10 and 4.5.16.1. 
\smallskip
{\bf 4.5.18.2. Remark.} One of the advantages of AS or RAS stratifications in comparison with Newton polygon stratifications is: using different determinants or norm maps, we can put different natural scheme structures on the strata, which are not a priori reduced (there is always logic in the choice of $m$ and of the quasi Artin--Schreier systems of equations involved). For instance, referring to Example 3, the third stratum of ${\rm Spec}(R)$ is logically defined by the equation $x_1^2=0$. 
\medskip
{\bf 4.5.18.3. Question.} {\it Is it true that locally (in the Zariski or the flat topology) a RAS stratification is obtained from a AS stratification by pull backs?}
\medskip
We hope to come back to the ideas of section 4.5.18 in a future paper. 
\medskip\smallskip
{\bf 4.6. Examples and main properties.}
Here we list the main properties of the theory of Shimura-ordinary types, of $G$-ordinary points and of their $G$-canonical lifts, and present some examples to illustrate them as well as the complexities which arise.
\medskip
{\bf P1} {\it $\tau$ is the (usual) ordinary type iff $k(v)=\FF_p$.}
\medskip
{\bf Proof:} From the description of $\tau$ given in 4.1, we deduce
(with the notations of 4.1.1.1) that $h^0$ has only the eigenvalues 0 and 1 iff
$\bar h_i=\bar h_1$, $\forall i\in S(1,d)$, and this is equivalent to $\sg(\bar h_1)=\bar h_1$, i.e. it is equivalent to $\mu$ being defined over $\ZZ_p$. So everything results from the paragraph of 4.1 referring to [Mi3, 4.6-7].
\medskip
{\bf Example 1.} In the case of a standard PEL situation $(f,L_{(p)},v,\Mb)$, if $(G^{\rm ad},X^{\rm ad})$ is of simple $C_\ell$ type, we have $E(G,X)=\QQ$; here either $\ell\ge 2$ or $\ell=1$ and all factors of $G^{\rm ad}_{\RR}$ are non-compact. So we get only usual ordinary types. Similarly, for the SHS's involving the classical Spin modular varieties of odd dimension and rank 2 considered in [Va2, 5.7.5], we get only usual ordinary types.
\medskip
{\bf P2} {\it If $k(v)=\FF_p$ and if $y$ and $z$ are as in 4.4.1 2), then $z^*(\Ma)$ is the (usual) canonical lift of (cf. P1) the ordinary abelian variety $y^*(\Ma)$.}
\medskip
{\bf Proof:} This is a direct consequence of c) of 4.4.1 3) and of the description of the $\sg$-crystal $(M,\vph)$ of 4.1.1 giving birth to $\tau$. 
\medskip
{\bf P3} {\it If $\tau$ is $e(1,1)$, i.e. if $\tau$ is the formal isogeny type associated to supersingular abelian varieties of dimension $e$, then $G$ is a torus.}
\medskip
{\bf Proof:} If $\tau$ is $e(1,1)$, then the Shimura Lie $\sg_k$-crystal attached to a $G$-ordinary point $y\in\Mn_{k(v)}(k)$ has only slopes $0$. From this and 4.4.1 1) we get that the deformation dimension of the Shimura $\sg_k$-crystal attached to $y$ is $0$. So $\dim_{\CC}(X)=0$, cf. 4.5.15.2.1 2). But $G$ is a torus iff $\dim_{\CC}(X)=0$, cf. axioms [Va2, (SV1-3) of 2.3]. This ends the proof.
\medskip
{\bf P4} {\it The slopes of Shimura $F$-crystals attached to $G$-ordinary points of $\Mn_{k(v)}$ are of the form ${i\over d}$, with $i\in S(0,d)$.}
\medskip
{\bf Proof:} This is a consequence of  4.2.1 a) and 4.1.2.2.
\medskip
{\bf Example 2.} If $k(v)=\FF_{p^2}$ and if $G$ is not a torus, then $\tau$ is of the form $r(1,0)+s(1,1)+r(0,1)$, with $r,s\in \NN$ (cf. P1, P3 and P4).
\medskip
{\bf Example 3.} We consider the case of Picard surfaces, with $f$ the map defined in [Go]. If $k(v)=\FF_{p^2}$, then $\tau$ is $2(1,0)+(1,1)+2(0,1)$.
\medskip
{\bf Example 4.} This is the generalization of Example 3. Let $E$ be an imaginary
quadratic field. Let $V$ be a $(m+n)$ dimensional $E$-vector space (with $m,n\in\NN$) and let $\Mj_V:V\times V\to E$ be a non-degenerate hermitian form on $V$ which has signature $(m,n)$ over $\RR$. We can assume $m\ge n$. Let ${\rm Sh}(G,X)$ be the Shimura variety defined as in [Go], for this situation. Similarly to [Go] it can be embedded in a Siegel modular variety $f:{\rm Sh}(G,X)\hookrightarrow {\rm Sh}\bigl({\rm GSp}(W,\psi),S\bigr)$ with
$\dim_\QQ(W)=2(m+n)$. We have $E(G,X)=E$ if $m>n$ and $E(G,X)=\QQ$ if $m=n$. Moreover $G^{\rm der}_\RR=SU(m,n)_{\RR}$. We have:
\medskip
{\bf Fact.} {\it If $v$ is a prime of $E$ such that $(v,2)=1$ and $k(v)=\FF_{p^2}$ and if $(f,L_{(p)},v)$ is a SHS, then $\tau$ is $2n(1,0)+(m-n)(1,1)+2n(0,1)$.}
\medskip
{\bf Proof:} To check this we use 4.1. $G_{W(\FF_{p^2})}$ is a split group. We choose a $W(\FF_{p^2})$-basis $\{e_1,...,e_{m+n},f_1,...,f_{m+n}\}$ of $M:=L_p^*\otimes_{\ZZ_p} W(\FF_{p^2})$ such that:
\medskip
-- $\tilde\psi(e_i,e_j)=\tilde\psi(f_i,f_j)=0$, $\forall i,j\in S(1,n+m)$, while $\tilde\psi(e_i,f_j)$ is $0$ or $1$ depending on $i$ being different or equal to $j$;
\smallskip
-- the $W(\FF_{p^2})$-submodule of $M$ generated by all $e_i$'s (resp. by all $f_i$'s) is $G_{W(\FF_{p^2})}$-invariant;
\smallskip
-- there is a Borel subgroup $B$ of $G_{\ZZ_p}$ normalizing the $W(\FF_{p^2})$-submodule of $M$ generated by $e_1$,..., $e_i$, $\forall i\in S(1,n+m)$, and having the property that the maximal torus of $B_{W(\FF_{p^2})}$ normalizing all $e_i$'s is obtained from a torus of $B$ by extension of scalars.
\medskip
We consider a cocharacter $\mu:\GG_m\to B_{W(\FF_{p^2})}$ which acts trivially on $e_{n+1}$,..., $e_{n+m}$, $f_1$,..., $f_n$ and as the inverse of the identical character on $e_1$,..., $e_n$, $f_{n+1}$,..., $f_{n+m}$. The natural action of ${\rm Gal}(\FF_{p^2}/\FF_p)$ on cocharacters of $G_{W(\FF_{p^2})}$ takes $\mu$ into the cocharacter which acts trivially on $e_{m+1}$,..., $e_{n+m}$, $f_1$,..., $f_m$ and as the inverse of the identical character on $e_1$,..., $e_m$, $f_{m+1}$,..., $f_{n+m}$. Using these and 4.1, the Fact follows.
\medskip
{\bf P5} {\it The slopes of ${\rm Lie}_G(\tau)$ are precisely $-1$, $0$ and $1$ iff $G$ is not a torus and $k(v^{\rm ad})=\FF_p$, where $v^{\rm ad}$ is the prime of $E(G^{\rm ad},X^{\rm ad})$ divided by $v$.}
\medskip
{\bf Proof:} If $G$ is not a torus and $k(v^{\rm ad})=\FF_p$, then the image of $\bar h_i$ in ${\rm Lie}(G^{\rm ad}_{B(k(v))})$ is non-zero and does not depend on $i\in S(1,d)$; so $h^0$ acts on ${\rm Lie}(G_{B(k(v))})$ as $\bar h_1$ does. So the ``if" part follows from b) of 4.4.1 2) and 4.1.1.2. If ${\rm Lie}_G(\tau)$ has precisely the slopes $-1$, $0$ and $1$, then $G$ is not a torus; the fact that $k(v^{\rm ad})=\FF_p$ is argued as in the proof of P1. This ends the proof.
\medskip
 Let ${\got g}:={\rm Lie}(G_{W(k)})$ and let $\bigl({\got g},\vph_0,F^0({\got g}),F^1({\got g})\bigr)$ be the Shimura filtered Lie $\sg_k$-crystal attached to a $G$-canonical lift ${\rm Spec}(W(k))\to\Mn_{k(v)}$. From P5 and F2 of 3.10.7 we get:
\medskip
{\bf P6} {\it ${\got p}_0:=W_0({\got g},\vph_0)$ is $F^0({\got g})$ iff $k(v^{\rm ad})=\FF_p$.}
\medskip
{\bf P7} {\it ${\got p}_0$ can be the Lie subalgebra of a Borel subgroup of $G_{W(k)}$.}
\medskip
This is equivalent to: the rank of $G$ is equal to the number of slopes $0$ of $({\got g},\vph_0)$. From the formulas of 3.10.6 we deduce that, if this is so, then all simple, adjoint factors of $G^{\rm ad}_{\CC}$ are of $A_\ell$ Lie type,  $\ell\in\NN$. What we actually need is: any cyclic adjoint factor of $({\got g},\vph_0,F^0({\got g}),F^1({\got g}))$ has maximal $A$-spreading (cf. F3 of 3.10.7). 
\medskip
{\bf Concrete example.} Let $F$ be a totally real number field which is a Galois extension of $\QQ$. Let $\tilde G$ be an absolutely simple, adjoint group over $F$ of $A_{\ell}$ Lie type, $\ell\ge 2$, such that denoting $G^{\rm ad}:={\rm Res}_{F/\QQ} \tilde G$, $G^{\rm ad}_{\RR}$ is the adjoint group of $\prod_{i=1}^{[{\ell+1\over 2}]} SU(i,\ell+1-i)_{\RR}$. So $[F:\QQ]=[{\ell+1\over 2}]$. We assume the existence of a prime $w_1$ of $F$ unramified over a rational prime $p>2$ and such that: 
\medskip
-- $\tilde G$ is unramified over $\QQ_p$, 
\smallskip
-- $w_1$ is the only prime of $F$ over $p$, and 
\smallskip
-- $\tilde G$ over the completion $F_{w_1}$ of $F$ w.r.t. $w_1$ is non-split. 
\medskip
We consider an arbitrary Shimura pair of the form $(G^{\rm ad},X^{\rm ad})$. Its reflex field $E(G^{\rm ad},X^{\rm ad})$ is a totally imaginary extension of $F$. Let $w$ be a prime of it dividing $w_1$. We have $[k(w):k(w_1)]=2$. Argument: $\tilde G$ splits over the unramified quadratic extension of $F_{w_1}$ (cf. Fact 1 of 4.3.6) and so this degree is at most $2$; on the other hand, due to the assumed non-trivial involution, it can not be 1. 
\smallskip
So, from F3 of 3.10.7 we get that for any SHS $(f,L_{(p)},v)$, with $v$ dividing $w$ and $f:(G,X)\hookrightarrow (GSp(W,\psi),S)$, P7 holds. Such a SHS does exist cf. [Va2, 6.4.2] and 3.2.6; for $p=3$ cf. also 2.3.5.1.
\medskip
{\bf Example 5.} P7 holds for the case considered in Example 3.
\medskip
{\bf Example 6.} If in Example 4 we have $m+n\ge 4$, then ${\got p}_0$ is not the Lie subalgebra of a Borel subgroup. If $m=n$, then ${\rm Lie}_G(\tau)$ is always an ordinary type.
\medskip
{\bf P8} {\it For any rational number $r=a/d$, with $a\in S(1,d)$, there are situations when the slopes of Shimura Lie $F$-crystals attached to $G$-ordinary points of $\Mn_{k(v)}$ are precisely $-r,0$ and $r$.}
\medskip
{\bf Proof:} We take $(G,X)$ to be of $B_\ell$ (or $C_\ell$ or  $D_\ell^\RR$) type and such that $G^{\rm ad}={\rm Res}_{F/\QQ} G^1$, with $F$ a totally real number field satisfying $[F:\QQ]=d$ and with $G^1$ an absolutely simple $F$-group. We assume that for precisely $a$-embeddings of $F$ into $\RR$, $G^1_{\RR}$ is non-compact. We also assume the existence of a prime $v$ of $E(G,X)=F$ (cf. [De2, 2.3.12]) such that $k(v)=\FF_{p^d}$, with $p>2$ a prime for which $G$ is unramified over $\QQ_p$. 
Now for a SHS $(f,L_{(p)},v)$ (cf. [Va2, 6.4.2] for $p\ge 5$; see \S 6 for $p=3$) defined by an injective map $f:(G,X)\hookrightarrow ({\rm GSp}(W,\psi),S)$, the slopes of the Shimura Lie $F$-crystals attached to $G$-ordinary points of $\Mn_{k(v)}$ are precisely
$-{a\over d}$, $0$ and ${a\over d}$ (cf. 3.10.6 i) and Case 1 of 3.10.6 iii); see also b) of 4.4.1 2)).
\smallskip
$\forall r\in\QQ\cap (0,1]$, we can choose a number field $F$ and a reductive group $G^1$ over $F$ such that the above assumptions are satisfied. For instance, if we are dealing with the $B_\ell$ type, fixing $F$ subject to the above requirements, we can take $G^1$ to be the $SO$-group of the quadratic form $a_1x_1^2+a_2x_2^2+x_3^2+x_4^2+ ...x_{2\ell+1}^2$ in $2\ell+1$ variables over $F$, with $a_1$, $a_2\in F$ such that (cf. approximation theory) for precisely $a$ (resp. $d-a$) embeddings $F\hookrightarrow\RR$, they are both negative (resp. positive).
Starting from 3.10.6 iii) and iv), we can construct examples for P8 with $(G,X)$ of $A_\ell$ or $D_\ell^\HH$ type; we leave this to the reader.
\medskip
{\bf P9} {\it If the slopes of Shimura $F$-crystals attached to $G$-ordinary points of $\Mn_{k(v)}$ are rational numbers of the interval $(0,1)$, then no cyclic factor of the Shimura adjoint Lie $\sg_{\bar k}$-crystal attached to a $\bar k$-valued $G$-ordinary point of $\Mn_{k(v)}$ is totally non-compact (see def. 3.10.5). The converse of this is not true.}
\medskip
{\bf Proof:} The first part is a direct consequence of F5 of 3.10.7. The second part can be seen through trivial examples. 
\smallskip
As a corollary we get:
\medskip
{\bf P10} {\it If there is a totally non-compact cyclic factor of the Shimura adjoint Lie $\bar\sg$-crystal attached to an $\FF$-valued $G$-ordinary point of $\Mn_{k(v)}$, then the abelian varieties obtained by pulling back $\Ma$ through $G$-ordinary points with values in algebraically closed fields, have non-zero $p$-ranks.}
\medskip
For future references we now include an extra elementary property. We assume that the representation $G_{\CC}\to GL(W\otimes_{\QQ} \CC)$ is irreducible. So $G^{\rm ad}_{\RR}$ has only 1 simple, non-compact factor. Let $n_f$ be the number of simple factors of $G^{\rm ad}_{\RR}$. It is well known that it is an odd number, that $G^{\rm ad}$ is a simple $\QQ$--group, and that $2e$ is an $n_f$-th power of a natural number (for instance, see [Pi, \S 4]). We refer to 4.3.1. $\Mh^{\rm nc}$ has only 1 element; we denote it by $i$. As $G^{\rm ab}=\GG_m$, we have $E(G^{\rm ad},X^{\rm ad})=E(G,X)$. As $G^{\rm ad}_{\RR}$ has only 1 simple, non-compact factor, $d_i|d$; so $d\in\{d_i,2d_i\}$ (cf. Fact 1 of 4.3.6).
\medskip
{\bf P11} {\it The abelian varieties obtained by pulling back $\Ma$ through $G$-ordinary points with values in algebraically closed fields, have non-zero $p$-ranks: in case $d=d_i$ (resp. in case $d=2d_i$), these $p$-ranks are equal to ${e\over {2^{d-1}}}$ (resp. are ${e\over {2^{d-1}}}\bigl({{m+1}\over {(2m+1)}}\bigr)^d$; here $m\in\NN$ is such that $G^{\rm ad}$ is of $A_{2m+1}$ Lie type). So, these $p$-ranks are 1 precisely in the case when $G^{\rm ad}_{\QQ_p}$ is a simple $\QQ_p$-group of $A_1$ Lie type.}
\medskip
{\bf Proof:} Based on 4.2.1 a), we can work in the context of 4.1.1. The faithful representation of $G^{\rm der}_{B(\FF)}$ on $M\otimes_{W(k(v))} B(\FF)$ is the tensor product representation of $n_f$ irreducible alternating representations $V_1$,..., $V_{n_f}$; here $V_i$'s are $B(\FF)$-vector spaces of the same dimension. So we can identify $M\otimes_{W(k(v))} B(\FF)$ with $V_1\otimes_{B(\FF)} V_2\otimes_{B(\FF)}... \otimes_{B(\FF)} V_{n_f}$. The extension $\mu_{B(\FF)}$ to $B(\FF)$ of a cocharacter $\mu:\GG_m\to G_{W(k(v))}$ as in 4.1.1, acts on $M\otimes_{W(k(v))} B(\FF)$ via one such irreducible representation $V_{i_{\mu}}$ (with $i_{\mu}\in S(1,n_f)$) and so we can identify $\mu$ with a cocharacter of $GL(V_{i_{\mu}})$. Moreover, ${\rm Gal}(\FF)$ acts on the image of $\mu_{B(\FF)}$ via its quotient ${\rm Gal}(k(v)/\FF_p)$. If $d=d_i$ (resp. if $d=2d_i$), then two different elements $\gamma_1,\gamma_2\in {\rm Gal}(k(v)/\FF_p)$ can not take this image to factor through the same irreducible representation (resp. can take this image to factor through the same irreducible representation $V_{\tilde i_{\mu}}$ iff $\gamma_1\gamma_2^{-1}$ is the generator $\tilde\tau$ of ${\rm Gal}(k(v)/\FF_{p^{d_i}})=\ZZ/2\ZZ$). It is convenient to identify the set $I_p(G^{\rm ad})$ with the set $SS:=S(1,n_f)$ in such a way that for any $j\in SS$, ${\rm Lie}(G^{\rm der}_{B(\FF)})$ acts on $V_j$ via ${\rm Lie}(G_j)$; here we use the notations of 4.3.1.1 with $k=\FF$.
\smallskip
We first treat the case $d_i=d$. We consider the natural direct sum decomposition
$$V_j=V_j^1\oplus V_j^0,$$
with $V_j^0$ (resp. with $V_j^1$) as the maximal $B(\FF)$-vector subspace of $V_j$ on which the extension to $B(\FF)$ of $\sg^s\mu({1\over p})\sg^{-s}$ acts trivially, $\forall s\in S(1,d)$ (resp. on which at least one of these $\sg^s\mu({1\over p})\sg^{-s}$'s does not act trivially). If $j\in\Mh_i$, then $\dim_{B(\FF)}(V_j^1)=\dim_{B(\FF)}(V_j^0)$. So the case $d_i=d$ follows by elementary computation, once we remark that the $B(\FF)$-subspace of $V_1\otimes_{B(\FF)} V_2\otimes_{B(\FF)}... \otimes_{B(\FF)} V_{n_f}$ corresponding to the slope $1$ for $(M\otimes_{W(k(v))} W(\FF),\vph\otimes 1)$ is 
$$V_1^{a_1}\otimes_{B(\FF)} V_2^{a_2}\otimes_{B(\FF)}... \otimes_{B(\FF)} V_{n_f}^{a_{n_f}},$$
where $a_j$ is $1$ or $0$ depending on the fact that $j$ is or is not in $\Mh_i$, $\forall j\in SS$.
\smallskip
If $d=2d_i$, then $G^{\rm ad}$ is of $A_{2m+1}$ Lie type, with $m\in\NN$, and the irreducible representations mentioned are associated to the minimal weight $\bar\om_{m+1}$ of the $A_{2m+1}$ Lie type (for instance, see loc. cit.; see [Bou2, planche I] for weights). Similarly we get a direct sum decomposition 
$$V_j=V_j^0\oplus V_j^{1\over 2}\oplus V_j^1,$$ 
with $V_j^0$ as above and with $V_j^{1\over 2}$ and $V_j^1$ normalized by all above mentioned cocharacters, such that $V_j^0\oplus V_j^{1\over 2}$ is the $B(\FF)$-vector subspace of $V_j$ perpendicular on $V_j^0$ w.r.t. the alternating form on it centralized by the natural image of ${\rm Lie}(G_{jB(\FF)})$ in ${\rm End}(V_j)$. So this case follows, via a similar trivial computation, once we remark that $\mu_{B(\FF)}$ (resp. that the extension to $B(\FF)$ of the cocharacter $\tilde\tau\circ\mu\circ\tilde\tau^{-1}$ of $G_{W(k(v))}$) is associated to the first (resp. to the last) node of the Dynkin diagram of the $A_{2m+1}$ Lie type: for $j\in\Mh_i$, $\dim_{B(\FF)}(V_j^1)=C_{2m}^m$; so the $p$-ranks in this case are $(C_{2m}^m)^{d_i}(C_{2m+2}^{m+1})^{n_f-d_i}={e\over {2^{d-1}}}\bigl({{m+1}\over {(2m+1)}}\bigr)^d$. 
\medskip
{\bf Remark.} If $d=d_i$, then using the fact that $V_j^1$ and $V_j^0$ have the same dimension, $\forall j\in\Mh_i$,  we get that the number of slopes ${l\over d}$ of $\tau$ is exactly ${{C_d^l}\over {2^{d-1}}}e$; here $l\in S(0,d)$, cf. P4. Similarly, if $d=2d_i$, we can express the number of slopes ${l\over d}$ of $\tau$ as a sum involving combinatorial numbers.  
\medskip
{\bf Example 7: The case of curves.} We assume that $\dim_{\CC}(X)=1$, i.e. that $G^{\rm ad}={\rm Res}_{F/\QQ} G^1$, with $F$ a totally real number
field and with $G^1$ an absolutely simple $F$-group of $A_1$ Lie type such that for exactly one embedding of $F$ into $\RR$, $G^1_\RR$ is a non-compact (so split) group (cf. [De2, 2.3.4]). We impose no restrictions on $(G^{\rm ab},X^{\rm ab})$. Let $(f,L_{(p)},v)$ be a SHS defined by an injective map $f:(G,X)\hookrightarrow ({\rm GSp}(W,\psi),S)$. Let $v^{\rm ad}$ be the prime of $E(G^{\rm ad},X^{\rm ad})=F$ (cf. [De2, 2.3.12]) divided by $v$.
\smallskip
Let $G^{\rm ad}_{\QQ_p}=\prod_{i\in\Mh} {\rm Res}_{F_i/\QQ_p} G^1_i$, where $F\otimes_\QQ\QQ_p=\prod_{i\in\Mh}F_i$, with $F_i$ an unramified finite field extension of $\QQ_p$, and where $G^1_i$ is an absolutely simple $F_i$-group. Let $i_0\in\Mh$ corresponding to the non-compact factor of $G^{\rm ad}_\RR$. Let $d_0:=[F_{i_0}:\QQ_p]\ge 1$. Then 
$$k(v^{\rm ad})=\FF_{p^{d_0}}$$ 
and the slopes of the Shimura adjoint Lie $F$-crystal attached to a $G$-ordinary point of $\Mn_{k(v)}$ are precisely $-{1\over d_0}$,  $0$ and ${1\over d_0}$; the multiplicity of $-1/d_0$ (or of $1/d_0$) is $d_0$, while the multiplicity of $0$ is $3\bigl(\dim_\QQ(F)-d_0\bigr)+d_0=3\dim_\QQ(F)-2d_0$ (cf. 3.10.6 i)). So if $d_0=1$ P10 applies.
\smallskip
We assume now that moreover the representation $G_{\CC}\to GL(W\otimes_{\QQ} \CC)$ is irreducible. So $G^{\rm ab}=\GG_m$. So $E(G,X)=E(G^{\rm ad},X^{\rm ad})$ is $F$ itself (cf. loc. cit.). Then the $p$-rank of any abelian variety obtained by pulling back $\Ma$ through a $G$-ordinary point with values in an algebraically closed field, is precisely $2^{\dim_{\QQ}(F)-d_0}$ (cf. P11). 
\medskip
{\bf 4.6.1. Complements. 1)} There are plenty of examples when the Newton polygon of the Shimura-ordinary type of a SHS $(f,L_{(p)},v)$ does not have integral slopes. Such examples can be constructed starting from [Va2, 6.5.1.1] and its proof. We briefly recall from loc. cit. how we can construct them, by restricting to very particular situations. 
\smallskip
Let $F$ be a totally real number field such that it has precisely one prime $v$ dividing $p$ and this prime is unramified over $p$. We assume $F\neq\QQ$. Let $E$ (resp. $K^0$) be a totally imaginary quadratic extension of $F$ in which $v$ splits (resp. which has precisely one prime $w$ unramified over $v$). Let (this is reviewed in 4.6.7.1 below) $\tilde G$ be a simply connected group over $F_{(p)}$ such that:
\medskip
-- it splits over $E_{(p)}$;
\smallskip
-- its adjoint is an absolutely simple $F_{(p)}$-group of $C_{l}$ Lie type;
\smallskip
-- all its extensions to $\RR$ are either compact or split and precisely one such embedding is compact.
\medskip
Let $(W_{E_{(p)}},\tilde\psi)$ be a symplectic space over $E_{(p)}$ such that $\tilde G_{E_{(p)}}=Sp(W_{E_{(p)}},\tilde\psi)$. Let $L_{(p)}:=W_{E_{(p)}}\otimes_{F_{(p)}} K^0_{(p)}$ and let $W:=L_{(p)}\otimes_{Z_{(p)}} \QQ$. In loc. cit. it is constructed a perfect alternating form $\psi:L_{(p)}\otimes_{\ZZ_{(p)}} L_{(p)}\to\ZZ_{(p)}$ and an injective map $f:(G,X)\hookrightarrow (GSp(W,\psi),S)$, such that:
\medskip
-- the Zariski closure of $G$ in $GL(L_{(p)})$ is a reductive group $G_{\ZZ_{(p)}}$, with $G_{\ZZ_{(p)}}^{\rm der}={\rm Res}_{F_{(p)}/\ZZ_{(p)}} \tilde G$, and \smallskip
-- under this identification, the representation of $G_{\ZZ_{(p)}}^{\rm der}$ on $L_{(p)}$ is obtained naturally via the tautological representation of $\tilde G_{E_{(p)}}$ on $W_{E_{(p)}}$, and
\smallskip
-- we have $E(G,X)=K^0$. 
\medskip
We assume that the triple $(f,L_{(p)},w)$ is a SHS (for instance, this is so if $p\not | 6(l+1)$, cf. [Va2, 6.5.1.1 vi)]). We have a direct sum decomposition 
$$L_{(p)}^*\otimes_{\ZZ_{(p)}} \ZZ_p=L_p^1\oplus L_p^2$$ 
of $G_{\ZZ_p}$-modules, corresponding to the $2$ primes of $E$ dividing $v$. Moreover, $L_p^s\otimes_{\ZZ_p} W(k(w))$ is a direct sum of $[k(w):\FF_p]$ absolutely irreducible $G_{W(k(w))}$-modules, whose simple factors are permuted transitively by $\sg_{k(w)}$, $s=\overline{1,2}$. From this and the fact that $\mu$ of 4.1 acts trivially precisely on one such simple factor (this can be read out --cf. loc. cit.-- from [De2, 2.3.9]), we get that the Shimura-ordinary type defined by $(f,L_{(p)},w)$ does not have integral slopes.
\medskip
The examples hinted at above are such that the irreducible representations of the faithful representation $G^{\rm der}_{\CC}\hookrightarrow GL(W\otimes_{\QQ} \CC)$ involve no tensors products. So, in some sense, we are at the opposite pole of the situation of P11.  
\smallskip
{\bf 2)} In [Va5] we will show that the results of 4.1-6 remain true (under proper formulation) for the integral canonical models to be constructed in \S 6 for quotients of Shimura varieties of Hodge type w.r.t. large classes of non-hyperspecial subgroups.
\smallskip
{\bf 3)} We reobtained the ``ordinary theory" for special fibres of  integral canonical models of Siegel modular varieties, w.r.t. primes $p\ge 2$ (for $p=2$, cf. 4.14.3 below). In particular, we obtained a
completely new proof of the density (see [Kob], [NO], [FC] and [EO]) of ordinary points in these special fibres.
\smallskip
{\bf 4)} The whole of 4.1 as well as properties P1 to P11 above remain true under proper formulation in the context of Shimura $F$-crystals over perfect fields (cf. 3.11.2). The only difference is: we have to change $k(v)=\FF_{p^d}$ (with $d\in\NN$) with the degree of definition $(=d)$, and  $k(v^{\rm ad})=\FF_{p^{d_0}}$, with the $A$-degree of definition $(=d_0)$; moreover, in P3 (resp. P5) we have to change the condition $G$ is (resp. is not) a torus to any one of the equivalent conditions 3.2.2 a) to d) (resp. with the condition that 3.2.2 a) to d) do not hold).
\smallskip
{\bf 5)} In particular, 4.1 and P1 to P11 can be applied to the context of abelian varieties with Hodge cycles without an a priori defined principal polarization (for instance, there is no problem if we have a polarization which is not principal), as long as we are in a context of Shimura $F$-crystals over perfect fields to which we can perform 4.2.3-4 (so we are in a reductive context even in the \'etale $\ZZ_p$-context; in \S 5 we will see that this is automatically so, cf. 1.15.1).
\smallskip
{\bf 6)} It can happen that we have two standard Hodge situations $(f,L_{(p)},v_i)$, $i=\overline{1,2}$, with $v_1$ and $v_2$ dividing the same rational prime $p\ge 3$, such that their attached Shimura-ordinary types are different. For instance, this is the case if $k(v_1)=\FF_p$, while $k(v_2)=\FF_{p^d}$, with $d\ge 2$, cf. P1. 
\medskip
{\bf 4.6.2. An application.} Let $E$ be a number field and let $A$ be an abelian variety over $E$. Let ${\rm Spec}(O^r_E)\to{\rm Spec}(O_E)$ be the maximal open subscheme of the spectrum of the ring of integers $O_E$ of $E$, over which $A$ extends to an abelian scheme $A^r$. Let $E(A)$ be the reflex field of the Shimura variety attached to $A$ (cf. 2.1).
\medskip
{\bf 4.6.2.1. The proof of Theorem 7 of 1.7.} We now prove Theorem 7 of 1.7. So we assume $E(A)\neq\QQ$. We consider an embedding $E\hookrightarrow\CC$; let 
$\tilde L_\ZZ:=H_1(A_{\CC},\ZZ)$ be the first group of the Betti homology of $A_\CC$. The Mumford--Tate group of $A$ (of $A_\CC$) is a reductive subgroup $G_A$ of $GL(\tilde L_\ZZ\otimes_{\ZZ} \QQ)$. There is $N_0(A)\in\NN$ such that for any prime $p\ge N_0(A)$, the Zariski closure of $G_A$ in $GL(\tilde L_\ZZ\otimes_{\ZZ} \ZZ_{(p)})$ is a reductive group over $\ZZ_{(p)}$. If also $p\ge\max\{\dim_E(A),5\}$, then (cf. [Va2, 5.8.6 and def. 5.8.1]) there is a $\ZZ_{(p)}$-well positioned family of tensors of $\Mt(\tilde L_\ZZ\otimes_{\ZZ} \ZZ_{(p)})$ for the group $G_A$. Let now $v$ be a prime of $O^r_E$ dividing a rational prime 
$$p\ge\max\{N_0(A),5,\dim_E(A)\}.$$ 
\indent
Let $V:=\widehat{O_{(v)}^{\rm sh}}$. So $V$ is a finite, flat extension of $V_0:=W(\FF)$. Let $A_V$ be the abelian variety over $V$, obtained from $A^r$ through the logical  morphism ${\rm Spec}(V)\to {\rm Spec}(O_E^r)$. All Hodge cycles of $A_V$ are defined over $V_1[{1\over p}]$, with $V_1$ a finite, flat, DVR extension of $V$ (cf. [De4]; in fact it is very easy to see that we can take $V_1=V$ but this is irrelevant here). We choose a family ${(w_\al)}_{\al\in\Mj(V_1)}$ of Hodge cycles of $A_{V_1}$ such that $G_A$ is the subgroup of $GL(\tilde L_\ZZ\otimes_{\ZZ} \QQ)$ fixing the Betti realizations of its members.
\smallskip
Let $\bar e:=[V_1:V_0]$. In what follows, to simplify the presentation, we use similar notations as before, though the situation is entirely unrelated to a fixed SHS. Repeating the arguments of [Va2, 5.2-3], we get an abelian scheme $A_{\tilde R{\bar e}}$ endowed with a family $({w_\al}^{\tilde R{\bar e}})_{\al\in\Mj(V_1)}$ of Hodge cycles over the $V_0$-algebra $\tilde R{\bar e}$ used in [Va2, 5.2.1], such that: 
\medskip
{\bf i)} through the canonical $V_0$-epimorphism $\tilde R{\bar e}\twoheadrightarrow V_1$ defined by a choice of a uniformizer of $V_1$ (see [Va2, 5.2]), it becomes $A_{V_1}$, together with ${(w_\al)}_{\al\in\Mj(V_1)}$;
\smallskip
{\bf ii)} through the canonical $V_0$-epimorphism $\tilde R{\bar e}\twoheadrightarrow V_0$, it becomes an abelian variety $A_{V_0}$ for which the quintuple $(M,F^1,\vph,G_{V_0},(t^0_{\al})_{\al\in\Mj(V_1)})$ is a Shimura filtered $\sg_{\FF}$-crystal; here $M:=H^1_{\rm crys}(A_{\FF}/V_0)$, $F^1$ is the Hodge filtration of $M$ defined by $A_{V_0}$, $\vph$ is the Frobenius endomorphism of $M$, ${(t^0_\al)}_{\al\in\Mj(V_1)}$ is the family of de Rham components of the family of Hodge cycles ${(w^0_\al)}_{\al\in\Mj(V_1)}$ of $A_{V_0}$ obtained naturally from $(w_{\al}^{\tilde R_{\bar e}})_{\al\in\Mj(V_1)}$, and $G_{V_0}$ is the Zariski closure in $GL(M)$ of the subgroup of $GL(M[{1\over p}])$ fixing $t^0_\al$, $\forall\al\in\Mj(V_1)$ (see [Va2, 5.3.4]).
\medskip
Let $\Mp$ be the Newton polygon of a (any) $G_{V_0}$-ordinary $\sg_{\FF}$-crystals produced by $(M,\vph,G_{V_0})$. From 3.1.0 a) we get:
\medskip
{\bf (4.6.2.1.1)} {\it The Newton polygon of the abelian variety over $k(v)$ obtained from $A$ by reduction w.r.t. $v$, is above $\Mp$.} 
\medskip
Let $E_1$ be the composite field of $E$ and $E(A)$ and let $w$ be a prime of $E_1$ dividing $v$. Let $w_A$ be the prime of $E(A)$ divided by $w$. We have:
\medskip
{\bf (4.6.2.1.2)} {\it If $k(w_A)\ne\FF_p$, then $\Mp$ is not the Newton polygon of an ordinary type.} 
\medskip
This results from 4.6 P1: the principal polarization existing in the context of a SHS played no role in obtaining 4.1 and P1 to P11 above (cf. 4.6.1 5)). It is also worth pointing out: as the Newton polygon is an isogeny invariant, [Mu, cor. 1 of p. 234] allows us to assume that $A$ is principally polarized.
\medskip
Now Theorem 7 of 1.7 is a consequence of (4.6.2.1.1-2): there are infinitely many primes of $E_1$ dividing primes of $E(A)$ not having residue fields with a prime number of elements.
\medskip
{\bf 4.6.2.2. Remark.} We have the following precise form of Theorem 7 of 1.7. Let $v$ be a prime of $O^r_E$ such that:
\medskip
a) it divides a prime $p\ge 3$ such that $G_A$ is unramified over $\QQ_p$;
\smallskip
b) there is a prime $w$ of $E_1$ dividing $v$ and dividing a prime $w_A$ of $E(A)$ for which
 $k(w_A)$ has more than $p$ elements.
\medskip 
Then $A$ has a non-ordinary reduction w.r.t. $v$ (and moreover there is a Newton polygon $\Mp$ of an abelian variety $B$ over $\FF$, with $\dim_{\FF}(B)=\dim_{\FF}(A)$, which is the Newton polygon of a Shimura-ordinary type --pertaining to $G_A$-- which is not an ordinary type and such that the Newton polygon of the reduction of $A$ w.r.t. $v$ is above $\Mp$). 
\medskip
This will become transparent in the light of \S 5-8, cf. 2.3.8 2). In [Va5] we will see that in a), in most situations, we can allow $p$ to be $2$ as well.
\medskip
{\bf 4.6.2.3. Questions.} {\bf 1)} Let $O^r_{E_1}$ be the normalization of $O^r_E$ in $E_1$ and let $A^r_{E_1}$ be the abelian scheme over ${\rm Spec}(O^r_{E_1})$, obtained from $A^r$ by pull back. We do not assume $E(A)$ is not $\QQ$. Is it true that for an infinite set of primes $\Ml$ of $O^r_{E_1}$, $A^r_{E_1}$ has a $G_A$-ordinary reduction which is not ordinary?
\smallskip
{\bf 2)} For $v$ as in 4.6.2.2, let $(M,F^1,\vph,G_{V_0})$ be a Shimura $\sg_{\FF}$-crystal constructed as in 4.6.2.1. How is its $CM$-deviation (or $SCM$-deviation) varying in terms of $v$?
\medskip
{\bf 4.6.2.4. Remark.} Property P1 supports the ordinary conjecture.
\medskip
{\bf 4.6.2.5. Remark.}
We regard the above non-ordinary reduction criterion as the natural generalization of the (well known) fact that any abelian variety over a number field $E$, having complex multiplication over $\bar E$, has a non-ordinary reduction w.r.t. an infinite set of primes of $E$.
\medskip
{\bf 4.6.2.6. Exercise.}
Let $(f,L_{(p)},v)$ be a SHS. Let $\tau$ be the Shimura-ordinary type attached to it. Show that there is an infinite number of standard Hodge situations $(f,L_{(q)},v_q)$ (with $v_q$ a prime of $E(G,X)$ dividing an odd rational prime q) which have $\tau$ as their attached Shimura-ordinary type. Hint: If the reflex field $E(G,X)$ is a Galois extension of $\QQ$, this is obvious; now consider the Galois extension $E_1(G,X)$ of $\QQ$ generated by $E(G,X)$, and use 4.1 and Tchebotarev's density theorem.
\smallskip
This simple fact makes 4.6.2.3 1) an interesting question and suggests that the ordinary conjecture should admit variants for Shimura-ordinary types, which are not ordinary. See also 3.13.4 7).
\medskip
{\bf 4.6.3. The case when the degree of definition is 1.} For the sake of future references we point out the following two particular cases of P1 and P2. 
\smallskip
{\bf A. The split case.} {\it If $(f,L_{(p)},v)$ is a SHS and $G_{\QQ_p}$ is a split reductive group, then the Shimura-ordinary type we obtain is the ordinary type and the abelian varieties over Witt rings of perfect fields, obtained from $\Ma$ by pull back through $G$-canonical lifts of $\Mn$, are (usual) canonical lifts.}
\medskip
{\bf B.} {\it If $(f,L_{(p)},v)$ is a SHS and if an $\FF$-valued $G_{W(\FF)}$-ordinary point of $\Mn_{k(v)}$ gives birth to a Shimura $\sg_{\FF}$-crystal whose  degree of definition is 1, then all abelian varieties over $W(\FF)$ obtained from $\Ma$ by pull back through $W(\FF)$-valued $G_{W(\FF)}$-canonical lifts of $\Mn$, are (usual) canonical lifts and so have complex multiplication.} 
\medskip
{\bf 4.6.4. An application to [Va2, 6.8.6].} We consider the situation and notations described in [Va2, 6.8.0 and 6.8.6] but we do not assume that we have a non-trivial involution (cf. end of AE.0). So $E(G^{\rm ad},X^{\rm ad})$ is a totally imaginary quadratic extension of a totally real number field (resp. is a totally real number field) if the involution is non-trivial (resp. is trivial). We assume we are dealing with a prime $v^{\rm ad}$ of $E(G^{\rm ad},X^{\rm ad})$ dividing a prime $p>2$. Also we do not assume the $\QQ$--rank of $G^{\rm ad}$ is positive: we just assume that all simple factors of $G^{\rm ad}_{\RR}$ are non-compact.
\medskip
{\bf A. Exercise.} Show that the criterion b) of [Va2, 6.8.2] is satisfied. Hint: the situation gets reduced to a standard PEL situation; so we can use 4.4.5, the density part of 4.2.1 and [Mi1, 4.11-2].
\medskip
{\bf Solution.} We use the notations of [Va2, 6.8]; so $V_0:=W(\FF)$. The idea is: any $V_0$-valued Shimura-canonical lift $z$ of $\Mm^\prime_{V_0}$ lifts to a $V_0$-valued point of $\Mn_{V_0}$. If the generic fibre of $z$ is a special point of the generic fibre of $\Mm^\prime_{V_0}$, then this is a consequence of [Mi1, 4.11-2], cf. [Va2, 3.2.3.1 5)] (for the passage from $K_0$-valued points to $V_0$-valued points) and cf. the following Lemma.
\medskip
{\bf Lemma 1.} {\it Let $f_1:(T_1,\{h_1\})\to (T,\{h\})$ be a cover between two $0$ dimensional Shimura pairs. Let $p$ be an arbitrary prime such that the subtorus $T_0$ of $T_1$ which is the kernel of $f_1$ is unramified over $\QQ_p$. Let $H_{1p}\subset T_1(\QQ_p)$ be a compact, open subgroup containing the hyperspecial subgroup $H_{0p}$ of $T_0(\QQ_p)$. Let $H_p$ be the image of $H_{1p}$ in $T(\QQ_p)$. Then any $K_0$-valued point of ${\rm Sh}_{H_p}((T,\{h\})$ lifts to a $K_0$-valued point of ${\rm Sh}_{H_{1p}}((T_1,\{h_1\})$.}
\medskip
{\bf Proof:} This is a consequence of the reciprocity map (for instance, see [Va2, 2.7]) and of class field theory, once we remark that we have a natural identification
$$T_1(\QQ_p)/H_{1p}T_1(\QQ)=T(\QQ_p)/H_pT(\QQ).$$
This identification is a consequence of the fact that $f_1$ is a cover and of the fact (see [Mi3, 4.10]) that $T_0(\QQ_p)=T_0(\QQ)H_{0p}$.
\medskip
We are left to show that the situation gets reduced to a situation in which we know that always the generic fibre of $z$ is a special point. In other words, based on 4.4.5, we are left to show that the situation gets reduced to a standard PEL situation. But this is the statement of [Va2, 6.8.6], cf. end of AE.0: in what follows we present two ways to argue this (and so implicitly two solutions of the Exercise), in some general forms convenient for future references. 
\medskip
{\bf First way.} In what follows we mainly just recall parts of [Va2, 6.5.1.1] to explain the referred statement of [Va2, 6.8.6]. Working with totally independent notations, we have the following variant of 2.3.5.1 and of its $p=2$ analogue.
\medskip
{\bf Lemma 2.} {\it Let $(G_0,X_0)$ be an adjoint Shimura pair whose simple factors are of $C_n$ or $D_n^{\HH}$ type for some $n\in\NN$. We assume $G_{0\RR}$ has no compact factors. Let $p\ge 3$ (resp. $p=2$) be a prime such that $G_0$ is unramified over $\QQ_p$. Then there is an injective map $f\colon (G,X)\hookrightarrow ({\rm GSp}(W,\psi),S)$, with $(G^{\rm ad},X^{\rm ad})=(G_0,X_0)$, and there is a $\ZZ_{(p)}$-lattice $L_{(p)}$ of $W$ good w.r.t. $f$, such that for any prime $v$ of $E(G,X)$ dividing $p$ we get a standard (resp. a $p=2$ standard) PEL situation $(f,L_{(p)},v,\Mb)$, for a suitable $\ZZ_{(p)}$-subalgebra $\Mb$ of ${\rm End}(L_{(p)})$.}
\medskip
{\bf  Proof:} As in [Va2, 6.5.1], we can assume $(G_0,X_0)$ is a simple Shimura pair. [Va2, i) to v) of 6.5.1.1] tells us how to construct an injective map $f\colon (G,X)\hookrightarrow ({\rm GSp}(W,\psi),S)$, with $(G^{\rm ad},X^{\rm ad})=(G_0,X_0)$ and how to choose a $\ZZ_{(p)}$-lattice $L_{(p)}$ of $W$ good w.r.t. $f$, such that, denoting by $G_{\ZZ_{(p)}}$ the Zariski closure of $G$ in $GL(L_{(p)})$ we have:
\medskip
-- the derived group $G^{\rm der}_{V_0}$ is a product of isomorphic semisimple groups having absolutely simple adjoints, and 
\smallskip
-- the faithful representation $G_{V_0}^{\rm der}\hookrightarrow GL(L_{(p)}\otimes_{\ZZ_{(p)}} V_0)$ is a direct sum of standard irreducible representations of rank $2n$ of the factors of this product. 
\medskip
Warning: only [Va2, 6.5.1.1 vi)] is worked out under the assumption that $p$ does not divide the number $B(G_0)$ introduced in [Va2, 5.7.2]; so, as we referred to [Va2, i) to v) of 6.5.1.1], we do not need to assume that $p$ does not divide $B(G_0)$. So in what follows we take $p\ge 2$. 
\smallskip
[Va2, 6.5.1.1] is a $\ZZ_{(p)}$-version of [De2, 2.3.10]. So, as $G_{0\RR}$ has no compact factors, [De2, 2.3.13] allows us to assume $G^{\rm ab}=\GG_m$. So the double centralizer $C$ (i.e. the centralizer of the centralizer) of $G^{\rm der}_{\ZZ_{(p)}}$ in $GL(L_{(p)})$ is such that its extension to $V_0$ is a product of as many copies of $GL(V_0^{2n})$ as factors of $G_{0\RR}$ (for $p=2$ cf. also the absolutely irreducible part of 2.3.18 B2); moreover, the representation of $C_{V_0}$ on $L_{(p)}\otimes_{\ZZ_{(p)}} V_0$ is a direct sum of standard irreducible representations of rank $2n$ of these factors. $C$ contains $G_{\ZZ_{(p)}}$ and ${\rm Lie}(C)$ is stable under the involution of ${\rm End}(L_{(p)})$ defined by $\psi$. From the classification of [Ko2, top of p. 375] we get that the maximal integral subgroup $G_{1\ZZ_{(p)}}$ of the intersection $C\cap GSp(L_{(p)},\psi)$ has $G^{\rm der}_{\ZZ_{(p)}}$ as a subgroup and its generic fibre is generated by $G^{\rm der}$ and by a subtori of $Z(C_{\QQ})$ (for $p=2$ cf. also 2.3.18 B2-3); so based on [Va2, 3.1.6 and 4.3.9] we get that $G_{1\ZZ_{(p)}}$ is reductive. So we just have to replace:
\medskip
-- $G$ by the generic fibre $G_1$ of $G_{1\ZZ_{(p)}}$;
\smallskip
-- $(G,X)$ by the Shimura pair $(G_1,X_1)$ having the property that we have a natural injective map $(G,X)\hookrightarrow (G_1,X_1)$;
\smallskip
-- to replace $f$ by the resulting injective map $(G_1,X_1)\hookrightarrow (GSp(W,\psi),S)$.
\medskip\noindent
Moreover, we have to take $\Mb$ as the $\ZZ_{(p)}$-subalgebra of ${\rm End}(L_{(p)})$ formed by all elements fixed by $G_1$. 
\smallskip
For $p\ge 3$ it is well known that the resulting quadruple $(f,L_{(p)},v,\Mb)$ is a standard PEL situation (for instance, cf. [Va2, 5.6.3] and AE.1; see also 2.3.6 and 2.3.8 4)). For $p=2$, it is 2.3.18 B which guarantees that the quadruple $(f,L_{(2)},v,\Mb)$ is a $p=2$ standard PEL situation. This ends the proof of Lemma 2 and so of the first solution of the Exercise.    
\medskip
{\bf Second way.} This second way is in essence just an elaboration of the first way, which from some points of view is more convenient (and precise). It is entirely detailed in [Va10, 3.0-1 and the paragraph after 3.1]; though [Va10] is a continuation of this paper, the mentioned part of it, is entirely independent of it. However, for the convenience of the reader we recall most of the details, in the following two Steps. 
\smallskip
{\bf Step 1.} Let $F$ be the totally real number field such that $G^{\rm ad}$ is the Weil restriction from $F$ to $\QQ$ of an absolutely simple $F$-group $\tilde G$. As $p$ is such that $G$ is unramified over $\QQ_p$, $F$ is unramified above $p$ (see [Va2, 6.5.1]). Following the pattern of [Va2, 6.6.2], based on [Va2, 4.3.16] and on 4.11.3 below, we can replace $F$ (resp. $\tilde G$) by any other totally real number field $F_1$ unramified above $p$ and containing $F$ (resp. by $\tilde G_{F_1}$). To avoid forwarding, we present a second reason why we can perform this replacement: we just have to apply the following general Lemma to the context (of Shimura triples) of [Va2, 6.2.1, 6.8.0 and 6.8.6]; its notations are entirely independent of the previous ones and it will be needed later on in connection to Shimura varieties of special type. For the meaning of EP, EEP and SEP we refer to [Va2, 3.2.3 3) and 4)].
\medskip 
{\bf Lemma 3.} {\it Let $f_1:(G,X,H)\hookrightarrow (G_1,X_1,H_1)$ be an injective map of Shimura triples. We assume that $H$ is a subgroup of $G(\QQ_p)$, with $p$ a prime which is relatively prime to $2$ times the order of the center of the simply connected semisimple group cover of $G$, that $E(G,X)=E(G^{\rm ad},X^{\rm ad})$ and that $G_{\CC}^{\rm der}$ and $G_{1\CC}^{\rm der}$ are products of isomorphic semisimple groups whose adjoints are simple. Let $f_2:(G_2,X_2,H_2)\to (G_1,X_1,H_1)$ be a cover (see defs. of [Va2, 2.4.0 and 3.2.6]), with $G^{\rm der}_2$ a simply connected semisimple group and with $E(G_2,X_2)=E(G_1,X_1)$ (cf. [Va2, 3.2.7 10)]). Let $(G_2^\prime,X_2^\prime,H_2^\prime)$ be an arbitrary Shimura triple having the same adjoint as $(G_2,X_2,H_2)$ and such that $G_2^{\prime\rm der}$ is a simply connected semisimple group. We also assume that $(G,X,H)$ and $(G_2^\prime,X_2^\prime,H_2^\prime)$ have integral canonical models $\Mn$ and respectively $\Mn_2^\prime$, that $\Mn_2^\prime$ has the EP and that $\Mn_{2W(\FF)}^\prime$ is quasi-projective and has the EEP. Then any Shimura triple having the same adjoint as $(G,X,H)$, has an integral canonical model.}
\medskip
{\bf Proof:} We consider, as in [Va2, 3.2.7 3)], the fibre product $f_{3}:(G_3,X_3,H_3)\to (G,X,H)$ and $f_{32}:(G_3,X_3,H_3)\to (G_2,X_2,H_2)$ of $f_1$ and $f_2$. We have 
$$G_3:=G\times_{G_1} G_2.$$
\noindent 
The kernel of $f_{3}:G_3\to G$ is the same as of $f_2$ and so $f_3$ is a cover. Moreover, $f_{32}$ is an injective map and $G_3^{\rm der}$ is a simply connected semisimple group. From the assumptions on $\Mn_2^\prime$ and [Va2, 6.2.3] we get that $(G_2,X_2,H_2)$ has an integral canonical model $\Mn_2$ having the EP. From the proof of loc. cit. and [Va2, 3.2.14-15] we get that the connected components of $\Mn_2$ and $\Mn_2^\prime$ are isomorphic over an \'etale cover of the normalization of $\ZZ_{(p)}$ in the composite field of $E(G_2,X_2)$ and $E(G_2^\prime,X_2^\prime)$. So $\Mn_{2W(\FF)}$ also has the EEP and is quasi-projective. From [Va2, 6.2.2 a)] we get that $(G_1,X_1,H_1)$ has an integral canonical model $\Mn_1$ and the natural morphism (see [Va2, 3.2.7 4)]) $q_2:\Mn_2\to\Mn_1$ is a pro-\'etale cover. 
\smallskip
We consider similarly the natural morphism $q:\Mn\to\Mn_1$. Let $q_{32}:\Mn_3^\prime\to\Mn_2$ be its pull back via $q_2$. The resulting morphism $q_3:\Mn_3^\prime\to\Mn$ is a pro-\'etale cover, and so $\Mn_3^\prime$ is a regular, formally smooth $\ZZ_{(p)}$-scheme. From [Va2, 3.2.14-15] we get that the generic fibre of $q_{32}$ is a closed embedding and that $q_{32}$ itself is pro-finite. $\Mn_{3W(\FF)}^\prime$ has the EEP, cf. [Va2, 3.2.3.1 5)]. Moreover, from [Va2, 3.2.14] we get that ${\rm Sh}_{H_3}(G_3,X_3)$ is naturally an open closed subscheme of the generic fibre of $\Mn_3^\prime$. So the Zariski closure $\Mn_3$ of ${\rm Sh}_{H_3}(G_3,X_3)$ in $\Mn_3^\prime$ is a pro-\'etale cover of $\Mn$ (as $f_3$ is a cover, the morphism
${\rm Sh}_{H_3}(G_3,X_3)\to {\rm Sh}_{H}(G,X)$ is --see [Mi1, 4.11-13]-- surjective). Again from [Va2, 3.2.3.1 5)] we get that $\Mn_{3W(\FF)}$ has the EEP and that $\Mn_3$ has the EP. So, from [Va2, 3.2.3.1 1)] we get that $\Mn_3$ has also the SEP; this implies the existence of a natural continuous right $G_3(\AA_f^p)$-action on $\Mn_3$ extending the canonical one on its generic fibre. We conclude: $\Mn_3$ is the integral canonical model of $(G_3,X_3,H_3)$.  
\smallskip
From [Va2, 6.2.2-3] we get that any triple having the same adjoint as $(G_2,X_3,H_3)$ (and so as $(G,X,H)$) has an integral canonical model. This proves Lemma 3.    
\medskip
{\bf Step 2.} We come back to our context of the Exercise. Based on the replacement part of Step 1 and on [Va10, 3.0-1 and the paragraph after 3.1], we deduce that we can assume we have an embedding
$f:(G,X)\hookrightarrow (GSp(W,\psi),S)$ such that:
\medskip
{\bf a)} the adjoint pair of $(G,X)$ is $(G^{\rm ad},X^{\rm ad})$;
\smallskip
{\bf b)} $f$ factors through an injective map $f_1:(G_1,X_1)\hookrightarrow (GSp(W,\psi),S)$, which is a PEL type embedding (see 2.1);
\smallskip
{\bf c)} there is a $\ZZ_{(p)}$-lattice of $W$ such that for any prime $v_1$ of $E(G_1,X_1)$ dividing $p$, we have a standard PEL situation $(f_1,L_{(p)},v_1,\Mb_1)$, for a suitable $\ZZ_{(p)}$-subalgebra $\Mb_1$ of endomorphisms of $L_{(p)}$;
\smallskip
{\bf d)} $G_1^{\rm ad}$ is the Weil restriction from $F$ to $\QQ$ of an absolutely simple $F$-group of $A_{2l}$ Lie type;
\smallskip
{\bf e)} the Zariski closure of $G$ in $GL(L_{(p)})$ is a reductive group.
\medskip
We recall that the construction of [Va10, 3.0-1 and the paragraph after 3.1] is based on four things:
\medskip
{\bf f)} standard Galois cohomology shows that by enlarging $F$ (to $F_1$) we can assume $\tilde G$ is obtained by pull back from a particular (see g) below) adjoint group $\bar G$ over ${\rm Spec}(\QQ)$;
\smallskip
{\bf g)} the part of [He, p. 445] referring to $SO^*(2n)$ (with $n:=l$ to fit the notations of [Va2, 6.8.6]) has a natural $\ZZ_{(p)}$-version and so we can take $\bar G$ to be the adjoint group of the generic fibre of the $\ZZ_{(p)}$-version $SO^*(2n)_{\ZZ_{(p)}}$ of $SO^*(2n)$;
\smallskip
{\bf h)} the $\ZZ_{(p)}$-version of [Sa, 3.3]: we have as well a natural $\ZZ_{(p)}$-version $SU(n,n)_{\ZZ_{(p)}}$ of $SU(n,n)_{\RR}$ and a natural monomorphism $SO^*(2n)_{\ZZ_{(p)}}\hookrightarrow SU(n,n)_{\ZZ_{(p)}}$ which is a $\ZZ_{(p)}$-version of the standard one over $\RR$ (described at the level of Lie algebras in loc. cit.);
\smallskip
{\bf i)} 2.3.5.1.
\medskip  
Coming back to properties a) to e), from [De2, 2.3.13], as $G^{\rm ad}_{\RR}$ has no compact factors, we deduce that we can assume $G^{\rm ab}=\QQ$. But, from the explicit construction of [Va10, 3.0 D)] and from 2.3.5.1, as in [Va2, end of 6.6.5.1] (see also the first way above), we can enlarge the center of $G$ such that we get a standard PEL situation $(f_0,L_{(p)},v,\Mb_0)$, with $f$ factoring through $f_0:(G_0,X_0)\hookrightarrow (GSp(W,\psi),S)$ in such a way that we have $(G_0^{\rm ad},X_0^{\rm ad})=(G^{\rm ad},X^{\rm ad})$ (we have $E(G_0,X_0)=E(G,X)=E(G^{\rm ad},X^{\rm ad})$). 
This ends the second solution of the Exercise.   
\medskip
{\bf B.} The types of Shimura varieties of preabelian type which are neither of abelian type nor of compact type and were postponed in [Va2] (cf. [Va2, 6.8.6] and AE.0), are among the ones treated in A above. This (cf. also [Va2, 1.4]) proves the Second Main Corollary of 1.14.3. 
\medskip
{\bf C. Remarks.} It is worth recalling (cf. [Va2, 6.2.7]) that any naive attempt of handling B just at the level of elementary arguments involving DVR's, is meaningless. Also, one could try to get a proof of Lemma 2 of A or to solve 2.3.5.1, starting from the classification of [Sh, \S 5]: in our opinion such an approach would be much more complicated (if successful at all).
\medskip
{\bf 4.6.5. Exercise.} Prove d) of 4.4.1 3) for the case $k(v)=\FF_p$. Hint: use P1, P2 and Fontaine's comparison theory. 
\medskip
{\bf 4.6.6. Remark.} Using the same ideas as in b) of 2.3.18 B3, d) of 4.4.1 3) follows immediately for the case of a standard PEL situation. We have no problems with the $D$ cases, cf. 4.2.3-4 and the following simple fact: referring to the standard PEL situation $(f,L_{(p)},v,\Mb)$ of 2.3.5.1, each connected component of the extension to $B(\FF)$ of the intersection of $GSp(W,\psi)$ with the centralizer of $\Mb$ in $GL(W)$, has a $K_0$-valued point which extends to a $V_0$-valued point of $GSp(L_{(p)}\otimes V_0,\psi)$.
\smallskip
Using this, one can obtain a much simpler and shorter proof of 4.2.1 a) for standard PEL situations: a logical analogue of 3.1.8.1, 3.6.6.0, and [dJ1] or [Fa2, th. 10] is all that is needed to complete a proof of 4.2.1 a) for standard PEL situations in a very small number of pages.
\medskip
{\bf 4.6.7. A digression.} We consider a Shimura group pair $(G_0,[\mu_0])$ over $\ZZ_p$, with $G_0$ a simple $\ZZ_p$-group. Let $q\in\NN$ be the smallest number for which there is a cocharacter $\mu_1:\GG_m\to G_{0W(\FF_{p^q})}$ such that $(G_0,[\mu_0])=(G_0,[\mu_1])$ and there is a Borel subgroup $B_0$ of $G_0$ with the property that ${\rm Lie}(B_{0W(\FF_{p^q})})$ is normalized by the image of $\mu_1$ and $\GG_m$ acts on it through $\mu_1$ via the trivial and the inverse of the identical cocharacters (cf. [Mi3, 4.6-7] and 2.2.6). We assume that the Shimura Lie $\sg_{\FF_{p^q}}$-crystal $({\rm Lie}(G_{0W(\FF_{p^q})}),\sg_{\FF_{p^q}}\mu_1({1\over p}))$ defined as in 4.1.1.2 is of $A_n$, $B_n$, $C_n$, $D_n^{\RR}$ or $D_n^{\HH}$ type (see 3.10.5). We have:
\medskip
{\bf Lemma.} {\it {\bf 1)} There is a Shimura quadruple $(G_1,X_1,H_1,v_1)$ of adjoint, abelian type such that its attached Shimura group pair (see end of 4.1) is (isomorphic to) $(G_0,[\mu_1])$.
\smallskip
{\bf 2)} There is a Shimura quadruple $(G_1,X_1,H_1,v_1)$ of adjoint, compact, abelian type such that its attached Shimura group pair has $(G_0,[\mu_1])$ as a product factor.}
\medskip
{\bf Proof:} Let $m$ be the number of simple factors of $G_{0W(\FF)}$. So $G_0$ is the Weil restriction from $W(\FF_{p^m})$ to $\ZZ_p$ of an absolutely simple $W(\FF_{p^m})$-group $G_0^1$. We start (cf. approximation theory involving unequivalent valuations) with a totally real number field $F_1$ such that $[F_1:\QQ]=m$ and there is a prime $v$ of it dividing $p$ and with $k(v)=\FF_{p^m}$; so $v$ is unramified above $p$. 
Let $F_2$ be a totally real quadratic extension of $\QQ$ in which $p$ splits. Let $F$ be the composite field of $F_1$ and $F_2$. We have $[F:\QQ]=2m$. For 1) (resp. for 2)) we take $G_1$ to be the Weil restriction from $F_1$ to $\QQ$ (resp. from $F$ to $\QQ$) of (see below) a suitable absolutely simple $F_1$-group (resp. $F$-group) $G_1^1$. 
\smallskip
Let $F_1^\prime$ be the smallest Galois extension of $\QQ$ containing $F_1$. We identify canonically the quotient set ${\rm Gal}(F_1^\prime/\QQ)/{\rm Gal}(F_1^\prime/F_1)$ with ${\rm Gal}(k(v)/\FF_p)$. So we also identify the set 
$$RE_1$$ 
of embeddings of $F_1$ into $\RR$ with the set of embeddings of $k(v)$ into $\FF$. Let $U$ be the set of primes of $F_1$ which are either archimedean or are $v$. For $w\in U$, we denote by $F_{1w}$ the completion of $F_1$ w.r.t. $w$. 
\smallskip
Let $\tilde G_1$ be the split, simple, adjoint group over $F_1$ of the same Lie type as $G_0$. Let $Aut_0(\tilde G_1)$ be the connected component  of the origin of $Aut(\tilde G_1)$. Let $Aut_1(\tilde G_1)$ be the subgroup of $Aut(\tilde G_1)$ defined as follows. If $G_0$ is not of $D_4$ Lie type, then $Aut_1(\tilde G_1)$ is $Aut(\tilde G_1)$ itself. If $G_0$ is of $D_4$ Lie type (so $Aut(\tilde G_1)$ is the semidirect product of $Aut_0(\tilde G_1)$ and of $S_3$), we take $Aut_1(\tilde G_1)$ to be the semidirect product of $Aut_0(\tilde G_1)$ by a $\mu_2$ subgroup of the symmetric group $S_3$ such that the element of $H^1({\rm Gal}(F_{1v}),Aut(\tilde G_1)_{F_{1v}}(\overline{F_{1w}}))$ defining $G_{0F_{1v}}^1$ is defined by an element of the image of $H^1({\rm Gal}(F_{1v}),Aut_1(\tilde G_1)_{F_{1v}}(\overline{F_{1w}}))$ in it (cf. Fact 2 of 4.3.6). Let 
$$H^1:=H^1({\rm Gal}(F_1),Aut_1(\tilde G_1)(\overline{\QQ})).$$ 
\indent
We first deal with 1). For each $\tau\in RE_1$, (in order to have 1)) the pull back $G_1^1(\tau)$ of $G_1^1$ to $\RR$ via it is uniquely determined by the pair 
$$(G^1_{0W(\FF)},[\mu_1(\tau)])$$ 
we naturally get from $(G_0,[\mu_1])$ via the embedding $W(k(v))\hookrightarrow W(\FF)$ corresponding to $\tau$. More precisely, $G_1^1(\tau)$ is compact if $[\mu_1(\tau)]$ is trivial, it is a split, simple, adjoint group of $C_n$ Lie type if $G_0$ is of $C_n$ Lie type and $[\mu_1(\tau)]$ is non-trivial, it is $SO(2,2n-1)_{\RR}$ if $G_0$ is of $B_n$ Lie type and $[\mu_1(\tau)]$ is non-trivial, etc. (see [He, p. 518]). 
\smallskip
It is known that the map
$$m(Aut_1(\tilde G_1),U):H^1\to \oplus_{w\in U}H^1({\rm Gal}(F_{1w}),Aut_1(\tilde G_1)_{F_{1w}}(\overline{F_{1w}}))$$
is surjective. 
If $Aut_1(\tilde G_1)$ is connected, then the surjectivity of $m(Aut_1(\tilde G_1),U)$ is a consequence of [Mi1, B.3 and B.24]. If $Aut_1(\tilde G_1)$ is the semidirect product of $Aut_0(\tilde G_1)$ and of $\mu_2$, then the surjectivity of $m(Aut_1(\tilde G_1),U)$ is a consequence of loc. cit. (applied to $Aut_0(\tilde G_1)$) and of the structure of the Brauer group over $F_1$ (i.e. of the fact that the similarly defined map $m(\mu_2,U)$ is surjective).
\smallskip
We deduce the existence of an absolutely simple $F_1$-group $G_1^1$ whose pull back to $F_{1v}=B(k(v))$ is $G_{0B(k(v))}^1$ and whose pull back to $\RR$ via $\tau$ is $G_1^1(\tau)$ (if $G_0$ is of $D_4$ Lie type we need to add that $S_3$ is the semidirect product of $\mu_2$ by a normal subgroup of it of order relatively prime to $[\CC:\RR]$; so the natural map $H^1({\rm Gal}(\CC/\RR),Aut_1(\tilde G_1)(\CC))\to H^1({\rm Gal}(\CC/\RR),Aut(\tilde G_1)(\CC))$ is surjective). Let $G_1:={\rm Res}_{F_1/\QQ} G_1^1$. 
\smallskip
If $({\rm Lie}(G_{0W(\FF_{p^q})}),\sg_{\FF_{p^q}}\mu_1({1\over p}))$ is of $C_n$, $B_n$ or $D_{n+4}^{\RR}$ type, $n\in\NN$, then there is a unique way to get a $G_1(\RR)$-conjugacy class $X_1$ of homomorphisms ${\rm Res}_{\CC/\RR} \GG_m\to G_{1\RR}$ so that $(G_1,X_1)$ is a Shimura pair. In general, for each $\tau\in RE_1$ there are at most two Shimura group pairs of the form $(G_{1\RR},[\mu_1^{\tau}])$ (here $G_{1\RR}$ is obtained from $G_1$ via $\tau$), cf. [De2, 1.2.7-8]. So we first choose 
$$X_1=\prod_{\tau\in RE_1} X_1^{\tau}$$ 
arbitrarily such that $(G_1,X_1)$ is a Shimura pair; here $X_1^{\tau}$ is the $G_{1\RR}(\RR)$-conjugacy class of homomorphisms ${\rm Res}_{\CC/\RR} \GG_m\to G_{1\RR}$ (defined via $\tau$ and) naturally associated to $[\mu_1^{\tau}]$. For each $\tau\in RE_1$ such that the pair $(G^1_{0W(\FF)},[\mu_{11}(\tau)])$ we naturally get from the triple $(G_1,X_1,\tau)$ (here we use the identification of $RE_1$ with ${\rm Gal}(k(v)/\FF_p)$) is $(G^1_{0W(\FF)},[(\mu_{1}(\tau)^{-1}])$, we replace $[\mu_1^{\tau}]$ by $[(\mu_1^{\tau})^{-1}]$. After all such replacements we get that 1) holds for the prime $v_1$ of $E(G_1,X_1)$ dividing $p$ and constructed as follows. We consider an embedding $B(\FF_{p^q})\hookrightarrow\CC$ and (via it) we view $E(G_1,X_1)$ naturally as a subfield of $B(\FF_{p^q})$; we take $v_1$ such that $W(\FF_{p^q})$ is a faithfully flat $O_{(v_1)}$-algebra. This takes care of 1) (we can take as $H_1$ the hyperspecial subgroup $G^1_0(W(\FF_{p^m}))$ of $G_1(\QQ_p)=G_1^1(F_1\otimes_{\QQ} \QQ_p)=G_1^1(F_{1v})=G^1_{0B(\FF_{p^m})}$). 
\smallskip
For 2) we proceed similarly: just one modification needs to be made. We fix an embedding $\tau_2$ of $F_2$ into $\RR$. Let $RE$ be the set of embeddings of $F$ into $\RR$. We identify $RE_1$ with the subset of $RE$ corresponding to embeddings whose restriction to $F_2$ is $\tau_2$. For $\tau\in RE\setminus RE_1$ (resp. for $\tau\in RE_1$) we take $G_1^1(\tau)$ to be compact (resp. we take it to be defined as above). $F\otimes_{\QQ} \QQ_p$ is identified naturally with two copies of $F_1\otimes_{\QQ} \QQ_p$: each such copy corresponds to primes of $F$ dividing $p$ and a fixed prime of $F_2$ dividing $p$; identifying the two primes of $F_2$ dividing $p$ with the two embeddings of $F_2$ in $\QQ_p$ and with the two emebeddings $\tau_2$ and $\overline{\tau_2}$ of $F_2$ and $\RR$, the mentioned two coppies of $F_1\otimes_{\QQ} \QQ_p$ corresspond to the two disjoint subsets $RE\setminus RE_1$ and $RE_1$ of $RE$. We get that the Shimura group pair attached to a Shimura quadruple $(G_1,X_1,H_1,v_1)$ constructed as above is the product of $(G_0,[\mu_1])$ with the pair $(G_0,[\mu_{{\rm triv}}])$ defined by the cocharacter $\mu_{{\rm triv}}:\GG_m\to G_0$ who has a trivial image. This ends the proof.
\medskip
{\bf Corollary 1.} {\it The extension of $({\rm Lie}(G_{0W(\FF_{p^q})}),\sg_{\FF_{p^q}}\mu_1({1\over p}))$ to $\FF$ is isomorphic to the Shimura adjoint Lie $\bar\sg$-crystal attached to a principally quasi-polarized Shimura $\bar\sg$-crystal $(N,\vph_N,G_N,p_N)$. More precisely, $(N,\vph_N,G_N,p_N)$ is attached to an $\FF$-valued Shimura-ordinary point of $\Mn$, where $\Mn$ is the integral canonical model of a particular SHS.}
\medskip
{\bf Proof:} This is a consequence of 1) and of the Existence Property of 1.10. For $p=3$, cf. also \S 6. 
\medskip
{\bf Corollary 2.} {\it All group actions $\TT^0$ of 3.13.7.1 have $\dim_k(P_k^0)-\dim_k(N_k^-)=\dim_k(P_k^0)-dd((M,\vph,G))$ dimensional orbits.}
\medskip
{\bf Proof:} Based on 2) and on Proposition of 4.5.15.2.1, this follows from the following Fact.
\medskip
{\bf Fact.} {\it For any SHS $(f,L_{(p)},v)$, with $(G,X)$ of compact type, the quasi-ultra stratification of $\Mn_{k(v)}$ has 0 dimensional strata.}
\medskip
As $\Mn/H_0$ is a projective scheme (see 2.3.3.1), the Fact follows from the quasi-affineness part of 4.5.15.2.5. From Lemma 2 of 4.6.4 A we get: 
\medskip
{\bf Corollary 3.} {\it We assume the Shimura Lie $\sg_{\FF_{p^q}}$-crystal $({\rm Lie}(G_{0W(\FF_{p^q})}),\sg_{\FF_{p^q}}\mu_1({1\over p}))$ is either of $A_n$ Lie type or is of totally non-compact $C_n$ or $D_n^{\HH}$ type. Then in the second part of Corollary 1 we can assume we have a standard PEL situation.}   
\medskip
{\bf 4.6.7.1. Back to 4.6.1 1).} We come back to 4.6.1 1). Let $U$ be as in the proof of the Lemma of 4.6.7. Based on the mentioned proof, we deduce the existence of an absolutely simple $F$-group $\tilde G^0_F$ of $C_n$ Lie types which splits over $F_v$ and whose pull backs to $\RR$ (via different embeddings of $F$) are as desired in 4.6.1 1). $\tilde G_F^0$ splits above all primes of $F$ but a finite set $FS$ not containing $p$. For $v_1\in SF$, $\tilde G_{F_{v_1}}^0$ splits over a quadratic extension $E_{v_1}$ of the completion $F_{v_1}$ of $F$ w.r.t. $v_1$. We now consider a totally imaginary quadratic extension $E$ of $F$ such that no element of $FS$ splits in it and moreover, for all $v_1$, the completion of $E$ w.r.t. its prime dividing $v_1$ is $E_{v_1}$ (cf. Krasner's lemma and the approximation theory). We get: the pull back $\tilde G_E^0$ of $\tilde G_F^0$ to $E$ splits at all places of $E$. So $\tilde G_E^0$ splits. So we can take (cf. [Va2, 3.1.3]) as $\tilde G$ a simply connected group over $O_{(v)}=F_{(p)}$ having $\tilde G_F^0$ as the adjoint group of its generic fibre.
\medskip
{\bf 4.6.8. Remark.} 
Based on 3.13.7.4 B, to check that the finiteness property of 4.5.15.1 holds if all simple factors of $(G^{\rm ad},X^{\rm ad})$ are of some $A_n$, $C_n$ or $D_n^{\HH}$ types, we can assume (cf. Corollary 3 of 4.6.7) that in 4.5.15.1 we are dealing with a standard PEL situation; in such a context it is easy to check it starting from [Oo3]. We will come back to this in \S 9-10.
\smallskip
So to check that the finiteness property of 4.5.15.1 holds in general, we can assume (cf. the above paragraph, 3.13.7.4 B and Corollary 1 of 4.6.7) that $G_{\RR}^{\rm ad}$ has no compact factors and that (cf. also Fact 1 of 2.3.5.2) $(G^{\rm ad},X^{\rm ad})$ has all its simple factors of some $B_n$ or $D_n^{\RR}$ type and we can proceed in many ways (like: by induction on $n$, or by working in the context of symmetric objects of $\Mm\Mf_{[0,2]}(W(\FF))$, etc.). However, the (presently abstract) context of $E_6$, $E_7$ and $D_{\ell}^{\rm mixed}$ types (see 3.13.7.7) can not be that easily studied and so we defer to \S 10 for more of a principal approaches (then the one above relying on the use of relative PEL situations in a sense similar to the one of [Va2, 4.3.16]).   
\medskip\smallskip
{\bf 4.7. Crystalline coordinates.}
\medskip
{\bf 4.7.0. $G$-deformations.} We consider an arbitrary point $y:{\rm Spec}(k)\to\Mn/H_0$. Let $(A_y,p_{A_y}):=y^*(\Ma_{H_0},\Mp_{\Ma_{H_0}})$. By a $G$-deformation of $(A_y,p_{A_y})$ over an artinian local $W(k)$-algebra $AL$ having $k$ as its residue field, we mean a deformation which is obtained from $(\Ma,\Mp_{\Ma})$ by pull back via a morphism ${\rm Spec}(AL)\to\Mn/H_0$ lifting $y$.
\medskip
{\bf 4.7.1. First properties.} Let now $y:{\rm Spec}(k)\to\Mu/H_0$ be a $G$-ordinary point. 
For simplifying the notations we assume $y$ has attached to it a principally quasi-polarized Shimura $\sg_k$-crystal of the form $\bigl(M_y,\vph_y,G_{W(k)}$, $(t_\al)_{\al\in\Mj^\prime},\psi_y\bigr)$ (cf. 2.3.10). Let $\mu_y:\GG_m\to G_{W(k)}$ be its canonical split (cf. 3.1.6). It produces a direct sum decomposition $M_y=F^1_y\oplus F^0_y$, with $\be\in\GG_m(W(k))$ acting through
$\mu_y$ on $F^i_y$ as the multiplication with $\be^{-i}$. Let $N$ be the subgroup of $G^{\rm der}_{W(k)}$ which acts trivially on $F^0_y$ and on $M_y/F_y^0$. Let $R$ be the $W(k)$-algebra of the formal completion of $N$ at the origin. We choose an isomorphism $\tilde f:R\tilde\to W(k)[[x_1,\ldots,x_{e_0}]]$, with $e_0=\dim_{\CC}(X)$ (see [Va2, 5.4.7]). We are mainly interested in the case when $k=\bar k$; but for a great part of what follows this is not needed. It is easy to see that the group $N$ is isomorphic to $\GG^{e_0}_a$ (i.e. it can be identified with the group scheme defined by its Lie algebra via the exponential map). We can choose the coordinates $x_1,\ldots,x_{e_0}$ (i.e. we can choose $\tilde f$) such that the $i$-th $\GG_a$ copy of $N$ is ${\rm Spec}(W(k)[x_i])$. So $N$ as a scheme is ${\rm Spec}(W(k)[x_1,\ldots,x_{e_0}])$. We consider the Frobenius lift of $N$ which at the $W(k)$-algebras level takes $x_i\to x_i^p$, $i=\overline{1,e_0}$. It induces by completion a Frobenius lift $\Phi_R$ of $R$.
\smallskip
As the case of a torus (of a $0$ dimensional Shimura variety of Hodge type) is trivial, we assume from now on that $G$ is not a torus.
Let $N_0$ be the smallest connected subgroup of $G_{W(k)}$ such that:
\medskip
\item{\bf i)} it contains $N$;
\smallskip
\item{\bf ii)} its Lie algebra is taken by $p\vph_y$ into itself. 
\medskip
For all that follows, we rely heavily on the terminology and formulas of 3.10.6. It is easy to see (to be compared also with 3.4.3.0 and 3.11.2 C) that:
\medskip
\item{\bf 0)} $N_0$ is a product indexed by the number of simple factors of $G_{\QQ_p}^{\rm ad}$, cf. 4.3.1.
\smallskip
\item{\bf 1)} $N_0$ is a unipotent, smooth subgroup of $G_{W(k)}$.
\smallskip
\item{\bf 2)} We get a Lie $\sg_k$-crystal $({\rm Lie}(N_0),\vph_y,F^0({\rm Lie}(N_0)),F^1({\rm Lie}(N_0)))$, with $F^i({\rm Lie}(N_0)):=F^i({\rm End}(M))\cap {\rm Lie}(N_0)$, $i=\overline{0,1}$. In fact $F^1({\rm Lie}(N_0))=\{0\}$, and ${\rm Lie}(N_0)$ is the Lie subalgebra of ${\rm Lie}(G_{W(k)})$ corresponding to negative slopes of $({\rm Lie}(G_{W(k)}),\vph_y)$ (this and 1) can be read out --over $\bar k$-- from 3.4.3.0, via 3.11.1 c); see also 3.11.2 C and 3.10.7).
\smallskip
\item{\bf 3)} $N=N_0$ iff the slopes of ${\rm Lie}_G(\tau)$ are precisely $-1$, $0$ and $1$, cf. 2) and 4.6 P5.
\smallskip
\item{\bf 4)} $N_0$ is an abelian group iff the spreading of any  cyclic adjoint factor attached to $y$ is 0 or 1, cf. def. in 3.10.4.4 and 4.3.1. This automatically happens if all simple factors of $(G^{\rm ad},X^{\rm ad})$ are of $B_l$, $C_l$ or $D_l^{\RR}$ type; but it is not necessarily true if there are simple factors of $(G^{\rm ad},X^{\rm ad})$ of $A_{\ell}$ or $D_{\ell}^{\HH}$ type ($\ell\in\NN$).
\smallskip
\item{\bf 5)} In general, $N_0$ is a nilpotent group of order of nilpotency not greater than some $m\in\NN$ iff the spreading of any cyclic adjoint factor attached to $y$ is not greater than $m$ (this and the first part of 4) are a consequence of the characteristic $0$ analogue of 3.5.4 (6), as it can be easily deduced from 3.10.7).
\smallskip
\item{\bf 6)} $N$ is a normal subgroup of $N_0$ (for instance, this follows from 3.1.4).
\medskip
So if $(G^{\rm ad},X^{\rm ad})$ has only simple factors of $D_{\ell}^{\HH}$ type ($\ell\in\NN$), then $N_0$ has the order of nilpotency at most 2. If $(G^{\rm ad},X^{\rm ad})$ is a simple, adjoint Shimura variety of $A_{\ell}$ type, then the order of nilpotency of $N_0$ is smaller or equal to $\ell$ and equality can hold (cf. 4.6 P7); in fact the order of nilpotency of $N_0$ can be any integer in $S(1,\ell)$, cf. 5).
\smallskip 
We also regard $y$ as a closed embedding $y:{\rm Spec}(k)\hookrightarrow\Mn_{W(k)}/H_0$. Let $\hat O_y$ be the completion of the local ring of $y$ in $\Mn_{W(k)}/H_0$.
Let $(A_{\hat y},p_{A_{\hat y}})$ be the principally polarized abelian scheme over ${\rm Spec}(\hat O_y)$ obtained from $(\Ma_{H_0},\Mp_{\Ma_{H_0}})$ through the canonical morphism ${\rm Spec}(\hat O_y)\to\Mn/H_0$. We consider the following principally quasi-polarized Shimura filtered $F$-crystal 
$$
_NM_{\hat y}:=(M_y,F^1_y,\vph_y,G_{W(k)},N,\tilde f,(t_\al)_{\al\in\Mj^\prime},\psi_y)
$$ 
over $R/pR$. There is a unique connection $_N\nabla_{\hat y}$ on $M_y\otimes_{W(k)} R$ which is integrable, nilpotent mod $p$ and such that the Frobenius endomorphism 
$$\Phi_y:=n(\vph_y\otimes 1)$$ 
of $M_y\otimes_{W(k)} R$ is $_N\nabla_{\hat y}$-parallel, cf. 2.2.10. Here $n\in N(R)$ is the universal element defined by the natural morphism ${\rm Spec}(R)\to N$; we view it as an $R$-linear endomorphism of $M_y\otimes_{W(k)} R$. 
\smallskip
Similarly we define $_{G^{\rm der}_{W(k)}}M_{\hat y}$,
starting from a Frobenius lift $\Phi_{R_1}$ of the $W(k)$-algebra $R_1$ of the completion of $G^{\rm der}_{W(k)}$ at the origin, which takes the kernel of the canonical $W(k)$-epimorphism $R_1\twoheadrightarrow R$ into itself, inducing $\Phi_R$ on $R$.
Let 
$$\tilde e_0:=\dim_{W(k)}(N_0).$$ 
As a scheme $N_0={\rm Spec}(W(k)[x_1,...,x_{\tilde e_0}])$, while as a group scheme $N_0$ is obtained from $\GG_a$ groups by a repeated process of forming short exact extensions defining semidirect products (cf. 1), or 5), or [SGA3, Vol. III, 4.1.2 of p. 172]). In other words, such an identification can be chosen so that each closed subscheme ${\rm Spec}(W(k)[x_i])$ of $N_0$, $i\in S(1,\tilde e_0)$, is a $\GG_a$ subgroup of $N_0$, to be referred as  the $i$-th $\GG_a$ copy of $N_0$. Let ${\rm Spec}(R_0)$ be the completion of $N_0$ in its origin. So $R_0=W(k)[[x_1,...,x_{\tilde e_0}]]$, with ${\rm Spec}(W(k)[[x_i]])$ as the completion of the $i$-th $\GG_a$ copy of $N_0$ in its origin, $i=\overline{1,\tilde e_0}$. We choose such an identification, so that (cf. also 6)) the natural epimorphism $R_0\twoheadrightarrow R$ is defined by $x_i$ goes to $x_i$, if $i\in S(1,e_0)$, and into $0$, if $i\in S(1+e_0,\tilde e_0)$. We assume the Frobenius lift $\Phi_{R_1}$ of $R_1$ is such that it also takes the kernel of the natural epimorphism $R_1\twoheadrightarrow R_0$ into itself, inducing a Frobenius lift of $R_0$ which takes $x_i$ into $x_i^p$, $i=\overline{1,\tilde e_0}$. 
\medskip
{\bf 4.7.2. Lemma.} {\it The principally quasi-polarized filtered $F$-crystal over ${\rm Spec}(\hat O_y/p\hat O_y)$ endowed with a family of crystalline tensors, defined by $(A_{\hat y},p_{A_{\hat y}})$ and de Rham components of Hodge cycles with which $A_{\hat y}$ is naturally endowed, is isomorphic to $_NM_{\hat y}$.}
\medskip
{\bf Proof:} This is a consequence of Fact 4 of 2.3.11 (via 2.2.21). For the sake of convenience we recall the details. As in 2.3.11, we start with a $W(k)$-morphism ${\rm Spec}(R_1)\buildrel{a_1}\over\lra {\rm Spec}(\hat O_y)$ lifting the natural morphism ${\rm Spec}(k)\to {\rm Spec}(\hat O_y)$ defined by $y$ and such that the principally quasi-polarized filtered  $F$-crystal over $R_1/pR_1$ endowed with tensors and associated to the pull back through $a_1$ of $(A_{\hat y},p_{A_{\hat y}})$ and  de Rham components of the mentioned Hodge cycles of $A_{\hat y}$, is isomorphic to $_{G^{\rm der}_{W(k)}}M_{\hat y}$. But the pull back of $_{G^{\rm der}_{W(k)}}M_{\hat y}$ through the canonical $W(k)$-monomorphism ${\rm Spec}(R)\buildrel{b_1}\over\hookrightarrow{\rm Spec}(R_1)$ is $_NM_{\hat y}$ (due to the compatibility of the Frobenius lifts of $R_1$ and $R$ with the $W(k)$-epimorphism  $R_1\twoheadrightarrow R$). But the $W(k)$-morphism $c_1:=b_1\circ a_1:{\rm Spec}(R)\to{\rm Spec}(\hat O_y)$  is
an isomorphism (this is the same as [Va2, 5.4.7-8]). This proves the Lemma.
\smallskip
We denote by $\bar\Om_{R/W(k)}$ the free $R$-module generated by $dx_i$, $i=\overline{1,e_0}$. $\bar\Om_{R/W(k)}$ is naturally a direct summand of the $R$-module $\Om_{R/W(k)}$. Let $\dl_N$ be the connection on $M_y\otimes_{W(k)} R$ annihilating $M_y$.
\medskip
{\bf 4.7.3. Lemma.} {\it The connection $_N\nabla_{\hat y}$ is of the form $\dl_N+\be_N$, where 
$$\be_N\in {\rm Lie}(N_0)\otimes_{W(k)} \Om_{N/W(k)}\otimes_{W(k)} R={\rm Lie}(N_0)\otimes_{W(k)} \bar\Om_{R/W(k)}.$$}
\indent
{\bf Proof:} Similarly to $_NM_{\hat y}$ we define $_{N_0}M_{\hat y}$. $_NM_{\hat y}$ is obtained from $_{N_0}M_{\hat y}$ by extension through the $W(k)$-epimorphism $R_0\twoheadrightarrow R$. As this $W(k)$-epimorphism is compatible with the chosen Frobenius lifts, the Lemma is a consequence of [Fa2, rm. ii) after th. 10] applied to $(M_y,F^1_y,\vph_y)$ and to the subgroup $N_0$ of $GL(M_y)$ (cf. also the formulas of 3.6.18.4 B) and of 3.6.18.7.1).
\medskip
{\bf 4.7.3.1. Notations.} We write 
$$\be_N=\sum_{j=1}^{\tilde e_0}n_j\otimes\vep_j, 
$$
with $\{n_j\}$ a $W(k)$-basis of the Lie algebra of the $j$-th $\GG_a$ copy of $N_0$ and with $\vep_j\in\bar\Om_{R/W(k)}$
(cf. 4.7.3). The elements $n_j$, $j\in S(1,\tilde e_0)$, are viewed as endomorphisms of $M_y\otimes_{W(k)} R$.
\medskip
{\bf 4.7.4. Levels of complexity.} Let $v^{\rm ad}$ be the prime of $E(G^{\rm ad},X^{\rm ad})$ divided by $v$. To construct (good) $G$-multiplicative coordinates and then to extend them to the general case of 4.7.14.1 below, we distinguish five different cases (levels of complexity). They are:
\medskip
\item{\bf A.} $N=N_0$.
\smallskip
\item{\bf B.} $N_0$ is an abelian group different from $N$.
\smallskip
\item{\bf C.} $N_0$ is a nilpotent group of order of nilpotency two. 
\smallskip
\item{\bf D.} $N_0$ is a nilpotent group of order of nilpotency greater than two.
\smallskip
\item{\bf E.} The general case referred to in 4.7.14.1 below.
\medskip
The condition $N=N_0$ is equivalent to $k(v^{\rm ad})=\FF_p$, cf. 4.7.1 2) and  4.6 P5.
First we describe the case A (so $\tilde e_0=e_0$); warning: 7.5-10 below refer just to this case. A quite general theory of crystalline coordinates is elaborated in 4.7.11. In 4.7.12-18 we deal with different aspects of it (in particular, we handle to a great extend cases B to E). 
\smallskip
From now on we assume $k=\bar k$.
\medskip
{\bf 4.7.5. Theorem.} {\it We assume $N=N_0$. We have:
\smallskip
{\bf i)} There is a $W(k)$-basis $\{a_i|1\le i\le e\}$ of $F^0_y$ and a $W(k)$-basis $\{b_i|1\le i\le e\}$ of $F^1_y$ and a permutation $\pi$ of $\{a_1,\ldots,a_{e},b_1,\ldots,b_{e}\}$ satisfying the conditions:
\medskip
\item{\bf (4.7.5.1)} $_N\nabla_{\hat y}(a_i)=0$, $\forall i\in S(1,e)$;
\medskip
\item{\bf (4.7.5.2)} $_N\nabla_{\hat y}(b_i)=\sum^{e_0}_{j=1}$ $n_j(b_i)\otimes\vep_j$, $\forall i\in S(1,e)$;
\medskip
\item{\bf (4.7.5.3)} $\Phi_y(a_i)=\pi(a_i)$ if $\pi(a_i)\in\{a_1,\ldots,a_e\}$;
\medskip
\item{\bf (4.7.5.4)} $\Phi_y(a_i)=n(\pi(a_i))$ if $\pi(a_i)\in\{b_1,\ldots,b_e\}$;
\medskip
\item{\bf (4.7.5.5)} $\Phi_y(b_i)=p\pi(b_i)$ if $\pi(b_i)\in\{a_1,\ldots,a_e\}$;
\medskip
\item{\bf (4.7.5.6)} $\Phi_y(b_i)=pn(\pi(b_i))$ if $\pi(b_i)\in\{b_1,\ldots,b_e\}$.
\medskip
Furthermore $\vep_i$'s are closed 1-forms verifying the following equations:
\medskip
\item{\bf (4.7.5.7)} If $i\in S(1,e)$ is such that $\pi(b_i)\in\{a_1,\ldots,a_e\}$,
then
$$0=\sum^{e_0}_{j=1}\vph(n_j)\bigl(\pi(b_i)\bigr)\otimes d\Phi_{R*}(\vep_j)\,;$$
\item{\bf (4.7.5.8)} If $i\in S(1,e)$ is such that $\pi(b_i)\in\{b_1,\ldots,b_e\}$, then 
$$dn\bigl(\pi(b_i)\bigr)+\sum_{j=1}^{e_0} \,n_j\bigl(\pi(b_i)\bigr)\otimes \vep_j=\sum^{e_0}_{j=1}\vph(\,n_j)\bigl(\pi(b_i)\bigr)\otimes d\Phi_{R*}(\vep_j);$$
\item{\bf (4.7.5.9)} If $i\in S(1,e)$ is such that $\pi(a_i)\in\{b_1,\ldots,b_e\}$, then 
$$dn\,\bigl(\pi(a_i)\bigr)+\sum^{e_0}_{j=1}n\bigl(\,n_j\bigl(\pi(a_i)\bigr)\bigr)\otimes\vep_j=0.$$
\smallskip
{\bf ii)} There are formal power series $\mu_j\in B(k)[[x_1,\ldots,x_{e_0}]]$, $j=\overline{1,e_0}$, such that:
\medskip
\item{\bf (4.7.5.10)} $\vep_j=d\mu_j$, $\forall j\in S(1,e)$;
\medskip
\item{\bf (4.7.5.11)} $\forall i\in S(1,e)$ such that $\pi(b_i)\in\{a_1,\ldots,a_e\}$ we have
$$\bigl(\sum^{e_0}_{j=1}\vph(n_j)\bigl(\pi(b_i)\bigr) d\Phi_{R*}(\mu_j)\bigr)(0)=0;$$
\item{\bf (4.7.5.12)} $\forall i\in S(1,e)$ such that $\pi(b_i)\in\{b_1,\ldots,b_e\}$ we have
$$\bigl((n-1_{M_y\otimes_{W(k)} R})(\pi(b_i))+p\sum_{j=1}^{e_0} n_j\bigl(\pi(b_i)\bigr)\mu_j\bigr)(0)=
\sum^{e_0}_{j=1}\bigl(\vph(p\,n_j)\bigl(\pi(b_i)\bigr) d\Phi_{R*}(\mu_j)\bigr)(0);$$
\item{\bf (4.7.5.13)} $\forall i\in S(1,e)$ such that $\pi(a_i)\in\{b_1,\ldots,b_e\}$ we have 
$$\bigl((\log\,n)(\pi(a_i))+\sum_{j=1}^{e_0} n_j(\pi(a_i))\mu_j\bigr)(0)=0.$$
\smallskip
Moreover $\mu_j(0,...,0)=0$ and the series $q_j=\exp(\mu_j)$ belongs to $R$ and satisfies $q_j(0,...,0)=1$, $\forall j\in S(1,\tilde e_0)$.}
\medskip
{\bf 4.7.6. Corollary.} {\it The elements $q_j-1\in R$, $j=\overline{1,e_0}$, form together with $p$ a regular system of parameters of $R$, i.e. the $W(k)$-homomorphism $R\to R$ sending $x_j$ to $q_j-1$, $j=\overline{1,e_0}$, is an isomorphism. Moreover, if we take the Frobenius lift of $R$ to be defined by $q_j\to q_j^p$, $j=\overline{1,e_0}$, then the Frobenius endomorphism of the underlying $R$-module $M\otimes_{W(k)} R$ of $_NM_{\hat y}$ takes $b_j$ into $p\pi(b_j)$ and takes $a_j$ into $\pi(a_j)$, $\forall j\in S(1,e)$.}
\medskip
{\bf 4.7.7. Definition.} The elements $q_j-1\in R$ are called $G$-canonical or Shimura-canonical coordinates of ${\rm Spec}(R)$ or of ${\rm Spf}(R)$.
\medskip
{\bf 4.7.8. Corollary.} {\it The moduli formal scheme $\Mm_G(y)$ of $G$-deformations of the principally polarized abelian variety $(A_y,p_{A_y})$, over (spectra of) artinian local $W(k)$-algebras having $k$ as their residue field, is canonically isomorphic to the formal torus $\Mt_{e_0}$ of the completion of $\GG_m^{e_0}$ in its origin; the principally polarized abelian variety over (the spectrum of) $W(k)$, obtained from $(\Ma_{H_0},\Mp_{\Ma_{H_0}})$ by pull back through the $G$-canonical lift of $y$, corresponds to the unit element of $\Mt_{e_0}(W(k))$.}
\medskip
{\bf 4.7.9. Corollary.} {\it We can assume that the permutation $\pi$ is trivial iff $k(v)=\FF_p$, i.e. iff $A_y$ is an ordinary abelian variety. If this is the case then the formal torus of $G$-deformations of $(A_y,p_{A_y})$ is a formal subtorus of the formal torus of deformations of $A_y$ or of $GSp$-deformations of $(A_y,p_{A_y})$
(i.e. of deformations of the principally polarized abelian variety $(A_y,p_{A_y})$).}
\medskip
{\bf 4.7.10. About the proofs of 4.7.5-9.}
The proofs of 4.7.5-6 are entirely analogous to the proofs of [De3, 1.4.2 and 1.4.7]. As this is the only part of our work which is similar to the classical theory of ordinary points of special fibres of integral canonical models of Siegel modular varieties, here
we just indicate which are the differences occurring. We do not bother to mention what happens when we choose different Frobenius lifts of $R$ (see [De3]).
\medskip
The first part of b) of 4.4.1 3) guarantees that the filtered $\sg_k$-crystal $(M_y,F^1_y,\vph_y)$ is cyclic diagonalizable. We deduce the existence of a $W(k)$-basis $\{a_i|1\le i\le e\}$ of $F^0_y$, of a $W(k)$-basis $\{b_i|1\le i\le e\}$ of $F^1_y$, and of a permutation $\pi$ of $BS:=\{a_1,\ldots,a_{e},b_1,\ldots,b_{e}\}$ such that $\vph_y(b_i)=p\pi(b_i)$ and $\vph_y(a_i)=\pi(a_i)$, $\forall i=\overline{1,e}$. We always make choices such that if a subset of $\{b_1,...,b_e\}$ or of $\{a_1,...,a_e\}$ is permuted cyclically by $\pi$, then this subset has only 1 element.  
\smallskip
The differences from [De3, 1.4] come from the fact that $\pi$ is not always the trivial permutation $1_{BS}$ of $BS$.
\smallskip
(4.7.5.1-2) just say what is the action of
$_N\nabla_{\hat y}$ on the $W(k)$-basis $BS$ of $M_y$, using the expression of $_N\nabla_{\hat y}$
given in 4.7.3.1. (4.7.5.3-6) just make explicit the action of $\Phi_y$ on elements of $BS$, starting from the above expression of $\vph_y$. The proof of the fact that $\vep_j$ are closed $1$-forms is as in [De3, 1.4.2], starting from the fact that the connection $_N\nabla_{\hat y}$ is integrable. (4.7.5.7-9) just write down what means that $\Phi_y$ is $_N\nabla_{\hat y}$-parallel, the Frobenius lift of $R$ being $\Phi_R$ (cf. [De3, 1.1.3.3 and 1.4.2.5]; see also ($E_1$) and ($E_2$) of 3.6.1.1.1 2)).
\medskip
To see 4.7.5 ii), we follow closely the proof of [De3, 1.4.2.5 ii)]. First the Poincar\'e lemma assures the existence of formal power series $\mu_j\in B(k)[[x_1,\ldots,x_{e_0}]]$ such that $d\mu_j=\vep_j$; they are unique up to a constant. Now
(4.7.5.7-9) can be rewritten as:
\medskip
\item{\bf (4.7.10.1)} $d\Bigl(\sum^{e_0}_{j=1}\vph(n_j)\bigl(\pi(b_i)\bigr)
d\Phi_{R*}(\mu_j)\Bigr)=0$;
\medskip
\item{\bf (4.7.10.2)} $d\Bigl(n\bigl(\pi(b_i)\bigr)+\sum_{j=1}^{e_0} n_j\bigl(\pi(b_i)\bigr)\mu_j-\vph(\,n_j)
\bigl(\pi(b_i)\bigr) d\Phi_{R*}(\mu_j)\Bigr)=0$;
\medskip
\item{\bf (4.7.10.3)} $d\Bigl(\log\,n\pi(a_i)+\sum^{e_0}_{j=1} n_j\bigl(\pi(a_i)\bigr)\mu_j\Bigr)=0$.
\medskip
$\vph(p\,n_j)\in M_y\otimes_{W(k)} R$ is non-zero mod $p$ and $N$ is an abelian group; so we have no problems in defining $\log n$. The equations (4.7.10.1-3) are of similar type: they are obtained by taking the $d$-operator of linear
combinations of functions; (4.7.25) is an equation similar to the one obtained by taking the $d$-operator of [De3, 1.4.2.7]. The values of $\mu_j$'s at 0 are determined by the equations (4.7.16-18). To see why this is so, let $t_j:=n_j\mu_j(0)$, $j\in S(1,e_0)$, and let $t:=\sum^{e_0}_{j=1}t_j$, We have $t\in {\rm Lie}(N)[{1\over p}]$. To determine $\mu_j(0)$, $j\in S(1,e_0)$, is equivalent to determine $t$. The equations (4.7.16-18) can be
rewritten respectively as:
\medskip
\item{\bf (4.7.10.4)} $\vph(t)\bigl(\pi(b_i)\bigr)=0$;
\medskip
\item{\bf (4.7.10.5)} $\vph(pt)\bigl(\pi(b_i)\bigr)=pt\bigl(\pi(b_i)\bigr)$;
\medskip
\item{\bf (4.7.10.6)} $t\bigl(\pi(a_i)\bigr)=0.$
\medskip
\medskip
Let $B_0:=\{a_1,\ldots,a_e\}\cup\bigl\{\pi(a_i),\ldots,\pi(a_e)\bigr\}$. $t$  acts trivially on the subspace of $M_y$ generated by $B_0$, i.e. on $F^0_y\cup\vph_y(F^0_y)$. It is easy to see that there is a unique $t\in {\rm End}(M_y)$ such that
$t(F^0_y)=\{0\}$, $t$ acts trivially on $M_y/F^0_y$ and 4.7.27-29 hold. Argument: $\vph(t)-t$ is uniquely determined by these requirements and so $\vph(t)=t$; as $t$ acts trivially on $F^0_y$ and $M_y/F^0_y$, we get $t=0$.
\smallskip 
There is another way to see that $t=0$: if $\pi$ is the trivial permutation this is obvious; but, as to be pointed out in 4.7.11 1) below, the situation gets reduced to such a context. 
\smallskip
The last fact of ii), to prove that $q_j\in R$ (obviously $q_j(0,...,0)=1$), is similar to  [De3, 1.4.5], starting from (4.7.10.1-3) and using the element
$$t_R:=\sum n_ju_j\in {\rm Lie}(N)[{1\over p}][[x_1,\ldots,x_{e_0}]].$$
\indent
The proof of 4.7.6 is entirely analogous to the proof of [De3, 1.4.7]: we are in a context of a moduli scheme ($\Mn_{W(k)}$) of principally polarized abelian schemes and so the Kodaira--Spencer map
$$gr_N\nabla_{\hat y}:\bar T_{R/W(k)}\to {\rm Hom}_R(F^1_y\otimes_{W(k)} R, F^0_y\otimes_{W(k)} R)$$
of $_N\nabla_{\hat y}$ (with $\bar T_{R/W(k)}:=\bar\Om_{R/W(k)}^\ast$) is injective, having  ${\rm Lie}(N)\otimes_{W(k)} R$ as its image (we have a canonical injective $W(k)$-linear map 
${\rm Lie}(N)\hookrightarrow {\rm Hom}_{W(k)}(F^1_y,F^0_y)$).
\medskip
Now 4.7.8 is obvious. The first statement of 4.7.9 is just a restatement of 4.6 P1. The second statement of 4.7.9 is a particular case of the principle expressed in 4.7.17 below: it is Fact 4 of 2.3.11 which allows us to be in a context involving just Shimura $p$-divisible groups; it reobtains in the case of a SHS $(f,L_{(p)},v)$, with $k(v)=\FF_p$, the main result of [No1].
\medskip
{\bf 4.7.11. On the general theory of crystalline coordinates for $p$-divisible groups over $k$.} We itemize the ideas and results. We recall $2$ is invertible in $k$.
\smallskip
{\bf 1)} Based on [De3, 1.4] and on the results pertaining to a SHS and presented till 4.7.4, the above part 4.7.5-29 is essentially trivial. What we did: we interpreted [De3, 1.4] as a property of deformations of Shimura adjoint Lie $\sg_k$-crystals. In other words, the $p$-divisible group over $k$ associated to $(M_y,\vph_y)$ (of 4.7.1) is a direct sum of two $p$-divisible groups $D^1_y$ and $D_y^2$ such that:
\medskip
-- the slopes of $D_y^1$ are precisely 0 and 1, while $D_y^2$ does not have integral slopes;
\smallskip
-- the $G$-deformations of $(A_y,p_{A_y})$ (see 4.7.0) correspond to deformations of $D_y^1$, while ``keeping" $D_y^2$ fixed.
\medskip
So, if 
$$M_y=M_y^1\oplus M_y^2$$
 is the direct sum decomposition corresponding to $D_y^1$ and $D_y^2$, then (cf. 4.7.1 2))
$${\rm Lie}(N)\subset {\rm End}(M_y^1).$$ 
The main goal of including 4.7.5-29 here is: to get some familiarity with the new type of equations (differential or not) showing up when we pass from (Shimura) ordinary $\sg_k$-crystals to cyclic diagonalizable Shimura $\sg_k$-crystals. The equations 4.7.5.4-5, 4.7.5.7, 4.7.5.9, 4.7.5.11, 4.7.5.13, 4.7.10.1 and 4.7.10.3 are typical to the cyclic diagonalizable situation: they do not show up in the classical setting of [De3, 1.4].  
\smallskip
{\bf 2)} Let $D$ be a $p$-divisible group over $k$.  Let $m$ be the product of the dimensions of $D$ and of its dual $D^t$. It is known (for instance, see [Il, 4.8]; see also 4.12.13 4) below) that $D$ has a universal deformation $\Md$ over 
$$R_m:=W(k)[[w_1,...,w_m]].$$ 
We assume $m\ge 1$. For any $W(k)$-epimorphism 
$$q:R_m\twoheadrightarrow R_n:=W(k)[[z_1,...,z_n]],$$ 
with $n\in\NN$, $n\le m$, we get naturally (from $\Md$) a $p$-divisible group $\Md_q$ over $R_n$. If $n=m$, then $z_i=w_i$, $\forall i\in S(1,m)$.
\medskip
{\bf The two goals.} The first goal of the general problem of crystalline coordinates (for $p$-divisible groups) is to find a good description (if possible canonical) of the $p$-divisible object $\DD(\Md_q)$ of $\Mm\Mf_{[0,1]}^\nabla(R_n)$, for epimorphisms $q$ of interest; we have freedom (as part of the good description) for the choice of a Frobenius lift $\Phi_{R_n}$ of $R_n$. If we have $\Phi_{R_n}(z_i)=z_i^p$ (resp. $\Phi_{R_n}(z_i+1)=(z_i+1)^p$), $\forall i\in S(1,n)$, we speak about an additive (resp. multiplicative) type situation. These are the most common (and handable) type of situations. The second goal of the problem of crystalline coordinates is to compare, in the case we do not have a canonical description, the different descriptions obtained in the process of achieving the first goal. In what follows we essentially concentrate on the first goal. 
\medskip
{\bf The extra Shimura structure.} Very often the $\sg_k$-crystal $(\tilde M,\tilde\vph)$ of $D$ gets naturally the structure of a Shimura $\sg_k$-crystal $(\tilde M,\tilde\vph,\tilde G,(\tilde t_\al)_{\al\in\tilde\Mj})$ with an emphasized family of tensors, and then we take $q$ to be the maximal formally smooth quotient of $R_m$ preserving (in the crystalline sense of 2.2.20 AX) the tensors involved (for concrete examples see [Va2, 5.4-5]); so $n=dd((\tilde M,\tilde\vph,\tilde G))$ (cf. 2.2.21 UP and 2.3.17.2). So we get a $p$-divisible object with tensors ${\got C}_q$ of $\Mm\Mf_{[0,1]}^\nabla(R_n)$, defined by $q$; forgetting the tensors it is $\DD(\Md_q)$. ${\got C}_q$ is uniquely associated to a uni plus versal Shimura $p$-divisible group over $R_n$. The Theorem 3.6.18.5 (as well as its relative counterparts of 3.6.18.7), in essence, achieves the first goal of the theory of crystalline coordinates. What we have to do is explained in the following two steps.
\medskip
{\bf a)} {\it First we choose, according to needs and desires, a Frobenius lift of $R_n$ taking the ideal $I_n:=(z_1,...,z_n)$ into itself, and a Frobenius endomorphism $\Phi$ of $\tilde M\otimes_{W(k)} R_n$ of the form $g_{R_n}(\tilde\vph\otimes 1)$, with $g_{R_n}\in\tilde G(R_n)$ which modulo $I_n$ is the identity. We choose also a filtration $F^1(\tilde M)$ of $\tilde M$ such that $(\tilde M,F^1(\tilde M),\tilde\vph,\tilde G)$ is a Shimura filtered $\sg_k$-crystal. We get, as in 3.6.18, a $p$-divisible object with tensors 
$${\got C}:=(\tilde M\otimes_{W(k)} R_n,F^1(\tilde M)\otimes_{W(k)} R_n,g_{R_n}(\tilde\vph\otimes 1),(\tilde t_\al)_{\al\in\tilde\Mj})$$ 
of $\Mm\Mf_{[0,1]}(R_n)$.}
\smallskip
{\bf b)} {\it We determine connections which respect the $\tilde G_{R_n}$-action and make ${\got C}/p{\got C}$ to be viewed as an object of $\Mm\Mf_{[0,1]}^\nabla(R_n)$. We look at the resulting Kodaira--Spencer maps (one for each such connection) and we pick up one such connection of whose Kodaira--Spencer map is injective (an isomorphism onto its image) modulo the maximal ideal $(I_n,p)$ of $R_n$; the choices in a)  should be made in such a way that such a connection does exist. Then using 3.6.18.7.0 we lift the picked up connection to a connection which respects the $\tilde G_{R_n}$-action and makes ${\got C}$ to be viewed as a $p$-divisible object with tensors of $\Mm\Mf_{[0,1]}^\nabla(R_n)$. So ${\got C}$ and ${\got C}_q$, as $p$-divisible objects with tensors of  $\Mm\Mf_{[0,1]}^\nabla(R_n)$, are isomorphic under an isomorphism which mod $I_n$ is the identity (cf. 2.2.21 UP). Here and in what follows, the respecting of the $\tilde G_{R_n}$-action is in the same sense as of 3.6.1.1.1.}
\medskip
Using these two steps a) and b), we regain (cf. 3.6.18.2) the results of [De3, 2.1], without any computational effort. After a brief digression on K3 surfaces, in 3) to 9) below we deal with cyclic diagonalizable contexts.
\medskip
{\bf 2A) Digression: the case of K3 surfaces.} A great part of the above ideas of 4.7.11 and of 3.6.18 can be adapted to the context of alternating or symmetric $p$-divisible objects of $\Mm\Mf_{[-1,1]}(R_n)$ (or of $\Mm\Mf_{[0,2]}(R_n)$) which are versal (cf. also 3.6.1.6 Ph). Below, as a sample, we briefly digress on the geometric context of K3 surfaces (see [Va6] for general abstract contexts; see [Va12] for other geometric context: of cubic fourfolds, etc.) 
\smallskip
The case of (the existence of) crystalline coordinates for the deformation space of an ordinary K3 surface $S$ over an arbitrary algebraically closed field $\tilde k$ of characteristic $p\ge 2$, is ``settled" by [Fa2, th. 10 and the remarks after] (see also [De3, 2.1 B] for the case $p\ge 3$); for $p=2$, ``settled" is used in the sense of the first goal mentioned above: however, the Fact of 4.14.3 E and 2B) below can be combined to achieve as well the second goal. [Fa2, th. 10] is not stated as such; however, starting from [De3, 2.1.6] it can be checked it applies. 
\smallskip
To explain this, let ${\got C}$ (resp. ${\got C}^1$) be the symmetric principally quasi-polarized, uni plus versal, filtered $F$-crystal over $\tilde k[[t_1,...,t_{20}]]$ of the universal deformation space of $S$ over ${\rm Spf}(\tilde k[[t_1,...,t_{20}]])$ (resp. over $\tilde k[[t_1,...,t_{20}]]$ constructed as the one of [De3, 2.1.7] having $H$ as its underlying $W(k)[[t_1,...,t_{20}]]$-module but for $p=2$). Warning: the Frobenius lift of $W(k)[[t_1,...,t_{20}]]$ we work with takes $t_i$ into $t_i^p$, $\forall i\in S(1,20)$. We assume the canonical lift of $S$ is defined by the $W(k)$-epimorphism $W(k)[[t_1,...,t_{20}]]\twoheadrightarrow W(k)$ annihilating all $t_i$'s. Forgetting the symmetric principally quasi-polarization, ${\got C}$ and ${\got C}^1$ are objects of $p-\Mm\Mf_{[0,2]}(W(k)[[t_1,...,t_{20}]])$. Moreover, they are ordinary in the sense of [De3, 1.3.3]. So let ${\got U}_0\hookrightarrow {\got U}_1\hookrightarrow {\got C}$ (resp. ${\got U}_0^1\hookrightarrow {\got U}_1^1\hookrightarrow {\got C}$) be the slope filtration in (filtered $F$-crystals) of ${\got C}$ (resp. of ${\got C}_1$) obtained as in [De3, 2.1.6] (see top of [De3, p. 109]). So ${\got U}$ and ${\got U}^1$ are uni plus versal objects of $p-\Mm\Mf_{[0,1]}(W(k)[[t_1,...,t_{20}]])$ and so we can apply [Fa2, th. 10] to them: we deduce the existence of a unique $W(k)$-isomorphism $I:{\rm Spec}(W(k)[[t_1,...,t_{20}]])\tilde\to {\rm Spec}(W(k)[[t_1,...,t_{20}]])$ which as a $W(k)$-isomorphism of $W(k)[[t_1,...,t_{20}]]$ leaves invariant the ideal $(t_1,...,t_{20})$ and which is such that $I^*({\got U}^1)$ is isomorphic to ${\got U}$ under an isomorphism $ISO$ which mod $(t_1,...t_{20})$ is the natural identification. Referring to [De3, 2.1.7], the element $c$ of it is uniquely determined by $a$, $b_1$,..., $b_{20}$. So, as we are in a symmetric principally quasi-polarized context, there is a unique way to extend $ISO$ to an isomorphism $I^*({\got C}^1)\tilde\to {\got C}$; modulo $(t_1,...,t_{20})$ is as well the natural identification. 
\smallskip
The case of polarized ordinary K3 surfaces over $\tilde k$ of degree relatively prime to $p$ (this excludes entirely the case $p=2$) can be deduced as well from a) and b)  (or from 4.7.5 performed abstractly, i.e. in the context of Shimura-ordinary $\sg_k$-crystals having the $A$-degree of definition equal to $1$; see Theorem of 4) below) and from 2.2.21 UP, by a natural process of moving (i.e. of lifting things) from the context of $SO(2,19)$-groups to the context of ${\rm GSpin}(2,19)$-groups. In other words, we can work entirely in the context of a Shimura filtered $\sg_{\tilde k}$-crystal $(\tilde M,F^1(\tilde M),\tilde\vph,\tilde G,(\tilde t_{\al})_{\al\in\tilde\Mj})$, such that:
\medskip
-- $(\tilde M,\tilde\vph)$ is ordinary and $(\tilde M,F^1(\tilde M),\tilde\vph)$ is diagonalizable;
\smallskip
-- $\tilde G^{\rm ad}$ is an absolutely simple, split $W(\tilde k)$-group of $B_{10}$ Lie type;
\smallskip
-- $\tilde G^{\rm ab}=\GG_m$;
\smallskip
-- the faithful representation $\tilde G^{\rm der}\hookrightarrow GL(\tilde M)$ is the spin representation of the simply connected group cover of $\tilde G^{\rm ad}$ (so $\dim_{W(k)}(\tilde M)=2^{10}$). 
\medskip
We consider a direct summand $\tilde M_1$ of $\Mt(M)$ normalized by $\tilde G$ and such that:
\medskip
a) $\tilde M_1[{1\over p}]$ is normalized by $\tilde\vph$;
\smallskip
b) the representation of $\tilde G$ on $\tilde M_1$ factors through $\tilde G^{\rm ad}$, giving birth to the symmetric faithful representation of $\tilde G^{\rm ad}$ of dimension $21$;
\smallskip
c) there is a symmetric perfect form $\tilde\psi_1$ on $\tilde M_1$ fixed by $\tilde G^{\rm ad}$ and $\tilde\vph$ at the same time.
\medskip
The existence of $\tilde M_1$ can be deduced easily by considering $\ZZ_p$-structures as in 2.2.9 8) and applying arguments as in 2.3.9 B and E.
The mentioned lifting process can be stated as (see \S7 and [Va12] for general forms of it):
\medskip
{\bf 2B) Proposition (the lifting process in a K3 context).} {\it We assume $\tilde k$ is just perfect of odd characteristic. Let $R$ be a local, henselian, regular, formally smooth $W(\tilde k)$-algebra having $\tilde k$ as its residue field and let $\Phi_R$ be a Frobenius lift of $R$. Then any symmetric $p$-divisible object $(M_R,(F^i(M_R))_{i\in\{0,1,2\}},\Phi_{M_R},\nabla,p_{M_R})$ of $\Mm\Mf_{[0,2]}^\nabla(R)$, with $F^2(M_R)$ of rank $1$ and with $M_R$ of rank $21$, is isomorphic to the symmetric $p$-divisible object of $\Mm\Mf_{[0,2]}^\nabla(R)$ obtained from a filtered $F$-crystal with tensors of the form 
$${\got C}^{\rm lift}=(\tilde M\otimes_{W(\tilde k)} R^\wedge,F^1(\tilde M)\otimes_{W(\tilde k)} R^\wedge,\tilde g_R(\tilde\vph\otimes 1),\nabla_0,(\tilde t_{\al})_{\al\in\tilde\Mj}),$$ 
with $\tilde g_R\in\tilde G(R^\wedge)$, by moving from $\tilde M$ to $(\tilde M_1(1),\tilde\psi_1(1))$ (in a way similar to the Fact 2 of 2.2.10 but involving the Tate twist $W(\tilde k)(1)$).}
\medskip
{\bf Proof:} We can assume $\Phi_R$ is of additive type in the maximal point of ${\rm Spec}(R)$. ``Ignoring" $\nabla$ and $\nabla_0$, the existence of ${\got C}^{\rm lift}$ is a consequence of the henselian assumption: any element of $\tilde G^{\rm ad}(R^\wedge)$ lifts to an element of $\tilde G(R^\wedge)$. The existence and uniqueness of $\nabla_0$ are implied by 3.6.18.4.4 1) and 3.6.18.8. But 3.6.18.8 also implies that under the process of moving from $\tilde M$ to $\tilde M_1(1)$, $\nabla_0$ ``becomes" $\nabla$. This ends the proof.
\medskip
{\bf 3)} We look closer at [De3] to understand the real principles underlying the main results of loc. cit. We see that there are precisely three main ideas. We formulate them in the language of this paper.
\medskip
{\item {\bf a)}} {\it We have a Shimura $\sg_k$-crystal which admits a ``canonical split": with the notations of 2), this is a cocharacter $\mu:\GG_m\to\tilde G$ with special properties; in particular, it produces a direct sum decomposition $\tilde M=\tilde F^1\oplus\tilde F^0$, with $\be\in\GG_m(W(k))$ acting through $\mu$ on $\tilde F^s$ as the multiplication with $\be^{-s}$, $s=\overline{0,1}$ (to be compared with 3.1.4).}
\smallskip  
{\item {\bf b)}} {\it The resulting Shimura filtered $\sg_k$-crystal $(\tilde M,\tilde F^1,\tilde\vph,\tilde G)$ is cyclic diagonalizable.}
\smallskip
{\item {\bf c)}} {\it We can construct good additive coordinates. By this we mean [De3, 1.4.2] (see also 5) below).}
\medskip
{\bf 3')} Using the cocharacter $\mu$, we get a natural $\ZZ_p$-structure of $(\tilde M,\tilde F^1,\tilde\vph,\tilde G)$, as in 4.4.7. So $\tilde M=M_{\ZZ_p}\otimes_{\ZZ_p} W(k)$, with $M_{\ZZ_p}$ a free $\ZZ_p$-module; under this identification, $\tilde G$ is obtained, by extension of scalars, from a reductive subgroup $\tilde G_{\ZZ_p}$ of $GL(M_{\ZZ_p})$. Moreover $\mu$ is defined over a finite unramified DVR extension $W(k_{\mu})$ of $\ZZ_p$; we always choose the smallest such extension. Let $T$ be the smallest torus of $\tilde G_{\ZZ_p}$ such that $\mu$ is obtained by pull back from a cocharacter of $\tilde G_{W(k_{\mu})}$ factoring through $T_{W(k_{\mu})}$ (cf. b) of 3) and 2.2.16.1). Let $k_{\mu}^{\rm ad}$ be the smallest subfield of $k_{\mu}$ such that the cocharacter of the image of $T$ in $\tilde G^{\rm ad}$ naturally defined by $\mu$, is defined over $W(k_{\mu}^{\rm ad})$. Let $C$ be the centralizer of $T$ in $\tilde G_{\ZZ_p}$. It is a reductive group over $\ZZ_p$. 
\smallskip
The principle behind the choice of $W(k)$-bases in [De3, 1.4.2] and in 4.7.5 is: 
\medskip\noindent
{\it We pick up a maximal torus $T_C$ of $C^{\rm der}$; only the $C(\ZZ_p)$-conjugacy class of $T_C$ matters. If $C^{\rm der}$ is split, then the most natural choice for $T_C$ is: a split maximal torus of $C^{\rm der}$.}
\medskip
When $C^{\rm der}$ is not split, still we have quite natural choices for $T_C$. For instance, we can take $T_C$ to contain a maximal split torus of $C^{\rm der}$ and to split over the smallest unramified extension of $\ZZ_p$ over which $C^{\rm der}$ splits (cf. [Ti2, 1.10]); more concretely, we can take $T_C$ as a maximal torus of a Borel subgroup of $C^{\rm der}$. For the most general applications it is important to take $T_C$ arbitrarily. Let $T_C^1$ be the maximal torus of $\tilde G_{\ZZ_p}$ generated by $T_C$ and by the maximal subtorus of $Z(C)$; as in the proof of 3.11.1 we get a Shimura filtered $\sg_k$-crystal $(\tilde M,\tilde F^1,\tilde\vph,T^1_C)$. We have 
$$\tilde\vph({\rm Lie}(T_C^1)\otimes_{\ZZ_p} W(k))={\rm Lie}(T_C^1\otimes_{\ZZ_p} W(k)).$$
Let $\{t_1,...,t_{\dim_{W(k)}(T_C^1)}\}$ be a $W(k)$-basis of ${\rm Lie}(T_C^1)\otimes_{\ZZ_p} W(k)$ formed by elements fixed by $\tilde\vph$. 
\smallskip
{\bf 4)} It is easy to see that choosing a multiplicative type Frobenius lift of $R_n$, in 2) we can not, for an arbitrary cyclic diagonalizable Shimura $\sg_k$-crystal $(\tilde M,\tilde\vph,\tilde G)$, choose $g_{R_n}$ to be the identity element of $\tilde G(R_n)$. Even for the case of the $p$-divisible group of a supersingular elliptic curve over $k$ (so $m=n=1$) we can not choose $g_{R_1}=1_{\tilde M\otimes_{W(k)} R_n}$. To get out from this deadlock, in general we have to either use the additive situation (see 5) below) or to work with the multiplicative situation, with a carefully chosen $g_{R_n}$ (see 6) below). In what follows we assume $dd((\tilde M,\tilde\vph,\tilde G))>0$. 
\smallskip
By $\tilde G$-canonical (multiplicative) coordinates for $(\tilde M,\tilde\vph,\tilde G)$ we mean that, referring to a) and b) of 2), we can choose simultaneously the following five things:
\medskip
{\bf i)} a lift $\tilde F^1$ of it of quasi CM type;
\smallskip
{\bf ii)} $n=dd((\tilde M,\tilde\vph,\tilde G))$;
\smallskip
{\bf iii)} a multiplicative type Frobenius lift of $R_n$  ($\Phi_{R_n}(z_i+1)=(z_i+1)^p$, $\forall i\in S(1,n)$);
\smallskip
{\bf iv)} $g_{R_n}=1_{\tilde M\otimes_{W(k)} R_n}$; 
\smallskip
{\bf v)} a connection which respects the $\tilde G$-action and of whose Kodaira--Spencer map modulo $(I_n,p)$ is injective.
\medskip
If the choices i) to v) are possible, we say $(\tilde M,\tilde\vph,\tilde G)$ has $\tilde G$-canonical (multiplicative) coordinates. We have:
\medskip
{\bf Theorem.} {\it We consider a cyclic diagonalizable Shimura  $\sg_k$-crystal $(\tilde M,\tilde\vph,\tilde G)$. The following three statements are equivalent:
\medskip
{\bf a)} it has $\tilde G$-canonical (multiplicative) coordinates;
\smallskip
{\bf b)} its attached Shimura adjoint Lie $\sg_k$-crystal ${\got L}$ is ordinary, with its multiplicity of the slope $-1$ being maximal, i.e. equal to $dd((\tilde M,\tilde\vph,\tilde G))$;
\smallskip
{\bf c)} it is a $\tilde G$-ordinary $\sg_k$-crystal and its $A$-degree of definition is 1.}
\medskip
{\bf Proof:} If b) holds, then the cyclic adjoint factors of any $\tilde G$-ordinary $\sg_k$-crystal ${\got C}_0$ produced by $(\tilde M,\tilde\vph,\tilde G)$ are either trivial or are totally non-compact of constant type (cf. 3.1.0 c) and 3.10.6). From this and 3.9.2 we get: $(\tilde M,\tilde\vph,\tilde G)$ is a $\tilde G$-ordinary $\sg_k$-crystal. As its attached Shimura adjoint Lie $\sg_k$-crystal is ordinary, we get (cf. F2 of 3.10.7 and Fact of 3.11.2 D): its $A$-degree of definition equal to 1. So b) implies c). Obviously c) implies b) (to be compared with 4.6 P4). 
\smallskip
From the proof of 3.6.18.4 B) (see P6 of it) we get (cf. also 2)), that b) implies a). In other words, choosing ii) to iv) to hold and choosing $\tilde F^1$ to be the $\tilde G$-canonical lift of $(\tilde M,\tilde\vph,\tilde G)$ (4.4.13.2 shows that in fact we have no other choice for $\tilde F^1$), based on 3.6.18.4 P6 (see also 3.6.18.2), we get that we can assume that v) holds as well.
\smallskip
We assume now that a) holds. From 3.6.18.4 P5 we get that subject to the choices i) to iv), the image modulo $(I_n,p)$ of the Kodaira--Spencer map of any $\tilde G$-invariant connection, is a $k$-vector space of dimension at most equal to the multiplicity of the slope $-1$  of the Shimura adjoint Lie $\sg_k$-crystal attached to $(\tilde M,\tilde\vph,\tilde G)$. So from ii) and v) we get that this multiplicity is equal to $dd((\tilde M,\tilde\vph,\tilde G))$. But this implies that all cyclic adjoint factors of ${\got C}_0$ are either trivial or are totally non-compact of constant type. As above, we get that $(\tilde M,\tilde\vph,\tilde G)$ is a $\tilde G$-ordinary $\sg_k$-crystal. From F2 of 3.10.7 we get that b) holds. This proves the Theorem. 
\medskip
In particular, we reobtain in a much faster manner the results of 4.7.6-7.
\smallskip
{\bf 5)} We come back to the general context of 3'). Let $\tilde N$ be the commutative, unipotent subgroup of $\tilde G$ acting trivially on $\tilde F^0$ and on $\tilde M/\tilde F^0$. We take $R_n$ to be the completion of $\tilde N$ in its origin (so $n=\dim_{W(k)}(\tilde N)=dd((\tilde M,\tilde\vph,\tilde G))$). $\tilde N$ has a natural $W(k_{\mu}^{\rm ad})$-structure; let $\tilde N_{W(k_{\mu}^{\rm ad})}$ be the subgroup of $\tilde G_{W(k_{\mu}^{\rm ad})}$ which over $W(k)$ is $\tilde N$. If $\Phi_{R_n}$ is of additive type, then we speak about additive coordinates. Due to the condition on the Kodaira--Spencer map of b) of 2), we do not have too many choices in connection to additive coordinates: the $p$-divisible object with tensors ${\got C}_q$ of $\Mm\Mf_{[0,1]}^\nabla(R_n)$ is a Shimura filtered $F$-crystal of the form $(\tilde M,\tilde F^1,\tilde\vph,\tilde G,\tilde N,\tilde f)$ (cf. 2.2.21 UP); so the only choice we have is in choosing (the automorphism) $\tilde f$ (of $R_n$). The choice of $\tilde f$ corresponds to a choice of an additive type Frobenius lift of $R_n$. [Fa2, rm. iii) after th. 10] tells us that it does not matter which choice we make. We propose one such choice, which to us looks convenient (at least from the point of view of 8) below).
\smallskip
The maximal subtorus $T_C^2$ of $T^1_C$ contained in $\tilde G^{\rm der}$ (resp. the image of $T$ in $\tilde G^{\rm ad}$) might not be split. It splits over an unramified finite DVR extension $W(k_2)$ (resp. $W(k_1)$) of $\ZZ_p$. As before, we choose the smallest such extension. We have natural $\ZZ_p$-monomorphisms $W(k_{\mu}^{\rm ad})\hookrightarrow W(k_1)\hookrightarrow W(k_2)$. Let $\tilde N_0$ be the smallest integral, closed subgroup of $\tilde G$ such that it contains $\tilde N$ and its Lie algebra is taken by $p\tilde\vph$ into itself. We consider two variants.
\medskip
{\bf Variant 1.} Let $\GG_a(i)$, $i\in S(1,n)$, be different $\GG_a$ subgroups of $\tilde N_{W(k_1)}$, whose Lie algebras are normalized by $T_{W(k_1)}$. Most common, they are not uniquely determined. 
\medskip
{\bf Variant 2.} Let $\GG_a(i)$, $i\in S(1,n)$, be the different $\GG_a$ subgroups of $\tilde N_{W(k_2)}$, whose Lie algebras are normalized by ${T^2_C}_{W(k_2)}$. They are uniquely determined.
\medskip
Under both these variants, we feel inclined of choosing $\Phi_{R_n}$ such that it takes the subring 
$$R(i):=W(k_2)[[z_i]]$$ 
of $R_n$ defined via (completion in the origin of) the projection (it is a group homomorphism) of $\tilde N_{W(k_2)}$ on $\GG_a(i)$ having $\prod_{j\in S(1,n)\setminus\{i\}} \GG_a(j)$ as its kernel, into itself; it is also natural to take 
$$\Phi_{R_n}(z_i)=z_i^p,$$ 
$\forall i\in S(1,n)$. Under both these variants, we can use $T_C$ (and so $T_C^2$) to determine uniquely $z_i$'s, through different considerations but this (presently) looks irrelevant to us.
\smallskip
{\bf 6)} We move now to the case when $\Phi_{R_n}$ is of multiplicative type (i.e. it is as in iii) of 4)). We work under one of the two variants of 5). Two natural questions arise: 
\medskip
{\bf Q} {\it What happens if the conditions of the Theorem of 4) are not satisfied? Is there a best choice of $g_{R_n}$ in such a case?} 
\medskip
We propose here one good choice of $g_{R_n}$ working out in all cases (i.e. in the case of the mentioned Theorem, provided we work under Variant 1, the element $g_{R_n}$ to be introduced below is the identity). We consider 
$$a_i\in {\rm Lie}(\GG_a(i))$$ 
generating this Lie algebra. Warning: we always assume that the $W(k)$-subbasis $\{a_1,...,a_{n},t_1,...,t_{\dim_{W(k)}(T_C^1)}\}$ of ${\rm Lie}(\tilde G)$ extends (by adding elements at the ``end") to a $W(k)$-basis $\{e_1,...,e_{\dim_{W(k)}(\tilde G)}\}$ of ${\rm Lie}(\tilde G)$ as in 2.2.12 c). Up to multiplying each $a_i$ by some invertible element of $W(k)$, this is automatic so under Variant 2; under Variant 1 we need to use the cyclic diagonalizability property in a way entirely similar to 3.4.3.0 (1) (but without having good ``control" on the exponents of $p$), in order to make suitable choices so that this assumption is satisfied. 
\smallskip
For $i\in S(1,n)$, let $n_i\in\NN\cup\{\infty\}$ be the infimum of the set of numbers $m_i\in\NN$ such that 
$\tilde\vph^{m_i}(a_i)\subset F^0({\rm Lie}(\tilde G))[{1\over p}]$. Let 
$$S(1,n)^0:=\{i\in S(1,n)|n_i\in\NN\}.$$ 
\indent
For each $i\in S(1,n)\setminus S(1,n)^0$ we consider an element $b_i\in {\rm Lie}(\tilde G)$ subject to the requirements:
\medskip
\item{\bf R1} $\tilde\vph(pb_i)=b_i$;
\smallskip
\item{\bf R2} the $W(k)$-span of $b_j$'s is the same as the $W(k)$-span of $a_j$'s, where $j\in S(1,n)\setminus S(1,n)^0$.
\medskip
Let $g_{R_n}\in\tilde N(R_n)$ be defined by the formula 
$$g_{R_n}=1_{\tilde M\otimes_{W(k)} R_n}+\sum_{i\in S(1,n)^0} a_iz_i.$$  
We have:
\medskip
{\bf Theorem.} {\it We can choose a connection on $\tilde M\otimes_{W(k)} R_n$ respecting the $\tilde G_{R_n}$-action, making 
$$(\tilde M\otimes_{W(k)} R_n,\tilde F^1\otimes_{W(k)} R_n,g_{R_n}(\tilde\vph\otimes 1))$$ 
to be viewed as a $p$-divisible object of $\Mm\Mf_{[0,1]}^\nabla(R_n)$ and whose Kodaira--Spencer map is an injection (modulo $(I_n,p)$).}
\medskip
{\bf Proof:} As $(\tilde M,\tilde F^1,\tilde\vph)$ is cyclic diagonalizable we choose a $W(k)$-basis of $\tilde M$ as in 2.2.1 d). Using it, we consider the convenient matrix form (as in 3.6.8.6 (CMF)) of the equations describing connections on $\tilde M\otimes_{W(k)} R_n/pR_n$ which respect the $\tilde G_{R_n}$-action and make $(\tilde M\otimes_{W(k)} R_n/pR_n,\tilde F^1\otimes_{W(k)} R_n/pR_n,g_{R_n}(\tilde\vph\otimes 1))$ to be viewed as an object of $\Mm\Mf_{[0,1]}^\nabla(R_n)$. We view (cf. 3.6.8.6.1) this convenient matrix form as an equality between two elements of $LIE/pLIE$, where $LIE:={\rm Lie}(\tilde G)\otimes_{W(k)} \bar\Om_{R_n/W(k)}$. 
\smallskip
We consider the $R_n$-basis $\{e_i\otimes dz_j|i\in S(1,\dim_{W(k)}(\tilde G)),\, j\in S(1,n)\}$ of $LIE$. Due to the cyclic diagonalizability assumption on $\{e_1,...,e_{\dim_{W(k)}(\tilde G)}\}$, identifying the coefficients w.r.t. its reduction mod $p$ of the mentioned two elements of $LIE/pLIE$ and taken them mod $(I_n,p)$, we come across systems of equations of first type very close in spirit to the ones of (27) of 3.6.18.4 B). More precisely, if $\pi$ is the permutation of $S(1,\dim_{W(k)}(\tilde G))$ obtained as in 2.2.12 c) (so $\tilde\vph(e_i)=p^{\vep_i}e_{\pi(i)}$ for some $\vep_i\in S(-1,1)$), the equations (with coefficients in $k$) are of the form
$$d_{\pi(i)j}+x_{\pi(i)j}=b_{ij}x_{ij}^p,\leqno (CYCLIC)$$
$i\in S(1,\dim_{W(k)}(\tilde G))$, $j\in S(1,n).$ This is in conformity with 3.6.8.6 (CMF) and 3.6.8.6.2 (so the right hand side of (CYCLIC) corresponds to the right hand side of 3.6.8.6 (CMF), etc.). We have:
\medskip
{\bf i)} $b_{ij}=0$ iff $i\notin S(1,n)$ (cf. the definition of $\pi(F^0)$ of 3.6.8.6 and the fact that $\Phi_{R_n}$ is of multiplicative type);
\smallskip
{\bf ii)} $d_{ij}\neq 0$ iff $i=j\in S(1,n)^0$ (cf. the shape of $g_{R_n}$ and the definition of $C_lE_G^{-1}$ in 3.6.8.6). 
\medskip
We fix $j_0\in S(1,n)^0$ and $i\in S(1,\dim_{W(k)}(\tilde G))$ and we look just at the equations involving $x_{ij_0}$, $x_{\pi(i)j_0}$,..., $x_{\pi^{u_i-1}(i)j_0}$, with $u_i\in\NN$ as the smallest number such that $\pi^{u_i}(i)=i$. If $i\notin S(1,n)\setminus S(1,n)^0$ and if $j_0\notin\{i,\pi(i),...,\pi^{u_i-1}(i)\}$ then from i) and ii) we get $x_{ij}=0$. Warning: if $i\in S(1,n)\setminus S(1,n)^0$, then $x_{ij}$ can be non-zero. 
\smallskip
So based on the previous paragraph and on ii), we easily get: $\forall i\in S(1,n)^0$, the endomorphism $\bar e_i$ of $M_y/pM_y$ defined by evaluating any such connection at ${\dl\over {\dl z_i}}$ and taking everything modulo $(I_n,p)$, is of the following ``$F^0$-truncated" form
$$\bigl(\sum_{s=0}^{n_i} b_s(i)\tilde a_s(i)\bigr)+\sum_{i\in S(1,n)\setminus S(1,n)^0} c_ia_i,\leqno (TRUNC)$$ 
where $b_s(i)\in k$, with $b_0(i)\neq 0$, where $c_i\in k$ (are not necessarily zero), and where $\tilde a_s(i)$ is the reduction mod $p$ of $e_{a(i,s)}$, with $a(i,s)\in S(1,\dim_{W(k)}(\tilde G))$ uniquely determined by the property that $\tilde\vph^{s}(a_i)\in B(k)e_{a(i,s)}$. We choose one such connection $\nabla+p$ with the property that for any $i\in S(1,n)\setminus S(1,n)^0$, the similarly defined endomorphism $\bar e_i$ of $M_y/pM_y$ is a non-zero multiple of $b_i/pb_i$ (cf. R1 and the proof of 3.6.18.2).  
\smallskip
The Kodaira--Spencer map of $\nabla_p$: its image modulo the ideal $(z_1,...,z_n)$ of $R_n/pR_n$ contains the $k$-span of $b_i/pb_i$'s, $i\in S(1,n)\setminus S(1,n)^0$, and by induction on the (increasing) possible values of $n_i$'s, with $i\in S(1,n)^0$, it contains (cf. (TRUNC)) $a_i/pa_i$. So the Theorem follows (cf. a) and b) of 2)) from the standard way (see 3.6.18.7) of lifting $\nabla_p$ to a connection on $\tilde M\otimes_{W(k)} R_n$ having all required properties. 
\medskip
{\bf Corollary.} {\it ${\rm Spf}(R_n)$ has a natural structure of a formal torus: it is defined by the coordinates $z_i+1$, $i\in S(1,n)$.}  
\medskip
In the case of the Theorem of 4), we have $k_1=\FF_p$ and so $\tilde N_0=\tilde N$. So $S(1,n)^0$ is the empty set and we have 
$$g_{R_n}=1_{\tilde M\otimes_{W(k)} R_n}.$$ 
Moreover, we can choose $a_i$'s such that $\tilde\vph(pa_i)=a_i$, $\forall i\in S(1,n)$. For such a choice of $a_i$'s we can choose the connection $\nabla$ on $\tilde M\otimes_{W(k)} R_n$ to have the following standard logarithmic form
$$\nabla=\dl+\sum_{i\in S(1,n)} a_id(ln(z_i+1)).\leqno (LN)$$
Here $\dl$ is the connection on $\tilde M\otimes_{W(k)} R_n$ that annihilates $\tilde M$. Obviously such a connection has all the properties mentioned in the Theorem. Accordingly, we view Theorem as the generalization of [Ka3, 3.7.1-3 and 4.3.1-2]: if $\tilde G=GL(\tilde M)$ and if $(\tilde M,\tilde F^1,\tilde\vph)$ is the filtered $\sg$-crystal associated to the canonical lift of an ordinary $p$-divisible group over $k$ having its dimension equal to half of its rank, then for a very particular choice of $a_i$'s --subject to the equalities $\tilde\vph(pa_i)=a_i$-- we get exactly loc. cit.  
\smallskip
In the general case it is still arguable which are the best choices for $a_i$'s and $z_i$'s. Of course the $a_i$'s and $z_i$'s are interrelated. We leave this problem of the best (i.e. the most natural) such choices to get matured on its due time. 
\smallskip
{\bf 7)} The ideas of 3.6 (for instance see 3.6.1.3 and 3.6.8.1.2) suggest that $W(k)$-algebras of the form $R_m$ and $R_n$ are not always the right ones for getting crystalline coordinates. Often they should be replaced by the $p$-adic completion of ind-\'etale algebras over smooth $W(k)$-algebras. For instance. we can work with the $p$-adic completion of the henselization of the localization of $W(k)[z_1,...,z_n]$ w.r.t. its maximal ideal $(p,z_1,...,z_n)$ (cf. also 3.6.20 5)). But these henselizations are still ``loosing" a lot of information.
\smallskip
Very often the $W(k)$-algebras we ``can get", have ``plenty" of $k$-valued points (i.e. their spectra are regular, $AG$ $k$-schemes). For instance, working with Shimura $\sg_k$-crystals which are ordinary, using the Fact of 3.6.18.9 (and of 3.6.19), we can construct canonical crystalline coordinates for $W(k)$-algebras whose spectra are the $p$-adic completion of $\NN$-pro-\'etale covers of the $p$-adic completion of open, affine subschemes of ${\rm Spec}(W(k)[z_1,...,z_n])$, like the one obtained by inverting all $z_i$'s; so we have uni plus versal deformations of Shimura $p$-divisible groups over these spectra which are ordinary (in the usual way): they are defined or associated to ``prescribed" (see 3.6.1.3 and 3.6.11) filtered $F$-crystals with tensors. 
\smallskip
More generally, referring to the Corollary of 6), instead of $W(k)[[z_1,...,z_n]]$ we can work with the $p$-adic completion of an $\NN$-pro-\'etale cover of the spectrum of 
$$R_n^{\rm al}:=W(k)[z_1,...,z_n][{1\over {\prod_{i\in S(1,n)} z_i+1}})],$$
cf. Fact of 3.6.18.9 applied to 
$${\got C}_q^{\rm al}:=(\tilde M\otimes_{W(k)} R_n^{\rm al},\tilde F^1\otimes_{W(k)} R_n^{\rm al},g_{R_n^{\rm al}}(\tilde\vph\otimes 1)),$$ 
with $g_{R_n^{\rm al}}$ having the same expression as $g_{R_n}$ but being ``viewed" as an element of $\tilde G(R_n^{\rm al})$.  
This opens a completely new horizon: variation of crystalline coordinates in families. 
\smallskip
When we are in the context of a principally polarized abelian variety, we do not know for which ``prescriptions" we still get (principally) polarized abelian schemes over the spectra of these $W(k)$-algebras (in some variants we might have to deal with polarizations which are not necessarily principal). We recall: we do need a polarization in order to get abelian schemes over such spectra (cf. [FC, 1.10 a)]). 
\smallskip
We do expect that the study of variations of crystalline coordinates will lead to nice modular properties. We hope to come back to this idea in a future paper.  
\smallskip
{\bf 8)} In 3) to 6) we had some choices: 
\medskip
-- of $\mu$ (and so of $\tilde F^1$ and of $\tilde N$);
\smallskip
-- of $T_C$;
\smallskip
-- and  of $z_i$ and $a_i$, $i\in S(1,n)$. 
\medskip
We do not stop to see if we can make all these choices to result in something canonical (up to isomorphisms of $(\tilde M,\tilde\vph,\tilde G$)). However, two samples are in order.
\medskip
{\bf U.} If $(\tilde M,\tilde F^1,\tilde\vph,\tilde G)$ is $U$-ordinary, then $\tilde F^1$ is uniquely determined and so $\mu$ too.
\smallskip
{\bf T.} If in fact $(\tilde M,\tilde F^1,\tilde\vph,\tilde G)$ is $T$-ordinary, then $T_C$ is as well uniquely determined. Implicitly, all $a_i$ and $z_i$'s are uniquely determined up to multiplication with elements of $\GG_m(W(k_2))$. 
\medskip
One can go even further on (see [Va9]) in the $T$-ordinary case, to get a canonical choice of $a_i$'s and $z_i$'s. Moreover, we have:
\medskip
{\bf Corollary.} {\it The $\tilde G$-canonical multiplicative coordinates of the Theorem of 4) (``defined" by $z_i+1$'s of iii) of 4)) are canonical (up to isomorphisms of the form: $z_i+1$ goes to $\prod_{j=1}^n (z_j+1)^{n_{ij}}$, where $n_{ij}\in\ZZ_p$ are the entries of an invertible $n\times n$ matrix).}
\medskip
{\bf Proof:} This can be checked in the same manner as in the case of an ordinary $p$-divisible group (i.e. of when we have $\tilde G=GL(\tilde M)$). 
\smallskip
First, the need to allow changes by isomorphisms is implied by the fact that in v) of 4), we did not specified which connection we choose: the set of such connections is in one-to-one correspondence to the set of invertible $n\times n$ matrices with coefficients in $\ZZ_p$ (this can be deduced from 3.6.18.2 --see equation (24) of it-- by remarking that, due to iv) of 4), --in our case-- all $b_{ijl}$'s are $0$; see also 3.6.18.7.0). Of course, we can always choose the connection as in (LN) of 6). However, for the sake of flexibility we think it is important not to fix a priori a connection subject to v) of 4). 
\smallskip
Second, there is a unique multiplicative type Frobenius lift of $R_n$ for which $g_{R_n}$ is the identity element (here a fixed choice in v) of 4) is implicit). This is entirely the same as [Ka4, A2.2 (2)]: using the direct sum decomposition $M_y=M_y^1\oplus M_y^2$ of 1), we can add extra variables (their number is $n_y^1-n$, where $n_y^1$ is the number of slopes $-1$ of $({\rm End}(M_y^1),\vph_y)$) so that the situation gets reduced to loc. cit. (see Fact of 4.14.3 E below for a more of a principle approach; the mentioned place is stated in terms of $p=2$ but a) of it holds for any prime). This ends  the proof.
\medskip
This Corollary motivates the terminology: $\tilde G$-canonical.
\smallskip
Coming back to the general case of 3), 6) gives us natural good structures of a formal torus of the moduli formal scheme $\M(\tilde M,\tilde\vph,\tilde G)$ of deformations of $(\tilde M,\tilde\vph,\tilde G)$ (see the proof of 3.12.1): we do not stop to see if there is a unique such group structure on $\M(\tilde M,\tilde\vph,\tilde G)$ (the origin being fixed; it is defined by the $F^1$-filtration $\tilde F^1$ of $\tilde M$). However, the part of the above Corollary and of its proof referring to isomorphisms and connections still applies to 6); we get that the coordinates of the Corollary of 6) are determined by the choice of $\tilde F^1$, $a_i$'s and $z_i$'s, $i\in S(1,n)^0$, up to suitable isomorphisms of $R_n$ which satisfy the following 2 properties (see also 9) below for the second one):
\medskip
{\bf P1.} $z_i$ is mapped into $z_i+\tilde z_i$, with $\tilde z_i\in (z_j)_{j\in S(1,n)\setminus S(1,n)^0}$, $\forall i\in S(1,n)^0$;
\smallskip
{\bf P2.} $z_i+1$ is mapped into $\tilde z_i+\prod_{j\in S(1,n)\setminus S(1,n)^0} (z_j+1)^{m_{ij}}$, with $\tilde z_i\in (z_j)_{j\in S(1,n)^0}$, $\forall i\in S(1,n)\setminus S(1,n)^0$; here $m_{ij}\in\ZZ_p$ are the entries of an invertible $\abs{S(1,n)\setminus S(1,n)^0}\times \abs{S(1,n)\setminus S(1,n)^0}$ matrix. 
\medskip
Due to the fact that $c_i$'s of the proof of 6) can be non-zero, we presently can not be more precise (than P1-2) on the structure of the isomorphisms of $R_n$ up to which the coordinates of the Corollary of 6) are determined.
\smallskip
We assume now that there is a reductive subgroup $\tilde G_0$ of $\tilde G$ such that $(\tilde M,\tilde\vph,\tilde G_0)$ is also a Shimura $\sg_k$-crystal. It is not difficult to see (in connection to 6)) that by working under Variant 1 of 5) we can make all choices to be functorial w.r.t. the monomorphism $\tilde G_0\hookrightarrow\tilde G$, so that the moduli formal scheme of deformations of $(\tilde M,\tilde\vph,\tilde G_0)$ is a subtorus of the formal torus $\M(\tilde M,\tilde\vph,\tilde G)$.
\smallskip
{\bf 9)} We come back to the general case of 6). Let $r$ be the Lie $p$-rank of $(\tilde M,\tilde\vph,\tilde G)$. We have $r\le dd((\tilde M,\tilde\vph,\tilde G))$. If equality holds then we are in the situation of the Theorem of 4). We assume now that $1\le r<dd((\tilde M,\tilde\vph,\tilde G))$. In this case there is a uniquely determined quotient 
$$q_{n,r}:R_n\twoheadrightarrow R_{r}:=W(k)[[x_1,...,x_r]]$$ 
such that (in the sense of 4)) we can take, for a suitable Frobenius lift $\Phi_{R_r}$ of $R_r$, $g_{R_r}$ to be the identity. The uniqueness is obtained immediately, by writing the connection on $\tilde M\otimes_{W(k)} R_r$ in a logarithm form
$\dl+\sum_{i\in S(1,r)} b_i(r)d(ln(x_i+1))$, with $b_i(r)\in {\rm Lie}(\tilde G)\otimes_{W(k)} R_r$; as we are dealing with injective Kodaira--Spencer maps and as $g_{R_r}$ is the identity, $\Phi_{R_r}$ has to be of essentially multiplicative type and so (cf. 3.6.18.1.1) we can assume that it takes each $x_i+1$ into $(x_i+1)^p$. We get immediately that all $b_i(r)$'s belong to ${\rm Lie}(\tilde G)$ and are fixed by $p\tilde\vph$. Based on the proof of 6), we conclude: 
$${\rm Ker}(q_{n,r})=(z_i)_{i\in S(1,n)^0}.$$ 
\indent
{\bf Corollary.} {\it Regardless of the choices made, there is a canonical subtorus of $\M(\tilde M,\tilde\vph,\tilde G)$ of dimension $r$.}
\medskip
{\bf 10)} All of 1) to 9) can be entirely adapted (limitations do apply) to the case when $k$ is just perfect. For instance, the part of b) of 2) involving liftings still makes sense if $k$ is just $\dim_{W(k)}(\tilde G)$-simply connected (to be compared with 3.6.18.5.4 3)).   
\medskip
{\bf 4.7.12. Definition.} We come back to the situation described in 4.7.5. The lifts of $y$ (to DVR's which are finite extensions of $W(k)$) corresponding to torsion points of the formal torus of $G$-deformations of $A_y$, are called $G$-quasi-canonical lifts of $y$.
\medskip
{\bf 4.7.13. Exercise.} We assume 4.4.6 and we refer to 4.7.12. If $k=\FF$, prove that the $G$-quasi-canonical lifts of a $G$-ordinary point $y:{\rm Spec}(k)\to\Mu$ are giving birth to special points of $\Mn$ (i.e. the abelian varieties obtained from $\Ma$ by pull back through these $G$-quasi-canonical lifts have complex multiplication). Show that they are the only lifts of $y$ to DVR's of mixed characteristic, which are special points of $\Mn$. Hint: if $k(v)=\FF_p$, this is well known; if $k(v)\ne\FF_p$ then use 4.4.6, 4.7.11 1) and just imitate the standard arguments of the case $k(v)=\FF_p$.
\medskip
{\bf 4.7.14. Remark.} 
We do not assume anymore $N=N_0$. 
From 4.7.11 3), 4) and 6) (resp. from 4.7.11 5)) we deduce that for any toric point $y_1:{\rm Spec}(k)\to\Mn_{k(v)}$ there are good $G$-multiplicative (resp. $G$-additive) coordinates of the moduli formal scheme of $G$-deformations of $(A_{y_1},p_{A_{y_1}}):=y_1^*(\Ma,\Mp_{\Ma})$. The same thing remains valid if we do not work in a principally polarized context involving (as well de Rham components of) Hodge cycles. Some particular situations are listed below. 
\medskip
{\bf a)} For any point $y_1:{\rm Spec}(k)\to\Mn_{k(v)}$ whose attached Shimura adjoint Lie $\sg_k$-crystal is isomorphic to the extension to $k$ of the Shimura adjoint Lie $\sg_k$-crystal of ${\got C}_\om$ of 4.1.5 (i.e. for any $k$-valued $G(\om)$-ordinary point $y_1$ of $\Mn_{k(v)}$), we get good $G(\om)$-multiplicative coordinates of the moduli formal scheme of $G$-deformations of $(A_{y_1},p_{A_{y_1}})$. 
\smallskip
{\bf b)} For any point $y:{\rm Spec}(k)\to\Mu$, we get good $GSp$-multiplicative coordinates (resp. good multiplicative coordinates) of the moduli formal scheme of deformations of the principally polarized abelian variety $(A_y,p_{A_y})$ (resp. of deformations of $A_y$).
\smallskip
{\bf c)} For any supersingular elliptic curve $E_s$ over $k$ we get a good additive coordinate $t$ of its moduli formal scheme $\M_1$ of deformations. This reobtains [Og, 3.15]. Warning: in loc. cit., in order to keep track of principal quasi-polarizations, the $\sg_k$-crystal $(M_1,\vph_1)$ of $E_s$ is described by a $W(k)$-basis $\{a,b\}$ of $M_1$ such that $\vph_1(a)=pb$ and $\vph_1(b)=-a$, and not by a $W(k)$-basis $\{a_1,b_1\}$ such that $\vph_1(a_1)=pb_1$ and $\vph_1(b_1)=a_1$; as we can take $a_1={\al}a$ and $b_1={\al}b$, with $\al\in W(k)$ such that $\sg_k^2(\al)=-1$, the passage from $\{a,b\}$ to $\{a_1,b_1\}$ is automatic.
\smallskip
This suggests, that in the case of perfect fields, we can work as well with variants of the notion of cyclic diagonalizability: for instance, referring to 2.2.1 d), by allowing $\vph(e_i)$ to be plus or minus $p^{\vep_i}e_{\pi(i)}$ we get the notion of almost cyclic (or circular) diagonalizability. 
\smallskip
4.7.11 6) provides us as well with a good multiplicative coordinate $t$ of $\M_1$. It is worth pointing out: in both additive and multiplicative cases, the Frobenius endomorphism of $M_1\otimes W(k)[[t]]$ takes $a$ into $pb$ and $b$ into $a+tb$ (of course the Frobenius of $W(k)[[t]]$ depends on which case we are). 
\smallskip
{\bf d)} For any abelian variety over $k$ which is a product of abelian varieties $A_{y_1}$ as in a) or b), we get good multiplicative coordinates of its moduli formal scheme of deformations.
\medskip
{\bf 4.7.14.1. Remark.} More generally than d) we have the Corollary of 1.8. It is a direct consequence of 2.2.18 and of 4.7.11 6).
\medskip
{\bf 4.7.14.2. Remarks.} Corollary of 1.8 is the most general abelian varieties context in which we presently  can prove the existence of useful (crystalline) (additive or multiplicative)  coordinates. From 2.2.19 we get that it does not handle the case of all abelian varieties over $k$. In the cases left, the hopes of obtaining any theory of multiplicative coordinates which has any significance at all (i.e. it is useful and motivated), are based on 3.13. So the problem 13 from [Oo1, p. 13] is still open for the cases left. We leave it to the reader to formulate 4.7.14.1 in the relative context (of abelian varieties over $k$ whose attached $\sg_k$-crystals have the extra structure of a cyclic diagonalizable Shimura $\sg_k$-crystal) (cf. 4.7.11 5) and 6)). 
\medskip
{\bf 4.7.14.3. Study problem.} In what follows we refer to 4.7.14 a). For any good $G(\om)$-multiplicative coordinates, we can define $G(\om)$-quasi-canonical lifts of $y_1$. It looks to us interesting to study them. In particular, if $y_1$ is a $G$-ordinary point and $N\neq N_0$, (as we had choices in 4.7.11) we do not know if (or when) these $G$-quasi-canonical lifts are uniquely determined or if (or when) the facts of 4.7.13 remain true. 
\medskip 
{\bf 4.7.15. Exercise.} Coming back to 4.7.3, prove that $\forall n\in\NN$, $_N\nabla_{\hat y}$ mod $p^n$ is algebraic, i.e. the
differential forms $\vep_j$ of 4.7.3.1, $j=\overline{1,e_0}$, are algebraic over $N_{W_n(k)}$. Hint: use 4.4.7 to reduce to the case $k=\overline{k(v)}$; then use the fact that everything, like the bases
${(a_i)}_{1\le i\le e}$ and ${(b_i)}_{1\le i\le e}$, $\mu_y$, etc., are definable over (the Witt ring of) a finite field containing $k(v)$; another way: cf. 3.6.1.3.
\medskip
{\bf 4.7.16. Problem.} Obtain the general form of the main result of [No1] and of [No2] and show that all of [No2, section 2] remains valid in the context of $G$-ordinary points. Hint: if $N=N_0$ use 4.7.5-29, 4.7.13 and 4.7.17 below and just imitate the arguments of [No1-2]; if $N\neq N_0$ then use 4.7.11 and work formally (i.e. abstractly).
\medskip
{\bf 4.7.17. Formal subtori.} Every formal subscheme of $\Mm_G(y)$ of 4.7.8 which can be \lq\lq defined\rq\rq\ starting from a smooth subgroup $N_1$ of $N$ whose Lie-algebra is taken by $p\vph_y$ into itself is a formal subtorus. This can be checked using 4.7.5-29 or 4.7.11 8). The same thing remains true in the general context of 4.7.11 4). Here, for the sake of future reference, we state this principle explicitly, just for the abstract context of Shimura $\sg_k$-crystals. We use the notations of 4.7.11 2) to 4) and assume that the equivalent statements of the Theorem of 4.7.11 4) hold. We have:
\medskip
{\bf Theorem.} {\it If there is a reductive subgroup $\tilde G_1$ of $\tilde G$ such that $(\tilde M,\tilde\vph,\tilde G_1)$ is also a Shimura-ordinary $\sg_k$-crystal, then the formal torus of dimension $n_1:=dd((\tilde M,\tilde\vph,\tilde G_1))$ of deformations of $(\tilde M,\tilde F^1,\tilde\vph,\tilde G_1)$ is canonically a formal subtorus of the formal torus of deformations of $(\tilde M,\tilde F^1,\tilde\vph,\tilde G)$.}
\medskip
{\bf Proof:} We consider a family of tensors $(\tilde t_{\al})_{\al\in\tilde\Mj_1}$, with $\Mj_1$ a set containing $\Mj$, such that $(\tilde M,\tilde\vph,\tilde G_1,(\tilde t_{\al})_{\al\in\tilde\Mj_1})$ is a Shimura $\sg_k$-crystal with an emphasized family of tensors. We first remark, that the $F^1$-filtration of $\tilde M$ defining the $\tilde G_1$-canonical lift of $(\tilde M,\tilde\vph,\tilde G_1)$ is $\tilde F^1$, cf. 4.4.13.2-3 and 3.11.1 a).
Based on 2.2.21 UP, it is enough to show that there is a filtered $F$-crystal over the special fibre ${\rm Spec}(k[[z_1,...,z_n]])$ of ${\rm Spec}(R_n)$, with $n$ as in 4.7.11 4), such that:
\medskip
-- it is a $p$-divisible object of $\Mm\Mf_{[0,1]}^\nabla(R_n)$ of the form $(\tilde M\otimes_{W(k)} R_n,\tilde F^1\otimes_{W(k)} R_n,\tilde\vph\otimes 1,\tilde\nabla)$, with the connection $\tilde\nabla$ on $\tilde M\otimes_{W(k)} R_n$ as in v) of 4.7.11 4);
\smallskip
-- the Frobenius lift of $R_n$ takes $z_i+1$ to $(z_i+1)^p$;
\smallskip
--  $\tilde\nabla$ modulo the ideal of $R_n$ generated by $z_{n_1+1},...,z_n$ respects the $\tilde G_1$-action (i.e. it annihilates $\tilde t_{\al}$, $\forall\al\in\tilde\Mj_1$).
\medskip
The construction of such a filtered $F$-crystal is trivial: ${\rm Lie}(\tilde N)$ has a $W(k)$-basis formed by elements fixed by $p\tilde\vph$ an such that a $W(k)$-subbasis of it is a $W(k)$-basis of ${\rm Lie}(\tilde N_1)$, where $\tilde N_1:=\tilde G_1\cap\tilde N$; so we just need to ``apply" (LN) of 4.7.11 6) to it. This proves the Theorem.
\medskip
{\bf 4.7.18. Problem.} State a theorem similar to 4.7.5 in the general case, i.e. when $k(v^{\rm ad})$ is not $\FF_p$. Hint: use 4.7.4 and 4.7.11 2) to 6).
\medskip\smallskip
{\bf 4.8. The Galois property of $G$-ordinary points.}
\medskip
{\bf 4.8.1. Notations.} We continue to assume $k=\bar k$. We start with a morphism $y:{\rm Spec}(k)\to\Mu/H_0$. Let ${\rm Spec}(V)\buildrel{z_V}\over\to\Mn/H_0$ be a lift of $y$, with $V$ a DVR which is a finite, flat extension of $W(k)$. Let $K:=V[{1\over p}]$. Let $(A_V,p_{A_V}):=z_V^*(\Ma_{H_0},\Mp_{\Ma_{H_0}})$ and let $(A_K,p_{A_K})$ be its generic fibre. Let 
$$\rho:{\rm Gal}(K)\to GL\bigl(H^1_{\acute et}(A_{\bar K},\QQ_p)\bigr)(\QQ_p)$$
be the $p$-adic Galois representation associated to $A_K$. Let $G_{\QQ_p}\hookrightarrow GSp\bigl(H^1_{\acute et}(A_{\bar K},\QQ_p),p_{A_K}\bigr)$ be the subgroup fixing the $p$-components of the \'etale components of the family of Hodge cycles $(w_\al)_{\al\in\Mj^\prime}$ with which $A_{\bar K}$ is naturally endowed.
\medskip
{\bf 4.8.2. Theorem.} {\it There is a finite field extension $K_1$ of $K$ such that the restriction $\rho_1:=\rho|{\rm Gal}(K_1)$ factors through the group of $\QQ_p$-valued points of a connected solvable subgroup of $G_{\QQ_p}$.}
\medskip
{\bf Proof:} We can assume $K_1$ is big enough so that the family $(w_\al)_{\al\in\Mj^\prime}$ is defined over $K_1$. In fact, as we assumed $k=\bar k$, for this part we can take $K_1=K$; argument: as $V$ is strictly henselian, $A_V$ has level-$m$ structure, $\forall m\in\NN$, $(m,p)=1$. So $\rho_1$ factors through $G_{\QQ_p}$. We can assume $G$ is not a torus. 
\smallskip
Let $(M_y,F_y^1,\vph_y,G_{W(k)})$ be the Shimura filtered $\sg_k$-crystal attached to the $G$-canonical lift $z:{\rm Spec}(W(k))\to\Mn/H_0$ of $y$. 
Let $\mu_y:\GG_m\to GL(M_y)$ be its canonical split.
Let $\bar h_i$, $i=\overline{1,d}$, be the elements of ${\rm Lie}(G_{W(k)})$ defined as in 4.1 but working with $\mu_y$ (cf. also 4.4.7). 
\smallskip
Let $N_1$ be the smallest subgroup of $G_{W(k)}$ such that ${\rm Lie}(N_1)$
contains ${\rm Lie}(N_0)$ ($N_0$ was def. in 4.7.1) and the Lie algebra generated by $\bar h_i$, $i=\overline{1,d}$. $N_1$ is an integral, solvable subgroup of $G_{W(k)}$. Let $B_{N_1}$ be a Borel subgroup of $G_{W(k)}$ containing $N_1$ and such that its Lie algebra is taken by $p\vph_y$ into itself. The solvability of $N_1$ and the existence of such Borel subgroups is a consequence of 4.1.4.1 and of b) of 4.4.1 3). Let $\tilde N_1^0$ be the Zariski closure in $G_{W(k)}$ of the connected component of the origin of the generic fibre of the intersection $\tilde N_1$ of all such Borel subgroups $B_{N_1}$ of $G_{W(k)}$ containing $N_1$. ${\rm Lie}(N_1)$ and ${\rm Lie}(\tilde N_1^0)$ are taken by $p\vph_y$ into themselves. 
\smallskip
With the notations of 4.7.1 we get that ${\rm Lie}(N_1)\otimes_{W(k)} R$ and ${\rm Lie}(B_{N_1})\otimes_{W(k)} R$ are taken by $n(\vph_y\otimes 1)$ into themselves. So for the morphism ${\rm Spec}(V_1)\to {\rm Spec}(R)$ (with $V_1$ the ring of integers of $K_1$), defined naturally by $z_V$, we get (cf. 4.7.2) that the Galois representation $\rho_1$ leaves invariant ${\rm Lie}(N_2)$, ${\rm Lie}(\tilde N_2)$ and ${\rm Lie}(B_{N_2})$, with 
$$N_2\subset \tilde N_2^0\subset B_{N_2}\subset G_{\QQ_p}$$ 
such that $B_{N_2}$ is a Borel subgroup of $G_{\QQ_p}$ and $N_2$ and $\tilde N_2^0$ are connected, solvable subgroups; here $N_2$, $\tilde N_2$ and $B_{N_2}$ are respectively the $\QQ_p$-\'etale correspondents (via Fontaine's comparison theory; see [Fa2, th. 5]) of the crystalline ones $N_1$, $\tilde N_1^0$ and $B_{N_1}$. In order to apply loc. cit. (mainly we need [Fa2, th. 7]), we just have to remark that the $W(k)$-homomorphism $R\to V_1$ lifts to a $W(k)$-homomorphism $R\to R_{V_1}$, where the $W(k)$-algebra $R_{V_1}$ is defined as in [Fa2, \S 2] starting from an arbitrary uniformizer of $V_1$.
\smallskip
But as $B_{N_2}$ is its own normalizer in $G_{\QQ_p}$, we get that $\rho_1$ factors through the group of $\QQ_p$-valued points of $B_{N_2}$. Repeating the argument for any Borel subgroup of $G_{W(k)}$ as above, we get that, by extending again $K_1$ by a finite cover (as $\tilde N_1^0$ might be different from $\tilde N_1$),  $\rho_1$ factors through the group of $\QQ_p$-valued points of
$\tilde N_2^0$. This proves the Theorem.
\medskip
{\bf 4.8.3. Remarks. a)} 4.8.2 solves a question of C.-L. Chai. It was R. Pink who pointed out the construction of $N_1$ (of the proof), provided 4.7.2, c) of 4.4.1 2) and a) and b) of 4.4.1 3) were
known. Always the unipotent radical of $\tilde N_1^0$ is the unipotent radical $N_0$ of $N_1$. Argument: from 4.7.1 2) we get that $N_0$ is the unipotent radical of the parabolic subgroup $P_{\le 0}$ of $G_{W(k)}$ whose Lie algebra is $W_0({\rm Lie}(G_{W(k)}),\vph_y)$; so, as the integral, closed subgroup $P_{=0}^0$ of $G_{W(k)}$ whose generic fibre has $W_0({\rm Lie}(G_{B(k)}),\vph_y)$ as its Lie algebra is reductive and has a natural $\ZZ_p$-structure $P_{=0\ZZ_p}^0$ (see 3.11.2 C), we just have to remark that:
\medskip
-- the intersection of all Borel subgroups of $P_{=0\ZZ_p}^{0\rm der}$ is a finite, flat group scheme over $\ZZ_p$;
\smallskip
-- $N_1$ is contained (cf. 3.11.2 C) in $Z(P_{=0})N_0$.
\medskip
Warning: not always we have ${\rm Lie}(\tilde N_1^0)={\rm Lie}(N_1)+{\rm Lie}(Z^0(G_{W(k)}))$, with $Z^0(G_{W(k)})$ as the maximal torus of $Z(G_{W(k)})$.
\smallskip
{\bf b)} Using [De4, 3.1 and 3.2 a)] and Fontaine's comparison theory (in a way similar to the proof of the Fact of 2.2.9 1)), we get that we do not need to use $\tilde N_1^0$ as an intermediary: $\rho_1$ factors through $N_2(\QQ_p)$, without having to replace $K$ by $K_1$. 
\smallskip
{\bf c)} There are plenty of examples in which $\tilde N_1^0$ is a Borel subgroup of $G_{W(k)}$ (for instance if we are in the context of 4.6 P7; if moreover, $G^{\rm ad}_{\ZZ_p}$ is a $\ZZ_p$-simple group then $N_1\cap G_{W(k)}^{\rm der}$ is a Borel subgroup $G_{W(k)}^{\rm der}$) and plenty of examples when it is not (for instance, if $G_{\ZZ_p}^{\rm ad}$ is a simple, split group of $C_{\ell}$ Lie type, with $\ell\ge 2$).
\smallskip
{\bf d)} A result similar to 4.8.2 can be proved for an arbitrary toric point ${\rm Spec}(k)\to\Mn_{k(v)}$. Warning: in such a general context $N_1$ is not necessarily a solvable group. So, if we still want $\tilde N_1$, then we have to define it using parabolic subgroups containing $N_1$. 
\smallskip
There are many cases when the image of $N_1$ in $G^{\rm ad}_{W(k)}$ is not contained in a Borel subgroup of $G^{\rm ad}_{W(k)}$; for instance, this is the case if we are in a context modeled on 4.4.13.3.1, with $\tilde G_0$ of the mentioned place of $A_1$ Lie type: in such a case the image of $N_1$ in $G^{\rm ad}_{W(k)}$ is $G^{\rm ad}_{W(k)}$ itself.   
\smallskip
{\bf e)} 4.8.2 and d) are properties of Shimura-ordinary $\sg_k$-crystals (instead of 4.7.2 we have to use 2.2.21). Moreover the condition $k=\bar k$ can be weaken: we just need the existence of a Borel subgroup $B_{N_1}$ of $G_{W(k)}$ of the type mentioned above; warning: in such a situation we can not (in general) replace $\tilde N_2$ by $N_2$. 
\smallskip
We assume now that $k$ is a finite field, that we are in the context of 4.8.1 and that 4.4.6 holds. We have:
\medskip
{\bf Theorem.} {\it Up to a passage to a finite field extension of $k$, the existence of such a Borel subgroup $B_{N_1}$ is automatic.}
\medskip
{\bf Proof:} We use the previous notations. Based on 4.4.6, by passage from $k$ to a finite field extension of it, we can assume that the generic fibre of $P_{\le 0}$ has a maximal torus whose Lie algebra is generated by elements fixed by $\vph_y$. Let $N_3$ be the connected subgroup of $G_{B(k)}$ generated by this maximal torus and by the generic fibre of the unipotent radical $P_{<0}$ of $P_{\le 0}$. Let $P_{<0}^{\acute et}\subset N_4\subset P_{\le 0}^{\acute et}\subset G_{\QQ_p}$ be the subgroups corresponding to $P_{< 0}$, $N_3$ and respectively $P_{\le 0}$ via Fontaine's comparison theory (applied to $z_V$). So the algebraic envelope $ALE$ of $\rho$ is a subgroup of $P_{\le 0}^{\acute et}$ such that its image in $P_{\le 0}^{\acute et}/P_{<0}^{\acute et}$ belongs to the image of $N_4$ in it. So $ALE$ is a subgroup of $N_4$. This ends the proof.
\medskip
Using this Theorem, we obtain (for instance, cf. 4.4.5 and 4.6.3 B) plenty of quasi-ordinary crystalline representations in the sense of [Pi, 2.10].
\smallskip
{\bf f)} We stated this section 4.8 in terms of $\QQ_p$-valued points; but everything can be restated in terms of $\ZZ_p$-valued points. 
\medskip\smallskip
{\bf 4.9. The passage from Shimura varieties of Hodge type to Shimura varieties of preabelian type.} Very often (like in the study of K3 surfaces, of cubic fourfolds, etc.), we have to consider integral aspects of Shimura varieties of preabelian type which are not of Hodge type. So here we extend a great part of 4.1-8 to the (often abstract) context of integral canonical models of Shimura quadruples of preabelian type; on the way we show that many properties pertaining to a SHS $(f,L_{(p)},v)$ are canonical, i.e. are properties just of $\Mn$ itself, not depending on who or how $f$ is.  
\medskip
{\bf 4.9.1. The initial setting.} We start with an arbitrary Shimura variety ${\rm Sh}(G_0,X_0)$ of adjoint, abelian type. Let $(G_0,X_0,H_0,v_0)$ be a Shimura quadruple with $v_0$ dividing a rational prime
$p\ge 3$. We fix a connected component $X_0^0$ of $X_0$ (it does not matter which one, cf. [Va2, 3.3.3]). Let $(G_i,X_i,H_i,v_i)$, $i=\overline{1,3}$, be Shimura quadruples having $(G_0,X_0,H_0,v_0)$ as their adjoint quadruples and such that $X_0^0\subset X_i$ (cf. [Va2, 2.4]). We assume the existence of two standard Hodge
situations $(f_i,{L_i}_{(p)},v_i)$ defined by $f_i:(G_i,X_i)\hookrightarrow
\bigl({\rm GSp}(W_i,\psi_i),S_i\bigr)$, $i=\overline{1,2}$, such that $H_i=G_i({L_i}_{(p)}\otimes_{\ZZ_{(p)}} \ZZ_p)$ (for $p\ge 5$ cf. [Va2, 6.4.2]; for $p=3$ cf. \S 6]). For them we use the standard notations of 2.3.1-3, except that we put a right lower index $i$ everywhere,
$i=\overline{1,2}$.
\medskip
We denote by $\Mn_0$ and $\Mn_3$ the integral canonical models of $(G_0,X_0,H_0,v_0)$ and respectively of $(G_3,X_3,H_3,v_3)$ (cf. [Va2, 6.4.1] for $p\ge 5$; for $p=3$ see \S 6). We get natural $O_{(v_0)}$-morphisms $\Mn_i\buildrel{j_i}\over\to\Mn_0$ (cf. [Va2, 3.2.7 4)]), $i=\overline{1,3}$.
\medskip
{\bf 4.9.1.1. The definition of $\Mc_i$'s.} Let $V_0:=W(\FF)$ and let 
$$\Mc_0:=\cap_{i_3\in FM} \, j_3({\Mn_3}_{V_0}),$$ 
where $FM$ is the set of all finite maps $i_3:(G_3,X_3,H_3,v_3)\to (G_0,X_0,H_0,v_0)$ such that $X_0^0\subset X_3$. The fact that $\Mc_0$ is non-empty is implied by the fact that $\Mc_0\times_{V_0}{}_j\CC$ is non-empty (with $j:V_0\hookrightarrow\CC$ an arbitrary $O_{(v_0)}$-monomorphism). The fact that $\Mc_0$ is an open closed subscheme of $\Mn_{0V_0}$ is implied for $p\ge 5$ by [Va2, 6.4.5.1] (for $p=3$ cf. \S 6). Similarly, working with the connected component $X_0^0$ of $X_i$, we define an open closed embedding $\Mc_i\hookrightarrow{\Mn_i}_{V_0}$, $i=\overline{1,3}$.
\medskip
{\bf 4.9.2. The new SHS.} Let $(G_4,X_4,H_4,v_4)$ be a Shimura quadruple which together with maps 
$h_i:(G_4,X_4,H_4,v_4)\to(G_i,X_i,H_i,v_i)$, $i=\overline{1,2}$, is a quasi fibre product of the below finite maps $q_1$ and $q_2$, i.e. we have a commutative diagram
$$
\def\mapright#1{\smash{
\mathop{\longrightarrow}\limits^{#1}}}
\def\mapdown#1{\Big\downarrow
\rlap{$\vcenter{\hbox{$\scriptstyle#1$}}$}}
\matrix{(G_4,X_4,H_4,v_4)&\mapright{h_1}&(G_1,X_1,H_1,v_1)\cr
\mapdown{h_2}&&\mapdown{q_1}\cr
(G_2,X_2,H_2,v_2)&\mapright{q_2}&(G_0,X_0,H_0,v)\times \bigl(\GG_m,\{x_0\},\GG_m(\ZZ_p),p\bigr)\cr}$$
in the category $qf-Sh$ defined in [Va2, 3.2.7 3)], with $G_4$ as the connected component of the origin of $G_2\times_{G_0\times\GG_m} G_1$; this diagram satisfies a quasi universal property as in [Va2, 2.4.0]. Here $q_i$ is defined by the adjoint map $(G_i,X_i,H_i,v_i)\to (G_0,X_0,H_0,v_0)$ and via the canonical map $h_i:(G_i,X_i)\to(GSp^{\rm ab},S^{\rm ab})$ defined via $f_i$, $i=\overline{1,2}$, where $\{x_0\}:=S^{\rm ab}$ and $GSp^{\rm ab}=\GG_m$. We assume $X_0^0\subset X_4$, cf. loc. cit. and [Va2, 3.2.7 3)]. 
\smallskip
We call $(G_4,X_4)$
as a Hodge twist of $(G_1,X_1)$ and $(G_2,X_2)$ and denote it by $(G_1,X_1)\tilde\times(G_2,X_2)$. 
\smallskip
We consider a Hodge quasi product $(G_1\times_{\Mh}G_2,X_1\times_{\Mh}X_2)$
of $(G_1,X_1)$ and $(G_2,X_2)$, and a Hodge quasi product $({\rm GSp}(W_1,\psi_1)\times_{\Mh} {\rm GSp}(W_2,\psi_2),S_1\times_{\Mh} S_2)$ of $({\rm GSp}(W_1,\psi_1),S_1)$ and $({\rm GSp}(W_2,\psi_2),S_2)$, cf. [Va2, Example 3 of 2.5]. We use similar notations for quadruples. We can  assume that 
$$X_4\subset X_1\times_{\Mh} X_2\subset S_1\times_{\Mh} S_2.$$ 
So we get (see 4.9.2.0 below) a new SHS 
$$(f_4,{L_1}_{(p)}\oplus {L_2}_{(p)},v_4)$$ 
defined by the composite injective map $f_4$ 
$$(G_4,X_4)\buildrel{\tilde i_4}\over\hookrightarrow \bigl({\rm GSp}(W_1,\psi_1)\times_H{\rm GSp}(W_2,\psi_2),S_1\times_H S_2\bigr)
\buildrel{i}\over\hookrightarrow ({\rm GSp}(W_1\oplus W_2,\psi_1\oplus \psi_2),S),$$ 
where $\tilde i_4:=(f_1\times_{\Mh} f_2)\circ i_4$, with 
$$i_4:(G_4,X_4)\hookrightarrow (G_1\times_H G_2,X_1\times_H X_2)$$ 
and with
$$f_1\times_{\Mh} f_2: (G_1\times_{\Mh} G_2,X_1\times_{\Mh} X_2)\hookrightarrow ({\rm GSp}(W_1,\psi_1)\times_{\Mh} {\rm GSp}(W_2,\psi_2),S_1\times_{\Mh} S_2)$$ 
as the natural inclusions, and where $i$ is a Segre embedding, as used in loc. cit. 
\smallskip
Let $\Mn_4$ be the integral canonical model of $(G_4,X_4,H_4,v_4)$ and let $\Mc_4$ be the open closed subscheme of ${\Mn_4}_{V_0}$ defined similarly to $\Mc_0$. We have a natural $O_{(v_0)}$-morphism $\Mn_4\to\Mn_1\times\Mn_2$. It is defined via the composite of $i_4$ with the natural injective map 
$$(G_1\times_H G_2,X_1\times_H X_2)\hookrightarrow (G_1\times G_2,X_1\times X_2),$$ 
cf. [Va2, 3.2.16 and 3.2.7 4)]. So we speak about points of $\Mn_4$ mapping into points of $\Mn_i$, $i\in\{1,2\}$.
\medskip
{\bf 4.9.2.0. The argument.} The fact that the triple $(f_4,{L_1}_{(p)}\oplus {L_2}_{(p)},v_4)$ is indeed a SHS, is related to the proof of Milne's conjecture (see 1.15.1, d) of 4.4.1 3) and 4.4.12). Let $\Mm_4$ be the integral canonical model of the Shimura quadruple associated naturally to the SHS $(1_{({\rm GSp}(W_1\oplus W_2),\psi_1\oplus\psi_2),S)},{L_1}_{(p)}\oplus {L_2}_{(p)},p)$. From constructions we have 
$$H_4=G_4(\QQ_p)\cap GSp({L_1}_{(p)}\oplus {L_2}_{(p)},\psi)(({L_1}_{(p)}\oplus {L_2}_{(p)})\otimes_{\ZZ_{(p)}} \ZZ_p).$$ 
So, from [Va2, 3.2.7 4)] we get a natural morphism 
$$\tilde m:\Mn_4\to\Mm_4.$$ 
We need to show that for any morphism $z:{\rm Spec}(W(k))\to\Mn_4$, working with $\tilde m$, the condition $(*)$ of 2.3.4 is satisfied. The idea behind checking this is rooted in [BLR, th. 1 of p. 109] and the density part of 4.2.1. To explain it, we index the main paragraphs by capital letters.
\medskip
{\bf A.} [Va2, 3.3.2] allows us to assume $z$ factors through $\Mc_4$. Moreover, 2.3.12-13 makes sense even without knowing that we are in a reductive context: using canonical splits of filtered $\sg_k$-crystals endowed with tensors, we can perform [Va2, 5.4-5] directly in the context of abelian unipotent groups. This goes as follows. Let $(M,F^1,\vph_4,(t_{\al})_{\al\in\Mj_4^\prime})$ be the $p$-divisible object with tensors of $\Mm\Mf_{[0,1]}(W(k))$ naturally associated to $z$ and $\tilde m$ (to be compared with 2.3.4; of course some choices --like of a family of tensors $(v_{\al})_{\al\in\Mj^\prime_4}$ of $\Mt((W_1\otimes_{\QQ} W_2)^*)$ such that $G_4$ is the subgroup of $GL(W_1\otimes_{\QQ} W_2)$ fixing them, like of a $\ZZ$-lattice of $W_1\otimes_{\QQ} W_2$, etc.-- are in order: they can be made as in 2.3.1-2). 
\smallskip
We consider the reductive subgroup $\tilde G_{4B(k)}$ of $GL(M[{1\over p}])$ fixing $t_{\al}$, $\forall\al\in\Mj_4^\prime$. Let $\tilde G_{4W(k)}$ be its Zariski closure in $GL(M)$. Using the canonical split cocharacter of $(M,F^1,\vph_4)$ (see also 2.2.1.2), as in 4.7.1 we define a subgroup $N$ of $\tilde G_{4W(k)}$. $N$ is always a smooth group over $W(k)$, identifiable with the affine group scheme defined by ${\rm Lie}(N)$, even if by chance $\tilde G_{4W(k)}$ is not a reductive subgroup of $GL(M)$. So, 2.2.10 can still be performed in the context of the $7$-tuple $(M,F^1,\vph_4,\tilde G_{W(k)},N,\tilde f,(t_{\al})_{\al\in\Mj_4^\prime})$, and so we can still perform 2.3.12-13 in this context (even if the maximal integral subgroup of $\tilde G_{W(k)}$ normalizing $F^1$ is not smooth). The last thing we need to add: we have a variant of the Theorem of 2.3.11 in the context of $(M,F^1,\vph_4,\tilde G_{W(k)},N,\tilde f,(t_{\al})_{\al\in\Mj_4^\prime})$, as the result [Va2, 4.1.5] (used in [Va2, 5.4.5]) involves arguments purely in characteristic $0$.   
\medskip
{\bf B.} Coming back to our concrete situation, from the above paragraph and the density part of 4.2.1, we deduce that it is enough to check condition $(*)$ of 2.3.4 for $z$ lifting a $G_4$-ordinary point of $\Mc_{4k(v_4)}$. The set of such points is the same, regardless of which morphism $\Mn_4\to\Mn_1$ or $\Mn_4\to\Mn_2$ we use to define them --by being mapped into a $G_i$-ordinary point of $\Mn_{ik(v_i)}$, $i\in\{1,2\}$--: the arguments of 4.9.3 below, at the level of isocrystals, apply, as the fact that a Shimura $F$-crystal over a perfect field is or is not Shimura-ordinary can be ``read out" (cf. 3.1.0 c)) from the Newton polygon of its attached Shimura adjoint Lie isocrystal. 4.9.2.1 to 4.9.3 below are stated in terms of Shimura adjoint Lie $\bar\sg$-crystals, starting from the fact that the triple $(f_4,{L_1}_{(p)}\oplus {L_2}_{(p)},v_4)$ is a SHS; however, without assuming this, it can be performed independently but in the context of adjoint Lie isocrystals over $k$. Of course, if we desire to simplify the things, we can define such $G_4$-ordinary points of $\Mn_{4k(v_4)}$, as points mapping into $G_i$-ordinary points of $\Mn_{ik(v_i)}$, $\forall i\in\{1,2\}$ (sort of intersection of subschemes of $G_{4k(v_4)}$: [Va2, 6.4.5.1] and 4.2.1 guarantee that this intersection is still an open, dense subscheme of $\Mn_{4k(v_4)}$).
\smallskip
Moreover, we can assume $z$ maps into a $G_i$-canonical lift of $\Mn_i$, $\forall i\in\{1,2\}$ (again cf. 4.9.3 below, performed with filtered isocrystals over $k$; b) of 4.4.1 2) implies: Shimura-canonical lifts can be ``read out" from attached Shimura adjoint filtered Lie isocrystals). This implies that the Galois representation defined by $z$ (working with $\tilde m$) factors through the $\ZZ_p$-valued points of a torus of $G_{\ZZ_p}$ (see also 4.5.11.2.1). This is the context in which we can prove directly Milne's conjecture referred to in 1.15 (for details see \S 5). So condition $(*)$ of 2.3.4 is satisfied for $z$.  
\medskip
{\bf C.} To avoid referring to \S 5, we now present an alternative (independent) proof of why $(*)$ of 2.3.4 is satisfied for $z$. We consider the Shimura filtered Lie isocrystal (definable as in 2.3.10) 
$$(LI_4[{1\over p}],\vph_4,F^0(LI_4[{1\over p}]),F^1(LI_4[{1\over p}]))$$ 
associated to $z$ (via $\tilde m$); here $LI_4:={\rm Lie}(\tilde G_{4B(k)})$. It is identifiable with a filtered Lie subisocrystal of the direct sum 
$$\bigl(LI_1[{1\over p}],\vph_1,F^0(LI_1[{1\over p}]),F^1(LI_1[{1\over p}]))\oplus (LI_2[{1\over p}],\vph_2,F^0(LI_2[{1\over p}]),F^1(LI_2[{1\over p}])\bigr),$$ 
where $(LI_i,\vph_i,F^0(LI_i),F^1(LI_i))$, with $LI_i:={\rm Lie}(G_{iW(k)})$, is the Shimura filtered Lie $\sg_k$-crystal associated to the morphism ${\rm Spec}(W(k))\to\Mn_i$ defined by $z$ (via $f_i$), $i\in\{1,2\}$.
Under this identification, using projections, $\bigl(\bigl[LI_4[{1\over p}],LI_4[{1\over p}]\bigr],\vph_4\bigr)$ can be identified with $\bigl(\bigl[LI_i[{1\over p}],LI_i[{1\over p}]\bigr],\vph_i\bigr)$, $i=\overline{1,2}$.  
We just need to show that 
$$LI_4:=LI[{1\over p}]\cap (LI_1\oplus LI_2)$$ 
is the Lie algebra of a reductive subgroup (it is $\tilde G_{4W(k)}$) of $G_{1W(k)}\times G_{2W(k)}\subset GL(M)$. [Va2, 4.3.13] takes care of the maximal torus of $Z(\tilde G_{4W(k)})$, while [Va2, 3.1.6] allows us to concentrate just on the ``der" part of $\tilde G_{4W(k)}$, and so on the ``adjoint" part of $\tilde G_{4W(k)}$. But for the adjoint context, as $G^{\rm ad}_{i\ZZ_{(p)}}$ does not depend on $i\in\{1,2,4\}$, we can appeal to 4.2.3-5. So, the Galois representation 
$$\Gamma_k\to G^{\rm ad}_{4\ZZ_p}(\ZZ_p)\hookrightarrow Aut({\rm Lie}(G^{\rm ad}_{4\ZZ_p}))(\ZZ_p)$$
defined naturally by the $W(k)$-valued point of $\Mm_1\times\Mm_2$ obtained naturally from $z$ (we recall that, cf. 2.3.7, $\Mm_i$ is the extension to $O_{(v_i)}$ of the integral canonical model of the SHS $(1_{({\rm GSp}(W_i,\psi_i),S_i)},{L_i}_{(p)},p)$), defines, (by inverting $p$ and) via Fontaine's comparison theory, a Shimura adjoint Lie isocrystal 
$$({\rm LIE},\tilde\vph):=({\rm Lie}(G_{4B(k)}^{\rm ad}),\vph_4)\subset (LI_4[{1\over p}],\vph_4)$$ 
over $k$. Moreover, using 2.3.10 for the $W(k)$-valued point of $\Mn_i$ defined by $z$ (via $f_i$), we deduce (via the above identifications) the existence of a $W(k)$-lattice ${\rm LIE}^i_W$ of ${\rm LIE}$, such that the quadruple $({\rm LIE}_W^i,\tilde\vph,F^0({\rm LIE}_W^i),F^1({\rm LIE}_W^i))$ is (naturally identifiable with) the Shimura adjoint filtered Lie $\sg_k$-crystal attached to $(LI_i,\vph_i,F^0(Li_i),F^1(LI_i))$, $i=\overline{1,2}$. Here $F^j(LIE_W^i):=LIE_W^i\cap F^j(LI_i[{1\over p}])$, $j\in\{0,1\}$ and $i\in\{1,2\}$. We have:
\medskip
{\bf D. Claim.} {\it ${\rm LIE}_W^1={\rm LIE}_W^2$.}
\medskip
{\bf Proof:} Let $B^+(W(k))$ and $\be_0$ have the same significance as in [Fa2, \S 4]: the notations to be used are as in 2.3.18.1 E (but this time $B^+(W((k))$ is defined by completing using the $PD$-topology). We have (cf. the logical passage in [Fa2, th. 7 of \S 6] to $End$'s) regardless of what odd prime $p$ we are dealing with, the following $B^+(W(k))$-monomorphisms:
$${\rm LIE}_W^i\otimes_{W(k)} \be_0B^+(W(k))\hookrightarrow {\rm Lie}(G^{\rm ad}_{4\ZZ_p})\otimes_{\ZZ_p} B^+(W(k))\hookrightarrow {\rm LIE}_W^i\otimes_{W(k)} \be^{-1}_0B^+(W(k)).\leqno (2)$$ 
Strictly speaking, the mentioned passage gives us such an inclusion at the level of Lie algebras of derived groups and not of adjoint groups: it is 3.1.8.1 (TOR) which allows us to move from the derived context to the adjoint context (if $T_G$ is a maximal torus of a semisimple group $\tilde G$ over $W(k)$, then ${\rm Lie}(\tilde G)$, as a $W(k)$-module, is the direct sum of ${\rm Lie}(T_G)$ and of a $W(k)$-submodule of ${\rm Lie}(\tilde G)$ which depends only on $G^{\rm ad}$, cf. the existence of Weyl's decompositions).
\smallskip  
So for $p\ge 5$, the Claim is a consequence of [Fa2, th. 7 of \S 6]. To see that it is also true for $p=3$ we rely on loc. cit., on 2.3.18.1 E and on 3.11.8.1. So, based on 3.11.8.1, we start with the \'etale slope type direct sum decomposition
$${\rm Lie}(G^{\rm ad}_{1\ZZ_p})={\rm Lie}(G^{\rm ad}_{2\ZZ_p})={\rm Lie}(G^{\rm ad}_{4\ZZ_p})=\oplus_{\al\in SSL} L_{\al},$$
where $SSL$ is the set of slopes of $(LIE,\tilde\vph)$, and with the slope type direct sum decompositions
$$(LIE^i_W,F^0(LIE_W^i),F^1(LIE_W^i),\tilde\vph)=\oplus_{\al\in SSL} {\got C}^i_{\al},$$
$i\in\{1,2\}$. From the adjoint form of 3.11.8.1 (ZERO) we get that $L_0\otimes_{\ZZ_p} W(\FF)$ can be naturally identified with $W(0)(LIE_W^i,\tilde\vph)$, $i\in\{1,2\}$. So we have ${\got C}^1_0={\got C}^2_0$. Similarly, as (2) is obtained from a similar $B^+(W(k))$-monomorphism associated to a $p$-divisible group over $W(k)$, by passage to the adjoint context, if $\{-1,1\}\subset SSL$ we get that $\be_0L_1$ (resp. $\be_0^{-1}L_{-1}$) is naturally identifiable with the Lie subalgebra of $LIE_W^i$ corresponding to slope $1$ (resp. $-1$) of $(LIE_W^i,\tilde\vph)$, $i\in\{1,2\}$. So ${\got C}_1^1={\got C}_1^2$ and ${\got C}_{-1}^1={\got C}_{-1}^2$. 
\smallskip
We recall from [Fa2, \S 4] that $\be_0^2$ has an image in $gr^1:=F^1(B^+(V))/F^2(B^+(W(k)))$ which is $3$ times a generator of this free $\overline{W(k)}^\wedge$-module of rank 1. So using entirely the same arguments as in 2.3.18.1 E, from (2) we get that ${\got C}^i_{\al}$ is determined by the $\Gamma_k$-module $L_{\al}$, $i\in\{1,2\}$, $\forall\al\in (0,1)\cup (-1,0)$. We conclude ${\got C}^1_{\al}={\got C}^2_{\al}$, $\forall\al\in SSL$, and so $LIE_W^1=LIE^2_W$.  
This proves the Claim.
\medskip
{\bf E.} $G_{iW(k)}$ is the connected component of the origin of the group of Lie automorphisms of $LIE_W^i$, $i\in\{1,2\}$ (cf. end of 2.2.13). From this and the Claim we get that the Zariski closure in $G_{1W(k)}^{\rm ad}\times G_{2W(k)}^{\rm ad}$ of the subgroup of $G_{1B(k)}^{\rm ad}\times G_{2B(k)}^{\rm ad}$ whose Lie algebra is ${\rm LIE}$, is an adjoint group $\tilde G^{\rm ad}_{4W(k)}$; its projections on $G_{1W(k)}^{\rm ad}$ and on $G_{2W(k)}^{\rm ad}$ are isomorphisms (for instance, cf. [Va2, 3.1.2.1 c)]). This takes care of the adjoint part of $\tilde G_{W(k)}$. The passage from adjoint groups to derived groups is trivial: as the normalization of $\tilde G^{\rm ad}_{4W(k)}$ in the field of fractions of $\tilde G^{\rm der}_{4B(k)}$ is a semisimple group, we get a natural homomorphism from it into $G_{1W(k)}^{\rm der}\times G_{2W(k)}^{\rm der}$, and again loc. cit. implies that it is a closed embedding. We conclude (based on [Va2, 3.1.6]): $\tilde G_{W(k)}$ is a reductive group. This ends the argument of why $(f_4,{L_1}_{(p)}\oplus {L_2}_{(p)},v_4)$ is a SHS.
\medskip
{\bf 4.9.2.1. Pro-\'etale covers.}
For the SHS $(f_4,{L_1}_{(p)}\oplus {L_2}_{(p)},v_4)$ we use the standard notations of 2.3.1-3, except that we put a right lower index 4 everywhere. We (naturally) get a commutative diagram of pro-\'etale covers:
$$
\def\mapright#1{\smash{
\mathop{\longrightarrow}\limits^{#1}}}
\def\mapdown#1{\Big\downarrow
\rlap{$\vcenter{\hbox{$\scriptstyle#1$}}$}}
\matrix{\Mc_4 &\mapright{m_2} &\Mc_1\cr
\mapdown{m_1}&&\mapdown{n_1}\cr
\Mc_2 &\mapright {n_2} &\Mc_0\cr}$$
(cf. the definition of $\Mc_i$, $i\in\{0,1,2,4\}$, and  [Va2, 3.2.7 4) and 6.4.5.1]).
\medskip
{\bf 4.9.3. Points.} Let $z_0\in\Mc_0(V_0)$. Let $z_i\in\Mc_i(V_0)$ be such that $n_i\circ z_i=z_0$, $i=\overline{1,2}$. Let $z_4\in\Mc_4(V_0)$ be such that $m_1\circ z_4=z_2$ and let $z_3:=m_2\circ z_4$. We denote by $y_i\in\Mc_{s_i}(\FF)$ the $\FF$-valued point defined by $z_i\in\Mc_{s_i}(V_0)$, with $i=\overline{0,4}$, with $s_i=i$ if $i\ne 3$ and with $s_3=1$. Let ${\got C}^{\rm ad}_{y_i}$ be the Shimura adjoint Lie $\bar\sg$-crystal attached to $y_i$, $i=\overline{1,4}$, obtained using the SHS $(f_{s_i},{L_{s_i}}_{(p)},v_{s_i})$, with ${L_4}_{(p)}:={L_1}_{(p)}\oplus {L_2}_{(p)}$. From the construction of $f_4$, we deduce the existence of a monomorphism ${\got C}_{y_4}^{\rm ad}\hookrightarrow {\got C}_{y_3}^{\rm ad}\oplus {\got C}_{y_2}^{\rm ad}$ between Shimura adjoint Lie $\bar\sg$-crystals, which produces (via its composite with the natural projections on the factors of ${\got C}_{y_3}^{\rm ad}\oplus {\got C}_{y_2}^{\rm ad}$; see 4.9.2.0 E) isomorphisms of ${\got C}_{y_4}^{\rm ad}$ with ${\got C}_{y_2}^{\rm ad}$ and with ${\got C}_{y_3}^{\rm ad}$.
\medskip
{\bf 4.9.4. Lemma.} {\it ${\got C}^{\rm ad}_{y_i}$, $i=\overline{1,4}$, are all isomorphic.}
\medskip
{\bf Proof:} We have just seen that ${\got C}^{\rm ad}_{y_4}$ is isomorphic to ${\got C}^{\rm ad}_{y_2}$ and ${\got C}^{\rm ad}_{y_3}$. The fact that ${\got C}^{\rm ad}_{y_1}$ and ${\got C}^{\rm ad}_{y_3}$ are isomorphic results from the fact that $z_1$ and $z_3$ are two lifts to $\Mc_1$ of $z_0$ and so $z_3=i_a(z_1)h$, with $h\in G_1(\AA^p_f)$ and with 
$$a:(G_1,X_1,H_1,v_1)\buildrel{\sim}\over\to(G_1,X_1,H_1,v_1)$$ 
an automorphism defined by an element of $G_0(\ZZ_{(p)}):=G_0(\QQ)\cap H_0$ (cf. [Va2, 3.3.1]); here as well as below we denote by $i_a$ the automorphism of $\Mn_1$ naturally defined by $a$ (cf. [Va2, 3.2.7 4)]). But the Shimura adjoint Lie $\bar\sg$-crystal attached to $y_1$ is isomorphic to the one attached to $i_a(y_1)$ (the argument of this is entirely analogous to the ones which produced isomorphisms of ${\got C}_{y_4}^{\rm ad}$ with ${\got C}_{y_2}^{\rm ad}$ and with ${\got C}_{y_3}^{\rm ad}$: we need to ``put together" the SHS $(f_1,{L_{1}}_{(p)},v_1)$ and the SHS $(f_1\circ a_1,{L_{1}}_{(p)},v_1)$, where we denote by $a_1$ the isomorphism $(G_1,X_1)\tilde\to (G_1,X_1)$ defined by $a$), and so (see Fact 6 of 2.3.11) it is isomorphic to the one attached to $y_3=i_a(y_1)h$.
\medskip
{\bf 4.9.5. Terminology.} We refer to ${\got C}^{\rm ad}_{y_1}$ as the Shimura adjoint Lie $\bar\sg$-crystal attached to $y_0$.
The above Lemma guarantees that it is well defined. For any algebraically closed field $k$ of characteristic $p$ and for every $y_{01}\in\Mc_0(k)$, we define similarly the Shimura adjoint Lie $\sg_k$-crystal attached to $y_{01}$.
\medskip
{\bf 4.9.5.1. Remark.} In all above operations, we can keep track of inner isomorphism classes as in 4.2.10: so in 4.9.4 (as well as in 4.9.7 below)
we can speak replace ``isomorphic" by inner isomorphic.
\medskip
{\bf 4.9.6. Translations.} If now $z_0^0\in\Mn_0(V_0)$ is the translation of $z_0$ by an element $t_0\in G_0(\AA^p_f)$ and we work with (the assumed to exist) $z_1^0$, $z_2^0$, $z_3^0$, $z_4^0$, and $y^0_1$, $y^0_2$, $y^0_3$, $y^0_4$ having a similar meaning, we get similarly that $z^0_3=i_{a^0}(z_1^0)h^0$, with
$h^0\in G_1(\AA^p_f)$ and with $a^0:(G_1,X_1,H_1,v_1)\buildrel{\sim}\over\to
(G_1,X_1,H_1,v_1)$ an automorphism defined by an element of $G_0(\ZZ_{(p)})$. We have:
\medskip
{\bf 4.9.7. Proposition.} {\it The Shimura adjoint Lie $\bar\sg$-crystals attached to $y^0_1$, $y^0_2$, $y^0_3$ and $y^0_4$ are all isomorphic to ${\got C}^{\rm ad}_{y_1}$.}
\medskip
{\bf Proof:} The argument that the Shimura adjoint Lie $\bar\sg$-crystals attached to $y^0_1$, $y^0_2$, $y^0_3$ and $y^0_4$ are isomorphic to each other is as in 4.9.4. The fact that the Shimura adjoint Lie $\bar\sg$-crystal attached to $y^0_1$ is isomorphic to ${\got C}^{\rm ad}_{y_1}$ can be checked, without any reference to $y^0_2$ or to $y^0_3$, as follows. Based on [Va2, 3.3.2] and Fact 6 of 2.3.11 (applied to the SHS $(f_1,{L_1}_{(p)},v_1)$) we can assume $z_0^1\in\Mc_1(V_0)$. We start working with a finite morphism $\Mn_1/H_{01}\to \Mn/H_{00}$, with $H_{01}$ as in 2.3.3 and with $H_{00}$ a compact, open subgroup
of $G_0(\AA_f^p)$ containing the image of $H_{01}$. As $G_0(\ZZ_{(p)})$ is dense in $G_0(\AA_f^p)$ (cf. the strong approximation theorem; see [Pr, Theorem A]), we can assume $t_0$ mod $H_{00}$ is naturally defined by an element $a\in G_0(\ZZ_{(p)})$. Choosing $H_{01}$ and $H_{00}$ to be small enough, we can assume $a$ leaves invariant $X_0^0$. So the part of 4.9.4 referring to automorphisms applies once more: we can assume the images of $z_0$ and $z_0^0$ in $\Mn_0/H_{00}(V_0)$ are the same. We can assume $H_{01}$ and $H_{00}$ are such that $\Mc_1/H_{010}\to\Mc_0/H_{000}$ is a Galois cover, with $H_{0i0}$ as the subgroup of $H_{0i}$ leaving $\Mc_i$ invariant, $i=\overline{0,1}$. So, as in [Va2, 6.2.2 F)] (see also AE.4) we can write $\Mc_0/H_{000}$ as the quotient of $\Mc_1/H_{010}$ through the action of a subgroup of automorphisms of $(G_1,X_1,H_1)$. So there is $a_1\in {\rm Aut}((G_1,X_1,H_1))$ such that the images of $z^0_1$ and of $i_{a_1}(z_1)$ in $\Mc_1/H_{010}$ are the same. So the Proposition follows from the proof of 4.9.4 and from Fact 6 of 2.3.11 (applied to SHS $(f_1,{L_1}_{(p)},v_1)$). 
\medskip
{\bf 4.9.7.1. Remark.} A second proof of 4.9.7 can be obtained by combining 2.3.5.6-7.
\medskip
{\bf 4.9.8. Theorem.} {\it For any integral canonical model $\Mn$ of a Shimura quadruple
$(G,X,H,v)$ of preabelian type, with $v$ diving a rational prime $p\ge 3$, there is a unique $G(\AA^p_f)$-invariant stratification (called the refined canonical Lie stratification) of $\Mn_{k(v)}$ in reduced, locally closed subschemes such that:
\medskip
a)  if $G$ is not a torus, its strata are indexed by a subset of the set $RLNP(G^{\rm ad},X^{\rm ad},v^{\rm ad})$ (defined in 4.5.7 and 4.5.8 1)); if $G$ is a torus, then we have only one stratum;
\smallskip
b) it is functorial w.r.t. finite maps between Shimura quadruples;
\smallskip
c) it has a unique open stratum $\Mu$ (so $\Mu$ is dense in $\Mn_{k(v)}$);
\smallskip
d) any point $y:{\rm Spec}(k)\to\Mu$, with $k$ a perfect field, has a uniquely determined $G$-canonical lift $z:{\rm Spec}(W(k))\to\Mn$ (here the uniqueness part is implied by e) and f) below; however, see also 4.9.17.1 below); 
\smallskip
e) the $G$-canonical lifts are functorial w.r.t. finite maps between Shimura quadruples;
\smallskip
f) if $\Mn$ and $(G,X,H,v)$ are related to a SHS $(f,L_{(p)},v)$ as in 4.1, then it is the refined canonical Lie stratification of 4.5.2 and we reobtain the $G$-canonical lifts of 4.4.2;
\smallskip 
g) it is invariant under isomorphisms, i.e. it is invariant under the natural action of the group ${\rm Aut}((G,X,H))$ (defined in [Va2, 3.2.7 9)]) on $\Mn_{k(v)}$.}
\medskip
{\bf Proof:} For the case $p=3$ we refer to \S 6. We assume now $p\ge 5$. If $G$ is a torus, then the Theorem is trivial (cf. [Va2, 3.2.8 and 3.3.3]). We assume now $G$ is not a torus. Let 
$$(G_0,X_0,H_0,v_0):=(G^{\rm ad},X^{\rm ad},H^{\rm ad},v^{\rm ad}).$$ 
We use the notations of 4.9.1-7. We recall: in order to define the set $RLNP(G_0,X_0,v_0)$, a bijection $f_{\tilde H}$ as in 4.3.2 is chosen. 
\smallskip
A point
$y_0\in\Mn(\FF)$ belongs to the stratum indexed by some $s\in RLNP(G_0,X_0,v_0)$ iff there is $t_0\in G_0(\AA_f^p)$ such that the right translate $y_0^t$ of $y_0$ through $t_0$ belongs to $\Mc_0(\FF)$ and the sequence of Newton polygons defined by cyclic factors of the Shimura adjoint Lie
$\bar\sg$-crystal attached to $y_0^t$, is the same as that of $s$; 4.9.4 and 4.9.7 guarantee that this is well defined. 
The fact that the things do not depend on the connected component $X_0^0$ of $X_0$ we chose is implied by [Va2, 3.3.3]. So, from very definitions the refined canonical stratification of $\Mn_{k(v)}$ is $G(\AA_f^p)$-invariant. 
\smallskip
To see that the partition of $\Mn_0(\FF)$ in disjoint subsets defines locally closed subschemes of $\Mn_{0k(v_0)}$ we just have to remark that the pull back (of these disjoint sets) to $\Mn_1(\FF)$ through the morphism $j_1$ (of 4.9.1) are forming the partition of $\Mn_1(\FF)$ in subsets corresponding to the refined canonical Lie stratification of ${\Mn_1}_{k(v_1)}$ defined (see 4.5.2) by the SHS $(f_1,L_{1(p)},v_1)$. So [Va2, 6.4.5.1] applies. This takes care of a) and f). 
\smallskip
The stratification of $\Mn_{k(v)}$ is obtained from the one of ${\Mn_0}_{k(v_0)}$ by the pull back operation through the natural morphism $\Mn\to\Mn_0$ (see [Va2, 3.2.7 4)]). This takes care of b), while c) results from 4.2.1, via [Va2, 6.4.5.1]. The existence and uniqueness of $G$-canonical lifts is a direct consequence of the following remark (it is just a particular case of the filtered versions of 4.9.4 and 4.9.7, cf. 4.4.1 1)):
\medskip
{\bf (4.9.8.1)} {\it If in 4.9.4 (resp. in 4.9.7) $z_2$ (resp. $z_2^0$) is a Shimura-canonical lift, then $z_1$ and $z_3$ (resp. $z_1^0$ and $z_3^0$) are also
Shimura-canonical lifts.}
\medskip
So it makes sense to say that $z_0$ is a $G_0$-canonical lift of $y_0$. For a point $y:{\rm Spec}(k)\to\Mn_{k(v)}$ of the open stratum, with $k$ just a perfect field, we first get a $G$-canonical lift $\bar z:{\rm Spec}(W(k_1))\to\Mn$ (with $k_1$ a pro-finite Galois extension of $k$) of the point $\bar y\in\Mn(k_1)$ factoring through $y$, which due to the uniqueness property asserted by c) and to (4.9.8.1), it is defined over $W(k)$, i.e. we get a $G$-canonical lift $z:{\rm Spec}(W(k))\to\Mn_{k(v)}$ of $y$. This takes care of d) and e).
\smallskip
To check that the refined canonical Lie stratification of $\Mn_{k(v)}$ is invariant under the action of ${\rm Aut}((G,X,H))$, we can assume $G=G_0$ (cf. the definition of ${\rm Aut}((G,X,H))$). But in this case the arguments of the proof of 4.9.4 apply entirely; we just need to:
\medskip
-- replace inner automorphisms $a_1$ of $(G_1,X_1)$ by isomorphisms $a_2:(G_1,X_1)\tilde\to (G_1,X_2)$ (with $X_2$ uniquely determined by $X_1$ and $a_2$), and to
\smallskip
-- point out that we can always choose $f_1$ such that 2.3.6 applies to it (cf. the Existence Property of 1.10); so 2.3.6 applies as well to $f\circ a_2^{-1}$. 
\medskip
This takes care of g) and ends the proof.
\medskip
{\bf 4.9.8.2. Corollary.} {\it We assume $\dim_{\CC}(X)>0$. Then the complement in $\Mn_{k(v)}$ of the open, dense stratum of the refined canonical Lie stratification is a $G(\AA_f^p)$-invariant closed subscheme of $\Mn_{k(v)}$, different from $\Mn_{k(v)}$ and of pure codimension 1.}
\medskip
{\bf Proof:} [Va2, 6.4.2 and 6.4.5.1] makes the reduction (cf. the proof of 4.9.8) of the situation to the case when we have a SHS. But this case is handled by 4.3.6. This ends the proof of the Corollary.
\medskip
This Corollary fulfills the promise of the end of the first paragraph of [Va2, 3.2.12].  
\medskip
{\bf 4.9.9. Variants.} A result entirely similar to 4.9.8 can be stated in the context of:
\medskip
i) canonical Lie stratifications (defined for a SHS in 4.5.1);
\smallskip
ii) (refined) Lie (non-) stable stratifications (defined for a SHS in 4.5.15);
\smallskip
iii) quasi-ultra stratifications (defined for a SHS in 4.5.15.2);
\smallskip
iv) Faltings--Shimura--Dieudonn\'e adjoint stratifications (defined for a SHS in 4.5.16). 
\medskip
The same proof applies (it is the simple Lemma 4.10.0 below which allows us to speak about the quasi-ultra stratification of $\Mn_{k(v)}$ and not only of $\Mn_{\FF}$). In \S10 we will see that this remains true, even in the context of the absolute and so in the context of the ultra stratifications: using 4.5.6.1 and 4.9.8 we get this just for the $C_n$ types (and for other particular situations pertaining to the $B_n$ types), $n\in\NN$. Here we just mention, that the argument (of the proof of 4.9.4) used to get 4.9.8 g), applies entirely to give us that the $\rho$'s stratifications (and so also the absolute, the pseudo-ultra, the ultra and the Faltings--Shimura--Dieudonn\'e adjoint or standard stratification) of 4.5 are ${\rm Aut}((G,X,H))$-invariant. 
\smallskip
Also 4.5.16 a) to c) and 4.5.16.1 extend automatically to all these Faltings--Shimura--Dieudonn\'e adjoint stratifications. 
\medskip
{\bf 4.9.9.1. Ultimate type of stratifications.} For what follows we refer to (the proof of) 4.9.8; see also 2.3.3.2 and 2.3.5.4. We consider the Faltings--Shimura--Dieudonn\'e adjoint stratification of $\Mn_{0k(v_0)}$. Its pull back to $\Mn_1$ is the Faltings--Shimura--Dieudonn\'e adjoint stratification as defined in 4.5.16. By the ultimate adjoint stratification of $\Mn_{k(v)}$ we mean the pull back to $\Mn_{k(v)}$ of the maximal (i.e. the most refined) $G_0(\AA_f^p)$-invariant, ${\rm Aut}(G_0,X_0,H_0)$-invariant stratification of $\Mn_{0k(v_0)}$ whose strata are open, closed subschemes of the Faltings--Shimura--Dieudonn\'e adjoint stratification of $\Mn_{0k(v_0)}$.
\smallskip
We now refer to 4.5.16. By the ultimate stratification of $\Mn_{k(v)}$ we mean the most refined $G(\AA_f^p)$-invariant, ${\rm Aut}(G,X,H)$-invariant stratification of $\Mn_{k(v)}$ whose strata are open, closed subschemes of the Faltings--Shimura--Dieudonn\'e stratification of $\Mn_{k(v)}$. 
\smallskip
Based on 4.5.16 and 4.5.16.1 and the last sentence of 4.9.9, 4.5.16 a) to c) and 4.5.16.1 extend automatically to all these ultimate (adjoint) stratifications.
\medskip
{\bf 4.9.9.2. Grothendieck's specialization category.} We refer to 4.5.16 (resp. to 4.9.8). By Grothendieck's specialization category (resp. adjoint category)
$$\GG\SS\CC(\Mn_{k(v)})\,\,({\rm resp.}\, \GG\SS\AA\CC(\Mn_{k(v)}))$$
of $\Mn_{k(v)}$ we mean the category:
\medskip
{\bf 1)} whose objects are the reduced, locally closed subschemes of $\Mn_{k(v)}$ which are strata of the ultimate (resp. ultimate adjoint) stratification of $\Mn_{k(v)}$;
\smallskip
{\bf 2)} whose set of morphisms from one object (stratum) $s_1$ to another object (stratum $s_2$) is empty or has precisely one element depending on the fact that each (it is enough one) connected component of $s_1$ specializes to some connected component of $s_2$. 
\medskip
$\GG\SS\AA\CC(\Mn_{k(v)})$ depends only on the integral canonical model of $(G^{\rm ad},X^{\rm ad},H^{\rm ad},v^{\rm ad})$. However, there are many very important ways to look at it (and its subcategories) and so it is relevant to think of it of being ``of $\Mn_{k(v)}$". We have a canonical functor
$$\GG\SS\CC(\Mn_{k(v)})\to\GG\SS\AA\CC(\Mn_{k(v)})$$
which takes a stratum of the ultimate stratification to the corresponding stratum of the ultimate stratification of $\Mn_{k(v)}$. We do not know when it is an isomorphism of categories. 
\smallskip
We have a variant of $\GG\SS\CC(\Mn_{k(v)})$ (resp. of $\GG\SS\AA\CC(\Mn_{k(v)})$): we allow the strata to be in the sense of 2.1; so the objects are elements of the class of the ultimate (resp. ultimate adjoint) stratification of $\Mn_{k(v)}$. Another variant, we fix an algebraically closed field $k$ containing $\FF$ and we define $\GG\SS\CC(\Mn_k)$ (resp. $\GG\SS\AA\CC(\Mn_k)$) using reduced, locally closed subschemes of $\Mn_k$ which are strata of the ultimate (resp. ultimate adjoint) stratification of $\Mn_{k(v)}$.  
\medskip
{\bf 4.9.9.3. Remark.} One might worry that the ultimate adjoint stratifications are too refined, i.e. that too many strata of them have very low dimensions. We think this is not relevant, as one of the main goals of them is to solve the Real Problem of 1.6.5; so the categories of 4.9.9.2 do ``capture" all the information one might think is lost by passing to too refined stratifications (providing --in the same time-- a lot more). On the other hand, from the Example of 3.15.7 H we get (cf. also the regularity part of 4.5.15.6; the passage from SHS's to the preabelian type context is achieved via the proof of 4.9.8):
\medskip
{\bf Corollary.} {\it The dimension of the stratum of the ultra stratification of $\Mn_{k(v)}$ to which an $\FF$-valued toric point $y$ of $\Mn_{k(v)}$ belongs is at least the negative $sp$-invariant of the Shimura adjoint Lie $\bar\sg$-crystal attached to (see 4.9.5 and 4.9.17 below) $y$.}
\medskip
{\bf 4.9.10. Definition.} The points with values in perfect fields of the open, dense stratum of $\Mn_{k(v)}/\tilde H_0$ (with $\tilde H_0$ a compact subgroup of $G(\AA_f^p)$) are called $G$-ordinary (or $G_{W(k(v))}$-ordinary or Shimura-ordinary) points.
\medskip
{\bf 4.9.11. Corollary.} {\it We have a natural functor from the full subcategory of qf-Sh of whose objects are Shimura quadruples $(G,X,H,v)$ of preabelian type, with $(v,2)=1$,
to the category of stratified schemes (in reduced, locally closed subschemes) whose morphisms are morphisms of schemes such that the pull back of any stratum is a stratum, taking such a quadruple into the scheme defined by the special fibre of its integral canonical model and endowed with the
refined  canonical Lie stratification (the action of the functor on morphisms being the logical one, cf. [Va2, 3.2.7 4)]).}
\medskip
This is a consequence of the proof of 4.9.8. We have variants of it:  cf. 4.9.9. 
\medskip
{\bf 4.9.12. Corollary.} {\it Under the right translation by an element of $G(\AA^p_f)$, a $G$-canonical lift of $\Mn$ is taken into a $G$-canonical lift.}
\medskip
{\bf 4.9.13. Remark.} Let $\tilde H_0$ be an arbitrary compact subgroup of $G(\AA_f^p)$. As in 4.4.2.1 we define uniquely the $G$-canonical lift ${\rm Spec}(W(k))\to\Mn/\tilde H_0$ of a $G$-ordinary point ${\rm Spec}(k)\to\Mn_{k(v)}/\tilde H_0$.  
 In \S 12 we will see (cf. 4.4.6) that a $G$-canonical lift ${\rm Spec}\bigl(W(\FF)\bigr)\to\Mn/\tilde H_0$ gives birth to a special point ${\rm Spec}(B(\FF))\to {\rm Sh}_{\tilde H_0}(G,X)$.
\medskip
{\bf 4.9.14. Refined canonical stratifications.} We start with a question. With the notations of 4.5, is it true that the $f$-canonical stratification of $\Mn_{k(v)}$ defined in 4.5.9 is indeed canonical, i.e. it does not depend on the SHS $(f,L_{(p)},v)$ producing it?
\smallskip
The expectation (cf. 4.5.9) of the existence of a refined $f$-canonical stratification of $\Mn_{k(v)}$, suggests (in general) a negative answer to this question. It seems to us, that the supposed existing refined $f$-canonical stratification of $\Mn_{k(v)}$ (cf. 4.5.9), still depends in general on the SHS producing it (and so that we can not stop mentioning ``$f$-"). However, one could proceed as in 4.9.8 to define a refined canonical stratification of $\Mn_{k(v)}$ which does not depend on $f$ (sort of intersection of all refined $f$-canonical stratifications of $\Mn_{k(v)}$); using this one could define the refined canonical stratification of the special fibre of any integral canonical model as in 4.9.8.
\medskip
{\bf 4.9.15. An expectation.} Even before one truly defines the refined canonical stratifications modeled on [Oo3] (i.e. in the sense of 4.5.9 and 4.9.14), we formulate here our feelings:
\medskip
{\bf Expectation.} {\it The quasi-ultra stratification of $\Mn_{k(v)}$ coincides with or is a refinement of the refined canonical stratification of $\Mn_{k(v)}$, whenever this last stratification can be satisfactorily defined.}
\medskip
For instance, it is easy to see that in the case of Siegel modular varieties they coincide (cf. [Oo3] and 1) of 3.13.7.4 D)).
\medskip
{\bf 4.9.16. Exercise.} {\bf a)} Show that the theory of good $G$-multiplicative (or $G$-additive) coordinates presented in 4.7.14 extends to the case of Shimura-ordinary points of an integral canonical model $\Mn$ as in 4.9.8. 
\smallskip
{\bf b)} If the prime $v^{\rm ad}$ of $E(G^{\rm ad},X^{\rm ad})$ divided by $v$ has $\FF_p$ as its residue field, then the $G$-multiplicative coordinates we get are canonical (i.e. are uniquely determined). 
\medskip
Hint: for the canonical part, use 4.7.11 4) and the idea of the proof of 4.9.4.
\medskip
{\bf 4.9.17. Some extensions.} Let $(G,X,H,v)$, $p$ and $\Mn$ be as in 4.9.8. Let $y\in\Mn(k)$, with $k$ an algebraically closed field of characteristic $p$. Till the end of 4.9.17.5, we assume $G$ is not a torus. Let $\Mn_0$ be the integral canonical model of the quadruple $(G_0,X_0,H_0,v_0):=(G^{\rm ad},X^{\rm ad},H^{\rm ad},v^{\rm ad})$ and let $y_0\in\Mn_0(k)$ be the composite of $y$ with the natural morphism $\Mn\to\Mn_0$. Let $\Mc_0$ be the open closed subscheme of $\Mn_0$ defined as in 4.9.1.1, starting from a fixed connected component of $X$. Let $h\in G^{\rm ad}(\AA_f^p)$ be such that the translation $y_{01}$ of $y_0$ by $h$ belongs to $\Mc_0(k)$ (cf. [Va2, 3.3.2]). The Shimura adjoint Lie $\sg_k$-crystal ${\got C}_{y_{01}}^{\rm ad}$ (cf. 4.9.5) is called the Shimura adjoint Lie $\sg_k$-crystal attached to $y$. It is attached to a Shimura $\sg_k$-crystal.
\smallskip
The reason why is this well defined is nothing else but 4.9.4, 4.9.7 and the proof of 4.9.8 (performed for $W(k)$ instead of $V_0$). In the case of a SHS $(f,L_{(p)},v)$ we recover the definition (see 2.3.10) of a Shimura adjoint Lie $\sg_k$-crystal attached to a point $y\in\Mn(k)$.
\medskip
{\bf 4.9.17.0. The quasi-affine, regular property.} We have:
\medskip
{\bf Corollary.} {\it Each stratum $\Ms_0$ of the quasi-ultra (or ultimate adjoint, or Faltings--Shimura--Dieudonn\'e adjoint) stratification of $\Mn_{k(v)}$ is regular and quasi-affine. Moreover, any stratum of the ultra stratification of $\Mn_{k(v)}$ is quasi-affine.}
\medskip
{\bf Proof:} The regularity part follows from 4.5.15.2.2, 4.5.16 and the proof of 4.9.8 (cf. also 4.9.9). The quasi-affineness part follows from 4.5.15.2.5 and the proof of 4.9.8 (cf. also 4.9.9). 
\medskip
{\bf 4.9.17.0.0. Remark.} Based on global forms of 2.3.9 and on 2.3.15.1 (or just based on the existence of actions $\TT$ as in 3.13.7.1; their orbits are always integral subschemes of reductive groups and so are quasi-affine), in \S 5 we will reobtain the quasi-affineness part of the above Corollary, without appealing (see 4.5.9) to [Oo3]. 
\medskip
{\bf 4.9.17.1. Remark.} In the situation of 4.9.17, we similarly define the Shimura adjoint filtered Lie $\sg_k$-crystal attached to a point ${\rm Spec}(W(k))\to\Mn$. So the canonical lifts ${\rm Spec}(W(k))\to\Mn$ can be recognized (see 4.4.1 2)) as follows: their Shimura adjoint filtered Lie $\sg_k$-crystals are (cyclic diagonalizable and) of parabolic (and Borel) type.
\medskip
{\bf 4.9.17.2. Definitions.} A point $y:{\rm Spec}(k_1)\to\Mn_{k(v)}$, with $k_1$ a perfect field, is called a toric point if the Shimura adjoint Lie $\sg_{\overline{k_1}}$-crystal attached to the $\overline{k_1}$-valued point of $\Mn_{k(v)}$ naturally defined by $y$ is isomorphic to the extension to $\overline{k_1}$ of the Shimura adjoint Lie $\bar\sg$-crystal attached to a Shimura $\bar\sg$-crystal ${\got C}_{\om}$ showing up in 4.1.5, for an adequately chosen SHS $(f,L_{(p)},v)$. If ${\got C}_{\om}$ is a $U$-ordinary (resp. $T$-ordinary) $\bar\sg$-crystal, then $y$ is called an $U$-ordinary (resp. $T$-ordinary) point. Similarly to 4.4.0 we define Borel and reductive points of $\Mn_{k(v)}$.
\medskip
{\bf 4.9.17.2.1. Some implications.} 3.11.6 B) implies that any $k_1$-valued Borel point is a Shimura-ordinary point and the converse holds if $k_1=\overline{k_1}$. Also 2.2.19.2 and the Criterion of 2.2.22 1) imply that any toric point is a reductive point.  
\medskip
{\bf 4.9.17.3. Extension of 4.5.11 and 4.5.13.} 4.9.4 and 4.9.7 show that the notion of toric point is intrinsic (canonical). In particular, if $(G,X,H,v)$ is obtained from a SHS $(f,L_{(p)},v)$, then we ``regain" the toric points defined in 4.5.11.2. Moreover, 4.12.12.6 below shows that there are plenty of toric points which are not $G$-ordinary. As in 4.5.11.2.1, the toric points can be defined in terms of Galois representations (for $p=3$ cf. also the Criterion of 2.2.22 1)). 
\smallskip
Moreover, due to the same reasons, we can speak about CM levels of non-toric points of $\Mn_{k(v)}$ with values in algebraically closed fields and so, by natural extension, with values in fields.
\medskip
{\bf 4.9.17.4. Remark.} Based on 4.9.17, it makes sense to speak, as in 3.6.16 1), about $\Mn$ having the completion property: we just have to rephrase 3.6.15 A in terms of Shimura adjoint Lie $\sg_k$-crystals attached to geometric points of $\Mn$.  
\medskip
{\bf 4.9.17.5. $U$-canonical lifts.} From 4.4.13.2 we get, as in the part of the proof of 4.9.8 referring to $G$-canonical lifts, that any $U$-ordinary point $y:{\rm Spec}(k_1)\to\Mn_{k(v)}$, with $k_1$ a perfect field, has a uniquely determined lift $z:{\rm Spec}(W(k_1))\to\Mn$ (called the $U$-canonical lift of $y$) with the property that the Shimura adjoint filtered Lie $\sg_{\overline{k_1}}$-crystal attached to the $W(\overline{k_1})$-valued point of $\Mn$ factoring through $z$, is of toric type. As in 4.9.8 b), $U$-canonical lifts are functorial w.r.t. finite maps between Shimura quadruples.
\medskip
{\bf 4.9.18. Comment.}
We think the terminology (refined) canonical Lie stratification is justified. The word Lie is meant to distinguish it from the other $\rho$-stratifications introduced in 4.5.4 (for a SHS), while the word canonical is meant to emphasize that in the case of a Shimura quadruple $(G,X,H,v)$ emerging from a SHS $(f,L_{(p)},v)$, the (refined) (canonical) Lie stratification we get for $\Mn_{k(v)}$ does not depend on the SHS $(f,L_{(p)},v)$ producing it. 
\medskip
{\bf 4.9.19. Exercise.} Let $(G_i,X_i,H_i,v_i)$, $i=\overline{1,2}$, be two Shimura quadruples of preabelian type, with $v_1$ and $v_2$ dividing the same prime $p\ge 3$. Let $w$ be a prime of $E(G_1\times G_2,X_1\times X_2)$ dividing $v_1$ and $v_2$. We get a Shimura quadruple 
$$(G_1\times G_2,X_1\times X_2,H_1\times H_2, w).$$ 
Let $\Mm_i$ be the extension to ${\rm Spec}(O_{(w)})$ of the integral canonical model $\Mn_i$ of $(G_i,X_i,H_i,v_i)$, $i=\overline{1,2}$. It is known that 
$$\Mm:=\Mm_1\times_{O_{(w)}}\Mm_2$$ 
is the integral canonical model of $(G_1\times G_2,X_1\times X_2,H_1\times H_2,w)$ (cf. [Va2, 3.2.16]). Show that the refined canonical Lie stratification of $\Mm_{k(w)}$ is the product stratification of the refined canonical Lie stratifications of $\Mm_{1k(v)}$ and of $\Mm_{2k(v)}$. Hint: first reduce to the case when $G_1\times G_2$ is an adjoint group and then use [Va2, 6.5.1].
\medskip
{\bf 4.9.20. Exercise.} With the notations of 4.9.17, let $H_0$ be a compact subgroup of $G(\AA_f^p)$ such that either $p\not | t(G^{\rm ad})$ and $H_0\times H$ is a subgroup of $G(\AA_f)$ smooth for $(G,X)$ or $p| t(G^{\rm ad})$ and $H_0\times H$ is $p$-smooth for $(G,X)$. Show that for any morphism $f_s:{\rm Spec}(R)\to\Mn/H_0$, with $R$ a faithfully flat, regular, formally smooth $W(k)$-algebra, we can define uniquely the Shimura adjoint filtered Lie $F$-crystal (over ${\rm Spec}(R/pR)$) attached to $f_s$.
\smallskip
A proof of this exercise will be given in \S 5. In what follows we need it only in the case when $R/pR$ is integral and the first fundamental group of ${\rm Spec}(R/pR)$ is trivial: this case can be treated entirely as in (the filtered versions of) 4.9.3-4 and 4.9.7 (where we had $R/pR=\FF$).
\medskip
{\bf 4.9.21. The relative Oort--Moonen problem.} Let $(G,X,H,v)$ be a Shimura quadruple of preabelian type, with $v$ dividing a prime $p>2$. Let $\Mn$ be its integral canonical model (cf. [Va2, 6.4.1] for $p>3$ and cf. \S 6 for $p=3$). Let $F$ be a number field containing $E(G,X)$. Let $H_0\subset G(\AA_f^p)$ be as in 4.9.20. Let $Z$ be an irreducible subvariety of ${\rm Sh}_{H_0\times H}(G,X)_F$. Let $Z_{(v)}$ be the Zariski closure of $Z$ in ${(\Mn/H_0)}_{F_{(v)}}$, with $F_{(v)}$ as the normalization of $O_{(v)}$ in $F$. 
\medskip
{\bf Problem.} We assume there is a prime $w$ of $F$ dividing $v$ and a $G$-ordinary point $y:{\rm Spec}(\overline{k(v)})\to\Mn_{k(w)}/H_0$ factoring through the special fibre of $Z_{(v)}$, such that $Z_{(v)}$ has (see below) formally quasi-linear components at $y$. Show that $Z$ is of preabelian type, i.e. for any embedding $F\hookrightarrow\CC$, every irreducible component of $Z\times_F\CC$ is a subvariety of $\Mn_{\CC}/H_0$ of preabelian type (in the sense of [Mo, 3.7]). 
\medskip
Due to the existence of $G$-multiplicative coordinates (see 4.9.16) the notion of having formally quasi-linear components at $y$ can be defined as in [Mo]. It seems to us that we can use 4.7.11 and 4.7.17, to prove that this notion is independent of all choices referred to in 4.7.11 8); if this turns out not to be the case, then we have to define this notion by working w.r.t. all --allowed, i.e. logical-- choices mentioned in 4.7.11 8).
\smallskip
 For not being too long, we simplify the things by just giving an example (without being particular about the precise choices defining $G$-multiplicative coordinates) in the case when $w$ is unramified over $v$. If the closed subscheme (defined naturally by $Z_{(v)}$) of the spectrum of the completion of the local ring of the point $y_1:{\rm Spec}(\overline{k(v)})\to\Mn_{W(\overline{k(v)})}/H_0$ defined by $y$, is a formal subtorus of the formal torus (of $G$-multiplicative coordinates) defined (see 4.9.16) by the $G$-ordinary point $y_1$, then (by very definitions) $Z_{(v)}$ has formally quasi-linear components at $y$. 
\medskip
{\bf 4.9.21.1. Exercise.} Using [Mo, 5.2] show that the relative Oort--Moonen problem is true if $k(v^{\rm ad})=\FF_p$ (here $v^{\rm ad}$ is the prime of $E(G^{\rm ad},X^{\rm ad})$ divided by $v$). Hint: use 4.7.11 4) and [Va2, 6.5.1.1] (referring to the Existence Property of 1.10, if $k(v_0)=\FF_p$, then based on arguments of loc. cit. we can assume $k(v^0)=\FF_p$).
\medskip
{\bf 4.9.22. Definition.} A simple, adjoint Shimura variety ${\rm Sh}(G,X)$ is said to be of $\ZZ$-End type if there is a principally polarized abelian variety $(A,p_A)$ over $\CC$, with ${\rm End}(A)=\ZZ$, whose attached (see 2.1) Shimura variety has ${\rm Sh}(G,X)$ as its adjoint variety.
\medskip
We refer to the situation of 4.9.2. Let $E$ be a number field. Let $H_{00}$ be the subgroup of $G_0(\AA_f^p)$ generated by the images of $H_{01}$ and $H_{02}$. We recall, cf. 2.3.3, that $H_{0i}$ is a compact, open subgroup of $G_i(\AA_f^p)$ subject to some conditions, $i\in\{1,2\}$. Accordingly, we assume $H_{00}$ is compact as well. We consider two morphisms 
$$a(i):{\rm Spec}(E)\to {\rm Sh}_{H_{0i}\times H_i}(G_i,X_i),$$ 
$i\in\{1,2\}$, giving birth to the same morphism ${\rm Spec}(E)\to {\rm Sh}_{H_{00}\times H_0}(G_0,X_0)$. Let $A_i$ be the abelian variety over $E$ naturally associated to $a(i)$ and the embedding $f_i$ of 4.9.1 (implicitly, a $\ZZ$-lattice $L_i$ of $W_i$ such that we get a perfect form $\psi_i:L_i\otimes_{\ZZ} L_i\to\ZZ$ is a priori chosen; see 2.3.2). We have:
\medskip
{\bf 4.9.23. Corollary.} {\it The ordinary conjecture is true for $A_1$ iff it is true for $A_2$.} 
\medskip
{\bf Proof:}
We can assume $E$ is big enough so that $G_{iE}^{\rm ab}$, $i=\overline{1,2}$, are split tori. We consider only those primes of $E$ such that:
\medskip
a) they are unramified over rational primes;
\smallskip
b) their residue fields have a prime number of elements;
\smallskip
c) $A_1\times A_2$ has good reduction w.r.t. them. 
\medskip
The set of such primes has Dirichlet density 1. Let $w$ be such a prime. Above, the role of $p$ is just to fix some notations (without introducing others); so, based on [Va2, 5.8.6] (to be compared with the proof of 4.6.2.1) we deduce that we can assume $k(w)=\FF_p$. From a) and c) we deduce the existence of a $O_{(v_i)}$-morphism $b(i):{\rm Spec}(W(\FF))\to\Mn_i$ whose generic fibre factors through $a(i)$. From 4.9.4 we get that the Shimura adjoint Lie $\bar\sg$-crystal attached to $b(1)$ is ordinary iff the Shimura adjoint Lie $\bar\sg$-crystal attached to $b(2)$ is so. From b) and 4.3.1.1 we get that at the level of $T$-degrees of definition we have $d_T(v_1)=d_T(v_2)=1$. So from end of 4.3.1.1 and 4.6 P1 we get: $A_1$ has an ordinary reduction w.r.t. $w$ iff $A_2$ has. The Corollary follows.   
\medskip
The ordinary reduction is true for a (finite) product of abelian varieties over a number field iff it is true for each member of it. Based on this and on [Va2, 6.5.1], from 4.9.23 and [Pi, 7.2] we get:
\medskip
{\bf 4.9.23.1. Corollary.} {\it Let $A$ be an abelian variety over a number field $E$. We assume the Shimura variety attached to it is such that all simple factors of its adjoint variety are of $\ZZ$-End type and are not of $C_n$ type with $n\ge 3$. Then the ordinary conjecture is true for $A$.}
\medskip
{\bf 4.9.23.2. Exercise.} Show that there are situations to which 4.9.23.1 applies, with $A$ an absolutely simple abelian variety such that ${\rm End}(A_{\bar E})\neq\ZZ$. Hint: use Mumford curves and 2.3.5.1; for more sophisticated examples, use [Va2, 6.4.2 and 6.5-6] and move from the $A_1$ Lie type to the $B_n$ Lie types, with $n\in\NN$ congruent to $1$ or $2$ mod $4$.    
\medskip\smallskip
{\bf 4.10. Hecke orbits of $G$-ordinary points.}
Let $(G,H,H,v)$ be a Shimura quadruple of preabelian type, with $(v,2)=1$, and let $\Mn$ be its integral canonical model over $O_{(v)}$. Let $\Mu\hookrightarrow\Mn_{k(v)}$ be the open, dense subscheme of $G$-ordinary points (cf. 4.9.8). We have:
\medskip
{\bf 4.10.0. Lemma.} {\it The field of fractions $FF(\Mc_0)$ of any connected component $\Mc_0$ of $\Mn_{k(v)}$ contains $\FF$.}
\medskip
{\bf Proof:} For $p=3$ we refer to \S 6. Let now $p\ge 5$. From 2.3.9 A we get that the Lemma holds for a SHS. We can assume $(G,H,H,v)$ is adjoint. So, from [Va2, 6.4.1-2, 6.2.2, 6.2.3.1 and 6.2.3] we get that $FF(\Mc_0)$ contains a subfield of $\FF$ whose Galois group is an $M$-torsion group for some $M\in\NN$. As such a subfield must be $\FF$ itself, the Lemma follows. 
\medskip
{\bf 4.10.1. Conjecture.} $\Mu$ {\it is the smallest non-empty open subscheme of $\Mn_{k(v)}$ which is $G(\AA^p_f)$-invariant.}
\medskip
{\bf 4.10.1.1. Reformulation.} Standard arguments (based on 4.10.0) show that another way to formulate 4.10.1 is: the Hecke orbit $o(y)$ of a $\bar k$-valued $G$-ordinary point $y$ (of $\Mu_{\bar k}$) with values in $\bar k$, i.e. the set of points obtained from $y$ through $G(\AA^p_f)$-translations, is dense
in $\Mn_{\bar k}$. In other words, $o(y)$ is never contained in a closed subscheme of $\Mn_{\bar k}$ which is not $\Mn_{\bar k}$ itself.
\medskip
{\bf 4.10.2. Motivations.} The expectations of the validity of 4.10.1 are based on:
\medskip
{\bf a)} 4.2.8.1;\par\noindent
\smallskip
{\bf b)} exercise 8) of 4.5.6;\par\noindent
\smallskip
{\bf c)} the part d) (or c)) of 4.4.1 3), and b) and c) of 4.4.1 2);\par\noindent
\smallskip
{\bf d)} the validity of 4.10.1 in the case when ${\rm Sh}(G,X)$ is a Siegel modular
variety [Ch2];\par\noindent
\smallskip
{\bf e)} 4.10.3 and 4.10.5.1 below.
\medskip
{\bf 4.10.3. Proposition.} {\it The above conjecture is true if $\dim_{\CC}(X)\le 1$.}
\medskip
{\bf Proof:} If $\dim_{\CC}(X)=0$, i.e. if $G$ is a torus, this is obvious (cf. [Va2, 3.2.8 and 3.3.1]). If $\dim_{\CC}(X)=1$, then 4.2.8.1 holds: the case of non-compact Shimura curves, i.e. of Shimura curves having the same adjoint curve as the elliptic modular curve, is a consequence of [FC, ch. IV, 5.10 and 6.12]. From 4.9.12 and from the fact the $G(\AA^p_f)$-orbit of the $B(\bar k)$-valued point of ${\rm Sh}(G,X)_{B(\bar k)}$ defined by the $G$-canonical lift of $y$ is dense in ${\rm Sh}(G,X)_{B(\bar k)}$, we get that $o(y)$ is dense in $\Mn_{\bar k}$. This ends the proof. 
\medskip
{\bf 4.10.4. Variants.} There are variants of 4.10.1. For instance by working not with the full group $G(\AA_f^p)$ of Hecke operators but just with $G(\QQ_q)$, where $q$ is a rational prime, different from $p$ and such that $G$ is unramified over $\QQ_q$ (cf. [Ch2]), or more generally with an adequate (i.e. of interest) subgroup of $G(\AA_f^p)$. However in such a variant we have to be careful: $G(\QQ_q)$ might not permute transitively the connected components of ${\rm Sh}_H(G,X)_{\CC}$ (for instance, this happens when $G_{\QQ_p}$ is anisotropic). 
\smallskip
Another variant consists in replacing $G$-ordinary points and $\Mu$ respectively by $G(\om)$-ordinary points and by the ${\got s}_{\om}$ (or by the ${\got s}^a_{\om}$) stratum (of 4.5.11.1). Again some precautions have to be taken: (theoretically) we might not get dense Hecke orbits. So in some of these variants 4.10.1 needs to be reformulated.  
\medskip
{\bf 4.10.5. Remark.} The gluing principle of 3.6.19 D can be interpreted as a local Hecke correspondence. Using this we will prove in [Va7] the following IOTA (if one then all) criterion:
\medskip
{\bf 4.10.5.1. Criterion.} {\it Let $(G,X,H,v)$ be a quadruple of preabelian type, with $k(v^{\rm ad})=\FF_p$. We assume that either $p\ge 3$ or $p=2$ and the Shimura quadruple $(G,X,H,v)$ is associated to a $p=2$ SHS. Let $\Mu$ be the $G$-ordinary locus of the integral canonical model $\Mn$ of $(G,X,H,v)$ (for $p=2$ cf. 4.14.3 A below). If the $G(\AA_f^p)$-orbit of a (one given) $\FF$-valued $G$-ordinary point of $\Mu$ is dense in $\Mn_{k(v)}$ (i.e. it is dense in $\Mu)$, then the $G(\AA_f^p)$-orbit of any $G$-ordinary point of $\Mu$ is dense in $\Mn_{k(v)}$.} 
\medskip
This criterion recovers all previous known cases of density of Hecke orbits in positive characteristic and provides plenty of new examples of such densities (like the case of classical Spin modular varieties, see [Va2, 5.7.5] and [Va4]).
The philosophy behind it can be roughly formulated as follows:
\medskip
{\bf Ph.} {\it We refer to 4.7.11 7) and we assume $S(1,n)^0$ is the empty set. Under the mentioned interpretation, in the $\NN$-pro-\'etale topology of ${\rm Spec}(R_n^{\rm al}/pR_n^{\rm al})$, the Hecke correspondences associated to ${\got C}_q^{\rm al}$ and two distinct $k$-valued points $y_1$ and $y_2$ of ${\rm Spec}(R_n^{\rm al}/pR_n^{\rm al})$, are ``very close" to the translation by the $k$-valued point ``$y_2-y_1$" of ${\rm Spec}(R_n^{\rm al}/pR_n^{\rm al})$ (w.r.t. its natural group structure as a torus). Here ``very close" refers to the fact that we need to change the initial Frobenius lift of (the $p$-adic completion of) $R_n^{\rm al}$ in order to get exactly the mentioned translation.}
\medskip\smallskip
{\bf 4.11. The functorial behavior of $G$-ordinary points.}
\medskip
{\bf 4.11.1. Two questions.} We consider now a SHS $(f,L_{(p)},v)$ defined by an injective map
$f:(G,X)\hookrightarrow\bigl({\rm GSp}(W,\psi),S\bigr)$ which factors as the composite $f_1\circ f_2$ of two injective maps, with $f_2:(G,X)\hookrightarrow(G_1,X_1)$ and $f_1:(G_1,X_1)\hookrightarrow\bigl({\rm GSp}(W,\psi),S\bigr)$. Let $v_1$ be the prime of $E(G_1,X_1)$ divided by $v$. We assume the triple $(f_1,L_{(p)},v_1)$ is as well a SHS. We still use for it the standard notations of 2.3.1-3 but everything having a right lower index 1, except $L_{(p)}$. We have a natural morphism $i_1:\Mn\to\Mn_1$ (cf. [Va2, 3.2.7 4)]). Two natural questions arise:
\medskip
{\bf Q1.} {\it When $i_1$ takes $G$-ordinary points of $\Mn_{k(v)}$ into $G_1$-ordinary points of ${\Mn_1}_{k(v_1)}$?}
\medskip
{\bf Q2.} {\it If $i_1$ takes  $G$-ordinary points of $\Mn_{k(v)}$ into $G_1$-ordinary points of ${\Mn_1}_{k(v_1)}$, is it true that a $G$-canonical lift of a $G$-ordinary point $y$ of $\Mn_{k(v)}$ is mapped by $i_1$ into a $G$-canonical lift of $i_1\circ y$?}
\medskip
$Q1$ can be reformulated (cf. 4.2.1): when the open, dense stratum of $\Mn_{k(v)}$ is the pull back through $i_1$ of the open, dense stratum of ${\Mn_1}_{k(v_1)}$?
\medskip
The case $k(v_1)=\FF_p$ is completely handled by 4.6 P1 and P2: the answer to $Q1$ is yes when $k(v)=\FF_p$, and then $Q2$ also has a positive answer. What about $k(v_1)=\FF_{p^{q_1}}$
with $q_1\in\NN$, $q_1\ge 2$? As $p\ge 3$, by combining 4.4.8 2), 2.3.17 and 3.1.1 c) we get:
\medskip
{\bf 4.11.1.1. Corollary.} {\it The answer to Q2 is always yes.}
\medskip
Related to Q1 we have:
\medskip
{\bf 4.11.2. Theorem.} {\it The morphism $i_1$ maps $G$-ordinary points of $\Mn_{k(v)}$ into $G_1$-ordinary points of $\Mn_{1k(v_1)}$ (and so $G$-canonical lifts of $\Mn$ into $G_1$-canonical lifts of $\Mn_1$) iff there is a triple $(T,\mu,B_1)$, with $T$ a torus of $G_{\ZZ_p}$, with $\mu$ an injective cocharacter of $T_{W(k(v))}$ and with $B_1$ a parabolic subgroup of ${G_1}_{\ZZ_p}$, such that:
\medskip
i) Under extension via an $O_{(v)}$-monomorphism $W(k(v))\hookrightarrow\CC$, $\mu$ becomes $G(\CC)$-conjugate to the cocharacters $\mu^\ast_x$, $x\in X$;
\smallskip
ii) Under the direct sum decomposition $L^\ast_p\otimes_{\ZZ_p} W(k(v))=F^1\oplus F^0$ produced by
$\mu$ ($\be\in\GG_m\bigl(W(k(v))\bigr)$ acts through $\mu$ on $F^i$ as the multiplication with $\be^{-i}$), $F^1$ is normalized by $B_{1W(k(v))}$ (i.e. ${\rm Lie}(B_1)\otimes_{\ZZ_p} W(k(v))\subset F^0({\rm Lie}(G_{1W(k(v))}))$);
\smallskip
iii) $B_1\cap G_{\ZZ_p}$ contains a Borel subgroup $B$ of $G_{\ZZ_p}$.
\medskip
Moreover, if such a triple $(T,\mu,B_1)$ does exist, then $k(v)=k(v_1)$.}
\medskip
{\bf Proof:} Let $(T,\mu,B_1)$ be a triple of the mentioned form, such that i) to iii) hold. The fact that $i_1$ takes $G$-ordinary points of $\Mn_{k(v)}$ into $G_1$-ordinary points of ${\Mn_1}_{k(v_1)}$ results from 4.1-2: the Shimura-ordinary type of the SHS $(f,L_{(p)},v)$ (resp. of the SHS $(f_1,L_{(p)},v_1)$) can be computed, cf. 4.1.1, starting from the data of the triple $(T,\mu,B)$ (resp. $(T,\mu,B_1)$); so these two Shimura-ordinary types are equal. 
\smallskip
We now assume $i_1$ takes $G$-ordinary points into $G_1$-ordinary points. Let $z:{\rm Spec}\bigl(W(\FF)\bigr)\to\Mn$ be a $G$-canonical lift of a $G$-ordinary point $y:{\rm Spec}(\FF)\to\Mn_{k(v)}$, and let ${\got C}:=(M_y,F^1_z,\vph_y,G_{W(\FF)})$ be the Shimura filtered $\bar\sg$-crystal attached to $z$. We get another Shimura filtered $\bar\sg$-crystal 
${\got C}_1:=(M_y,F^1_z,\vph_y,{G_1}_{W(\FF)})$; it is the one attached to $y_2\circ z$ and (cf. 4.11.1.1) is a Shimura-canonical lift.
\smallskip
The canonical split cocharacter $\mu_y:\GG_m\to G_{W(k)}$ of ${\got C}$ when composed with the monomorphism $G_{W(k)}\hookrightarrow {G_1}_{W(k)}$, is the canonical split of ${\got C}_1$ (cf. 3.1.5). So the degree of definition of ${\got C}$ is equal to the degree of definition of ${\got C}_1$ and so (cf. 4.4.7) $k(v)=k(v_1)$.
\smallskip
Writing $\vph_y=a\mu_y({1\over p})$, $a$ becomes a $\sg_k$-linear automorphism of $M_y$. Then $M^a_y:=\bigl\{x\in M_y|a(x)=x\bigr\}$ is a $\ZZ_p$-structure on $M_y$. We get monomorphisms 
$\tilde G_{\ZZ_p}\buildrel{\tilde f_2}\over\hookrightarrow
\tilde G_{1\;\ZZ_p}\hookrightarrow GL(M^a_y)$ and a cocharacter $\mu_y:\GG_m\hookrightarrow\tilde G_{W(k(v))}$. All these are as in 4.4.7. Let $\tilde B$ be a Borel subgroup of $\tilde G_{\ZZ_p}$ normalizing $F^1_z$ and such that $\mu_y$ factors through $\tilde B_{W(k(v))}$. Let $\tilde B_1$ be a parabolic subgroup of $\tilde G_{1\;\ZZ_p}$ normalizing $F^1_z$ and containing $\tilde B$. The existence of $\tilde B$ is a direct consequence of b) of 4.4.1 2), while we can take $\tilde B_1$ such that the extension to $W(\FF)$ of its Lie algebra is the Lie subalgebra of ${\rm Lie}(G_{1W(\FF)})$ corresponding to non-negative slopes of $({\rm Lie}(G_{1W(\FF)}),\vph_y)$, cf. the parabolic property of 4.4.1 1) and cf. the existence of the $\ZZ_p$-structures of 3.11.2 C. Let $\tilde T$ be a maximal torus of $\tilde B$ such that $\mu_y$ factors through $\tilde T_{W(k(v))}$, cf. 3.11.2 B. 
\smallskip
$\tilde G_{\ZZ_p}$ is an inner form of $G_{\ZZ_p}$ (argument: Fontaine's comparison theory implies $\tilde G_{\QQ_p}$ is an inner form of $G_{\QQ_p}$). So it makes sense to speak about inner isomorphisms $\tilde G_{\ZZ_p}\tilde\to G_{\ZZ_p}$ or $\tilde G_{\ZZ_p}^{\rm ad}\tilde\to G_{\ZZ_p}^{\rm ad}$. The same applies in the context of $G_{1\ZZ_p}$ and $\tilde G_{1\ZZ_p}$.  
\medskip 
{\bf Claim*.} {\it We have inner isomorphisms $\tilde G_{\ZZ_p}\tilde\to G_{\ZZ_p}$ and $\tilde G_{1\ZZ_p}\tilde\to {G}_{1\ZZ_p}$, through which $\tilde f_2$ becomes the natural monomorphism ${\bar f}_2:G_{\ZZ_p}\hookrightarrow{G_1}_{\ZZ_p}$ as subgroups of $GL(L^\ast_p)$.} 
\medskip
So i) to iii) hold for the triple $(T,\mu,B_1)$ corresponding to the triple $(\tilde T,\mu_y,\tilde B_1)$ under these isomorphisms. The Claim is implied by d) of 4.4.1 3) (applied to $(f,L_{(p)},v)$ and $z$), via Lang's theorem; as we postponed the proof of d) of 4.4.1 3), we present another argument for why a triple $(T,\mu,B_1)$ of the mentioned form exists such that i) to iii) hold. 
\smallskip
The fact that in c) of 4.4.1 3) we can take $n=1$ or not for the SHS $(f,L_{(p)},v)$ and for the point $y$, is expressed in terms of a torsor of $\tilde G_{\QQ_p}$ being trivial or not. But the resulting torsor of $\tilde G^{\rm ad}_{\QQ_p}$ (i.e. of $G^{\rm ad}_{\QQ_p}$) is trivial, cf. b) of 4.4.1 2). So it corresponds to a torsor of the group $Z(\tilde G_{\QQ_p})=Z(G_{\QQ_p})$. As this last group is included in $\tilde T_{\QQ_p}$ (and so in $\tilde B_{\QQ_p}$ and $\tilde B_{1\QQ_p}$), we can consider the twist $(T_{\QQ_p},B_{\QQ_p},B_{1\QQ_p})$ of the triple $(\tilde T_{\QQ_p},\tilde B_{\QQ_p},\tilde B_{1\QQ_p})$. We view $T_{\QQ_p}$, $B_{\QQ_p}$ as subgroups of $G_{\QQ_p}$ and $B_{1\QQ_p}$ as a subgroup of $G_{1\QQ_p}$. Let $B$ (resp. $B_1$) be the Zariski closure of $B_{\QQ_p}$ (resp. of $B_{1\QQ_p}$) in $G_{\ZZ_p}$ (resp. in $G_{1\ZZ_p}$). It is a Borel (resp. a parabolic) subgroup of $G_{\ZZ_p}$ (resp. of $G_{1\ZZ_p}$), cf. Fact of 2.2.3 3). Let $T$ be the Zariski closure of $T_{\QQ_p}$ in $G_{\ZZ_p}$. It is a maximal torus (based on [Va2, 4.3.9] this can be read out from the adjoint context). But there is a unique cocharacter $\mu$ of $T_{W(k(v))}$ with the property that i) to iii) hold for $(T,\mu,B_1)$ (this can be checked easily in the adjoint context). This ends the proof of the Theorem.
\medskip
{\bf 4.11.3. Example.} We consider the case when $f_2$ satisfies the conditions:
\medskip
a) The connected components of the origin of $Z(G)$ and of $Z(G_1)$ coincide;
\smallskip
b) $G^{\rm ad}={\rm Res}_{F/\QQ} G^s$, with $F$ a totally real number field and with $G^s$ an absolutely simple $F$-group;
\smallskip
c) $f_2$ induces an injective map $f^{\rm ad}_2:(G^{\rm ad},X^{\rm ad})\hookrightarrow(G_1^{\rm ad},X_1^{\rm ad})$ which is the injective map $f_{2\;F_1}:{\rm Sh}(G^{\rm ad},X^{\rm ad})\hookrightarrow {\rm Sh}^{F_1}(G^{\rm ad},X^{\rm ad})$ defined in [Va2, 6.6.1] (with $F_1$ a totally real number field
containing $F$) (i.e. $G_1^{\rm ad}={\rm Res}_{F_1/\QQ} G^s_{F_1}$ and $f^{\rm ad}_2:G^{\rm ad}\hookrightarrow G_1^{\rm ad}$
is the canonical inclusion).
\medskip
Then $E(G_1,X_1)=E(G,X)$ and 4.11.2 i) to iii) hold for the triple $(T,\mu,B_1)$, where $(T,\mu,B)$ is a triple constructed as in 4.1 for the SHS $(f,L_{(p)},v)$ and such that $B$ contains $T$, and where $B_1$ is the unique Borel subgroup of $G_{1\ZZ_p}$ containing $B$. We get (cf. also 4.11.1.1): 
\medskip
{\bf Fact.} {\it $i_1$ takes $G$-ordinary points of $\Mn_{k(v)}$ into $G_1$-ordinary points of $\Mn_{1k(v_1)}$ and $G$-canonical lifts of $\Mn$ into $G_1$-canonical lifts of $\Mn_1$.}
\medskip
This Fact will be used in \S 6 and \S 12 (in connection to 1.15.8).
\medskip
{\bf 4.11.4. A third question.} There is an extra natural question related to 4.11.2:
\medskip
\item{\bf Q3.} {\it Does $i_1$ automatically take $G$-ordinary points into $G_1$-ordinary points if we know only that $k(v)=k(v_1)$?}
\medskip
{\bf 4.11.5. Remark.} If $k(v)=k(v_1)\ne\FF_p$ the answer to Q3 is (in general) no. There are many examples with $G$ a torus which provide a negative answer to Q3. If $k(v)=k(v_1)=\FF_p$ then a triple $(T,\mu,B_1)$ such that 4.11.2 i) to iii) hold, automatically exists and this motivates the form of 4.6 P1 and P2.
\medskip
{\bf 4.11.6. Example.} Here we concentrate just on Shimura $\sg$-crystals (4.12.12.6 below tells that everything below can be realized geometrically). What follows is very much related to the part of the proof of 4.3.6 involving the completion property. We assume $G^{\rm ab}_{\ZZ_p}$ is a split torus and $G^{\rm ad}_{\ZZ_p}$ is an absolutely simple, non-split group of $A_\ell$ Lie type (so $\ell\ge 2$); so $G_{\ZZ_p}^{\rm ad}$ splits over $W(\FF_{p^2})$. Let $B$ be a Borel subgroup of $G_{\ZZ_p}$. Let $T$ be a maximal torus of it. As $G_{\ZZ_p}$ splits over $W(\FF_{p^2})$, so does $B$ and $T$. We assume we have a SHS $(f,L_{(p)},v)$ such that $G^{\rm ad}_{\RR}=SU(a,\ell+1-a)_{\RR}^{\rm ad}$, for some $a\in S(1,\ell)\setminus\{{{\ell+1}\over 2}\}$. We also assume that the reflex field $E(G,X)$ (it is a quadratic totally imaginary extension of $\QQ$) has precisely one unramified prime $v$ above $p$; so $k(v)=\FF_{p^2}$. Let $\mu:\GG_m\to T_{W(\FF_{p^2})}$ be a cocharacter such that the conditions of 4.1 are satisfied (for $(f,L_{(p)},v)$). 
\smallskip
Let $\om\in G_{\ZZ_p}\bigl(W(\FF_{p^2})\bigr)$ normalizing $T_{W(\FF_{p^2})}$ and taking
$B_{W(\FF_{p^2})}$ into its opposite $B_{W(\FF_{p^2})}^{\rm opp}$ (see 4.1.4.1). Let
$\sg_1:=\om\bar\sg$; we view it as a $\bar\sg$-linear automorphism of $L_{(p)}^*\otimes_{\ZZ_{(p)}} W(\FF)$. Let $T_1\hookrightarrow\tilde G_{\ZZ_p}$ be an inclusion of reductive groups which over $W(\FF)$ becomes the inclusion $T_{W(\FF)}\hookrightarrow G_{W(\FF)}$ and which is defined at the level of Lie algebras by 
$${\rm Lie}(T_1):={\rm Lie}\bigl(T_{W(\FF)}\bigr)^{\sg_1}\hookrightarrow {\rm Lie}\bigl(G_{W(\FF)}\bigr)^{\sg_1}=:{\rm Lie}(\tilde G_{\ZZ_p})$$ 
(the upper right index $\sg_1$ means we take the elements fixed by $\sg_1$). We get:
\medskip
-- {\it the Shimura $\bar\sg$-crystals $\bigl(L^\ast_p\otimes_{\ZZ_p} W(\FF),\sg_1\mu_{W(\FF)}({1\over p}),T_{1\;W(\FF)}\bigr)$ and $\bigl(L^\ast_p\otimes_{\ZZ_p} W(\FF),\break\sg\mu({1\over p})\otimes 1,G_{W(\FF)}\bigr)$ are both Shimura-ordinary, with the degree of definition equal to 2.} 
\medskip
But the Newton polygon of $\bigl(L^\ast_p\otimes_{\ZZ_p} W(\FF),\sg_1\mu_{W(\FF)}({1\over p}),T_{1\;W(\FF)}\bigr)$ is different from the Newton polygon of $\bigl(L^\ast_p\otimes_{\ZZ_p} W(\FF),\sg\mu({1\over p})\otimes 1,G_{1\;W(\FF)}\bigr)$. Argument (cf. 3.9.2 and 3.1.0 a)): the Lie stable $p$-ranks of $\bigl(L^\ast_p\otimes_{\ZZ_p} W(\FF),\sg\mu({1\over p})\otimes 1,G_{W(\FF)}\bigr)$ and of $\bigl(L^\ast_p\otimes_{\ZZ_p} W(\FF),\sg_1\mu_{W(\FF)}({1\over p}),G_{W(\FF)}\bigr)$ are different, the second one being $0$.
\medskip
{\bf 4.11.6.1. Remark.} Always $i_1$ takes $G$-ordinary points into toric points, cf. 4.5.11.2.1. It is quite common, even if $k(v)\neq k(v_1)$, that $i_1$ takes $G$-ordinary points into $U$-ordinary points; see \S 6 for examples.
\medskip
{\bf 4.11.7. Example.} We assume that, as in 4.3.1, we have a product decomposition $G_{1\ZZ_p}^{\rm ad}=\times_{i\in\Mh} G_{1i}$ indexed by the same set $\Mh$ as for $G$. We also assume that the image of ${\rm Lie}(G_i)$ in ${\rm Lie}(G_{1\ZZ_p}^{\rm ad})$ lies in the direct summand ${\rm Lie}(G_{1i})$. As in 4.3.1.1 we define, for $i\in\Mh^{\rm nc}$, the $i$-th $A$-degree $d_{1A}(i)$ of definition of $(f_1,L_{(p)},v_1)$. We can have $k(v^{\rm ad})=k(v_1^{\rm ad})$, without having $d_A(i)=d_{1A}(i)$, $\forall i\in\Mh^{\rm nc}$; easy examples can be obtained with $(G^{\rm ad},X^{\rm ad})$ a simple, adjoint Shimura pair of $D_n^{\HH}$ type with non-trivial involution and with $(G_1^{\rm ad},X^{\rm ad}_1)$ a simple, adjoint Shimura pair of $A_{2n}$ type. We have the following general result:
\medskip
{\bf Fact.} {\it If there is $i\in\Mh^{\rm nc}$ such that $d_A(i)\neq d_{1A}(i)$, then a $G$-ordinary point of $\Mn_{k(v)}$ is not mapped by $i_1$ into a $G_1$-ordinary point of ${\Mn_1}_{k(v_1)}$.} 
\medskip
This is implied by the cyclic diagonalizability part of b) of 4.4.1 3) and by 4.4.13.3.
\medskip
{\bf 4.11.8. Exercise.} We assume that $k=\bar k$, that $i_1$ takes $G$-ordinary points of $\Mn_{k(v)}$ into $G_1$-ordinary points of $\Mn_{1k(v_1)}$, that neither $G$ nor $G_1$ are tori and that the $A$-degree of definition of a (any) $G_1$-canonical lift ${\rm Spec}(W(\FF))\to\Mn_1$ is 1. Show that: if $y:{\rm Spec}(k)\to\Mn_{k(v)}$ is a $G$-ordinary point, then the formal torus (of the moduli formal scheme) of $G$-deformations of $(A_y,p_{A_y})$ is a subtorus of the formal torus of $G_1$-deformation of $(A_y,p_{A_y})$; here we use the language and notations of 4.7.0 and 4.7.14. Hint: use 4.11.2, 4.7.17 and 4.7.11 4).
\medskip
This is a generalization, in the case of an arbitrary SHS whose $A$-degree of definition is 1, of the  main result of [No1]. There are variants (cf. 4.7.11 8)) of this exercise for an arbitrary $A$-degree of definition; warning: we do not know yet how useful they are (as they are not a priori canonical). 
\medskip\smallskip
{\bf 4.12. Integral and generalized Manin problems and applications.} Let $e\in\NN$ and let $p$ be a rational prime. Let $k$ be an algebraically closed field of characteristic $p$. 
\medskip
The (original) Manin problem asked for the proof of the predicted classification of isocrystals over $k$ attached to (principally polarized) abelian varieties over $k$ of dimension $e$. This problem was positively solved by the Honda--Serre--Tate's theory (cf. [Ta, p. 98] and [Mu, cor. 1 of p. 234]). For a second, more recent proof, see [Oo2]; see also 4.12.12.4. 
\medskip
{\bf 4.12.0. Overview of 4.12.} The main goal of this section is to formulate (and solve in some cases: see 4.12.4, the proof of 4.12.12, and 4.12.13 2) and 3); see \S 11 for the general case) different integral Manin problems for (special fibres of) integral canonical models of Shimura varieties of preabelian type w.r.t. primes dividing an odd rational prime. For the solution of these integral problems in the cases in which we know the completion property, we rely heavily on 3.6.15 A and 4.4.1. 
\smallskip
It turns out that the combination of the proof of 4.12.12 (we have in mind its application 4.12.12.6) with 3.9.7.2 and with a refinement of 4.3.6, implies the strong form of the generalized Manin problem for the integral canonical models of 4.9.8 (see 4.12.12.6.4; see 4.12.12.6.3.4 for the abstract extension to the generalized Shimura context); it is one of the main four ingredients needed to solve all these integral Manin problems (see 4.12.12.7 i) to iv)). 
\smallskip
Though we first treat the case $p\ge 3$, as a roundabout we show that (see 4.12.12.1 and 4.14.3 J below) the same ideas apply for $p=2$. 
\smallskip
These integral Manin problems have many applications: below we include relatively few of them; see \S 9-14 for many more. 
\smallskip
From now on, except for a special reference, $p\ge 3$.
\medskip
{\bf 4.12.1. The integral Manin problem for Siegel modular varieties.} Determine all  isomorphism classes of principally quasi-polarized $\sg_k$-crystals attached to principally polarized abelian varieties over $k$ of dimension $e$.
\medskip
{\bf 4.12.2. The integral Manin problem for a SHS $(f,L_{(p)},v)$.} Determine all isomorphism classes of principally quasi-polarized Shimura $\sg_k$-crystals and of Shimura Lie $\sg_k$-crystals over $k$ attached to $k$-valued points of $\Mn_{k(v)}$.
\medskip  
{\bf 4.12.3. The integral Manin problem for a Shimura quadruple $(G,X,H,v)$ of preabelian type, with $v$ dividing a rational prime $p\ge 3$.} Determine all isomorphism classes of Shimura adjoint Lie $\sg_k$-crystals attached (cf. 4.9.17) to $k$-valued points  of the special fibre $\Mn_{k(v)}$ of the integral canonical model $\Mn$ of the quadruple $(G,X,H,v)$.
\medskip
{\bf 4.12.4. Solution of 4.12.1.} Let $({\rm GSp}(W,\psi),S)$ define a Siegel modular variety, with $\dim_{\QQ}(W)=2e$. Let $L_{(p)}$ be a $\ZZ_{(p)}$-lattice of $W$ such that $\psi: L_{(p)}\otimes_{\ZZ_{(p)}} L_{(p)}\to\ZZ_{(p)}$ is a perfect form. We get a SHS $(1_{({\rm GSp}(W,\psi),S)},L_{(p)},p)$. Let $\{w_i|i\in S(1,2e)\}$ be a $\ZZ_{(p)}$-basis of $L_{(p)}^*$ such that for $i,j\in S(1,2e)$ satisfying $j\ge i$, $\tilde\psi(w_i,w_j)$ is equal to $1$ if $j-i=e$ and is equal to $0$ otherwise. Let $M_0:=L_{(p)}^*\otimes_{\ZZ_{(p)}} W(k)$ and let $F^1\subset M_0$ be the $W(k)$-submodule of $M_0$ generated by $w_i$, $i\in S(1,e)$. Let $\vph$ be the $\sg_k$-linear endomorphism of $M_0$ defined by: $\vph(w_i)=pw_i$, if $i\in S(1,e)$, and $\vph(w_i)=w_i$, if $i\in S(e+1,2e)$. We have:
\medskip
{\bf Answer:} {\it The set $\Ms_e(k)$ of isomorphism classes of $\sg_k$-crystals attached to principally polarized abelian varieties over $k$ of dimension $e$ is precisely the set of isomorphism classes of principally quasi-polarized $\sg_k$-crystals defined by elements of the set 
$$\{(M_0,g\vph,\tilde\psi)|g\in Sp(M_0,\tilde\psi)(W(k))\}.$$}
\indent
To argue this answer, we recall the following obvious isogeny property: 
\medskip
{\bf The isogeny property.} {\it If $(M_0^\prime,g\vph,\tilde\psi)$, with $M_0^\prime$ a $W(k)$-lattice of $M_0[{1\over p}]$, is the principally quasi-polarized $\sg_k$-crystal attached to a principally polarized abelian variety $(A^\prime,p_{A^\prime})$ over $k$, then there is a principally polarized abelian variety $(A,p_A)$ over $k$, $\ZZ[{1\over p}]$-isogeneous to $(A^\prime,p_{A^\prime})$ and of whose associated principally quasi-polarized $\sg_k$-crystal, under this $\ZZ[{1\over p}]$-isogeny, is identifiable with $(M_0,g\vph,\tilde\psi)$.}
\medskip
From [Ta, p. 98] and [Mu, cor. 1 of p. 234] we deduce the existence of an abelian variety $\tilde A$ over $k$ whose attached isocrystal is $(M_0[{1\over p}],g\vph)$ and which has a principal polarization $\tilde p_{\tilde A}$. But for any two principal quasi-polarizations $\psi_1$ and $\psi_2$ of $(M_0[{1\over p}],g\vph)$, the triples $(M_0[{1\over p}],g\vph,\psi_1)$ and $(M_0[{1\over p}],g\vph,\psi_2)$ are isomorphic. Argument: using Dieudonn\'e's classification of isocrystals and duals, the situation gets reduced to the case when $(M_0[{1\over p}],g\vph)$ is of pure slope ${1\over 2}$; in such a case, using the standard argument which shows that any two non-degenerate alternating forms on an even dimensional $B(k)$-vector space are isomorphic, we can write $(M_0[{1\over p}],g\vph,\psi_1)$ as a direct sum of principally quasi-polarized isocrystals of rank $2$ (for them the statement is trivial). So we can assume $\tilde\psi$ is defined by $\tilde p_A$. So the above answer follows from the mentioned isogeny property.
\medskip
If we assume 3.6.15 A, then the answers to 4.12.2-3 are provided by 4.4.1 as follows. 
\medskip
{\bf 4.12.5. Solution of 4.12.2 in the case $\Mn$ has the completion property.} Let $M_0:=M\otimes_{W(k(v))} W(k)$, $F_0^1:=F^1\otimes_{W(k(v))} W(k)$ and let $\vph_0:=\vph\otimes 1$ be the Frobenius endomorphism of $M_0$, where $M$, $F^1$ and $\vph$ have the same significance as in 4.1. Using 4.1.6, we have:
\medskip
{\bf Answer*:} {\it The set $\Ms_{(f,L_{(p)},v)}(k)$ of isomorphism classes of principally quasi-polarized Shimura $\sg_k$-crystals attached to $k$-valued points of $\Mn$ is precisely the set of isomorphism classes of principally quasi-polarized Shimura $\sg_k$-crystals defined by elements of the set 
$$\{(M_0,g\vph_0,G_{W(k)},(v_{\al})_{\al\in\Mj^\prime},\tilde\psi)|g\in G^0_{W(k)}(W(k))\}.$$}
\indent
{\bf Answer:} {\it The set ${\rm Lie}(\Ms_{(f,L_{(p)},v)})(k)$ of isomorphism classes of Shimura Lie $\sg_k$-crystals attached to $k$-valued points of $\Mn$ is precisely the set defined by elements of the set $\{({\rm Lie}(G_{W(k)}),g\vph_0)|g\in G^0_{W(k)}(W(k))\}$.}
\medskip
The proofs of these two answers are a direct consequence of 3.6.15 A and of d) of 4.4.1 3)  (resp. and of c) of 4.4.1 2)) for the first answer (resp. for the second one), cf. also 4.4.12. As we postponed the proof of d) of 4.4.1 3), we invite the reader to formulate a version of the first answer, inspired from 3.6.15 A, for $k$-valued points of a fixed connected component $\Mc$ of $\Mn_k$ (in fact 4.2.8.1 takes care of all connected components of $\Mn_k$ at once); we can also use 4.4.7 and 2.3.13.1: we get the above starred answer (for $\Mc$) but with $(M,F^1,\vph,G_{W(k(v))},(v_{\al})_{\al\in\Mj^\prime})$ replaced by another Shimura-canonical lift of a Shimura-ordinary $\sg_k$-crystal with an emphasized family of tensors indexed by $\Mj^\prime$ (so we do not know a priori that it is $1_{\Mj^\prime}$-isomorphic to $(M,F^1,\vph,G_{W(k(v))},(v_{\al})_{\al\in\Mj^\prime})$).
\medskip
{\bf 4.12.6. Solution of 4.12.3 in the case the integral canonical model $\Mn^{\rm ad}$ of the adjoint quadruple of $(G,X,H,v)$ has the completion property.} We refer to 4.9.17.4 for terminology. Let $G_{\ZZ_{(p)}}$ be the reductive group over $\ZZ_{(p)}$ whose generic fibre is $G$ and whose group of $\ZZ_{(p)}$-valued points is $H$, cf. [Ti2, 3.4.1] and [Va2, 3.1.3]. Let $v^{\rm ad}$ be the prime of $E(G^{\rm ad},X^{\rm ad})$ divided by $v$. Let $\mu^{\rm ad}:\GG_m\to G_{W(k(v^{\rm ad}))}^{\rm ad}$ be a cocharacter which under an $O_{(v^{\rm ad})}$-monomorphism $W(k(v^{\rm ad}))\hookrightarrow\CC$, is $G^{\rm ad}(\CC)$-conjugate to the cocharacters $\mu_x:\GG_m\hookrightarrow G^{\rm ad}_{\CC}$, with $x\in X^{\rm ad}$, defined in [Va2, 2.2]. We have: 
\medskip
{\bf Answer:} {\it The set $\Ms_{(G,X,H,v)}(k)$ of isomorphism classes of Shimura adjoint Lie $\sg_k$-crystals attached to $k$-valued points of $\Mn$ is precisely the set of isomorphism classes of Shimura adjoint Lie $\sg_k$-crystals defined by elements of the set 
$$\{({\rm Lie}(G_{W(k)}^{\rm ad}),g(\sg\mu^{\rm ad}({1\over p})\otimes 1)|g\in G_{W(k)}^{\rm ad}(W(k))\}.$$} 
\indent
To argue this, we consider a SHS $(f_1,L_{1(p)},v_1)$ such that $(G_1^{\rm ad},X_1^{\rm ad},H_1^{\rm ad},v_1^1)=(G^{\rm ad},X^{\rm ad},H^{\rm ad},v^{\rm ad})$, cf. [Va2, 6.4.2] and 2.3.6 (for $p=3$ cf. \S 6). The connected components of $\Mn^{\rm ad}_{W(\FF)}$ are permuted by $G^{\rm ad}(\AA_f^p)$, cf. [Va2, 3.3.2]. So the answer follows from 4.9.17.4 and b) of 4.4.1 2), once we remark that each connected component of $\Mn_{W(\FF)}$ is naturally a pro-\'etale cover of a connected component of $\Mn^{\rm ad}_{W(\FF)}$ (cf. [Va2, 6.4.5.1]). 
\medskip
{\bf 4.12.6.1. Fact.} {\it  $\Mn^{\rm ad}$ has the completion property, provided 3.6.15 A is true for a standard Hodge situation, involving a Shimura quadruple having the same adjoint as $(G,X,H,v)$.}
\medskip 
{\bf Proof:} We first assume $G^{\rm ad}$ is a $\QQ$--simple, adjoint group. Let $o(G^{\rm ad})$ be the order of the center of the simply connected group cover of $G_{\CC}^{\rm ad}$. The proof of the Fact in the case when $p\not |o(G^{\rm ad})$ is a direct consequence of 4.9.17 and of b) of 4.4.1 2): if $p\not |o(G^{\rm ad})$ then the natural homomorphism $G_{W(k)}^{\rm der}(W(k))\to G_{W(k)}^{\rm ad}(W(k))$ is surjective. 
\smallskip 
As we assumed $p\ge 3$, $p$ divides $o(G^{\rm ad})$ only in the case when  $(G^{\rm ad},X^{\rm ad})$ is of $A_l$ Lie type, with $p$ dividing $l+1$. But in such a case there is an injective map $f: (G_1,X_1)\hookrightarrow ({\rm GSp}(W,\psi),S)$, with $(G_1^{\rm ad},X_1^{\rm ad})=(G^{\rm ad},X^{\rm ad})$ and with $G_1^{\rm der}$ a simply connected, semisimple group, such that we have a standard PEL situation $(f,L_{(p)},v_1,\Mb)$, with $v_1$ dividing $v^{\rm ad}$ and with the homomorphism ${G_1}_{W(k)}(W(k))\to G_{W(k)}^{\rm ad}(W(k))$ surjective, cf. 2.3.5.1. The surjectivity part results from the fact that for any simple factor of $G^{\rm ad}_{W(k)}$ there is a general linear group over $W(k)$ which is included in ${G_1}_{W(k)}$ and which has the given factor as its adjoint group, cf. 2.3.5.1 and [Ko2, (A) of p. 395]. The completion property in the context of $\Mn_1$ is stated (see 3.6.15 A) in terms of elements of $G_{1W(k)}^0(W(k))$; but from the proof of 2.3.16 we deduce that provided we ignore the principal quasi-polarizations, it can be stated as well in terms of elements of $G_{1W(k)}(W(k))$. So, again directly from def. 4.9.17 we get the validity of the above statement in this case, provided $\Mn_1$ has the completion property. This proves the weaker form of the Fact obtained by replacing ``for a standard Hodge situation" by ``for a certain standard Hodge situation" (it is [Va2, 6.5.1] and 2.3.5.2 which allow us to eliminate the assumption that $G^{\rm ad}$ is $\QQ$--simple). 
\smallskip
We now prove the Fact in its general form, without making any distinction on how $p$ is and without assuming $G^{\rm ad}$ is $\QQ$--simple. Let $\tilde T$ be the image of $T$ of 4.1 in $G^{\rm ad}_{\ZZ_p}$. We have (cf. (DER) of 2.2.6 1); see 4.2.7 for the def. of $m_{\rm qu}$) 
$$G^{\rm ad}_{W(\FF)}(W(k))={m_{\rm qu}}_{W(k)}(G_{W(\FF)}(W(k))))\tilde T(W(k)).$$ 
From this and the adjoint variant 3.11.7 of 3.11.1 c) (see also 3.5.5) we get (the notations are as of 4.12.6): 
\medskip
{\bf Subfact.} {\it Any Shimura adjoint Lie $\sg_k$-crystal $({\rm Lie}(G^{\rm ad}_{W(k)}),g(\sg\mu^{\rm ad}({1\over p})\otimes 1))$, with $g\in G^{\rm ad}_{W(k)}(W(k))$, is inner isomorphic to $({\rm Lie}(G^{\rm ad}_{W(k)}),g^1(\sg\mu^{\rm ad}({1\over p})\otimes 1))$ for some $g^1\in G^{\rm der}_{W(k)}(W(k))$.} 
\medskip
From this Subfact the Fact follows.
\medskip
{\bf 4.12.7. The Manin problem for a SHS $(f,L_{(p)},v)$.} Classify all  
principally quasi-polarized Shimura isocrystals (resp. all Shimura adjoint Lie isocrystals) over $k$ defined by principally quasi-polarized Shimura $\sg_k$-crystals (resp. by Shimura adjoint Lie $\sg_k$-crystals) attached to $k$-valued points of $\Mn$.  
\medskip
{\bf 4.12.8. The Manin problem for a Shimura quadruple $(G,X,H,v)$ of preabelian type, with $v$ dividing a rational prime $p\ge 3$.} Classify all adjoint Lie isocrystals over $k$ defined by Shimura adjoint Lie $\sg_k$-crystals attached to $k$-valued points of the integral canonical model of $(G,X,H,v)$.
\medskip
{\bf 4.12.9. Remarks. 1)} To our knowledge the result (answer) of 4.12.4 has not been stated as such for Siegel modular varieties of dimension greater than $3$. It says: 
\medskip
{\bf Fact.} {\it Any principally quasi-polarized $\sg_k$-crystal (i.e. any principally quasi-polarized $p$-divisible group) over $k$, is associated to a principally polarized abelian variety over $k$, and any principally quasi-polarized filtered $\sg_k$-crystal (i.e. any principally quasi-polarized $p$-divisible group over $W(k)$), is associated to a principally polarized abelian scheme over $W(k)$.}
\medskip
{\bf 2)} 4.12.5-6 form implicitly a solution of the problems 4.12.7-8, provided we assume the completion property. But it is very much desirable, to have a more concrete (and finite) answer to 4.12.7-8, similar to the answer of the (original) Manin problem in terms of formal isogeny types. From numerical point of view, we will deal with this in \S10 (see 4.12.10 below for a sample); from the point of view of toric points, solutions of 4.12.7-8 can be obtained starting from 4.12.12.6.4 below.
\medskip 
{\bf 4.12.10. Example.} Let $n\ge m$ be natural numbers and let $l:=m+n$. Let $(f,L_{(p)},v)$ be a SHS such that $G_{\RR}^{\rm der}=SU(n,m)_{\RR}$, $G_{\ZZ_p}^{\rm der}=SL(N)$, with $N$ a free $\ZZ_p$-module of rank $l$, and $\Mn$ has the completion property. Then the set of isomorphism classes of Shimura adjoint Lie isocrystals over $k$ defined by Shimura adjoint Lie $\sg_k$-crystals attached to $k$-valued points of $\Mn$, is the set obtained by taking ${\got s}{\got l}(*)$ of isomorphism classes $(*)$ of isocrystals $(\Mi,\vph_{\Mi})$ over $k$ having only  slopes which are rational numbers of the interval $[0,1]$, whose total sum (involving multiplicities) is $m$, and with $\Mi$ a free $B(k)$-vector space of dimension $l$. This is a consequence of 4.12.6 and of Subfact of 4.12.6.1 (cf. also 4.9.17).
\medskip
For more examples see \S10 and [Va4].
\medskip
{\bf 4.12.11. Problems. 1)} Let $(G,X,H,v)$ be a Shimura quadruple of preabelian type, with $v$ dividing a rational prime $p\ge 3$. Let $k_1(v)$ be a finite field extension of $k(v)$ and let $H_0$ be a compact, open subgroup of $G(\AA_f^p)$ as in 4.9.20. Let $\sg_1$ be the Frobenius automorphism of $W(k_1(v))$. Determine all Shimura adjoint Lie $\sg_1$-crystals over $k_1(v)$ attached (see 4.9.20) to $k_1(v)$-valued points of $\Mn/H_0$, where $\Mn$ is the integral canonical model of the quadruple $(G,X,H,v)$. In particular, determine conditions on $H_0$ and $k_1(v)$, which assures $\Mn/H_0(k_1(v))$ is not empty. 
\medskip
{\bf 2)} For a fixed subgroup $H_0$ of $G(\AA_f^p)$ as in 1), determine the smallest algebraic field extension  $k_1(v)$ of $k(v)$ with the property that any isomorphism class of Shimura adjoint Lie $\bar\sg$-crystals (resp. isocrystals) over $\FF$ attached to $\FF$-valued points of $\Mn$, is defined (by extension to $\FF$ and by inverting $p$) by a Shimura adjoint Lie $\sg_1$-crystal (resp. isocrystal) attached to a $k_1(v)$-valued point of $\Mn/H_0$.
\medskip
{\bf 3)} The same problems 1) and 2) but for a SHS $(f,L_{(p)},v)$ (in this case the reference to 4.9.20 can be substituted by 2.3.10).
\medskip
{\bf 4.12.11.1. Remark.} If there is a SHS $(f,L_{(p)},v)$, with $H=G(L_{(p)}\otimes_{\ZZ_{(p)}} \ZZ_p)$, and if either all Hodge cycles of $\Ma$ are as well Hodge cycles of $\Ma_{H_0}$ or $k_1(v)$ is big enough and we have restrictions on the Lie types of the simple factors of $G^{\rm ad}$ similar in nature to the ones of 2.3.9, then in \S 14 we will see that the answer to problem 4.12.11 1) is of the same form as the one given in 4.12.6 but with $k$ replaced by $k_1(v)$ and with $g$ running through a subset $\Mg\bigl(k_1(v),\Mn/H_0\bigl)$ of $G_{W(k_1(v))}^{\rm ad}(W(k_1(v)))$. Also in \S 14 we will deal with the problem of determining this subset and of what can be said when no assumption on $H_0$ is made.
\medskip
{\bf 4.12.12. Theorem.} {\it Any $p$-divisible group $D$ over $k$ has a uni plus versal deformation (to a $p$-divisible group) over a smooth, affine $W(k)$-scheme $\Mm_D$.}
\medskip
{\bf Proof:} Let $(M,\vph)$ be the $\sg_k$-crystal defined by $D$ and let $(M,F^1,\vph)$ be the filtered $\sg_k$-crystal associated to an arbitrary lift of $D$ to a $p$-divisible group over $W(k)$. If $F^1$ is not a proper summand of $M$ we can take $\Mm_D={\rm Spec}(W(k))$. So we can assume we have a Shimura filtered $\sg_k$-crystal $(M,F^1,\vph,GL(M))$. Let $a:=\dim_{W(k)}(F^1)$ and $b:=\dim_{W(k)}(M/F^1)$. We have $ab\neq 0$. We can assume $\vph=g\vph_0$, with $g\in SL(M)(W(k))$ and with $\vph_0$ a Frobenius endomorphism of $M$ such that $(M,F^1,\vph_0,GL(M))$ is a $GL(M)$-canonical lift. So $F^1$ is left invariant by $\vph_0$ and the slopes of the isocrystal defined by $(M,\vph_0)$ are 0 and 1 with multiplicities $b$ and respectively $a$. If $ab=1$, then $D$ is the $p$-divisible group of an elliptic curve and so the Theorem is well known. We assume now that $ab\ge 2$. 
\smallskip
We consider a SHS $(f,L_{(p)},v)$ such that:
\medskip
{\bf a)} $G^{\rm der}_{\RR}$ is $SU(a,b)_{\RR}$;
\smallskip
{\bf b)} $G$ splits over $\QQ_p$; so $k(v)=\FF_p$ and $G^{\rm der}_{\ZZ_p}=SL(\ZZ_p^{a+b})$;
\smallskip
{\bf c)} there is a torus $T$ of $G$ of dimension $0$ or $1$ depending on the fact that $a$ is or is not equal to be $b$, such that $T_{\RR}$ is compact and the monomorphisms ${\rm Res}_{\CC/\RR} \GG_m\to G_{\RR}$ defining elements of $X$ factor through $G_{1\RR}$, where $G_1$ is the subgroup of $G$ generated by $G^{\rm der}$, by $T$ and by the $1$ dimensional split torus of scalar automorphisms of $W$;
\smallskip
{\bf d)} $E(G,X)$ is $\QQ$ or a totally imaginary quadratic extension of $\QQ$ depending on the fact that $a$ is or is not equal to $b$;
\smallskip
{\bf e)} there is $m\in\NN$ such that the representation $G_{\ZZ_p}\to GL(L_{(p)}^*\otimes_{\ZZ_{(p)}} \ZZ_p)$ is a direct sum of $2m$ irreducible representations $\rho_1$, $\rho_2$,..., $\rho_{2m}$ such that the restriction of $\rho_i$ to $G^{\rm der}_{\ZZ_p}$ is the representation associated to the minimal weight $\bar\om_1$ (resp. $\bar\om_{a+b-1}$) if $i\in S(1,m)$ (resp. if $i\in S(m+1,2m)$): here we use the standard notations for weights (see [Bou2, planche I]); in other words, we have a direct sum decomposition as $G_{\ZZ_p}$-modules $L_{(p)}^*\otimes_{\ZZ_{(p)}} \ZZ_p=\oplus_{i=1}^{2m} L^i$, with $L^i=\ZZ_p^{a+b}$, such that $A\in G_{\ZZ_p}^{\rm der}(\ZZ_p)=SL(\ZZ_p^{a+b})(\ZZ_p)$ acts as $A$ on $L^1$,..., $L^m$ and as the inverse of the transpose of $A$ on $L^{m+1}$,..., $L^{2m}$; 
\smallskip
{\bf e')} moreover, $\rho_i$ is the dual of $\rho_{i+m}$ w.r.t. $\psi$, i.e., for $i$, $j\in S(1,2m)$, with $\abs{j-i}\neq m$, $L^i$ is perpendicular on $L^j$ w.r.t. $\psi$;
\smallskip
{\bf f)} $(G,X)$ is of compact type.
\medskip
{\bf Argument.} The existence of such a SHS is just a particular case of 2.3.5.1: in fact we get a standard PEL situation $(f,L_{(p)},v,\Mb)$, with $\Mb$ as the $\ZZ_{(p)}$-subalgebra of ${\rm End}(L_{(p)})$ formed by endomorphisms fixed by $G_{\ZZ_{(p)}}$. We just have to remark a couple of things. To explain them, we start with an arbitrary adjoint Shimura pair $(G_0,X_0)$ of compact type such that $G_0$ splits over $\QQ_p$ and $G_{0\RR}$ is $SU(a,b)_{\RR}^{\rm ad}$. It can be constructed as follows. Let $l$ be a prime different from $p$ and let $K_1$ be a totally imaginary quadratic extension of $\QQ$ splitting at $p$ and $l$. Let $(G^1,X^1)$ be a Shimura pair constructed as in [Go] starting from a non-degenerate hermitian form on a $K_1$-vector space $V$ of dimension $a+b$ which has signature $(a,b)$ over $\RR$ (i.e. when extended to $V\otimes_{\QQ} \RR$). $G_1$ splits over $\QQ_p$ and over $\QQ_l$. We twist $G^{1\rm ad}$ by an inner form, corresponding to a central division algebra over $\QQ$ of dimension $(a+b)^2$ which splits at $p$ and at $\infty$ and which is ramified at $l$: we get $G_0$. $G_0$ splits over $\QQ_p$, has rank $0$ over $\QQ_l$ (for instance, cf. [Bo2, 23.1]), and so over $\QQ$; moreover, as $G_{0\RR}=G^{1\rm ad}_{\RR}$, we get a Shimura pair $(G_0,X_0)$, with $X_0:=X^{1\rm ad}$.  
\smallskip 
The reflex field $E(G_0,X_0)$ is $\QQ$ if $a=b$ and a totally imaginary quadratic extension $K_0$ of $\QQ$ splitting above $p$ if $a\neq b$ (cf. [De2, 2.3.6] and the fact that $G_0$ is absolutely simple). We apply 2.3.5.1 to $(G_0,X_0)$: we get a standard PEL situation $(f,L_{(p)},v,\Mb)$, with $f:(G,X)\hookrightarrow (GSp(W,\psi),S)$ an injective map for which we have $(G^{\rm ad},X^{\rm ad})=(G_0,X_0)$. This takes care of a) and f) (cf. [BHC] for the f) part). As 2.3.5.1 is just a restatement of [Va2, 6.5.1.1 and Case 2 of 6.6.5.1] (which refers to the case $p=3$ as well), we just need to point that we can arrange the things so that b), c), d), e) and e') hold as well. In loc. cit. we considered a totally imaginary quadratic extension of $\QQ$: we take it to split above $p$ and in case $a\neq b$ (in order to get c)) we need it to be $K_0$. 
\smallskip
We first treat the case $a=b$. In this case c) and d) follow from [De2, 2.3.13]: the centralizer $C_1$ in $GSp(L_{(p)},\psi)$ of the Zariski closure $G_{1\ZZ_{(p)}}$ of $G_1$ in $GL(L_{(p)})$ is a reductive group over $\ZZ_{(p)}$ which splits over $\ZZ_p$ (we need to choose $\tilde\psi$ in [Va2, 6.5.1.1 v)] in the similar way to the one mentioned in [Va2, Case 2 of 6.6.5.1], so that e') and e) hold mod $p$); so we take $G$ to be generated by $G_1$ and (cf. [Har, 5.5.3]) by the generic fibre of a maximal torus of $C_1$ which splits over $\ZZ_p$ (as in the proof of Lemma 2 of 4.6.4 we argue that $G_{\ZZ_p}$ is reductive). So $G_{\ZZ_p}$ is split and so b) follows. As $G_{\ZZ_p}$ is split, e) and e') follow directly from their mod $p$ versions. 
\smallskip
From now on, till the end of the argument, we assume $a\neq b$. d) is a consequence of [Va2, 6.5.1.1 c)]: the $\QQ$--algebra $K_S$ of the proof of [De2, 2.3.10] and used in [Va2, 6.5.1.1 and Case 2 of 6.6.5.1] is $K_0$ itself; so, cf. [De1, 2.3.9], $E(G^{\rm ab},X^{\rm ab})=K_0$ and so $E(G,X)=K_0$. c) follows from [Va2, 6.5.1.1 iv)]: we take $T$ to be the subtorus of dimension 1 of the torus defined by invertible elements of $K_S\otimes_{\QQ} K_0$ which over $\RR$ is compact and c) holds (the fact that we can take $T$ to be of dimension 1 and not 2 is implied by our choice $K_S=K_0$); obviously $G_1$ splits over $\QQ_p$. To get b) and e), we just have to take $G$ to be defined as above, via $C_1$. This ends the argument of the existence of a SHS subject to the requirements a) to f).
\medskip
{\bf The completion property.} We assume $H_0$ is small enough so that all Hodge cycles of $\Ma$ are Hodge cycles of $\Ma_{H_0}$ as well. So $\Mb\otimes_{\ZZ_{(p)}} \ZZ_p$ acts on the $p$-divisible group of the extension of $\Ma_{H_0}$ to $\Mn_{W(k)}/H_0$. We deduce (cf. e)) that this $p$-divisible group is a direct sum of $2m$ $p$-divisible groups $\Md_1$,..., $\Md_{2m}$ of rank $a+b$: $\Md_i$ ``corresponds" to $L^i$. Each one of them is a uni plus versal deformation in each $k$-valued point of $\Mn_{W(k)}/H_0$ and has dimension $a$ or $b$; moreover, $\forall i\in S(1,m)$, $\Md_i^t$ is $\Md_{i+m}$ (cf. e')). We can assume $\Md_1$ has dimension $a$ in each $k$-valued point of $\Mn_{W(k)}/H_0$. 
\smallskip
Let $f_1$ be the restriction of $f$ to $G_1$. The triple $(f_1,L_{(p)},v)$ is as well a SHS (the argument is the same as the one of b) of Fact 1 of 2.3.5.2). Let $\Mn_1$ be the integral canonical model associated to it. 
\smallskip
We start checking that $\Mn$ has the completion property. $\Mn_1$ is naturally an open closed subscheme of $\Mn$ (cf. [Va2, 3.2.15]). From this, [Va2, 3.3.2], and the fact (see 2.3.3.1) that $\Mn/H_0$ is a projective $O_{(v)}$-scheme we get: it is enough to check that $\Mn_1$ has the completion property. So based on e) and e') we can treat the situation as if $m=1$: the representation of $G_{1\ZZ_p}$ on $L^1$,..., $L^m$ are isomorphic, cf. above constructions (this can be checked over $\CC$ and so we just have to trace $K_S$ and $K_0$ in the proof of [De2, 2.3.10]). 
 Accordingly, in what follows, not to complicate the notations we keep in mind that we can treat situation as if $m=1$ but we still deal with the checking that $\Mn$ itself has the completion property. We follow the ideas of 4.12.4. First, the analogous isogeny property can be easily checked in this context. With the notations of 4.12.4, we obviously have:
\medskip
{\bf The isogeny property for $(f,L_{(p)},v,\Mb)$.} {\it We assume we have a triple $(A^\prime,p_{A^{\prime}},\Mb)$ which is obtained from $(\Ma,\Mp_{\Ma})$ and its natural family of $\ZZ_{(p)}$-endomorphisms (still denoted by $\Mb$) by pull back via a $k$-valued point $y^\prime$ of $\Mn$. If any invertible element of $\Mb$, when viewed naturally (via de Rham components) as an automorphism of $M_0^\prime$, is also an automorphism of $M_0$, then all these elements (assumed to be de Rham components of $\ZZ_{(p)}$-automorphisms of $A$) are de Rham components of $\ZZ_{(p)}$-automorphisms of $A$.}
\medskip
We consider a lift $(A_{W(k)},p_{A_{W(k)}},\Mb)$ to $W(k)$ of the triple $(A,p_A,\Mb)$ we get (its existence is argued using Grothendieck--Messing's theory of [Me, ch. 4-5] and standard arguments on representations of products of matrix $W(k)$-algebras; see also [Va2, 4.3.11], AE.1 and e) of Exercise of 2.3.18 B3). For $s\in\NN$, with $(s,p)=1$, we endow $(A_{W(k)},p_{A_{W(k)}})$ with the level-$s$ symplectic similitude structure obtained naturally from the one (see 2.3.2) of $(A^\prime,p_{A^\prime})$ via the mentioned $\ZZ[{1\over p}]$-isogeny. The fact that the triple $(A_{W(k)},p_{A_{W(k)}},\Mb)$ together with the mentioned level-$s$ symplectic similitude structures is obtained from $(\Ma,\Mp_{\Ma},\Mb)$ and its natural level-$s$ symplectic similitude structures by pull back through a $W(k)$-valued point of $\Mn$, is a consequence of the following four things:
\medskip
-- the Hasse principle holds for $G$ even if $a+b$ is even, see the first two paragraphs of [Ko2, ch. 7] (referring to the notations of loc. cit., in our case $F_0=\QQ$ and so $D=F^{\times}$; our $K_0$ is nothing else but $F$ of loc. cit.);
\smallskip
-- the triples $(H^1_{\rm crys}(A^\prime/W(k)),p_{A^\prime},\Mb)$ and $(H^1_{\rm crys}(A/W(k)),p_{A},\Mb)$ are isomorphic (this is just the simplest case of the equivalent for $p$ of b) of 2.3.18 B3); so (cf. also 4.6.8 and Lang's theorem) the triples (see 2.3.1 and 4.1 for the first one) $(L_p^,\tilde\psi,\Mb)$ and $(H^1_{\acute et}(A_{\overline{B(k)}},\ZZ_p),p_{A_{W(k)}},\Mb)$ are isomorphic: this takes care of the $\QQ_p$-context;
\smallskip
-- [Ko2, 4.3] takes care of the $\RR$-context;
\smallskip
-- the choice of the level-$s$ symplectic similitude structures takes care of the $\QQ_l$-context, for any prime $l$ different from $p$.
\medskip 
Conclusion: the triple $(A,p_A,\Mb)$ is as well obtained from $(\Ma,\Mp_{\Ma},\Mb)$ by pull back via a $k$-valued point $y$ of $\Mn$. 
\smallskip
Dieudonn\'e's classification of isocrystals over $k$ implies that any isocrystal over $k$ of a $p$-divisible group $E$ is isomorphic to the isocrystal of a cyclic diagonalizable filtered $\sg_k$-crystal over $k$ whose $F^1$-filtration has the same rank as the dimension of $E$. So, from this and from the above isogeny property, as we can treat the situation as if $m=1$, we get: to show that $\Mn$ has the completion property, we just need to check that for any cyclic diagonalizable Shimura $\sg_k$-crystal ${\got C}_N=(N,F^1_N,\vph_N,GL(N))$, with $N$ a free $W(k)$-module of rank $a+b$ and with $F^1_N$ a direct summand of it of rank $a$, there is a $k$-valued $y_N$ point of $\Mn_{W(k)}/H_0$ such that the Shimura $\sg_k$-crystal of $y_N^*(\Md_1)$ (see 2.2.9 9)), is isomorphic to $(N,\vph_N,GL(N))$. It is 2.3.17, which allows us to pass on to a Shimura filtered context: there is a $W(k)$-valued point $z_N$ of $\Mn$ lifting $y_N$ and such that the Shimura filtered $\sg_k$-crystal of $z_N^*(\Md_1)$ is isomorphic to $(N,F^1_N,\vph_N,GL(N))$.
\smallskip
Let $T_N$ be a maximal torus of $GL(N)$ such that $(N,F^1_N,\vph_N,T_N)$ is a Shimura $\sg_k$-crystal (cf. 2.2.16). We consider the Galois representation associated to it; it can be viewed as a homomorphism
$$\rho_N:\Gamma_k\to\tilde T_N(\ZZ_p),$$
where $\tilde T_N$ is a suitable maximal torus of the image $GL^1$ of $G_{\ZZ_p}$ in the $GL$-group $GL^1$ of the first $\ZZ_p^{a+b}$-summand $L^1$ of $L_{(p)}^*\otimes_{\ZZ_{(p)}} \ZZ_p$. $\tilde T_{NW(\FF)}$ is naturally equipped with a cocharacter $\mu_N$: identifying $T_N$ with $\tilde T_{NW(k)}$ (cf. 2.2.16.2.1), it is the cocharacter defining the canonical split of $(N,F^1_N,\vph_N,T_N)$; it acts trivially on a $W(k)$-submodule of $L^1\otimes_{\ZZ_p} W(k)$ of rank $b$. Let (cf. [Har, 5.5.3]) $T_0$ be a maximal $\ZZ_{(p)}$-torus of $G_{1\ZZ_{(p)}}$ such that:
\medskip
g) over $\RR$ it is the extension of a compact torus by a $1$ dimensional split torus;
\smallskip
h) its extension to $\ZZ_p$ has an image $T_1$ in $GL^1$ which is $GL^1(\ZZ_p)$-conjugate to $\tilde T_N$.
\medskip
Let $g\in GL^1(\ZZ_p)$ be such that $gT_Ng^{-1}=T_1$. 
We consider a cocharacter 
$$\mu(T_0):\GG_m\to T_{0\CC}$$
such that:
\medskip
i) viewed as a cocharacter of $G_{\CC}$ it is $G(\CC)$-conjugate to cocharacters $\mu_x$, $x\in X$;
\smallskip
j) viewing $T_{1\CC}$ as a quotient of $T_{0\CC}$, its natural cocharacter defined by $\mu(T_0)$ is obtained from the cocharacter $g\mu_Ng^{-1}$ of $T_{1W(\FF)}$ by extension of scalars (under an $O_{(v)}$-monomorphism $W(\FF)\hookrightarrow\CC$).
\medskip
$\mu(T_0)$ is uniquely determined: if $a=b$ then $T_1=T_{0\ZZ_p}$; if $a\neq b$ then this is a consequence of the structure of $G_{1\ZZ_p}$ (it is naturally isomorphic to $GL^1\times\GG_m$). 
\smallskip
As $T_0$ is the extension of a $1$ dimensional split torus by a torus which is compact over $\RR$, from i) we deduce that the image of $\mu(T_0)$ is defined over $\RR$. We consider the subgroup of $T_{0\RR}$ generated by this image and by the split subtorus of $T_{0\RR}$. It is nothing else but ${\rm Res}_{\CC/\RR} \GG_m$, cf. i). So we get naturally a monomorphism  
$$h_0:{\rm Res}_{\CC/\RR} \GG_m\hookrightarrow {T_0}_{\RR}.$$
From [De2, 1.2.8] (see also the general result of [Ko2, 4.3]) we get: the $G(\RR)$-conjugacy class of $h_0$ (viewed as a cocharacter of $G_{\RR}$) is $X$ itself; in fact the map 
$$(G_1,X_1)\to (G^{\rm ad},X^{\rm ad})=(G_1^{\rm ad},X_1^{\rm ad})$$ 
is a cover (cf. its construction; the kernel of $G_1\to G_1^{\rm ad}$ is ${\rm Res}_{K_0/\QQ} \GG_m$) and so $X_1=X=X^{\rm ad}$ (see [Mi1, 4.11]).
We deduce the existence of a special quadruple (see [Va2, 3.2.10])
$$({T_0}_{\QQ},\{h_0\},T_0(\ZZ_p),v_p)\hookrightarrow (G,X,H,v),$$ 
with $v_p$ a prime of $E(T_0,\{h_0\})$ dividing $v$, which has the property that the cocharacter of ${T_0}_{W(k(v_p))}$ we naturally get (as in 4.1) is modeled on $\mu_N$ (i.e. its extension to $\CC$ via an inclusion $W(k(v_p))\subset\CC$ extending the natural inclusion $O_{(v_p)}\subset\CC$, is $\mu(T_0)$ itself). Now we can easily check (cf. 4.1.1 and 4.2.3.1) that any $k$-valued point of $\Mn_{W(k)}/H_0$ factoring through the integral canonical model of this special quadruple, can be chosen as the desired $k$-valued point $y_N$. We conclude: $\Mn$ (and so --as we can treat the situation as if $n=1$-- also $\Mn_1$) has the completion property.
\smallskip
So there is $y\in\Mn_{W(k)}/H_0$ such that $y^*(\Md_1)$ is isomorphic to $D$ (cf. 4.12.5); as $k(v)=\FF_p$, to get 4.12.5 for our SHS's $(f,L_{(p)},v)$ and $(f_1,L_{(p)},v)$ we can refer to 4.6.5 instead of d) of 4.4.1 3. This ends the proof.   
\medskip
{\bf 4.12.12.0. Variants and comments.} If we do not require 4.12.12 f) to hold, we can take $m=1$, i.e. we can work out the things in the context of $(G^1,X^1)$ itself, without performing any inner twist.
\smallskip
What was important in 4.12.12?  Answer: that over $\ZZ_p$ we get a ``right" direct sum decomposition of $L_{(p)}^*\otimes_{\ZZ_{(p)}} \ZZ_p$ into irreducible $G_{\ZZ_p}$-representations so that we can ``single out" a $p$-divisible group $\Md_1$ as above which is ``related to" $D$. In fact we do not need a direct sum decomposition: we just need an adequate irreducible subrepresentation; however, if we do have such a ``right" direct sum decomposition of $L_{(p)}^*\otimes_{\ZZ_{(p)}} \ZZ_p$, then it is much easier to check that an isogeny property (analogue to the of 4.12.4) holds. So 4.12.12 e) can be replaced by more general versions, inspired from [Va2, 6.5-6].
\smallskip
In future we will use the above paragraph to show that 4.12.12 remains true (under proper formulation) for any generalized Shimura $p$-divisible object $(M_0,\vph_0,G_0)$ over $k$: for the case when $G_0$ together with its filtration class is related to Shimura varieties of preabelian type see \S 7; for the case when the situation is related to Shimura varieties of special type see [Va6] (cf. also 3.15.6). Moreover it seems to us that the condition $k=\bar k$ is not (really) needed: in other words, (at least) in many situations it is enough to work with a perfect field (as the combination of Theorem 2 of 3.15.1, 3.15.3 1) and 4.12.12 suggests).  
\medskip 
{\bf 4.12.12.1. The case $p=2$.} The natural question arises: what about the case $p=2$? Defining the completion property for $p=2$ as in 3.6.16 1) and 4.9.17.4, we have:
\medskip
{\bf Answer:} {\it 4.12.4 remains true for $p=2$. Also, 4.12.5 remains true for the case of a $p=2$ SHS $(f,L_{(2)},v)$ whose integral canonical model $\Mn$ has the completion property.}
\medskip
Moreover, the general case of the proof of the Fact of 4.12.6.1 involving a Subfact, holds as well for $p=2$. The same proofs (for $p\ge 3$) apply to get these results for $p=2$. We also have:
\medskip
{\bf 4.12.12.2. Corollary.} {\it The result of 4.12.12 remains true for $p=2$: any $2$-divisible group over an algebraically closed field $\tilde k$ of characteristic 2 has a uni plus versal deformation (to a $2$-divisible group) over a smooth, affine $W(\tilde k)$-scheme.}
\medskip
{\bf Proof:} This is a consequence of 4.12.12 and 2.3.18. We just need to remark that we can still define $\rho_N$ so that it can be used in the same manner, cf. the part of 2.3.18.2 referring to ``isolating aside" the slopes $0$ and $1$ (in 4.14 A below it is pointed out that we can still appeal to 4.1.1 and 4.2.3.1 for $p=2$).  
\medskip
Similarly, the variant of 4.12.12 mentioned in 4.12.12.0 and involving $m=1$, holds as well for $p=2$. The surjectivity part of the Claim of 2.4.1 implies (cf. also 2.3.18.1 D) that 4.12.9 1) holds for $p=2$ (even for principally quasi-polarized filtered $\sg_k$-crystals).
\medskip
{\bf 4.12.12.3. Tate deformation Hopf algebras and Shimura $p$-divisible groups.} Let $\tilde k$ be a perfect field of characteristic $p\ge 2$. Let $(\tilde M,\tilde F^1,\tilde\vph,\tilde G,(\tilde t_{\al})_{\al\in\Mj})$ be the not necessary quasi-split Shimura filtered $\sg_{\tilde k}$-crystal associated to a Shimura $p$-divisible group $\tilde D$ over $W(\tilde k)$. Let $\tilde d:=dd((\tilde M,\tilde\vph,\tilde G))$. Let $\tilde R$ be the $p$-adic completion of the henselization $\tilde R(h)$ of the localization of $W(\tilde k)[x_1,...,x_{\tilde d}]$ w.r.t. its prime ideal $(x_1,...,x_{\tilde d})$. There is (cf. Theorem 2 of 3.15.1) a uni plus versal Shimura $p$-divisible group $D_{\rm def}(\tilde D)$ over ${\rm Spec}(\tilde R)$, which in a $W({\tilde k})$-valued point $\tilde z$ of ${\rm Spec}(\tilde R)$ is $\tilde D$. Warning: $D_{\rm def}(\tilde D)$ is not uniquely determined; here uniqueness is up to automorphisms of $D_{\rm def}(\tilde D)$ and of $\tilde R$.
\smallskip 
We refer to $D_{\rm def}(\tilde D)$ as a Tate deformation Shimura $p$-divisible group of $\tilde D$. Writing $D_{\rm def}[p^n]={\rm Spec}(R_n(\tilde D))$, we refer to $R_n(\tilde D)$ as a level-$n$ Tate deformation Hopf algebra of $\tilde D$. We think it is an interesting (and extremely important) problem to study (classify) these finite, flat, Hopf $\tilde R$-algebras. 
\smallskip
We have a similar problem (much harder) to understand (classify) the similarly defined Hopf $\tilde R(h)$-algebras (for this problem we assume that at least one such $D_{\rm def}(\tilde D)$ is definable over $\tilde R(h)$ itself).
\medskip
{\bf 4.12.12.4. Exercise.} Use the proof of 4.12.12 to get a proof of the original Manin problem. Hint: using products of principally polarized abelian varieties, the problem gets reduced to the case when we are dealing with a Newton polygon $NP$ having only the slopes $\al$ and $1-\al$, with $\al\in (0,{1\over 2})$; writing $\al={a\over {a+b}}$, with $(a,b)=1$, $a,b\in\NN$, $b>a$, we can assume the multiplicities of $\al$ and $1-\al$ for $NP$ are both $a+b$ and so we can use the SHS's mentioned in the proof of 4.12.12 (or in 4.12.12.0).
\smallskip
It is worth pointing out that this proof does not appeal to [Mu, cor. 1 of p. 234]. 
\medskip
{\bf 4.12.12.5. Corollary.} {\it Let $k_1$ be an algebraically closed field of arbitrary positive characteristic $p$. Let $R_1$ be a complete, noetherian, local $W(k_1)$-algebra having $k_1$ as its residue field. Let $D_1$ be a $p$-divisible group over ${\rm Spec}(R_1)$. Then $D_1$ is a direct summand of the $p$-divisible group of a principally polarized abelian scheme $(A_1,p_{A_1})$ over ${\rm Spec}(R_1)$ of dimension not greater than the rank of $D_1$.}
\medskip
{\bf Proof:} If $D_1$ is \'etale or of multiplicative type, then we can take $(A_1,p_{A_1})$ to be the pull back (via the morphism ${\rm Spec}(R_1)\to {\rm Spec}(\ZZ_p)$) of the canonical lift of an ordinary abelian variety over $\FF_p$ which has a principal polarization and is of dimension equal to the rank of $D_1$. The other cases are a consequence of the part of 4.12.12.0 and of 4.12.12.2 referring to $m=1$; so based on Serre--Tate's deformation theory we can take $(A_1,p_{A_1})$ such that:
\medskip
a) its principally quasi-polarized $p$-divisible group is isomorphic to the direct sum of $D_1$ and of its dual, together with its natural principal polarization;
\smallskip
b) the projector which achieves the direct sum decomposition of a) is a $\ZZ_p$-linear combination of endomorphisms of $A_1$ which preserve this decomposition.
\medskip
{\bf 4.12.12.6. The existence of toric points.} The proof of 4.12.12 can be adapted to show the existence of ``as many types of toric points as expected". Let $(f,L_{(p)},v)$ be a SHS. Let $\om$ and ${\got C}_{\om}$ be as in 4.1.5. Let $(G_1,X_1,H_1,v_1)$ be a Shimura quadruple having the same adjoint as $(G,X,H,v)$. We have:
\medskip
{\bf Corollary.} {\it Any connected component $\Mc_1$ of the extension to $W(\FF)$ of the integral canonical model $\Mn_1$ of $(G_1,X_1,H_1,v_1)$ has $\FF$-valued points whose attached Shimura adjoint Lie $\bar\sg$-crystals are isomorphic to the Shimura adjoint Lie $\bar\sg$-crystal attached to ${\got C}_{\om}$.}
\medskip
{\bf Proof:}
As the connected components of $\Mn_{1W(\FF)}$ are permuted transitively by $G_1(\AA_f^p)$ (see [Va2, 3.3.2]), based on 4.9.2.1 and the proof of 4.9.8, we can assume $G_1$ is an adjoint group. We can assume $(f,L_{(p)},v)$ is such that $f$ is a good embedding w.r.t. $p$, cf. 2.3.6 (see [Va2, 6.4.2] for $p\ge 5$ and \S6 for its $p=3$ analogue). Starting from the pair $(T^{\om},\mu_{\om})$ of 4.1.5.4, we construct as in the proof of 4.12.12, a Shimura quadruple 
$$(T^{\om}_{\QQ},h_{\om}^{\rm poss},T_{\QQ_p}^{\om}(\ZZ_p),v_p(\om)^{\rm poss}),$$ 
with $T^{\om}_{\QQ}$ a maximal torus of $G$ which is a $\QQ$--form of the generic fibre of $T^{\om}$ and whose Zariski closure in $G_{\ZZ_{(p)}}$ is a torus, and with $T_{\QQ_p}^{\om}(\ZZ_p)$ as the hyperspecial subgroup of $T^{\om}_{\QQ_p}(\QQ_p)$. Let $X_{\om}^{\rm poss}$ be the $G(\RR)$-conjugacy class of $h_{\om}^{\rm poss}$. We get a Shimura pair $(G,X_{\om}^{\rm poss})$. The resulting injective map (cf. [De1, 1.6])
$$(G,X_{\om}^{\rm poss})\hookrightarrow (GSp(W,\psi),S)$$ 
is denoted by $f_{\om}^{\rm poss}$. Denoting by $v_{\om}^{\rm poss}$ the prime of $E(G,X_{\om}^{\rm poss})$ divided by $v_p(\om)^{\rm poss}$, we get another triple 
$$(f_{\om}^{\rm poss},L_{(p)},v_{\om}^{\rm poss});$$ 
2.3.6 and our assumption on $f$ imply that it is still a SHS. Let $\Mn_{\om}^{\rm poss}$ be its integral canonical model. 
\smallskip
$X_{\om}^{\rm poss}$ might not be $X$ itself; however, from [De2, 1.2.7-8], as $\mu_{\om}$ and $\mu$ are $G_{\ZZ_{(p)}}(W(\FF))$-conjugate, we get that $(G^{\rm ad},X_{\om}^{{\rm poss}\, {\rm ad}})=(G^{\rm ad},X^{\rm ad})$. So the Corollary follows, once we show that any $W(\FF)$-valued point $z_0$  of $\Mn_{\om}^{\rm poss}$ factoring through the integral canonical model of $(T^{\om}_{\QQ},h_{\om}^{\rm poss},T^{\om}(\ZZ_p),v_p(\om)^{\rm poss})$ (see [Va2, 3.2.8]), has attached to it a Shimura adjoint Lie $\bar\sg$-crystal isomorphic to the Shimura adjoint Lie $\bar\sg$-crystal attached to ${\got C}_{\om}$.
\smallskip
As this is a statement which can be reformulated purely in the \'etale $\ZZ_p$-context, in order to benefit from previous notations, we argue it by working in the context of the SHS $(f,L_{(p)},v)$ itself. So we use the notations of 4.2.1-4, with $z_0$ factoring through the integral canonical model $\Mt$ of a special quadruple $(T,\{h\},H_T,v_T)\hookrightarrow (G,X,H,v)$ having the property that the Zariski closure of $T$ in $G_{\ZZ_{(p)}}$ is a maximal torus. So $V_0:=W(\FF)$. Denoting by $f|T$ the restriction of $f$ to $(T,\{h\})$, the triple $(f|T,L_{(p)},v_T)$ is a SHS (cf. 2.3.6 and [Va2, 4.3.13]). We apply 4.1.1 and 4.2.3-4 to it. Let $(v_{\al})_{\Mj_T}$, with $\Mj_T$ a set containing $\Mj^\prime$, be the set of tensors of $\Mt(W^*)$ formed by adding the elements of ${\rm Lie}(T)$ to $\Mj^\prime$. Working in the \'etale context we construct abstractly (as in 4.1.1) starting from $(f|T,L_{(p)},v_T)$, a Shimura filtered $\bar\sg$-crystal 
$$(M_T\otimes_{W(k(v_T))} V_0,F^1_T\otimes_{W(k(v_T))} V_0,\vph_T,T_{V_0},(v_{\al})_{\al\in\Mj_T}).$$ 
Working in the crystalline cohomology context, we consider (as in 4.2.3-4) the Shimura filtered $\bar\sg$-crystal $(M_0,F^1_0,\vph_0,T_{V_0},(t_{\al})_{\al\in\Mj_T})$ attached to the $V_0$-valued point (defined by) $z_0$ of $\Mt$. Using the fact that two maximal tori of $G_{V_0}^{\rm ad}$ are $G_{V_0}^{\rm ad}(V_0)$-conjugate and that $N_{T_1}(K_0)=N_{T_1}(V_0)T_1(K_0)$, with $T_1$ as the image of $T_{V_0}$ in $G_{V_0}^{\rm ad}$ and with $N_{T_1}$ as its normalizer in $G_{V_0}^{\rm ad}$, we deduce that we can assume that, in the present situation of $(f,L_{(p)},v)$ and $z_0$,  $g_0$ of 4.2.4 is a $K_0$-valued point of $T_1$. So the Shimura adjoint filtered Lie $\bar\sg$-crystal $({\rm Lie}(G_{V_0}^{\rm ad}),\vph_0)$ defined by the natural inclusion of $G_{V_0}$ in $GL(M_0)$ is inner isomorphic to the Shimura adjoint filtered Lie $\bar\sg$-crystal $({\rm Lie}(G_{V_0}^{\rm ad}),g\vph_T)$ defined by the natural inclusion of $G_{V_0}$ in $GL(M_T\otimes_{W(k(v_T))} V_0)$, with $g\in T_1(V_0)$. Using $\ZZ_p$-structures (as in 2.2.9 8) and as in the proof of b) of 4.4.1 2)), we can assume $g$ is the identity element of $T_1(V_0)$ (variant: the proof of 3.11.1 c) can be entirely adapted to get that $({\rm Lie}(G_{V_0}^{\rm ad}),\vph_T)$ and $({\rm Lie}(G_{V_0}^{\rm ad}),g\vph_T)$ are isomorphic under an isomorphism defined by an element of $T_1(V_0)$). This ends the proof.
\medskip
{\bf 4.12.12.6.0. Remarks.} {\bf 1)} If any homomorphism $h_s:{\rm Res}_{\CC/\RR} \GG_m\to GSp(W\otimes_{\QQ} {\RR},\psi)$ defined by an element of $S$ (we recall that $f:(G,X)\hookrightarrow (GSp(W,\psi),S)$) and factoring through $G_{\RR}$, is defined by an element of $X$, then $X_{\om}^{\rm poss}=X$, $\forall\om\in W_G$. 
\smallskip
[Ko2, 4.3] points out: this is always the case if $f$ is a PEL type embedding, with all simple factors of $G^{\rm ad}$ of some $A_n$ or $C_n$ Lie type. Using [De2, 1.2.8] we can construct many more examples. For instance, the above paragraph applies in the context of [Va2, 5.7.5]. Warning: [Va2, 2.5.1] points out that not always $h_s$ is defined by an element of $X$.
\smallskip
{\bf 2)} [Va2, 3.2.7.1 f) and 3.3.3] imply that replacing the injective map
$$(T^{\om}_{\QQ},h_{\om}^{\rm poss},T^{\om}(\ZZ_p),v_p(\om)^{\rm poss})\hookrightarrow (G,X_{\om}^{\rm poss},H,v_{\om}^{\rm poss})$$ 
by its composite with an isomorphism 
$$(G,X_{\om}^{\rm poss},H,v_{\om}^{\rm poss})\tilde\to (G,X,H,v)$$ 
defined by an element of $G^{\rm ad}_{\ZZ_{(p)}}(\ZZ_{(p)})$, we can always assume $X_{\om}^{\rm poss}=X$ and so $v_{\om}^{\rm poss}=v$. So in the above proof of 4.12.12.6, references to 4.9.8 can be entirely avoided, as long as we are just dealing with $(f,L_{(p)},v)$.
\medskip
The operations of the proof of 4.12.12.6 keep track of inner automorphisms of ${\rm Lie}(G^{\rm ad}_{W(k)})$ (cf. the def. of $T_{\mu_{\om}}$ in 4.1.5.4 and 4.2.10) and so we conclude (cf. 4.1.5.5):
\medskip
{\bf 4.12.12.6.1. Corollary.} {\it The number of strata of the quasi-ultra stratification of $\Mn_{1k(v_1)}$ is at least equal to the number of elements of the quotient set $W_G/R_G(v)$.}
\medskip
{\bf 4.12.12.6.2. Examples.} {\bf A.} We first consider the case when all simple factors of $G^{\rm ad}$ are of $A_1$ Lie type. So the number of non-compact simple factors of $G^{\rm ad}_{\RR}$ is $\dim_{\CC}(X)=\dim_{\CC}(X_1)$. From 3.13.7.6.0 and 4.12.12.6.1 we get that the number of strata of the quasi-ultra stratification of $\Mn_{1k(v_1)}$ is precisely $2^{\dim_{\CC}(X_1)}$. 
\smallskip
These strata are in one-to-one correspondence to functions 
$$f_s: S(1,\dim_{\CC}(X_1))\to \{0,1\}$$
 in such a way that the dimension of the stratum ${\got s}$ corresponding to $f_s$ is precisely 
$$\dim_{\CC}(X_1)-\sum_{l=1}^{\dim_{\CC}(X_1)} f_s(l)$$ 
(cf. 3.13.7.6.0 and 4.5.15.2.1). A stratum ${\got s}_1$ corresponding to $f_{s_1}$ specializes to another stratum ${\got s}_2$ corresponding to $f_{s_2}$ iff $f_{s_2}(l)\ge f_{s_1}(l)$, $\forall l\in S(1,\dim_{\CC}(X_1))$ (cf. 3.13.7.6.0). In particular there is only one open (resp. closed stratum: it is of dimension $0$). In fact this closed stratum has a non-empty intersection with each connected component of $\Mn_{1k(v_1)}$. To see this it is enough to show that the connected components of $\Mn_{1k(v_1)}$ are permuted transitively by $G_1(\AA_f^p)$. It is enough to show --cf. [Va2. 3.3.2]-- that each connected component of $\Mn_{1W(\FF)}$ has a connected special fibre.
\smallskip
We can assume $(G^{\rm ad},X^{\rm ad})$ is (cf. the proof of 4.9.8) simple and is not (cf. [Va2, 6.4.11]) of compact type. So [BHC, 11.4 and 11.6] implies: $(G^{\rm ad},X^{\rm ad})$ is the adjoint variety of a Hilbert--Blumenthal variety. Based on the proof of 4.9.8 we can assume $(f,L_{(p)},v)$ is a SHS associated to a Hilbert--Blumenthal variety. But this case follows from [Ch1]. As a conclusion: for the $A_1$ Lie type, the star of 4.2.8.1 can be removed. 
\smallskip
This example represents the logical generalization of [GO]: loc. cit. deals only with Hilbert--Blumenthal varieties.
\medskip
{\bf B.} We assume all simple factors of $(G^{\rm ad},X^{\rm ad})$ are such that:
\medskip
\item{{\bf ASS}} {\it $\forall i\in\Mh^{\rm nc}$, the $i$-th cyclic factor of a Shimura-ordinary point of $\Mn_{k(v)}$ is one among the list of 3.13.7.6.2.} 
\medskip
For instance this is the case if:
\medskip
\item{{\bf EX}} All simple factors of $(G^{\rm ad},X^{\rm ad})$ are of $B_n$, $C_n$ or $D_n^{\RR}$ type, or of $A_{2l+1}$ Lie type such that all non-compact, simple factors of $G^{\rm ad}_{\RR}$ are isomorphic to $SU(l+1,l+1)_{\RR}^{\rm ad}$, or of $D_{2l}^{\HH}$ type which involve groups which split over $\QQ_p$, or of $A_{2l}$ or $D_{2l+1}^{\HH}$ type which involve groups which over $\QQ_p$ are absolutely simple but do not split ($l\in\NN$). 
\medskip
From 3.13.7.6.1 C, 4.5.15.2.1 and 4.12.12.6.1 we get (under ASS; cf. also the list of 3.13.7.6.2 and the cyclic diagonalizability property of 3.13.7.6.3):
\medskip
{\bf Corollary.} {\it There is a reduced, closed subscheme ${\got s}_1^0$ of $\Mn_{1k(v_1)}$ of dimension $0$ such that any stratum of the quasi-ultra stratification of $\Mn_{1k(v_1)}$ of dimension $0$ is a scalar extension of ${\got s}_1^0$. All points of ${\got s}_1^0$ are toric. Moreover, the Shimura adjoint Lie $\sg_{\bar k}$-crystals attached to $\bar k$-valued points of ${\got s}_0$ are all inner isomorphic. 
\smallskip
Similarly, all Shimura $\sg_{\bar k}$-crystals endowed with natural families of tensors (indexed by the set $\Mj^\prime$ of 2.3.1) attached to $\bar k$-valued points of ${\got s}^0$ belonging to a fixed connected component of $\Mn_{k(v)}$ are all $1_{\Mj^\prime}$-isomorphic; here ${\got s}^0$ has the same significance for $\Mn_{k(v)}$ as ${\got s}_1^0$ for $\Mn_{1k(v_1)}$.}
\medskip
{\bf C.} We briefly go through two examples of pivotal points which are not among those of B, just to illustrate what type of computations we will have to perform in \S 9-10 in order to handle the cases not covered by B. We consider a SHS $(f,L_{(p)},v)$. We assume $G^{\rm ad}_{\RR}$ is $SU(a,b)^{\rm ad}$ and that $G^{\rm ad}$ splits over $\ZZ_p$. 
\smallskip
If $c:={b\over a}\in\NN$, then any $W(k)$-valued toric point of $\Mn_{k(v)}$ whose Shimura adjoint Lie $\sg_k$-crystal is isomorphic to $({\rm End}(M_1),\vph_1)$, where $(M_1,\vph_1)$ is the direct sum of $a$ copies of the $\sg_k$-crystal $(M_2,\vph_2)$ which has a lift $(M_2,F^1_2,\vph_2)$ which is (isomorphic to) the circular diagonalizable filtered $\sg_k$-crystal having (see B of 2.2.22 3)) the type $(1,0,0,...,0)$ (the number of $0$'s is $c$), is pivotal. This is a consequence of the Proposition of 4.5.15.2.1. We just need to remark that:
\medskip
-- the group of automorphisms $GA_2$ of the truncation mod $p$ of $(M_2,\vph_2)$ has (this is an easy exercise) dimension $c$ (and so the group of automorphisms of $(M_1,\vph_1)$ has dimension $ca^2=ab=\dim_{\CC}(X)$);
\smallskip
-- the group of automorphisms $GA_2^\prime$ of the Faltings--Shimura--Hasse--Witt map of $(M_2,\vph_2,GL(M_2))$ has the same dimension as $GA_2$ (the natural homomorphism $GA_2\to GA_2^\prime$ has finite kernel and from 3.13.7.1.1 and 3.13.7.2 we get that it is an epimorphism). 
\medskip
Similarly, if $(a,b)=(2,3)$ one can check that starting from the type $(1,0,1,0,0)$ and taking $(M_1,\vph_1)=(M_2,\vph_2)$ we get pivotal points.  
\medskip
{\bf 4.12.12.6.3. The null strata.} We refer to 4.9.8. The stratum of the refined Lie canonical stratification of $\Mn_{k(v)}$ whose geometric points are such that the Newton polygons of their attached Shimura adjoint Lie $F$-crystals are having only the slope $0$, is called the null stratum of $\Mn_{k(v)}$; its points with values in fields are called null points (of $\Mn_{k(v)}$). It is the natural generalization of the usual supersingular strata of special fibres of integral canonical models of Siegel modular varieties. We now prove 1.12.1 A) (in some sense we just refine a little bit the proof of 4.3.6).
\medskip 
Based on the proof of 4.9.8, we can assume we are dealing with the integral canonical model $\Mn_{k(v)}$ of a SHS $(f,L_{(p)},v)$. So we refer to 4.1 and 4.3.1 with $k=\FF$; we assume $B$ contains $T$. Based on 4.12.12.6.0, it is enough to show the existence of $\om_0\in N(W(\FF))$ (see 4.1.5) such that the Shimura Lie $\sg_{\FF}$-crystal attached to ${\got C}_{\om_0}$ has all slopes equal to $0$. First, we consider the examples provided by EX of 4.12.12.6.2 B. For each $i\in\Mh^{\rm nc}$, let $\om_0^i\in G_{iW(\FF)}(W(\FF))$ normalizing the image of $T_{W(\FF)}$ in $G_{iW(\FF)}$ and such that:
\medskip
-- looking at the simple factors of $G_{iW(\FF)}$ only one component in them of $\om_0^i$ is non-trivial; let $G_{i_0}$ be this factor corresponding to the non-trivial component $\om_{00}^i$ of $\om_0^i$ (it is defined by picking up an element ${i_0}$ of the set $\Mh_i$ of 4.3.1.1) (warning: we do not assume $G_{i_0}$ is such that the image of $\mu_{W(\FF)}$ in it is non-zero);
\smallskip
-- $\om_{00}^i$ takes the image of $B_{W(\FF)}$ in $G_{i_0}$ into its opposite w.r.t. the image of $T_{W(\FF)}$ in $G_{i_0}$.
\medskip
We take $\om_0\in N(W(\FF))$ such that its image in $G_{iW(\FF)}$ is $\om_0^i$, $\forall i\in\Mh^{\rm nc}$. The Newton polygon of the Shimura Lie $\sg_{\FF}$-crystal attached to ${\got C}_{\om_0}$ has only one slope: it is $0$; this is an immediate consequence of the choice of $\om_0$ (see the numbering in 3.4.3.2 and the item XI in [Bou2, planche I to IV]). This takes care of the $B_n$, $C_n$ and $D_n^{\RR}$ types entirely. To take care of the other types, we can assume $(G^{\rm ad},X^{\rm ad})$ is simple of $A_n$ ($n\ge 2$) or $D_n^{\HH}$ ($n\ge 4$) type. It is enough to deal with a fixed $i\in\Mh^{\rm nc}$. 
\smallskip
Second, we assume that $G^i$ of 4.3.1.1 is either non-split (i.e. we are in the case when the $i$-th cyclic adjoint factor associated to a $G$-ordinary point of $\Mn_{k(v)}$ is with involution) of $A_n$ or $D_{2n+1}$ Lie type or it is split of $D_{2n}$ Lie type, with $n\ge 2$. In this case the same choice of $\om_i$ (and the same argument) makes the thing work as in the case when we were working with the examples provided by EX of 4.12.12.6.2 B. 
\smallskip
So we are left with the cases: $G^i$ is split of $A_n$ or $D_{2n+1}$ Lie type or is non-split of $D_{2n}$ Lie type, with $n\ge 2$. We use the same type of Weyl element $\om_0^i$: we assume ${i_0}$ has been fixed and we still denote by $\om_{00}^i\in G_{i_0}(W(\FF))$ the (unique) non-trivial component of $\om_0^i$. Let $j_i\in\NN$ be the smallest number such that the cocharacter $\om_0^i\mu_{iW(\FF)}({\om_0^{i}})^{-1}$ of $T_{W(\FF)}$ is fixed by $\bar\sg^{j_i}$. It is enough to show that we can choose $\om_{00}^i$ such that the product 
$$PDIT:=\prod_{s=1}^{j_i} \bar\sg^s\om_0^i\mu_{iW(\FF)}({\om_0^{i}})^{-1}\bar\sg^{-s}$$ 
is trivial (one can view this additively as well). Warning: here, in order to benefit from 4.1.2 (IT), we choose a Frobenius endomorphism of $M\otimes_{W(k(v))} W(\FF)$ of the form $\bar\sg\om_0^i\mu({1\over p})$ and not of the usual form $\om_0^i\bar\sg\mu({1\over p})$.
\smallskip
Third, we assume $G^i$ is split of $A_n$ Lie type. We identify $G_{i_0}$ with the adjoint group of the $SL_{n+1}$-group $G_{i_0}^{\rm sc}$ of a free $W(\FF)$-module $M_{i_0}$ of rank $n+1$. We consider a $W(\FF)$-basis $\{e_1,...,e_{n+1}\}$ of $M_{i_0}$ such that:
\medskip
-- each element of it is normalized by the maximal torus of $G_{i_0}^{\rm sc}$ which has the same image in $G_{i_0}$ as $T_{W(\FF)}$;
\smallskip
-- $\forall i\in S(1,n+1)$, $<e_1,...,e_i>$ is normalized by the Borel subgroup of $G_{i_0}^{\rm sc}$ which has the same image in $G_{i_0}$ as $B_{W(\FF)}$.
\medskip
We take $\om_{00}^i$ to be the image in $G_{i_0}(W(\FF))$ of any element $\om_0^{i\rm sc}$ of $G_{i_0}^{\rm sc}(W(\FF))$ which takes $<e_i>$ into $<e_{i+1}>$, $i=\overline{1,n+1}$, with $e_{n+2}:=e_1$. The fact that $PDIT$ is trivial is a consequence of the circular property of $\om_0^{i\rm sc}$ as well as of 3.4.2.1-2.
\smallskip
Last, we assume $G_i$ is split of $D_{2n+1}$ Lie type or non-split of $D_{2n}$ Lie type, with $n\ge 2$. The proof of 4.3.6 (via the proof of Fact 2 of it) allows us to shift our attention to the non-split $A_3$ Lie type and so the above part applies. This ends the proof of 1.12.1 A).
\medskip
{\bf 4.12.12.6.4. The strong form of the generalized Manin problem.} We now prove 1.12.1 B) and C). The arguments of 4.12.12.6.3 were presented in the abstract form. So they apply in the general context to give us that in each class of Shimura $\sg_{k}$-crystals there is an element defined by a Shimura $\sg_k$-crystal whose attached Shimura Lie $\sg_k$-crystal has all slopes $0$. We have:
\medskip
{\bf Corollary.} {\it We refer to the notations of 3.9.7.2 (with $k=\bar k$). There is $w_0\in R_0(W(k))$ normalizing $T$ which depends only on the Shimura adjoint Lie $\sg_k$-crystal attached to $(M,g\vph,R_0)$ such that the Newton polygons of $(M,\vph)$ and of $(M,hg\vph)$ are the same as of the cyclic diagonalizable Shimura $\sg_k$-crystal $(M,w_0g\vph,T)$.}
\medskip
{\bf Proof:} Let $w_0\in R_0(W(k))$ normalizing $T$ be such that the Shimura adjoint Lie $\sg_k$-crystal attached to the Shimura $\sg_k$-crystal $(M,w_0g\vph,R_0)$ has all slopes $0$. The Corollary follows from the Corollary of 3.9.7.2 and from the class invariant part of 2.2.24.1 applied to $Cl(M,g\vph,R_0)$. 
\medskip
Corollary can be entirely adapted to the generalized Shimura context provided by the $\rho$-representations of 4.5.4 (cf. 3.9.7.2.1 1)). So 1.12.1 B) follows from it and from 4.12.12.6.0. To detail this, as in 4.12.12.6.3 we can assume we are dealing with the integral canonical model $\Mn_{k(v)}$ of a SHS $(f,L_{(p)},v)$. We refer to 4.5.4. We can assume $y$ is not a null point (cf. 1.12.1 A)). We can assume that the parabolic subgroup $P_y$ of $G_{W(k)}$ whose Lie algebra is $W_0({\rm Lie}(G_{W(k)}),\vph_y)$ has a maximal torus $T_y$ normalizing $F^1_z$ (cf. 2.3.17 and 3.9.7.1). We can assume $\mu_y$ factors through $T_y$. Let $R_y$ be the Levi subgroup of $P_y$ containing $T_y$. Let $UR_y$ be the unipotent radical of $P_y$. Let $g_y\in P_y(W(k))$ be such that the quadruple $(M_y,F^1_z,g_y\vph_y,T_y)$ is a Shimura filtered $\sg_k$-crystal. Let $w_y\in R_y(W(k))$ normalizing $T_y$ and such that all slopes of $({\rm Lie}(R_y),w_yg_y\vph_y)$ are 0. Writing $\vph_y=a_y\mu_y({1\over p})$, we consider a $\ZZ_p$-structure $(N_y,G_{\ZZ_p},(t_{\al})_{\al\in\Mj^\prime})$ as in 2.2.9 8) (so $N_y$ is the $\ZZ_p$-submodule of $M_y$ formed by elements fixed by $a_y$); w.r.t. it we can write $\sg_k$ instead of $a_y$. 
\smallskip
For what follows, (in connection to the $\rho$-stratifications of 4.5.4) it is more convenient to work directly with the $\ZZ_p$-structures provided by 2.2.9 8), without appealing to d) of 4.4.1 3). $_{\rho}{\got C}_z$ (of 4.5.4) is isomorphic to $(N_y\otimes_{\ZZ_p} W(k),\sg_k\rho_{B(k)}(\mu_y({1\over p})))$. So its Newton polygon is the same as the Newton polygon of the following generalized Shimura $p$-divisible object
$$_{\rho}{\got C}_y^{\rm mod}:=\Bigl(N_y\otimes_{\ZZ_p} W(k),\rho_{W(k)}(w_yg_y)\sg_k\rho_{B(k)}(\mu_y({1\over p})),{\rm Im}(\rho_{W(k)})\Bigr)$$ 
(cf. the mentioned adaptations of the Corollary; so this is a consequence of the $P_{>0}(B(k))$-conjugacy part of 3.9.7.2 --$P_{>0}$ of 3.9.7.2 is denoted here by $UP_y$-- and of 2.2.24.1 applied to $\bigl(N_y\otimes_{\ZZ_p} W(k),\rho_{W(k)}(*)\sg_k\rho_{B(k)}\bigl(\mu_y({1\over p})\bigr),R_y\bigr)$, with $*\in\{g_y,w_yg_y\}$), for any $\rho$ as in 4.5.4. The upper index ``mod" stands for modification (of the non-filtered version of $_{\rho}{\got C}_z$).
\smallskip
From 4.12.12.6 we deduce the existence of a $k$-valued toric point $y_1$ of $\Mn_{k(v)}$ such that its attached Shimura adjoint Lie $\sg_k$-crystal is inner isomorphic to the one attached to $(M_y,F^1_z,w_yg_y\vph_y,G_{W(k)})$ (the $\ZZ_p$-structures of 3.11.2 B applied in the context of $(M_y,F^1_z,w_yg_y\vph_y,T_y)$ allow us to work over $k$ and not just over $\FF$; of course, here we can assume that $k=\FF$). But it is an easy exercise to see that the Newton polygon of $_{\rho}{\got C}_{z_1}$, with $z_1$ a $W(k)$-valued point of $\Mn$ lifting $y_1$ and such that its attached Shimura filtered $\sg_k$-crystal is cyclic diagonalizable, is the same as the Newton polygon of $_{\rho}{\got C}_y^{\rm mod}$, for any $\rho$ as in 4.5.4. It can be solved as follows. If we assume d) of 4.4.1 3), then we just need to use arguments with semisimple elements as in 2.2.24.1 and as in the proof of 4.2.7, in the context of the class $Cl(_{\rho}{\got C}_y^{\rm mod})$. If we do not assume d) of 4.4.1 3), then we just need to add that the rational variant of 4.4.12 obtained by combining c) of 4.4.1 3) with 2.3.13.1, allows us still to appeal to the mentioned arguments. We conclude: $y$ and $y_1$ are $k$-valued points of the same stratum of the absolute stratification of $\Mn_{k(v)}$.
\smallskip
1.12.1 C) follows from 1.12.1 B) and the above paragraph (cf. also 4.1.5-6, 4.2.10 and 4.5.11). As in the elimination process of 3.13.7.3.1 B, we get that in fact in 1.12.1 C) we can replace the Weyl group of $G^1_{\CC}$ by the Weyl group of the extension to $\CC$ of the product of the simple, non-compact factors of $G^{1\rm ad}_{\RR}$.
\medskip
{\bf Exercise.} Show that the map $f_G$ of 4.5.11 a) is surjective. Hint: combine the Corollary with 4.5.15.2.5 or with the arguments of the proof of 4.4.1 referring to the independence part of b) of 4.4.1 3).
\medskip
{\bf 4.12.12.6.4.1. Examples.} The estimate of 1.12.1 C) is very gross in general. However, there are examples in which its bound is attained (for instance, this happens when all factors of $G^{\rm ad}_{\ZZ_p}$ are split of $A_1$ Lie type and $\Mh^{\rm c}$ is the empty set). We assume now that we have a SHS $(f,L_{(v)},p)$, with $(G^{\rm ad},X^{\rm ad})$ simple of $B_n$ Lie type and with $G^{\rm ad}_{\RR}$ having at most $2$ non-compact factors. We have:
\medskip
{\bf Corollary.} {\it The absolute stratification of $\Mn_{k(v)}$ has:
\smallskip
a) $n+1$ strata if $G^{\rm ad}_{\RR}$ has only one non-compact factor;
\smallskip
b) $(n+1)^2$ strata if $G^{\rm ad}_{\RR}$ has two non-compact factors and $\Mh^{\rm nc}$ has $2$ elements;
\smallskip
c) $n([{n\over 2}]+2)-[{n\over 2}]([{n\over 2}]+1)$ strata if $G^{\rm ad}_{\RR}$ has two non-compact factors and $\Mh^{\rm nc}$ has $1$ element.}
\medskip
{\bf Proof:} We use 4.5.6.2 A and 4.12.12.6.4. Let $m$ (resp. $l$) be the number of simple (resp. of non-compact, simple) factors of $G^{\rm ad}_{\RR}$. We refer to [Bou2, planche II] for the description of the Weyl group of the $B_n$ Lie type; below we use it freely. Everything boils down to the following situation. For $s\in S(1,m)$, we consider a free $W(\FF)$-module $V_s$ of rank $2n+1$ and a $W(\FF)$-basis $\Mb_s=\{e_1^s,e_2^s,...,e_{2n+1}^s\}$ of it. We consider a cocharacter $\mu$ of $GL(V_1\oplus_{W(\FF)}...\oplus_{W(\FF)} V_m)$ which:
\medskip
-- acts trivially on $e_j^s$, if $j\ge 3$ or if $s>l$;
\smallskip
-- acts as the multiplication with $p^j$ on $e_{{3-j}\over 2}^s$, if $s\in S(1,l)$ and $j\in S(1,2)$.
\medskip
 We have a permutation $\pi$ of $S(1,m)$ and we look at the Newton polygons of $h\mu(p)$ where 
$$h\in GL(V_1\oplus_{W(\FF)}...\oplus_{W(\FF)} V_m)(W(\FF))$$ 
fixes $e^s_{2n+1}$ and takes each pair $(e_{2q_1-1}^s,e_{2q_1}^s)$ into a pair of one of the following two forms $(e_{2q-1}^{\pi(s)},e_{2q}^{\pi(s)})$ or $(e_{2q}^{\pi(s)},e_{2q-1}^{\pi(s)})$, with $q$, $q_1\in S(1,n)$, $\forall s\in S(1,m)$. Here $h\mu(p)$ is viewed as a $B(\FF)$-valued point of $GL(V_1\oplus_{W(\FF)}...\oplus_{W(\FF)} V_m)$ and Newton polygons are obtained as usual (using ${1\over q}$ times valuations of the eigenvalues of $h^q$, with $q\in\NN$ big enough). For each non-empty subset $C$ of $S(1,m)$ permuted transitively by $\pi$ we want to first count the number $n(C)$ of Newton polygons of the restriction $r(C)$ of $h\mu(p)$ to $V(C):=(\oplus_{s\in C} V_s)[{1\over p}]$, with $h$ varying as described, and second we want to take the product of all these numbers $n(C)$. All Newton polygons we get are $0$-symmetric and so it is enough to deal with non-negative slopes. 
\smallskip
If $C\cap S(1,l)$ is the empty set, then $n(C)=1$. If $C\cap S(1,l)$ has only one element $s_0$, then the Newton polygon of $r(C)$ has either only slope $0$ or the only positive slope of it is ${1\over {\abs{C}s}}$ with multiplicity $\abs{C}s$, where $s$ is the number of pairs of the form $(e_{2q-1}^{s_0},e_{2q}^{s_0})$, $q\in S(1,n)$, or of the form $(e_{2q}^{s_0},e_{2q-1}^{s_0})$, $q\in S(2,n)$, which are in the orbit or $(e_1^{s_0},e_{2}^{s_0})$. So $n(C)=n+1$. This takes care of a) and b). 
\smallskip
If $l=2$ and $S(1,2)\subset C$, then we split the computation of $n(C)$ as follows. We count:
\medskip
i) $1$ for the Newton polygon having all slopes $0$;
\smallskip
ii) $n$ for Newton polygons which have only one positive slope $2\over {\abs{C}a}$ with multiplicity $a\abs{C}$, $a\in S(1,n)$; they correspond to the case when $e_1^1$ and $e_1^2$ are in the same orbit under $h$ but this orbit does not contain $e_2^1$ (and so does not contain $e_2^2$);
\smallskip
iii) $n-1$ for Newton polygons which have only one positive slope $1\over {\abs{C}a}$ with multiplicity $a\abs{C}$, $a\in S(1,n-1)$; they correspond to the case when the orbit of $e_j^1$ under $h$ does not contain $e_{2-j}^1$ or $e_j^2$ or $e_{2-j}^2$ but $e_{2-j}^2$ and $e_{j}^2$ are in the same orbit under $h$ for some fixed $j\in S(1,2)$;
\smallskip
iv) $1$ for each pair $(a,b)\in S(1,n)\times S(1,n)$, with $a+b\le n$ and $a<b$. 
\medskip
For each pair as in iv) we get a Newton polygon which has $2$ positive slopes ${1\over {\abs{C}a}}$ and ${1\over {\abs{C}b}}$ with multiplicities $a\abs{C}$ and respectively $b\abs{C}$; it corresponds to the cases when the orbit of $e_1^j$ under $h$ does not contain $e^j_2$ or $e_1^{3-j}$ or $e^{3-j}_2$ but has $a_j$ elements, $j=\overline{1,2}$, with $(a_j,b_j)$ equal to $(a,b)$ or with $(b,a)$. Warning: the case corresponding to $a=b$ gives birth to Newton polygons which were already counted (see ii)). We end up 
with
$$1+n+(n-1)+\sum_{j=1}^{[{n\over 2}]} (n-2j)$$
distinct Newton polygons. This sum is exactly the number of c). This ends the proof.
\smallskip
If $n=1$ (or $n=2$ or $n=3$ or $n=4$ or $n=5$) then in c) we get 2 (or 4 or 7 or 10 or 14) strata. 
\medskip
{\bf 4.12.12.6.4.2. Remark.} We refer to 4.5.8 1). $N(G^{\rm ad},X^{\rm ad},v^{\rm ad})$ is the number of strata of the canonical Lie stratification of $\Mn_{k(v)}$ and $N_1(G^{\rm ad},X^{\rm ad},v^{\rm ad})$ is the number of strata of the refined canonical Lie stratification of $\Mn_{k(v)}$. This follows from 1.12.1 B). So, in 4.9.8 a) we can replace ``indexed by a subset of $RLNP(G^{\rm ad},X^{\rm ad},v^{\rm ad})$" by: indexed by $RLNP(G^{\rm ad},X^{\rm ad},v^{\rm ad})$.
\medskip
{\bf 4.12.12.6.4.3. Remark.} We consider a $p$-divisible group with a reductive structure $(\tilde M,\tilde\vph,\tilde G)$ of $p-\Mm\Mf_{[a,b]}(W(k))$. Then there is $\tilde g\in \tilde G(W(k))$ such that $(\tilde M,\tilde g\tilde\vph,\tilde G)$ has a lift which is cyclic diagonalizable and all the slopes of its attached $b-a$-Lie $\sg_k$-crystal are zero. The proof of this is entirely the same as 4.12.12.6.3. We can assume $\tilde G_{B(k)}$ is adjoint and the Lie algebras of its simple factors are permuted transitively by $\tilde\vph$. As in 3.2.3, we can assume there is a maximal torus $\tilde T$ of $\tilde G$ such that its Lie algebra is normalized by $\tilde\vph$; so we look for $\tilde g$'s normalizing $\tilde T$. If $\tilde G$ is of $G_2$, $F_4$, $E_7$, $E_8$, $B_n$ or $C_n$ Lie type, then (cf. item IX in [Bou2, planches II, III, and VI to IX]) we can use the same type of Weyl elements as in the part of 4.12.12.6.3 referring to 4.12.12.6.2 EX. If $\tilde G$ is of $A_n$ or $D_n$ Lie type, we can use the same type of Weyl elements as in the cases of 4.12.12.6.3 corresponding to the same Lie types. 
\smallskip
If $\tilde G$ is of $E_6$ Lie type we refer to [Va6] and [Va8]. Here we will just mention that in many cases, the situation gets reduced to the $A_5$ Lie type and so we can apply the previous paragraph (the ``embedding" of $A_5$ into $E_6$ is provided by ``removing" the second node of the Dynkin's diagram of the $E_6$ Lie type). In particular, this reduction is possible if $({\rm Lie}(\tilde G),\tilde\vph)$ is a cyclic Shimura adjoint Lie $\sg_k$-crystal of $E_6$ type.  So we conclude: 
\medskip
{\bf Corollary (the abstract strong --Weyl-- form of the generalized Manin problem).} {\it The Corollary of 4.12.12.6.4 holds for the generalized Shimura context.}
\medskip
{\bf 4.12.12.6.5. Extra types of points.} We refer to 4.12.12.6. A point $y$ of $\Mn_{1k(v_1)}$ with values in an algebraically closed field is called:
\medskip
-- non-integral, if its Lie $p$-rank (defined as in 4.3.6, via 4.9.17) is $0$;
\smallskip
-- quasi-pivotal, if it is a toric and a non-integral point at the same time;
\smallskip
-- super-null, if it is a toric and a null point and the same time;  
\smallskip
-- pivotal, if it is super-null and belongs to a stratum of the quasi-ultra stratification of dimension $0$;
\smallskip
-- non-null, if it is not a null point;
\smallskip
-- partially null, if there is $i\in\Mh^{\rm nc}$ such that the $i$-th cyclic adjoint factor attached to $y$ has all slopes $0$;
\smallskip
-- totally non-null, if it is not partially null;
\smallskip
-- partially $U$-ordinary (resp. partially Shimura-ordinary), if there $i\in\Mh^{\rm nc}$ such that the $i$-th cyclic adjoint factor attached to $y$ is $U$-ordinary (resp. has a lift which is of parabolic type); 
\smallskip
-- quasi $U$-ordinary, if it is a toric and totally non-null point at the same time;
\smallskip
-- $isom$-canonical, if its attached Shimura adjoint $F$-crystal is $isom$-canonical;
\smallskip
-- quasi-final, if its attached Shimura adjoint $F$-crystal is quasi-final.
\medskip
By natural passage, we use the same terminology for points of $\Mn_{1k(v_1)}$ with values in a field. Each null point is non-integral and so each super-null point is quasi-pivotal. The toric points of the Corollary of 4.12.12.6.2 B are pivotal.
\medskip
{\bf Definition.} A lift $z\in\Mn_{1W(k)}$ of a partially $U$-ordinary point $y\in\Mn_{1k}$ is called a partial $U$-canonical lift, if $\forall i\in\Mh^{\rm nc}$ such that the $i$-th cyclic adjoint factor attached to $y$ is $U$-ordinary, the lift of it defined naturally by $z$ is a $U$-canonical lift. 
\medskip
{\bf 4.12.12.6.6. More on quasi $U$-ordinariness and on the $A_1$ Lie type case.} 
\smallskip
{\bf 1)} 4.4.13.2 implies: any $U$-ordinary point of $\Mn_{1k(v_1)}$ is quasi $U$-ordinary. This generalizes 4.6 P3.
\smallskip
{\bf 2)} If all factors of $G^{\rm ad}$ are of $A_1$ Lie type, then the converse of 1) holds. More precisely, we have:
\medskip
{\bf Fact.} {\it Any quasi $U$-ordinary point of $\Mn_{1k(v_1)}$ is a $T$-ordinary point.}
\medskip
{\bf Proof:} We assume $k=\bar k$. As we are dealing with the $A_1$ Lie types, for any $k$-valued quasi $U$-ordinary point $y$ of $\Mn_{1k(v_1)}$, the Lie subalgebra ${\got p}_{=0}$ of ${\rm Lie}(G^{\rm ad}_{W(k)})$ corresponding to the slope $0$ of its attached Shimura adjoint Lie $\sg_k$-crystal $({\rm Lie}(G^{\rm ad}_{W(k)}),\vph_y)$, is (cf. also 2.2.19.2 and the Criterion of 2.2.22 1)) the Lie algebra of a maximal torus of $G^{\rm ad}_{W(k)}$. So $({\rm Lie}(G^{\rm ad}_{W(k)}),\vph_y)$ has lifts which are of toric type. So the Fact follows from 4.4.13.2 and def. 4.4.13.1 c).   
\smallskip
{\bf 3)} We refer to 4.12.12.6. We denote by 
$$NSU(\Mn_{1k(v_1)})$$ 
the number of strata of the quasi-ultra stratification of $\Mn_{1k(v_1)}$ having $U$-ordinary points. If $G^{\rm ad}$ has all simple factors of $A_1$ Lie type, based on 2) it can be computed as follows (see [Va9] for general computations). We can refer just to the context of a SHS $(f,L_{(p)},v)$, with all factors of $G^{\rm ad}$ of $A_1$ Lie type (cf. 4.12.12.6.2 A). We use the notations of 4.3.1.1, with $k=\FF$. 
\smallskip
We fix $i\in\Mh^{\rm nc}$. Let $d(i)$ be as in 4.5.15.2.1: it is the number of elements $j\in\Mh_i$ such that the image of the composite of the cocharacter $\mu_{W(k)}$ of $G_{W(k)}$ with the natural projection of $G_{W(k)}$ on $G_j$ is non-trivial; let $\Ms_i$ be the subset of $\Mh_i$ formed by such $j$'s. We consider the $i$-th cyclic adjoint filtered Lie $\bar\sg$-crystal $({\rm Lie}(G_{iW(\FF)}),\vph_z,F^0({\rm Lie}(G_{iW(\FF)})),F^1({\rm Lie}(G_{iW(\FF)}))$ attached to a Shimura-canonical lift $z:{\rm Spec}(W(\FF))\to\Mn$ (see 4.2.4).  Based on b) of 4.4.1 2), we can assume $\vph_z$ is $\vph\otimes 1$ (see 4.1.1). We work under 4.1.4.1. 
\smallskip
Based on 2) and 4.12.12.6.2 A, we need to determine the number of elements $\om_i$ of the Weyl group $W_G(i)$ of $\prod_{j\in\Ms_i} G_j$ such that, denoting by $g_{\om_i}\in\prod_{j\in\Ms_i} G_j(W(\FF))$ an arbitrary representative of it normalizing the image $T(i)$ of $T_{W(\FF)}$ in $G_{iW(\FF)}$, 
$${\got L}(\om_i):=({\rm Lie}(G_{iW(\FF)}),g_{\om_i}(\vph\otimes 1),F^0({\rm Lie}(G_{iW(\FF)})),F^1({\rm Lie}(G_{iW(\FF)}))$$ 
is of toric type. We consider an arbitrary $j\in\Ms_i$. Let $H_j$ be the $\GG_a$ subgroup of $G_j$ which is normalized by $T(i)$ and such that ${\rm Lie}(H_j)\subset F^0({\rm Lie}(G_{iW(\FF)}))$. Let $x$ be an arbitrary generator of ${\rm Lie}(H_j)$. $\bigl(g_{\om_i}(\vph\otimes 1)\bigr)^{2d(i)}$ takes $x$ into a multiple of $p^{n(\om_i)}x$ by an invertible element of $W(\FF)$; here 
$$n(\om_i)\in S(-d(i),d(i)).$$
\indent
Based on the circular property expressed by 3.4.2.1 and on the $0$-symmetric property of 2.2.3 1), we have: ${\got L}(\om_i)$ is $T$-ordinary iff $n(\om_i)\neq 0$. The number of elements $\om_i\in W_G(i)$ such that $n(\om_i)\neq 0$ is precisely the number of $2d(i)$-tuples $(x_1,...,x_{2d(i)})$ formed by elements of the set $\tilde I:=\{-1,1\}$ and having the following three properties:
\medskip
{\bf P1.} The product $x_sx_{s+d(i)}$ does not depend on $s\in S(1,d(i))$;
\smallskip
{\bf P2.} The sum of its entries is non-zero;
\smallskip
{\bf P3.} Its first entry is $1$.
\medskip
P1 just expresses that $g_{\om_i}$ has just $d(i)$ non-zero components in the simple factors of $G_{iW(\FF)}$, while P3 pays attention to the choice of $H_j$. So we need to find the number $n_{d(i)}$ of $d(i)-1$-tuples formed by elements of the set $\tilde I$ such that the sum of their entries is different from $-1$. If $d(i)$ is odd, we have:
$$n_{d(i)}=2^{d(i)-1}.$$
If $d(i)$ is even, we have:
$$n_{d(i)}=2^{d(i)-1}-C_{d(i)-1}^{{d(i)}\over 2}.$$ 
\indent
Allowing $i\in\Mh^{\rm nc}$ to vary, we conclude:
$$NSU(\Mn_{k(v)})=\prod_{i\in\Mh^{\rm nc}} n_{d(i)}.\leqno NSU(A_1)$$
\indent
{\bf 4)} Referring to 3), the Newton polygon of ${\got L}(\om_i)$ is uniquely determined by the sum $\sum_{i=1}^{2d(i)} x_i$. So to count the number of strata of the refined canonical Lie stratification of $\Mn_{k(v)}$, we need to count the number of such sums which are non-negative, for $2d(i)$-tuples $(x_1,...,x_{2d(i)})$ which are subject to P1 and have entries in $\tilde I$. As this number is obviously $[{{d(i)+3}\over 2}]$, from 4.12.12.6.4.2, 4.5.6.1 and the proof of 4.9.8 we get that the number of strata of the absolute stratification of $\Mn_{1k(v_1)}$ is 
$$N_1(G^{\rm ad},X^{\rm ad},v^{\rm ad})=\prod_{i\in\Mh^{\rm nc}} \bigl[{{d(i)+3}\over 2}\bigr].\leqno NP(A_1)$$
\indent
{\bf 4.12.12.6.7.} ${\rm piv}${\bf-invariants.} We refer to 4.5.15.2.1. From its proof we get that locally in the Zariski topology of $\Mn_{\FF}/H_0$, the $\FF$-valued pivotal points define a locally closed subscheme. As the Weyl group $W_G$ of 4.1.5 is finite, the number of isomorphism classes of Shimura adjoint Lie $\bar\sg$-crystals attached to $\FF$-valued pivotal points of $\Mn_{\FF}$ is finite (4.12.12.6.2 B already points out many situations in which it is 1). We deduce that the number of $\FF$-valued pivotal points of $\Mn_{\FF}/H_0$ is finite. Based on 4.9.9, this extends to $\Mn_{1\FF}/H_{01}$, with $H_{01}$ an arbitrary compact, open subgroup of $G_1(\AA_f^p)$. We denote by 
$${\rm piv}_{(G_1,X_1,v_1)}(H_{01})$$ 
the number of $\FF$-valued pivotal points of $\Mn_{1\FF}/H_{01}$ and refer to it as the ${\rm piv}$-invariant of $H_{01}$ w.r.t. the triple $(G_1,X_1,v_1)$. Its independence on $H_1$ is checked as in [Va2, 6.4.6 3)]. 
\medskip
{\bf 4.12.12.6.7.1. The case of curves.} We assume $\dim_{\CC}(X)=1$. So ${\rm Sh}(G,X)$ is a Shimura curve; in particular, $G^{\rm ad}$ is a simple $\QQ$--group of $A_1$ Lie type. From 4.12.12.6.2 A we get that the quasi-ultra stratification of $\Mn_{k(v)}$ has precisely $2$ strata: the Shimura-ordinary locus and its complement, to be referred as the pivotal locus. Based on 4.12.12.6.2 B), this matches 4.12.12.6.5; we get that all points of $\Mn_{k(v)}$ with values in fields are toric. We also get that ${\rm piv}_{(G_1,X_1,v_1)}(H_{01})$ counts the number of $\FF$-valued of $\Mn_{1\FF}/H_{01}$ which are not Shimura-ordinary. See [Ih, (6) of p. 17] for the computation of ${\rm piv}_{(G_1,X_1,v_1)}(H_{01})$ in many particular cases.  
\medskip
{\bf 4.12.12.7. A second approach to (integral) Manin problems.} Besides the above approach (see 4.12.4 and 4.12.12) to solving integral Manin problems and based on:
\medskip
{\bf i)} the isogeny property stated in [Va2, 1.7] (see also 1.15.7), on 
\smallskip
{\bf ii)} Theorem 14 of 1.12.1, on
\smallskip
{\bf iii)} d) of 4.4.1 3), and on
\smallskip
{\bf iv)} the passage from the Hodge type to the preabelian type (see 4.9.4, 4.9.7-9 and 4.12.5-6),
\medskip\noindent
there is a second approach. As there is nothing else to be done to replace iii) or iv), we need to refer just to a second approach to the completion property of 3.6.15 A. It relies on: 
\medskip
{\bf v)} the possibility of determining which strata of the ultimate adjoint stratification of $\Mn_{k(v)}$ are closed; 
\smallskip
{\bf vi)} the assumed existence of pivotal points (based on 4.12.12.6 and on 3.13.7.3, in \S 10 we will include computations pointing out that --at least in majority of cases-- 4.12.12.6.2 B holds in general);
\smallskip
{\bf vii)} the Expectation of 3.13.7.1. 
\medskip
The use of pivotal points is an idea first developed (in connection to Newton polygons) for Siegel modular varieties in [Oo2]. Part of the philosophy behind v) to vii) can be formulated as the following question:
\medskip
{\bf Q} {\it When the category $\GG\SS\CC(\Mn_{k(v)})$ has a final object?}
\medskip
The density part of 4.2.1 and b) of 4.4.1 2) can be interpreted as: 
\medskip
{\bf Corollary.} {\it $\GG\SS\CC(\Mn_{k(v)})$ has an initial object.}
\medskip
{\bf 4.12.12.8. Dieudonn\'e's invariants.} Let $(G_1,X_1,H_1,v_1)$ and $H_{01}$ have the same significance as in 4.12.12.6 and 4.12.12.6.7. Let $q_1\in\NN$ be such that $k(v_1)=\FF_{p^{q_1}}$. Let $q\in\NN$. The number 
$$D_q(G_1,X_1,H_1,v_1,H_{01})$$ 
of strata of the ultimate adjoint stratification of $\Mn_{1k(v_1)}$ containing $\FF$-valued points lifting $\FF_{p^{qq_1}}$-valued points of $\Mn_{1k(v_1)}/H_{01}$, is referred as the level-$q$ Dieudonn\'e invariant of $H_{01}$ w.r.t. the triple $(G_1,X_1,v_1)$.    
\medskip
{\bf 4.12.13. Variants.} We take $p\ge 2$.
\medskip
{\bf 1)} Theorem 2 of 3.15.1 (or 3.6.7.1 and 3.11.1 for $p\ge 3$) can be interpreted as the solution of an integral Manin problem: Any Shimura filtered $\sg_k$-crystal $(M,F^1,\vph,\tilde G)$ can be connected (through the $p$-adic completion of an $\NN$-pro-\'etale scheme over a smooth, affine $W(k)$-scheme, having a special fibre which is a connected, $AG$ $k$-scheme) with a $\tilde G$-canonical lift $(M,F^1,\vph_1,\tilde G)$. Moreover, Theorem 2 of 3.15.1 and 3.6.18.8.3 point out that often there are variants of these over any perfect field of characteristic $p$.
\smallskip
We have as well variants for quasi (even pseudo) Shimura filtered $F$-crystals over perfect fields (cf. 3.6.1.6). 
\smallskip
{\bf 2)} The proofs of 4.12.12 and 4.12.12.2 show that (cf. also 2.3.3.1), in 4.12.12 and 4.12.12.2, if $D$ and its dual are not \'etale, then instead of an affine $W(k)$-scheme $\Mm_D$ we can work with a connected, projective, smooth scheme $\Mp_D$ over $W(k)$, such that the $p$-divisible group over it is uni plus versal in each $k$-valued point of $\Mp_D$. In other words, in the case when $D$ is the $p$-divisible group of an elliptic curve, we can still use the first paragraph of 4.12.12.0 in the context of a ($p=2$) SHS $(f,L_{(p)},v)$, with $G^{\rm ad}$ an absolutely simple $\QQ$--group of $A_1$ Lie type whose $\QQ$--rank is $0$ and which splits over $\QQ_p$. 
\smallskip
This usage goes as follows. It is 2.3.5.1 (and 2.3.18 A for $p=2$) which guarantees the existence of such a ($p=2$) SHS. Moreover, as its construction (see [Va2, 6.5.1.1 and 6.6.5]) is based on [De2, 2.3.10], from [De2, 2.3.13] we get that we can assume $G^{\rm ab}=\GG_m$. So $G_{\ZZ_p}$ is the $GL$-group of a free $\ZZ_p$-module of rank $2$. We consider a non-trivial subrepresentation of it of its representation on $L_{(p)}^*\otimes_{\ZZ_{(p)}} \ZZ_p$ which is defined by a direct summand $DS$ of rank $2$. It is 4.3.6 (or 4.12.12.6.1) which guarantees that there are $k$-valued points of $\Mn$, such that the $p$-divisible groups associated to them via Fontaine's comparison theory and this $DS$ (see 4.12.12.0), are isomorphic to the $p$-divisible group of an arbitrary a priori given elliptic curve over $k$.  
\smallskip
Moreover, choosing $H_0$ small enough, we can assume $\Mp_D$ is of general type (for instance, see [Mi4, 1.2 of \S 2]). 
\smallskip
{\bf 3)} We have a variant of Theorem 13 of 1.12 in the principally quasi-polarized context: its proof is entirely the same (though slightly easier due to the fact that we do not have to split the discussion --see the proof of 4.12.12-- in terms of some reflex field being $\QQ$ or not; we always can assume it is $\QQ$).
\smallskip
{\bf 4)} There are simple proofs of Serre--Tate's deformation theory: for instance, see Drinfeld's proof in [Ka3, ch. 1]. Also in [MFK, ch. 6.3] a simple theory of deformations of abelian varieties is presented. So, based on this last two loc. cit. and 4.12.12.5 we reobtain the local deformation theory of $p$-divisible groups over a perfect field $k_1$, i.e. we reobtain [Il, 4.8] (standard Galois descent allows us to descend from $\overline{k_1}$ to $k_1$). Warning: this approach avoids entirely the use of [Me], [BBM], [BM] or of [Il]. 
\medskip
Let $m\in\NN\cup\{0\}$. Let $R:=W(k)[[x_1,...,x_m]]$ and let $\Phi_R$ be a Frobenius lift of it.  
\medskip
{\bf 4.12.14. The general local integral Manin problem for Siegel modular varieties.} We consider a SHS of the form $(1_{({\rm GSp}(W,\psi),S)},L_{(p)},p)$. So, with the standard notations, $\Mm=\Mn$. The problem is to determine all principally quasi-polarized $p$-divisible objects of $\Mm\Mf_{[0,1]}(R)$ associated to principally quasi-polarized $p$-divisible groups over $R$, obtained from the one of $(\Ma,\Mp_{\Ma})$ by pull back through morphisms ${\rm Spec}(R)\to\Mn$. 
\medskip
{\bf 4.12.15. The general local integral Manin problem for a SHS $(f,L_{(p)},v)$.}
Let $(f,L_{(p)},v)$ be a SHS. The problem is: which principally quasi-polarized $p$-divisible objects with tensors of $\Mm\Mf_{[0,1]}(R)$ are associated to the pull back of the principally quasi-polarized $p$-divisible group of $(\Ma,\Mp_{\Ma})$ and of (de Rham components of) Hodge cycles $(w_{\al}^{\Ma})_{\al\in\Mj^\prime}$ of $\Ma$ through morphisms ${\rm Spec}(R)\to\Mn$?
\medskip
{\bf 4.12.16. The general local integral Manin problem for a Shimura quadruple $(G,X,H,v)$ of preabelian type, with $(v,2)=1$.} The problem is to determine which Lie $p$-divisible objects of $\Mm\Mf_{[-1,1]}(R)$ are isomorphic to Shimura adjoint filtered Lie $F$-crystals over $R/pR$  attached to (cf. 4.9.20) morphisms ${\rm Spec}(R)\to\Mn$. 
\medskip
{\bf 4.12.17. The solution of 4.12.14.} The answer is: any principally quasi-polarized $p$-divisible object of $\Mm\Mf_{[0,1]}(R)$. Argument: 3.6.18.7.1 a) (applied in the context of a principally quasi-polarized pseudo Shimura $\sg_k$-crystal $(M_0,g\vph,Sp(M_0,\psi),\psi)$), allows us to assume we are dealing with a principally quasi-polarized $p$-divisible object of $\Mm\Mf_{[0,1]}^\nabla(R)$. So the answer follows from 4.12.4, the deformation theory of abelian varieties and 2.2.21 UP.
\medskip
{\bf 4.12.18. The solution of 4.12.15 in the case $\Mn$ has the completion property.} We use the notations of 4.12.5. The answer is: any such $p$-divisible object of the form 
$$(M_0\otimes_{W(k)} R,F^1_0\otimes_{W(k)} R,g_R(\vph_0\otimes 1),(v_{\al})_{\al\in\Mj^\prime},\tilde\psi),$$ 
with $g_R$ an arbitrary element of $G^0_{\ZZ_{(p)}}(R)$. As in 4.12.17, this is an easy consequence of 3.6.18.7.1 a) (applied in the context of a principally quasi-polarized pseudo Shimura $\sg_k$-crystal $(M_0,g\vph_0,G^0_{W(k)},(v_{\al})_{\al\in\Mj^\prime},\tilde\psi)$), 4.12.5 and 2.2.21 UP (cf. also Fact 4 of 2.3.11).
\medskip
{\bf 4.12.19. The solution of 4.12.16 in case $\Mn^{\rm ad}$ (of 4.12.6) has the completion property.} We use the notations of 4.12.6. The answer is: any such $p$-divisible object of the form 
$$
({\rm Lie}(G^{\rm ad}_R),g_R^{\rm ad}(\sg\mu^{\rm ad}({1\over p})\otimes 1),F^0({\rm Lie}(G_R^{\rm ad})),F^1({\rm Lie}(G_R^{\rm ad}))), 
$$
with $g_R^{\rm ad}\in G^{\rm ad}_{\ZZ_{(p)}}(R)$. Here $F^0({\rm Lie}(G_R^{\rm ad}))$ and $F^1({\rm Lie}(G_R^{\rm ad}))$ are obtained as usual from (the extension to $R$ of) the cocharacter $\mu^{\rm ad}$ of 4.12.6. 
\smallskip
To argue this, we consider a SHS $(f_1,L_{1(p)},v_1)$ as in 4.12.6. Using the same argument as the one of 4.12.6.1 handling the general form of its Fact, we can assume $g_R^{\rm ad}$ is of the form $g_1g_R^{2\rm ad}$, with $g_R^{2\rm ad}$ liftable to an $R$-valued element $g_R$ of $G_{1R}^{\rm der}$ congruent to the identity modulo the ideal $(x_1,...,x_m)$ of $R$ and with $g_1\in G^{\rm ad}_{W(k)}(W(k))$. Using this we deduce the existence of a principally quasi-polarized $p$-divisible object ${\got C}=(M_1\otimes_{W(k)} R,F^1_1\otimes_{W(k)} R,g_R(\vph_1\otimes 1),(v_{1\al})_{\al\in\Mj^\prime},\tilde\psi)$ with tensors of $\Mm\Mf_{[0,1]}(R)$ as in 4.12.18 (the index $0$ being replaced by 1), and whose attached Shimura adjoint Lie $F$-crystal is isomorphic to the one of the above answer. We can assume the Shimura filtered $\sg_k$-crystal $(M_1,F_1^1,\vph_1,G_{1W(k)},(v_{1\al})_{\al\in\Mj^\prime})$ is attached to a $W(k)$-valued point $z_1$ of $\Mn_1$, cf. 4.12.6. As in 4.12.18, 3.6.18.7.1 a) allows us to assume the existence of a $G_{1W(k)}$-invariant connection on $M_1\otimes_{W(k)} R$ making ${\got C}$ to be viewed as a principally quasi-polarized $p$-divisible object with tensors of $\Mm\Mf_{[0,1]}^\nabla(R)$. So the answer follows by passage to associated Shimura adjoint $F$-crystals, once we remark (cf. 2.2.21 UP and the Fact 4 of 2.3.11 applied to $(f_1,L_{1(p)},v_1)$) that ${\got C}$ is obtained via pull back from a uni plus versal Shimura filtered $F$-crystal of the form $(M_1,F_1^1,\vph_1,G_{1W(k)},N_1,\tilde f,(v_{1\al})_{\al\in\Mj^\prime})$.
\medskip
4.12.16-18 admit a non-filtered version (cf. 2.2.1.6). Not to be too long this is not going to be presented here.   
\medskip
{\bf 4.12.20. Exercise.} State and prove the analogues of 4.12.14-15 and 4.12.17-18 for $p=2$. Hint: just use 2.3.18, 3.14 B and 4.12.12.1.
\medskip
{\bf 4.12.21. General (integral) Manin problems.} Given an arbitrary $\ZZ_{(p)}$-scheme $Y$ we would like to understand which principally quasi-polarized $p$-divisible groups (resp. which isogeny classes of principally quasi-polarized $p$-divisible groups) over $Y$ are coming from principally polarized abelian schemes (resp. from isogeny classes of principally polarized abelian schemes) over $Y$ of relative dimension $e$. This forms the general integral Manin problem (resp. the general Manin problem) for the Siegel modular variety of 4.12.4. Similarly we can define the general (integral) Manin problem for an arbitrary SHS $(f,L_{(p)},v)$: we need to replace $p$-divisible groups by Shimura $p$-divisible groups (see 2.2.20; in \S5-10 we will gradually extend the referred definition: we will define Shimura $p$-divisible groups over any $O_{(v)}$-scheme). We can also define the general (integral) Manin problem in the context of 4.12.3: again 4.9.20 can be extended to any $O_{(v)}$-scheme (see \S 5 for details).  
\smallskip
A general (integral) Manin problem looks very hard in such a generality. In the cases when we are dealing with (the $p$-adic completion of) a (regular, formally) smooth, separated scheme $Y$ (over DVR's of index of ramification 1 over $\ZZ_{(p)}$ or $O_{(v)}$) we hope that some partial solutions can be obtained based on 3.6.
\medskip\smallskip
{\bf 4.13. Local forms of the invariance principle.} There are two main local forms of this principle: one in the adjoint context and one in the context of Shimura $p$-divisible groups; their proofs are very much the same and so they are combined. Let $p\ge 3$ be a rational prime. 
\smallskip
We start with the first form. Let $(G_i,X_i,H_i,v_i)$, $i=\overline{1,2}$, be two Shimura quadruples of preabelian type, with $v_i$ dividing $p$. Let $\Mn_i$ be the integral canonical model of $(G_i,X_i,H_i,v_i)$. Let $H_{0i}\subset G_i(\AA_f^p)$ be a subgroup such that either $p| t(G_i^{\rm ad})$ and $H_{0i}\times H_i$ is a subgroup of $G_i(\AA_f)$ $p$-smooth for $(G_i,X_i)$ or $p$ does not divide $t(G_i^{\rm ad})$ and $H_{0i}\times H_i$ is smooth for $(G_i,X_i)$ (cf. AE.4.2 a)). Let $y_i:{\rm Spec}(k)\to\Mn_{ik}/H_{0i}$, $i=\overline{1,2}$, be such that the non-trivial part of the Shimura adjoint Lie $\sg_k$-crystal attached to $y_1$ is isomorphic to the non-trivial part of the Shimura adjoint Lie $\sg_k$-crystal attached to $y_2$. We fix such an isomorphism $\Mi$. Let $O_{y_i}$ be the local ring of $y_i$ (in $\Mn_{iW(k)}/H_{0i}$) and let $O_{y_i}^{{\rm h}\wedge}$ be the $p$-adic completion of its henselization. We assume $k=\bar k$. Let ${\got C}_{y_i}^{\rm ad}$ be the Shimura adjoint filtered Lie $F$-crystal  over $O(p)_{y_i}^{\rm h}:=O_{y_i}^{{\rm h}\wedge}/pO_{y_i}^{{\rm h}\wedge}$ attached (see 4.9.20) to the canonical morphism ${\rm Spec}(O_{y_i}^{{\rm h}\wedge})\to\Mn_{iW(k)}/H_{0i}$. Let 
$${\got C}_{y_i}^{\rm nt}$$ 
be its non-trivial part (it is defined as in 3.10.1; it is a direct factor of ${\got C}_{y_i}^{\rm ad}$). The refined canonical Lie stratification of $\Mn_{ik(v_i)}$ makes ${\rm Spec}(O(p)_{y_i}^{\rm h})$ to be a stratified scheme.
\smallskip
Let $n\in\NN$. We have:
\medskip
{\bf 4.13.1. Theorem.} {\it There is an isomorphism 
$$\rho:{\rm Spec}(O_{y_1}^{{\rm h}\wedge})\tilde\to {\rm Spec}(O_{y_2}^{{\rm h}\wedge})$$
inducing an isomorphism 
$$\rho(p):{\rm Spec}(O(p)_{y_1}^{\rm h})\tilde\to {\rm Spec}(O(p)_{y_2}^{\rm h})$$
of stratified schemes and such that $\rho^*({\got C}_{y_2}^{\rm nt}/p^n{\got C}_{y_2}^{\rm nt})$ is isomorphic to ${\got C}_{y_1}^{\rm nt}/p^n{\got C}_{y_1}^{\rm nt}$, through an isomorphism respecting the Lie structures and lifting $\Mi$ mod $p^n$.}
\medskip
We present now the second form. Let $(f^i,L_{(p)}^i,v^i)$, $i=\overline{1,2}$, be two standard Hodge situations. Using the standard notations of 2.3.1-3 for the SHS $(f^i,L_{(p)}^i,v^i)$ but putting everywhere an upper right index $i$, let $y_i:{\rm Spec}(k)\to\Mn^i_k/H_0^i$ be a morphism. Here $k$ is an arbitrary perfect field containing $k(v^1)$ and $k(v^2)$. We assume $H_0^i$ is small enough so that all Hodge cycles of $\Ma^i$ are as well Hodge cycles of $\Ma_{H_0^i}^i$. We have:
\medskip
{\bf 4.13.2. Theorem.} {\it We assume the principally quasi-polarized  Shimura $\sg_k$-crystal attached to $y_1$ is isomorphic to the principally quasi-polarized Shimura $\sg_k$-crystal attached to $y_2$. We fix such an isomorphism $\Mi$. Then there is an isomorphism
$$\rho:{\rm Spec}(O_{y_1}^{{\rm h}\wedge})\tilde\to {\rm Spec}(O_{y_2}^{{\rm h}\wedge})$$
inducing an isomorphism
$$\rho(p):{\rm Spec}(O(p)_{y_1}^{\rm h})\tilde\to {\rm Spec}(O(p)_{y_2}^{\rm h})$$ 
of stratified schemes and such that the pull back through $\rho$, of the crystalline counterpart of the kernel of the multiplication by $p^n$ of the principally quasi-polarized Shimura $p$-divisible group $\Md_{y_2}$ over ${\rm Spec}(O_{y_2}^{\rm h})$ obtained from the principally quasi-polarized Shimura $p$-divisible group over $\Mn^2/H_0^2$ (see Fact 3 of 2.3.11)  by pull back through the canonical morphism ${\rm Spec}(O_{y_2}^{\rm h})\to\Mn^2/H_0^2$, is isomorphic to the crystalline counterpart of the kernel of the multiplication by $p^n$ of the principally quasi-polarized Shimura $p$-divisible group $\Md_{y_1}$ over ${\rm Spec}(O_{y_1}^{\rm h})$, defined similarly as $\Md_{y_2}$ (here $O_{y_i}^{\rm h}$, $O(p)_{y_i}^{\rm h}$ have the same meaning as in 4.13.1), through an isomorphism lifting $\Mi$ mod $p^n$.}
\medskip 
{\bf Proofs:} In connection to 4.13.1 we refer to definitions 4.9.17 and 4.9.20 (in the case of a SHS, 4.9.20 is trivial, just introducing some terminology). For the parts involving stratifications, we just need to take $n$ big enough and to apply 3.15.7 BP1. So 4.13.1-2 are a direct consequence of (the second variant of) 3.6.14.4 (for 4.13.1 the passage from Lie algebras of reductive groups to Lie algebras of their adjoints, in the context of crystalline counterparts, is entirely trivial). 
\medskip
{\bf 4.13.3. Remarks.} {\bf 1)} The form of 4.13.2, with $n=\infty$, where instead of taking henselization we take the completion, was obtained previously: the Serre--Tate's deformation theory implies this form for a standard PEL situation; the general case can be obtained by just copying [Va2, 5.4].
\smallskip
{\bf 2)} The condition $k=\bar k$ in 4.13.1 can be weaken: see \S 5 (cf. 4.9.20).
\smallskip
{\bf 3)} We have variants of 4.13.1-2, where instead of refined canonical Lie stratifications, we use its variants listed in 4.9.9.
\medskip\smallskip
{\bf 4.14. Final remarks.} Here we include some remarks which are sort of overview of great parts of [Va2] and of this paper.
\medskip 
{\bf 4.14.1. The generic approach.} Many results presented here or in [Va2] are trivial or very easy, from the generic point of view, i.e. working with a prime $p$ big enough. But the main inconvenience of this simple approach is that we do not have at all effective upper bounds (of how big $p$ has to be); so, from the point of view of concreteness, it is void. We now list some results which generically ``are very easy" and motivate why this is so.
\medskip
{\bf a)} [Va2, 3.4.1 and 3.4.7] is all we need to get generically the existence of integral canonical models of Shimura varieties of Hodge type.
\smallskip
{\bf b)} [Va2, 6.2.2 a), 6.2.3 and 6.8] is all we need to extend this to Shimura varieties of preabelian type.
\smallskip
{\bf c)} [Va2, 6.4.4 and 6.4.11] generically are obvious.  
\smallskip
{\bf d)} The density part of 4.2.1 is obvious generically (due to the existence of smooth toroidal compactifications of canonical models of Shimura varieties and of c)).
\smallskip
{\bf e)} Using [Va2, 3.2.11] and d) we get (via 3.1.0 a) and c)) the generic form of 4.2.1 itself.
\smallskip
{\bf f)} Using [Va2, 3.2.11] and d), 4.6 P1 is obvious generically.
\medskip
For e) and f) cf. also the proof of 4.12.12 or cf. 4.12.12.6.
\medskip
{\bf 4.14.2. The special type case.} There are plenty of generalized Shimura (filtered) $F$-crystals and of Shimura (filtered) Lie $F$-crystals involving the  $E_6$, $E_7$ and $D_n^{\rm mixed}$ types (see 3.6.1.6 and 3.10.5 for terminology). So, even before the proof of the existence of integral canonical models of Shimura varieties of special type and of the validity of the expectation of [Va2, a) and b) of 3.2.7 8)] (in loc. cit., strictly speaking just a variant of [Va2, 5.6.5 f)] is satisfied, cf. 2.2.8 3) and the beginning paragraph of 3.4; see also AE.0), we can speak about stratifications (as in 4.5) of schemes of moduli of Shimura (adjoint) filtered Lie $F$-crystals (often attached to generalized Shimura filtered $F$-crystals) of special type. For the study of these stratifications we refer to [Va6]. Here we just mention two things:
\medskip
{\bf 1)} All formulas of 4.5.15.2 hold for the generalized Shimura context, whenever we have global deformations which are uni plus versal; based on 3.15.6 D (such deformations do exist and moreover) the arguments of 4.5.15.2.5 need no modification at all. 
\smallskip
{\bf 2)} All of 4.7.11 and of 4.7.12-17 (as well as their $p=2$ version; see 4.14.3 E below) can be entirely adapted to the generalized Shimura context (again, based on 3.15.6 D and E no modification of arguments are needed; 3.15.6 B shows that even the part of 4.7.11 6) referring to [Ka4, A2.2 (2)] can be entirely adapted). In particular, the third paragraph of 4.7.11 2A) can be performed without any reference (via some lifting process) to $p$-divisible groups.
\medskip
{\bf 4.14.3. The case of a $p=2$ SHS.} From 2.3.18 and 3.14 we get that most of the results of 4.1-13 remain true for a $p=2$ SHS $(f,L_{(2)},v)$ (see also [Va5]). However some of them need some reformulation: below we refer to such a $p=2$ SHS. Warning: below we do not always state explicitly at each step the implicit references to 3.14.
\medskip
{\bf A.} All results of 4.1-3 remain true: no extra reformulation is needed (cf. 3.14 C). Warning: 4.2.8.1* remains true (the argument of [Va2, 6.4.1.1 2)] applies for $p=2$ as well) in the compact type case but we can not say presently anything very precise about non-compact type cases. Also, the quadratic form in two variables of the proof of Fact 2 of 4.3.6, has to be replaced by $x_1^2+x_1x_2+{\al}x_2^2$, where $\al\in\tilde k$ is such that the equation $x^2+x+\al$ has no roots in $\tilde k$.
\medskip
{\bf B.} 4.4.1 1) and 3) as well as b), c) and d) of 4.4.1 2) remain true without any modification. We do expect that d) of 4.4.1 3) and 4.4.12 also hold. But b) of 4.4.1 2) and its version for $U$-ordinary points (see 4.4.13.1 d)) as well as 4.4.8 1) and 3) need substantial reformulation.
\smallskip
As $G$-canonical lifts for $p\ge 3$ were defined at the level of filtrations in the crystalline cohomology context (cf. a) of 4.4.1 2)), for $p=2$, for their (desired) uniqueness (and the parts of 4.4-11 involving it, like 4.4.5-6) we have to proceed very carefully, cf. 2.3.18.1. 
\medskip
{\bf B1.} As a first warning, using the existence of (to be mentioned in E below) $G$-crystalline (additive or multiplicative) coordinates we get:
\medskip
{\bf Warning.} {\it For any $k$-valued $G$-ordinary point $y$ of $\Mn_{k(v)}$, denoting its attached Shimura $\sg_k$-crystal by $(M_y,\vph_y,G_{W(k)})$, a lift $z:{\rm Spec}(W(k))\to\Mm_{W(k(v))}$ of $y$ (viewed as a $k$-valued point of $\Mm_{k(v)}$) such that $(M_y,F^1_0,\vph_y,G_{W(k)})$, with $(M_y,F^1_0,\vph_y,p_{M_y})$ as the principally quasi-polarized filtered $\sg_k$-crystal attached to $z$, is a Shimura filtered $\sg_k$-crystal, does not necessarily factor through $\Mn_{W(k(v))}$.}
\medskip
Easy examples can be obtained with $G$ a torus or starting from 4.14.3.1 below.
\smallskip
So even if $k(v)=\FF_2$ (i.e. even when we can use 4.6 P1 cf. D below), we can not get (at least theoretically) the uniqueness part by specifying that we choose the usual canonical lift. 
\smallskip
We explain the situation, in the slightly more general case when $k(v^{\rm ad})=\FF_2$. We use the notations and terminology of 4.7.0, 4.7.8 and 4.7.11 1). The formal moduli scheme of $G$-deformations of $(A_y,p_{A_y})$ can be naturally identified with a closed formal subscheme of the formal torus (see E below) of deformations of the $2$-divisible group $D_y^1$ (of 4.7.11 1)). This formal subscheme is the translation of a formal subtorus by a $2$-torsion point $P_{2-{\rm tor}}$: this is nothing else but the $p=2$ analogue of 4.7.17 (cf. also E below). This formal subtorus is uniquely determined but not $P_{2-{\rm tor}}$, provided it is non-zero. So if $P_{2-{\rm tor}}$ is the identity element, then we can define the $G$-canonical lift of $y$ uniquely but not otherwise. We do not know when such a $2$-torsion point is not the identity; we just point out that this problem is very much related to variants of Milne's conjecture of 1.15.1. In \S5 and \S6 we will prove, using such variants of Milne's conjecture, that in many situations such a torsion point is the identity. 
\smallskip
So, in general, we have to allow either a finite number of such $G$-canonical lifts or to try to use something else. The most logical ``something else" would be to use 3.14 I; however, as pointed out in the above paragraph, without quite a lot of work, we can not check that we can apply 3.14 I in our geometric contexts of $p=2$ SHS's. 
\smallskip 
As a second warning we have: some $G$-canonical lifts might be defined only up to passage from $k$ to an abelian extension whose Galois group is a subgroup of $(\ZZ/2\ZZ)^m$, with $m\in\NN$, $m\le\dim_{\CC}(S)$, cf. the proof of 2.3.18.1 B and its logical variant in the principally quasi-polarized context. It can checked that $m\le\dim_{\CC}(X)$ (the case $k(v^{\rm ad})=\FF_p$ can be deduced from 2.3.18.1.1 and the above part on $2$-torsion points; the general case can be deduced by adapting the proof of B2 below: see also the Exercise 4.14.3.1 below).
\smallskip
Presently we do not know what in general would be the most practical ``something else".  For instance, we could try to choose one such $G$-canonical lift to which a ``bigger" (or specified) part of the $\ZZ$-Lie algebra of $G$-endomorphisms (or of the $\ZZ$-algebra of endomorphisms or of the group of $G$-automorphisms or of automorphisms) of $A_y$, lifts, or (this seems the easiest way) in terms of $2$-divisible groups (i.e. of Galois representations with $\ZZ_2$-coefficients); in case $k$ is not $1$-simply connected, we could try to specify that the $G$-canonical lift has to be defined over $W(k)$, etc.  
\smallskip
However we do have:
\medskip
{\bf B2. Proposition.}  {\it If the Lie $p$-rank $p-{\rm Lie}_G(y)$ (see 4.3.8 1)) is $0$ (i.e. if each point of $\Mn_{k(v)}$ with values in a field is non-integral in the sense of 4.12.12.6.5), then we can define a unique $G$-canonical lift of $y$.}
\medskip
{\bf Proof:} We refer again to the direct sum decomposition $M_y=M_y^1\oplus M_y^2$ of 4.7.11 1); so in what follows $N$ and $N_0$ are defined as in 4.7.1. We have a direct sum decomposition 
$${\rm End}(M_y)={\rm End}(M_y^1)\oplus {\rm Hom}(M_y^1,M_y^2)\oplus {\rm Hom}(M_y^2,M_y^1)\oplus {\rm End}(M_y^2).$$ 
From 3.14 I we get that the $F^1$-filtration $F^1$ of $M_y$ defined by any $G$-canonical lift $z$ of $y$ can be similarly decomposed as $F^1=F^1_1\oplus F^1_2$, with $F^1_j\subset M_y^j$, $j=\overline{1,2}$. So we have a direct sum decomposition of $F^0({\rm End}(M))$ in four direct summands similar to the one of ${\rm End}(M)$. Accordingly, the moduli scheme ${\rm Spec}(R)$ of deformations of $D_y^1\oplus D_y^2$ is ``decomposed into four pieces". Using the language of homomorphisms of $W(k)$-algebras, by this we mean: there are natural $W(k)$-epimorphisms $q_i:R\twoheadrightarrow R_i$, $i=\overline{1,4}$, with $R$, $R_1$,..., $R_4$ as $W(k)$-algebras of formal power series, such that the tangent space of ${\rm Spec}(R)$ is a direct sum of the images (via these $W(k)$-epimorphisms) of the tangent spaces of ${\rm Spec}(R_i)$, $i=\overline{1,4}$; $q_1$ is defined by the moduli scheme of deformations of $D_y^1$, $q_4$ is defined by the moduli scheme of deformations of $D_y^2$, etc. Here, as well as in what follows, the tangent spaces are defined w.r.t. the $W(k)$-epimorphism $R\twoheadrightarrow W(k)$ defining the trivial deformation of $D_z^1\oplus D_z^2$ and factoring through $q_i$, $i=\overline{1,4}$; $D_z^i$ is a fixed lift of $D_y^i$, $i=\overline{1,2}$, such that the $p$-divisible group of $z^*(\Ma)$ is naturally decomposed as a direct sum $D_z^1\oplus D_z^2$ (cf. also 2.3.18.1.1).
\smallskip
We can view these four pieces also as projections. To exemplify it, we concentrate on the first factor corresponding to $q_1$. We first consider a complete, local $W(k)$-subalgebra $R^4$ of $R$ generated by $W(k)$ and by a maximal set of regular parameters of $R$ belonging to the kernel of $q_4$, then we consider a complete, local $W(k)$-subalgebra $R^{34}$ of $R^4$ generated by $W(k)$ and by a maximal set of regular parameters of $R^4$ belonging to the kernel of the natural $W(k)$-epimorphism $R^4\twoheadrightarrow R_3$ defined by $q_3$, and then we consider a complete, local $W(k)$-subalgebra $R^{234}$ of $R^{34}$ generated by $W(k)$ and by a maximal set of regular parameters of $R^{34}$ belonging to the kernel of the natural $W(k)$-epimorphism $R^{34}\twoheadrightarrow R_2$ defined by $q_2$. The composite $q^1:R^{234}\to R_1$ of the $W(k)$-monomorphism $R^{234}\hookrightarrow R$ with $q_1$, is a $W(k)$-isomorphism.
\smallskip
Let $O_y$ be the local ring of $y$, viewed as a $k$-valued point of $\Mn_{W(k)}$; we view it as a $W(k)$-algebra. The hypothesis $p-{\rm Lie}_G(y)=0$ gets translated in: 
${\rm Lie}(N)$ has a trivial projection on ${\rm End}(M_y^1)$; so also ${\rm Lie}(N_0)$ has a trivial projection on ${\rm End}(M_y^1)$. We first consider the particular case when ${\rm Lie}(N)\subset {\rm End}(M_y^2)$. In this case, as ${\rm Spec}(\widehat{O_y})$ is integral, we are just deforming $D_z^2$, while fixing $D_z^1$ and so 2.3.18.1 C applies directly. 
\smallskip
We now come back to the general case. We assume the existence of another $W(k)$-valued point $z_1$ (different from $z$) of $\Mn_{W(k)}$ lifting $y$ and whose attached Shimura filtered $\sg_k$-crystal is $(M_y,F^1,\vph_y,G_{W(k)})$. From 2.3.18.1.1 (cf. also E below) we get that it corresponds to (i.e. it is defined by) a $2$-torsion point of the formal torus of deformations of $D_y^1$. So, from the point of view of the moduli scheme of deformations of $D_y^1\oplus D_y^2$, it is defined by a $W(k)$-valued point of ${\rm Spec}(R)$ factoring through ${\rm Spec}(R_1)$. So, as the $2$-torsion points are ``recorded" by the tangent space of ${\rm Spec}(R_1)$, the tangent map of the $W(k)$-morphism obtained by composing the natural $W(k)$-morphism $m_y:{\rm Spec}(\widehat{O_y})\to {\rm Spec}(R)$ with the natural $W(k)$-morphism ${\rm Spec}(R)\to {\rm Spec}(R^{234})$, is non-trivial. But this contradicts the mentioned trivial projection (the image of the Kodaira--Spencer map naturally associated to $m_y$ can be identified with a direct summand of ${\rm Lie}(N_0)$, cf. 4.7.2 and its logical $p=2$ version). This proves the Proposition.
\medskip
{\bf B3.} B2 applies entirely to $U$-canonical lifts. So in 4.9.17.5 we do not get the uniqueness of $U$-canonical lifts of $U$-ordinary points (in G below we do not repeat this). So 4.4.8 1) and 3) have to be modified as well. However, the above part on the uniqueness of $G$-canonical lifts applies as well to non-integral, $U$-ordinary points; so Proposition of B2 points out the following thing: for $p=2$ it is easier to understand $U$-canonical lifts of $U$-ordinary points which in some sense are far from being ordinary. 
\smallskip
Moreover, we have:
\medskip
{\bf Corollary.} {\it We consider a quasi-pivotal point $y\in\Mn_{k(v)}(k)$. Then any $W(k)$-valued point $z$ of $\Mn$ lifting $y$ is uniquely determined by the Hodge filtration $F^1_z$ of $M_y$ defined by $z^*(\Ma)$.}
\medskip
{\bf Proof:} The proof of B2 takes care of the case when the Shimura filtered $\sg_k$-crystal $(M_y,F^1_z,\vph_y,G_{W(k)})$ attached to $z$ is cyclic diagonalizable. It is the Claim of 2.4.1 which allows us to pass the uniqueness to any mentioned $z$. This ends the proof.
\medskip
{\bf B4.} Also 4.4.4 has to be modified, in case a $G$-canonical lift is not uniquely determined by its filtration in the crystalline cohomology context. The easiest way:
\medskip
{\bf Corollary.} {\it Referring to the notations of 4.4.4 we have:
$$4{\rm End}_G(A_y)\subset {\rm End}_G(A_z)$$
and
the elements of ${\rm Aut}_G(A_y)$ which restricted to $A_y[4]$ are identity, are also elements of ${\rm Aut}_G(A_z)$.}
\medskip
{\bf Proof:} This is a consequence of Serre--Tate's deformation theory.
\medskip
Accordingly, 4.4.5-6 remain true. The rest, i.e. 4.4.7 and 4.4.9-12, does not need modifications.
\medskip
{\bf C.}  4.5 involves just Newton polygons, parabolic groups and Shimura (adjoint Lie) $F$-crystals; so it needs no modification at all besides the fact that the split group of $B_n$ (resp. of $D_n$) Lie type of 4.5.6.1 has to be defined using the quadratic form $x_0^2+\sum_{i=1}^n x_{2i-1}x_{2i}$ (resp. $\sum_{i=1}^n x_{2i-1}x_{2i}$), to be compared with 2.3.18 B2. For the fact that all of 4.5.15.2 remains true for $p=2$, cf. 2.3.18 and 3.14 C and J. It seems to us that the automorphism invariants of 4.5.15.3 have ``jumps" in the case $p=2$ provided we are dealing with the $B_n$ and the $C_n$ Lie types with $n\ge 2$ (cf. the special cases in [BT, 4.2.3]; see also 3.5.4 and 3.14 C and end of J). Here by ``jumps" we mean: irregular behavior (unpredictable by looking just at odd primes). 
\medskip
{\bf D.} Modulo some modifications, the part of 4.6 till 4.6.3 inclusive, except the parts of 4.6 P2 and 4.6.3 A referring to canonical lifts, remains true without modifications. In 4.6.1 1) (resp. in connection to 4.6 P8), for $p=2$ we presently need to restrict to the $C_1$ (resp. to the $A_{\ell}$) Lie type. Out of 4.6.2, only 4.6.2.2 and 4.6.2.6 make any sense for $p=2$. Moreover, P2 and 4.6.3 A can be modified as follows: any abelian variety over $W(\bar k)$ obtained via a $W(\bar k)$-valued $G$-canonical lift of $\Mn$, is isogeneous to the (usual) canonical lift of its special fibre (see 2.3.18.1 E). 
\smallskip
Warning: Exercise 4.6.4 A remains true (cf. its first proof) for $p=2$; but (presently) it does not lead to the construction of integral canonical models of (some) Shimura varieties of preabelian type which are not of abelian type (cf. the limitations of [Va2, 6.8.1]). 4.6.4 B makes no sense for $p=2$. 4.6.5 can be checked in many situations (see \S5, \S6 and [Va5]) but for its complete $p=2$ analogue we have to postpone to [Va3]. 4.6.6 holds for $p=2$ standard PEL situations (cf. b) of 2.3.18 B3: the argument pertaining to $B(\FF)$- and $W(\FF)$-valued points of 4.6.6 holds for $p=2$ as well); so d) of 4.4.1 3) holds for all $p=2$ standard PEL situations of 2.3.18 B. 
\smallskip
Lemma of 4.6.7 holds for $p=2$. So 3.14 J implies that 4.6.7-8 applies for $p=2$ as far as 4.14.4.1.1 and 4.14.4.1.2 2) below allow. More precisely, the first part of Corollary 1 (resp. of Corollary 3) of 4.6.7 does hold for $p=2$, as [Va2, 6.5-6] can be performed abstractly, with no reference to $p=2$ SHS's (resp. as we have 2.3.18 B and Lemma 2 of 4.6.4; see also 4.14.4.1.1); moreover, Corollary 2 of 4.6.7 holds for $p=2$ for the $A_n$, and the totally non-compact $C_n$ and $D_n^{\HH}$ types (of 3.10.5); for the totally non-compact $C_n$ and $D_n^{\HH}$ types, cf. also 4.14.4.1.2 2).  
\smallskip
In connection to the $B_n$ and $D_n^{\RR}$ types, 4.6.8 for $p=2$ has to be interpreted just in terms of principally quasi-polarized Shimura $F$-crystals, with no reference to a SHS.
\medskip
{\bf E.} Related to 4.7 there are just three things to be modified. First, part of [De3, 1.4.2] and of [Ka4, A2.2] are dealing just with odd primes; so to get 4.7.5 we have to proceed as explained in 4.7.11 (see 2) and 4) of it especially: they handle the case $p=2$ as well). Second, in 4.7.8 the unit element can correspond to any one of the $G$-canonical lifts of $y$. Accordingly, in 4.11.8 (as well as in geometric variants of the end of 4.7.11 8) and of 4.7.17)  we have to speak not about formal subtori but about the translate of a formal subtorus by a $2$-torsion point. Moreover, in 4.7.11 2) for $p=2$, we define $q:R_m\to R_n$ to be a maximal formally smooth quotient of $R_m$ preserving the tensors involved (i.e. we are in a universal context as of 2.2.21). Similarly, $q_{n,r}$ of 4.7.11 9) is uniquely determined up to translation with a $2$-torsion point.    
\smallskip
Third, [Og, 3.14] implies that in the part of 4.7.11 8) referring to Theorem of 4.7.11 4), with $p=2$, the Frobenius lift $\Phi_{R_n}$ of $R_n$ which makes $g_{R_n}$ to be the identity element, is not unique; we recall $n=\dim_{W(k)}(\tilde N)$. We have:
\medskip
{\bf Fact.} {\it The number of such Frobenius lifts is precisely $2^n$, and so it is finite.}
\medskip
{\bf Proof:} Referring to the last paragraph of 4.7.11 6), we choose (as in 3.6.18.2 in the context where all $b_{ijl}$'s are $0$) a connection $\nabla$ on $\tilde M\otimes_{W(k)} R_n$ which is of the form
$$\dl+\sum_{i=1}^n x_ia_idz_i,$$
with $\dl$ the connection that annihilates $\tilde M$ and with $x_i\in R_n$ such that the Kodaira--Spencer map of $\nabla$ is injective modulo the maximal ideal of $R_n$; for the convenience of the computations we assume $x_i={1\over {z_i+1}}$ (cf. (24) of 3.6.18.2 and (LN) of 4.7.11 6)) and that $\tilde\vph(pa_i)=a_i$, $\forall i\in S(1,n)$, cf. end of 4.7.11 6). So $x_idz_i=d({\rm ln}(z_i+1))$. 
Such a connection allows us to treat the situation as if $n=1$: as the product of any two elements of ${\rm Lie}(\tilde N)$, viewed as endomorphisms of $\tilde M$, is $0$, based on [Og, 3.14], we get immediately that the number of Frobenius lifts of $R_n$ which makes $g_{R_n}$ to be the identity element is at least $2^n$. All such lifts are obtained (cf. loc. cit.) by mapping for $i$ running through an arbitrary subset of $S(1,n)$, $z_i+1$ into $-(z_i+1)^2$ (and not into $(z_i+1)^2$).
\smallskip
To check that the number of Frobenius lifts of $R_n$ which makes $g_{R_n}$ to be the identity element is precisely $2^n$ we just need to remark the following two things:
\medskip
{\bf a)} $\Phi_{R_n}$ is uniquely determined by its Teichm\"uller lift ${\rm Spec}(W(k))\hookrightarrow {\rm Spec}(R_n)$, i.e. there is no Frobenius lift of $R_n$ different from $\Phi_{R_n}$ which makes $g_{R_n}$ to be the identity element and whose Teichm\"uller lift ${\rm Spec}(W(k))\hookrightarrow {\rm Spec}(R_n)$ is the same as of $\Phi_{R_n}$; 
\smallskip
{\bf b)} any Teichm\"uller lift as in a) defines a $2$-torsion point of the $n$ dimensional formal torus ${\rm Spf}(R_n)$. 
\medskip
We first argue b). Using a direct sum decomposition $\tilde M=\tilde M_1\oplus\tilde M_2$ as in 4.7.11 1), we can replace $\tilde M$ by $\tilde M_1$. So, using the $p=2$ variant of 4.7.17, we can assume $\tilde G=GL(\tilde M)$ and this last case is well known (see 2.3.18.1 B and 2.3.18.1.1). The argument for a) goes as follows. Due to the injectivity of the Kodaira--Spencer map, any such Frobenius lift $\Phi_{R_n}^1$ has to be of essentially multiplicative type and so (cf. 3.6.18.1.1) of multiplicative type. So, as in the part of 4.7.11 8) referring to ``up to an isomorphism", we can assume the existence of $\tilde z_i\in R_n$, $\forall i\in S(1,n)$, such that:
\medskip
-- $\Phi_{R_n}^1(\tilde z_i+1)=(\tilde z_i+1)^2$, $\forall i\in S(1,n)$;
\smallskip
-- $R_n=W(k)[[\tilde z_1,...,\tilde z_n]]$;
\smallskip
-- $\tilde z_i\in(z_1,...,z_n)$, $\forall i\in S(1,n)$ (here is the place where the part on Teichm\"uller lifts is used);
\smallskip
-- w.r.t. $\tilde z_i$'s, $\nabla$ gets the same logarithmic form as above.
\medskip
From 2.2.21 UP we get that the $W(k)$-isomorphism of $R_n$ that takes $z_i$ into $\tilde z_i$ is the identity isomorphism, and so $\Phi_{R_n}^1$ is not different from $\Phi_{R_n}$. This proves the Fact. 
\medskip
In another order of ideas, the uniqueness part of [Og, 3.15] for $p=2$ is a consequence of 2.3.18.1 C. 
\medskip
{\bf F.} The whole of 4.8 remains true for $p=2$, as Fontaine's comparison theory with rational coefficients is still true in the case of good reduction and mixed characteristic $(0,2)$, see [Fa1-2]. This is the reason we stated 4.8 in terms of $\QQ_p$-valued points (however, see 4.8.3 f) and the involved 2.3.18.1 E).
\medskip 
{\bf G.} Also 4.9 can be adapted to a great extend to the case of a $p=2$ SHS $(f,L_{(2)},v)$. The limitations are as follows. First, in 4.14.3.2 below we include just $p=2$ versions of [Va2, 6.2.2] and not of [Va2, 6.1.1]. Second, we keep quiet (i.e. we postpone a precise statement, till B1 is ``sorted out" fully) regarding the fact that $(f_4,{L_1}_{(2)}\oplus {L_2}_{(2)},v_4)$ is or is not a $p=2$ SHS; accordingly we keep quiet about a $p=2$ analogue of 4.9.20. 
\smallskip
So we get 4.9.8 using just: integral canonical models of Shimura quadruples having the same adjoint as the Shimura quadruple $(G,X,H,v)$ (underlying $(f,L_{(2)},v)$), involving reductive groups over $\QQ$ whose derived subgroups are quotients of $G^{\rm der}$, and whose existence is guaranteed by 4.14.3.2.1 below; warning: 4.9.8 can be worked out just in a rational context (i.e. working with Shimura isocrystals attached to points). Unfortunately, for the variants of 4.9.8 listed in 4.9.9, the rational approach is not enough (except for the case of the canonical Lie stratifications and, when appropriate, of $\rho$-stratifications). Also, in connection to 4.9.10-13 and 4.9.19-21 we have to deal only with such types of integral canonical models. Moreover, in connection to 4.9.20-21, due to limitations of 4.14.4.1 b) below, often we need to state explicitly that we assume $\Mn$ is a pro-\'etale cover of $\Mn/H_0$.
\smallskip
4.9.14-6 and 4.9.17.0 hold for the case of a $p=2$ SHS. Warning: based on reasons explained in B1-2, the toric points of 4.9.17.3 can not be always defined for $p=2$ in terms of Galois representations. The notions of 4.9.21 make sense for $p=2$. 4.9.21.1 holds for $p=2$ as far as allowed by 4.14.4.1.1 and by 4.14.4.1.2 2) below. 4.9.23 does not pertain to a fixed prime. 
\medskip
{\bf H.} 4.10 holds for $p=2$, provided we use integral canonical models $\Mn$ of the type mentioned in G (i.e. are constructed via 4.14.3.2.1 below, etc.).
\medskip
{\bf I.} No modifications are needed to 4.11 (except referring in 4.11.1 to 2.2.1.5.1 instead of [Va2, 3.2.7 4)]); due to reasons explained in 4.14.3 B1, presently 4.11.1.1 and 4.11.2-4 have to be interpreted in a sense which allows a finite number of $G$-canonical lifts of a $G$-ordinary point.
\medskip
{\bf J.} Related to 4.12 see 4.12.12.1-2 and 4.12.20. In connection to 4.12.12.6, 4.12.12.6.1, 4.12.12.6.2 B, 4.12.12.6.3, 4.12.12.6.4 and 4.12.12.6.4.1-2 we presently assume (cf. the limitations of the rational approach of paragraph G) that $(G,X,H,v)=(G_1,X_1,H_1,v_1)$ and that we are in a context where $X_{\om}^{\rm poss}=X$, $\forall\om\in W_G$: it is 4.12.12.6.0 2) which allows us to always assume this second part.
\smallskip
In particular, based on 4.12.12.6.0 2), the following sections 4.12.12.6.2, 4.12.12.6.2 B, 4.12.12.6.3, 4.12.12.6.4 and 4.12.12.6.4.1-2 hold as well for a $p=2$ standard PEL situation. Moreover, the Corollary of 4.12.12.6.4 and 4.12.12.6.4.3 hold under no restriction for $p=2$, as their arguments involve just Weyl elements and Newton polygons. 4.12.15 and 4.12.18 hold for a $p=2$ SHS (as the Theorem of 2.3.11 does).
\medskip
{\bf K.} The second variant of 3.6.14.4 can be entirely adapted for $p=2$, by just working in the context of Shimura $F$-crystals. So 4.13.1-2 remain true without any modification except the one (see end of G): often we need to state explicitly that we assume $\Mn_i$ is a pro-\'etale cover of $\Mn_i/H_{0i}$, $i=\overline{1,2}$. Accordingly, in connection to 4.13.2 we do not need to assume that there is a $W(k)$-valued point $z_i$ of $\Mn_i$ lifting $y_i$, $i=\overline{1,2}$, such that the principally quasi-polarized Shimura $2$-divisible groups over $W(k)$ obtained by pull backs via $z_1$ and $z_2$, are isomorphic: based on the surjectivity part of the Claim of 2.4.1, it is enough to start with $z_1$ and $z_2$ such their attached principally quasi-polarized Shimura filtered $\sg_k$-crystals are isomorphic. 
\medskip
{\bf 4.14.3.1. Exercise.} We refer to 2.4.1. Show that the number of points of each fibre of $m_D^G(k)$ is at least $2^{s(-1)}$ and at the most $2^{dd((M,\vph,G))}$. Here $s(-1)$ is the number of slopes $-1$ of the Shimura Lie $\sg$-crystal attached to $(M,\vph,G)$. Hints: for the at least part, use the $p=2$ variant of 4.7.11 9); for the at most part use first a similar specialization argument as in 2.3.18.1 B to reduce the situation to a Shimura-ordinary $F$-crystal over $k$ and then use 2.3.18.1.1 and arguments at the level of dimensions of $k$-vector spaces.  
\medskip
{\bf 4.14.3.2. The passage from the Hodge type to the abelian type.} We follow very closely [Va2, 6.2.2] and its corrections of AE.4 and AE.4.1. We start with the situation of [Va2, 6.2.1]. Let $f:(G,X)\to (G_1,X_1)$ be a cover such that $E(G,X)=E(G_1,X_1)$. We consider a map $(G,X,H,v)\to (G_1,X_1,H_1,v)$ defined by $f$ between two Shimura quadruples, with $v$ a prime of $E(G,X)$ dividing a rational prime $p=2$. As usual $\FF=\overline{k(v))}$.
\smallskip
Due to 2.2.1.4, the proof of [Va2, 6.2.3] holds for $p=2$: the only place where we used the fact that we were dealing with odd primes in the mentioned proof, was in the construction of a Galois-descent datum; but 2.2.1.4 guarantees that this datum can still be constructed for $p=2$. 
\smallskip
[Va2, 6.2.2 a)] for $p=2$ does apply only to the $E_6$ and $A_{2n}$ Lie types, with $n\in\NN$, cf. also [Va2, 6.4.6 6)] and AE.4.2. Accordingly, in Lemma 3 of 4.6.4, the part ``2 times" can be removed. We have the following $p=2$ analogue of [Va2, 6.2.2 b)]:
\medskip
{\bf 4.14.3.2.1. Theorem.} {\it We assume the existence of a $p=2$ SHS $(f_2,L_{(2)},v_2)$, such that $(G_2^{\rm ad},X_2^{\rm ad},H_2^{\rm ad},v_2^{\rm ad})=(G_1^{\rm ad},X_1^{\rm ad},H_1^{\rm ad},v_1^{\rm ad})$, $G^{\rm der}=G_2^{\rm der}$ and the centralizer of $G_{\ZZ_{(2)}}$ in $GL(L_{(2)})$ is a reductive group (here we use the standard notations for $(f_2,L_{(2)},v_2)$ except that we put a lower right index $2$ everywhere except for $L_{(2)}$). We also assume $\Mn_{2k(v_2)}$ has enough quasi-pivotal points, i.e. we assume each connected component of $\Mn_{2k(v_2)}$ has $\FF$-valued quasi-pivotal points. Then $(G,X,H,v)$ and $(G_1,X_1,H_1,v_1)$ have integral canonical models $\Mn$ and respectively $\Mn_1$; moreover, the natural morphism (cf. 2.2.1.5.1) $\Mn\to\Mn_1$ is a pro-\'etale cover.}
\medskip
{\bf Proof:} The existence of $\Mn$ has been already argued. As $\Mn_2$ is a quasi-projective integral model (in the sense of [Va2, 3.2.3 2')]; the argument is as in the proof of [Va2, 6.4.2]), [Va2, 6.2.3] implies $\Mn$ is a quasi-projective integral model. To show that $\Mn_1$ exists, we follow the proof of [Va2, 6.2.2 b)]: we construct $\Mn_1$ as the quotient of $\Mn$ via a group action of a pro-finite, abelian group which is an $M$-torsion group, for some $M\in\NN$; here we just mention that $2$ does not divide $M$ iff $G^{\rm ad}$ has all its simple factors of some $A_{2n}$ Lie type, $n\in\NN$ (cf. also [Va2, 6.2.3.1]).
\smallskip
As any connected component $\Mc^0$ of $\Mn_{W(\FF)}$ can be identified with a connected component $\Mc_2$ of $\Mn_{2W(\FF)}$, the fact that $\Mn_1$ has the EP, follows once we show that the mentioned group action is free (cf. [Va2, 6.2.2 A) to D)]). We choose a connected component $X^0$ of $X$. We can assume it is a connected component of $X_2$ as well (cf. [Va2, 3.3.3]). We can choose $\Mc^0=\Mc_2$ such its set of complex points contains those defined by equivalence classes of the form $[x,1]$, $x\in X^0$. We follow entirely the pattern of AE.4.1. Let $(G_2^\prime,X_2^\prime)$, $h$, $y$, $GA$, $(A_y,p_{A_y})$, $(A_y^\prime,p_{A_y^\prime})$, $M$, $a$ and $a_M$ have the same significance as in AE.4.1. We assume $h$ does not fix all geometric points of $\Mc_2$. If $y$ is a quasi-pivotal point, then we reach a contradiction as in AE.4.1 (cf. 4.14.3 B3). So, based on our hypothesis on the existence of quasi-pivotal points, we get that $h$ does not act trivially on any connected component of $\Mc_{2\FF}$. If $y$ is not a quasi-pivotal point, then as in AE.4.1 we get that $a_M$ leaves invariant the $F^1$-filtration of $M$ defined by any $W(\FF)$-valued point $z$ of $\Mc_2$ lifting $y$. Based on 2.3.5.6.1 B (INTR) (see also end of 2.3.18 A) and on the proof of 2.4.2, we get that $h$ acts trivially on the connected component of $\Mc_{2\FF}$ to which $y$ belongs. Contradiction. This ends the proof.
\medskip
For future references, we also mention here explicitly the following $p=2$ version of the corrected version of [Va2, 6.2.2.1] mentioned in the Proposition of AE.4. Let 
$$g\in {\rm Aut}(G_2,X_2,H_2)\subset Aut(G_2)(\QQ)$$ 
such that it takes $X^0$ into $X^0$ and normalizes a compact, open subgroup $H_{02}$ of $G_2(\AA_f^2)$. So $g$ gives birth (via 2.2.1.5.1) to an automorphism $a_g$ of the image $\Mc_{2H_{02}}$ of $\Mc_2$ into $\Mn_{2W(\FF)}/H_{02}$. 
\medskip
{\bf 4.14.3.2.2. Proposition.} {\it We assume the order of $a_g$ is finite. If $a_g$ acts freely on $\Mc_{2H_{02}}$ but fixes an $\FF$-valued point $y$ of it, then the power of $p$ dividing the order of $a_g$ is bounded above independently of $H_{02}$ and $g$.} 
\medskip
{\bf Proof:} Let $g_0$ have the same significance as in [Va2, p. 495]; so it is an automorphism of the non-trivial factor of the Shimura adjoint Lie $\sg_{\FF}$-crystal attached to $y$. Let $z\in\Mc_{2H_{02}}(W(\FF))$ be a lift of $y$ and let ${\got C}_z:=(M_y,F^1_z,\vph,G_{2W(\FF)})$ be its attached Shimura filtered $\sg_{\FF}$-crystal. Let ${\rm Spec}(R)$ be the completion of the local ring of $y$ and let ${\rm Spec}(R_1)$ be obtained as in 2.4, starting from ${\got C}_z$. From the Claim 2.4.1 and from 4.14.3.1 we get that the $W(\FF)$-morphism $FIL_0:{\rm Spec}(R^{PD})\to {\rm Spec}(R_1^{PD})$ is such that each $W(\FF)$-valued point of ${\rm Spec}(R_1^{PD})$ lifts to precisely $2^m$ $W(\FF)$-valued points of ${\rm Spec}(R^{PD})$, where $m\in\NN\cup\{0\}$ is at most $dd(({\got C}_z))$. 
\smallskip
We can assume $g_0$ has an order a power of $p$. As AE.6.1 points out, we can not always argue immediately that $g_0$ is an inner automorphism even if $g$ is. However, as the group of outer automorphisms of any adjoint group over $W(\FF)$ is finite, we get that there is $b(G^{\rm ad}_{2W(\FF)})\in\NN$ which does not depend on $g$ or on $H_{02}$ such that $g_0^{b(G^{\rm ad}_{2W(\FF)})}$ is an inner automorphism. So $g_0^{b(G^{\rm ad}_{2W(\FF)})t(G^{\rm ad}_{2W(\FF)})}$ is the identity automorphism. But this together with the existence of $m$ implies that the order $o$ of the $W(\FF)$-automorphism of ${\rm Spec}(R)$ defined by $a_g$ is bounded above independently of $H_{02}$ and $g$.  We have (cf. 2.4.2): 
$$o\le 2b(G^{\rm ad}_{2W(\FF)})t(G^{\rm ad}_{2W(\FF)}).$$
This ends the proof.
\medskip
{\bf 4.14.3.2.3. Back to 2.3.18 A.} We use the notations of the beginning of 2.3.18. We assume the centralizer of $G_{\ZZ_{(2)}}$ in $GL(L_{(2)})$ is a reductive group. We also assume $\Mn_{k(v)}$ has enough quasi-pivotal points; for instance, this is so if the connected components of $\Mn_{\FF}$ are permuted transitively by $G(\AA_f^p)$, cf. 4.14.3 J (in particular, this applies if $(G,X)$ is of compact type). From 4.14.3.2.1 we get that $\Mn^{\rm ad}$ exists and moreover the $O_{(v)}$-morphism $q_{\Mn}$ (resp. $q_{\Mn_{1}}$) obtained as in 2.3.3.2 (resp. as in 2.3.5.5) exists and is a pro-\'etale cover of its image (resp. is a pro-\'etale cover). 
\smallskip
Similarly, Fact 2 of 2.3.5.2 holds for $p=2$ due to the mentioned $p=2$ version of [Va2, 6.2.3] (see 4.14.3.2). 
\medskip
{\bf 4.14.3.3. References for 1.13.4.} 1.13.4 a) follows from 4.14.3 A and B. For 1.13.4 b) cf. 4.14.3.2.1 (or 4.14.4.1.1 below). 1.13.4 c) (resp. d) follows from 4.14.3 C (resp. J). 1.13.4 e) follows from 4.14.3 B and D.
\medskip
{\bf 4.14.4. Generalized Serre lemma.} 4.14.3.2.1-2 point out that [Va2, 6.4.6 6)] has a version in mixed characteristic $(0,2)$. To state it as a Corollary, we need to adapt AE.4.2 and 4.4.13.5 to the context of 4.14.3.2.1. 
\smallskip
We refer to the mentioned context. Due to the limitations of 4.14.3 G, here we define a non-integral, $U$-ordinary point of $\Mn_{1k(v_1)}$ to be a point with values in a field $k_1$  such that a $G_1(\AA_f^2)$-translate of it lifts, after potentially replacing $k_1$ by $\overline{k_1}$ and identifying connected components of $\Mn_{W(\FF)}$ with ones of $\Mn_{2W(\FF)}$, to a non-integral, $U$-ordinary point of $\Mn_{2\FF}$. We have:
\medskip
{\bf 4.14.4.1. Corollary.} {\it Let $H_{01}$ be a compact subgroup of $G_1(\AA_f^2)$ such that $H_{01}\times H_1$ is smooth for $(G_1,X_1)$. 
\medskip
{\bf a)} Then $\Mn_1$ is a pro-\'etale cover of $\Mn_1/H_{01}\times H_1$ above all points of $\Mn_{1k(v_1)}$ which are non-integral, $U$-ordinary points.
\smallskip
{\bf b)} We assume the centralizer of $G_{2\ZZ_{(2)}}$ in $GL(L_{(2)})$ is a reductive group over $\ZZ_{(2)}$ and $\Mn_{2k(v_2)}$ has enough quasi-pivotal points. If $H_{01}\times H_1$ is $2$-smooth for $(G_1,X_1)$, then $\Mn_1$ is a pro-\'etale cover of $\Mn_1/H_{01}\times H_1$.}
\medskip
{\bf Proof:} The proof of a) is entirely the same as 4.4.13.5, cf. 4.14.3 B3. b) is obtained from 4.14.3.2.1 in the same way (see 2.3.8.5) AE.4.2 b) was deduced from AE.4.1.
\medskip
{\bf 4.14.4.1.1. A conclusion.} The concrete conclusion of the methods developed above in connection to the $p=2$ theory of Shimura varieties, can be gathered in the following particular case of 4.14.4.1 b):
\medskip
{\bf Corollary.} {\it We consider a Shimura quadruple $(G_1,X_1,H_1,v_1)$ of abelian type, with $v_1$ dividing $2$. We assume each simple factor of $(G_1^{\rm ad},X_1^{\rm ad})$ is either of some $A_n$ Lie type or is of some $C_n$ or $D_n^{\HH}$ type defined by an adjoint group whose extension to $\RR$ has no simple, compact factors. We consider a $p=2$ standard PEL situation $(f,L_{(2)},v,\Mb))$, with $(G,X,H,v)=(G_1^{\rm ad},X_1^{\rm ad},H_1^{\rm ad},v_1^{\rm ad})$ (cf. Lemma 2 of 4.6.4 and the $p=2$ version --see 2.3.18 A-- of 2.3.5.1). We also assume $\Mn_{k(v)}$ has enough quasi-pivotal points (for instance this holds if $(G,X)$ is of compact type). Then $(G_1,X_1,H_1,v_1)$ has an integral canonical model which is quasi-strongly smooth.}
\medskip
{\bf Proof:} As we are dealing with a $p=2$ standard PEL situation  $(f,L_{(2)},v,\Mb))$, the centralizer of $G_{\ZZ_{(2)}}^{\rm der}$ in $GL(L_{(2)})$ is a reductive group and so, as in Fact 2 of 2.3.5.5, we deduce the existence of a standard cover defined by $(f,L_{(2)})$. So the Corollary follows from 4.14.3.2.1.
\medskip
{\bf 4.14.4.1.2. Remarks.} {\bf 1)} The extra condition of 4.14.3.2.1 and 4.14.4.1 b) referring to centralizers, can be eliminated, cf. 2.3.5.6.1 A. However, in connection to the $p=2$ theory of Shimura varieties to be fully developed in \S6 and [Va5], this is irrelevant (as the $p=2$ SHS's to be ``constructed" will follow the pattern of [Va2, 6.5-6]).   
\smallskip
{\bf 2)} We take $p\ge 2$. We consider a Shimura quadruple $(G_0,X_0,H_0,v_0)$ of adjoint, abelian type, with $v_0$ dividing $p$, such that any simple factor $(G_1,X_1)$ of $(G_0,X_0)$ is of $C_n$ or $D_n^{\HH}$ type and satisfies the following property:
\medskip
{\bf The essentially non-compact property.} {\it Every simple factor $SF$ of $G_{1\QQ_p}$ is such that the Shimura group pair $(G_{0\ZZ_p},[\mu_0])$ defined by $(G_0,X_0,H_0,v_0)$ is such that $\mu_{0B(\FF)}$ has either a trivial image in $SF_{B(\FF)}$ or has a non-trivial image in all simple factors of $SF_{B(\FF)}$.}
\medskip
Then Lemma 2 of 4.6.4 can be entirely adapted to it: we just have to replace ($p=2$) standard PEL situations by ($p=2$) SHS's. The only two modifications in its proof needed to be made are as follows. First, $C_{V_0}$ is a product of $GL(V_0^{2n})$ copies whose number may exceed (in general) the number of simple factors of $G^{\rm ad}_{\RR}$. Second, we can work everything in the context of $G^1_{\ZZ_p}$ instead of $G_{\ZZ_p}$, where $G^1_{\ZZ_p}$ is the subgroup of $G_{\ZZ_p}$ generated by the maximal torus of $GSp(L_{(p)}\otimes_{\ZZ_{(p)}} \ZZ_p,\psi)$ whose extension to $V_0$ is contained in the center of $C_{V_0}$ and by each factor of $G^{\rm der}_{\ZZ_p}$ having a $\ZZ_p$-simple adjoint whose generic fibre is an $SF$ as above. 
So the alternating representation of $G^1_{\ZZ_p}$ on $L_{(p)}\otimes_{\ZZ_{(p)}} \ZZ_p$ gets decomposed into alternating subrepresentations, each one of them being a $\ZZ_p$-version of the $(p=2$) standard PEL situations; Warning: some of them are related to $\ZZ_p$-groups which are tori (i.e. are not necessarily of $C_n$ or $D_n$ Lie type). So [Va2, 5.6.3] (cf. also AE.1) for $p\ge 3$ and 2.3.18 B for $p=2$ can be applied entirely (using the language of abelian schemes endowed with $\ZZ_p$-endomorphisms and not with $Z_{(p)}$-endomorphisms) to get that $(f,L_{(p)},v)$ is a ($p=2$) SHS. 
\smallskip
Moreover, we have a version of 4.14.4.1.1 in the context of abelian quadruples having $(G_0,X_0,H_0,v_0)$ as their adjoints.  
\smallskip
{\bf 3)} [Va2, 6.4.11 D] still holds for the compact type cases of 2) or of 4.14.4.1.1, cf. the part of 2.2.5.1 referring to [Va2, 3.2.12] and the fact that the arguments of [Va2, 6.4.1.1 2)] hold as well for $p=2$.
\medskip
{\bf 4.14.4.2. Serre lemma via the $p$-adic approach.} We now show that 4.14.4.1 b) and AE.4.2 imply the classical Serre lemma of [Mu, p. 207]. Let $(A,p_A)$ be a polarized abelian variety over an algebraically closed field $\tilde k$ of arbitrary characteristic. Let $N\in\NN\setminus\{1,2\}$. Let $a$ be an automorphism of $(A,p_A)$ acting trivially on $A[N]$. We need to show that $a$ is the trivial automorphism. 
\smallskip
We can assume $p_A$ is a principal polarization (cf. Zarhin's trick; see [M-B, p. 205]). If $\tilde k$ is of characteristic $0$, then we can assume it is $\CC$ and so the fact that $a$ is trivial can be easily checked starting from the classical theorem of Riemann. But, with the notations of [Va2, 4.1], this can be restated as: the subgroup $K(N)$ of 2.3.3 is smooth for $(GSp(W,\psi),S)$ (for instance, this can be read out from Artin's method; see [FC, p. 23-4]: these two pages have an analogue which considers level-$N$ symplectic similitude structures). Moreover, as the image of $K(N)$ in $PSp(W,\psi)(\QQ_l)$ is an open subgroup having no torsion, it is $p$-smooth for $(GSp(W,\psi),S)$ as well; here $l$ is either an odd prime dividing $N$ or is $2$, depending on the fact that $4$ does not or does divide $N$. So assuming that the characteristic of $\tilde k$ is a prime $p$ not dividing $N$, from AE.4.2 (resp. from 4.14.4.2 b)) for $p\ge 3$ (resp. for $p=2$) applied to Siegel modular varieties, we get $\Mm/K^p(N)$ is smooth over $\ZZ_{(p)}$ and $\Mm$ is a pro-\'etale cover of it; here $\Mm$ is as in [Va2, 3.2.9] and $K^p(N)$ is as in 2.3.3. But this implies (for instance, via [FC, p. 23-4]) that $a$ is the trivial automorphism.
\smallskip
We are left with the cases: $N=p$ is odd or $N=4$ and $p=2$. These cases are trivial: as $a$ has finite order, the automorphism of $H^1_{\rm crys}(A/W(\tilde k))$ defined by $a$ is trivial and so $a$ acts trivially on $A(N^m)$, $\forall m\in\NN$ (cf. also 2.3.18.1 D for $N=4$). So $a$ is trivial.
\smallskip
There is one extra thing worth pointing out. The mentioned Artin's method can be used to get that, provided $p$ does not divide $N$, there is $l\in\NN$ such that $\Mm/K(N^l)$ is smooth over $\ZZ_{(p)}$ and $\Mm$ is a pro-\'etale cover of it: there is $q(a,A)\in\NN$ such that $a$ does not act trivially on $A[N^{q_A}]$; the fact that the stack of principally polarized abelian schemes of dimension equal to $\dim_{\tilde k}(A)$ is of finite type over ${\rm Spec}(\ZZ)$ (see [FC, p. 23-4]) implies that we can choose $q(a,A)$ to be independent of $a$ and $A$. Using this we get the existence of integral canonical models of Siegel modular varieties, without appealing to Serre lemma. This also shows that any automorphism of a polarized abelian variety over a field has finite order.
\medskip
{\bf 4.14.4.3. Comment.} In the classical situation of Siegel modular varieties, most common one works with some level-$N$ structure (not necessarily a symplectic or symplectic similitude one), see [FC], [MFK], etc., or a very slight variant of it (see the classical case of elliptic modular curves). However this has considerable limitations and disadvantages in the general context of Shimura varieties of preabelian type. This forced us to introduce in [Va2, 2.11] the abstract notion of smooth subgroup for a given Shimura pair; Definition 1 of AE.4 is just a logical continuation of loc. cit. Some of the advantages we gain working in this abstract, general context can be read out from the easiness of stating different results: see [Va2, 6.4.4 and 6.4.6 1) and 2)], AE.4.1-2, 4.14.4.1, etc. 
\medskip
{\bf 4.14.5. Some extensions.} Let $k$ be a perfect field. The greatest part of the terminology used in all previous sections can be extended to (filtered) $\sg_k$-$\Ms$-crystals. In particular we get $\sg_k$-$\Ms$-crystals which are $U$-ordinary, $T$-ordinary, $G$-ordinary, potentially (or quasi or strongly) cyclic diagonalizable (of level $n$), etc. Similarly, in a geometric context we get $U$-ordinary, $T$-ordinary, $G$-ordinary, null points, etc; see [Va9] for an extensive study. Also, notions like deviations (see 3.13 and 4.5.15.0), $a$-Lie $\sg_k$-crystals attached to $\sg_k$-$\Ms$-crystals and $[a,b]$-filtered Lie $\sg_k$-crystals attached to filtered $\sg_k$-$\Ms$-crystals (see 2.2.3 2)), CM levels, NP variation functions, etc., do extend as well. These extensions from the context of $\Mm\Mf_{[0,1]}(W(k))$ to the contexts of other categories $\Mm\Mf_{[a,b]}(W(k))$ are automatic. 
\medskip
{\bf 4.14.6. The complete list of invariants.} Let $(f,L_{(v)},v)$ be a ($p=2)$ SHS. In 2.3 and 4.1-5 (cf. also 4.14.3 for $p=2$) we introduced many invariants attached to it. Here, for the sake of future references, we list them all. These invariants are: rational numbers, $n$-tuples of non-negative integers, sets, stratifications, $F$-crystals over perfect fields, groups, etc. Grouping them in some logical way, they are:
\medskip
a) $E(G,X)$ and $e$ (see 2.3.1);
\smallskip
b) $\Mf(G^{\rm ad})$ and $q_{G^{\rm ad}}$ (see 2.3.5.2);
\smallskip
c) $d(f)$, $d_0(f)$, $d_1(f)$,..., $d_6(f)$, $d(f_p)$, $d_0(f_p)$, $d_3(f_p)$ and $d_4(f_p)$ (see 2.3.6.1);
\smallskip
d) $\Ml$ and $\Mg$ (see 2.3.11);
\smallskip
e) $\mu$, $T_{\mu}$, ${\got C}_{(f,v)}$, $\tau$, $r(f,v)$, $FL(f,v)$, $IL(f,v)$, $SKL(f,v)$, $CL(f,v)$, $SDL(f,v)$, ${\rm Lie}_G(\tau)$, $r(v^{\rm ad})$ and $SS(v^{\rm ad})$ (see 4.1, 4.1.1 and 4.1.1.2);
\smallskip
f) $d=d(v)$ (see 4.1.1.1);
\smallskip
g) $\Ml$, $\bar\vph$, $I_{\bar\vph}$, $d_l$'s and $m_{(f,v)}(\tau_q)$'s (see 4.1.1.4 and 4.1.2.2; here $l\in I_{\bar\vph}$, $q\in\NN$ divides $d$ and $\tau_q\in\Mt_q$);
\smallskip
h) $W_G$, $R_G(v)$, $\mu_{\om}$, $T_{\mu_{\om}}$, ${\got C}_{\om}$, $\Mp_{\om}$, ${\rm Lie}_G(\Mp_{\om})$ and $\tau_{\om}$'s (see 4.1.5 and 4.1.5.4-5; here $\om\in W_G$);
\smallskip
i) $k(v^{\rm sp}_G)$, $d_{\om}$'s, $SDD(f,v)$ and $\Mm\Mm$ (see 4.1.5.2-3 and 4.1.6; here $\om\in W_G$);
\smallskip
j) $\Mh$, $\Mh^{\rm nc}$, $\Mh^{\rm c}$, $\Mh_i$'s, $d_i$'s, $\gamma_p(G^{\rm ad})$ and $I_p(G^{\rm ad})$ (see 4.3.1 and 4.3.1.1; here $i\in\Mh$); 
\smallskip
k) $d_A(i)$'s, $\vep_i(v^{\rm ad})$'s, $\vep_i^{\rm c}(v^{\rm ad})$'s, $d_A(v)$ and $d_T(v)$ (see 4.3.1.1-2); here $i\in\Mh^{\rm nc}$);
\smallskip
l) $R(\tau)$, $R^r(\tau)$, $a_i(v^{\rm ad})$, $a(v^{\rm ad})$, $a^r(v^{\rm ad})$, $D_i^+(v^{\rm ad})$'s and $D_i^-(v^{\rm ad})$'s (see 4.3.1.2; here $i\in\Mh^{\rm nc}$);
\smallskip
m) $SLSp-rk(\Mn_{k(v)})$ (see 4.3.3);
\smallskip
n) $NP(\Mn_{k(v)})$, $LNP(\Mn_{k(v)})$, $RLNP(\Mn_{k(v)})$, $N(G^{\rm ad},X^{\rm ad},v^{\rm ad})$ and $N_1(G^{\rm ad},X^{\rm ad},v^{\rm ad})$ (see 4.5.7);
\smallskip
o) the standard stratifications of 4.5 and 4.9.9.1: the canonical Lie, the refined canonical Lie, the $\rho$-stratification, the absolute, the $f$-canonical, the Lie stable, the refined Lie stable, the Lie non-stable, the pseudo-ultra, the quasi-ultra (i.e. the Faltings--Shimura--Hasse--Witt), the ultra, the Faltings--Shimura--Dieudonn\'e (adjoint, principal or standard), the ultimate and the ultimate adjoint stratifications of $\Mn_{k(v)}$ and the level-$n$ Faltings--Shimura--Dieudonn\'e equivalence relation of $\Mn_{k(v)}(\FF)$ (here $n\in\NN$);
\smallskip
p) $f_G$ (see 4.5.11 a));
\smallskip
q) $d(i)$'s and $d^i(\om)$'s (see 4.5.15.2.3; here $i\in\Mh$);
\smallskip
r) $\GG\SS\CC(\Mn_{k})$ (see 4.9.9.2; here $k$ is either $k(v)$ or an algebraically closed field containing $k(v)$);
\smallskip
s) ${\rm aut}_i(\om)$'s (see 4.5.15.3; here $i\in\Mh$ and $\om\in W_G$);
\smallskip
t) $NSU(\Mn_{k(v)})$ and $W_G(i)$ (see 4.12.12.6.6 3); here $i\in\Mh^{\rm nc}$);
\smallskip
u) $c_p(\tilde H_0)$, $c_p(\tilde H_0,v)$ and ${\rm piv}_{(G,X,v)}(\tilde H_0)$ (see 2.3.3 and 4.12.12.6.7; here $\tilde H_0$ is as in 2.3.3);
\smallskip
v) $D_q(G,X,H,v,\tilde H_0)$ (see 4.12.12.8).
\medskip
{\bf 4.14.6.1. Class invariants.}
Let $\FF(\dim_{\CC}(X))$ be the algebraically closed field of transcendental degree $\dim_{\CC}(X)$ over $\FF$. So any point of $\Mn_{k(v)}$ with values in an algebraically closed field, factors through a point of $\Mn_{k(v)}$ with values in a subfield of $\FF(\dim_{\CC}(X))$. Another reason to introduce $\FF(\dim_{\CC}(X))$ is: we do not stop to check that all class invariants to be defined below are preserved by extensions (to algebraically closed fields). We consider the class 
$$CL_{(f,v)}:=Cl(M\otimes_{W(k(v))} W(\FF(\dim_{\CC}(X))),\vph\otimes 1,G_{W(\FF(\dim_{\CC}(X))}),$$ 
see 2.2.22 2) and 4.1.1. We assume $G$ is not a torus; so $CL_{(f,v)}$ has at least $2$ elements (cf. 4.3.6 and its proof). 
\smallskip
The maximal CM level of $CL_{(f,v)}$ is denoted by $CML(CL_{(f,v)})$ (see 3.13.7). The smallest $n\in\NN$ such that:
\medskip
-- the isomorphism deviation of any representative of an arbitrary element of $CL_{(f,v)}$ is less or equal to $n$, is denoted by $ISOM(CL_{(f,v)})$ (cf. 3.15.7 BP1);
\smallskip
-- the Newton polygon (resp. the Newton polygon Lie) deviation of any representative of an arbitrary element of $CL_{(f,v)}$ is less or equal to $n$, will be denoted by $NP(CL_{(f,v)})$ (resp. by $NPL(CL_{(f,v)})$ (its existence is implied by [Ka2, 1.4.4] or by the existence of $ISOM(CL_{(f,v)})$);
\smallskip
-- the number of slopes of any representative of an arbitrary element of $CL_{(f,v)}$ is less or equal to $n$, is denoted by $NSL(CL_{(f,v)})$;
\smallskip
-- the denominator of any representative of an arbitrary element of $CL_{(f,v)}$ is less or equal to $n$, is denoted by $DEN(CL_{(f,v)})$;
\smallskip
-- the cyclic rank (see 3.9.1.0) of any representative of an arbitrary element of $CL_{(f,v)}$ is less or equal to $n$, is denoted by $CR(CL_{(f,v)})$.
\medskip
The non-negative integers $CML(CL_{(f,v)})$, $CR(CL_{(f,v)})$, $NP(CL_{(f,v)})$, $NPL(CL_{(f,v)})$, $NSL(CL_{(f,v)})$, $DEN(CL_{(f,v)})$ and $ISOM(CL_{(f,v)})$ are referred as the standard class invariants of $(f,L_{(p)},v)$.
\medskip
{\bf 4.14.6.2. Types of points.} For most useful types of points (with values in fields, in some cases --like the one of Borel points-- assumed to be perfect) of special fibres of integral canonical models of Shimura varieties of preabelian type we refer to 4.9.17.2-3 and to 4.12.12.6.3-4; cf. also 4.14.3 for the case of special fibres in characteristic $2$.
\medskip
{\bf 4.14.6.3. Remark.} For each point $w$ of $\Mn/H_0$ with values in some field $F$, we can attach to $w^*(\Ma_{H_0},\Mp_{\Ma_{H_0}})$ an incredible amount of invariants. See [Va4] and [Va10] for samples.
\vfill\eject
\centerline{}
\bigskip
\bigskip
\centerline{\bigsll {\bf AE Addenda and errata to [Va2]}}
\bigskip\bigskip
\medskip
We have the following addenda and errata to [Va2]; all references below till the beginning of AE.1 are made to loc. cit. We use ``$\rightsquigarrow$" as a substitute for ``should be replaced by". AE.0 is a list of replacements. In AE.1-2 and AE.5 we mention three inexactitudes. Two mistakes are handled in AE.3-4. AE.6 consists of some extra specifications meant to help the reader less familiar with (group) schemes.
\medskip
{\bf AE.0.} Always by ``the closure" in some scheme we meant ``the Zariski closure". On p. 404 l. 17 (resp. l. 23) ``prove" $\rightsquigarrow$ ``proof" (resp. ``connected component" $\rightsquigarrow$ ``maximal subtorus"). On p. 410 l. 7 ``DMF" $\rightsquigarrow$ ``DMS". In 1.4.1 (resp. in 2.7), ``4.6.10" $\rightsquigarrow$ ``6.4.10" (resp. ``$\overline{T(Q)}$" $\rightsquigarrow$ ``$\overline{T(\QQ)}$"). In 2.4.3 and of 3.2.7 9), ``type" $\rightsquigarrow$ ``type and all its simple factors are such that [De2, 1.2.8 (ii)] applies" and ``(of finite index)" $\rightsquigarrow$ ``(it is of finite index if $G$ is an adjoint group)". At the bottom of p. 415 (resp. in 2.6), ``an" $\rightsquigarrow$ ``a suitable" (resp. ``action" $\rightsquigarrow$ ``action (via conjugation)"). 
On p. 423 l. 15 and l. 29 ``$\widehat{Y_2}$" $\rightsquigarrow$ ``$\widehat{Y}_2$". In 3.1.4, ``its centralizer" $\rightsquigarrow$ ``the connected component of the origin $C$ of its centralizer".
\smallskip
In 3.2.1.1 3) (resp. on. p. 427 l. 3), ``3.2.1 3)" $\rightsquigarrow$ ``3.2.1 4)" (resp. ``almost" $\rightsquigarrow$ ``almost healthy"). 
In 3.2.2 1) (resp. in 3.2.2.2), ``this" $\rightsquigarrow$ ``this and its proof" (resp. ``of the limit" $\rightsquigarrow$ ``of the projective system"). In 3.2.2.4 b) (resp. in 3.2.3.1 1)), ``3.2.3.2.1" $\rightsquigarrow$ ``3.2.2.3.1" (resp. ``3.2.2 1)" $\rightsquigarrow$ ``3.2.2 3)"). 
\smallskip
On p. 435, ``for 2) above" $\rightsquigarrow$ ``as in 3.2.3.1 2)". In b) of 3.2.7 8), ``for Shimura" $\rightsquigarrow$ ``for generalized Shimura",  ``, f) and g)" $\rightsquigarrow$ ``, and a variant of f) and g)", and ``$,F^1$" should be deleted. On p. 451 l. 1, ``[M," $\rightsquigarrow$ ``[Ma,". On p. 454 l. 13 (resp. l. 40), ``$g^{-1}\circ k$ of $H_1(A,\QQ)\otimes\AA_f$" $\rightsquigarrow$ ``$g^{-1}\circ k:H_1(A,\QQ)\otimes\AA_f\tilde\to W\otimes\AA_f$" (resp. ``isomorphism" $\rightsquigarrow$ ``isomorphisms").  In 4.1.3, ``$F^0(H^1_{dR}(A/Z)\otimes H^1_{dR}(A/Z)^*)^{\otimes m(\al)}$" $\rightsquigarrow$ ``$F^0((H^1_{dR}(A/Z)\otimes H^1_{dR}(A/Z)^*)^{\otimes m(\al)})$". In 4.2, ``the elements of their centers are semisimple endomorphism of $W$" $\rightsquigarrow$ ``their centers are generated by semisimple endomorphisms of $W$ and the trace forms on their centers are perfect" and ``(cf. [Bou1, ch. 1.1, th.4])" should be deleted.
\smallskip
On p. 461 l. 34, ``3.4.6 3)" $\rightsquigarrow$ ``4.3.6 3)" and ``to its" $\rightsquigarrow$ ``by its". On p. 463 l. 19 (resp. 22), ``a natural" $\rightsquigarrow$ ``naturally an" (resp. ``$[[[$" $\rightsquigarrow$ ``$[[$"). On p. 464 l. 1 and l. 24, ``component" $\rightsquigarrow$ ``components". On p. 464 l. 13, ``automorphism" $\rightsquigarrow$ ``automorphisms". On p. 464 l. 34 (resp. l. 41), ``subgroup" $\rightsquigarrow$ ``$\GG_a$ subgroup" (resp. ``isomorphism:" $\rightsquigarrow$ ``isomorphism onto its image $\GG_{a,\al}$:"). On p. 468 l. 24, ``4.3.4" $\rightsquigarrow$ ``4.3.2". In 4.3.14, ``the closure of $G_0$ in $GSp(L,\psi)$" $\rightsquigarrow$ ``a reductive group having the Zariski closure of $G_0$ in $GSp(L,\psi)$ as a reductive subgroup containing its derived subgroup" and ``4.3.13" $\rightsquigarrow$ ``4.3.13 (cf. also 4.3.9)". On p. 471 l. 19, ``situation" $\rightsquigarrow$ ``PEL situation". In 4.3.17, ``group $\tilde H$" $\rightsquigarrow$ ``subgroup $\tilde H$".
\smallskip
On p. 472 l. 7, ``3.2.7" $\rightsquigarrow$ ``3.2.7.1". In 5.2.1, ``$n\in\NN$." $\rightsquigarrow$ ``$n\in\NN\cup\{0\}$.". On p. 475 l. 2 and 4, in order to be match [Fa2], ``p-adic" should be removed. In (5.2.7), ``$p$" $\rightsquigarrow$ ``$p\be$".  In 5.2.3, ``$\xi:=f$" $\rightsquigarrow$ ``$\xi:=f_e$". In 5.2.13, ``$\nabla$-parallel" $\rightsquigarrow$ ``annihilated by $\nabla$". On p. 476 l. 12 (resp. l. 33), ``theory" $\rightsquigarrow$ ``comparison theory" (resp. ``have used" $\rightsquigarrow$ ``used"). On p. 477 l. 37, ``positive" $\rightsquigarrow$ ``non-negative". On p. 479 l. 1 (resp. l. 19), ``as defined" $\rightsquigarrow$ ``defined" (resp. ``ring homomorphism" $\rightsquigarrow$ ``$V_0$-homomorphisms"). On p. 480 l. 11 (resp. l. 40), ``$\nu\vert Y:=i_0^*\nu$" $\rightsquigarrow$ ``$\nu\vert Y_0:=i_0^*(\nu)$" (resp. ``$pi^{e-1}T$" $\rightsquigarrow$ ``$\pi^{e-1}T$").
\smallskip
On p. 488 l. 9 (resp. l. 24), ``groups" $\rightsquigarrow$ ``groups of classical Lie type" (resp. ``(cf. [De2, 1.3.6])" should be deleted). On p. 489, ``4.3.10.1)" $\rightsquigarrow$ ``4.3.10.1 1)", ``$\Ma(G,X,W)=1$") $\rightsquigarrow$ ``$\Ma(G,X,W)$ is 1 or 2 depending on the fact that $l$ mod $4$ is or is not $1$ or $2$", and ``$G_2^{\rm der}=G_1^{\rm der}$, but $G_2^{\rm ab}$ is a torus of dimension 2 (the representation ${G_1^{\rm der}}_{\CC}\to GL(W_{\CC})$ is not irreducible, just the representation ${G_1^{\rm der}}_{\RR}\to GL(W_{\RR})$ is irreducible)" $\rightsquigarrow$ ``$G_2^{\rm der}$ is isogeneous to $G_1^{\rm der}$ times an $SL_2$-group, and $G_2^{\rm ab}$ is a torus of dimension 1 (the representation ${G_1^{\rm der}}_{\RR}\to GL(W_{\RR})$ is not irreducible; see [Sa, p. 458])". The last sentence of 5.8.1 should be deleted. On p. 491 l. 15 (resp. l. 32), ``$A_{p+1}$" $\rightsquigarrow$ ``$A_{p-1}$" (resp. ``this" $\rightsquigarrow$ ``the corollary").
\smallskip
In 6.2.2 D), ``and so $S_0$ is a pro-\'etale cover of $S_1$. This implies that" $\rightsquigarrow$ ``and so, provided the morphism $\Mm\to\Mm_1$ is a pro-\'etale cover, $S_0$ is a pro-\'etale cover of $S_1$ and". In 6.2.2 E) ``complex points are" $\rightsquigarrow$ ``set of complex points contains those".  On p. 498 l. 32, ``3.2. 7 10)" $\rightsquigarrow$ ``3.2.7 10)". In 6.2.6 2), ``6.2.4.1" $\rightsquigarrow$ ``6.2.4". In 6.4.1.1 2), ``[FC, ii)" $\rightsquigarrow$ ``[FC, iv)". On p. 503, ``6.4.1.1 2) and of [Mi4, 4.13]" $\rightsquigarrow$ ``[Mi4, 4.13], via the same argument used in the first key fact of the proof of 3.2.3.2 b)". On p. 504 l. 6 (resp. l. 12), ``isomorphism" $\rightsquigarrow$ `isomorphisms" (resp. ``preabelian." $\rightsquigarrow$ ``preabelian type."). On p. 504 l. 36, ``6.2.2." $\rightsquigarrow$ ``6.2.2 a).". In 6.4.9, ``healthy" $\rightsquigarrow$ ``any healthy". On p. 506 l. 15, ``simple" $\rightsquigarrow$ ``simply".
\smallskip
On p. 507 l. 21, ``forms" $\rightsquigarrow$ ``alternating forms" and ``$\tilde G_{3\ZZ_{(p)}}$" $\rightsquigarrow$ ``$\tilde G^{c+\rm der}_{3\ZZ_{(p)}}$". On p. 507 l. 30, ``6.6.5.1" $\rightsquigarrow$ ``i)" and ``component" $\rightsquigarrow$ ``component of the origin". On p. 508 l. 1, ``i" $\rightsquigarrow$ ``$i$". In 6.6.2 d) and in 6.6.4, ``is the centralizer" $\rightsquigarrow$ ``is contained and has the same derived subgroup as the centralizer". In b) of p. 510, ``$B(\tilde G^{\rm ad}_{E_{(p)}})$" $\rightsquigarrow$ ``$B(\tilde G^{\rm ad}_{E})$". In e) of p. 511 ``is the centralizer of a torus $\tilde T$ of $\tilde G_{E_{(p)}}$" $\rightsquigarrow$ ``is the derived subgroup of the centralizer of a torus $\tilde T$ of $\tilde G_{E_{(p)}}$ in $\tilde G_{E_{(p)}}$". In 6.6.5.1, ``4.3.6.2)" $\rightsquigarrow$ ``4.3.6 2)" (resp. ``keeping the" $\rightsquigarrow$ ``keeping") and ``. In fact, referring" $\rightsquigarrow$ ``: referring". On p. 509, ``so that $G^s$ splits over the completion of $F_1$ with respect to any finite prime of the ring of integers of $F_1$" should be deleted. In 6.6.6, ``group" $\rightsquigarrow$ ``group of classical Lie type". In 6.7.1, ``If there is a prime $l$ such that $G$ is ramified over $\QQ_l$ (for instance if $F$ or $E(G,X)$ is different from $\QQ$)" $\rightsquigarrow$ ``If there is a prime $l$ which mod $4$ is 2 or 3 or if there are two distinct primes $l$ such that $G_{\QQ_l}$ is unramified,". In 6.8.2, ``condition" $\rightsquigarrow$ ``conditions".
\smallskip
In 6.8.6, ``full prove" $\rightsquigarrow$ ``complete proof", ``$E(G^{\rm ad},X^{\rm ad})=\QQ$" $\rightsquigarrow$ ``$E(G^{\rm ad},X^{\rm ad})$ is a totally real number field", ``$E(G,X)=\QQ$" $\rightsquigarrow$ ``$E(G,X)=E(G^{\rm ad},X^{\rm ad})$, ``It is an easy" $\rightsquigarrow$ ``If $E(G,X)=\QQ$ it is an easy", ``of $\QQ$" $\rightsquigarrow$ ``of a totally real number field", ``the situation can not be reduced" $\rightsquigarrow$ ``the situation still gets reduced", ``moreover the ideas" $\rightsquigarrow$ ``On the other hand, the ideas" and ``this second case" $\rightsquigarrow$ ``these two cases".
\medskip
{\bf AE.1.} In [Va2, 4.3.11], as we worked with the natural trace form on the Lie algebra ${\got g}{\got s}{\got p}$ of a group of symplectic similitudes, the condition $p$ does not divide half of the rank $r_L$ of $L$ needs to be added in order to have it perfect (i.e. in order to have the fourth paragraph of [Va2, p. 469] applying). If $p$ is odd and divides $r_L$, then, provided we work with $Sp$'s groups instead of $GSp$'s groups, it still applies entirely, cf. [Va2, 3.1.6]. So, its only application (for $p\ge 3$), i.e. [Va2, 5.6.3], does not need to be modified. The Claim of [Va2, 4.3.11] is still true, under a very slight (compulsory) restriction, even if $p=2$ and $4| r_L$, see a) and c) of 2.3.18 B3.
\medskip
{\bf AE.2.} The argument of the first sentence of [Va2, 6.4.6 5)] ought to be enlarged. As in [Va2, 3.2.7 11)], we consider two injective maps $(G_1,X_1,H_1)\hookrightarrow (G,X,H)$ and $(G_2,X_2,H_2)\hookrightarrow (G,X,H)$, such that all simple factors of $(G_1^{\rm ad},X_1^{\rm ad})$ and $(G_2^{\rm ad},X_2^{\rm ad})$ are of preabelian and respectively of special type and we have a natural identification $G^{\rm ad}=G_1^{\rm ad}\times G_2^{\rm ad}$. [Va2, 6.4.1] (resp. [Va2, 3.4.1]) points out that $(G_1,X_1,H_1)$ (resp. that $(G_2,X_2,H_2)$) has an integral canonical model $\Mn_1$ (resp. has a quasi-projective normal integral model $\Mn_2$ having the EP). As in [Va2, 3.2.16] we get that $(G_1\times G_2,X_1\times X_2,H_1\times H_2)$ has a quasi-projective normal integral model $\Mn_{12}$ having the EP. Using the fact that the intersection $G_1^{\rm der}\cap G_2^{\rm der}$ (taken inside $G^{\rm der}$) is a finite group scheme of order relatively prime to $p$, from [Va2, 6.2.3] and [Va2, A) to E) of 6.2.2 b)] we get that $(G,X,H)$ has a normal integral model $\Mn$ over the normalization of $\ZZ_{(p)}$ in $E(G_2\times G_1,X_2\times X_1)$, which is a quotient of $\Mn_{12}$ through a free (see [Va2, 3.4.5.1]) action. From [Va2, 3.2.12] we get: if $(G,X,H)$ has an integral canonical model, then $\Mn$ is smooth; so also $\Mn_{12}$ and $\Mn_2$ are smooth. So we can replace $(G,X,H)$ by $(G_1\times G_2,X_1\times X_2,H_1\times H_2)$. So [Va2, 6.4.6 5)] follows from [Va2, 6.2.2-3 and 6.4.1] or from [Va2, 6.2.4 and 6.4.1]. 
\medskip
{\bf AE.3.} The condition of [Va2, 4.3.6 3)] on the existence of ${\rm Lie}$ is not enough; the reason is: it misses data needed to relate ${H^\prime_1}_{R/I_1+I_2}$ with ${H^\prime_2}_{R/I_1+I_2}$. Loc. cit. was thought as taking aside part of the argument of [Va2, 4.3.7 4)]; so in loc. cit. we had in mind the context of [Va2, 4.3.7 4)] (see the Proposition below). It turns out that [Va2, 4.3.7 4) and 5)] as well needs some mild assumptions (they can be read out from the Proposition below and so they will not be restated explicitly). 
\smallskip
The impact of this is (cf. what follows below): the sentence (between parantheses) containing ``it is instructive not to do so" (used before [Va2, Case a) of the proof of 4.3.10]) should be deleted (however, i) below and combinations of i) and iii) below allow its usage under some restrictions). Warning: this deletion does not necessarily need to be performed in the context of the similar sentence (used before [Va2, Part 1) of the proof of 4.3.10]) or in connection to [Va2, 4.3.10.2], cf. iii) below. However, the proof of Proposition below is very much the same as parts of the proof of [Va2, 4.3.10 b)] and so (cf. also AE3.1 3) below) one might still prefer to perform such a deletion.
So what follows is just of group theoretical interest; the reader interested just in applications to Shimura varieties can skip to AE.4.
\smallskip
Coming back to [Va2, 4.3.6 3)], its conclusion holds if the $O$-algebras $R/(I_{j}\cap_{i\in SI} I_i)$, for any $j\in I$ and for every non-empty subset $SI$ of $I\setminus\{j\}$, are still faithfully flat, or if $O$ is of equal characteristic $0$. In both situations we get that ${H^\prime_1}_{R/I_1+I_2}$ and ${H^\prime_2}_{R/I_1+I_2}$ coincide (for the second one cf. also the part below referring to ii)). Let $p$ be an arbitrary rational prime. We start with a definition.
\medskip
{\bf Definition.} A faithful representation $\rho_T:T\hookrightarrow GL(V)$ of a split torus over an arbitrary field $\tilde k$ is called $p$ Lie recoverable, if $T$ is generated by a finite family $(T_j)_{j\in J}$ of subtori of it such that for any $j\in J$, the representation of $T_j$ on $V$ does not involve two distinct characters of $T_j$ whose difference inside the group of characters of $T_j$ (viewed additively) is divisible by $p$. If $J$ has just $1$ element, then we also say $\rho_T$ is elementary $p$ Lie recoverable.
\medskip
We situate ourselves in the context of the first paragraph of [Va2, 4.3.6 3)]. We assume that there is a noetherian, local, integral $O$-subalgebra $R_0$ of $R$ such that $R$ is the strict henselization of $R_0$, that $M$ is obtained by extension of scalars from a free $R_0$-module $M_0$ and that $H^\prime$ is obtained by pull back from a reductive subgroup $H^\prime_0$ of $GL(M_0)_{R_0[{1\over \pi}]}$. We also assume the existence of a projector $\Pi$ of ${\rm End}(M)$ on ${\rm Lie}$ fixed by $H^\prime$. So we have a direct sum decomposition
${\rm End}(M)={\rm Lie}\oplus {\rm Ker}(\Pi)$ preserved by the action of $H^\prime$ on ${\rm End}(M[{1\over {\pi}}])$. Let $l$ (resp. $l_R$) be the residue field of $O$ (resp. of $R$). Let $\tilde k$ be the algebraic closure of the field of fractions $\tilde k_0$ of $R_0$. 
\smallskip
If $H_0^\prime$ is non-trivial and semisimple, then $s({H^\prime}_{0\tilde k},M_0\otimes_{R_0} \tilde k)\in\NN$ is defined by: it is the maximal dimension of indecomposable $ST$-submodules of $M_0\otimes_{R_0} \tilde k$, with $ST$-running through all standard semisimple subgroups of $H^\prime_{0\tilde k}$ whose adjoints are $PSL_2$'s groups. Here ``standard" is used in the same manner as in [Va2, 4.3.3] but at the level of subgroups and not of Lie subalgebras of ${\rm Lie}({H^\prime}_{0\tilde k})$. If $O$ is of characteristic $0$, then $s({H^\prime}_{0\tilde k},M_0\otimes_{R_0} \tilde k)=s({\rm Lie}({H^\prime}_{0\tilde k}),M_0\otimes_{R_0} \tilde k)$ (cf. Weyl's complete reductibility theorem and the standard connections between subrepresentations of ${\rm Lie}({H^\prime}_{0\tilde k})$ and of $H^\prime_{0\tilde k}$). In [Va2, 4.3.10.2], by $s({\got g}_0,W)$ in case ${\got g}_0$ is over a field of positive characteristic, we meant $s(G_0,W)$.
\medskip
{\bf Proposition.} {\it The Zariski closure of $H_0^\prime$ in $GL(M_0)$ is a reductive group if moreover any one of the following three conditions hold:
\medskip
i) $l$ is of characteristic $p$, $H^\prime$ is a torus and the representation $\rho_{\tilde k}$ of ${H^\prime}_{0\tilde k}$ on $M_0\otimes_{R_0} \tilde k$ is $p$ Lie recoverable;
\smallskip
ii) $l$ is of characteristic $0$;
\smallskip
iii) $H_0^\prime$ is non-trivial and semisimple, the characteristic of $l$ is an odd prime $p\ge s({H^\prime}_{0\tilde k},M_0\otimes_{R_0} \tilde k)$ and moreover the semisimplification $AUT^{\rm ss}$ of the connected component $AUT^0$ of the origin of $Aut({\rm Lie}\otimes_R l_R)_{\rm red}$ is smooth of dimension equal to the relative dimension of $H^\prime_0$.} 
\medskip
{\bf Proof:} It is enough to show $H_R^\prime$ is a reductive subgroup of $GL(M)$. As in [Va2, 4.3.6 3)] we can concentrate most of the time just on $1$, $2\in I$. 
\smallskip
The first step is to check that the pull backs of ${H^\prime_1}_{R/I_1+I_2}$ and ${H^\prime_2}_{R/I_1+I_2}$ to ${{\rm Spec}(R/I_1+I_2)}_{\rm red}$ coincide. For ii) this is a consequence of [Bo2, 7.1]. For i) this can be checked as follows. The generic fibre $GF$ of $H^\prime_0$ splits over a Galois field extension of $\tilde k_0$ contained in the residue field of any point of ${\rm Spec}(R/I_i)$ of codimension $0$, $\forall i\in I$. So we can assume $GF$ is a split torus. So, based on arguments similar to the ones of [Va2, 3.1.6] pertaining to $T_R^2$ of loc. cit., we can assume $\rho_{\tilde k}$ is elementary $p$ Lie recoverable; but this case follows from the very definitions. In other words: we can recover the maximal direct sum decomposition of $M\otimes_R R/I_i$ normalized by $H^\prime_i$ from the representation of ${\rm Lie}\otimes_R R/I_i$ on $M\otimes_R R_i$; moreover, the character of the representation of $H^\prime_i$ on any fixed member of it can be read out from the representation of $GF$ on $M_0\otimes_{R_0} \tilde k_0$ ($i=\overline{1,2}$).
\smallskip
For the time being we assume that the pull backs of ${H^\prime_1}_{R/I_1+I_2}$ and ${H^\prime_2}_{R/I_1+I_2}$ to ${{\rm Spec}(R/I_1+I_2)}_{\rm red}$ coincide even for iii); so let ${H^\prime_{12}}_{\rm red}$ be the reductive subgroup of $GL(M\otimes_R {(R/I_1+I_2)}_{\rm red})$ with which these two pull backs coincide. Let $R_{12}:=R/I_1+I_2$ and let $N$ be its ideal of nilpotent elements. Let $n\in\NN$ be such that $N^n=\{0\}$. We consider a maximal tori $T_i$ of ${H^\prime_i}_{R_{12}}$, $i=\overline{1,2}$. We can assume that $T_1$ and $T_2$ lift the same maximal torus $T_{12}$ of ${H^\prime_{12}}_{\rm red}$, cf. [SGA3, Vol. II, 3.6 of p. 48]. Loc. cit. implies as well the existence of $g\in GL(M\otimes_R R_{12})$ which mod $N$ is the identity and $gT_1g^{-1}=T_2$ as tori of $GL(M\otimes_R R_{12})$.
\smallskip
The second step is to show, by induction on $j\in\{1,...,n\}$, that we can assume ${\rm Lie}(T_1)$ and ${\rm Lie}(T_2)$, as $R_{12}$-submodules of ${\rm End}(M\otimes_R R_{12})$, coincide modulo $N^j$. The passage from $j$ to $j+1$ goes as follows. We can assume $g$ mod $N^j$ is the identity. The existence of $g$ implies the existence of $A\in N^j{\rm End}(M\otimes_R R_{12})$ such that under the automorphism $1_{{\rm End}(M\otimes_R R_{12})}+{\rm ad}(A)$ of ${\rm End}(M\otimes_R R_{12})$, ${\rm Lie}(T_1)$ mod $N^{j+1}$ is mapped onto ${\rm Lie}(T_2)$ mod $N^{j+1}$. But due to the existence of $\Pi$, modulo $N^{j+1}$ we can assume $A\in N^j{\rm Lie}\otimes_R R_{12}$. As the Lie algebra of $H_1^\prime$ is ${\rm Lie}\otimes_R R/I_1$, we can assume $g$ mod $N^{j+1}$ is an $R_{12}/N^{j+1}$-valued point of $H_1^\prime$. So, by replacing $T_1$ by an $H_1^\prime(R_{12})$-conjugate of it, we can assume that the restrictions mod $N^{j+1}$ of ${\rm Lie}(T_1)$ and ${\rm Lie}(T_2)$ are the same. So, by induction, we can assume ${\rm Lie}(T_1)={\rm Lie}(T_2)$. 
\smallskip
The third step is to check that ${\rm Lie}(T_1)={\rm Lie}(T_2)$ implies $T_1=T_2$. For i) this is a consequence of Definition: as above, we can recover the maximal direct sum decomposition of $M\otimes_R R_{12}$ normalized by $T_1$ as well as the characters through which $T_1$ acts on its members from ${\rm Lie}(T_1)$ (and GF). The same applies for ii). 
\smallskip
For iii) we have to proceed as in [Va2, p. 465]: 
\medskip
-- if $p>2s({H^\prime}_{0\tilde k},M_0\otimes_{R_0} \tilde k)$, then we can apply the part of the previous paragraph referring to i) (argument: we need to consider as in loc. cit. the standard $1$ dimensional split subtori of $T_1$ in order to get that we are in the context of elementary $p$ Lie recoverable representations of their fibres over $\tilde k$; so the representation of $T_{1\tilde k}$ on $M_0\otimes_{R_0} \tilde k$ is $p$ Lie recoverable); 
\smallskip
-- if $s({H^\prime}_{0\tilde k},M_0\otimes_{R_0} \tilde k)\le p\le 2s({H^\prime}_{0\tilde k},M_0\otimes_{R_0} \tilde k)$, then we have to use as in loc. cit. standard ${\got s}{\got l}_2$ Lie subalgebras of ${\rm Lie}\otimes_R R_{12}$ in order to be able to recover the maximal direct sum decomposition of $M\otimes_R R_{12}$ normalized by $T_1$, from the representations of ${\rm Lie}(T_1)$ and of ${\rm Lie}\otimes_R R_{12}$ on it.
\medskip
What follows next for ii) and iii) is very much the same as [Va2, p. 464 and p. 466]: using exponential maps (for iii) this is allowed as $p\ge s({H^\prime}_{0\tilde k},M_0\otimes_{R_0} \tilde k)$), we can recover from $T_1$ and ${\rm Lie}\otimes_R R_{12}$ the $\GG_a$ subgroups of ${H^\prime_i}_{R_{12}}$ normalized by $T_1$ and so, based on [SGA3, Vol. III, 4.1.2 of p. 172] we get that an open subscheme of ${H^\prime_1}_{R_{12}}$ coincides with an open subscheme of ${H^\prime_2}_{R_{12}}$ and so ${H^\prime_1}_{R_{12}}={H^\prime_2}_{R_{12}}$. The rest of the argument that $H_R^\prime$ is a reductive group is based on amalgamated sums and is as in [Va2, 4.3.6 3)]. This paragraph forms the fourth step. 
\smallskip
So to end the proof, we need to argue the first step for iii). The natural homomorphism from ${H^\prime_i}_{l_R}$ to $AUT^{\rm ss}$ is an isogeny (it is nothing else but the natural central isogeny ${H^\prime_i}_{l_R}\to {H^{\prime{\rm ad}}_i}_{l_R}$; as $p$ is odd this can be checked easily using standard ${\got s}{\got l}_2$ Lie subalgebras of ${\rm Lie}\otimes_R l_R$). We get: ${\rm Lie}(T_1)\otimes_{R/I_1} l_R$ is the Lie algebra of a maximal torus of ${H_2^\prime}_{l_R}$ and so we can assume it is the Lie algebra of ${\rm Lie}(T_2)\otimes_R {l_R}$. As in the above part referring to the subcases corresponding to the fact that $p$ is or is not greater than $2s({H^\prime}_{0\tilde k},M_0\otimes_{R_0} \tilde k)$, we get ${T_1}_{l_R}={T_2}_{l_R}$. As in the previous paragraph, we get: ${H^\prime_1}_{l_R}$ and ${H^\prime_2}_{l_R}$ coincide. But the role of $l_R$ is just to fix notations: the same applies to any point ${\rm Spec}(\tilde l_R)\to {\rm Spec}(R_{12})$, with $\tilde l_R$ an algebraically closed field, lifting a point of the special fibre of ${\rm Spec}(O)$; we just need to mention that ${\rm Lie}\otimes_R \tilde l_R={\rm Lie}\otimes_R l_R\otimes_{l_R} \tilde l_R$ (cf. the uniqueness theorem of [SGA3, Vol. III, p. 313-4] and the fact that $H_0^\prime$ is over $R_0[{1\over {\pi}}]$).
\medskip
{\bf AE.3.1. Remarks.} {\bf 1)} In practice we use i) and iii) intermingled. The last part of the above paragraph points out that the conditions of iii) are such that they can be checked over $O$, using an arbitrary split semisimple group over $O$ whose extension to $R/I_1$ is isomorphic to $H^\prime_1$ (it does not matter which element $1\in I$ we use); so all three parts of the Proposition can be used to get right versions of [Va2, 4.3.7 4) and 5)]. Moreover, it is easy to see (based on [Bo2, 7.1] and on arguments involving non-degenerate trace forms) that the two assumptions of the second paragraph before the Proposition, are not needed in connection to ii). 
\smallskip
{\bf 2)} In [Va2, 4.3.6 3)] we could add: the fact that $H^\prime$ is split can be obtained via [SGA3, Vol. III, 1.5 of p. 329] applied over $R/I_1+I_2$ (using the fact that we are in a local, strictly henselian context) and then making $\pi$ invertible; so the sentence of [Va2, 4.3.6 3)] referring to it should have been placed after the argument that $H^\prime_{12}$ exists (i.e. after showing that we have ${H^\prime_1}_{R_{12}}={H^\prime_2}_{R_{12}}$). 
\smallskip
{\bf 3)} The part of iii) referring to $AUT^0$ automatically holds if the Killing form on ${\rm Lie}(H_1)$ is perfect: in such a case $AUT^0$ is semisimple (cf. [Va2, Part 1 of p. 463-4]).
\medskip
{\bf AE.4.} Let $p$ be a rational prime. A great part of the proof of [Va2, Lemma of p. 496] is wrong and in fact there are contraexamples to it (for instance, with $H_O$ of $A_{p-1}$ Lie type). There are two mistakes in its proof. First, in the case $p=2$ and $\Mg=\mu_2$, we need to assume that the special fibre of $\Mg$ is a subscheme of $P_{O/pO}$; as mentioned in loc. cit. this case is trivial: we can assume $O$ is complete and so [SGA3, Vol. II, 3.6 of p. 48] implies $\Mg$ is contained in $gP_Og^{-1}$, with $g\in H_O(O)$ congruent to the identity mod $2$. Second, the last part of it referring to an action of $\mu_p$ on ${\rm Spec}(\Mo)_{O_2}$ makes no sense. We did not have time to rethink the whole issue (i.e. to list all cases when it holds), so here we just mention that [Va2, 3.4.5.4 3)] implies that [Va2, Lemma of p. 496] holds at least if $p-1$ is greater than the relative dimension of each factor of $H_O/P_O$ defined by a simple factor of $H_O$ (here we can assume $O$ is strictly henselian).
\smallskip
However, [Va2, 6.2.2.1] still holds in different forms (some weaker and some slightly stronger). We take a longer root, in order to fully exploit what can be done about it. We start with an arbitrary Shimura pair $(G,X)$. The ``root of the problem" is in [Va2, 2.11]. In loc. cit. we defined subgroups of $G(\AA_f)$ smooth for $(G,X)$: such a subgroup $H$ is said to be smooth for $(G,X)$ if it is compact and ${\rm Sh}(G,X)$ is a pro-\'etale cover of ${\rm Sh}_H(G,X)$. Warning: in order to avoid always taking quotients (see [De2, 2.1.12]), we allow $H$ to have elements which act trivially on ${\rm Sh}(G,X)$. For instance any compact, neat subgroup of $G(\AA_f)$ is smooth for $(G,X)$. However, smooth subgroups are far more general than compact, neat subgroups. Here we mention just two examples.
\medskip
{\bf Example 1.} If $G$ is a torus, then any compact subgroup of $G(\AA_f)$ is smooth for $(G,X)$ (cf. [Va2, 3.2.8]).
\medskip
{\bf Example 2.} We assume $G$ is a $\QQ$--simple, adjoint group. We assume $G_{\QQ_p}$ has two non-trivial factors $G_1$ and $G_2$ and there is a non-trivial, finite subgroup $H$ of $G_1(\QQ_p)$. Then $H$ is smooth for $(G,X)$. This is a consequence of the structure of the set ${\rm Sh}(G,X)(\CC)$ (see [De2, 2.1.1]): as $G(\QQ)$ has trivial intersection with any $G(\AA_f)$-conjugate of $H$, any $g\in H$ acts freely on this set and so on ${\rm Sh}(G,X)$. We get: there are compact, open subgroup of $G(\AA_f)$ which contain $H$, are smooth for $(G,X)$ but are not neat.
\medskip
We propose some extra classes of compact subgroups of $G(\AA_f^p)$ which are ``between" the neat ones and the ones smooth for $(G,X)$. Let $H$ be a compact subgroup of $G(\AA_f)$.
\medskip
{\bf Definition 1.} {\bf a)} $H$ is called S-smooth for $(G,X)$ (here S stands for strongly) if it is smooth and for each connected component $\Mc$ of ${\rm Sh}_H(G,X)_{\CC}$ there is a compact, open subgroup $H_1$ of $G(\AA_f^p)$ containing the subgroup $H(\Mc)$ of $H$ leaving invariant $\Mc$ and such that the set of complex points of the connected component of ${\rm Sh}_{H_1}(G,X)_{\CC}$ dominated naturally by $\Mc$ can be identified (see [De2, p. 266-7]) with $\Sigma\setminus X^0$, with $X^0$ a connected component of $X$ and with $\Sigma$ an arithmetic subgroup of $G^{\rm ad}(\QQ)$ which does not have $2$-torsion. 
\smallskip
{\bf b)} We say $H$ is $p$-smooth if it is smooth and for each connected component $\Mc$ of ${\rm Sh}_H(G,X)_{\CC}$ there is a rational prime $l(\Mc)$ different from $p$ such that the image of any  ${\rm pro}$-$p$ subgroup $H_p(\Mc)$ of $H(\Mc)$ (with $H(\Mc)$ as in a)) in $G^{\rm ad}(\QQ_{l(\Mc)})$ is trivial. If $\Ms$ is a set of rational primes having at least two elements, we say $H$ is $\Ms$-smooth, if it $p$-smooth for any $p\in\Ms$.
\medskip
Let $v$ be a prime of $E(G,X)$ dividing $p$. Let $O$ be a faithfully flat, DVR $O_{(v)}$-algebra. We consider a smooth (resp. normal) integral model $\Mn$ of ${\rm Sh}_H(G,X)$ over $O$. b) suggests the following variant of [Va2, def. 3.4.8].
\medskip
{\bf Definition 2.} {\bf a)} $\Mn$ is said to be quasi-strongly smooth (resp. quasi-strongly normal) if for any compact, open subgroup $H_0$ of $G(\AA_f^p)$ such the subgroup $H_0\times H$ of $G(\AA_f)$ is $p$-smooth for $(G,X)$, $\Mn$ is a pro-\'etale cover of $\Mn/H_0$. 
\smallskip
{\bf b)} If $p=2$ we say $\Mn$ is $A$-strongly smooth (resp. $A$-strongly normal), if $\Mn$ is a pro-\'etale cover of $\Mn/H_0$, for any subgroup $H_0$ of $G(\AA_f^p)$ such the subgroup $H_0\times H$ of $G(\AA_f)$ is S-smooth for $(G,X)$ (here $A$ stands for almost).
\medskip
Referring to [Va2, Lemma of p. 496], there is $n\in\NN$ easy computable in terms of $H_O$, such that $H_O(O)$ has no elements of order $p^n$. This is an immediate consequence of the following two facts, applied in the context of any faithful representation of $H_O$ (for instance, the adjoint one):
\medskip
-- $H_{\FF}(\FF)$ has no semisimple elements of order $p$ and all its unipotent elements have an order dividing a fixed power of $p$;
\smallskip
-- [Va2, 3.4.5.2] and its variant for $p=2$ (if $p=2$, any $O$-endomorphism of a free $O$-module of finite rank of order a power of $2$ and which mod $4$ is congruent to the identity, is the trivial endomorphism). 
\medskip
We now start correcting the parts of [Va2] relying on [Va2, Lemma of p. 496]. In what follows, $(G,X,H,v)$, $(G_1,X_1,H_1,v)$, $(G_2,X_2,H_2)$, $X^0$, $\Mm_1$, $\Mm_2$, $p>2$, $L_p$, $(W,\psi)$, $V_0$ and $\Mc^0=\Mc_2$ are as in [Va2, 6.2.2 b) and its proof]. To conform with the notations of this paper, we denote $L_p$ by $L_{(p)}$. Let $G_{2\ZZ_{(p)}}$ be the reductive group over $\ZZ_{(p)}$ obtained by taking the Zariski closure of $G_2$ in $GL(L_{(p)})$. We take $O:=V_0$ and $H_O:=G_{2V_0}^{\rm ad}$. We consider $g\in {\rm Aut}(G_2,X_2,H_2)$; we do not assume it is inner (i.e. is defined by an element of $G^{\rm ad}_{2\ZZ_{(p)}}(\ZZ_{(p)})$). Moreover, $H_2^p$ is a compact, open subgroup of $G_2(\AA_f^p)$ normalized by $g$ and such that the morphism $\Mm_2\to\Mm_2/H_2^p$ is a pro-\'etale cover. On [Va2, p. 494-7], ``$\Mm_2/H_2^p\times H_2$" $\rightsquigarrow$ ``$\Mm_2/H_2^p$" and ``$\Mm_{2V_0}/H_2^p\times H_2$" $\rightsquigarrow$ ``$\Mm_{2V_0}/H_2^p$". The $\ZZ$-lattice $L$ of [Va2, p. 494] is as in [Va2, 5.1.2].
\smallskip
We assume $g$ is such that the automorphism $a_g$ of $\Mn_{2V_0}/H_2^p$ it defines naturally leaves invariant the image $\Mc(H_2^p)$ in $\Mn_{2V_0}/H_2^p$ of the (chosen) connected component $\Mc_2$ of $\Mn_{2V_0}$. Warning: this does not necessarily imply that $g$ leaves invariant $X^0$. We also assume $a_g$ acts freely on the generic fibre of $\Mc(H_2^p)$. We have four different possibilities (which can be combined) in connection to [Va2, 6.2.2.1]:
\medskip
\item{\bf a1)} we require $X_2=X_2^{\rm ad}$ and that in fact $a_g$ acts freely on the generic fibre of $\Mn_{2V_0}/H_2^p$; 
\smallskip
\item{\bf a2)} we require the order $o(a_g)$ of $a_g$ to be divisible by $p^{n+n_1}$, where $p^{n_1}$ is the greatest power of $p$ dividing the order of the image of $g$ into the group of outer automorphisms of the product of all simple, non-compact factors of $G^{\rm ad}_{\RR}$;  
\smallskip
\item{\bf a3)} we require that $g$ over $\RR$ normalizes $X^0$ and each simple factor of $G^{\rm ad}_{\RR}$;
\smallskip
\item{\bf a4)} $g$ is defined by an element of $G_{2\ZZ_{(p)}}^{\rm ad}(\ZZ_{(p)})$ which is contained in a compact subgroup $H_2^{p{\rm ad}}(g)$ of $G_2^{\rm ad}(\AA_f^p)$ containing the image of $H_2^p$ in $G_2^{\rm ad}(\AA_f^p)$ and having the property that $H_2^{p{\rm ad}}(g)\times H_2^{\rm ad}$ is $p$-smooth for $(G_2^{\rm ad},X_2^{\rm ad})$. 
\medskip
If we work under a3), then there is $m\in\NN$ relatively prime to $p$ such that $g^m$ is an inner automorphism of $G^{\rm ad}_2$. This is a consequence of the structure of the group of outer automorphisms of an absolutely simple group of classical Lie type and (in case $p=3$ and $G^{\rm ad}$ has simple factors of $D_4$ Lie type) of the fact that $g$ leaves invariant $X^0$. 
\smallskip
If $g$ is inner and we work under a1), then either $g$ is not an element of finite order or $a_g$ does not fix any $\FF$-valued point of $\Mc(H_2^p)$. Argument: if it is of finite order, then we can assume (cf. [Va2, 3.4.5.1]) that this order is a power of $p$ and so $g$ over $\RR$ (as $p\ge 3$ there is no invertible element of $\RR$ of order $p$) is contained in a maximal compact subgroup of $G_2^{\rm ad}(\RR)$; so it fixes an element of $X_2=X_2^{\rm ad}$ (cf. also [De2, 1.2.8]) and so $a_g$ does not act freely on the generic fibre of $\Mn_{2V_0}/H_2^p$. Contradiction. 
\medskip
{\bf Proposition.} {\it If we work under a2), then $a_g$ acts freely on $\Mc(H_2^p)$. In particular, this applies if $g$ is inner and $G^{\rm ad}_{V_0}(V_0)$ has no element of order $p$.}  
\medskip
{\bf Proof:} If $a_g$ fixes an $\FF$-valued point $y$ of $\Mc(H_2^p)$, then we can define $g_0$ as in [Va2, p. 495]; its order is the same as the order of $g$. As in loc. cit. we can assume $o(a_g)$ is a power of $p$ and that $g_0^{p^{n_1}}$ is inner. As $p^n$ divides $o(g_0^{p^{n_1}})$ we reach a contradiction to the definition of $n$. 
\medskip
{\bf Remarks.} For another form of this Proposition, see 4.4.13.5; see \S13 for extra improvements of the Proposition (similar to the one of 4.4.13.5) obtained based on the extra fact that in [Va2, p. 495], $g_0$ was defining an automorphism of a Shimura adjoint $\sg$-crystal over $\FF$ (see def. 2.2.11 4) and end of 2.2.22 2)). This Proposition does not suffice in applications to [Va2, 6.2.2 b)]; so to considerable strengthen it, we need to bring into the picture some extra tools.
\medskip
{\bf AE.4.1.} We can assume $G_1$ is adjoint (the action of a subgroup of a group acting freely is free); so $(G_1,X_1)=(G_2^{\rm ad},X_2^{\rm ad})$. Based on [Va2, 6.2.3 and 6.4.2] we can assume $G_1$ is $\QQ$--simple. The only cases of [Va2, 6.2.2 b)] not covered by [Va2, 6.2.2 a)] are those in which $p$ is odd and $G_2^{\rm ad}$ has factors of $A_{pm-1}$ Lie type. So based on [Va2, 6.5.1.1 and 6.6.5.1], we can assume we are in the context of a standard PEL situation $(f_2,L_{(p)},v_2,\Mb_2)$ (as defined in 2.3.5, cf. also 2.3.5.1), where $v_2$ is a prime of $E(G_2,X_2)$ dividing the prime $v^{\rm ad}$ of $E(G^{\rm ad},X^{\rm ad})=E(G_2^{\rm ad},X_2^{\rm ad})$ divided by $v$. 
\smallskip
We consider an injective map $(G_2,X_2)\hookrightarrow (G_2^\prime,X_2^\prime)$, with $G_2^\prime$ as the subgroup of $GL(W)$ generated by $G_2$ and by the center of the centralizer of $G_2$ in $GL(W)$. The Zariski closure of $G_2^\prime$ in $GL(L_{(p)})$ is a reductive group $G^\prime_{2\ZZ_{(p)}}$ and moreover we get a cover 
$$q_2^\prime:(G_2^\prime,X_2^\prime)\to (G_2^{\rm ad},X_2^{\rm ad});$$ 
the first part can be seen immediately inside $GL(L_{(p)}\otimes_{\ZZ_{(p)}} V_0)$, starting from the classification provided by [Ko2, top of p. 375 and p. 395], while the second part is a consequence of the fact that the center $Z(G_2^\prime)$ of $G_2^\prime$ is the group scheme defined by invertible elements of an \'etale $\QQ$--algebra $AL$ of endomorphisms of $W$. $\Mm_2$ is an open closed subscheme of the integral canonical model $\Mm_2^\prime$ of the Shimura quadruple $(G_2^\prime,X_2^\prime,G^\prime_{2\ZZ_{(p)}}(\ZZ_p),v_2)$ (cf. [Va2, 6.2.3 and 3.2.15]).
\smallskip
So, based on [Mi1, 4.13], we can assume that a connected component of $\Mm_{1V_0}$ is the quotient of $\Mc^0=\Mc_2$ (see [Va2, p. 494]) by a group of automorphisms $GA$ of $\Mc_2$ which are defined by right translation by elements of a subgroup of the group $GR$ of $\AA_f^p$-valued points of the center $Z(G_2^\prime)$ of $G_2^\prime$; based on [De2, 2.1.12], in fact we can replace $GR$ by the familiar group (of the class field theory)
$$Z(G_2^\prime)(\AA_f^p)/Z(G_2^\prime)(L_{(p)}),$$
cf. also [Va2, 3.3.1]. We can assume $AL\subset\Mb_2\otimes_{\ZZ_{(p)}}\QQ$.
\smallskip 
We consider an element $h\in Z(G_2^\prime)$ defining an element of $GA$. We assume it fixes a point $y\in\Mc_2(\FF)$. $y$ gives birth to a quadruple 
$$Q_y=(A_y,p_{A_y},\Mb_2,(k_N)_{N\in\NN,\,(N,p)=1}),$$ 
where $(A_y,p_{A_y})$ is a principally polarized abelian variety over $\FF$ of dimension ${1\over 2}\dim_{\QQ}(W)$, endowed with a family of $\ZZ_{(p)}$-endomorphisms (still denoted by $\Mb_2$) and having (in a compatible way; see [Va2, top of p. 455]) level-$N$ symplectic similitude structure $k_N$, $N\in\NN$ with $(N,p)=1$, such that some axioms are satisfied, cf. the standard interpretation of $\Mn_2$ as a moduli scheme (to be compared also with [Ko2, ch. 5]). We have a similar modular interpretation for $\FF$-valued and $V_0$-valued points of $\Mn_2^\prime$, provided we work in a $\ZZ_{(p)}$-context; in such a context we speak about principally $\ZZ_{(p)}$-polarized abelian varieties, $\ZZ_{(p)}$-isomorphisms (i.e. isomorphisms up to $\ZZ_{(p)}$-isogenies as defined in 2.1), etc. So the translation of $y$ by $h$ gives birth to a similar quadruple 
$Q_y^\prime=(A_y^\prime,p_{A^\prime_y},\Mb_2,(k_N^\prime)_{N\in\NN,\,(N,p)=1}),$
with $(A_y^\prime,p_{A^\prime_y})$ a principally $\ZZ_{(p)}$-polarized abelian variety over $\FF$ which is $\ZZ_{(p)}$-isogeneous to $(A_y,p_{A_y})$. Due to this $\ZZ_{(p)}$-isogeny, we can identify $H^1_{\rm crys}(A_y^\prime/V_0)$ with $M:=H^1_{\rm crys}(A_y/V_0)$.
\smallskip
The fact that $h$ fixes $y$ means that these quadruples are isomorphic, under a $\ZZ_{(p)}$-isomorphism $a:A_y\tilde\to A_y^\prime$. The automorphism $a_M$ of $M$ we get (via $a$ and the mentioned identification), as an element of ${\rm End}(M)$, belongs to the $\QQ$--vector space generated by the crystalline realizations of the $\QQ$--endomorphisms of $A_y$ naturally defined by elements of ${\rm Lie}(Z(G_2^\prime))$: this can be read out from the \'etale context with $\QQ_l$-coefficients, where $l$ is an arbitrary prime different from $p$. So $a_M$ leaves invariant any $F^1$-filtration of $M$ defined by a $V_0$-valued point lift $z$ of $\Mc_2$ lifting $y$. Such a lift is determined by the mentioned filtration, cf. the deformation theory (see [Me, ch. 4-5]) of polarized abelian varieties endowed with endomorphisms; see also [Va2, 5.6.4]. So from the modular interpretation of $\Mn_2$ we get: $h$ fixes all these lifts. So $h$ acts trivially on $\Mc_2$. We conclude: $\Mc_2$ is a pro-\'etale cover of $\Mc_2/GA$. This completes the corrected proof of [Va2, 6.2.2 b)].
\medskip
{\bf AE.4.1.1.} The use of $\ZZ_{(p)}$-isogenies in AE.4.1 can be entirely avoided. This goes as follows. Let $G_3:=G_2^{\rm ad}\times\GG_m$; identifying $\GG_m$ with the quotient of $G_2$ by its subgroup $G_2^0$ acting trivially on $\psi$, we get naturally an epimorphism $q_2:G_2\twoheadrightarrow G_3$. Let $(G_3,X_3)$ be the Shimura pair such that $q_2$ defines a finite map 
$$q_2:(G_2,X_2)\to (G_3,X_3).$$  
\indent
{\bf Fact.} {\it Let $A:={\rm Ker}(q_2)$. For any field $l$ of characteristic $0$, the group $H^1({\rm Gal}(l),A(\bar l))$ is a $2$-torsion group (which in general --like for $l=\RR$-- is non-trivial).}
\medskip
{\bf Proof:} From the structure of $G_2^0$ (for instance, see [Ko2, ch. 5 and 7]) we get that $A$ is a product of Weil restriction of scalars from some totally real number fields to $\QQ$ of $1$ dimensional tori which over $\RR$ are compact (if $\Mb\otimes_{\ZZ_{(p)}} \QQ$ is a simple $\QQ$--algebra, then the mentioned product has only one factor). The Fact follows.
\medskip
So $q_2$ is not a cover but from the point of view of free actions (see [Va2, 3.4.5.1]) it is ``close enough". In other words, the image of $\Mc_{2\CC}$ in the quotient of ${\rm Sh}_{H_2}(G_2,X_2)_{\CC}$ by $A(\AA_f^p)$ is a (potentially infinite) Galois cover of the image of $\Mc_{2\CC}$ in ${\rm Sh}_{H_3}(G_3,X_3)_{\CC}$ whose Galois group is a $2$-torsion group, cf. also [Va2, 3.3.1]; here $H_3:=H_2^{\rm ad}\times\GG_m(\ZZ_p)$ and $\Mc_{2\CC}$ is obtained via extension of scalars through an arbitrary $O_{(v_2)}$-monomorphism $V_0\hookrightarrow\CC$. So, referring to AE.4.1, eventually by replacing $h$ by $h^2$, we can assume $h\in G_2(\AA_f^p)$. So AE.4.1 can be performed ``in terms" of $q_2$ and not of $q_2^\prime$. 
\smallskip
However, the context of $q_2^\prime$ is more convenient for generalizations (see the whole of 2.3.5.5-6) and moreover it can be adapted (see the proof of 4.14.3.2.1) for $p=2$.  
\medskip
{\bf AE.4.2.} [Va2, 6.2.2.1] was used in [Va2, 6.2.4.1] and so implicitly in [Va2, 6.4.4]. For the meaning of ${\rm mid}(X)$ and of the torsion number $t(G^{\rm ad})$ we refer to 2.1. We have:
\medskip
{\bf Exercise.} Let $(\tilde G,\tilde X)$ be an adjoint Shimura pair, with all simple factors of $\tilde G$ of classical Lie type. We assume $\tilde G_{\QQ_p}$ is unramified. Then the inequality ${\rm mid}(X)+3\le p$ implies: $p$ does not divide $t(G^{\rm ad})$. If ${\rm mid}(X)+2=p$, then it can happen that $p| t(G^{\rm ad})$. Hints: for the first part use the classification of [He, p. 518] and the fact that a $GL_n$-group over the Witt ring $WR$ of a perfect field of characteristic $p$, with $n+1<p$, can not have $WR$-valued points of order $p$; for the second part, use $\QQ$--forms of $SU(1,p-2)^{\rm ad}_{\RR}$ or of $SO(2,p-2)_{\RR}$.
\medskip
So, in connection to the part of the first paragraph of AE.4 referring to [Va2, 3.5.4 3)], it is shorter and more convenient to work with $t(G^{\rm ad})$ rather than ${\rm mid}(X)$; accordingly, in all that follows we state our results in terms of $t(G^{\rm ad})$'s. However, it is useful to keep in mind that for Shimura varieties whose adjoint varieties are products of simple, adjoint varieties which are either of $B_{{p-1}\over 2}$ type or of $A_{p-2}$ type involving $\QQ$--simple groups which over $\RR$ have only simple factors of the form $SU(a,p-1-a)_{\RR}$, with $a\in \{0,1\}$, one can do a little bit better (with respect to the odd prime $p$). We have the following correction of [Va2, 6.2.4.1 and 6.4.4]:
\medskip
{\bf  Correction.} {\it {\bf a)} In [Va2, 6.2.4.1] we need to either assume $p\not | t(G^{\rm ad})$ or to replace strongly smooth by quasi-strongly smooth.
\smallskip
{\bf b)} In [Va2, 6.4.4] we need to assume $H_{\Ms}\times H^{\Ms}$ is $\Ms_1$-smooth for $(G,X)$, where $\Ms_1$ is the set of odd rational primes belonging to $\Ms$ and dividing $t(G^{\rm ad})$.}
\medskip
b) is a consequence of a), cf. the proof of [Va2, 6.4.4]. The proof of a) is very much the same as AE.4.1 (cf. also Proposition of AE.4); the only difference: fixing a connected component $\Mc$ of ${\rm Sh}_H(G,X)_{\CC}$, we have ``control" only on one prime $l(\Mc)$ different from $p$ (cf. a) of Definition 1). As the part of AE.4.1 involving $G_2^\prime$ has not been previously used (i.e. known outside the context of standard PEL situations), a) is detailed (as an application) in 2.3.5.8. From a) we get:
\medskip
{\bf Corollary.} {\it We refer to AE.4. If we work under a4), then $a_g$ acts freely on $\Mc(H_2^p)$.}
\medskip
For the sake of completeness we repeat here the philosophy of [Va2, 6.4.6 6)] in a slightly corrected form as suggested by Correction.
\medskip
{\bf Ph.} {\it To generalize Serre's lemma [Mu1, p. 207] to the context of an integral canonical models $\Mn$ of a Shimura quadruple $(G,X,H,v)$ of preabelian type (with $v$ dividing a rational prime $p\ge 5$) and of a compact subgroup $H_0$ of $G(\AA_f^p)$, it is enough to check things in characteristic zero (i.e. to check that ${\rm Sh}(G,X)$ is a pro-\'etale cover of ${\rm Sh}_{H_0\times H}(G,X)$) as well as, in case $p| t(G^{\rm ad})$, that  $H_0\times H$ is in fact $p$-smooth (and not just smooth) for $(G,X)$.}  
\medskip
{\bf AE.5.} For a detailed (slightly corrected) version of [Va2, 5.5.1] we refer to 2.4.
\medskip\smallskip
{\bf AE.6.} In Step 6 of [Va2, 3.1.3.1] we could add that $Y_2^\prime$ is obtained from $\widehat{Y}_2$ as in Step 4 of loc. cit. 
In the proof of [Va2, 3.1.6], we could add first that the kernel of $\rho$ is the extension of a finite, flat group scheme by a torus, before concluding that it is a finite, flat group scheme (i.e. that this torus is trivial). Also we could add that the fastest way to see that $G_R^1$ is a reductive group, is to remark that we have a short exact sequence 
$$0\to G^{0\rm der}_R\to G_R^1\to T_R/C_R\to 0.$$ 
In [Va2, 3.2.1.1 15) and 16) and the proof of 3.2.2.1], we do not need to add quasi-separatedness to the quasi-compactness assumptions, despite the fact that in [EGA IV, \S 8] this is always done. It is [FC, 2.7 of p. 9] which implies that in connection to abelian schemes over normal schemes, quasi-separatedness condition can be dropped. Similarly, as $Y_1$ is normal and geometrically connected over the field $K_V$ of characteristic 0, in [Va2, top of p. 431] we can still appeal to the fundamental exact sequence of [SGA1, p. 253], even if $Y_{1K_V}$ is not quasi-compact.   
\smallskip
It is worth pointing out that in [Va2, proof of 3.2.2.1], the existence of $O_1^\prime$ in the case when $Y$ is not quasi-compact, can be seen as well using Weil divisors, starting from the embedding of $\Ma_{d(A_{U_1}),1,N}$ in the projective scheme $\bar\Ma_{d(A_{U_1}),1,N}$ as of [FC, \S 5 of ch. V] and using [FC, iv) of 10.1 of p. 88]. Before the fourth sentence of [Va2, proof of 3.2.2.3] we could add that we can assume $Z$ is local, and using descent that it is strictly henselian; so the last paragraph of [Va2, proof of 3.2.2.1] can be imitated in the context of $Z$ as well (to bring the things back to the level of $O$-schemes of finite type).
\smallskip
We could add to [Va2, 4.2] that $[{\got g},{\got g}]$ is algebraizable (see [Bo2, 7.6]) and so (as we are in characteristic $0$ and ${\got g}$ is reductive) the restriction of ${\rm Tr}$ to it is non-degenerate.  We could add in the beginning of the proof of the Claim 1 of [Va2, 4.3.10 a)], that ${\got g}(M)$ is free (as $\pi({\got g})$ is enveloped by $M$). In [Va2, 4.3.10.2], we could add that any strictly henselian DVR is a Nagata ring (easy consequence of [Ma, Cor. 2 of p. 234]). Related to [Va2, p. 463], we could add that ${\got g}(M)$ is naturally isomorphic to its dual and so it is irrelevant if we use $R[[{\got g}(M)]]$ or $R[[{\got g}(M)^*]]$.
\smallskip
Related to [Va2, end of 6.2.2 C)], we could add that the forwarding to [Va2, 6.2.3] refers to the logical version of loc. cit. over $\CC$. Related to the unramified part of the proof of [Va2, 6.4.5], we could add the reference to [Va2, 3.4.5.2] (as $p>2$). It is not stated explicitly, but on [Va2, p. 495] the argument that $g_0$ is not an outer automorphism did use the hypothesis $p>2$: there are no outer automorphisms of odd order of a simple, adjoint group $G$ over $\CC$ leaving invariant a parabolic subgroup $P$ of $G$ with the property that $G/P$ is the dual of a hermitian symmetric domain; however if $G/P$ is the dual of the hermitian symmetric domain associated to $SU(n,n)_{\RR}$ or to $SO(2,2n-2)_{\RR}$, then we can get such outer automorphisms of order $2$. In the Fact of [Va2, 6.8.0] we could add that the pro-finite group is as well abelian.
\medskip
{\bf AE.6.1.} We could add to the theory of [Va2, 4.3], that in applications to Shimura varieties of Hodge type, we need to consider only DVR's of mixed characteristic which are complete and have a residue field $k$ which is the algebraic closure of a finite field and we need to check just that condition [Va2, 4.3.5] holds (in a suitable context) for any $O$-algebra $R$ which is local, integral, strictly henselian and has $k$ as its residue field. This is so due to the fact that the rings $Re$ and $Re^1$ of [Va2, 5.2.17.1] have all these properties (cf. [Va2, 5.2.1.1] for the strictly henselian part). Moreover, based on [Va2, 4.3.7 5') and 6)] we can assume $R$ is as well complete and excellent. Warning: the way the proof of [Va2, 4.3.10 b)] is organized (we assumed that $G_0(M)_{R^0}$ is already semisimple), we still need to work in a reduced context (i.e. we still need to replace above integral by reduced). 
\medskip
\vfill\eject
\centerline{}
\vfill\eject
\centerline{}
\bigskip
\bigskip
\centerline{\bigsll {\bf Appendix}}
\bigskip\bigskip
\medskip
It is desirable to make the connection between the present work and [RR]. There are three results here related to [RR]. They were obtained independently and their main part precedes [RR], cf. [Va1] and 1.15.10. These results are: 
\medskip
{\bf A)} The part of 3.1.0 expressed in 3.6.6.1 2). It is related to the main results of [RR, ch. 3-4]. Loc. cit. obtains 3.6.6.1 2) (as well as its variant in a generalized Shimura context) in more general contexts. However our ideas work (without any change) to even more general contexts (see below).
\smallskip
{\bf B)} 3.2.10 (FORM). Its idea is very much the same as of the general result [RR, 1.17].
\smallskip
{\bf C)} 4.5.4 and 4.5.6.1. The idea of 4.5.4 is the same underlying [RR, def. 3.3]. The Proposition 4.5.6.1 is related to the main result 3.6 of [RR]. It can be used to restate this main result in a simpler way, provided we restrict the class of groups. We use the notations of [RR, ch. 2-3]. We have: 
\medskip
{\bf The Lie criterion (preliminary form).} {\it If all factors of $G^{\rm ad}$ are of $B_3$, $C_n$ ($n\in\NN$), $F_4$ or $G_2$ Lie type, then for two points $s_1$, $s_2$ of $S$, we have 
$\bar b(s_1)=\bar b(s_2)$ iff the sequence of Newton polygons of the cyclic factors of the Lie isocrystal over $\bar s_1$ (defined by ${\rm Lie}(G)$ or, in the sense of [RR, 3.3], corresponding to the adjoint representation of $G$ on its Lie algebra) is equal to
the sequence of Newton polygons of the cyclic factors of the Lie isocrystal over $\bar s_2$.}
\medskip
The proof of this is entirely the same as of 4.5.6.1: it is 2.2.24.1 (or [RR, 3.4 (i)]) which allows us to use semisimple elements as in the proof of 4.5.6.1. The ``specialty" of $F_4$ and $G_2$ Lie types (from which the part of the Criterion pertaining to them follows) is: if we order the positive roots w.r.t. a chosen basis of the root system of any of these Lie types, on ``the top part" we have a total ordering. To be quick, we explain what we mean by ``the top part" directly in terms of the positive roots listed in [Bou2, planche VIII and IX]. We denote by $>$ the ordering we referred to. For the $G_2$ Lie type we have:
$$3\al_1+2\al_2>3\al_1+\al_2>2\al_1+\al_2>\al_1+\al_2.$$ 
For the $F_4$ Lie type we have:
$$2\al_1+3\al_2+4\al_3+2\al_4>\al_1+3\al_2+4\al_3+2\al_4>\al_1+2\al_2+4\al_3+2\al_4>\al_1+2\al_2+3\al_3+2\al_4.$$
All the above roots are $>$ than all positive roots which are not mentioned. 
\smallskip
We do not know when this criterion is true without its assumption on the Lie types of the simple factors of $G^{\rm ad}$; 4.5.6.2 and the Fact of the proof of 4.5.6.1 point out that there are simple variants (the easiest ones are for the $B_n$ Lie types; see 4.5.6.2 A) of this criterion in many cases involving factors of $G^{\rm ad}$ of $A_n$, $B_n$ or $D_n$ Lie types ($n\ge 2$). In practice it is easier to work with the above criterion (or variants of it) than with the partial ordering on $B(G)$ defined in [RR, 2.3]. Moreover, Exercise 4.5.6 4) points out that, from the point of view of stratifications, [RR, 3.6] brings nothing new to the specialization theorem. Loc. cit. is just a restatement of the specialization theorem as provided by [RR, 2.2] and by the Newton map as defined in [RR] (in particular it does not represent a new proof of this theorem). [RR, ch. 1, 2.2, 2.4, 3.6-12] contains lots of useful theoretical tools (for instance, the language is practical and the proof of [RR, 4.2] is short). Unfortunately, [RR] is of extremely little use w.r.t. the Real Problem of 1.6.5.
\medskip
As a conclusion: the work of \S 1-4 is integral and so more refined (handling deformation aspects and subtler things mod $p$, like 3.9 and 3.13) and more restricted (we work in a context of good reduction and of generalized Shimura $p$-divisible objects), while [RR] works rationally and so less refined but more general. 
\smallskip
However, the greatest part of the ideas and results presented in 3.6 can be used (i.e. extended) to even more general contexts. What follows (except Theorem 2 below) is an afterthought, suggested by [RR, 4.2]. Though 3.9.9 and 1) of 3.15.7 B point out that a lot can be done in the non-smooth context, for the time being let's refer just to the smooth one. We start with an arbitrary filtered $\sg_k$-$\Ms$-crystal
$$
{\got C}=(\tilde M,(F^i(\tilde M))_{i\in S(\tilde a+1,\tilde b)},\tilde\vph,\tilde G);
$$ 
we do not assume that the $W$-condition holds for it or that $\tilde G$ is a reductive group (see 3.6.1.5). We work with $p\ge 2$. There are four steps in achieving this extension.
\medskip
{\bf Step 1. Part A.} If the $W$-condition holds for ${\got C}$ and $\tilde G$ is a quasi-split reductive group, then 3.2.3, 3.6.6, 3.6.6.0, 3.6.6.2 and the proof of 3.11.1 (relying on 3.4.11) remain valid, with no essential change (besides the logical one in connection to the part of 3.13.7.8 involving $N_j$'s, with $j<0$).
So the obstructions of proving 3.1.0, 3.4-5, 3.6.1.2-3 and 3.11.1 in this context come just from:
\medskip
-- 3.6.8.9 (cf. the discussion in 3.6.1.6),
\smallskip
-- from 3.5.6,
\smallskip
-- and from 3.13.7.9 (cf. its end part postponing the general definition of Lie stable $p$-ranks). 
\medskip
The obstruction coming from 3.5.6 is not a very serious one (as we expect 3.5.6 to be a known result and), as it can be substituted (cf. 3.6.6.0) by the expectation of 3.6.1.6. Moreover, the obstruction coming from end of 3.13.7.9 is just for the time being. So the serious obstruction is in achieving suitable global deformations of ${\got C}$, as predicted by 3.6.1.6. We hope to come back to this in a future paper, as we think it is one of the most urgent thing needed to be ``sorted out" in connection to Fontaine categories. 
\smallskip
However, from the point of Newton polygons, the ``lack of deformations" can be substituted by the specialization theorem (or by loc. cit.). Let $g_0\in G(W(k))$ be such that $\tilde\vph_0:=g_0\tilde\vph$ takes the Lie algebra of a Borel subgroup of $\tilde G$ into itself (this is as in 3.2.3). As in 3.6.6 (cf. also 3.6.6.2) and 3.6.6.0 we get that there is an affine, open, dense $U={\rm Spec}(R)$ of $\tilde G$ through which the origin of $\tilde G$ factors, such that for any $u\in U(W(\bar k))$, $u(\tilde\vph_0\otimes 1)$ and $\tilde\vph_0\otimes 1$ are isomorphic under an isomorphism defined by an element of $\tilde G(W(\bar k))$ which has the analogue shape of the element $h$ of 3.13.7.8. We would like to add just four extra things. 
\medskip
{\bf T1.} First, the fact that we use inner conjugations just by elements of $G(W(k))$ of the mentioned shape, is irrelevant from the point of view of the proof of 3.6.6 (see the end of the paragraph before the Expectation of 3.13.7.8).
\smallskip
{\bf T2.} Second, it is irrelevant what shape the exponential maps have (cf. also end of 3.14 A); in connection to the generalization of 3.6.6 the only important think is that the property IND of the proof  of 3.6.6 still holds. 
\smallskip
{\bf T3.} Third, in the proof of 3.6.6 the particular shape of parabolic subgroups involved played no role at all. 
\smallskip
{\bf T4.} Fourth, the case when we have (see the proof of 3.6.6) $\abs{I_0}=1$ and $\abs{\tilde I_0}=3$ (so we are dealing with the $D_4$ Lie type) is, as mentioned in 3.6.6.2, entirely the same as the case when we have $\abs{I_0}=1$ and $\abs{\tilde I_0}=2$. The explicit checking of this is left as an exercise for the diligent reader. We just add: the case $\abs{I_0}=1$ and $\abs{\tilde I_0}=3$ can be handled easily as well by computers. 
\smallskip
This goes as follows. We need to check that, for $\abs{I_0}=1$ and $\abs{\tilde I_0}=3$, a suitable morphism $m_{U}:\tilde U\to U$ is such that its image has the right dimension. $m_U$ is either constructed as the morphism $m_H$ of the proof of 3.6.6 or (which is even better) as the orbit map of $\TT_H$ of 3.13.7.8 defined by the identity element of $H(k)$. $m_{U}$ is the extension to $k$ of a morphism $m_{U,3}:\tilde U_{\FF_{p^3}}\to U_{\FF_{p^3}}$ over $\FF_{p^3}$, and $\tilde U$ is a connected, smooth group scheme of dimension less than $2(24+4)=56$ (see [Bou, planche IV] for the fact that the dimension of an absolutely simple adjoint group of $D_4$ Lie type is 28); if we construct $m_U$ via $\TT_H$ of 3.13.7.8, then we can replace 56 by 52 as well as below we can replace 55 by 51. Moreover, the number of possibilities for such morphisms $m_U$ is finite (as they depend only on the choice of some parabolic subgroups of an absolutely simple, split, adjoint group over $\FF_{p^3}$ of $D_4$ Lie type) and each such morphism $m_U$ is defined by an orbit of a group action. Let $c(m_U)\in\NN$ be the number of connected components of the maximal reduced subgroup $\tilde U^{\rm st}$ of $\tilde U$ mapped by $m_U$ into a point. Let $c\in\NN$ be the maximum of all such possible $c(m_U)$'s. 
\smallskip
Let $q\in\NN\setminus\{1,2\}$ be such that:
\medskip
{\bf a)} any reductive group over $\FF_{p^3}$ of dimension at most 55 splits over $\FF_{p^q}$ and 
\smallskip
{\bf b)} any \'etale, finite, flat group scheme $GS$ over $\FF_{p^3}$ of rank at most $c$, over $\FF_{p^q}$ is defined by an abstract group. 
\medskip
If any such morphism $m_U$ is not dominant, then the number $NR(q)$ of $\FF_{p^q}$-valued points of $\tilde U_{\FF_{p^3}}^{\rm st}$ is 
$$c(m_U)P(p^q)\leqno (POINTS)$$ 
where $P(x)$ is a polynomial with coefficients in $\ZZ$ which belongs to a finite list of polynomials of order of degree at least equal to $\dim_k(\tilde U)-\dim_k(U)+1$. The finite list depends only on the number $55$ and is obtained easily by combining [Bo2, 21. 15-17] (see also [GLS, Table 2.2 of p. 39]) with the fact that the number of split, reductive groups over $\FF_{p^q}$ of dimension at most 55 is finite, and with the fact that each reduced, connected, affine, algebraic group over $\FF_{p^q}$ is the semidirect product of a reductive group with a unipotent group which is obtained from $\GG_a$ via short exact sequences. By choosing $q$ big enough (depending only on $c$ and on $55$), it can be easily checked by computers that none of these possible expressions for $NR(q)$ is the right one. 
\smallskip
Lower bounds for such $q$'s are effectively computable provided $c$ is effectively computable. In practice is not very easy to compute $c(m_U)$'s; however, (SES) of 3.13.7.1 can be easily used to substitute b) by two other conditions b') and b''). For b') we just need to add replace in b) ``of rank at most $c$" by the requirement that $GS$ has a normal series of length at most 55 whose factors are direct sums of $\ZZ/p\ZZ$. For b'') we just need to substitute in b) $c$ by $c^\prime$ which is the maximum of effectively computable numbers $c^\prime(m_U)$. 
\medskip
{\bf Part B.} Either 3.15.7 BP2 or the estimates of [Ka2, 1.4-5] (to be compared as well with 3.15.8 or with the proof of [Ka2, 2.3.1-2]) imply that there is $n\in\NN$, not depending on $\tilde G$, such that the Newton polygon of $(\tilde M\otimes_{W(k)} W(\bar k),\tilde h(\tilde\vph_0\otimes 1))$ depends only on the reduction mod $p^n$ of $\tilde h\in\tilde G(W(\bar k))$. We can assume $U$ is small enough so that we have an isomorphism 
$$m_U:{\rm Spec}(W_n(R/pR))\to U_{W_n(k)}$$ 
inducing the natural identification on special fibres. We fix some $\tilde h\in G(W(k))$ and we consider the following object of $p-\Mm(W(R/pR))$ 
$$(\tilde M\otimes_{W(k)} W(R/pR),h_{\rm univ}\tilde h(\tilde\vph_0\otimes 1)),$$
with $h_{\rm univ}\in \tilde G(W(R/pR))$ an element which mod $p^n$ is the universal element defined naturally via $m_U$. From the proof (of [Ka2, p. 143]) of the specialization theorem worked just mod $p^n$ and from the existence of $U$ we get directly:
\medskip
{\bf Theorem 1.} {\it The Newton polygon of $(\tilde M\otimes_{W(k)},\tilde h\tilde\vph_0))$ is above the Newton polygon $NP_0$ of $(\tilde M,\tilde\vph_0)$.}
\medskip
Theorem 1 is equivalent to [RR, 4.2]. For future references we state here the following refinement of Theorem 1 for the generalized Shimura context. We assume ${\got C}=(\tilde M,(F^i(\tilde M))_{i\in S(\tilde a+1,\tilde b)},\tilde\vph,\tilde G)$ is a generalized Shimura $p$-divisible object of $\Mm\Mf_{[\tilde a,\tilde b]}(W(k))$. Let $A_{\tilde G}$, $B_{\tilde G}$, $C_{\tilde G}$ and $D_{\tilde G}$ be the subsets of $\tilde G(W(k))$ defined as in 3.1.3 and 3.2.10 (either Theorem 1 or the combination of the arguments of 3.7 with 3.15.6 D allows us to define $A_{\tilde G}$ and $B_{\tilde G}$). 
\medskip
{\bf Theorem 2.} {\it {\bf 1)} We have: $A_{\tilde G}=B_{\tilde G}=C_{\tilde G}=D_{\tilde G}g_0$. 
\smallskip
{\bf 2)} $g\in G(W(k)$ belongs to these sets iff the (refined) Lie stable $p$-rank of $(\tilde M,g\tilde\vph_0,\tilde G)$ is equal to the one of $(\tilde M,\tilde\vph_0,\tilde G)$.
\smallskip
{\bf 3)} If $g\in A_{\tilde G}$, then there is a unique lift of $(\tilde M,g\tilde\vph_0,\tilde G)$ such that the Shimura adjoint filtered Lie $\sg$-crystal attached to it is of parabolic type (i.e. is Shimura-ordinary).}
\medskip
The proof of this is entirely the same as the implications of 3.2.4. We just need to point out: either the deformation part of 3.15.6 or Theorem 1 implies that we can still use the specialization part of 3.2.5-6 as for the context of Shimura $\sg$-crystals (cf. also 3.7). No doubt, once the Faltings--Shimura--Hasse--Witt invariants are defined (see 3.13.7.9), in connection to Theorem 2, the extra hypothesis of being in a generalized Shimura context, will be automatically removed. 
\smallskip
One extra thing: if $p\ge\tilde b-\tilde a+2$, then the Galois property of 4.8.2 holds in the generalized Shimura context of ${\got C}$: no arguments need to be changed, cf. the deformation theory of 3.15.6 and [Fa2, th. 5]. In fact [CF] can be used to eliminate the assumption $p\ge\tilde b-\tilde a+2$ (we recall that 4.8.2 is stated with $\QQ_p$-coefficients). 
\medskip
{\bf Step 2.} If the $W$-condition holds for ${\got C}$ and $\tilde G$ is a reductive group which is not quasi-split, then 3.8 still applies. So, whenever we can define $\tilde G$-canonical lifts over $\bar k$, they are in fact defined over $k$. Moreover, there is $g\in G(W(k))$ such that the Newton polygon of $(\tilde M,g\tilde\vph)$ is the same as the Newton polygon $NP_0$ obtained as in Theorem 1 but for the pull back of ${\got C}$ to $\bar k$.
\medskip
{\bf Step 3.} If the $W$-condition holds for ${\got C}$ but $\tilde G$ is not a reductive group, we have to proceed as follows. For most problems (like isomorphisms classes, etc.) it matters only how $\tilde G^\wedge$ looks like. This, most common, reduces the problem to the case when $\tilde G$ is the extension of a reductive group by a nilpotent group. So to work out 3.4-5 and 3.6.6 in this situation we have to combine the case of a reductive group with the case of a nilpotent group. For how this is accomplished see [Va5] (cf. 4.6.1 2)). Warning: from the point of view of Newton polygons, by giving up the condition that $\tilde G$ is smooth, we can assume $\tilde G_{B(k)}$ is reductive; however, from the point of subtler things (mod powers of $p$) we can not make this assumption (to be compared with 3.9.7.2-3). 
\medskip
{\bf Step 4.} The last step is to deal with the case when the $W$-condition does not hold for ${\got C}$ (to be compared with 2.2.9 1) and 1')). This, at a first glance, looks as an obstacle in getting 3.4-5 and 3.6.6 in this context. However, if $\tilde G$ is reductive, Theorem 1 should be provable in this context, either by using global deformations or by shifting the focus from ${\got C}$ to 
$${\rm Lie}({\got C}):=({\rm Lie}(\tilde G),(F^j({\rm Lie}(\tilde G))_{j\in S(L\tilde a,L\tilde b)},\tilde\vph,\tilde G/Z(\tilde G)).$$ 
Here $L\tilde a,L\tilde b\in\ZZ$, $L\tilde b\ge L\tilde a$, are the logical values; the quotient group $\tilde G/Z(\tilde G)$ acts on ${\rm Lie}(\tilde G)$ via inner conjugation. How ``successful" we are (i.e. of what we achieved) in this shift of focus, can be read out from 2.2.16.5; for instance, if $\tilde G/Z(\tilde G)$ is an adjoint group and if $({\rm Lie}(Z(\tilde G)),\tilde\vph)$ has only slopes $0$, then we are in the context of 1) or 2).
We hope to come back to this step in a future paper.
\medskip
As an afterthought, we point out: using [RR] and unpublished work of T. Zink, in [We] it is also checked the density part of 4.2.1 involving its a) part, in the cases mentioned in 4.6.6 as well as in the case $p=2$ of the $A$ and $C$ cases of [Ko2].
\vfill\eject
\centerline{}
\bigskip
\bigskip
\centerline{\bigsll {\bf References}}
\bigskip\bigskip
\medskip
\item{[BBM]} P. Berthelot, L. Breen and W. Messing, {\it Th\'eorie de Dieudonn\'e crystalline}, LNM {\bf 930} (1982), Springer--Verlag.
\medskip
\item{[Be]} P. Berthelot, {\it Comologie cristalline des sch\'emas de caract\'eristique $p>0$}, LNM {\bf 407} (1974), Springer--Verlag.
\medskip
\item{[BHC]} A. Borel and Harish-Chandra, {\it Arithmetic subgroups of algebraic groups}, Ann. of Math. (3) {\bf 75} (1962), p. 485--535.
\medskip
\item{[Bl]} D. Blasius, {\it A p-adic property of Hodge cycles on abelian varieties}, Proc. Symp. Pure Math. {\bf 55}, Part 2 (1994), p. 293--308.
\medskip
\item{[BLR]} S. Bosch, W. L\"utkebohmert and M. Raynaud, {\it N\'eron models}, Springer--Verlag (1990).
\medskip
\item{[BM]}  P. Berthelot and W. Messing, {\it Th\'eorie de Dieudonn\'e crystalline III}, in the Grothendieck Festschrift I, Birkh\"auser, Progr. in Math. {\bf 86} (1990), p. 173--247. 
\medskip
\item{[Bo1]} A. Borel, {\it Properties and linear representations of Chevalley groups}, LNM {\bf 131} (1970), Springer--Verlag, p. 1--55.
\medskip
\item{[Bo2]} A. Borel, {\it Linear algebraic groups}, Grad. Text Math. {\bf 126} (1991), Springer--Verlag.
\medskip
\item{[Bou1]} N. Bourbaki, {\it Lie groups and Lie algebras}, Chapters {\bf 1--3}, Springer--Verlag (1989).
\medskip
\item{[Bou2]} N. Bourbaki, {\it Groupes et alg\`ebres de Lie}, Chapitre {\bf 4--6}, Hermann (1968).
\medskip
\item{[Bou3]} N. Bourbaki, {\it Groupes et alg\`ebres de Lie}, Chapitre {\bf 7--8}, Hermann (1975).
\medskip
\item{[BT]} A. Borel and J. Tits, {\it Homomorphismes \lq\lq abstraits\rq\rq\ de groupes alg\`ebriques simples}, Ann. of Math. (3) {\bf 97} (1973), p. 499--571.
\medskip
\item{[Bu1]} A. Buium, {\it Geometry of $p$-jets}, Duke Math. J., Vol. {\bf 82}, no. 2 (1996), p. 349--367.
\medskip
\item{[Bu2]} A. Buium, {\it Geometry of Fermat Ad\`eles}, preprint (1999), p. 1--65.
\medskip
\item{[CF]} P. Colmez and J.-M. Fontaine, {\it Construction des représentations $p$-adiques semi-stables}, Inv. Math. {\bf 140}, no. 1 (2000), p. 1--43.
\medskip
\item{[Ch1]} C.-L. Chai, {\it Arithmetic Compactifications of Hilbert modular spaces}, Appendix to A. Wiles ``The Iwasawa conjecture for totally real number fields", Ann. of Math. (3) {\bf 131} (1990), p. 541--554.
\medskip
\item{[Ch2]} C.-L. Chai, {\it Every ordinary symplectic isogeny class in positive characteristic is dense in the moduli}, Inv. Math. {\bf 121}, no. 3 (1995), p. 439--479.
\medskip
\item{[Ch3]} C.-L. Chai, {\it Newton polygons as lattice points}, Amer. J. Math. {\bf 122}, no. 5 (2000), 967--990.
\medskip
\item{[Cr]} J. E. Cremona, {\it Algorithms for modular elliptic curves}, Cambridge Univ. Press, Cambridge (1992).
\medskip
\item{[De1]} P. Deligne, {\it Travaux de Shimura}, S\'em. Bourbaki, Expos\'e 389, LNM {\bf 244} (1971), Springer--Verlag.
\medskip
\item{[De2]} P. Deligne, {\it Vari\'et\'es de Shimura: Interpretation modulaire, et techniques de construction de mod\`eles canoniques}, Proc. Symp. Pure Math. {\bf 33}, Part 2 (1979), p. 247--290.
\medskip
\item{[De3]} P. Deligne, {\it Cristaux ordinaires et coordon\'ees canoniques}, LNM {\bf 868} (1981), Springer--Verlag, p. 80--137.
\medskip
\item{[De4]} P. Deligne, {\it Hodge cycles on abelian varieties}, Hodge cycles, motives, and Shimura varieties,
LNM {\bf 900} (1982), Springer--Verlag, p. 9--100.
\medskip
\item{[De5]} P. Deligne, {\it Cat\'egories tannakiennes}, The Grothendieck Festschrift, Progr. Math. {\bf 87}, Vol. II (1990), Birkh\"auser,  p. 111--196.
\medskip
\item{[De6]} P. Deligne, {\it Le groupe fondamental de la droite projective moins trois points}, Galois groups over $\QQ$, Math. Sci. Res. Inst. Publ. {\bf 16} (1989), p. 79--297. 
\medskip
\item{[Dem]} M. Demazure, {\it Lectures on $p$-divisible groups}, LNM {\bf 302} (1972), Springer--Verlag.
\medskip
\item{[Di]} J. Dieudonn\'e, {\it Groupes de Lie et hyperalg\`ebres de Lie sur un corps de caract\'erisque $p>0$ (VII)}, Math. Annalen {\bf 134} (1957), p. 114--133. 
\medskip
\item{[EGA III]} A. Grothendieck et al., {\it Etude locale des sch\'emas et des morphismes de sch\'ema. I.}, Publ. Math. IHES {\bf 11} (1961).
\medskip
\item{[EGA IV]} A. Grothendieck et al., {\it Etude locale des sch\'emas et des morphismes de sch\'ema}, Publ. Math. IHES {\bf 20} (1964), {\bf 24} (1965), {\bf 28} (1966), {\bf 32} (1967).
\medskip
\item{[EO]} T. Ekedahl and F. Oort, {\it Connected subsets of a moduli space of abelian
varieties}, preprint (1994).
\medskip
\item{[Fa1]} G. Faltings, {\it Crystalline cohomology and $p$-adic Galois representations},
Algebraic Analysis, Geometry, and Number Theory, Johns Hopkins Univ. Press (1990), p. 25--79.
\medskip
\item{[Fa2]} G. Faltings, {\it Integral crystalline cohomology over very ramified valuation rings}, J. of Am. Math. Soc., Vol. {\bf 12}, no. 1 (1999), p. 117--144.
\medskip
\item{[Fa3]} G. Faltings, {\it Almost \'etale extensions}, preprint series 1998 (3), MPI, Bonn.
\medskip
\item{[FC]} G. Faltings and C.-L. Chai, {\it Degeneration of abelian varieties},  A series of modern surveys in math., Vol. {\bf 22} (1990), Springer--Verlag.
\medskip
\item{[FI]} J. -M. Fontaine and L. Illusie, {\it $p$-adic periods: a survey}, Bombay (1989), Pr\'epublic. of Univ. Paris-Sud 90--53.
\medskip
\item{[FL]} J. -M. Fontaine and G. Laffaille, {\it Construction de repr\'esentations p-adiques}, Ann. Sci. \'Ecole Norm. Sup. {\bf 15}, no. 4 (1982), p. 547--608.
\medskip
\item{[Fo]} J.-M. Fontaine, {\it Groupes $p$-divisibles sur les corps locaux}, J. Ast\'erisque {\bf 47--48} (1977).
\medskip
\item{[GLS]} D. Gorenstein, R. Lyons and R. Soloman, {\it The classification of the finite simple groups, Number 3}, Surv. and Monog., Vol. {\bf 40} (1998), No. 3, A. M. S.   
\medskip
\item{[Go]} B.B. Gordon, {\it Canonical models of Picard modular surfaces}, The Zeta function of Picard modular surfaces, Les Publications CRM, Montreal (1992), p. 1--28.
\medskip
\item{[GO]} E. Z. Goren and F. Oort, {\it Stratifications of Hilbert modular varieties}, J. Alg. Geom. {\bf 9}, no. 1 (2000), p. 111--154.
\medskip
\item{[Gr]} A. Grothendieck, {\it Groupes de Barsotti--Tate et cristaux de Dieudonn\'e}, S\'em. Math. Sup. {\bf 45}, Pr. de l'Universit\'e de Montr\'eal (1970).
\medskip
\item{[Ha]} R. Hartshorne, {\it Algebraic geometry}, Grad. Text Math. {\bf 52} (1977), Springer--Verlag.
\medskip
\item{[Har]} G. Harder, {\it \"Uber die Galoiskohomologie halbeinfacher Matrizengruppen II}, Math. Z. {\bf 92} (1966), p. 396--415.
\medskip
\item{[He]} S. Helgason, {\it Differential geometry, Lie groups, and symmetric spaces}, Academic Press, New--York (1978).
\medskip
\item{[Hr]} E. Hrushovski, {\it The First Order Theory of the Frobenius}, manuscript (4/5/1999).
\medskip
\item{[Hu1]} J. E. Humphreys, {\it Introduction to Lie algebras and representation theory}, Grad. Text. Math. {\bf 9} (1997), seventh edition, Springer--Verlag.
\medskip
\item{[Hu2]} J. E. Humphreys, {\it Conjugacy classes in semisimple algebraic groups}, Math. Surv. and Monog., Vol. {\bf 43} (1995), A. M. S. 
\medskip
\item{[HK]} R. Hain and M. Kim, {\it The crystalline fundamental group and the De Rham--Witt complex}, preprint (4/2001).
\medskip
\item{[Ih]} Y. Ihara, {\it Shimura curves over finite fields and their rational points}, Applications of Curves over finite fields, Cont. Math. {\bf 245} (1999), p. 15--23.
\medskip
\item{[Il]} L. Illusie, {\it D\'eformations des groupes de Barsotti--Tate}, J. Ast\'erisque {\bf 127} (1985), p. 151--198.
\medskip
\item{[IKO]} T. Ibukiyama, T. Katsura and F. Oort, {\it Supersingular curves of genus two and class numbers}, Comp. Math. {\bf 57} (1986), p. 127--152.
\medskip
\item{[Ja]} J.C. Jantzen, {\it Representations of algebraic groups}, Academic Press  (1987).
\medskip
\item{[dJO]} J. de Jong and F. Oort, {\it Purity of the stratification by Newton polygons}, J. of Amer. Math. Soc. {\bf 13}, no. 1 (2000), p. 209--241.
\medskip
\item{[dJ1]} J. de Jong, {\it Finite locally free group schemes in characteristic $p$ and Dieudonn\'e modules}, Inv. Math. {\bf 114}, no. 1 (1993), p. 89--137.
\medskip
\item{[dJ2]} J. de Jong, {\it Crystalline Dieudonn\'e module theory via formal and rigid geometry}, Publ. Math. IHES {\bf 82} (1995), p. 5--96.
\medskip
\item{[dJ3]} J. de Jong, {\it Homomorphisms of Barsotti--Tate groups and crystals in positive characteristic}, Inv. Math. {\bf 134}, no. 2 (1998), p. 301--333. 
\medskip
\item{[Ka1]} N. Katz, {\it Travaux de Dwork}, S\'eminaire Bourbaki, Exp. 409, LNM {\bf 383} (1983), Springer--Verlag.
\medskip
\item{[Ka2]} N. Katz, {\it Slope filtration of $F$-crystals}, Journ\'ees de G\'eom\'etrie alg\'ebrique de Rennes 1978, J. Ast\'erisque {\bf 63}, Part 1 (1979), p. 113--163.
\medskip
\item{[Ka3]} N. Katz, {\it Serre--Tate local moduli}, LNM {\bf 868} (1981), Springer--Verlag, p. 138--202.
\medskip
\item{[Ka4]} N. Katz, {\it Appendix to expose V}, LNM {\bf 868} (1981), Springer--Verlag, p. 128--137.
\medskip
\item{[KMRT]} M.-A. Knus, A. Merkurjev, M. Rost and J.-P. Tignol, {\it The book of involutions}, Coll. Public., Vol. {\bf 44} (1998), A. M. S. 
\medskip
\item{[Kob]} N. Koblitz, {\it p-adic valuations of zeta function over families of varieties defined over finite fields}, Comp. Math. {\bf 31}, no. 2 (1975), p. 119--218. 
\medskip
\item{[Ko1]} R.E. Kottwitz, {\it Isocrystals with additional structure}, Comp. Math. {\bf 56}, no. 2 (1985), p. 201--220.
\medskip
\item{[Ko2]} R.E. Kottwitz, {\it Points on some Shimura varieties over finite fields}, J. of Am. Math. Soc., Vol. {\bf 5}, no. 2 (1992), p. 373--444.
\medskip
\item{[La]} G. Laffaille, {\it Groupes p-divisibles et modules filtr\'es: le cas peu ramifi\'e}, Bull. Soc. Math. de France {\bf 108} (1980), p. 187--206.
\medskip
\item{[LO]} K.-Z. Li and F. Oort, {\it Moduli of supersingular abelian varieties}, LNM {\bf 1680} (1998), Springer--Verlag.
\medskip
\item{[LR]} R. Langlands and M. Rapoport, {\it Shimuravarietaeten und Gerben}, J. reine angew. Math. {\bf 378} (1987), p. 113--220.
\medskip
\item{[Ma]} H. Matsumura, {\it Commutative algebra}, The Benjamin/Cummings Publishing Co., Inc. (1980).
\medskip
\item{[Man]} J.I. Manin, {\it The theory of formal commutative groups in finite characteristic}, Russian Math. Surv. {\bf 18}, no. 6 (1963), p. 1--83. 
\medskip
\item{[M-B]} L. Moret-Bailly, {\it Pinceaux de vari\'et\'es ab\'eliennes}, J. Ast\'erisque {\bf 129} (1985).
\medskip
\item{[Me]} W. Messing, {\it The crystals associated to Barsotti--Tate groups, with applications to abelian schemes}, LNM {\bf 264} (1972), Springer--Verlag.
\medskip
\item{[MFK]} D. Mumford, J. Fogarty and F. Kirwan, {\it Geometric invariant theory}, A series of modern surveys in math., Vol. {\bf 34} (1994), third edition, Springer--Verlag.
\medskip
\item{[Mi1]} J.S. Milne, {\it The points on a Shimura variety modulo a prime of good reduction}, The Zeta function of Picard modular surfaces, Les Publications CRM, Montreal (1992), p. 153--255.
\medskip
\item{[Mi2]} J.S. Milne, {\it On the conjecture of Langlands and Rapoport}, manuscript (9/1995).
\medskip
\item{[Mi3]} J.S. Milne, {\it Shimura varieties and motives}, Proc. Symp. Pure Math. {\bf 55}
(1994), Part 2, p. 447--523.
\medskip
\item{[Mi4]} J.S. Milne, {\it Canonical models of (mixed) Shimura varieties and automorphic vector bundles}, in Automorphic forms, Shimura varieties and L-functions, Vol. I, Perspectives in Math. {\bf 10} (1990), Academic Press, p. 283--414.
\medskip
\item{[MM]} T. Matsusaka and D. Mumford, {\it Two fundamental theorems on deformations of polarized varieties}, Am. J. Math. {\bf 86} (1964), p. 668--683.
\medskip
\item{[Mo]} B.J.J. Moonen, {\it Special points and linearity properties of Shimura varieties}, Ph.D. thesis, Univ. of Utrecht (1995). 
\medskip
\item{[Mu]} D. Mumford, {\it Abelian varieties}, Tata Inst. of Math. Research (1988).
\medskip
\item{[No1]} R. Noot, {\it Hodge classes, Tate classes, and local moduli of abelian varieties}, Thesis Rijksuniv. Utrecht (1992).
\medskip
\item{[No2]} R. Noot, {\it Models of Shimura varieties in mixed characteristc}, J. Alg. Geom. {\bf 5}, no. 1 (1996), p. 187--207.
\medskip
\item{[Nor]} P. Norman, {\it Lifting abelian varieties}, Inv. Math. {\bf 64} no. 3 (1981), 431--443.
\medskip
{\item{[NO]} P. Norman and F. Oort, {\it Moduli of abelian varieties}, Ann. of Math. (3) {\bf 112} (1980), p. 413--439.
\medskip
\item{[Og]} A. Ogus, {\it Elliptic crystals and modular motives},
manuscript (3/8/2000), to appear in Adv. in Math.
\medskip
\item{[Oo1]} F. Oort, {\it Some questions in algebraic geometry}, preprint Utrecht Univ. (1995).
\medskip
\item{[Oo2]} F. Oort, {\it Newton polygons and formal groups: conjectures by Manin and Grothendieck}, Ann. of Math. (2) {\bf 152} (2000), p. 183--206.
\medskip
\item{[Oo3]} F. Oort, {\it A stratification of a moduli space of polarized abelian varieties in positive characteristic}, Asp. of Math. E {\bf 33} (1999), p. 47--64.
\medskip
\item{[Pi]} R. Pink, {\it $l$-adic algebraic monodromy grpups, cocharacters, and the Mumford--Tate conjecture}, J. reine angew. Math. {\bf 495} (1998), p. 187--237.
\medskip
\item{[Pr]} G. Prasad, {\it Strong approximation for semi-simple groups over function fields}, Ann. of Math. (2) {\bf 105} (1977), p. 553--572.
\medskip
\item{[Ra]} M. Raynaud, {\it Anneaux Locaux Hens\'eliens}, LNM {\bf 169} (1970), Springer--Verlag.
\medskip
\item{[Ro]} J. J. Rotman, {\it An introduction to the theory of groups}, Grad. Text Math. {\bf 52} (1995), fourth edition, Springer--Verlag.
\medskip
\item{[RR]} M. Rapoport and M. Richartz, {\it On the classification and specialization of $F$-isocrystals with additional structure}, Comp. Math. {\bf 103} (1996), p. 153--181.
\medskip
\item{[Sa]} I. Satake, {\it Holomorphic embeddings of symmetric domains into a Siegel space}, Am. J. Math. {\bf 87} (1965), p. 425--462.
\medskip
\item{[Sam]} P. Samuel, {\it Lectures on unique factorization domains}, Tata Inst. of Math. Research (1964).
\medskip
\item{[Se1]} J.-P. Serre, {\it Sur les groupes de Galois attach\'es aux groupes $p$-divisibles}, Proc. Conf. Local Fields, Springer--Verlag, Heidelberg (1967), p. 118--131.
\medskip
\item{[Se2]} J.-P. Serre, {\it Groupes alg\'ebriques associ\'es aux modules de Hodge--Tate}, Journ\'ees de G\'eom. Alg. de Rennes, J. Ast\'erisque {\bf 65} (1979), p. 155--188.
\medskip
\item{[SGA1]} A. Grothendieck, et al. {\it Rev\^etements \'etales et groupe fondamental}, LNM {\bf 224} (1971), Springer--Verlag.
\medskip
\item{[SGA3]} M. Demazure, A. Grothendieck, et al. {\it Schemes en groupes}, Vol. I-III,
LNM {\bf 151--153} (1970), Springer--Verlag.
\medskip
\item{[SGA4]} M. Artin, A. Grothendieck, et al. {\it Th\'eorie des topos et cohomologie \'etale des sch\'emas}, Vol. I, LNM {\bf 269} (1972), Springer--Verlag.
\medskip
\item{[Sh]} G. Shimura, {\it Moduli of abelian varieties and number theory}, Proc. Sympos. Pure Math., Vol. {\bf 9} (1966), p. 312--332, A. M. S., Providence, R. I.  
\medskip
\item{[Shi1]} A. Shiho, {\it Crystalline fundamental groups. I. Isocrystals on log crystalline site
and log convergent site}, J. Math. Sci. Univ. Tokyo {\bf 7}, no. 4 (2000), p. 509--656.
\medskip
\item{[Shi2]} A. Shiho, {\it Crystalline fundamental groups. II}, preprint (1/2001).
\medskip
\item{[Ta]} J. Tate, {\it Classes d'isog\'enie des vari\'et\'es sur un corps fini (d'apr\`es J. Honda)}, S\'eminaire Bourbaki 1968/69, Exp. 352, LNM {\bf 179} (1971), p. 95--110.
\medskip 
\item{[Ti1]} J. Tits, {\it Classification of algebraic semisimple groups}, Proc. Sympos. Pure Math., Vol. {\bf 9} (1966), p. 33--62, A. M. S., Providence, R. I.
\medskip
\item{[Ti2]} J. Tits, {\it Reductive groups over local fields}, Proc. Sympos. Pure Math. {\bf 33} (1979), Part 1, p. 29--69.
\medskip
\item{[Va1]} A. Vasiu, {\it Integral canonical models for Shimura varieties of Hodge type}, thesis, Princeton Univ. (1994).
\medskip
\item{[Va2]} A. Vasiu, {\it Integral canonical models of Shimura varieties of preabelian
 type}, Asian J. Math., Vol. {\bf 3}, no. 2 (1999), p. 401--518.
\medskip
\item{[Va3]} A. Vasiu, {\it Toroidal smooth compactifications of integral canonical
models of Shimura varieties of preabelian type}, to be submitted for publication.
\medskip
\item{[Va4]} A. Vasiu, {\it Examples of Shimura varieties and of their integral canonical models}, to be submitted for publication.
\medskip
\item{[Va5]} A. Vasiu, {\it The $p=2$ and $p=3$ theory of Shimura varieties of preabelian type and the existence of integral canonical models of quotients of Shimura
varieties of Hodge type with respect to non-hyperspecial subgroups}, to be submitted for publication.
\medskip
\item{[Va6]} A. Vasiu, {\it Generalized Shimura $p$-divisible objects of special type}, to be submitted for publication.
\medskip
\item{[Va7]} A. Vasiu, {\it Density of Hecke orbits in positive characteristic}, to be submitted for publication.
\medskip
\item{[Va8]} A. Vasiu, {\it On the classification of cyclic diagonalizable $p$-divisible objects}, to be submitted for publication.
\medskip
\item{[Va9]} A. Vasiu, {\it On the $U$-ordinariness in geometric situations}, to be submitted for publication.
\medskip
\item{[Va10]} A. Vasiu, {\it Shimura varieties and the Mumford--Tate conjecture, part I}, version 1/15/2000, to be submitted for publication.
\medskip
\item{[Va11]} A. Vasiu, {\it On the classification of versal generalized Shimura $p$-divisible objects}, to be submitted for publication.
\medskip
\item{[Va12]} A. Vasiu, {\it Standard arithmetics of classes of polarized projective varieties: examples, results and problems}, to be submitted for publication. 
\medskip
\item{[Va13]} A. Vasiu, {\it Points of integral canonical models of Shimura varieties of preabelian type and the Langlands-Rapoport conjecture}, preprint, ETH, FIM, Z\"urich, Sep. 1996, p. 1--76.
\medskip
\item{[Vo]} P. Vojta, {\it Nagata's embedding theorem}, preprint (appendix to planned book; 1997).
\medskip
\item{[We]} T. Wedhorn, {\it Ordinariness in good reductions of Shimura varieties of PEL-type}, Ann. Sci. \'Ecole Norm. Sup. (4) {\bf 32} (1999), p. 575--618.
\medskip
\item{[Wi]} J.-P. Wintenberger, {\it Un scindage de la filtration de Hodge pour certaines
variet\'es algebriques sur les corps locaux}, Ann. of Math. (2) {\bf 119}, (1984), p. 511--548.
\vskip 0.5 in
\line{\hfill\vbox{\hbox{\it address:}
\medskip
\hbox{Adrian Vasiu}
\hbox{University of Arizona}
\hbox{Department of Mathematics}
\hbox{617 N. Santa Rita}
\hbox{P.O. Box 210089}
\hbox{Tucson, AZ-85721-0089, U.S.A.}
\medskip
\hbox{e-mail: adrian@math.arizona.edu\ \ \ \ \ }}}

\end